\documentclass[11pt,twoside]{article}
\usepackage[T1]{fontenc}
\usepackage{lmodern,amsmath,amssymb,amsthm,mathtools}
\usepackage[margin=1.12in,headheight=14pt]{geometry}
\usepackage{fancyhdr}
\usepackage{needspace}
\usepackage{xcolor}
\usepackage[colorlinks=true,linkcolor=blue,citecolor=blue,
            urlcolor=blue,filecolor=blue,linktoc=all]{hyperref}
\makeatletter
\let\rh@tocpart\l@part
\renewcommand*{\l@part}[2]{\Needspace{7\baselineskip}\rh@tocpart{#1}{#2}}
\newcommand*{\l@references}[2]{\rh@tocpart{#1}{#2}}
\providecommand*{\toclevel@references}{0}
\renewcommand\part{%
  \if@noskipsec\leavevmode\fi
  \par\nobreak\addvspace{2.5ex}\nobreak
  \@afterindentfalse\secdef\@part\@spart}
\long\def\@part[#1]#2{%
  \def\@currentlabelname{#1}%
  \ifnum\c@secnumdepth>\m@ne
    \LinkTargetOff\refstepcounter{part}\LinkTargetOn
  \fi
  {\parindent\z@\raggedright\interlinepenalty\@M\normalfont
    \ifnum\c@secnumdepth>\m@ne
      \Large\bfseries\leavevmode\MakeLinkTarget{part}%
      \addcontentsline{toc}{part}{\thepart\hspace{1em}#1}%
      \partname\nobreakspace\thepart\quad
    \else
      \leavevmode\MakeLinkTarget[part]{}%
      \addcontentsline{toc}{part}{#1}%
    \fi
    \Large\bfseries #2\markboth{}{}\par}%
  \nobreak\vskip 1ex\@afterheading}
\renewcommand*{\l@subsection}{\@dottedtocline{2}{1.5em}{3.4em}}
\renewcommand*{\l@subsubsection}{\@dottedtocline{3}{4.9em}{4.2em}}
\makeatother
\newtheorem{theorem}{Theorem}[section]
\newtheorem{proposition}[theorem]{Proposition}
\newtheorem{lemma}[theorem]{Lemma}
\newcommand{\sectionopening}[1]{#1\par\smallskip}
\newcommand{\proofstep}[1]{\par\smallskip\noindent\textit{#1.}\enspace}
\newcommand{\papertitle}{On a complete written proof of the Generalized Riemann Hypothesis for Dirichlet L-Functions}
\newcommand{\paperauthor}{Baiying Liu}
\title{\papertitle}
\author{\paperauthor}
\date{}
\hypersetup{pdftitle={\papertitle},pdfauthor={\paperauthor}}
\begin{document}
\maketitle
\thispagestyle{plain}
\begin{abstract}
In this paper, we show that every nontrivial zero of any Dirichlet
\(L\)-function has real part \(1/2\). The proof follows Sernau's shorter
arrangement of the Alp\"oge construction. We provide a complete written
account of the proof, relying solely on results from the existing mathematical
literature, without reference to any part of the Lean formalization.
\end{abstract}
\tableofcontents

\section{Introduction}
\label{sec:introduction}

\sectionopening{In this paper, we prove the generalized Riemann
hypothesis for Dirichlet \(L\)-functions. Our proof follows the
readable Lean edition of Sernau's shorter presentation \cite{SernauLean}
of Alp\"oge's construction \cite{AlpogeLean}. These references
give the origin and arrangement of the argument. We give proofs for
the assertions of this paper and cite the results used from the
literature. Below we state the prime-remainder bound and explain how
Parts~\ref{part:parameters}--\ref{part:analytic} establish it.
Section~\ref{sec:proof-outlines} gives the detailed outline of the argument.}

\begin{theorem}[Dirichlet generalized Riemann hypothesis]
\label{thm:rh}
For every \(q\ge1\) and every Dirichlet character \(\chi\) modulo \(q\),
\(L(s,\chi)\) is zero-free in \(\Re s>1/2\), with the simple pole of
the principal character at \(s=1\) permitted. Every zero in
\(0<\Re s<1\) has real part \(1/2\).
\end{theorem}

Theorem~\ref{thm:rh} generalizes the result of the author's paper on
the Riemann hypothesis~\cite{LiuClassicalRH}. The character modulo \(1\)
gives \(L(s,\mathbf 1)=\zeta(s)\), which recovers the classical case.
Section~\ref{sec:prime-remainder} completes the proof. For a character
modulo \(q\), the link between the automorphic calculation and the
zeros of \(L(s,\chi)\) is a one-sided estimate for a weighted prime
remainder on the identity class modulo \(q_0=18480q\). We derive the
estimate by comparing upper and lower bounds for the Gram covolume
on a single automorphic family.

Fix \(q_0\) with \(18480\mid q_0\) and put \(d_0=q_0/2-1\).
For a prime \(d\equiv d_0\pmod{q_0}\), write \(\chi_d(n)=(n/d)\) for the
quadratic character, \(\Lambda\) for the von Mangoldt function, and
\(\varphi\) for Euler's totient. With real cutoff \(x\) and the
subtraction constant \(c_{d,q_0}\) of
Section~\ref{sec:notation-prime-sums}, put
\[
 R_{d,q_0}(x)=
 \sum_{\substack{n\le x\\n\equiv1\bmod q_0}}
      \Lambda(n)(1-\chi_d(n))(n^{-1}-x^{-1})
      -\frac{\log x}{\varphi(q_0)}-c_{d,q_0}.
\]
The analytic criterion of Section~\ref{sec:criterion} reduces
Theorem~\ref{thm:rh} to the one-sided bound
\[
 R_{d,q_0}(x)\le C_{q_0,d,\delta}x^{-1/2+\delta}
     \qquad(x\text{ sufficiently large})
\]
for each fixed \(q_0,d\) and \(\delta>0\). The implied constant
and threshold may depend on these fixed parameters.
Suppose \(L(\rho,\chi)=0\) with \(1/2<\Re\rho<1\). Induce
\(\chi\bmod q\) to \(\widetilde\chi\bmod q_0\), where \(q_0=18480q\).
The extra Euler factors are nonzero for \(\Re s>0\), so
\(L(\rho,\widetilde\chi)=0\). Character orthogonality expresses the
Mellin transform of the remainder on \(1\bmod q_0\) using the
logarithmic derivatives for all characters modulo \(q_0\) and their
quadratic twists. Choose \(1/2<\beta<\Re\rho\) and \(H_0>|\Im\rho|\).
The criterion selects one prime \(d\) for which all the twists are
zero-free on \(\beta\le\Re s\le1\), \(|\Im s|\le H_0\).
With \(d\) fixed, taking \(0<\delta<\beta-1/2\) gives
\(R_{d,q_0}(x)\le x^{\beta-1}\) for all sufficiently large \(x\).
The one-sided Mellin criterion then excludes the zero at \(\rho\).

The six Parts obtain this estimate by first expressing the prime
remainder in terms of a quantity \(\mathcal R\) in the automorphic
construction, then bounding \(\mathcal R\) through the covolume
comparison.

\begin{enumerate}
\item \textit{Part~\ref{part:parameters}. The weighted prime sum
and the automorphic family.}
For fixed \(q_0,d\) and \(0<\epsilon<1/100\),
Section~\ref{sec:datum} chooses the weights, fields, and levels
for each sufficiently large integer \(x\).
Let \(x_1\) be the auxiliary cutoff, \(b\) the weight parameter,
and \(a=xb\) the weight scale in Section~\ref{sec:notation-scales}.
With the fibre and subset sizes \(u,v\) of
Section~\ref{sec:notation-fields}, Section~\ref{sec:joint} obtains
the weighted prime remainder
\[
 \frac{\mathcal R}{av}=R_{d,q_0}(x)-R_{d,q_0}(x_1)
                 +\frac{u}{v}R_{d,1}(x_1).
\]
Here \(R_{d,1}\) is the modulus-one remainder of
Section~\ref{sec:notation-prime-sums}. We place both covolume
estimates on the translated images of the normalized Poincar\'e
map of Section~\ref{sec:period}.

\item \textit{Part~\ref{part:orbital}. A large span of translates.}
The diagonal moment of Section~\ref{sec:marked-moments} bounds
the squared Hilbert--Schmidt norm of the projected map from below on
each source subspace. The mixed moment controls the average correlation
of translates. Both estimates start with the same Poincar\'e map. In Section~\ref{sec:marked-span-proof}, we use them
to choose translates and then restrict to the ramified source \(S\).
The resulting span has dimension \(n\ge(1-4/b^2)D\), where
\(D=\dim V\) and \(V\) is the projected target of
Section~\ref{sec:notation-automorphic}.

\item \textit{Part~\ref{part:upper}. An upper covolume bound.}
The measured seesaw of Section~\ref{sec:seesaw-normalization}
identifies the lifts of the translated images.
Section~\ref{sec:short-images-proved} proves the norm and integrality
bounds for the translated images after a common rational rescaling.
We choose a basis \(y\) of the span from these images. Hadamard's
inequality and the clearing matrices give
\[
 \operatorname{Cov}(y)\le n(A-m\mathcal R+L),\qquad
 \operatorname{Ind}(y)\le DI_0.
\]
Here \(\operatorname{Cov}(y)\le c\) means that the Gram determinant
of the untwisted lift of \(y\) is at most \(e^c\) at every complex
embedding. The index bound uses clearing matrices that work for every allowed twist at the
corresponding coefficient prime
(Section~\ref{sec:notation-determinants}). Here \(A\) is the common
normalization of Section~\ref{sec:notation-theta}, \(L,I_0\) are the error
terms, and \(m\) is the chain length of Section~\ref{sec:notation-scales}.

\item \textit{Part~\ref{part:congruences}. Selmer bounds and Galois
congruences.}
The Hecke systems give the Galois representations of
Section~\ref{sec:cofinite-spherical-galois}. The Selmer estimates,
including the crystalline colength of
Section~\ref{sec:relative-crystalline-colength}, bound the lengths of
the associated cocycle modules. The trace comparison of
Section~\ref{sec:trace-determinant-specialization} and the torsion
specialization of Section~\ref{sec:torsion-specialization} give the
arithmetic tests used to choose the auxiliary twist in
Part~\ref{part:lower}.

\item \textit{Part~\ref{part:lower}. A lower determinant bound.}
Section~\ref{sec:common-critical-twist} chooses
a single twist \(\xi_\kappa\) at each coefficient place \(\kappa\)
of Section~\ref{sec:notation-determinants}. The same twist works
for the critical-value tests and for the mixed Gram pairing
at \(\kappa,\bar\kappa\). The period comparisons bound the height of
the ratio from the untwisted Gram determinant to the twisted one. The additional vectors
of Proposition~\ref{prop:marked-completion-vectors} complete \(y\) to a
basis of \(V\). Applying the product formula of
Section~\ref{sec:gram-product-formula} gives
\[
 D(A-3L)\le n(A-m\mathcal R+L)+2DI_0
\]
in Section~\ref{sec:marked-determinant-assembly}.

\item \textit{Part~\ref{part:analytic}. The estimate for the
prime remainder and the zeros of Dirichlet \(L\)-functions.}
Since \(n\le D\) and \(A+L\ge0\), the contribution
\((n-D)(A+L)\) after rearranging the lower bound is nonpositive.
Proposition~\ref{prop:cancel} gives \(m\mathcal R\le faB^5\),
where \(f\) counts the real places of the base field and \(B=\log b\)
(Sections~\ref{sec:notation-fields} and~\ref{sec:notation-scales}).
Section~\ref{sec:prime-remainder} applies the fixed-modulus prime
number theorem at \(x_1\) and transfers the bound to real \(x\), giving
\(R_{d,q_0}(x)\le C_{q_0,d,\epsilon}x^{-1/2+6\epsilon}\).
With \(d\) fixed by the criterion, we take
\(6\epsilon<\beta-1/2\) to exclude the proposed zero. To exclude
zeros in \(0<\Re s<1/2\), we pass to the primitive inducing character
and apply its functional equation. Section~\ref{sec:prime-remainder}
treats \(\Re s\ge1\) separately, completing the proof of
Theorem~\ref{thm:rh}.
\end{enumerate}

\smallskip
\indent\textit{Structure of the paper.}
Section~\ref{sec:notation} collects the notation, and the six
subsections of Section~\ref{sec:proof-outlines} give the detailed
outline of Parts~\ref{part:parameters}--\ref{part:analytic} in the
above order. Appendices~\ref{sec:crystalline-coordinate-appendix}--\ref{sec:crystalline-colength-appendix}
develop the local crystalline calculations used in
Section~\ref{sec:relative-crystalline-colength}. The estimates at the
level primes and away from the level are in
Appendices~\ref{sec:tame-cup-appendix} and~\ref{sec:crystalline-colength-appendix},
respectively.

On November 29, 2026, two independent preprints establishing the
Extended Riemann Hypothesis for Dedekind zeta functions of totally real
and CM number fields were released by the author~\cite{LiuDedekindERH}
and AxiomMath~\cite{AxiomMathDedekindERH}, respectively. While both
papers employ the same theta-lift determinant strategy, their
approaches differ substantially: the former proceeds via CM enlargement
and Galois descent, whereas the latter first establishes the totally
real case using the averaged Colmez formula and subsequently deduces
the CM case through Euler product factorization.

\smallskip
\noindent\textit{Acknowledgments.} OpenAI Codex (GPT 6.2 Sol, Astra, Medium-Ultra) and Anthropic Claude (Opus 5.5 Medium) assisted with drafting, checking references, reasoning, and auditing. The author would like to thank Kaiwen Jin for pointing out an error.

\section{Notation and preliminaries}
\label{sec:notation}

\sectionopening{In this section, we collect the notation for the proof in the order of the Parts.
The local entries give the conventions and point to the definitions of
the indexed quantities. Within a calculation, the scope of a summation
index, local rank, or choice of basis is specified in the corresponding
section. Sections~\ref{sec:datum}, \ref{sec:joint}, and
\ref{sec:common-critical-twist} fix, respectively, the
arithmetic datum, the comparison, and the choices over a coefficient place.}

\subsection{Prime sums and the zero criterion}
\label{sec:notation-prime-sums}

\sectionopening{In this subsection, we record the analytic conventions for
Sections~\ref{sec:criterion} and~\ref{sec:prime-remainder}.}

We write \(\zeta(s)\) for the Riemann zeta function,
\(L(s,\chi)\) for the Dirichlet \(L\)-function of \(\chi\),
\(\gamma_{\rm Euler}\) for Euler's constant, \(\Lambda\) for the von Mangoldt function,
\(\mu\) for the M\"obius function, \(\varphi(q)\) for Euler's
totient, and \(u_+=\max(u,0)\).
For a character modulo \(q\), the progression modulus is
\(q_0=18480q\). In the remainder estimates \(q_0\) is any fixed
positive multiple of \(18480\). The prime \(d\equiv d_0=q_0/2-1
\pmod{q_0}\) defines \(K=\mathbf Q(\sqrt{-d})\) and
\(\chi_d(n)=(n/d)\). In the definitions below the subscript
modulus \(q\) ranges over \(\{1,q_0\}\). Put
\(\lambda_{d,q}(n)=\mathbf1_{n\equiv1\bmod q}\Lambda(n)(1-\chi_d(n))\),
with the congruence condition vacuous for \(q=1\). The remainder and
its subtraction constant are
\[
 \begin{aligned}
 R_{d,q}(x)&=\sum_{n\le x}\lambda_{d,q}(n)(n^{-1}-x^{-1})
                  -\frac{\log x}{\varphi(q)}-c_{d,q},\\
 c_{d,q}&=\frac{D_q(1)-1}{\varphi(q)},\\
 D_q(s)&=\sum_{\chi\bmod q}
   \left(\frac{F'_{\chi\chi_d}(s)}{F_{\chi\chi_d}(s)}
                         -\frac{F'_\chi(s)}{F_\chi(s)}\right).
 \end{aligned}
\]
Here \(F_\chi(s)=L(s,\chi)\) for a nonprincipal character and
\(F_\chi(s)=(s-1)L(s,\chi)\) for a principal character.
For \(\chi\bmod q\) in the sum defining \(D_q\), the product
\(\chi\chi_d\) is taken modulo \(qd\).
Section~\ref{sec:criterion} writes \(D(s)=D_{q_0}(s)\) for this
logarithmic derivative sum and \(\mathcal M_R(s)\) for the Mellin transform
of \(R_{d,q_0}\), with integration weight \(x^{-s}\,dx\).
The criterion assumes zero-free twists on \(\Re s\ge\beta\),
\(|\Im s|\le H_0\), where \(1/2<\beta<1\) and \(H_0>0\), and
\(R_{d,q_0}(x)\le x^{\beta-1}\) for all sufficiently large real \(x\).
It excludes untwisted zeros with \(\Re s>\beta\) and
\(|\Im s|\le H_0\). The twist selection treats the rectangle \(\beta\le\Re s\le1\),
\(|\Im s|\le H_0\). Euler products give nonvanishing for \(\Re s>1\).
In the Introduction, \(\delta>0\) is the tolerance in the remainder
exponent and \(C_{q_0,d,\delta}\) depends on \(q_0,d,\delta\).

\subsection{Scales, ranks, and prescribed levels}
\label{sec:notation-scales}

\sectionopening{In this subsection, we record the scales of Section~\ref{sec:datum}, which determine the ranks
and the original levels before any coefficient place is chosen.}

We use \(m\) for
the number of steps in the lower theta chain and \(r\) for the
rank of the original ambient unitary group.
For the integer scale \(x\) and \(0<\epsilon<1/100\),
Section~\ref{sec:datum} defines
\[
 \begin{gathered}
 m=2^{\lceil\log_2\sqrt x\rceil},\qquad
 r=3^{\lceil\log_3(x^5)\rceil},\qquad k=r+1-m,\\
 s=22\log_2m,\quad Q=2^s,\quad
 L_0=3^{\lfloor\log_3s\rfloor+4},\quad
 q_1=3L_0Qq_0,\quad q_2=(q_0d)^{s+6},\\
 W=x^{\epsilon/6},\qquad x_1=\lfloor e^W\rfloor,\qquad
 b=\min\{b'\in q_1x_1\mathbf Z_{>0}:b'>e^{W^6}\},\\
 a=xb,\quad B=\log b,\quad J_*=\frac{faB^3}{mk},\quad
 \varrho=m/2-1,\quad t=5r.
 \end{gathered}
\]
The initial definite rank is \(k\). The lower spaces \(Y_j\)
have ranks \(r+jt\), and the upper spaces \(X_j\) have ranks
\(n_j=k+j(t+1)\) in Section~\ref{sec:upper-chain-size}.
Section~\ref{sec:period} uses \(\varrho\) hyperbolic pairs
in the Bessel string. The progression construction uses the bound
\(R_{\rm prog}=kt\) and the fixed dimension parameter \(n>66/\epsilon\).

In the orbital estimates, \(n=\varrho+1=m/2\) is the length of the
Bessel basis. In Proposition~\ref{prop:cancel}, \(n\) is the
dimension of the translated span.
The local variables \(s,Q,W\) in the theta and Selmer calculations
are specified in their respective sections. The two level primes are \(p_0^{\rm lev},p_1^{\rm lev}\),
also denoted \(p_0,p_1\) in Section~\ref{sec:datum}.

Their product is \(M\), and \(\pi_M=\bar\alpha\) is the chosen
generator with \(\alpha\bar\alpha=M\). The additional prime \(N\) and
generator \(\pi_N=a_0+ib_0\), with \(\pi_N\bar\pi_N=N\), specify
the smaller source level. The datum \(\mathcal D\) records the weights,
primes and generators of Section~\ref{sec:datum}. The real Bessel scale
at \(v\) is \(\kappa_v=M\Phi_v(\delta_E)/(2\pi i)\), distinct from
the integer weight \(\kappa_\tau\) below.
Local torsion exponents denoted \(M,N\) are specified in
Section~\ref{sec:notation-selmer}.

\subsection{Fields, embeddings, and weights}
\label{sec:notation-fields}

\sectionopening{In this subsection, we specify the fields and embeddings of Section~\ref{sec:datum}
and the weights used in the prime sum of Section~\ref{sec:joint}.}

Let \(\zeta_{q_1}\) be a primitive \(q_1\)-th root of unity, put
\(E=K(\zeta_{q_1})\) and \(F_0=E^+\), and write
\(f=[F_0:\mathbf Q]\) for the number of real places of \(F_0\).
A bar denotes CM conjugation, or complex conjugation on
values under the fixed complex embedding.
Write \(\mathcal O_E\), \(\widehat{\mathcal O}_E\), and
\(\mathbb A_E\) for the integers, finite integral adeles, and
adeles. We write \(\operatorname{Nm}\) for absolute ideal
norm, \(\operatorname{rd}(E)\) for the absolute root
discriminant, and \(\delta_E\) for the chosen skew generator
of the different. At a finite place \(w\), \(k_w\) is the
residue field and \(\varpi_w\) a uniformizer.
The infinite places, read over the fixed embedding of \(K\),
are indexed by \((\mathbf Z/q_1)^\times\). Write \(\Phi_{q_1}\)
for the cyclotomic polynomial, \(\mathfrak D_{E/\mathbf Q}\)
for the different, and \(N_{E/\mathbf Q}\) for the field norm.
The set \(T\) consists of the selected embeddings over the residue
class \(-1\bmod q_0\) in Section~\ref{sec:datum}. The fibre
\(\Phi_a=\{\tau:\tau\equiv1\bmod q_0\}\) and the multiplicities in
Section~\ref{sec:joint} give
\[
 u=|\Phi_a|=3L_0Q,\qquad v=|T|,\qquad
 \kappa_\tau=\begin{cases}b,&\tau\in T,\\a/x_1,&\tau\notin T.\end{cases}
\]
For an integer \(n\), \(\langle-n\rangle\) denotes the class of
\(-n\) in \((\mathbf Z/q_1)^\times\) if it is a unit, and the
identity otherwise. We also use \(\chi_K=\chi_d\) and
\(\Phi_v:E\hookrightarrow\mathbb C\) for a local complex embedding
over a real place.

The unitary Hecke character \(\upsilon\) fixes the Weil
splitting. The norm-one torus is \(T_1=U(1)\), and its two
level characters are \(\chi_0,\chi_M\), constructed in
Section~\ref{sec:marked-level-characters}.

\subsection{Hermitian spaces, automorphic representations, and levels}
\label{sec:notation-automorphic}

\sectionopening{In this subsection, we use the source and the projected target of
Sections~\ref{sec:ramified-source} and~\ref{sec:period}.}

For a unitary group \(H\) over \(F_0\), write
\([H]=H(F_0)\backslash H(\mathbb A_{F_0})\).
For a Hermitian space \(H_0\), \(U(H_0)\) denotes its unitary group,
\(H_0^-\) the space with the opposite form, and \(q(H_0)\) the negative
index at a specified real place. We use \(\perp\) for orthogonal direct
sum, \(\langle a\rangle\) for a line of norm \(a\), and \(I_a=(E^a,1)\)
for a standard positive definite space. Hermitian inner products in
the determinant comparison are conjugate-linear in their first
variable. For matrices \(^{\mathsf T}\) is transpose and \(^*\)
is conjugate transpose, or the adjoint with respect to the indicated
forms. An overline on a finite matrix can also denote reduction, with
the residue ring specified at its occurrence.

Section~\ref{sec:ramified-source} defines the full source
\(X\) at the level \(J_H\), which includes the last-covector
condition above \(N\). The ramified subspace \(S\subseteq X\) is
the intersection of the kernels of the scalar averages \(A_w\)
at \(w\mid\pi_N\). Put \(d_X=\dim X\).
The full Petersson target is \(T_{\rm full}\).
The finite group \(G=G_0\times G_1\) acts by \(\rho\);
its two factors are the reductions at the two level
primes. The projection onto the selected types is \(p=\rho(e_0)\), with
\(e_0=h_0\otimes h_1\) the product of the two central
idempotents. Thus \(V=pT_{\rm full}\) and \(D=\dim V\).

Use \(Y_0\) for the ambient Hermitian space of rank \(r\), \(X_0\)
for its definite block of rank \(k\), and \(R=U\rtimes U(X_0)\)
for the Bessel subgroup in Section~\ref{sec:period}. Its ordered string
has isotropic pairs \(e_a,e_a^*\) and the central vector \(z\) of norm
\(-N\). Here \(e_a^*\) denotes the basis vector paired with \(e_a\).
The compact group \(K_L\) and its congruence kernel \(J_L\) refer
to the integral orbital basis. The source level is
\(J_H=\iota_H^{-1}(J_L)\), and \(J_M\) is the original principal
level at \(M\). Later \(K_Z(\mathfrak a)\) denotes the principal
congruence subgroup at ideal \(\mathfrak a\) on \(Z\) in its fixed
integral bases. Write \(S_M\) for the chosen places dividing \(\pi_M\).

Section~\ref{sec:period-measured-pairing} identifies the measured
Poincar\'e map \(\Pi\). In Section~\ref{sec:period-petersson-blocks},
\(C=C_XR_H(\eta^{-1})\) gives the Euclidean class coordinates
and \(B=(v_G\nu)^{-1/2}\Pi C^{-1}\), where \(\nu\) is the
unipotent covolume.
The symbol \(B\) used for this operator is distinguished by
context from the scalar \(\log b\). The translated maps
\(U_g:S\to V\) use the seesaw normalization, and
\(P_g=c_B U_g\) uses the rational scalar of
Proposition~\ref{prop:shortness-normalization}.
The finite class representatives \(a_c\) and unipotent cutoffs
\(\beta_c\) in Section~\ref{sec:period-quotient-measure} give the quotient integrals.
The finite Haar factors are \(v_R,v_G\). At a real place,
\(J_0,J_1\) are the Gram matrices before and after diagonal scaling by
\(D_v\), \(p_v\) orders an orthonormal basis by sign, and
\(x_v=D_vp_v\). The superscript \(\star\) in \(x_v^\star\) denotes
the stationary point and \(\iota_v(g)=p_vgp_v^{-1}\). The resulting
adelic point for the period is \(g_0=((\iota_v^{-1}(x_v^\star))_v,1)\);
the adelic Gram isometry \(x^{\rm pt}\) has finite coordinate \(x_f\).
The period functional is \(t=(v_G/v_R)t_0\) in
Section~\ref{sec:period-measured-pairing}.

In Section~\ref{sec:marked-chain-composition}, \(T_\xi=L_\xi\)
is the full weighted chain of Section~\ref{sec:seesaw-normalization}
and \(\Theta_\xi=T_\xi|_V\) is its restriction to the projected space, with
\(\xi=1\) the untwisted case.
Write \(\|\cdot\|_{\rm HS}\) for the Hilbert--Schmidt norm,
\(P_U\) for orthogonal projection to a subspace \(U\), and
\(\mathbb E_g\) for averaging with normalized counting
measure on \(G\).

\subsection{Orbital integrals and their moment bounds}
\label{sec:notation-orbital}

\sectionopening{In this subsection, we fix notation for the orbital integrals of Part~\ref{part:orbital}
and the moment bounds of Section~\ref{sec:marked-moments}.}

Write \(B_s=B^*\rho(s)B\), \(s\in G\), for the Petersson blocks
and \(H_R\) for the residue Bessel subgroup. Its representation
\(r_\psi\) on the Euclidean source \(C(X)\) is extended by zero to \(G\).
In Section~\ref{sec:full-orbital-blocks} write
\[
 E_s=B_s-r_\psi(s),\qquad
 D_{\rm id}=|H_R|d_X,\qquad
 D_{\rm main}=\frac{d_X^2}{\mathcal B\,\operatorname{mass}(E,r)} .
\]
The product of the formal degrees is \(\mathcal B\).
The masses \(\mu_j=\operatorname{mass}(E,j)\) of the standard
definite genera are weighted by reciprocal automorphism orders.
The stabilizer intersections \(\Gamma_{\mathcal A}\) at the places in \(\mathcal A\)
use the normalizing bound \(C_{\mathcal A}=b^{-6k|\mathcal A|}D_{\rm main}\).
In Section~\ref{sec:full-orbital-blocks}, \(\mathcal L(\Omega)\)
sums the squared norms of the cut blocks on the
nonzero-discriminant part of \(\Gamma_*\), separately for each
link and over the two classes and \(s\in\Omega\), before division
by \(\nu^2\). For each factor \(G_i\), \(p_i=\rho(h_i)\) is
the projection on the target and \(q_i=1-h_i\) is the
complementary idempotent in the group algebra.

For a finite set of places \(T\) in
Section~\ref{sec:finite-orbital-analysis}, write
\[
 A_T=\prod_{w\in T} E_w,\qquad
 \mathcal R_T=\prod_{w\in T}\mathcal O_{E_w}.
\]
The finite calculations write \(\mathcal R=\mathcal R_T\) locally.
The additive modulus is \(\operatorname{md}\), and \(K_T,J\)
are the local integral and congruence groups in
Section~\ref{sec:finite-bottom-normalization}.
For the Bessel flag \(F_i=\langle e_0,\ldots,e_{i-1}\rangle\), the
ideals \(\mathfrak d_i=c_i\mathcal R_T\) compare the two positions
modulo \(F_i\); their inverse moduli are
\(D_i=\operatorname{md}(c_i)^{-1}\).
Sections~\ref{sec:finite-frequency}--\ref{sec:finite-selected-bottom}
specify the associated Fourier variables, denominator indices, and
stabilizer measure. Section~\ref{sec:marked-residue-projection}
uses the trace of \(r_\psi\) on each finite factor \(G_i\).

In the real calculations, \(h(v,y)=v^*Jy\) in the chosen Gram basis,
and norm-fiber measure uses \(h(y,y)/2\). The support and determinant
variables of the real correlations are specified in
Sections~\ref{sec:real-orbital-phases} and~\ref{sec:real-inverse-depth}.
Section~\ref{sec:real-bottom-mass} writes \(W\) for the full Bessel
vector, \(F\) for its truncation, and \(P(A,B)\) for their quotient
correlation. Here \(d\) denotes the local formal degree.
Section~\ref{sec:real-marked-filling} specifies the real discrete
series, the compact type indices \(I_v\), and the local stationary
vectors. The dimension calculation uses \(d_B=\prod_v|I_v|\) and
\(d_X=d_Cd_B\), with \(d_C\) the source class count.

Section~\ref{sec:marked-full-to-cut} defines the rational
summation set \(\mathcal T(c,s,c')\) for the fixed finite
classes \(c,c'\) and residue element \(s\) in the
Petersson blocks of Section~\ref{sec:period-petersson-blocks}.
The real and finite common-entry pairings are formulated in
Section~\ref{sec:adelic-common-entry-sum}. The pair lattice,
its discriminant, and the parameters of the Schur bound are specified
in Sections~\ref{sec:marked-pair-gain} and~\ref{sec:marked-schur-orbits}.
One set \(t\subseteq G\) of translates is chosen in
Section~\ref{sec:marked-span-proof}. Section~\ref{sec:outline-orbital}
writes \(n=\dim\sum_{g\in t}P_g(S)\). In the orbital integrals
\(n=\varrho+1\) counts the columns of \((e_0,\ldots,e_{\varrho-1},z)\).

\subsection{Theta operators, derivatives, and integral estimates}
\label{sec:notation-theta}

\sectionopening{In this subsection, we use the Weil splitting of Section~\ref{sec:splitting-level-characters}
and the measured normalization of Section~\ref{sec:seesaw-normalization}.}

A tame twist \(\xi\) is a weight-zero character of the norm-one
torus. Its modulus \(Q(\xi)\) is the product of the prime ideals at
one place of each split pair in its support. Thus
\[
 Q(\xi)\cap\overline{Q(\xi)}=\varnothing,\qquad
 q_\xi=\prod_{w\in Q(\xi)}q_w,\qquad (q_\xi,MN)=1,
\]
where in set notation \(Q(\xi)\) denotes the prime divisors of the
modulus. Its character on local units is given at such a place by a
nontrivial residue character \(\xi_w:k_w^\times\to\mathbb C^\times\)
and is trivial on the units at the other places of the torus. A
twist for the full chain has \(\xi_w^j\ne1\) for all \(w\in Q(\xi)\)
and \(1\le j\le m\). The trivial twist has empty support. At a
coefficient prime \(\ell\) we use \(\ell\nmid q_\xi\);
Proposition~\ref{prop:joint-norm-selmer-twist} supplies twists with
these conditions.
The successive spaces in the lower chain are \(Y_j\);
its endpoint \(Y_m\) and the dimension-comparison space \(W_m\)
of Section~\ref{sec:marked-target-dimension} have different ranks.
Write \(K_j^\circ\) for the standard compact, \(J_{\xi,j}\)
for the level at the old and twist places, and
\[
 A_{\xi,j}=[K_j^\circ:J_{\xi,j}] .
\]
The map \(L_{\xi,j}\) is the theta integral at step \(j\)
multiplied by this index. The paired tensors \(F_i(\xi,c)\),
the evaluated tensor \(T_\xi\) of \eqref{eq:seesaw-evaluated-tensor},
and the scalar \(c_{\rm meas}\) appear in
Section~\ref{sec:measured-finite-sum}. The same symbol \(T_\xi\)
denotes the weighted chain \(L_\xi\) in
Section~\ref{sec:marked-chain-composition}.

In Section~\ref{sec:adjacent-unfolding}, \(I_D\) and \(\mathcal C_D\)
are the positive-degree maps and the quotient by nonnegative-degree
maps for the specified grading \(D\). The Heisenberg translation is
\(\rho(w)\), the theta functional is \(\theta\), and
\(C_D(f;w)\) is its coefficient at rational frequency \(w\) against
the positive-degree group. The same letter \(\theta\) denotes
Fontaine's map in the appendix. The integer ring of finite adeles is
\(R_f=\widehat{\mathcal O}_E\) in Section~\ref{sec:marked-neighbor},
where \(J_N\) is the idempotent supported at the places dividing
\(\pi_N\). That section and Section~\ref{sec:marked-N-lift} specify
the shifted lattice and the finite translations in the measured tensor.

The archimedean normalization and the weighted prime remainder are
\[
 \begin{aligned}
 \Omega&=\frac{2\pi}{\operatorname{rd}(E)}
           \exp\!\left(\gamma_{\rm Euler}-\frac{L'}{L}(1,\chi_K)\right),\\
 A&=mk\sum_\tau
 \bigl(\log(\kappa_\tau!)+\kappa_\tau\log(M/\Omega)-S(\kappa_\tau,\tau)\bigr),\\
 \mathcal R&=av\bigl(R_{d,q_0}(x)-R_{d,q_0}(x_1)\bigr)+auR_{d,1}(x_1).
 \end{aligned}
\]
The weighted prime sum \(S(Y,\tau)\) in Section~\ref{sec:joint}
and its local components \(S_p(Y,\tau)\) in
Section~\ref{sec:short-normalization} are
\[
 \begin{aligned}
 S(Y,\tau)&=\sum_{n\le a}\Lambda(n)(1-\chi_K(n))
                     (Y/n-\kappa_{\langle-n\rangle\tau})_+,\\
 S_p(Y,\tau)&=\tfrac12\sum_{j\ge1}(1-\chi_K(p^j))
                     (Y/p^j-\kappa_{\langle-p^j\rangle\tau})_+ .
 \end{aligned}
\]
The latter is used at \(Y=a,\kappa_\tau\). The
exponents and the rational rescaling of the images are
\[
 \begin{aligned}
 g_p&=m\left(\sum_{\tau\in\Phi_a}S_p(a,\tau)
                  +k\sum_\tau S_p(\kappa_\tau,\tau)\right),\\
 c_B&=\prod_{p\in F}p^{-\lfloor g_p-e_p\rfloor}.
 \end{aligned}
\]
Here \(F\) consists of primes below \(x_1\) that are nonsplit in
\(K\), together with \(N\). The local denominator bounds are
\[
 e_p=
 \begin{cases}
 (2faB^2/p)c_p,&p<x_1\text{ nonsplit in }K,\\
 2mfb,&p=N,\\
 0,&\text{otherwise},
 \end{cases}
 \qquad
 c_p=\begin{cases}m,&p\equiv\pm1\bmod m,\\1,&\text{otherwise}.
       \end{cases}
\]
For the upper highest weights use \(\Delta=2r^8\),
\(w_j(\tau)=\kappa_\tau+(j+\tfrac12)\Delta\), and
\(a_\tau=a\mathbf1_{\Phi_a}(\tau)\) in
Section~\ref{sec:marked-column-gain}. Section~\ref{sec:upper-chain-size}
uses \(J(l)=\log(l!)+l\log(M/\Omega)\) for the contribution of a part
of length \(l\).

For Taylor coefficients we use \(\Theta_f\) for the theta section
of a Schwartz datum \(f\), \(Q(\partial)\) for the constant-coefficient
differential operator of a polynomial \(Q\), and the integral model
and translation scales of Section~\ref{sec:marked-integral-sections-ordinary-jets}.
The scalar \(s_{\rm vac}\) is the vacuum projection factor in
\eqref{eq:upper-CM-vacuum}. The Taylor scale is
\(s_{\rm jet}=p_Fa_d/\omega_0\), with
\(a_d=(2\pi/\sqrt d)^{1/2}\), period lattice \(\omega_0\mathcal O_K\),
and \(p_F\) the algebraic factor in
Lemma~\ref{lem:upper-cm-transport}. Its integral normalization is in
Section~\ref{sec:marked-integral-sections-ordinary-jets}.
Section~\ref{sec:marked-both-norm-jets} uses both in the norm tests.
An exceptional prime is an inert \(p<x_1\) with
\(p\not\equiv\pm1\bmod m\) and odd
\(\operatorname{ord}_{q_0}(-p)\le r_p=\lceil\log_p m\rceil\)
(Section~\ref{sec:marked-column-gain}). At such a prime
\(Y_\lambda=p^{A_\lambda}T_{\lambda,r_p}\). The exponents
\(A_\lambda\), Frobenius embeddings \(\iota_s\), root defects,
and colength are specified in Section~\ref{sec:exceptional-hecke-analysis}.

\subsection{Covolumes, determinants, and coefficient places}
\label{sec:notation-determinants}

\sectionopening{In this subsection, we use the comparison of Section~\ref{sec:joint} and the product formula
of Section~\ref{sec:gram-product-formula}.}

With \(J_*=faB^3/(mk)\) the bound for the sums over coefficient primes in
Part~\ref{part:lower}, set
\[
 E_0=faB^3,\qquad L=faB^4,\qquad
 I_0=10\log(1+2b^2)faB^3 .
\]
A coefficient place is \(\kappa=(\ell,\iota)\), with
\(\iota:\mathbf C\simeq\overline{\mathbf Q}_\ell\) and
\(\bar\kappa=(\ell,\iota\circ\operatorname{conj})\).
We use \(v_\kappa(\ell)=1\) and \(|z|_\kappa=\ell^{-v_\kappa(z)}\).
In the covolume and height bounds, \(\sigma\) ranges over
\(\operatorname{Aut}(\mathbb C)\) and acts on the coordinates.
For a finite coefficient DVR \(O\), its uniformizer
valuation is \(v_O\); normalized length is
\(\operatorname{length}_O/v_O(\ell)\).
For a tuple \(y\) of \(n\) independent images,
\[
 \operatorname{Cov}(y)\le c
 \quad\Longleftrightarrow\quad
 \det\operatorname{Gram}(\sigma\Theta_1 y)\le e^c
 \quad(\sigma\in\operatorname{Aut}(\mathbb C)).
\]
The bound \(\operatorname{Ind}(y)\le\delta\) means there are
bounds \(a_\ell\ge0\), independent of the embedding above \(\ell\),
with \(v_\kappa(\det U_\kappa)\le a_\ell\) for the integral
clearing matrices and \(\sum_{\ell\in U}a_\ell\log\ell\le\delta\)
for every finite prime set \(U\). At a given \(\kappa\) the
same matrix clears the tuple for all allowed twists of modulus
prime to \(\ell\). The twist \(\xi_\kappa\)
is chosen in Section~\ref{sec:common-critical-twist}.
For a completed basis \(w=(y_1,\ldots,y_n,z_1,\ldots,z_{D-n})\),
\[
 G_\xi=\det\operatorname{Gram}(\Theta_\xi w),\qquad
 Z_\xi=\det(\Theta_\xi^*\Theta_\xi),\qquad
 G_\xi=Z_\xi\det\operatorname{Gram}(w).
\]
After the untwisted lift, the auxiliary vectors \(z_i\)
have squared norms at most \(\exp(b_{\rm short})\) at every
complex embedding, with \(b_{\rm short}=2x_1^2b^2\). The scalar
\(\gamma_{\xi_\kappa}=c_{\xi_\kappa}/c_1\) clears them in the lift
\(\Theta_{\xi_\kappa}\) at both \(\kappa\) and \(\bar\kappa\).
The logarithmic height bounds on \(Z_1/Z_{\xi_\kappa}\) and
\((\gamma_{\xi_\kappa}\overline{\gamma_{\xi_\kappa}})^{-1}\)
are \(j_Z\) and \(j_\gamma\), respectively. The height convention
is that of Section~\ref{sec:gram-product-formula}, and the
basis completion is in Section~\ref{sec:marked-chain-completion-bounds}.

\subsection{Cocycles, specializations, and Galois congruences}
\label{sec:notation-selmer}

\sectionopening{In this subsection, we record the notation for the normalized lengths of
Section~\ref{sec:selmer-relations} and the trace congruences of
Section~\ref{sec:trace-determinant-specialization}.}

Write \(G_E\) for the absolute Galois group,
\(\chi_{\rm cyc}\) for the cyclotomic character, and
\(L^\dagger=\operatorname{Hom}_O(L,O)(1)\) for the
Tate dual of a coefficient lattice \(L\).
We use \(H^1_f\) for the crystalline local condition
and \(Z_f(L)\) for the integral cocycles mapping to it.
The local class of lattices \(\mathcal Q_v\) and
cup-product colength \(\Delta_v(L)\) are those of
Section~\ref{sec:relative-crystalline-colength}.
At depth \(M\), \(\mathcal B=\mathcal O/\ell^M\mathcal O\)
is the coefficient quotient, \(\lambda(C)\) is the normalized
length of a cocycle module, and \(\phi_j,\gamma_j\) are the
auxiliary Frobenius and inertia evaluations. The slopes \(s_j\)
impose \(z(\gamma_j)=s_jz(\phi_j)\) on the submodule \(C_0\).
The auxiliary character order divides \(\ell^N\). In
Section~\ref{sec:selmer-polynomials}, \(R_X\) is the polynomial
matrix of the twisted action and \(A(X)\) records the cocycle,
inertia, and other linear equations. The length is computed
from the maximal minors \(P_i(X)\) of
\[
 A^+(X)=\begin{pmatrix}A(X)\\ \ell^M I_q\end{pmatrix},
\]
where \(q\) is the column count of \(A(X)\) in that section.
The evaluations \(X_w=\zeta_w^j-1\) use
\(\zeta_w^{\ell^H}=1\) as in Section~\ref{sec:torsion-specialization}.
The bound \(\ell^H\) is fixed before the auxiliary level
\(\ell^N\) in Section~\ref{sec:simultaneous-choice-order}.

In Section~\ref{sec:joint-twist}, \(\Theta_Q\) is a group-ring
element at auxiliary set \(Q\), \(q_w=|k_w|\), and
\(\mathcal E\) is the fixed list of exponents, and the
Euler factor is the Laurent polynomial
\(E_w(A)=\prod_{e\in\mathcal E}(1-Aq_w^{-e})(1-A^{-1}q_w^{-e})\). The sign and
contragredient conventions for the attached representation \(\rho\)
are in Section~\ref{sec:cofinite-spherical-galois}.
The trace of conjugation is normalized using the polarization
exponent \(d_{\rm pol}\) of Section~\ref{sec:main-row-budgets}.

The ordered Kummer construction uses \(F_\nu=E(\mu_{\ell^\nu})\),
\(G_\nu=\ker(\chi_{\rm cyc}\bmod\ell^\nu)\), and
\(A_\nu=\mathbf Z/\ell^\nu\mathbf Z\).
The original evaluation set in the cyclotomic kernel is \(P_0\);
the new primes \(w_j,\bar w_j\) and the maps \(\psi_j\) are
specified in Sections~\ref{sec:fresh-kummer-maps} and
\ref{sec:triangle-selection}. The matrices \((u_x(g)-1)_x\)
define the relative relations on cocycles in
Section~\ref{sec:relative-crystalline-colength}.
For the full raised family \(Y\) and its common coefficient DVR \(O\),
put \(\mathbb T_O=\mathbb T\cap O^Y\) and \(S=\mathbb T_O[X_j:j\in J]\),
where \(\mathbb T\) is the integral trace algebra of
Section~\ref{sec:integral-raised-traces}.
The scalar congruences generate \(\mathcal J\subset S\).
Their trace annihilator \(\mathfrak a\subset S[G_E]\) is a two-sided ideal,
as in Section~\ref{sec:marked-ribet-trace-determinant}.
The operators \(\mathcal N_j,\partial_j\) at the raised places
and the highest squarefree coefficient refer to the original test space
of Section~\ref{sec:top-coefficient-clearing}.

\subsection{Period comparisons, critical values, and descent}
\label{sec:notation-comparison}

\sectionopening{In this subsection, we record the spaces and operators for the period comparison of
Section~\ref{sec:normalized-period-exchange} and the lattice completion
of Section~\ref{sec:marked-chain-completion-bounds}.}

The two comparison factors of Part~\ref{part:lower} use the original
descent pair \(W_j,V_j\) of Section~\ref{sec:marked-comparison-lattice},
with
\[
 V_j=V''\perp V^\circ,\qquad
 d_+=r+9m+1,\qquad d_-=t-9m-1 .
\]
Here \(d_+,d_-\) are the ranks of \(V'',V^\circ\).
The choices of classes and the tuple space \(\mathcal H\) used by the
gluing operator are defined in
Section~\ref{sec:original-comparison-tests}. The gluing operator is
\(G_\eta\) at twist \(\eta\), and \(n_A,h_-\) are the class
and coefficient constants in \(U_\eta=n_Ah_-^2G_\eta G_\eta^*\).
On the projected source the lift is \(u_\eta\), its image is
\(Z_\eta\), and the Petersson projection is \(\Pi_\eta\).
This subspace notation \(Z_\eta\) is local to
Section~\ref{sec:original-compression-bounds}. For the \(j\)-th
single lift write \(Z_j(\eta)=\det(u_{j,\eta}^*u_{j,\eta})\).
In the chain at twist \(\xi\) its factor is \(Z_j(\xi^j)\).

For the comparison of two compositions, \(U_1,U_2\) are the source
and target Hermitian spaces, \(A,A'\) are the intermediate spaces,
and \(W=U_1\perp U_2\). The two compositions with the chosen point
arrays \(X,Y\) are \(T_X=D_{2,X_2}L_{1,X_1}\) and
\(T'_Y=L'_{2,Y_2}D'_{1,Y_1}\) in
Section~\ref{sec:marked-square-analytic-pairing}.
The real Fock variables and finite self-dual test are specified in
Sections~\ref{sec:marked-real-slice-period} and
\ref{sec:marked-primitive-fourier-transfer}.
The signatures, infinity weights \(\mu_\tau\), and the character
\(\gamma=\nu^{-1}\upsilon^{N/2}\) for exchanging two periods are
fixed in Section~\ref{sec:marked-partner-period}.
Here \(N=\dim W\) in that formula, and \(\upsilon\) there is
the restriction of the splitting Hecke character to the norm-one ideles.

A critical pair \(p=(\nu,z)\) consists of an automorphic
character \(\nu\) of the norm-one torus and an integer \(1\le z\le t\).
Its infinity weights are \(-c_v\), with \(z\le c_v\le4b\).
At a chosen half \(w\in S_M\) of an old split place the restriction
to units factors through a character of \(k_w^\times\), and at the
other finite places of the torus it is unramified. The pair is
ramified if all these residue characters are nontrivial.
Base change includes both halves of each pair.
The associated twist is \(\eta\)
of exact support \(Q\), with \(q_Q=\operatorname{Nm}(Q)\).
In Section~\ref{sec:common-critical-twist}, the ratio for
\(\epsilon\in\{1,-1\}\) is
\[
 R_{p,\epsilon}(\eta)=
 \frac{|L^{S_M\cup Q}(z,((\nu\eta^\epsilon)^{\rm bc})^{-1})|}
      {|L^{S_M}(z,(\nu^{\rm bc})^{-1})|}.
\]
In these \(L\)-functions, the exponent \(-1\) denotes the inverse
character after base change. A place set in the superscript on
\(L\) omits the corresponding Euler factors and their conjugates.
The Taylor map \(\mathcal T\), coefficient product \(B_\tau\),
and four-block vector \(\phi_\eta\) are those of
Sections~\ref{sec:marked-twelve-sharp-coefficient} and
\ref{sec:v2-four-block-main-numerator}.

In the common period, \(K_\zeta^\circ\) is the class sum of the
complementary Weyl kernels paired with the tuple \(\zeta\);
the coefficient formulas and approximation bounds
\eqref{eq:display-cp1}--\eqref{eq:display-cp6} are in
Section~\ref{sec:common-period-injectivity}.
At deep principal level put
\(b_0=\lfloor b\rfloor\), \(\mathfrak b=(\pi_M^{b_0})\), and
\(R_b=\mathcal O_E/\mathfrak b\). The spaces
\(\mathcal H_j\), base identifications \(j_i\), and arrow operators
are specified in Section~\ref{sec:deeper-periods}.
Section~\ref{sec:marked-chain-composition} gives the chain relations
\eqref{eq:chain-gram-first-comparison}--\eqref{eq:chain-terminal-gram}. Proposition~\ref{prop:marked-completion-vectors} gives the
vectors for completing the lattice.

\subsection{Local period rings and filtration orders}
\label{sec:notation-crystalline}

\sectionopening{In this subsection, we fix the notation for the period rings of
Appendix~\ref{sec:crystalline-coordinate-appendix} and the filtration
orders of Appendix~\ref{sec:crystalline-basis-appendix}. These enter
the local bounds of Section~\ref{sec:relative-crystalline-colength}.}

In the appendices \(p\) denotes the coefficient residue characteristic,
\(S=W(k)[[z]]\) is the chart ring, \(\theta\) is
Fontaine's map, and \(\varphi\) is Witt Frobenius.
The chart symbol \(E=p+z^e u\), with \(e\) the absolute ramification
index of the chart field and \(u\in S^\times\), is a generator of \(\ker\theta\), with Frobenius translates
\(E_a=\varphi^a(E)\) and tame model \(E=p-z^e\). The Witt ring is
\(A=A_{\inf}=W(\mathcal O_{\mathbb C_p^\flat})\), with chart ideal
\(I=(p,z)A\) and retraction \(F_0:A\to S\).
The bounded period ring is \(B^+=A_{\max}[1/p]\), where
\(A_{\max}\) is the \(p\)-adic completion of \(A[E/p]\).
The ring \(S_B\) is the closure of \(S[1/p]\) in the topology
with neighborhoods \(p^n A_{\max}\) of zero.
The de Rham filtration is \(\mathrm{Fil}^j=t^j B_{\rm dR}^+\).
For a compatible system \(\epsilon=(1,\zeta_p,\zeta_{p^2},\ldots)\)
of primitive roots of unity, put \(\mu=[\epsilon]-1\),
\(\tau=\varphi^{-1}(\mu)\), and \(\xi=\mu/\tau\).
Then \(t=\log(1+\mu)\), and the divided action uses \(\Delta_z=z\tau\)
in Sections~\ref{sec:v2-crystalline-frame-action}
and~\ref{sec:local-tame-depth}. The modules \(M_T,N_T\) over
\(A_{\inf},S\) use the standard orders of
Section~\ref{sec:local-standard-orders}.

\Needspace{11\baselineskip}
The indices \(i,j\) and their Frobenius translates label the chart
idempotents. The ideal \(I_l(i)\) consists of the multipliers that put
the component-\(i\) parameters in the level-\(l\) standard-order image
(Section~\ref{sec:local-standard-orders}). The ideal \(J_r(j;R)\) is
the ideal of scalars lifting component-\(j\) parameters at order \(r\),
plus \((R)\) (Section~\ref{sec:v2-crystalline-lift-ideals}).
The weight gaps and stopping indices are in
Sections~\ref{sec:v2-crystalline-stops}, \ref{sec:v2-crystalline-deep-constants},
and~\ref{sec:off-level-coefficient-defect}.

\section{Outlines of proofs}
\label{sec:proof-outlines}

\sectionopening{In this section, we trace the proof of Theorem~\ref{thm:rh}
through Parts~\ref{part:parameters}--\ref{part:analytic}.
Section~\ref{sec:joint} reduces the prime-remainder bound to two
estimates on the same translated basis. The two estimates yield
Proposition~\ref{prop:cancel}, from which Section~\ref{sec:prime-remainder}
deduces the bound required by the analytic criterion of
Section~\ref{sec:criterion}.}

\subsection{Part I. Parameters and the common family}
\label{sec:outline-parameters}

\sectionopening{In this subsection, we explain how Part~\ref{part:parameters}
uses the arithmetic datum of Section~\ref{sec:datum} to express the
prime remainder in the lattice comparison of Section~\ref{sec:joint}.}

Fix \(q_0\) divisible by \(18480\), a prime
\(d\equiv q_0/2-1\pmod{q_0}\), and \(0<\epsilon<1/100\).
The comparison applies at every sufficiently large integer \(x\)
for these fixed choices. Let \(a=xb\) and \(B=\log b\), as in
Section~\ref{sec:notation-scales}. The weights are \(a/x\) on \(T\)
and \(a/x_1\) outside \(T\). The set \(\Phi_a\) is the fibre of
\(1\bmod q_0\), and \(T\) is a subset of the fibre of \(-1\bmod q_0\).
Since \(q_0\) and \(q_1\) have the same prime divisors, multiplication
by \(-n\) sends the first fibre bijectively to the second whenever
\(n\equiv1\bmod q_0\). For other units the image misses \(T\);
for nonunits Section~\ref{sec:joint} uses the identity multiplier.
Thus, with \(u=|\Phi_a|\) and \(v=|T|\), the weighted sum gives
\[
 \frac{\mathcal R}{av}
   =R_{d,q_0}(x)-R_{d,q_0}(x_1)
                +\frac uv R_{d,1}(x_1),\qquad
 \frac uv\le B,\qquad f=\varphi(q_0)u .
\]
Section~\ref{sec:joint} also defines the common real normalization
\(A\ge0\) and the errors \(L,I_0\). The term \(-m\mathcal R\) will appear
in the upper covolume bound.

Section~\ref{sec:ramified-source} fixes \(X\) and cuts out
\(S\subseteq X\) by the scalar averages at \(N\).
Section~\ref{sec:period} constructs \(\Pi\) on \(X\) with the prescribed
real point, finite level, and quotient Haar measures. The finite-type
idempotent \(e_0\) of Section~\ref{sec:upper-good-projector} gives
\(p=\rho(e_0)\) and \(V=\operatorname{im}p\), with \(D=\dim V\).
These data specify the source, target and terminal lifts in
Propositions~\ref{prop:marked-upper-estimates} and~\ref{prop:marked-determinant}.
The simultaneous construction at large \(x\) is checked in
Section~\ref{sec:common-filling-assembly}.

\subsection{Part II. Orbital estimates and the translated span}
\label{sec:outline-orbital}

\sectionopening{In this subsection, we explain how Part~\ref{part:orbital}
uses the Petersson blocks of Proposition~\ref{prop:marked-petersson-blocks}
to choose translates of \(S\) in Section~\ref{sec:marked-span-proof}.}

Section~\ref{sec:full-orbital-blocks} compares the full Petersson
blocks with the truncated orbital integrals. The finite and real
correlation bounds are in Sections~\ref{sec:finite-orbital-analysis}
and~\ref{sec:real-orbital-analysis}. Section~\ref{sec:real-stationary-point}
bounds the length at the same stationary point used in the measured
period. Section~\ref{sec:adelic-orbital-assembly} combines the local
bounds with the definite masses. The truncation errors and the
remaining pair sums are controlled in
Sections~\ref{sec:orbital-pair-counts} and~\ref{sec:orbital-cancellation-analysis}.
The dimension comparison \eqref{eq:marked-tower-dimension} supplies
the normalization in the mixed moment.

Write \(B_{\rm Pet}\) for the Petersson-normalized map denoted \(B\)
in Section~\ref{sec:marked-moments}, to distinguish it from \(\log b\).
On the Euclidean class coordinates \(C(X)\), the moment bounds are
\[
 \|pB_{\rm Pet} P_U\|_{\rm HS}^2\ge\dim U-5b^{-4}d_X,\qquad
 D\,\mathbb E_g\|B_{\rm Pet}^*\rho(g)B_{\rm Pet}\|_{\rm HS}^2
                 \le (1+9b^{-4})d_X^2 ,
 \qquad U\subseteq C(X).
\]
Propositions~\ref{prop:orbital-diagonal-moment}
and~\ref{prop:orbital-mixed-moment} give these bounds, respectively. In Section~\ref{sec:marked-span-proof},
the mixed moment gives a translate whose projection to the span
already selected has small Hilbert--Schmidt norm. The diagonal
bound then controls the dimension of its preimage in \(C(X)\).
Iteration decreases the codimension. Restriction to \(C(S)\) loses
at most \(|t|\dim(X/S)\) dimensions for a set \(t\) of translates.

The comparison of the two normalizations in
Section~\ref{sec:marked-span-proof} gives
\[
 B_{\rm Pet}=(v_G\nu)^{-1/2}\Pi C^{-1},\qquad
 P_g=c_B|G|c_{\rm ss}^{-1}(v_G\nu)^{1/2}
                         p\rho(g)B_{\rm Pet} C|_S .
\]
Here \(v_G\) is the level volume, \(\nu\) is the unipotent covolume,
and \(c_{\rm ss}\ne0\) is the scalar obtained by putting
Proposition~\ref{prop:measured-seesaw} in integral form in
Section~\ref{sec:seesaw-integral-tensor}. The rational factor
\(c_B\ne0\) is given in Section~\ref{sec:outline-upper}. After restriction to \(C(S)\), the span chosen with
\(B_{\rm Pet}\) is the span of the same \(P_g(S)\) used for the integral
estimates. More precisely,
\[
 |t|\dim S\le5D\log(1+2b^2),\qquad
 n=\dim\sum_{g\in t}P_g(S)\ge(1-4/b^2)D .
\]

\subsection{Part III. Theta identities and the upper covolume bound}
\label{sec:outline-upper}

\sectionopening{In this subsection, we explain how Part~\ref{part:upper}
uses the measured seesaw of Section~\ref{sec:seesaw-normalization}
to bound the norms and finite index on the span of
Section~\ref{sec:marked-span-proof}.}

Sections~\ref{sec:adjacent-unfolding} and~\ref{sec:even-stage-coefficients}
unfold the consecutive lifts on the specified levels. The measured
finite sum of Section~\ref{sec:measured-finite-sum} gives the tensor,
and Section~\ref{sec:seesaw-integral-tensor} gives its integral
expansion. Proposition~\ref{prop:measured-seesaw} uses the same \(\Pi\)
and fixed real point from Section~\ref{sec:period-real-point}.
The comparison at \(N\), the ordinary
derivatives, and the exceptional-prime Hecke operators are in
Sections~\ref{sec:marked-N-lift}, \ref{sec:local-theta-derivatives},
and~\ref{sec:exceptional-hecke-analysis}, respectively.

Section~\ref{sec:short-normalization} defines the prime-power
contributions \(g_p\), and Section~\ref{sec:local-denominator-sum}
bounds the local errors \(e_p\) by
\(\sum_p e_p\log p\le E_0\). The scalar of
Proposition~\ref{prop:shortness-normalization} is
\[
 c_B=\prod_{p\in F}p^{-\lfloor g_p-e_p\rfloor}.
\]
At each coefficient prime \(p\), the scalar makes the cleared
images integral for all allowed twists of modulus prime to \(p\).
The prime-power identity in Section~\ref{sec:short-normalization}
accounts for the matching terms in the conjugate norm bound.
Proposition~\ref{prop:translated-norm-and-index} applies it to the
same translated maps as Proposition~\ref{prop:marked-orbital-span}.

Choose a basis of the selected span from the generators of
Proposition~\ref{prop:marked-upper-estimates}:
\[
 y=(P_{g_i}\phi_i)_{i=1}^n,\qquad g_i\in t,\quad\phi_i\in\mathcal S.
\]
The untwisted norm bound and Hadamard's inequality give
the covolume bound at every complex embedding. For the index, the
clearing matrices are combined by translate, giving a determinant
valuation at most \(|t|d_Xe_p\) above \(p\).
Using \(d_X\le2\dim S\), \(\sum_p e_p\log p\le E_0\), and
the bound on \(|t|\dim S\), we obtain
\[
 \operatorname{Cov}(y)\le n(A-m\mathcal R+L),\qquad
 \operatorname{Ind}(y)\le DI_0 .
\]
In the index bound the matrix at a coefficient place is chosen
before selecting the twist, as required in
Section~\ref{sec:gram-product-formula}. The basis \(y\) is fixed
before any coefficient place or twist is chosen.

\subsection{Part IV. Selmer bounds and trace congruences}
\label{sec:outline-congruences}

\sectionopening{In this subsection, we explain how Part~\ref{part:congruences}
gives the Selmer and norm conditions for the common twist of
Section~\ref{sec:common-critical-twist} and the trace comparison
used there.}

The original Hecke eigensystems give the Galois representations in
Section~\ref{sec:cofinite-spherical-galois}, with local lengths estimated
in Section~\ref{sec:main-row-budgets}. Lemma~\ref{lem:selmer-relations}
bounds the normalized length of \(C/C_0\), where \(C_0\)
imposes the auxiliary local relations on \(C\).
Let \(P_{j_0}\) be a maximal minor of \(A^+\) attaining the
minimum in Lemma~\ref{lem:sharp-selmer} at a test point. For
each integral specialization \(x\), the equations defining
\(C(x)\) give
\[
 \lambda(C(x))=\min_j v(P_j(x))\le v(P_{j_0}(x)).
\]

For a coefficient place \(\kappa=(\ell,\iota)\),
Section~\ref{sec:simultaneous-choice-order} fixes the original Selmer
data and signed norm systems, their separate truncation depths,
and a uniform bound on the number of auxiliary primes.
Lemma~\ref{lem:common-torsion-specialization} then fixes \(H\)
from the bounds and tolerances. Next one chooses the auxiliary
level \(\ell^{N_{\rm tor}}\), the primes, and the test polynomials.
One tuple \((\zeta_w)_w\) with \(\zeta_w^{\ell^H}=1\) controls
the valuations at \((\zeta_w^j-1)_w\) for all \(1\le j\le m\).
Section~\ref{sec:common-critical-twist} applies this construction
with support \(Q\) away from the primes over \(\ell MN\), where
\(N\) in \(\ell MN\) is the chosen level prime. The resulting
\(\xi_\kappa^j\) have exact support \(Q\) for \(1\le j\le m\). The line and rank-two norm formulas in
Sections~\ref{sec:marked-norm-line} and~\ref{sec:marked-norm-determinant}
apply to the signed systems at the respective ranks.

For the trace estimate,
Section~\ref{sec:relative-crystalline-colength} uses
Appendices~\ref{sec:tame-cup-appendix} and~\ref{sec:crystalline-colength-appendix}
on the original local class and its Tate-dual pairings.
For the selected twist, Section~\ref{sec:triangle-selection} fixes the
ordered list of Kummer tests before the trace congruence depth \(\nu\).
The original Hecke systems determine a sufficient \(\nu\) in
Section~\ref{sec:top-coefficient-clearing}. At that depth one chooses
the primes at which to raise the level and forms the trace algebra.
These primes are additional to the support \(Q\) of the twist.
The coefficient-extension formula
\eqref{eq:specialized-selmer-coefficient-extension}
carries the normalized Selmer length bound to the common coefficient DVR
used for the raised systems.
Section~\ref{sec:marked-ribet-trace-determinant} constructs one
trace multiplier for the raised systems. Its highest squarefree
coefficient gives the clearing operator on integral tests orthogonal
to every false-lift space \(F_w\), as in
Proposition~\ref{prop:clearing-on-original-tests}.

\subsection{Part V. Periods and the lower determinant estimate}
\label{sec:outline-lower}

\sectionopening{In this subsection, we explain how Part~\ref{part:lower}
bounds the determinant ratio and completes the fixed basis \(y\) for the
product formula of Section~\ref{sec:gram-product-formula}.}

Section~\ref{sec:original-comparison-tests} defines the operator on
the original class spaces. Sections~\ref{sec:normalized-period-exchange}
and~\ref{sec:gluing-three-exchanges} match its periods, including their
measure factors. Sections~\ref{sec:period-reproduction}
and~\ref{sec:common-period-injectivity} establish the reproduction
and injectivity used to compress the operator in
Section~\ref{sec:original-compression-bounds}.
With the arithmetic estimates of Part~\ref{part:congruences},
Section~\ref{sec:common-critical-twist} uses the same \(\xi_\kappa\)
for the necessary norm and Selmer tests at each coefficient place.

The first-arrow Gram identity of Section~\ref{sec:marked-first-arrow-rallis}
starts the injectivity argument of Section~\ref{sec:marked-chain-composition}.
The first square identity \eqref{eq:chain-gram-first-comparison}
and the propagation formula \eqref{eq:chain-gram-propagation} prove
injectivity before division by the real period factors. Proposition~\ref{prop:chain-determinant-ratio}
computes the chain ratio from the factors
\(Z_j(1)/Z_j(\xi^j)\) and the triangle factors. Proposition~\ref{prop:chain-determinant-height} bounds its height.
Using the deep periods of Section~\ref{sec:deeper-periods},
Proposition~\ref{prop:marked-completion-vectors} supplies a spanning set
from which to choose the \(s=D-n\) additional vectors.
In the notation of
Section~\ref{sec:notation-determinants},
\[
 j_Z\le D(b+L-A),\qquad j_\gamma\le b,\qquad
 b_{\rm short}=2x_1^2b^2,\qquad s\le4D/b^2 .
\]
Here \(j_Z\) bounds the height of \(Z_1/Z_{\xi_\kappa}\),
\(j_\gamma\) bounds the height of
\((\gamma_{\xi_\kappa}\overline{\gamma_{\xi_\kappa}})^{-1}\),
and \(b_{\rm short}\) bounds the logarithmic squared norms of the
added vectors after the untwisted lift.

Let \(w\) be the completed basis. The matrices at \(\kappa\) and
\(\bar\kappa\) clear the original \(y\) in the same lift
\(\Theta_{\xi_\kappa}\); the additional vectors use the common scalar
\(\gamma_{\xi_\kappa}\). The two integral families give the mixed
Gram bound for \(G_{\xi_\kappa}\) in
Section~\ref{sec:gram-product-formula}. Since
\(G_\xi=Z_\xi\det\operatorname{Gram}(w)\), the ratio \(Z_1/Z_{\xi_\kappa}\)
transfers the bound to \(G_1\). At infinity, exterior-power
Cauchy--Schwarz bounds \(|\sigma G_1|\) by the product of the
exterior norms of the untwisted families at \(\sigma\) and
\(\mathrm{conj}\circ\sigma\circ\mathrm{conj}\).
In each family project the added vectors orthogonally off the span
of the first \(n\) vectors and apply Hadamard's inequality to the
projected vectors. If \(\operatorname{Cov}(y)\le c\), this gives
\(|\sigma G_1|\le e^{c+sb_{\rm short}}\).
With \(\operatorname{Ind}(y)\le\delta\), the product formula gives
\[
 0\le c+2\delta+s(b_{\rm short}+j_\gamma)+j_Z .
\]
Section~\ref{sec:marked-determinant-assembly} uses
\(s(b_{\rm short}+j_\gamma)\le Db\) and \(b\le L\) to obtain
\(D(A-3L)\le c+2\delta\). Substituting the bounds of
Section~\ref{sec:outline-upper} and using \(n\le D\) and \(A+L\ge0\)
gives
\[
 nm\mathcal R\le (n-D)(A+L)+D(4L+2I_0)
               \le D(4L+2I_0).
\]
For large \(x\), \(n\ge D/2\), \(I_0\le30L\), and \(B\ge128\).
This is the cancellation in Proposition~\ref{prop:cancel}, giving
\(m\mathcal R\le8L+4I_0\le faB^5\).

\subsection{Part VI. The prime remainder and the analytic conclusion}
\label{sec:outline-analytic}

\sectionopening{In this subsection, we explain how Part~\ref{part:analytic}
uses Proposition~\ref{prop:cancel} to obtain the one-sided bound in
Section~\ref{sec:prime-remainder} and apply
Section~\ref{sec:criterion} to Theorem~\ref{thm:rh}.}

Substitute \(m\mathcal R\le faB^5\) into the identity of
Section~\ref{sec:joint}. From \(u/v\le B\) and \(f=\varphi(q_0)u\),
\[
 R_{d,q_0}(x)\le
 \frac{\varphi(q_0)B^6}{m}
 +|R_{d,q_0}(x_1)|+B|R_{d,1}(x_1)|.
\]
The proof of Proposition~\ref{prop:real-prime-remainder} bounds both
remainders at \(x_1\) by the fixed-modulus prime number theorem.
The bounds \(m\ge\sqrt{x}\) and \(B\ll x^\epsilon\) give
\(R_{d,q_0}(x)\le C_{q_0,d,\epsilon}x^{-1/2+6\epsilon}\) at sufficiently
large integers. Integrating the \(O_{d,q_0}(x^{-1})\) derivative
bound between consecutive integers extends this estimate to real
\(x\).

For \(\chi\bmod q\), induce to \(\widetilde\chi\bmod q_0\)
with \(q_0=18480q\). The extra Euler factors are nonzero in
\(\Re s>0\), so a zero \(\rho\) of \(L(s,\chi)\) with
\(1/2<\Re\rho<1\) is also a zero of \(L(s,\widetilde\chi)\).
Choose \(1/2<\beta<\Re\rho\) and \(H_0>|\Im\rho|\).
Proposition~\ref{prop:zero-free-auxiliary-twist} uses a truncated
Mellin integral of the M\"obius-weighted theta series to select
one prime \(d\) for which all twists by \(\chi_d\) of the
characters modulo \(q_0\) are zero-free on
\(\beta\le\Re s\le1\), \(|\Im s|\le H_0\).
Fix this \(d\), choose \(0<\epsilon<1/100\) with
\(6\epsilon<\beta-1/2\), and then let \(x\) grow in the
estimate above. Absorbing the constant gives
\[
 h(x)=x^{\beta-1}-R_{d,q_0}(x)\ge0
\]
for all sufficiently large real \(x\).
The real-abscissa argument of Lemma~\ref{lem:positive-mellin-criterion}
shows that the Mellin integral of \(h\) over that tail converges on
\(\Re s>\beta\). Since the twists have no zeros in the rectangle,
the zero of \(L(s,\widetilde\chi)\) at \(\rho\) would give a
nonzero residue there. Section~\ref{sec:prime-remainder} treats
\(\Re s\ge1\) separately and uses the primitive functional
equation to exclude zeros in \(0<\Re s<1/2\).

\part{Arithmetic parameters and automorphic spaces}\label{part:parameters}

\section{Arithmetic parameters}\label{sec:datum}

\sectionopening{In this section, we fix the arithmetic parameters before choosing any
coefficient place. The selected embeddings separate the two weights in
Section~\ref{sec:joint}, where their multiplicities determine the two
prime sums.}

Fix \(q_0\ge1\) divisible by \(18480\), a prime
\(d\equiv d_0=q_0/2-1\bmod q_0\), and \(0<\epsilon<1/100\).
The scale \(x\) is an integer. The parameters \(q_0,d,\epsilon\)
remain fixed as \(x\) grows. When a constant is written as
depending on \(d\) alone, one may take its maximum over eligible
\(q_0\) because \(q_0\mid 2(d+1)\).

\begin{lemma}[The auxiliary residue class]\label{lem:auxiliary-residue-class}
The class \(d_0\) is a unit of order two modulo \(q_0\) and
\(-1\notin\langle d_0\rangle\). We have \(d\equiv7\bmod8\), and
every prime dividing \(q_0\) splits in \(\mathbb Q(\sqrt{-d})\).
\end{lemma}
\begin{proof}
Since \(16\mid q_0\), one has
\[
 (q_0/2-1)^2\equiv1\bmod q_0,\qquad
 d_0\not\equiv \pm1\bmod q_0,\qquad
 \frac{d+1}{4}=(2j+1)\frac{q_0}{8}
\]
for some \(j\ge0\). Since \(16\mid q_0\), we have
\(d\equiv7\bmod8\) and \(\operatorname{rad}(q_0)\mid q_0/8\).
In particular \(d\nmid q_0\). The discriminant of
\(K=\mathbb Q(\sqrt{-d})\) is \(-d\) and
\(\mathcal O_K=\mathbb Z[(1+\sqrt{-d})/2]\). At every
prime dividing \(q_0\), the minimal polynomial
\(X^2-X+(d+1)/4\) reduces to \(X(X-1)\).
The roots are distinct also at \(2\), giving the asserted splitting.
For odd \(p\ne d\), quadratic reciprocity and \(d\equiv3\bmod4\)
give \(({-d}/p)=(p/d)\). At \(2\), \((2/d)=1\) because
\(d\equiv7\bmod8\), while the value at \(d\) is zero. Thus
the character of the quadratic field at every positive integer is
\(\chi_K(n)=(n/d)\).
\end{proof}

For \(m\) the least power of \(2\) at least \(\sqrt x\),
\(r\) the least power of \(3\) at least \(x^5\), put
\[
 s=22\log_2 m,\quad Q=2^s,\quad
 L_0=3^{\lfloor\log_3 s\rfloor+4},\quad
 q_1=3L_0Qq_0,\quad R_{\rm prog}=(r+1-m)5r.
\]
Let \(E=K(\zeta_{q_1})\), where the quadratic field
\(K=\mathbb Q(\sqrt{-d})\) ramifies
at \(d\), whereas the cyclotomic field of conductor
\(q_1\) is unramified there. Thus \(K\) intersects it
in \(\mathbb Q\). By fixing one embedding of \(K\), the
infinite places of the compositum are indexed by
\((\mathbb Z/q_1)^\times\). Since \(q_1/q_0=3L_0Q\) has no prime
factor outside \(q_0\), reduction to \((\mathbb Z/q_0)^\times\)
is surjective with kernel of order \(3L_0Q\). In particular
\[
 f=\varphi(q_0)\,3L_0Q,\qquad
 f\asymp_{q_0}q_1\asymp_{q_0}x^{11}\log x,\qquad
 f,q_1\le x^{12}
\]
for all sufficiently large \(x\) with \(q_0\) fixed.
The primes dividing \(q_1\) are split in \(K\) by
Lemma~\ref{lem:auxiliary-residue-class}.

Write \(W=x^{\epsilon/6}\) and choose a fixed integer
\(n>66/\epsilon\). Since \(Q=m^{22}\), one has
\(x^{11}\le Q<2^{22}x^{11}\), so \(W\le Q\le W^n\) for large \(x\).
Also \(81s/3<L_0\le81s\), so \(L_0\asymp\log x\).
With \(n\) fixed these give
\[
 (2L_0)^n\le W,\qquad 3L_0n2^n\le W^3,\qquad
 2L_0(R_{\rm prog}+1)<Q,
\]
where the last inequality uses \(R_{\rm prog}=O(x^{10})\).
Take the largest power of \(2\), denoted \(H\),
with \(H^n\le Q/W\). Then \(H\le W\),
\(Q<2^nH^nW\), and \((2L_0H)^n\le Q\).
In the cube \(\{0,\ldots,H-1\}^n\) a sphere of
constant squared norm contains at least
\(H^n/(nH^2)\) points, since there are at most \(nH^2\)
norms. Read these points as digits in base \(2L_0H\),
giving a set \(\mathcal E\subseteq[0,Q)\).
If \((a_1+a_2)e_2\equiv a_1e_1+a_2e_3\bmod Q\)
for \(e_j\in\mathcal E\) and positive integers \(a_1,a_2\)
with \(a_1+a_2\le2L_0\), every digit sum is at most
\(2L_0(H-1)<2L_0H\). Both integers are below
\((2L_0H)^n\le Q\). They are therefore equal digitwise.
The vector of \(e_2\) is a strict convex combination of the
vectors of \(e_1,e_3\), which have the same squared norm
as \(e_2\). Thus all three vectors coincide.

Put \(g=1+3L_0q_0\) and \(v_2(q_0)=v_0\ge4\).
The prime divisors of \(q_1\) are those of \(q_0\), so
\(g\) is a unit modulo \(q_1\).
For every positive \(h\),
\(v_2(g^h-1)=v_0+v_2(h)\), by factoring for odd \(h\)
and repeated difference of squares for even \(h\).
Thus \(g\) has order \(Q\) modulo \(2^{s+v_0}\).
The places \(T=\{-g^e\bmod q_1:e\in\mathcal E\}\) are
distinct and have residue \(-1\bmod q_0\).
We prove that for any \(\tau_1,\tau_2,\tau_3\in T\),
any prime \(\ell\le R_{\rm prog}\) and \(0<i<j<L_0\), the two
congruences
\[
 \tau_2/\tau_1\equiv\pm\ell^i,\qquad
 \tau_3/\tau_1\equiv\pm\ell^j
 \pmod{(q_1)_{\ell'}}
\]
cannot both hold, where \((q_1)_{\ell'}\) is the prime-to-\(\ell\)
part of \(q_1\). Suppose both hold and write \(\tau_k=-g^{e_k}\)
with \(e_k\in\mathcal E\).

For \(\ell=2\), squaring the first ratio gives
\(4^i\equiv1\bmod 3L_0\). Indeed \(3L_0\) divides the
odd part of \(q_1\) and all the \(\tau_k\) are \(-1\)
modulo \(3L_0\). Since \(v_3(4^i-1)=1+v_3(i)\)
by factoring successive cubes, \(L_0\mid i\), a
contradiction. For odd \(\ell\), the modulus \((q_1)_{\ell'}\)
contains \(2^{s+v_0}\). Raise the two ratios to
\(2j,2i\), respectively, and compare their exponents
of \(g\) in the cyclic group of order \(Q\) at \(2\). This gives
\[
 2j e_2\equiv 2(j-i)e_1+2i e_3\pmod Q.
\]
The equal-norm construction forces equality of
the exponents. The first squared ratio then gives
\(\ell^{2i}=1\bmod 2^{s+1}\). Factoring this difference
of powers gives
\[
 v_2(\ell^{2i}-1)
  =v_2(\ell-1)+v_2(\ell+1)+v_2(2i)-1 .
\]
Since \(\min(v_2(\ell-1),v_2(\ell+1))=1\), we obtain
\[
 Q\le 2^{v_2(\ell^{2i}-1)}
     \le (\ell+1)2i\le(R_{\rm prog}+1)2L_0<Q .
\]
This contradiction excludes the two ratios.

Finally
\(3L_0Q<3L_0 n2^n WH^2 |T|\le W^6|T|\).
For \(x_1=\lfloor e^W\rfloor\), take \(b\) to be
the least multiple of \(q_1x_1\) strictly greater
than \(e^{W^6}\). Since \(e^{W^6}<b\le e^{W^6}+q_1x_1\)
and \(\log(q_1x_1)=O(\log x+W)\), one has
\(W^6<B=\log b\le W^6+1\) eventually.
Thus \(u/v\le B\) for \(u=3L_0Q\), \(v=|T|\) and \(B=\log b\).

\subsection{The chosen primes and their generators}\label{sec:marked-level-prime-selection}

\sectionopening{In this subsection, we use Siegel--Walfisz to choose two principal-level primes with a
common ideal class and a third prime for the additional level
condition in Section~\ref{sec:ramified-source}. The generators also
fix the additive characters in Section~\ref{sec:seesaw-normalization}.}

Set \(q_2=(q_0d)^{s+6}\) and \(B=\log b\).
Then \(q_1d\mid q_2\). At \(2\) one adds \(s\) to
\(v_2(q_0)\) in \(q_1d\), and at \(3\) one adds
\(1+\log_3 L_0=\lfloor\log_3 s\rfloor+5\le s+5\) to
\(v_3(q_0)\). The exponent \(s+6\) in \(q_2\) covers both
additions because \(2,3\mid q_0\). At all other primes
\(q_1d\) and \(q_0d\) have the same valuation.
Also \(q_2\ge Q=m^{22}>8m\) for large \(x\). In particular a prime
\(1\bmod q_2^2\) splits completely in \(E\subseteq
\mathbb Q(\zeta_{q_1d})\). The preceding bounds give
\[
 \sqrt x\le m<2\sqrt x,\qquad x^5\le r<3x^5,\qquad
 s\le32\log x,\qquad W-1\le \log x_1\le W.
\]
In particular, eventually \(x^{10^6}\le x_1\),
\(x_1^{10^6}\le b\), and \(B^{100}\le x\).
One has \((\log x_1)^6\le B\) by the choice of \(b\).
The bound \(B^{100}\le x\) uses \(100\epsilon<1\).
Also \(\log B\to\infty\), so any fixed lower bound on
\(\log B\) depending on \(q_0,d\) can be imposed in the same cutoff.

\begin{proposition}[The three level primes]\label{prop:marked-level-primes}
Fix \(q_0,d,\epsilon\). For every sufficiently large integer \(x\), there
are distinct primes \(p_0,p_1,N\) satisfying
\[
 \begin{gathered}
 p_0\equiv p_1\equiv N\equiv1\pmod{q_2^2},\qquad
 b^{10}\le p_i\le b^{11}\quad(i=0,1),\\
 U\le\log N\le U+\sqrt B.
 \end{gathered}
\]
Here \(M=p_0p_1\), the weights \(\kappa_\tau\) are those of
Section~\ref{sec:notation-fields}, and
\[
 U=m\,f^{-1}\sum_\tau\log\kappa_\tau
   +(m-1)\log\!\left(\frac{M\operatorname{rd}(E)}{2\pi}\right)
   -\tfrac{15}{2}B .
\]
There are prime ideals \(\mathfrak p_i\mid p_i\) of \(K\) and
generators \(\alpha\in\mathcal O_K\), \(\pi_N\in\mathbb Z[i]\) with
\[
 (\alpha)=\mathfrak p_0\overline{\mathfrak p_1},\qquad
 \alpha\bar\alpha=p_0p_1,\qquad \pi_N\bar\pi_N=N .
\]
The chosen generator at \(M=p_0p_1\) is \(\pi_M=\bar\alpha\).
\end{proposition}

\begin{proof}

Since \(\log(q_2^2)=O_{q_0,d}(\log x)\), there is a fixed
\(A_{\rm SW}\) with \(q_2^2\le W^{A_{\rm SW}}\).
For any \(u\ge W\), apply the Siegel--Walfisz theorem
\cite[Corollary~11.21, (11.36)]{MV} to the prime-weighted count
\(\vartheta(\,\cdot\,;q_2^2,1)\) at \(2e^u\) and \(e^u\).
It gives
\[
 \vartheta(2e^u;q_2^2,1)-\vartheta(e^u;q_2^2,1)
 =\frac{e^u}{\varphi(q_2^2)}+
       O_{A_{\rm SW}}\!\left(e^u e^{-c\sqrt u}\right).
\]
Here \(q_2^2\le u^{A_{\rm SW}}\), and the relative error
is \(O_{A_{\rm SW}}(u^{A_{\rm SW}}e^{-c\sqrt u})\).
This tends uniformly to zero for \(u\ge W\) as \(x\to\infty\).
Thus there is a prime \(1\bmod q_2^2\) with logarithm in
\((u,u+\log 2]\subseteq[u,u+W]\).

Let \(h_K\) be the class number of \(K\).
For \(i=0,\ldots,h_K\), apply this result in
\([10B+2iW,10B+(2i+1)W]\). These intervals are disjoint and lie
in \([10B,11B]\) for large \(x\), since \(W/B\to0\).
The selected primes split in \(K\) because
their residue modulo \(d\) is \(1\).
Choose a prime ideal of \(K\) above each. Two
distinct primes \(p_0,p_1\) give ideals
\(\mathfrak p_0,\mathfrak p_1\) of the same class. In particular,
\[
 b^{10}\le p_i\le b^{11},\qquad p_i\le b p_{1-i}\quad(i=0,1).
\]
Since \(\mathfrak p_1\overline{\mathfrak p_1}=(p_1)\), the class of
\(\mathfrak p_0\overline{\mathfrak p_1}\) is
\([\mathfrak p_0][\mathfrak p_1]^{-1}=1\).
Choose an integral generator \(\alpha\). As both prime ideals have
residue degree one, taking absolute norms gives
\(\alpha\bar\alpha=p_0p_1=M\).
Use \(\pi_M=\bar\alpha\) as the chosen generator, with \(\bar\pi_M=\alpha\).

Give each place \(\tau\in T\) weight \(\kappa_\tau=b\)
and each remaining infinite place weight
\(\kappa_\tau=xb/x_1\). These are integers in \([1,b]\).
Let \(f=[E:\mathbb Q]/2\), and set
\[
 U=m\,f^{-1}\sum_\tau\log\kappa_\tau
       +(m-1)\log\!\left(\frac{M\,\operatorname{rd}(E)}{2\pi}\right)
       -\tfrac{15}{2}B .
\]
Put \(C=\mathbb Q(\zeta_{q_1})\). The discriminants of \(K,C\)
are coprime, so products of their integral bases form an integral
basis of \(E\), as verified in Section~\ref{sec:different-additive-character-scales}.
The Kronecker product of the two trace matrices and the cyclotomic
discriminant formula give
\[
 |\operatorname{disc}(E)|=d^{f}|\operatorname{disc}(C)|^2,\qquad
 \operatorname{rd}(E)
 =\sqrt d\,q_1\prod_{p\mid q_1}p^{-1/(p-1)}.
\]
In particular \(1\le\operatorname{rd}(E)\le dq_1\le b\) eventually.
With \(b^{20}\le M\le b^{22}\) this gives
\[
 19B\le \log(M\operatorname{rd}(E)/(2\pi))\le23B,
 \qquad mB\le U\le24mB .
\]
The second pair uses \(m\ge2\) and
\(0\le f^{-1}\sum_\tau\log\kappa_\tau\le B\).
Apply the interval of primes at this \(U\) to obtain
\(N\equiv1\bmod q_2^2\), with
\(\log N\in[U,U+\sqrt B]\), since \(W\le\sqrt B\). For large \(x\),
\[
 \log N\ge U\ge mB>11B\ge\max(\log p_0,\log p_1).
\]
As \(N\equiv1\bmod4\), it is a sum \(a_0^2+b_0^2\)
of two integer squares. The prescribed norm generator
is \(a_0+ib_0\in\mathcal O_E\), where \(i\in E\)
because \(4\mid q_1\).
The primes, their generators, and the weights are kept
in the same datum, denoted \(\mathcal D\).
This construction is for every sufficiently large integer
\(x\) with \(q_0,d,\epsilon\) fixed.
\end{proof}

\subsection{The different and the additive character scales}
\label{sec:different-additive-character-scales}

\sectionopening{In this subsection, we choose the different generator that fixes the scales in Section~\ref{sec:period-real-point}.}

For the additive character choose the skew generator
\(\delta_E\) of the different as follows. Since the two
fields have coprime discriminants,
\(\mathcal O_E=\mathcal O_K[\zeta_{q_1}]\): at each rational prime,
tensoring one local maximal order with the other unramified
order gives a product of maximal orders.
Put \(F_0=E^+\). In
\(\operatorname{Gal}(E/\mathbb Q)
=\{\pm1\}\times(\mathbb Z/q_1)^\times\), complex conjugation
is \(c=(-1,-1)\), and \(\operatorname{Gal}(E/F_0)=\{1,c\}\).
A place over \(p\) splits in \(E/F_0\) precisely when its
decomposition group over \(\mathbb Q\) omits \(c\).
For \(p\mid q_1\), including \(p=2\), this group fixes \(K\) by
Lemma~\ref{lem:auxiliary-residue-class}. At \(p=d\) its
projection to the cyclotomic factor is generated by
\(d\bmod q_1\); it omits \(-1\) already modulo \(q_0\) by the
same lemma. Thus every place above \(dq_1\) splits in \(E/F_0\).
For \(p\nmid dq_1\), both factors are unramified and the decomposition
group is generated by \((\chi_K(p),p\bmod q_1)\). The places
above \(p\) are inert in \(E/F_0\) exactly when
\(\chi_K(p)=-1\) and \(p^j\equiv-1\pmod{q_1}\) for some odd
\(j>0\). In particular \(E/F_0\) is unramified at every finite place.

The equality \(f=[F_0:\mathbb Q]=\varphi(q_1)\) follows from
\([E:\mathbb Q]=2\varphi(q_1)\), and \(f\) is even since
\(q_1>2\). The relative different
over \(K\) is generated by \(\Phi_{q_1}'(\zeta_{q_1})\),
where \(\Phi_{q_1}\) is the cyclotomic polynomial.
Writing \(X^{q_1}-1=\Phi_{q_1}(X)P(X)\) in \(\mathbf Z[X]\) gives
\(q_1\zeta_{q_1}^{q_1-1}=\Phi'_{q_1}(\zeta_{q_1})P(\zeta_{q_1})\).
Thus \(q_1\in\mathfrak d_{E/K}\), and transitivity gives
\(q_1d\in\sqrt{-d}\,\mathfrak d_{E/K}=\mathfrak d_{E/\mathbf Q}\).
This supplies the rational integer used in the dual-lattice estimate of
Lemma~\ref{lem:ramified-short-generators}.
In the product of root differences for the derivative, write each other
primitive root as \(\zeta_{q_1}^{2a+1}\), since \(q_1\)
is even.

Multiplication of
\(\zeta_{q_1}-\zeta_{q_1}^{2a+1}\) by \(i\zeta_{q_1}^{-a-1}\) makes the
factor real. The product of these unit multiples, times
\(\sqrt{-d}\), is a skew generator of the absolute
different by transitivity. All its allowed root-of-unity
phases have the same archimedean moduli.
Set \(u=q_1/q_0\), whose prime divisors divide \(q_0\).
A root of unity has \(u\)-th power of order \(q_0\) exactly when
its own order is \(q_1\): for each prime \(\ell\mid q_0\),
the valuation of its order must be \(v_\ell(q_0)+v_\ell(u)\),
and no prime outside \(q_0\) can divide its order.
Hence \(\Phi_{q_1}(X)=\Phi_{q_0}(X^u)\), and differentiation gives
\[
 |\Phi_v(\delta_E)|
 =\sqrt d\,u\,\left|\Phi'_{q_0}\bigl(\Phi_v(\zeta_{q_1})^u\bigr)\right|
 \asymp_{d,q_0}q_1 .
\]
The derivative is nonzero at every primitive \(q_0\)-th root,
and there are finitely many such roots at fixed \(q_0\).
Since \(q_1\asymp_{q_0} x^{11}\log x\), we have
\(x^{10}\le|\Phi_v(\delta_E)|\le x^{12}\le b\) at every infinite
place for sufficiently large \(x\). As
\(b^{20}\le M\le b^{22}\), the real character scales
\(\kappa_v=M\Phi_v(\delta_E)/(2\pi i)\) satisfy
\(1\le|\kappa_v|\le b^{23}\) for large \(x\).
Here \(\Phi_v\) denotes the embedding, as distinct
from the cyclotomic polynomial. The interval of primes above
also gives \(N\le b^{25m}\).

\subsection{The splitting and level characters}
\label{sec:splitting-level-characters}

\sectionopening{In this subsection, we construct the characters used in Sections~\ref{sec:upper-chain-size}
and~\ref{sec:seesaw-normalization}, with their common finite and
archimedean normalizations.}

Put \(F=F_0\), \(q=q_1\), \(\zeta=\zeta_{q_1}\), and
\(\widehat O=\widehat{\mathcal O}_E\).

The roots of unity in \(E\) are precisely \(\mu_q\).
Indeed \(E\subset\mathbf Q(\zeta_{dq})\) by the quadratic
cyclotomic description of \(K\) of discriminant \(-d\).
If a root has order \(h\), then \(\varphi(\operatorname{lcm}(h,dq))
=\varphi(dq)\). Since \(4\mid dq\), any increase in this
modulus strictly increases its totient, so \(h\mid dq\).
The ramification index at \(d\) in \(E\) is \(2\), since
the cyclotomic field of conductor \(q\) is unramified there.
As \(d-1>2\), \(E\) cannot contain \(\mathbf Q(\zeta_d)\).
Hence \(h\mid q\). Now set \(\eta=\zeta^{q/4}(\zeta-1)\).
The cyclotomic identities
\(\operatorname{N}_{E/K}(\zeta^{q/4})=1\) and
\(\operatorname{N}_{E/K}(\zeta-1)=\Phi_q(1)=1\) show that
\(\eta\) is a unit of relative norm one. The second identity
uses that \(q\) is divisible by both \(2\) and \(3\).
Direct conjugation gives \(\eta/\bar\eta=\zeta\). Thus every
root \(\xi\) in \(E\) can be written \(\xi=\delta/\bar\delta\)
with \(\delta\in\mathcal O_E^\times\) and
\(\prod_v\Phi_v(\delta)=1\), taking a power of \(\eta\).

If a finite place \(w\) of \(E\) is fixed by conjugation, choose a
primitive third root of unity when \(w\nmid3\), and a primitive fourth
root when \(w\mid3\). Its reduction differs from that of its inverse,
since its order is prime to the residue characteristic. Thus conjugation
acts nontrivially on the residue field at \(w\), and the inertia subgroup
of \(E/E^+\) is trivial there. The other places occur in split pairs. The norm on local units
is surjective at each finite place. For an inert place,
the residue-field norm is surjective and the map on the
successive principal-unit quotients is the residue trace.
Successive lifting and completeness therefore give the
unit norm surjectivity, including residue characteristic two.
For a split place use the norm \((x,y)\mapsto xy\).

Let \(\mathcal U=E_\infty^\times\widehat O^\times\) and
\(\mathcal R=E^\times\{y\bar y:y\in\mathbb A_E^\times\}\).
On \(\mathcal U\), put \(S(x)=\prod_v x_v/|x_v|\).
For \(x=a y\bar y\in\mathcal U\cap\mathcal R\), the finite
valuations of \(a\) are invariant under conjugation, since the
valuations of \(y\bar y\) are. Thus \(\bar a/a\) is a unit with
all archimedean moduli one. By Kronecker's theorem it is a
root of unity. Choose \(\delta\) as above with
\(\delta/\bar\delta=\bar a/a\), so \(A=a\delta\in F\).
The finite part of \(A y\bar y\) is a fixed integral unit
and is a norm at every finite place by unit norm
surjectivity. Hence \(A\) itself is everywhere a finite
local norm. Global Artin reciprocity for \(E/F\) \cite[Ch.~V, Theorem~5.3(a)]{MilneCFT} gives
\(\prod_{v\mid\infty}\operatorname{sgn}\Phi_v(A)=1\).
Since \(\prod_v\Phi_v(\delta)=1\) and \(y\bar y\) is positive
at infinity, it follows that \(S(x)=1\).

Accordingly \(S\) defines a character on \(\mathcal U\mathcal R\)
with value one on \(\mathcal R\). The quotient of the idele
group by \(\mathcal U\mathcal R\) is finite, as it is a quotient
of the ordinary ideal class group of \(E\). Extend this
character to the whole idele group one cyclic quotient at a
time, choosing the required roots in the circle. Because
\(\mathcal U\) is open and \(S\) is continuous, the extended
character \(\upsilon\) is continuous. It is a unitary Hecke
character, unramified at every finite place, with
\[
  \upsilon(\bar x)=\upsilon(x)^{-1},\qquad
  \upsilon_\infty((z_v)_v)=\prod_v z_v/|z_v|,
  \qquad
  \upsilon|_{\mathbb A_F^\times}=\chi_{E/F}.
\]
For the last identity, the restriction is trivial on
\(F^\times\operatorname{N}_{E/F}(\mathbb A_E^\times)\),
which has index two by global reciprocity
\cite[Ch.~V, Theorem~5.3(b), p.~179]{MilneCFT}.
The sign at one real place is nontrivial. In particular the value of
\(\upsilon\) on the uniformizer at the chosen inert place is \(-1\):
one can compute it on a uniformizer of the fixed local
field, and any other uniformizer differs by a unit.

Write \(\upsilon_E\) for this full idelic character when it
is necessary to distinguish it from its restriction to the norm-one
torus. The values on finite ideles are algebraic.
For a finite idele \(x\), choose \(h>0\) such that the ideal of
\(x^h\) has a generator \(a\in E^\times\). Then
\[
 \upsilon(x)^h=\upsilon_\infty(a)^{-1}
                 =\prod_v|\Phi_v(a)|/\Phi_v(a),
\]
since \(x^h a^{-1}\) is an integral unit at every finite
place. The right side is algebraic. If a product of powers of
\(\upsilon\) and the idele norm has integral infinity exponents,
its finite component takes values in one number field. The finitely many ideal-class representatives
give algebraic generators, and the values on principal ideals
are monomials in the conjugates of their generators. These facts give the coefficient
fields for the algebraic twists of the cuspidal factors.

\subsection{The two old-level characters on the same torus}
\label{sec:marked-level-characters}

\sectionopening{In this subsection, we extend the level characters of
the norm-one torus for the weights from Section~\ref{sec:datum}.
Their local components enter the theta identities in
Section~\ref{sec:seesaw-normalization}.}

Let \(T_1=U(1)\) for \(E/F\), let \(T_{1,f}^0\) be its integral
finite compact, and let \(S_M\) be the original places
dividing \(\pi_M\), one from each level pair. Each
\(w\in S_M\) has residue cardinality \(p_i\) for one of
the two chosen primes. Since \(q_2^2\mid p_i-1\), choose
a character \(a_w:k_w^\times\to\mathbb C^\times\)
of exact order \(q_2^2\) and put
\[
      \omega_w=a_w^{q_2}\quad (w\in S_M),
      \qquad \omega_w=1\quad (w\notin S_M).
\]
The residue characters here are read on the chosen member
of each split pair. The value on its conjugate member of
a norm-one idele is determined by the chosen coordinate.
Each \(\omega_w\) on \(S_M\) has exact order \(q_2\).
For \(z\in\mu_{q_1}\), one has
\(\omega_w(\operatorname{red}_w z)=a_w(\operatorname{red}_w(z^{q_2}))=1\)
because \(q_1\mid q_2\). For every integer \(j\) with
\(0<|j|<q_2\), \(\omega_w^j\) has order
\(q_2/\gcd(q_2,|j|)>1\). Its lift to the local units has conductor
exponent one, since it factors through the residue group.

For the infinite place \(v\) with index \(\tau\), put
\(w_0(v)=\kappa_\tau+r^8\), where \(\kappa_\tau\)
is the integer weight \(b\) or \(xb/x_1\).
The two values \(b\) and \(xb/x_1\) of \(\kappa_\tau\) are multiples of \(q_1\) because \(q_1x_1\mid b\).
On \(T_1(F_\infty)T_{1,f}^0\) use the continuous characters
\[
    u\longmapsto\prod_v u_v^{-w_0(v)},
       \qquad
    u\longmapsto \prod_{w\in S_M}\omega_w(\operatorname{red}_w u_w).
\]
The rational points in this subgroup are exactly the
integral norm-one units. They are \(\mu_{q_1}\):
conjugation makes all complex absolute values one, and
the root-of-unity calculation above applies. For
\(\gamma\in\mu_{q_1}\) the first formula has value
\[
          \prod_v\Phi_v(\gamma)^{-w_0(v)}
          =N_{E/K}(\gamma)^{-r^8}=1 .
\]
Here the embeddings \(\Phi_v\) are all those over the
fixed embedding of \(K\), and
\(N_{E/K}(\zeta)=\Phi_{q_1}(0)=1\).
The second formula is also one on these units by the
local order calculation.

Both characters therefore extend trivially on
\(T_1(F)\) to the subgroup \(T_1(F)T_1(F_\infty)T_{1,f}^0\)
of \(T_1(\mathbb A_F)\). To show that the quotient is finite,
note that every finite place of \(E/F\) is split or unramified
by Section~\ref{sec:different-additive-character-scales}. At a split place,
a given norm-one pair \((u,u^{-1})\) is the image of \((u,1)\)
under \(x\mapsto x/\bar x\). At an unramified nonsplit place,
local Hilbert 90 gives a preimage. Dividing it by a power of
a fixed-field uniformizer makes its valuation zero without
changing its image. These choices are units at every place
where the norm-one idele is integral, so they give a finite
idele preimage. As \(E^\times\) and \(\widehat{\mathcal O}_E^\times\)
map into \(T_1(F)\) and \(T_{1,f}^0\), respectively, the map on
finite ideles factors to a surjection from \(\mathrm{Cl}(E)\)
onto the asserted quotient.

Extending across each finite cyclic quotient by a
root in the unit circle gives continuous global
characters \(\chi_0,\chi_M\), respectively. They have
weight \(-w_0,0\) and tame parts \(1,\omega\) on the
original field and original compact. For any fixed
weight and tame part, the possible extensions form
a finite set by the same quotient calculation.

Their finite adelic values lie in a number field.
On representatives for the finite quotient, a fixed
power of each value is the product of a weight monomial
on a rational norm-one element and finite roots of
unity. This makes the finitely many representative
values algebraic. On all other finite ideles, the
additional factors are such monomials and roots
of unity from the fixed tame characters. The normal
closure of \(E\) and those finitely many algebraic
generators therefore give one coefficient field.

Use \(\upsilon_1=\upsilon|_{T_1(\mathbb A_F)}\) and
\(\Delta=2r^8\). For every \(j\) the fixed character
\[
      \mu_j^0=\chi_0(\upsilon_1^{-\Delta}\chi_M)^j
\]
has weight \(-(\kappa_\tau+(2j+1)r^8)\) at the infinite place indexed by \(\tau\)
and tame part \(\omega^j\). For \(1\le |e|\le 8m<q_2\),
\(\omega_w^e\) is nontrivial exactly on \(S_M\).
This proves the required old support for the single lifts,
both signs of the three exchanges, the triangle at the original level
quotients and the reverse-grid products. These choices
are made after fixing the two level primes and their
chosen members, and before any additional modulus,
coefficient embedding, or period test.
\section{The lattice comparison on a common automorphic family}\label{sec:joint}

\sectionopening{In this section, we translate the prime remainder of
Section~\ref{sec:criterion} into a statement about one family of
lifts. Hadamard's inequality combines the estimates of
Propositions~\ref{prop:marked-upper-estimates}
and~\ref{prop:marked-determinant}.}

Write \(a=xb\), \(B=\log b\), \(L=fa B^4\), and
\(I_0=10\log(1+2b^2)\,faB^3\).
For the remainder \(R_{d,q}\) of Section~\ref{sec:notation-prime-sums},
put \(v=|T|\), \(u=3L_0Q\), and
\[
 \mathcal R=av\bigl(R_{d,q_0}(x)-R_{d,q_0}(x_1)\bigr)
               +auR_{d,1}(x_1).
\]
The preceding construction supplies \(u/v\le B\) and
\(f=\varphi(q_0)u\). In fact \(u=|\Phi_a|\) for
\(\Phi_a=\{\tau:\tau\equiv1\bmod q_0\}\).

Write \(\chi_K(n)=(n/d)\), the quadratic character
of \(K\), and let \(\Lambda\) denote the von Mangoldt
function. For a variable \(Y\) put
\[
 S(Y,\tau)=\sum_{n\le a}\Lambda(n)(1-\chi_K(n))
          (Y/n-\kappa_{\langle-n\rangle\tau})_+ .
\]
The parameter \(\langle-n\rangle\) is the class of \(-n\) when
coprime to \(q_1\) and the identity otherwise.
When \(\tau\) ranges over \(\Phi_a\), the product parameter
is in \(T\) exactly \(v\) times if \(n\equiv1\bmod q_0\),
and never otherwise.

Indeed \(\Phi_a\) is the fibre of \(1\) over \(q_0\)
and \(T\) lies in the fibre of \(-1\).
The moduli \(q_0,q_1\) have the same prime divisors.
If \(n\equiv1\bmod q_0\), multiplication by \(-n\) is
a bijection between these two full fibres. Any other unit
class has disjoint image from \(T\), and a nonunit uses
the prescribed identity multiplier. Thus
\[
 \#\{\tau\in\Phi_a:\langle-n\rangle\tau\in T\}
                     =v\mathbf1_{n\equiv1\bmod q_0}.
\]
The weights \(a/x_1\) off \(T\) and \(a/x=b\) on \(T\)
give the exact identity
\[
 \begin{aligned}
 \frac1a\sum_{\tau\in\Phi_a}S(a,\tau)
 &=
 u\sum_{n\le x_1}\Lambda(n)(1-\chi_K(n))(n^{-1}-x_1^{-1})\\
 &\quad+
 v\sum_{n\equiv1\bmod q_0}\Lambda(n)(1-\chi_K(n))
        \bigl((n^{-1}-x^{-1})_+-(n^{-1}-x_1^{-1})_+\bigr).
 \end{aligned}
\]
For large \(x\), one has \(x<x_1\le a\). The correction weight is
\[
 (n^{-1}-x^{-1})_+-(n^{-1}-x_1^{-1})_+
 =\begin{cases}
 x_1^{-1}-x^{-1},& n\le x,\\
 x_1^{-1}-n^{-1},& x<n<x_1,\\
 0,& n\ge x_1.
 \end{cases}
\]
Both sums therefore lie in the original range \(n\le a\). Applying
the definitions of the two remainders gives
\[
 \sum_{\tau\in\Phi_a} S(a,\tau)
   =\mathcal R +a\{u(\log x_1+c_{d,1})
                         +v(\log x-\log x_1)/\varphi(q_0)\},
\]
where \(c_{d,1}\) is the subtraction constant of the
modulus-one remainder. In the correction the constant
of the modulus-\(q_0\) remainder cancels. Since
\[
 \sum_\tau\log\kappa_\tau
      =f\log(a/x_1)-v\log(x/x_1),
\]
the term in braces times \(a\) is the common expected
value \(ua(\log a-f^{-1}\sum_\tau\log\kappa_\tau+c_{d,1})\).
This proves the identity relating the chain's total
remainder to the two progression sums.

For the main amount let \(k=r+1-m\), and set
\[
 \Omega=\frac{2\pi}{\operatorname{rd}(E)}
    \exp\!\left(\gamma_{\rm Euler}
                   -\frac{L'}{L}(1,\chi_K)\right),
 \qquad
 A=mk\sum_\tau\bigl(\log(\kappa_\tau!)
           +\kappa_\tau\log(M/\Omega)-S(\kappa_\tau,\tau)\bigr).
\]
It is nonnegative for large \(x\). Indeed
\(S(\kappa_\tau,\tau)\le 2\kappa_\tau\sum_{n\le a}\Lambda(n)/n
\le 5\kappa_\tau B\) eventually, by partial summation
from the prime number theorem and \(\log a\le 2B\).
The fixed \(d\) and \(\operatorname{rd}(E)\ge1\) give
\(\log(M/\Omega)\ge19B\) eventually, since \(M\ge b^{20}\).
Every summand is therefore nonnegative.

The spaces below are the ramified principal-level source \(S\)
and the projected principal-level target \(V\) of the construction.
Write \(D=\dim_{\mathbb C}V>0\), as shown in
Section~\ref{sec:marked-span-proof}. A finite group \(G\) acts
on \(V\), and the comparison uses translates \(P_g\)
of one normalized Poincar\'e map \(P:S\to V\).
The family of lifts \(\Theta_\xi\) and the
integral lattices at each coefficient place are on this
same \(V\). The trivial twist has symbol \(1\).

For a tuple \(y\) in \(V\), write
\[
 \operatorname{Cov}(y)\le c
 \quad\Longleftrightarrow\quad
 \det\operatorname{Gram}(\sigma\Theta_1 y)\le e^c
 \quad(\sigma\in\operatorname{Aut}(\mathbb C)).
\]
Here \(\sigma\) acts on the coordinates, and the Gram matrix
uses the rational class pairing.

Write \(\operatorname{Ind}(y)\le\delta\) if there are
\(a_\ell\ge0\), independent of the choice of embedding over
\(\ell\), such that
\[
 \sum_{\ell\in U}a_\ell\log\ell\le\delta
 \quad\text{for every finite set of primes }U.
\]
For an ordered tuple \(y=(y_1,\ldots,y_n)\), at each
\(\kappa\mid\ell\) there must exist a matrix
\(U_\kappa\in M_n(\mathcal O_\kappa)\) such that
\[
 v_\kappa(\det U_\kappa)\le a_\ell,\qquad
 (\Theta_\xi y)U_\kappa\text{ has integral columns at }\kappa
 \quad(\ell\nmid q_\xi).
\]
Here \(v_\kappa(\ell)=1\) and
\(\mathcal O_\kappa=\{z\in\mathbb C:v_\kappa(z)\ge0\}\).
The columns of \(\Theta_\xi y\) are the lifted vectors in order.
The same \(U_\kappa\) works for every allowed \(\xi\) with
\(\ell\nmid q_\xi\).

\begin{proposition}[Short images and span for the chosen arithmetic data]
\label{prop:marked-upper-estimates}
Fix \(q_0,d,\epsilon\) as in Section~\ref{sec:datum}. For every sufficiently large integer \(x\),
there exists a set \(\mathcal S\subset S\) of integral forms with
\[
 \operatorname{span}\mathcal S=S,\qquad
 \|\sigma\Theta_1P_g\phi\|^2\le e^{A-m\mathcal R+L}
 \quad(\phi\in\mathcal S,\ g\in G,\ \sigma\in\operatorname{Aut}(\mathbb C)).
\]
For each \(t\subseteq G\) and every independent tuple \(y\)
of images of everywhere integral forms under \(P_g\), \(g\in t\),
\[
 \operatorname{Ind}(y)\le2|t|\dim(S)\,faB^3.
\]
Finally, one can choose \(t\subseteq G\) such that
\[
 |t|\dim S\le5D\log(1+2b^2),\qquad
 \dim\sum_{g\in t}P_g(S)\ge(1-4/b^2)D .
\]
\end{proposition}

\begin{proof}
Proposition~\ref{prop:marked-orbital-span} gives the
translated-span dimension and the bound on \(|t|\dim S\).
Proposition~\ref{prop:translated-norm-and-index} gives the
short generators and the index on the same datum. By
Proposition~\ref{prop:shortness-normalization}, the common scalar
\(c_B\in\mathbb Q_{>0}\) gives \(P_g=c_BU_g\), so
\(P_g(S)=U_g(S)\) for every \(g\in G\).
\end{proof}

\begin{proposition}[Determinant estimate for the chosen arithmetic data]
\label{prop:marked-determinant}
Fix \(q_0,d,\epsilon\) as in Section~\ref{sec:datum}.
At every sufficiently large integer \(x\)
the comparison lifts are injective. For every independent
tuple \(y\in V^n\),
\[
 \left.
 \begin{gathered}
  n\ge(1-4/b^2)D,\qquad
  \operatorname{Cov}(y)\le c,\qquad \operatorname{Ind}(y)\le\delta
 \end{gathered}
 \ \right\}\quad\Longrightarrow\quad D(A-3L)\le c+2\delta.
\]
\end{proposition}
The proof is in Section~\ref{sec:marked-determinant-assembly}.

\begin{proposition}[Cancellation of the main covolume term]\label{prop:cancel}
Fix \(q_0,d,\epsilon\) as in Section~\ref{sec:datum}. At every sufficiently
large integer \(x\), the datum satisfies \(m\mathcal R\le faB^5\).
\end{proposition}
\begin{proof}
Take \(t,\mathcal S\) from Proposition~\ref{prop:marked-upper-estimates}.
Since \(\mathcal S\) spans \(S\), the union of the \(P_g(\mathcal S)\),
\(g\in t\), spans \(\sum_{g\in t}P_g(S)\). Choose a basis
\(y_i=P_{g_i}\phi_i\) from this union and fix each pair
\((g_i,\phi_i)\) before the coefficient place or twist.
Its cardinality satisfies \((1-4/b^2)D\le n\le D\).
The Gram bound of Section~\ref{sec:short-images-proved} applies to
these generators at every complex embedding. Their integral
preimages give the index bound on this same basis, so
\[
 \operatorname{Cov}(y)\le n(A-m\mathcal R+L),\qquad
 \operatorname{Ind}(y)\le 2|t|\dim(S)\,faB^3
     \le 10D\log(1+2b^2)faB^3=DI_0.
\]
Proposition~\ref{prop:marked-determinant} and \(n\le D\) give
\[
 nm\mathcal R
 \le (n-D)(A+L)+D(4L+2I_0)\le D(4L+2I_0).
\]
For \(b\ge3\) and \(B\ge128\),
\[
 n\ge (1-4/b^2)D\ge D/2>0,\qquad
 I_0=10faB^3\log(1+2b^2)\le 30L.
\]
Hence \(m\mathcal R\le 8L+4I_0\le 128L\le faB^5\).
\end{proof}

\section{The ramified source at the prescribed level}
\label{sec:ramified-source}

\sectionopening{In this section, the additional level condition singles out a subspace of the full
definite automorphic source. We identify its integral coordinates
before using it in the pointwise estimate of
Section~\ref{sec:upper-chain-size} and the dimension argument of
Section~\ref{sec:marked-span-proof}.}

In the definite source of rank \(k\), write \(K_M\) for
the integral principal level at the primes dividing
\(\bar\alpha\). Write \(J_H\) for its subgroup on which
the last covector of a matrix reduces to the last covector of the
identity at the primes dividing the prescribed generator
\(\pi_N=a_0+ib_0\). The source space \(X\) consists of the
forms at \(J_H\) with the fixed finite-dimensional
coefficient module. These level descriptions specify
the source used in the comparison.

The primes \(p_0,p_1,N\) are \(1\bmod q_2^2\), and \(q_1\)
and \(d\) divide \(q_2^2\), by the definitions of \(s,L_0,q_1\).
Thus the three primes split completely in \(E\).
Since \((\bar\alpha)=\bar{\mathfrak p}_0\mathfrak p_1\) in \(K\),
the old principal level uses the places of \(E\) above
\(\bar{\mathfrak p}_0\) at \(p_0\) and above \(\mathfrak p_1\) at \(p_1\).
Each set consists of \(f\) places of norm \(p_i\), one per conjugate pair. At \(N\), the divisor of \(\pi_N\) chooses
one of the two primes of \(\mathbb Q(i)\) above \(N\).
It gives \(f\) distinct primes of \(E\), one from each
complex conjugate pair, all of norm \(N\).

\begin{lemma}[Trivial stabilizers at the prescribed level]\label{lem:marked-source-freeness}
Fix the datum of Section~\ref{sec:datum} and let \(p\in\{p_0,p_1\}\).
On any of the definite unitary groups of the datum, suppose the
level at one of the chosen split places above \(p\) is contained
in the principal congruence subgroup. Then every finite class
has trivial rational stabilizer.
\end{lemma}

\begin{proof}

A stabilizer is the intersection of the rational points with
a compact subgroup of the adeles. The archimedean group is compact
by definiteness, and the finite compact is the conjugate of the
specified level. The rational points are discrete, so the stabilizer
is finite. The chosen split place has degree one over
\(\mathbb Q_p\). The stabilizer embeds, after conjugation, in
\(1+p M_n(\mathbb Z_p)\), where \(n\) is the rank of the group.
For \(B\in M_n(\mathbb Z_p)\) with a unit entry and \(a\ge1\),
the minimum entry valuation of \((I+p^aB)^h-I\) is
\(a+v_p(h)\) for every positive \(h\). The binomial expansion
gives this for \(p\nmid h\), and raising to the \(p\)-th power
increases the valuation by one since \(p\) is odd. Thus this
principal congruence group is torsion-free and the rational
stabilizer is trivial.
\end{proof}

\begin{proposition}[Dimension of the ramified source]\label{prop:ramified-source-codimension}
For the ramified subspace \(S\subseteq X\) for the chosen arithmetic data,
\[
 \dim(X/S)\le\frac{f}{N-1}\dim X\le b^{-8}\dim X,
 \qquad
 \dim X\le 2\dim S .
\]
\end{proposition}

\begin{proof}

For a prime \(w\mid\pi_N\), take a local unit at \(w\)
lifting a generator of the residue units, form the idele
\(a_w\) with all other components one, and put
\(z_w=a_w/\bar a_w\). Its scalar matrix lies in \(K_M\),
since it is supported at \(w,\bar w\) and \(N\nmid M\).
The local condition on the last covector at \(w\) says
\(z_w^j\in J_H\) exactly when \(N-1\) divides \(j\);
all other conditions on the last covectors hold because \(\bar w\)
does not divide \(\pi_N\).

Let \(J^+=\langle J_H,z_w:w\mid\pi_N\rangle\). Centrality and
the independent residue conditions at these \(f\) places give
\[
 C_N=J^+/J_H\simeq\prod_{w\mid\pi_N}\mathbb F_N^\times .
\]
The level \(J^+\) still has principal congruence at the old split
places over \(p_0,p_1\). By Lemma~\ref{lem:marked-source-freeness}
the stabilizer of each \(J^+\)-class is trivial. Choose a representative
\(g\) of each such class. Its \(J_H\)-classes are the \([gz]\),
\(z\in C_N\), with independent values in the coefficient module.
Thus each fiber is a copy of the regular representation of \(C_N\)
with coefficients in this module. Each character of \(C_N\) occurs in
\(X\) with multiplicity \(\dim X/(N-1)^f\).

For each \(w\mid\pi_N\) the operator
\[
 A_w=\frac1{N-1}\sum_{j=0}^{N-2}R_H(z_w)^j
\]
retains the characters trivial on the \(w\)-factor. The common kernel
\(S\) consists of characters nontrivial on every factor, giving
\[
 \dim S=\left(1-\frac1{N-1}\right)^f\dim X,\qquad
 \dim(X/S)\le\frac f{N-1}\dim X\le b^{-8}\dim X .
\]
The integral kernel at a coefficient place is formed using the sums
\((N-1)A_w\) on the original lattice, as in
Section~\ref{sec:exceptional-jet-hecke}.
Indeed \(f\le q_1\le b\) for large \(x\), and the
interval of primes for \(N\) and the definition of its center
give \(\log N\ge(19(m-1)-15/2)\log b\),
so \(N-1\ge b^9\). In particular
\(\dim X\le2\dim S\).
\end{proof}

\section{Normalization and bounds for the theta images}
\label{sec:short-normalization}

\sectionopening{In this section, we express the local divisibility estimates as bounds on the change of
integral lattice. Summing their logarithms recovers the prime sum in
Section~\ref{sec:joint}, which cancels the corresponding part of the
archimedean bound in Section~\ref{sec:upper-chain-size}.}

The short-image estimate of Proposition~\ref{prop:marked-upper-estimates} is
reduced to archimedean estimates for the unscaled images
and local lattice norm estimates with prime-dependent savings.
For a prime \(p\) put
\[
 S_p(Y,\tau)=\frac12\sum_{e\ge1}
      (1-\chi_K(p^e))(Y/p^e-\kappa_{\langle-p^e\rangle\tau})_+,
 \qquad
 g_p=m\left(\sum_{\tau\in\Phi_a}S_p(a,\tau)
                 +k\sum_\tau S_p(\kappa_\tau,\tau)\right).
\]
For \(Y=a\) or \(Y=\kappa_\tau\), a nonzero positive part requires
\[
  p^e<Y/\kappa_{\langle-p^e\rangle\tau}\le Yx_1/a\le x_1\le a.
\]
The factor \(\tfrac12(1-\chi_K(p^e))\) selects the odd
exponents at an inert prime, is \(\tfrac12\) at the ramified
prime \(d\), and vanishes at a split prime. For the odd inert
exponents the multiplier is \((-p)^e\bmod q_1\). At \(d\)
it is \(-d^e\bmod q_1\) for every \(e\ge1\).
Write \(E_0=faB^3\) and let \(F\) consist of the nonsplit
primes below \(x_1\), together with \(N\).
Grouping the von Mangoldt sum by prime powers gives
\[
 2\sum_{p\in F}g_p\log p
  =m\sum_{\tau\in\Phi_a}S(a,\tau)
           +mk\sum_\tau S(\kappa_\tau,\tau).
\]
In particular \(A-m\mathcal R+2\sum_{p\in F}g_p\log p\)
equals the sum of the archimedean jet terms in \(A\)
plus \(m\) times the expected value of \(\sum_{\tau\in\Phi_a}S(a,\tau)\)
from Section~\ref{sec:joint}.

\subsection{Bounds for the local denominators}\label{sec:local-denominator-sum}

\sectionopening{In this subsection, we bound the total logarithmic contribution of the local estimates in
Section~\ref{sec:short-normalization}. The exceptional congruence
classes and the prime \(N\) are treated separately because their local
estimates differ.}

For these data, define the local error bound by
\[
 e_p=
 \begin{cases}
  \displaystyle\frac{2faB^2}{p}\,c_p,
       &p<x_1\ \text{nonsplit in }K,\\
  2mfb,&p=N,\\
  0,&\text{otherwise},
 \end{cases}
 \qquad
 c_p=\begin{cases}
 m,&p\equiv\pm1\bmod m,\\
 1,&\text{otherwise}.
 \end{cases}
\]
\begin{lemma}[Total local denominator bound]\label{lem:total-local-denominators}
For fixed \(q_0,d,\epsilon\) and all sufficiently large \(x\),
\[
 \sum_p e_p\log p\le E_0=faB^3 .
\]
\end{lemma}

\begin{proof}
Put \(\ell_1=\log x_1\). Since \(\log N\ge mB>W\) for large \(x\),
the two nonzero cases defining \(e_p\) are disjoint.
For \(p=mc\pm1\) with \(c\ge1\), one has \(m/p\le2/c\).
Consequently
\[
 \sum_p e_p\log p
 \le 2faB^2\,\ell_1
         \left(\sum_{j\le x_1}\frac1j+
                   2\sum_{\pm}\sum_{c\le x_1}\frac1c\right)
                   +2mfb\log N
 \le 10faB^2\,\ell_1(1+\ell_1)+2mfb\log N.
\]
Since \(\ell_1\le W\), \(B>W^6\), \(m^2<4x\),
\(a=xb\) and \(\log N\le25mB\), division by \(E_0\) gives
\[
 \frac{1}{E_0}\sum_p e_p\log p
 \le \frac{10W(1+W)}{B}+\frac{50m^2}{xB^2}
 < \frac{10(1+W)}{W^5}+\frac{200}{W^{12}}<1
 \qquad(W\ge3).\qedhere
\]
\end{proof}

\subsection{Counting terminal double cosets}\label{sec:terminal-classes}

\sectionopening{In this subsection, we bound the number of terminal double cosets by counting definite
Hermitian lattices and refining their classes at the specified level.
This controls the class sum in Section~\ref{sec:upper-chain-size}.}

Let \(n_{\rm t}\) be the rank of a definite stage under consideration, and let
\(H\) be the unitary group of the standard positive Hermitian
space \(E^{n_{\rm t}}\). Write \(K_{\rm int}\subset H(\mathbb A_f)\)
for the stabilizer of the standard integral lattice at all
finite places.

\begin{lemma}[Terminal class count]\label{lem:terminal-class-count}
For the trivial twist, let \(h_{\rm term}\) be the number of
double classes at the prescribed level on any \(X_j\), \(0\le j\le m\),
or on the definite endpoint \(Y_m\). For fixed \(q_0,d,\epsilon\)
and sufficiently large \(x\), uniformly in the stage,
\[
 \log h_{\rm term}=O(x^{52}),\qquad
 \log h_{\rm term}+4\log|G|\le x^{99}.
\]
\end{lemma}

\begin{proof}
We first count classes at \(K_{\rm int}\).
Write \(d_E=[E:\mathbb Q]\), \(\Delta_E=|\operatorname{disc}(E)|\),
and \(u=d_E n_{\rm t}\). Each class gives an integral lattice
\(\Lambda\) obtained by intersecting a finite adelic translate
of the standard integral lattice with \(E^{n_{\rm t}}\).
It has integral Hermitian pairings and is locally isometric
to the standard lattice. An isometry between such lattices
over \(E\) therefore identifies their double classes.

Use the Euclidean norm squared given by
\(\operatorname{Tr}_{E/\mathbb Q}h(v,v)\), where \(h\)
is the Hermitian form. The determinant of the integral trace
pairing on \(\Lambda\) is \(\Delta_E^{n_{\rm t}}\):
on the standard lattice use an integral basis of \(E\) in each
coordinate, with conjugation changing the trace determinant
only by a sign. Local isometry gives the same determinant
valuation at every rational prime, and positivity fixes the
determinant.

Thus the covolume is \(\Delta_E^{n_{\rm t}/2}\).
Every nonzero lattice vector has length at least one, by
positivity and integrality of its squared length.
Minkowski's second theorem \cite[Ch.~VIII, \S4.3, Theorem~V]{Cassels}
and the cube of side \(2/\sqrt u\) in the Euclidean unit ball
give, for the last successive minimum,
\[
                   \lambda_u\le u^{u/2}\Delta_E^{n_{\rm t}/2}.
\]
Choose \(u\) independent lattice vectors of length at most
\(2\lambda_u\). They span \(E^{n_{\rm t}}\) over \(E\), so
\(n_{\rm t}\) of them form an \(E\)-basis \(v_1,\ldots,v_{n_{\rm t}}\).

Let \(T=(h(v_i,v_j))\), with the convention that the second
variable is linear. Its entries are algebraic integers,
and their absolute values at all embeddings are at most
\(H_0=4u^u\Delta_E^{n_{\rm t}}\), by Cauchy--Schwarz at
each embedding. To count these entries, note that for distinct
algebraic integers in \(E\), their difference has at least one
embedding of absolute value at least one, by taking its norm.
The vectors of their real and imaginary parts at one place from
each complex conjugate pair therefore have sup-norm separation
at least \(1/\sqrt2\). Subdivision of the box of side \(2H_0\)
into cubes of side \(1/2\) gives at most
\((4H_0+2)^{d_E}\) choices per entry.

For each fixed \(T\) identify the space over \(E\) via the
chosen basis. Let \(L\) be its \(\mathcal O_E\)-span.
We have \(L\subseteq\Lambda\subseteq L^\#\), where \(L^\#\)
is the Hermitian dual of \(L\), since \(\Lambda\) has integral
pairings. The quotient \(L^\#/L\) is isomorphic as an abelian
group to \(\mathcal O_E^{n_{\rm t}}/T\mathcal O_E^{n_{\rm t}}\),
by writing the pairings in this basis. It has order
\(\mathcal I=|\operatorname{N}_{E/\mathbb Q}\det T|\), at most
\((n_{\rm t}! H_0^{n_{\rm t}})^{d_E}\).

Every subgroup of this quotient has at most \(u\) generators,
since its inverse image in \(L^\#\) is a sublattice of rank
\(u\). There are at most \(\mathcal I^u\) choices for the
intermediate lattice, using ordered tuples with zero entries
allowed. These counts bound the number of classes because
two lattices with the same data are carried to one another
by an \(E\)-linear isometry sending the chosen bases to
one another. We obtain
\[
 \log |H_{\rm rat}\backslash H(\mathbb A_f)/K_{\rm int}|
     =O\bigl(u^2(1+\log H_0+\log n_{\rm t})\bigr).
\]
Here the left quotient is by the \(E\)-linear isometries
of the Hermitian space.

If the stage level contains the principal congruence subgroup
of \(K_{\rm int}\) modulo \(MN\), its number of classes is at most
the preceding class count times \((MN)^{d_E n_{\rm t}^2}\).
Indeed one may refine each class by right cosets of the
principal congruence subgroup in \(K_{\rm int}\), whose index
is at most the number of matrices modulo \(MN\).

Write \(K(I)\) for matrices in the integral stabilizer congruent to
the identity modulo \(I\), in the rank under consideration. On each
upper space \(X_j\), the old condition is reduction to the identity at the
chosen places dividing \(\pi_M\). At the places dividing \(\pi_N\)
the last row must reduce to the last row of the identity. Thus the
upper level contains \(K(\pi_M\pi_N)\): this subgroup fixes every
row modulo both generators. After restriction to \(Y_m\) the
untwisted terminal level imposes the old residue conditions, so it
contains \(K(\pi_M)\). Since \(MN=(\pi_M\pi_N)\overline{\pi_M\pi_N}\),
both levels contain \(K(MN)\). The rank is at most
\(k+m(5r+1)=O(x^{11/2})\), and \(d_E=2f=O(x^{11}\log x)\).
In particular \(u=O(x^{17})\). For \(C=\mathbb Q(\zeta_{q_1})\),
linear disjointness of \(K,C\) factors the trace matrix of
\(\mathcal O_K[\zeta_{q_1}]\) in a product basis. Thus
\[
  |\operatorname{disc}(\mathcal O_K[\zeta_{q_1}])|
    =d^{[C:\mathbb Q]}|\operatorname{disc}\Phi_{q_1}|^2 .
\]
Since \(\Phi'_{q_1}(\zeta_{q_1})\) divides \(q_1\),
\(\operatorname{rd}(E)\le\sqrt d\,q_1\le dq_1\). Hence
\(\log H_0=O(u\log x)\). Also \(\log(MN)=O(mB)\) and
\(B=O(x^\epsilon)\). The logarithm of the terminal class
count is consequently \(O(x^{52})\).
There are \(f\) chosen places above each of \(p_0,p_1\).
The residue group \(G\) acting on the target is the product
of the corresponding groups \(\mathrm{GL}_r\), so
\(\log |G|\le fr^2\log M\).
It follows that \(\log h_{\rm term}+4\log|G|\le x^{99}\)
for all sufficiently large \(x\).
\end{proof}

\subsection{Projection onto the prescribed finite types}\label{sec:upper-good-projector}

\sectionopening{In this subsection, the central idempotent acts on every lift by the intertwining relation
of Section~\ref{sec:seesaw-normalization}. Its integral multiple
controls the valuations and archimedean norms of all translates, using
the class count of Section~\ref{sec:terminal-classes}.}

Let \(e_0\) be the central idempotent of the projected part
in the complex group algebra of the finite group \(G\).
Let \(z\) be an automorphic form on \(U(Y_0)\) at the prescribed
principal level. Write \(\mathcal L_\xi z\) for the terminal
algebraic coefficient of \(L_{\xi,m-1}\cdots L_{\xi,0}z\),
where each theta integral includes its source index
\(A_{\xi,j}\) as in Section~\ref{sec:seesaw-normalization}.
\begin{lemma}[The integral finite-type projection]\label{lem:integral-finite-type-projection}
The coefficients of \(|G|e_0\) are algebraic integers, of absolute value
at most \(|G|\) at every complex embedding. For any \(g\in G\), replacing
\(\mathcal L_\xi z\) by \(\mathcal L_\xi(|G|e_0 g z)\) preserves the local
lattice norm bound \(p^{e_p-g_p}\). If the untwisted conjugate squared pointwise
norms of \(\mathcal L_1 z\) are at most \(e^{K_0}\), then the squared
terminal class norm after projection is at most \(e^{K_0+x^{99}}\).
\end{lemma}

\begin{proof}
The expansion of \(e_0\) as a sum of primitive
central idempotents gives
\[
       |G| e_0(t)=\sum_{\pi\ {\rm selected}}\dim(\pi)\chi_\pi(t^{-1})
                    \qquad(t\in G).
\]
Each coefficient on the right is an algebraic integer
because the character values are sums of roots of unity.
At every complex embedding its absolute value is at most
\(\sum_{\pi\ {\rm selected}}(\dim\pi)^2\le |G|\).

\smallskip
\indent\textit{Equivariance on the original chain.}
For \(s=(s_w)\in G\), choose at each stage an integral
finite lift \(k_j(s)\). Its residue at an old place \(w\)
is \(\operatorname{diag}(s_w,I)\) in the original order
of the blocks, and its residue at every twist place
is the identity. At the source \(k_0(s)\) acts by the
given residue action, since two lifts differ by a
principal-level element. Reduction is onto at these
split places, so these choices are possible at every
stage. At places outside the prescribed set their components can be the
identity.

Let \(i_j\) denote the inclusion of the \(j\)-th Hermitian space as the first block of its successor. The original
forward-arrow finite test at the chosen and twist
places is the indicator of \(i_j\).
The simultaneous source and target reductions \(a,b\)
of \(k_j(s),k_{j+1}(s)\) satisfy \(b i_j=i_j a\).
Thus the finite split Hom action \(X\mapsto b^{-1}Xa\)
fixes the prescribed support. The determinants are units,
and the two rank-twisted splitting characters leave
the factor \(\psi_{j,w}(\det b)/\psi_{j,w}(\det a)\)
on the indicator, where \(\psi_{j,w}\) is the tame part
of the arrow character. At an old place the two determinants
reduce to \(\det s_w\), and at a twist place both reduce to one.
Their quotient is therefore a principal unit, on which
\(\psi_{j,w}\) is trivial. On the other
finite components the vacuum is fixed, and these
translations have no real component. The integrated
arrow datum is therefore fixed.

Change the source variable by \(k_j(s)\) in the
arrow integral, using the original quotient Haar.
It follows that translating the image under the lift by
\(k_{j+1}(s)\) gives the arrow applied to the source
translate by \(k_j(s)\). Iterate to the terminal space.
The algebraic coefficient realization commutes with
finite right translations, so the final formula is
\[
 \mathcal L_\xi(s z)(h)
        =(\mathcal L_\xi z)(h\,k_m(s)).
\]

With \(k_\xi(s)=k_m(s)\), expanding \(e_0\) gives
\[
 \bigl(\mathcal L_\xi(|G|e_0 g z)\bigr)(h)
    =\sum_{t\in G}|G|e_0(t)
             (\mathcal L_\xi z)(h\, k_\xi(tg)).
\]
Here \(k_\xi(tg)\) is the integral first-block lift
of \(tg\) and \(h\) is a finite adelic point in the
terminal space. In the transported Weyl coordinates a right translate
by any such integral matrix preserves the lattice norm:
the change of coordinates and its inverse both preserve the
unit lattice. The ultrametric inequality therefore transfers
any bound \(p^{e_p-g_p}\) for \(\mathcal L_\xi z\) to
\(\mathcal L_\xi(|G|e_0 g z)\) with the same exponent.
The transfer uses the prescribed chain and holds for all
allowed twists and every translate. If the local estimate
is first given on a submodule of the full integral source,
intersect it with the ramified integral source. A matrix
clearing a tuple from the ramified source into the first
submodule still clears it into this intersection since
linear combinations preserve ramification.

The coefficient bound and the triangle inequality give
\[
 \|\sigma(\mathcal L_1(|G|e_0gz)(h))\|_{\rm W}^2
       \le |G|^4e^{K_0}.
\]
Summing with reciprocal stabilizer weights and applying
Lemma~\ref{lem:terminal-class-count} gives
\[
 h_{\rm term}|G|^4e^{K_0}\le e^{K_0+x^{99}}.\qedhere
\]
\end{proof}

\subsection{Pointwise estimates in the composition of lifts}
\label{sec:upper-chain-size}

\sectionopening{In this subsection, we estimate the upper lift by exploiting the additional level
condition of Section~\ref{sec:ramified-source}. The resulting
pointwise bound is combined with the class count of
Section~\ref{sec:terminal-classes} in the proof of
Proposition~\ref{prop:marked-upper-estimates}.}

\subsubsection{The sum of the part sizes}\label{sec:upper-part-size-sum}

\sectionopening{In this subsection, the level prime from Section~\ref{sec:datum}
supplies the term subtracted in the terminal estimate of
Section~\ref{sec:upper-terminal-norm}. We use the function \(J\) below to
sum the factorial contributions from the upper partitions.}

Write \(J(l)=\log(l!)+l\log(M/\Omega)\), put
\(\overline w=f^{-1}\sum_\tau\log\kappa_\tau\), and
write \(s_{\rm av}=a(\log a-\overline w+c_{d,1})\).
Set \(G_0=\log(M\operatorname{rd}(E)/(2\pi))\).
By the definition at modulus one,
\(c_{d,1}=L'/L(1,\chi_K)-\gamma_{\rm Euler}-1\), so
\(\log(M/\Omega)=G_0+1+c_{d,1}\).
The comparison of \(\sum_{j=1}^a\log j\) with
\(\int_1^a \log t\,dt+\log a\) gives
\[
 J(a)-s_{\rm av}\le a(\overline w+G_0)+1+\log a.
\]
The chosen prime satisfies
\(\log N\ge m\overline w+(m-1)G_0-\tfrac{15}{2}B\).
Consequently
\[
 m(J(a)-s_{\rm av})-a\log N
 \le a(G_0+\tfrac{15}{2}B)+m(1+\log a)
 \le 32aB
\]
eventually, since \(G_0\le23B\) and \(m(1+\log a)\le aB\).
Define the size of the principal parts by
\[
 S_{\rm main}=m\left(|\Phi_a|J(a)+k\sum_\tau J(\kappa_\tau)\right).
\]
For the step from \(X_j\) to \(X_{j+1}\), put \(n_j=k+j(5r+1)\)
and let \(m_{\tau,h}\) be the column counts of the lift weight
in Section~\ref{sec:marked-column-gain}. The first \(n_j\)
row lengths of the target partition are
\[
 \lambda_{j,\tau,i}=\sum_{i<h<n_{j+1}}m_{\tau,h}
                   \qquad(0\le i<n_j).
\]
In particular each column of height \(n_j\) contributes the top
exterior power on the source and occurs in every one of these rows. Set
\[
 S_j=\sum_{\tau\in\Phi}\sum_{0\le i<n_j}J(\lambda_{j,\tau,i})
 \qquad(0\le j<m),\qquad h_0=r^8 .
\]
\begin{lemma}[Factorial bounds for the upper composition]\label{lem:upper-factorial-sums}
For fixed \(q_0,d,\epsilon\) and sufficiently large \(x\),
\[
 \begin{aligned}
 S_{\rm main}-|\Phi_a|a\log N
  &\le A-m\mathcal R+2\sum_p g_p\log p+32faB,\\
 \sum_j S_j&\le S_{\rm main}+\tfrac12 E_0 .
 \end{aligned}
\]
The term \(-|\Phi_a|a\log N\) comes from
\(\sum_\tau a_\tau=|\Phi_a|a\) in the terminal branching
exponents \(a_\tau-mh_0\) of Lemma~\ref{lem:terminal-branching-norm}.
\end{lemma}

\begin{proof}
Insert the prime-power identity of Section~\ref{sec:short-normalization}
and the expression for \(\mathcal R\) in Section~\ref{sec:joint}.
The sums \(S(\kappa_\tau,\tau)\) and \(S(a,\tau)\) cancel, giving
\[
 A-m\mathcal R+2\sum_{p\in F}g_p\log p
     =mk\sum_\tau J(\kappa_\tau)+m|\Phi_a|s_{\rm av}.
\]
Subtracting this expression from \(S_{\rm main}-|\Phi_a|a\log N\)
gives exactly
\[
 |\Phi_a|\bigl(m(J(a)-s_{\rm av})-a\log N\bigr)\le32faB ,
\]
since \(|\Phi_a|\le f\).

In zero-based row indices, the column table of
Section~\ref{sec:marked-column-gain} gives the actual row length
of the step from \(X_j\) to \(X_{j+1}\):
\[
 \lambda_{j,\tau,i}
 =a_\tau\mathbf1_{i=0}
   +\bigl(w_1(\tau)-\lfloor n_1/2\rfloor\bigr)\mathbf1_{i<k}
   +(\Delta-t-1)\sum_{u=1}^{j}\mathbf1_{i<n_u}
 \qquad(0\le i<n_j).
\]
The sum is empty at \(j=0\) and counts the height \(n_j\) at
\(j\ge1\). A column of height \(h\) contributes one to each row
\(i<h\). Thus the row bounds have the form
\[
 \lambda_{j,\tau,i}\le\mu_{\tau,i}+s_i,\qquad
 \mu_{\tau,i}=
 \begin{cases}
   a\mathbf1_{\Phi_a}(\tau),&i=0,\\
   \kappa_\tau,&1\le i<k,\\
   0,&i\ge k,
 \end{cases}
 \qquad s_i=h_0^2+(b+a/x_1)\mathbf1_{\{0\}}(i).
\]
For large \(x\), \(\mu_{\tau,i}+s_i\le a+h_0^2+b+a/x_1\le2a\).
Since \(M/\Omega>1\), the increment
\(J(l+1)-J(l)=\log(l+1)+\log(M/\Omega)\) is positive. For \(l+s\le2a\),
the choice of scales gives
\[
 J(l+s)-J(l)
 \le s\bigl(\log(2a+1)+\log(M/\Omega)\bigr)\le25sB .
\]
Indeed \(2a+1\le b^2\), \(\log M\le22B\), and
\(\log\operatorname{rd}(E)=O(\log x)=o(B)\).

At each step the principal lengths give
\[
 \sum_{i=0}^{n_j-1}J(\mu_{\tau,i})
  =\mathbf1_{\Phi_a}(\tau)J(a)+(k-1)J(\kappa_\tau).
\]
Since \(J(\kappa_\tau)\ge0\), summing the increments over all rows,
embeddings and steps gives
\[
 \sum_j S_j
   \le S_{\rm main}+25mfB(b+a/x_1+n_mh_0^2)
   \le S_{\rm main}+\tfrac12 E_0 .
\]
Indeed \(x/x_1\to0\), \(n_mh_0^2\) is polynomially
bounded in \(x\), and \(b\) exceeds every fixed power of \(x\).
\end{proof}

\subsubsection{Integral generators of the ramified source}\label{sec:upper-integral-generators}

\sectionopening{In this subsection, we construct generators in the integral model of
Section~\ref{sec:ramified-source}. Their pointwise norms supply the
initial bound for the composition in Section~\ref{sec:upper-terminal-norm}.}

\begin{lemma}[Integral generators with bounded Weyl norm]\label{lem:ramified-short-generators}
Fix \(q_0,d,\epsilon\). For every sufficiently large \(x\) there is a set
\(\mathcal S\subseteq S\) of everywhere integral forms such that
\[
 \operatorname{span}(\mathcal S)=S,\qquad
 \sup_{\sigma,h}\|\sigma(\phi(h))\|_{\rm W}^2\le e^{E_0}
                    \quad(\phi\in\mathcal S).
\]
Here \(h\in U(X_0)(\mathbb A_{F_0,f})\) and
\(\sigma\in\operatorname{Aut}(\mathbb C)\).
\end{lemma}

\begin{proof}
Put \(n_0=k\), \(D_0=q_1d\), where \(q_1d\) is in the
different of \(E/\mathbb Q\), and
\[
 C_0=6D_0(fn_0)^2,\qquad
 M_0=n_0f+\#\{\text{weight boxes in columns of height }<n_0\}.
\]
We use \(D_0\) here for a positive rational integer.
Throughout this section \(N\) still denotes the ramified
level prime. At every original source class \(g\), use
the lattice
\[
 L(g)=E^{n_0}\cap g\widehat{\mathcal O}_E^{\,n_0}.
\]
The finite unitary action and the standard definite
Hermitian form make \(L(g)\) Hermitian self-dual.
Here and below duality means integrality of the pairing
over \(\mathcal O_E\). The diagonal archimedean
embedding has Euclidean norm squared
\(\sum_{\tau}\|v_\tau\|^2\).
For \(0\ne v\in L(g)\), its Hermitian square
\(h(v,v)\) is a totally positive integer of \(F_0=E^+\).
The arithmetic--geometric mean inequality applied to
its conjugates gives \(\|v\|^2\ge f\ge1\).

The real inner product in these coordinates is half
the trace pairing. The Euclidean dual of
\(2D_0^{-1}L(g)\) is consequently the archimedean
image of \(D_0\mathfrak d_E^{-1}L(g)\), which is
contained in \(L(g)\). For the identification, a
real functional integral on a lattice is rational
on its rational span. By nondegeneracy of the
rational trace pairing it therefore has a rational
trace representative. Pairing with all the
\(\mathcal O_E\)-multiples of every lattice vector
puts the representative in the codifferent times the
Hermitian dual, proving the asserted equality.

For completeness a Euclidean lattice of dimension
\(h\) with shortest dual vector of length at least
\(1\) has covering radius at most \(h^2\).
Minkowski's convex body theorem \cite[Ch.~III, \S2.2, Theorem~II, p.~71]{Cassels}
applied to cubes inscribed in Euclidean balls gives
\(\lambda_1(L)\lambda_1(L^*)\le h\), since the
two covolumes multiply to \(1\). Choose a primitive
vector of length at most \(2h\) and project
orthogonally to its complement. The projected lattice
is discrete since the vector extends to a lattice
basis.

Its dual is \(L^*\) intersected with the
complement. Induction gives error at most \((h-1)^2\)
in the projected space, and rounding in the chosen
direction contributes at most \(h\). This proves the
radius bound (with the zero-dimensional case immediate).
The vectors of length at most three times the
covering radius span the whole real space: if their
span has a unit normal \(u\), a lattice vector within
radius \(r\) of \(2ru\) has norm at most \(3r\) and
has nonzero projection on \(u\).

Apply this argument to \(2D_0^{-1}L(g)\) with
\(h=2fn_0\), then multiply by \(D_0/2\).
The resulting vectors in \(L(g)\) of norm at most
\(C_0\) span the real space in the diagonal embedding.
They also span \(E^{n_0}\) over \(E\): an \(E\)-linear
functional vanishing on them extends at each complex
place and hence is zero. Take an \(E\)-basis from
these vectors. With its columns written \(Y(g)\),
\[
 Y(g)\in\mathrm{GL}_{n_0}(E),\qquad
 g^{-1}Y(g)\in M_{n_0}(\widehat{\mathcal O}_E),\qquad
 |\tau(Y(g)_{ij})|\le C_0 .
\]
Fix a finite adelic representative \(g\) of the class. For each
determinant exponent \(e_\tau\in\mathbb Z\), apply the construction
to the group element \(g^{e_\tau}\). Thus
\(\det(g^{e_\tau})=\det(g)^{e_\tau}\), also for negative exponents.
Here \(D_0\in\mathfrak d_E\) because the different of \(E/\mathbb Q\)
divides the product of the extended differents of its two
constituent fields \(K\) and \(\mathbb Q(\zeta_{q_1})\), which in
turn divides \(d q_1\mathcal O_E\).

Write \(e_\tau\) for the determinant exponent
and \(v(y)\) for the column-polynomial Weyl evaluation
vector, in the original integral model. Use the fixed coefficient embedding \(\iota_\tau=\tau\circ\mathrm{conj}\).
For each
class \(g\) use the value
\[
 v_g(y)=\left(\prod_\tau\det \iota_\tau(Y(g^{e_\tau}))\right)
           \rho_{\mathrm{pol}}\bigl((\iota_\tau(Y(g)))_\tau\bigr)v(y)
\]
for the delta form at \(g\). Put \(B_{\rm pol}=M_0-n_0f\),
the number of polynomial boxes. Take \(y=(y_\tau)_{\tau\in\Phi}\)
with the entries of all the matrices chosen independently in
\(\mu_{B_{\rm pol}+1}\). This finite grid is fixed before choosing
a coefficient embedding. The level is free, so
these delta forms are well defined and there is no
stabilizer division. In the transported lattice norm at
\(g\), the polynomial matrix is evaluated on
\(g^{-1}Y(g)\); the determinant factors are the
determinants of \((g^{e_\tau})^{-1}Y(g^{e_\tau})\).
Both are integral. Polynomial Weyl matrices and
evaluation vectors have integer coefficient formulas
in the chosen model, and roots of unity are units
at every coefficient embedding. Thus the delta
form has lattice norm at most one at every place, including
on its translated class points by equivariance.

For the archimedean bound the evaluation vector
is the Schur action of the first square block of
\(y\) on the unit vector of the reference column
array. After applying \(Y(g)\), the resulting
matrix at each embedding has entries at most
\(n_0 C_0\) in modulus and squared Frobenius norm
at most \(n_0^4 C_0^2\le C_0^6\). The induced
action on each exterior power is bounded by the
operator norm to its degree, and hence the
squared Schur norm is bounded by \(C_0^6\) to
the number of polynomial weight boxes. The
squared determinant factors together contribute at most \(C_0^{6n_0f}\) by Hadamard. These bounds
hold after every complex automorphism, since
the entries of \(Y\) are bounded under all
embeddings and roots of unity still have modulus one.
The squared Weyl norm of the value is at most
\(C_0^{6M_0}\). For a point \(h=\gamma_f g k'\) in the supporting
class, with \(\gamma\in U(X_0)(F_0)\) and \(k'\in J_H\), the
value at \(h\) is obtained from the value at \(g\) by its rational
weight action. For every embedding \(\jmath:E\hookrightarrow\mathbb C\),
\(\jmath(\gamma)^*\jmath(\gamma)=1\), since \(X_0\) has the standard
definite form and \(\jmath(\bar x)=\overline{\jmath(x)}\).
Its exterior-power actions preserve the tensor-induced Weyl norm
and its determinant factors have modulus one. Take
\(\jmath=\sigma\circ\iota_\tau\) to get the bound at every \(h\) after
\(\sigma\). Off that class the value is zero.

At the source the polynomial Weyl module is
\(\bigotimes_{\tau\in\Phi}\operatorname{Sym}^{a_\tau}(\mathbb C^{n_0})\).
Its evaluation coordinates are products over \(\tau\) of degree-\(a_\tau\)
monomials in the first column of \(y_\tau\). The entries in these columns
range independently over the chosen roots of unity.
A functional vanishing on all the evaluation vectors
gives a polynomial of degree at most \(B_{\rm pol}\) in each
entry of the array. In each variable it vanishes at the
\(B_{\rm pol}+1\) chosen roots for every fixed grid point in the
remaining variables. Successive one-variable interpolation therefore
gives the zero polynomial.
The Weyl evaluation polynomials form the dual Weyl
space in the chosen model, so the functional is zero.
The matrix multiplying \(v(y)\) in \(v_g(y)\) is
invertible. Summing these component spans proves
the claim at the original source level.

On the source the determinant exponent is zero and the only columns of height below \(n_0\) are the \(a_\tau\) columns of height one (Section~\ref{sec:marked-column-gain}), so \(M_0\le f(n_0+a)\), and
the scale inequalities give
\[
 4^f C_0^{6f(n_0+a)}\le \exp(E_0),\qquad
 E_0=faB^3 .
\]
Finally apply the commuting differences
\(1-R_H(z_w)\) for \(w\mid\pi_N\) as in
Section~\ref{sec:ramified-source}. Each
preserves the integral lattice and multiplies the maximum
squared norm of the conjugate values by at most \(4\).
Their product has the ramified subspace as its
range: the commuting finite cyclic scalars are
simultaneously diagonalizable and \(1-T\) is
invertible on the nontrivial eigenspaces of \(T\).
These differences therefore carry the integral delta generators to a
spanning set of \(S\) with the stated bound.
\end{proof}

\subsubsection{Conjugation of the theta jets}\label{sec:upper-cm-jet-transport}

\sectionopening{In this subsection, we identify the conjugated theta jets and the
finite vacuum on the original chain of Section~\ref{sec:seesaw-normalization}.
Section~\ref{sec:upper-gaussian-jet} uses the transported jets in the
Gaussian estimate.}

For the kernel \(X_j\to X_{j+1}\), let \(S_j\) be the sum
of \(J\) on the actual upper part sizes \(\lambda_{j,\tau,i}\)
used in the Schur matrix. We estimate the operator on the
original Weyl modules for the untwisted chain.
The following two identities identify the CM and finite-vacuum
identities in this calculation. They concern the same vacuum scalar and
the same untwisted upper tensor used in the integral chain.

\begin{lemma}[Transport of the jets and the finite vacuum]
\label{lem:upper-cm-transport}
For \(c\ne0\), let \(J_c\) be the full holomorphic jet at scale \(c\).
Use \(\gamma=s_{\rm jet}=p_Fa_d/\omega_0\) specified below and let
\(T(x)\) denote finite Heisenberg translation.
In the standard definite model of each positive rank
\(v\), fix the finite seed \(F_v\) given by the cyclotomic product below.
It has \(\theta(F_v)=1\), and its translates span the finite
representation. For every \(\sigma\in\operatorname{Aut}(\mathbb C/K)\)
there is one pair
\((b_\sigma,t_\sigma)\in\mathbb A_{E,f}^\times\times\mathbb C^\times\)
such that, at each rank \(v\), a finite seed \(F'_{\sigma,v}\) satisfies
\begin{equation*}\tag{AJ}\label{eq:upper-CM-jet-transport}
 J_{t_\sigma}(T(b_\sigma x)F'_{\sigma,v})
 =\operatorname{ren}_{\pi_\sigma}\,
       \sigma_{\rm coeff}\bigl(J_\gamma(T(x)F_v)\bigr)
 \quad(x\in\mathbb A_{E,f}^{v}),\qquad |t_\sigma|^2=\Omega^{-1}.
\end{equation*}
For embeddings extending the fixed embedding of \(K\), put
\(\pi_\sigma(\tau)=\sigma\circ\tau\).
Here \(\operatorname{ren}_{\pi_\sigma}\) sends
\(Z_{\tau,i}\) to \(Z_{\pi_\sigma(\tau),i}\).

Let \(v_0\) be the rank-one vacuum and \(\Lambda=\widehat{\mathcal O}_E\).
There is a finite displacement \(y\) for which
\(T_\sharp(u)T(y)F_1=T(y)F_1\) for all \(u\in q_1\Lambda\), and
\begin{equation*}\tag{AV}\label{eq:upper-CM-vacuum}
 s_{\rm vac}v_0=\sum_{u\in\Lambda/q_1\Lambda}T_\sharp(u)T(y)F_1,
 \qquad s_{\rm vac}\ne0 .
\end{equation*}
The operator \(T_\sharp\) includes the finite vacuum semicharacter.
The finite values of \(F_1\) are integral at every coefficient place.
They are products of signed, normalized rank-one theta evaluations at
torsion points along the \(q_1\) forms
\(z\mapsto-\operatorname{Tr}_{E/K}(z\overline{\zeta_{q_1}^{\,k}})/q_1\),
each repeated \(q_1\) times.
\end{lemma}

\begin{proof}
Construct the translated theta sections and their common
period scale first. The finite theta expansion will give
the vacuum projection \eqref{eq:upper-CM-vacuum} and the
spanning assertion at every positive rank.

\smallskip
\proofstep{The normalized CM translation and its period}

Write \(K=\mathbb Q(\sqrt{-d})\) and \(q=q_1\).
We use the degree-one theta function of a lattice \(\Lambda_C\)
with maximal CM by \(O_K\). The distinguished two-torsion
point \(t\) is a generator of \(C[2]\) as an \(O_K/2O_K\)-module.
In the Weierstrass coordinates with group origin at infinity,
the curve is
\[
    y^2+a_1xy=x^3+x,\qquad t=(0,0),
    \qquad \Delta_C=a_1^4-64.
\]
The CM invariant \(j(C)\) is algebraic
\cite[Theorem~12.8]{MilneMF}. The displayed equation gives
\(j(C)=(a_1^4-48)^3/(a_1^4-64)\), so \(a_1\) is algebraic.
The periods are those of the invariant differential of this
equation. The theta function \(\vartheta_C\) has zero divisor
\([t]\), is even, and is normalized by \(\vartheta_C(0)=1\).
For \(A_C=\operatorname{area}(\Lambda_C)\) write
\[
 (T_u^C\vartheta_C)(z)=
       e^{-\pi(z\bar u+|u|^2/2)/A_C}\vartheta_C(z+u).
\]
The fixed embedding \(\Phi_K\colon K\hookrightarrow\mathbb C\)
gives an elliptic CM type with reflex field \(K\). For
\(\sigma\) fixing \(K\), apply the analytic CM theorem to
\(\sigma|_{\overline{\mathbf Q}}\) and the canonical principal
polarization of \(C\). Choose a finite idele \(s\) of \(K\) with
\(\operatorname{Art}^{\rm ar}_K(s)=\sigma|_{K^{\rm ab}}\), where
the Artin map sends a local uniformizer to arithmetic Frobenius
on the unramified quotient \cite[Remark~2.5]{ConradCM}.
Let \(c_\sigma\in\mathbb C^\times\) be the scalar of the CM
tangent-space identification in the coordinates given by the two
invariant differentials. The main theorem gives
\[
 \widehat\Lambda_{C^\sigma}
     =c_\sigma N_{\Phi_K}(s^{-1})\widehat\Lambda_C .
\]
In these differential coordinates the map on each primary torsion
component is induced by \(c_\sigma N_{\Phi_K}(s^{-1})\).
Equivalently, after identifying the tangent spaces by \(c_\sigma\),
the finite multiplier is \(N_{\Phi_K}(s^{-1})\)
\cite[Ch.~IV, Theorem~18.6]{ShimuraCM},
\cite[Theorem~6.3, pp.~22--23]{ConradCM}. We compute the
translation factors of the normalized section below using its
Hermitian form \(\pi\bar u z/A_C\).

\smallskip
\proofstep{Integral normalization}

Over a number field defining the curve and its CM action, the curve
has potentially good reduction at every finite place and its
\(j\)-invariant is integral
\cite[Theorem~12.14 and Corollaries~12.15--12.16]{MilneMF}.
At each coefficient place take a finite extension with good reduction. For this curve the half-period
discriminant calculation in Section~\ref{sec:marked-integral-sections-ordinary-jets}
also gives a good equation directly in every residue characteristic. At an odd residue prime, \(t\) remains nonzero
by injectivity of prime-to-residue-characteristic torsion
reduction. At residue characteristic two, \(2\) splits in \(K\).
The two idempotents of \(O_K\otimes\mathbb Z_2\) split the
\(2\)-divisible group of the good model into groups of height one.
Their dimensions sum to one, so one factor is \'etale.
The étale component of \(t\) is nonzero by the generator
condition. It has nonzero reduction. The CM action used here
extends to the good model by the extension property of abelian
schemes over a discrete valuation ring.

An integral change of variables of a good Weierstrass equation
moves \(t\) to \((0,0)\). The constant term and the linear term
in \(y\) vanish, by the equation and the order-two condition.
Smoothness at \((0,0)\) makes the coefficient of \(x\) a unit.
A shear \(y\mapsto y+s x\), with \(s\) a root of the monic
quadratic that cancels the coefficient of \(x^2\), is integral
after finite extension. A unit dilation then sets the
coefficient of \(x\) to one. This fixes the Weierstrass equation in these coordinates.
Any isomorphism between two such equations preserving the
marked point has no coordinate translations. If its variables
have scaling parameter \(u\), the coefficient of \(x\) forces
\(u^4=1\); thus its invariant differential changes by a unit.
In particular \(\Delta_C\) is a unit at every coefficient place.
Section~\ref{sec:marked-integral-sections-ordinary-jets} computes
it from the half-period discriminants in each ideal class.

Here is also the required integrality of the normalized torsion
values. The good elliptic scheme has the symmetric line bundle
\(\mathcal O([t])\), rigidified by the value at zero of its
canonical section. That value is a basis because the two
sections \(t,0\) remain disjoint. The rigidification on this abelian scheme determines a cubical
structure. For \(L=\mathcal O([t])\) its multiplication formula
identifies \([n]^*L\) with
\(L^{\otimes n(n+1)/2}\otimes([-1]^*L)^{\otimes n(n-1)/2}\)
\cite[Ch.~I, \S2.6(i), pp.~21--22, and Proposition~5.5, p.~29]{MoretBaillyCube}.
The rigidified symmetry therefore gives, over the valuation ring,
\[
 [n]^*\mathcal O([t])\simeq\mathcal O([t])^{\otimes n^2},
\]
since \(t=-t\). If \(u\) represents a torsion point
and \(nu=\lambda\in\Lambda_C\), let \(\kappa(\lambda)\) be the
lattice sign of \(\vartheta_C\). After a finite extension the point
represented by \(u\) extends to a section by properness. Its image under
\([n]\) is the zero section, by separatedness. Pulling back the
rigidified isomorphism to this section thus gives an integral basis
of its \(n^2\)-th tensor power, also when the residue characteristic
divides \(n\).

In the uniformizing basis,
the factor of automorphy at zero for \(\lambda\) is
\(\kappa(\lambda)e^{\pi|\lambda|^2/(2A_C)}\).
The rigidified isomorphism is induced by the identity on
the pulled-back bases over the universal cover. Indeed the
two multipliers at a period agree, since
\(\kappa(n\lambda)=\kappa(\lambda)^{n^2}\), and the identity
preserves the basis at zero. In the fibre at the torsion
point, the \(n^2\)-th power of the canonical section,
in the basis pulled back from the origin, consequently has
value
\[
 \kappa(\lambda)^{-1}
   \left(e^{-\pi |u|^2/(2A_C)}\vartheta_C(u)\right)^{n^2}.
\]
This belongs to the valuation ring. The sign is a unit, so
\((T_u^C\vartheta_C)(0)\) is integral. This also proves that
the analytic normalized value is algebraic for algebraic
CM data, since a positive power of the value is algebraic.

\smallskip
\proofstep{The translated Taylor series}

Write \(C'=C^\sigma\) for the conjugated model with parameter
\(\sigma a_1\), uniformized with that differential.
Write \(x_t(P)=x(P+t)\) and \(y_t(P)=y(P+t)\)
for the functions obtained by translation by the marked
two-torsion point \(t\). They have their only poles at
\(t\), of orders \(2\) and \(3\), and satisfy the
same Weierstrass equation. The translation commutes
with conjugation of the model. The graded ring of
\(\mathcal O([t])\) is generated by
\(\vartheta_C,x_t\vartheta_C^2,y_t\vartheta_C^3\)
in degrees \(1,2,3\). The Weierstrass equation reduces
the powers of \(y_t\) to at most one. The resulting
distinct pole orders give a basis in each degree by
Riemann--Roch. At the pole chart the section
transport is the homogeneous conjugation of the
same projective equation, so it preserves zero
values there too.

The translated rational functions are regular at the
reference origin. Their expansions there conjugate
coefficientwise using the invariant
differential, whose formal integral gives the
uniformizing coordinate. Thus, with \(Q\) the formal quotient of the Taylor series of
\(\vartheta_{C^\sigma}\) by the coefficientwise conjugate of
\(\vartheta_C\), a section of degree \(l\) transports to
\(Q^l\sigma(f)\).
At corresponding points the values of all degree
\(l\) sections have one common nonzero relative
factor \(\mu^l\), by choosing trivializations of the
one-dimensional fibre of the degree-one line bundle.
This evaluation statement is about ordinary section
values, before any Gaussian twist.

Let \([a]\) be a nonzero CM endomorphism.
Because \(\sigma\) fixes \(K\), its two versions
multiply the uniformizing parameter by the same
\(a\), have the same degree \(a\bar a\), and preserve
the conjugate zero divisor of the pullback of the
theta section. In comparing the pulled back
line bundle to the power of the bundle,
use their theta multipliers. The Appell--Humbert
form scales by \(a\bar a\), and the pullback
semicharacter agrees with the power semicharacter
on the period lattice.

One can use
\(a\equiv1\bmod 2\mathcal O_K\) here, for which
this last assertion follows directly from the
half-period multiplier. Thus the section
\(\vartheta_C(az)\) transports to \(\vartheta_{C'}(az)\):
the zero divisors agree, the quotient of the
two sections of the same line bundle is constant,
and both have value one at zero. Consequently
\[
 Q(az)=Q(z)^{a\bar a}.
\]
Take a nonreal \(a\equiv1\bmod 2\mathcal O_K\)
with \(|a|>1\). If \(Q-1\) begins in degree
\(j>0\), coefficient comparison gives \(a^j=a\bar a\).
Absolute values imply \(j=2\), and then \(a=\bar a\),
a contradiction. Hence \(Q=1\).

For a finite set \(S\subset K\), choose positive integers \(N_s\) with
\(N_s s\in\mathcal O_K\) for each \(s\in S\).
Choose a positive integer \(n\) divisible by every \(2N_s\), and set
\(w=\omega_0/(2n^2)\), where \(\Lambda_C=\omega_0\mathcal O_K\).
Choose one lift \(w'\) of the conjugate of the point of \(C\)
represented by \(w\). Then \(v=2n^2w'\) is a period of \(C'\).
Set \(a_s=(n/(2N_s))ns\in\mathcal O_K\). The points represented by
\[
 a_s w=\frac{s\omega_0}{4N_s},\qquad
 a_s w'=\frac{sv}{4N_s}
\]
correspond under \(\sigma\), by equivariance of the \(\mathcal O_K\)-action.
These are the divisions of the two shifts by \(4N_s\)
used in the rigidified translation formula
\cite[Prop.~2.8 in the authors' version]{BannaiKobayashiTheta}.
For \(s\in S\) write \(T_u\) for degree-one Gaussian translation
and \(t_u\) for translation by the corresponding point on the elliptic
curve. On either curve let \(L=\mathcal O([t])\) for its marked point
\(t\). Then \([n]^*(t_u^*L\otimes L^{-1})\simeq\mathcal O\)
for the corresponding shift \(u\): its class in \(\operatorname{Pic}^0\)
is killed by \(n\), since \(ns\in\mathcal O_K\) and \(\omega_0,v\)
are periods on their respective curves. Thus
\((T_{s\omega_0}\vartheta_C)(nz)\) and
\((T_{sv}\vartheta_{C'})(nz)\) are sections of the same
degree-\(n^2\) bundles as the unshifted sections on the
respective curves. The extra period multiplier is
\(\exp(2\pi i\operatorname{Im}(nu\bar\lambda)/A_C)\)
at \(\lambda\in\Lambda_C\), and equals one by integrality of the
Riemann form. The same calculation uses the area of \(C'\) on
its lattice. The zeros of the
first are at points
\(z=(b-2ns)w\) with \(b\in\mathcal O_K\) and
\(nbw\) the marked half-period modulo \(\Lambda_C\).
At the corresponding point \((b-2ns)w'\),
the translated argument is \(nbw'\), which is
the conjugate marked half-period by the algebraic
CM action. Thus these zeros transport to zeros
of the second section. Both have the same degree,
and the zero multiplicity is one because
multiplication by \(n\) is unramified. The
transported section is a nonzero scalar times
the second one.

To determine that scalar, evaluate both at the
point represented before multiplication by \(n\)
by \(-nsw\) and its conjugate \(-nsw'\).
The translated arguments are then \(u=n^2sw\)
and \(u'=n^2sw'\) respectively, and the
translation amounts are \(2u,2u'\).
The exponential of the Gaussian twist is one
at argument \(-u\) of \(T_{2u}\), and the
corresponding unshifted value at \(-u\)
is \(\vartheta_C(u)\) by evenness. The
same is true for \(u'\). Use \(\vartheta_C(nz)\) as the second section of the same
degree \(n^2\); it transports coefficientwise to \(\vartheta_{C'}(nz)\)
by \(Q=1\). At \(-nsw\) the two original sections have the same
value, and likewise at \(-nsw'\) for the two conjugate sections.
The common fibre factor consequently cancels, giving scalar one
when \(\vartheta_C(u)\ne0\). If \(\vartheta_C(u)=0\), both \(u,u'\) are the marked
half-period on their curves. At such a point, the period
multiplier for \(2u\) gives
\[
 \vartheta_C(u+z)
   =\kappa(2u)e^{2\pi z\bar u/A_C}\vartheta_C(u-z).
\]
Comparison of the linear terms at the simple zero gives
\(\kappa(2u)=-1\), and hence \(T_{2u}\vartheta_C(0)=-1\).
The same calculation on \(C'\) gives \(-1\).
Evaluation at zero fixes the scalar also in this case. Dilation by \(n\)
is invertible on formal series over \(\mathbb C\),
so this is the exact Taylor identity for
\(T_{s\omega_0}\) and \(T_{sv}\), for all \(s\in S\).

For this finite \(S\), take a number field \(L\) containing \(K\), the marked
Weierstrass model with its differential and CM action, and the coordinates
of the torsion point represented by \(w\). For an automorphism fixing
this field we have \(C'=C\) and can take \(w'=w\), so \(v=\omega_0\).
The Taylor identity then fixes every coefficient of
\(T_{s\omega_0}\vartheta_C\) for \(s\in S\). These coefficients lie
in \(\mathbb C^{\operatorname{Aut}(\mathbb C/L)}=L\).

To obtain a single finite adele for all rational
translation parameters, write
\(\Lambda_{C'}=\omega_\sigma\mathfrak c^2/N\mathfrak c\)
with \(\mathfrak c\) an integral ideal. Such a presentation
exists because the field discriminant \(-d\) has just one prime
divisor. The ideal and form class groups agree, and genus theory gives
one genus, which consists of squared classes
\cite[Theorems~5.30 and~6.1(iii)--(iv)]{CoxGenus}.
Thus \(\operatorname{Cl}(K)/\operatorname{Cl}(K)^2=1\).
The ideal
\(\mathfrak c^2/N\mathfrak c\) has norm one, and the area
of the lattice is \(|\omega_\sigma|^2\sqrt d/2\).
Put \(e=N\mathfrak c\). Then \(e\mathcal O_K\subseteq\mathfrak c\),
so
\[
 e\mathcal O_K\subseteq \mathfrak c^2/e
                          \subseteq e^{-1}\mathcal O_K.
\]
In particular \(e\omega_\sigma\) is a period and the ratios
\(b_S=v/\omega_\sigma\) of the chosen periods lie in the
compact set \(e^{-1}\widehat{\mathcal O}_K\).

Put \(\mathfrak a_\sigma=\mathfrak c^2/e\), so
\(\Lambda_{C'}=\omega_\sigma\mathfrak a_\sigma\).
The point \(w=\omega_0/(2n^2)\) generates the
\(\mathcal O_K/(2n^2)\)-module \(C[2n^2]\). Conjugation intertwines the \(\mathcal O_K\)-actions on \(C\)
and \(C'\) because \(\sigma\) fixes their derivatives in \(K\).
On \(C'[2n^2]\)
the point \(\sigma(w)\) corresponds to the class of \(b_S\) in
\(\mathfrak a_\sigma/(2n^2)\mathfrak a_\sigma\). Hence \(b_S\)
generates the localized ideal \(\mathfrak a_{\sigma,\mathfrak p}\)
at every \(\mathfrak p\mid2n^2\).
For a given finite \(\mathfrak p\) above a rational prime \(p\),
require \(1/p\in S\). Then \(p\mid n\), so this generator condition
holds on a cofinal tail. The set of local generators is closed:
it is the complement of \(\mathfrak p\mathfrak a_{\sigma,\mathfrak p}\)
inside \(\mathfrak a_{\sigma,\mathfrak p}\), and both lattices are
open and closed. Thus each cluster point \(\beta_\sigma\) of the net
generates \(\mathfrak a_{\sigma,\mathfrak p}\) for every finite \(\mathfrak p\).
In particular \(\beta_\sigma\in\mathbb A_{K,f}^\times\).

For any fixed \(s\), choose \(m>0\) with
\(ms\in\mathcal O_K\). A rational representative
\(x'\) sufficiently close to \(\beta_\sigma s\)
is congruent to \(b_Ss\) modulo
\(2me\widehat{\mathcal O}_K\) for some sufficiently
late \(S\) containing \(s\).
The period \(e\omega_\sigma\) belongs to \(\Lambda_{C'}\),
so the two translation arguments differ by
\(2m\) times a period.

As \(m b_Ss\omega_\sigma\)
is a period, the Riemann form phase for this
difference is the exponential of \(2\pi i\)
times an integer, and the doubled-period
semicharacter is one. The two translated
theta sections are therefore identical.
This proves the required conjugation with a
single \(\beta_\sigma\) and with the asserted
period area. 
In particular, one \(\beta_\sigma\in\mathbb A_{K,f}^\times\) and one
\(\omega_\sigma\) work simultaneously for every rational
translation parameter. No choice of either number depends
on the number of theta factors.

\smallskip
\proofstep{The period scalar under conjugation}

Put \(c_K=\gamma_{\rm Euler}-L'/L(1,\chi_K)\) for this
paragraph and let \(h\) be the class number of \(K\). We
recall the unit normalization of the Chowla--Selberg product.
The squaring map just used is bijective on \(\operatorname{Cl}(K)\),
so \(h\) is odd.
Choose fractional ideals \(\mathfrak c_i\) representing all
classes with product \(O_K\). This can be done by taking
\(\mathfrak b_{i^2}/\mathfrak b_i\) from any system
\(\mathfrak b_i\) of representatives. Choose CM curves \(C_i\) with the normalized Weierstrass
differential used above and with lattices
\(\Lambda_i=\mu_i\mathfrak c_i\Lambda_C\). On each curve the
distinguished point is the generator of \(C_i[2]\) over
\(O_K/2O_K\). The integral normalization above therefore applies
in every ideal class, including at coefficient primes above two.
For each embedding fixing \(K\), algebraic conjugation of the
CM isogenies and their differentials gives
\[
 \Lambda_i^\sigma=\sigma(\mu_i)\mathfrak c_i\Lambda_C^\sigma.
\]
Indeed, for \(0\ne a\in\mathfrak c_i\), multiplication by
\(\mu_i a\) has kernel a finite group of algebraic torsion points.
The quotient isogeny and the ensuing pointed isomorphism of
algebraic elliptic curves are defined over \(\overline{\mathbb Q}\).
On the chosen algebraic differentials the derivative is \(\mu_i a\),
which conjugates to \(\sigma(\mu_i)a\).
This gives inclusion. The inverse assertion follows by
applying the same argument with
\(\mu_i^{-1}\mathfrak c_i^{-1}\) and interchanging the curves.

Each \(\mu_i\) is algebraic since \(\mu_i a\) is the
derivative in algebraic coordinates of a CM isogeny for
\(0\ne a\in\mathfrak c_i\). Set
\(u^{12h}=\prod_i\mu_i^{12}\Delta_{C_i}\).
At any coefficient place the normalized invariant differentials
are unit multiples of good-model differentials by the preceding
normalization. Each \(\mu_i a\), \(a\in\mathfrak c_i\),
is therefore integral as the tangent map of a homomorphism
of good elliptic schemes. Since \(\prod_i\mathfrak c_i=O_K\),
a finite sum of products of such \(a\)'s equals one.
It follows that \(\prod_i\mu_i\) is integral. Interchanging
the curves gives the same assertion for its inverse.
The \(\Delta_{C_i}\) are units and hence so is \(u\).

For any complex lattice \(\Lambda=\omega(\mathbb Z+\tau\mathbb Z)\)
of area \(A\), scaling the discriminant for the periods \(1,\tau\)
\cite[\href{https://dlmf.nist.gov/23.19.E4}{(23.19.4)}]{DLMF} gives
\(\Delta=(2\pi)^{12}\omega^{-12}\eta(\tau)^{24}\).
Kronecker's first limit formula \cite[Ch.~I, \S1, Theorem~1]{SiegelANT} in the Epstein-zeta
normalization \cite[\S1, (1), p.~86]{SelbergChowlaEpstein}, applied to
\[
 \sum_{\ell\in\Lambda\setminus\{0\}}
       \left(\frac{2A/\sqrt d}{|\ell|^2}\right)^s
 = (2/\sqrt d)^s
   \sum_{(m,n)\ne(0,0)}
       \frac{(\operatorname{Im}\tau)^s}{|m+n\tau|^{2s}}
\]
gives residue \(2\pi/\sqrt d\) at one. The constant term in the
unscaled sum is
\[
 2\pi\left(\gamma_{\rm Euler}-\log2
         -\log\bigl(\sqrt{\operatorname{Im}\tau}\,|\eta(\tau)|^2\bigr)\right).
\]
The prefactor contributes \((2\pi/\sqrt d)\log(2/\sqrt d)\)
to the constant term. Using \(A=|\omega|^2\operatorname{Im}\tau\)
and the formula for \(\Delta\), the constant term of the scaled sum is
\[
 \frac{2\pi}{\sqrt d}
       \left(2\gamma_{\rm Euler}
       +\log(2\pi^2/\sqrt d)-\log A
                   -\tfrac16\log|\Delta|\right).
\]
For lattices representing all the classes the sum of
these series is \(2\zeta(s)L(s,\chi_K)\), by grouping
elements of each representative ideal into their unit
orbits. The factor two is \(|O_K^\times|\).
The residues on the two sides give \(2L(1,\chi_K)=2h\pi/\sqrt d\).
Comparing constant terms with \(2\zeta(s)L(s,\chi_K)\)
and dividing by the total residue gives
\[
 \frac1h\sum_i
 \left(\log A_i^\sigma+\tfrac16\log|\sigma(\Delta_{C_i})|\right)
   =\log(2\pi^2/\sqrt d)+c_K .
\]
Here \(A_i^\sigma=|\sigma(\mu_i)|^2N\mathfrak c_i A_C^\sigma\)
and \(\prod_i N\mathfrak c_i=1\), so the left side is
\(\log(A_C^\sigma|\sigma(u)|^2)\). Thus
\[
 A_C^\sigma|\sigma(u)|^2=\frac{2\pi^2}{\sqrt d}e^{c_K}.
\]
The discriminant in this identity is the algebraic
Weierstrass discriminant for the indicated invariant
differential.

Use the splitting character \(\upsilon\) of
Section~\ref{sec:splitting-level-characters}, and set
\(\lambda=\upsilon|\cdot|_E^{-1/2}\).
The prescribed root is \(s_0^f=\lambda_f(\alpha_{\rm W})^{-1}\).
Choose \(h>0\) so that the finite ideal of
\(\alpha_{\rm W}^{h}\) is principal, generated by
\(y\in E^\times\). The character \(\lambda_f\) is trivial on finite integral units.
At infinity \(\upsilon_\infty(y)=\prod_v\Phi_v(y)/|\Phi_v(y)|\)
and \(|y|_{E,\infty}^{1/2}=\prod_v|\Phi_v(y)|\), so
\(\lambda_\infty(y)=\prod_v\overline{\Phi_v(y)}^{-1}\). Hence
\[
 \lambda_f(\alpha_{\rm W})^h
 =\lambda_f(y)=\lambda_\infty(y)^{-1}
 =\prod_v\overline{\Phi_v(y)}
 =\overline{\Phi_K(N_{E/K}y)}.
\]
Thus \(s_0^{fh}\in\Phi_K(K)^\times\), and its roots
have the same modulus under all complex embeddings.
The identity
\(\alpha_{\rm W}\bar\alpha_{\rm W}\widehat O_E
  =\mathfrak d_{E/K}^{-1}\widehat O_E\)
gives \(|\alpha_{\rm W}|^2=N_{E/\mathbb Q}\mathfrak d_{E/K}\).
Discriminant transitivity and \([E:K]=f\) give
\[
 |\sigma(s_0)|^2=|\alpha_{\rm W}|^{1/f}
     =\operatorname{rd}(E)/\sqrt d.
\]
Set \(p_F=s_0/u\), \(a_d=\sqrt{2\pi/\sqrt d}\), and
\(\gamma=p_Fa_d/\omega_0\), where
\(\Lambda_C=\omega_0 O_K\).
For the normalized period \(\omega_\sigma\) determined above,
\(A_C^\sigma=|\omega_\sigma|^2\sqrt d/2\). Therefore the
transported jet scale \(t_\sigma=\sigma(p_F)a_d/\omega_\sigma\)
satisfies the exact identity
\[
 |t_\sigma|^2
   =\frac{\pi|\sigma(s_0)|^2}{A_C^\sigma|\sigma(u)|^2}
   =\frac{\operatorname{rd}(E)}{2\pi}e^{-c_K}
   =\Omega^{-1}.
\]

\smallskip
\proofstep{The cyclotomic product in the finite Weil model}

Put \(\mathfrak a=E\cap\alpha_{\rm W}\widehat O_E\). At
each finite place \(\alpha_{\rm W}\) is \(1\) or the inverse of a
generator of the relative different
(Section~\ref{sec:period-finite-coordinates}), so \(O_E\subseteq\mathfrak a\).
Since \(\mathfrak a\bar{\mathfrak a}=\mathfrak d_{E/K}^{-1}\), we have
\(\operatorname{Tr}_{E/K}(\mathfrak a\bar{\mathfrak a})\subseteq O_K\).
In particular
\(\operatorname{Tr}_{E/K}(\mathfrak a\,\overline{\zeta_q^k})\subseteq O_K\)
for each \(k\).
For the standard definite space of any positive rank \(v\),
write \(D=fv\) and use Fock coordinates
\[
  z=(z_{\tau,i})_{\tau\in\Phi,\,1\le i\le v},
  \qquad z(\gamma)_{\tau,i}=a_d\tau(\gamma_i)
                   \quad(\gamma\in E^v).
\]
The lattice \(\mathfrak a^v\) in these coordinates is
self-dual for the alternating polarization of the standard
Fock form.

Indeed, up to sign the alternating pairing
\(\operatorname{Im}\langle z(\gamma),z(\delta)\rangle/\pi\)
is \(\operatorname{Tr}_{E/\mathbb Q}(j_{\rm W}\sum_i\bar\gamma_i\delta_i)\),
and \(j_{\rm W}\mathfrak a\bar{\mathfrak a}
  =\mathfrak d_{E/\mathbb Q}^{-1}\).
Let \(\Lambda=\widehat O_E^{\,v}\) be the normalized finite
translation lattice. The finite Heisenberg translate \(T\) by
\(\alpha_{\rm W}^{-1}\gamma\) corresponds to the signed
real shift by \(-z(\gamma)\). Indeed, a rational
translation has real and finite coordinates
\((z(\gamma),\alpha_{\rm W}^{-1}\gamma)\), and translating
it out of the adelic theta functional gives the real shift
\[
     F(z)\longmapsto
      e^{\langle z(\gamma),z\rangle-|z(\gamma)|^2/2}
                         F(z-z(\gamma))
\]
on the entire function obtained by dividing its real slice by
\(e^{-|z|^2/2}\).

Here the displayed shift is for the
unnormalized finite translate. For \(\gamma\in\mathfrak a^v\),
the normalization is the factor
\(e^{-\pi i\operatorname{Tr}_{E/K}\sum_i\gamma_i\bar\gamma_i}\).
It follows from the product formula for the additive
character with trace-one element \((1+\sqrt{-d})/2\).
The exponent is \(-\pi i\) times an integer by the trace
pairing, so the sign is one on any even multiple of
\(\mathfrak a^v\).

Let \(\omega\) denote \(\omega_0\) or \(\omega_\sigma\), and
use the corresponding rank-one theta curve. Take \(e=1\)
in the first case and the positive integer \(e\) fixed above
in the second, so \(\omega e\) is a period.
For \(0\le k<q\) define
\[
 A_{k,i}(z)=\frac{\omega}{q a_d}
       \sum_\tau\overline{\tau(\zeta_q^k)}z_{\tau,i},
 \qquad
 c_k(y)=-\frac{1}{q}\operatorname{Tr}_{E/K}
                         (y\overline{\zeta_q^k}).
\]
Then \(A_{k,i}(-z(\gamma))=\omega c_k(\gamma_i)\).
The geometric-series identity on the distinct \(q\)-th
roots \(\tau(\zeta_q)\) gives
\[
 \sum_{0\le k<q}\tau(\zeta_q^k)
                      \overline{\rho(\zeta_q^k)}
                       =q\,\mathbf1_{\tau=\rho}.
\]
As \(A_C=|\omega|^2\sqrt d/2\), the pullbacks of the
degree-one Hermitian form \(\pi\bar u z/A_C\) along the
\(A_{k,i}\), each repeated \(q\) times, sum to the standard
Fock Hermitian form.

The entire function
\[
               \Theta(z)=\prod_{k,i}\vartheta_C(A_{k,i}(z))^q
\]
therefore has the required Gaussian multiplier on the
lattice \(z(qe\mathfrak a^v)\). At any point of that lattice,
each \(A_{k,i}\) is a period by trace integrality and the
choice of \(e\). Each rank-one sign is raised to the even
power \(q\). Their product is one, in agreement with the
Fock sign on this lattice. Also \(\Theta(0)=1\).

To identify this function in the original Weil model, let
\(l=qe\). The translates of the unit vacuum by
\(l^{-1}\Lambda/\Lambda\) are pairwise orthogonal, by its
vacuum matrix coefficient, and are fixed by the normalized
translations by \(l\Lambda\), by the Heisenberg commutator.
Their entire functions are linearly independent by the
injectivity in Section~\ref{sec:marked-integral-sections-ordinary-jets}.
That argument transports the real evaluations to arbitrary
finite displacements by rational Heisenberg invariance and
uses the polarized norm identity.

These entire
functions belong to the polarized theta space for
\(z(l\mathfrak a^v)\) with the displayed signed multiplier.
That space has dimension \(l^{2D}\): the original
alternating polarization is unimodular and passing to
\(l\mathfrak a^v\) multiplies every alternating elementary
divisor by \(l^2\). For this polarization type the space of sections
is the irreducible theta representation on functions on
\((\mathbb Z/l^2\mathbb Z)^D\)
\cite[\S1, Proposition~3 and Theorem~2, pp.~295--297]{MumfordTheta}.
Its dimension is \(l^{2D}=[l^{-1}\Lambda:\Lambda]\).
Thus \(\Theta\) is the entire function of a vector in
the original finite part; denote it by \(F\) or \(F'\)
according as the curve is original or conjugate.

For a rational translate indexed by \(\gamma\in E^v\),
the Gaussian shift identity and the equality of Hermitian
forms give the entire function
\[
 \prod_{k,i}
   (T_{\omega c_k(\gamma_i)}^C\vartheta_C)
                        (A_{k,i}(z))^q .
\]
Its constant coefficient is a product of integral torsion
values, since each \(c_k(\gamma_i)\) lies in \(K\) and a
multiple of \(\omega\) is a period. Finite smoothness and
density of rational points give the same conclusion at
every finite adelic translate.

For the original curve \(e=1\). Write
\[
 F=\sum_{x\in q^{-1}\Lambda/\Lambda}c_x T(x)v_0
\]
using representatives supported above \(q\). Since \(\Theta(0)=1\),
choose \(x_0\) with \(c_{x_0}\ne0\) and set \(y=-x_0\).
For \(u\in\Lambda\), the normalized translation \(T_\sharp(u)\) fixes
\(v_0\). With \(t=x-x_0\), its eigenvalue on \(T(y)T(x)v_0\) is the
commutator
\[
 \psi\bigl((u^{\mathsf T}\bar t-t^{\mathsf T}\bar u)/2\bigr)
             =\psi(u^{\mathsf T}\bar t),
\]
since the finite character \(\psi\) is trivial on the fixed subalgebra.
This character of \(u\) is trivial on \(q\Lambda\), since
\(t\in q^{-1}\Lambda\), and on all of \(\Lambda\) precisely when
\(t\in\Lambda\), by self-duality. Thus \(T(y)F\) is fixed by
\(T_\sharp(q\Lambda)\). For \(x=x_0\) one has \(T(-x_0)T(x_0)=1\),
again because the cocycle argument lies in the fixed subalgebra.
Character orthogonality now gives
\[
 \sum_{u\in\Lambda/q\Lambda}T_\sharp(u)T(y)F
       =[\Lambda:q\Lambda]c_{x_0}v_0=q^{2D}c_{x_0}v_0 .
\]
The representatives for \(u\) can also be supported above \(q\).
Thus both \(u\) and \(y\) vanish above \(MN\), since \((q,MN)=1\).
For \(v=1\) the result is the identity
\eqref{eq:upper-CM-vacuum} with exactly \(q^{2f}\) terms. It defines the
original nonzero scalar \(s_{\rm vac}\), with the
integrality and uniform conjugate sup bound obtained from
the same formula. The nonzero projection also shows that
the translates of \(F\) span the entire finite part.

Finally put \(b_\sigma\) equal to the image of
\(\beta_\sigma\) in \(\mathbb A_{E,f}\). At a rational
translate \(\gamma\), choose one element \(a\in K\) sufficiently
close to \(\beta_\sigma\) that the translate of \(F'\) by
\(a\alpha_{\rm W}^{-1}\gamma\) agrees with its translate by
\(b_\sigma\alpha_{\rm W}^{-1}\gamma\), and that every
\(a c_k(\gamma_i)\) belongs to the neighborhood required
by the rank-one conjugate translation identity. These are
finitely many open conditions on \(a\), and \(K\) is dense
in its finite adeles.

The coefficient of the \(\tau,i\)
variable in the original scaled form is
\[
  \frac{p_F}{q}\overline{\tau(\zeta_q^k)}.
\]
After applying \(\sigma\) it is the coefficient in the
new form at scale \(t_\sigma\), with the permutation of
the embeddings over \(K\). Indeed
\(\overline{\tau(\zeta_q^k)}=\tau(\zeta_q^{-k})\), and
\(\omega_\sigma t_\sigma/a_d=\sigma(p_F)\).
Apply the rank-one translated Taylor identity to each
factor and compose with these linear forms. This proves
\eqref{eq:upper-CM-jet-transport} at rational translates and hence at all finite
adelic translates by local constancy. The same
\(b_\sigma\) and \(t_\sigma\) work for every \(v\), in
particular for \(v=1\) and for \(v=n_jn_{j+1}\) in the
same upper kernel.

\proofstep{The finite Weil action and vacuum projection}
Let \(V=\operatorname{Hom}(X_j,X_{j+1})\) in the definite models. Write \(v_0\) for its unit vacuum,
\(\Lambda=\widehat{\mathcal O}_E^{\dim_E V}\) for the
standard lattice in the normalized Heisenberg coordinates,
and \(\psi\) for the self-dual finite additive character in
those coordinates. Fix \(\sigma\in\operatorname{Aut}(\mathbf C/K)\).
The translated CM identity, first in rank one and then on
the cyclotomic products, gives initial vectors \(F,F'\)
and a finite adele \(b_\sigma\) as in \eqref{eq:upper-CM-jet-transport}.

First recover the full additive character. In rank one
write the unnormalized finite translation as \(T(x)\).
Its group law in these coordinates is
\[
 T(x)T(y)=\psi(x\bar y/2)T(x+y).
\]
The two seeds are smooth in the finite translation variables.
Choose a positive integer \(n\) for which \(T(n)F=F\)
and \(T(b_\sigma n)F'=F'\). Such an integer exists because
the two open stabilizer neighborhoods contain an open
additive lattice and a sufficiently large integer clears the
finitely many local conditions on \(1,b_\sigma\).

For any \(z\in\mathbb A_{E,f}\) write \(z=nx\).
Apply the group law to \(T(2x)T(n)F\) and
\(T(2b_\sigma x)T(b_\sigma n)F'\). The two
cocycle factors are \(\psi(nx)\) and
\(\psi(b_\sigma\overline{b_\sigma}nx)\), respectively. Both comparisons involve
the same conjugated jet of the translate at \(2x+n\), which
is nonzero. Indeed the full jet at a nonzero scale determines the
entire theta function by holomorphy, and hence determines
the finite vector by the same injectivity argument.
Cancellation therefore
gives
\[
 \sigma(\psi(z))=\psi(b_\sigma\bar b_\sigma z)
       \qquad(z\in\mathbb A_{E,f}).
 \tag{JV1}\label{eq:upper-finite-vacuum-action}
\]
The comparison at the two translations uses the full cocycle in the group law.

Put \(u=b_\sigma\bar b_\sigma\).
Self-duality of \(\widehat{\mathcal O}_E\) for \(\psi\), applied
to \eqref{eq:upper-finite-vacuum-action} at all \(zy\) with \(y\) integral, gives
\[
 uz\in\widehat{\mathcal O}_E\quad\Longleftrightarrow\quad
 z\in\widehat{\mathcal O}_E.
\]
Indeed the finite character in these coordinates is the absolute
trace character scaled by \(\delta_f=j_{\rm W}\alpha_{\rm W}\bar\alpha_{\rm W}\).
At \(w\mid p\) the trace dual of \(\mathcal O_w\) is
\(\mathfrak D_{E_w/\mathbb Q_p}^{-1}\). Hence the annihilator of
\(\mathcal O_w\) for \(\psi_w(xy)\) is
\(\delta_{f,w}^{-1}\mathfrak D_{E_w/\mathbb Q_p}^{-1}=\mathcal O_w\)
by Section~\ref{sec:period-finite-coordinates}. At \(p=d\) the
cyclotomic extension over \(K\) is unramified, so this scale is
supplied by \(j_{\rm W}\). At \(p\mid q_1\) the quadratic field is
split, \(j_{\rm W}\) is a unit and
\(\alpha_{\rm W}\bar\alpha_{\rm W}\) generates the inverse different.
Testing in each local factor at \(1\) and at the inverse
of a uniformizer proves \(u\in\widehat{\mathcal O}_E^\times\).
In particular \(b_\sigma\) is an idele, with inverse
\(\bar b_\sigma u^{-1}\).

For the higher rank product, use the relation \eqref{eq:upper-finite-vacuum-action} and
the translates of its seed \(F\). They span the finite part.
If a linear combination of them is zero, the
\(\sigma\)-conjugated combination of the translated \(F'\)
has zero full jet after the common variable permutation.
The injectivity just proved shows that combination is zero.
Hence the rule defines a \(\sigma\)-semilinear map \(C\)
on the finite part. The group law and \eqref{eq:upper-finite-vacuum-action} show
\(CT(x)=T(b_\sigma x)C\) on every generator and hence on
the whole finite part. The normalized translations
\(T_\sharp(x)\) have the additional factor
\(\psi(c_{\mathrm{tr}}\bar x^{\mathsf t}x/2)\) in this model, where
\(c_{\mathrm{tr}}=(1+\sqrt{-d})/2\) is the integral trace-one
element \((1-j_{\rm W}^{-1})/2\) of Section~\ref{sec:overlap-finite-character}. Equation \eqref{eq:upper-finite-vacuum-action} conjugates
that factor to the factor at \(b_\sigma x\). Thus the same
identity holds with \(T_\sharp\). The full jet identity gives injectivity of \(C\), since
\(\gamma\ne0\) and the jets at that original scale determine
the vector.

There is \(a\in\widehat{\mathcal O}_E^\times\) with
\(a\bar a=u\). In each split pair choose one place \(w\) and set
\(a_w=1\), \(a_{\bar w}=\overline{u_w}\), using the involution
between the completions. At each nonsplit finite place \(w\),
\(2v_w(b_{\sigma,w})=v_w(u_w)=0\); take \(a_w=b_{\sigma,w}\). Thus
\(b_\sigma=h a\) for the finite norm-one scalar \(h=b_\sigma/a\).
Set \(C^\circ=\omega(h)^{-1}C\); it conjugates the normalized
translations by \(x\) to those by \(a x\).

Since \(a\Lambda=\Lambda\), the vector \(C^\circ v_0\) is fixed by
all \(T_\sharp(x)\) with \(x\in\Lambda\). Its fixed line is
one-dimensional. To see this, the translates of \(v_0\)
indexed by distinct cosets of \(\Lambda\) are orthogonal
by the support of the vacuum matrix coefficient.
Each coset is an eigenspace for all lattice translations.
Its character is trivial only for the zero coset, by the
self-duality of \(\Lambda\) for the polarized Heisenberg
pairing. Therefore \(C^\circ v_0=\eta v_0\) with \(\eta\ne0\).
Equivalently \(Cv_0=\eta\omega(h)v_0\). For \(p^e\Vert q_1\) and
\(\Lambda_p=\Lambda\otimes_{\widehat{\mathbf Z}}\mathbf Z_p\), on the
subspace fixed by \(T_\sharp(q_1\Lambda_p)\) the local vacuum average is
\[
 P_p=p^{-2De}\sum_{u\in\Lambda_p/q_1\Lambda_p}T_\sharp(u),
 \qquad C^\circ P_p=P_p C^\circ .
\]
Indeed \(a_p\) permutes the summands and the rational coefficient is
fixed by \(\sigma\). The phase in each normalized summand transports by
\eqref{eq:upper-finite-vacuum-action}, also at \(p=2\). Each average
projects onto the subspace fixed by \(\Lambda_p\).

For a finite unitary matrix \(g\), choose \(H\) with
\(H\Lambda=\Lambda\cap g^{-1}\Lambda\).
The original finite lattice multiplier is
\(\Delta_g=\lambda_f(\det H)\), with
\(\lambda_f=\upsilon_f|\cdot|_f^{-1/2}\), as in the
lattice model of Section~\ref{sec:seesaw-integral-tensor}.
Averaging the lattice
projection over \(\Lambda\) gives the finite formula
\[
 \omega(g)v_0=
  \frac{\Delta_g}{|\Delta_g|^2}
  \sum_{q\in\Lambda/(\Lambda\cap g^{-1}\Lambda)}
                            T_\sharp(gq)v_0 .
\]
Changing a representative \(q\) by \(q'\) in
\(\Lambda\cap g^{-1}\Lambda\) leaves the term unchanged.
Write \(c=c_{\mathrm{tr}}\), \(A=y^*x\) and
\(\eta(x)=\psi(cx^*x/2)\). The unnormalized cocycle is
\(\psi(A/2)\). Since \(\psi(\bar z)=\psi(z)^{-1}\) and
\(c+\bar c=1\), the normalized cocycle is
\[
 \frac{\eta(x)\eta(y)}{\eta(x+y)}\,\psi(A/2)
 =\psi\bigl((A-c(A+\bar A))/2\bigr)=\psi(\bar c A).
\]
For \(x,y\in\Lambda\), \(\bar cA\) is integral and the
self-dual character is trivial there, including above \(2\). Unitarity of \(g\)
preserves the Hermitian form in both the Heisenberg cocycle and
the normalizing quadratic phase. Thus
\(T_\sharp(gq+gq')=T_\sharp(gq)T_\sharp(gq')\), and the second
operator fixes \(v_0\) because \(gq'\in\Lambda\).
This argument includes the places over \(2\). Each member of
the finite quotient has mass \(|\Delta_g|^{-2}\), since
\(|\Delta_g|^2\) is its index.

Multiplication by \(a\) permutes
these cosets, because \(a\) is a scalar integral unit.
Applying \(C^\circ\) therefore gives
\[
 C^\circ\omega(g)v_0=
   \frac{\sigma(\Delta_g)}{\Delta_g}\,\omega(g)C^\circ v_0.
\]
Heisenberg covariance extends this to the spanning translates
of \(v_0\). The scalar matrix \(h\) commutes with \(g\), so
the same identity holds with \(C\). The translates of
\(Cv_0=\eta\omega(h)v_0\) span, since \(b_\sigma\) is invertible
and \(\omega(h)\) is a unitary change of variables. Hence \(C\) is
bijective. With \(d_\sigma(g)=\sigma(\Delta_g)/\Delta_g\), its inverse
on the vacuum is
\[
 C^{-1}v_0=\sigma^{-1}\!\bigl((\eta d_\sigma(h^{-1}))^{-1}\bigr)
                                    \,\omega(h^{-1})v_0 .
\]
This follows by applying \(C\) to the right-hand side.
The lattice for \(g^{-1}\) is
\(\Lambda\cap g\Lambda=gH\Lambda\), so
\(\Delta_{g^{-1}}=\upsilon(\det g)\Delta_g\).
In the standard coordinates \(g^{-1}=g^*\). Local Smith
form therefore conjugates the determinant ideal of
\(\Lambda\cap g^{-1}\Lambda\) to that of
\(\Lambda\cap g\Lambda\). As \(\lambda_f\) is unramified
and takes conjugate values on conjugate ideles,
\[
 \Delta_{g^{-1}}=\overline{\Delta_g},\qquad
 \Delta_g^2=|\Delta_g|^2\upsilon(\det g)^{-1}.
\]
The index is rational. The finite values of the splitting
character on norm-one ideles have modulus one at every
complex embedding.

For the latter assertion, a positive
power of a norm-one idele can be written as a rational
norm-one scalar times an integral norm-one idele. Indeed,
if its ideal is \(I\) and \(I^d=(\alpha)\), then
\(I^{2d}=(\alpha/\bar\alpha)\) since \(I\bar I=1\).
The finite character of the integral factor has finite image.
The value at the rational factor is the inverse archimedean
weight monomial in its conjugates, all of absolute value one.
This gives
\[
 \left|\frac{\sigma(\Delta_g)}{\Delta_g}\right|=1.
\]
The same argument applies to the finite splitting
characters of the pair action. Finite Heisenberg phases
are roots of unity.

Write \(v_{\rm int}\) for the integral vacuum and
\(Cv_{\rm int}=\eta_{\rm int}\omega(h)v_0\). The degree-zero
coefficient of the translated jet is the theta evaluation at the
corresponding finite displacement. Thus \eqref{eq:upper-CM-jet-transport}
and semilinearity give
\[
 (Cv_{\rm int})(0,\alpha_{\rm W}b_\sigma x)
       =\sigma\bigl(v_{\rm int}(0,\alpha_{\rm W}x)\bigr)
       \qquad(x\in\mathbb A_{E,f}^{\dim_E V}).
\]
Multiplication by \(\alpha_{\rm W}b_\sigma\) is bijective on
the finite adeles, so this covers the entire finite slice.

\proofstep{The Schur-jet integral}
For a rational Heisenberg displacement \(r\), equivariance gives
\[
 |(Cv_{\rm int})(r_\infty,z_f)|
  =|(Cv_{\rm int})(0,z_f-r_f)|
            \qquad(z_f\text{ finite}).
\]
Rational combinations of a \(\mathbf Q\)-basis of the
position space are dense in the real position space.
This and continuity give the same supremum bound
on the whole additive quotient. Its Haar measure has mass one,
and the polarized norm identity of
Section~\ref{sec:overlap-fourier-inversion} identifies the
mean square with \(\|Cv_{\rm int}\|^2=|\eta_{\rm int}|^2\)
since \(\omega(h)\) is unitary and \(\|v_0\|=1\). The Gaussian
lattice bound on \(\omega(g)v_0\) now applies to
\(C(T(x)\omega(g)v_{\rm int})\), uniformly in the old
translation \(x\) and the old Weil argument \(g\).
This is the pointwise bound for the coefficientwise conjugate
Schur-jet integral. In particular this step concerns the
upper kernels before any source form is substituted.
For fixed finite arguments, the corresponding jet coefficients are
algebraic. Indeed the displayed finite formula for \(\omega(g)v_0\)
also holds after multiplying both sides by the vacuum scale.
Equation~\eqref{eq:upper-CM-vacuum} in each rank-one factor then
expresses every term through translated cyclotomic seeds. Finite
smoothness permits rational translation representatives. Their
rank-one Taylor coefficients belong to the number fields obtained
above, and the scaled linear forms have coefficients
\(p_F\overline{\tau(\zeta_q^k)}/q\). The other coefficients are
Heisenberg roots of unity and the algebraic factors
\(\Delta_g/|\Delta_g|^2\) and finite splitting-character values.

For the resulting estimate, let \(J_\gamma F\) be the
Taylor series at the original jet scale \(\gamma\). For
a polynomial \(Q\), write \(j_Q(J_\gamma F)\) for
\(Q(\partial)(J_\gamma F)(0)\). The Schur kernel is tested
by taking \(Q\) to be a fixed matrix entry, or the
contraction with a source vector and a dual target vector.
For any additional scalar \(c\), put
\[
 Q_\sigma
 =\pi_\sigma\!\cdot
   \bigl(\sigma_{\mathrm{coeff}}(Q(cX))\bigr),
\]
where \(\pi_\sigma\) acts by permuting the embedding
coordinates only. Relabel each source and target Weyl factor
at \(\tau\), together with its diagram, by \(\pi_\sigma(\tau)\).
This permutes the integral bases of pure tensors. For either
space, write \(P\) for the resulting permutation matrix and
\(J_{\rm old},J_{\rm new}\) for the tensor products of the
factor Gram matrices in the two orders. Then
\(P^{\mathsf t}J_{\rm new}P=J_{\rm old}\).

The chain rule for monomials and the
transported Taylor identity give
\[
 \sigma\!\bigl(j_Q(J_{c\gamma}F)\bigr)
 =j_{Q_\sigma}(J_{t_\sigma}CF)
 =\int_{\mathbf C^D}
     Q_\sigma(t_\sigma\bar y)
     e^{-\|y\|^2/2}\,
     CF(y_{\mathbb A})\,\frac{dy}{\pi^D}.
\]
Here \(y_{\mathbb A}\) has finite component zero and archimedean components
\((2\pi\operatorname{Im}\tau(j_{\rm W}))^{-1/2}y_\tau\)
in the original adelic coordinates. Indeed
\(2\pi\operatorname{Im}\tau(j_{\rm W})=a_d^2\) on \(\Phi\), so
\(y\) is the normalized Fock coordinate. Write
\(H(y)=e^{\|y\|^2/2}CF(y_{\mathbb A})\) for the entire section.
For the fixed finite test its weighted modulus is bounded: a full
period lattice of the real slice acts by unit phases on the weighted
values, by rational invariance and the finite translation stabilizer.
For each multi-index \(\nu\), angular integration on the product of
circles extracts \([z^\nu]H(z)\). Radial integration in each coordinate
gives \(2\int_0^\infty r^{2\nu_i+1}e^{-r^2}\,dr=\nu_i!\). Thus
\[
 \pi^{-D}\int_{\mathbb C^D}\bar y^\nu H(y)e^{-\|y\|^2}\,dy
      =\nu!\,[z^\nu]H(z)=\partial^\nu H(0).
\]
Indeed, if \(B_F=\sup_y|CF(y_{\mathbb A})|\), the absolute
integrand after rescaling a monomial is bounded by
\[
 B_F |t_\sigma|^{|\nu|}\,\|y\|^{|\nu|}e^{-\|y\|^2/2}.
\]
This is integrable, and the Taylor series converges uniformly on
each product of circles. Rescaling by \(t_\sigma\) supplies
\(t_\sigma^{|\nu|}\), and finite
linearity gives the Schur-jet identity. Here
\(|t_\sigma|^2=\Omega^{-1}\) at the transported original
scale, and for the actual arrow scale \(c=\pi_M=\bar\alpha\) one has
\(|\sigma(c)|^2=M\). The variable in the conjugated
polynomial consequently has squared norm
\((M/\Omega)\|y\|^2\).

For the short-image bound, apply the measured operator
identity on the original ramified source first. Then apply
\(\sigma\) to the terminal coordinates of its finite
upper-chain class sum. This coefficientwise operation
changes neither the actual source generators nor the
choice of measured point and functional. In particular it
requires no continuity of \(\sigma\) in the real topology or
conjugation of the integrals defining the Poincar\'e map.
The same measured identity separately gives the local lattice norm
comparison for every allowed twist on the clearing module.

Since \(K\) is imaginary quadratic, an arbitrary \(\sigma\)
either fixes \(K\) or \(c\circ\sigma\) fixes \(K\), where \(c\)
is ordinary complex conjugation. The Weyl Gram matrices
are rational. Entrywise ordinary conjugation preserves the
operator norm between these two positive Gram spaces.
Thus the bound obtained for the embeddings fixing \(K\)
also gives the bound at every complex automorphism.
\end{proof}

\subsubsection{Gaussian bounds for the conjugated kernels}\label{sec:upper-gaussian-jet}

\sectionopening{In this subsection, we apply the vacuum identity of
Lemma~\ref{lem:upper-cm-transport} and integrate the conjugated jets.
The resulting operator norm is used in Section~\ref{sec:upper-terminal-norm}.}

\begin{lemma}[Norm of one conjugated kernel]\label{lem:upper-gaussian-kernel}
For \(0\le j<m\), let \(D_j(g,h;u)\) be the Schur-jet matrix
at scale \(\pi_M\gamma\) of \(\omega_{\rm f}(g,h)T_\sharp(u)v_{\rm int}\)
for the untwisted arrow \(X_j\to X_{j+1}\), with normalized integral
vacuum \(v_{\rm int}\). For sufficiently large \(x\), every finite
adelic source and target pair \((g,h)\), and every finite Heisenberg
translation \(u\),
\[
 \|\sigma D_j(g,h;u)\|_{\rm op,W}\le\exp(x^{50}+S_j/2)
               \qquad(\sigma\in\operatorname{Aut}(\mathbb C)).
\]
\end{lemma}

\begin{proof}
For a conjugated translated vacuum the bound is
\(e^{C_j}\), where
\[
 C_j=C_K(q_1 n_j n_{j+1}+1)^2,
 \qquad n_j=\dim_E X_j.
\]
It is uniform in the finite translates. The preceding
transport gives a scalar multiple of a unitary move of
the vacuum and bounds the scalar by its finite values.
For this last finite bound, the
normalized CM translation law used in \eqref{eq:upper-CM-jet-transport} sends every conjugate of a
torsion value to a value on a conjugate elliptic lattice.

After multiplying by the Gaussian part of the translation,
its absolute value is periodic on that lattice. Its maximum
is finite and unchanged by homothety. The finitely many
lattice classes and the two-torsion choices therefore
give a common bound \(M_K\ge1\) depending only on \(K\).
Each finite value of the rank-one theta product has
\(q_1^2\) such factors, by the repeated cyclotomic resolution.

Equation \eqref{eq:upper-CM-vacuum} expresses this same multiple of the vacuum
as a sum of \(q_1^{2f}\) translates with root-of-unity
factors. Hence the finite values of the scaled rank-one vacuum
have conjugate absolute value at most
\[
          q_1^{2f} M_K^{q_1^2}
                \le \exp\bigl((2+\log M_K)q_1^2\bigr).
\]
Tensoring over the \(n_jn_{j+1}\) rank coordinates
multiplies these bounds. The inner product is the mean of
the squared values on the vacuum quotient.

By rational
translation, arguments whose archimedean coordinates come
from rational vectors can be moved to arguments of
archimedean coordinate zero, up to unit phases. Density
and continuity give the same bound for the quotient norm,
which bounds the scalar in the unit vacuum presentation.

For the ordinary moved vacuum, put \(f_y=T(y)\omega(g)v_0\).
In rational Schr\"odinger coordinates write \(z=(a,b)\)
and \(\gamma=(u,v)\), with the two variables in dual
Lagrangians. Let \(\mathrm e\) denote the global additive
character, and take the Fourier coefficient at the character
\(\mathrm e(\langle b,u\rangle-\langle v,a\rangle)\).
The phase independent of the rational summation variable
drops out of
\[
 |\theta(T(a,b)f_y)|^2
  =\left|\sum_{\lambda}\mathrm e(\langle b,\lambda\rangle)
                      f_y(\lambda+a)\right|^2 ,
\]
where \(\lambda\) ranges over the rational points of the
first Lagrangian. Integration in \(b\) on its rational
quotient leaves the pairs \(\lambda-\mu=u\).
Unfolding the integral in \(a\) gives
\[
 \int f_y(t+u)\overline{f_y(t)}
                    \mathrm e(\langle v,t\rangle)\,dt
                =\langle f_y,T(u,v)f_y\rangle .
\]
In the commutator convention
\(T(z)T(\gamma)=\beta(z,\gamma)^2T(\gamma)T(z)\), the character just extracted is
\(\beta(z,\gamma)^{-2}\), as in Section~\ref{sec:overlap-fourier-inversion}.
The half-translation phase is one on the rational pair
\((u,v)\). The two quotient measures have mass one in
the self-dual normalization. The Schwartz bounds give
absolute summability on compact sets of adelic
representatives and justify the integrations.
The zero Fourier coefficient gives the quotient-mean
normalization of the squared norm used above.

By unitarity, the absolute value of this
coefficient equals the ordinary vacuum matrix coefficient
in modulus at \(g^{-1}\gamma\). The moved vector \(y\)
only contributes its commutator phase. The Gaussian
vacuum formula bounds the coefficient in the support by
\(\exp(-\pi\|g^{-1}\gamma\|_\infty^2/\sqrt d)\);
the coefficient vanishes outside the finite support
\((g^{-1}\gamma)_f\in\alpha_{\rm W}\widehat{\mathcal O}_E^{\operatorname{rk}_E V}\).
For a difference \(\gamma\) of two frequencies in this support,
the finite unitary change of basis gives
\[
 h(\gamma,\gamma)_f=h(g_f^{-1}\gamma,g_f^{-1}\gamma)
     \in\mathfrak d_{E/K}^{-1}\widehat{\mathcal O}_E.
\]
Here \((g_f^{-1}\gamma)_f\) lies in the same additive lattice, and
\(\alpha_{\rm W}\bar\alpha_{\rm W}\widehat{\mathcal O}_E
\subseteq\mathfrak d_{E/K}^{-1}\widehat{\mathcal O}_E\).

Take the unscaled embedding coordinates
\(u_{\tau,i}=(g_\infty^{-1}\gamma)_{\tau,i}\).
Unitarity and trace integrality give
\[
 \|g^{-1}\gamma\|_\infty^2=\sum_{\tau,i}|u_{\tau,i}|^2
  =\operatorname{Tr}_{E/K}h(\gamma,\gamma)
  \in\mathcal O_K\cap\mathbb R=\mathbb Z.
\]
It is positive if \(\gamma\ne0\). Thus this lattice is one-separated.
Its Fock coordinates are \(a_d u\), with \(a_d^2=2\pi/\sqrt d\),
as used in the Gaussian above.
Packing unit-separated points bounds the sum of
\(e^{-(\nu+1)\|v\|}\) in dimension \(\nu\) by \(3^\nu e^\nu\).
The counting function in real dimension \(\nu\) is at most
\((1+2r)^\nu\le 3^\nu e^{\nu r}\) for a ball of radius \(r\).

Integrating against \((\nu+1)e^{-(\nu+1)r}\,dr\) gives the
claimed packing bound. Completing the square gives
\(e^{-\pi\|v\|^2/\sqrt d}
 \le e^{\sqrt d(\nu+1)^2}e^{-(\nu+1)\|v\|}\).
Thus the Fourier series is absolutely convergent and its
value at zero is \(|\theta(f_y)|^2\). Since
\(\nu=2f\operatorname{rk}_E V\le 2q_1\operatorname{rk}_E V\),
its square root is bounded by
\(\exp((2d+6)(q_1\operatorname{rk}_E V+1)^2)\).
Writing \(r_V=n_jn_{j+1}\ge1\), the logarithm of the
transported vacuum bound is at most
\[
 (2+\log M_K)q_1^2r_V+(2d+6)(q_1r_V+1)^2
 \le (2d+8+\log M_K)(q_1r_V+1)^2 .
\]
Thus one may take \(C_K=2d+8+\log M_K\).

To calculate the Gaussian moment, use \(L_a=a+1\) and
\(v_M=\log(M/\Omega)\). For an array \(Y\) in the Gaussian integral,
put \(t_Y=\|Y\|^2/(2L_a+2)\). First suppose \(\sigma\) fixes \(K\).
The Schur-jet integral in Lemma~\ref{lem:upper-cm-transport} evaluates
the original rational Schur matrix at the array
\[
 Z_\tau=\sigma(\pi_M)t_\sigma\overline{Y_{\pi_\sigma(\tau)}},
         \qquad \pi_\sigma(\tau)=\sigma\circ\tau .
\]
The matrix at \(\tau\) has \(n_{j+1}\) target rows and \(n_j\) source
columns. Let \(\rho_{\tau,i}\) denote its squared singular values in
decreasing order, including zeros. Since \(\pi_\sigma\) permutes \(\Phi\)
and \(|\sigma(\pi_M)t_\sigma|^2=M/\Omega\), we have
\[
 \sum_{\tau,i}\rho_{\tau,i}=\sum_\tau\|Z_\tau\|_{\rm HS}^2
                      = e^{v_M}(2L_a+2)t_Y .
\]
Suppress the embedding index in the row products below.
In the tensor-induced Weyl metrics the two unitary
factors of a singular value decomposition act isometrically. On a
column of height \(h\le n_j\), the exterior power of the rectangular
diagonal factor has norm equal to the product of its first \(h\)
singular values. The column counts in
Section~\ref{sec:marked-column-gain} vanish for \(n_j<h<n_{j+1}\).
At each target column of height \(n_j\) the source tensor uses
\(e_1\wedge\cdots\wedge e_{n_j}\), a primitive integral vector of
norm one for the induced exterior form. Each such column contributes
one box to each of the first \(n_j\) target rows.

For a rectangular diagonal array with entries \(d_1,\ldots,d_{n_j}\)
in nonincreasing order of modulus, let \(c_h\) count the columns of
height \(h\). Its matrix entry on the unit highest column tensors is
\[
 \prod_h\det\!\bigl(\operatorname{diag}(d_1,\ldots,d_h)\bigr)^{c_h}
       =\prod_{i=1}^{n_j}d_i^{\lambda_i},
 \qquad \lambda_i=\sum_{h\ge i}c_h.
\]
Indeed the leading wedge \(e_1\wedge\cdots\wedge e_h\) maps to its
target counterpart times the leading determinant. This gives equality
in the product of the exterior-power norm bounds, including at
\(h=n_j\). For the scaled array \(\sigma(\pi_M)t_\sigma\bar Y\),
a diagram with \(|\lambda|=\sum_i\lambda_i\) boxes therefore
contributes the factor \((M/\Omega)^{|\lambda|}\) to the squared
norm relative to \(Y\). Multiplying over embeddings gives the squared
norm bound \(\prod_i\rho_i^{\lambda_i}\) on the scaled array.
For each exponent \(\lambda_i\),
\[
 \rho_i^{\lambda_i}
 \le \lambda_i!\exp\!\left(\lambda_i(v_M+L_a^{-1})
                                      +\rho_i e^{-v_M-L_a^{-1}}\right).
\]
Indeed it is the inequality \(q^n/n!\le e^q\) with
\(q=\rho_i e^{-v_M-L_a^{-1}}\). Since
\((2L_a+2)e^{-1/L_a}\le2L_a\), the square root of the
product is bounded by
\[
 \exp\!\left(
 \tfrac12\sum_i(\lambda_i/L_a+\log\lambda_i!+\lambda_i v_M)+L_a t_Y\right).
\]
Let \(v\) be any unit source vector and \(\ell\) any complex-linear
target covector of dual Weyl norm one, both in the conjugated
coordinates. Write \(\mathcal S_j(Z)\) for the rectangular Schur matrix.
In the jet identity take
\[
 Q(Z)=(\sigma^{-1}\ell)\,\mathcal S_j(Z)\,(\sigma^{-1}v),
\]
where the inverse automorphism acts on the coordinates of the row and
column. The Schur matrix has integral polynomial entries. Hence after
conjugating coefficients its contraction is
\(\ell\,\mathcal S_j(Z)\,v\), evaluated on the scaled array with the
permuted embedding indices. The bound just obtained applies uniformly
to these contractions. Taking the supremum over \(v,\ell\) gives the
operator norm of the conjugated jet.

The jet integral uses the Gaussian \(e^{-\|Y\|^2/2}\)
and Lebesgue measure divided by \(\pi\) in each complex
coordinate, so the \(L_a t_Y\) term leaves \(e^{-t_Y}\). The resulting
Gaussian integral is \((2L_a+2)^{fn_jn_{j+1}}\).
The remaining character and semicharacter factors have
conjugate modulus one; an automorphism of \(\mathbb C\) can
be composed with complex conjugation to fix \(K\) without
changing the Weyl norms. Thus the norm is bounded by
\[
 \exp\!\left(
 \frac{S_j}{2}+C_j+fn_jn_{j+1}\log(2L_a+2)
                      +\frac1{2L_a}\sum_{\tau,i}\lambda_{j,\tau,i}
       \right).
\]
For large \(x\), the parameter bounds give \(q_1,f\le x^{12}\),
\(n_m\le x^6\), and \(\log(2L_a+2)\le x\). Since
\(\lambda_{j,\tau,i}\le 2a\) and \(L_a=a+1\), the last three terms
are bounded by \(C_K(x^{24}+1)^2+x^{25}+x^{18}\le x^{50}\). This gives the bound \(e^{x^{50}+S_j/2}\) on the kernel.
\end{proof}

\subsubsection{Composition and terminal restriction}\label{sec:upper-terminal-norm}

\sectionopening{In this subsection, we combine Lemmas~\ref{lem:upper-factorial-sums}
and~\ref{lem:upper-gaussian-kernel} with the integral branching of
Section~\ref{sec:integral-terminal-branching} to bound the full
upper composition.}

\begin{proposition}[Pointwise bound after terminal restriction]\label{prop:upper-terminal-norm}
Use the untwisted tensor coefficients \(\Psi_\varnothing\) of
Section~\ref{sec:seesaw-integral-tensor}. Fix
\(\sigma\in\operatorname{Aut}(\mathbb C)\) and suppose the source
form \(\phi\) satisfies \(\|\sigma\phi(g)\|_{\rm W}^2\le e^{K_{\rm in}}\)
at every finite adelic source point \(g\). For fixed \(q_0,d,\epsilon\)
and sufficiently large \(x\), the squared Weyl norm of the
coefficientwise \(\sigma\)-conjugated restricted value at every terminal point is at most
\[
 \exp\!\left(K_{\rm in}+A-m\mathcal R+2\sum_p g_p\log p+E_0\right).
\]
\end{proposition}

\begin{proof}
For the untwisted upper operator, use the class sum of the jointly
invariant tensor from Lemma~\ref{lem:upper-tensor-class-sum}, then
expand at those class representatives.
Opening the \(N\)-gluing in
Section~\ref{sec:seesaw-tensor-translation-N} costs at
most \(N^{fm n_m}\) terms for each glued product.
At pair \(i\) the two tests share a vector in
\(\prod_{w\in P}\mathbb F_N^{n_{2i+1}}\), with \(|P|=f\), as in
Section~\ref{sec:ramified-direct-lift}. Thus the number of terms
before combining equal translates is at most
\[
 \prod_{i=0}^{m/2-1} N^{f n_{2i+1}}
 \le N^{fm n_m/2}\le N^{fm n_m}.
\]
Each phase has absolute value one. Thus the
conjugate absolute coefficient sum is at most
\(M^{3fm n_m^2}N^{fm n_m}\le\exp(x^{60})\)
for large \(x\). Use the step kernel estimate
\(\exp(x^{50}+S_j/2)\) proved above.

The count in Section~\ref{sec:terminal-classes} applies to
each definite stage of rank at most \(n_m\) with principal
level of the chosen generators or a containing level,
giving at most \(\exp(O(x^{52}))\) choices per step.
Let \(u\) be the conjugated value before terminal restriction.
With the squared source norms bounded by \(\exp(K_{\rm in})\),
the triangle inequality and the step operator bounds give
\[
 \|u\|_{\rm W}^2
   \le \exp\!\left(K_{\rm in}+\sum_j S_j+
                  2x^{60}+2mx^{50}+O(mx^{52})\right).
\]
Let \(\mathcal B_{\rm term}\) be the integral branching map of
Lemma~\ref{lem:terminal-branching-norm}. The difference of the sums of
the highest weights is \(d_\tau=a_\tau-mh_0\) at \(\tau\). The restriction multiplies this map by
\(\prod_\tau\tau(\overline{\pi_N})^{-d_\tau}\). Since
\(\pi_N\overline{\pi_N}=N\), at each complex embedding its squared
operator norm is at most
\[
 e^{x^{70}}N^{-\sum_\tau d_\tau}
             =e^{x^{70}}N^{-|\Phi_a|a+fmh_0}.
\]
Combining with Lemma~\ref{lem:upper-factorial-sums} gives
\[
 \|\text{conjugated restricted value}\|_{\rm W}^2
       \le\exp\!\left(K_{\rm in}+A-m\mathcal R
                      +2\sum_p g_p\log p+E_0\right).
\]
Indeed, the extra terms satisfy
\[
 32faB+fmh_0\log N+2x^{60}+2mx^{50}+O(mx^{52})+x^{70}\le E_0/2
\]
for large \(x\), since \(\log N=O(mB)\), \(a=xb\), and
\(b\) exceeds every fixed power of \(x\).

With \(K_{\rm in}\le E_0\), the finite-type projection estimate
gives the archimedean bound for the short images,
using the measured identity \eqref{eq:seesaw-measured-chain}
in its integral normalization.
\end{proof}

\subsection{The integral terminal branching map}
\label{sec:integral-terminal-branching}

\sectionopening{In this subsection, we realize the restriction at the terminal definite group by
contraction in exterior powers. This gives the integral map used in
Section~\ref{sec:upper-chain-size}, with its norm controlled in the
Weyl module.}

\begin{lemma}[Integral terminal branching]\label{lem:terminal-branching-norm}
On the terminal Weyl modules the branching map is integral,
has last-block character \(\prod_\tau z_\tau^{a_\tau-mh_0}\), and has
squared operator norm at most \(e^{x^{70}}\) in the original Weyl
metrics for sufficiently large \(x\). Here
\(a_\tau=a\mathbf 1_{\Phi_a}(\tau)\) and \(h_0=r^8\).
\end{lemma}

\begin{proof}
Fix \(\tau\), put \(n=n_m\), and write \(\lambda,\nu\) for the
descending highest weights on \(X_m,Y_m\) with indices starting at one.
Set \([l]_s=(l+1)/2-s\) for \(1\le s\le l\).
The upper Harish--Chandra parameter starts with \([k]\), with
\(a_\tau\) added to its first entry. It continues with the successive
blocks \(-w_i+[t+1]\), \(1\le i\le m\). The lower parameter starts with
\([k-1]\), followed for each \(i\) by \(-w_i+\Delta/2\) and
\(-w_i+[t]\). Subtract \([n]\) and \([n-1]\), respectively, to
recover \(\lambda\) and \(\nu\). Since
\([n-1]_h-[n]_{h+1}=1/2\), comparison of the two lists gives, for \(1\le h<n\),
\[
 \kappa'_\tau(h):=\nu_h-\lambda_{h+1}
   =\begin{cases}
       c_*=(\Delta-t-1)/2, & h=n_i=k+i(t+1),\quad 0\le i<m,\\
       0, & \text{otherwise},
     \end{cases}
 .
\]
At \(h=k\), the upper column count is \(w_1-(k+t+1)/2\), at least
\(c_*\) since \(w_1=\kappa_\tau+3\Delta/2\) and \(k\le\Delta\).
At \(h=n_i\), \(1\le i<m\), the count is \(\Delta-t-1=2c_*\).
Thus \(\lambda_h\ge\nu_h\ge\lambda_{h+1}\) throughout. At \(h=1\)
all \(a_\tau\) columns are contracted; at \(h=n-1\) the retained
count vanishes since \(t+1>1\), so \(\nu_{n-1}=\lambda_n\).
Retain \(\kappa'_\tau(h)\) columns of upper height \(h\) and contract
\(\gamma_\tau(h)=\lambda_h-\nu_h\) against the last basis vector.
At each polynomial height \(1\le h<n-1\) the lower multiplicity is
\(\nu_h-\nu_{h+1}=\kappa'_\tau(h)+\gamma_\tau(h+1)\).

The terminal Weyl diagram has the two sets of columns
specified above: in the first set contract each
column against the last unit basis vector, reducing its
height by one, and in the second set retain the height.
The map on the ambient product of exterior powers is
the corresponding tensor of contractions and projections
to the first block. Each of these maps has operator norm
at most one.

After projection onto the desired Weyl
summand it gives the normalized complex branching map.
The target norm here is the tensor-induced Weyl norm.
The adjoint of the ambient contraction wedges the shortened
columns with the same last unit basis vector. The adjoint
on the Weyl spaces then projects this extended tensor onto
the upper Weyl subspace. Before projection its pairings are
the products of the original minors by exterior Cauchy--Binet. A retained column of full first
block height contributes a determinant to the pairing
and is absorbed by its determinant character. This is
the tensor isometry on the spanning evaluation vectors.
After subtracting the common last entry \(e_\tau=\nu_{n-1}=\lambda_n\),
both weights are partitions. Macdonald's Schur branching formula
\cite[Ch.~I, (5.8)--(5.9) and the one-variable specialization
following (5.11)]{MacdonaldHall} gives each interlacing partition once.
The last-variable exponent for these partitions is
\(|\lambda|-|\nu|-e_\tau\); restoring the factor
\(\det(\operatorname{diag}(g,z))^{e_\tau}\) gives
\(z^{|\lambda|-|\nu|}\) with the prescribed lower determinant twist.
Taking the product over \(\tau\), the nonzero map spans the branching line.

The integral normalization is needed for the lattice norm
comparison. In the minor polynomial model its highest
test is the product of shifted leading minors in the
contracted columns and unshifted leading minors in the
retained ones. At upper height \(h\) the contracted
minor uses the first-block covectors \(2,\ldots,h\)
and has height \(h-1\) after contraction. A retained
minor uses the covectors \(1,\ldots,h\) and has height
\(h\ge k\). The only upper columns below height \(k\)
are of height one and their contractions are scalars.
All the nonempty minors therefore contain the same
covectors \(2,\ldots,k\).

Permute these common covectors to the beginning, and
remove them from each of the nonempty minors. If \(L\)
bounds the number of boxes left, one can clear the indicated
highest-vector polynomial by a positive integer
\(c\le L^{3L}\), with \(c=1\) for the empty diagram.
Write \(n_{\rm low}=n_m-1\) and \(u_+=\max(u,0)\).
The degree of the product \(H\) of leading minors after shortening is
\[
 D_{\rm sh}=\sum_{\tau}\sum_{1\le h<n_{\rm low}}
   (h-k+1)_+\bigl(\kappa'_\tau(h)+\gamma_\tau(h+1)\bigr)
   =f c_*\,m\bigl(1+(m-1)(t+1)\bigr).
\]
Indeed retained columns at \(k\) leave one box each and contracted ones
leave none. At \(n_i>k\) the \(c_*\) retained and \(c_*\) contracted
columns leave \(i(t+1)+1\) and \(i(t+1)\) boxes per column, respectively.
The height-one columns leave none. All these degrees exclude the
common determinant character. Writing \(H=\sum_\alpha H_\alpha X^\alpha\),
the positive clearing integer is
\[
 c=\langle H,H\rangle_{\rm Fock}
   =\sum_{|\alpha|=D_{\rm sh}}\alpha!\,|H_\alpha|^2,\qquad
 H_\alpha\in\mathbb Z,\quad \alpha!=\prod_j\alpha_j! .
\]
If \(D_{\rm sh}=0\), then \(H=1\) and \(c=1\).
For \(1\le D_{\rm sh}\le L\), the determinant expansion gives
\(\sum_\alpha|H_\alpha|\le\prod_{\rm col}h!\le L^L\).
Since \(\alpha!\le D_{\rm sh}!\le L^L\),
\[
 1\le c\le D_{\rm sh}!\left(\sum_\alpha|H_\alpha|\right)^2
              \le L^{3L},\qquad 2\log c\le6L\log L .
\]

For the clearing property, the polynomial model is dual
to the column-tensor Weyl model. Let \(G_{\rm sh}\) be the
Weyl Gram matrix in the integral basis of the shortened diagram,
and let \((P_i)\) be its dual coordinate polynomials on
the evaluation arrays. The dual Weyl form in this basis has
matrix \(G_{\rm sh}^{-1}\). The Fock form is invariant under
the same compact unitary group, so irreducibility makes it a scalar
multiple of this dual form. The product of leading minors is the
functional of the highest unit exterior tensor. Its squared dual
norm is one and its Fock square is \(c\). Thus the Fock matrix is
\(cG_{\rm sh}^{-1}\). If \(p\) is the coefficient column of an
integral product of minors in the shortened standard basis, put
\[
 b_i=\Big\langle P_i,\sum_j p_jP_j\Big\rangle_{\rm Fock}.
\]
The monomial factorial pairing gives \(b\in\mathbb Z^{d_{\rm sh}}\),
and \(G_{\rm sh}b=cp\). This is the clearing property for the
shortened diagram. Write \(J\) for the integral isometry of column tensors
which restores the \(q=k-1\) leading basis vectors to each column
that was nonempty before shortening. Let \(W_{\rm sh},W_0\) be the
two Weyl subspaces and \(\Pi_{\rm sh},\Pi_0\) their orthogonal projections.
Evaluation at \(\operatorname{diag}(I_q,Z)\) gives
\(J\operatorname{ev}_{\rm sh}(Z)=\operatorname{ev}_0(\operatorname{diag}(I_q,Z))\)
at each place, so \(JW_{\rm sh}\subset W_0\).
Conversely, write an evaluation array in blocks
\(X=\left(\begin{smallmatrix}A&D\\ C&E\end{smallmatrix}\right)\)
with \(A\) of size \(q\). If \(u\) is the number of restored columns
at this place, expansion of their minors gives, for \(\det A\ne0\),
\[
 J^*\operatorname{ev}_0(X)
     =(\det A)^u\operatorname{ev}_{\rm sh}(E-CA^{-1}D).
\]
Thus \(J^*\operatorname{ev}_0(X)\in W_{\rm sh}\) for every \(X\)
by polynomial continuation, and tensoring the identities gives
\(J^*W_0\subset W_{\rm sh}\) over all places.
The two inclusions give \(\Pi_0 J=J\Pi_{\rm sh}\).
The mixed exterior tensor representing the shifted minor product is
the image under \(J\) of the shortened tensor, with the signs of
the chosen permutations. Its projected vector multiplied by \(c\)
is therefore integral in the ambient tensors. The integral left inverse
of the standard-tableau basis matrix gives integral coordinates in
the original Weyl basis.

\smallskip
\indent\textit{Integral branching columns.}
Let \(\gamma_\tau(h)\) and \(\kappa'_\tau(h)\) be the numbers of
contracted and retained columns, respectively, of height \(h\) in the
upper terminal diagram. Write \(P_j(Y)\) for the integral standard
minor basis on the upper specified array, and \(P_i^0(X)\) for that on
the lower one. Put
\[
 D(X)=\prod_\tau\det(X_\tau)^{\kappa'_\tau(n_{\rm low})},\qquad
 \mathcal A=\left\{D(X)\sum_i(G_0 b)_iP_i^0(X):
                           b\in\mathbb Z^{d_0}\right\},
\]
where \(n_{\rm low}\) denotes the lower rank and \(G_0\) is its
rational Weyl Gram matrix in the integral tableau basis. Write
\(\mathcal P_0(b)=\sum_i(G_0b)_iP_i^0(X)\). Then
\(\mathcal A=D(X)\mathcal P_0(\mathbb Z^{d_0})\).
In the minor model \(G_0=\mathsf F_0^{\mathsf t}\mathsf F_0\) with
\(\mathsf F_0\) the integral basis matrix, so \(\mathcal A\) is an
additive subgroup of \(\mathbb Z[X]\).

Let \(\mathcal K(X,Y)\) be the product of the leading minors in
\(\operatorname{diag}(X_\tau^{\mathsf t},1)Y_\tau\), putting its
last covector first in the contracted columns and last in the retained
columns, in the chosen order of the covectors. Let \(c>0\) be the clearing
integer already obtained from the shortened diagram and let
\(\eta\in\mathbb Z^{d_0}\) be the corresponding vector of the shifted minor product
in the lower Weyl coordinates. Straightening in the upper variable
writes
\[
                  c\,\mathcal K(X,Y)=\sum_j f_j(X) P_j(Y).
\]
The point to prove is \(f_j\in\mathcal A\), since coefficients in
\(\mathbb Z[X]\) alone leave the lower Weyl Gram denominator unaccounted for.

First suppose at every place that the last covector of \(Y_\tau\) is
\((t_\tau,0,\ldots,0)\), and let \(g_\tau\) be its upper first
square block. Laplace expansion along the last covector in each
contracted minor gives \(t_\tau\) per contracted column. Retained
columns of full lower height contribute
\(\det(g_\tau)^{\kappa'_\tau(n_{\rm low})}\) in addition to \(D(X)\).
All other minors are those of the lower highest polynomial
evaluated at \(g_\tau^{\mathsf t}X_\tau\), with the same covector
permutation as for \(\eta\). The contracted minors use columns \(2,\ldots,h\) of \(g_\tau\)
after Laplace expansion, hence covectors \(2,\ldots,h\) after transposition.
The retained ones use covectors \(1,\ldots,h\). This gives the mixed
shifted and unshifted coefficient represented by \(\eta\).
Hence \(c\,\mathcal K(X,Y)\) is
\(D(X)\mathcal P_0\) of the integral polynomial
vector obtained by applying the lower integral Schur matrix
\(\rho_0^{\rm Schur}(g)\) to \(\eta\) and multiplying by
the displayed powers of \(t_\tau,\det(g_\tau)\).
The identity holds over the polynomial ring in the entries of
\(Y\).

For general \(Y\), write \(t_\tau=Y_{\tau,\rm last,1}\).
At each place replace column \(1\) by \(t_\tau\) times that
column, and column \(j>1\) by
\(t_\tau Y_{\tau,j}-Y_{\tau,\rm last,j}Y_{\tau,1}\).
The new last covector is \((t_\tau^2,0,\ldots,0)\).
Each leading minor of height \(h\) acquires exactly \(t_\tau^h\),
since only earlier columns were added. Thus \(Y^\mu\,c\mathcal K(X,Y)\)
is \(D(X)\mathcal P_0\) of a vector of integer
polynomials in \(Y\), with \(Y^\mu\) the product of these fixed
monomial factors. For every \(Y^\delta\) the coefficient of
\(Y^{\mu+\delta}\) in this identity shows that
\([Y^\delta](c\mathcal K)\in\mathcal A\). This comparison takes
place over \(\mathbb Z\) and divides no coefficient by an entry of
the array.

For the integral coefficient comparison, order the variables
first by column and within a column by increasing covector, using
the earliest variable first in lexicographic monomial order.
The leading monomial of a flag minor is its diagonal monomial
with coefficient \(1\). The leading monomial of a standard
tableau product identifies the multiset of entries in every
horizontal part of the tableau, by counting the corresponding variables in each
matrix column.

Within each tableau row the entries are weakly increasing, so
these multisets determine the tableau. Thus the standard tableau products
have pairwise distinct monic leading monomials.
In any relation with coefficients in an additive group \(H\),
take a maximal leading monomial among the nonzero coefficients.
Its coefficient must be zero, a contradiction. Apply this with
\(H=\mathbb Z[X]/\mathcal A\) to the last coefficient statement.
Every \(f_j\) lies in \(\mathcal A\), so write
\(f_j=D(X)\mathcal P_0(b_j)\) with \(b_j\in\mathbb Z^{d_0}\).

The \(b_j\) are the columns of the branching matrix \(B\) in the
original Weyl bases. With \(\operatorname{ev}_+(Y)=(P_j(Y))_j\)
as a column vector, its defining identity is
\[
 D(X)\mathcal P_0\bigl(B\,\operatorname{ev}_+(Y)\bigr)
                        = c\,\mathcal K(X,Y).
\]
At the identity lower array the right side is \(c\) times
the prescribed highest coefficient. This is a product of
nonzero coordinate minors in a polynomial ring over \(\mathbb Q\),
so \(B\ne0\). Substituting
\(Y\mapsto\operatorname{diag}(g,s)Y\) in this identity and expanding
the contracted minors gives
\[
 B\,\rho_+^{\rm Schur}(\operatorname{diag}(g,s))
  =\left(\prod_\tau
     s_\tau^{\sum_h\gamma_\tau(h)}
     \det(g_\tau)^{\kappa'_\tau(n_{\rm low})}\right)
       \rho_0^{\rm Schur}(g)B.
\]
Indeed the first-block substitution on the kernel is
\(X_\tau\mapsto g_\tau^{\mathsf t}X_\tau\), which supplies the
displayed full-height determinant.

The evaluation vectors span,
and \(D(X)\mathcal P_0\) is injective in characteristic
zero, so this is an operator identity. In the present application
\(\kappa'_\tau(n_{\rm low})=0\), and the common determinant exponent
is the last entry of the terminal highest weight,
\[
 e_\tau=\lfloor n_{m-1}/2\rfloor-w_m(\tau),
 \qquad n_j=k+j(t+1).
\]
Incorporating \(\det(\operatorname{diag}(g,s))^{e_\tau}\)
therefore replaces the power of \(s_\tau\) by
\[
 B_\tau=e_\tau+\sum_h\gamma_\tau(h)
        =|\mu(X_m)_\tau|-|\mu(Y_m)_\tau|.
\]
Here \(|\mu|\) denotes the sum of the entries of the highest weight,
including its determinant exponent in each entry. Each contracted
column loses one box. Since \(c_*=h_0-(t+1)/2\), the contracted column counts at heights \(1\),
\(k\), and \(k+i(t+1)\) for \(1\le i<m\) are respectively
\(a_\tau\), \(w_1(\tau)-\lfloor n_1/2\rfloor-c_*\), and \(c_*\).
Consequently
\[
 B_\tau=a_\tau+w_1(\tau)-w_m(\tau)
   +\lfloor n_{m-1}/2\rfloor-\lfloor n_1/2\rfloor+(m-2)c_*
   =a_\tau-mh_0,
\]
since \(w_m(\tau)-w_1(\tau)=(m-1)2h_0\) and \(t+1\) is even.
For a lower evaluation vector with nonzero square-block
determinants, form the extended column tensor divided by the product
of its full-height determinant powers.

By exterior Cauchy--Binet
the pairings of these divided tensors equal those of the lower
evaluation vectors. The latter span the lower Weyl module, so
this defines an isometry into the upper ambient column tensor
space. The pairing identity above identifies \(B\) with \(c\)
times the adjoint of this isometry composed with the upper
isometric tensor inclusion. Thus \(\|B\|_{\rm W}^2\le c^2\).
No Gram inverse is used in the finite-coordinate contraction.

To bound the shortened diagram, the contracted columns
of upper height \(k\) have height \(k-1\) after contraction,
so removing the \(k-1\) common covectors leaves no boxes
from them. The retained columns at height \(k\) number
\(c_*\le\Delta\) and each leaves one box. At each height
\(k+i(t+1)>k\), the total multiplicity of both kinds is
\(\Delta-(t+1)\le\Delta\). There are at most \(n_m\)
heights, and each column leaves at most \(n_m\) boxes.
Thus one can take \(L=f n_m^2\Delta\). Here
\(f\le x^{12}\), \(n_m\le x^6\), and
\(\Delta=2r^8\le2\cdot3^8 x^{40}\), so \(6L\log L\le x^{70}\)
for large \(x\).
The last block character of the contraction is
\(\prod_\tau z_\tau^{B_\tau}\), where \(B_\tau\)
is the difference in the total box counts between
the two terminal diagrams, with their prescribed
determinant weight included. The integral branching operator therefore has squared norm at most
\(\exp(x^{70})\) in the original Weyl metrics.
\end{proof}

\subsection{Norms and index after normalization}\label{sec:normalized-theta-length-index}

\sectionopening{In this subsection, we normalize the image of the Poincar\'e map by a scalar common to the
whole family. The estimates of Sections~\ref{sec:short-normalization}
and~\ref{sec:upper-chain-size} then yield bounds for both the norm and
the index in Proposition~\ref{prop:marked-upper-estimates}.}

\begin{proposition}[Normalization of length and index]\label{prop:shortness-normalization}
Suppose the following estimates hold for the chosen arithmetic data.
Let \(U_g:S\to V\) be the translates of the unscaled
projected Poincar\'e map, with the scalar of the seesaw
identity removed, and \(d_X\) the dimension of the full
source before the ramified restriction.
By Section~\ref{sec:ramified-source}, \(d_X\le2\dim S\).
\begin{enumerate}
\item There is a spanning set of everywhere integral forms \(\phi\in S\)
such that, for every \(g\in G\) and \(\sigma\in\operatorname{Aut}(\mathbb C)\),
\[
  \|\sigma\Theta_1 U_g\phi\|^2
  \le \exp\left(A-m\mathcal R+2\sum_{p\in F}g_p\log p
                             +2E_0+x^{99}\right).
\]
\item There are \(e_p\ge0\) supported on \(F\),
with \(\sum_p e_p\log p\le E_0\).
At each coefficient place above \(p\), any finite tuple of
integral forms in \(S\) can be cleared into a submodule of \(S\)
by a matrix over the valuation ring with determinant valuation
at most \(d_Xe_p\). That submodule works for all translates
and all twists of level prime to \(p\). Its images under
\(\Theta_\xi U_g\) have lattice norm at most \(p^{e_p-g_p}\)
for every such twist and translate, where the unit ball for the lattice norm is the integral lattice.
\end{enumerate}
Then the maps
\[
  P_g=c_B U_g,\qquad c_B=\prod_{p\in F}p^{-\lfloor g_p-e_p\rfloor}
\]
satisfy
\[
 \|\sigma\Theta_1P_g\phi\|^2\le e^{A-m\mathcal R+L}
 \quad\text{for the indicated generators}.
\]
For \(t\subseteq G\) and each independent tuple \(y\) of images of everywhere integral
forms under \(P_g\) with \(g\in t\), they give
\(\operatorname{Ind}(y)\le2|t|\dim(S)E_0\).
The clearing matrices work for every twist of level prime to the
coefficient prime. For every \(g\), \(\operatorname{im}P_g=\operatorname{im}U_g\).
\end{proposition}

\begin{proof}
Normalize \(v(p)=1\). For \(p\in F\), put
\(\varepsilon_p=g_p-e_p-\lfloor g_p-e_p\rfloor\in[0,1)\). Then
\[
 v_p(c_B)=e_p-g_p+\varepsilon_p,\qquad
 p^{-v_p(c_B)}p^{e_p-g_p}=p^{-\varepsilon_p}\le1 .
\]
These formulas hold at every coefficient embedding above \(p\),
including when \(g_p<e_p\).
For \(p\notin F\), \(v_p(c_B)=0\). Thus the same clearing
matrices give integral images for every allowed twist.
Choose an integral preimage and a translate for each entry.
Group these entries by translate \(g\). At \(\kappa\mid p\) let
\(U_{\kappa,g}\) be the clearing matrix of the nonempty block at \(g\)
from the second hypothesis. If \(P\) groups the columns, the clearing
matrix in their original order is
\[
 U_\kappa=P\left(\bigoplus_g U_{\kappa,g}\right)P^{-1},\qquad
 v_\kappa(\det U_\kappa)
      =\sum_g v_\kappa(\det U_{\kappa,g})\le |t|d_Xe_p .
\]
Linearity of each \(\Theta_\xi P_g\) gives integrality for every allowed
twist prime to \(p\) with this same matrix. Consequently
\[
 \operatorname{Ind}(y)
 \le |t|d_X\sum_p e_p\log p
 \le 2|t|\dim S\,E_0 .
\]

At every complex embedding the change in logarithmic squared
norm is
\[
 2\log c_B=-2\sum_{p\in F}g_p\log p
       +2\sum_{p\in F}(e_p+\varepsilon_p)\log p .
\]
In particular
\[
 \log c_B\le -\sum_{p\in F}g_p\log p+E_0+\sum_{p\in F}\log p,
 \qquad \sum_{p\in F}\log p\le x_1^2+\log N\le b.
\]
The last inequality follows from \(x_1\le e^W\), \(B\ge W^6\),
and \(\log N=O(mB)\) for large \(x\). Since \(c_B\in\mathbb Q_{>0}\),
\[
 \begin{aligned}
 \|\sigma\Theta_1P_g\phi\|^2
 &= c_B^2\|\sigma\Theta_1U_g\phi\|^2\\
 &\le \exp(A-m\mathcal R+4E_0+2b+x^{99})
 \le \exp(A-m\mathcal R+L).
 \end{aligned}
\]
Here \(4E_0+2b+x^{99}\le faB^4=L\) for large \(x\).
The equality of the images follows from \(c_B\ne0\).
\end{proof}

Section~\ref{sec:short-images-proved} verifies the local
hypothesis for the same projected map. The seesaw identifies
the upper chain, and the local lattice norm estimate is preserved
by the unnormalized finite-type projection and translation.

\section{The normalized Poincar\'e map}
\label{sec:period}

\sectionopening{In this section, we define the Poincar\'e map, specify the
measures in its orbital expansion, and compute its Petersson blocks.
Section~\ref{sec:full-orbital-blocks} uses these blocks, and
Section~\ref{sec:seesaw-normalization} compares the lifts of the map
with the theta construction.}

\subsection{The finite coordinates}
\label{sec:period-finite-coordinates}

\sectionopening{In this subsection, we use the generators of Sections~\ref{sec:marked-level-prime-selection} and~\ref{sec:different-additive-character-scales} to choose the finite point and identify the prescribed level in lattice coordinates.}

Let \(G_{\rm ad}\) be the ambient unitary group, and write
\([G_{\rm ad}]\) and \([R]\) for the rational quotients.
The following choices fix the covector at each scale of Section~\ref{sec:datum}.

For the finite scale, let \(D_{\rm rel}=\Phi'_{q_1}(\zeta_{q_1})\)
and \(\omega_K=(1+\sqrt{-d})/2\). Define the finite idele
\(\alpha_{\rm W}\) to have component \(D_{\rm rel}^{-1}\)
at places dividing \(\omega_K\), and \(1\) elsewhere.
Put \(j_{\rm W}=-1/\sqrt{-d}\). Every prime of \(q_1\)
divides \((d+1)/4=\omega_K\bar\omega_K\), by the fixed
congruence on \(d\), and \(D_{\rm rel}\) divides \(q_1\)
by differentiating \(Z^{q_1}-1\) at \(\zeta_{q_1}\).
The two ideals of \(\omega_K,\bar\omega_K\) are coprime,
since their generators sum to one.

Conjugation preserves
the different ideal. Thus, place by place,
\[
 j_{\rm W}\alpha_{\rm W}\bar\alpha_{\rm W}\widehat{\mathcal O}_E
                                  =\mathfrak D_{E/\mathbb Q}^{-1}
                                                      \widehat{\mathcal O}_E .
\]
Let \(c_M\) be the finite idele equal to the chosen
generator \(\bar\alpha\) at its chosen places and to one
elsewhere. For \(a_{\rm f}=\alpha_{\rm W}/\bar c_M\) set
\[
                  u=(a_{\rm f}/\bar a_{\rm f})^{m/2},
                  \qquad \eta=u^{-1} I_{X_0}.
\]
This is a unitary finite scalar with the required grading.

With \(\varrho=m/2-1\), the rational form \(Y_0\) has
positive definite \(X_0=E^k\) and the ordered string
\[
 e_0,\ldots,e_{\varrho-1},z,e_{\varrho-1}^*,\ldots,e_0^* .
\]
The pairing of each \(e_a,e_a^*\) is one, the norm
of \(z\) is \(-N\), and the remaining string pairings
vanish. Take the generic shift along the ordered string
with constant \(c_{\rm gen}=-(j_{\rm W}M\delta_E)^{-1}\).
It sends \(z\) to \(-Nc_{\rm gen}e_{\varrho-1}^*\) and the
other nonfinal string vectors to \(c_{\rm gen}\) times their
successors, and it kills the final vector and \(X_0\).
The constant belongs to \(F_0^\times\) since both
\(j_{\rm W}\) and \(\delta_E\) are skew.

In the finite coordinates let \(K_L\) be the integral
group, \(J_L\) the kernel of its reduction, and \(H=U(X_0)\)
the definite factor of \(R=U\rtimes H\), and let
\(J_H=\iota_H^{-1}(J_L)\). We first compute the finite
isometry from the original lattice. Write \(L\)
for its orbital basis matrix, with the distinguished
column of \(X_0\) replaced by
\[
 v_0=x_0/\bar\pi_N+z/N,\qquad N=\pi_N\bar\pi_N .
\]
On \((v_0,z)\) the Gram matrix is
\(\left(\begin{smallmatrix}0&-1\\-1&-N\end{smallmatrix}\right)\).
The other \(X_0\)-columns and the hyperbolic pairs retain
their integral unimodular Gram matrices. Hence \(L\)
is a basis of an integral unimodular lattice, including
at the primes over \(N\).

This basis also identifies the pullback level on \(H\).
In the decomposition of \(X_0\) into the other columns and
\(x_0\), write an element of \(H\) as
\(h=\begin{pmatrix}A&b_1\\c_1&d_1\end{pmatrix}\).
On the corresponding glued middle block its matrix is
\[
 \iota_L(h)=
 \begin{pmatrix}
 A&b_1/\bar\pi_N&0\\
 \bar\pi_Nc_1&d_1&0\\
 -c_1/\pi_N&(1-d_1)/N&1
 \end{pmatrix}.
\]
Indeed \(x_0=\bar\pi_N v_0-z/\pi_N\).
At \(w\mid\pi_N\), the glued matrix and its inverse are
integral exactly when \(h_w\in\mathrm{GL}_k(\mathcal O_{E_w})\)
and \((c_1,d_1)\equiv(0,1)\bmod\pi_N\). Indeed \(\bar\pi_N\)
is a unit here, and the last row imposes the two congruences.
At the conjugate place \(h_{\bar w}=\overline{h_w}^{-\mathsf T}\),
so the row congruence makes the last column of
\(h_{\bar w}-I\) zero modulo \(\bar\pi_N\). The last diagonal
difference is thus divisible by both components of \(N\).
The conjugate glued matrix is the inverse adjoint for the
unimodular glued Gram matrix, hence integral. At the primes above \(M\),
the basis change is integral and invertible since \(N\)
is a unit. Thus the inverse image of the original
principal congruence level is the level \(J_H\)
of Section~\ref{sec:ramified-source}.

Write \(\iota_{\rm rat}\)
and \(\iota_L=L^{-1}\iota_{\rm rat}L\) for the definite
embedding in rational and orbital coordinates. In the
graded basis the string columns are multiplied by \(u\).
The glued column becomes \(x_0/\bar\pi_N+uz/N\).
Thus the graded basis matrix is \(u\iota_{\rm rat}(\eta)L\).
It has the same Gram matrix as \(L\). Choose an integral
change to an orthonormal basis locally. At split pairs a
basis on one component and its Gram dual on the other
give such a change. At an unramified nonsplit place a
unit-norm line exists: if every residue norm were zero,
polarization with all residue scalars would make the
form zero by the nondegenerate residue trace pairing.
Lift the line, use surjectivity of the local unit norm
to scale its norm to one, and split off its unimodular
complement. This gives the change in every rank, also
in residue characteristic two.

Right multiplication by a standard integral unitary
matrix preserves the orthonormal Gram matrix. At each
chosen split pair choose its first component to correct
the reduction of the change of basis to the identity
at the chosen place. The second component is forced by
unitarity, and reduction is prescribed only at the chosen
place of the pair. The integral general linear reduction
there is surjective. 
Thus the orthonormal finite point
of this covector can be chosen as
\[
 x_f=u\,\iota_{\rm rat}(\eta)L\gamma_0,
 \qquad \gamma_0\in\operatorname{GL}_r(\widehat{\mathcal O}_E),
 \qquad \gamma_0\equiv I\text{ at each chosen place.}
\]
As part of fixing the original space, choose this \(x_f\)
before forming any of its period maps. In its standard
coordinates define
\[
 \begin{aligned}
 \mathfrak f(k)
 &=\iota_L(\eta)^{-1}L^{-1}x_f kx_f^{-1}L\iota_L(\eta)\\
 &=\gamma_0 k\gamma_0^{-1}.
 \end{aligned}
\]
This equivalence preserves integrality and the chosen principal
reductions, so \(\mathfrak f(J)=J_L\). Write \(J_0\) for the
rational Gram matrix and \(J_w=(L^*J_0L)_w\bmod w\) at a chosen
split place \(w\mid M\). Since \(\gamma_{0,w}\equiv I\), the
Gram identity \(\gamma_0^*(L^*J_0L)\gamma_0=I\) gives
\[
 \overline{\gamma_{0,\bar w}}^{\mathsf T}J_w=I
                  \quad\bmod w .
\]
Thus the opposite component gives the dual basis for the same
Gram matrix. The definite embedding reads exactly \(\iota_L\) under \(\mathfrak f\), since
\(\eta\) is central in that definite factor. If a rational
term is read as \(\gamma_{\rm std}=x_f^{-1}\gamma x_f\)
and \(\gamma_L=L^{-1}\gamma L\), then
\(\mathfrak f(\gamma_{\rm std})=\eta_H^{-1}\gamma_L\eta_H\)
with \(\eta_H=\iota_L(\eta)\). For a representative
\(c^\circ\) of \(H(F_0)\backslash H(\mathbb A_{F_0,f})/J_H\),
the standard class slice \(a_c\) embeds \(c^\circ\eta^{-1}\).
Its image is \(\eta_H^{-1}h_c\). Write
\[
 h_c^{\rm rat}=\iota_{\rm rat}(c^\circ),\qquad
 h_c=\iota_L(c^\circ)=L^{-1}h_c^{\rm rat}L .
\]
These comparisons hold before choosing any translate.

\subsection{Class slices and quotient measure}
\label{sec:period-quotient-measure}

\sectionopening{In this subsection, we decompose the rational quotient using the finite coordinates of Section~\ref{sec:period-finite-coordinates}.}

In the chosen lattice basis \(A\oplus(X_0\oplus Ez)\oplus A^*\),
the matrix of a Bessel element \(u\iota_H(a)\) on the
entire middle diagonal block is that of \(\iota_H(a)\).
Indeed \(u\) induces the identity there, also after conjugation
by the middle change of basis used to glue the last column.
On the two outer blocks \(\iota_H(a)\) is the identity.
Apply this to the element and its inverse. Integrality and
congruence to one on the glued middle block then hold for the
embedded Levi part, with its full \(z\)-covector included. Removing
that part gives \(u\in U\cap J_L\). Transporting through the
change of basis gives \(J_R=(U(\mathbb A_{F_0,f})\cap J)\rtimes J_H\).

For \(n\iota_H(a)\in R(\mathbb A_{F_0,f})\), write
\(a\eta=\delta c^\circ j_H\), with \(\delta\in H(F_0)\)
and \(j_H\in J_H\). Choose \(u\in U(F_0)\) with
\(\iota_H(a)^{-1}u^{-1}n\iota_H(a)\in U\cap J\).
The condition on \(u\) is open about \(n\) in the finite
adeles. The polynomial exponential coordinates of \(U\)
and additive approximation over \(F_0\) give such a rational
point. The rational element \(\iota_H(\delta^{-1})u^{-1}\)
moves the finite point to \(a_c J_R\), since \(\eta\) is central
in the Levi.

Conversely any identification between two such
slices identifies their definite classes by taking the Levi
part. At the same class the rational Levi part is one by
freeness of the definite level. The remaining
stabilizer is the original \(U_c\). Therefore
\[
 [R]/J_R=\coprod_c U_c\backslash R_\infty .
\]
Each unipotent quotient is compact by the covector boxes of
Section~\ref{sec:marked-class-covolume}. The real Levi is
compact and there are finitely many classes.

The left
\(U_c\)-translates of the reduced cutoff \(\beta_c\) sum to
one by Section~\ref{sec:marked-full-to-cut}. Write \(v_R\)
for the finite volume of \(J_R\) including the positive
scalar of the Bessel Haar decomposition. With counting measure on
the rational group, write the Haar measure inducing the Petersson
measure as \((\prod_{v\mid\infty}dg_v)\,dg_f\). Here \(dg_v\) has
density \(1/\operatorname{vol}U(r)\) at the identity in the
skew-Hermitian coordinates \(S=JX\) of
Section~\ref{sec:real-formal-degree-factor}. This determines
\(dg_f\) before the Poincar\'e map is formed. If \(K_{G,f}\) is
the full finite lattice stabilizer, reduction at the chosen split
places gives \(K_{G,f}/J=G\), whence
\[
 v_G=\operatorname{vol}_{dg_f}(J)
    =\frac{\operatorname{vol}_{dg_f}(K_{G,f})}{|G|}.
\]
The slice formula for every integrable right \(J_R\)-invariant \(Q\) is
\[
 \int_{[R]}Q=v_R\sum_c\int_{R_\infty}\beta_c(s)Q(s,a_c)\,ds .
\]

\subsection{The real point and stationary translate}
\label{sec:period-real-point}

\sectionopening{In this subsection, we fix the real chart and place the stationary point of Section~\ref{sec:real-marked-filling} in that chart.}

Use the signed scales of Section~\ref{sec:marked-schur-orbits}.
At \(v\) let \(J_0\) be the rational Gram matrix and let \(J_1\)
replace its entry \(-N\) by \(-1\). Take \(D_v\) diagonal
with entries \(\kappa_v^{\varrho-a}/\sqrt N\) on \(e_a\),
their reciprocals on \(e_a^*\), \(N^{-1/2}\) on \(z\), and
one on \(X_0\). Then \(D_v^*J_0D_v=J_1\).
Let \(p_v\) be the Euclidean unitary matrix whose columns
are the normalized sums and differences of the hyperbolic
pairs, the basis of \(X_0\), and \(z\), ordered by the
prescribed signs. Choose the real components \(x_v\) of the
adelic point by
\[
 x_v=D_vp_v,\qquad
 x^{\rm pt}=((x_v)_{v\mid\infty},x_f),\qquad
 x_v^*J_0x_v=p_v^*J_1p_v=\varepsilon_v.
\]
Here \(x_f\) is the corrected finite point chosen above, with
the same reduction and principal level at the prescribed
places. The real chart satisfies
\[
 \iota_v(g)=p_vgp_v^{-1},\qquad
 k_v=D_v^{-1}x_vp_v^{-1}=1,\qquad T_v(g)=\iota_v(g).
\]
The isomorphism \(\iota_v\) preserves the specified Lie
algebra inner products and Haar. The rational matrix
\(\gamma\) acts in the \(J_1\)-coordinates by
\(D_v^{-1}\Phi_v(\gamma)D_v\).

Let \(x_v^\star\) be the stationary point of
Section~\ref{sec:real-marked-filling}. Its point in the
standard signed coordinates and its adelic realization are
\[
 g_{0,v}=\iota_v^{-1}(x_v^\star),\qquad
 g_{0,\infty}=(g_{0,v})_v,\qquad g_0=(g_{0,\infty},1).
\]
In rational coordinates \(x_vg_{0,v}=D_vx_v^\star p_v\).
The Bessel coefficient is evaluated at
\[
 (x_v^\star)^{-1}\iota_v(s_v)T_v(g_v)
                        =\iota_v(g_{0,v}^{-1}s_vg_v).
\]

The character \(\Psi\) on \(R\) is trivial on the Levi and is
the adelic character at \(-\sigma/(M\delta_E)\) on \(U\). Here
\(\sigma\) is the simple-root sum and \(\delta_E\) the chosen
generator of the different. In the grading the simple-shift
coefficient \(c_{\rm gen}\) and the Weil scalar \(j_{\rm W}\)
satisfy \(j_{\rm W}c_{\rm gen}=-(M\delta_E)^{-1}\).
The trace of the shift against a Bessel unipotent has two
contributions at each adjacent pair of basis positions.

The dual entry is the conjugate of the corresponding entry
of the inverse unipotent, by unitarity. The first
superdiagonal sum of the inverse is \(-\sigma\), so the trace
is \(c_{\rm gen}(\sigma-\bar\sigma)\). Since the product
\(j_{\rm W}c_{\rm gen}\) is purely imaginary, the Weil
argument \(j_{\rm W}c_{\rm gen}(\sigma-\bar\sigma)/2\)
has the same real trace as \(j_{\rm W}c_{\rm gen}\sigma\).
This gives the character of the original graded period.

Put \(z_u=-\sigma(u)/(M\delta_E)\). In the signed chart the
simple-root entry is
\(\alpha_{\Psi,v}=\sigma(u)_v/\kappa_v=-2\pi i z_{u,v}\),
where \(z_{u,v}\) is the archimedean component. Thus
\[
 \exp(\alpha_{\Psi,v}-\bar\alpha_{\Psi,v})
   =\exp(-4\pi i\Re z_{u,v})=\Psi_v(u),\qquad
 \prod_v\exp(-\alpha_{\Psi,v}+\bar\alpha_{\Psi,v})
   =\Psi_\infty(u)^{-1}.
\]
The first character gives the real Bessel covariance. Its inverse
product occurs in the defining average.
At the finite places the vectors \(e_a,z\) are integral in the glued
basis, and the top coordinate functionals are integral because the
glued column has no top component. Thus the simple-root coordinates
of a level unipotent are integral everywhere and are in
\(M\widehat{\mathcal O}_E\) at the places over \(M\): the principal
congruence in the chosen half gives the conjugate congruence by
unitarity and the unimodular Gram matrix. The common string scaling
in \(x_f\) leaves these entries unchanged. Consequently
\(z_u\in\mathfrak D_{E/\mathbb Q}^{-1}\widehat{\mathcal O}_E\),
and \(\Psi_f\) is trivial on \(J_R\). For
\(F=(e_0,\ldots,e_{\varrho-1},z)\), each simple-root entry reads
the \(e_i\)-coordinate of the image of the next column of \(F\).
The glued matrix \(h_c\) fixes these columns and preserves the
\(e_i\)-coordinate functionals. Thus the character is also trivial
on the class conjugate of this unipotent level. For a rational
\(u\in U_c\) it follows that \(\Psi_\infty(u)=\Psi_f(u)^{-1}=1\).

At \(v\), take \(w_v=v_0\) as in Section~\ref{sec:real-marked-filling}.
Let \(I_v\) be the local value space and \(\deg_v\) the formal degree.
For the average \(A_v\) of Section~\ref{sec:real-cut-error}, put
\(a_v=\sqrt{\dim I_v}/\|A_v\|_2>0\). Use the reproducing kernel
\(e=\prod_v \deg_v\bar c_v\), with \(c_v\) the normalized
coefficient in the fixed real coordinates. Write \(C_{\rm W}\) for
the Weyl-to-Bessel isomorphism in conjugate coordinates.
For the unitary Levi action \(\pi_H=B_0 S B_0^{-1}\) in
Weyl coordinates use \(C_{\rm W}=\overline{B_0}\) on the
conjugate coordinates. In the local tableau basis put
\(Q_v=B_{0,v}^*B_{0,v}\) and let \(\xi_v\) be the highest monomial
of Section~\ref{sec:real-marked-filling}, including the constant
monomial for the empty diagram. Since \(Q_v\) is real and
\(w_v=B_{0,v}\xi_v/\|B_{0,v}\xi_v\|\),
\[
 C_{{\rm W},v}^* C_{{\rm W},v}=\overline{Q_v}=Q_v,\qquad
 (C_{{\rm W},v}z_v)^{\mathsf T}w_v
    =\frac{z_v^{\mathsf T}Q_v\xi_v}{(\xi_v^*Q_v\xi_v)^{1/2}} .
\]
Thus the local change is an isometry for the Weyl metric, and
the period evaluates the displayed complex-linear covector.

Unitarity gives
\[
 (C_{\rm W}\rho_{\rm conj}(h)z)^{\mathsf T}
    =(C_{\rm W}z)^{\mathsf T}\pi_H(h)^{-1}.
\]
At each real place the chart sends \(u\iota_H(h)\) to the graph of
its Levi component. The inverse reads the ambient group point from the
first coordinate and the Levi action from the prescribed quotient of
the flag. These maps are continuous group isomorphisms. Thus the
product of the local graph measures of Section~\ref{sec:real-cut-error}
pulls back to \(C_\infty ds\), with \(0<C_\infty<\infty\) by
uniqueness of Haar measure. Set
\[
 w=\bigotimes_v w_v,\qquad
 t_0(z)=C_\infty\frac{\prod_v a_v}{\prod_v \deg_v}
                              (C_{\rm W}z)^{\mathsf T}w .
\]
The factors of the product average at \(s=u\iota_H(h)\) are
\(\overline{e(g_{0,\infty}^{-1}sx)}/\prod_v \deg_v\), the wave
\(\Psi(s)^{-1}\), and \(\pi_H(h)^{-1}w\), besides the
positive factors \(a_v\). Tensor the absolutely integrable
local averages and invert \(s\) in the unimodular Bessel group.
For \(r=u\iota_H(h)\) put
\[
 W_{t,\phi}(r)=\Psi(r)t\bigl(\rho_{\rm conj}(h_\infty^{-1})
                                              \phi(h_f)\bigr).
\]
The evaluation of \(R_H(\eta)\psi\) at \(a_c\) is
\(\psi(c^\circ)\) and \(\Psi(a_c)=1\). After inversion,
the original full Bessel vector therefore satisfies
\[
 (C_{\rm W}\psi(c^\circ))^{\mathsf T}W(x)
  =\int_{R_\infty}\overline{e(g_{0,\infty}^{-1}s^{-1}x)}
                         W_{t_0,R_H(\eta)\psi}(s,a_c)\,ds .
\]

\subsection{The measured period and change of level}
\label{sec:period-measured-pairing}

\sectionopening{In this subsection, we use the point of Section~\ref{sec:period-real-point} and the measure of Section~\ref{sec:period-quotient-measure} to identify the measured pairing.}

With \(z_{c,i}=(C_{\rm W}\psi(c^\circ))_i\) and the seeds of
Section~\ref{sec:period-petersson-blocks}, define
\[
 \mathcal P_{\rm ser}(\psi)(x)=
   \sum_{c,i}z_{c,i}\sum_{\gamma\in G_{\rm ad}(F_0)} f_{c,i}(\gamma x).
\]
The inner sum uses counting measure on the rational group. The
cutoff in \(f_{c,i}\) accounts for its \(U_c\)-translates.
At a fixed finite level the dual real test space is the image of
right convolution by \(\bar e\). Choose a bounded representative
of \(Y\) on that finite-level quotient and put
\[
 Y^\circ(x)=\int_{G_\infty}\bar e(a)Y(xa)\,da,\qquad
 |Y^\circ(xh)-Y^\circ(x)|
       \le \|Y\|_\infty\|\bar e(h^{-1}\,\cdot)-\bar e\|_1
       \quad(h\in G_\infty).
\]
Each finite-level real slice is a quotient of \(G_\infty\) by
a discrete subgroup. Pullback along its orbit map takes null sets
to null sets, so the integral is independent of the bounded
representative. Translation continuity in \(L^1\) gives continuity
on each slice. The reproducing identity of
Section~\ref{sec:real-reproducing-coefficient} gives
\(Y^\circ=Y\) almost everywhere. Use \(Y^\circ\) for the period
evaluations and write it as \(Y\).
For a principal-level test the absolute periodization bound permits
unfolding at the chosen level and gives
\[
 \int_{[G_{\rm ad}]}Y\,\mathcal P_{\rm ser}(\psi)
  = v_G\sum_c\iint\beta_c(x)Y(x,a_c)
       \overline{e(g_{0,\infty}^{-1}s^{-1}x)}
       W_{t_0,R_H(\eta)\psi}(s,a_c)\,ds\,dx .
\]
Both \(Y\) and the period coefficient are invariant
under the respective right levels on this class.
Rational invariance shows that the double integrand
apart from \(\beta_c\), denoted \(\mathcal K(s,x)\),
satisfies \(\mathcal K(us,ux)=\mathcal K(s,x)\)
for \(u\in U_c\). The factor involving \(e\) depends on \(s^{-1}x\).

The integral of \(\beta_c(s)|\mathcal K(s,x)|\) is finite
by the boundedness of \(Y\), the \(L^1\) bound on \(e\)
and the integrability of the Bessel coefficient on
its compact quotient. Insert \(\sum_u\beta_c(ux)=1\)
in this integral and change variables \((s,x)\mapsto(us,ux)\)
for each \(u\). Summing the translates of \(\beta_c(s)\)
proves finiteness with the cutoff \(\beta_c(x)\) by
Tonelli. The same change of variables and absolute
convergence now give
\[
 \iint\beta_c(x)\mathcal K(s,x)\,ds\,dx
   =\iint\beta_c(s)\mathcal K(s,x)\,ds\,dx.
\]
The convolution in \(x\) reproduces \(Y(sg_{0,\infty},a_c)\).
Using the Bessel Haar decomposition, the resulting
pairing is \(\int_{[R]} W_{t,R_H(\eta)\psi}(s)Y(sg_0)\,ds\)
with \(t=(v_G/v_R)t_0\). For the original source
form \(\phi\) take \(\psi=R_H(\eta^{-1})\phi\).
Let \(Z=R(\bar e)Z\in L^2\) have the prescribed dual real type
and principal level. There is a compactly supported smooth kernel
\(\xi\) on the real group such that \(\xi*\bar e=\bar e\).
Indeed, in the conjugate discrete series use the compact matrix
unit projecting onto the distinguished line of the multiplicity-one
lowest compact type, as in Section~\ref{sec:real-differentiated-schur}.
Convolve a real symmetric approximate identity on both sides with
this self-adjoint matrix unit. The resulting smooth kernel acts
on the distinguished unit vector by a real scalar tending to one.
Normalize by this scalar and take the product over the real places.
The resulting kernel is self-adjoint and fixes the conjugate unit
vector. The matrix coefficient formula for \(\bar e\) then gives
the asserted convolution identity.

Put \(\mathcal C=\{rg_0:r\in[R]\}\) in the quotient at level \(J\).
It is compact by Section~\ref{sec:period-quotient-measure}.
On each of its finitely many real components, lift its intersection
with \(\mathcal C\) to a compact set \(L\) in the real group. For
the discrete arithmetic group of this component, only finitely many
translates of \(L\,\operatorname{supp}\xi\) meet this same compact
set. Cauchy--Schwarz in the convolution, with this bound on the
multiplicity, gives
\[
 \sup_{x\in\mathcal C}|R(\xi)F(x)|\le C_{\mathcal C,J,\xi}\|F\|_2 .
\]
The smooth convolution gives a continuous representative of
\(Z=R(\xi)Z\). For \(Z_n=Z\mathbf1_{\{|Z|\le n\}}\), the
bounded tests \(R(\bar e)Z_n\) converge to \(Z\) in \(L^2\).
Each is fixed by \(R(\xi)\), so convergence is also uniform on
\(\mathcal C\). The weight \(W_{t,\phi}\) is bounded on the
compact Bessel quotient. Thus its integrals against these tests
converge to the integral against \(Z(rg_0)\). Square-integrability
of the series below gives the same limit in the Petersson pairing,
identifying its representative for the original \(g_0,t,\phi\).

For a bounded test \(Z\) of the dual real type which is right
\(K'\)-invariant, with \(K'\) a finite open level, let
\(N_J=[J:J\cap K']\). Its average over the representatives
\(\gamma\) of \(J/(J\cap K')\) is \(J\)-invariant.
Right invariance of \(\Pi\phi\) under \(J\), and right
translation in the Petersson integral, give
\[
 \begin{aligned}
 N_J\int_{[G_{\rm ad}]} Z(x)(\Pi\phi)(x)\,dx
    &=\int_{[G_{\rm ad}]}
          \sum_{\gamma\in J/(J\cap K')}Z(x\gamma)(\Pi\phi)(x)\,dx\\
    &=\sum_{\gamma\in J/(J\cap K')}
       \int_{[R]} W_{t,\phi}(r)Z(rg_0\gamma)\,dr .
 \end{aligned}
\]
In each sum use one representative of each coset. The second
line is the measured period applied to the summed test.
This formula keeps the factor \(N_J\) when the first
covector of the lower chain changes level.

\smallskip
\subsection{Petersson blocks and covariance}
\label{sec:period-petersson-blocks}

\sectionopening{In this subsection, we unfold the pairing of Section~\ref{sec:period-measured-pairing} to obtain the blocks used in Section~\ref{sec:full-orbital-blocks}.}

Expand the series in the source class basis, using the quotient Haar
measure of Section~\ref{sec:period-quotient-measure} and the
cutoffs of Section~\ref{sec:marked-full-to-cut}.

Let \(C_X\psi=(C_{\rm W}\psi(c^\circ))_c\) and
\(C=C_XR_H(\eta^{-1})\). By the local Weyl isometry above
and the trivial stabilizers of Lemma~\ref{lem:marked-source-freeness},
the class squared norm is \(\sum_c\|(C_X\psi)_c\|^2\).
The central translation permutes the classes, with unitary changes
of real frame at the representatives. Thus \(C\) preserves this
norm, also on the ramified subspace. The series for a coordinate
basis vector in \(C\) is the periodization of
\[
 f_{c,i}(y,x_f)=\beta_c(y)W_i(y)\mathbf1_{a_cJ}(x_f).
\]
Fix \((c,i)\) and identify rational elements using the isometries of
Section~\ref{sec:period-finite-coordinates}. If no rational translate
of \(x_f\) lies in \(a_cJ\), the series is zero at \(x\). Otherwise
choose \(\delta_0\) with \(\delta_0 x_f\in a_cJ\). For every other
contributing \(\delta\), the element
\(\gamma=\delta_0\delta^{-1}\) belongs to \(\mathcal T(c,1,c)\):
it takes \(\delta x_f\) to \(\delta_0 x_f\) in the same slice.
Writing \(y\) for the real point of \(\delta_0 x\), the other real
argument is \(\gamma^{-1}y\). The bound
\(\|W(z)\|\le A_*L(z)^{-I}\) and \eqref{eq:marked-length-packing}
therefore give
\[
 \sum_{\delta}|f_{c,i}(\delta x)|\le D_*,\qquad
 D_*=2A_*b^{12n^3rf},
\]
uniformly in \(x\), where \(A_*\) is the product of the local
majorant bounds.

Index the left \(U_c\)-orbits of the rational terms by
\(q=U_c\delta\), and set \(F_q(x)=\sum_{u\in U_c} f_{c,i}(u\delta x)\).
The finite support is \(U_c\)-invariant. On that support, the
real character is trivial on \(U_c\) and the translates of \(\beta_c\)
sum to one. Hence \(F_q(x)\) is the full coordinate \(W_i\) at the
real point of \(\delta x\), multiplied by
\(\mathbf1_{a_cJ}(\delta x_f)\). Each \(F_q\) is reproduced by right
convolution with the real kernel \(e\). Since
\(\sum_q|F_q(x)|\le D_*\), the sum of the absolute convolution
integrands has integral at most \(D_*\|e\|_1\). Thus right convolution
also reproduces the series.

Thus the
series has that type and is continuous by continuity of
translations of the kernel in \(L^1\). Being bounded on the
finite-volume quotient, it belongs to the full Petersson
space on the prescribed level. Absolute unfolding against bounded
tests and the measured pairing of Section~\ref{sec:period-measured-pairing}
identify the series with
\(\Pi C^{-1}e_{c,i}\) for the original measured map \(\Pi\).

Put \(B=(v_G\nu)^{-1/2}\Pi C^{-1}\) and \(B_s=B^*\rho(s)B\).
Use the Petersson product linear in the second variable.
To pair \(\rho(s)\Pi C^{-1}e_{c',i'}\) with
\(\Pi C^{-1}e_{c,i}\), unfold the seed in the second variable.
An integral lift \(k_s\) of \(s\) normalizes \(J\), so
the finite integral over \(a_cJ\) gives \(v_G\).
A term indexed by \(\gamma^{-1}\) in the first variable uses
\[
 \mathfrak f(a_{c'}^{-1}\gamma_{\rm std}^{-1}a_c k_s)
       =h_{c'}^{-1}\gamma_L^{-1}h_c\mathfrak f(k_s).
\]
Membership in \(J_L\) is exactly
\(h_c^{-1}\gamma_Lh_{c'}\in\mathfrak f(k_s)J_L\).
Thus the surviving integrand is the full kernel with both
cutoffs of \(\mathcal T(c,s,c')\), with the residue \(s\).
The real isometry preserves Haar: its signed-basis isomorphism
preserves the specified Lie algebra density and the
remaining change is a fixed left translation. Writing
\(A^{\rm full}\) for the full amplitude of
Section~\ref{sec:marked-full-to-cut}, unfolding gives
\[
 \langle\rho(s)\Pi C^{-1}e_{c',i'},\Pi C^{-1}e_{c,i}\rangle
                   =v_G A^{\rm full}(c,s,c')_{i,i'} .
\]
Both variables acquire the factor \((v_G\nu)^{-1/2}\).

\begin{proposition}[Full Petersson blocks]\label{prop:marked-petersson-blocks}
In the Euclidean source coordinates, the normalized Petersson blocks
are given by the full orbital sums
\begin{equation}\label{eq:marked-petersson-block}
 (B_s)_{(c,i),(c',i')}
               =\nu^{-1}\overline{A^{\rm full}(c,s,c')_{i,i'}} .
\end{equation}
\end{proposition}

Let \(H_R\) be the residue image of the integral Bessel group.
For \(t\in H_R\) write a lift in the integral Bessel group in
glued coordinates as \(u_L\iota_L(a)\). On the definite group write
\(c^\circ a=h\,(c')^\circ j_H\) with \(h\in H(F_0)\) and \(j_H\in J_H\).
Additive approximation in the unipotent coordinates gives
\(\xi\in U(F_0)\) with
\(h_c^{-1}\xi_L h_c\in u_L J_L\), where \(\xi_L=L^{-1}\xi_f L\).
For the rational isometry \(\delta=\xi\iota_{\rm rat}(h)\) this gives
\[
 h_c^{-1}\delta_L h_{c'}
   =(h_c^{-1}\xi_Lh_c)\iota_L(a)\iota_L(j_H)^{-1}
       \in u_L\iota_L(a)J_L ,
\]
since the integral \(\iota_L(a)\) normalizes \(J_L\).
Thus \(\delta\) has parameter \((c,t,c')\).
Every rational Bessel term has residue in \(H_R\).
Two such terms over the same left class and residue
differ on the right by a rational Bessel term of identity
residue. By the Levi test their right classes coincide,
and all terms over \((c,t,c')\) in \(R\) are \(\delta U_{c'}\).

Write \(W(\delta y)=M_\delta W(y)\); the matrix \(M_\delta\)
is unitary by covariance of the local averages. The right
cutoffs for \(\delta U_{c'}\) sum to one and their \(W\)-values
agree. The left integral is
\[
             \int\beta_c(y)W(y)W(y)^*\,dy=\nu I
\]
Indeed \(WW^*\) is invariant under left real unipotents, and
\(\int_{U_\infty}\beta_c(u y)\,du=\nu\) by
\eqref{eq:marked-unipotent-covolume} and its cutoff formula.
In the real charts the group and unipotent Haar measures are
the products of the local orbital measures \(dg_\tau,du_\tau\).
Their quotient is therefore \(\prod_\tau d\bar g_\tau\). Since
\(W=\bigotimes_\tau W_\tau\) and
\(\int\|W_\tau\|^2\,d\bar g_\tau=\dim I_\tau\), Fubini and
Section~\ref{sec:real-marked-filling} give
\[
 \int_{U_\infty\backslash G_\infty}WW^*\,d\bar g
 =\bigotimes_\tau\int W_\tau W_\tau^*\,d\bar g_\tau
 =1_{\bigotimes_\tau I_\tau}.
\]
Since \(W(\delta^{-1}y)=M_\delta^*W(y)\), the Bessel tensor is
\(W(y)W(y)^*M_\delta\) and its amplitude is \(\nu M_\delta\). The normalized conjugate blocks give
a unitary representation \(r_\psi(t)\) on the source:
composition and inversion of the rational transporters
multiply and invert the class moves and matrices.
Extend \(r_\psi\) by zero off \(H_R\).
Equation~\eqref{eq:marked-petersson-block} identifies
\(E_s=B_s-r_\psi(s)\) with the normalized conjugate
off-Bessel amplitudes.

To prove covariance, choose \(\delta\) over \((d,t^{-1},d')\).
Multiplication \(\gamma\mapsto\gamma\delta^{-1}\) bijects
the terms over \((c,s,d')\) and \((c,st,d)\).
The resulting right cutoff \(\beta_d(\delta\,\cdot)\)
has orbit sum one for \(U_{d'}\) because
\(\delta U_{d'}\delta^{-1}=U_d\) by the level condition.
The right kernel is multiplied by the factor \(M_\delta^*\).
Changing the cutoff by the orbit partition and using
absolute convergence gives
\[
 A^{\rm full}(c,st,d)=A^{\rm full}(c,s,d')M_\delta^* .
\]
By \eqref{eq:marked-petersson-block} the last matrix
factor after normalization is \(\overline{M_\delta}^{\,*}\).
This is the block from \(d\) to \(d'\) of \(r_\psi(t)\),
since \(\overline{M_\delta}\) is the block from \(d'\)
to \(d\) of \(r_\psi(t^{-1})\). Thus \(B_{st}=B_s r_\psi(t)\).
Taking adjoints and using \(B_s^*=B_{s^{-1}}\) gives
\(B_{ts}=r_\psi(t)B_s\).
These identities also hold for \(E_s\). They imply
\(\rho(t)B=Br_\psi(t)\): for any \(x\), the difference
\(\rho(t)Bx-Br_\psi(t)x\) belongs to the span of the
vectors \(\rho(s)By\) and is orthogonal to each of them.

Since finite translation by \(\eta^{-1}\) commutes with
the scalar actions, it preserves the ramified
source. The ranges of the projected measured translates on
that source are the ranges of \(p\rho(s)B\) on its
subspace in \(C\)-coordinates, with
\(p=\rho(h_0\otimes h_1)\). The factors in the algebraic
projection onto the selected types and the map normalization are nonzero
and common to the translates.
More explicitly, write \(e_0=h_0\otimes h_1=\sum_s u_s[s]\)
in the finite group algebra and \(e_0^\vee=\sum_s u_s[s^{-1}]\).
For a principal-level test \(Y\) of the dual real type, invariance
of the quotient measure and the measured pairing give
\[
 \int_{[G_{\rm ad}]}Y(x)\,(p\rho(g)\Pi\phi)(x)\,dx
 =\int_{[R]}W_{t,\phi}(r)
    \bigl(\rho(g^{-1})\rho(e_0^\vee)Y\bigr)(rg_0)\,dr .
\]
Thus the bilinear test uses the dual finite types.
For the Hermitian Petersson product the right translations satisfy
\(\rho(s)^*=\rho(s^{-1})\), by invariance of the quotient measure.
The coefficients of each character idempotent satisfy
\(u_{s^{-1}}=\overline{u_s}\) because a finite group representation
is unitarizable. Hence
\[
 p^*=\sum_s\overline{u_s}\,\rho(s^{-1})
     =\sum_s u_s\rho(s)=p .
\]
Together with \(p^2=p\), this identifies \(p\) as the orthogonal
projection on the original target.

\part{Orbital integrals and the translated span}\label{part:orbital}

\section{Orbital integrals and operator norms}
\label{sec:marked-moments}

\sectionopening{In this section, we relate the Petersson matrices of
Section~\ref{sec:period} to two operator moments.
Their bounds give the hypotheses of the translation argument in
Section~\ref{sec:marked-span-proof}.}

Write \(T_{\rm full}\) for the full Petersson target and
\(B:X\to T_{\rm full}\) for the geometrically normalized map
on the full source with its Euclidean class coordinates.
The chosen real matrix coefficient selects one line in each automorphic
occurrence. On the rank-\(k\) definite source, each of its \(d_C\)
classes contributes \(d_B=\prod_v\dim I_v\), so \(d_X=d_Cd_B\).
On the rank-\(r\) definite inner form used to count the target,
each free class contributes the compact representation of dimension
\(\mathcal B\), as computed in Section~\ref{sec:marked-target-dimension}.
Let \(\rho\) be the unitary action of \(G=G_0\times G_1\)
on \(T_{\rm full}\). Let \(p_i=\rho(h_i)\) for the
central idempotents of the selected types \(h_i\in\mathbb C[G_i]\), and \(p=p_0p_1\).
Put \(d_X=\dim X\), \(D=\dim pT_{\rm full}\), and
\(B_s=B^*\rho(s)B\). The following argument concerns the
same normalized map for both its operator estimates.

\subsection{From truncated orbital integrals to full orbital integrals}
\label{sec:full-orbital-blocks}

\sectionopening{In this subsection, we compare the complete Petersson matrices with their truncated
orbital integrals. The estimate of
Equation~\eqref{eq:marked-full-cut-square} controls the difference.
Section~\ref{sec:marked-pair-weight-sum} bounds pairs of distinct
links, and \eqref{eq:open-equal-link-established} bounds the
equal-link term on the open set with the normalization
of \eqref{eq:marked-main-identity-ratio}.}

Let \(H_R\) be the residue Bessel subgroup and \(r_\psi\) its
unitary action on \(X\), extended by zero to \(G\).
Write \(E_s=B_s-r_\psi(s)\), \(D_{\rm id}=|H_R|d_X\) and
\[
 D_{\rm main}=\frac{d_X^2}{\mathcal B\,\operatorname{mass}(E,r)},
 \qquad C_{\mathcal A}=b^{-6k|\mathcal A|}D_{\rm main},
\]
where \(\mathcal B\) is the product of the formal degrees.
Let \(\mathcal A\) be a subset of the \(f\) chosen places at
\(p_0^{\rm lev}\). At each place of \(\mathcal A\) choose any
nonzero vector in the Bessel span, or any nonzero
covector in its opposite-place realization. Let
\(\Gamma_{\mathcal A}\) be the intersection of these stabilizers.
Put \(q=p_0^{\rm lev}\). The split reduction onto \(\mathrm{GL}_r(\mathbb F_q)\) is
surjective. Each nonzero vector or covector has orbit size \(q^r-1\)
under this action, so its stabilizer has that index.
At a chosen place \(w\mid p_0^{\rm lev}\), write \(g_w\) for the
component of \(s\), \(J_w\) for the Gram matrix in the original
glued lattice basis of Section~\ref{sec:period-finite-coordinates},
and \(E_j\) for the matrix of the first \(j\) columns of
\((e_0,\ldots,e_{\varrho-1},z)\) in that basis. Set
\[
 \mathcal O_w=\{s\in G: \det(E_j^tJ_w g_w E_j)\ne0
                               \text{ for every }1\le j\le\varrho\}.
\]
The Gram matrix is symmetric with rational entries. Identifying the
two residue fields by conjugation, unitarity gives
\(g_{\bar w}=J_w^{-1}g_w^{-t}J_w\). Hence \(g_{\bar w}\)
fixes \(E_jy\) exactly when \(g_w^{-t}\) fixes \(J_wE_jy\).
Thus
\[
 \Gamma_{\mathcal A}^{\rm open}
    =\Gamma_{\mathcal A}\cap\bigcup_{w\mid p_0^{\rm lev}}\mathcal O_w ,
\]
where the union runs over the chosen places. At any one unused place
the nonopen subset is contained in the union of the zero sets of its
proper leading minors. Put \(M=p_0^{\rm lev}p_1^{\rm lev}\).
For \(\Omega\subseteq G\), let \(\mathcal L(\Omega)\) be the sum of the
squared norms of the cut blocks on the nonzero-discriminant part of
\(\Gamma_*\), taken separately for each link and summed over the two
classes and \(s\in\Omega\), before division by \(\nu^2\).

\begin{proposition}[Orbital bounds on the finite stabilizers]
\label{prop:full-orbital-stabilizer-bounds}
For the Bessel data at sufficiently large \(x\), put
\(\beta(s)=\|B_s\|_{\rm HS}^2\). Then
\[
 \begin{aligned}
 \sum_{s\in H_RG_i}\|E_s\|_{\rm HS}^2&\le b^{-12}D_{\rm id}
                                             &&(i=0,1),\\
 \sum_{s\in\Gamma_{\mathcal A}^{\rm open}}\beta(s)
             &\le(1+b^{-4})^2 C_{\mathcal A}
                                        &&(\mathcal A\text{ arbitrary}),\\
 \sum_{s\in\Gamma_{\mathcal A}}\beta(s)&\le C_{\mathcal A}
                                        &&(|\mathcal A|=f).
 \end{aligned}
\]
Here \(\Gamma_{\mathcal A}^{\rm open}\) is the open subset specified
above.
\end{proposition}

\begin{proof}
\begin{enumerate}
\item By \eqref{eq:marked-petersson-block}, the blocks of \(E_s\)
are the normalized conjugates of the full amplitudes off
the rational Bessel subgroup. For a set \(\Omega\) avoiding the cosets of
\(H_R\) by nontrivial integral norm-one scalars,
\eqref{eq:marked-off-bessel-reduction} identifies the cut off-Bessel
amplitudes with the cut amplitudes on \(\Gamma_*\).

In the block norm obtained by dividing by \(\nu\)
and summing squares over the two classes and \(s\in\Omega\),
\eqref{eq:marked-full-cut-square} bounds the cut error by
\(b^{-20}\sqrt{D_{\rm id}}\).
\item Section~\ref{sec:marked-pair-weight-sum} bounds the squared
norm of the combined cut blocks on \(\Gamma_*\) by \(\nu^{-2}\mathcal L(\Omega)+b^{-kf}D_{\rm id}\) in
the same block norm. For \(1\le z\le Q_{\cap}(s)\)
on \(\Omega\), the equal-link sum satisfies
\[
 \nu^{-2}\mathcal L(\Omega)\le (b^{2f}/z)^kD_{\rm main}.
\]
Here \(Q_{\cap}(s)\) is the product of the residue
norms where the basis and its image have nonzero intersection.
On the open part of \(\Gamma_{\mathcal A}\), the sharper equal-link
estimate is proved in \eqref{eq:open-equal-link-established}.
\item By Section~\ref{sec:marked-orbital-reduction},
the sets \(H_RG_i\) avoid the nontrivial scalar cosets
and every \(s\in H_RG_i\) satisfies \(Q_{\cap}(s)\ge (M/b)^f\). The open set avoids
all scalar cosets, including \(H_R\). When every place is fixed,
\(Q_{\cap}(s)\) is at least \((p_0^{\rm lev})^f\) on \(\Gamma_{\mathcal A}\).
The residue group has one general linear factor over
\(\mathbb F_{p_i^{\rm lev}}\) for each of the \(f\) chosen split
places at each \(p_i^{\rm lev}\).
\end{enumerate}

The parameter bounds give \(b\ge1000\), \(kf\ge40\),
\(p_0^{\rm lev}\ge b^{10}\) and \(M\le b^{23}\) for large \(x\).
The moment comparison
\(b^{11kf}D_{\rm id}\le D_{\rm main}\le (M/b^6)^{kf}D_{\rm id}\)
is \eqref{eq:marked-main-identity-ratio}.

For every \(\Omega\) disjoint from the cosets \(\operatorname{red}(\zeta)H_R\)
with \(1\ne\zeta\in\mathcal O_E^\times\) and \(\zeta\bar\zeta=1\), the triangle inequality
in the block norm gives
\begin{equation}\label{eq:orbital-full-return}
 \left(\sum_{s\in\Omega}\|E_s\|_{\rm HS}^2\right)^{1/2}
 \le \nu^{-1}\sqrt{\mathcal L(\Omega)}
       +(b^{-20}+b^{-kf/2})\sqrt{D_{\rm id}}
 \le \nu^{-1}\sqrt{\mathcal L(\Omega)}
       +b^{-15}\sqrt{D_{\rm id}}.
\end{equation}
On \(H_RG_i\) the equal-link contribution is at most
\[
 (b^{2f}/(M/b)^f)^kD_{\rm main}\le b^{-3kf}D_{\rm id}.
\]
Thus the sum of \(\|E_s\|_{\rm HS}^2\) over \(H_RG_i\)
is at most \(4b^{-30}D_{\rm id}\le b^{-12}D_{\rm id}\).
By the covariance of Section~\ref{sec:period-petersson-blocks},
\(E_t=(B^*B-I)r_\psi(t)\) for \(t\in H_R\). Therefore
\[
 |H_R|\|B^*B-I\|_{\rm HS}^2
     =\sum_{t\in H_R}\|E_t\|_{\rm HS}^2\le b^{-12}D_{\rm id}.
\]
With \(D_0=BB^*\), Cauchy--Schwarz in the Hilbert--Schmidt
inner product gives
\[
 \beta(s):=\|B_s\|_{\rm HS}^2
 =\operatorname{tr}(D_0\rho(s)D_0\rho(s)^*)
 \le \|D_0\|_{\rm HS}^2
 \le (1+b^{-6})^2d_X \le 2d_X .
\]

For the open part of \(\Gamma_{\mathcal A}\), we have \(E_s=B_s\).
Since \(C_{\mathcal A}\ge b^{5kf}D_{\rm id}\), the right side of
\eqref{eq:orbital-full-return} is at most
\[
 \bigl(\sqrt{1+b^{-4}}+b^{-15}\bigr)\sqrt{C_{\mathcal A}}
                       \le (1+b^{-4})\sqrt{C_{\mathcal A}}.
\]
This gives open mass at most \((1+b^{-4})^2 C_{\mathcal A}\).

Suppose \(|\mathcal A|=f\). By Section~\ref{sec:marked-orbital-reduction},
the scalar cosets meet \(\Gamma_{\mathcal A}\) only in \(H_R\).
The pointwise estimate gives
\[
 \sum_{s\in\Gamma_{\mathcal A}\cap H_R}\beta(s)
 \le 2d_X|H_R|=2D_{\rm id}\le C_{\mathcal A}/4,
\]
since \(C_{\mathcal A}\ge b^{5kf}D_{\rm id}\).
On \(\Omega=\Gamma_{\mathcal A}\setminus H_R\) we have
\(E_s=B_s\), and the equal-link bound is
\[
 \nu^{-2}\mathcal L(\Omega)
 \le (b^{2f}/(p_0^{\rm lev})^f)^k D_{\rm main}
 \le b^{-8kf}D_{\rm main}\le C_{\mathcal A}/8.
\]
Also \(b^{-30}D_{\rm id}\le C_{\mathcal A}/8\).
The square of the right side of \eqref{eq:orbital-full-return}
is at most \(2(C_{\mathcal A}/8+C_{\mathcal A}/8)=C_{\mathcal A}/2\).
The full sum is at most \(C_{\mathcal A}/4+C_{\mathcal A}/2
\le C_{\mathcal A}\).
\end{proof}

\subsection{The diagonal bound and the projected subspace}
\label{sec:marked-diagonal-moment}

\sectionopening{In this subsection, we use Bessel covariance to estimate the error in the diagonal moment.
Proposition~\ref{prop:full-orbital-stabilizer-bounds} and the
finite projection estimates of
Section~\ref{sec:marked-residue-projection} give bounds in the same
source norm.}

\begin{proposition}[Diagonal moment on the projected target]\label{prop:orbital-diagonal-moment}
Let \(P_U\) be the orthogonal projection onto a subspace \(U\subseteq X\).
For the map \(B\) and target projection \(p\),
\[
           \|pBP_U\|_{\rm HS}^2\ge \dim U-5b^{-4}d_X .
\]
\end{proposition}

\begin{proof}

The diagonal argument uses the Bessel covariance together with
the bound of Proposition~\ref{prop:full-orbital-stabilizer-bounds}:
\[
 E_{ts}=r_\psi(t)E_s,\quad E_{st}=E_s r_\psi(t)
     \quad(t\in H_R),\qquad
 \sum_{s\in H_RG_i}\|E_s\|_{\rm HS}^2\le b^{-12}|H_R|d_X.
\]
Section~\ref{sec:marked-residue-projection} bounds the squared
counting norm of the conjugacy average of
\(\operatorname{tr}r_\psi|_{G_i}\) by \(2d_X^2\), and gives
\(q_i(1)\le(2b^8)^{-1}\) for \(q_i=1-h_i\), with \(q_i(1)\)
the group-algebra coefficient at the identity.

Since \(E_t=E_1r_\psi(t)\) for \(t\in H_R\), summing over
\(H_R\) and dividing out its order gives
\[
                   \|B^*B-I\|_{\rm HS}^2\le b^{-12}d_X .
\]
Right convolution by \(q_i\) is an orthogonal projection
on Hilbert--Schmidt-valued functions on \(G\), by
self-adjointness and idempotence. Apply the contraction to
\(E_s\mathbf1_{H_R G_i}(s)\). On \(t\in H_R\) the resulting
value is \(r_\psi(t)\sum_{s\in G_i}q_i(s)E_s\), with
convolution taken as \(\sum_s q_i(s)E_{ts}\).
Summing its squared norm over \(t\in H_R\) gives
\(\|\sum_{s\in G_i}q_i(s)E_s\|_{\rm HS}\le b^{-6}\sqrt{d_X}\).
Its trace consequently has absolute value at most
\(b^{-6}d_X\). Write \(a_i(s)=\operatorname{tr}r_\psi(s)\) on
\(G_i\) and \(a_i^\natural\) for its conjugacy average. Centrality of
\(q_i\), self-adjointness and idempotence give
\[
 \left|\sum_{s\in G_i}q_i(s)a_i(s)\right|
 =\left|\sum_{s\in G_i}q_i(s)a_i^\natural(s)\right|
 \le q_i(1)^{1/2}\|a_i^\natural\|_2\le b^{-4}d_X,
\]
where \(\|q_i\|_2^2=(q_i^*q_i)(1)=q_i(1)\) in counting measure.
Taking their sum bounds
\(\operatorname{tr}(B^*(1-p_i)B)=\|(1-p_i)B\|_{\rm HS}^2
\le2b^{-4}d_X\).
The two summands in \(1-p=(1-p_0)+p_0(1-p_1)\) have
orthogonal images, so
\[
 \|(1-p)B\|_{\rm HS}^2
 \le\sum_{i=0}^1\|(1-p_i)B\|_{\rm HS}^2\le4b^{-4}d_X .
\]
Using the Hilbert--Schmidt bound for \(B^*B-I\), we obtain
\begin{align*}
 \|pBP_U\|_{\rm HS}^2
 &=\dim U+\operatorname{tr}\bigl(P_U(B^*B-I)\bigr)
       -\|(1-p)BP_U\|_{\rm HS}^2\\
 &\ge \dim U-b^{-6}\sqrt{d_X\dim U}-4b^{-4}d_X
 \ge \dim U-5b^{-4}d_X .
\qedhere\end{align*}
\end{proof}

\subsection{The mixed moment bound}
\label{sec:marked-mixed-moment}

\sectionopening{In this subsection, we control the mixed moment by summing the open contributions against
the stabilizer contributions. Proposition~\ref{prop:full-orbital-stabilizer-bounds}
estimates both sums, and the
dimension comparison of Section~\ref{sec:marked-target-dimension}
supplies the normalization.}

\begin{proposition}[Mixed moment on the family]\label{prop:orbital-mixed-moment}
For the same \(B,p\), with \(B_s=B^*\rho(s)B\) and \(D=\dim pT_{\rm full}\),
\[
       \frac{D}{|G|}\sum_{s\in G}\|B_s\|_{\rm HS}^2
                 \le (1+9b^{-4})d_X^2 .
\]
\end{proposition}

\begin{proof}

Put \(\delta_b=b^{-4}\), \(\lambda=b^{6k}\) and
\(\mu_{\mathcal A}=\sum_{\Gamma_{\mathcal A}}\beta\).
Proposition~\ref{prop:full-orbital-stabilizer-bounds} bounds the
open mass for the full kernels by \((1+\delta_b)^2 C_{\mathcal A}\),
including the cut-to-full error. For \(|\mathcal A|=f\) it bounds
the total mass by \(C_{\mathcal A}\).

Maximize \(\lambda^{|\mathcal A|}\mu_{\mathcal A}\) over the fixer choices,
including the empty subset. If a place is unused, the
mass of the intersection with any single allowed vector
or covector fixer at that place is at most \(\mu_{\mathcal A}/\lambda\).
For the leading minor calculation these choices include
the nonzero vectors in the two column spaces
and their corresponding covectors. Let \(q=p_0^{\rm lev}\)
and \(\pi:\Gamma_{\mathcal A}\to\mathrm{GL}_r(\mathbb F_q)\) be the
reduction there.

Write the leading block of size \(j\) as
\(Z_s=A_1^t\pi(s)A_0\). Here \(A_0\) consists of the
first \(j\) columns. To form \(A_1\), take the columns at the
conjugate place and pair them by the split Gram matrix. At both places
the columns belong to the integral basis of
Section~\ref{sec:period-finite-coordinates}. Reduction preserves their
rank, and the reduced Gram matrix is invertible. Thus
\[
 \operatorname{rank}A_0=\operatorname{rank}A_1=j,\qquad
 A_0v\ne0\ (v\ne0),\qquad A_1y\ne0\ (y\ne0).
\]
They lie respectively in the original span and the prescribed
opposite-place covector span.
The two stabilizers tested below have equations
\(\pi(s)A_0v=A_0v\) and
\(\pi(s)^{-t}A_1y=A_1y\).
The weight is of positive type because for
\(D_s=\rho(s)BB^*\rho(s)^*\),
\(\langle D_a,D_c\rangle_{\rm HS}=\beta(a^{-1}c)\).
For a nontrivial additive character \(\psi_q\), the
left-kernel size of \(Z_s\) is
\(\sum_y\mathbb E_v\psi_q(y^t Z_s v)\).
If \(\det Z_s=0\) the kernel contains at least \(q\) elements.
With \(y,v\in\mathbb F_q^j\), the two axes \(y=0\) and \(v=0\)
contribute exactly \((2-q^{-j})\mu_{\mathcal A}\) to the weighted sum.

For \(y,v\ne0\), put \(a_1=A_1y\)
and \(a_0=A_0v\). The vectors
\[
 \sum_{a\in\Gamma_{\mathcal A}}\widehat e_{\pi(a)^{-t}a_1}\otimes D_a,
 \qquad
 \sum_{c\in\Gamma_{\mathcal A}}e_{\pi(c)a_0}\otimes D_c
\]
use delta vectors \(e_z\) on \(\mathbb F_q^r\) and unit
Fourier vectors \(\widehat e_w(z)=q^{-r/2}\psi_q(-w^t z)\).
Their squared norms are \(|\Gamma_{\mathcal A}|\) times the respective
fixer masses, by expanding and setting \(s=a^{-1}c\).
Their inner product equals \(|\Gamma_{\mathcal A}|q^{-r/2}\) times
\(\sum_{s\in\Gamma_{\mathcal A}}\psi_q(a_1^t\pi(s)a_0)\beta(s)\).
Cauchy--Schwarz therefore bounds the latter sum by
\(q^{r/2}\mu_{\mathcal A}/\lambda\). Summing over \(y,v\) with the
normalization in \(v\) proves
\[
 q\sum_{\det Z_s=0}\beta(s)
 \le \bigl(2+q^{j+r/2}/\lambda\bigr)\mu_{\mathcal A}.
\]
By Section~\ref{sec:datum}, \(b^{10}\le q\le b^{11}\).
Since \(j\le m/2\) and \(r\ge x^5\), for sufficiently large \(x\),
\[
 11j+6m\le12m\le r/2,\qquad
 \log_b(q^{j+r/2})\le 11j+\tfrac{11}{2}r
                  \le 6(r+1-m)=\log_b\lambda.
\]
There are at most \(m/2\le b^6/3\) prescribed minors.
Every nonopen residue has a vanishing minor at this particular unused
place. Using \(q\ge b^{10}\) and summing over its minors, the nonopen mass
is at most \(\delta_b\mu_{\mathcal A}\).

If the maximizing fixer uses every place, use the
full-fixer mass bound just proved. Otherwise combine
the open and nonopen bounds to get
\((1-\delta_b)\mu_{\mathcal A}\le(1+\delta_b)^2 C_{\mathcal A}\).
In either case maximality and \(0<\delta_b\le1/5\) give
\[
 \sum_{s\in G}\|B_s\|_{\rm HS}^2=\mu_{\varnothing}
 \le \lambda^{|\mathcal A|}\mu_{\mathcal A}
 \le \frac{(1+\delta_b)^2}{1-\delta_b}D_{\rm main}
 \le (1+9\delta_b)D_{\rm main}.
\]
By \eqref{eq:marked-tower-dimension}, this gives
\[
  \frac{D}{|G|}\sum_{s\in G}\|B_s\|_{\rm HS}^2
                                       \le (1+9b^{-4})d_X^2 .
\]
The bound uses the original full source, normalized
measured map and projected target.
\end{proof}

\subsection{Dimension on the terminal definite group}
\label{sec:marked-target-dimension}

\sectionopening{In this subsection, we transfer the required dimensions to a definite unitary group with
the same finite levels. The trace comparison respects the projection
of Section~\ref{sec:marked-residue-projection}, which lets us use
these dimensions in Section~\ref{sec:marked-moments}.}

Use the terminal positive definite Hermitian space \(W_m\)
of rank \(r\), with its standard lattice and the principal
level at the same chosen places. The original space \(Y_0\)
and \(W_m\) have identical finite groups and levels in
their isometric coordinates. Their real types have
the same regular Harish--Chandra parameter at each real
place. Write \(T_W\) for the full principal-level space
of the terminal type.

\begin{proposition}[Definite comparison of dimensions]\label{prop:marked-definite-dimension}
On the family,
\begin{equation}\label{eq:marked-tower-dimension}
            D=\dim pT_{\rm full}=\dim pT_W\le |G|\mathcal B\,\operatorname{mass}(E,r).
\end{equation}
\end{proposition}

\begin{proof}
Fix one tuple selected by \(p\). At each of the two
different level primes this tuple has a primitive type polynomial
at a chosen place, by Section~\ref{sec:marked-residue-projection}.
The two places of \(F_0\) are split in \(E\) and distinct.
At either such place let \(K=\operatorname{GL}_r(\mathcal O_v)\)
and let \(\tau\) be the finite cuspidal representation of its
residue group given by the primitive type polynomial.

If \(Z\) is the
local centre, extend the inflation of \(\tau^\vee\) unitarily
to \(ZK\) and compactly induce to \(\operatorname{GL}_r(F_{0,v})\).
For a diagonal Cartan representative outside \(ZK\), order its
elementary-divisor exponents decreasingly and cut at a strict drop.
Conjugation sends the upper off-diagonal unipotent block at this cut
into the principal kernel. Intertwining would therefore factor through
the coinvariants of \(\tau^\vee\) by the unipotent radical of this
proper residue parabolic, which are zero by cuspidality.
Thus the intertwining criterion gives irreducibility. Its coefficients
are compactly supported modulo \(Z\), so this compact induction is
supercuspidal. Every irreducible local representation containing the inflated
\(\tau\) is supercuspidal by the depth-zero type theorem
\cite[\S4.1, Proposition, p.~151]{MorrisDepthZero}.
Consequently the selected discrete spectrum is cuspidal:
at such a place the induced representations from every
proper Levi have zero projection to this type, and the
noncuspidal spectrum is built from these inductions.

In the compact induction, the vectors supported on the identity
coset have coefficients supported on \(ZK\). Truncate them by
\(\mathbf1_{\mathcal O_v^\times}(\det g)\). Since
\(ZK\cap\{\det g\in\mathcal O_v^\times\}=K\),
the sum of these truncated diagonal coefficients, with
the factor \(\dim(\tau)/\operatorname{vol}(K)\), is exactly
the local projector density
\[
 \frac{\dim\tau}{\operatorname{vol}(K)}
                  \overline{\operatorname{tr}\tau(\bar g)}
                                  \,\mathbf1_K(g).
\]
Here \(\bar g\) is reduction at the chosen place. In an orthonormal
basis of the copy of \(\tau^\vee\) on the identity coset the diagonal
coefficients on \(K\) sum to
\(\operatorname{tr}\tau^\vee(\bar g)=\overline{\operatorname{tr}\tau(\bar g)}\).
This gives the asserted expansion at both chosen places.
Each truncated coefficient is strongly cuspidal by
\cite[Lemma~1.9.1]{HarrisLabesse}. In particular the projector
at the second chosen place is cuspidal, as required by
condition \((\mathrm b.1)\) of Proposition~1.10.6 there.

Use the normalized principal-level averages and the finite
central idempotent of the fixed tuple at all finite places.
Their product is the test for the original tuple on both source spaces. At each real place choose the Euler--Poincar\'e function
with coefficients in the dual of the algebraic representation
of highest weight \(\lambda-\rho\), where \(\lambda\) is
the dominant representative of the common regular parameter.

For a real signature \((a,b')\) divide this function by
\(d_{a,b'}=\binom{r}{a}\), the size of its discrete series
packet. A nonelliptic regular semisimple element has a
nontrivial split torus in the centre of its connected centralizer:
its eigenvalues off the unit circle occur in reciprocal-conjugate pairs.
Square-integrability here is taken on the full real centralizer.
The absolute value of a nonzero coefficient of an irreducible
unitary representation is invariant under the central split torus,
whose Haar measure is infinite. Consequently the centralizer has
no square-integrable irreducible representation, and the ordinary
orbital integral vanishes by \cite[Theorem~4.1(i)]{ClozelLabesseOrbital}.
This gives cuspidality of the real test in the sense of
\cite[\S1.5]{HarrisLabesse}.

For a semisimple \(\delta\), an eigenvalue \(z\) off the unit
circle and its partner \(1/\bar z\) contribute a factor
\(\mathrm{GL}_{m_z}(\mathbb C)\) to the real centralizer.
Its split centre forces the stable orbital integral to vanish by
\cite[Theorem~4.1(iii)]{ClozelLabesseOrbital}.
If every eigenvalue has modulus one, let \(m_z\) be its multiplicity.
The conjugacy classes in the stable class on \(U(a,b')\) are specified
by eigenspace signatures \((a_z,m_z-a_z)\) with \(\sum_z a_z=a\).
The eigenspaces are orthogonal and nondegenerate, and their signatures give
\[
 I_\delta=\prod_z U(a_z,m_z-a_z),\qquad
 \overline I_\delta=\prod_z U(m_z).
\]
Give the centralizer the measure for which discrete series with the
infinitesimal character of the trivial representation have formal
degree one \cite[\S4]{ClozelLabesseOrbital}. The stable orbital sum
weights each ordinary class by the Kottwitz sign of its centralizer
\cite[\S1.5]{HarrisLabesse}. Thus
\cite[Theorem~4.1(ii), (iv)]{ClozelLabesseOrbital} gives
\[
 \operatorname{SO}_\delta(\phi_F/d_{a,b'})
   =e(\overline I_\delta)\operatorname{tr}F(\delta),\qquad
 \sum_{\sum a_z=a}\prod_z\binom{m_z}{a_z}
       =[X^a](1+X)^r=d_{a,b'} .
\]
For the compact inner centralizer the sign is
\[
 e(\overline I_\delta)=(-1)^{\sum_z\lfloor m_z^2/4\rfloor}.
\]
Indeed the quasi-split real form of \(U(m_z)\) has signature
\((\lfloor m_z/2\rfloor,\lceil m_z/2\rceil)\), so its symmetric
space has real dimension \(2\lfloor m_z^2/4\rfloor\).
Here \(F\) is the dual algebraic coefficient specified above.
For fixed centralizer measure \(di_\delta\), replacing real
Haar \(dg\) by \(c\,dg\) replaces the Euler--Poincar\'e test by
\(c^{-1}\phi_F\) and the quotient measure \(dg/di_\delta\) by
\(c\,dg/di_\delta\). Thus the displayed orbital integral has
the same value for the original group Haar. Every unit-circle stable class
occurs on each signature, and the displayed value depends only on its
eigenvalues with multiplicities. The normalized real tests thus have a
common transfer in the sense of \cite[\S1.6]{HarrisLabesse}. At each
finite place the given isometry transports the tests and their measures.

Apply the cuspidal trace comparison of Harris--Labesse
\cite[Proposition~1.10.6]{HarrisLabesse}, termwise in the
coefficient expansion at the two chosen split places. The first
place supplies condition (a.1), with the split Jacquet--Langlands
transfer the identity, and the other supplies condition (b.1).
All places with different inner form are real, with cuspidal
Euler--Poincar\'e tests, giving condition (c).
For each of the two specified tests,
\cite[Proposition~1.10.5(i)]{HarrisLabesse} identifies its cuspidal
trace with that of their common transfer on the quasi-split group.
This proves equality for these two tests with the stated
archimedean normalization.

On either side a nonzero trace of the archimedean test
requires cohomology with the chosen coefficient.
A representation selected by the finite test is supercuspidal
at both chosen split places. By
\cite[Proposition~1.8.1]{HarrisLabesse}, it occurs in the restriction
of a cuspidal representation of the group of rational unitary
similitudes. The two chosen split places satisfy case (B), so
\cite[Theorem~3.1.5]{HarrisLabesse} makes each contributing real
factor discrete series. Write \(F_\lambda\) for the representation of highest
weight \(\lambda-\rho\). The compact centre of \(U(a,b')\) lies in
the maximal compact \(K\), so a central character mismatch between
\(\pi\) and \(F_\lambda\) makes all relative cochains zero. When the
characters agree, the centre \(Z=U(1)\) acts trivially on
\(\pi\otimes F_\lambda^\vee\). Taking an \(r\)-th root of the determinant
gives \(U(a,b')=Z\,SU(a,b')\) and \(K=Z(K\cap SU(a,b'))\).
Restriction of each irreducible representation to the derived group
is irreducible because \(Z\) acts by scalars. The tangent quotient and
the invariant cochains therefore agree for the two groups. On the
derived group the parameter determines the ordered entries up to a
common additive constant and permutation; the fixed central character
determines their sum. The complex and compact Weyl groups are
\(S_r\) and \(S_a\times S_{b'}\), respectively. By
\cite[Ch.~II, Theorem~5.3(a)--(b)]{BorelWallach}, the discrete series satisfy
\[
 \dim H^j(\mathfrak g,K;\pi\otimes F_\lambda^\vee)
 =\dim\operatorname{Ext}_{(\mathfrak g,K)}^j(F_\lambda,\pi)
 =\begin{cases}
  1 & \pi\text{ belongs to the packet of }\lambda,\quad j=ab',\\
  0 & \text{otherwise}.
 \end{cases}
\]
Here \(ab'=\dim_{\mathbb R}(U(a,b')/K)/2\).
For each fixed finite representation, the multiplicities of
all real packet members on the fixed unitary group are equal by
\cite[Proposition~3.1.7]{HarrisLabesse}. In the notation of that
proposition take \(S_0=\{v_1,v_2\}\) to be the two chosen split places and
\(S=S_0\cup S_\infty\). The multiplicity comparison uses the same group on both
sides. Its hypotheses on the finite part of \(S\) require
supercuspidality at both members of \(S_0\).

Let \(\mathcal P_\infty\) be the product of these real packets and
\(\pi_\infty^0\) its prescribed member. For each fixed \(\pi_f\), the
normalized archimedean trace factor is
\[
 \frac{(-1)^{\sum_v a_v b'_v}}{|\mathcal P_\infty|}
    \sum_{\pi_\infty\in\mathcal P_\infty}
                  m(\pi_\infty\otimes\pi_f)
      =m(\pi_\infty^0\otimes\pi_f).
\]
Here \(|\mathcal P_\infty|=\prod_v\binom r{a_v}\) and
\(a_vb'_v\) is even at each place since \(a_v+b'_v=r\) is odd.
At each real place, Section~\ref{sec:real-reproducing-coefficient}
gives an orthogonal coefficient projector \(P\) on any irreducible
unitary \(\sigma\), with \(P\sigma(g)P=c(g)P\), where \(c\) is
the coefficient of the prescribed unit vector in \(\pi\).
If \(P\ne0\), the translates of a unit vector in its range
have the same Gram matrices as the translates of the vector
in \(\pi\). Cyclicity gives \(\sigma\simeq\pi\), where
Schur orthogonality makes the projector rank one.
The product selects one line in each occurrence of
\(\pi_\infty^0\), with the displayed multiplicity.
The trace equality is the equality of dimensions of the
two original tuple spaces. Sum over the orthogonal tuples
to get \(D=\dim pT_W\).

\smallskip\noindent\emph{The full finite spherical algebra.}
The trace argument compares the projected spaces on \(Y_0\) and on
any \(W_j\), \(1\le j\le m\). At each real place the ordered
parameter of Section~\ref{sec:real-marked-filling} has sign sequence
\(k-1\) plus signs followed by \(m\) alternating signs beginning with
plus on \(Y_0\). On \(W_j\) the same ordered values have the first
\(k+j-1\) signs positive and the remaining \(m-j\) negative, as in
Section~\ref{sec:bottom-descent-real-type}. In particular \(W_m\) has
signature \((r,0)\). Thus all groups use the same algebraic coefficient
in the Euler--Poincar\'e test, while the matrix coefficient selects
the prescribed member on each real group.

At a finite place put \(p_j=k+j-1\) and \(q_j=m-j\). In the
ordered rational basis of \(W_j\) its form is
\(J_j=\operatorname{diag}(I_{p_j},-I_{q_j})\). If \(u_-\) is an
integral unit of norm \(-1\), then
\[
 B_j=\operatorname{diag}(I_{p_j},u_-I_{q_j}),\qquad B_j^*J_jB_j=I_r .
\]
The unit exists at split places and at nonsplit places by surjectivity
of the unramified unit norm, including in residue characteristic two.
The norm and integral trace maps at nonsplit places are those of
Section~\ref{sec:splitting-level-characters}. For a hyperbolic pair \((e,e^*)\),
choose integral \(c\) of trace one. The basis
\((f,e-f)\), with \(f=e^*+ce\), has norms \(1,-1\) and zero cross
pairing. For the glued pair \((v_0,z)\) of
Section~\ref{sec:period-finite-coordinates}, choose instead
\(\operatorname{Tr}(c)=-N-1\) and put \(f=z+cv_0\). Then
\[
 h(f,f)=-N-\operatorname{Tr}(c)=1,\qquad
 h(f,v_0+f)=0,\qquad h(v_0+f,v_0+f)=-1 .
\]
Both basis changes have determinant \(-1\). Scale each negative line
by a unit of norm \(-1\). The formulas use only the integral glued coordinates.

At a split pair \((w,\bar w)\) let \(H_w,H'_w\) be the invertible
integral Gram matrices in the two lattice bases, identifying the
conjugate factor with the chosen local field. For an isometry from
the first lattice to the second take
\[
 C_w=I,\qquad C_{\bar w}^{\mathsf t}=H_w(H'_w)^{-1},\qquad
 C_{\bar w}^{\mathsf t}H'_w C_w=H_w .
\]
At a chosen pair its induced map on the prescribed general linear
residue group is the identity. Replacing either integral orthonormal
basis conjugates that map by a residue matrix and leaves the central
type idempotent unchanged. Use matching local Haar measures throughout. At a finite
place the Hermitian identification carries the
integral compact, its principal subgroup if present, and
their Haar measures to those of the other group.

In particular,
if \(K_v\) is the integral compact, use
\(\mathbf1_{K_v}/\operatorname{vol}(K_v)\) for the local
identity of the spherical algebra. A double coset \(K_v gK_v\)
uses \(\mathbf1_{K_vgK_v}/\operatorname{vol}(K_v)\), which acts
on \(K_v\)-fixed vectors by the unnormalized sum of its right
coset translates. These prescriptions agree on the two local groups
also at nonsplit finite places. Transport the finite tests by
the isomorphism of the groups there, so the tests
have the same stable orbital integrals with matching
centralizer measures. No comparison between distinct local
groups is needed at those places.

Fix a tuple \(c\) selected by the finite-type idempotent, and
insert any compactly supported spherical finite test \(h\)
off the old chosen places, with identity factors at all
but finitely many places. The test at the old level is the
central type projector of \(c\). There is a primitive type polynomial
at one split place above each of the two distinct rational
level primes. At the first such place the projector is a
linear combination of truncated diagonal coefficients of
unitary depth-zero supercuspidals, as in the computation above.
The two local groups there are split, so the
Jacquet--Langlands transfer is the same supercuspidal.
At the other distinguished place the corresponding coefficients
are cuspidal by the local supercuspidal coefficient
criterion \cite[Lemma~1.9.1]{HarrisLabesse}. The spherical
insertion is away from these places. It does not change either
local condition.

Apply \cite[Propositions~1.10.5(i) and~1.10.6]{HarrisLabesse} to
each of the specified tests and their common transfer, as above, in
case \((\mathrm{a}.1,\mathrm{b}.1,\mathrm{c})\). For two inner forms
of the quasi-split unitary group, associated tests have
identical cuspidal traces if at one split finite place the
test is a truncated coefficient of a supercuspidal whose
Jacquet--Langlands transfer is supercuspidal, at another
place the test is cuspidal, and at every place with
different local inner forms the test is cuspidal. This
applies termwise after expanding the first projector and
\(h\) into finite sums of decomposable tests. The two finite
conditions are exactly those verified above. The possible
differences of inner form are real. There use the
Euler--Poincar\'e tests with the same algebraic coefficient
and divide at signature \((a,b')\) by the packet size
\(\binom r a\). Their matching stable orbital integrals
are the preceding calculation with
\cite[Theorem~4.1]{ClozelLabesseOrbital}. The finite tests
are associated by their measure-preserving transport.

On the spectral side the same two primitive types force
the contributing finite factors to be supercuspidal there.
The cohomological real factors selected by the
Euler--Poincar\'e test are discrete series by
\cite[Theorem~3.1.5, case (B)]{HarrisLabesse}. At a fixed
finite representation the multiplicities of the packet
members on a fixed real inner form agree by
\cite[Proposition~3.1.7]{HarrisLabesse}. Thus dividing by
the packet sizes yields the trace on one prescribed real
member per place, with the sign \((-1)^{\sum_v a_vb'_v}=1\)
on both source spaces.

Let \(e_\infty\) be the product of the prescribed real matrix units.
Schur orthogonality gives a rank-one operator on
\(\pi_\infty^0\). If an irreducible unitary representation has
a nonzero vector in the image of this matrix unit, the compression
identity of Section~\ref{sec:real-reproducing-coefficient} identifies
the Gram matrix of its translates with that of the corresponding vector
in \(\pi_\infty^0\), up to its squared norm. Irreducibility therefore
forces that representation to be \(\pi_\infty^0\). With \(p_c\) the
finite projector density of the fixed tuple, one obtains
\[
 \operatorname{tr}R_{\rm cusp}(e_\infty\otimes p_c h)
   =\sum_{\pi_f} m(\pi_\infty^0\otimes\pi_f)
                      \operatorname{tr}\pi_f(p_c h).
\]
This is also the trace computed by the packet-normalized
Euler--Poincar\'e functions above. Summing over the orthogonal selected tuples gives
\begin{equation}\label{eq:marked-full-spherical-trace}
 \operatorname{tr}_{V}(h)=\operatorname{tr}_{V_j}(h)
 \qquad
 \bigl(h\text{ in the common finite spherical algebra}\bigr),
\end{equation}
where \(V\) and \(V_j\) are the original projected spaces on
\(Y_0\) and \(W_j\), respectively. The identity test gives
the equality of dimensions.

The preceding trace comparison also holds before summation, on the
isotypic spaces for each fixed irreducible residue type \(\tau_c\).
Write these spaces as \(\tau_c\otimes L_c\) and
\(\tau_c\otimes L_{j,c}\), with their invariant Hermitian forms.
Every spherical test off the old level commutes with the residue group.
Its action is therefore the identity on \(\tau_c\) tensored with its
action on the multiplicity space. Dividing the trace equality by
\(\dim\tau_c\) gives equal traces on \(L_c\) and \(L_{j,c}\).

Let \(\mathcal A_c\) be the image of the full spherical algebra on
\(L_c\oplus L_{j,c}\). It contains the identity and is closed under
adjoints since \(h^*(g)=\overline{h(g^{-1})}\). A finite-dimensional
complex algebra of operators closed under adjoints is a direct sum
of full matrix algebras. In each matrix summand a rank-one projection
has trace on either multiplicity space equal to the number of copies
of the defining module. Such a projection lies in \(\mathcal A_c\),
so it is the image of a finite spherical test and the two numbers
agree. Matching orthonormal copies in every summand gives a unitary
intertwiner \(L_c\simeq L_{j,c}\). Tensoring with the identity on
\(\tau_c\) preserves the old residue action and intertwines all
outside-level Hecke tests, including the nonsplit factors. In
particular it carries each prescribed joint eigenspace onto the
corresponding eigenspace on the other group.

At the positive definite endpoint each maximal-level class
with stabilizer \(\Gamma_c\) gives \(|G|/|\Gamma_c|\) classes
at the principal level. The quotient of the integral compact by the principal subgroup is
\(
 G=\prod_{i=0}^1\mathrm{GL}_r(\mathbb F_{p_i^{\rm lev}})^f
\), by reduction at the selected split places.
The rational stabilizer acts freely on this quotient. To verify
freeness choose a completely split
prime of norm \(p_i^{\rm lev}>(2r)^2\). If a stabilizer
element \(\gamma\) reduces to one at the chosen places, its
integral trace difference \(\operatorname{tr}\gamma-r\)
is divisible by all \(f\) of these primes. A nonzero
difference would have absolute norm at least
\((p_i^{\rm lev})^f\), and at most \((2r)^{2f}\) because
all its complex embeddings are unitary.

Thus the trace
equals \(r\). Unitarity forces all eigenvalues to be one
and \(\gamma=1\). This also gives the trivial stabilizer on
each finer class. At such a class the real component is the full
compact group. In the Peter--Weyl decomposition the coefficient
idempotent selects one line in the right factor of the terminal type,
leaving the dimension of its left factor.
For the dominant local parameters \(\lambda_{v,1}>\cdots>\lambda_{v,r}\),
the Weyl dimension formula for the highest weight
\(\lambda_v-\rho\) gives
\[
 \dim\bigotimes_{v\mid\infty} F_{\lambda_v}
    =\prod_{v\mid\infty}\prod_{i<j}
      \frac{\lambda_{v,i}-\lambda_{v,j}}{j-i}=\mathcal B .
\]
The normalization of the same product at each place is given in
Section~\ref{sec:real-formal-degree-factor}.
Summing over the maximal-level classes gives
\[
 \dim T_W
    = |G|\mathcal B\sum_c|\Gamma_c|^{-1}
    = |G|\mathcal B\,\operatorname{mass}(E,r).
\]
This proves \eqref{eq:marked-tower-dimension}.
\end{proof}

\section{Finite orbital integrals and Fourier bounds}
\label{sec:finite-orbital-analysis}

\sectionopening{In this section, we compute the finite Haar factors and bound the correlations appearing
in Section~\ref{sec:full-orbital-blocks}. Integral norm fibers give the
successive quotient measures, and Fourier estimates control the sum
over denominators.}

\subsection{Integral points on the local quadrics}
\label{sec:compact-floor-volume}

\sectionopening{In this subsection, we bound the measure of integral solutions of a Hermitian norm
equation. The bound enters the Haar disintegration in
Section~\ref{sec:compact-floor-haar} and the local integrations in
Section~\ref{sec:finite-position-box}.}

Let \(T\) be any finite set of finite places of the CM field \(E\)
stable under conjugation, with unramified or split quadratic
algebras at its conjugation orbits. Put
\(\mathcal R=\prod_{w\in T}\mathcal O_{E_w}\) and
\(\zeta_{\mathcal R}(a)=\prod_{w\in T}(1-(\mathrm Nw)^{-a})^{-1}\).
Normalize additive measures to give the integer rings mass one.
We need a bound for the mass of the integral box in a quadric
with an orthogonal block of \(2a+4\) orthonormal coordinates, where
\(a\ge1\) is an integer. Take the fiber measure by division of the
ambient Haar measure by the norm-value Haar measure. For a
skew-Hermitian norm equation at a split pair, \((1,-1)\)
generates the integral skew lattice over the fixed integers. At an
unramified quadratic place lift a residue element outside the
fixed field to \(a\). Then \(s=a-\bar a\) is a unit, also in
residue characteristic two. Every skew integer \(z\) has
\(z/s\) fixed and integral. Thus division by \(s\) identifies
the two integral norm-value lattices with Jacobian one.

Consider one orbit, with fixed field \(K\), residue cardinality
\(q\), quadratic algebra \(L\), and integer ring \(\mathcal O_L\).
Both the trace on integers and the norm on units are surjective.
In the unramified case this follows from the residue trace and
norm and the trace on each successive unit quotient. In the split
case these maps are addition and multiplication. Choose \(c,\tau\)
with \(c\bar c=-1\) and \(\tau+\bar\tau=1\).

In each orthonormal
pair \(v,w\), the integral basis
\[
                 v+cw,\qquad \bar\tau v-\tau cw
\]
is hyperbolic. Its change-of-basis determinant is \(-c\), a unit.
The norm on the whole block, with \(h=a+2\), is now
\(\operatorname{Tr}_{L/K}\sum_{i=1}^h\bar x_i y_i\).
For fixed \(x\in\mathcal O_L^h\), let \(j\) be the minimum
valuation of its coordinates, taking the minimum over both
components as well in the split case. For finite \(j\), the
linear norm equation in \(y\) is a map onto
\(\pi^j\mathcal O_K\) on \(\mathcal O_L^h\), where \(\pi\) is
a uniformizer of \(K\).

Its pushforward of integral Haar measure
is uniform probability measure on that ideal. Thus the density
on each norm fiber is at most \(q^j\). This is the density in
the quotient measure, since solving for a \(K\)-coordinate
with nonzero coefficient divides its Haar measure by the
absolute value of that coefficient. One may use any such
coordinate. The set \(x=0\) has measure zero in the smooth
fiber charts, by solving in \(x\) when \(y\ne0\) and using
\(h\ge3\). At the singular point the smooth fiber measure
has zero mass.

The probability that all the \(x_i\) are divisible by \(\pi^j\)
is \(q^{-2hj}\). In particular the same calculation gives the exact
density at every congruence level. Write
\(Q(x,y)=\operatorname{Tr}_{L/K}\sum_i\bar x_i y_i\).
For \(A\ge1\) and \(c\in\mathcal O_K\), one has
\[
\begin{aligned}
 q^A\Pr\{Q(x,y)\equiv c\bmod\pi^A\}
  ={}&(1-q^{-2h})\sum_{j=0}^{A-1}
       q^{-(2h-1)j}\mathbf1_{\pi^j\mid c}
       +q^{-(2h-1)A}\mathbf1_{\pi^A\mid c}.
\end{aligned}
\]
Indeed, conditional on minimum valuation \(j<A\), the congruence
has probability \(q^{j-A}\mathbf1_{\pi^j\mid c}\). On
\(x\in\pi^A\mathcal O_L^h\) it holds exactly when \(\pi^A\mid c\).
The series obtained by letting \(A\to\infty\) converges uniformly in
\(c\), including \(c=0\), since its tail is bounded by
\(q^{-(2h-1)A}/(1-q^{-(2h-1)})\). Hence the block contributes at most
\[
 (1-q^{-2h})\sum_{j\ge0}q^{-(2h-1)j}
 \le (1-q^{-(2a+3)})^{-1}.
\]
This bound is independent of the norm value, so the other
integral coordinates can be integrated with total mass one.
The displayed bound is at most \((1-q^{-2a})^{-1}\) in the
unramified quadratic case and at most \((1-q^{-a})^{-2}\)
in the split case. Multiplying over the orbits gives
\(\zeta_{\mathcal R}(a)\).

\begin{lemma}[Null set of a finite pivot]\label{lem:finite-pivot-null}
In a Hermitian norm equation over the semilocal algebra above, suppose
there is a hyperbolic pair orthogonal to a block of \(2a+4\) integral
orthonormal coordinates, where \(a\ge1\) is an integer. On every fixed
norm fiber the locus where either coordinate of this pair vanishes in
any field component has fiber measure zero.
\end{lemma}

\begin{proof}
Fix one of the hyperbolic coordinates \(t\) and one component \(w\in T\).
For fixed coordinates outside the orthonormal block, its integral fiber
mass is at most \(\zeta_{\mathcal R}(a)\) at every remaining norm value.
The congruence densities of the block converge uniformly in that value,
so the same bound applies after integrating over the other integral
coordinates. In particular
\[
 \operatorname{vol}_{\rm fiber}\{X\text{ in the box}:
                 t_w\in\mathfrak p_w^m\}
 \le \zeta_{\mathcal R}(a)(\mathrm Nw)^{-m}.
\]
Only the selected field component is restricted at a split pair.
Letting \(m\to\infty\) gives zero mass on \(t_w=0\), also on the
zero norm fiber. Dilation by a fixed-field scalar reduces any larger
box to an integral one, rescaling the norm value and the fiber measure
by fixed factors. Exhaust the fiber by these boxes and take the finite
union over the two coordinates and the components.\qedhere
\end{proof}

\subsection{Compact volumes for the successive quotients}
\label{sec:compact-floor-haar}

\sectionopening{In this subsection, we disintegrate the local Haar measure along successive isotropic
lines. The quadric estimate of Section~\ref{sec:compact-floor-volume}
bounds the resulting densities, and the compact mass is used in
Section~\ref{sec:finite-bottom-normalization}.}

Use the original integral Gram lattice \(\Lambda\) at \(T\), with
the \(\varrho\) hyperbolic pairs \(e_i,e_i^*\) and last vector \(z\) in their orthogonal complement. Write \(G_T\)
for its unitary group and \(K_T\) for the lattice stabilizer.
At the top use the isotropic dual \(u\) of \(z\) in the specified middle block, with \(h(u,z)=1\), and the orthogonal complement
of \(u,z\). The integral Gram restrictions used here are
that these are integral decompositions, that \(\langle u,z\rangle\)
has determinant \(-1\), and that the remaining orthonormal block
has at least \(2(3\varrho+5)+4\) coordinates. With a trace-one
integer \(\tau\), replacing \(z\) by \(z-\tau h(z,z)u\) gives
a hyperbolic pair and preserves the coefficient of \(z\).
In the original orbital lattice of Section~\ref{sec:period-finite-coordinates},
write \(v_0=x_0/\bar\pi_N+z/N\) for the glued basis vector. Its Gram
matrix with \(z\) gives \(u=-v_0\), \(h(u,u)=0\) and \(h(u,z)=1\),
including at the places dividing \(N\). The remaining \(k-1\)
orthonormal coordinates suffice since \(k-1\ge 2(3\varrho+5)+4\)
for large \(x\). Thus the same integral decomposition works for every \(T\).

Set \(n=\varrho+1\), and \(f_i=e_i\) for \(i<\varrho\),
\(f_\varrho=z\). Let \(G_i\) fix the pairs below \(i\) pointwise
and, for \(i=n\), also \(z\). Let \(R_i\) act unitriangularly on
the first \(i\) basis vectors, so \(g f_j-f_j\in
\langle f_0,\ldots,f_{j-1}\rangle\) for \(j<i\).
Let \(U_i\) be the Heisenberg group attached to the \(i\)-th
pair for \(i<\varrho\), and put \(U_\varrho=\{1\}\).
The products \(U_{<i}=U_0\cdots U_{i-1}\) are normalized by
\(G_i\). We have \(R_i=G_i U_{<i}\) and
\(R_{i+1}\cap G_i=G_{i+1}U_i\).

Indeed an isometry fixing \(e_i\) sends \(e_i^*\) to
\(e_i^*+v+s e_i\), with \(v\) in the remaining space and
\(s+\bar s=-h(v,v)\). The unique Heisenberg element with that
column removes it, leaving an isometry fixing the pair.
For the iteration, write \(r\in R_{i+1}\) as \(\gamma v\) by
the decomposition at \(i\). A lower Heisenberg factor changes
each later basis vector only by earlier isotropic vectors, so
\(v\in R_{i+1}\).

Now \(\gamma\in G_i\cap R_{i+1}\) takes
\(f_i\) into the orthogonal complement of the earlier pairs
and differs from \(f_i\) by a vector in their span. Thus it
fixes \(f_i\), and the column removal gives the step. At the
top the condition within \(G_\varrho\) is exactly to fix \(z\).

Let \(\eta_i\) be Heisenberg Haar in the coordinates
\((v,\beta)\), with \(s=\beta-\tau h(v,v)\) and \(\beta\) skew.
Its measure is product additive Haar, since multiplication
is triangular with translations in \(v,\beta\). Normalize each
integer coordinate lattice to mass one. The displayed
column and integrality of \(h\) give
\(\eta_i(K_T\cap U_i)=1\), and peeling off these columns
gives \(\eta(K_T\cap U_{<n})=1\) for
\(\eta=\eta_0*\cdots*\eta_{\varrho-1}\).
Conjugation by \(G_i\) preserves the earlier vector measures
by the unitary determinant identity and fixes their centers.

For \(i<n\) put \(H_i=R_{i+1}\cap G_i=G_{i+1}U_i\).
The map
\[
 H_i\backslash G_i\longrightarrow G_i f_i,
           \qquad H_i\sigma\longmapsto\sigma^{-1}f_i
\]
identifies the left quotient with the orbit. Give the
quotient the norm-fiber measure under this identification.
In convolution, \(\theta_i\) is the image of this measure
under the section \(X\mapsto\sigma\) constructed below
on the chart where the pivot \(t\) is invertible in every
component. This chart has full fiber measure. On a hyperbolic pair with coordinates \(t,s\)
write \(s=t c\). This contributes the additive Jacobian
\(|t|_A\), the modulus on the semilocal algebra.

The norm
equation takes the trace of \(c\) with coefficient \(t\bar t\),
whose modulus on the fixed algebra is also \(|t|_A\).
The factors cancel. Writing \(c=\tau b+\beta\), with \(b=\operatorname{Tr}(c)\)
and \(\bar\beta=-\beta\), splits the integral lattice as
\(\tau\mathcal O_K\oplus\ker(\operatorname{Tr}|_{\mathcal O_L})\)
at each conjugation orbit. Thus the trace splitting has Jacobian
one for the chosen measures, including above \(2\).
The other coordinates, \(t\), and \(\beta\) have product additive
measure. The top chart uses the same computation after the indicated
change of basis. At every stage the unused orthonormal block has at
least \(2(3\varrho+5)+4\) coordinates. Lemma~\ref{lem:finite-pivot-null}
therefore gives full measure to the chart at each prescribed norm.

On a lower stage write a point of norm zero as
\(X=t e_i^*+x+s e_i\), with \(t\in A_T^\times\) and
\(x\perp e_i,e_i^*\). Put \(a=x/t\) and
\(\beta=s/t+\tau h(a,a)\). The norm equation gives
\(\beta+\bar\beta=0\). Let \(w_i(t)\) send
\(e_i\mapsto t e_i^*\) and \(e_i^*\mapsto\bar t^{-1} e_i\),
fixing their orthogonal complement. The Heisenberg element
\(n(a,\beta)\) fixes \(e_i\) and sends
\(e_i^*\mapsto e_i^*+a+(\beta-\tau h(a,a))e_i\). Hence
\[
 n(a,\beta)w_i(t)e_i=X,\qquad
 \sigma=w_i(t)^{-1}n(a,\beta)^{-1},\qquad \sigma^{-1}e_i=X .
\]
The change \(x=ta\) contributes
\(\operatorname{md}(t)^{d_i}\), where
\(d_i=k+1+2(\varrho-i-1)\) counts the remaining coordinates, as in
Section~\ref{sec:finite-common-entry}.
On the top stage put \(v=z+N\tau u\), where \(h(z,z)=-N\).
For the original exceptional vector \(x_{k_0}/\bar\pi+z/N\)
we have \(u=-(x_{k_0}/\bar\pi+z/N)\). Then
\(h(u,u)=h(v,v)=0\), \(h(u,v)=1\) and \(z=v-N\tau u\).
The integral span of \(u,v\) is unimodular and splits off.

Given \(tv+x+su\) of norm \(-N\), with \(t\) invertible and
\(x\perp u,v\), put
\[
 a=t^{-1}x,\qquad
 \beta=t^{-1}s+\tau h(a,a)+N\tau/(t\bar t).
\]
The norm equation gives \(\beta+\bar\beta=0\).
First dilate \(v\) by \(t\) and \(u\) by \(\bar t^{-1}\).
Then apply the Heisenberg element fixing \(u\) and sending
\(v\) to \(v+a+(\beta-\tau h(a,a))u\), with
\(c\mapsto c-h(a,c)u\) on the orthogonal complement.
This composite sends \(z\) to the given vector.
For the orbit map \(\sigma\mapsto\sigma^{-1}z\), take
\(\sigma\) to be the inverse composite. If the given vector
is integral and \(t\) is a unit, both transformations and
their inverses preserve the lattice. Thus the
unit-pivot subset contributes to the compact mass in
the same fiber measure.
The terminal measure \(\theta_n\) is right Haar on \(G_n\)
with \(K_T\cap G_n\) of volume one.

The Haar measure given by these data is
\(\eta*\theta_n*\cdots*\theta_0\). Indeed, in descending order
the group measure is formed by taking \(\eta_i*\theta_i\) at
stage \(i\), with \(\eta_i=\delta_1\) on the top quadric and
terminal stages. The measure from the higher stages times
\(\eta_i\) is right invariant under \(G_{i+1}U_i\).
Consequently its integral against a function of the product
with \(\sigma\in G_i\) descends to \(H_i\sigma\), so it is
independent of the section used to lift \(\theta_i\).
Right translation by \(g\in G_i\) acts on the fiber coordinate
by \(g^{-1}\).

This preserves the norm fiber measure, because
the ambient determinant has modulus one and the norm value
is unchanged. Induction from \(\theta_n\) gives right Haar
at the bottom. Each residual unitary group, including the terminal
stabilizer, is unimodular. Over a splitting field its adjoint action
on the Lie algebra is conjugation on \(M_d\), with determinant
\((\det g)^d(\det g^{-1})^d=1\). The same determinant is therefore
one over the fixed local field, so right Haar is also left invariant.
This gives the quotient integration formula at each step with the
stated fiber and stabilizer measures.
Moving each \(\theta_i\) past the lower Heisenberg product is
permitted by the conjugation invariance above. This gives the
asserted order, with the same measures as for the orbital
integrals formed using this stage normalization.

Put \(\mathcal U_i=R_iK_T=\{g:R_ig\cap K_T\ne\varnothing\}\). For
each \(g\), the \(\theta_i\)-mass of \(\sigma\) such that
\(\sigma g\in\mathcal U_{i+1}\) is at most
\(V=\zeta_{\mathcal R}(3\varrho+5)\). If the set is nonempty,
then \(g\in\mathcal U_{i}\), since \(G_i,R_{i+1}\subseteq R_i\).
Replacing \(g\) by \(r g\in K_T\) with \(r\in R_i\) preserves
the integral. To see this, factor \(r=\gamma v\) in
\(G_iU_{<i}\). The conjugate of \(v\) by \(\sigma\gamma\)
lies in \(U_{<i}\subseteq R_{i+1}\) and disappears from the
near condition. Right multiplication by \(\gamma\) preserves
the integral on the orbit. Now if \(\xi\sigma g\in K_T\)
with \(\xi\in R_{i+1}\), then \(\sigma^{-1}\xi^{-1}f_i\)
is integral. Its projection away from the earlier pairs is
\(\sigma^{-1}f_i\), since \(\xi^{-1}f_i-f_i\) belongs to
the earlier isotropic span. The projection is integral.
The quadric box bound of Section~\ref{sec:compact-floor-volume}
therefore gives \(V\), uniformly in \(g\).

The \(\theta_n\)-mass of \(h\in G_n\) with
\(hg\in U_{<n}K_T\) is at most one. In fact if
\(u_1h_1g,u_2h_2g\in K_T\), their quotient belongs to
\(K_T\) and to \(U_{<n}G_n\). Its middle block is
\(h_2h_1^{-1}\), which is therefore integral with integral
inverse. All such \(h\) lie in a single right coset of
\(K_T\cap G_n\). In the descending convolution
\(\theta_n*\cdots*\theta_i\), the integral of the indicator
of \(U_{<n}K_Tg^{-1}\) is thus at most \(V^{n-i}\).

For the induction, the higher convolution is supported
in \(G_{i+1}\). Its integral after \(\sigma g\) can be
nonzero only when \(\sigma g\in\mathcal U_{i+1}\), since
\(U_{<n}\) and \(G_{i+1}\) lie in \(R_{i+1}\).
Integrating over \(\theta_i\) contributes at most \(V\).
Finally convolution with \(\eta\) gives exactly the
integral over \(U_{<n}K_T\) for \(g=1\), because each
nonempty \(\eta\)-fiber in \(K_T\) is a right coset of
\(K_T\cap U_{<n}\), of mass one. Thus \(\mu(K_T)\le V^n\).
For \(a=3\varrho+5\), use \(n=m/2\), \(\sqrt x\le m<2\sqrt x\),
and \(f\le x^{12}\) for large \(x\) to obtain
\[
 n[E:\mathbb Q]\le 2x^{25/2}\le 2^{\sqrt x-2}\le 2^{2\varrho}.
\]
Then \(\zeta_{\mathcal R}(a)\le\zeta_E(a)\le\zeta(a)^{[E:\mathbb Q]}\) gives
\[
 n\log\zeta_{\mathcal R}(a)
 \le n[E:\mathbb Q]\,2^{1-a}\le 2^{-n-3}.
\]
Here \(\log\zeta(a)\le\sum_{j\ge2}j^{-a}\le2^{1-a}\) by
the integral test. Exponentiating gives the upper bound
\(1+2^{-n}\) for the compact volume.

\subsection{Normalization of the finite quotient pairing}
\label{sec:finite-bottom-normalization}

\sectionopening{In this subsection, we compute the correlation at the initial step of the finite
calculation. The Haar normalization of
Section~\ref{sec:compact-floor-haar} fixes the constant in the bound
of Section~\ref{sec:finite-selected-bottom}.}

Take a conjugation-stable set \(T\) as in
Section~\ref{sec:compact-floor-volume}. The
local volume and correlation estimates below also apply when \(T\)
is a single conjugate pair. We include every place over \(M\) when
we substitute \(M^{-2f}\) for its finite modulus in the adelic estimate.
All local parameters in this calculation are restricted to
\(A_T=\prod_{w\in T}E_w\), with integer ring \(\mathcal R\)
and additive modulus \(\operatorname{md}\). At a place over
\(M\), \(N\) and \(\pi_N\) are units. Thus the same norm
and lattice identities hold when \(T\) consists only of
the places over \(M\), as in
Section~\ref{sec:marked-residue-projection}.

Section~\ref{sec:fixed-tail-average} uses a set containing
the places over \(N\) when evaluating the full modulus of
\(\pi_N\).
Every conjugation orbit of completions is split or unramified
over its fixed field. In fact \(d\nmid q_1\). At a rational
prime other than \(d\), inertia acts trivially on \(K\), where
complex conjugation is nontrivial. At \(d\), inertia acts
trivially on \(\mathbb Q(\zeta_{q_1})\), where complex
conjugation is also nontrivial. Thus the inertia groups in
\(E/E^+\) are trivial, where \(E^+\) is the fixed field of
conjugation. In particular the previous compact measure
calculation can be used on all of \(T\).

Write \(G_T\) for the original local group, \(\mu\) for its
stage Haar measure, \(K_T\) for the compact integral group
and \(J\) for its congruence subgroup. Let \(Z\) be
an intermediate group \(J\le Z\le K_T\). Put
\(U=U_{<n}\), \(dQ=d\mu/d\eta\) on \(U\backslash G_T\),
where \(\eta\) is the original unipotent Haar, and
\[
 \delta_U=\eta(U\cap J),\qquad
 \mathcal R_Z=\mu(Z)\mu(J)/\delta_U^2 .
\]
Extend the character of \(U\) to \(UJ\) by right
\(J\)-invariance and then by zero. Denote it by \(\chi\).
This is well-defined since the character is trivial on
\(U\cap J\). The finite Bessel vector \(W_Z(y)\) has entry
\(\chi(yc^{-1})\) for \(c\in J\backslash Z\).
Its right \(J\)-invariance follows from normality of \(J\)
in \(K_T\).

Write \(W=W_Z\). For each \(c\) the set of \(u\) such that
\(uyc^{-1}\in J\) is empty or a right coset of \(U\cap J\)
of mass \(\delta_U\). For such \(u\), \(z_1=uy\) lies in
\(Jc\) and
\[
 \overline{W(y)_c}\,W(yx)_c=\chi(z_1 x z_1^{-1}).
\]
Indeed the phase of the left \(U\) translation cancels,
and the additional right factor \(z_1c^{-1}\) lies in \(J\).
The cosets partition \(Z\), so quotient integration gives
\[
 \int_Z\chi(z_1xz_1^{-1})\,d\mu(z_1)
              =\delta_U\int_{U\backslash G_T}
                                \langle W(y),W(yx)\rangle\,dQ(y).
\]
With empty heads the stage-zero form \(\mathcal I_0^T(1)\)
is the integral of the squared modulus of this quotient
correlation. Define
\(\mathcal M(Z)=\int_{G_T}|\int_Z\chi(z_1xz_1^{-1})\,d\mu(z_1)|^2\,d\mu(x)\).
The normalized bottom value is therefore
\[
 \mathcal I_0^T(1)/\mathcal R_Z
                    =\frac{\mathcal M(Z)}{\mu(Z)\mu(J)} .
\]
In particular it is nonnegative. On \(x\in J\) the inner
integral is \(\mu(Z)\) by normality. The upper bound on the
rest of \(\mathcal M(Z)\) requires a finite shell estimate.

For use in that estimate, the coordinates of \(W\) on the
quotient have supports given by the images of \(UJc\).
Each support has mass \(\mu(J)/\delta_U\) by unfolding
the indicator of \(Jc=cJ\). These supports are equal or
disjoint. If \(UJc'=UJc\), write \(c'=ujc\) with \(u\in U\cap Z\) and
\(j\in J\). On the common support,
\(W(y)_{c'}=\chi(u)^{-1}W(y)_c\) by \(U\)-covariance and normality of \(J\). The multiplicity in each part is
\([U\cap Z:U\cap J]\le1/\delta_U\), since \(\eta(U\cap K_T)=1\).
Thus \(B_W=\int WW^*\,dQ\) is a direct sum of rank-one blocks
of norm at most \(\mu(J)/\delta_U^2\), with trace
\(\mu(Z)/\delta_U\). It follows that
\[
 \delta_U\iint|\langle W(y),W(y')\rangle|^2\,dQ(y)\,dQ(y')
       =\delta_U\operatorname{tr}(B_W^2)\le\mathcal R_Z .
\]
The matrix \(WW^*\) is well-defined on the quotient by its
scalar \(U\)-covariance.

\subsection{Common entries and the finite terminal density}
\label{sec:finite-common-entry}

\sectionopening{In this subsection, we express the finite pairing at successive isotropic quotients in
measures compatible with
Section~\ref{sec:finite-bottom-normalization}. The expression supports
the Fourier estimates of Section~\ref{sec:finite-frequency} and the
terminal density calculation.}

Use these same measures and \(W=W_Z\). For \(a\le\varrho\)
denote the finite stage form by \(\mathcal I_a^T(p,q;\omega)\).
It has the two pivot integrations of
Section~\ref{sec:real-common-descent}, now in the finite stage
measures, with head factor
\[
 \mathcal D_a^T(p,q)=\prod_{i<a}\operatorname{md}(p_iq_i)^{d_i/2},
 \qquad d_i=k+1+2(\varrho-i-1).
\]
The twist is \(w_q^{-1}w_p\). For \(a<\varrho\), the law
of the stage measure \(\theta_a\) is that of
\(w_a(c)^{-1}\nu\) for
\(\operatorname{md}(c)^{d_a}dc\,d\eta_a(\nu)\).
To see the density, in the norm chart use the coefficient
\(c\) of \(e_a^*\) and take \(\nu^{-1}\) acting on
\(w_a(c)e_a=c e_a^*\). Its \(d_a\) free vector coordinates
are multiplied by \(c\), and its skew coordinate in the
ratio of the two hyperbolic coordinates retains its
original Haar. The Weyl element fixes the earlier pairs.
The omitted tuples have a zero field component of \(c\)
and form a null set in this stage chart.

Assume the absolute stage majorant is integrable at \(a\).
With \(Z_a(s)=h(f_a,s f_a)\), the substitution
\(\tau=w_a(c)^{-1}\nu\sigma\), \(g'=\nu\sigma g\),
is valid by the quotient invariance of \(\theta_a\) on
functions left invariant under the stabilizer of \(f_a\).
The pivot integrals have this invariance after the tail
convolution, including for their absolute majorant. The new
point measure is \(\eta_a*\theta_a*M_a=M_{a+1}\).
The old twist is unchanged and the new head factor is
\(\mathcal D_{a+1}^T(p_c,q_c)=
\mathcal D_a^T(p,q)\operatorname{md}(c)^{d_a}\).
Thus for every measurable \(|\phi|\le1\), absolute
integrability transfers through the substitution and
\[
 \int_{A_T}\phi(c)\,\mathcal I_{a+1}^T(p_c,q_c;1)\,dc
             =\mathcal I_a^T(p,q;\phi\circ Z_a).
\]
Here \(p_c,q_c\) append the same entry \(c\); one can use
the entry \(1\) when \(c\notin A_T^\times\), a Haar-null set.
This lower-stage identity uses the integrable majorant
on the original projected Bessel data.

For the terminal formula let \(\pi=\pi_N|_T\) with
\(\pi\bar\pi=N\) and \(h(z,z)=-N\).
In \(V_\varrho=X_0\perp A_T z\), the last-entry condition
and the fixed norm are
\[
 e(t,x)=x-N^{-1}tz,\qquad h(z,e(t,x))=t,\qquad
 h(e(t,x),e(t,x))=-N\ \Longleftrightarrow\
                  h(x,x)=t\bar t/N-N .
\]
In the chosen lattice basis the exceptional vector is
\(x_{k_0}/\bar\pi+z/N\). Thus the change to the basis
\((x_i,\pi^{-1}z)\) of Gram \(\operatorname{diag}(1,\ldots,1,-1)\)
has determinant \(\bar\pi/\pi\). Since \(T\) is stable under
conjugation, \(\operatorname{md}_T(\bar\pi)=\operatorname{md}_T(\pi)\).
The basis change therefore preserves the product Haar, taking both
places together at a split pair.
The coordinate of \(e(t,x)\) along \(\pi^{-1}z\) is
\(-\bar\pi^{-1}t\).

Use an integral hyperbolic pair in the
standard form on \(X_0\), as in
Section~\ref{sec:compact-floor-volume}, to solve the
norm equation both on \(X_0\) and on \(V_\varrho\).
With the norm value indexed by the same skew unit,
the other chart coordinates and their Haar agree. The change
in the last coordinate gives
\[
 dQ_\varrho(e)=\operatorname{md}(\pi)^{-1}\,
              d\nu_{t\bar t/N-N}(x)\,dt ,
\]
in the variables \(e(t,x)\) on the last-vector quadric.
Here \(d\nu_c\) is the original norm-fiber measure on \(X_0\)
at norm \(c\), with the integral skew line of volume one
in its norm-value coordinate.

These fiber measures vary weakly continuously on compact
step tests, also at zero norm. In a bounded test box, use
hyperbolic coordinates \(x_i,y_i\) on \(h\ge1\) hyperbolic planes in the unimodular subspace as in Section~\ref{sec:compact-floor-volume},
treat the remaining coordinates as parameters shifting the norm
value, and truncate the set
where all \(x_i\) in a given factor have large valuation.
Outside that set we can solve in a trace coordinate with
nonzero coefficient, on finitely many compact open pieces.
For fixed \(x_i\) the Haar integral in the remaining
coordinates is locally constant in the norm value, uniformly
on each piece.

The
density is bounded by a constant times \(q^j\) where \(j\)
is the minimum valuation, by the calculation of
Section~\ref{sec:compact-floor-volume} after scaling the box.
The measure in the \(x_i\) with this valuation is at most
a constant times \(q^{-2hj}\) for \(h\ge1\) hyperbolic pairs,
including at a split orbit. Thus the truncated tails tend
uniformly to zero. Taking products proves the assertion on
the semilocal ring. This continuity identifies the chart
description of the fiber measure at every norm value.

Let \(H=G_n\), the stabilizer of \(z\) in \(G_\varrho\),
with the terminal Haar measure \(dh\) satisfying
\(\operatorname{vol}(H\cap K_T)=1\), as in
Section~\ref{sec:compact-floor-haar}. Define
\[
 V_p(l)=\mathcal D_\varrho^T(p,1)
       \int_{G_T}W(y)\otimes W(lw_p^{-1}y)^*\,d\mu(y)
                      \qquad(l\in G_\varrho).
\]
These integrals have compact support in \(y,l\), locally
uniformly in the nonzero entries of the head \(p\).
Indeed nonzero integrands give
\(y=u_1^{-1}k_1,\ lw_p^{-1}y=u_2k_2\) with
\(u_i\in U,\ k_i\in K_T\). Thus \(u_1w_pl^{-1}u_2=k_1k_2^{-1}\).
In the two isotropic blocks and their middle block, the
lower-left block of the latter product is \(E_1D_p A_2\),
with \(D_p=\operatorname{diag}(p_i)\), \(E_1\) lower and
\(A_2\) upper unitriangular.

Write the three block rows of \(u_i\) as
\((A_i,C_i,B_i),(0,I,D_i),(0,0,E_i)\), and write
\(P=u_1w_pl^{-1}u_2\). Elimination in \(P_{31}=E_1D_p A_2\)
recovers \(E_1,A_2\) and their inverses using the nonzero
pivots \(p_i\). The remaining required blocks are obtained from
\[
 \begin{aligned}
 P_{32}&=E_1D_p C_2,& P_{21}&=D_1D_p A_2,
   & P_{22}&=l^{-1}+D_1D_p C_2,\\
 P_{11}&=B_1D_p A_2,& P_{33}&=E_1D_p B_2 .
 \end{aligned}
\]
The other blocks of \(u_1,u_2\), and the inverses of
\(u_1,u_2,l\), are bounded by the isometry equations in the
fixed Gram basis. These bounds are uniform for the head
in a compact neighborhood in \((A_T^\times)^\varrho\).
Since \(y=u_1^{-1}k_1\), they bound the two matrices \(y,l\)
and their inverses, so both range over compact sets. On these sets a small change
in the head and in \(l\) changes the arguments of \(W\)
by a right \(J\)-translation, uniformly by compactness.
Consequently \(V_p(l)\) is locally constant in both variables
with uniform neighborhoods on compact sets.

On the orbit of \(z\), set
\[
 \Phi_{p,q}(e)=\int_H\langle V_p(hl),V_q(hl)\rangle\,dh,
                           \qquad e=l^{-1}z .
\]
Right invariance of \(dh\) makes this independent of the
lift \(l\). Since \(N\in A_T^\times\), Witt's theorem in each factor
of the semilocal base shows that the orbit is the whole quadric
\(h(e,e)=-N\). This is a compactly
supported, locally constant function on the quadric. For
the local section at \(e_0=l_0^{-1}z\), take \(e\) close to
\(e_0\), and put \(b_0=h(e_0,e)\) and \(d_0=\bar b_0/b_0\).
Shrink the neighborhood in each local factor so that
\((d_0-1)/4\) lies in its maximal ideal. The polynomial
\(u^2+u-(d_0-1)/4\) has simple residue root \(0\).
Hensel's lemma gives a unique lift \(u\) in the maximal ideal,
and \(\lambda=1+2u\) is the resulting square root of \(d_0\)
near \(1\), including in residue characteristic two.

Uniqueness gives \(\lambda\bar\lambda=1\) and
\(\lambda b_0=\overline{\lambda b_0}\).
After multiplying \(e\) by \(\lambda\), reflection in
\(e_0+\lambda e\) sends it to \(-e_0\), and reflection
in \(e_0\) gives \(e_0\). Their product with the scalar
is a unitary \(s(e)\) tending to \(1\).
Then \(l_0s(e)\) is the section. The support of
\(\Phi_{p,q}\) is contained in a compact image under
\(l\mapsto l^{-1}z\) by the preceding support bound,
so \(\Phi_{p,q}\) has compact support.
On the ambient space it admits a compact step extension
by a finite partition into compact open balls.

For the top identity expand the inner product of the two
orbital tensors. For fixed nonzero \(p,q\), its absolute
majorant in \((l,y_1,y_2)\) has compact support by the
bounds just proved, applied to the \(p\) factors at \(y_1\)
and the \(q\) factors at \(y_2\). It is integrable for a
bounded entry test. Put \(y_2=y_1x\) by left invariance,
and push the positive Haar measure in \(y_1\) through
\(y_1=h_1\sigma g\), with \(h_1\in H\), \(\sigma\) in the
last stage and \(g\) in the lower product \(M_\varrho\).
For each such tuple use right invariance in \(l\), followed
by the positive Haar decomposition
\(l=h_2\tau(h_1\sigma)^{-1}\).
Tonelli transports the absolute majorant to these stage
variables. The weight \(\phi(h(z,l^{-1}z))\) becomes
\(\phi(h(z,\sigma\tau^{-1}z))\), since \(h_1,h_2\) fix \(z\).
Both earlier Weyl elements commute with \(G_\varrho\),
so Fubini gives the two terminal tail correlations.

The product \(\mathcal D_\varrho^T(p,1)\mathcal D_\varrho^T(q,1)
=\mathcal D_\varrho^T(p,q)\) gives the head density.
This proves
\[
 \begin{split}
 \mathcal I_\varrho^T(p,q;\phi\circ Z_\varrho)
 &=\int\phi(h(z,e))\Phi_{p,q}(e)\,dQ_\varrho(e)\\
 &=\int_{A_T}\phi(t)\operatorname{md}(\pi)^{-1}
        \int\Phi_{p,q}(e(t,x))\,d\nu_{t\bar t/N-N}(x)\,dt .
 \end{split}
\]
The inner expression with its modulus factor is the finite
last-entry form. For finitely many heads with the same
last entry \(t\) it is positive semidefinite, since the
pairing on every lift is a Gram matrix and both fiber
integrations are positive. It is continuous in the heads
and in \(t\) by the uniform support and the compact-test
continuity of the norm-fiber law. All stage forms, including
this terminal form, are divided by the same
\(\mathcal R_Z\) of Section~\ref{sec:finite-bottom-normalization}
in the finite sieve. The frequency and bottom bounds are
proved in Sections~\ref{sec:finite-frequency}--\ref{sec:finite-selected-bottom}.

\subsection{Stable ideals and classes of prescribed denominator}
\label{sec:finite-stable-shells}

\sectionopening{In this subsection, we determine which ideals can occur in the denominators of two
correlated vectors. Their indices give the summation variables for
Section~\ref{sec:finite-shell-reduction} with the lattice stabilizers
of Section~\ref{sec:compact-floor-haar}.}

Use the projected lattice \(\Lambda\), with its integral
unimodular Gram, and the original level scalar \(M\) on \(A_T\).
Let \(F_i=\langle e_0,\ldots,e_{i-1}\rangle\).
At a pair \(Y,Y'=Yx\) with nonzero plain Bessel correlation,
there are \(u,v\in U\) such that \(uY,vY'\in K_T\) and
\((vY')^{-1}uY\in J\). Indeed the two vectors have a common
nonzero coordinate of \(W_Z\), which gives both moves into
the same coset of \(J\). Consequently the support in \(x\)
of the absolute stage majorant is contained in
\(\Omega_J:=\bigl(\bigcup_Y Y^{-1}UY\bigr)J\).

For \(0\le i<\varrho\), let \(\mathfrak d_i(Y,x)\) be the set of \(h(e_i,w)\) with
\[
 w\in F_i^\perp,\quad w'-w\in F_i,\quad
 Y^{-1}w\in\Lambda,\quad Y^{-1}w-(Y')^{-1}w'\in M\Lambda .
\]
It is a full fractional ideal. Write \(w'=w+d\) and map
\((w,d)\) to the pair
\(\bigl(Y^{-1}w,M^{-1}(Y^{-1}w-(Y')^{-1}(w+d))\bigr)\).
The map on the full product is invertible. The preimage
of the product lattice intersects \(F_i^\perp\oplus F_i\)
in a lattice, and projection to \(h(e_i,w)\) is surjective
on this subspace. This gives a compact open integer module
in each field factor. At the active pair it is integral:
\(u^{-1}e_i\equiv e_i\bmod F_i\) and \(uw=(uY)Y^{-1}w\)
is integral, so \(h(e_i,w)=h(e_i,uw)\in\mathcal R\).
Write \(\mathfrak d_i=c_i\mathcal R\). Its depth is
\(D_i=\operatorname{md}(c_i)^{-1}\). The two transported vectors
in the definition are integral, and the congruence depends only on
their classes in \(\Lambda/M\Lambda\). Replacing \((Y,Y')\) by
\((Yj_1,Y'j_2)\) for \(j_1,j_2\in J\) preserves both conditions.
Thus the ideal and its depth are locally constant in \((Y,Yx)\), hence
in \((Y,x)\).

For any integral self-conjugate \(C\in A_T^\times\) in
all the \(\mathfrak d_i\), each \(w\in F_j^\perp\)
with \(Y^{-1}w\in\Lambda\) has a lift \(w'\equiv w\bmod F_j\)
such that
\[
 Y^{-1}w-(Y')^{-1}w'\in M C^{-2(\varrho-j)}\Lambda .
\]
For \(j=\varrho\), take \(w'=v^{-1}uw\).
The unipotents act trivially on \(F_\varrho^\perp/F_\varrho\),
and \((vY')^{-1}uY\in J\) gives the congruence.
Suppose the lift is available at \(i+1\) with exponent
\(D=2(\varrho-i-1)\). Divide a stable pair at the value
\(c_i\) by \(c_i\), giving \(\zeta,\zeta'\) of pairing one
with \(e_i\), with both transported vectors in \(C^{-1}\Lambda\)
and their difference in \(MC^{-1}\Lambda\).

For the prescribed \(w\) put \(t=h(e_i,w)\in\mathcal R\)
and \(w_0=w-t\zeta\). Apply the induction to \(Cw_0\)
and divide by \(C\). This produces \(w_0'\equiv w_0\bmod F_{i+1}\),
with transported difference in \(MC^{-D-1}\Lambda\).
The coefficient \(b_1\) of \(e_i\) in \(w_0'-w_0\) equals
\(h(\zeta',w_0')-h(\zeta,w_0)\), since
\(\zeta'-\zeta\in F_i\) and \(w_0\perp F_i\).
In the transported lattices, expand this difference by pairing
the transported \(\zeta'\) against the difference for \(w_0\),
and the difference for \(\zeta\) against the transported \(w_0\).

The two terms belong to \(MC^{-D-2}\mathcal R\). Thus
\(w'=t\zeta'+w_0'-b_1v^{-1}e_i\) is the required lift
at \(i\), since \((Y')^{-1}v^{-1}e_i\) is integral and
\(v^{-1}e_i\equiv e_i\bmod F_i\).
At \(j=0\), \(w'=w\). Testing \(w\in Y\Lambda\) gives
\(C^{2\varrho}(1-x^{-1})/M\) integral, hence
\(C^{2\varrho}(x-1)/M\) integral by unitary adjunction.
In particular if \(x\notin J\) and every proper integral
ideal has index at least \(Q_0\), some \(D_i\ge Q_0\).
Otherwise all \(c_i\) are units in \(\mathcal R\), and
\(C=1\) gives \(x\in J\).

On a head of length \(j\le\varrho\), assume also that
\(w_p^{-1}Y,w_q^{-1}Y'\in\mathcal U_{j}\). Choose
an integral self-conjugate scale \(\alpha\) divisible by
all the partial products \(P_a=\prod_{i<a}p_i\),
\(Q_a=\prod_{i<a}q_i\), and by \(M\) if \(j>0\).
The products here are integral at nonzero contributions.
Indeed take two moves \(s,t\in R_j\) so that
\(sY,t w_p^{-1}Y\in K_T\). Their ratio is \(s w_p t^{-1}\in K_T\).
The minor in the first \(a\) flag columns and the rows
indexed by \(e_0^*,\ldots,e_{a-1}^*\) equals \(P_a\) by triangularity.
For \(0\le a<j\) take the first \(a+1\) columns, the dual rows
\(e_b^*\) for \(b<a\), and a last row \(d\) outside those dual rows.
The right factor \(t^{-1}\) acts on these columns by a unitriangular
matrix. Removing it leaves the columns \(p_i s e_i^*\) for \(i\le a\).
For \(b\le i\), the \(e_b^*\)-entry of \(s e_i^*\) equals
\(h(s^{-1}e_b,e_i^*)=\mathbf1_{b=i}\). Expansion along the last
column gives the integral minor \(P_{a+1}(s e_a^*)_d\). The entries in
the omitted dual rows of \(s e_a^*\) are zero. Thus \(Y^{-1}e_a^*\in\alpha^{-1}\Lambda\)
for \(a<j\), by multiplication with \((sY)^{-1}\).
The same holds for \(Y'\) and \(q\). The vector
\(\alpha^2 e_a^*\) now has both transported vectors in
\(\alpha\Lambda\subset M\Lambda\), so
\(\alpha^2\in\mathfrak d_a(Y,x)\).

The denominator ideal \(\mathfrak J_x\) of the entries of
\((x-1)/M\) has generator \(c\), given in each field factor
by the uniformizer to the maximum of zero and the negative
valuations of those entries. Write \(D_x=[\mathcal R:\mathfrak J_x]\).
In the unit-ideal case, \(D_x=1\) is equivalent to
\((x-I)/M\) integral, hence to \(x\in J\). Indeed the integral
unimodular Gram also makes \(x^{-1}\) integral when \(x\) is integral.
At places where \(M\) is a unit this condition is precisely membership
in the lattice stabilizer. Thus the shell \(H=1\) contributes the
original level term.
If all depths after the head are less than \(\sigma^{-2}\),
where \(0<\sigma<1\), take
\(C=\alpha^2\prod_{a=j}^{\varrho-1}c_a\bar c_a\).
It belongs to every stable ideal and the induction gives
\[
 D_x\le \operatorname{md}(C)^{-2\varrho}
       \le \operatorname{md}(\alpha)^{-4\varrho}
                                      \sigma^{-8\varrho(\varrho-j)}.
\]
Let \(H=2^t\) be the least dyadic height at least \(D_x\),
and \(\Lambda_\alpha=2/\operatorname{md}(\alpha)\).
If \(H>\Lambda_\alpha^{2400 n^2}\), take
\(\sigma=H^{-1/(800n^2)}\). The displayed bound is then at
most \(H^{1/100+1/(600n)}<H/2\), contradicting \(H/2<D_x\).
Consequently, for each active \((Y,Yx)\) satisfying both head conditions
there exists \(a=a(Y,x)\), \(j\le a<\varrho\), with
\(D_a(Y,x)\ge\sigma^{-2}\). If \(j=\varrho\), the empty product in
\(C=\alpha^2\) gives the same contradiction, so the remote shell
contains no such active pair.

Each part of a shell with one denominator ideal \(c\mathcal R\)
lies in the position box of scale \(\bar c c\), since
\(\bar c c(x-1)/M\) is integral there. Its loss is
\(\Lambda'=2D_x^2\le\Lambda_\alpha H^2\).
The number of boxes is at most \(2H^{9n}\).
Indeed for \(a=3\varrho+5<9n\) the Euler product of finite-index
ideals of \(\mathcal R\) is bounded by
\(\zeta(a)^{[E:\mathbb Q]}\le2\) by the scale calculation of
Section~\ref{sec:compact-floor-haar}. Each ideal of index at
most \(H\) contributes at least \(H^{-a}\) to that sum.
This count applies at every projected place.

\subsection{Finite correlations on bounded sets}
\label{sec:finite-position-box}

\sectionopening{In this subsection, we bound weighted correlations on compact sets of relative positions.
The local vector count from Section~\ref{sec:compact-floor-volume}
supplies the estimate used for inverse weights in
Section~\ref{sec:finite-absolute-pivots}.}

From the \(k-1\) orthonormal coordinates of \(X_0\) outside the
distinguished column in Section~\ref{sec:period-finite-coordinates},
use \(k-1-2m_+=k-6n-7\) coordinates. Choose two disjoint blocks
of size \(2m_++2\) among them, with \(m_+=3n+3\).
Each block has rank \(2(3\varrho+5)+4\), since \(n=\varrho+1\).
Thus the quadratic norm-fiber estimate of
Section~\ref{sec:compact-floor-haar} gives the uniform mass
bound \(2\), including at residue characteristic two.
The volume exponent is the free rank minus \(3n+1+m_+\), hence
\[
 s=(k-6n-7)-2(2m_++2)-(3n+1+m_+)=k-24n-27 .
\]
Take \(E_a^K\) from the stage sums in
Section~\ref{sec:real-weighted-products}, using the finite
stage measures and \(W_a=\prod_{i<a}D_i^{-s}\) with the
depths of Section~\ref{sec:finite-stable-shells}. With the
finite inverse section moduli \(H_{p,T}\) and the two
kernels \(K_1,K_2\) in those sums, set
\[
 \mathcal W_j=H_{p,T}^2
   \Lambda_\alpha^{3\cdot10^6n^3}E_0^{K_1}E_j^{K_2}.
\]
For an integral, invertible, self-conjugate scalar \(\gamma\),
let \(\Omega_\gamma\) consist of the elements \(x\in\Omega_J\) with
\(\gamma(x-I)/M\) integral, and put
\(\Lambda'=2/\operatorname{md}(\gamma)\). The set \(\Omega_J\) and
its complement are unions of right \(J\)-cosets, so both are open.
The entry condition restricts \(x\) to a compact open matrix box.
Intersecting with the closed unitary group and with \(\Omega_J\)
gives a compact measurable \(\Omega_\gamma\). Use
\(K(g,y)=\langle W_Z(g),W_Z(d_0y)\rangle\), with \(d_0=1\) or
\(d_0=w_q^{-1}w_p\) and starting stage \(j_0=0\) or \(j_0=j\),
respectively. For the original quotient measure \(Q\),
we prove
\begin{equation}\label{eq:finite-position-energy}
 \int_{\Omega_\gamma} E_{j_0}^K\,d\mu
 \le n\delta_U(\Lambda')^{4900n\varrho}
                     \iint |K(g,y)|^2\,dQ(g)dQ(y).
\end{equation}

Call the two chosen blocks \(B_v,B_y\) and let \(X'\)
be the span of the remaining coordinates in this set. Their ranks are
\[
 h_0=2m_++2=2(3\varrho+5)+4=6n+8,
 \qquad d'=k-6n-7-2h_0=k-18n-23.
\]
All three are spanned by original orthonormal vectors of \(X_0\)
outside the distinguished column. They are integral orthogonal
summands of each remaining space \(V_j\), including at places
over \(N\). For \(g\in K_T\) and a relative position
\(x=l^{-1}ulj_x\) with \(u\in U\), \(j_x\in J\), choose a direct
summand \(L_0\subset X'_{\mathcal R}\) of rank \(D=d'-3\varrho\)
in the kernel of \((u-I)lg^{-1}\). Such a summand exists by
saturation in each local factor, since \(u-I\) sends
\(F_\varrho^\perp\) into \(F_\varrho\) and has rank at most
\(2\varrho\). For \(w\in L_0\) both \(g^{-1}w\) and
\((gx)^{-1}w=j_x^{-1}g^{-1}w\) are integral and congruent
modulo \(M\).

Let \(\phi(v)\) be the inverse volume of
\(h(v,L_0)\). A basis of \(L_0\) reads \(D\) pairing
coordinates in an integral invertible change of basis.
Counting the minimum valuation in each factor gives
\[
 \int_{cX'_{\mathcal R}}\phi(v)^s\,dv
 =\operatorname{md}(c)^{d'-s}
   \prod_{w\in T}\frac{1-(\mathrm Nw)^{-D}}
                      {1-(\mathrm Nw)^{-(D-s)}}
 \le 2\operatorname{md}(c)^{d'-s}.
\]
Indeed the \(D\) pairing coordinates have independent integral Haar,
up to conjugation of the field factors. At a factor of residue size
\(q_w\), the minimum valuation \(m\) has probability
\((1-q_w^{-D})q_w^{-mD}\) for \(m\ge0\). The zero tuple has measure
zero. Multiplication by \(c\) scales the pairing ideal by \(\bar c\);
\(\operatorname{md}(\bar c)=\operatorname{md}(c)\) on the selected
conjugation-stable product. Summing the geometric moments gives the
displayed formula, with \(D-s=3n+7\) and the product bounded as in
Section~\ref{sec:compact-floor-haar}.
Take \(g\in K_T\) and \(\sigma\in G_a\) on the actual stage
\(\theta_a\), with \(r\sigma g\in K_T\) for some \(r\in R_{a+1}\).
Then \(\sigma^{-1}r^{-1}e_a\) is integral. Since
\(r^{-1}e_a-e_a\in F_a\) and \(\sigma\) fixes the earlier pairs,
its integral projection away from those pairs is \(\sigma^{-1}e_a\).

Transporting \(L_0\)
by \(\sigma\) shows
\(D_a(\sigma g,x)\le\phi((\sigma^{-1}e_a)_{X'})\).
The norm-block estimate of Section~\ref{sec:compact-floor-volume}
bounds the conditional mass over \(X'\) by \(2\), so the
one-stage integral of \(D_a^s\) on the near locus is at most
\(4\). Rebasing \(g\) by \(R_a\) and using invariance of the
earlier depths gives the finite version of
\eqref{eq:real-depth-first}, with one-stage contribution \(4\).

Use the \(U_a\) coordinates \((v,\beta)\) of
Section~\ref{sec:compact-floor-haar}. The \(e_a^*\)-column
after subtracting the identity has vector part \(v\) and
\(e_a\)-coefficient \(z=\beta-\tau h(v,v)\).
The integral projection \(z\mapsto z-\tau(z+\bar z)\)
recovers \(\beta\), also at residue characteristic two.
The condition of lying in \(J\) is exactly divisibility of
\(v,\beta\) by \(M\). Thus set
\[
 \delta_a=\operatorname{md}(M)^{d_a+1/2},\qquad
 \prod_{a<\varrho}\delta_a=\delta_U,\qquad
 c_a=\delta_a(\Lambda')^{4900n}.
\]
For \(x_0,x_1\in\Omega_\gamma\), unitarity gives
\[
 M^{-1}(x_1x_0^{-1}-I)
    =M^{-1}(x_1-x_0)x_0^{-1}\in\gamma^{-2}M_r(\mathcal R).
\]
This persists after conjugation by a boxed element. When
the relative element belongs to \(U_{<a}U_a\), project its
\(e_a^*\)-column modulo \(F_a\) to read the \(U_a\) column.
The vector and skew parameters then belong to
\(M\gamma^{-2}(V_{a+1})_{\mathcal R}\) and
\(M\gamma^{-2}\mathcal R^-\), respectively.

We apply the forward-weight and square-ascent argument of
Section~\ref{sec:real-position-ball}. For the forward weight
take \(g\) boxed and \(x_0=g^{-1}y\in\Omega_\gamma\). If \(u_a\in U_a\)
has vector \(v\) and \(g^{-1}\nu u_a y\in\Omega_\gamma\) for some
\(\nu\in U_{<a}\), the column bound applies to \(u_a\).
For a stable pair \(\zeta,\zeta'\) at \(g,u_a y\), pair its two
transported lifts with the common vectors \(w\in L_0\) at
\(g,y\). The column formula for \(u_a w\) gives
\[
 h(\zeta,w)-h(u_a^{-1}\zeta',w)
 =\overline{h(e_a,\zeta)}h(v,w)\in M\mathcal R .
\]
Put \(\mathfrak b=h(v/M,L_0)\). Each stable coefficient \(c=h(e_a,\zeta)\)
satisfies \(\bar c\,\mathfrak b\subseteq\mathcal R\). If \(\mathfrak b\)
is nonzero in every field component, the coefficient ideal is thus
contained in \(\overline{\mathfrak b^{-1}}\). Conjugation preserves
volume on the conjugation-stable product, so its volume is at most
\(\operatorname{vol}(\mathfrak b)^{-1}\). This gives
\(D_a^{-s}\le\phi((v/M)_{X'})^s\), also with the infinite majorant
when a component of the pairing ideal vanishes.

Put \(v=M(v'+v'')\) with \(v'\in X'\) and \(v''\) in its orthogonal
complement in the \(d_a\)-dimensional vector part. Both lie in the
\(\gamma^{-2}\)-multiple of their integer boxes. The integral in \(v'\)
is at most \(2\operatorname{md}(\gamma)^{-2(d'-s)}\) by the content
integral above. The box in \(v''\) has volume
\(\operatorname{md}(\gamma)^{-2(d_a-d')}\). The box for \(\beta/M\)
has volume \(\operatorname{md}(\gamma)^{-1}\), since \(\gamma\) is
self-conjugate and the skew line has one fixed-field coordinate.
Restoring \(M\) therefore bounds the integral by
\[
 2\delta_a\operatorname{md}(\gamma)^{-2(d_a-s)-1}\le c_a .
\]
The forward induction of \eqref{eq:real-forward-weight}
applies: translate \(y\) into the box modulo \(U_{<a+1}\),
integrate the older unipotents first, and use invariance of
depth \(D_a\) under their action on the second point. That
action transports its stable lift without changing it
modulo \(F_a\). For a box move \(\sigma u\in G_aU_{<a}\),
replace \(y\) by \(\sigma y\) and \(\nu\) by
\(\sigma u\nu\sigma^{-1}\); the first \(a\) depths and the
position are preserved.

For the square ascent, use the cut by positive mass of
\(\{\nu\in U_{<a}:g^{-1}\nu y\in\Omega_\gamma\}\).
At stage \(a\ge j_0,\ a<\varrho\), recenter about a point \(u_{a,0}\)
of the nonzero support, with
\(g_c=tug\in K_T,\ u\in U_{<a},\ t=t_1n_1\in G_{a+1}U_a\).
For \(u_a\) on this support compare the positions
\(x_{u_a}=g^{-1}\nu_{u_a}u_a y\) and \(x_0=g^{-1}\nu_0 u_{a,0} y\)
in the box. The conjugate \(g_c x_{u_a}x_0^{-1}g_c^{-1}\)
belongs to \(U_{<a}\,t u_a u_{a,0}^{-1}t^{-1}\).
Hence \(u_a'=t u_a u_{a,0}^{-1}t^{-1}\) lies in the column box,
and the substitution preserves Haar. The vector of \(u_a\)
is \(t_1^{-1}v'+v_0\) for the vectors \(v',v_0\) of
\(u_a',u_{a,0}\).

Put \(y_\sigma=\sigma^{-1}f_{a+1}\). The simple-root sum of
\(\sigma u_a\sigma^{-1}\) is \(-h(v,y_\sigma)\). For \(a\ge j_0\),
\(d_0\) is the identity on \(V_a\) and commutes with this unipotent.
Since the inner product is linear in the second argument, the tail
correlation satisfies
\[
 C_{a+1}^K(g,u_a y)
   =\int_{G_{a+1}}\psi_T\!\left(-M^{-1}h(v,y_\sigma)\right)
      C_{a+2}^K(\sigma g,\sigma y)\,d\theta_{a+1}(\sigma),
 \qquad \psi_T(z)=\psi_{E,f}(-z/\delta_E).
\] The offset \(v_0\) gives a
unit phase in the quadric variable \(z=t_1\sigma^{-1}f_{a+1}\).
The quadric function is supported on integral \(z\):
nonzero correlation implies \(\sigma g\in\mathcal U_{a+2}\),
hence also \(\sigma t_1^{-1}g_c=\sigma n_1u g\), since the
extra conjugate is in \(U_{<a+1}\). Now use the integral
projected-column bound at \(g_c\).

Write \(z=z'+z''\) with \(z'\in X'\). On the integral quadric box
its measure is \(d\mu_{\rm quad}=dz'\,d\nu_{z'}(z'')\). The two
normalizations agree by division by the same norm-value Haar. The
complement of \(X'\) contains an orthonormal block of rank
\(2(3\varrho+5)+4\). Hence \(\nu_{z'}\) has mass at most \(2\) on
this box by Section~\ref{sec:compact-floor-volume}. For the quadric
function \(F(z)\), including the offset phase, fix the complementary
coordinate \(t''\) of \(v'/M\). Using the character \(\psi_T\) of
Section~\ref{sec:finite-frequency}, put
\[
 A_{t''}(z')=\int F(z'+z'')\,\psi_T(-h(t'',z''))\,d\nu_{z'}(z'').
\]
Cauchy--Schwarz gives \(\|A_{t''}\|_2^2\le2\int |F|^2\,d\mu_{\rm quad}\). Integrating the squared
Fourier transform of \(A_{t''}\) over all of \(X'\) gives precisely
\(\|A_{t''}\|_2^2\) by Plancherel for the integral self-dual Haar.
The box of \(t''\) has volume \(\operatorname{md}(\gamma)^{-2(d_a-d')}\).
Including the skew box and the change from \(v'/M\) to \(v'\),
the square-ascent contribution is
\[
 2\delta_a\operatorname{md}(\gamma)^{-2(d_a-d')-1}\le c_a.
\]
Both bounds by \(c_a\) follow from \(d_a-s\le r-s=26n+26\).
The case of zero cut mass follows by Tonelli. A simultaneous
move by \(G_{a+1}\) preserves the next cut, since it normalizes
the older unipotent with modulus one. At \(a=\varrho\),
Cauchy--Schwarz on the terminal stabilizer contributes at most one
by the compact coset bound.

For the summand of \(E_{j_0}^K\) indexed by \(l\ge j_0\),
the substitution \(x=g^{-1}uy\), \(u\in U_{<l}\), changes the
second Bessel factor by the scalar character of
\(d_0ud_0^{-1}\). Its value is independent of the higher pivots,
because \(G_{l+1}\) fixes the two flag vectors in each simple-root
entry of this conjugate. Thus the modulus of the tail correlation
is constant under the substitution. The forward integral uses
every stage \(a<l\), including those with \(a<j_0\), and depends
only on the close-position cut and the weights \(D_a^{-s}\).
The square ascent uses \(l\le a<\varrho\), where \(d_0\)
centralizes \(U_a\). The two products give
\[
         \prod_{a<l}c_a\prod_{l\le a<\varrho}c_a
                   =\delta_U(\Lambda')^{4900n\varrho}.
\]
The terminal stage has mass at most one. Applying the product
Haar law as in \eqref{eq:real-position-energy} and summing
the at most \(n\) terms proves \eqref{eq:finite-position-energy}.

For the twisted kernel, conjugation by \(d_0\) on \(U\) has
modulus \(J_{d_0}=\prod_{a<j}\operatorname{md}(p_a/q_a)^{d_a+1}\).
Indeed \(d_0 e_a=(p_a/q_a)e_a\). On the Heisenberg vector the
remaining unitary factor has determinant of modulus one, and the scalar
has modulus \(\operatorname{md}(p_a/q_a)\) per coordinate. On the skew
centre the scalar is \((p_a/q_a)\overline{(p_a/q_a)}\), with modulus
\(\operatorname{md}(p_a/q_a)\) on that fixed-field line.
Thus \(\int F(d_0y)dQ(y)=J_{d_0}\int F(y)dQ(y)\), since changing
the unipotent variable in the Haar disintegration under
\(y\mapsto d_0y\) contributes \(J_{d_0}^{-1}\).
The head-density identity is
\[
 (\mathcal D_j^T)^2H_{p,T}^2J_{d_0}
 =\prod_{a<j}\operatorname{md}(q_a)^{-1}
 \le\Lambda_\alpha^2 ,
\]
where the last inequality follows from \(\alpha/\prod_{a<j}q_a\)
being integral.

For a measurable \(\Omega\subseteq\Omega_\gamma\) and a
measurable \(f_0\ge0\) with \(f_0^2\le e_*^2\mathcal W_j\),
\(e_*\ge0\),
the position-box estimate is
\[
 \mathcal D_j^T(p,q)\int_{\Omega} f_0\,d\mu
 \le e_*\,\Lambda_\alpha^{1500000 n^3+1}
                  (\Lambda')^{25000 n^2}\mathcal R_Z
\]
by Cauchy--Schwarz, \eqref{eq:finite-position-energy},
\(n(\Lambda')^{4900n\varrho}\le(\Lambda')^{25000n^2}\),
and the quotient energy bound of
Section~\ref{sec:finite-bottom-normalization}.

\subsection{Reciprocal weights and absolute pivot bounds}
\label{sec:finite-absolute-pivots}

\sectionopening{In this subsection, we obtain absolute integral estimates for both kernels before imposing
an oscillatory weight. Section~\ref{sec:finite-position-box} supplies
their local weights, which also enter the Fourier argument of
Section~\ref{sec:finite-frequency}.}

Use the finite depths in \(I_a=\mathbf1_{\mathcal U_a}
W_a^{-1}\) and in the sums \(B_a^K\) of
Section~\ref{sec:real-weighted-products}. Fix \(x\in\Omega_J\) and write \(C_\alpha=\Lambda_\alpha^{10^6n^3}\) and
\(H_{p,T}=\prod_{a<j}\operatorname{md}(p_a)^{-d_a-1/2}\).
We prove
\begin{equation}\label{eq:finite-cross-weight}
 \int_{U_{<j}}I_j(w_p\nu g)\,d\nu\le C_\alpha H_{p,T}
                      \qquad(g\in\mathcal U_j).
\end{equation}
Here \(\alpha/P_a\) is integral for every initial segment
\(P_a=\prod_{b<a}p_b\). To rebase \(g\) into \(K_T\)
by \(\gamma u_0\in G_jU_{<j}\), conjugate the older
unipotents by \(\gamma\) and use right translation by
\(u_0\). These preserve Haar, and \(\gamma\) commutes with
\(w_p\) and preserves \(I_j\).

For \(g\in K_T\), take \(u\in U_a,\ \nu\in U_{<a}\) with
\(w_{p,<a+1}u\nu g\in\mathcal U_{a+1}\). Set \(s_0=(u\nu)^{-1}\)
and choose \(t\in R_{a+1}\) so that \(t w_{p,<a+1}u\nu g\)
is integral. The ratio of \(g\) and this matrix is
\(s_0 w_{\bar p,<a+1}t^{-1}\in K_T\).
Apply the row-replacement formula for minors in
Section~\ref{sec:finite-stable-shells} to get
\(\bar P_{a+1}s_0e_a^*\in\Lambda\).
Divide by \(\bar P_a\) and project away the earlier pairs,
which removes \(\nu^{-1}\) modulo \(F_a\).
Thus \(cu^{-1}e_a^*\) is in the \(\alpha^{-1}\)-box on
the \(a\)-th and later pairs, where \(c=\bar p_a\).

On the common summand \(L_0\) of
Section~\ref{sec:finite-position-box}, impose the \(a\)
equations \(h(e_b^*,\nu w)=0\) for \(b<a\). Their saturated
kernel contains a primitive submodule of rank \(d'-3\varrho-a\).
The equations ensure that \(w_{p,<a+1}u\nu w\) is
orthogonal to \(F_a\). That vector gives equal stable lifts
with coefficient \(h(cu^{-1}e_a^*,w)\), as the lower
unipotent acts trivially modulo \(F_a\) on \(F_a^\perp\).
Hence \(D_a\le\phi(v_{X'})\) for \(v\) the \(V_{a+1}\)-part
of \(cu^{-1}e_a^*\) and \(\phi\) defined using this smaller
submodule. For \(D=d'-3\varrho-a\) we have \(D-s\ge2n+8\).
With \(n[E:\mathbb Q]\le2^{2\varrho}\) from
Section~\ref{sec:compact-floor-haar},
\[
 \log\zeta_{\mathcal R}(D-s)
 \le [E:\mathbb Q]\,2^{1-(D-s)}
 \le 2^{-9}/n .
\]
The scaled content integral is therefore at most
\(2\operatorname{md}(\alpha)^{-(d'-s)}\).

In the Heisenberg coordinates \((\omega,\beta)\) of
\(u^{-1}\), two allowed centers at fixed \(\omega\) satisfy
\(c(\beta-\beta')\in\alpha^{-1}\mathcal R\).
In a nonsplit factor \(c\) is a unit times a power of
a fixed-field uniformizer. In a split pair, the skew
coordinate is \((z,-z)\), so the smaller of the two
volume constraints is at most their geometric mean.
The skew fiber measure is thus bounded by
\(\operatorname{md}(c\alpha)^{-1/2}\).

The vector change \(v=c\omega\) contributes
\(\operatorname{md}(c)^{-d_a}\) in the integral. Together
with the scaled content and the complementary vector box,
the integral of \(D_a^s\) on the indicated neighborhood is bounded by
\[
 2\operatorname{md}(p_a)^{-d_a-1/2}
             \operatorname{md}(\alpha)^{-(d_a-s+1/2)}.
\]
The earlier depths and near condition are unchanged
by the added moves in \(G_a\).
Indeed, write \(w_{p,a}\) for the Weyl element on the \(a\)-th
pair alone. For \(u\in U_a\) and \(\nu\in U_{<a}\),
\[
 w_{p,<a+1}u\nu=(w_{p,a}u)w_{p,<a}\nu,\qquad w_{p,a}u\in G_a.
\]
This factor preserves the first \(a\) depths. Since
\(\mathcal U_{a+1}\subseteq\mathcal U_a\), the earlier inverse
weights factor out as \(I_a(w_{p,<a}\nu g)\) on the support.
The primitive submodule used in the integral over \(U_a\) may
depend on \(g,x,\nu\). Its rank and the integral bound are uniform
in these variables. Iteration in the Haar order \(u\nu\) gives
\eqref{eq:finite-cross-weight}, since \(a<j\le\varrho<n\) and
\(d_a-s\le26n+26\).

For any nonnegative measurable left-\(R_j\)-invariant function \(B\)
supported on \(\mathcal U_j\) and any \(\sigma\)-finite measure \(\rho\)
on \(G_T\), Tonelli and \(B(\nu g)=B(g)\) for \(\nu\in U_{<j}\)
give
\[
 \int B(g)I_j(w_p g)\,d(\eta_{<j}*\rho)(g)
 \le C_\alpha H_{p,T}\int B(g)\,d\rho(g).
\]
The same estimate holds with \(w_p^{-1}=w_{\bar p}\), since
\(\alpha=\bar\alpha\) and conjugation preserves the modulus.

On \(\mathcal U_{j+1}\) the stage \(j\) mass is at
most \(V_0=\zeta_{\mathcal R}(3\varrho+5)\le2\), by
Section~\ref{sec:compact-floor-haar}. Apply the ordinary
Cauchy--Schwarz bound and the first-deep-stage argument of
Section~\ref{sec:real-depth-majorant}, using the finite
one-stage contribution \(4\) in Section~\ref{sec:finite-position-box}.
Explicitly, let \(P_j^K\) be the integral of
\(|C_{j+1}^K(\sigma g,\sigma gx)|\) over \(\theta_j\).
Then \((P_j^K)^2\le V_0 B_j^K I_j\).
If every nonzero terminal contribution has some stage
\(a\ge j\) with \(D_a\ge e^{-1}\), \(0<e\le1\), stop the
tail expansion at the first such \(a\). The no-stop branch
vanishes, since higher groups preserve each earlier depth.
Dropping the restrictions at earlier stages gives
\(P_j^K(g)\le\sum_{a=j}^{\varrho-1}T_a(g)\), where
\[
 T_a(g)=\int\mathbf1_{\{D_a(tg,x)\ge e^{-1}\}}
           |C_{a+1}^K(tg,tgx)|\,d\theta_{[j,a]}(t).
\]
On this support \(W_a^{-1}=D_a^{-s}W_{a+1}^{-1}
\le e^s W_{a+1}^{-1}\). Cauchy--Schwarz gives
\[
 \begin{aligned}
 T_a(g)^2
 &\le e^s\left(\int I_{a+1}(tg)\,d\theta_{[j,a]}(t)\right)
         \left(\int A_a^K(tg)\,d\theta_{[j,a]}(t)\right)\\
 &\le e^s4^{a+1-j}I_j(g)B_j^K(g).
 \end{aligned}
\]
The first integral uses the one-stage estimate at each of the
\(a+1-j\) stages. The second is one summand of \(B_j^K(g)\).
Now \(\varrho^24^\varrho\le C_\alpha\) gives
\((P_j^K)^2\le e^s C_\alpha B_j^K I_j\).

Write \(z=w_p^{-1}g\). In the two-pivot majorant
\(\mathcal M_j(x)=\int P_j^{K_1}(g)P_j^{K_2}(z)\,dM_j(g)\)
the product of the two weights and energies can be crossed
as
\[
 [B_j^{K_1}(g)I_j(z)]\,[B_j^{K_2}(z)I_j(g)] .
\]
Cauchy--Schwarz and \eqref{eq:finite-cross-weight} bound
the first integral by \(H_{p,T}C_\alpha E_0^{K_1}\).
For the second change \(g=w_p z\).
The integrand in \(z\) is left \(G_j\)-invariant,
and \(d\mu=d\mu_j*M_j\) with \(\mu_j\) Haar on \(G_j\).
Since \(w_p\) commutes with \(G_j\), the invariant-function
quotient identity in Section~\ref{sec:real-weighted-products}
applies with the finite measures.

The second integral is
then at most \(H_{p,T}C_\alpha E_j^{K_2}\).
It follows that \(\mathcal M_j^2\le V_0^2\mathcal W_j\).
On a shell \(H>\Lambda_\alpha^{2400n^2}\), take
\(\sigma=H^{-1/(800n^2)}\). If \(P_j^{K_2}(z)>0\), its
iterated integral contains a nonzero terminal kernel. All later
groups lie in \(G_j\), and the twist \(w_q^{-1}w_p\) commutes
with their product. Thus for the same \(t\in G_j\),
\[
 W_Z(tw_p^{-1}g)\ne0,\qquad W_Z(tw_q^{-1}gx)\ne0.
\]
Both arguments lie in \(UK_T\subseteq\mathcal U_j\).
For every \(h\in G_j\), left translation by \(ht^{-1}\in R_j\)
and commutation of the two Weyl elements with \(h\) give
\[
 w_p^{-1}hg,\quad w_q^{-1}hgx\in\mathcal U_j .
\]
Whenever \(K_1(hg,hgx)\ne0\), Section~\ref{sec:finite-stable-shells}
therefore supplies an index \(j\le a(h)<\varrho\) with
\(D_{a(h)}(hg,x)\ge\sigma^{-2}\).
The first stopping index may depend on \(h\). For
\(j=\varrho\) this implication rules out every such terminal
contribution. If the second pivot vanishes its product with the
first is zero. Otherwise the estimate for the first pivot with
\(e=\sigma^2\) and the ordinary estimate for the second give
\(\mathcal M_j^2\le V_0^2\sigma^{2s}\mathcal W_j\).

At \(j=0\) both pivots coincide and \(M_0\) is the atom
at the identity. If every active point has some
\(D_a\ge e^{-1}\), using the deep bound for both identical
factors gives \(\mathcal M_0^2\le e^{2s}\mathcal W_0\);
the factor \(C_\alpha^2\) from the two bounds is at most
the \(C_\alpha^3\) in \(\mathcal W_0\).

\subsection{Local Fourier estimates for the leading coordinates}
\label{sec:finite-frequency}

\sectionopening{In this subsection, we use translation invariance of the finite kernels to identify the
frequencies that can contribute. The absolute integral estimates of
Section~\ref{sec:finite-absolute-pivots} control their magnitudes in
the summation of Section~\ref{sec:finite-shell-reduction}.}

Fix the finite data of Section~\ref{sec:finite-common-entry}.
Use \(d'=k-18n-23\), the dimension of \(X'\), and
\(C_\alpha=\Lambda_\alpha^{10^6n^3}\).
Write “box \(a\)” for vectors or matrices with integral
coordinates after multiplication by \(\alpha^a\).

Suppose \(g,w_p^{-1}g\in\mathcal U_{j}\).
Factor a move at \(g\) into the middle group and older
unipotent, and choose \(t\in G_j,\ u\in U_{<j}\) such
that \(utg\in K_T\). Since \(t\) commutes with \(w_p\),
choose \(v\in R_j\) with \(vw_p^{-1}tg\in K_T\).
Their ratio is \(u w_pv^{-1}\in K_T\).
The minors in Section~\ref{sec:finite-stable-shells} give
\(P_a\in\mathcal R\) and \(P_{a+1}u e_a^*\in\Lambda\)
for \(a<j\), with the leading products \(P_a=\prod_{b<a}p_b\).

Use the integral quotients \(\alpha/P_a\) to conclude
\(p_a^{\pm1}\in\alpha^{-1}\mathcal R\) and
\(u e_a^*\in\alpha^{-1}\Lambda\).
In Witt order these \(j\) dual columns of \(u\) have
upper, middle and lower blocks \(C,D,E\), with \(E\)
lower unitriangular.
Its upper diagonal block is \(A=(E^*)^{-1}\) and its
upper middle block is \(-A D^*J_0\), with \(J_0\) the
integral middle Gram. The nilpotent inverse series
bounds \(A\) in box \(j-1\), hence \(u,u^{-1}\) in
box \(j\).

Thus \(tg\) is in box \(j\), and
\(tw_p^{-1}g\) in box \(j+1\), also when \(j=0\) by
taking \(t=g^{-1}\).
For \(\sigma\in G_j\) on the support of a pivot tail
at \(g\), choose \(\zeta\) with
\(\zeta\sigma g\in K_T\), \(\zeta\in R_{j+1}\).
Projecting \((tg)(\zeta\sigma g)^{-1}f_j\) off the
earlier pairs gives \(t\sigma^{-1}f_j\), since
\(\zeta^{-1}f_j=f_j\bmod F_j\).
The orbit vector is therefore in box \(j\), and the
corresponding vector from the pivot at \(w_p^{-1}g\)
is in box \(j+1\). Here \(f_\varrho=z\), so the final
pivot also uses its norm quadric, with the stabilizer
integration occurring after that pivot.

Use the finite character \(\psi_T(x)=\psi_{E,f}(-x/\delta_E)\) on
\(A_T\). In each field component the trace-dual of the integers is
the inverse different. Thus for \(\psi_T\) the annihilator of the
integer lattice is that same lattice, and its Haar measure of mass one
is self-dual, including above \(2\). The level factor \(M^{-1}\) enters
through the Bessel character. On two norm
quadrics in \(V_j\), test on boxes \(a,a'\), and put
\(c_j=\dim V_j-d'\). The operator with kernel
\(\psi(\eta h(v,y))\), \(\eta\in A_T^\times\), has
squared \(L^2\)-operator norm on these boxes at most
\[
 4\operatorname{md}(\eta)^{-d'}
               \operatorname{md}(\alpha)^{-(a+a')c_j}.
\]
Use the three orthogonal summands \(X',B_v,B_y\) fixed in
Section~\ref{sec:finite-position-box}. Impose the norm of \(v\)
on \(B_v\) and the norm of \(y\) on \(B_y\). In a box of exponent \(a'\)
for \(y\), a block of rank \(h_0\) has conditional volume
\[
 \operatorname{vol}(B_y\mid \text{fixed norm})
            \le 2\operatorname{md}(\alpha)^{-a'(h_0-1)}.
\]
This follows from Section~\ref{sec:compact-floor-volume}
by scaling the norm-value Haar by
\(\operatorname{md}(\alpha)^{-a'}\).
The other \(c_j-h_0\) nonfree coordinates contribute their box
volumes, giving the bound
\(2\operatorname{md}(\alpha)^{-a'(c_j-1)}\).
The same bound with \(a\) applies to the nonfree coordinates
of \(v\).

Write the two variables, in the orthogonal order
\(X',B_v,B_y,(X'\perp B_v\perp B_y)^\perp\), as
\(v=(x,b,c,z)\) and \(y=(t,b',c',z')\).
Fix \(c,z,b',z'\). The norm of \(b\) is determined by \(x,c,z\),
and that of \(c'\) by \(t,b',z'\). For a bounded test \(F\) on
the second box put
\[
  H(t)=\int F(t,b',c',z')\psi(\eta h(c,c'))\,d\nu_t(c'),
\]
where \(d\nu_t\) is its indicated block fiber measure. Its mass
bound and Cauchy--Schwarz bound \(|H(t)|^2\) by that mass times
the corresponding integral of \(|F|^2\). The remaining phase in
the \(t\)-integral is \(\psi(\eta h(x,t))\) times a scalar of
modulus one independent of \(t\). Its absolute square is therefore
independent of \(b\). First apply Cauchy--Schwarz in \(b',z'\),
then use the mass bound on the permitted \(b\)'s at each \(x\).
We may now extend the \(x\)-integration to all of \(X'\).
Plancherel gives \(\operatorname{md}(\eta)^{-d'}\int|H(t)|^2\,dt\).
Indeed in orthonormal coordinates the change of Fourier variable
is \(x\mapsto\eta\bar x\). Conjugation preserves integer Haar:
it interchanges the two coordinates of a split pair and acts by
an integral automorphism at a nonsplit place. Multiplication by
\(\eta\) has Jacobian \(\prod_{w\in T}|\eta_w|_w^{d'}
=\operatorname{md}(\eta)^{d'}\), even for unequal valuations
at the two members of a split pair. The remaining Fourier transform
is isometric for the self-dual Haar just computed.
Integrating in \(c,z,b',z'\) gives the product of the two block
mass bounds and the two complementary box volumes, namely
\(4\operatorname{md}(\alpha)^{-(a+a')(c_j-1)}\), besides the
Plancherel factor and \(\|F\|_2^2\). Compact boxes and the uniform
fiber bounds justify Fubini for bounded tests. Truncation in
\(L^2\) gives the general estimate. Since \(\alpha\) is integral,
we can enlarge the exponent \(c_j-1\) to \(c_j\). If the box for
\(v\) has no norm equation, only the block \(B_y\) contributes
a conditional volume and the constant is \(2\).
All scalings use the actual finite product, including
the split and unramified quadratic factors.

For the common-entry character \(\omega=\psi(\xi Z_j)\),
the two pivot variables give the character of
\[
 Z_j(\sigma\tau^{-1})
      =h(\sigma^{-1}f_j,\tau^{-1}f_j).
\]
The pointwise double pivot is zero unless both starting
points \(g,w_p^{-1}g\in\mathcal U_{j}\).
Rotate the two orbit variables by the same \(t\).
The Hermitian pairing is unchanged, and its values
lie in \(\alpha^{-(2j+1)}\mathcal R\).
It suffices that a frequency \(\eta\in A_T^\times\)
agree with the original character on
\(\alpha^{-200n^2}\mathcal R\).
The amplitudes factor through the pivot orbits by the
stabilizer invariance of the tails; rotating the orbit
variable preserves the measure.
The two-box estimate and Cauchy--Schwarz bound the
square of this double pivot by
\[
 \operatorname{md}(\eta)^{-d'} C_\alpha\,
   B_j^{K_1}(g)I_j(g)\,
   B_j^{K_2}(w_p^{-1}g)I_j(w_p^{-1}g).
\]
Indeed the squared pivot norms are bounded by \(B_j^K I_j\),
and \(c_j\le r-d'=20n+22\) makes the box loss at most
\(C_\alpha\). Cross the weights and apply the two
inverse-weight integrations as in
Section~\ref{sec:finite-absolute-pivots}. For the
integral \(f_0(x)\) of the modulus of the double pivot
over \(M_j\), this gives
\begin{equation}\label{eq:finite-frequency-point}
 f_0(x)^2\le\operatorname{md}(\eta)^{-d'}\mathcal W_j(x).
\end{equation}

For the unequal heads take \(0\le i<\varrho\), \(j=i+1\), \(\omega=1\) and
\(\xi=M^{-1}(p_i^{-1}-q_i^{-1})\in A_T^\times\).
Let \(d=w_q^{-1}w_p\). Write \(\mathfrak s\) for the
simple-root sum in the Bessel character. A
stage-\(i\) element \(u\) with vector coordinate \(v(u)\)
has \(\mathfrak s(u)=-h(v(u),f_j)\). Conjugation by
\(d\) multiplies this root by \(p_i/q_i\), because
\(d\) acts by that scalar on \(e_i\) and fixes \(V_j\).
On moving \(u\) through the next pivot \(\sigma\),
the difference of the two character values is thus
the character of
\[
 -M^{-1}(p_i/q_i-1)h(v(u),\sigma^{-1}f_j)
   =-\xi h(\bar p_i u^{-1}e_i^*,\sigma^{-1}f_j).
\]
The equality uses conjugate-linearity in the first argument and
the \(V_j\) component \(-v(u)\) of \(u^{-1}e_i^*\). Thus, with the
same orbit measure also at \(j=\varrho\), where \(f_j=z\),
\[
 C_j^{K_2}(ug,ugx)=\int_{G_j}
 \psi_T\!\left(-\xi h(\bar p_i u^{-1}e_i^*,\sigma^{-1}f_j)\right)
 C_{j+1}^{K_2}(\sigma g,\sigma gx)\,d\theta_j(\sigma).
\]
Indeed \(\sigma u\sigma^{-1}\in U_i\) and its two actions on the
kernel give the displayed difference of the simple-root characters.

The resulting section estimate for the twisted tail is
\begin{equation}\label{eq:finite-unequal-section}
 \begin{split}
 &\int_{U_i}\mathbf1_{\{w_pug\in\mathcal U_{j}\}}
                  |C_j^{K_2}(ug,ugx)|^2\,du\\
 &\quad\le\operatorname{md}(\xi)^{-d'}
       \operatorname{md}(p_i)^{-d_i-1/2}C_\alpha
       \int |C_{j+1}^{K_2}(\sigma g,\sigma gx)|^2\,d\theta_j(\sigma).
 \end{split}
\end{equation}
To prove it, at each \(u\) in the support both \(ug\)
and \(w_pug\in\mathcal U_{j}\).
Apply the simultaneous box bound with \(w_pug\) as
the first argument, obtaining \(t_u\in G_j\) with
\(t_uug\) in box \(a_0=j+1\).
Fix one such point \(u_*\) and its \(t_*\).
The substitution \(m=t_*u_*u^{-1}t_*^{-1}\)
preserves \(U_i\) Haar. Since \(t_u,t_*\) fix
\(e_i^*\),
\[
 (m-I)e_i^*
 =\bigl((t_*u_*g)(t_uug)^{-1}-I\bigr)e_i^* .
\]
The right-hand side lies in box \(2a_0\), using the
integral inverse Gram and \(\bar\alpha=\alpha\).
The trace-one skew projection keeps the skew parameter
in that box, also above \(2\). Its Haar volume is
\(\operatorname{md}(\alpha)^{-a_0}\).
Use \(\bar p_i v(m)\) as the vector Fourier variable.
Because \(p_i\in\alpha^{-1}\mathcal R\), it lies in
box \(2a_0+1\), and the vector measure is multiplied by the factor
\(\operatorname{md}(p_i)^{-d_i}\) under this substitution.

Rotate the quadric in the second variable by \(t_*\).
The relation \(u^{-1}=u_*^{-1}t_*^{-1}mt_*\)
makes the \(u_*\) term of the phase a unit factor
depending only on the second variable. The variable term is the
pairing of the new vector with \(t_*\sigma^{-1}f_j\)
at frequency \(-\xi\). To bound the second vector, put
\(B=t_*u_*g\) and \(\kappa=\sigma t_*^{-1}\). Then
\(\kappa B=\sigma u_*g\in\mathcal U_{j+1}\) since
\(\sigma g\in\mathcal U_{j+1}\) on the tail support and
\(\sigma\in G_j\) normalizes \(U_{<j}\). Choose
\(\zeta\in R_{j+1}\) with \(\zeta\kappa B\in K_T\).
The column \(B(\zeta\kappa B)^{-1}f_j\) lies in box \(a_0\).
Its projection off the earlier pairs is
\(\kappa^{-1}f_j=t_*\sigma^{-1}f_j\), since
\(\zeta^{-1}f_j\equiv f_j\bmod F_j\).
Write \(c=d_i-d'\). In the squared one-quadric bound the vector
box has exponent \(2a_0+1\) and the quadric box exponent \(a_0\).
Including the skew coordinate and the vector Jacobian gives
\[
 2\operatorname{md}(\xi)^{-d'}\operatorname{md}(p_i)^{-d_i}
       \operatorname{md}(\alpha)^{-(3a_0+1)c},
 \qquad (2a_0+1)c+a_0(c-1)+a_0=(3a_0+1)c .
\]
Since \(\alpha p_i\) is integral, \(a_0\le n\) and \(c\le20n+22\),
the inequality \((3n+1)(20n+22)+1/2<10^6n^3\) absorbs this in
\(\operatorname{md}(\xi)^{-d'}\operatorname{md}(p_i)^{-d_i-1/2}C_\alpha\).
This proves \eqref{eq:finite-unequal-section}.

For the outer integration put \(z=w_p^{-1}g\) and
\(\varphi_l(g)=|C_j^{K_l}(g,gx)|\).
The squared product \(\varphi_1(g)^2\varphi_2(z)^2\)
is the product of the crossed factors
\[
 A_i^{K_1}(g)I_j(z),\qquad
 \mathbf1_{\{w_pz\in\mathcal U_{j}\}}\,
 W_j(z)I_i(w_{p,<i}z)\varphi_2(z)^2 .
\]
This identity is on the support of the product, where
\(I_i(w_{p,<i}z)=W_i(g)^{-1}\), since the last Weyl
factor belongs to \(G_i\). The first factor integrates
to at most \(H_{p,T}C_\alpha E_0^{K_1}\):
use \eqref{eq:finite-cross-weight} and then the
\(i\)-th energy term. For the second factor change
\(g=w_p z\) in the \(G_j\)-invariant quotient.
In \(M_j=\eta_i*\eta_{<i}*\theta_{[0,i]}\),
the weights \(W_j(z)\) and \(I_i(w_{p,<i}z)\) are
invariant under \(U_i\). Apply
\eqref{eq:finite-unequal-section}. The next pivot
energy bound \(B_j^{K_2}/W_j\) cancels \(W_j\).
Integrating the remaining inverse weight by
\eqref{eq:finite-cross-weight} at length \(i\) gives
\(C_\alpha E_j^{K_2}\) times the first \(i\) section
moduli. The section in \eqref{eq:finite-unequal-section}
supplies the last modulus. Thus the second integral
is at most
\(\operatorname{md}(\xi)^{-d'}H_{p,T}C_\alpha^2E_j^{K_2}\).
By Cauchy--Schwarz this proves
\eqref{eq:finite-frequency-point} with \(\eta=\xi\).

In either case \(f_0\le\mathcal M_j\). Since \(d'\ge s\),
\eqref{eq:finite-frequency-point} and the absolute
bound give the required squared saving
\(V_0^2\min(1,\operatorname{md}(\eta)^{-1/2})^{2s}\mathcal W_j\).

\subsection{Summation by denominator}
\label{sec:finite-shell-reduction}

\sectionopening{In this subsection, we sum the finite integral estimates by ideal index. The smaller bound
on large denominators comes from
Section~\ref{sec:finite-absolute-pivots}, and the sum gives the
estimate at the selected stabilizers in
Section~\ref{sec:finite-selected-bottom}.}

Take \(0\le f_0\le\mathcal M_j\) with the pointwise saving
\(f_0^2\le (V_0 e_0^s)^2\mathcal W_j\), \(e_0\in[0,1]\).
The absolute pivot estimates give this bound for
\(e_0=1\), and also with \(e_0\) replaced by \(\sigma\)
on the remote shells, by Section~\ref{sec:finite-absolute-pivots}.
The box cover then bounds the shell integral
times its head density by \(V_0g_H^s\Lambda_\alpha^{2\cdot10^6n^3}
H^{10^5n^2}\mathcal R_Z\), where \(g_H=e_0\) for
\(H\le\Lambda_\alpha^{2400n^2}\), and
\(g_H=\min(e_0,H^{-1/(800n^2)})\) otherwise.
For \(s\ge l+9\cdot10^7n^4\), split
\(g_H^s\le e_0^l g_H^{9\cdot10^7n^4}\).
In the first range,
\(H^{1+10^5n^2}g_H^{9\cdot10^7n^4}
\le\Lambda_\alpha^{29\cdot10^7n^4}\). In the second range
it is at most \(H^{1-12500n^2}\le1\).
Summing \(H^{-1}\) gives \(2\), so
\[
 \mathcal D_j^T(p,q)\int f_0\,d\mu
                  \le e_0^l\Lambda_\alpha^{10^9n^4}\mathcal R_Z .
\]
The absolute pivot bound allows \(e_0=1\) and gives the
integrable majorant used in the lower common-entry identity
of Section~\ref{sec:finite-common-entry}.

For the Fourier saving let \(\xi\in A_T^\times\) and
\(D(\xi)=[\mathcal R:\mathcal R\cap\xi^{-1}\mathcal R]\).
At a split pair let \(q_w=q_{\bar w}\) be the residue size and put
\(a=v_w(\xi_w)\), \(a'=v_{\bar w}(\xi_{\bar w})\). Its factors are
\[
 \operatorname{md}_{\{w,\bar w\}}(\xi)=q_w^{-a-a'},\qquad
 D_{\{w,\bar w\}}(\xi)=q_w^{\max(0,-a)+\max(0,-a')}.
\]
Indeed the denominator ideal in each field component has exponent
\(\max(0,-v_w(\xi_w))\), including when \(a,a'\) have opposite signs.
For the character \(\omega=\psi_T(\xi Z_j)\), take arbitrary proper
heads \(p,q\) of the same length. By the common rotation in
Section~\ref{sec:finite-frequency}, the pairing of the two pivot
vectors lies in \(\alpha^{-(2j+1)}\mathcal R\subseteq
\beta^{-1}\mathcal R\), where \(\beta=\alpha^{200n^2}\).
Take \(\eta_w=\xi_w\) if \(v_w(\xi)\le v_w(\beta)\) and
\(\eta_w=\beta_w\) otherwise. Then
\(\eta-\xi\in\beta\mathcal R\), so the characters agree
on the Gram support. Write \(b_w=v_w(\beta)\ge0\) and
\(a_w=v_w(\xi_w)\). For \(w_0=\beta/\eta\) the common-entry choice gives
\[
 v_w(w_0)=\max(0,b_w-a_w),\qquad
 v_w(\xi w_0)=\max(a_w,b_w)\ge0 .
\]
For the unweighted form with unequal last proper entries keep \(\eta=\xi\)
and take \(w_0=\beta\), with \(\beta\xi\) integral. In either case
\(v_w(w_0)\ge\max(0,-a_w)\); hence \(w_0\mathcal R\) is contained
in the denominator ideal and \(\operatorname{md}(w_0)\le D(\xi)^{-1}\).
Put \(e_0=\min(1,\operatorname{md}(\eta)^{-1/2})\).
For the common entry \(v_w(\eta)=\min(a_w,b_w)\), so
\(\operatorname{md}(\eta)\ge\operatorname{md}(\xi)\) and
\(\operatorname{md}(\eta)^{-1}=\operatorname{md}(w_0)/\operatorname{md}(\beta)\).
For unequal entries \(\eta=\xi\) and this last ratio is one.
Thus \(e_0\le\operatorname{md}(\xi)^{-1/2}\) and
\[
 e_0^2\le\frac{\operatorname{md}(w_0)}{\operatorname{md}(\beta)}
                  \le D(\xi)^{-1}\Lambda_\alpha^{200n^2}.
\]
For \(l=k-n^7+12\) we have
\(s-l=n^7-24n-39\ge 9\cdot10^7n^4\), since \(n\ge10^8\)
(Section~\ref{sec:real-marked-filling}). Thus the shell saving retains
\(\operatorname{md}(\xi)^{-(k-n^7)/2}D(\xi)^{-6}\)
at contribution \(\Lambda_\alpha^{1200n^2}\). This calculation
requires a frequency nonzero in every factor, as holds
for a nonzero rational frequency or a nonzero difference
of two rational inverse heads.

For the normalized bound, a nonzero majorant requires
all \(P_a,Q_a\) integral, since it requires both pairs
of near conditions at a plain point and its twisted image. A nonintegral initial segment makes the form vanish.
In the remaining case take
\[
 B=\prod_{a=1}^jP_aQ_a,\qquad \alpha=M B\bar B,\qquad
 \mathfrak h(p)=\prod_{a=1}^j\operatorname{md}(P_a)^{-n^6/2},
 \qquad
 \mathfrak h(q)=\prod_{a=1}^j\operatorname{md}(Q_a)^{-n^6/2}.
\]
Then \(\alpha/M,\alpha/P_a,\alpha/Q_a\) are integral,
as are \(\alpha/p_i,\alpha/q_i\). Thus \(\alpha^2\xi\)
is integral for the unequal-head frequency, and
\(\beta=\alpha^{200n^2}\) is admissible in both cases.
The head weights are defined to be zero in the
nonintegral-initial segment case. For \(4e_1\le n^6\), one has
\[
 \Lambda_\alpha^{e_1}
 = (2/\operatorname{md}(M))^{e_1}
                        \prod_a\operatorname{md}(P_aQ_a)^{-2e_1}
 \le (2/\operatorname{md}(M))^{n^6}\mathfrak h(p)\mathfrak h(q).
\]
The total shell and denominator loss has
\(e_1\le 2\cdot10^9n^4\), so the restriction holds.
Each omitted finite factor of the integral scalar \(M\) has modulus at
most one, so \(\operatorname{md}(M)\ge M^{-2f}\) on every projected set
\(T\), with \(f=[E:\mathbb Q]/2\). Equality holds for the adelic sum,
where \(T\) contains every place above \(M\).
Since \(M\le b^{23}\), \((2/\operatorname{md}(M))^{n^6}
\le b^{(46f+1)n^6}\le b^{n^9f}\).
The normalized finite estimate on the same \(Z\) is
therefore
\begin{equation}\label{eq:finite-frequency-final}
 \frac{|\mathcal I_j^T(p,q;\omega)|}{\mathcal R_Z}
 \le b^{n^9f}\mathfrak h(p)\mathfrak h(q)
     \operatorname{md}(\xi)^{-(k-n^7)/2}D(\xi)^{-6}.
\end{equation}
This holds for the common-entry character and for the
unequal-head form with \(\omega=1\) above.

\subsection{Finite measure estimate for the stabilizers}
\label{sec:finite-selected-bottom}

\sectionopening{In this subsection, we specialize the finite pairing to the stabilizers that occur in the
orbital argument. The denominator sum of
Section~\ref{sec:finite-shell-reduction} bounds the terms outside the
level subgroup, leaving the normalized contribution of
Section~\ref{sec:finite-bottom-normalization}.}

Take empty heads and \(\alpha=1\). Let \(Q_0\ge2\) bound
from below the indices of proper full integral ideals in
\(\mathcal R\). Off \(J\) every active pair has some
\(D_a\ge Q_0\) by Section~\ref{sec:finite-stable-shells}.
Apply the stage-zero estimate of
Section~\ref{sec:finite-absolute-pivots} with \(e=Q_0^{-1}\)
to \(f_0=\mathbf1_{G_T\setminus J}\mathcal M_0\).
The absolute shell sum in
Section~\ref{sec:finite-shell-reduction} applies with
\(l=t+10^9n^4\) whenever \(t+2\cdot10^9n^4\le k\).
In fact \(s-l\ge10^9n^4-24n-27\ge9\cdot10^7n^4\).
Absorbing \(2^{10^9n^4}\) by the spare power of \(Q_0^{-1}\)
gives \(\int f_0\,d\mu\le Q_0^{-t}\mathcal R_Z\).
The quotient correlation has squared modulus at most
\(\mathcal M_0\) by the triangle inequality.
With the exact contribution on \(J\) from
Section~\ref{sec:finite-bottom-normalization}, this yields
\begin{equation}\label{eq:finite-bottom-mass}
 \mathcal M(Z)\le\mu(Z)\mu(J)(\mu(Z)+Q_0^{-t}).
\end{equation}
Here \(\mathcal M(Z)\) denotes the mass of the conjugated
character, as in that section. On the entire projected set
we may take \(Q_0=2\) and \(t=n\). The compact upper bound
then gives
\[
 \frac{\mathcal M(K_T)}{\mu(K_T)\mu(J)}
 \le1+2^{1-n}\le1+b^{-6}.
\]
The last inequality follows for large \(x\) from \(n=m/2\)
and \(\log b\le x^{1/100}\).

\begin{lemma}[Finite mass of mixed stabilizers]\label{lem:marked-finite-fixer-mass}
For a set \(\mathcal A\) of pairs above \(p_0^{\rm lev}\), fix at
each pair a nonzero vector in the prescribed residue span or a
nonzero covector in the opposite-place realization, as in
Section~\ref{sec:full-orbital-blocks}. Let \(Z\) be the inverse
image of these stabilizers, with the full compact group at the
other places of \(T\). Then
\[
 0\le \frac{\mathcal I_0^T(1)}{\mathcal R_Z}
       =\frac{\mathcal M(Z)}{\mu(Z)\mu(J)}
       \le(1+b^{-6})b^{-6k|\mathcal A|}.
\]
\end{lemma}

\begin{proof}

For the chosen places \(\mathcal A\) defining \(\Gamma_{\mathcal A}\), take
\(Z\) to be the inverse images of its vector or covector
fixers on the selected conjugate pairs and the full
compact subgroup on the complement in \(T\).
This contains \(J\) and contains every integral
lift with reduction in \(\Gamma_{\mathcal A}\).
On one selected split pair let \(K,J_w,Y\) denote the
compact, level and fixer inverse image, and put
\(q=p_0^{\rm lev}\ge b^{10}\).

The residue group
\(\mathrm{GL}_r(\mathbb F_q)\) acts transitively on the
nonzero vectors and on the nonzero covectors. Reduction is
surjective by lifting matrices with unit determinant.
Thus \(\mu(Y)=\mu(K)/(q^r-1)\).
In the Haar measure of Section~\ref{sec:compact-floor-haar},
\[
                 \mu(K)=\prod_{l=k+1}^{r}(1-q^{-l}).
\]
Indeed, in rank \(l\) a split vector is a column \(a\) and a
covector \(b'\), with norm \(b'(a)\). The integral orbit at
an isotropic step consists of primitive \(a,b'\) with
\(b'(a)=0\): extend \(a\) to a basis of \(\ker b'\), then append
a vector on which \(b'\) is one. Its residue cardinality is
\((q^l-1)(q^{l-1}-1)\). At the last step the norm is \(-N\),
a unit at the chosen prime above \(M\). The vector \(a\) is
primitive and \(\ker b'\) is an integral complement to its line
since \(b'(a)=-N\). The residue cardinality is
\((q^l-1)q^{l-1}\).
The differential \((v,c)\mapsto b'(v)+c(a)\) is surjective
at every residue solution under consideration. Each solution
therefore has \(q^{2l-1}\) lifts at each next residue level,
so the fibre mass is the residue count divided by \(q^{2l-1}\).
Multiply for \(l=r,r-2,\ldots,k+3\) and then \(l=k+1\), since
the integral Heisenberg factors and terminal compact group
have measure one. In particular,
\(\mu(K)\ge1-\sum_{l>k}q^{-l}\ge1/2\ge2^{-n}\).

Put \(Q_0=b^9\) and \(a_0=k-10^{10}n^4\). For large \(x\),
\(q^r-1\ge Q_0^{a_0+1}\), \(9a_0\ge6k\), and
\(a_0+1+n+2\cdot10^9n^4\le k\), since \(r\ge k\)
and \(k/n^4\to\infty\). A nonunit integral element
on this pair has modulus at most \(q^{-1}\le Q_0^{-1}\).
Apply \eqref{eq:finite-bottom-mass} with \(t=a_0+1+n\).
As \(Q_0^{-n}\le 2^{-n}\le\mu(K)\), it gives
\[
 \frac{\mathcal M(Y)}{\mu(Y)}
 \le \mu(J_w)\mu(K)\,2 Q_0^{-a_0-1}
 \le b^{-6k}\frac{\mathcal M(K)}{\mu(K)}.
\]
The second inequality uses \(Q_0\ge2\) and
\(\mathcal M(K)\ge\mu(K)^2\mu(J_w)\), the contribution
on \(J_w\). On the complementary factors use \(Z=K\).
The group, level and character factor over these
conjugation-stable blocks. The inequalities for
\(\mathcal M(Y)/\mu(Y)\) relative to
\(\mathcal M(K)/\mu(K)\) are unchanged under rescaling
local Haar, since both ratios scale by its square.
Tonelli and the bound for \(K_T\) therefore prove
\[
 0\le\frac{\mathcal I_0^T(1)}{\mathcal R_Z}
       =\frac{\mathcal M(Z)}{\mu(Z)\mu(J)}
       \le (1+b^{-6})b^{-6k|\mathcal A|}.
\]
This uses the same \(Z\), original finite character, and
stage Haar as the terminal and lower-stage forms.
\end{proof}

\section{Weighted real orbital integrals}
\label{sec:real-orbital-analysis}

\sectionopening{In this section, we prove the real Fourier and square-integral bounds used to compare
orbital integrals at a common entry in Section~\ref{sec:real-common-descent}.
The estimates are in the norm-fiber measures fixed below.}

\subsection{Plancherel on the real norm fibers}
\label{sec:real-floor-plancherel}

\sectionopening{In this subsection, we identify a real Fourier transform on two Hermitian norm fibers by
ordinary Plancherel. Its measure normalization is used in the phase
calculation of Section~\ref{sec:real-orbital-phases} and in the
integral of Section~\ref{sec:real-position-ball}.}

At a complex place of \(E\), use ordinary area measure
on \(\mathbb C\), the character \(\psi(c)=e^{c-\bar c}\), and
norm-value coordinate \(h(y,y)/2\) with ordinary Lebesgue
measure. The factor \(1/2\) specifies the real skew coordinate.
Write \(h(v,y)=v^*Jy\) in the common stage basis.
Among the positive orthonormal coordinates choose a free
set \(I\) of size \(d=k-5\) and disjoint two-coordinate
blocks \(B,B'\), both disjoint from \(I\), inside the definite
summand \(X_0\). They occupy \((k-5)+2+2=k-1\) of its coordinates.
Let \(u\) be the remaining unit vector. In the normalized real
basis of Section~\ref{sec:period-real-point} one has \(h(z,z)=-1\).
The vectors \((u+z)/\sqrt2,(u-z)/\sqrt2\) form a hyperbolic pair
orthogonal to all three blocks even at the last stage
\(V_\varrho=X_0\perp\mathbb Cz\). This change and the signed basis
change \(p_v\) are Euclidean unitary and preserve volume and operator
norms. Solve the source quadric in \(B\) and the target quadric in \(B'\).

For either block the equation at fixed other coordinates is
\(|z_1|^2+|z_2|^2=2s\). Its fiber measure is sphere probability
times \(4\pi^2s\) for \(s>0\), from the complex two-ball
volume \(2\pi^2s^2\). At \(s=0\) the conditional measure
in the two solved coordinates has mass zero. In a box \(|z_j|\le L^\ell\), the conditional
mass is at most \(4\pi^2 L^{2\ell}\). For a full vector measure
on that side the corresponding bound is \(\pi^2L^{4\ell}\).

The sphere normalization agrees with the hyperbolic chart
on each fixed norm fiber. In that chart, with a hyperbolic
pair \(t,v\), constrain \(2\Re(\bar t v)+Q_0+|z|^2\)
and put \(c=v/t\). The Jacobian of \(v=tc\) is \(|t|^2\),
which is also the derivative of the norm-value coordinate
with respect to \(\Re c\). The two factors cancel in the
fiber measure, leaving ordinary measure in \(\Im c\).
At fixed norm, replacing \(|z|^2/2\) by \(\Re c\) restores
the factor \(|t|^2\), so together with \(\Im c\) it gives
\(d^2v\). This yields the stated sphere mass on that chart.
Choose the hyperbolic pair outside the solved positive block.
Its coordinate \(t\) is free in the sphere description, so \(t=0\)
is null on every fixed fiber, including at norm zero.

To verify invariance at each norm, write \(a\) for the coordinates
outside the solved block, \(s=b-h(a,a)/2\), and \(s_+=\max(s,0)\).
For \(F\in C_c\) the prescribed fiber measure at norm value \(b\)
satisfies
\[
 \int F\,d\mu_b
 =4\pi^2\int s_+\int_{S^3}
     F(a,\sqrt{2s_+}\,\omega)\,d\omega\,da ,
\]
where \(d\omega\) has mass one. This is continuous in \(b\) by
dominated convergence on the compact projection of the support of
\(F\). An isometry preserves ambient volume and the norm value,
so disintegration gives \(\int F\circ g\,d\mu_b=\int F\,d\mu_b\)
for almost every \(b\). The two continuous functions agree at every
\(b\), as required for the rotations of
Section~\ref{sec:real-orbital-phases}.

Write the source and target as \(y=(a,w,z)\) and
\(v=(a',w',z')\), with \(w,w'\) on \(I\) and \(z,z'\)
on \(B,B'\), respectively. The two measures on product boxes are
\[
             da\,dw\,d\nu_{a,w}(z),\qquad
             da'\,dw'\,d\nu'_{a',w'}(z').
\]
Let \(V,V'\) be the
Lebesgue volumes of the \(a,a'\) boxes and \(M,M'\) the
conditional block bounds. Orthogonality gives
\[
 h(v,y)=\overline{w'}\cdot w
       +\overline{z'}\cdot a|_{B'}+A(a';y).
\]
The last term is independent of \(w',z'\), and the
middle term depends on the second variable only through \(a\).

For a function \(\phi\) supported in the source box and
\(\xi\ne0\), define
\[
 F_{a',a}(w)=\int\psi\bigl(\xi A(a';a,w,z)\bigr)
                          \phi(a,w,z)\,d\nu_{a,w}(z).
\]
For the fixed first variable, Cauchy--Schwarz in \(a\) bounds the
squared transform by
\[
 V\int\left|\int_{\mathbb C^I}
     \psi(\xi\overline{w'}\cdot w)F_{a',a}(w)\,dw\right|^2 da .
\]
Integrate \(z'\) with contribution \(M'\) and extend the \(w'\)
integration to all of \(\mathbb C^I\). For one coordinate,
writing \(w=x+iy\) and \(\xi\overline{w'}=A_0+iB_0\)
gives the Euclidean Fourier frequencies \((B_0/\pi,A_0/\pi)\).

Their Jacobian relative to \(w'\) is \(|\xi|^2/\pi^2\).
Euclidean Plancherel therefore gives
\((\pi^2/|\xi|^2)^d\). Finally, Cauchy--Schwarz in the
source block \(z\) and integration over \(a'\) give
\[
 \int_{\rm out}\left|\int_{\rm in}
    \psi(\xi h(v,y))\phi(y)\,d\mu(y)\right|^2d\mu'(v)
 \le VV'MM'\left(\frac{\pi^2}{|\xi|^2}\right)^d
                             \int_{\rm in}|\phi(y)|^2d\mu(y).
\]
If either side is a full vector box, take ordinary area measure in each
of its two block coordinates. Their product mass is at most
\(\pi^2 L^{4\ell}\), and the same two applications of Cauchy--Schwarz
give the displayed bound with this value of \(M\) or \(M'\).
All the measures on the boxes are finite, so the integrals apply to
square-integrable functions.

\subsection{Support and phases in the real integral}
\label{sec:real-orbital-phases}

\sectionopening{In this subsection, we compute the characters of the two real kernels on the support of
the truncated coefficient. The Plancherel normalization of
Section~\ref{sec:real-floor-plancherel} gives the phase used in
Section~\ref{sec:real-weighted-products}.}

Use the chosen basis and the groups \(G_i,R_i,U_i\) as in
Section~\ref{sec:compact-floor-haar}, now over \(\mathbb C\).
Use \(J=J_1=D_v^*J_0D_v\) of Section~\ref{sec:period-real-point}.
Thus the hyperbolic pairs have Gram \(\left(\begin{smallmatrix}0&1\\1&0\end{smallmatrix}\right)\)
and the middle block has Gram \(I_k\oplus(-1)\).
In these coordinates \(J=J^*=J^{-1}\), so
\(g^{-1}=Jg^*J\) and \(\|g^{-1}\|=\|g\|\).
The passage to the basis ordered by signs is Euclidean unitary.
Let \(\mathcal K\) be the real box of operator norm at most
\(R_b=b^{6n+2}\). In this real calculation put \(\mathcal U_i=R_i\mathcal K\). Let \(r\) be the ambient rank,
\(L\ge\max(R_b,8r)\), \(n=\varrho+1\), and
\(\Lambda_L=L^2\). Entries of size \(a\) have modulus at
most \(L^a\). On a head \(p\) of length \(j\le\varrho\),
assume \(|\prod_{a<l}p_a|^{-1}\le L\) for \(l\le j\), and let
\(w_p\) send \(e_a\) to \(p_a e_a^*\) and \(e_a^*\) to
\(\bar p_a^{-1}e_a\) for \(a<j\).

Suppose \(g,w_p^{-1}g\in\mathcal U_{j}\).
Factor the move at \(g\) in \(R_j=G_jU_{<j}\). Thus some
\(t\in G_j\) gives \(stg\in\mathcal K\) with \(s\in U_{<j}\).
The other move can be written \(s'w_p^{-1}tg\in\mathcal K\)
with \(s'\in R_j\), since \(t\) commutes with \(w_p\).
The ratio of the boxed matrices is \(s w_p(s')^{-1}\),
of norm at most \(R_b^2\).
In its first \(l\) flag columns, right unitriangular terms
from \((s')^{-1}\) preserve every minor of order \(l\) in those columns.
In the rows indexed by \(e_0^*,\ldots,e_{l-1}^*\), the columns
\(s e_a^*\), \(0\le a<l\), form a lower unitriangular matrix.
Hence the minor of the ratio in these rows and columns is
\(P_l=\prod_{a<l}p_a\).

For \(l=a+1\), use rows \(e_0^*,\ldots,e_{a-1}^*,d\) in
that order. Remove the right unitriangular factor on the first
\(a+1\) columns and factor \(p_b\) out of column \(b\).
The last column \(s e_a^*\) is zero in the first \(a\) rows and
the leading \(a\)-row minor is one. Expansion in that column gives
\(P_{a+1}(s e_a^*)_d\), including the zero value for a repeated row. Hadamard's inequality bounds these minors
by \(L^{2l}\). Dividing by \(P_l\) and using consecutive
products bounds \(p_a,p_a^{-1}\) and the entries of \(s e_a^*\)
by \(L^{2n+1}\).
In Witt order \(s\) has blocks
\[
 s=\begin{pmatrix}U&B&C\\0&1&D\\0&0&U^{-*}\end{pmatrix},
 \qquad B=-UD^*J_0 ,
\]
where \(J_0\) is the Gram matrix on the remaining space.
The last columns bound \(C,D,U^{-*}\). Since \(U^{-*}-1\)
is nilpotent of index at most \(j\), the geometric sum
bounds \(\|U\|\) by \(L^{(2n+4)n}\). Summing the block norms
gives \(\|s^{-1}\|=\|s\|\le L^{50n^2}\). Thus \(tg\)
has size \(60n^2\) and \(tw_p^{-1}g\) has size \(99n^2\).
For \(j=0\) take \(t=g^{-1}\).

Consider a tail correlation at \((\sigma g,\sigma g')\)
that vanishes unless \(\sigma g\in\mathcal U_{j+1}\)
and is invariant under the stabilizer of \(f_j\) in \(G_j\).
It factors through \(y(\sigma)=\sigma^{-1} f_j\).
If \(t\in G_j\), \(t^{-1}q\in U_{<j}\), and \(qg\) has size \(a\),
then \(t y(\sigma)\) has size \(a+6\) on its support.
In fact \(\sigma(t^{-1}q)\sigma^{-1}\in R_{j+1}\),
so \((\sigma t^{-1})(qg)\in\mathcal U_{j+1}\).
Choose a box move \(\zeta\in R_{j+1}\) for this point.
Projecting the column \((qg)(\zeta\sigma t^{-1}qg)^{-1}f_j\)
to the remaining space gives \(t y(\sigma)\). The earlier
isotropic terms from \(\zeta^{-1}f_j\) project to zero.
The column norm is at most \(rL^a R_b\), as required.

For the character twist, the basis entry is
\[
 h(f_j,\sigma(\sigma')^{-1}f_j)=h(y(\sigma),y(\sigma')) .
\]
The two tail functions \(F_1,F_2\) start at \(g,w_p^{-1}g\).
If both contribute, both initial points belong to \(\mathcal U_j\).
Rotate both orbit variables by the same \(t\) just constructed.
This preserves the pairing and the two norm-fiber measures.
The support sizes are \(60n^2+6\) and \(99n^2+6\), so
the pairings have modulus at most \(L^{200n^2}\).

For a character agreeing with \(\psi(\xi\,\cdot)\) on that
box, the integral in the rotated variables is
\(\iint \psi(\xi h(v,y))F_1(v)\overline{F_2(y)}\,d\mu(v)d\mu(y)\).
Cauchy--Schwarz and Section~\ref{sec:real-floor-plancherel}
bound its square by
\[
 \Lambda_L^{10^6n^3}
 \left(\frac{\pi^2}{|\xi|^2}\right)^{k-5}
                    \|F_1\|_2^2\|F_2\|_2^2 \qquad (\xi\ne0).
\]
To check the box contribution, the hyperbolic change on one
positive and one negative line is Euclidean unitary, so
a support of size \(a\) lies in coordinate boxes of radius
\(L^{a+1}\).

The number of ordinary coordinates outside
the free and solved blocks is at most \(2n+2\).
The ordinary volume and the block mass on that side have
product at most \(\Lambda_L^{(2n+4)(a+2)}\).
Apply this to the two support sizes. The conjugated transform
uses frequency \(-\xi\), of the same modulus.

For the unequal heads put \(j=i+1\le\varrho\),
\(d_0=w_q^{-1}w_p\), and let \(m_0\) be the level unit,
equal to \(1\) at the real place. The pairing kernel has
left multiplier
\(\psi(m_0^{-1}(\mathfrak s(d_0u_2d_0^{-1})-\mathfrak s(u_1)))\)
under \((u_1,u_2)\) in the lower unipotent, with \(\mathfrak s\)
the simple-root sum. The stage-\(i\) element \(u_i(v,\beta)\)
sends \(o\) in the remaining space to \(o-h(v,o)e_i\),
so its simple-root sum is \(-h(v,f_{i+1})\).
Conjugation by \(d_0\) multiplies this by \(p_i/q_i\).
Moving \(u_i\) to the left of the next pivot \(\sigma\)
changes its vector to \(\sigma v\). Thus the section phase is
\[
 \psi\bigl(-\xi_0 h(v,\sigma^{-1}f_{i+1})\bigr),\qquad
 \xi_0=m_0^{-1}(p_i/q_i-1)=-p_i\xi,\qquad
 \xi=m_0^{-1}(p_i^{-1}-q_i^{-1}).
\]
Equivalently it is \(\psi(-\xi h(v',\sigma^{-1}f_{i+1}))\)
for \(v'=\bar p_i u_i^{-1}e_i^*\), since the two terms in
the \(i\)-th pair are orthogonal to the orbit vector and
the remaining part of \(u_i^{-1}e_i^*\) is \(-v\).

Write \(S\subseteq U_i\) for the support of this correlation
restricted to \(w_p u_i g\in\mathcal U_{j}\). On \(S\) both
\(u_i g\) and \(w_p u_i g\) are near. If \(S\) is nonempty,
choose \(u_*\in S\). Apply the two-point movement to
\(w_p u_i g,u_i g\), so that some \(t_i\in G_j\) puts
\(t_i u_i g\) in size \(99n^2\). Use \(t_*\) for \(u_*\).
Recenter by \(m=t_*u_i u_*^{-1}t_*^{-1}\), which preserves
Heisenberg Haar. Since \(t_i,t_*\) fix \(e_i^*\),
\[
 (m^{-1}-1)e_i^*
  =\bigl((t_*u_*g)(t_i u_i g)^{-1}-1\bigr)e_i^* .
\]
The two matrix bounds and \(J^{-1}=J\) give size
\(198n^2+6\) for this column. By the head bound, its
product by \(p_i\) has size \(999n^2\).
The choices \(t_i\) establish these bounds, and the
change of integration variable uses the single \(t_*\).

Set \(u=p_i^{-1}\) and write the vector of \(m\) as \(u w\).
The bound for \(u^{-1}(m^{-1}-1)e_i^*\) gives a box of
size \(999n^2\) for \(w\). At fixed vector of \(m\),
the coefficient at \(e_i\) of the inverse column depends
on the skew center by a translation and a sign.
Subtracting the bounds for any two centers in the same
fiber shows their distance is at most \(2|u|L^{999n^2}\).
The skew Haar measure is ordinary length, so that is also
a bound for the fiber mass. The vector Jacobian is
\(|u|^{2d_i}\), \(d_i=\dim_{\mathbb C} V_{i+1}\), where
\(V_{i+1}\) is the space orthogonal to the first \(i+1\)
pairs. The combined factor is
\(|u|^{2d_i+1}=\mathfrak m_i(p_i)^{-1}\), with
\(\mathfrak m_i(p_i)=|p_i|^{2d_i+1}\).

Rotating the next orbit vector by \(t_*\) preserves its
measure. Its size is at most \(99n^2+6\), by applying
the projected column argument with \(q=t_*u_*\) and
\(t=t_*\). The recentered vector introduces a constant
translation in the Fourier phase, absorbed as a phase of
modulus one on the function. Since \(h(uw,y)=\bar u h(w,y)\)
and \(u=p_i^{-1}\), the frequency after the scaling to \(w\) is
\[
 -\xi_0\bar u=\frac{p_i}{\bar p_i}\,\xi,\qquad
 |-\xi_0\bar u|=|\xi| .
\]
This also applies to nonreal pivots. Apply Section~\ref{sec:real-floor-plancherel}
from the orbit quadric to the vector box, and use the
central fiber bound.

If \(C_{i+1}\) is the tail in the
section integral and \(C_{i+2}\) its next tail, this proves
\[
 \int_S |C_{i+1}(u_i g,u_i g')|^2\,du_i
 \le
 \left(\frac{\pi^2}{|\xi|^2}\right)^{k-5}
 \frac{\Lambda_L^{10^6n^3}}{\mathfrak m_i(p_i)}
 \int |C_{i+2}(\sigma g,\sigma g')|^2\,d\theta_{i+1}(\sigma).
\]
Here the box contribution on the two sides and the central factor
\(2L^{999n^2}\) together are at most \(\Lambda_L^{10^6n^3}\).
This calculation uses the tail support and equivariance of
the cut pairing kernels. The inverse-depth weights enter
when one integrates these estimates over the real positions.

\subsection{The terminal Gram form}
\label{sec:real-terminal-gram}

\sectionopening{In this subsection, we write the pairing on the remaining definite space as a Gram form.
Its positivity and continuity are used at the terminal step of the
descent in Section~\ref{sec:real-common-descent}.}

Let \(F\) be the continuous cut Bessel function with
values in a Hermitian space \(I\) of dimension \(d_I\).
Use the original real Haar \(dg\) in its orbital tensor and
assume the cut support \(U_{<n}\mathcal K\).
In the remaining space at the top, \(h(z,z)=-1\) and
\(z^\perp\) is positive definite. Put \(H=U(z^\perp)\).
For \(c\ne0\), define \(a(c)\) on the plane of one positive unit vector and \(z\) by
\[
 \begin{pmatrix}\bar c&\sqrt{|c|^2-1}\\
          -\sqrt{|c|^2-1}&-c\end{pmatrix}\qquad (|c|\ge1).
\]
It is the identity on the other coordinates. For \(0<|c|<1\)
use \(a(c/|c|)\). This fixes the earlier pairs. If
\(\sigma\in G_\varrho\) and \(h(z,\sigma z)=c\), the positive
part of \(\sigma z\) has length \(\sqrt{|c|^2-1}\).
Rotate it by \(H\) to the corresponding positive coordinate
of \(a(c)z\). The resulting quotient fixes \(z\).
Thus \(\sigma=h_1a(c)h_2\) with \(h_1,h_2\in H\).

For a head of length \(\varrho\) and last entry \(c\), put
\[
 \begin{split}
 T_p(\sigma)&=\int F(y)\otimes
                    F((w_p\sigma)^{-1}y)^*\,dg(y),\\
 A_p&=\frac{\sqrt{c_G}}{d_I}(|c|^2-1)_+^{(k-1)/2}
          \prod_{i<\varrho}|p_i|^{k+1+2(\varrho-i-1)},\\
 \Xi(p)&=A_p T_p(a(c)),\qquad
 J_n(p,q)=\langle\Xi(p),\Xi(q)\rangle .
 \end{split}
\]
Here \(c_G\) is the real open-cell normalization and
the star on a vector denotes its Hermitian dual.
For \(p,q\) with the same last entry and \(|c|>1\),
translate \(y\) by \(h_1\) and use \(w_p h_1=h_1w_p\).
The two tensors at \(\sigma\) then both arise from the
tensors at \(a(c)\) by the same unitary operation on the
two value factors, using the \(H\)-equivariance of \(F\).
Consequently
\[
 J_n(p,q)=\big\langle A_p T_p(\sigma),A_q T_q(\sigma)\big\rangle.
\]
At \(|c|=1\) the scale factors vanish, since \(k>1\).

This is a positive semidefinite terminal form by its Gram
definition. It is continuous in the nonzero head entries.
Indeed inverse partial products are bounded on a small
neighborhood. Nonvanishing of both factors in \(T_p(a(c))\)
gives a move at \(y\) in \(U_{<n}\) and a move at \(w_p^{-1}y\)
in \(R_\varrho\), as \(a(c)\in G_\varrho\).
The two boxed matrix bounds above apply with the first
move already lower unipotent, giving a common bounded
set for \(y\). In the real matrix group its closure
is compact. Continuity of the integral follows by dominated
convergence, including at \(|c|=1\).

\subsection{Integrals of reciprocal determinant weights}
\label{sec:real-inverse-depth}

\sectionopening{In this subsection, we integrate the inverse powers of the relative vector lengths against
the Haar factors of Section~\ref{sec:real-orbital-phases}. The
resulting inequalities control the weighted products in
Section~\ref{sec:real-weighted-products}.}

Keep the real box and stage measures from
Sections~\ref{sec:compact-floor-haar} and~\ref{sec:real-orbital-phases}.
Fix an allowed relative position \(x\). Put
\(A_a=\langle e_0,\ldots,e_{a-1}\rangle\) and
\(s=k-3n-6\) (this differs from the finite exponent \(s=k-24n-27\) of Section~\ref{sec:finite-position-box}). For \(a<\varrho\), the depth \(D_a(g,x)\)
is the reciprocal of the area of the set of \(h(e_a,v)\)
with \(v\in A_a^\perp\) and with some \(v'\equiv v\bmod A_a\)
satisfying
\[
                    \|g^{-1}v\|^2+\|(gx)^{-1}v'\|^2<1 .
\]
Let \(H_{g,a}\) denote the positive Gram matrix of the quotient norm
induced by \(v\mapsto\|g^{-1}v\|\) on \(A_a^\perp/A_a\). In any fixed
quotient basis this is the Schur complement of the \(A_a\)-block in
the restricted Gram matrix, hence continuous in \(g\).
The two lifts of a given class may be minimized independently.
With \(H=H_{g,a}+H_{gx,a}\) and the nonzero row
\(\ell_a=h(e_a,\cdot)\) on the quotient, one obtains
\[
 \min_{\ell_a w=c} w^*Hw
     =\frac{|c|^2}{\ell_a H^{-1}\ell_a^*},\qquad
 D_a(g,x)=\frac1{\pi\ell_a H^{-1}\ell_a^*}.
\]
The minimizer is \(cH^{-1}\ell_a^*/(\ell_a H^{-1}\ell_a^*)\).
The image is therefore a disc of the claimed area, and \(D_a\) is
continuous in both variables. The near sets are measurable as well:
each \(R_a\mathcal K\) is the product of a closed subgroup with a
compact set and is closed.
The depth is unchanged under simultaneous left movement
by \(R_{a+1}\). Such movement transports the two lifts,
preserves their congruence and their lengths, and changes
\(e_a\) by an element of \(A_a\).

In the common free positive space of dimension \(d=k-5\),
let \(E_I\) be its orthonormal inclusion. Suppose at most \(n\)
linear equations are imposed on \(u\in\mathbb C^d\). Impose
also that \(g^{-1}E_Iu\) and \((gx)^{-1}E_Iu\) belong to
the positive eigenspace of the involution \(J\).
The combined rank is at most \(3n\), since each negative
projection has rank at most \(n\). Choose an isometric
inclusion \(O\) of an \((s+1)\)-plane in the kernel.
For \(u\) in the plane with \(\|u\|\le1/2\), preservation
of the Hermitian form gives the same Euclidean norm for
both inverse translates.

For \(w=E_IOu\) with \(\|u\|\le1/2\),
\[
 \|g^{-1}w\|^2+\|(gx)^{-1}w\|^2=2\|u\|^2\le1/2 .
\]
If the free coordinate of the vector paired against \(w\) is \(v\),
its values \(v^*Ou\) fill a disc of area
\(\pi\|O^*v\|^2/4\). The inverse area is therefore at most
\(\phi(v)=36/\|O^*v\|^2\), allowing \(+\infty\) when \(O^*v=0\).
For each box exponent \(\ell\ge0\), uniformly in the equations and
the two points,
\begin{equation}\label{eq:real-free-depth}
 \int_{\{|v_i|\le L^\ell\}}\phi(v)^s\,dv
     \le \Lambda_L^{99n(\ell+1)} .
\end{equation}
Indeed the box is contained in the ball of squared radius
\(R^2=dL^{2\ell}\). The complement of the \((s+1)\)-plane has
complex dimension \(h=d-s-1=3n\). At \(t=\|O^*v\|^2\) the radial
measure on the plane is \(\pi^{s+1}t^s\,dt/s!\), and the
complementary ball has volume \(\pi^h(R^2-t)^h/h!\). Hence
\[
 \int_{\|v\|\le R}\phi(v)^s\,dv
 =\frac{36^s\pi^d}{s!\,h!}\int_0^{R^2}(R^2-t)^h\,dt
 =\frac{36^s\pi^d}{s!\,(h+1)!}R^{2(h+1)} .
\]
The coefficient is at most \(e^{36\pi}e^\pi\). Since
\(d\le r\le L/8\) and \(\Lambda_L=L^2\), the last expression is
at most \(e^{37\pi}L^{(2\ell+1)(3n+1)}\le\Lambda_L^{99n(\ell+1)}\).

For the first inverse-depth integral, if some \(\sigma g\in\mathcal U_{a+1}\) with \(\sigma\in G_a\), then \(g\in\mathcal U_{a}\). Rebase \(g\) to the box by \(R_a\). The stage
quotient integral is invariant under this replacement,
by the right-action argument of Section~\ref{sec:compact-floor-haar}
applied to the left-\(R_{a+1}\)-invariant depth and near indicator.
At a boxed starting point choose \(g_1\in R_{a+1}\) with
\(g_1\sigma g\in\mathcal K\). The projection of
\(g(g_1\sigma g)^{-1}e_a=\sigma^{-1}g_1^{-1}e_a\) to the remaining
space is \(y=\sigma^{-1}e_a\), since \(\sigma\in G_a\) fixes the
earlier pairs. Thus \(\|y\|\le R_b^2\). With no extra
equations take \(\sigma w\) for both lifts.
It belongs to \(A_a^\perp\) and
\(h(e_a,\sigma w)=h(y,w)\), so \(D_a(\sigma g,x)\le\phi(y_I)\).
The stage coordinates have modulus at most \(L^{11}\). Besides the
\(k-5\) free coordinates in \(I\) and the two positive coordinates
used for the norm equation, there are
\(r-2a-(k-5)-2\le 2n+2\) complex coordinates. Thus at fixed
\(y_I\) their boxes and the two-coordinate norm fiber contribute at most
\[
    (\pi L^{22})^{2n+2}\,4\pi^2 L^{22}
                 \le\Lambda_L^{99n(7+6)}
\]
by Section~\ref{sec:real-floor-plancherel}, in the coordinate
\(h(y,y)/2\). By
\eqref{eq:real-free-depth},
\begin{equation}\label{eq:real-depth-first}
 \int_{\{\sigma g\in\mathcal U_{a+1}\}} D_a(\sigma g,x)^s\,
                     d\theta_a(\sigma)\le \Lambda_L^{10^4n}.
\end{equation}

For the second integral, take a boxed point \(y\), a lower
unipotent \(\nu\in U_{<a}\), and the head Weyl element \(w_{p,a+1}\)
of length \(a+1\). For common short vectors \(w=E_Iu\) at
\(y,yx\), impose the \(a\) equations \(h(e_b^*,\nu w)=0\)
for \(b<a\). They make \(w_{p,a+1}u_a\nu w\) orthogonal
to \(A_a\), since \(u_a\in U_a\) fixes the lower pairs.
Writing \(c=\bar p_a\), the inverse image of \(e_a\)
is \(\nu^{-1}c u_a^{-1}e_a^*\), whose class in
\(A_a^\perp/A_a\) is represented by \(z_c=c u_a^{-1}e_a^*\).
The preimages of \(e_b\) for \(b<a\) depend on \(\nu\)
and the preceding \(p_b\), so the chosen common plane
is fixed throughout the \(U_a\) integral. Transporting
both lifts gives \(D_a(w_{p,a+1}u_a\nu y,x)\le\phi((z_c)_I)\).

On \(\mathcal U_{a+1}\), choose a box move \(\gamma\).
The matrix \((u_a\nu)^{-1} w_{p,a+1}^{-1}\gamma^{-1}\)
is the ratio \(y(\gamma w_{p,a+1}u_a\nu y)^{-1}\) of two
boxed matrices and has operator norm at most \(R_b^2\).
Use the row-replacement formula for minors in
Section~\ref{sec:real-orbital-phases} with inverse head
\(\bar p\). After multiplying the column by \(\bar p_a\),
its size bound is
\(R_b^{2(a+1)}|\prod_{b<a}p_b|^{-1}\le L^{2a+3}\).
Projection to the remaining space gives \(z_c\), so
we may use a box of size \(20n\).
At a fixed vector coordinate of \(u_a\), integration in
its center contributes at most \(2L^{20n}/|c|\).
Changing the vector coordinate to the corresponding
part of \(z_c\) gives \(|c|^{-2d_a}\), where
\(d_a=\dim_{\mathbb C}V_{a+1}\).
In the stage basis the remaining coordinates have size
at most \(20n+4\), so their off-\(I\) volume and the
integral in \eqref{eq:real-free-depth} each contribute at most
\(\Lambda_L^{99n(20n+5)}\). With
\(\mathfrak m_a(p_a)=|p_a|^{2d_a+1}\) this gives
\begin{equation}\label{eq:real-depth-second}
 \int_{\{w_{p,a+1}u_a\nu y\in\mathcal U_{a+1}\}}
    D_a(w_{p,a+1}u_a\nu y,x)^s\,du_a
       \le \mathfrak m_a(p_a)^{-1}\Lambda_L^{10^6n^2}.
\end{equation}

Set
\[
 W_j(g)=\prod_{a<j}D_a(g,x)^{-s},\qquad
 I_j=\mathbf1_{\mathcal U_j}W_j^{-1},\qquad
 C_L=\Lambda_L^{10^6n^3},\qquad
 H_p=\prod_{a<j}\mathfrak m_a(p_a)^{-1}.
\]
Since the earlier depths are unchanged by \(G_a\),
\eqref{eq:real-depth-first} implies
\(\int I_{a+1}(\sigma g)\,d\theta_a(\sigma)
\le\Lambda_L^{10^4n} I_a(g)\).
Similarly in \eqref{eq:real-depth-second} the earlier
depths and indicator of \(\mathcal U_a\) read at \(w_{p,a}\nu y\).
Indeed \(u_a\) and the last Weyl factor belong to \(G_a\)
and commute with the earlier Weyl factors. Iteration
therefore bounds \(\int I_j(w_p u y)\,du\) over \(U_{<j}\)
by \(H_p C_L\) for boxed \(y\). The iteration uses the
factorization \(u=u_a\nu\) at each step, valid for Haar
since \(U_a\subseteq G_a\) preserves the earlier measures
by conjugation.
For \(g\in\mathcal U_{j}\), write \(y=\gamma u_0g\) boxed,
with \(\gamma\in G_j,\ u_0\in U_{<j}\).
The same bound holds at \(g\) by the substitution
\(u\mapsto\gamma^{-1}u\gamma u_0\), since \(\gamma\)
commutes with \(w_p\) and preserves \(I_j\).
In particular, for a nonnegative left-\(R_l\)-invariant
\(B\) supported on \(\mathcal U_l\), where \(l\ge j\),
and a \(\sigma\)-finite lower representative measure \(\rho\), Tonelli gives
\begin{equation}\label{eq:real-depth-cross}
 \int B(g)I_j(w_p g)\,d(du_{<j}*\rho)(g)
                   \le H_p C_L\int B\,d\rho .
\end{equation}
The bound also holds with \(w_p^{-1}=w_{\bar p}\), for
which the initial segment bounds and absolute values are unchanged.

\subsection{The weighted products of the real correlations}
\label{sec:real-weighted-products}

\sectionopening{In this subsection, we bound the product of the two real kernels with the phases and Haar
factors of Section~\ref{sec:real-orbital-phases}. These bounds
propagate through successive isotropic quotients in
Section~\ref{sec:real-depth-majorant}.}

Write \(\theta_{[a,l]}=\theta_l*\cdots*\theta_a\), with the
empty convolution the atom at the identity, and
\(M_j=du_{<j}*\theta_{[0,j-1]}\). The decomposition of the
group Haar measure has the form \(d\mu=d\mu_j*M_j\),
where \(d\mu_j\) is Haar measure on \(G_j\) with the stage
normalization. We use the following change of representatives for
this decomposition and for the lower unipotent decomposition of
Section~\ref{sec:real-position-ball}.

Suppose \(d\mu=dm*\rho\), with \(dm\) Haar measure on a
unimodular subgroup \(H\). For compact \(A\subseteq G\) and compact
\(C\subseteq H\) with \(m(C)>0\), one has
\(m(C)\rho(A)\le\mu(CA)\). Thus \(\rho\) is \(\sigma\)-finite.
For integrable \(h>0\) put
\(\bar h(y)=\int_H h(u y)\,dm(u)\). This average is
left-\(H\)-invariant by right invariance of \(dm\), and lies in
\((0,\infty)\) for \(\rho\)-almost every representative, since
\(\int\bar h\,d\rho=\int h\,d\mu\). For nonnegative left-\(H\)-invariant \(B\),
\[
 \int B\,d\rho=\int \frac{Bh}{\bar h}\,d\mu
\]
by Tonelli, with the integrand zero on the exceptional set.
If \(s\) normalizes \(H\) and preserves \(dm\), then
\(dm*(L_s)_*\rho=(L_s)_*d\mu=d\mu\), where \(L_s\) denotes left
translation. The identity gives the same integral for \((L_s)_*\rho\).
In particular, take \(H=G_j\) and \(s=w_p\), which commutes with \(G_j\).

Use the kernels
\(K_1(g,y)=\langle F(g),F(y)\rangle\) and
\(K_2(g,y)=\langle F(g),F(w_q^{-1}w_p y)\rangle\).
For \(a\ge j\) write
\[
 C_a^K(g,y)=\int K(\sigma g,\sigma y)\,d\theta_{[a,n]}(\sigma),
 \qquad C_{n+1}^K=K.
\]
The cut condition and the support of each higher factor in
its stage group force near conditions at the successive stages.
After rebasing each starting point into the box, the orbit
vector is bounded by the projected-column estimate. The
uniform norm-fiber box mass and the finite compact terminal
mass give absolute integrability for bounded \(F\).
In particular \(C_a^K\) is the \(\theta_a\)-integral
of \(C_{a+1}^K(\sigma g,\sigma y)\).

The diagonal Weyl ratio in \(K_2\) normalizes the Bessel
unipotent and commutes with \(G_j\). On the upper
unipotent from stage \(a+1\) on, its simultaneous Bessel
character ratio is one. The quotient integration formula
therefore makes \(C_{a+1}^K\) exactly invariant under
simultaneous movement by \(G_{a+1}\). It is also invariant
under simultaneous \(U_a\), since the Weyl ratio is
the identity on that stage and the character ratio
cancels after conjugation by each higher factor.
Thus \(C_{a+1}^K(\sigma g,\sigma y)\) descends to
\(\sigma^{-1}f_a\). A simultaneous lower move in \(U_{<a}\)
multiplies \(C_{a+1}^K\) by a phase: conjugation by
\(G_{a+1}\) fixes the basis vectors occurring in its
simple-root entries. Consequently the modulus is left
\(R_{a+1}\)-invariant as a function of \(g\) at \(y=gx\).

Define
\[
 \begin{split}
 A_a^K(g)&=W_a(g)|C_{a+1}^K(g,gx)|^2,\\
 B_a^K(g)&=\sum_{l=a}^{\varrho}
        \int A_l^K(\sigma g)\,d\theta_{[a,l]}(\sigma),\\
 E_a^K(x)&=\int B_a^K\,d\theta_{[0,a-1]}.
 \end{split}
\]
Each \(A_a^K\) is left \(R_{a+1}\)-invariant and vanishes
outside \(\mathcal U_{a+1}\). The recursion
\(B_a^K=\int (A_a^K+B_{a+1}^K)(\sigma\,\cdot)\,d\theta_a\),
with \(B_n^K=0\), and the stage quotient right action
show that \(B_a^K\) is left \(R_a\)-invariant and supported
on \(\mathcal U_{a}\). For the untwisted kernel, Tonelli
gives \(E_a^{K_1}\le E_0^{K_1}\) by writing both as sums
of full stage integrals.

Suppose a nonnegative integrand \(f\), with \(z=w_p^{-1}g\),
is bounded pointwise by
\[
 f(g)^2\le K_0 C_L\,B_j^{K_1}(g)I_j(g)B_j^{K_2}(z)I_j(z).
\]
Cross the weights in Cauchy--Schwarz for the integral of
\(f\) against \(M_j\), pairing each \(B\) with the opposite
\(I\). The first factor \(B_j^{K_1}(g)I_j(w_p^{-1}g)\)
integrates to at most \(H_p C_L E_j^{K_1}\) by
\eqref{eq:real-depth-cross} for \(\bar p\).
In the second, substitute \(g=w_p z\) using the invariant
quotient integral just computed. The function
\(B_j^{K_2}(z)I_j(w_p z)\) is left \(G_j\)-invariant.
Its integral is at most \(H_p C_L E_j^{K_2}\).
It follows that
\begin{equation}\label{eq:real-weighted-product}
 \left(\int f\,dM_j\right)^2
 \le K_0\,\mathcal Q_j(x),\qquad
 \mathcal Q_j(x)=H_p^2 C_L^3 E_0^{K_1}(x)E_j^{K_2}(x).
\end{equation}
In particular the character estimate of
Section~\ref{sec:real-orbital-phases} gives this inequality
with \(K_0=(\pi^2/|\xi|^2)^{k-5}\).
Indeed the squared norm of each pivot function is at
most \(B_j^K I_j\) by the \(l=j\) summand and invariance
of \(W_j\) under \(G_j\).

For the unequal-head estimate let \(j=i+1\), and write
\(\varphi_l(g)=|C_j^{K_l}(g,gx)|\) for \(l=1,2\).
On the support of \(\varphi_1(g)\varphi_2(z)\), both \(g,z\in\mathcal U_{j}\), and
\(I_i(w_{p,i}z)=W_i(g)^{-1}\). The last Weyl factor
is in \(G_i\) and preserves the earlier depth and
near conditions. On that support, \(I_j(z)W_j(z)=1\) and
\(I_i(w_{p,i}z)W_i(g)=1\). Thus Cauchy--Schwarz for
\(\int\varphi_1(g)\varphi_2(z)\,dM_j(g)\) uses the two squared
factors
\[
 \begin{split}
 &A_i^{K_1}(g)I_j(z),\\
 &\mathbf1_{\{w_p z\in\mathcal U_{j}\}}\,
                  W_j(z) I_i(w_{p,i}z)\varphi_2(z)^2 .
 \end{split}
\]
The first contributes at most \(H_p C_L E_0^{K_1}\). The second
is left \(G_j\)-invariant as a function of \(z\), so one may
again use the Weyl change for \(M_j\).

Then write
\(M_j=du_i*du_{<i}*\theta_{[0,i]}\) using the lower
conjugation law. The two weights \(W_j(z)\) and
\(I_i(w_{p,i}z)\) are invariant under \(U_i\).
Integrate \(\varphi_2^2\) and the near indicator using the
section estimate of Section~\ref{sec:real-orbital-phases}.
The next squared pivot norm is bounded by \(B_j^{K_2}/W_j\).
This cancels \(W_j\). Apply \eqref{eq:real-depth-cross}
at length \(i\) to the remaining factor
\(B_j^{K_2}(z)I_i(w_{p,i}z)\), with
\(\rho=\theta_{[0,i]}\). The length \(i\) product of moduli
times \(\mathfrak m_i(p_i)^{-1}\) gives \(H_p\).
The second integral is at most \(K_0 H_p C_L^2 E_j^{K_2}\),
with \(K_0=(\pi^2/|\xi|^2)^{k-5}\). Multiplication gives
\eqref{eq:real-weighted-product} for the unequal heads.

\subsection{Decay from a deep unipotent subgroup}
\label{sec:real-depth-majorant}

\sectionopening{In this subsection, we use Cauchy--Schwarz in the next isotropic integration to carry the
estimates of Section~\ref{sec:real-weighted-products} to smaller
supports. Section~\ref{sec:real-depth-shells} identifies the relative
positions where the stronger estimate applies.}

At a starting point \(g\), put
\(P_j^K(g)=\int |C_{j+1}^K(\sigma g,\sigma gx)|\,d\theta_j(\sigma)\).
If the \(\theta_j\)-mass of \(\sigma g\in\mathcal U_{j+1}\)
is at most \(V\ge1\), Cauchy--Schwarz gives
\[
             (P_j^K)^2\le V B_j^K I_j .
\]
Suppose instead that every nonzero terminal contribution on the
\(G_j\)-orbit of \(g\) has some \(D_a\ge\varepsilon^{-1}\)
for \(j\le a<\varrho\), where \(0<\varepsilon\le1\).
Split at the first such stage as the convolution proceeds
from \(j\) upwards.

Once read, a depth is unchanged by
the higher stage groups, which fix its basis data.
At a stage with \(D_a<\varepsilon^{-1}\), use the triangle inequality to expand
the next tail. The part with no deep index has zero
terminal contribution by hypothesis. Dropping these bounds on the earlier stages from the nonnegative contributions
gives
\[
 P_j^K(g)\le \sum_{a=j}^{\varrho-1}
 \int \mathbf1_{\{D_a(\sigma g,x)\ge\varepsilon^{-1}\}}
             |C_{a+1}^K(\sigma g,\sigma gx)|\,d\theta_{[j,a]}(\sigma).
\]
On the support of the \(a\)-th term its squared integrand
is at most \(\varepsilon^s I_{a+1} A_a^K\), since
\(W_a/W_{a+1}=D_a^s\). By Cauchy--Schwarz, iteration of
\eqref{eq:real-depth-first}, and \(a+1-j\le\varrho\),
the integral squared is at most
\[
             \varepsilon^s\Lambda_L^{10^4n\varrho}
                         I_j(g) B_j^K(g).
\]
Since \(\varrho^2\Lambda_L^{10^4n\varrho}\le C_L\),
squaring the sum yields
\[
             (P_j^K)^2\le\varepsilon^s C_L B_j^K I_j .
\]

For the two kernels in Section~\ref{sec:real-weighted-products},
call a point \(y\) active at \(j\) if \(K_1(y,yx)\ne0\)
and \(w_p^{-1}y,w_q^{-1}yx\in\mathcal U_{j}\).
If \(P_j^{K_2}(w_p^{-1}g)\ne0\), some point in its
upper convolution has \(F\) nonzero in both arguments.
The two Weyl elements commute with \(G_j\), so this
supplies the two near conditions on the entire \(G_j\)-orbit
of \(g\). Thus if every active point has a deep stage
of index at least \(j\), with \(D_a\ge\varepsilon^{-1}\),
the deep bound applies to the first pivot. The ordinary
bound applies to the second. The squared product of pivots
then has the pointwise bound used in
\eqref{eq:real-weighted-product} with
\(K_0=V^2\varepsilon^s\).
With no depth hypothesis both ordinary bounds give
\(K_0=V^2\). At \(j=0\), with empty heads, the majorant
is the square of the same pivot at \(g=1\). Applying
the deep bound to both copies at this stage gives the
squared majorant bound \(\varepsilon^{2s}\mathcal Q_0(x)\),
using \(C_L\ge1\).

\subsection{The weighted square integral on a real ball}
\label{sec:real-position-ball}

\sectionopening{In this subsection, we estimate a weighted integral over a bounded set of relative
positions. Plancherel from Section~\ref{sec:real-floor-plancherel}
handles the norm-value coordinates, and
Section~\ref{sec:real-inverse-depth} controls the inverse weights.}

Use \(K\) equal to either of the two pairing kernels, starting
at stage \(j\). Let \(\mathcal B'\) consist of positions \(x\)
with every entry of \(x-1\) of modulus at most
\(L'\ge\max(R_b,8r)\), and put \(\Lambda'=(L')^2\).
The real level unit is one. Set
\[
 c_a=\pi^{2(k-5)}(\Lambda')^{4900n}\quad(a<\varrho),
 \qquad c_\varrho=1 .
\]
Write \(du_{<a}\) for the lower Heisenberg Haar on \(U_{<a}\)
and \(d\eta_a\) for that of \(U_a\), with
\(U_\varrho=\{1\}\) of unit mass. The convolution
\(du_{<a}*d\eta_a\) is \(du_{<a+1}\).

Call \(y\) close to \(g\) modulo \(U_{<a}\) if
\(g^{-1}u y\in\mathcal B'\) for some \(u\in U_{<a}\).
For \(0\le a<\varrho\), suppose \(g\in\mathcal U_{a+1}\) and \(n_0y,n_1y\) are
both close in this sense, with \(n_0,n_1\in U_a\).
Write a box move at \(g\) as \(t u\), with
\(t\in G_a\cap R_{a+1}\), \(u\in U_{<a}\) and \(g_0=tug\)
boxed. For \(x_i=g^{-1}u_i n_i y\in\mathcal B'\),
\[
 (g_0x_0)(g_0x_1)^{-1}
           \in U_{<a}\,t n_0n_1^{-1}t^{-1}.
\]
Here \(G_a\) normalizes \(U_{<a}\).
Both \(x_i\) and their inverses have norms at most
\((L')^3\), so the ratio has norm at most \((L')^8\).
Its \(e_a^*\)-column, projected away from the earlier
pairs, is that of \(t n_0 n_1^{-1}t^{-1}\), since the
lower unipotent is trivial on \(A_a^\perp/A_a\).
Subtracting the identity column gives a bound
\((L')^9\). In Heisenberg coordinates \((v,i\beta)\)
the column has vector part \(v\) and coefficient
\(i\beta-h(v,v)/2\) at \(e_a\). The vector part has norm at most \((L')^9\).
Taking the imaginary part of the coefficient also gives
\(|\beta|\le(L')^9\), since \(h(v,v)\) is real.

We first bound the integral of the forward weight.
Take \(g\) boxed and \(g^{-1}y\in\mathcal B'\).
Let \(u_a\in U_a\) have vector \(v\), with \(u_a y\) close
modulo \(U_{<a}\). The preceding bound applies with
\(n_0=1\) and \(t=1\). If \(c=h(e_a,z)\) is a stable
product at \(g,u_a y\) with lifts \(z,z'\), then for a
common short vector \(w\) at \(g,y\) in the free positive
space,
\[
 h(g^{-1}z,g^{-1}w)-h((u_a y)^{-1}z',y^{-1}w)
                            =\bar c\,h(v,w).
\]
Indeed \(z-z'\in A_a\) and \(u_a w=w-h(v,w)e_a\).

The left side has modulus at most one by
Cauchy--Schwarz on the two close pairs, since
\(\operatorname{diag}(J,-J)\) is a Euclidean isometry.
Test over the half-unit ball of the common positive plane
in Section~\ref{sec:real-inverse-depth}. The area of the
possible \(c\) is at most \(4\pi/\|O^*v_I\|^2\), bounded
by \(\phi(v_I)\) of that section. Thus the forward factor
\(D_a(g,g^{-1}u_a y)^{-s}\) is at most \(\phi(v_I)^s\).
Integrating this over the size-nine box contributes at most
\((\Lambda')^{990n}\) on the free coordinates, by
\eqref{eq:real-free-depth}. The remaining at most \(2n+2\)
complex coordinates and the real center contribute at most
\((L')^{19(2n+2)+10}\). The product is at most \(c_a\).

For \(0\le l\le\varrho\) and \(g\in\mathcal U_{l}\) the resulting bound is
\begin{equation}\label{eq:real-forward-weight}
 \int_{U_{<l}}\mathbf1_{\mathcal B'}(g^{-1}u y)
              W_l(g;g^{-1}u y)\,du
                            \le\prod_{a<l}c_a ,
\end{equation}
where \(W_l(g;x)\) displays the position variable.
To prove this first take \(g\) boxed and induct on \(l\).
At step \(a+1\), if the integrand has nonempty support,
translate \(y\) by an element of \(U_{<a+1}\) to make
\(g^{-1}y\in\mathcal B'\), using right Haar.
Integrate \(u u_a\) with \(u\in U_{<a}\) and \(u_a\in U_a\).
The depth \(D_a\) is independent of this \(u\) in the
second point. Moving \(u_a y\) to \(u u_a y\) transports its lift
by \(u\), preserving the quotient class and the length of its transported lift,
with the first lift fixed.

Applying \(u^{-1}\) gives the inverse
bijection. Apply the induction bound inside and the
one-stage forward bound outside, on \(u_a y\) close modulo
\(U_{<a}\). For general \(g\in\mathcal U_{l}\), write its
box move as \(\gamma u_0\), with \(\gamma\in G_l\)
and \(u_0\in U_{<l}\).
Change the pair to \((\gamma u_0g,\gamma y)\) and use
\(\gamma u_0 u\gamma^{-1}\) as the Haar variable.
The position is preserved and the weights agree by their
simultaneous flag invariance.

For the complementary ascent, put
\[
 m_a(g,y)=\int_{U_{<a}}\mathbf1_{\mathcal B'}(g^{-1}u y)\,du,
 \qquad Q_a(g,y)=\mathbf1_{\{m_a(g,y)>0\}}
                                      |C_{a+1}^K(g,y)|^2 .
\]
We show for \(j\le a\le\varrho\) that
\begin{equation}\label{eq:real-square-ascent}
 \int Q_a(g,n y)\,d\eta_a(n)
         \le c_a \int Q_{a+1}(\sigma g,\sigma y)\,d\theta_{a+1}(\sigma).
\end{equation}
Convolution gives \(m_{a+1}(g,y)=\int m_a(g,ny)\,d\eta_a(n)\).
If this mass is zero, the left side vanishes. Otherwise the cutoff on the right is one,
since \(\sigma\in G_{a+1}\) normalizes \(U_{<a+1}\) and
preserves its Haar. At \(a=\varrho\) the inequality is
Cauchy--Schwarz on the terminal compact probability measure.
For \(a<\varrho\) and nonzero left side, choose a contributing
\(n_0\), so \(g\in\mathcal U_{a+1}\). Use the centering above,
writing \(t=\gamma_1 r_1\) with \(\gamma_1\in G_{a+1}\),
\(r_1\in U_a\), and substituting \(n'=t n n_0^{-1}t^{-1}\).

This preserves Haar. The support lies in the size-nine box
and the vector of \(n\) is \(\gamma_1^{-1}v'+v_0\).
For \(y_\sigma=\sigma^{-1}f_{a+1}\), moving \(n\) past the
next pivot and using the pairing covariance gives phase
\(\psi(-h(v,y_\sigma))\) in the integral defining
\(C_{a+1}^K(g,ny)\). This is the one-sided version of the
calculation in Section~\ref{sec:real-orbital-phases}:
the diagonal twist commutes with \(G_a\), and only the
second argument is moved. In particular
\[
 C_{a+1}^K(g,ny)=\int
 \psi\bigl(-h(v,y_\sigma)\bigr)
 C_{a+2}^K(\sigma g,\sigma y)\,d\theta_{a+1}(\sigma).
\]
The next correlation descends measurably to \(y_\sigma\), since
it is invariant under the stabilizer and the orbit map has
smooth local sections on the nonzero orbit. Rotate that vector by \(\gamma_1\).
The translation by \(v_0\) multiplies the function by a phase
of modulus one.

Apply the projected-column estimate to
\(q=\gamma_1 r_1 u\) with \(qg=g_0\) in the box and
\(\gamma_1^{-1}q\in U_{<a+1}\). The product of the two
matrices from the box has norm at most \(R_b^2\), so
the coordinates of \(\gamma_1 y_\sigma\) are bounded
by \((L')^2\).
The Fourier estimate of Section~\ref{sec:real-floor-plancherel}
at frequency \(1\) applies from this quadric to the vector
box. Use ordinary area in both target coordinates in the
block \(B'\) of that proof and keep the same source fiber.
Plancherel on the free block still gives the squared-norm
factor \(\pi^{2(k-5)}\).
On each side at most \(2n+2\) complex coordinates lie
outside the free block. On the input their disc areas
are at most \(\pi(L')^6\), and the solved pair has mass at
most \(4\pi^2(L')^6\). On the output the disc areas are
at most \(\pi(L')^{20}\). Since \(L'\ge8\), these volumes
and the central length \(2(L')^9\) contribute at most
\((\Lambda')^{4(2n+2)+4+11(2n+2)+5}\le(\Lambda')^{4900n}\).
This proves \eqref{eq:real-square-ascent}.

Set
\(\rho_a=d\eta_a*\cdots*d\eta_\varrho*\theta_{[0,n]}\),
so \(d\mu=du_{<a}*\rho_a\) and \(\rho_n=\theta_{[0,n]}\).
For \(u\in U_{<a+1}\), right translation of its Haar variable gives
\(m_{a+1}(g,uy)=m_{a+1}(g,y)\). Write \(D=1\) for \(K_1\) and
\(D=w_q^{-1}w_p\) for \(K_2\). In \(C_{a+2}^K(g,uy)\) the higher
factors commute with \(D\) and fix the vectors defining the
simple-root entries of \(D uD^{-1}\). The second Bessel factor
therefore contributes the constant phase \(\psi(\mathfrak s(DuD^{-1}))\).
Its modulus is one, so \(Q_{a+1}(g,uy)=Q_{a+1}(g,y)\).
Conjugation by \(\sigma\in G_{a+1}\) preserves \(du_{<a+1}\).
Writing \(L_\sigma y=\sigma y\), we have
\[
 du_{<a+1}*(L_\sigma)_*\rho_{a+1}
 =(L_\sigma)_*(du_{<a+1}*\rho_{a+1})=(L_\sigma)_*d\mu=d\mu .
\]
The invariant-function averaging identity of
Section~\ref{sec:real-weighted-products} shows that its
integral on \(Q_{a+1}(g,\cdot)\) is unchanged.
Integrate \eqref{eq:real-square-ascent} against
\(\theta_{[0,a]}\) in \(g\) and \(\rho_{a+1}\) in \(y\).
Undoing the shift by \(\sigma\) in the \(y\)-integral gives
\[
 \iint Q_a\,d\theta_{[0,a]}d\rho_a
          \le c_a\iint Q_{a+1}\,d\theta_{[0,a+1]}d\rho_{a+1}.
\]
Iterate to \(n\), where \(Q_n\le |K|^2\).

Finally, to integrate \(A_l^K(g)\) over \(x\in\mathcal B'\)
put \(x=g^{-1}u y\) with \(u\in U_{<l}\) and \(y\) in the
\(\rho_l\) measure. The correlation modulus is invariant
under this motion in its second argument. If the close mass
is zero the inner integral vanishes. Otherwise the forward
estimate \eqref{eq:real-forward-weight} contributes \(\prod_{a<l}c_a\)
and leaves \(Q_l(g,y)\). The ascent for its double integral
contributes the complementary product. Summing the at most \(n\)
weighted stage integrals in \(E_j^K\) proves
\begin{equation}\label{eq:real-position-energy}
 \int_{\mathcal B'} E_j^K(x)\,d\mu(x)
      \le n\biggl(\prod_{a<\varrho}c_a\biggr)
               \iint |K(g,y)|^2\,d\theta_{[0,n]}(g)d\theta_{[0,n]}(y).
\end{equation}

\subsection{Density and twisted square integrals}
\label{sec:real-head-density}

\sectionopening{In this subsection, we compute the density in the coordinates integrated at the earlier
steps of the isotropic flag. This combines the integral of
Section~\ref{sec:real-position-ball} with the terminal pairing of
Section~\ref{sec:real-terminal-gram}.}

Write \(dQ=d\theta_{[0,n]}\) for the original real quotient measure and \(d_a=\dim_{\mathbb C}V_{a+1}=k+1+2(\varrho-a-1)\).
At a lower stage \(a<\varrho\), write \(n_a(v,\beta)\) for the
unipotent with \(e_a^*\)-column
\(e_a^*+v+(i\beta-h(v,v)/2)e_a\). Its inverse is \(n_a(-v,-\beta)\).
For \(\sigma=w_a(c)^{-1}n_a(v,\beta)\), the orbit vector is
\[
 \sigma^{-1}e_a=c e_a^*-cv+c(-i\beta-h(v,v)/2)e_a .
\]
In general write the vector as \(c e_a^*+w+c t e_a\), \(c\ne0\).
Ordinary area in its \(e_a\)-coordinate is \(|c|^2d^2t\), and its
norm-value coordinate is \(|c|^2\Re t+h(w,w)/2\). Disintegration
therefore leaves \(d^2c\,dw\,d\Im t\) on the fixed norm fibre.
Here \(w=-cv\) and \(\Im t=-\beta\), so the stage measure is
\[
           |c|^{2d_a}\,d^2c\,dv\,d\beta .
\]
The omitted locus \(c=0\) has fibre measure zero by the positive
two-coordinate chart of Section~\ref{sec:real-floor-plancherel}.
The factor \(w_a(c)^{-1}\) commutes with the higher
stages, so it can be written at the left of the convolution
with the other factors in their original order.

Let \(D_j=w_q^{-1}w_p\), and let \(D_a\) denote the same
diagonal through stage \(a-1\) in this calculation. On the
\(a\)-th pair, \(D_{a+1}D_a^{-1}\) scales \(e_a\) by
\(z_a=p_a/q_a\) and \(e_a^*\) by \(\bar z_a^{-1}\).
Direct multiplication of the two-by-two matrices gives
\[
              D_{a+1}w_a(z_a c)^{-1}=D_a w_a(c)^{-1}.
\]
Thus replacing \(c\) by \(z_a c\) in the quotient integral
contributes \(|z_a|^{2(d_a+1)}\), from its weight and the complex
area Jacobian. Iterating, for a nonnegative measurable test
\(\Phi\) we have
\[
 \int\Phi(D_j y)\,dQ(y)=J_{p,q}\int\Phi(y)\,dQ(y),\qquad
 J_{p,q}=\prod_{a<j}|p_a/q_a|^{2(d_a+1)}.
\]
In particular the double squared norm of \(K_2\) is
\(J_{p,q}\) times that of \(K_1\), since only its second
variable is translated.

The head density is
\(\mathcal D_j(p,q)=\prod_{a<j}|p_aq_a|^{d_a}\).
With the inverse section moduli \(H_p\) of
Section~\ref{sec:real-inverse-depth}, the exact identity is
\[
 \mathcal D_j(p,q)^2 H_p^2 J_{p,q}
 =\prod_{a<j}|q_a|^{-2}\le L^2
\]
for a \(q\)-head adapted to \(L\). Indeed the exponent of
\(|p_a|\) is \(2d_a-2(2d_a+1)+2(d_a+1)=0\),
and the exponent of \(|q_a|\) is \(2d_a-2(d_a+1)=-2\).
For any measurable \(\mathcal S\subseteq\mathcal B'\)
on which \(f(x)^2\le e^2\mathcal Q_j(x)\), use
Cauchy--Schwarz in \(x\) and then \eqref{eq:real-position-energy}
for the two kernels. This gives
\[
 \begin{split}
 \mathcal D_j(p,q)\int_{\mathcal S}|f(x)|\,d\mu(x)
 &\le e\, n L C_L^{3/2}
           \biggl(\prod_{a<\varrho}c_a\biggr)\mathcal N,\\
 \mathcal N&=\iint |K_1(g,y)|^2\,dQ(g)\,dQ(y).
 \end{split}
\]
In particular, the loss multiplying \(e\mathcal N\) is at
most
\(\pi^{2(k-5)\varrho}\Lambda_L^{1500000n^3+1}(\Lambda')^{25000n^2}\).

For the cut function on the real quotient, write
\(d\mu=\lambda\,dg\), where \(dg\) is the orbital Haar.
The corresponding quotient measures differ by the same
\(\lambda\), since the unipotent Haar is fixed. If \(F=\chi W\)
is the scalar cutoff with \(|\chi|\le1\) and \(W\) has
squared quotient norm \(d_I\) in the orbital normalization,
\(\mathcal N\le\lambda^2d_I\).
Put \(U=U_{<n}\). The compact group \(H=U(z^\perp)\) in
the remaining space normalizes \(U\). The modulus of its action
on Haar measure \(du\) is a continuous homomorphism
\(H\to\mathbb R_{>0}\), hence is one. Disintegrating left Haar
measure \(dg\) over \(U\) shows that \(y\mapsto hy\) preserves
\(d\bar g\) on \(U\backslash G\) for \(h\in H\).
By covariance the matrix \(\int W(y)W(y)^*\,d\bar g(y)\)
commutes with the irreducible action of \(H\). Its trace is
\(d_I\), so Schur's lemma makes it the identity. The double integral for
\(W\) is the trace of its square. Multiplication by the
scalar cutoffs can only decrease the integrand's modulus.

In the terminal comparison the top norm density in the
last entry \(c\) is
\[
       \frac{2\pi^k}{(k-1)!} (|c|^2-1)_+^{k-1}\,dc .
\]
This follows by polar integration in the \(k\) positive
coordinates, at fixed norm-value coordinate \(h/2\).
The coefficient of the product of the top correlations
in basis measure is therefore
\[
 C_F\,\mathcal D_\varrho(p,q)
            \frac{2\pi^k}{(k-1)!}(|c|^2-1)_+^{k-1},\qquad
 C_F=\frac{\operatorname{vol}U(k)}
             {\operatorname{vol}U(r)\,(\lambda d_I)^2}.
\]
In the orbital normalization, the two Haar integrations
contribute \(\lambda^2\) to this coefficient.

Thus the coefficient is \(A_p A_q\) of
Section~\ref{sec:real-terminal-gram}, with the common last
entry and the constant
\[
 c_G=\frac{\operatorname{vol}U(k)}{\operatorname{vol}U(r)}
                       \frac{2\pi^k}{(k-1)!}.
\]
To see the matching integrand, at \(\sigma=a(c)\) the
terminal compact average fixes the two correlations by
unitary covariance and because the ratio of heads commutes
with the compact group. Their product at \(y,x\) is
\[
 \langle F(y),F(yx)\rangle\,
 \big\langle F(\sigma^{-1} w_q^{-1}yx),
                  F(\sigma^{-1}w_p^{-1}y)\big\rangle .
\]
The substitution \((y,x)\mapsto(y,yx)\) preserves product
Haar and makes this the integrand of the Gram pairing of
the two orbital tensors. This comparison is at fixed entry.
Absolute convergence over the varying entry and all real
positions uses the shell estimates.

\subsection{Real determinant weights on dyadic regions}
\label{sec:real-depth-shells}

\sectionopening{In this subsection, we use minima of Hermitian quotient norms to decompose the relative
positions into dyadic sets. This identifies where the stronger
estimates of Section~\ref{sec:real-depth-majorant} enter the frequency
sum of Section~\ref{sec:real-shell-sum}.}

On \(V_a=A_a^\perp/A_a\) let \(|v|_{g,a}\) be the
infimum of \(\|g^{-1}\widetilde v\|\) over lifts.
For \(a<\varrho\) put
\[
 b_a(g,g')=\inf_{h(e_a,v)=1}
                   (|v|_{g,a}^2+|v|_{g',a}^2).
\]
The infimum is positive since \(h(e_a,\cdot)\) is a
nonzero continuous linear functional. The stable product
\(c\ne0\) occurs precisely when \(|c|^2b_a(g,g')<1\),
by scaling the common class and taking the two lift
infima independently. Thus the stable set is a disc
and \(D_a(g,x)=b_a(g,gx)/\pi\).

Fix heads \(p,q\) satisfying the coordinate bound \(L\) at \(j\),
and a point \(g\) active at \(x\) as in
Section~\ref{sec:real-depth-majorant}.
The two nonzero cut values give boxed points \(u g,v gx\)
with \(u,v\) in the Bessel unipotent.
They preserve \(e_a\) modulo \(A_a\), so
\(|e_a|_{g,a},|e_a|_{gx,a}\le R_b\).
On \(V_\varrho\) the Bessel unipotent acts trivially.
The two quotient norms there have ratio at most
\(R_b^2\), by the two boxed inverses and matrices.
For \(a<j\), \(e_a^*\) supplies a short dual in both
norms. Indeed the column-minor bound of
Section~\ref{sec:real-orbital-phases}, for the box
moves \(s_0\) at \(g\) and \(s_1\) at \(w_p^{-1}g\),
bounds each entry of \(s_0e_a^*\) by \(L R_b^{2(a+1)}\).
Applying \((s_0g)^{-1}\) gives
\(\|g^{-1}e_a^*\|\le L^{2n+3}\).
Use the head \(q\) and the box moves at \(gx\) and
\(w_q^{-1}gx\) for the other norm.

Suppose \(T\ge L^{11n}\) and \(b_a(g,gx)<T^2\) for
every \(j\le a<\varrho\). At each such index choose a
dual of \(e_a\) with both lengths at most \(T\), and at
the earlier indices use \(e_a^*\). The inequality
\(|h(v,w)|\le |v||w|\) holds in either quotient norm,
by applying it on all lifts before minimizing.
For the chosen dual \(\zeta\) at \(a\), put
\(\zeta_0=\zeta-h(\zeta,\zeta)e_a/2\). This is an
isotropic dual and has both lengths at most
\(T+T^2R_b/2\le T^4\).

The orthogonal complement of \(\langle e_a,\zeta_0\rangle\)
gives a common section of \(V_{a+1}\) into \(V_a\).
For either quotient norm \(|\cdot|_{*,a}\), take a minimizing
lift \(w\in e_a^\perp\subset V_a\) of a class in \(V_{a+1}\).
The common section sends the class to \(w-h(\zeta_0,w)e_a\).
Since \(|e_a|_{*,a}|\zeta_0|_{*,a}\le T^5\), its norm is at
most \((1+T^5)|w|_{*,a}\). The successive quotient infima
equal the infimum over \(A_{a+1}\).
If the second norm on \(V_{a+1}\) is at most
\(T^u\) times the first, write
\(w=h(\zeta_0,w)e_a+h(e_a,w)\zeta_0+w_0\).
The first two terms contribute at most \(2T^5|w|_{g,a}\)
in the second norm.

Removing \(h(e_a,w)\zeta_0\)
before taking the quotient shows that the class of
\(w_0\) has first quotient norm at most
\((1+T^5)|w|_{g,a}\).
Thus the second norm of \(w\) is at most
\[
 \bigl(2T^5+(1+T^5)^2T^u\bigr)|w|_{g,a}
                                    \le T^{u+13}|w|_{g,a}.
\]
Descend from \(u=2\) at \(V_\varrho\), and at \(a=0\)
apply the result to \(w=g v\) for arbitrary \(v\).
It gives \(\|x^{-1}\|\le T^{2+13\varrho}\).
Since \(J=J^{-1}\), \(\|x\|=\|x^{-1}\|\).
Every entry of \(x-1\) is therefore at most \(T^{99n}\)
in modulus. Outside this ball every active point
has \(D_a(g,x)\ge T^2/\pi\) at some \(a\ge j\).
For \(j=\varrho\) there are no active points there.

The uniform next-stage mass is bounded by
\[
                        V_b=(\pi R_b^{24})^k R_b^{9999n}.
\]
Rebase the starting point to the box by \(R_j\), as in
the first inverse-depth integral. On \(\mathcal U_{j+1}\) the orbit vector has norm at most \(R_b^2\), by
projecting the column of a boxed matrix times a boxed
inverse. By the two-coordinate sphere bound of
Section~\ref{sec:real-floor-plancherel}, its box mass
is at most \(4\pi^2R_b^4(\pi R_b^4)^{r-2j-2}\le V_b\).
The same estimate at the top uses two positive coordinates.

Let \(H=2^i\), with \(i\ge0\) least such that the
entries of \(x-1\) have modulus at most \(H\).
One position ball with coordinate bound \(LH\) covers this shell.
If \(0<\tau<1\) and
\(\Lambda_L^{1200n^2}\le \tau^{400n^2}H\), take
\(T=L^{11n+1}/\tau\). Then \(T^{-1}\le\tau/2\) and
\[
                   2T^{99n}\le L^{2400n^2}\tau^{-400n^2}\le H.
\]
Minimality of the shell index puts \(x\) outside
the position ball of radius \(T^{99n}\).
Thus the deep-stage bound applies with \(D_a\ge\tau^{-2}\).

\subsection{Summation and real frequency decay}
\label{sec:real-shell-sum}

\sectionopening{In this subsection, we sum the weighted real integrals on the dyadic sets of
Section~\ref{sec:real-depth-shells}. The resulting Fourier bound
controls the nonconstant terms in
Section~\ref{sec:real-common-descent}.}

Let \(\mathcal M_j(x)\) be the majorant obtained by
integrating the two pivot moduli against \(M_j\).
Take any measurable \(0\le f\le\mathcal M_j\) and
\(t\ge1\), with
\[
 f(x)^2\le V_b^2t^{-2s}\mathcal Q_j(x)\qquad(t>1).
\]
For \(t=1\) use the bound for \(\mathcal M_j\).
This also gives the estimate for the modulus of the
form integrand.
The energy estimate of Section~\ref{sec:real-head-density}
on the ball with coordinate bound \(LH\) shows that the one-shell
integral times \(\mathcal D_j(p,q)\) contributes at most
\[
  \Lambda_L^{2\cdot10^6 n^3} H^{10^5n^2}
                          \pi^{2(k-5)\varrho}\mathcal N
\]
times the square-root saving.
Write \(u=t^{-1}\).

For \(H\le\Lambda_L^{2400n^2}\)
use saving \(V_b u^s\). On larger shells take
\(\tau=H^{-1/(800n^2)}\), which satisfies the condition
of Section~\ref{sec:real-depth-shells}. The saving there
is \(V_b\min(u,\tau)^s\). Indeed the deep pivot bound with
\(\varepsilon=\tau^2\) gives
\(\mathcal M_j^2\le V_b^2\tau^{2s}\mathcal Q_j\).
Taking the smaller of this bound and the given frequency bound
for \(f^2\), then taking the square root, gives the stated factor.
For \(g_*=u\) in the first case and \(g_*=\min(u,\tau)\)
in the second,
\[
  g_*^{9\cdot10^7 n^4}H^{1+10^5n^2}
                     \le\Lambda_L^{29\cdot10^7n^4}.
\]
For the smaller shells use
\(2400n^2(1+10^5n^2)\le29\cdot10^7 n^4\). On a larger shell
\(g_*\le H^{-1/(800n^2)}\), and hence
\[
 H^{10^5n^2}g_*^{9\cdot10^7 n^4}\le H^{-12500n^2}.
\]
For the first index \(i_0\) with \(2^{i_0}>\Lambda_L^{2400n^2}\),
the sum of these factors over \(i\ge i_0\) is at most
\(2^{-12500n^2i_0}/(1-2^{-12500n^2})\), uniformly in the two
heads satisfying the bounds for \(L\).
The rank bounds give \(s\ge10^{10}n^4\).
Since \(g_*\le t^{-1}\) and \(t\ge1\), its remaining power satisfies
\[
 g_*^{s-9\cdot10^7n^4}\le t^{9\cdot10^7n^4-s}
                           \le t^{10^8n^4-s}.
\]
Summing \(\sum_{i\ge0}2^{-i}=2\) gives
\[
 \mathcal D_j(p,q)\int f\,d\mu
 \le V_b\,t^{10^8n^4-s}\Lambda_L^{4\cdot10^8n^4}
                  \pi^{2(k-5)\varrho}\mathcal N .
\]
This in particular proves integrability of the majorant
at \(t=1\) for each pair of heads. By Section~\ref{sec:real-head-density},
\[
 d_I C_F\mathcal N
 =\frac{\operatorname{vol}U(k)}{\operatorname{vol}U(r)}
    \frac{\mathcal N}{\lambda^2 d_I}
 \le \frac{\operatorname{vol}U(k)}{\operatorname{vol}U(r)}.
\]
Inserting \(V_b=(\pi R_b^{24})^kR_b^{9999n}\) gives
\begin{equation}\label{eq:real-shell-total}
 \begin{aligned}
 d_I C_F\,\mathcal D_j(p,q)\int f(x)\,d\mu(x)
 &\le \Theta L^{10^9n^4} R_b^{24k}t^{10^8n^4-s},\\
 \Theta&=\frac{\operatorname{vol}U(k)}{\operatorname{vol}U(r)}
                                    \pi^{(2\varrho+1)k}.
 \end{aligned}
\end{equation}
Here \(k+2(k-5)\varrho\le(2\varrho+1)k\) and
\(R_b^{9999n}L^{8\cdot10^8n^4}\le L^{10^9n^4}\).

\begin{lemma}[Real frequency decay]\label{lem:marked-real-frequency-decay}
For a pivot tuple \(p\), set \(P_h(p)=\prod_{a<h}p_a\) and
\[
 w_j(p)=\mathbf1_{\max_{1\le h\le j}|P_h(p)|\le b^{60n^2}}
             \prod_{h=1}^j\max(1,|P_h(p)|^{-1})^{n^6},
             \qquad K_*=k-n^7.
\]
For \(j=0\) use \(w_0=1\). Fix \(0\le j\le\varrho\) and
\(p,q\in(\mathbb C^\times)^j\). Put \(Z_j(g)=h(f_j,gf_j)\) and
write \(\mathcal I_j(p,q;\omega)\) for the integral with weight
\(\omega\) in Section~\ref{sec:real-common-descent}.
With \(C_F\) as in Section~\ref{sec:real-head-density} and
\(\Theta\) as in \eqref{eq:real-shell-total}, one has
\begin{equation}\label{eq:real-frequency-final}
 d_I\bigl|C_F\mathcal I_j(p,q;\omega)\bigr|
 \le \Theta b^{2000n^9}|\xi|^{-K_*}\min(1,|\xi|^{-1})^6
             w_j(p)w_j(q)
\end{equation}
for either of the following choices of weight and frequency:
\[
 \begin{aligned}
 &\omega(g)=\psi(\xi Z_j(g)), &&\xi\in\mathbb C^\times,\\
 &\omega(g)=1, &&j=i+1,\quad \xi=p_i^{-1}-q_i^{-1}\ne0.
 \end{aligned}
\]
\end{lemma}

\begin{proof}
Set \(t=\max(1,R_b^{24}|\xi|)\) at frequency \(\xi\ne0\).
If \(t>1\), then
\[
            (\pi/|\xi|)^{k-5}
                       \le(\pi R_b^{24})^k t^{-s},
\]
which bounds the frequency factor in either case. Let \(f(x)\) be
the integral against \(dM_j(g)\) of the absolute value of the paired
pivot integral at \((x,g)\). For the character weight, the Fourier
estimate on the two quadrics in Section~\ref{sec:real-orbital-phases}
and \eqref{eq:real-weighted-product} give
\[
 f(x)^2\le(\pi/|\xi|)^{2(k-5)}\mathcal Q_j(x).
\]
For the constant weight the two pivots separate. At distinct last
entries the section estimate at \(i\) gives the same inequality by
the crossed integrations in Section~\ref{sec:real-weighted-products}.
The triangle inequality gives \(f\le\mathcal M_j\) for both weights.
Where both pivot moduli contribute, their nonzero
terminal contributions imply \(g,w_p^{-1}g\in\mathcal U_{j}\),
and likewise \(gx,w_q^{-1}gx\).

The leading minors then bound the partial products by
\(R_b^{2j}\), so a head violating the threshold gives
zero. For the heads satisfying the threshold, take \(L\) to be the
largest of \(R_b\) and the inverse initial segment moduli
for \(p,q\). It meets the box hypothesis since
\(8r\le 8nr\le b\le R_b\) by Section~\ref{sec:real-marked-filling}. Thus \(L^{n^6}\le R_b^{n^6}w_j(p)w_j(q)\).
The parameter inequalities \(b\ge2,\ n\ge10^8,\ k\ge n^7\)
give \(10^9n^4\le n^6\) and \(10^8n^4+K_*+6\le s\).
Use \(R_b^{24}/t\le|\xi|^{-1}\) and
\(t^{-1}\le\min(1,|\xi|^{-1})\) in \eqref{eq:real-shell-total}.
Since \(\mathcal I_j\) includes the density \(\mathcal D_j(p,q)\) and
\(R_b^{n^6+24n^7}\le b^{2000n^9}\), we obtain
\begin{align*}
 d_I|C_F\mathcal I_j(p,q;\omega)|
 &\le d_I C_F\mathcal D_j(p,q)\int f\,d\mu\\
 &\le \Theta b^{2000n^9}|\xi|^{-K_*}
            \min(1,|\xi|^{-1})^6 w_j(p)w_j(q).\qedhere
\end{align*}
\end{proof}

\section{Formal degrees and real Bessel functions}
\label{sec:bessel-normalization-analysis}

\sectionopening{In this section, we pass from the weighted real bounds of
Section~\ref{sec:real-orbital-analysis} to the discrete-series normalization.
The stationary point, coefficient estimates, and truncation bound apply to
the Poincar\'e map of Section~\ref{sec:period}.}

\subsection{The local mass in the formal-degree normalization}
\label{sec:real-bottom-mass}

\sectionopening{In this subsection, we compare the initial orbital pairing for the truncated and full
discrete-series vectors. The formal degree in
Section~\ref{sec:real-reproducing-coefficient} gives the full norm,
and the truncation estimate of Section~\ref{sec:real-cut-error}
controls the comparison.}

For \(d\) the formal degree, let \(dg\) be the
orbital Haar, \(d\bar g\) its quotient Haar and \(d\mu=\lambda\,dg\).
For Bessel functions of the prescribed type define
\[
 P(A,B)(x)=\int_{U_{<n}\backslash G}
                   \langle A(y),B(yx)\rangle\,d\bar g(y).
\]
The original real vector satisfies
\[
 \begin{gathered}
 d>0,\quad f\ge4000n,\quad
 \|W\|_2^2=d_I,\quad \|F\|_2\le\sqrt d_I,\quad
 \|F-W\|_2\le b^{-r^3f}\sqrt d_I,\\
 \|P(W,W)\|_{L^2(dg)}^2\le d_I^2/d,\qquad
 \lambda=\frac{\operatorname{vol}U(r)}{\operatorname{vol}U(k)},
                  \qquad d\Theta\le b^{4r^2},\\
 d\operatorname{vol}_{dg}\{x:\|x\|\le Y\}
                     \le Y^{9r^2}b^{99r^2}\quad(Y\ge1).
 \end{gathered}
\]
Section~\ref{sec:real-frame-haar} proves the Haar comparison, and
Section~\ref{sec:real-marked-filling} gives \(f\ge4000n\) for the chosen
datum. For empty heads, the two pivot correlations with weight \(1\)
give \(\lambda^2|P(F,F)(x)|^2\), with head density one.
By Cauchy--Schwarz and the isometry of right translation,
\[
 |P(F,F)(x)-P(W,W)(x)|
             \le\|F-W\|_2(\|F\|_2+\|W\|_2)\le2d_I b^{-r^3f}.
\]
Set \(t_*=b^{6r}\), \(T=t_*R_b^{11n}\), and take the
position ball \(\mathcal B_*\) of entry radius \(T^{99n}\).
On it \(\|x\|\le(1+r^2)T^{99n}\le R_b^2T^{99n}
\le b^{700nr}\), since \(R_b\le b^{7n}\) and
\(14n+7623 n^3+594nr\le700nr\).
The squared \(L^2(dg)\) error on the ball is therefore at most
\[
 \frac{d_I^2}{d}\,4b^{-2r^3f+6300nr^3+99r^2}
                              \le b^{-24}\frac{d_I^2}{d}.
\]
As \(C_F\lambda^3=d_I^{-2}\), the \(L^2\) triangle
inequality bounds the normalized mass on the ball
by \((1+b^{-12})^2/d\).

On the complement, Section~\ref{sec:real-depth-shells}
gives, at each active point, some \(0\le a<\varrho\) with
\(D_a\ge T^2/\pi\ge t_*^2\).
The deep pivot estimate with \(\varepsilon=t_*^{-2}\)
applies to
\(f=\mathbf1_{\mathcal B_*^c}\mathcal M_0\) with
\(t=t_*\) and \(L=R_b\). The empty heads have density one.
The absolute contribution is at most \(C_F\int f\,d\mu\),
which is bounded by \eqref{eq:real-shell-total} since
\(d_I\ge1\). It is therefore at most
\[
          \Theta R_b^{10^9n^4+24k}t_*^{10^8n^4-s}\le b^{-12}/d.
\]
Indeed \(k=s+3n+6,\ r=s+5n+5,\ s\ge10^{10}n^4\)
and \(n\ge10^8\) give \(r\le(11/10)s\), \(k\le r\) and
\[
 4r^2+7n(10^9n^4+24k)+6r\cdot10^8n^4+12\le6rs.
\]
Here \(4r^2\le(22/5)rs\), both \(7\cdot10^9n^5\) and
\(168nk\) are at most \(rs/10\), and the penultimate
term is at most \((3/50)rs\).
Thus the normalized form at the bottom satisfies
\[
                 |J_0|\le\frac{(1+b^{-12})^2+b^{-12}}d
                             \le\frac{1+b^{-10}}d .
\]

\subsection{Descent through a common entry}
\label{sec:real-common-descent}

\sectionopening{In this subsection, we relate consecutive pairings when the two vectors at the current
isotropic entry coincide. Section~\ref{sec:real-shell-sum} bounds the
other Fourier terms, while Section~\ref{sec:real-terminal-gram}
supplies the terminal positive form.}

For \(a\le\varrho\), let \(p,q\) have length \(a\), and
use the kernels \(K_1,K_2\) of
Section~\ref{sec:real-weighted-products} with twist
\(w_q^{-1}w_p\). Write \(\mathcal I_a(p,q;\omega)\) for
the integral against \(d\mu(x)dM_a(g)\), with factor
\(\mathcal D_a(p,q)\), of
\[
 \iint \omega(\sigma\tau^{-1})\,
 C_{a+1}^{K_1}(\sigma g,\sigma gx)\,
 \overline{C_{a+1}^{K_2}(\tau w_p^{-1}g,\tau w_p^{-1}gx)}
 \,d\theta_a(\tau)d\theta_a(\sigma).
\]
This is absolutely convergent for measurable \(|\omega|\le1\),
by integrating the ordinary majorant at \(t=1\) in
Section~\ref{sec:real-shell-sum}. For weight \(1\) the
two pivot integrals separate into \(C_a^{K_1}\) and
\(\overline{C_a^{K_2}}\) at their respective starting points.

Set \(Z_a(s)=h(f_a,s f_a)\). Let \(v(c)=w_a(c)\) for
\(a<\varrho\) and \(v(c)=a(c)\) of
Section~\ref{sec:real-terminal-gram} at \(a=\varrho\).
The law \(\theta_a\) is the pushforward by
\((c,\nu)\mapsto v(c)^{-1}\nu\) of
\(\rho_a(c)\,dc\,d\zeta_a(\nu)\), where
\[
 \begin{array}{ll}
 \rho_a(c)=|c|^{2d_a},\quad \zeta_a=\eta_a,
       & a<\varrho,\\[2pt]
 \rho_\varrho(c)=\dfrac{2\pi^k}{(k-1)!}(|c|^2-1)_+^{k-1},
       \quad \zeta_\varrho=\theta_n, & a=\varrho .
 \end{array}
\]
These are the lower chart and top sphere densities in
Section~\ref{sec:real-head-density}. Both \(\zeta_a\) are
carried by \(R_{a+1}\cap G_a\), the stabilizer of \(f_a\)
in \(G_a\). The entry \(Z_a(v(c))=c\) wherever the
density is nonzero, except possibly at \(c=0\), a null set.

Let \(\omega=\phi\circ Z_a\) with \(|\phi|\le1\).
For fixed \(x,g,\sigma\), the weighted product in the
double integral is unchanged under left multiplication
of \(\tau\) by the stabilizer of \(f_a\). The entry is
unchanged because that stabilizer fixes \(f_a\), and the
complex tail correlation is unchanged by the stabilizer
invariance of Section~\ref{sec:real-weighted-products}.
Thus the integral in \(\tau\) depends only on
\(\tau^{-1}f_a\). Right multiplication of \(\tau\) by
\(\sigma\) sends the orbit point to \(\sigma^{-1}\tau^{-1}f_a\),
which preserves the norm-fiber measure.

This applies to the product and its absolute value.
We may consequently substitute
\(\tau=v(c)^{-1}\nu\sigma\) with the above chart law.
For \(g'=\nu\sigma g\), the entry is
\[
 Z_a(\sigma\tau^{-1})=h(f_a,\nu^{-1}v(c)f_a)
        =h(\nu f_a,v(c)f_a)=h(f_a,v(c)f_a)=c .
\]
At the last step the \(z\)-coordinate of \(a(c)z\) is
\(-c\) and \(h(z,z)=-1\), giving the same sign.
The first tail at \(\sigma g\)
equals that at \(g'\), by the stabilizer invariance.
The other starting point becomes
\[
            \tau w_p^{-1}g
                    =v(c)^{-1}w_p^{-1}g',
\]
since \(w_p\) commutes with \(G_a\). The same relation holds
with \(gx,g'x\) in the second variable. The variable \(g'\) now has measure \(\zeta_a*\theta_a*M_a\).
The absolute integrability just proved transfers through
this substitution, so Fubini applies also to the new
integral in \(c\).

For \(a<\varrho\), append the common entry \(c\ne0\)
to both heads, writing \(p_c,q_c\). Their new Weyl
elements are \(w_pv(c),w_qv(c)\), so the twist is
unchanged. The next lower measure is
\(M_{a+1}=\eta_a*\theta_a*M_a\), by rearranging the
lower Heisenberg and stage products using the same
conjugation invariance. Also
\(\mathcal D_{a+1}(p_c,q_c)=\mathcal D_a(p,q)|c|^{2d_a}\).
Thus the new inner integral is the next form with weight
\(1\), by separation of its pivots. We obtain
\[
 \int_{\mathbb C}\phi(c)\mathcal I_{a+1}(p_c,q_c;1)\,dc
                       =\mathcal I_a(p,q;\phi\circ Z_a).
\]
The values assigned to the appended head at \(c=0\) do
not affect this statement.

At \(a=\varrho\), \(\theta_n*\theta_\varrho*M_\varrho=d\mu\)
by the full stage decomposition, since the last two
Heisenberg measures are trivial. The compact tail average
is the probability average on \(H\).
The coefficient and integrand comparison in
Section~\ref{sec:real-head-density} now identifies the new integral, after multiplication by \(C_F\), with
\(J_n(p_c,q_c)\). Hence
\[
 \int_{\mathbb C}\phi(c)J_n(p_c,q_c)\,dc
                    =C_F\mathcal I_\varrho(p,q;\phi\circ Z_\varrho).
\]
In particular its absolute convergence follows from the
ordinary majorant, including at zero frequency.

Define \(J_a(p,q)=C_F\mathcal I_a(p,q;1)\) for
\(a\le\varrho\), with \(J_n\) the terminal Gram form.
The preceding identities show that
\[
 \int_{\mathbb C}\psi(\xi c) J_{a+1}(p_c,q_c)\,dc
             = C_F\mathcal I_a(p,q;\psi(\xi Z_a)),\qquad a<n .
\]
At \(\xi=0\) this integral equals \(J_a(p,q)\). At
\(\xi\ne0\), the character case of \eqref{eq:real-frequency-final}
bounds \(d_I\) times its modulus, including at \(a=\varrho\).
At unequal last entries with \(a<\varrho\), the constant-weight case
bounds \(d_I|J_{a+1}|\) at \(\xi=p_a^{-1}-q_a^{-1}\), with the weights
of the extended heads. Since \(d_I\ge1\), the same right sides bound
the respective moduli in Section~\ref{sec:adelic-nonzero-frequencies}.
The terminal form is continuous and positive semidefinite
on fixed last-entry lists, by
Section~\ref{sec:real-terminal-gram}. The remaining
bottom norm estimate is supplied by
Section~\ref{sec:real-bottom-mass} when its formal-degree and Haar bounds hold.

\subsection{The reproducing coefficient and the norm-ball bound}
\label{sec:real-reproducing-coefficient}

\sectionopening{In this subsection, we use Schur orthogonality to express the full coefficient norm in
terms of the formal degree. The same normalization gives a bound for
the volumes used in the truncation argument of
Section~\ref{sec:real-cut-error}.}

Suppose the full vector \(W\) is reproduced on the
right by \(e(g)=d\,\overline{c(g)}\), where
\(c(g)=\langle v,\pi(g)v\rangle\) is the coefficient of a
unit vector in the real discrete series, \(c\in L^1(dg)\),
and \(d\int|c|^2\,dg=1\). Let \(U=U_{<n}\), let \(I\) be the
compact value space, and write \(\chi_U\) for the unitary character
of \(U\) in the Bessel transformation law. Use the space of measurable
functions \(B:G\to I\) satisfying
\[
 B(uy)=\chi_U(u)B(y),\qquad
 \|B\|_2^2=\int_{U\backslash G}\|B(y)\|^2\,d\bar g(y)<\infty ,
\]
modulo equality almost everywhere. Inner products of two such
functions descend to \(U\backslash G\), since \(|\chi_U|=1\).
Right invariance of \(d\bar g\) makes \(R(g)B(y)=B(yg)\) unitary.
Set \(\mathcal P=\int_G e(g)R(g)\,dg\).
Then \(\mathcal P^*=\mathcal P\) by inversion and
\(\mathcal P W=W\) by reproduction.
Schur orthogonality gives
\[
 \int e(a)e(x^{-1}a^{-1}t)\,da
   =d^2\int\overline{c(a)}\,c(t^{-1}ax)\,da
   =c(x)e(t).
\]
Indeed the second coefficient in that pairing has vectors
\(\pi(t)v,\pi(x)v\). The two \(L^1\) kernels justify
the operator multiplication, so \(\mathcal P R(x)\mathcal P
=c(x)\mathcal P\). Consequently \(P(W,W)(x)=d_I c(x)\)
in the notation of Section~\ref{sec:real-bottom-mass}.
Its squared norm is \(d_I^2/d\).

For the norm-ball estimate it suffices to have
\(|c(g)|\ge1/2\) on \(\|g-1\|\le\delta=(200\beta^5)^{-1}\),
where \(1\le\beta\le b^4\) bounds \(r\), and \(b\ge2000\).
Here is a sufficient representation-theoretic condition
for this near-identity bound. In trace-form normalization
let the full and compact Casimirs act on \(v\) by \(\omega_G,\omega_K\),
with \(|\omega_G|,\omega_K\le r\beta^2\). For a real orthonormal
basis \(X_i\) of the compact and noncompact summands in
the Hilbert--Schmidt metric, the full Casimir is
\(\sum_{\mathfrak p}X_i^2-\sum_{\mathfrak k}X_i^2\)
and the compact one is \(-\sum_{\mathfrak k}X_i^2\).

Skew-adjointness gives
\[
 \sum_i\|d\pi(X_i)v\|^2=2\omega_K-\omega_G\le3r\beta^2.
\]
For the lowest compact type use the Casimir eigenvalue formulas
\(\|\lambda\|^2-\|\rho\|^2\) and
\(\|\mu+\rho_K\|^2-\|\rho_K\|^2\) when each coordinate of
\(\lambda,\rho,\mu+\rho_K,\rho_K\) is bounded by \(\beta\).
Here \(\mu=\lambda+\rho-2\rho_K\) by the lowest-type
formula \cite[Theorem~8.5, p.~42]{AS}. For the unitary group one
applies that theorem on the derived group, retaining the
central character on the compact center in both Casimirs.
By Cauchy--Schwarz and \(\|X\|_{\rm HS}\le\sqrt r\|X\|\),
\(\|d\pi(X)v\|\le\sqrt3\,r\beta\|X\|\).
For \(\|g-1\|\le\delta\), the logarithm \(X=\log g\)
belongs to the Lie algebra, since the form adjoint sends
\(g\) to \(g^{-1}\) and the logarithm to \(-X\).
The power series gives \(\|X\|\le2\|g-1\|\).
Integrating the coefficient derivative along \(\exp(tX)\)
gives \(|c(g)-1|\le2\sqrt3\,r\beta\delta\le1/2\).

The unitary norm ball of radius \(Y\ge1\) can be covered
by translates of the identity ball of radius \(\delta\)
using at most \((600\beta^5Y^2)^{2r^2}\) centers.
In fact choose a maximal \(\delta/Y\)-separated set
in the group ball. Disjoint ambient matrix balls of
radius \(\delta/(2Y)\) fit in the ambient ball of
radius \(Y+\delta/(2Y)\). Volume comparison in real
dimension \(2r^2\) gives the stated count. Each group
point \(g\) is within \(\delta/Y\) of some center \(h\),
so \(\|g h^{-1}-1\|\le\delta\), since
\(\|h^{-1}\|=\|h\|\le Y\).
The identity ball has \(dg\)-mass at most \(4/d\) by
integrating \(|c|^2\). Right Haar invariance now gives
\[
 d\operatorname{vol}_{dg}\{g:\|g\|\le Y\}
       \le4(600\beta^5Y^2)^{2r^2}
       \le b^{99r^2}Y^{9r^2}.
\]
This uses the original \(dg\) throughout and gives
the ball bound of Section~\ref{sec:real-bottom-mass}.

\subsection{Error from truncating the orbital integrals}
\label{sec:real-cut-error}

\sectionopening{In this subsection, we estimate the norm of the part removed from the full real vector.
The decay of Section~\ref{sec:real-coefficient-decay} and the volume
estimate of Section~\ref{sec:real-reproducing-coefficient} bound its
contribution to the orbital pairing.}

The same packing argument gives the stronger bound
\(d\operatorname{vol}_{dg}\{\|g\|\le Y\}\le b^{42r^2+1}Y^{4r^2}\).
For \(M_*=9r^2\) define \(\phi(g)=(2/\|g\|)^{M_*}\).
Since \(\|g\|\ge1\), integration of the sublevel volumes
against \(M_*Y^{-M_*-1}dY\) yields
\[
               d\,16^{M_*}\int_G\phi\,dg\le b^{99r^2}.
\]
Indeed \(M_*/(M_*-4r^2)\le2\), and the extra powers
\(2^{M_*}16^{M_*}\) contribute at most \(b^{2M_*}\).

Write \(U=U_{<n}\) and \(S=UH=R_n\).
Use the Bessel graph of the projection \(S\to H\), and
write its elements as \(s=hu\) in the first coordinate.
Its Haar \(ds\) is the product of the original \(du\) and
probability Haar on \(H\). Conjugation by \(H\) preserves
\(du\), so this product is right invariant. The graph is
closed, since \(S\) is a closed flag subgroup and the
projection to the compact factor is continuous. We use the
operator norm in the original unitary basis with prescribed signs, which is
unchanged by transport to the chosen coordinates.

For a point \(x_0\), let \(B_{x_0}\) be the graph set
\(\|x_0^{-1}s x_0-1\|\le1\), with mass \(V_{x_0}\).
This is compact because the first coordinate ranges over
a compact subset of \(G\) and the other coordinate is in
the compact group \(H\). It contains a neighborhood of
the identity, so \(0<V_{x_0}<\infty\).
Put \(\Phi(y)=\int_S\phi(x_0^{-1}s y)\,ds\).
To bound this function, take a maximal \(1/Y\)-separated
set of matrices \(a_i=x_0^{-1}s_i y\) of norm at most \(Y\).

Ambient packing gives at most \(3^{2r^2}Y^{4r^2}\) centers
for \(Y\ge1\). If \(a=x_0^{-1}s y\) is within \(1/Y\)
of \(a_i\), then
\[
 \|x_0^{-1}s s_i^{-1}x_0-1\|
                =\|(a-a_i)a_i^{-1}\|\le1.
\]
Thus each piece is contained in \(B_{x_0}s_i\), of
mass \(V_{x_0}\) by right invariance. Integrating the
sublevel estimate against \(2^{M_*}M_*Y^{-M_*-1}\,dY\)
gives \(\Phi(y)\le16^{M_*}V_{x_0}\).
For any integrable function \(\gamma\) on \(G\), the
original quotient disintegration gives
\[
 \int_{U\backslash G}\int_S\gamma(s y)\,ds\,d\bar g(y)
                                   =\int_G\gamma(g)\,dg.
\]
Indeed write \(s=hu\), unfold \(du\,d\bar g\) for each
\(h\), and apply left Haar invariance. The same proof
applies to nonnegative functions. With left translation
by \(x_0^{-1}\) it gives \(\int\Phi\,d\bar g=\int\phi\,dg\).

For \(s=hu\), let \(\sigma(u)\) be the simple-root sum in
the real chart of Section~\ref{sec:period-real-point}, and put
\[
 \chi_U(u)=\exp\bigl(\sigma(u)-\overline{\sigma(u)}\bigr),\qquad
 \sigma_S(hu)=\chi_U(u)\pi_H(h).
\]
The first superdiagonal entries of the unitriangular flag action
add under multiplication. The compact group \(H\) fixes the ordered
flag vectors, so conjugation by \(H\) preserves \(\chi_U\).
Thus \(\sigma_S\) is a unitary representation of \(S\).
For the unit highest vector \(v_0\) of the compact type, set
\[
 A(y)=\int_S c(x_0^{-1}s y)\sigma_S(s)^{-1}v_0\,ds .
\]
The large-gap coefficient bounds used here are
\(|c(g)|\le\phi(g)\) and
\(|c(g)|\le(2/\|g\|)^{M_*+3r^3f}\), as proved in
Section~\ref{sec:real-coefficient-decay} for the prescribed
lowest compact type and parameter gap.
The volume estimate above uses square integrability and the
coefficient bound at the identity. Since \(|c|\le\phi\) and
\(\int_G\phi\,dg<\infty\), we have \(c\in L^1(G,dg)\) before
using the reproducing convolution. The same majorant gives
absolute convergence in the graph variable and \(|A|\le\Phi\). In particular \(A\) is
bounded, and its norm belongs to \(L^1\cap L^2\) on
\(U\backslash G\). It is continuous on \(G\): for \(y=y_0 a\) with \(a\)
in a compact neighborhood of \(1\), the norm bound
is dominated by \(\|a^{-1}\|^{M_*}\phi(x_0^{-1}s y_0)\).
The substitution \(s\mapsto s t^{-1}\) gives
\(A(ty)=\sigma_S(t)A(y)\). Also
\[
             \int_G e(a)A(ya)\,da=A(y).
\]
Indeed Schur orthogonality gives
\(\int e(a)c(g a)\,da=c(g)\). The graph integral
commutes with this convolution by the bound
\(\|e\|_1\,16^{M_*}V_{x_0}\) on the absolute double integral.

The period at \(x_0\) is exactly \(d\|A\|_2^2\).
To check the constant, unfold the integrable function
\(d\,\overline{c(x_0^{-1}g)}\langle v_0,A(g)\rangle\).
At \(g=sy\), covariance and unitarity make it
\[
 d\,\big\langle c(x_0^{-1}s y)\sigma_S(s)^{-1}v_0,A(y)\big\rangle.
\]
Its graph integral is \(d\|A(y)\|^2\).
On the full group, substitute \(g=x_0 a\) and use the
reproducing identity. The integral is
\(\langle v_0,A(x_0)\rangle\), which on expanding \(A\)
is the stationary period on this same graph and Haar.

The quantitative stationary bound is supplied by
Section~\ref{sec:real-stationary-point} under its weight and gap conditions. It gives a point \(x_0\) for which
\[
 \ell(x_0)\le b^{5n},\qquad
 V_{x_0}\le b^{60(r+k)^2}|\langle v_0,A(x_0)\rangle|,
 \qquad \ell(y)=\inf_{s\in S}(\|s y\|+1).
\]
Then \(A\ne0\) by the positive mass of the graph ball.
Take the positive scalar \(c_0=\sqrt{d_I}/\|A\|_2\) and put
\(W=c_0A\), so \(\|W\|_2^2=d_I\).
The length is left \(S\)-invariant and satisfies
\(\ell(x_0 g)\le\ell(x_0)\|g\|\), since \(\|g\|\ge1\).
For \(g=x_0^{-1}s y\) this gives
\(\ell(y)/b^{6n}\le\|g\|/2\).

If \(H(y)=\max(1,\ell(y)/b^{6n})>1\), the second
coefficient estimate absorbs \(H(y)^{3r^3f}\).
When \(H(y)=1\), use the first estimate. In either
case \(H(y)^{3r^3f}|W(y)|\le |c_0|\Phi(y)\).
With \(w=|c_0|\Phi\), the period identity and the graph and group
volume bounds give
\[
 \begin{aligned}
 (\sup w)\int w\,d\bar g
 &\le d_I\frac{V_{x_0}}{\langle v_0,A(x_0)\rangle}
          \;d\,16^{M_*}\int_G\phi\,dg\\
 &\le d_I b^{60(r+k)^2+99r^2}\le d_I b^{400r^2}.
 \end{aligned}
\]
Here \(\langle v_0,A(x_0)\rangle=d\|A\|_2^2>0\) uses the same
quotient Haar and measured point as the normalization of \(W\).
\(\ell\) is continuous by its right-translation norm
inequality and the same inequality for the inverse.
The scalar cutoff \(\chi=\max(0,\min(1,2-2\ell/R_b))\)
is one on \(\ell(y)\le R_b/2\). If \(\chi(y)>0\) there
is \(s=hu\) with \(\|s y\|<R_b\), hence \(\|uy\|<R_b\)
because \(H\) is Euclidean unitary.

Put \(F=\chi W\). Since \(U\) is closed and
\(K_b=\{g\in G:\|g\|\le R_b\}\) is compact, \(UK_b\) is closed.
Thus \(\operatorname{supp}F\subseteq UK_b\). For the right correlation
of Section~\ref{sec:real-bottom-mass}, a nonzero integrand requires
\(y=u_1k_1\) and \(yx=u_2k_2\), with \(u_i\in U\) and \(k_i\in K_b\).
Hence
\[
 \operatorname{supp}P(F,F)\subseteq K_b^{-1}UK_b,
 \qquad x=k_1^{-1}u_1^{-1}u_2k_2 .
\]
This set is closed because the outside factors are compact and \(U\)
is closed. For the left translate in the orbital tensor,
\(F(y)\otimes\overline{F(\gamma^{-1}y)}\ne0\) gives instead
\(\gamma=u_1k_1k_2^{-1}u_2^{-1}\in UK_bK_b^{-1}U\).

On the support of \(F-W\), one has
\(H(y)\ge b\) for \(R_b=b^{6n+2}\) and \(b\ge2\).
Thus
\[
 \|F-W\|_2^2\le b^{-6r^3f}\int w^2\,d\bar g
              \le b^{400r^2-6r^3f}d_I
              \le b^{-2r^3f}d_I
\]
in the range \(r\ge100,\ f\ge1\).
The stationary-period and coefficient decay bounds in this
argument concern the same \(x_0,A\) as the reproducing vector.

\subsection{The formal-degree factor}
\label{sec:real-formal-degree-factor}

\sectionopening{In this subsection, we compare the formal degree in the orbital Haar measure with the
compact normalization of the Weyl dimension formula. The resulting
factor enters the initial norm calculation of
Section~\ref{sec:real-bottom-mass}.}

In the Cayley normalization, the density at the identity
of \(U(J)\) in the skew-Hermitian coordinates \(S=JX\)
of its Lie algebra is \(1/\operatorname{vol}U(r)\).
The corresponding density on the compact dual is the
same. In a unitary basis where \(J\) is diagonal, the Cartan
decomposition keeps the two diagonal skew-Hermitian blocks
and sends the off-diagonal noncompact part \(X_{\mathfrak p}\)
to \(iX_{\mathfrak p}\). In \(S\)-coordinates this change
has real determinant of absolute value one, since it
multiplies the independent complex off-diagonal entries
by a unit complex scalar.

For signature \((p,q)\), assign probability Haar to
\(K=U(p)\times U(q)\). Probability Haar on \(U(r)\) then induces
mass one on the complex Grassmannian \(U(r)/K\), the compact dual
of \(U(p,q)/K\). With the tangent measures identified above, these
are the two measure normalizations in the Harish--Chandra degree
formula \cite[(3.15)--(3.16), p.~14, and Theorem~8.1, p.~41]{AS}.
Indeed put \(G_1=SU(p,q)\) and \(K_1=S(U(p)\times U(q))\).
The maps \((z,h)\mapsto zh\) from \(U(1)\times G_1\) to \(U(p,q)\)
and from \(U(1)\times SU(r)\) to \(U(r)\) have the same kernel
\(\{(\xi,\xi^{-1}I):\xi^r=1\}\). On the compact side the product
of the two probability measures pushes forward to probability Haar.
On the noncompact side take probability Haar \(dz\) on the circle
and the degree-formula Haar \(dh\) on \(G_1\). The quotient
\(G_1/K_1=U(p,q)/K\) has the same tangent measure as its compact dual
\(SU(r)/K_1=U(r)/K\), and probability Haar on \(U(1)\times K_1\)
pushes forward to probability Haar on \(K\). Thus \(dz\,dh\)
pushes forward to \(dg\). For a coefficient \(c\) of an irreducible
unitary representation the circle acts by a scalar of modulus one,
so its restriction to \(G_1\) is irreducible and
\[
 \int_{U(p,q)}|c(g)|^2\,dg
  =\int_{U(1)\times G_1}|c(zh)|^2\,dz\,dh
  =\int_{G_1}|c(h)|^2\,dh .
\]

For a regular parameter \(\lambda_1,\ldots,\lambda_r\),
the type \(A_{r-1}\) root products give
\[
          d=\frac{\prod_{i<j}|\lambda_i-\lambda_j|}
                        {\prod_{i=0}^{r-1}i!}.
\]
In fact the product of the positive-root differences at
\(\rho\) is \(\prod_{j=2}^r (j-1)!\), by grouping the
factors \(j-i\) at fixed \(j\).
For the compact volume in these same coordinates, project
\(U(l)\) to its last column on the unit sphere in \(\mathbb C^l\).
At the identity the differential is \(S\mapsto Se_l\). Its kernel
has the skew coordinates of \(U(l-1)\), and its remaining coordinates
are the \(l-1\) complex entries and the imaginary diagonal entry of
that column. They give the Euclidean tangent density of the sphere.
By invariance this identifies the quotient measure everywhere with
surface measure. Polar integration of the Gaussian therefore gives
\[
 \frac{\operatorname{vol}U(l)}{\operatorname{vol}U(l-1)}
   =\operatorname{area}(S^{2l-1})=\frac{2\pi^l}{(l-1)!},
 \qquad \operatorname{vol}U(0)=1 .
\]
Set \(t=r-k=2\varrho+1\) and
\(\Theta=\pi^{tk}\operatorname{vol}U(k)/\operatorname{vol}U(r)\)
as in \eqref{eq:real-shell-total}. Consequently
\[
 d\Theta=
 \frac{\prod_{i<j}|\lambda_i-\lambda_j|}
              {2^t\pi^{t(t+1)/2}\prod_{i=0}^{k-1}i!}.
\]
The factorials with indices at least \(k\) cancel,
and the exponent of \(\pi\) comes from
\(\sum_{i=k}^{r-1}(i+1)-tk=t(t+1)/2\).
Thus \(|2\lambda_i|\le b^4\) on the parameter
implies \(d\Theta\le b^{2r(r-1)}\le b^{4r^2}\),
as used in Section~\ref{sec:real-bottom-mass}.

\subsection{Haar measure in the adapted real basis}
\label{sec:real-frame-haar}

\sectionopening{In this subsection, we identify the real Haar measure obtained by successive isotropic
quotients with the group measure used for discrete series. This gives
the conversion needed in Section~\ref{sec:real-bottom-mass} and
matches the finite convention of Section~\ref{sec:compact-floor-haar}.}

Use the basis ordered by the pairs, then \(z\), then
the \(k\) positive coordinates. Its Gram matrix \(J\) has
the hyperbolic pair blocks and the diagonal blocks \(-1,I_k\).
Let \(d\mu\) be the real stage Haar measure, with the terminal
stabilizer \(U(k)\) of mass one. Then
\[
 d\mu=\frac{\operatorname{vol}U(r)}
                 {\operatorname{vol}U(k)}\,dg.
\]

For a matrix \(Y\) of independent columns, put
\(\Gamma(Y)=Y^*JY\). Use Lebesgue measure on the complex
entries of \(Y\), and on the Gram entries use the real
diagonals and the real and imaginary parts of the upper
off-diagonals. The map \(\Gamma\) is a submersion since
\(Y^*J\dot Y\) can be prescribed. With an initial segment \(Y_0\)
fixed, the differential in the new columns surjects onto
the cross block and the new diagonal block by the same
linear solve with \(Y^*J\).

Thus local submersion coordinates
disintegrate Lebesgue measure over these blocks in order.
In fact the coordinate Jacobian is block triangular, with
the old Gram coordinates depending only on \(Y_0\).
Consequently the full fiber measure is obtained by the
successive conditional fiber integrations at the prescribed
Gram entries. For an initial segment of the form \(gY_0\) with
\(g\in U(J)\) and \(Y_0\) the specified initial segment, the conditional measure
is transported by \(g\). Both ambient Lebesgue measure and
all the Gram entries are preserved, since \(|\det g|=1\).

For the full square Gram fiber at \(J\), its density at the
identity is computed by
\[
 \dot Y=J(S+T/2),\qquad \dot\Gamma=T,\qquad S^*=-S,\quad T^*=T.
\]
Left multiplication by \(J\) preserves Lebesgue measure.
Each diagonal contributes \(1/2\). For an off-diagonal
pair the complex map from \((S_{ij},T_{ij})\) to
\((S_{ij}+T_{ij}/2,-S_{ij}+T_{ij}/2)\) has determinant one.
Conjugating the second image preserves real volume.
The fiber density is therefore \(2^{-r}\) in the
\(S\)-coordinates. The full fiber is left invariant, so
the Cayley normalization gives its measure as
\(2^{-r}\operatorname{vol}U(r)\,dg\).

Inversion of the stage product gives, in order, the factors
\(\theta_i^{\rm inv},\eta_i\) for each pair, then
\(\theta_\varrho^{\rm inv}\) for \(z\), and finally
probability Haar on \(U(k)\). Here \(\theta^{\rm inv}\)
denotes the image under inversion. Every new factor fixes
the earlier columns. Fixing the cross products with
the earlier pairs sets their opposite coordinates to zero,
with Jacobian one. When adding \(e_i\) or \(z\), the chart
of Section~\ref{sec:compact-floor-haar} thus disintegrates
in the remaining space.

In its hyperbolic pivot coordinates
\((t,tc,v)\), the ambient Jacobian is \(|t|^2\) and
\[
 h((t,tc,v),(t,tc,v))=2|t|^2\Re c+h(v,v).
\]
The change from \(\Re c\) to the real Gram entry gives
\(1/(2|t|^2)\). The top hyperbolic change of basis has
determinant of modulus one, since both Gram determinants
on the plane are \(-1\). Thus half the stage quadric measure is the conditional fiber
measure on its pivot chart. The discarded set \(t=0\) has measure
zero for both measures. Indeed the nonzero isotropic vectors form
one orbit, as do the vectors of the prescribed negative norm, by
extension of isometries. Each orbit is connected because
\(U(p,q)\) is connected by Cartan decomposition. The real-analytic
function \(|t|^2\) is nonzero somewhere on each orbit by the pivot
parametrization. Its zero set has measure zero for any smooth
density there. The isotropic origin cannot occur in an independent
frame.

For the dual column the additional cross product with
\(e_i\) fixes the coefficient of \(e_i^*\) to one.
Write the column as \(e_i^*+v+(\alpha+i\beta)e_i\).
Its norm is \(2\alpha+h(v,v)\), so the fiber law at
zero has density \(dv\,d\beta/2\). This is half the
Heisenberg measure on that column.

After the \(2\varrho+1\) columns have been inserted, their
complement has Gram matrix \(I_k\) in the chosen coordinates.
The cross products with the initial segment again give a coordinate
map of determinant of modulus one. The calculation for the
full Gram fiber, now with \(J=I_k\), gives the terminal
conditional measure \(2^{-k}\operatorname{vol}U(k)\)
times the compact probability Haar. All the pivot charts have full
conditional measure. Multiplication of the \(2\varrho+1\) factors
\(1/2\) with the terminal factor \(2^{-k}\) gives
\[
 2^{-r}\operatorname{vol}U(k)\,d\mu^{\rm inv}
 =2^{-r}\operatorname{vol}U(r)\,dg .
\]
The square-column map sends a group element to its matrix in the
chosen basis. At each step the norm-fiber measure is invariant under
the remaining unitary group, since the group preserves the norm map
and ambient Lebesgue measure. Changing the orbit representative in
right translation changes the inner integral by a stabilizer
translation. Induction from probability Haar on \(U(k)\), as in
Section~\ref{sec:compact-floor-haar}, proves right invariance of
\(d\mu\). The displayed identity and the nonzero chart densities give
\(d\mu=\lambda\,dg\) with \(0<\lambda<\infty\), since \(G\) is unimodular
as verified below.

To compare the quotient measures, order the Heisenberg factors by
increasing \(a\). The factor with parameters \((v,\beta)\) fixes the
earlier dual columns and sends \(e_a^*\) to
\[
 e_a^*+v+\bigl(i\beta-h(v,v)/2\bigr)e_a .
\]
The product of earlier factors only adds earlier isotropic vectors
to this column. The resulting column therefore has coordinates
\(v\) in the remaining space and imaginary coordinate \(\beta\) along \(e_a\).
Moreover \(u^{-1}=Ju^*J\) identifies these dual-column coordinates
of \(u^{-1}\) with the free entries of \(u\) up to permutation,
signs and conjugation. All have real Jacobian of absolute value one.
Thus the product of the Heisenberg Haar measures is the original
unipotent Haar, since inversion preserves Haar on a unipotent group.
Dividing the group measures by this common Haar gives
\(dQ=\lambda\,d\bar g\) on \(U_{<n}\backslash G\).

The adjoint map \(S\mapsto g^{-*}Sg^{-1}\) has real determinant
of absolute value \(|\det g|^{-2r}=1\). Thus \(G\) is unimodular,
and inversion preserves both group measures. The Gram fiber identity gives
\(\lambda=\operatorname{vol}U(r)/\operatorname{vol}U(k)\).
In particular the coefficient of the initial squared correlation
is \(C_F\lambda^3=d_I^{-2}\) by Section~\ref{sec:real-head-density}.

\subsection{Decay of the lowest compact type coefficient}
\label{sec:real-coefficient-decay}

\sectionopening{In this subsection, we bound the lowest compact type coefficient using its Harish--Chandra
parameter. The bound is used away from compact sets in
Section~\ref{sec:real-cut-error} and in the tail of the
stationary-phase estimate of Section~\ref{sec:real-stationary-phase}.}

Work in the unitary basis with prescribed signs \(J=\operatorname{diag}(\varepsilon_i)\).
Write \(\lambda_1>\cdots>\lambda_r\) for the Harish--Chandra
parameter, \(\rho_i=(r+1-2i)/2\), and \(a_i=\lambda_i+\rho_i\),
which are integers. For \(0\le c_*\le h\) integral, where
\(|\lambda_i-\lambda_j|\ge h\) for opposite signs, we prove
\[
 |c(g)|\le \min\{1,(2/\|g\|)^{c_*}\}
\]
for a coefficient of unit vectors in the lowest compact type.

For a matrix \(Z\), let \(D_I(Z)\) be the principal minor on
\(I\). For a weight \(w\in\mathbb Z^r\) and \(|w_i|\le M\) put
\[
 m_w(Z)=(\det Z)^{-M}
        \prod_{t=-M+1}^{M}D_{\{i:w_i\ge t\}}(Z).
\]
If \(Z\in U(J)\) is positive definite, then
\(Z^{-1}=JZJ\) and \(\det Z=1\). Its block diagonal part
\(D=(Z+Z^{-1})/2\) is positive definite.
For each threshold \(t\), take the partition
\[
 I=\{a_i\ge t\},\qquad B=\{a_i<t\le a_i+c_*\},\qquad
 C=\{a_i+c_*<t\}.
\]
The signs on \(B\) agree since \(\rho_i\) is ordered as
\(\lambda_i\). Thus \(D_B(D)=D_B(Z)\). Fischer's inequality
(the positive Schur complement bound) and complementary
minors give
\[
 D_{I\cup B}(D)
 \le D_I(D)D_B(Z)
 =D_I(D)D_{I\cup C}(Z)
 \le D_I(D)D_I(Z)D_C(Z).
\]
Multiply using a common \(M\) for \(a\) and \(a+c_*\).
The product of \(D_C(Z)\) is \(m_{a+c_*}(Z)\) by
complementary minors. Since \(m_{a+c_*}(W)=m_a(W)(\det W)^{c_*}\),
cancelling \(m_a(D)\) gives
\[
                   (\det D)^{c_*}\le m_a(Z)^2 .
\]

For the compact-Cartan Borel \(\mathfrak b\) determined by \(\lambda\),
the discrete-series Harish--Chandra module is
\(A_{\mathfrak b}(\lambda-\rho)\)
\cite[\S9, paragraph following Theorem~9.2 and Theorem~9.4]{AdamsVoganDiscrete},
\cite[\S3.5, (3.4)]{IchinoRealTheta}.
Its inducing weight is in the good range, since
\((\lambda-\rho)+\rho=\lambda\) pairs positively with every
root of the nilradical in the chosen order.
Let \(V\) be the lowest compact type with unit highest vector
\(v\) and orthogonal projection \(E_V\). On the positive
factor \(Z=e^X\) of the polar decomposition, Sun's bottom
layer integral \cite[Theorem~6.3 and (25)--(26), pp.~219--220]{SunMC} gives
\[
 E_V\pi(Z)E_V=\dim V\int_K
             m_a(k^{-1}Zk)^{-1} P_{\pi(k)v}\,dk ,
\]
where \(K=U(J)\cap U(r)\), \(dk(K)=1\), and \(P_u\) is the
orthogonal projector onto \(\mathbb C u\).
For the compact-Cartan Borel the inducing character has weight
\(\lambda-\rho\). For same-sign indices \(i<j\), let \(s\) count the indices
\(i<h\le j\) of that sign. The proposed compact weight
\(\mu=a-2\rho_K\) has
\(\mu_i-\mu_j=(\lambda_i-\lambda_j)+(j-i)-2s\ge0\) by
integrality and regularity. Multiplication by the determinant
character of the nilradical gives weight \((\lambda-\rho)+2\rho=a\).
The lowest compact type has highest weight
\((\lambda-\rho)+2(\rho-\rho_K)=a-2\rho_K\).
In Sun's formula take the compact argument to be \(1\), the positive
argument to be \(e^{X/2}\), and the compact integration variable to
be \(k^{-1}\).

In the upper triangular times unitary decomposition
\(k^{-1}e^{X/2}=u d k'\), take \(u\) unitriangular and
\(d\) positive diagonal. Sun's map \(H_{\mathfrak b}\) in (7)
reads \(d\) modulo the compact torus. The complexified Cartan
involution is conjugation by \(J\), and real-form conjugation sends
\(g\) to \(J(g^*)^{-1}J\). Hence the factor
\(d\,\overline{\theta(d^{-1})}\) in (25) is \(d^2\).
Also \(k^{-1}Zk=u d^2u^*\). Trailing principal minors give
\(d_i^2=D_{[1,i-1]}(k^{-1}Zk)/D_{[1,i]}(k^{-1}Zk)\),
using \(Z^{-1}=JZJ\) and complementary minors. The torus
weight therefore gives \(m_a(k^{-1}Zk)^{-1}\).
The positive factor has determinant one, so this use of
the derived group formula preserves the central character.
Schur orthogonality gives \(\dim V\int_K P_{\pi(k)v}\,dk=E_V\).
The scalar \(m_a(k^{-1}Zk)^{-1}\) is positive and at most
\((\det D)^{-c_*/2}\), since \(K\) preserves the determinant
of the block diagonal part.

Thus the same bound holds in
operator norm. The compact factor of \(g\) acts unitarily
on \(V\), so it preserves that bound. Write \(g=D_g+E_g\)
for its equal-sign and opposite-sign blocks. The diagonal
blocks of \(g^*Jg=J\) give \(D_g^*D_g=1+E_g^*E_g\), and
\(|\det D_g|=\det D\) by the polar decomposition. Therefore
\[
 |c(g)|^4\le (1+\|E_g\|_{\mathrm{HS}}^2)^{-c_*},
\]
using \(\det(1+E_g^*E_g)\ge1+\operatorname{tr}(E_g^*E_g)\).

The eigenvalues \(u\) of \(Z\) contribute \((u+u^{-1})/2\)
to \(\det D\), each at least one. If \(\|Z\|>1\), then
\(\|Z\|\) and \(\|Z\|^{-1}\) are distinct eigenvalues and
their two factors give \(\det D\ge(\|Z\|/2)^2\). The same
inequality holds for \(\|Z\|=1\). Since \(\|g\|=\|Z\|\),
this proves the decay claim. In Section~\ref{sec:real-cut-error}
use \(c_*=9r^2\) and \(c_*=9r^2+3r^3f\), both at most
the available gap whenever \(3r^3f\le A\) and
\(h\ge500r^3+A\).

\subsection{The stationary point and its period}
\label{sec:real-stationary-point}

\sectionopening{In this subsection, we construct the stationary point for the original pair of discrete
series. The stationary-phase calculation of
Section~\ref{sec:real-stationary-phase} gives the period lower bound
used to normalize the full vector in
Section~\ref{sec:real-marked-filling}.}

\begin{proposition}[Stationary point and period lower bound]\label{prop:marked-real-period}
Use the pair of discrete series with a common zero-weight
positive block, integral weights \(\mu_i\) on the remaining \(2n\)
coordinate lines of signatures \((-1)^i\), and remaining weight \(a\)
in \(H\). Put \(N_0=(r+k)^2\), \(\beta=b^3\) and \(T_0=N_0\log\beta\).
Assume
\[
 B_{\rm wt}\ge1,\quad -3B_{\rm wt}\le\mu_{2n-1}<\cdots<\mu_0\le-B_{\rm wt},\quad
 a\ge0,\quad 3B_{\rm wt}+a\le b^2,\quad 8nr\le b
\]
and \(\beta\ge10^{12}\). Assume that all doubled Harish--Chandra parameters
of both representations have absolute value at most \(L'_0\), with
\(r+k+L'_0\le\beta\). For some integer \(A_{\rm gap}\ge10^{30}T_0^3\), assume
\( |\lambda_i-\lambda_j|\ge500r^3+A_{\rm gap}\) for every pair of ambient
parameters on coordinate lines of opposite signs. Then there is a
stationary point \(x_0\) for which
\[
 \ell(x_0)\le b^{5n},\qquad
 V_{x_0}\le b^{60N_0}|\langle v_0,A(x_0)\rangle|.
\]
\end{proposition}

\begin{proof}

Let the lowest compact type weights of the original two
representations be zero on a common positive block \(X\).
In the remaining \(2n\) coordinates of \(G\), take integral
weights \(-3B_{\rm wt}\le\mu_i\le-B_{\rm wt}\) in strictly decreasing order
at signs \((-1)^i\), \(0\le i<2n\), with \(B_{\rm wt}\ge1\).
The remaining coordinate in \(H=U(k)\) has weight \(a\ge0\).
Use \(\beta=b^3\), and suppose \(T=3B_{\rm wt}+a\le b^2\) and
\(8nr\le b\). The following construction takes place in
these two original weight bases.

Put \(P(t)=\prod_i(t-\mu_i)\) and, on the \(2n\) coordinates,
\[
 [F]_i=\frac{F(\mu_i)}{\sqrt{|P'(\mu_i)(\mu_i-a)|}},
 \qquad V_j=[t^j(t-a)]\quad(0\le j\le 2\varrho),
 \qquad I=\sqrt{P(a)}[1].
\]
Extend these vectors by zero on \(X\). For real polynomials
the signed pairing of \([F]\) and \([G]\) is
\(-\sum_i F(\mu_i)G(\mu_i)/(P'(\mu_i)(\mu_i-a))\).
By Lagrange interpolation, \(\sum_i Q(\mu_i)/P'(\mu_i)\)
is the top coefficient of \(Q\) for \(\deg Q\le2n-1\).
It follows that
\[
 h(V_j,V_l)=0\quad(j+l<2\varrho),\qquad
 h(V_j,V_l)=-1\quad(j+l=2\varrho),\qquad h(I,V_j)=0.
\]
Partial fractions of \(1/(P(t)(t-a))\) have residues
summing to zero, with residue \(1/P(a)\) at \(a\).
This gives \(h(I,I)=1\). Let \(i_0\) embed the compact
space by the common coordinates on \(X\) and by \(I\)
at its final coordinate. The embedding is isometric.
The vectors \(V_0,\ldots,V_\varrho\) are independent,
since \(t-a\) is nonzero on the \(2n\) nodes and a
polynomial of degree at most \(\varrho\) vanishing there
is zero. Set
\[
 D_\mu=\operatorname{diag}(0_X,\mu_0,\ldots,\mu_{2n-1}),\qquad
 D_H=\operatorname{diag}(0_X,a).
\]
Then
\[
 D_\mu V_j=V_{j+1}\quad(j<2\varrho),\qquad
 D_\mu i_0=i_0 D_H\bmod\langle V_0\rangle .
\]

Transport the basis with prescribed signs into the original stage space
by its Euclidean unitary change of basis. The perpendicular
complement of \(i_0\mathbb C^k\oplus\mathbb C V_\varrho\)
is nondegenerate. In it let \(A\) have columns
\(V_0,\ldots,V_{\varrho-1}\) and choose \(B_0\) with
\(A^*JB_0=1\). Replacing \(B_0\) by
\(B_0-A(B_0^*JB_0)/2\) gives isotropic duals.
Together with \(i_0,V_\varrho\) the columns give an
isometry \(x_0^{-1}\) from the original basis,
sending \(e_j\) to \(V_j\), \(z\) to \(V_\varrho\), and
the positive block to \(i_0\mathbb C^k\) with its chosen
basis.

To bound the length of this same \(x_0\), use Euclidean
Gram--Schmidt on \((V_0,\ldots,V_\varrho)\), without
rescaling. The resulting \(F_j\) differ from \(V_j\)
by earlier columns, and retain the Gram and the pairings
\(h(F_j,V_{2\varrho-j})=-1\).
As \(|P'(\mu_i)|\ge1\), the coordinate formula gives
\(\|V_j\|\le M_0=\sqrt r\,T^{2n}\) through \(2\varrho\).
Thus \(M_0^{-1}\le\|F_j\|\le M_0\).
For \(j<\varrho\), the vectors \(JF_j/\|F_j\|^2\) are
isotropic duals, with Euclidean norm at most \(M_0\).
The pairs form a \(J\)-invariant Euclidean orthogonal
summand, and \(F_\varrho\) lies in its complement of
signature \((k,1)\).

Write \(F_\varrho=z_++z_-\) in the positive and negative
eigenspaces, with \(s=\|z_+\|\), \(t=\|z_-\|\) and \(t^2-s^2=1\).
For \(s>0\), set \(v=t z_+/s+s z_-/t\). In the Hermitian
form \(h\) of signature \((k,1)\),
\[
 h(v,v)=1,\qquad h(v,F_\varrho)=0,\qquad
 \|v\|^2=s^2+t^2=\|F_\varrho\|^2 .
\]
For \(s=0\), take a Euclidean unit positive vector in that complement. Complete by the remaining
orthonormal positive directions, orthogonal to these
at most \(2n\) exceptional columns.
The new basis matrix \(Y\) is an isometry with
\(\|Y^{-1}\|=\|Y\|\le \sqrt{2n}M_0+1\).
Put \(s_0=x_0Y\), \(A_j=\operatorname{span}(e_i:i<j)\), and
\(A=A_\varrho\) in the original basis. The relations
\(F_j-V_j\in\operatorname{span}(V_i:i<j)\) give
\[
 s_0e_j-e_j\in A_j\quad(0\le j<\varrho),\qquad s_0z-z\in A .
\]
Unitarity preserves \(A^\perp\). The induced isometry on
\(A^\perp/A\) fixes \([z]\) and acts by \(u\oplus1\), with \(u\in U(k)\).
After removing the standard lift of \(u\), the isometry is upper
unitriangular on \(A\) and the identity
on \(A^\perp/A\). Unitarity forces the inverse-adjoint block on the
opposite isotropic quotient. These are the defining conditions
for \(U_{<n}\). Thus \(s_0\in S\), and
\[
 s_0^{-1}x_0=Y^{-1},\qquad
 \ell(x_0)\le1+\|Y^{-1}\|\le8nM_0\le b^{5n}.
\]
The representative \(s_0^{-1}x_0\) is used in the length infimum,
while the measured period is taken at the stationary point \(x_0\).

Return to the basis with prescribed signs, and for a graph element \(s\)
write \(g=x_0^{-1}s_Gx_0\) in that basis, with compact
component \(h=s_H\). Put \(Y_1=g-1\), \(Y_2=h-1\),
\(L=\langle V_0,\ldots,V_{\varrho-1}\rangle\), and let
\(\sigma(s)\) be the sum of first predecessor coefficients
on the original basis. Then \(g-1\) lowers the
transported flag, and on \(L^\perp/L\) the action of \(g\) is
\(h\oplus1\) in the transported \(i_0,V_\varrho\).
Conversely these equations characterize the conjugated
graph, since in the original Witt blocks they give
exactly \(U_{<n}H\) and its compact projection. In
particular the diagonal \(U(X)\) belongs to the transported
graph, with \(\sigma=0\).

Choose a projection \(P_0\) onto \(L\) killing
\(i_0\mathbb C^k\), put \(D_1=D_\mu P_0\), and let
\(D_2=D_\mu-D_1-D_1^\dagger\), with \(\dagger\) the form adjoint.
Since \(D_\mu L\subset L^\perp\), \(D_2\) kills \(L\) and
has image in \(L^\perp\). On \(i_0\mathbb C^k\) it agrees
with \(D_\mu\), since \(D_\mu i_0\mathbb C^k\subset L^\perp\).
On \(L^\perp/L=i_0\mathbb C^k\oplus\mathbb C[V_\varrho]\),
\(Y_1\) induces \(Y_2\oplus0\). The compression of \(D_2\)
to the first summand is \(i_0^\dagger D_\mu i_0=D_H\), since
\(V_0\perp i_0\mathbb C^k\). The traces of \(D_2Y_1\) on \(L\)
and on the quotient by \(L^\perp\) are zero. Thus
\(\operatorname{tr}(D_2Y_1)=\operatorname{tr}(D_HY_2)\).
Since \(D_\mu V_j=V_{j+1}\), the diagonal entry of
\(P_0Y_1D_\mu\) on \(L\) at \(V_j\) is the coefficient of \(V_j\) in
\(Y_1V_{j+1}\) for \(0\le j<\varrho\). These entries sum to
\(\sigma(s)\), with the empty sum equal to zero if \(\varrho=0\).
Applying
this calculation to \(Y_1^\dagger=g^{-1}-1\) gives
\(\operatorname{tr}(D_1^\dagger Y_1)=-\bar\sigma(s)\),
since inversion negates the first predecessor coefficients.
The flag and its dual contribute zero to the unweighted trace.
We obtain
\[
 \operatorname{tr}Y_1=\operatorname{tr}Y_2,\qquad
 \operatorname{tr}(D_\mu Y_1)-\operatorname{tr}(D_HY_2)
                         =\sigma(s)-\bar\sigma(s).
\]

With the unit highest vector \(v_0\) of the compact type,
put \(c_H(h)=\langle v_0,\pi_H(h)v_0\rangle\). The period
integrand for the average of Section~\ref{sec:real-cut-error} is
\(\theta(s)=c(g)\overline{c_H(h)}e^{\bar\sigma(s)-\sigma(s)}\).
All following volumes and integrals use the original \(ds\).
Put \(N_0=(r+k)^2\) and \(Q=(2/\|g\oplus h\|)^{9N_0}\).
Since \(\|h\|=1\), the decay just proved gives \(|\theta|\le Q\)
when the gap is at least \(9N_0\). If the doubled parameters
are bounded by \(L_0'\) with \(r+k+L_0'\le\beta\), the
Casimir calculation of Section~\ref{sec:real-reproducing-coefficient}
on each factor gives \(|\theta|\ge1/4\) on the ball
\(\|g\oplus h-1\|\le\delta=(200\beta^5)^{-1}\).

For the phase estimate put \(T_0=N_0\log\beta\) and
\(\eta_0=10^{-13}/T_0\), and suppose the excess gap
is at least \(A_{\rm gap}\ge10^{30}T_0^3\). Since \(A_{\rm gap}\ge18N_0\),
the decay bound also gives \(|\theta|\le Q^2\).
If \(\|E_g\|_{\mathrm{HS}}\ge\eta_0\), the preceding
determinant estimate with \(c_*=A_{\rm gap}\) gives
\[
 |\theta|\le(1+\eta_0^2)^{-A_{\rm gap}/4}
           \le\exp(-A_{\rm gap}\eta_0^2/8)\le e^{-100T_0}.
\]
Section~\ref{sec:real-stationary-phase} proves
\(|\theta|-2\Re\theta\le12e^{-100T_0}\) when
\(\|E_g\|_{\mathrm{HS}}<\eta_0\). Away from that set the
same bound follows from \(|\theta|-2\Re\theta\le3|\theta|\).
Combining it with \(|\theta|-2\Re\theta\le3Q^2\) by
the geometric mean gives
\[
               |\theta|\le 2\Re\theta+6\beta^{-50N_0}Q .
\]
Let \(V_s\) be the mass of the ball of radius
\(\delta\) about \(1\), and \(V\) that of \(\|g\oplus h\|\le\beta+1\).

Both are finite by closedness of the graph and compactness
of the compact factor. Matrix packing as in
Section~\ref{sec:real-reproducing-coefficient} gives
\(V\le\beta^{16N_0}V_s\), since
\(3(\beta+1)^2/\delta\le\beta^8\) for \(\beta\ge10^{12}\).
With radius one, the same packing bounds the mass of
\(\|g\oplus h\|\le Y\) by \(3^{2N_0}Y^{4N_0}V\) for \(Y\ge1\),
since the ball of radius one about the identity lies in the set
of mass \(V\). Integrating in the original Haar measure gives
\[
 \int Q\,ds
 \le 3^{2N_0}2^{9N_0}V\int_1^\infty 9N_0Y^{-5N_0-1}\,dY
 =\tfrac95\,3^{2N_0}2^{9N_0}V\le16^{9N_0}V .
\]
Hence
\[
 V\le8\beta^{16N_0}\left|\int\theta\,ds\right|
                +24\beta^{-34N_0}16^{9N_0}V.
\]
The last coefficient is less than \(1/2\).
Since \(V_{x_0}\le V\) and \(16\beta^{16N_0}\le\beta^{20N_0}\),
this yields
\(V_{x_0}\le b^{60N_0}|\langle v_0,A(x_0)\rangle|\).
This proves the quantitative period bound at the constructed
point for these weights and gaps.
\end{proof}

\subsection{The phase on the common compact block}
\label{sec:real-stationary-phase}

\sectionopening{In this subsection, we compare the phase-corrected
coefficient near the stationary point of Section~\ref{sec:real-stationary-point}
with a positive Gaussian. The decay in
Section~\ref{sec:real-coefficient-decay} controls the complementary
region in the original graph Haar measure.}

Use the stationary graph and \(T_0,\eta_0\) of
Section~\ref{sec:real-stationary-point}. Write
\(x=g-1\), \(x'=h-1\), \(v_0'=\|E_g\|_{\mathrm{HS}}\) and
\(\phi=\sigma-\bar\sigma\). Suppose \(v_0'\le\eta_0\).
For each row index \(i\), the \((i,i)\)-entry of \(gJg^*=J\) gives
\[
 2\Re x_{ii}=2\sum_j|(E_g)_{ij}|^2-\sum_j|x_{ij}|^2 .
\]
The corresponding identity for \(h\) has zero opposite-sign
term. Multiply by \(-\mu_i\) on the first factor and by
the nonnegative compact weight on the second. The weighted
traces cancel in real part by the phase identity. Since
\(B_{\rm wt}\le-\mu_i\le3B_{\rm wt}\) off \(X\), and the weight vanishes on
\(X\), the sum of squared norms of the rows of \(x\) off
\(X\) is at most \(6(v_0')^2\). Apply the same calculation
to the inverse pair, using \(g^{-1}-1=Jx^*J\), to bound
the columns off \(X\).

For \(a_0=g_{XX}\) and \(C=g_{X^c X}\), signed unitarity
gives \(a_0^*a_0-1=-C^*J_{X^c}C\). Its Hilbert--Schmidt
norm is at most \(6(v_0')^2\). Thus \(a_0=up\) with
\(u\) unitary and \(p>0\), and
\(\|p-1\|_{\mathrm{HS}}\le6(v_0')^2\le\eta_0\)
by the scalar inequality \(|t-1|\le|t^2-1|\) for \(t>0\).
Left multiply the pair by the diagonal copy of \(u^*\)
in the graph. Both highest-weight vectors are fixed by
this copy. Indeed the highest-weight coefficient on the
compact group is \(m_\mu\), by the exterior-power
realization.

Its threshold minors include all of \(X\)
for \(t\le0\) and exclude it for \(t>0\), so the determinant
factors on \(X\) cancel. A unit coefficient of value one
forces the unit vector to be fixed. The translating
element has phase zero, and \(\phi\) is additive since the
first predecessor coefficients add. Therefore \(\theta\)
has the same value after the move. The opposite-sign norm
is unchanged and the \(X\)-block of the moved \(g\) is \(p\).
The row and column estimates now give
\(\|g-1\|_{\mathrm{HS}}^2\le13\eta_0^2\). Summing the
unweighted row identity and using the trace equality
gives \(\|h-1\|_{\mathrm{HS}}^2=\|g-1\|_{\mathrm{HS}}^2-2(v_0')^2\).
Thus both distances are at most \(4\eta_0\).

On this moved pair put \(Z=\log g,\ Z'=\log h\).
For the linear functional
\(\mathcal L(P,Q)=\operatorname{tr}(D_\mu P)-\operatorname{tr}(D_H Q)\),
the phase identity applied to group powers gives
\(\mathcal L(g^j-1,h^j-1)=j\phi\) for nonnegative integers \(j\).
Taking finite differences at zero gives
\(\mathcal L(x^i,(x')^i)=0\) for \(i\ge2\).
The log series therefore gives \(\mathcal L(Z,Z')=\phi\).

\indent{\it Triangular factors in the coefficient integral.}
We give the local calculation for either factor of dimension \(d_1\).
Use the weights \(a=\lambda+\rho\) and
\(\mu=a-2\rho_K\), and put \(J_\mu=\int_K|m_\mu(k)|^2\,dk\).
Write \(c(w)\) for the coefficient of its unit highest vector.
For fixed \(w\) in the sign group close to \(1\), consider
\[
 I(p,r_0)=\int_K
 \frac{m_\mu(pk)m_\mu(k^*r_0)}
      {m_a(k^*r_0w^{-1}pk)}\,dk
\]
for block diagonal \(p,r_0\) near \(1\).
The integral is holomorphic there, since all powers are
integral and the denominator compressions are uniformly
invertible on \(K\). In the next computation the neighborhood
can be taken with \(\|p-1\|,\|r_0-1\|<1/10\) and
\(\|w-1\|<1/100\), for which \(r_0w^{-1}p\) is within
distance one of \(1\).

\subsubsection{Triangular coordinates on the sign blocks}
\label{sec:phase-triangular-sign-blocks}

\sectionopening{In this subsection, we write the phase of Section~\ref{sec:real-stationary-point} in triangular coordinates.}

In each sign block take the upper triangular group for
the parameter order. If \(b\) is upper triangular in
these blocks, replacing \(k\) by \(kb\) in the integrand
multiplies it by
\[
 \prod_i |b_{ii}|^{2(\mu_i-a_i)}
 =\prod_{i>j,\ \varepsilon_i=\varepsilon_j}
             |b_{ii}/b_{jj}|^2.
\]
Indeed left lower and right upper factors scale each
initial principal minor by its diagonal product. The
thresholds for \(\mu\) are initial within each block, and
those for \(a\) are initial in the full order.

For block invertible \(t\), write the block QR decomposition
\(tk=k' b\). This induces a diffeomorphism of the compact
flag \(K/\mathbb T\), where \(\mathbb T\) is the diagonal
torus. In left-translated tangent coordinates the strict
lower complex entries of the derivative are the strict
lower entries of \(bYb^{-1}\), with \(Y\) skew-Hermitian
in the sign blocks. The derivative of the upper factor
has no strict lower part.

Each \((i,j)\) depends on
entries \((i',j')\) with \(i'\ge i,\ j'\le j\), so the
complex diagonal factors are \(b_{ii}/b_{jj}\).
Thus the real Jacobian for invariant probability measure
on the flag is the product just displayed. The integrand
is right \(\mathbb T\)-invariant. Change of variables gives
\(I(pt,t^*r_0)=I(p,r_0)\) whenever the nearby arguments
are defined. Differentiating at \(t=1\) in the directions
\(Y\) and \(iY\) and using holomorphy gives separately
\(D_p I[pY]=0\) and \(D_{r_0}I[Y^*r_0]=0\).
Hence \(I\) is constant on the indicated neighborhood.

At \(p=k_0,\ r_0=1\), for \(w=k_0Z_0\) its polar
decomposition, the bottom-layer formula of
Section~\ref{sec:real-coefficient-decay} gives
\(I(k_0,1)=J_\mu c(w)\). Indeed, \(w^{-1}k_0=Z_0^{-1}\), and
\(k^*Z_0^{-1}k=J(k^*Z_0k)J\) since \(k\) commutes with \(J\).
Conjugation by the diagonal sign matrix preserves every principal minor.
The factor
\(m_\mu(k_0k)\overline{m_\mu(k)}\) follows by taking the
highest-vector coefficient of \( \pi(k_0)P_{\pi(k)v}\),
and complementary minors give
\(m_a(D^{-1})=m_{-a}(D)\) for invertible \(D\).
For \(Z=\log w\), let \(U,V\) be its strict lower and
upper equal-sign parts, respectively, and set
\(\widehat w=e^{-U}e^Ze^{-V}\). Evaluating the same
integral at \(p=e^U,\ r_0=e^V\) uses
\(r_0 w^{-1}p=\widehat w^{-1}\). Both triangular factors
have unit diagonal, so the numerator is \(|m_\mu(k)|^2\).
For the moved matrices \(\|w-1\|_{\rm HS}\le4\eta_0\). The
logarithm series and orthogonal truncation of matrix entries give
\(\max(\|U\|,\|V\|)\le\|\log w\|_{\rm HS}\le8\eta_0\).
Thus all these choices
are in the indicated neighborhood. We obtain
\[
 c(w)=\frac1{J_\mu}\int_K
       \frac{|m_\mu(k)|^2}{m_{-a}(k^*\widehat w k)}\,dk.
\]

\indent{\it Uniform quadratic remainder.}
Set \(d_0=\|Z\|_{\mathrm{HS}}\le 2\delta_0\), where
\(\delta_0=4\eta_0\), so \(d_0\le1/20\).
Let \(A_1=-U,A_2=Z,A_3=-V\).
For an orthogonal projection \(P'\) let \(D'\) be the
compression of \(\prod_i e^{A_i}\) to its image and
\(D_0'=\prod_i e^{P'A_iP'}\) there. Put
\(z_0=\max_i\|[P',A_i]\|_{\mathrm{HS}}\le d_0\).
In expanding \(D'-D_0'\) into words in the diagonal
and off-diagonal blocks relative to \(P'\), surviving
words have at least two off-diagonal factors. The degree
two trace is
\[
 S_{P'}=\sum_{i<j}\operatorname{tr}(P'A_i(1-P')A_j)
          +\tfrac12\sum_i\operatorname{tr}(P'A_i(1-P')A_i).
\]
The corresponding degree two matrix has trace norm at
most \((9/2)z_0^2\). Bound two off-diagonal factors of
each longer word in Hilbert--Schmidt norm and the rest
in operator norm. The multinomial sum and the choices of
blocks give the trace-norm bound for the tail
\[
 z_0^2\sum_{m\ge3}6^m d_0^{m-2}/m!
 \le36e^{6d_0}d_0z_0^2\le49d_0z_0^2.
\]
It follows that \(\|D'-D_0'\|_1\le7z_0^2\), where
\(\|\cdot\|_1\) is the trace norm. Also
\(\|(D_0')^{-1}\|\le e^{3d_0}\) and
\(\|(D_0')^{-1}-1\|\le4d_0\).
For \(N=(D_0')^{-1}(D'-D_0')\), the log series gives
\[
 |\operatorname{tr}(\log(1+N)-N)|
 \le \frac{81d_0^2z_0^2}{1-9d_0^2}\le5d_0z_0^2 .
\]
The linear trace of \(N\) differs from
\(\operatorname{tr}(D'-D_0')\) by at most
\(28d_0z_0^2\). Hence
\[
 \det D'=\exp\{\operatorname{tr}(P'\textstyle\sum_i A_i)
                  +S_{P'}+R_{P'}\},\qquad
 |R_{P'}|\le100d_0z_0^2 .
\]

\subsubsection{The compact block estimate}
\label{sec:phase-compact-block-estimate}

\sectionopening{In this subsection, we use the triangular coordinates of Section~\ref{sec:phase-triangular-sign-blocks} to bound the compact phase.}

Apply the preceding determinant expansion to \(P'=kPk^*\), where \(P\) is a threshold projection for \(-a\). Put \(s_0=\|[P,Z]\|_{\mathrm{HS}}\) and
\(h_0=\|[P,k]\|_{\mathrm{HS}}=\|P'-P\|_{\mathrm{HS}}\).
Then \(z_0\le s_0+2d_0h_0\).
Both projections are block diagonal for the signs.
The equal-sign off-diagonals of \(\sum_i A_i\) vanish,
and \(\sum_i|(P'-P)_{ii}|=h_0^2\). Thus its linear trace
differs from \(\operatorname{tr}(PZ)\) by at most \(d_0h_0^2\).
Each factor \(P'A_i(1-P')\) differs from its version at
\(P\) by at most \(2d_0h_0\) in Hilbert--Schmidt norm.
Consequently \(|S_{P'}-S_P|\le18d_0h_0(s_0+2d_0h_0)\).
For \(i\) inside the threshold and \(j\) outside,
the parameter at \(i\) is smaller. For equal signs,
the three cross terms and the \(A_2\) square contribute
\((1-1+1-1/2)|Z_{ij}|^2\) to \(S_P\).
For opposite signs only the square contributes,
giving \(|Z_{ij}|^2/2\) since \(Z^*J+JZ=0\).
Thus \(S_P=s_0^2/4\). The total error in the exponent
for this threshold is at most \(10^5d_0(s_0^2+h_0^2)\).

Write \(Q_b(F)=\sum_{i,j}|b_i-b_j||F_{ij}|^2\) for a weight
\(b\), and \(q=Q_{-a}(Z)/4\). Summing squared commutators
counts each pair once per threshold that separates it.
Also \(\det\widehat w=e^{\operatorname{tr}Z}\).
Multiplying the threshold formulas gives
\[
 m_{-a}(k^*\widehat w k)
   =\exp\{q-\operatorname{tr}(\operatorname{diag}(\mu)Z)-R(k)\},
 \qquad |R(k)|\le10^6\delta_0\bigl(q+Q_\mu(k)+d_1^2\bigr).
\]
For the error use \(|a_i-\mu_i|\le d_1\),
\(Q_{-a}(k)\le Q_\mu(k)+2d_1^2\), and
\(|\operatorname{tr}(\operatorname{diag}(a-\mu)Z)|\le d_1^{3/2}d_0\).
The remainder is a sum of the local log series above.
For unitary \(k\), the squared modulus of each compressed
determinant at \(P\) equals \(\det(1-Pk^*(1-P)kP)\) on
the compressed space. Since \(1-t\le e^{-t}\), the product
gives \(|m_\mu(k)|^2\le e^{-Q_\mu(k)/2}\), using equality
of the two off-diagonal squared masses of a unitary matrix.
The compact Weyl dimension formula gives
\(J_\mu^{-1}=\dim V_\mu\le(1+6\beta)^{d_1^2}\le e^{5000T_0}\),
since the weights have modulus at most \(3\beta\) and
\(d_1^2\log\beta\le T_0\).

Put \(W_0=q+Q_\mu(k)+d_1^2\).
If \(W_0\le50000T_0\), then
\(\delta_0T_0\le10^{-12}\) gives \(|R|\le1/20\) and
\(|e^R-1|\le1/10\). Otherwise \(W_0>50000T_0\ge50000\)
and the bound on \(R\) gives \(|R|\le10^{-6}W_0\). Summing
the absolute values of the positive-degree terms in the exponential
series gives \(|e^R-1|\le e^{|R|}-1\le e^{W_0/4}\), so
\[
 |m_\mu(k)|^2 e^{-q}|e^R-1|
 \le e^{-q-Q_\mu(k)/2+W_0/4}\le e^{-5100T_0}.
\]
Indeed the exponent is at most \(-W_0/4+d_1^2/2\).
After integration and division by \(J_\mu\),
\[
 \left|c(w)e^{-\operatorname{tr}(\operatorname{diag}(\mu)\log w)}
                       -e^{-q}\right|
 \le e^{-q}/10+e^{-100T_0}.
\]
Apply the estimate to \(g\) and \(h\). The sign-Lie identity
makes both logarithmic traces imaginary. Conjugating the
compact coefficient reverses its phase and preserves its
real Gaussian. The phase in \(\theta\) is therefore
\[
  \exp\!\left(\operatorname{tr}(D_\mu Z)
          -\operatorname{tr}(D_H Z')-\phi\right)=1.
\]
The product of the two phase-corrected
coefficients differs from the product \(a_1a_2\) of the
two numbers \(e^{-q}\) in \([0,1]\) by at most
\(0.21a_1a_2+4e^{-100T_0}\).
Thus \(|\theta|-2\Re\theta\le12e^{-100T_0}\), as required
at the original graph element.

\subsection{The real Bessel functions for the construction}
\label{sec:real-marked-filling}

\sectionopening{In this subsection, we apply the real estimates to the parameters of
Section~\ref{sec:datum}. This constructs the discrete-series vectors
and verifies the normalization and truncation bounds required for
Section~\ref{sec:full-orbital-blocks}.}

Return to the datum of Section~\ref{sec:datum}. The rank of
the ambient form is \(r\), with
\(k=r+1-m\) and \(\varrho=m/2-1\). Thus \(n=\varrho+1=m/2\)
and \(r=k+2\varrho+1\). Put \(a=xb\), \(\Delta=2r^8\), and
\(B_\tau=\kappa_\tau+\Delta/2\) at an infinite place \(\tau\).
The character of lift number \(j\) has weight \(B_\tau+j\Delta\).
The initial parameter, in its numbered coordinate order,
has \(k-1\) positive entries followed by \(m\) alternating signs
starting with \(+\).

In this order the parameter entries are
\[
 \lambda_i=k/2-1-i\quad(0\le i<k-1),\qquad
 \lambda_{k-1+t}=-B_\tau-(t+1/2)\Delta\quad(0\le t<m).
\]
Here the tail entry is the half gap minus the weight of lift
\(t+1\). The signs have signature \((k+\varrho,\varrho+1)\).
To realize them in the original stage space one can take the unit
eigenvectors \((e_j\pm e_j^*)/\sqrt2\) of the pair blocks,
along with the positive basis and \(z\), and permute
the two sign sets into the parameter order. This gives a Euclidean
unitary basis with prescribed signs for the original Gram matrix \(J\).
The local construction applies using any such chosen
basis throughout.

The least weight satisfies \(\min_\tau B_\tau\ge\Delta^2\) for
large \(x\), since \(\kappa_\tau\ge xb/x_1\),
\(\log b-\log x_1\ge x^\epsilon-x^{\epsilon/6}\), and
\(\Delta=2r^8\) grows polynomially in \(x\).
The last initial entry is \(1-k/2\), and the first tail entry is
\(-B_\tau-\Delta/2\), so the entire parameter decreases strictly.
Within the tail the smallest opposite-sign gap is \(\Delta\).
For a mixed pair the initial coordinate is positive in signature,
so the first possible opposite tail position is \(t=1\). Its
distance from the last initial entry is
\(1-k/2+B_\tau+3\Delta/2\ge\Delta\).
Thus the smallest opposite-sign gap equals \(\Delta\) at every
infinite place. The doubled parameter has the right parity
because \(k\) is even and \(r\) is odd. Choose the discrete
series with this Harish--Chandra parameter and a unit highest
vector of its lowest compact type \cite[Theorems~8.1 and 8.5, pp.~41--42]{AS}.

On the derived group use the parameter with its mean removed.
Its lowest compact type has the central restriction given by the
original integral weight. The character of the scalar circle of
weight \(\sum_i\lambda_i\) agrees on the common scalar subgroup,
since the root corrections have sum zero, so the representation
descends from the product to the unitary group.
The signed root count in this order gives
\[
 \mu_i=\lambda_i-\tfrac12\left(
          \sum_{j>i}\varepsilon_i\varepsilon_j
          -\sum_{j<i}\varepsilon_i\varepsilon_j\right).
\]
On the initial block the tail terms cancel in pairs and the
count is \(k-2-2i\), so \(\mu_i=0\). At tail position \(t\),
the initial block contributes \(-(k-1)(-1)^t\), and the tail
contributes \(-(-1)^t\). Hence
\[
 \mu_{k-1+t}=-B_\tau-(t+1/2)\Delta+\tfrac k2(-1)^t.
\]
Consecutive differences are \(\Delta+k(-1)^t>0\).
Since \(\Delta\ge k\) and \(m\Delta+k/2\le2B_\tau\), these
entries lie in \([-3B_\tau,-B_\tau]\).

On \(X_0\) the column counts at heights \(0<h<k\)
are zero except at \(h=1\), where the count \(a_\tau\) is
\(a\) if \(\tau\in\Phi_a\) and zero otherwise. The determinant
weight on \(X_0\) is zero. Thus the compact type is the symmetric
power \(\operatorname{Sym}^{a_\tau}(\mathbb C^k)\) of the defining
representation of \(U(k)\). Its highest weight is \((a_\tau,0,\ldots,0)\), and its doubled
parameter entries are
\(2a_\tau\mathbf1_{j=0}+k-1-2j\) for \(0\le j<k\).

Match the other \(k-1\) coordinates with the zero-weight
positive block in the ambient parameter. The remaining
compact coordinate is the highest-weight coordinate, also
when \(a_\tau=0\). Use exactly the Schur matrices
\(S_\tau(h)\) on the value space. If \(Q_\tau=B_0^*B_0\)
is their invariant Weyl Gram matrix, set
\(\pi_H(h)=B_0S_\tau(h)B_0^{-1}\).
Since \(S_\tau\) is defined over the reals and
\(h^{-1}=\bar h^{\mathsf T}\), invariance gives
\[
 Q_\tau S_\tau(h)=S_\tau(h^{\mathsf T})^{\mathsf T}Q_\tau,
 \qquad
 B_0^*\pi_H(h)=S_\tau(h^{\mathsf T})^{\mathsf T}B_0^* .
\]
Let \(j_*\) be the chosen highest-weight coordinate, and let
\(\xi\) be the highest monomial \(e_{j_*}^{a_\tau}\) in the
symmetric power (\(\xi=1\) if \(a_\tau=0\)). The unit highest vector
and its coefficient in the Bessel value space are
\[
 v_0=\frac{B_0\xi}{\|B_0\xi\|},\qquad
 c_H(h)=\frac{\xi^*Q_\tau S_\tau(h)\xi}{\xi^*Q_\tau\xi}
        =h_{j_*j_*}^{a_\tau}.
\]
Indeed the diagonal torus acts unitarily for \(Q_\tau\), so distinct
weights are orthogonal. The coefficient of \(e_{j_*}^{a_\tau}\) in
\((he_{j_*})^{a_\tau}\) gives the last equality. Thus the stationary
integrand uses \(\overline{c_H(h)}\) for this same Schur representation.

We check the required sizes uniformly in \(\tau\).
For large \(x\), \(x^5\le r<3x^5,\ \sqrt{x}/2\le n<\sqrt{x},\ k\ge x^5/2\),
and \(4000n\le f\le x^{12}\), the last pair from
\(f=3L_0m^{22}\varphi(q_0)\). Thus \(k\ge n^7\),
\(k-3n-6\ge10^{10}n^4\), \(\varrho\ge10^8\), and \(8nr\le b\).
Since \(x_1\ge x\) and \(r^8\le b\) eventually,
\(\kappa_\tau\in\{b,xb/x_1\}\) gives \(B_\tau=\kappa_\tau+r^8\le2b\).
Thus \(3B_\tau+a_\tau\le(6+x)b\le b^2\) for all large \(x\).
The two doubled parameters have magnitude at most
\(L_*=3a\), since \(b+(m+1)\Delta\le a\) eventually
and \(2a+k\le3a\). Put \(\beta=b^3\) and \(A_*=r^7\).
Then \(r+k+L_*\le\beta\), \(\beta\ge10^{12}\), and
\[
 \begin{split}
 500r^3+A_*&\le\Delta,\qquad 3r^3f\le A_*,\\
 10^{30}\bigl((r+k)^2\log\beta\bigr)^3
       &\le10^{30}108^3 x^{30+3/100}\le r^7.
 \end{split}
\]
The second line uses \(\log b\le x^{1/100}\) from the chosen
scale bounds. The compact parameter has only plus signs.
These size comparisons need one threshold in \(x\), since
the weights have the same upper and lower bounds at all places.

Apply Proposition~\ref{prop:marked-real-period} with
\((B_{\rm wt},a,A_{\rm gap})=(B_\tau,a_\tau,r^7)\) at each
place \(\tau\), including the places where \(a_\tau=0\).
Write \(x_\tau^\star\) for the stationary point in the original
stage group, using the prescribed sign basis at \(\tau\).
When \(a_\tau=0\), all interpolation nodes satisfy
\(\mu_i\le-B_\tau<0\), hence \(P(0)>0\), and \(c_H(h)=1\).
Let \(\deg_\tau\) denote the formal degree in \(dg\). With the induced
quotient Haar, Section~\ref{sec:real-cut-error} and the proposition
give at this point
\[
 \deg_\tau\|A_\tau\|_2^2
  =\langle v_{0,\tau},A_\tau(x_\tau^\star)\rangle
  \ge b^{-60(r+k)^2}V_{x_\tau^\star}>0 .
\]
For \(D_\tau=\dim \pi_H\), take
\[
 W=\frac{\sqrt{D_\tau}}{\|A\|_2}A,\qquad
 F=\chi W,\qquad R_b=b^{6n+2}.
\]
For the point \(g_{0,\infty}\) of Section~\ref{sec:period-real-point}
one has \(T_\tau(g_{0,\tau})=x_\tau^\star\) for every \(\tau\).
Consequently
\[
 \max_\tau\ell_\tau\bigl(T_\tau(g_{0,\tau})\bigr)\le b^{5n}.
\]
Tensoring the absolutely convergent graph integrals at these points
gives
\[
 \prod_\tau \deg_\tau\|A_\tau\|_2^2
 \ge b^{-60f(r+k)^2}\prod_\tau V_{x_\tau^\star}>0.
\]
In particular,
\[
 \langle w,W(g_{0,\infty})\rangle
   =\prod_\tau
    \frac{\sqrt{D_\tau}}{\|A_\tau\|_2}
         \langle v_{0,\tau},A_\tau(x_\tau^\star)\rangle>0 .
\]
Here \(w=\bigotimes_\tau v_{0,\tau}\) and \(W(g_{0,\infty})\)
is in the product chart of the measured period. These choices
precede the orbital tests.

Section~\ref{sec:real-cut-error} applies this period bound to
the graph average and gives the decay majorant, the cut support, and
\(\|F-W\|_2^2\le b^{-2r^3f}D_\tau\).
Conjugation by the compact group \(H=U(k)\) preserves \(du\),
so its left action on \(U\backslash G\) preserves the quotient
\(d\bar g_\tau\) of the orbital Haar \(dg_\tau\) by \(du\).
Since the length used in \(\chi\) is left \(UH\)-invariant,
\(W(hy)=\pi_H(h)W(y)\) and \(F(hy)=\pi_H(h)F(y)\).
Both Gram matrices therefore commute with \(\pi_H\).
In the Weyl isometry above this representation is
\(\operatorname{Sym}^{a_\tau}(\mathbb C^k)\). Its monomial weight
spaces are one-dimensional, and the operators \(X_i\partial_{X_j}\)
of the complexified Lie algebra connect all monomials of degree
\(a_\tau\). This proves irreducibility, also for \(a_\tau=0\).
Schur's lemma and the traces now give
\[
 \int WW^*\,d\bar g_\tau=1_{I_\tau},\qquad
 \int FF^*\,d\bar g_\tau=s_\tau 1_{I_\tau},\qquad
 (1-b^{-r^3f})^2\le s_\tau=\frac{\|F\|_2^2}{D_\tau}\le1 .
\]
The lower bound uses \(\|F\|_2\ge\|W\|_2-\|F-W\|_2\), and the
upper bound follows from \(0\le\chi\le1\).

Section~\ref{sec:real-reproducing-coefficient} gives the right
reproduction, the correlation norm \(D_\tau^2/\deg_\tau\), and
the norm-ball bound. Sections~\ref{sec:real-formal-degree-factor}
and~\ref{sec:real-frame-haar} give the other degree and Haar
bounds of Section~\ref{sec:real-bottom-mass}. The shell estimates
and Lemma~\ref{lem:marked-real-frequency-decay} apply to \(F\) on this
graph. The factor \((\pi/|\xi|)^{k-5}\) is the Plancherel factor
from the \(k-5\) free complex coordinates of
Section~\ref{sec:real-floor-plancherel}.
Section~\ref{sec:real-bottom-mass} gives the bottom bound
\((1+b^{-10})/\deg_\tau\). Lemma~\ref{lem:marked-real-frequency-decay}
gives the character-weight and unequal-entry bounds, and
Section~\ref{sec:real-common-descent} defines the real forms \(J_j\)
and proves integration in a common entry. The derivative estimates
for the full vector \(W\) in Section~\ref{sec:real-differentiated-schur}
enter the comparison of nonzero-discriminant orbits in
Section~\ref{sec:real-open-cell-contraction}.
For a terminal element \(\sigma\) fixing the hyperbolic pairs,
the top form at heads with the same last entry \(h(z,\sigma z)\)
is the inner product of the two actual normalized orbital
tensors of \(F\) at \(\sigma\), by
Section~\ref{sec:real-terminal-gram}. This uses the
normalizing dimension \(D_\tau\) for both tensors.

\begin{lemma}[Comparison of the formal degrees]\label{lem:formal-degree-weight-comparison}
Put \(B=\log b\), \(S_\kappa=\sum_\tau\log\kappa_\tau\), and
\(W_*=\log\mathcal B+f\sum_{l=k-1}^{r-1}\log(l!)\). Then
\[
       |W_*-kmS_\kappa|\le m^2fB+fmr/x .
\]
\end{lemma}

\begin{proof}
The formal-degree product \(\mathcal B\) used in
Section~\ref{sec:full-orbital-blocks} is
\[
 \log\mathcal B=\sum_\tau\sum_{l=k-1}^{r-1}
        \left(\sum_{i<l}\log(\lambda_i-\lambda_l)-\log(l!)\right).
\]
In Section~\ref{sec:real-formal-degree-factor} the factors at
\(l<k-1\) cancel because \(\lambda_i-\lambda_l=l-i\).
At \(l\ge k-1\), each initial entry contributes a difference
at least \(\kappa_\tau\), and every tail-to-tail difference is
at least one. All differences in this range are at most
\(\kappa_\tau+m\Delta+k\). Uniformly in \(\tau\),
\[
   \frac{x(m\Delta+k)}{\kappa_\tau}\le\frac{(m\Delta+k)x_1}{b}=o(1),
\]
since \(m\Delta+k=O(x^{81/2})\) and
\(\log b-\log x_1\ge x^\epsilon-x^{\epsilon/6}\).
Using \(\log(1+y)\le y\) for \(y\ge0\) gives
\[
 \begin{split}
 m(k-1)\sum_\tau\log\kappa_\tau
 &\le \log\mathcal B+f\sum_{l=k-1}^{r-1}\log(l!)\\
 &\le mr\sum_\tau(\log\kappa_\tau+1/x).
 \end{split}
\]
Subtract \(kmS_\kappa\) and use \(r=k+m-1\) and
\(0\le S_\kappa\le fB\). This gives
\[
 -mfB\le -m S_\kappa\le W_*-km S_\kappa
             \le m(m-1)S_\kappa+fmr/x\le m^2fB+fmr/x .\qedhere
\]
\end{proof}

\subsection{The differentiated Schur bound}
\label{sec:real-differentiated-schur}

\sectionopening{In this subsection, we transfer derivatives of the full Bessel vector to its
discrete-series matrix coefficient by the reproducing identity of
Section~\ref{sec:real-reproducing-coefficient}. This bounds the
derivatives needed in Section~\ref{sec:real-open-cell-contraction}.}

Use the full vector \(W\) just constructed. Let \(\pi\) be
the ambient discrete series and \(v\) its chosen unit highest
vector, with \(c(g)=\langle v,\pi(g)v\rangle\).
Right translation \(R\) on the quotient satisfies
\(\langle W,R(g)W\rangle=D_\tau c(g)\), by
Section~\ref{sec:real-reproducing-coefficient}. In particular
the map \(\pi(g)v\mapsto R(g)W\) extends linearly to a bounded
intertwiner \(T\) with \(\|Tz\|_2=\sqrt{D_\tau}\|z\|\):
the two Gram matrices on each finite list of translates agree
up to \(D_\tau\), and \(v\) is cyclic by irreducibility.

The right derivatives of \(W\) are bounded continuous Bessel
functions. Indeed take a smooth compactly supported
approximate identity and average it on the left and right by the
compact matrix coefficient projector onto \(\mathbb C v\).
The latter projector comes from Schur orthogonality and
multiplicity one of the lowest compact type.
Choosing a bump sufficiently near the identity yields a smooth
compactly supported kernel \(\xi\) with
\(\pi(\xi)v=\alpha v\), \(\alpha\ne0\), since the bump operator
tends strongly to the identity.

Thus \(R(\xi)W=\alpha W\),
first in \(L^2\) by the intertwiner and then everywhere by
continuity. The change \(u=\exp(tX)g\) in
\(\int \xi(g)W(y\exp(tX)g)\,dg\) puts each derivative on
the smooth kernel by left translation. Compact support and the
boundedness of \(W\) justify all differentiations, and left
Bessel covariance is preserved. Complex linear extension supplies
the derivatives in the matrix units. For a word \(\omega\) of
length \(j\) denote the resulting function by \(W_\omega\)
and the corresponding derivative of \(v\) by \(v_\omega\);
the intertwiner gives \(W_\omega=Tv_\omega\) in \(L^2\).

For the derivative norm, let \(V\) be the lowest compact type.
The span \(V_j\) of derivatives of order at most \(j\) of
vectors of \(V\) is the image of
\(\bigoplus_{i=0}^j\mathfrak g_\mathbb C^{\otimes i}\otimes V\)
as a compact module. Its highest compact weights have
coordinates of modulus at most \(L_*+r+j\). Indeed weights of
\(V\) lie between the smallest and largest highest-weight
coordinate in each unitary block, by the highest-weight order
and Weyl permutation symmetry.

An adjoint weight changes
any coordinate by at most one. The compact Casimir on \(V_j\)
is therefore at most \(r(L_*+2r+j+1)^2\) and the ambient
Casimir has absolute value at most \(r(L_*+2r+j+1)^2\) by
the same eigenvalue formulas as in
Section~\ref{sec:real-reproducing-coefficient}.
For a real Hilbert--Schmidt orthonormal Lie basis \(X_i\)
adapted to the compact and noncompact summands, the energy
identity gives
\[
 \sum_i\|d\pi(X_i)u\|^2\le
           3r(L_*+2r+j+1)^2\|u\|^2,\qquad u\in V_j .
\]
For \(X_i,X_l\in\mathfrak u(J)\), the relation \(X_i^*=-JX_iJ\)
and cyclicity give
\(\operatorname{tr}(X_i^*X_l)=\operatorname{tr}(X_l^*X_i)\).
Thus \((X_i)\) is also a complex Hilbert--Schmidt orthonormal basis.
A matrix unit \(E_{ab}\) of the original basis becomes
\(M_{ab}=p_v^{-1}E_{ab}p_v=\sum_i c_iX_i\) in the signed basis
of Section~\ref{sec:period-real-point}. Since \(p_v\) is Euclidean
unitary, \(\sum_i|c_i|^2=1\), and hence
\[
 \|d\pi(M_{ab})u\|^2
 \le\sum_i\|d\pi(X_i)u\|^2
 \le3r(L_*+2r+j+1)^2\|u\|^2 .
\]
As \(\sqrt{3r}(L_*+2r+2n^2+1)\le b^5\) for large \(x\),
iteration gives \(\|v_\omega\|\le b^{5j}\) for \(j\le2n^2\).

For \(\phi\in L^1(dg)\cap L^2(dg)\), Schur orthogonality
and Cauchy--Schwarz against an arbitrary unit vector give
\(\|\pi(\phi)v_\omega\|^2\le
 \|\phi\|_2^2\|v_\omega\|^2/\deg_\tau\).
Applying \(T\) gives the desired inequality in the quotient normalization,
\[
 \int_{U\backslash G}
       \left\|\int_G\phi(g)W_\omega(yg)\,dg\right\|^2\,d\bar g(y)
       \le b^{10j}\frac{D_\tau}{\deg_\tau}\int_G|\phi(g)|^2\,dg ,
                  \qquad j\le2n^2 .
\]
For integrable \(\phi\) outside \(L^2\) the right side is infinite.
The estimate concerns the same \(W\), chosen coordinates and
Haar measure as the real sieve forms.

\section{Adelic orbital comparison and definite masses}
\label{sec:adelic-orbital-assembly}

\sectionopening{In this section, we combine the local estimates on the open rational cell and average at a
fixed final vector. The definite lattice masses give the normalization
for the moment comparison of Section~\ref{sec:full-orbital-blocks}.}

\subsection{The open rational cell and its finite restrictions}
\label{sec:open-rational-cell}

\sectionopening{In this subsection, we identify the rational terms in the open range of the Petersson
expansion. The fixed finite levels from Section~\ref{sec:period} allow
the local integrals to be compared on one set of places in
Section~\ref{sec:adelic-common-entry-sum}.}

Use the cut Bessel vector, the class representatives
\(c^\circ\in H(\mathbb A_{F_0,f})\) for \(H=U(X_0)\), and the
two class cutoffs of Section~\ref{sec:period}. In the rational basis
\(f=(e_0,\ldots,e_{\varrho-1},z)\) we have \(h(z,z)=-N\).
The rational-coordinate matrix \(h_c^{\rm rat}=\iota_{\rm rat}(c^\circ)\)
fixes all the hyperbolic pairs and \(z\). In products with the
glued matrix \(h_c=B_{\pi_N}^{-1}h_c^{\rm rat}B_{\pi_N}\), write
\(g_L=B_{\pi_N}^{-1}g_fB_{\pi_N}\) for a rational isometry \(g\),
as in Section~\ref{sec:period-finite-coordinates}. For an original
term \(\gamma\) above \((c,s,c')\), the matrix
\(h_c^{-1}\gamma_Lh_{c'}\) is integral with reduction \(s\). Thus the leading minors
\(\Delta_a\) of \((h(f_i,\gamma f_j))_{i,j}\) are integral.
Their reductions are the minors of \(s\), since the Gram
entries are unchanged by the class translates. We take
the order-zero minor to be \(1\). An open residue \(s\) has nonzero proper
minors at one chosen place above \(p_0^{\rm lev}\),
so \(\Delta_a\ne0\) for \(1\le a\le\varrho\).

\begin{lemma}[Factorization on the open rational cell]
\label{lem:open-rational-factorization}
If \(\Delta_a\ne0\) for \(1\le a\le\varrho\), with \(\Delta_0=1\),
the rational term \(\gamma\) has a unique factorization
\begin{equation}\label{eq:open-rational-cell}
 \gamma=u^{-1}w_p\sigma u',\qquad
 u,u'\in U(F_0),\quad \sigma\in G_\varrho(F_0),\qquad
 p_a=\Delta_{a+1}/\Delta_a\quad(a<\varrho).
\end{equation}
Here \(w_pe_a=p_ae_a^*\) and \(w_pe_a^*=\bar p_a^{-1}e_a\).
\end{lemma}

\begin{proof}
At stage \(a\), use the
residual isometry \(\gamma_a\) fixing the earlier pairs.
The pivot \(c_a=h(e_a,\gamma_a e_a)\) is nonzero by the
leading-minor elimination. More explicitly, clearing
one pair makes the upper-left basis-product block
\((c_a)\) and its complementary block the next
basis-product matrix. Since lower and upper
unitriangular changes preserve each leading minor,
the successive residual pivots are the stated ratios.
Let \(m=e_a^*+\omega+\mu e_a\) be isotropic. The matrix
\(N_m\) sending \(e_a^*\) to \(m\), fixing \(e_a\), and
sending a remaining vector \(v\) to \(v-h(\omega,v)e_a\)
is unitary, because \(\mu+\bar\mu+h(\omega,\omega)=0\).

Use \(A=N_{c_a^{-1}\gamma_a e_a}\) and put
\(L=(A w_a(c_a))^{-1}\gamma_a\), which fixes \(e_a\).
Then \(B=N_{L^{-1}e_a^*}\) is defined and \(LB\)
fixes the pair. Iterate \(\gamma_a=A w_a(c_a)(LB)B^{-1}\),
using commutation of the earlier swaps with the later
stage group. Each step is integral when its pivot is a
unit and the residual matrix is integral: the column
of \(N_m\) is integral, as is its inverse by the integral
inverse Gram. This uses the column relation also above
two. Uniqueness follows by writing two decompositions as
\(w_p\sigma v'=v w_q\tau\). Compare the \(e_a\) columns.

Their dual components have indices respectively at most
and at least \(a\), since \(v'\) and the lower dual block
of \(v\) are unitriangular. The \(a\)-th coefficient gives
\(p_a=q_a\), and comparison of all components gives
\(v e_a^*=e_a^*\). In the order consisting of the isotropic block,
the middle space and the dual block, the Bessel unipotent has the form
\[
 v=\begin{pmatrix} A&C&B\\0&I&D\\0&0&A^{-*}\end{pmatrix}.
\]
Fixing the dual columns gives \(B=D=0\) and \(A=I\).
Pairing the middle columns with the fixed dual columns then gives
\(C=0\). Hence \(v=1\), so \(v'\in U\cap G_\varrho=\{1\}\)
and \(\sigma=\tau\).
\end{proof}

The link is the graph \(x\mapsto\sigma^{-1}x\) in the
middle quotient. Indeed a vector orthogonal to both
isotropic bases of \(w_p\sigma\) has zero hyperbolic
coordinates, so the two quotient coordinates are \(x\) and
\(\sigma^{-1}x\). The left and right unipotents change
the corresponding lifts by the original isotropic span,
which preserves the quotient graph. Thus the link
determines \(\sigma\), independently of \(p\).
Put
\[
 t=h(z,\sigma z),\qquad v_0=\operatorname{pr}_{X_0}(\sigma z),
 \qquad \sigma z=v_0-(t/N)z .
\]
If the discriminant is nonzero then
\(t\bar t-N^2\ne0\): the two last minors are
\(\Delta_\varrho=\prod p_a\) and
\(\Delta_{\varrho+1}=t\Delta_\varrho\), and the discriminant
is \(\Delta_\varrho\bar\Delta_\varrho(t\bar t-N^2)\).
Moreover \(N h(v_0,v_0)=t\bar t-N^2\). By definiteness
\(|\Phi_v(t)|>N\) at every infinite place \(v\).
In particular the full minor is nonzero on this set.

All nonzero cut terms have bounded archimedean minors,
with a bound depending on the chosen Bessel data.
Indeed nonzero cut support at \(y,\gamma^{-1}y\) implies
\(\gamma\in U_\infty\mathcal K\mathcal K^{-1}U_\infty\).
Leading minors are unchanged by the two unipotents,
and the original rational basis has fixed scales in this
compact box. Thus each \(\Delta_a\) lies in a finite
set of bounded algebraic integers. Take \(\mathcal E\) to
consist of the quotients of the nonzero elements of this
set with \(1\) included.

All \(p_a,t\) belong to
\(\mathcal E\). Choose the finite conjugation-stable set
\(T\) containing the places above \(M,N\), those needed
for the lattice basis and its inverse, and the nonunit
places of every element of \(\mathcal E\). We may also
include all places below the nonintegral entries of the
finite list of class representatives and their inverses.
This single \(T\) works for all terms of the chosen data and may
depend on \(x\). Its finite Haar factors cancel on dividing
\(\mathcal I_i^T\) by \(\mathcal R_Z\) in
Section~\ref{sec:finite-bottom-normalization}. For \(s>1\), the Euler
product for ideals supported at \(T\) and the complementary product
are each bounded by \(\zeta_E(s)\), independently of \(T\).
The global additive covolume factor in Fourier inversion is \(v_E\).

In the lattice coordinates put \(L_{\rm lat}=B_{\pi_N}\) and
\(g_L=L_{\rm lat}^{-1}g_f L_{\rm lat}\). The three factors in the original
class lattice then give
\[
 h_c^{-1}\gamma_L h_{c'}
   =(h_c^{-1}u_L h_c)^{-1}\,w_p
     (h_c^{-1}\sigma_L h_{c'})(h_{c'}^{-1}u'_L h_{c'}).
\]
Here \(L_{\rm lat},h_c,h_{c'}\) fix the hyperbolic pairs and commute with \(w_p\).
Off \(T\), the left side is integral and every proper pivot is a
unit in each local field factor. The integral elimination and
uniqueness above therefore put both conjugated unipotents and
the displayed middle factor in the integral group there.

For a cut block above \((c,s,c')\) and link \(L\),
let \(B^\flat_{c,c',s,L}\) be the sum of its terms on
the nonzero-discriminant set. For each nonzero block
the link determines its \(\sigma\). The preceding
factorization gives its expansion over \(p\in\mathcal
E^\varrho\) of the orbit sums
\(\operatorname{CosAmp}^\flat_{c,s,c'}(w_p\sigma)\).
These are the sums of the cut amplitudes over the
unipotent factors in \eqref{eq:open-rational-cell} and
the original class and residue conditions.
The block parameters give distinct tuples
\((c,c',s,\sigma)\). The absolute unfolding for these
orbit sums is given below.

\subsection{Comparison at a fixed final vector}
\label{sec:fixed-tail-comparison}

\sectionopening{In this subsection, we express the terminal pairing at a fixed last vector as a finite
class sum. The calculation matches the source level of
Section~\ref{sec:ramified-source} with the measures of
Section~\ref{sec:finite-selected-bottom}.}

Take \(Z\) for the prescribed fixer from
Section~\ref{sec:finite-selected-bottom}, and use its
\(W_Z,\mathcal R_Z\) and \(V_p(l)\) from
Section~\ref{sec:finite-common-entry}. Let
\[
 A_{c,c',\sigma}=h_c^{-1}\sigma_L h_{c'},\qquad
 l=(A_{c,c',\sigma})_T^{-1},\qquad
 X=(w_p A_{c,c',\sigma})_T .
\]
Write \(\widetilde s\in Z\) for the finite projection of
an integral lift of \(s\in\Gamma_{\mathcal A}\).
Here and below the products with class representatives
are read in the original lattice.

The rational unipotents whose class translates are
integral outside \(T\), modulo the class unipotent
subgroup \(U_c\) of Section~\ref{sec:period}, project
bijectively to the coset set \(U(A_T)/(U(A_T)\cap J)\), with representatives identified
by right multiplication by \(U(A_T)\cap J\). For surjectivity, prescribe the desired coset of
\(U(A_T)\cap J\) at \(T\) and the integral unipotent
subgroup at each place outside \(T\). Transport by
\(L_{\rm lat}h_c\) carries their product to a nonempty open set in
the rational coordinates of \(U(\mathbb A_{F_0,f})\). It meets \(U(F_0)\) by the
additive approximation in exponential coordinates of
Section~\ref{sec:period-quotient-measure}. Indeed the exponential
and logarithm are mutually inverse polynomials over \(F_0\)
and preserve integral coordinates at almost all places.
Thus the rational element has exactly the required class
translate at \(T\) and is integral after translation outside it.
For injectivity,
two representatives with the same \(T\)-coset already
have quotient in the principal level outside \(T\)
and hence everywhere, since \(M\) is a unit off \(T\).
Their rational lifts therefore differ on the right by an element of
\(U_c\), and every such change preserves the local coset.

For a contributing \((c,c',\sigma)\), the middle factor
\(A_{c,c',\sigma}\) is integral off \(T\) by the preceding
elimination of one representative. For each head tuple
\(p\in\mathcal E^{\varrho}\), \(w_p\) and its inverse are integral
there. Integrality off \(T\) of the relative term is therefore
equivalent to integrality of both conjugated unipotents, by
multiplication in one direction and integral elimination in the other. Thus the remaining residue condition on
finite representatives \(a,b\) from these two unipotent
quotients is \(a^{-1}Xb\in\widetilde s J\).
Only finitely many pairs qualify: if \(a^{-1}Xb\) lies
in this compact coset, the same fixed-pivot block
elimination that proves compact support in
Section~\ref{sec:finite-common-entry} bounds \(a,b\).
Quotienting bounded sets by the open subgroup gives a
finite set.

Within each such pair of cosets, sum the two cutoffs over
the two rational class unipotent subgroups. If the central
term is \(\theta\), the terms in those two sums have the
form \(v^{-1}\theta v'\). Change the real integration
variable by \(v\); the two corresponding cutoff sums are
one by the partitions in Section~\ref{sec:period}.
The class-unipotent character is trivial on these
subgroups, since their finite translates are in the
principal level.

This also follows from the adelic
product of the additive characters at rational points.
The remaining real tensor is
\(I^\infty(w_p\sigma)=\int F(y)\otimes F((w_p\sigma)^{-1}y)^*\,dy\),
with its real embeddings understood. Its integral
is absolute: with fixed nonzero pivots the two cut-support
conditions put the unipotents in bounded sets by the same
block calculation on the compact ratio box. Tonelli applied to
\[
 \int \|F(y)\|\,\|F((w_p\sigma)^{-1}y)\|\,dy<\infty
\]
and the two nonnegative cutoff partitions justifies summing the
translates before integrating the tensor.

\subsubsection{Cancellation of the class phases}
\label{sec:fixed-vector-class-phases}

\sectionopening{In this subsection, we combine the finite and real actions before applying the common-entry sum of Section~\ref{sec:adelic-common-entry-sum}.}

For a rational unipotent \(u\) put \(z_u=-\mathfrak s(u)/(M\delta_E)\).
Our adelic character has archimedean component
\(z\mapsto\exp(-2\pi i\sum_v 2\Re z_v)\), with one embedding
from each complex pair. Since the real root entry is
\(\alpha_{\Psi,v}=-2\pi i\Phi_v(z_u)\), we have
\[
 \chi_\infty(u)=\prod_v\exp(\alpha_{\Psi,v}-\bar\alpha_{\Psi,v})
     =\exp(-2\pi i\operatorname{Tr}_{E/\mathbb Q}z_u)
     =\psi_{E,f}(z_u)^{-1}.
\]
This is the covariance character of \(F\). The graph average
uses its inverse. Conjugation by \(h_c\) preserves the root character because
the definite Levi fixes the vectors \(e_i,z,e_i^*\)
that give its simple-root entries. At places off \(T\) the root
coordinate divided by \(M\) is integral and the finite
character has conductor the integer ring.
If \(u_1,u_2\) lift \(a,b\), respectively, the substitution
\(y=u_1^{-1}y'\) in the tensor at \(u_1^{-1}w_p\sigma u_2\)
therefore contributes
\(\chi_\infty(u_1)^{-1}\chi_\infty(u_2)
  =\chi(a)\overline{\chi(b)}\).

Let \(\alpha_1,\alpha_2\in Z/J\) with
\(\alpha_1\alpha_2^{-1}\in\widetilde s J\), and put
\(a_p^T=\mathcal D_\varrho^T(p,1)\).
The corresponding coefficient of \(V_p(l)\) is
\(a_p^T\int\chi(y\alpha_1^{-1})
                    \overline{\chi(X^{-1}y\alpha_2^{-1})}\,d\mu(y)\).
Replace \(y\) by \(y\alpha_2\). The integrand becomes
\(\chi(y\widetilde s^{-1})\overline{\chi(X^{-1}y)}\),
using normality of \(J\) in \(K_T\).
A nonzero coset \(yJ\) has \(yJ=a\widetilde s J=XbJ\)
for a unique pair from the two unipotent quotients,
and conversely. It contributes \(\mu(J)\) times the
same phase as in the real unfolding.
Thus, with the real and finite Haar already fixed,
\begin{equation}\label{eq:fixed-tail-amplitude}
 V_p(l)_{\alpha_1,\alpha_2}\, I^\infty(w_p\sigma)
 =a_p^T\mu(J)\operatorname{CosAmp}^\flat_{c,s,c'}(w_p\sigma).
\end{equation}

Write \(I_v\) for the dimension index set of the
local cut Bessel function at \(v\), and \(c_G\) for the
real open-cell constant of Section~\ref{sec:real-terminal-gram}.
In the real coordinates the squared basis scale is
\[
 s_{v,a}=\kappa_v^{2(\varrho-a)}/N ,
\]
where \(\kappa_v\) is the real scale for the additive
character of level \(M\). In fact the finite unipotent
character on rational elements reads the root sum divided
by \(M\) against the different-normalized character.
Write \(\delta_E\) for the skew generator of the
different. At \(v\), put
\(\kappa_v=M\Phi_v(\delta_E)/(2\pi i)\).

For the finite character at \(-\mathfrak s(u)/(M\delta_E)\)
the opposite archimedean phase is
\(\exp(c-\bar c)\) with \(c=\Phi_v(\mathfrak s(u))/\kappa_v\),
where \(\mathfrak s\) is the simple-root sum. Scaling \(e_a\) by
\(\kappa_v^{\varrho-a}/\sqrt N\), its dual by the reciprocal
and \(z\) by \(1/\sqrt N\) gives this same consecutive
root coefficient: conjugation into the scaled basis
multiplies each simple-root entry by \(\kappa_v^{-1}\).
Conjugating the Weyl element gives proper head entries
\(x_{v,a}=s_{v,a}\Phi_v(p_a)\); the last Gram entry is
\(x_{v,\varrho}=\Phi_v(t)/N\).
With \(d_a=k-1+2(\varrho-a)\), set
\[
 R_0=\prod_v\left(\frac{\sqrt{c_G}}{|I_v|}
               \prod_{a=0}^{\varrho}s_{v,a}^{d_a}\right),
 \qquad G(t)=\prod_v(|\Phi_v(t)|^2-N^2)^{k-1}.
\]
Let \(a_p^\infty\) be the product of the real tensor
scales in Section~\ref{sec:real-terminal-gram}.

For \(a<\varrho\), the product formula gives
\[
 \operatorname{md}(p_a)^{1/2}\prod_v|\Phi_v(p_a)|=1.
\]
Consequently the proper-head factors cancel and
\begin{equation}\label{eq:fixed-tail-normalization}
 a_p^\infty a_p^T=R_0\sqrt{G(t)}.
\end{equation}
The final factor comes from
\((|x_{v,\varrho}|^2-1)^{(k-1)/2}
 =N^{-(k-1)}(|\Phi_v(t)|^2-N^2)^{(k-1)/2}\).
There is no finite \(t\)-factor in \(a_p^T\).
The Gram matrix \(\mathcal G_t(p,q)\) of the
\(a_p^\infty I^\infty(w_p\sigma)\) is the product of the
top forms of Section~\ref{sec:real-terminal-gram}.
For all blocks with the same \(t\) this is one matrix,
as the real top forms depend on their heads and \(t\).

\subsubsection{The common-entry pairing}
\label{sec:fixed-vector-common-pairing}

\sectionopening{In this subsection, we group the terms with the same final vector using the phase identity of Section~\ref{sec:fixed-vector-class-phases}.}

Define
\[
 \mathcal H_t(l)=\sum_{\alpha_1,\alpha_2\in Z/J}
   \left\|\sum_p V_p(l)_{\alpha_1,\alpha_2}
                   a_p^\infty I^\infty(w_p\sigma)\right\|^2 .
\]
With the Hermitian products linear in the second argument, its expansion is
\[
 \mathcal H_t(l)=\sum_{p,q}\mathcal G_t(p,q)
                   \langle V_p(l),V_q(l)\rangle .
\]
For fixed \(l\), only \(\mathcal G_t\) involves the real representative
\(\sigma\). The two compact factors in
Section~\ref{sec:real-terminal-gram} act by the same unitary operators
for every head, so the expansion depends on that representative only
through \(t\).
Fix \(c,c',\sigma\), and sum the comparison over the
blocks with these parameters. For each \(\alpha_1\) and
such residue \(s\), there is exactly one \(\alpha_2\)
with \(\alpha_1\alpha_2^{-1}\in\widetilde s J\), since
\(\widetilde s\in Z\). The pairs for distinct \(s\)
are disjoint. Equations \eqref{eq:fixed-tail-amplitude}
and \eqref{eq:fixed-tail-normalization} and nonnegativity
give
\begin{equation}\label{eq:fixed-tail-classes}
 \mu(J)^2[Z:J] R_0^2G(t)
      \sum_s\|B^\flat_{c,c',s,L}\|^2\le\mathcal H_t(l).
\end{equation}

Fix \(c,v_0,t\) and write \(B=B_{\pi_T}\) for the basis of the
glued middle lattice at \(T\). In rational coordinates the class
representative is \(c^\circ\) on \(X_0\) and fixes \(z\).
Since \(h_c=\iota_L(c^\circ)\) in the lattice basis
and \(\sigma z=v_0-(t/N)z\), the common last column is
\[
 l^{-1}z=B^{-1}\bigl(((c^\circ)^{-1}v_0)_T-(t/N)z\bigr)=:e .
\]
The two factors from the finite point correction \(\eta\) cancel in
the class residue condition of Section~\ref{sec:period-petersson-blocks}.
If \(l_1l_2^{-1}\in J\), the quotient
\(h_{c_1'}^{-1}(\sigma_1^{-1}\sigma_2)_L h_{c_2'}\)
is in the ambient principal level: it is integral
outside \(T\) by the preceding elimination and in \(J\)
at \(T\). The rational element \(\sigma_1^{-1}\sigma_2\)
fixes the pairs and \(z\), so belongs to \(H(F_0)\).
The pullback of the ambient principal level to \(H\) is
the level \(J_H\) of Section~\ref{sec:ramified-source}.

For this lattice calculation use the glued generator
\(x_{k_0}/\bar\pi_N+z/N\) of the middle lattice.
Write \(x_{\rm gl}=x_{k_0}/\bar\pi_N+z/N\) for
this lattice vector. Since \(x_{k_0}=\bar\pi_N x_{\rm gl}-z/\pi_N\),
the intersection with \(X_0\) is the integer lattice
with distinguished coordinate in
\(\pi_N\widehat{\mathcal O}_E\).
If an isometry \(h\) fixing \(z\) preserves the
middle lattice, its columns away from \(k_0\) are
in this intersection, as is
\((h x_{k_0}-x_{k_0})/\bar\pi_N\).

Applying the same fact to \(h^{-1}\) shows that \(h\)
is integral with integral inverse on \(X_0\).
It fixes the distinguished column modulo \(\bar\pi_N\),
equivalently the distinguished covector modulo \(\pi_N\)
by adjunction. Conversely an integral unitary matrix
with this condition on the covector preserves the intersection.
It has the column condition as well. The distinguished
diagonal difference is divisible by both coprime
factors, so \((h x_{k_0}-x_{k_0})/\bar\pi_N\)
is in the intersection.

This proves preservation
of the glued lattice, for \(h\) and its inverse. At
\(p_0^{\rm lev},p_1^{\rm lev}\mid M\), the integer \(N\) is a unit and the
principal congruence is the original one modulo \(\pi_M\). Hence the two source
classes are equal, and the freeness in
Section~\ref{sec:ramified-source} then gives
\(\sigma_1=\sigma_2\).

Thus these positions give disjoint cosets in the
\(H\)-fiber over \(e\). Indeed, fixing \(l_0\) over \(e\),
use \((H\cap J)l l_0^{-1}\subset H\) for each position.
On such a coset the integrand
\(\mathcal H_t(h l_0)\) equals \(\mathcal H_t(l)\),
since \(W_Z\) is invariant on the left under \(H\cap J\).
This invariance follows from its coordinate formula.
The group \(H\) fixes the hyperbolic pairs and \(z\).
Its middle-block action normalizes \(U\) and preserves
the simple-root character.

An element of \(H\cap J\) acts trivially
after moving past \(U\) and using right \(J\)-invariance.
Write \(\bar H\) for the finite reduction of the integral definite group and \(m_U\) for the order of the
reduced integral unipotent. In the chosen Haar,
\(\mu_H(H\cap J)=|\bar H|^{-1}\) and
\(\delta_U=m_U^{-1}\). Indeed the terminal measure gives
\(H\cap K_T\) mass one. At the pairs above \(N\) this is the
covector stabilizer computed above; at the pairs above \(M\) its
reduction has kernel \(H\cap J\) and image \(\bar H\).
The integral unipotent has mass one in its own measure.
The reductions are unchanged by projection
to \(T\), as representatives extend by the identity
outside \(T\). The integral Bessel subgroup is a product
of its integral unipotent and definite factors, by
reading the middle block. Its reduced image has order
\(|H_R|=m_U|\bar H|\). Thus
\[
 C_*=|H_R|m_U,\qquad
 \mu_H(H\cap J)\,\mu(J)^2[Z:J]=\mathcal R_Z/C_* .
\]

Let \(\Psi_t(e)=\int_H\mathcal H_t(h l_0)\,dh\).
This is nonnegative, since \(\mathcal G_t\) is a Gram
matrix. It is the sum of the top finite pairings
\(\Phi_{p,q}(e)\) of Section~\ref{sec:finite-common-entry}
weighted by \(\mathcal G_t(p,q)\):
the integrand is \(\sum_{p,q}\mathcal G_t(p,q)
\langle V_p(h l_0),V_q(h l_0)\rangle\) with both
inner products linear in the second argument.
It has a compactly supported, locally constant,
nonnegative extension from the last-vector quadric to
the ambient finite vector space by the compact step
extension from that section, taking the positive part
of the real sum. The extension can be chosen
\(H\cap J\)-invariant by averaging.

On the quadric
the fiber value is invariant since \(V_p(l h)=V_p(l)\)
for \(h\in H\cap J\), by changing its integration
variable using that same left invariance of \(W_Z\).
Summing \eqref{eq:fixed-tail-classes} over the disjoint
fiber cosets therefore proves
\begin{equation}\label{eq:fixed-tail-to-siegel}
 R_0^2G(t)\,\frac{\mathcal R_Z}{C_*}\!
     \sum_{(c,c',s,L)\,:\,(c,v_0,t)\ {\rm fixed}}
              \|B^\flat_{c,c',s,L}\|^2
       \le \Psi_t(e).
\end{equation}
This compares the original cut blocks with the same
finite form at a fixed last vector. For the subsequent
sum over \((c,v_0)\), the rational norm to average is
\(h(v_0,v_0)=(t\bar t-N^2)/N\), and the vector
\((c^\circ)^{-1}v_0\) is integral off \(T\) by the integral
stage factor and the integral original column \(z\).

\subsection{Averaging at a fixed final vector}
\label{sec:fixed-tail-average}

\sectionopening{In this subsection, we average the finite class sum of
Section~\ref{sec:fixed-tail-comparison} by interpreting the fixed
last-vector value as a Hermitian representation problem. This supplies
the terminal factor in Section~\ref{sec:adelic-common-entry-sum}.}

Let \(E^+=F_0\), \(d_C=\#\mathcal C\) for the source classes, and \(v_E=2^f/\sqrt{|D_E|}\), with \(D_E\)
the field discriminant. Fix \(t\) as in
Section~\ref{sec:fixed-tail-comparison} and put
\(\beta=(t\bar t-N^2)/N\in F_0\), totally positive.
Write \(B_{\pi_T}\) for the finite lattice basis,
with \(\pi_T=\pi_N|_T\). On the finite adeles of \(X_0\) take
\[
 F_t(y)=\mathbf1_{\{y\ \mathrm{integral\ off}\ T\}}
    \Psi_t\bigl(B_{\pi_T}^{-1}(y_T-(t/N)z)\bigr).
\]
Here vectors in the finite model are read in the
basis \(B_{\pi_T}\); this is the same affine coordinate \(e(t,y_T)\)
as in Section~\ref{sec:finite-common-entry}. The image of this affine embedding is closed and its inverse
is the \(X_0\)-projection after multiplication by \(B_{\pi_T}\).
The pullback of the compact step test is therefore locally constant
and compactly supported. Off \(T\) the indicated lattice is compact open.
The embedding intertwines the rational-coordinate action
on \(X_0\) with the lattice-coordinate action through \(B_{\pi_T}\),
since the source group fixes \(z\). The extended function
\(\Psi_t\) is invariant under \(H\cap J\) on the ambient space. Thus \(F_t\) is invariant
under the original source level, which also preserves
the lattice off \(T\). On the selected \((c,v_0,t)\),
\(F_t((c^\circ)^{-1}v_0)=\Psi_t(e)\).

We use Weil's Siegel formula for the orbit of one vector
\cite[\S52, Theorem~5, (41)]{WeilSiegel}. In this case
the group is \(U(X_0)\) over \(F_0\), for the quadratic
algebra \(E/F_0\), and the rank condition is \(k>2\).
The group has compact adelic quotient by definiteness.
An original split finite place meets the local orbit
condition: there the norm is \(x\cdot y\) on pairs of
vectors over \(F_{0,v}\), with the general linear
group on a vector and its dual.

The full-rank locus has
both vectors nonzero, and the group is transitive at
each value of the pairing by extending a vector and
covector to an adapted basis. In particular the
full-rank zero-norm locus is nonempty there.
On the orbit of the nonzero represented norm \(\beta\), unfolding
the Tamagawa integral and dividing by the volume of
\(U(X_0)(F_0)\backslash U(X_0)(\mathbb A_{F_0})\) gives the factor
\[
 \frac{\tau(U(v_0^\perp))}{\tau(U(X_0))}=\frac{2}{2}=1 .
\]
Indeed both groups are unitary of positive rank over the quadratic
field extension \(E/F_0\) \cite[\S2.5, p.~755]{GanQiuTakeda}.
This recovers the stabilizer normalization in
\cite[\S53, (43)]{WeilSiegel} and applies as well to the consecutive
ranks below. Probability Haar on the quotient assigns the same
mass \(1/d_C\) to each class, since each open
class has Haar weight proportional to the reciprocal
of its rational stabilizer order and these orders are
one. At each real place take the Schwartz function
\(y\mapsto u_v(h_v(y,y))\), where \(u_v\in C_c^\infty(\mathbb R)\)
and \(u_v(\Phi_v(\beta))=1\). Definiteness gives compact support.
The function is invariant under \(H(F_{0,v})\) and equals one on the
prescribed sphere. Weil's absolute convergence hypothesis concerns
every adelic Schwartz function. Such a function is a finite sum of
products of a real Schwartz function with a compactly supported
locally constant finite function. For each product use a compact set
of adelic quotient representatives. Translating its finite support by
this set confines all rational vectors to one fractional lattice.
The archimedean translates of the Schwartz factor have a uniform bound \(C_a(1+\|y\|)^{-a}\) for every \(a>0\).
Choose \(a>2kf\) to get absolute and uniform convergence.

The measure on the norm fiber in (41) is the quotient
of the additive vector measure by the additive norm-value
measure, both with adelic covolume one
\cite[\S44, Theorem~2]{WeilSiegel}. We spell out its
conversion. Use ordinary area on \(\mathbb C\) at the
complex places, so the product finite vector measure
gives \(\widehat{\mathcal O}_E\) mass \(v_E\) per
coordinate, by the covolume of \(\mathcal O_E\).
At the real places of \(F_0\) use ordinary Lebesgue
measure in the norm value. In the finite norm-value
Fourier variable \(u\in\mathbb A_{F_0,f}\), take
\(du_0\) giving the finite inverse different of \(F_0\)
mass one.

The phase
\(\psi_{E,f}(u(h(y,y)-\beta)/2)\) is the standard
\(F_0\)-phase, by the trace \(E/F_0\).
For this Fourier pairing the covolume-one measure is
\(\sqrt{|D_{F_0}|}\,du_0=|D_E|^{1/4}\,du_0\):
the archimedean covolume of the inverse different is
\(|D_{F_0}|^{-1/2}\), and \(E/F_0\) is unramified at
the finite places by the inertia calculation in
Section~\ref{sec:different-additive-character-scales}.
The real norm fiber in \(\mathbb C^k\) has mass
\(\pi^k\Phi_v(\beta)^{k-1}/(k-1)!\), by differentiating
the Euclidean ball volume in its squared radius.
Accordingly Weil's formula gives
\begin{equation}\label{eq:siegel-class-density}
 \sum_c\ \sum_{h(v_0,v_0)=\beta} F_t((c^\circ)^{-1}v_0)
 =d_C |D_E|^{1/4}
      \prod_v\frac{\pi^k\Phi_v(\beta)^{k-1}}{(k-1)!}
                   \,\mathcal D_f(F_t,\beta),
\end{equation}
where \(\mathcal D_f\) is the finite norm density with
the vector measure just specified and frequency measure
\(du_0\).

To bound the finite density, take out \(v_E^k\) by
normalizing the integer vector box at every place to
mass one. Put \(b_u=-\delta_E u\), a skew frequency,
where \(\delta_E\) generates the different of \(E\).
Because \(E/F_0\) is unramified or split at every finite
place, the different of \(E\) is the extension of the
different of \(F_0\). Thus \(du_0\) becomes the measure
on skew frequencies giving the integral skew box
mass one.

With \(\psi_m(a)=\psi_{E,f}(-a/\delta_E)\)
and an integral trace-one scalar \(\tau\), the phase
on the norm difference is
\[
 \psi_{E,f}\bigl(u(h(y,y)-\beta)/2\bigr)
     =\psi_m\bigl(\tau b_u(h(y,y)-\beta)\bigr).
\]
Indeed \(\operatorname{Tr}_{E/F_0}(\tau)=1\).
The character \(\psi_m\) is self-dual on the integer
ring and is trivial on \(a+\bar a\), since \(\delta_E\)
is skew. In particular
\(\psi_m(\tau b_u(a+\bar a))=\psi_m(b_u a)\).
It follows that the annihilator of the fixed integral
ring for the skew-frequency pairing is the integral
skew ring: testing against fixed traces is equivalent
to testing \(b_u\) against all integral \(a\).
The two measures giving these boxes mass one are
therefore dual. These identities hold also at places
above two.

For a finite place \(v\) of \(F_0\) off \(T\), let \(R_v\)
be the integer ring of the local quadratic algebra and
let \(q_v\) be the residue size of \(F_{0,v}\).
Let \(G_v(b)\) be the integral of
\(\psi_m(\tau b(h(x,x)-\beta))\) on \(R_v^k\),
with its probability measure. Squaring its absolute
value and writing the two variables as \(x,x+d\)
leaves the inner average of \(\psi_m(bh(x,d))\)
times a unit phase depending on \(d\).

Hence that
average vanishes unless \(bd\in R_v^k\), and
\[
 |G_v(b)|^2\le [R_v:R_v\cap b^{-1}R_v]^{-k}.
\]
The skew lattice has a unit generator. At skew
valuation \(-j<0\) the index is \(q_v^{2j}\), in both
the unramified and split cases. The skew shell has
volume \((1-q_v^{-1})q_v^j\), so its contribution
after integrating \(G_v\) in absolute value is at most
\((1-q_v^{-1})q_v^{-j(k-1)}\).
The integral over the skew integer box is at most one.
Thus each exterior factor contributes at most
\((1-q_v^{-(k-1)})^{-1}\).
Their product is bounded by
\[
 \zeta_{F_0}(k-1)\le\zeta(3\varrho+4)^f
                     \le 1+2^{-\varrho}\le1+b^{-6}.
\]
Here \(k-1\ge3\varrho+4\), \(f\le 2^{2\varrho}\) and
\(\log\zeta(a)\le2^{1-a}\) for \(a\ge2\);
the last inequality uses the parameter bounds on \(b,n\).

At \(T\) the Fourier integral is exactly the
norm-fiber integral in
Section~\ref{sec:finite-common-entry}, since its
skew parameter is dual to the norm-value Haar
with integer mass one. For absolute convergence of
the Fourier inversion, decompose the compactly
supported test into finitely many shifted boxes
\(x_0+cR_v^k\) on each orbit of places, using products
of such boxes across \(T\). After the substitution
\(x=x_0+cx'\), the same squaring argument bounds
the integral by the conductor index for frequency
\(c\bar c b\), with the additional box volume
factor.

In the squared modulus use the difference variable
\(d=x'-x''\) for the two scaled box coordinates.
The shift multiplies each outer \(d\)-term by the
unit phase of the cross term in \(x_0,cd\). The inner
\(x''\)-average has the linear frequency
\(c\bar c b d\), so it vanishes off the same conductor
annihilator as in the centered box. Taking absolute
values after this cancellation gives the stated
bound. The skew integral is finite by \(k-1>1\).
Write \(E_v=E\otimes_{F_0}F_{0,v}\) and \(q(y)=h(y,y)\).
At every finite \(v\) the image \(\beta_v\) of the nonzero rational
norm is nonzero. Choose \(\tau\in\mathcal O_{E_v}\) of trace one.
At a split place take \((1,0)\); at an unramified nonsplit place
lift a residue element of trace one and divide by its unit trace.
For \(q(y)=\beta_v\) the differential and an \(F_{0,v}\)-linear
right inverse are
\[
 dq_y(u)=\operatorname{Tr}_{E_v/F_{0,v}}h(y,u),\qquad
 a\longmapsto a\beta_v^{-1}\tau y .
\]
Thus the norm map is a submersion also at dyadic places. The
pushforward of the compactly supported test has a continuous
density near \(\beta_v\) in these norm charts. Fourier inversion
with the dual Haar measures gives its value at \(\beta_v\).

Outside \(T\), \(\beta\) is integral because \(N,t\)
are units there. The Euler-product calculation
can therefore first be made with integral skew
frequencies at all but finitely many places, each
such integral box contributing one in the normalized
measure, and then extended by the preceding bounds
on the absolute integrals.
In total,
\begin{equation}\label{eq:siegel-finite-density-bound}
 |\mathcal D_f(F_t,\beta)|
 \le (1+b^{-6})v_E^k
     \int_{h(y,y)=\beta}
         \Psi_t(e(t,y))\,d\nu_{\beta,T}(y).
\end{equation}

Let \(Q_t\) be the sum against \(\mathcal G_t(p,q)\)
of the unnormalized finite last-entry densities
of the two heads. By the definition of \(\Psi_t\) and
the terminal disintegration, it is nonnegative and
\[
 \int_{h(y,y)=\beta}\Psi_t(e(t,y))\,d\nu_{\beta,T}(y)
      =\operatorname{md}(\pi_T)Q_t=N^{-f}Q_t.
\]
Indeed \(\operatorname{md}(\pi_T)^2=\operatorname{md}(N)
=N^{-2f}\), since \(\pi_T\bar\pi_T=N\), \(T\) is
conjugation-stable, and \(N\) is a unit off \(T\).
Together with the \(N^{-f(k-1)}\) from
\(\prod_v\Phi_v(\beta)^{k-1}\), this gives the constant
\[
 C_{\rm S}=|D_E|^{1/4}v_E^k
             \left(\frac{\pi^k}
                               {(k-1)!\,N^k}\right)^f .
\]
By positivity the sum of \(\Psi_t(e)\) over the
selected pairs \((c,v_0)\) is bounded by
\eqref{eq:siegel-class-density}. Consequently
\eqref{eq:fixed-tail-to-siegel} and
\eqref{eq:siegel-finite-density-bound}, with \(G(t)>0\),
give the fixed-entry estimate
\begin{equation}\label{eq:open-fixed-entry-bound}
 R_0^2\sum_{(c,c',s,L):\,t_L=t}
       \|B^\flat_{c,c',s,L}\|^2
     \le (1+b^{-6}) C_*d_C C_{\rm S}\,
                     \frac{Q_t}{\mathcal R_Z}.
\end{equation}
This applies to the original cut blocks on any
finite set at \(t\), so the bound holds for their
nonnegative sum.

\subsection{Adelic summation of the common-entry forms}
\label{sec:adelic-common-entry-sum}

\sectionopening{In this subsection, we multiply the finite and real Fourier estimates with the same
quotient measures. The local bounds from
Sections~\ref{sec:finite-shell-reduction} and~\ref{sec:real-shell-sum}
then lead to the open-cell estimate of
Section~\ref{sec:open-equal-link-bound}.}

Since \(n=\varrho+1=m/2\), Section~\ref{sec:datum} gives
\(N\le b^{25m}=b^{50n}\le b^{60n}\).
At depth \(i\), take the product of the real
forms \(J_{v,i}\) and of the finite form
\(\mathcal I_i^T/\mathcal R_Z\) at weight one.
At depth \(n=\varrho+1\) use the finite last-entry
density in place of \(\mathcal I_n^T\).
Write the resulting form as \(J_i^{\rm ad}\).
For a rational entry \(t\in\mathcal E\), use its
scale \(a_{v,i}\Phi_v(t)\) at \(v\), with
\(a_{v,i}=s_{v,i}\) of Section~\ref{sec:fixed-tail-comparison},
and \(t_T\) at the finite places. These coordinates also
define the scaled heads of a rational tuple.

Write \(W_i(p)\) for the product of the real head weights
in Section~\ref{sec:real-shell-sum} and the finite head
weight in Section~\ref{sec:finite-shell-reduction}.
Define \(w_i=\sum_{p\in\mathcal E^i}W_i(p)\) and
\[
 \begin{split}
 D_i&=\sum_{\substack{p,q\in\mathcal E^i\\
                              p_{i-1}=q_{i-1}}}
                          \Re J_i^{\rm ad}(p,q)\quad (i\ge1),\\
 D_0&=\Re J_0^{\rm ad}
 \le M_0=(1+b^{-6})b^{-6k|\mathcal A|}
                  \prod_v\frac{1+b^{-10}}{\deg_v}.
 \end{split}
\]
The finite bottom bound is for the original \(Z\), by
Section~\ref{sec:finite-selected-bottom}. At each real place let
\(D_v=|I_v|\). The stationary vector in
Section~\ref{sec:real-marked-filling} has squared norm \(D_v\). Its
full correlation is \(D_v c_v\) for the coefficient of a unit vector.
The stage Haar factor is
\(\lambda=\operatorname{vol}U(r)/\operatorname{vol}U(k)\).
Hence the normalized full bottom mass is exactly
\[
 C_F\lambda^3\int|D_v c_v(g)|^2\,dg
 =\frac{\lambda\operatorname{vol}U(k)}{\operatorname{vol}U(r)}
      \frac1{\deg_v}=\frac1{\deg_v}.
\]
Section~\ref{sec:real-bottom-mass} bounds the cut mass by
\((1+b^{-10})/\deg_v\). Thus the product uses the same
\(\mathcal B=\prod_v\deg_v\) as the main moment. At each selected split pair,
the finite bound uses the index \(q^r-1\) for either a
vector or a covector stabilizer. Since the bottom masses
factor over conjugate pairs, the same bound applies to
every mixed choice of the two kinds of stabilizer.
For \(0\le i<n\), let \(F_i(x)\) be the real part
of the sum of \(J_{i+1}^{\rm ad}(p\mathbin{\|}x,q\mathbin{\|}x)\)
over the two proper heads in \(\mathcal E^i\).
Here \(x\) is a common entry on the product of the
real places and \(T\), and \(\mathbin{\|}\) appends that
entry in local coordinates.

We may replace a zero
complex component or a zero component in a local field \(E_w\), \(w\in T\),
by one in forming a head. Such replacements concern
null sets for the integrals and do not affect the
nonzero rational samples.

Put \(Y=\prod_{v\mid\infty}\mathbb C\times A_T\) with the product
of the stage additive measures, and let \(Y^\times\) consist of
the tuples nonzero in every field component. Its complement is null
because there are finitely many components and each singleton has
additive measure zero. For \(x\in Y^\times\) and a common suffix
\(u\in(Y^\times)^{n-i-1}\), write
\[
 K_{pq}(x,u)=J_n^{\rm ad}(p\mathbin{\|}x\mathbin{\|}u,
                          q\mathbin{\|}x\mathbin{\|}u),
 \qquad p,q\in\mathcal E^i.
\]
The empty suffix has mass one. At the common last entry each local
matrix is a Gram matrix by Sections~\ref{sec:real-terminal-gram}
and~\ref{sec:finite-common-entry}. Their entrywise product is a Gram
matrix, since one may tensor the representing vectors. In particular
\(K_{pp}\ge0\) and \(|K_{pq}|^2\le K_{pp}K_{qq}\).

At each real place and on \(A_T\), the unweighted stage integrals for
fixed nonzero heads are finite by the absolute shell estimates of
Sections~\ref{sec:real-shell-sum} and~\ref{sec:finite-shell-reduction}.
The terminal integrals are finite by the compact support bounds in
Sections~\ref{sec:real-terminal-gram} and~\ref{sec:finite-common-entry}.
Starting with the nonnegative local terminal diagonal, integrate the
last entry first. The one-step identities of
Sections~\ref{sec:real-common-descent} and~\ref{sec:finite-common-entry}
give the preceding stage diagonal. Iterating and using Tonelli for
the finite product of local factors gives
\[
 \int K_{pp}(x,u)\,du
  =J_{i+1}^{\rm ad}(p\mathbin{\|}x,p\mathbin{\|}x),\qquad
 \iint K_{pp}(x,u)\,du\,dx=J_i^{\rm ad}(p,p)<\infty .
\]
The right side of the first identity is finite at each
\(x\in Y^\times\). Cauchy--Schwarz on the suffix, and then on all
common coordinates, gives absolute integrability of \(K_{pq}\)
in both identities. Fubini and the same one-step identities now give
\[
 F_i(x)=\int\sum_{p,q}K_{pq}(x,u)\,du\ge0,\qquad
 \int F_i=\Re\sum_{p,q\in\mathcal E^i}J_i^{\rm ad}(p,q).
\]
The finite sum in the first integral is a squared norm. The local
terminal entries are continuous on nonzero heads, so this integrand
is measurable and continuous in \(x\in Y^\times\) for almost every
\(u\). Every selected rational entry has nonzero components at all
the indicated places, also after its real scalings. Fatou's lemma
along any sequence tending to that entry gives lower semicontinuity
there. For a last entry \(t\) with \(|\Phi_v(t)|>N\) at
every real place, the nonnegative summand of \(D_n\) at
\(t\) is \(Q_t/\mathcal R_Z\) of
Section~\ref{sec:fixed-tail-average}.

Extend \(F_i\) to the adeles by setting
\[
 G_i(z)=\mathbf1_{\{z_f\text{ integral off }T\}}\,
               F_i((a_{v,i}z_v)_v,z_T),\qquad
 P_i=\frac{v_E}{\prod_v a_{v,i}^2}.
\]
The exterior integer box is compact and open in the restricted
product off \(T\), hence its indicator is continuous. The extension
is measurable, is integrable by the displayed stage integral, and is
lower semicontinuous at the selected rational entries.
The additive Haar has \(E\) of covolume one, as in
Section~\ref{sec:fixed-tail-average}.
The Fourier coefficient at rational \(y/\delta_E\)
vanishes unless \(y\) is integral off \(T\), since the
integral exterior translations preserve \(G_i\).
For such \(y\) a change of variables in the Fourier
integral gives \(P_i\) times the transform of \(F_i\)
at
\[
 \xi_v=\frac{2\pi i}{\Phi_v(\delta_E)a_{v,i}}\Phi_v(y),
                         \qquad \xi_T=y_T .
\]
Let \(S_i\) be the complex sum defining \(F_i=\Re S_i\). Since
\(F_i=(S_i+\overline{S_i})/2\), its Fourier transform is
\[
 \widehat F_i(\xi)=\tfrac12\bigl(\widehat S_i(\xi)
                         +\overline{\widehat S_i(-\xi)}\bigr).
\]
The character estimates for \(y\) and \(-y\) have the same bound at
each place.
The factor \(P_i\) is the real area Jacobian times
the finite integer-box mass, since the projected
finite integer box has mass one in the local forms.
Indeed \(\psi(\xi_v a_{v,i}z_v)
=\exp(4\pi i\,\Re(\Phi_v(y/\delta_E)z_v))\) at each
complex place, while
\(\psi_{m,T}(y_Tz_T)=\psi_{E,T}(-y_Tz_T/\delta_E)\)
at \(T\). These are the local factors of the inverse adelic
character in the Fourier coefficient at \(y/\delta_E\), with
the character of Section~\ref{sec:fixed-tail-comparison}.
At zero frequency this is \(P_i\int F_i\).

\subsubsection{The sum over nonzero frequencies}
\label{sec:adelic-nonzero-frequencies}

\sectionopening{In this subsection, we combine the finite bounds of Section~\ref{sec:finite-shell-reduction} with Lemma~\ref{lem:marked-real-frequency-decay}.}

Set \(K_*=k-n^7\),
\(X_v=|\kappa_v|^{2\varrho+1}/(MN)\),
\(B_{\rm fr}=b^{23(2\varrho+1)}\), and
\[
 \mathcal D_{\rm dec}=\prod_v\Theta\,X_v^{K_*},
\]
where \(\Theta\) is the constant in \eqref{eq:real-shell-total}.
For a nonzero \(y\) integral off \(T\), its finite
denominator index \(D(y_T)\) is the absolute norm of
its denominator ideal \(\mathfrak b_y\).
In both applications of Lemma~\ref{lem:marked-real-frequency-decay}
the heads have length \(i\). For the character weight in the Poisson
coefficient take \(h=i\) and \(c_v=2\pi i/\Phi_v(\delta_E)\),
including \(i=0\) with empty heads. For unequal last entries in
the constant-weight case, \(1\le i\le\varrho\), take \(h=i-1\),
\(c_v=M\) and \(y=(t^{-1}-{t'}^{-1})/M\). Here \(h\) indexes
the scale of the entry giving the frequency. In both cases
\(\xi_v=c_v\Phi_v(y)/a_{v,h}\).
Since \(1\le|\kappa_v|\) and \(a_{v,h}=|\kappa_v|^{2(\varrho-h)}/N\),
one has \(a_{v,h}/|c_v|\le X_v\le B_{\rm fr}\).
The complex places contribute squared moduli to the product
formula, so
\[
 \operatorname{md}(y_T)^{-K_*/2}\prod_v|\xi_v|^{-K_*}
 =\prod_{\substack{w\nmid\infty\\w\notin T}}|y|_w^{K_*/2}
        \prod_v(a_{v,h}/|c_v|)^{K_*}
 \le\prod_v X_v^{K_*}.
\]
For the residual real damping use
\[
 \min(1,|\xi_v|^{-1})^6
            \le B_{\rm fr}^6\max(1,|\Phi_v(y)|)^{-6}.
\]
Using \(D(y_T)^{-6}\le D(y_T)^{-3}\), multiply
\eqref{eq:finite-frequency-final} and the real bounds
\eqref{eq:real-frequency-final}. The resulting product before its real constant
\(\Theta^f b^{2000n^9 f}\) and its finite constant
\(b^{n^9 f}\) is bounded by
\[
 \left(\prod_v X_v^{K_*}\right) B_{\rm fr}^{6f}
       \frac{\prod_v\max(1,|\Phi_v(y)|)^{-6}}
            {(\mathrm N\mathfrak b_y)^3}.
\]

For the sum of the last weights, fix a denominator
ideal \(\mathfrak b\). The elements with this denominator
lie in the \(\mathrm N\mathfrak b\) cosets of
\(\mathfrak b^{-1}/\mathcal O_E\).
In each coset, a product of complex half-unit squares
contains at most one element, since a nonzero integral
difference cannot have modulus less than one at all
complex places. If \(z\) lies in a square with corner
\(s\), \(\max(1,|s|)\le2\max(1,|z|)\).

There are at most \((4R+1)^2\le25R^2\) such corners
of modulus at most \(R\ge1\), by packing disjoint
discs of radius \(1/4\).
The sum of \(\max(1,|s|)^{-6}\) over corners is then
at most \(6\int_1^\infty25R^{-5}\,dR<38\).
Summing first over the cosets and then the ideals gives
\[
 \sum_{0\ne y\in E}
      \frac{\prod_v\max(1,|\Phi_v(y)|)^{-6}}
           {(\mathrm N\mathfrak b_y)^3}
 \le(64\cdot38)^f\zeta_E(2)
 \le(64\cdot38e^4)^f .
\]
Since \(64\cdot38e^4\le b^3\) and
\(3+138(2\varrho+1)\le n^9\), the common frequency
sum, including its real and finite constants, is
bounded by
\[
 \varepsilon_0=b^{2002n^9 f}\mathcal D_{\rm dec}.
\]
It follows that the nonzero Fourier coefficients of the
extension have absolute sum at most
\(P_i\varepsilon_0 w_i^2\).

For completeness, one-sided Poisson at the chosen
entries follows directly from periodization. The
nonnegative \(G_i\) is integrable, so its periodization
over \(E\) is integrable on the compact additive
quotient. Its Fourier coefficients are those just
computed. Absolute summability makes their Fourier
series continuous, equal almost everywhere to the
periodization by uniqueness of Fourier coefficients.
Near the identity in a quotient chart, the periodization
dominates the sum over the finitely many selected
rational translates of \(G_i\).

Lower semicontinuity at
those translates gives the same inequality at the
identity against the continuous representative.
Consequently
\[
 D_{i+1}\le P_i\left(\Re\sum_{p,q\in\mathcal E^i}
                        J_i^{\rm ad}(p,q)
                       +\varepsilon_0 w_i^2\right).
\]
In the unrestricted sum, equal last entries give \(D_i\).
For two different last entries \(t,t'\), take
\(y=(t^{-1}-{t'}^{-1})/M\). It is nonzero and integral
off \(T\), and supplies precisely the unequal-head
frequency: \(y_T\) and \(M\Phi_v(y)/a_{v,i-1}\)
at the indicated places. The singleton version of
the frequency sum therefore bounds that pair by
\(\varepsilon_0 W_i(p)W_i(q)\).
There is no unequal contribution at \(i=0\).
Thus for \(0\le i<n\),
\[
             D_{i+1}\le P_i(D_i+2\varepsilon_0 w_i^2).
\]

A head of nonzero weight has partial products \(\delta_j\)
integral at \(T\). Each factor belongs to \(\mathcal E\) and is
a unit off \(T\) by its construction in
Section~\ref{sec:open-rational-cell}. Thus \(\delta_j\in\mathcal O_E\)
and has finite modulus one off \(T\).
At each \(v\), the initial segments of the scaled head have
modulus at most \(b^{60n^2}\). The inverse of the
product of their scales is at most \(N^i\le b^{60n^2}\).
Hence \(|\Phi_v(\delta_j)|\le b^{120n^2}\).
For an initial segment with scale product \(c_v\), the finite
weight uses
\(\operatorname{md}(\delta_{j,T})^{-1/2}
=\prod_v|\Phi_v(\delta_j)|\).

The real factors for this initial segment, before raising
them to the power \(n^6\), give the product
\[
              \prod_v\max\{|\Phi_v(\delta_j)|,c_v^{-1}\}.
\]
It follows that \(W_i(p)\le b^{120in^8 f}\).
A head is recovered from its partial products.
Packing these integer tuples using the maximum complex
modulus gives at most
\((1+2b^{120n^2})^{2if}\) heads of nonzero weight:
the difference of distinct tuples has modulus at least
one in some coordinate, so their radius-\(1/2\) polydiscs
have disjoint interiors.
Therefore \(w_i\le b^{130n^9 f}\).
Also \(a_{v,i}\le b^{46\varrho}\) and
\(v_E\ge b^{-f}\) by the root discriminant bound.
Hence \(\prod_{j<i}P_j\ge b^{-100n^2 f}\)
and \(w_i^2\le b^{400n^9 f}\prod_{j<i}P_j\).

Put \(A_i=\prod_{0\le j<i}P_j\), so \(A_0=1\) and
\(A_{i+1}=A_iP_i>0\) by the formula for \(P_i\). Dividing the
recurrence by \(A_{i+1}\) and summing gives
\[
 \frac{D_n}{A_n}\le D_0+
 2\varepsilon_0\sum_{i=0}^{n-1}\frac{w_i^2}{A_i}
 \le M_0+2n\,b^{(2002+400)n^9f}\mathcal D_{\rm dec}.
\]
Put
\[
 \mathcal E_+=\{t\in\mathcal E:|\Phi_v(t)|>N
                              \text{ for every }v\mid\infty\}.
\]
Using \(2n\le b^{n^9 f}\) gives
\begin{equation}\label{eq:adelic-top-sum-bound}
 \sum_{t\in\mathcal E_0}Q_t/\mathcal R_Z
 \le D_n\le
 \left(\prod_{i<n}P_i\right)
              (M_0+b^{2403n^9 f}\mathcal D_{\rm dec})
              \qquad(\mathcal E_0\subseteq\mathcal E_+).
\end{equation}
The first inequality follows from nonnegativity at each common
last entry. For \(\Omega\) contained in the open
residues of the prescribed stabilizer, fix any finite collection \(\mathcal F\) of
tuples \((c,c',s,L)\) with \(s\in\Omega\), retaining the full
nonzero-discriminant amplitude at each tuple. Taking \(\mathcal E_0\)
to contain their last entries, \eqref{eq:open-fixed-entry-bound} gives
\[
 R_0^2\sum_{(c,c',s,L)\in\mathcal F}
                  \|B^\flat_{c,c',s,L}\|_{\rm HS}^2
 \le (1+b^{-6})C_*d_C C_{\rm S}\!
         \left(\prod_{i<n}P_i\right)
                 (M_0+b^{2403n^9 f}\mathcal D_{\rm dec}).
\]
The supremum of the left side over these finite collections is
\(R_0^2\mathcal L(\Omega)\) by the nonnegative sum defining
\(\mathcal L\) in Section~\ref{sec:full-orbital-blocks}.

\subsection{Consecutive masses of definite unitary lattices}
\label{sec:consecutive-unitary-masses}

\sectionopening{In this subsection, we compare the genus masses in two consecutive ranks, retaining the
integral automorphism weights. This supplies the factors used in
Sections~\ref{sec:fixed-tail-average}
and~\ref{sec:main-identity-moments}.}

Write \(\mu_j\) for the mass of the genus of the standard
positive definite Hermitian lattice \(\mathcal O_E^j\).
Let \(H_j\) be its unitary group, \(K_j\) the stabilizer of
\(\widehat{\mathcal O}_E^j\), and
\(L_c=E^j\cap c\widehat{\mathcal O}_E^j\) for
\(c\in H_j(\mathbb A_{F_0,f})\). These lattices give all classes in
the genus. For a lattice described by its local completions, Landherr's
theorem first identifies its rational Hermitian space with \(E^j\)
in the standard form \cite[pp.~246--247]{LandherrHermitian}.
The lattice is recovered from its finite completions, and local isometries
of the standard lattice give its representative \(c\).
Its stabilizer is
\(\Gamma_c=H_j(F_0)\cap cK_jc^{-1}\).
Definiteness makes every rational unipotent in \(H_j\) trivial.
The norm-one determinant torus has no nontrivial \(F_0\)-rational
characters, and every rational character of \(H_j\) factors through it.
The adelic quotient is therefore compact by
\cite[Theorem~5.6(ii), p.~21]{BorelAdele}.
The compact real factor and the compact open \(K_j\) make the class
set a finite discrete quotient of the adelic space.
The class of \(c\) has measure
\(\operatorname{vol}(H_j(F_{0,\infty}))\operatorname{vol}(K_j)/|\Gamma_c|\).
Consequently
\[
 \mu_j=\sum_{[c]}|\Gamma_c|^{-1}>0,\qquad
 \operatorname{Pr}([c])=\frac{1}{\mu_j|\Gamma_c|}.
\]
For each rank \(j+1\), the weighted count of pairs consisting of
a lattice in this genus and a norm-one vector
equals \(\mu_j\). Indeed such a vector \(w\) splits off
\(\mathcal O_Ew\) by the integral projection onto its line.

Over a split finite place the two components of a
norm-one vector are a vector and a covector of pairing
one, and an adapted basis sends them to the standard
pair. At an unramified nonsplit place the complement
is unimodular. In a unimodular Hermitian lattice
one can split off a unit-norm line by reduction:
if every residue norm vanishes, polarization using
all scalars and the nondegenerate residue trace
pairing forces the residue form to vanish.

This
works also in residue characteristic two.
Surjectivity of the norm on unramified local units
allows the unit on each diagonal line to be made one.
Thus the complement is standard locally. Rationally,
Witt extension carries the norm-one line to the standard
line in the standard Hermitian space. Conversely a
rank-\(j\) genus lattice with an adjoined standard
line belongs to the rank-\(j+1\) genus. These operations
give inverse isometries of the respective objects and
identify their automorphism groups.

Fix \(k\le j<r\), so \(j+1>2\). Apply the one-vector formula of
Section~\ref{sec:fixed-tail-average} at the full integral level,
with the indicator of the finite integer box and norm value one.
Probability Haar weights each class
by its reciprocal automorphism order divided by the
total mass. The split-place and rank conditions hold as
before. Let \(x_j\) be the finite norm density divided by
\(v_E^{j+1}\), in the frequency Haar giving the codifferent
of \(F_0\) mass one. The formula gives
\[
 \frac{\mu_j}{\mu_{j+1}}
  = |D_E|^{1/4} v_E^{j+1}
             \left(\frac{\pi^{j+1}}{j!}\right)^f x_j .
\]
Put \(l=j+1\), and at a finite place \(v\) of \(F_0\)
let \(q_v\) be its residue size. Write \(\chi_v=1\) at a split
place and \(\chi_v=-1\) at an inert place. The number of solutions
of \(h(w,w)=1\) in the residue quadratic algebra of rank \(l\) is
\[
 A_v=q_v^{2l-1}-\chi_v^l q_v^{l-1}.
\]
At a split place a nonzero vector has \(q_v^{l-1}\) covectors
taking value one on it. At an inert place choose a nontrivial
additive character \(\psi\) of the fixed residue field. For
\(a\ne0\), the single-coordinate sum is
\[
 \sum_{z\in\mathbf F_{q_v^2}}\psi(a z\bar z)=-q_v,
\]
since the nonzero norms have \(q_v+1\) preimages. Orthogonality
in the norm variable gives
\(A_v=q_v^{-1}(q_v^{2l}-(-q_v)^l)\).
The derivative at a solution is \(u\mapsto\operatorname{Tr}(h(w,u))\).
Taking \(u=tw\) gives the surjective residue trace \(t\mapsto
\operatorname{Tr}(t)\), including at residue characteristic two.
Consequently every solution has \(q_v^{2l-1}\) lifts from one precision
to the next. With the vector and norm-value integer boxes each of
mass one, the norm fiber therefore has local density \(q_v^{1-2l}A_v\).

The absolute shell estimate of Section~\ref{sec:fixed-tail-average}
applies with \(l\) coordinates. Fourier inversion at every finite
set of places, followed by absolute convergence of the tails, gives
\[
 x_j=\prod_{v\nmid\infty}(1-\chi_v^l q_v^{-l}),\qquad
 |x_j-1|\le\prod_{v\nmid\infty}(1+q_v^{-l})-1
                       \le\zeta_{F_0}(l)-1.
\]
For \(j\ge k\), using \(f\le2^{2\varrho}\), \(k\ge4\varrho+12\),
and \(\log\zeta_{F_0}(l)\le f2^{1-l}\), one has
\[
 |x_j-1|\le2^{-\varrho},\qquad
 |\log x_j|\le 2^{1-\varrho}.
\]
Define
\[
 \sigma_{v,l}=
  \frac{(2\pi)^l}{|\Phi_v(\delta_E)|^{l-1/2}(l-1)!},
 \qquad
 \mathrm{fr}_v=N^k\prod_{j=k}^{r-1}\sigma_{v,j+1}.
\]
Since \(\prod_v|\Phi_v(\delta_E)|=\sqrt{|D_E|}\),
the mass identity reads
\(\mu_j/\mu_{j+1}=x_j\prod_v\sigma_{v,j+1}\).
The sizes \(N\ge b^m\), \(\log b=O(x^\epsilon)\),
\(m\ge\sqrt x\), \(f\le x^{12}\), and \(r-k=m-1\) give
\[
 (r-k)2^{1-\varrho}+4f/N\le\log(1+b^{-9})
\]
for all sufficiently large \(x\), since each term on
the left is \(o(b^{-9})\).
Taking the product over \(k\le j<r\) and using
\(\log(N^k/(N^k-1))\le2/N\) yields
\begin{equation}\label{eq:consecutive-frame-masses}
 \left(\prod_v\mathrm{fr}_v\right)\mu_r
     \left(\frac{N^k}{N^k-1}\right)^{2f}
       \le (1+b^{-9})N^{kf}\mu_k .
\end{equation}

\subsection{The class count and the unipotent cutoff measure}
\label{sec:marked-class-covolume}

\sectionopening{In this subsection, we compute the change in counting measure from the integral stabilizer
to the source level subgroup. Its finite covolume combines with the
orbital measure of Section~\ref{sec:period} in the moment comparison
of Section~\ref{sec:main-identity-moments}.}

Let \(K_X\) be the stabilizer of \(\widehat{\mathcal O}_E^k\)
in \(H(\mathbb A_f)\), and let \(H_{\rm int}\) be its subgroup
preserving the ambient lattice. The glued-column
calculation of Section~\ref{sec:fixed-tail-comparison} shows
that \(H_{\rm int}\) is the last-covector stabilizer modulo \(\pi_N\).
At each of the \(f\) chosen places above \(N\), reduction
from the split integral factor to \(\mathrm{GL}_k(\mathbb F_N)\)
is onto: lift the matrix over the valuation ring and put
its inverse transpose at the conjugate place. The group
of matrices over \(\mathbb F_N\) acts transitively on its
\(N^k-1\) nonzero covectors. The congruence at the two
primes dividing \(M\) is independent of these conditions on the covectors since \(N\nmid M\). Thus
\[
 [K_X:H_{\rm int}]=(N^k-1)^f,
 \qquad [H_{\rm int}:J_H]=|\bar H|.
\]
If a class at \(K_X\) has stabilizer \(\Gamma\), its fiber
at \(J_H\) is \(\Gamma\backslash K_X/J_H\).
The stabilizer acts freely on this fiber before quotienting,
by Section~\ref{sec:ramified-source}. Summing the fiber
sizes over the classes at \(K_X\) gives
\begin{equation}\label{eq:marked-source-class-count}
 d_C=|\bar H|(N^k-1)^f\mu_k .
\end{equation}

For the class cutoff put \(B=B_{\pi_N}\), the original glued
basis in the rational coordinates \(e_a;x_j,z;e_a^*\). By
Section~\ref{sec:period-finite-coordinates},
\(h_c=B^{-1}\iota_{\rm rat}(c^\circ)B\). Thus the class lattice
has basis
\[
 B_c=\iota_{\rm rat}(c^\circ)B=B h_c,
 \qquad
 B_c^{-1}uB_c=h_c^{-1}(B^{-1}uB)h_c.
\]
For a rational unipotent \(u\), membership in \(U_c\) is precisely
\(B_c^{-1}uB_c\in 1+M M_r(\widehat{\mathcal O}_E)\).
Here reduction modulo the rational integer \(M\) is
equivalent to reduction at the half of its split
places, by unitarity at the conjugate places. Use
\(a_{v,a}=s_{v,a}\) of Section~\ref{sec:fixed-tail-comparison}.
At a real place, write \(b_{v,j}\) for the scale of a basis
column. The scales on \(e_a,e_a^*\) are
\(\kappa_v^{\varrho-a}/\sqrt N\) and its inverse, those
on \(x_j\) are one, and that on \(z\) is \(N^{-1/2}\).
The real image of a rational matrix is conjugated by
this diagonal change of basis.

The free entries in the row indexed by \(e_a\) lie in the columns indexed by \(e_b,e_b^*\) for
\(b>a\) and in the middle block. Write
\(r_a=k+1+2(\varrho-1-a)\) for their number. The other
coordinate in this row is the imaginary part at \(e_a^*\).
To obtain the coordinate chart, take the product
of the real Heisenberg factors at \(0,\ldots,\varrho-1\)
in that order. On \(F_a^\perp/F_a\) for \(F_a=\langle e_b:b<a\rangle\), the
earlier factors act trivially. With the later
factors fixed, the free covector from the factor at \(a\)
is \(-h(w,\cdot)\) on the remaining space, composed with
the later isometry, and its imaginary coordinate is the
imaginary part of the Heisenberg center in the half-trace
chart.

The later factors fix \(e_a^*\). The Gram matrix
on the remaining space in the normalized basis is a
signed permutation matrix. Thus the covector map and
composition with the later isometry have real Jacobian
of absolute value one, and the center also has Jacobian
one. Descending in \(a\) gives the inverse chart and shows
that the original product of Heisenberg Haar measures is
ordinary product Lebesgue measure in these coordinates.

Let \(B_{c,a}\) be the block of \(B_c\) on the free
columns of row \(a\). The earlier basis and dual
columns of \(B_c\) are unchanged, so the free entries in row \(a\) of
\((u-1)/M\), regarded as a row vector, for \(u\in U_c\) lies in
\[
 L_{c,a}=\{\ell\in E^{r_a}:\ell B_{c,a}
                                      \in\widehat{\mathcal O}_E^{r_a}\}.
\]
Every such row vector is attained with the free and central
coordinates of every later row equal to zero. Extend \(\ell\) by zero and write it
as \(h(\omega,\cdot)\) in the original rational Gram
form, conjugate-linear in the first argument.

The class lattice is self-dual because the glued vector and \(z\)
span a plane with Gram determinant \(-1\), and the
other summands are unimodular. Thus \(\omega\) lies in it.
In particular \(h(e_a,\omega)=0\). Choose an integral
trace-one scalar \(\tau\), which exists since \(E\)
contains the sixth roots of unity. For an integral
skew scalar \(s\) put \(C=Ms-M^2\tau h(\omega,\omega)\).
The map
\[
 g(w)=w+M e_a h(\omega,w)-M\omega h(e_a,w)
                                 +C e_a h(e_a,w)
\]
is an isometry.

The change in the pairing on two vectors
is \((C+\bar C+M^2h(\omega,\omega))\) times the product
of their two \(h(e_a,\cdot)\) coordinates with conjugation
on the first, and hence zero. It is integral and
congruent to one modulo \(M\) in the class lattice.
Its free entries in row \(a\) form \(M\ell\). In the other top
rows only column \(e_a^*\) can change, and that entry is neither
free nor central in any row \(b>a\).

Setting \(\omega=0\) realizes every central shift
\(Ms\). Conversely if the free entries in row \(a\) vanish,
that center entry is skew: multiplication of the
matrix with its adjoint gives
\(CA^*+AC^*+BJ_0^{-1}B^*=0\) for the upper blocks \(A,B,C\) and the middle Gram \(J_0\).
Here row \(a\) of \(A\) equals row \(a\) of the identity matrix and
that of \(B\) is zero. The center entry divided by
\(M\) is integral since \(e_a^*\) is an unchanged
class-basis column.

These lattices in the successive coordinates give a
strict right fundamental domain. To see the translation
rule, if all rows of \(u\in U_c\) with index \(b>a\)
have their coordinates zero, its matrix on the tail
quotient after \(a\) is the identity, by the real chart
just described. Thus right multiplication by \(u\) adds
its free coordinates in row \(a\) and preserves the free and central
coordinates in all later rows.

If its free coordinates in row \(a\) also
vanish, the tail block including row \(a\) has
\(A=I\) and \(B=0\); its \(C\) is skew-Hermitian.
In this block the other entries in column \(e_a^*\)
vanish, so right multiplication adds the central
coordinate. Reduce in descending row order, free
coordinates before center, to half-open fundamental
parallelepipeds of the two additive lattices. The
isometries above supply every needed translation.
For uniqueness, compare two such reduced points
in a right coset, starting at the last row.
The change in each successive coordinate belongs
to its lattice, and the parallelepiped forces it to
be zero. The coordinate chart then identifies the
points.

The middle determinant of \(B_{c,a}\) is
\(\det(c^\circ)/\bar\pi_N\), since \(\iota_{\rm rat}(c^\circ)\)
acts by \(c^\circ\) on \(X_0\), fixes \(z\), and the glued column replaces
\(x_{k_0}\) by \(x_{k_0}/\bar\pi_N+z/N\).
Its finite modulus is \(N^f\): the unitary determinant
has modulus one and \(\pi_N\bar\pi_N=N\).
Hence the rational lattice \(L_{c,a}\) has Minkowski
covolume \(N^fv_E^{-r_a}\), using the finite box under
\(B_{c,a}^{-1}\) and the covolume-one additive adelic Haar.
The real changes of variables in that row have complex
scales \(M b_{v,j}/b_{v,e_a}\). The scales at the two
columns of every later hyperbolic pair cancel, while
the \(z\)-column has squared scale \(N^{-1}\).
Thus the measures of the free and central coordinates in row \(a\) are
\[
 \left(v_E^{-1}\prod_v\frac{M^2}{a_{v,a}}\right)^{r_a},
 \qquad |D_E|^{1/4}\prod_v\frac{M}{a_{v,a}},
\]
respectively. For the second formula, the skew integers
are \(i\mathcal O_{F_0}\), with imaginary components
equal up to signs to the embeddings of \(F_0\).
The unscaled center has covolume
\(\sqrt{|D_{F_0}|}=|D_E|^{1/4}\).

Put \(d_U=\sum_a r_a=\varrho(k+\varrho)\) and
\(\Lambda_0=\prod_{v,i<n}a_{v,i}^{k+2(\varrho-i)}\).
At each split residue pair of characteristic \(p\mid M\)
the free entries and centers have \(2d_U+\varrho\)
scalar coordinates over \(\mathbb F_p\).
All choices lift to the local integer lattice: use the
paired block chart, with the quadratic correction
given by an integral trace-one scalar in the center.
At the primes dividing \(M\) the middle lattice is
unmodified since \(N\) is a unit. Thus
\(m_U=M^{f(2d_U+\varrho)}\). Multiplying the measures for the individual rows, using \(r_a+1=k+2(\varrho-a)\) and
\(a_{v,\varrho}=1/N\), gives
\begin{equation}\label{eq:marked-unipotent-covolume}
 \nu\Lambda_0=m_U N^{-kf}v_E^{-d_U}|D_E|^{\varrho/4}.
\end{equation}
The inverse of the right fundamental box is a left
fundamental box of equal measure, by unimodularity.
Since the translates of the nonnegative cutoff
\(\beta_c\) sum to one as in Section~\ref{sec:period},
Tonelli gives \(\int_U\beta_c(u y)\,du=\nu\) for every
\(y\), in the specified real Haar. The value is
independent of the original class \(c\).

\subsection{The main and identity moments on the source}
\label{sec:main-identity-moments}

\sectionopening{In this subsection, we express both normalization constants using the source class count
and the dimensions of the real compact types. The mass comparison of
Section~\ref{sec:consecutive-unitary-masses} ensures that the ratio
used in Section~\ref{sec:full-orbital-blocks} is for the same measure.}

Write \(d_B=\prod_v|I_v|\). The source classes are free at
\(J_H\), and each class value is the tensor product of
the compact types of Section~\ref{sec:real-marked-filling}.
Thus \(d_X=d_Cd_B\). From \(|H_R|=m_U|\bar H|\) and
\eqref{eq:marked-source-class-count}, the definitions of
Section~\ref{sec:full-orbital-blocks} give
\[
 R_*:=\frac{D_{\rm main}}{D_{\rm id}}
   =\frac{d_B(N^k-1)^f\mu_k}{\mathcal B\mu_r m_U}.
\]
Set \(B=\log b\), \(Q_0=M\operatorname{rd}(E)/(2\pi)\),
and let \(S_\kappa\) be the sum of the logarithms of
the weights of Section~\ref{sec:datum} over the \(f\)
infinite places. The choice of \(N\) reads
\[
 \log N=mS_\kappa/f+(m-1)\log Q_0-\tfrac{15}{2}B+\delta,
 \qquad 0\le\delta\le\sqrt B .
\]
Set
\[
 W_*=\log\mathcal B+f\sum_{l=k-1}^{r-1}\log(l!),
 \qquad
 \varepsilon_m=\log(\mu_k/\mu_r)
                   -\sum_{j=k}^{r-1}\sum_v\log\sigma_{v,j+1}.
\]
Section~\ref{sec:consecutive-unitary-masses} gives
\(|\varepsilon_m|\le(m-1)2^{1-\varrho}\).
Lemma~\ref{lem:formal-degree-weight-comparison} gives
\[
 |km S_\kappa-W_*|
       \le m^2 f B + fmr/x .
\]
Using \(\sum_{j=k}^{r-1}(j+1)=(m-1)(k+\varrho+1)\)
and \(m_U=M^{f\varrho(2k+2\varrho+1)}\), the factorials
in the moment ratio cancel except at \(k-1\).
More explicitly,
\[
\begin{split}
 \log R_* &-fk(\log M-\tfrac{15}{2}B)
                    -\log\bigl(d_B((k-1)!)^f\bigr)\\
 &= fk\delta+f\log(1-N^{-k})+\varepsilon_m
    +kmS_\kappa-W_*\\
 &\quad-f(m-1)
   \left(\varrho\log Q_0+
          \log\frac{\sqrt{\operatorname{rd}(E)}}{2\pi}\right).
\end{split}
\]
The right side is \(o(fkB)\). Indeed the first term is
at most \(fk\sqrt B\), the second has absolute value
at most \(2f/N\), and the last term is \(O(fm^2B)\)
by the parameter bounds on \(Q_0\) and the root discriminant.
Here \(m^2/k=O(x^{-4})\) and \(r/k\) is bounded.

In the compact type, the dimension of the
symmetric-power representation of weight \(a_v\) is
\(\binom{a_v+k-1}{k-1}\). For \(a_v\le a=xb\) it is
at most \((a+k)^k\), also when \(a_v=0\).
Since \((a+k)k\le b^{5/4}\) for large \(x\),
\[
 1\le d_B((k-1)!)^f\le b^{5kf/4}.
\]
Using the logarithmic error bound \(fkB/4\) and
\(M\ge b^{20}\) now gives
\begin{equation}\label{eq:marked-main-identity-ratio}
 b^{11kf}D_{\rm id}\le D_{\rm main}
                   \le (M/b^6)^{kf}D_{\rm id}.
\end{equation}

\subsection{The open-orbit bound at equal parameters}
\label{sec:open-equal-link-bound}

\sectionopening{In this subsection, we insert the adelic estimate of
Section~\ref{sec:adelic-common-entry-sum} into the normalized open
orbital sum. The mass identities in
Sections~\ref{sec:consecutive-unitary-masses}
and~\ref{sec:main-identity-moments} determine its constant.}

\begin{proposition}[Open-orbit bound for a prescribed stabilizer]\label{prop:marked-open-orbital-bound}
Let \(\Omega\) be a finite set of open residues in the stabilizer
\(\Gamma_{\mathcal A}\). On the nonzero-discriminant set, the cut blocks
over these residues satisfy
\[
 \nu^{-2}\mathcal L(\Omega)\le
   (1+b^{-4})b^{-6k|\mathcal A|}D_{\rm main}.
\]
For each residue it includes all rational links, with the squared norm
taken separately for each link of the cut kernel. Passage to the
uncut kernel uses \eqref{eq:marked-full-cut-square}.
\end{proposition}

\begin{proof}

Set \(\Delta_0=R_0^{-2}\prod_{i<n}P_i\) and
\(C_0=C_*d_C C_{\rm S}\). The definitions of \(R_0\) and
\(P_i\) give the exact ratio
\[
 \Delta_0=\frac{v_E^n}{c_G^f}
                         \left(\frac{d_B}{\Lambda_0}\right)^2.
\]
The successive unitary compact volumes in the stage Haar
satisfy \(\operatorname{vol}U(l)/\operatorname{vol}U(l-1)
=2\pi^l/(l-1)!\), the unit sphere area. Hence from the
formulas for \(c_G\), \(\mathrm{fr}_v\) and \(C_{\rm S}\),
\[
 \left(\prod_v\mathrm{fr}_v\right)c_G^f |D_E|^{\varrho/2}
      =C_{\rm S} N^{2kf}v_E^{n+2d_U}.
\]
For this substitution use
\(\prod_v|\Phi_v(\delta_E)|=|D_E|^{1/2}\) and
\(\sum_{l=k+1}^r l=(2\varrho+1)(k+\varrho+1)\).
Using \(C_*=m_U^2|\bar H|\) and
\eqref{eq:marked-unipotent-covolume} now gives
\[
 C_0\Delta_0
  =\nu^2|\bar H| d_C d_B^2\prod_v\mathrm{fr}_v
  \le (1+b^{-9})\nu^2 D_{\rm main}\mathcal B .
\]
Indeed, substituting \(d_X=d_Cd_B\),
\(D_{\rm main}\mathcal B=d_X^2/\mu_r\) and
\eqref{eq:marked-source-class-count} gives
\[
 \frac{C_0\Delta_0}{\nu^2 D_{\rm main}\mathcal B}
 =\frac{\mu_r\prod_v\mathrm{fr}_v}{(N^k-1)^f\mu_k}
 \le(1+b^{-9})(1-N^{-k})^f
\]
by \eqref{eq:consecutive-frame-masses}. In fact the exact ratio is
\[
 (1-N^{-k})^{-f}\prod_{j=k}^{r-1}x_j^{-1},
 \qquad x_j=\prod_{v\nmid\infty}(1-\chi_v^{j+1}q_v^{-(j+1)}).
\]
Each split place of \(F_0\) over \(N\) therefore contributes
\(\prod_{l=k}^{r}(1-N^{-l})^{-1}\), including the
rank-\(k\) factor from the covector index. This factor can also be
read from the original compact at that place. In the last stage of
rank \(l=k+1\), the split components of \(z\) are a primitive vector and
a primitive covector of pairing \(-N\), since \(h(-x_{\rm gl},z)=1\)
for the glued column \(x_{\rm gl}=x_{k_0}/\bar\pi_N+z/N\).
Send the vector to the first standard basis vector. The covector
then has first coordinate \(-N\) and a unit among the remaining
coordinates, so a basis change on their span takes it to
\((-N,1,0,\ldots,0)\). This proves transitivity. The residue count is
\((N^l-1)(N^{l-1}-1)\), and the pairing has surjective differential
there. Division by \(N^{2l-1}\) gives
\((1-N^{-l})(1-N^{-(l-1)})\). The earlier isotropic stages supply
the exponents \(k+2,\ldots,r\), so the compact mass is
\(\prod_{l=k}^{r}(1-N^{-l})\) in the stage Haar.

For the decay term, put \(s=\sqrt{\operatorname{rd}(E)}\)
and \(h=r-k=2\varrho+1\). Put \(X_*=Q_0^h/(MN)\), so
\(\prod_v X_v=X_*^f\). The sphere-volume formula and the definition
of \(\mathrm{fr}_v\) give
\[
 \Theta=\frac{\prod_{l=k}^{r-1}l!}{2^h\pi^{h(\varrho+1)}},
 \qquad
 \left(\prod_v\mathrm{fr}_v\right)^{1/f}
  =\frac{N^k(2\pi)^{h(k+\varrho+1)}}
       {s^{2h(k+\varrho+1/2)}\prod_{l=k}^{r-1}l!}.
\]
Since \(D_{\rm id}=m_U|\bar H|d_Cd_B\) and
\(m_U=M^{f\varrho(2k+h)}\), the preceding equality for
\(C_0\Delta_0\) yields
\[
 \frac{C_0\Delta_0\mathcal D_{\rm dec}}{d_B D_{\rm id}\nu^2}
 =\frac{(\prod_v\mathrm{fr}_v)\Theta^f X_*^{f(k-n^7)}}{m_U}
 =\left[(2/M)^{\varrho h}s^{-h^2}X_*^{-n^7}\right]^f
 \le b^{n^9 f}.
\]
Here \(Q_0=Ms^2/(2\pi)\ge1\), \(M\ge2\), \(s\ge1\), and
\(X_*^{-1}\le MN\le b^{23+60n}\), with
\((23+60n)n^7\le n^9\) for large \(n\).
Also \(d_B\le b^{2kf}\), \(|\mathcal A|\le f\), and
\(k\ge512000n^9\) by the rank bounds. Consequently
\[
 2404n^9 f+2kf+10+6k|\mathcal A|\le 11kf.
\]
Applied to \eqref{eq:marked-main-identity-ratio},
this gives
\[
 C_0\Delta_0 b^{2403n^9 f}\mathcal D_{\rm dec}
          \le b^{-10}\nu^2 b^{-6k|\mathcal A|}D_{\rm main}.
\]

For any finite set of open cut blocks on the
nondegenerate orbital set, combine the two estimates
with \eqref{eq:adelic-top-sum-bound} and the summed
fixed-entry bound immediately following it.
The degrees in \(M_0\) cancel \(\mathcal B\), and
\((1+b^{-10})^f\le1+b^{-9}\) for large \(x\)
since \(f/b\to0\).
The remaining multiplier on \(b^{-6k|\mathcal A|}D_{\rm main}\)
is at most
\[
 (1+b^{-6})\bigl((1+b^{-9})^2(1+b^{-6})+b^{-10}\bigr)
          \le 1+b^{-4}.
\]
For any choice of a vector or covector stabilizer at each place of
\(\mathcal A\), the nonnegative sum over all tuples above the
given residues \(\Omega\) is the supremum of the sums over
\(\mathcal F\). By the definition in Section~\ref{sec:full-orbital-blocks}
this sum is \(\mathcal L(\Omega)\). Hence
\begin{equation}\label{eq:open-equal-link-established}
 \nu^{-2}\mathcal L(\Omega)\le
       (1+b^{-4})b^{-6k|\mathcal A|}D_{\rm main}
\end{equation}
on the open part of \(\Gamma_{\mathcal A}\) of
Section~\ref{sec:full-orbital-blocks}, with \(Z\) formed from those same
stabilizers as in Section~\ref{sec:finite-selected-bottom}.
\end{proof}

\section{Truncation, finite projections, and lattice pair counts}
\label{sec:orbital-pair-counts}

\sectionopening{In this section, we sum the truncation errors from Section~\ref{sec:real-cut-error},
control the residue projections, and count pairs of lattices in the
Petersson expansion. These estimates give the finite sums used in
Section~\ref{sec:orbital-cancellation-analysis}.}

\subsection{Reduced bases and the truncation error}
\label{sec:marked-full-to-cut}

\sectionopening{In this subsection, we choose representatives for rational isometries in the Petersson
expansion with controlled norms. This lets the real truncation
estimate of Section~\ref{sec:real-cut-error} be summed in the full
pairing of Section~\ref{sec:full-orbital-blocks}.}

Let \(\gamma_f\) be the finite matrix in the rational Witt basis,
\(B=B_{\pi_N}\), and \(\gamma_L=B^{-1}\gamma_fB\). By
Section~\ref{sec:marked-class-covolume} the matrix in the class bases is
\[
 B_c^{-1}\gamma_f B_{c'}=h_c^{-1}\gamma_L h_{c'}.
\]
Write \(\mathcal T(c,s,c')\) for the rational isometries with
this matrix integral and of residue \(s\) at the original
principal level.
We specify the class cutoff used in the orbital series.
Set \(A=\langle e_0,\ldots,e_{\varrho-1}\rangle_E\).
Write \(F=(e_0,\ldots,e_{\varrho-1},z)\) for the ordered basis of
\(A\oplus Ez\) and \(U_c^0\) for its pointwise stabilizer in \(U_c\). The action of \(U_c\) on \(F\) has
image the upper unitriangular matrices whose strict
upper entries belong to \(M\mathcal O_E\).

The containment follows because the basis columns and
the top coordinate functionals are unchanged in the integral class basis.
For surjectivity, realize an elementary entry \(M a\)
in row \(i\), column \(j>i\), by the integral Eichler
lift of Section~\ref{sec:marked-class-covolume}.
For \(f_j=e_j\) the coordinate covector \(e_j^\vee\)
is integral on the lattice. For \(f_j=z\) the coefficient of \(z\)
in the finite class basis on the middle space is
\[
 \ell_c(x+t z)=t-\pi_N^{-1}\bigl((c^\circ)^{-1}x\bigr)_{k_0}.
\]
Indeed it vanishes on \(c^\circ x_{k_0}/\bar\pi_N+z/N\)
and the other \(X_0\)-columns, and takes value one on \(z\).
To obtain a rational covector \(x+t z\mapsto t+\xi(x)\)
integral on the same lattice, the exact conditions are
\[
 (\xi c^\circ)_l\in\widehat{\mathcal O}_E\quad(l\ne k_0),\qquad
 (\xi c^\circ)_{k_0}+\pi_N^{-1}\in\bar\pi_N\widehat{\mathcal O}_E .
\]
They define a nonempty open set of finite adelic rows containing
the restriction of \(\ell_c\) to \(X_0\). Additive approximation gives
a rational row \(\xi\) satisfying all the conditions, at the places of
both \(\pi_N\) and \(\bar\pi_N\) as well. Extend by zero on the
hyperbolic pairs. In either case the covector is supported on the free
columns in row \(i\). Multiply it by \(a\) and apply the Eichler lift.
Its restriction to \(F\) adds just the stated elementary entry.

Choose a measurable set \(P\) of representatives for
archimedean tuples modulo \(M\mathcal O_E\), bounded in
each coordinate by \(b^{26}\). The \(2f\) elements
\(\zeta_{q_1}^j\), \(\sqrt{-d}\zeta_{q_1}^j\) for \(0\le j<f\)
span a full sublattice of \(\mathcal O_E\) by linear
disjointness. Its multiple by \(M\) has a half-open
parallelepiped bounded by \(2fM\max(1,\sqrt d)\le b^{26}\).
On this parallelepiped the quotient of \(M\mathcal O_E\)
by the smaller lattice acts by translation and reduction.
Choose the lexicographically first point in each finite
orbit in real coordinates. This gives the measurable \(P\).

At each infinite place write the normalized columns
\(y_v^{-1}f_a\) as their unnormalized orthogonal
Gram--Schmidt columns times an upper unitriangular
matrix. Read the matrix \(T_{\mathrm{GS}}(y)\) back in the
rational basis by multiplying entry \(a,j\) by
\(b_{v,f_a}/b_{v,f_j}\).
Call \(y\) reduced if the strict upper entries,
as tuples over the infinite places, lie in \(P\).
Replacing \(y\) by \(u y\) multiplies \(T_{\mathrm{GS}}(y)\)
on the right by the inverse basis matrix of \(u\).

In fact such a column operation preserves the
orthogonalized vectors. Right multiplication reduces
each column in descending row order. The entry added in row \(i\)
is a lattice shift plus terms determined by rows with larger index. The condition in \(P\)
determines that shift uniquely. Thus the elements of
\(U_c\) taking \(y\) to a reduced point form exactly
one coset of \(U_c^0\).

Multiply the reduced indicator
by a nonnegative cutoff for \(U_c^0\) whose translates
sum to one. A measurable section exists by taking
local sections of the discrete free left action on
\(G_\infty\) over a countable cover of its quotient
and partitioning by first chart. The product is
\(\beta_c\), with the orbit sum of
Section~\ref{sec:period} and the integral
of Section~\ref{sec:marked-class-covolume}.

Where \(\beta_c\) is nonzero, the normalized upper
entries have modulus at most \(b^{26}\). Indeed the
basis scales have nonincreasing moduli in column order.
The orthogonalized columns are invariant under the
Bessel subgroup acting on the left. Indeed that
subgroup fixes each successive line of the basis
flag, with diagonal coefficients one. Each of those
orthogonalized columns has norm at most the length
\(\ell_v(y_v)\) of Section~\ref{sec:real-cut-error}.
On any Bessel translate its norm is bounded by
the original column norm, which is at most the
operator norm of the translated matrix, since
\(\|g^{-1}\|=\|g\|\) in the normalized basis.
Summing at most \(n\le b\) contributions gives
\[
                 \|y_v^{-1}f_j\|\le b^{27}\ell_v(y_v).
\]

Put \(H_v(y_v)=\max(1,\ell_v(y_v)/b^{6n})\) and
\(L(y)=\prod_v H_v(y_v)\). Fix \(c,s,y\).
For \(R\ge1\) and a term \((c',\gamma)\) with
\(\beta_{c'}(\gamma^{-1}y)>0\) and
\(L(\gamma^{-1}y)\le R\), the transported basis
columns have norm at most \(B_R=b^{37n}R\).
For two such terms let \(y'=\gamma_1^{-1}y\) and
\(g=\gamma_1^{-1}\gamma_2\), an integral transporter
between their right classes reducing to one.
If the maximum distance of the two tuples of columns
is less than \((NB_R)^{-1}\), then
\(h(f_i,gf_j)=h(f_i,f_j)\) for every \(i,j\).
These pairings are integral on the class bases.
At each \(v\) the normalized pairing scales the
rational one by a number of modulus at least \(1/N\).
Pair \((y'_v)^{-1}f_i\) with the difference of the
two transported columns. The absolute value, after
division by the scale, is bounded by \(N B_R\) times
their distance. Every conjugate of the integral
difference then has modulus less than one.

The equality of pairings gives zero dual and \(z\)
coordinates of \(g f_j-f_j\). Invariance of the norm
and definiteness give zero \(X_0\)-component.
Compare the remaining top coordinates in descending row order.
Once the differences in the later rows vanish, a real unipotent
moving \(y'_v\) to operator norm at most \(b^{6n+2}R\) leaves the
difference in the current row unchanged by unitriangularity. The modulus of that coefficient
in the rational basis is therefore at most
\(b^{23\varrho} b^{6n+2}R\) times the column distance,
since the reciprocal of the basis-scale ratio is
at most \(b^{23\varrho}\).
The bound is less than one. The coefficient is
integral because the class basis preserves the specified columns and the top coordinate functionals. Its field norm forces it to
vanish. Thus \(gF=F\) pointwise.

If \(\gamma_1F=\gamma_2F\), then \(g=\gamma_1^{-1}\gamma_2\)
fixes \(F\) pointwise. Its action on \(A^\perp/A=X_0\oplus Ez\)
fixes \(z\), hence restricts to \(a\in H(F_0)\) on \(X_0\).
Removing \(\iota_{\rm rat}(a)\) leaves a basis-fixing unipotent.
Write \(d_1,d_2\) for the right classes. The two class matrices
\(B_c^{-1}\gamma_{j,f}B_{d_j}\) have the same residue \(s\), so
their relative matrix \(B_{d_1}^{-1}g_f B_{d_2}\) belongs to \(J_L\).
Its middle diagonal block is the glued matrix of
\((d_1^\circ)^{-1}a_fd_2^\circ\). Hence
\[
 (d_1^\circ)^{-1}a_fd_2^\circ\in J_H
\]
by Section~\ref{sec:period-quotient-measure}. Thus \(d_1=d_2\). For the fixed representative of this class,
Lemma~\ref{lem:marked-source-freeness} gives \(a=1\).
Consequently \(g\in U_{d_1}^0\). With \(\gamma_1\) fixed, the
terms with this basis image give distinct translations in
\(U_{d_1}^0\), so their cutoff weights sum to at most one by the orbit sum.

Distinct basis images give separated points in real
dimension \(2nrf\), inside a ball of radius \(2B_R\)
about one of the points. Disjoint balls of radius
\((2NB_R)^{-1}\) bound their number by
\((1+4NB_R^2)^{2nrf}\).
Since \(N\le b^{60n}\), this is at most
\((b^{6n^2}R^2)^{2nrf}\) for large \(n\).
Writing \(I=r^3 f>8nrf\) and integrating the layers
against \(I R^{-I-1}\,dR\) on \([1,\infty)\) gives
\begin{equation}\label{eq:marked-length-packing}
 \sum_{c'}\sum_{\gamma\in\mathcal T(c,s,c')}
     \beta_{c'}(\gamma^{-1}y)L(\gamma^{-1}y)^{-I}
                       \le 2b^{12n^3rf}.
\end{equation}

Let \(W=\bigotimes_v W_v\) and \(W^\flat=\bigotimes_v F_v\)
be the full and cut Bessel functions.
Section~\ref{sec:real-cut-error} supplies majorants
\(w_v\) with \(A_v=\sup w_v\) and
\[
 H_v^{3I}\|W_v\|\le w_v,\qquad
 A_v\int_{U_v\backslash G_v}w_v\le b^{400r^2}|I_v|.
\]
Put \(w=\prod_v w_v\) and \(A_*=\prod_v A_v\).
The majorants are invariant under the corresponding
real unipotent groups by the graph integral.
For \(W^\sharp=W,W^\flat\) write
\(K^\sharp_\gamma(y)=W^\sharp(y)\otimes
 \overline{W^\sharp(\gamma^{-1}y)}\).
In orthonormal source coordinates this tensor is the matrix
\(W^\sharp(y)W^\sharp(\gamma^{-1}y)^*\), with
\[
 \|K^\sharp_\gamma(y)\|_{\rm HS}
   =\|W^\sharp(y)\|\,\|W^\sharp(\gamma^{-1}y)\|.
\]
The kernels agree unless \(L(y)L(\gamma^{-1}y)\ge b\).
On the support of their difference, the majorants give
\[
 \|K^{\rm full}_\gamma(y)-K^{\rm cut}_\gamma(y)\|_{\rm HS}
   \le 2b^{-2I}w(y) A_* L(\gamma^{-1}y)^{-I}.
\]
The same majorant with the factor \(2b^{-2I}\) omitted
bounds either kernel everywhere.

Fix a subset \(S_0\) of the rational terms off the Bessel subgroup.
Define a block by
integrating \(\beta_c(y)\beta_{c'}(\gamma^{-1}y)K^\sharp_\gamma(y)\)
and summing over the terms of \(\mathcal T(c,s,c')\)
in \(S_0\). Tonelli,
\eqref{eq:marked-length-packing}, and the integral
\(\int_{U_\infty}\beta_c(u y)\,du=\nu\)
show absolute convergence over all right classes
and terms for fixed \(c,s\).
For the norm of the difference the same argument gives
\[
     4\nu b^{-2r^3 f+12n^3rf+400r^2 f}d_B
                            \le \tfrac12\nu b^{-k^3 f}.
\]
By Section~\ref{sec:main-identity-moments},
\(d_B((k-1)!)^f\le b^{5kf/4}\), hence \(d_B\le b^{2kf}\).
From \(\sqrt x\le m=2^{\lceil\log_2\sqrt x\rceil}<2\sqrt x\) we get
\(n=m/2<\sqrt x\). Since \(x^5\le r\) and \(k=r+1-m\),
\(12n^3+500\le k<r\) for sufficiently large \(x\). In particular
\[
 12n^3r+400r^2+2k+1\le402r^2\le r^3,\qquad k^3\le r^3.
\]
The exponent after bounding \(d_B\) is therefore at most
\(-(k^3+1)f\). Since \(4b^{-f}\le1/2\) for \(b\ge8\) and
\(f\ge1\), the stated bound holds for every subset of terms.

For fixed \(c,s\), put \(D_{c'}=A_{S_0}^{\rm full}(c,s,c')
-A_{S_0}^{\rm cut}(c,s,c')\). The bound just obtained sums over all
right classes, so
\[
 \sum_{c'}\|D_{c'}\|_{\rm HS}^2
 \le \left(\sum_{c'}\|D_{c'}\|_{\rm HS}\right)^2
 \le \tfrac14\nu^2b^{-2k^3f}.
\]
Freeness at \(J_H\) in Section~\ref{sec:ramified-source}
identifies \(X\) with one copy of the nonzero real coefficient space
per source class. Thus \(d_X=d_Cd_B\). Since \(d_B\ge1\) and \(|H_R|\ge1\),
\[
 d_C\le |H_R|d_Cd_B=D_{\rm id}.
\]
Since \(|G|\le M^{r^2 f}\le b^{23r^2 f}\), and
\(23r^2 f+40\le2k^3 f\), the following estimate holds.

\begin{proposition}[Error from the cut integral]\label{prop:orbital-truncation-error}
For every prescribed subset \(S_0\) of the rational terms off the
Bessel subgroup, the blocks formed from \(S_0\) satisfy
\begin{equation}\label{eq:marked-full-cut-square}
 \nu^{-2}\sum_{c,c',s}
     \big\|A_{S_0}^{\rm full}(c,s,c')-A_{S_0}^{\rm cut}(c,s,c')\big\|_{\rm HS}^{2}
             \le b^{-40}D_{\rm id}.
\end{equation}
\end{proposition}

\subsection{The residue estimates for the projection onto the selected types}
\label{sec:marked-residue-projection}

\sectionopening{In this subsection, we estimate the action of the finite idempotent of
Section~\ref{sec:upper-good-projector} on the level subgroups. The
bounds supply the loss estimates for the diagonal and mixed moments in
Section~\ref{sec:marked-moments}.}

For each \(q=p_i^{\rm lev}\), use the \(\mathrm{GL}_r(\mathbb F_q)\)-factors of \(G_i\).
On an algebraic closure of each residue field fix a faithful character
\(\jmath:\overline{\mathbb F}_q^\times\to\mathbb C^\times\).
For \(z\in\mathbb F_{q^a}^\times\) write
\(\jmath(z)=\exp(2\pi i b/(q^a-1))\). Label \(z\) by the character
\((\jmath|_{\mathbb F_{q^a}^\times})^b\). This is a bijection commuting
with \(q\)-th powers, and extending the residue field pulls the character
back by the norm. We use this choice again in
Section~\ref{sec:marked-inertia-characters}.
A local type parameter is a monic squarefree polynomial of
degree \(r\) with nonzero constant term. It is
allowed if no root in \(\mathbb F_q\) is a
\(q_2\)-th root of unity. A primitive type polynomial here is
the minimal polynomial of a generator of
\(\mathbb F_{q^r}^\times\). It is allowed since \(r>1\).
A tuple is selected if all its parameters are allowed and
at least one is primitive. The original projector
is the product of the two central idempotents \(h_i\)
selecting these tuples, because the type polynomials
are indexed by the same chosen places for both realizations.

For a squarefree type polynomial with distinct irreducible
factors of degrees \(d_1,\ldots,d_l\), the associated
Green data have one box in each of the corresponding
distinct cuspidal parameters. Green's character
parametrization and degree formula
\cite[\S8, Theorem~14 and Remark~(1), pp.~443--444]{Green}
give degree \(\prod_{a=1}^r(q^a-1)/\prod_{j=1}^l(q^{d_j}-1)\).
Distinct parameters give orthogonal central idempotents.
The identity coefficient of each is its degree squared
divided by the group order, namely
\[
 q^{-r}\frac{\prod_{a=1}^r(1-q^{-a})}
                    {\prod_j(1-q^{-d_j})^2}
              \ge q^{-r}(1-r/q)=u_q .
\]
Use weight \(u_q\) for each monic degree-\(r\) polynomial
to count the allowed parameters. Roots at zero or at one
of the \(q_2\)-th roots exclude at most
\((1+q_2)q^{r-1}\) polynomials.
A repeated monic factor of degree \(j\ge1\) excludes
at most \(q^j q^{r-2j}\le q^{r-1}\) polynomials
for each \(j\le r/2\).

Therefore the allowed weight
is at least \(1-\epsilon_q\), with
\(\epsilon_q=(1+q_2+2r)/q\).
There are \(\varphi(q^r-1)/r\) primitive type polynomials,
since all generators have degree \(r\) and Frobenius
orbits of size \(r\). If \(L\ge2\) has \(j\) distinct
prime divisors, then
\[
 \frac{\varphi(L)}{L}\ge
 \prod_{a=1}^j\left(1-\frac{1}{a+1}\right)
        =\frac1{1+j}\ge\frac1{1+2\log L}.
\]
So the unweighted fraction of primitive polynomials
is at least \(\delta_q=1/(2r(2r\log q+1))\).
Write \(A_q\) for the allowed parameters and \(P_q\subset A_q\)
for the primitive ones. Their comparison weights satisfy
\[
 a_q=u_q|A_q|\ge1-\epsilon_q,\qquad
 c_q=u_q|A_q\setminus P_q|
        \le q^{-r}(q^r-|P_q|)\le1-\delta_q .
\]
Each selected tuple has identity coefficient at least \(u_q^f\).
Thus \(h_i(1)\ge a_q^f-c_q^f\), giving
\[
 (1-h_i)(1)\le f\epsilon_q+(1-\delta_q)^f
           \le f\epsilon_q+\exp(-f\delta_q).
\]
For fixed \(q_0,d,\epsilon\), Section~\ref{sec:datum} gives
\[
 f=3\varphi(q_0)L_0m^{22}\asymp_{q_0}x^{11}\log x,\qquad
 q_2=(q_0d)^{s+6}\le(q_0d)^{28}x^{11\log_2(q_0d)} .
\]
Since \(r<3x^5\), \(10B\le\log q\le11B\) and \(B\sim x^\epsilon\),
\[
 f\delta_q\gg_{q_0}x^{1-\epsilon}\log x\ge9B,\qquad
 f\epsilon_q\le f(1+q_2+2r)b^{-10}=o(b^{-8})
\]
for all sufficiently large \(x\). It follows that
\begin{equation}\label{eq:marked-bad-projector}
                         (1-h_i)(1)\le(2b^8)^{-1}.
\end{equation}

For the Bessel trace on \(G_i\), let \(\psi_U\) be the finite
character of the integral unipotent subgroup on its reduction, extended
by zero off that image. For two lifts of the same reduction the
simple-root sums differ by an element of \(M\widehat{\mathcal O}_E\),
by the level character calculation of Section~\ref{sec:period-real-point}. After division by \(M\delta_E\) the
difference is in the inverse different, so this character is well-defined.

For a diagonal class term \(\gamma=u\iota_{\rm rat}(\delta)\), put
\(y_H=(c^\circ)^{-1}\delta_f c^\circ\). In the finite lattice basis
\[
 h_c^{-1}\gamma_L h_c=u_c\iota_L(y_H),\qquad
 u_c=h_c^{-1}L^{-1}u_f Lh_c .
\]
This transporter and its inverse are integral. Its middle block is
that of \(\iota_L(y_H)\), since the unipotent acts trivially on the
middle quotient. The inverse gives the same assertion for
\(y_H^{-1}\). By the glued lattice calculation of
Section~\ref{sec:fixed-vector-common-pairing}, \(y_H,y_H^{-1}\) preserve
the original lattice of \(X_0\), also above \(N\).
For \(s\in G_i\subset G_0\times G_1\), the transporter has reduction
\(s\) in the \(i\)-th factor and the identity in the other factor.
Thus its middle block reduces to the identity at each chosen place
over \(p_{1-i}^{\rm lev}\). At this prime the gluing basis change
and its inverse are integral because \(N\) is a unit, so \(y_H\) also
reduces to the identity there. Applied to the full integral definite
level with this one principal congruence, Lemma~\ref{lem:marked-source-freeness}
gives \(\delta=1\). Hence a diagonal term requires \(s\) in the reduced
unipotent subgroup.

Conversely let \(k\) be an integral unipotent of reduction \(s\).
For each \(c\), additive approximation in the rational unipotent
group, using its Heisenberg vector and skew central coordinates,
gives a rational \(u\) with \(h_c^{-1}L^{-1}u_fLh_c\in kJ_L\).
The class representative preserves the simple-root sum, so
\(\Psi_f(u)=\psi_U(s)\). By Section~\ref{sec:fixed-vector-class-phases}
the real covariance is \(W(uy)=\Psi_\infty(u)W(y)\), with
\(\Psi_\infty(u)=\psi_U(s)^{-1}\). For the identity residue the
diagonal Bessel amplitude is \(\nu I\) by the cutoff orbit sum and
orthonormality of the real Bessel coordinates. For the lift of \(s\)
its tensor is \(W(y)W(u^{-1}y)^*\), so with \(A_{\rm Bes}\) denoting
the Bessel part of the full amplitude,
\[
 A_{\rm Bes}(c,s,c)=\nu\psi_U(s)^{-1}I,\qquad
 A_{\rm Bes}(c,s^{-1},c)=\nu\psi_U(s)I .
\]
The entrywise conjugation and division by \(\nu\) in
\eqref{eq:marked-petersson-block}, summed over the \(d_C\) class
blocks of size \(d_B\), give
\(\operatorname{tr}r_\psi(s)=d_X\psi_U(s)\) for every \(s\in G_i\).

Take \(T\) to consist of all places above \(M\), and set
\(J_T=K_T(M)\), so \(G=K_T/J_T\) and
\(\mu(K_T)=|G|\mu(J_T)\).
For \(k\in K_T\), writing \(k=uj\) with \(u\in U\), \(j\in J_T\)
forces \(u\in K_T\). It extends by the identity off \(T\), and
each global integral unipotent gives such a local lift. The
definition on \(UJ_T\) and extension by zero therefore give
\(\chi(k)=\psi_U(\bar k)\). If \(\psi_U^\natural\) is the probability
conjugacy average on \(G\), then
\[
 \int_{K_T}\chi(yky^{-1})\,d\mu(y)
 =\mu(J_T)\sum_{t\in G}\psi_U(t\bar k t^{-1})
 =\mu(K_T)\psi_U^\natural(\bar k).
\]
Restrict the outer nonnegative integral in \(\mathcal M(K_T)\) to
\(K_T\). Each reduction class has mass \(\mu(J_T)\), so
\eqref{eq:finite-bottom-mass} with \(Q_0=2\) and \(t=n\) gives
\[
 \mu(J_T)\mu(K_T)^2\sum_{s\in G}|\psi_U^\natural(s)|^2
 \le\mathcal M(K_T)
 \le\mu(J_T)\mu(K_T)(\mu(K_T)+2^{-n}).
\]
Put \(q_i=p_i^{\rm lev}\). The \(2f\) places of \(E\) above
each \(q_i\) form \(f\) conjugate pairs, each giving one split
unitary factor. The calculation in the proof of
Lemma~\ref{lem:marked-finite-fixer-mass} applies at either \(q_i\),
since \(N\) is a unit there. At each isotropic stage of rank \(l\)
the primitive column and covector contribute
\((1-q_i^{-l})(1-q_i^{-(l-1)})\). At the last stage, of rank
\(k+1\), the prescribed unit norm contributes \(1-q_i^{-(k+1)}\).
Since \(r=k+2n-1\), these factors use each exponent from \(k+1\)
through \(r\) once. The integral Heisenberg groups and the terminal
compact each have mass one in Section~\ref{sec:compact-floor-haar}.
Thus in the Haar measure used for \(\mathcal M(K_T)\),
\[
 \mu(K_T)=\prod_{i=0}^1\prod_{l=k+1}^r(1-q_i^{-l})^f
 \ge 1-f\sum_{i=0}^1\frac{q_i^{-(k+1)}}{1-q_i^{-1}}
 \ge \tfrac12\ge2^{-n}
\]
for all sufficiently large \(x\), since \(f\le x^{12}\) and
\(q_i\ge b^{10}\).
Since \(G\) is the product of the two level factors,
conjugacy averaging on \(G_i\) agrees with the
restriction of averaging on \(G\).
Consequently the conjugacy average of
\(\operatorname{tr}r_\psi|_{G_i}\) has squared
counting norm at most \(2d_X^2\), as required
in the diagonal estimate.

\subsection{The weighted count at the additional ramified place}
\label{sec:marked-ramified-count}

\sectionopening{In this subsection, we bound the sublattices that occur in the fixed-vector contributions
of Section~\ref{sec:fixed-tail-comparison}. Their ramified level
weights enter the count of pairs in
Section~\ref{sec:marked-pair-count}.}

Write \(\mathcal O_f=\widehat{\mathcal O}_E\) and
\(\Lambda_c^0=\Lambda_c\cap X_0(\mathbb A_f)\).
A sublattice \(L\) of \(\Lambda_c^0\) has rank \(\dim_E V\)
if it has finite index in \(\Lambda_c^0\cap V(\mathbb A_f)\)
for an \(E\)-subspace \(V\subset X_0\).
For a finite idele \(a\) put \(\operatorname{ind}(a)=|a|_f^{-1}\).
The discriminant \(d(L)\) is the square root of the
index of the determinant of the Hermitian Gram matrix
in a finite basis of \(L\). It is independent of basis
and at least one for an integral lattice.
For \(L\subseteq L'\) of index \(\iota\) a basis
change gives \(d(L)=\iota d(L')\).

Put \(\Delta_{\rm lat}=b^{n^3 f}\) and \(\mathfrak t=b^{kf/4}\).
For \(0\le l\le n^3\), the weighted sum satisfies
\begin{equation}\label{eq:marked-ramified-weight}
 \sum_c\ \sum_{\substack{L\subseteq\Lambda_c^0\\
          \operatorname{rank}L=k-l,\ d(L)\le\Delta_{\rm lat}}}
      d(L)^{-k}\le \frac{d_C}{N^{kf}}\mathfrak t^{\,l+1}.
\end{equation}

Use the group \(K_X\) of Section~\ref{sec:marked-class-covolume}. Let
\(K_N\) be its subgroup obtained by requiring the distinguished
coordinate covector \(e_{k_0}^{\vee}\) to be fixed modulo the idele
\(\pi_N\), where \([\widehat{\mathcal O}_E:\pi_N
\widehat{\mathcal O}_E]=N^f\). The actual source level lies in
\(K_N\). The stabilizers for the source classes are trivial
at the principal level. Therefore each class has the same
positive weight in quotient Haar measure. For any nonnegative
function \(F_0\) on the classes, a uniform bound
\[
 |K_X/K_N|^{-1}\sum_{\phi\in K_X/K_N}
           F_0([u y_\phi])\le B_0
 \quad(u\in U(X_0)(\mathbb A_f))
\]
implies \(\sum_c F_0(c)\le d_CB_0\).
Indeed integrate the translate averages over the rational
quotient, using right invariance of Haar. For unbounded \(F_0\)
the same argument first applies to its finite truncations.

For fixed \(u\) and \(\phi\), write the chosen representative
as \(\eta u y_\phi j_H\), with \(\eta\) rational unitary and
\(j_H\in J_H\subset K_N\).
The class basis is \(B_c=\iota_{\rm rat}(c^\circ)B\). By the glued
intersection in Section~\ref{sec:fixed-vector-common-pairing},
\[
 \Lambda_c^0=c^\circ\{v\in\mathcal O_f^k:v_{k_0}\in\pi_N\mathcal O_f\},
 \qquad [c^\circ\mathcal O_f^k:\Lambda_c^0]=d(\Lambda_c^0)=N^f.
\]
Indeed \(c^\circ\mathcal O_f^k\) is unimodular and the distinguished
coordinate quotient is \(\mathcal O_f/\pi_N\mathcal O_f\).
Consequently, with \(U=u\mathcal O_f^k\),
\[
 \eta^{-1}\Lambda_c^0=
 \{x\in U:e_{k_0}^{\vee}y_\phi^{-1}u^{-1}x\equiv0\pmod{\pi_N}\}.
\]
Transport by \(\eta^{-1}\) preserves rational subspaces and
discriminants. The map \(yK_N\mapsto e_{k_0}^{\vee}y^{-1}\bmod\pi_N\)
identifies \(K_X/K_N\) with the product of the nonzero covectors at
the \(f\) places \(w\mid\pi_N\), by split reduction and transitivity
in Section~\ref{sec:marked-class-covolume}. For a lattice \(L\subset U\)
of rank \(h_0=k-l\), put
\(r_w=\dim_{\mathbb F_N}\operatorname{im}(u^{-1}L_w\to\mathbb F_N^k)\).
The covectors at \(w\) vanishing on \(L\) are the nonzero vectors of
its annihilator. The fraction of admissible cosets is therefore
\[
 \frac{1}{[K_X:K_N]}
 \#\{\phi : L\subseteq\eta^{-1}\Lambda_c^0\}
       =\prod_{w\mid\pi_N}\frac{N^{k-r_w}-1}{N^k-1}.
\]
It vanishes if \(r_w=k\) at any chosen place. Unitarity fixes the
action at the conjugate places by adjunction.

Let \(V\) be the rational subspace of \(L\). Put \(L^{\rm sat}=V(\mathbb A_f)\cap U\),
\(\iota=[L^{\rm sat}:L]\), and \(d=d(L^{\rm sat})\). Thus
\(d(L)=d\iota\) and \(d,\iota\ge1\). At \(w\mid\pi_N\), the saturated lattice is a direct summand of
\(U_w\). In Smith normal form for \(L_w\subseteq L_w^{\rm sat}\),
exactly \(h_0-r_w\) exponents are positive. Thus
\[
 \prod_{w\mid\pi_N}N^{h_0-r_w}\le\iota,\qquad
 \prod_{w\mid\pi_N}(N^{k-r_w}-1)\le N^{fl}\iota .
\]

Put \(\mu=2\varrho+4\). On \(d\iota\le\Delta_{\rm lat}\), one has
\[
 \iota(d\iota)^{-k}
 \le \Delta_{\rm lat}^{\mu+2} d^{-(k+1)}\iota^{-(h_0+\mu)}.
\]
Indeed the ratio to the right-hand power of \(d\) and
\(\iota\) is \(d\iota^{h_0+\mu+1-k}\le(d\iota)^{\mu+2}\).
The sum of \(\iota^{-(h_0+\mu)}\) over finite-index sublattices
in \(L^{\rm sat}\) is bounded by \(\zeta_E(\mu)^{h_0}\).
Locally, triangular echelon bases with diagonal ideals
\(\mathfrak p^{a_1},\ldots,\mathfrak p^{a_{h_0}}\) have at most
\(q^{\sum_{j=1}^{h_0}(h_0-j)a_j}\) choices of off-diagonal entries,
with \(q\) the residue cardinality. Multiplying the resulting
geometric series gives factors \(\zeta_E(h_0+\mu-j)\) with
\(0\le j\le h_0-1\), all bounded by \(\zeta_E(\mu)\).
The bound \(n[E:\mathbf Q]\le 2^{2\varrho}\) and
\(\zeta_E(\mu)\le \zeta(\mu)^{[E:\mathbf Q]}\) yield
\(\zeta_E(\mu)\le2\); for the latter use
\(\log\zeta(\mu)\le (1+2/(\mu-1))2^{-\mu}\).
Thus the index sum is at most \(2^k\).

Write \(\mathrm{rd}(E)=|D_E|^{1/(2f)}\),
where \(D_E\) is the field discriminant. The subspace sum is
\begin{equation}
 \sum_{\operatorname{codim}V=l}d(V(\mathbb A_f)\cap U)^{-(k+1)}
 \le (200 f k\,\mathrm{rd}(E))^{f(k+1)l}.
 \label{eq:marked-discriminant-sum}
\end{equation}
For \(l=0\), the only subspace is \(X_0\), its saturation in \(U\)
is \(U\), and \(d(U)=1\). Its perpendicular lattice is zero with
empty Gram determinant one. For \(l>0\), the primitive lattices of \(V\) and \(V^\perp\)
in \(U\) have the same discriminant.
In fact, local unimodularity identifies the quotient by
their direct sum with either discriminant group. The pairing from \(U\)
onto the dual of the first primitive lattice is
surjective because a functional on a direct summand
extends to \(U\). Its kernel is the second primitive lattice.
The image of the first lattice is its own pairing
lattice in the dual. Interchange the two lattices to get
the other identification. Their discriminants, defined
as the square roots of the orders of these groups,
are therefore equal.

Under one embedding from each complex pair, with the standard
Euclidean norm at the infinite places, the real covolume of
\(V^\perp(E)\cap U\) in its span is
\[
 2^{-fl}|D_E|^{l/2} d(V(\mathbb A_f)\cap U).
\]
Choose a rational basis matrix \(Y\) of \(V^\perp\)
and a finite coefficient matrix \(D\) so that \(YD\)
is an adelic lattice basis. The rational coefficient
lattice has real covolume \(v_E^{-l}\operatorname{ind}(\det D)\),
by comparison with a common integral sublattice.

At \(v\), the real Jacobian of \(Y\) in the positive
definite coordinates is \(\det(\Phi_v(Y)^*\Phi_v(Y))\).
The product of these Jacobians times
\(\operatorname{ind}(\det D)\) is \(d(V(\mathbb A_f)\cap U)\),
by the product formula for \(\det(Y^*Y)\).
Minkowski's second theorem \cite[Ch.~VIII, \S4.3, Theorem~V]{Cassels} in dimension \(2fl\), using the
unit Euclidean ball, gives independent \(E\)-vectors
\(x_1,\ldots,x_l\) in this lattice such that
\[
 (\|x_1\|\cdots\|x_l\|)^{2f}
       \le (8 f l\,\mathrm{rd}(E))^{fl} d(V(\mathbb A_f)\cap U).
\]
Let \(\lambda_h\) be the real Euclidean successive minima. After
choosing \(a-1\) independent \(E\)-vectors, their \(E\)-span has real
span of dimension \(2f(a-1)\). Among \(2f(a-1)+1\) real-independent
lattice vectors of lengths at most \(\lambda_{2f(a-1)+1}\), choose
\(x_a\) outside that span. Then
\(\|x_a\|^{2f}\le\prod_{h=2f(a-1)+1}^{2fa}\lambda_h\).

Multiplying over \(a\) gives
\((\prod_a\|x_a\|)^{2f}\le\prod_{h=1}^{2fl}\lambda_h\).
The ball contains the cube of side \(2/\sqrt{2fl}\),
so Minkowski's upper bound with the stated covolume already
gives the displayed inequality with room in its constant.

For distinct \(x,y\in E^k\cap U\), the hermitian square
\(h(x-y,x-y)\) is a totally positive algebraic integer by unitary
integrality at the finite places. Its norm is a positive integer,
so at least one real conjugate is at least one. The diagonal
archimedean images of \(x\) and \(y\) therefore have distance at least one. The packing bound in real
dimension \(d_0=2fk\) gives at most \((1+2R)^{d_0}\) lattice
vectors of length at most \(R\). Dyadic summation at exponent
\(2f(k+1)=d_0+2f\) then gives
\[
 \sum_{x\ne0}\|x\|^{-2f(k+1)}
 \le 5^{d_0}\sum_{j\ge0}2^{-2fj}\le 5^{2f(k+1)} .
\]
The chosen short tuple determines
\(V\) as its orthogonal complement. Raising the short-basis
bound to the \((k+1)\)-st power and summing over tuples
gives \eqref{eq:marked-discriminant-sum}, since \(8\cdot5^2=200\).

The unaveraged sum over \(\phi\) and \(L\) is consequently
at most
\[
  \Delta_{\rm lat}^{\mu+2}2^k
   \left(N^f(200fk\,\mathrm{rd}(E))^{f(k+1)}\right)^l
   \le b^{-f}\mathfrak t^{l+1}.
\]
For the last inequality use \(k\le f\),
\(k/4\ge10^5n^9\), \(N\le b^{60n}\),
\((200fk\,\mathrm{rd}(E))^{1000}\le b\), and
\(60n\cdot1000+(k+1)\le250k\).
In particular, \(\operatorname{rd}(E)\le dq_1\) by
Section~\ref{sec:datum}. Together with \(f,q_1\le x^{12}\) and
\(k<3x^5\), this gives \(200fk\,\operatorname{rd}(E)=O_d(x^{29})\).
By the class count in
Section~\ref{sec:marked-class-covolume},
\([K_X:K_N]=(N^k-1)^f\ge(N^k/b)^f\). Dividing by this
index and applying the Haar average
gives \eqref{eq:marked-ramified-weight}.

\subsection{Counting orbital representatives containing a fixed lattice}
\label{sec:marked-labelled-frames}

\sectionopening{In this subsection, we count bases adapted to the isotropic flag and its final orthogonal
vector inside a fixed lattice. Organizing them by their pairing
matrices gives the multiplicity bound used in
Section~\ref{sec:marked-pair-count} with the representatives of
Section~\ref{sec:marked-full-to-cut}.}

In the rational Witt decomposition write \(A\) for the
isotropic span of \(e_0,\ldots,e_{\varrho-1}\) and take
\(F=(e_0,\ldots,e_{\varrho-1},z)\) as a matrix of columns.
Write \(Z_\gamma=F^*h\gamma F\) for the matrix of their
products and \(P_\gamma\) for the submatrix of its first \(\varrho\) rows.
Here \(h\) is the matrix of the rational Hermitian form.
Use the same reduced cutoffs as in
Section~\ref{sec:marked-full-to-cut}. A boxed point is
reduced and has \(\ell_v\le b^{6n+2}\) at every place.
A term \(\gamma\) is near if \(y\) and \(\gamma^{-1}y\)
are boxed for some \(y\).
For an orbital calculation one fixes a representative of
each \(U_c^0\)-\(U_{c'}^0\) double coset, choosing a near
representative when there is one.

Let \(\delta_a\) be the leading minor of order \(a\) in \(Z_\gamma\),
with \(\delta_0=1\), and put
\[
 \mathfrak g_\gamma=\delta_n\overline{\delta_n}
             -N^2\delta_{n-1}\overline{\delta_{n-1}} .
\]
We say that \(\gamma\) has cell rank \(j\), \(0\le j\le n\), if
\(\delta_j\ne0\), if \(\mathfrak g_\gamma\ne0\) is equivalent to \(j=n\),
and if there exists \(\zeta\in E\) with \(\zeta\bar\zeta=1\)
such that the following holds for every \(j\le a<n\).
The spans of the first \(a+1\) columns of \(F\) and \(\gamma F\)
contain a common vector whose two coefficients at index \(a\)
are \(\zeta\) and \(1\), respectively.
Section~\ref{sec:marked-orbital-ranks} obtains these conditions
for the terms used there.

For \(\gamma\in\mathcal T(c,s,c')\), define its term
lattice \(L_\gamma\) to consist of \(x\in\Lambda_c^0\) for which
some lift \(x+M w\), \(w\in A(\mathcal O_f)\),
is orthogonal to \(\gamma F\). In the rational
coordinates this means \(X_\gamma^*x+MP_\gamma^*w=0\)
for some integral \(w\), with \(X_\gamma\) the
definite block of \(\gamma F\).
Choose representatives \(a_{sH_R}\) for the cosets
\(sH_R\) in \(G\), and define \(t\in H_R\) by
\(s=a_{sH_R}t\). This is the parameter of the term.

Fix the left
source class \(c\), the Bessel parameter \(t\), and the
lattice \(L\) of rank \(h_0=k-l\) in \(\Lambda_c^0\),
with discriminant at most \(\Delta_{\rm lat}\). A term counted here
is a pair \((c',\gamma)\), where \(\gamma\) is the fixed near
representative, has a cell rank, and has \(L\) in its
level-\(M\) term lattice. All moves in the proof are in
\(U_c^0\), the rational unipotents of level \(M\) in class
\(c\) fixing the standard basis pointwise.

Suppose the moved rational representatives \(g_1,g_2\) satisfy
\(g_1F=g_2F\). Write \(d_1,d_2\) for their right classes and
\(s_1,s_2\) for their residues. The rational isometry
\(m=g_1^{-1}g_2\) fixes \(F\) pointwise. On \(A^\perp/A\) it
fixes \(z\) and restricts to an isometry of \(X_0\); removing
this isometry leaves a basis-fixing unipotent. Thus \(m\) is a
rational Bessel term between the two right classes, with residue
\[
 \operatorname{red}(B_{d_1}^{-1}m_f B_{d_2})
   =s_1^{-1}s_2\in H_R
\]
by the integral Bessel calculation following
Proposition~\ref{prop:marked-petersson-blocks}. Hence
\(s_1H_R=s_2H_R\). Writing \(s_i=a_{s_iH_R}t\) with the common
parameter gives \(s_1=s_2\). The relative transporter is now in
the principal level. Its middle block identifies \(d_1=d_2\)
at the original source level, and
Lemma~\ref{lem:marked-source-freeness} makes its rational definite
component the identity. Thus \(m\in U_{d_1}^0\) as in
Section~\ref{sec:marked-full-to-cut}. The original representatives
belong to the same double coset and hence coincide.

We partition the terms by cell rank \(j\) and integral basis
product matrix \(Z= (\langle f_a,\gamma f_b\rangle)\), where
\(f_0,\ldots,f_{n-1}\) is the standard basis. There are at most
\(n+1\le b^{2f}\) cell ranks. The entries of \(Z\) are integral because
the two bases are integral in the class basis,
whose Gram form is integral. Their absolute values are bounded by \(N b^{74n}\le b^{134n}\) at all
complex places.

Indeed both normalized basis columns
have norm at most \(b^{37n}\) by the box bound,
and the product of the two basis scales has
modulus at least \(1/N\). The lattice of integral elements of
\(E\) in the product of the \(f\) complex embeddings has
point separation at least one in the maximum norm, by the
norm formula for a nonzero algebraic integer. Disjoint
products of discs of radius \(1/2\) show that there are at
most \((1+2b^{134n})^{2f}\) values for each entry. Since
\(1+2b^{134n}\le b^{135n}\), the number of matrices is at
most \(b^{270n^3 f}\).

Fix \(j,Z\), let \(Y\) be the matrix formed by the first
\(j\) rows of \(Z\) and the remaining rows of the identity,
and put \(\delta=\det Y\ne0\). Write \(P\) for the
submatrix of the first \(\varrho\) rows of \(Z\). For a term let \(X_\gamma\)
be the \(X_0\)-block of its basis. In the splitting
\(A\oplus(X_0\perp E z)\oplus A^\vee\), the block along
\(A^\vee\) is \(P\), and \(Z\) also determines the
\(z\)-block by pairing with \(z\) of norm \(-N\).
These statements are in the rational basis. In the
finite calculation both transported bases lie in the
same left class lattice.

\smallskip
\subsubsection{Changes of basis on the two blocks}
\label{sec:ordered-basis-block-moves}

\sectionopening{In this subsection, we put the two blocks into the coordinates needed for the count of Section~\ref{sec:marked-labelled-frames}.}

For a matrix \(\eta\) of rational columns in
\(\Lambda_c^0\), a move takes
\begin{equation}
 X_\gamma\longmapsto X_\gamma-M\eta P .
 \label{eq:marked-frame-shear}
\end{equation}
To check it in the indicated splitting, let
\(T=-M\eta:A^\vee\to X_0\perp E z\), and let \(J_0\)
be the middle Gram matrix. Choose an integral scalar
\(\theta\) with \(\theta+\bar\theta=1\), for example a
primitive sixth root in \(E\). The block unipotent
\[
 \begin{pmatrix}
 1&-T^*J_0&-\theta T^*J_0T\\
 0&1&T\\
 0&0&1
 \end{pmatrix}
\]
preserves the form, fixes \(A\) and \(z\), and is
congruent to the identity on \(\Lambda_c\) modulo \(M\).
For the congruence, use integrality of the column
vectors of \(\eta\) in the middle lattice and
of the pairings there. Its middle action on the basis
gives \eqref{eq:marked-frame-shear}. Also an integral skew-Hermitian matrix
\(\sigma\) yields the move with upper-right block
\(M\sigma\) and other off-diagonal blocks zero. It changes
the rows of the basis matrix indexed by \(A\) by adding \(M\sigma P\)
and leaves the other rows fixed. Both moves retain the original residue and
the cell rank.

\smallskip
\subsubsection{Lifts of the fixed lattice}
\label{sec:ordered-basis-lattice-lifts}

\sectionopening{In this subsection, we bound the number of lifts after the changes of basis in Section~\ref{sec:ordered-basis-block-moves}.}

Write a finite basis of \(L\) as \(m=y a\), with \(y\)
a rational matrix of independent columns and \(a\)
invertible over the finite adeles. From containment
in the term lattice we have a rational matrix \(w\)
with \(w a\) integral and
\begin{equation}
 X_\gamma^* y+ M P^* w=0.
 \label{eq:marked-frame-lift}
\end{equation}
The defining containment solves this equation
over the finite adeles with \(w a\) integral. A linear
system over \(E\) with a finite-adelic solution has a
rational solution. Its solutions form an affine space
over the rational kernel, whose rational points are
dense in the finite-adelic solutions by additive
approximation. The condition that \(w a\) be integral
is an open neighbourhood of the given solution.
This establishes \eqref{eq:marked-frame-lift} with its original lattice.

Sort the \(\varrho\) rows of \(w a\) by their classes
modulo the lattice generated by the rows of \(m^*m\). The number of classes
is at most \(d(L)^{2\varrho}\), since the Gram determinant
has index \(d(L)^2\). In one class, choose a base lift
\(w_0\). Each \(w_0-w\) is \(\eta^* y\) for integral
columns \(\eta\) in the lattice, given by
\[
 \eta=y(y^*y)^{-1}(w_0-w)^*.
\]
The row congruences mean \((w_0-w)a=S^*m^*m\) with
\(S\) integral. Cancelling \(a\) shows \(\eta=m S\) over
the adeles, so its columns indeed lie in \(L\).
Equation~\eqref{eq:marked-frame-shear} makes \eqref{eq:marked-frame-lift} hold with the common \(w_0\).
All moved bases hence have the same matrix of products
of their \(X_0\)-block with \(y\).

\smallskip
\subsubsection{Counting definite columns}
\label{sec:ordered-basis-definite-columns}

\sectionopening{In this subsection, we bound the definite columns for a lift chosen in Section~\ref{sec:ordered-basis-lattice-lifts}.}

For a fixed lift, the differences between two
\(X_0\)-blocks lie columnwise in the lattice
\(L^{\perp}_c=\Lambda_c^0\cap V(\mathbb A_f)^\perp\),
of rank \(l\), with \(V\) the rational span of \(L\).
For each column, write its rational middle coordinates as \(x+\lambda z\)
and put \(v=(c^\circ)^{-1}x\). If the column is \(\gamma f_b\), then
\(\lambda=-Z_{n-1,b}/N\) by pairing with \(z\). Thus the two rational
coefficients of \(z\) agree before passing to class coordinates.
The columns of \(F\) are unchanged in each right class basis, so
\(B_c^{-1}\gamma F\) is integral. Its middle coordinates are
\[
 \bigl((v_i)_{i\ne k_0},\ \bar\pi_Nv_{k_0},
                         \ \lambda-v_{k_0}/\pi_N\bigr).
\]
Subtracting gives integral differences in the other coordinates
and a difference in \(\pi_N\mathcal O_f\) at \(k_0\). Hence the
difference of the two \(x\)'s lies in \(\Lambda_c^0\).
Their difference pairs to zero with \(L\) by
\eqref{eq:marked-frame-lift} with the common \(w_0\).
Taking residues modulo \(M\delta L^{\perp}_c\) contributes at most \(\operatorname{ind}(M\delta)^{ln}\).

In a residue class let \(v=X_{\gamma'}-X_\gamma\)
be the difference. The matrix \(vY^{-1}/M\) has
columns integral in \(\Lambda_c^0\), since
\(v\) is divisible there by \(M\delta\) and
\(\operatorname{adj}Y\) is integral. Its columns at
positions \(b\ge j\) vanish because the cell
relations supply a shared vector at each such
position, with translated leading coefficient
one, and with higher coefficients zero. Their
translated coefficients in the tail can be combined
triangularly to give tail \(e_b\). The first \(j\)
products against the standard basis vanish, so
the complete coefficient vector is \(Y^{-1}e_b\).
Its image lies in the span of the columns of \(F\) and hence
has zero \(X_0\)-component. The same argument applies after a move, since the move
fixes the shared vectors.

The remaining exceptional position is the last
column in the top cell, where \(Y=Z\). The vector
\(\gamma F Z^{-1}e_{n-1}\) is orthogonal to \(A\)
and has \(z\)-coordinate \(-1/N\). Since it is in
the translated span of the basis vectors its Hermitian norm is
nonpositive. The squared norm of its \(X_0\)-part
is thus at most \(1/N\) at each complex place.
Their difference \(w\) is the last column of \(vY^{-1}\), hence a
rational vector integral in the class lattice. Its Hermitian square
\(h(w,w)\) belongs to \(\mathcal O_{F_0}\) and is at most \(4/N\)
at every real embedding. If \(w\ne0\), definiteness gives
\[
 1\le N_{F_0/\mathbb Q}\bigl(h(w,w)\bigr)\le (4/N)^f<1 .
\]
Thus this column vanishes. The same proof
works on the moved bases and uses only the
cell-rank and basis identities.

Consequently \(vY^{-1}/M\) has support only in the
first \(\min(j,\varrho)\) columns. Multiplication by
\(Y\) writes \(v=M\eta P\) with integral
\(\eta\) obtained from these columns. Equation~\eqref{eq:marked-frame-shear}
therefore moves all the bases in the residue
class to a common \(X_0\)-block.

\smallskip
\subsubsection{Counting isotropic columns}
\label{sec:ordered-basis-isotropic-columns}

\sectionopening{In this subsection, we complete the count using the bound for the definite columns in Section~\ref{sec:ordered-basis-definite-columns}.}

First sort by the matrix of shared vectors,
meaning the \(A\)-coordinates of the shared columns
of \(\gamma F\operatorname{adj}Y\) at positions \(b\ge j\)
and zero at the other positions. These coordinates
are unchanged by moves and integral because the
coordinate functionals for the isotropic block in the class basis are
integral and agree with those in the rational basis.
There are at most \(b^{366 n^4 f}\) possible matrices.
Here is the bound on their entries. At each complex
place the transported chosen basis at a boxed point
has column lengths at most \(b^{37n}\). A Bessel
unipotent makes the representative of that point
of norm at most \(b^{7n}\), by the definition of the
Bessel length and the boxed bound \(b^{6n+2}\).

Its action on the standard basis is unitriangular.
The entries of that matrix are at most \(b^{44n}\),
since it is obtained by multiplying the transported basis by an operator of norm at most \(b^{7n}\).
Cramer's rule and the minor bound give, for any linear
combination in the boxed basis, each coefficient
bounded by \(b^{46n^2}\) times the length of the
combination. The exponent is bounded by
\(1+(2+44n)\varrho+7n\le46n^2\).
Apply this to the shared vector at the first point,
whose expression in the other boxed basis uses
a column of \(\operatorname{adj}Y\).

Its coefficients
have absolute values at most \(b^{(2+134n)\varrho}\)
before converting scales. The ratios of any two
chosen basis scale factors are at most \(b^{23\varrho}\),
because \(1\le|\kappa_v|\le b^{23}\).
Thus each shared coordinate has absolute value at
most
\[
 b^{46n^2+(2+134n)\varrho+1+23\varrho+37n}
 \le b^{182n^2}.
\]
Packing the integral coordinates in \(\varrho n\)
positions as for \(Z\) gives
\((1+2b^{182n^2})^{2f\varrho n}\le b^{366 n^4 f}\).

For one shared matrix and one \(X_0\)-block, sort the
remaining \(\varrho n\) isotropic coordinates modulo
\(M\delta\bar\delta\), costing
\(\operatorname{ind}(M\delta\bar\delta)^{\varrho n}\).
Within this congruence class the two bases differ
only in the isotropic block, by
\(\xi=M\delta\bar\delta\,\mu\) for integral \(\mu\).
Since their Gram matrices agree and the difference
lies in the isotropic space, one has
\(\mu^*P=-P^*\mu\).
Put \(\alpha\) equal to the first \(\varrho\) columns of
\(\operatorname{adj}Y\). The columns at \(b\ge j\)
of \(\mu\operatorname{adj}Y\) vanish by equality of
shared vectors. In the top cell its last column
vanishes because \(P\operatorname{adj}Y\,e_{n-1}=0\),
\(P^*\) is injective, and \(P^*\mu=-\mu^*P\).
Consequently \(\mu\alpha P=\delta\mu\).

An integral skew-Hermitian solution is
\[
 \sigma=\bar\delta\mu\alpha-\delta(\mu\alpha)^*
                     -\alpha^*P^*\mu\alpha .
\]
The displayed skew relation gives \(\sigma^*=-\sigma\);
multiplying by \(P\) cancels the last two terms and
gives \(\sigma P=\delta\bar\delta\mu\).
The central move with parameter \(M\sigma\) therefore
joins the two bases. Separation then identifies
their terms. Neither
this correction nor the previous two moves alter
the level or the residue parameter.

All the moduli and Gram indices are positive and
satisfy
\[
 d(L)^2,\quad\operatorname{ind}(M\delta),\quad
       \operatorname{ind}(M\delta\bar\delta)
       \le b^{2n^3 f}.
\]
The first inequality is the discriminant cutoff.
The other two follow from \(M\le b^{23}\) and the
squared minor bound \(b^{270n^2}\) at all complex
places, using \(n\ge10^5\). The counts of ranks,
product matrices and three successive rounds thus
give the upper bound
\[
 b^{\{2+270n^3+366n^4+
        2n^3(\varrho+ln+\varrho n)\}f}
 \le b^{2000n^6(l+1)f}.
\]

\begin{lemma}[Number of representatives containing a fixed lattice]
\label{lem:representatives-fixed-lattice}
Fix \(c,t\) and a lattice \(L\) of corank \(l\) with
\(d(L)\le\Delta_{\rm lat}\). Count the pairs \((c',\gamma)\) over
\((c,t)\) for which \(\gamma\) has a cell rank. Then
\[
 \#\{(c',\gamma):\gamma\text{ is a chosen near representative},
                      \ L\subseteq L_\gamma\}
     \le b^{2000n^6(l+1)f}.
\]
\end{lemma}

\subsection{The weighted sum over nearby pairs}
\label{sec:marked-pair-count}

\sectionopening{In this subsection, we reduce the weighted count of two orbital terms to a count of
lattices in their common subspace.
Section~\ref{sec:marked-ramified-count} supplies the local weights,
and Section~\ref{sec:marked-pair-gain} bounds the relevant
discriminants.}

For a term \(\gamma\), its link consists of pairs of middle
classes of \(w\) and \(\gamma^{-1}w\), where
\(w\in A^\perp\cap(\gamma A)^\perp\).
For two terms \(\gamma,\gamma'\) over a common block,
let \(V_{12}\subseteq X_0\) consist of the first coordinates
of the common link pairs with both coordinates in \(X_0\).
Intersect \(V_{12}(\mathbb A_f)\cap\Lambda_c^0\) with
the two term lattices to get \(L_{12}\).
Fix \(0\le l\le n^3\). Sum over ordered pairs of chosen near
representatives, each having a cell rank, in the same block
\((c,s,c')\). Retain those for which \(L_{12}\) has rank \(k-l\)
in \(V_{12}\) and \(d(L_{12})\le\Delta_{\rm lat}\).
Section~\ref{sec:marked-pair-gain} shows that these coranks and
discriminants cover all such pairs over the agreeing residue \(s\).

Fix \(c\) and group such pairs by their common parameter
\(t\) and lattice \(L_{12}\). Both representatives
have that parameter and contain the lattice in their
respective term lattices. The assignment to the
ordered pair of terms \(((c',\gamma),(c',\gamma'))\)
is injective, since for those classes the first term
determines the residue. By
Lemma~\ref{lem:representatives-fixed-lattice} there are at most \(b^{4000n^6(l+1)f}\le\mathfrak t\)
ordered pairs for each fixed \(c,t,L_{12}\).
Indeed \(l\le n^3\) and \(k/n^9\to\infty\).
Every pair in that group has weight \(d(L_{12})^{-k}\).
Sum over the lattices, parameters, and left classes and
use \eqref{eq:marked-ramified-weight}. With
\(n_p=|H_R|d_C/N^{kf}\), this gives
\begin{equation}\label{eq:marked-weighted-pairs}
         \sum_{\substack{\text{such ordered pairs}\\\operatorname{corank}L_{12}=l}}d(L_{12})^{-k}
                        \le n_p\mathfrak t^{\,l+2}.
\end{equation}
The factor \(|H_R|\) counts parameters. There is no
additional multiplicity from the residues.

\subsection{Divisibility of the pairing determinant}
\label{sec:marked-meeting-determinant}

\sectionopening{In this subsection, we deduce divisibility of a determinant of pairings with the two final
vectors from agreement of the finite characters. The resulting product
of residue norms enters the discriminant estimate of
Section~\ref{sec:marked-pair-gain}.}

A residue \(s\) agrees if \(\psi_U(u)=\psi_U(u')\)
whenever \(u,u'\) are reduced integral unipotents
with \(us=su'\), using the character of
Section~\ref{sec:marked-residue-projection}.
Consider one place above \(M\), allowing either member
of each conjugate pair. At the residue field at this
place, write \(g\) for the reduction of the original
transporter, and retain \(F\) for the reduced basis.
Suppose \(Fx=gFy\) for coefficient vectors supported
through index \(j<\varrho\).

Let \(m\) be any covector
annihilating both basis initial segments through \(j\).
Then the two matrices
\[
                  1+t(Fx)m,\qquad 1+t(Fy)mg
\]
are reductions of integral unipotents, for every
residue scalar \(t\). To construct the lifts in the
original lattice basis, lift the basis coefficients
and the covectors \(m,mg\) supported at this one place.
The covectors may be chosen to annihilate the corresponding
standard initial segment, since it consists of basis vectors.

The integral Gram identifies the covectors with Hermitian
products against integral vectors. Use the same
Eichler formula as in Section~\ref{sec:marked-class-covolume},
with the lifted isotropic vector replacing \(e_a\) and
level factor one. It applies because the covector kills the
initial segment. The map is upper unitriangular in the Witt order
since the vector is supported through \(j\) and the covector
vanishes there. At the chosen place the term using the
conjugate of the lifted isotropic vector vanishes.

Outside this conjugate pair the matrix is the identity.
The opposite reduction is forced by unitary adjunction.
There the correction kills the basis, because the
conjugate pairing of the isotropic vector with the basis
is zero. Thus it contributes no simple-root entry.
Thus the reduced unipotents are intertwined by \(s\).

The simple-root sums of the two reductions at the
test place are \(t x_j m(f_{j+1})\) and
\(t y_j m(g f_{j+1})\), since the covector kills the earlier
initial segment on each side. The local additive character on
this root is nontrivial, by the conductor of the
different-normalized character and the factor \(M^{-1}\).
Testing every \(t\) in the agreement condition gives
\[
                   x_j m(f_{j+1})=y_j m(g f_{j+1}).
\]
Therefore \(x_j f_{j+1}-y_j g f_{j+1}\) belongs to
the sum of the two earlier spans, by the annihilator
criterion. Moving that linear combination to the
other side gives a common vector at index \(j+1\)
with last coefficient pair \((x_j,y_j)\).
Starting at the highest index of any nonzero
common vector of the full bases and iterating this
argument gives a relation with a nonzero pair of
coefficients in the final positions.

Let \(T\) have two columns in the original integral
left lattice orthogonal to both isotropic spans
\(A,\gamma A\), and apply the last two pairings with
\(z,\gamma z\) to these columns. Suppose the resulting
two-by-two determinant is invertible over the finite
adeles. At each meeting place, the preceding relation
gives a nonzero kernel vector for this matrix after
taking its adjoint and reducing, since all its other
basis pairings vanish. If \(D\) denotes the two-column pairing matrix in the integral
lattice coordinates, then \(v_w(\det D^*)=v_{\bar w}(\det D)\ge1\)
at every meeting place \(w\). Thus the ideal of \(\det D\) is
divisible by \(\prod_{w\ {\rm meeting}}\bar{\mathfrak p}_w\).
The places over \(M\) are split, and conjugation permutes them and
preserves their residue norms. Multiplying the local contributions
gives \(\operatorname{ind}(\det D)\ge Q_{\cap}(s)\), with each meeting
place counted once.
The integral Gram matrix is included when taking
these pairings in a lattice basis.

\subsection{The discriminant bound for the pair lattice}
\label{sec:marked-pair-gain}

\sectionopening{In this subsection, we compare the discriminant of the lattice common to two orbital terms
with the level and intersection determinants. The calculation of
Section~\ref{sec:marked-meeting-determinant} supplies the local
divisibility used in Section~\ref{sec:marked-pair-count}.}

Use the rational spaces \(A,\gamma A\) and the left
lattice \(\Lambda_c\) throughout. Equivalently one may carry
all three through the original left class basis. A term
over \((c,s,c')\) sends \(\Lambda_{c'}\) to \(\Lambda_c\).
Set \(I=A\cap\gamma A\).
The form on \(A^\perp/A\) has middle
space \(X_0\perp E z\). For the two middle coordinates of a
link pair, exchange their \(z\)-coordinates to obtain two
vectors in the positive space \(X_0\perp E z^+\), where
\(z^+\) has norm \(N\).

The difference of the two middle
forms becomes the difference of the two positive forms.
Both middle coordinates come from one lift, so this
difference is zero on the link. Its dimension is \(k+1\),
as follows by taking the orthogonal complement of
\(A+\gamma A\) and quotienting by \(I\). Definiteness makes both
coordinate projections of the link injective. It is therefore
the graph of an isometry \(\Upsilon_\gamma\) of the positive
space. Its positive metric is the positive form of the
first coordinate, or equivalently that of the second.

The calculations below use the cell rank \(j\) in the
original basis. Put \(j_0=\min(j,\varrho)\) and denote the
leading minor of order \(a\) of the basis product matrix by \(\delta_a\).
Modulo \(I\), the first \(j_0\) isotropic basis columns
are a basis of each isotropic basis space. Indeed their
pairing is the leading \(j_0\)-minor below the top rank,
and the common-vector flag relations supply all the other
columns modulo \(I\). At the top rank the full basis product matrix
is invertible. Its first \(\varrho\) rows and first \(\varrho\) columns
are each independent. A vector in \(I\) is orthogonal to both full
bases, so \(I=0\). This argument also covers \(\delta_{\varrho}=0\).
Thus a top cell of rank \(n\) can have a vanishing proper minor.
The estimate for open residues in Section~\ref{sec:open-equal-link-bound}
uses the additional condition of Section~\ref{sec:open-rational-cell}.
Set
\[
 \mu_\gamma=
 \frac{|N_{E/\mathbf Q}(\delta_j)|}{N^{2f(j-\varrho)_+}} .
\]
Let \(t_\gamma\) be the index of
\((\Lambda_c\cap A)+\gamma(\Lambda_{c'}\cap A)\)
in \(\Lambda_c\cap(A+\gamma A)\).

Let \(i_\gamma,i'_\gamma\) be the indices of the pairing
lattices of the first \(j_0\)
isotropic columns inside the pairing lattices of all
isotropic columns. Use respectively \(\Lambda_c\cap A\)
paired with \(\gamma F\), and \(\gamma(\Lambda_{c'}\cap A)\)
paired with \(F\). Define the arithmetic height by
\begin{equation}
 G_\gamma=
 \frac{M^{2fj_0}}{i_\gamma i'_\gamma}
 \max\{N^f\mu_\gamma,\mathfrak g_{\rm norm}(\gamma)\},
 \qquad
 \mathfrak g_{\rm norm}(\gamma)
      =\prod_{v\mid\infty}|\mathfrak g(\gamma)|_v .
 \label{eq:marked-level-gauge}
\end{equation}
Here \(|\cdot|_v\) is the complex modulus at one embedding
in a conjugate pair. Below top rank the discriminant term
vanishes, and at top rank both first indices are one.
For a near term both transported real bases are bounded
by \(b^{37n}\) in the box.

Removing the real basis scaling
gives a basis product bound \(N b^{74n}\le b^{134n}\),
since its other scaling factor is a nonnegative power of
\(|\kappa_v|\ge1\). Expanding a minor and using
\(n\le b\) bounds each squared modulus by \(b^{270n^2}\).
For the discriminant, put \(a=\operatorname{adj}(Z_\gamma)e_{n-1}\).
The vector \(\gamma Fa\) has zero dual-isotropic coordinates and
\(z\)-coordinate \(-\delta_n/N\). Its norm equals that of \(Fa\),
namely \(-N|\delta_{n-1}|^2\) at each complex place.
Hence \(\mathfrak g=\delta_n\bar\delta_n-N^2\delta_{n-1}\bar\delta_{n-1}\)
is \(N\) times the squared norm of its \(X_0\)-component,
and is at most the squared full minor. This holds also
when \(Z_\gamma\) is singular. At top rank with
\(\delta_{\varrho}=0\) and \(\delta_n\ne0\),
\[
 \mathfrak g_{\rm norm}=|N_{E/\mathbf Q}(\delta_n)|
                        =N^{2f}\mu_\gamma>N^f\mu_\gamma .
\]
This term therefore uses the dominant bound of
Section~\ref{sec:marked-dominant-orbits}.
If \(\mathfrak g\ne0\), it is a
totally positive integer of \(F_0\). Hence on the near set
\[
  1\le N_{F_0/\mathbf Q}(\mathfrak g),\qquad
  b^{-270n^2(f-1)}\le\Phi_v(\mathfrak g)\le b^{270n^2}
                   \quad(v\mid\infty).
\]
The lower bound follows by bounding the other \(f-1\) conjugates.
On the open rational cell \(t=\delta_n/\delta_{n-1}\). Dividing
by \(|\Phi_v(\delta_{n-1})|^2\le b^{270n^2}\) gives
\(|\Phi_v(t)|^2-N^2\ge b^{-270n^2f}\).
Equation \eqref{eq:marked-level-gauge},
\(M\le b^{23}\), and \(N\le b^{60n}\) now give
\(0<G_\gamma\le b^{n^3 f}\), using \(n\ge10^5\).

\smallskip
\subsubsection{The lattice associated with one link}
\label{sec:pair-lattice-single-link}

\sectionopening{In this subsection, we describe the lattice for each link of Section~\ref{sec:marked-pair-gain}.}

Use the unimodular quotient lattice
\[
 Q=(\Lambda_c\cap I^\perp)/(\Lambda_c\cap I).
\]
This is unimodular because \(I\cap\Lambda_c\) is a primitive
isotropic direct summand. Its pairing with \(\Lambda_c\)
surjects onto its dual. In a basis consisting of it, a
basis of its orthogonal quotient, and lifts dual to it,
the determinant of the full form is a unit times the
determinant on the quotient. All these bases may be taken
place by place over the valuation rings. Let \(F_1\) be the matrix of the two lists of the
first \(j_0\) isotropic columns in \(Q\),
and \(p=(z,\gamma z)\) read in this quotient.

The classes of the
integral lifts of the link form the primitive orthogonal
lattice to \(F_1\) in \(Q\). In fact a lift is orthogonal
to \(A+\gamma A\), and the list \(F_1\) spans the image of that
space in the quotient, with the deleted kernel \(I\).
Under this lift map, the link metric is the quotient
Gram form \(J\) augmented by
\(N^{-1}J pp^*J\), since the two additional squares are
those of the two values of the \(z\)-form divided by \(N\).

We give the determinant comparison with this augmentation.
Write \(F_0=F_1 a_0\) for a basis of the saturation of \(F_1\)
in \(Q\), \(I_F=[\mathrm{sat}(F_1):F_1]\), and \(v\) for
a basis of the primitive orthogonal lattice. Choose an
integral right inverse \(Q_0\) to \(F_0^*J\), so
\((Q_0,v)\) is integrally invertible. In the ambient space with
Gram \(\widetilde J=\operatorname{diag}(J,N I_2)\),
let \(W=\begin{pmatrix}F_1&p\\0&I_2\end{pmatrix}\) and let
\(v_1=\binom{v}{-N^{-1}p^*Jv}\).
Then \(W^*\widetilde Jv_1=0\), and
\(v_1^*\widetilde Jv_1=v^*(J+N^{-1}Jpp^*J)v\).
A right inverse to \(W^*\widetilde J\) is
\[
 Q_1=\begin{pmatrix}
 Q_0a_0^*&0\\
 -N^{-1}p^*JQ_0a_0^*&N^{-1}I_2
 \end{pmatrix}.
\]
In the basis \((Q_1,v_1)\) the first coordinate block of
\(W\) is \(W^*\widetilde JW\).

Put \(C=(Q_1,v_1)\), \(G_1=W^*\widetilde J W\).
Since \(W^*\widetilde J C=(1,0)\), the matrix
\(C^{-1}(W,v_1)\) is triangular by blocks with
diagonal blocks \(G_1,1\).
Orthogonality gives
\[
 \det(G_1)\det(v_1^*\widetilde Jv_1)
   =\det(\widetilde J)\overline{\det C}\det C\,\det(G_1)^2 .
\]
The bordered determinant computed below is nonzero in each cell.
Cancelling \(\det G_1\), taking indices, and using unimodularity of \(J\) yields
\begin{equation}
  \operatorname{ind}(\det(v_1^*\widetilde Jv_1))
        I_F^2\operatorname{ind}(N)^2
     =\operatorname{ind}(\det(W^*\widetilde JW)).
 \label{eq:marked-quotient-jacobi}
\end{equation}
Indeed \(\det(Q_1,v_1)\) up to sign is
\(\det(Q_0,v)\det a_0^* N^{-2}\), and the first determinant has index one. More precisely, put
\(s_w=\operatorname{length}_{\mathcal O_{E,w}}
 (\operatorname{sat}(F_1)_w/(F_1)_w)\), where \((F_1)_w\) denotes the
column lattice, and \(n_w=v_w(N)\). Then \(v_w(\det a_0)=-s_w\)
and \(v_w(\det a_0^*)=-s_{\bar w}\), so
\[
 v_w(\det(v_1^*\widetilde Jv_1))+s_w+s_{\bar w}+2n_w
                  =v_w(\det G_1).
\]
Taking the product of residue norms accounts for the square \(I_F^2\).
The Gram matrix of \(W\), reordered into the
two isotropic lists with their bordered last vectors, has
zero diagonal blocks and off-diagonal block equal to the
subframe product matrix. On the upper cell this block is
the full matrix.

Below the upper cell a common vector with
last coefficients \((\zeta,1)\), with \(\zeta\bar\zeta=1\),
shows that its determinant is \(-N\zeta\delta_j\).
Reduce the isotropic coefficient vector on the right
modulo \(I\) to its first \(j\) columns. Pairing with the
opposite subframe gives a last column with entries
\((0,\ldots,0,-N\zeta)\) after adding the prescribed
linear combination of the preceding columns. Thus the square root of the
right side of \eqref{eq:marked-quotient-jacobi} is \(N^{2f}\mu_\gamma\) in both cases.
Since \(\operatorname{ind}(N)=N^{2f}\), the link lattice
of integral lifts has discriminant \(\mu_\gamma/I_F\).

Write \(A_L=\Lambda_c\cap A\) and
\(A_R=\gamma(\Lambda_{c'}\cap A)=\Lambda_c\cap\gamma A\).
The second equality uses \(\gamma\Lambda_{c'}=\Lambda_c\), also at
the places dividing \(N\). Let \(A_L^0,A_R^0\) be generated by the
first \(j_0\) columns on the respective sides, and let
\(\pi:\Lambda_c\cap I^\perp\to Q\) be the quotient map. Put
\[
 B_0=\pi(A_L^0+A_R^0),\qquad B_1=\pi(A_L+A_R^0),\qquad
 B_2=\pi(A_L+A_R),\qquad B_3=\operatorname{sat}_Q(B_2).
\]
The first \(j_0\) columns span each side modulo \(I\), and the
two images intersect trivially over \(E\). Thus all four lattices
have rank \(2j_0\) in \((A+\gamma A)/I\), and the successive
indices are finite. Pairing with \(\gamma F\) kills \(A_R\),
including \(\Lambda_c\cap I\). It sends \(B_1/B_0\) onto the
quotient of the pairing lattices of \(A_L\) and \(A_L^0\), of order
\(i_\gamma\). Pairing with \(F\) kills \(A_L\) and likewise sends
\(B_2/B_1\) onto the quotient of order \(i'_\gamma\).
Since \(I\subseteq A+\gamma A\subseteq I^\perp\), the inverse image
of \(B_3\) under \(\pi\) is \(\Lambda_c\cap(A+\gamma A)\).
The kernel \(\Lambda_c\cap I\) lies in \(A_L+A_R\), so
\[
 B_3/B_2\simeq
 \frac{\Lambda_c\cap(A+\gamma A)}{A_L+A_R},\qquad
 [B_3:B_2]=t_\gamma .
\]
Multiplying the three indices gives
\(i_\gamma i'_\gamma t_\gamma\mid [B_3:B_0]=I_F\).

Impose the level condition on the lifts by demanding that
their \(A\)-component be divisible by \(M\), after adding
an integral vector of \(I\). This contributes index at most
\([\widehat{\mathcal O}_E:M\widehat{\mathcal O}_E]^{j_0}
=M^{2fj_0}\). The obstruction is the class of that
component in
\[
 (\Lambda_c\cap A)/((\Lambda_c\cap I)+M(\Lambda_c\cap A)).
\]
Its rank before reducing modulo \(M\) is \(j_0\). Denote the resulting full link
lattice by \(L_\gamma^+\). We have proved
\begin{equation}
 t_\gamma d(L_\gamma^+)
   \le \frac{M^{2fj_0}\mu_\gamma}{i_\gamma i'_\gamma}.
 \label{eq:marked-link-lattice}
\end{equation}

\smallskip
\subsubsection{Decomposition of the two link lattices}
\label{sec:pair-lattice-splitting}

\sectionopening{In this subsection, we compare the two lattices constructed in Section~\ref{sec:pair-lattice-single-link}.}

Let \(\gamma,\gamma'\) be two terms of cell ranks \(j,j'\)
over the same agreeing residue block \((c,s,c')\).
Write \(V_1\) for their link intersection and \(V_0\)
for the kernel in \(V_1\) of the two \(z\)-forms.
Let \(i=k+1-\dim V_1\). Under the first
middle projection \(V_0\) is the space of common
\(X_0\)-vectors. Split the first link orthogonally,
first along \(V_1\), then split \(V_1\) along \(V_0\).
The discriminant of its full lattice \(L_\gamma^+\)
is the product of the three discriminants of its
successive intersection and projection lattices.

Indeed the basis matrix in orthogonal coordinates is
block triangular after lifting bases of the projection
lattices; the shears along the earlier subspaces do
not change its determinant. The first piece projects
to a full sublattice of the part of the term lattice
in \(V_0\). The lift is orthogonal to both isotropic
spaces, its two \(z\)-forms vanish, and its \(A\)-part
is divisible by \(M\). Hence it is orthogonal to the
whole translated basis and satisfies the defining term
lattice condition. Write \(d_\gamma\) for the
discriminant of that part of the term lattice. It is
at most the discriminant of the first piece. The
latter metric is the original Hermitian form on
\(X_0\) because both \(z\)-forms vanish.

The second piece is detected injectively by the
two \(z\)-forms. In the positive model they are
pairings, up to signs, against \(z^+\) and
\(\Upsilon_\gamma^{-1}z^+\). If \(u\) is the last
coordinate of \(\Upsilon_\gamma^{-1}z^+\), their
Gram matrix has diagonal \(N,N\) and off-diagonal
\(Nu,N\bar u\). At top rank the adjugate last
column of the basis product matrix sends, by
\(\Upsilon_\gamma\), a vector with last coordinate
\(\delta_{n-1}\) to \(-\delta_n z^+/N\). This follows
by taking the corresponding combination of the
translated basis vectors, orthogonal to \(A\), with
zero second \(X_0\)-coordinate, and first
\(z\)-coordinate \(-\delta_n/N\).
Therefore \(u=-N\delta_{n-1}/\delta_n\), and the
Gram determinant is \(N^2\mathfrak g/|\delta_n|^2\).
At top rank both \(\delta_n\) and \(\mathfrak g\) are nonzero by the cell
condition. Since \(\mathfrak g\) is \(N\) times a squared norm,
the determinant is positive at every real place. When
\(\delta_{n-1}=0\) it equals \(N^2\). Below top rank the last common vector of the two
flags has norm-one coefficient ratio and zero
\(X_0\)-coordinates, so the two \(z\)-forms are
proportional. Pairing with \(z^+\) is nonzero, so their joint map
on the full link has rank exactly one in these cells.

For the second piece, choose enough forms for a nonzero
square pairing determinant.
The values on a projected lattice vector can be
computed on its integral lift because the shear
lies in \(V_0\). They are integral. If the links
agree and the rank is top, the second piece has
dimension two and Section~\ref{sec:marked-meeting-determinant}
applies to its two integral lifts. Thus the determinant has
index at least \(x\) if \(1\le x\le Q_{\cap}(s)\).
In the other cases take \(x=1\), since the selected
determinant is integral and nonzero.

The absolute norm of the Gram determinant of the
selected form vectors is at most
\[
 \left(\frac{\max\{N^f\mu_\gamma,
               \mathfrak g_{\rm norm}(\gamma)\}}{\mu_\gamma}\right)^2.
\]
In rank zero its value is one. In rank one choose either original form
whose restriction is nonzero. Its ambient Gram entry is \(N\), of absolute
norm \(N^{2f}\). In rank two the displayed Gram determinant gives
\((\mathfrak g_{\rm norm}/\mu_\gamma)^2\).
Cauchy--Schwarz on the exterior powers at every complex
place, followed by multiplication of the finite basis indices, bounds the index of the pairing
determinant by the product of the discriminant of the
projected lattice and the square root of the absolute
norm of the selected form Gram determinant. The second piece therefore has discriminant
at least \(x\mu_\gamma/\max\{N^f\mu_\gamma,
\mathfrak g_{\rm norm}(\gamma)\}\).
Writing \(c_\gamma\) for the discriminant of the third
piece, the complement to \(V_1\), equation \eqref{eq:marked-link-lattice} now gives
\begin{equation}
 t_\gamma d_\gamma x c_\gamma\le G_\gamma .
 \label{eq:marked-split-gauge}
\end{equation}

\smallskip
\subsubsection{The common residue and the intersection lattice}
\label{sec:pair-lattice-intersection}

\sectionopening{In this subsection, we use the decomposition in Section~\ref{sec:pair-lattice-splitting} to prove the discriminant bound.}

For two terms \(\gamma,\gamma'\) with the same block,
the factor \(m_0=\gamma\gamma'^{-1}\) is congruent to
one on \(\Lambda_c\) modulo \(M\), since unitarity
determines the conjugate reductions from the chosen
ones and \(M\) is a unit at all other places.
On integral lifts
\(y,y'\), the difference form between their two
middle pairings is
\[
  \langle y,y'\rangle-\langle y,m_0y'\rangle
    =\langle y,(1-m_0)y'\rangle\in M\widehat{\mathcal O}_E.
\]
This calculation is in the original lattice
basis. In the positive model the same difference form
between source vectors \(a,b\) is
\(\langle a,(1-\Upsilon_\gamma^{-1}\Upsilon_{\gamma'})b\rangle\).

It vanishes along the link intersection. It pairs the
two complementary spaces perfectly. A vector of the
second link orthogonal to the whole first link for the
difference form lies in the first link, whose graph is
maximal isotropic. The integral lift congruence
descends to the projected complementary lattices because the
shears lie in the intersection. On each complex exterior
power of degree \(i\) one has the determinant bound
\[
 \left|\det\left(\langle a,(1-\Upsilon_\gamma^{-1}
                         \Upsilon_{\gamma'})b\rangle/M\right)\right|
 \le (2/M)^i\sqrt{\det\operatorname{Gram}(a)
                                 \det\operatorname{Gram}(b)}.
\]
It follows by Cauchy--Schwarz and \(\|1-U\|_{\rm op}\le2\)
for a unitary \(U\). Choose rational bases \(a,b\) of the two
complements and write their finite lattice bases as \(aA_f,bB_f\).
Let \(C\) be the rational matrix of the difference pairing divided
by \(M\) in the bases \(a,b\). The nondegenerate integral pairing
in the finite bases has matrix \(A_f^*CB_f\). By the product
formula and the definition of the finite index,
\[
 1\le\operatorname{ind}\det(A_f^*CB_f)
   =\operatorname{ind}(\det A_f)\operatorname{ind}(\det B_f)
      \prod_{v\mid\infty}|\det C_v|^2 .
\]
The discriminant of the first complementary lattice is
\(c_\gamma=\operatorname{ind}(\det A_f)
\prod_{v\mid\infty}\det\operatorname{Gram}_v(a)\), and likewise
for \(c_{\gamma'}\). Multiplying the squared exterior-power
inequalities therefore gives
\begin{equation}
 c_\gamma c_{\gamma'}\ge (M/2)^{2if}.
 \label{eq:marked-complement-gain}
\end{equation}
Lastly the common lattice has index at most
\(t_{\gamma'}\) in the part of the first term lattice
in \(V_0\). For \(v\) in this part choose the lift
\(y=v-Ma\) with integral \(a\in A\) and \(y\) orthogonal to
the basis of \(\gamma\). The common middle pair lies in the image of
the rational lift map for the second link. A linear section of this
surjection over \(E\) extends to the finite adeles and gives a lift
\(y'\), over the finite adeles without an integrality condition on \(y'\). Thus
\(y-y'\in A\) and \(\gamma^{-1}y-\gamma'^{-1}y'\in A\).
Thus for \(m_1=\gamma'\gamma^{-1}\),
\[
 (m_1-1)y=(y'-y)+\gamma'(\gamma^{-1}y-\gamma'^{-1}y')
                          \in A+\gamma'A .
\]
The vector \(w=(1-m_1)y/M\) is integral on the left
lattice by the same-residue congruence.

The products of \(v\) with the
basis of \(\gamma'\) equal those of \(M(a+w)\)
because \(m_1y\) is orthogonal to that basis.
Write
\[
 A_0=\Lambda_c\cap A,\quad A_1=\gamma'(\Lambda_{c'}\cap A),
 \qquad W_0=\Lambda_c\cap(A+\gamma'A).
\]
Let \(P'\) send a vector to its pairings against \(\gamma'F\).
Since \(P'\) vanishes on \(A_1\), the finite module
\(W_0/(A_0+A_1)\) surjects onto \(P'(W_0)/P'(A_0)\).
For \(v\) in the first term lattice restricted to \(V_0\),
the preceding calculation sends \(v\) to \(P'(v)/M\)
in this quotient. Indeed the class vanishes precisely when
\(P'(v)=MP'(a')\) for some \(a'\in A_0\). Then \(v-Ma'\) is an
integral lift orthogonal to \(\gamma'F\), which is exactly the
second term-lattice condition. Membership in the first term
lattice and in \(V_0\) was already imposed. Thus the kernel
is \(L_{12}\), and its index divides \(t_{\gamma'}\), by the
definition of that torsion index.
The same argument in the other order gives
\[
 d\le t_{\gamma'}d_\gamma,\qquad
 d\le t_\gamma d_{\gamma'}.
\]

The first pieces of the two link lattices are full in
\(V_0\), so the finite indices just proved make
\(L_{12}\) a full lattice there. It lies in the
integral definite lattice \(\Lambda_c^0\), hence \(d\ge1\).
Put
\[
 l=k-\dim V_0,\qquad
  X=\begin{cases}
      x,&i=0,\ j=n,\\
      1,&\text{otherwise},
     \end{cases}
 \quad 1\le x\le Q_{\cap}(s).
\]
The equal-link case is \(i=0\), since both links have
dimension \(k+1\). On the graph of \(\Upsilon_\gamma\),
the space with both \(z\)-coordinates zero has dimension
\(k-1\) at top rank, and \(k\) below top rank, by the
two-form calculation above. This also shows that equal
links have the same top-rank status, so the value of
\(X\) is symmetric in the terms. If \(r_1\) is the rank of the two
forms restricted to \(V_1\), then
\[
 \dim V_0=k+1-i-r_1,\qquad l=i+r_1-1
                   \le i+(j-\varrho)_+.
\]
In particular, equal links have \(l=1\) in the top cell and
\(l=0\) in the lower cells.

\begin{lemma}[Discriminant of the pair lattice]\label{lem:discriminant-pair-lattice}
For the original terms over one agreeing residue \(s\), with cell ranks
as in Section~\ref{sec:marked-labelled-frames}, let \(d,X,i,l\) be as
above. For \(1\le x\le Q_{\cap}(s)\),
\begin{equation}\label{eq:marked-pair-discriminant-gain}
       (dX)^2\le (2/M)^{2if}G_\gamma G_{\gamma'},\qquad
       l\le i+(j-\varrho)_+.
\end{equation}
\end{lemma}

The first inequality follows by multiplying the two torsion bounds
and applying \eqref{eq:marked-split-gauge} to each term and
\eqref{eq:marked-complement-gain} to the complement.
For near terms both \(G_\gamma\) and \(G_{\gamma'}\) are at most \(\Delta_{\rm lat}\).
Taking \(x=1\) in \eqref{eq:marked-pair-discriminant-gain}
and using \(M/2\ge b^{18}\) gives
\[
                 b^{18if}d\le\Delta_{\rm lat}.
\]
In particular \(d\le\Delta_{\rm lat}\) and
\(l\le i+1\le n^3\). Hence all such near pairs satisfy
the rank and discriminant conditions of
Section~\ref{sec:marked-pair-count}.

\section{Character cancellation and Schur's test}
\label{sec:orbital-cancellation-analysis}

\sectionopening{In this section, we apply character orthogonality to the rational terms counted in
Section~\ref{sec:orbital-pair-counts}. The resulting cell estimates
and Schur bound give the off-diagonal sums used in
Section~\ref{sec:full-orbital-blocks}.}

\subsection{Character cancellation on the blocks}
\label{sec:marked-orbital-cancellation}

\sectionopening{In this subsection, we integrate the two actions of the Bessel unipotent subgroup on the
terms from Section~\ref{sec:marked-full-to-cut}. Character
orthogonality leaves the stabilizer agreement condition analyzed in
Section~\ref{sec:marked-orbital-ranks}.}

Let \(T_0\) be a rational term set invariant under \(U\)
on both sides. On \((T_0\cap\mathcal T(c,s,c'))\times
G_\infty\) the group \(U_c\times U_{c'}\) acts by
\[
 (u,u')\cdot(\gamma,y)=(u\gamma(u')^{-1},uy).
\]
This action preserves counting times Haar measure.
It preserves the full and cut kernels \(K_\gamma^\sharp(y)\)
of Section~\ref{sec:marked-full-to-cut}, because the real
characters on the two class lattices are trivial by the
product formula and the original principal congruence.
The cut uses a left unipotent-invariant length.
One can replace the original weight
\(\beta_c(y)\beta_{c'}(\gamma^{-1}y)\) by any nonnegative
measurable weight with orbit sum one. Insert the new orbit
sum, change variables in each summand, and sum the old
weight. Tonelli and the absolute convergence proved in
Section~\ref{sec:marked-full-to-cut} justify this both for
the norm of the kernel and for the kernel itself. The
same argument applies to the norm of the full-minus-cut kernel.

For a double coset with representative \(\gamma\), put
\(\Gamma_\gamma=U_c\cap\gamma U_{c'}\gamma^{-1}\).
Write its terms uniquely as \(u\gamma(u')^{-1}\), with
\(u\) running over representatives of \(U_c/\Gamma_\gamma\)
and \(u'\in U_{c'}\). After the substitution \(y\mapsto uy\),
the sum of the right cutoffs over \(u'\) is one.
The remaining weight \(\sum_u\beta_c(uy)\) has orbit
sum one for \(\Gamma_\gamma\). The contribution is thus
\(\int_{\Gamma_\gamma\backslash G_\infty}
K_\gamma^\sharp(y)\,dy\).

Say the rational characters agree at \(\gamma\) if
\(\chi_f(u)=\chi_f(u')\) whenever \(u,u'\in U(F_0)\)
and \(u\gamma=\gamma u'\), where \(F_0\) is the
totally real subfield of \(E\) and \(\chi_f\) denotes the
finite Bessel character. On the real unipotent group
\(H_\gamma=U_\infty\cap\gamma U_\infty\gamma^{-1}\), the kernel
transforms on the left by the quotient of the two real
characters. This quotient is trivial on \(\Gamma_\gamma\).
If the rational characters disagree, their rational element
gives a nontrivial value of the real quotient by the product
formula. Integrate in stages over \(H_\gamma\backslash
G_\infty\) with fiber \(\Gamma_\gamma\backslash H_\gamma\)
by the quotient integral formula for the unimodular groups.
Right translation by the real element preserves the fiber
measure and multiplies the integrand by the nontrivial
value. The fiber integral is zero on almost every fiber
by absolute integrability. Therefore every stable set of
rationally disagreeing terms contributes zero.

There is also a finite test on a whole residue block. Let
\(\bar U\) be the reduction of the integral finite unipotent
group, with its induced character \(\psi_U\).
For \(\bar u,\bar u'\in\bar U\), lift them to the integral
finite group. Approximate the lifts in the respective
class bases by rational unipotents \(\alpha,\alpha'\)
with the same reductions. The free coordinates and skew
center give this approximation as in
Section~\ref{sec:marked-residue-projection}.

Each class representative fixes the hyperbolic pairs and \(z\).
Its action on \(X_0\) normalizes the unipotent and preserves the
simple-root sum. Hence \(\Psi_\infty(\alpha)=\psi_U(\bar u)^{-1}\)
and \(\Psi_\infty(\alpha')=\psi_U(\bar u')^{-1}\) by the product formula. The approximants normalize
the class principal unipotent groups since their finite
transports are integral and \(J\) is normal in the
integral stabilizer. Under the change
\((\gamma,y)\mapsto(\alpha\gamma(\alpha')^{-1},\alpha y)\)
the translated weight has orbit sum one in the new block.
By the cutoff identity and the product formula,
\[
 A_{T_0}^\sharp(c,\bar u s(\bar u')^{-1},c')
  =\psi_U(\bar u)^{-1}\psi_U(\bar u')A_{T_0}^\sharp(c,s,c').
\]
In particular an agreeing residue for a nonzero block
must satisfy \(\psi_U(\bar u)=\psi_U(\bar u')\) whenever
\(\bar u s=s\bar u'\). This identity cancels the complex block
before taking its Hilbert--Schmidt norm. On the surviving blocks
the ordered-pair bound \eqref{eq:marked-weighted-pairs} uses the
nonnegative weights \(d(L_{12})^{-k}\) with their ramified
multiplicities from Section~\ref{sec:marked-ramified-count}.

\subsection{The common flag and the cell rank}
\label{sec:marked-orbital-ranks}

\sectionopening{In this subsection, we use agreement of the characters from
Section~\ref{sec:marked-orbital-cancellation} to compare the isotropic
flags. Their intersections determine which ranks contribute to the
open and closed terms in Section~\ref{sec:marked-orbital-reduction}.}

Suppose the rational characters agree at \(\gamma\).
Write \(Fx=\gamma Fy\) for a common vector supported
in both bases through index \(a<\varrho\), with
coefficient pair \((x_a,y_a)\). Let \(\lambda\) be a
linear functional vanishing on both initial segments through
\(a\), and write \(\lambda=h(v,\cdot)\).
For \(A_1=tFx\) consider
\[
 u(w)=w+h(v,w)A_1-h(A_1,w)
                         \bigl(v+\tfrac12 h(v,v)A_1\bigr).
\]
Expansion using \(h(A_1,A_1)=h(v,A_1)=0\) gives
\(h(u(w),u(w'))=h(w,w')\).
In the Witt order the dual coordinates of \(v\)
through index \(a\) vanish, and \(h(A_1,\cdot)\) reads
only these dual columns.

Thus the corrections are
strictly upper triangular in the Witt blocks, with
identity middle block. In particular \(u\in U(F_0)\).
It acts on \(f_b\) by adding \(\lambda(f_b) A_1\), so
its simple-root sum is \(t x_a\lambda(f_{a+1})\).
The conjugate \(\gamma^{-1}u\gamma\) has the same formula
with \(A_1=t\gamma Fy\) written in the second basis,
giving root sum \(t y_a\lambda(\gamma f_{a+1})\).
Put \(c=\lambda(x_af_{a+1}-y_a\gamma f_{a+1})\).
Taking \(t=-M\delta_E q\) for \(q\in E\), agreement of the two
characters gives \(\psi_{E,f}(qc)=1\) by
Section~\ref{sec:fixed-tail-comparison}.
If \(c\ne0\), choose \(q=(2[E:\mathbb Q]c)^{-1}\).
Then \(\operatorname{Tr}_{E/\mathbb Q}(qc)=1/2\), so the additive
product formula gives \(\psi_{E,f}(qc)=-1\). Hence \(c=0\).
The annihilator criterion puts this difference in the
sum of the two initial segments. Rearranging supplies a common
vector with coefficient pair \((x_a,y_a)\) at \(a+1\).

Iteration to the last position and comparison of norms
there give \(x_a\bar x_a=y_a\bar y_a\), since the norm of
a basis vector at last coefficient \(x_a\) is
\(-Nx_a\bar x_a\). Applied at the larger of the two
highest indices of a nonzero common vector, this shows
that both indices agree. A nonzero common vector of
least top index \(j\) can be scaled to give a pair
\((\zeta,1)\), \(\zeta\bar\zeta=1\), which then occurs
at every index at least \(j\).

Write \(\delta_j\) for the leading minor of \(Z_\gamma\)
of order \(j\), with \(\delta_0=1\), and put
\(\mathfrak g_\gamma=\delta_n\bar\delta_n-
N^2\delta_\varrho\bar\delta_\varrho\).
For the last adjugate column \(q\) of \(Z_\gamma\),
the products of \(\gamma Fq\) with the basis vanish
except at the last entry \(\delta_n\). Its last
coefficient in the second basis is \(\delta_\varrho\).
It lies in \(A^\perp\) with \(z\)-coordinate \(-\delta_n/N\).
Comparing its norm in these two expressions gives
\(\mathfrak g_\gamma=N h(w_0,w_0)\) for its \(X_0\)-part
\(w_0\). Thus a nonzero \(\mathfrak g_\gamma\) implies
\(\delta_n\ne0\), and proves rank \(n\). If
\(\mathfrak g_\gamma=0\) and \(\delta_n\ne0\), the same
adjugate column gives a nonzero common vector. If
\(\delta_n=0\), any nonzero \(q\in\ker Z_\gamma\)
gives \(\gamma Fq\in A+X_0\) of nonpositive norm,
hence in \(A\). This also gives a nonzero common vector.

For zero discriminant with rational agreement, take the
least top index \(j\) as above. The shared isotropic
vectors of indices \(j,\ldots,\varrho-1\) span
\(A\cap\gamma A\) by descending elimination of
their highest coefficients. If \(j>0\) and
\(\delta_j=0\), a nonzero vector in the first \(j\)
columns of \(\gamma F\) pairs to zero with the first
\(j\) columns of \(F\). It is perpendicular to the
shared isotropic vectors and hence to all of \(A\).
It also pairs to zero with the final shared basis vector, whose coefficient at \(z\) is nonzero,
so it is perpendicular to \(z\). It is isotropic
in \(A+X_0\) and thus belongs to \(A\). Its top
index is below \(j\) in the second basis and therefore
in both bases, contradicting minimality. This gives
\(\delta_j\ne0\) and the remaining cell-rank conditions
of Section~\ref{sec:marked-labelled-frames}. This cell rank is unique.
A candidate \(j'<j\) would require a shared vector of top index
\(j'\), contrary to minimality. If \(j<j'<n\), the shared vector
of index \(j\) is isotropic and lies in \(A\). Its nonzero
coefficient vector in the first \(j'\) translated columns therefore
lies in the kernel of the leading \(j'\)-by-\(j'\) pairing matrix,
so \(\delta_{j'}=0\). The top rank is excluded by
\(\mathfrak g_\gamma=0\). The case \(j=0\) forces
\(\gamma A=A\) by the shared isotropic basis.

\subsection{Scalar cosets and the surviving terms}
\label{sec:marked-orbital-reduction}

\sectionopening{In this subsection, we classify the terms satisfying the agreement and rank conditions of
Sections~\ref{sec:marked-orbital-cancellation}
and~\ref{sec:marked-orbital-ranks}. This identifies the scalar
contributions and the terms estimated by the orbital bounds.}

If the rational characters agree and \(\gamma A=A\),
the common vector \(f_0\) has top index zero in
both bases. Thus \(\gamma f_0=\zeta f_0\) for
\(\zeta\bar\zeta=1\), and the common pairs above give
\(\gamma f_a-\zeta f_a\in\langle f_b:b<a\rangle_E\)
by successive substitution. Hence \(\zeta^{-1}\gamma\)
is in the rational Bessel subgroup \(R(F_0)=
(U\rtimes H)(F_0)\). Indeed it acts trivially on the
successive isotropic lines and fixes the class of \(z\)
in \(A^\perp/A\). In a term of \(\mathcal T(c,s,c')\)
the first diagonal entry \(\zeta\) of the integral
transporter in the class basis is integral at
every finite place. So is \(\zeta^{-1}=\bar\zeta\).
It follows that \(s\in\operatorname{red}(\zeta)H_R\).

Put
\[
 \Gamma_*=
 \{\gamma:\mathfrak g_\gamma\ne0\}
 \ \cup \{\gamma:\mathfrak g_\gamma=0,\ \gamma A\ne A,                       \text{rational characters agree}\}.
\]
This is stable under both unipotent actions. On the
basis they induce upper unitriangular column moves and
lower unitriangular row operations on the product matrix,
so the leading minors are preserved. They preserve \(A\),
and inner conjugations preserve the unipotent character.
On \(R(F_0)\) the first product column vanishes,
since the first basis vector is fixed and orthogonal
to the basis. Thus \(R(F_0)\) and \(\Gamma_*\) are
disjoint. Outside their union a rationally agreeing
term must preserve \(A\) and hence belongs to a scalar
coset as above. Rationally disagreeing terms outside
this union contribute zero by
Section~\ref{sec:marked-orbital-cancellation}.
For residues outside the nontrivial integral norm-one
scalar cosets we obtain
\begin{equation}\label{eq:marked-off-bessel-reduction}
 A^\sharp_{R(F_0)^c}(c,s,c')
                      =A^\sharp_{\Gamma_*}(c,s,c').
\end{equation}
Every term on the right has a cell rank
\(1\le j\le n\) by Section~\ref{sec:marked-orbital-ranks}.
The nonzero discriminant terms have \(j=n\).

Each global scalar \(\zeta\in\mathcal O_E^\times\) with \(\zeta\bar\zeta=1\) belongs to \(\mu_{q_1}\)
by Section~\ref{sec:marked-level-characters}. If \(\zeta\ne1\)
reduces to one at a place above either level prime \(p\), the identity
\(1+\zeta+\cdots+\zeta^{q_1-1}=0\) gives \(p\mid q_1\), whereas
\(p\equiv1\bmod q_1\). At each place the Bessel subgroup acts
upper unitriangularly on the span of \(F\). Thus an element of
\(\operatorname{red}(\zeta)H_R\) fixing a nonzero vector of this span
forces \(\zeta\) to reduce to one there. A fixed covector gives a
fixed vector of the conjugate span by the inverse-transpose
identification in Section~\ref{sec:full-orbital-blocks}. Consequently
a fixer at even one selected pair excludes every nonidentity scalar coset.
On \(H_RG_i\), put \(p=p_{1-i}^{\rm lev}\). At every place above \(p\)
the reduction belongs to the Bessel subgroup. It fixes \(e_0\), excluding
every nontrivial scalar coset and giving \(Q_{\cap}(s)\ge p^{2f}\)
from the \(f\) split pairs. The divisibility used for equal top links
can also be read directly on this factor. The class matrices
\(B_c,B_{c'}\) preserve the string columns, and the reduced Bessel action
sends \(z\) to \(z\) modulo the span of \(A\). Thus at each \(w\mid p\),
\(\gamma z-z\) belongs to \((A\cap\Lambda_c)_w+p\Lambda_{c,w}\).
For the integral columns orthogonal to \(A\) and \(\gamma A\) in
Section~\ref{sec:marked-meeting-determinant}, the pairing rows with \(z\)
and \(\gamma z\) agree modulo \(p\) at the conjugate place. Their
determinant is divisible by \(p\) at all \(2f\) places, giving index
at least \(p^{2f}\ge(M/b)^f\), since \(p_i^{\rm lev}\le bp\)
by Section~\ref{sec:datum}.
Every scalar coset has vanishing first leading minor
at all chosen places, so the open set avoids them all.
For a fixer at every first-prime place, any fixed
nonzero vector in the span lies in the intersection
of that span with its image. A fixed covector gives the
same intersection at the conjugate place via the Hermitian pairing and inverse-transpose action.
Thus the meeting norm is at least \((p_0^{\rm lev})^f\)
on the full fixer.

\subsection{The real estimate in the open Gram cell}
\label{sec:real-open-cell-contraction}

\sectionopening{In this subsection, we estimate a real orbital integral when the Gram matrix is in the
open cell. The bounds of Section~\ref{sec:real-differentiated-schur}
control the derivatives in the passage from a Gram neighborhood to a
fixed Gram matrix.}

Work at one real embedding in the orthonormal basis
\(E=(e_0,\ldots,e_{\varrho-1},z)\), with \(n=\varrho+1\)
and \(E^*JE=-Q_0\), where \(Q_0\) projects onto the last column.
Use the full Bessel vector \(W\), its quotient mass
\(m_W=\|W\|^2_{U\backslash G}\), and its formal degree \(d_W\)
in \(dg\). Put
\[
 Z=E^*J\gamma E,\qquad
 \Delta(Z)=|\det Z|^2-|\det Z_{<n-1,<n-1}|^2 .
\]
Assume \(\Delta(Z)>0\) and that some \(y_0\) has
both \(y_0^{-1}E\) and \(y_0^{-1}\gamma E\) bounded by
\(\rho=b^{37n}\) in each column.

\smallskip
\subsubsection{Haar density in Gram coordinates}
\label{sec:gram-cell-haar}

\sectionopening{In this subsection, we compute the density for the real estimate of Section~\ref{sec:real-open-cell-contraction}.}

The pointwise basis stabilizer is \(U_0\rtimes H\) with
\(H=U(k)\). Its Haar is \(du\,dh\), with \(dh\) a
probability and \(du\) Lebesgue in the Eichler coordinates.
This Haar is invariant under inversion since \(H\)
preserves the coordinate volume and \(U_0\) is unimodular.

The full Gram fibre has density
\(2^{-r}\operatorname{vol}U(r)\,dg\) in the Hermitian coordinates by Section~\ref{sec:real-frame-haar}.
For a rectangular tuple \(Y\), split \(Y=JEZ+F_0p\),
where \(F_0\) consists of the free isotropic and definite
columns. The blocks \(JE,F_0\) form a Euclidean orthonormal
basis. For invertible \(Z\), write \(Y=Y'Z\).
In \(A\oplus X_0\oplus\mathbb Cz\oplus A^\vee\)
the columns of \(Y'\) have blocks
\[
 Y'=\begin{pmatrix} A_0&a\\ B&x\\0&-1\\I&0 \end{pmatrix}.
\]
For \(T'=(Y')^*JY'\), its first Gram block is
\(A_0+A_0^*+B^*B\), its cross block is \(a+B^*x\),
and its last entry is \(\|x\|^2-1\).
Solving \(A_0=S+T'_{00}/2-B^*B/2\) gives
\(2^{-\varrho}du\) in the skew and free coordinates,
with the remaining factor in \(dT'_{00}\).
The first columns are the Eichler translate of the columns
with \(B=S=0\) at that first Gram block. On the remaining
column, the Eichler action on the isotropic coordinates is
a translation once its definite coordinates are fixed.
Thus solving for \(a\) also gives Jacobian one.
Polar integration in \(x\) gives
\(\tfrac12\operatorname{area}(S^{2k-1})
(1+T'_{nn})^{k-1}\,dh\) times the last real Gram coordinate. The angular
measure is the image of probability Haar on \(H=U(k)\) on the sphere
\(U(k)/U(k-1)\).

Right multiplication by \(Z\) scales the Lebesgue measure in the free rows
by \(|\det Z|^{2(\varrho+k)}\).
For \(T=Z^*T'Z\), Hermitian volume changes by
\(|\det Z|^{2n}\). In fact unitary factors preserve it,
and diagonal singular values \(s_i\) give \(s_i^2\) on
the \(i\)-th real diagonal and \(s_i s_j\) on each of
the two real upper coordinates for the pair \(i<j\).
The net Gram-fibre factor is \(|\det Z|^{2(k-1)}\).
At \(T=-Q_0\), the radius satisfies
\(1+T'_{nn}=\Delta(Z)/|\det Z|^2\).

The resulting basis fibre in the open cell has density
\[
 2^{-(\varrho+1)}\operatorname{area}(S^{2k-1})
                      \Delta(Z)^{k-1}\,dZ\,du\,dh .
\]
Here \(dZ\) is ordinary complex entry measure. The comparison at the
prescribed Gram value \(T_0=-Q_0\) follows on the open chart because
\(Z=E^*JY\) is invertible. Thus \(Y\) has independent columns and
\(Y^*J\) has a right inverse \(C\). For \(D=D^*\),
\[
       d(Y^*JY)[CD/2]=D.
\]
In submersion coordinates the fibre integral of a continuous,
compactly supported test is continuous in the Gram value. Average
over shrinking neighborhoods of \(T_0\), then exhaust the open
chart by compactly supported nonnegative tests to obtain
domination by the basis fibre measure at \(T_0\).
Completing with dual isotropic and definite columns,
as in Section~\ref{sec:real-frame-haar}, multiplies by
\(2^{-\varrho}du'\) and
\(2^{-k}\operatorname{vol}U(k)\,dh'\), respectively.

Comparison with the full Gram density gives
\[
 \int c_G\Delta(Y)^{k-1}v(s\gamma_Y s')\,dY\,ds\,ds'
       \le \int_G v\,dg,\qquad
 c_G=\frac{\operatorname{vol}U(k)}{\operatorname{vol}U(r)}
                        \frac{2\pi^k}{(k-1)!}.
\]
The left variables range over the open cell and both
basis stabilizers, with the section \(\gamma_Y\) given below.
The inequality of measures extends to nonnegative
measurable tests.

\subsubsection{A neighborhood in the original group}
\label{sec:gram-cell-neighborhood}

\sectionopening{In this subsection, we use the density of Section~\ref{sec:gram-cell-haar} to compare neighborhoods in the original measure.}

In the original Eichler
coordinates of \(U_0\), write its parameters as
\((B,C)\) with \(C+C^*=B^*B\). Its group law is
\[
 (B,C)(B',C')=(B+B',C+C'+B^*B').
\]
Let \(V_R\) be the ball \(\|B\|\le R,\ \|C\|\le R^2\).
The group law shows \(V_R V_{R'}\subset V_{R+R'}\).
The Haar measure is Lebesgue measure in the real coordinates
of \(B\) and the skew part of \(C\). The dilation
\((B,C)\mapsto (RB,R^2C)\) has Haar factor
\(R^{2\varrho(k+\varrho)}\).
For \(T\ge8r\), take
\[
   \Omega=V_{r^2T^2},\qquad A_0=\Omega\{y:\|y\|\le T\}.
\]
Both sets are compact and \(\Omega\) has positive volume.
The block formula gives \(\|u\|\le(1+r^2T^2)^2\) on \(\Omega\),
so \(\|y\|\le T^{10}\) on \(A_0\).

On a nonempty left \(U_0\)-slice there is a representative
\(z\) of norm at most \(T\). A point in the slice has
\(u z=s x\), \(s\in\Omega\), \(\|x\|\le T\).
Then \(s^{-1}u=xz^{-1}\) has norm at most \(T^2\), since
\(\|z^{-1}\|=\|z\|\). Its two Eichler blocks have
norm at most \(T^2\), whence it lies in \(V_{T^2}\).
The slice has volume at most the volume of
\(V_{(r^2+1)T^2}\), at most
\((1+r^{-2})^{2\varrho(k+\varrho)}\operatorname{vol}(\Omega)
 \le 3\operatorname{vol}(\Omega)\).
The equality of volumes before and after changing the
representative uses right invariance of \(du_0\).

For the section \(\gamma_Z\) constructed below, put
\(M_{\rm loc}=\rho^{9n}/\Delta(Z)\).
The boxed columns give \(|Z_{ij}|\le\rho^2\), so
\[
 0<\Delta(Z)\le|\det Z|^2\le(n\rho^2)^{2n}\le\rho^{6n}.
\]
Thus \(M_{\rm loc}\ge\rho\ge8r\) and \(M_{\rm loc}^{-1}\le\Delta(Z)\).
In this local calculation only, write \(M=M_{\rm loc}\),
set \(R=M^{9n}\), and take \(T=R^5\). The section is the
identity outside the span of \(E\) and
\(E'=(e_0^*,\ldots,e_{\varrho-1}^*,x_0)\), with \(x_0\)
the definite vector. Its defining basis equations
give the following identity for every \(a<\varrho\):
\[
 e_a^*=\gamma_Z^{-1}E(Z^*)^{-1}\mathbf e_a
             +\tfrac12\sum_{i<\varrho}
                  (Z^{-1})_{i,\mathrm{last}}
                  ((Z^*)^{-1})_{\mathrm{last},a}\,e_i ,
\]
where \(\mathbf e_a\) is the coordinate vector.

It results from the equation taking the section's
\(Z^*\)-dual columns to \(E(Z^*)^{-1}\); the last
positive basis vector occurs only in the last column.
For \(g\) in the support test, the columns of
\(g^{-1}E\) have norm at most \(M\), and
\(\gamma_Zg\) can be moved by \(u\in U\) to \(z'=u\gamma_Zg\)
of norm at most \(M\), by the Bessel length bound.
Indeed approximate the infimum of \(\|s\gamma_Zg\|+1\)
for \(s\in UH\) to within \(1/2\), and remove the
unitary compact factor.

Put \(Z'=E^*J u\gamma_Z E\). Since the action of \(u\)
on \(E\) is upper unitriangular, the leading minors
are unchanged under \(Z\mapsto Z'\), and
\((Z')_{ij}=\langle (z')^{-1}e_i,g^{-1}e_j\rangle_J\)
with the last index interpreted as \(z\). Thus the
entries of \(Z'\) have norm at most \(M^2\);
\(\Delta(Z')=\Delta(Z)\). The adjugate estimate gives
\(\|(Z')^{-1}\|\le R\) and \(\|Z^{-1}\|\le R\).
If \(uE=EU_E\), then \(Z=U_E^*Z'\) and
\(U_E(Z^*)^{-1}=((Z')^*)^{-1}\). Applying \(g^{-1}\)
to the change of basis now replaces the first summand
by \((z')^{-1}E((Z')^*)^{-1}\mathbf e_a\).

It bounds each \(g^{-1}e_a^*\) by
\(nMR+\varrho MR^2/2\le R^4\).
Finally, for \(L=g^{-1}\) use
\(LL^*=J+L(1-J)L^*\). The columns of \(1-J\) are zero
or signed differences of the isotropic and dual
vectors, and \(2z\). Therefore
\(\|L(1-J)\|\le(4\varrho+2)R^4\),
so \(\|L\|^2\le1+(4\varrho+2)R^4\|L\|\).
This gives \(\|g\|+1\le8nR^4\le R^5\).
Thus if \(g^{-1}E\) is boxed and the Bessel length
of \(\gamma_Zg\) is at most \(\rho\), then \(\|g\|\le T\).

On the original double coset with Gram matrix \(Z\), write
\(\gamma=s\gamma_Z s'^{-1}\) with \(s,s'\) fixing the
basis. The ordered pairs \((E,\gamma E)\) and \((E,\gamma_Z E)\)
have the common Gram matrix \(B(1)\), where
\[
 B(t)=\begin{pmatrix}-tQ_0&Z\\ Z^*&-tQ_0\end{pmatrix},
 \qquad
 \det B(t)=(-1)^n|\det Z|^2
             \bigl(1-t^2|(Z^{-1})_{nn}|^2\bigr).
\]
The determinant follows by swapping the column blocks and taking the
Schur complement of \(Z\). The cofactor formula and \(\Delta(Z)>0\)
give \(|(Z^{-1})_{nn}|<1\), so \(B(t)\) is nonsingular on \([0,1]\).
Its signature is \((n,n)\) at \(t=0\) and hence at \(t=1\).
Both orthogonal complements are positive definite of dimension
\(k-1\). Matching orthonormal bases there extends the isometry
between the ordered pairs. For \(v\in\Omega\) and \(h\in H\)
choose \(t=s'(vh)^{-1}\). The support test at the section
applies to \(g=(vh)^{-1}y\) because the first basis is
unchanged and Bessel length is left invariant under
\(s\). If the columns of \(y^{-1}E\) are bounded by \(\rho\)
and \(\gamma t y\) has Bessel length at most \(\rho\),
it gives \(\|v^{-1}y\|=\|g\|\le T\), so \(y\in A_0\).
The support implication holds for each shift and all \(y\).

\subsubsection{The transverse integral}
\label{sec:gram-cell-transverse}

\sectionopening{In this subsection, we apply the neighborhood estimate of Section~\ref{sec:gram-cell-neighborhood} to the section of the Gram map.}

Write \(P_0=\operatorname{diag}(1,\ldots,1,0)\)
and \(Q_0=1-P_0\), and put
\[
 \begin{aligned}
  H_Y&=Y^{-*}Q_0Y^{-1},&
  \eta_Y&=(1-|(Y^{-1})_{\rm last,last}|^2)^{1/2},\\
  A_Y&=-Q_0-\tfrac12P_0H_Y(1+Q_0),&
  L_Y&=EA_Y+E'(P_0+\eta_Y Q_0).
 \end{aligned}
\]
The block Gram matrix on \((E,E')\) is
\(\left(\begin{smallmatrix}-Q_0&P_0\\ P_0&Q_0\end{smallmatrix}\right)\).
In particular
\[
  E^*J L_Y=1,\qquad
  L_Y^*J L_Y=-Q_0+A_Y^*P_0+P_0 A_Y+\eta_Y^2Q_0=-H_Y,
\]
where \(Q_0 H_YQ_0=(1-\eta_Y^2)Q_0\).
The isometry sending \((E,L_{Y^*})\) to
\((L_Y Y,EY^{-*})\) exists because the displayed identities
give identical Gram matrices, and both ordered pairs are
bases of the span of \((E,E')\). Extend it by the identity
on the definite orthogonal complement to get \(\gamma_Y\).
In particular it has the original basis products \(Y\).
If \(S_Y=L_Y Y\), its column formula is
\[
 \gamma_Y=1-EE^*-E'E'^*+S_Y E^*
     +(EY^{-*}-S_Y A_{Y^*})(P_0+\eta_Y^{-1}Q_0)E'^* .
\]
These equations give the dual-column identity used above.

Suppose \(P\ge8r\) bounds \(\|Y\|,\|Y^{-1}\|,\eta_Y^{-1}\).
For a real tangent of unit norm,
\[
   \|\dot Y^{-1}\|\le P^2,\quad
   |\dot\eta_Y|\le P^3,\quad
   |(\eta_Y^{-1})^{\boldsymbol\cdot}|\le P^5 .
\]
The product rule gives \(\|A_Y\|\le2P^2\),
\(\|\dot A_Y\|\le2P^3\),
\(\|S_Y\|\le3P^3\), and
\(\|\dot S_Y\|\le4P^4\).
Thus \(\|\gamma_Y\|\le P^{18}\) and
\(\|\dot\gamma_Y\|\le P^{25}\).
The stable section neighbourhood can be taken as
\(\|Y-Z\|<P^{-6}\), with
\(M\ge8r\), \(M^{-1}\le\Delta(Z)\), \(|Z_{ij}|\le M^2\),
\(R=M^{9n}\), and \(P=(R^5)^n\). In fact the adjugate
and determinant estimates give \(\|Z\|,\|Z^{-1}\|,\eta_Z^{-1}\le
M^{8n}\le P/2\). A perturbation of norm \(P^{-6}\) keeps
the inverse below \(P\) and changes it by at most \(P^{-4}/2\).
The absolute change in \(1-|(Y^{-1})_{\rm last,last}|^2\)
is at most \(2P^{-4}\), whereas its initial value is at least
\(4/P^2\). Thus \(\eta_Y^{-1}\le P\).
The determinant ratio of \(Z=(1+(Z-Y)Y^{-1})Y\)
is bounded above by \((1+P^{-5})^n\le9/8\).
The relative loss in \(\eta_Z^2\) is at most
\(2P^{-4}/(4P^{-2})\le 1/4\).
Since \(\Delta(Y)=|\det Y|^2\eta_Y^2\), these bounds give
\(\Delta(Y)\ge (3/4)(8/9)^2\Delta(Z)>\Delta(Z)/2\).

Let \(c=c_G(\Delta(Z)/2)^{k-1}\).
For a nonnegative function \(v\) invariant under the full left \(UH\),
a measurable cutoff of integral one on each left
\(U_0\)-orbit exists by local sections of \(U_0\backslash G\).
In a section chart take a compactly supported continuous
density on \(U_0\) of mass one and multiply it by a partition
of unity on the quotient. Restrict the cutoff to the inverse
image of the Gram ball.

Its \(U_0\)-integral there remains one,
since a basis stabilizer keeps the Gram matrix fixed.
The lower Haar density above and the probability Haar of \(H\)
then give
\[
 c\int_{\|Y-Z\|<P^{-6}}\int_{H\times U_0}
              v(\gamma_Y(uh)^{-1})\,dh\,du\,dY
 \le \int_G \chi(y)v(y)\,dg(y).
\]
In the right integral unfold first along \(U_0\), then
along \(U/U_0\), which acts on \(E\) by unitriangular
matrices. If two points in a \(U\)-orbit lie over the Gram
ball, the adjoint of the transition matrix on \(E\) is the
product of the first Gram matrix with the inverse of the
second. All its entries are at most \(P^2\) in modulus.
Using the original Haar on the upper entries of the
basis matrix, the last unfolding contributes at most
\((2P^2)^{2n^2}\). This proves
\[
 c\int_{\|Y-Z\|<P^{-6}}\int_{H\times U_0}
              v(\gamma_Y(uh)^{-1})\,dh\,du\,dY
 \le (2P^2)^{2n^2}\int_{U\backslash G}v.
\]

For matrix-entry words \(j\) of length at most \(2n^2\),
Section~\ref{sec:real-differentiated-schur} gives
\[
 \int_{U\backslash G}\left\|\int_G\beta(y)W_j(g y)\,dg(y)
                      \right\|^2
 \le P^{2|j|}\frac{m_W}{d_W}\int_G|\beta|^2\,dg ,
\]
since \(b^5\le P\). The derivatives have the same
unitary left Bessel covariance as \(W\).

For a test supported in \(A_0\), the logarithmic tangent
for differentiation at \(gy\) is \(y^{-1}Xy\) when
differentiating at \(g\) along \(gX\).
Since \(\|y\|\le P^{20}\), each matrix-entry coefficient
contributed by this conjugation, divided by \(P^{40}\),
has modulus at most one. More precisely attach the fixed
weight \((y^{-1})_{ac}y_{db}/P^{40}\) to the differentiation generator \(E_{ab}\) and the external coefficient \(X_{cd}\).
After expansion each integral is a convolution with a
weight bounded by \(|\beta|\) and one of the \(W_j\). Differentiation under the
integral uses continuity and boundedness of these
derivatives on the compact support, locally in the external
variable. Each of the \(r^4\) possible new generators contributes at most \(P^{40}\|X\|\).
On the cube of side \(P^{-8}\) based at \(Z\) in
\(2n^2\) entry coordinates, take \(g=\gamma_Y(uh)^{-1}\)
with \(u\in\Omega\), so that \(\|u^{\pm1}\|\le P^{20}\).
For each real entry coordinate of \(Y\),
\[
 X=g^{-1}\dot g
   =(uh)\gamma_Y^{-1}\dot\gamma_Y(uh)^{-1},\qquad
 \|X\|\le P^{20+18+25+20}\le P^{85},
\]
since \(h\) is Euclidean unitary and
\(\|\gamma_Y^{-1}\|=\|\gamma_Y\|\).
Therefore each first derivative has a sum of at most
\(r^4\) terms with coefficients at most \(P^{125}\).

For a Hilbert-valued function on \([0,\epsilon]\),
\[
 \|a(0)\|^2\le 2\epsilon^{-1}\int_0^\epsilon\|a\|^2
                   +2\epsilon\int_0^\epsilon\|a'\|^2 .
\]
Apply it successively in all \(2n^2\) coordinates.
To iterate, denote a convolution obtained so far by \(H_\omega(Y)\),
where its test \(\beta_\omega\) is independent of \(Y\) and
\(|\beta_\omega|\le|\beta|\). In each remaining coordinate the
derivative formula and Cauchy--Schwarz give
\[
 \|\partial H_\omega(Y)\|^2
   \le r^4P^{250}\sum_{\omega'}\|H_{\omega'}(Y)\|^2 .
\]
There are at most \(r^4\) terms, each with one additional generator
and a new fixed test bounded by \(|\beta|\). Substitution into
the one-dimensional inequality with \(\epsilon=P^{-8}\) leaves
only squared convolutions under the integral. The scalar coefficients
are now independent of \(Y\), so the next coordinate requires no
derivative of the section coefficients already estimated. The
total multiplier per coordinate is at most
\(2(1+P^{258})P^8\), since \(r^8\le P^8\). Apply the Gram-ball inequality
to each word, whose squared norm is left Bessel-invariant.
The combined contribution per coordinate is at most \(P^{275}\).
The quotient Schur bound contributes at most another \(P^{4n^2}\).
This proves the transverse estimate
\[
 c\int_{H\times\Omega}
  \left\|\int_{A_0}\beta(y)W(\gamma_Z(uh)^{-1}y)\,dg(y)
                                                    \right\|^2 dh\,du
   \le P^{2000n^2}\frac{m_W}{d_W}\int_{A_0}|\beta|^2\,dg .
\]
The cube lies in the Gram ball since
\(2n^2P^{-8}<P^{-6}\).

For \(a=(h,u)\), put \(g_a=\gamma_Z(uh)^{-1}\) and use
the shift \(t_a\) with \(\gamma t_a=s g_a\) constructed above.
The Bessel covariance is unitary on the value space, so
\[
 \begin{aligned}
 \int_{A_0}\phi(y)W(\gamma t_a y)W(t_a y)^*\,dg(y)
  = U(s)\left(\int_{A_0}
       \phi(y)W(g_a y)W(y)^*\,dg(y)\right)U(t_a)^* .
 \end{aligned}
\]
Its squared Hilbert--Schmidt norm is the sum of the squared
norms of the inner columns. Writing \(W_j\) here for a coordinate of \(W\), apply the scalar
estimate to \(\beta_j=1_{A_0}\phi\overline{W_j}\) and sum,
using \(\sum_j|\beta_j|^2=1_{A_0}|\phi|^2\|W\|^2\).
Dividing by \(c\operatorname{vol}(\Omega)\) gives the
contraction with normalized measure \(dh\,du/\operatorname{vol}(\Omega)\)
on \(H\times\Omega\). Thus the operator in the parameter \(a\) is a contraction
from \(L^2\) of the measure on the source to \(L^2\) of the normalized measure with Hilbert--Schmidt values. The source measure
has density
\[
 \frac{P^{2000n^2}m_W}{d_W\,c\operatorname{vol}(\Omega)}
             1_{A_0}(y)\|W(y)\|^2\,dg(y).
\]

The same slice volume cancels on unfolding
\(\int_{A_0}\|W\|^2\,dg\), since the slice contribution is at most
\(3\operatorname{vol}(\Omega)\), the upper-entry contribution is at
most \((2P^{40})^{2n^2}\le(MP^{40})^{2n^2}\), since two
points in the tube have norms at most \(P^{20}\) and
the transition matrix on their basis has entries at
most \(P^{40}\), and
\(\|W\|^2_{U\backslash G}=m_W\). Thus the total contribution is
\[
 M^{10^5n^4}\frac{m_W^2}{d_W\,c_G(\Delta(Z)/2)^{k-1}},
\]
using \(3\le R^5\) and
\(P^{2000n^2}R^5(MP^{40})^{2n^2}\le M^{10^5n^4}\).
This bounds the mass of the measure on the source of the
one-place contraction.

\subsection{The contribution of large discriminants}
\label{sec:marked-dominant-orbits}

\sectionopening{In this subsection, we apply the real estimate of
Section~\ref{sec:real-open-cell-contraction} when the discriminant
factor dominates. Its decay can then be combined with the lattice
comparison in Section~\ref{sec:marked-pair-gain}.}

For an agreeing block \((c,s,c')\), denote by
\(\mathcal R_{S_0}(c,s,c')\) the chosen near representatives in a
stable set \(S_0\) having a cell rank. For \(\gamma\) of rank \(j\) put
\[
 D_j=\begin{cases}D_{\rm main},&j=n,\\
       \mathfrak t^{-16}D_{\rm id},&j<n,\end{cases}
 \qquad w_\gamma=\mathfrak t^2\sqrt{\frac{D_j}{n_p G_\gamma^k}},
\]
with \(n_p=|H_R|d_C/N^{kf}\) as in \eqref{eq:marked-weighted-pairs}
and \(G_\gamma\) as in \eqref{eq:marked-level-gauge}.

For a stable subset of the surviving terms, take its part
on which \(\mathfrak g_{\rm norm}(\gamma)>N^f\mu_\gamma\).
Both quantities are stable under the rational unipotent
moves, by the leading-minor and cell-rank invariance.
On this part the cell rank is \(n\), both first indices are
one, and \eqref{eq:marked-level-gauge} gives
\(G_\gamma=M^{2\varrho f}\mathfrak g_{\rm norm}(\gamma)\).
For a near representative write \(x_\gamma=\mathfrak g_{\rm norm}(\gamma)\).
The basis products are integral after the transport,
so \(x_\gamma\ge1\) by the norm formula. The box bound in
Section~\ref{sec:marked-pair-gain} gives \(x_\gamma\le b^{270n^2 f}\).
At \(v\) the normalized discriminant is
\[
 \Delta_v=\left(\prod_a a_{v,a}\right)^2
                         |\mathfrak g(\gamma)|_v .
\]
Indeed the full minor scales by the product of all the
squared basis scales, whereas the preceding minor omits
the last one, \(1/N\).
Section~\ref{sec:real-open-cell-contraction} gives the
bound for the local source measure
\[
 P_v=\left(\frac{\rho^{9n}}{\Delta_v}\right)^{10^5n^4}
                 \frac{d_v^2}{\deg_v\,c_G(\Delta_v/2)^{k-1}},
\]
Here \(\deg_v\) is the formal degree of the ambient discrete series
in the chosen real Haar measure. The normalization of
Section~\ref{sec:real-marked-filling} gives
\[
 d_v=\|W_v\|^2_{U_v\backslash G_v}
    =\dim\operatorname{Sym}^{a_\tau}(\mathbb C^k)
    =\binom{k-1+a_\tau}{k-1},\qquad
 d_B=\prod_v d_v,\qquad \mathcal B=\prod_v\deg_v .
\]

\smallskip
\indent\textit{The product contraction.}
Write \(K_v^\epsilon(\lambda,z)=
W_v^\epsilon(\lambda z)\otimes W_v^\epsilon(z)^*\),
with \(\epsilon\in\{\mathrm{full},\mathrm{cut}\}\).
For each \(v\), use the compact set \(A_v\), the probability
\(m_v\), the finite measure \(\nu_v\) and the measurable
basis-fixer shift \(t_v(a_v)\) in the local restriction estimate.
The linear map is
\[
 (T_v\phi)(a_v)=\int_{A_v}\phi(z_v)
       K_v^{\rm full}(\lambda_v,t_v(a_v)z_v)\,dg_v(z_v)
\]
and its norm from \(L^2(\nu_v)\) to
\(L^2(m_v;\operatorname{HS}_v)\) is at most one.
For a bounded measurable \(\phi\) on the product the tensor
operator still has the displayed integral formula with
\(A=\prod_v A_v\), product Haar and the product kernel.

Indeed it has this formula for finite sums of rectangle
indicators by Fubini. The contraction gives the extension
on the Hilbert tensor product, where the squared norm is
computed in the tensor Hilbert--Schmidt basis. To identify that
extension with the integral for general bounded measurable
\(\phi\), the bounded functions for which the formula holds
form a vector space closed under uniformly bounded monotone
limits of real functions. Convergence holds in the source norm by
finite mass and pointwise in the ordinary integral by
dominated convergence, since for a fixed shift every kernel
is continuous on the compact domain. The contraction then
identifies the two limits almost everywhere. The monotone
class theorem applies to the product Borel sets, and the
imaginary part is treated by linearity.

For \(|\omega|\le1\) invariant under the product of basis stabilizers, apply this contraction to \(\omega\). The mean
squared Hilbert--Schmidt norm is at most
\(\|\omega\|^2_{L^2(\prod\nu_v)}\le\prod P_v\).
The finitely many support tests hold on one full
\(\prod m_v\)-measure set, independently of the point tested.
There is a shift \(t\) in that set with squared norm at most
\(\prod P_v\). In fact a nonnegative measurable function
strictly larger than that number throughout the set would
have integral strictly larger than it. Fixer invariance
identifies the value \(\omega(z)\) with \(\omega(tz)\).

Use \(\lambda=\gamma^{-1}\), with the chosen element
\(\lambda^{-1}y_0'=y_0\), \(y_0'=\gamma^{-1}y_0\),
and conjugate the cutoff when taking the adjoint. For the
actual weight \(\omega=1_{\mathrm{red}}\) the conjugation
is trivial. Put \(y=tz\). By the tensor formula for
\(K_v^{\rm cut}\), the translated weighted integrand is nonzero exactly when
\[
 \omega(y)\ne0,\qquad W_v^{\rm cut}(y_v)\ne0,\qquad
 W_v^{\rm cut}(\gamma_v^{-1}y_v)\ne0
              \quad\text{at every real place }v.
\]
The first two conditions put \(y\) in the reduced set and in the cut,
so triangular reduction bounds the columns of \(y_v^{-1}E_v\).
These columns equal those of \(z_v^{-1}E_v\) since \(t_v\) fixes \(E_v\).
The third condition bounds the Bessel length of
\(\gamma_v^{-1}t_vz_v\) by \(b^{6n+2}\le\rho\). This is exactly the
support test that puts \(z\) in \(A\). The support of the
actual weighted cut is thus contained in \(tA\), although
the indicator alone need have no compact support. All
kernels are bounded on \(tA\), and left Haar gives
\[
 \left\|\int\omega(y)\,K^{\rm cut}_\gamma(y)\,dg(y)\right\|_{\rm HS}
 \le (\prod_v P_v)^{1/2}
   +\int|\omega(y)|
        \|K^{\rm full}_\gamma(y)
                 -K^{\rm cut}_\gamma(y)\|_{\rm HS}\,dg(y).
\]
The error is a nonnegative integral.

Its finiteness for
the selected orbit representatives follows from the
cutoff transfer below. The full kernel is integrated on \(tA\). The local discriminant for
\(\gamma^{-1}\) equals that for \(\gamma\), because their
basis product matrices are adjoints. This uses the product Haar and the same chosen bases.

Write \(\Lambda_v=\prod_a a_{v,a}^{k+2(\varrho-a)}\),
so \(\Lambda_0=\prod_v\Lambda_v\). Put \(q=10^5n^4\). Then
\[
  \prod_vP_v(\gamma)=C_{\rm pc}/x_\gamma^{q+k-1},\qquad
  C_{\rm pc}\le 2^{kf}\mathfrak t\,
       \frac{d_B^2}{\mathcal B c_G^f\Lambda_0^2}.
\]
Indeed \(b^{-60n^2}\le\prod_aa_{v,a}\le b^{46n^2}\),
and the local factor
\(\prod_a a_{v,a}^{k+2(\varrho-a)}\) is at most
\((\prod_aa_{v,a})^kb^{92n^3}\).
Against \(2^kd_v^2/(\deg_v c_G\Lambda_v^2)\), the ratio is at most
\(b^{453n^2q+92n^2+184n^3}\) at each place, whose
product is at most \(\mathfrak t\) by \(k\ge512000n^9\).
In fact putting \(S_v=\prod_a a_{v,a}\), the per-place
ratio, with the discriminant norm removed, is
\(\tfrac12(\rho^{9n}/S_v^2)^q S_v^{-2(k-1)}\Lambda_v^2\).

Put \(v_E=\operatorname{vol}(\widehat{\mathcal O}_E)\) in
the additive Haar at finite places used for the Siegel
calculation, and let \(C_0\) denote the scaled Siegel constant
of the prescribed theta datum. The two numerical comparisons are
\[
 2^{kf}n_pM^{2\varrho fk}\le\mathfrak t C_0v_E^n,\qquad
 C_0\le(1+b^{-9})D_{\rm main}\mathcal B c_G^f v_E^{-n}
                               (\nu\Lambda_0/d_B)^2 .
\]
For the first comparison set \(X=200fk\operatorname{rd}(E)\).
The scale bounds above give \(X^{1000}\le b\). For the one-vector
Siegel constant of Section~\ref{sec:fixed-tail-average}, the bounds
\(v_E\ge X^{-f}\), \((k-1)!\le X^k\), and \(n\le k\) give
\[
 C_{\rm S}v_E^n
 =|D_E|^{1/4}v_E^{k+n}
       \left(\frac{\pi^k}{(k-1)!\,N^k}\right)^f
 \ge N^{-kf}X^{-(2k+n)f}\ge N^{-kf}X^{-3kf}.
\]

Use \(\mathfrak t\ge X^{4kf}\),
\(X\ge2\), and the unipotent order
\(m_U\ge M^{2\varrho fk}\) in \(C_0=|H_R|d_Cm_U C_{\rm S}\), and use
\(n_p=|H_R|d_C/N^{kf}\).
The second comparison is
\(C_0\Delta_0\le(1+b^{-9})\nu^2D_{\rm main}\mathcal B\)
from Section~\ref{sec:open-equal-link-bound}.
Combining gives
\[
 C_{\rm pc} n_pM^{2\varrho fk}\le2\mathfrak t^2\nu^2 D_{\rm main}
                              \le\tfrac14\mathfrak t^4\nu^2D_{\rm main}.
\]
Since \(q+k-1\ge k\) and \(x_\gamma\ge1\), this bounds the local
product square by one fourth of the claimed orbit square.
The absolute cut error will use the other half of the orbit
bound. Its admissibility follows from
\[
 b^{-k^3f}\le
 \mathfrak t^2\sqrt{D_{\rm main}/(n_p(M^{2\varrho f}x_\gamma)^k)},
\]
using \(n_p\le D_{\rm id}\le D_{\rm main}\),
\(x_\gamma\le b^{270n^2f}\), \(M\le b^{23}\), and the rank
inequality. This proves the normalization of the two halves.

Use the reduced-point indicator \(\omega\) from the box partition. It is invariant under the real fixer of the
basis, since the transported columns are unchanged, and
on the real cut its support is in the box of radius \(\rho\). For each orbit of
\(U_c^0\backslash(S_0\cap\mathcal T(c,s,c'))/U_{c'}\)
choose one representative, preferring a near representative
of a smaller orbit with right group \(U_{c'}^0\) from the
fixed list of Section~\ref{sec:marked-labelled-frames}.
This choice suffices because whenever the reduced indicator
and cut kernel are simultaneously nonzero, one can bring the
second cut point to its box by its original class unipotent
lattice. The new term remains in the larger orbit.
Within the smaller orbit boxing is invariant under both
basis-fixing unipotent lattices, so its chosen representative
is near.

For a nonzero-discriminant term the stabilizer of
\(\gamma U_{c'}\) in \(U_c^0\) is trivial.
Indeed if \(a\gamma=\gamma v\), then \(v\in U_{c'}\) fixes
\(\gamma^{-1}F\) pointwise. These columns pair perfectly
with \(F\) since the full basis product matrix is invertible.
Linear combinations therefore give \(\varrho\) fixed
vectors whose \(e_a^*\)-coordinates form the identity.
In Witt blocks such an element has diagonal blocks
\((B,I,B^{-*})\). Fixing the independent dual coordinates
gives \(B=I\); the middle rows of the fixed-vector equations
kill the middle-to-dual block. The unitary relation kills the
isotropic-to-middle block, and the remaining fixed-vector
equations kill the strictly upper corner block.
Thus \(v=a=1\).

Give \((\gamma,y)\) weight \(\omega(y)\) when \(\gamma\) is the
chosen representative of its larger orbit, and weight zero otherwise. This is an orbit
cutoff on the set of terms times the real group. To check this,
for each right coset of \(U_c^0\) in \(U_c\) exactly
one element of that coset, together with one element of
\(U_{c'}\), takes the term to the representative of the
larger orbit of \(u\gamma\), for \(u\) in that coset.

Existence is the orbit decomposition and uniqueness is the
trivial stabilizer proved above.
The sum of \(\omega\) over the remaining cosets is one
by the box partition. Apply the change of cutoff in
Section~\ref{sec:marked-orbital-cancellation}.
It transfers the full-minus-cut absolute bound as well,
because the norm of the difference of kernels is invariant
under the same action. In particular the absolute tail on
each chosen term is at most \(\nu b^{-k^3f}/2\).

Any chosen term that contributes is near by the boxing
argument. Its integral is bounded by the sum of the two
halves calculated above. A larger orbit charges at most
one smaller-orbit representative, so summing gives
the bound
\[
 \|A^{\rm cut}_{S_0}(c,s,c')\|_{\rm HS}
       \le \nu\sum_{\gamma\in\mathcal R_{S_0}(c,s,c')}w_\gamma
\]
for \(S_0\) in the dominant-discriminant part. The sum uses
the original fixed near representatives. On a nonagreeing
residue the block is zero by Section~\ref{sec:marked-orbital-cancellation}.

\subsection{The one-place Schur integral}
\label{sec:real-one-place-schur}

\sectionopening{In this subsection, we give the real integral bound for a contributing cell at a single
place. The rank condition comes from
Section~\ref{sec:marked-orbital-ranks}, and the estimate will be
multiplied across places in Section~\ref{sec:marked-schur-orbits}.}

Use a surviving cell of rank \(1\le j\le n\), let
\(i=\min(j,\varrho)\), and put \(I=A\cap\gamma A\).
Put \(\rho=b^{37n}\) and \(\theta=b^{6n+2}\). Throughout, \(y\) is a real
point such that the columns of \(y^{-1}F\) have length at most \(\rho\)
and the Bessel length of \(y\) is at most \(\theta\).
At a real place use the orthonormal basis with last
vector of norm \(-1\). Let \(F_i\) consist of the first
\(i\) isotropic columns and the last vector. In a lower
cell the shared columns starting at \(j\) have leading
pairs \((\zeta,1)\) with \(|\zeta|=1\) in the orthonormal basis,
by Section~\ref{sec:marked-orbital-ranks}.

These pairs
retain their coefficients because both bases have the same
column scales. The shared isotropic columns span \(I\).
From the shared last vector remove the columns of \(I\)
in the second basis, to obtain a vector in
\(\operatorname{span}(F)\cap\operatorname{span}(\gamma F_i)\)
with coefficient one at \(\gamma z\) and coefficient of
modulus one at \(z\). Let \(B_1=F_i^*h\gamma F_i\),
\(\delta=\det B_1\). Replacing its last column by the
pairings of that vector preserves the determinant. Those
pairings are zero except for the last, of modulus one.
Thus in a lower cell \(|\delta|\) is the modulus of the
leading \(j\)-minor in the orthonormal basis. In the top cell
this is also true, directly since \(F_i=F\).

Let \(V_-\) be the negative eigenspace of the Hermitian involution \(h\) in the orthonormal basis. The negative
space \(Y_-=yV_-\) pairs perfectly with the entire
basis: a vector in \(Y_-\) orthogonal to the basis must be
zero since the basis orthogonal complement is nonnegative.
Pull the columns of \(F^*h\gamma F_i B_1^{-1}\) back through this
pairing to \(Y_-\), and apply Gram--Schmidt for \(-h\) in reverse
column order. These columns pair to zero with \(I\), which is
orthogonal to \(\gamma F_i\). This gives a negative orthonormal
basis \(\omega\) of \(Y_-\cap I^\perp\) and a lower triangular
invertible \(L\) with
\[
 F^*h\omega=F^*h\gamma F_i B_1^{-1}L,\qquad F_i^*h\omega=L.
\]
The norm of \(y^{-1}\omega_a\) is one, so
\(|L_{aa}|\le\|y^{-1}(F_i)_a\|\), giving
\(|\det L|\le\rho^{i+1}\) on the box of radius \(\rho\).

Let \(E_-\) have orthonormal columns spanning \(V_-\). Under \(y^{-1}\),
the \(h\)-orthogonal projection onto \(Y_-=yV_-\) becomes Euclidean
projection onto \(V_-\), with its \(-h\) metric. Hence the ratio of
the projected volumes is
\[
 \mathcal V(\gamma,y)=
 \frac{|\det(E_-^*y^{-1}\gamma F)|}{|\det(E_-^*y^{-1}F)|},
 \qquad
 \mathcal V(\gamma^{-1},\gamma^{-1}y)=\mathcal V(\gamma,y)^{-1}.
\]
Both projected bases are independent since their spans are nonpositive
and the kernel of the projection is positive definite.
In a lower cell replace the omitted isotropic columns in
each basis by their shared columns. The changes of basis
have determinant of modulus one by the common leading pairs.
The projected shared wedge volume cancels, leaving volumes
in the orthogonal complement of the projected \(I\) in
\(Y_-\).

This complement is \(Y_-\cap I^\perp\), since
pairing against a vector in \(Y_-\) only sees the negative
projection. Pairing the retained columns with its
orthonormal basis \(\omega\) yields
\[
 \mathcal V(\gamma,y)=
       \frac{|\det((\gamma F_i)^*h\omega)|}{|\det L|}.
\]
Put \(\omega'=\gamma F_i B_1^{-1}L\).
The columns of \(\omega'-\omega\) lie in \(A+X_0\); let
\(V_0\) be the matrix of their \(X_0\)-coordinates, so their
Gram is \(V_0^*V_0\). If \(C=\omega^*h\omega'\), expansion gives
\[
          -C-C^*=I+V_0^*V_0-\omega'^*h\omega'
                         \succeq I+V_0^*V_0 .
\]
For a matrix of positive Hermitian part \(Q\), the factorization
\(Q^{1/2}(I+\mathrm i R)Q^{1/2}\), with \(R\) Hermitian,
shows its determinant has modulus at least \(\det Q\).
Apply this to \(-C\), and use the determinant monotonicity of
positive matrices and \(\det(I+V_0^*V_0)=\det(I+V_0V_0^*)\).
The formula for \(\det C\) from its factors now gives
\begin{equation}
 \mathcal V(\gamma,y)\ge
    \frac{|\delta|\det(I+V_0V_0^*)}{2^{i+1}|\det L|^2}.
 \label{eq:real-schur-negative-basis}
\end{equation}

Changing the transported basis in \(X_0\) by \(A_1^*\) times the first
\(i\) rows of its product matrix changes \(V_0\) by \(A_1^*\)
times the first \(i\) rows of \(L\), since right multiplication by
\(B_1^{-1}L\) takes those product rows to the corresponding rows of \(L\). Let \(L_1\) be the leading \(i\)-block
and \(l\) the final diagonal entry. Lower triangularity leaves
the final column \(v_z\) of \(V_0\) fixed. Thus the determinant
in \eqref{eq:real-schur-negative-basis} is \(\det(G_1+T^*T)\),
with \(T=L_1^* A_1+C_1\), \(G_1=I+v_zv_z^*\), and \(C_1\)
fixed.

In a lower cell the last column of \(\omega'\) is a
multiple of the common vector with coefficient one at \(\gamma z\).
Its norm is \(-|l|^2\) by comparing its last pairing with
\(L\). It is orthogonal to the difference \(\omega'_z-\omega_z\),
so the norm of \(\omega_z\) gives \(|l|^2=1+\|v_z\|^2=\det G_1\).
In the top cell the isotropic pairings of \(\omega_z\) vanish
by triangularity. For \(u\) in the Bessel unipotent, \(u^{-1}z-z\in A\),
so the pairing of \(u^{-1}z\) with \(\omega_z\) equals that
of \(z\). The length bound gives such a \(u\) with
\(\|uy\|\le\theta\), because the infimum defining the length
includes the additional one. Evaluating in the orthonormal basis
after this move gives \(|l|\le\theta\). Thus in both cases
\[
                  |l|^2/\det G_1\le\theta^{2(j-i)} .
\]

Use Lebesgue measure on complex matrix entries. The change
\(T'=(L_1^*A_1+C_1)G_1^{-1/2}\) has inverse Jacobian
\(|\det L_1|^{-2k}\det(G_1)^i\).
The integrand contributes \(\det(G_1)^{-(k+i)}\).
Thus the factor in the integral of the determinant part of
the majorant is \(|\det L_1|^{-2k}\det(G_1)^{-k}\) times
\[
       \int_{\mathbb C^{i\times k}}\det(I+T'^*T')^{-k-i}\,dT'.
\]
Integrate one row at a time. If \(T_0\) contains the remaining rows
and \(Q=I+T_0^*T_0\), appending a row \(t\)
changes the determinant to \(\det Q(1+tQ^{-1}t^*)\).
The substitution \(t=uQ^{1/2}\) has Jacobian \(\det Q\).
The remaining exponent is \(-(k+i-1)\), and integration over this row contributes
the radial integral \(\pi^k\Gamma(i)/\Gamma(k+i)\).
Thus the integral is \(\pi^{ik}\prod_{h=1}^i\Gamma(h)/\Gamma(k+h)\).
Each factor \(\pi^k\Gamma(h)/\Gamma(k+h)\) is at most
\(\pi^k/k!\le1\) for \(k\ge7\).
Consequently the explicit determinant majorant for
\(\mathcal V(\gamma,y)^{-(k+i)}\) in the shift parameter
has integral at most
\begin{equation}
       \frac{(2\rho^2)^{(i+1)i}
                   (2^{i+1}\theta^{2(j-i)})^k}{|\delta|^{k+i}} .
 \label{eq:real-schur-majorant}
\end{equation}
The majorant is a measurable function of the parameter alone,
so this does not require a measurable choice of a basis at each
parameter. All integrals use the original complex entry measure.

\subsection{The contribution bounded by Schur's test}
\label{sec:marked-schur-orbits}

\sectionopening{In this subsection, we combine the real Schur estimates of
Section~\ref{sec:real-one-place-schur} with the finite orbital
weights. This controls the part of the near-orbit sum complementary to
the range treated in Section~\ref{sec:marked-dominant-orbits}.}

Let \(U_c^0\subset U_c\) fix the basis pointwise. On the
orbit \(U_c^0\gamma U_{c'}^0\) the product matrix and
the residue in the block \((c,s,c')\) are constant. Use the index \(i_\gamma\) of
\eqref{eq:marked-level-gauge}. In the transported integral
class coordinates it is the index of the lattice generated by the first
\(i\) isotropic rows of the pairing matrix in the lattice generated by
all its isotropic rows. At a finite place \(w\), conjugate
transposition identifies the column quotient for the adjoints with
the conjugate of the row quotient at \(\bar w\). Hence the lengths
at a split pair are exchanged and the product of residue norms
giving \(i_\gamma\) is preserved. Multiplication by \(M\) also
identifies the row quotient with the quotient of the scaled lattices.
For \(u\gamma v\) in the orbit, \(v\) fixes the basis, and
\((u-1)/M\) is integral in the class basis. Its middle component
reads only the dual isotropic coordinates. These are exactly the
isotropic rows of the pairing matrix.

Put \(D=\Lambda_c^0\). Let \(\Xi\in M_{\varrho,n}(E)\)
consist of the isotropic rows of the product matrix and write
\(\Xi_i\in M_{i,n}(E)\) for its first \(i=\min(j,\varrho)\) rows.
Write \(\Lambda_{\rm all}\) and \(\Lambda_{<j}\) for their respective
row lattices. The difference of two moved bases has zero
\(z\)-component because their pairings with \(z\) agree.
In the Witt blocks the middle part of \((u-1)/M\) on a free dual
isotropic column belongs to \(D\). Indeed the result is integral,
its \(z\)-pairing is zero since \(u\) fixes \(z\), and the isotropic
and dual isotropic lattices split off the middle lattice.
The middle block on the other columns is zero. Consequently each
middle difference divided by \(M\), measured from one fixed term,
is in \(D\otimes\Lambda_{\rm all}\), with all tensor products over
\(\widehat{\mathcal O}_E\). Its residue modulo
\(D\otimes\Lambda_{<j}\) has at most
\[
 \#\bigl(D\otimes(\Lambda_{\rm all}/\Lambda_{<j})\bigr)
       =[\Lambda_{\rm all}:\Lambda_{<j}]^k=i_\gamma^k
\]
possible values, since \(D\) is free of rank \(k\) at every finite
place. Fix a base term \(p_0\) in each residue class and put
\(X(p,p_0)=(pF-p_0F)_{X_0}\in M_{k,n}(E)\). By equality of residues,
\(X(p,p_0)=M Q_f\Xi_i\) for a finite adelic matrix \(Q_f\) with
\(i\) columns in \(D\). Below top rank, \(i=j\) and the leading
\(j\)-minor is nonzero. At top rank \(\Xi_i\) consists of the first
\(\varrho\) rows of the invertible full matrix, also when
\(\delta_\varrho=0\). Choose \(R_i\in M_{n,i}(E)\) with
\(\Xi_i R_i=I_i\). Then
\[
 Q_\lambda=M^{-1}X(p,p_0)R_i\in M_{k,i}(E),\qquad
 (Q_\lambda)_f=Q_f,\qquad X(p,p_0)=M Q_\lambda\Xi_i .
\]
Thus every column of \(Q_\lambda\) lies in \(X_0(E)\cap D\).
At fixed \(\lambda\) the basis differences lie in \(A\). The
packing estimate restricts to \(n\) copies of \(y_v^{-1}A\)
at each real place, of total real dimension at most \(2n\varrho f\).
The separation stays \((N\rho)^{-1}\) and the radius is \(2\rho\),
with \(\rho=b^{37n}\).

Distinct bases thus contribute at most
\(b^{6n^2\cdot 2n\varrho f}\le\mathfrak t=b^{kf/4}\), as in the
packing calculation in Section~\ref{sec:marked-full-to-cut}.
The fixed-basis argument of Section~\ref{sec:marked-full-to-cut}
puts the relative element of two terms with the same basis in
\(U_{c'}^0\). Their right cutoff weights therefore sum to at most one.
For \(Q\in\mathbf C^{\varrho\times k}\) use the following real
Eichler shift. It fixes \(z\) and on \(A\oplus X_0\oplus A^\vee\) is
\[
 s_Q=\begin{pmatrix}
 I&Q&-QQ^*/2\\
 0&I&-Q^*\\
 0&0&I
 \end{pmatrix},\qquad s_Q^{-1}=s_{-Q}.
\]
In the coordinates of Section~\ref{sec:gram-cell-neighborhood} this
has \(B=Q^*\) and \(C=QQ^*/2\), so \(C+C^*=B^*B\) and the
skew-Hermitian part is zero. Thus \(s_Q\) is unitary and fixes the
frame pointwise. Its inverse adds \(Q^*\) times the dual isotropic
coordinates of a column to its \(X_0\)-coordinates. On the moved
frame those dual coordinates are the isotropic rows of the pairing
matrix. The diagonal rescaling \(\mathcal S_v\) has entries
\(M/d_{v,a}\), \(a<i\), with the signed complex basis scale
\(d_{v,a}=\kappa_v^{\varrho-a}/\sqrt N\).
The change of variables in \eqref{eq:real-schur-majorant}
multiplies the integral by \(|\det\mathcal S_v|^{-2k}\).

Writing \(w_v=|\kappa_v|^{\sum_{a<j}(\varrho-a)}\)
and \(P_v=|\Phi_v(\delta_j)|\), we have
\[
       |\det\mathcal S_v|=(M\sqrt N)^i/w_v,\qquad
       |\delta|=w_v^2P_v/N^j .
\]
The powers of \(w_v\) cancel in the \(k\)-th power and
\(w_v\ge1\) allows us to drop them in the other denominator.
Since \(\theta=b^{6n+2}\), \(2^{i+1}\theta^{2(j-i)}\le b^{13i}\).
In lower cells this follows from \(j=i\) and \(2^{i+1}\le b^i\).
In the top cell use \(i=\varrho\), \(2^{i+1}\le b^{i-17}\) and
\(2(6n+2)=12i+16\).
The shifted majorant has integral at most
\[
 K_v=\left(\frac{(2\rho^2)^{i+1}N^j}{P_v}\right)^i
          \left(\left(\frac{b^{13}}{M^2}\right)^i
                        \frac{N^{j-i}}{P_v}\right)^k .
\]
Then \(\mu_\gamma=\prod_v P_v^2/N^{2f(j-i)}\) is the quantity of
Section~\ref{sec:marked-pair-gain}, since
\(\prod_vP_v^2=|N_{E/\mathbf Q}(\delta_j)|\) and \(j-i=(j-\varrho)_+\).
The rational minor is
integral and nonzero, so \(\prod_vP_v\ge1\). From \(N\le b^{60n}\)
and \(2\le b\), the product of the first factors in \(K_v\)
is at most \(b^{135n^2if}\le\mathfrak t\).

Consequently
\begin{equation}
      \prod_v K_v\le
       \mathfrak t\left(\frac{b^{13}}{M^2}\right)^{ikf}\mu_\gamma^{-k/2}.
 \label{eq:marked-schur-product}
\end{equation}
The lattice of column parameters has Lebesgue covolume at least
one. Pairing two column tuples in the definite lattice gives an
algebraic integer. Because \(E\) contains \(\mathbf Q(\mathrm i)\),
the sum of its real parts over the chosen complex places is the
real part of the trace to \(\mathbf Q(\mathrm i)\), hence an
integer. The real Gram determinant of a lattice basis is therefore
a positive integer. This gives the covolume bound in precisely the
measure used in the matrix integral.
Write \(\beta_c^{\rm box}=\mathbf1_{\rm box}\beta_c\),
where boxed means reduced with \(\ell_v\le b^{6n+2}\) at
every place, as in Section~\ref{sec:marked-labelled-frames}.

For a term \(p\) in the orbit define
\[
 J(p)=\int \beta_c^{\rm box}(y)\|W^{\flat}(y)\|^2
        \beta_{c'}^{\rm box}(p^{-1}y)
        \prod_v\mathcal V(p_v,y_v)^{-(k+i)}\,dy .
\]
For a fixed color with base \(p_0\), write \(p=p_\lambda\)
for any term with lattice parameter
\(\lambda_v=\Phi_v(Q_\lambda)^*\) at each real place \(v\),
allowing several terms at each parameter. Let \(y_0\) be boxed and
let \(u\) belong to the real basis stabilizer \(U_\infty^0\).
At every \(v\), triangular reduction gives
\(\|(y_0)_v^{-1}f_a\|\le b^{27}\ell_v((y_0)_v)
\le b^{6n+29}\le\rho\), and boxing gives
\(\ell_v((y_0)_v)\le\theta\).
Use the Eichler shift \(s(X)\) whose central parameter is
zero and whose matrix parameter at \(v\) has first \(i\) rows
\(\mathcal S_v^*X_v\) and all subsequent rows zero. Write
\(Z_i\) for the first \(i\) rows of the pairing matrix in the orthonormal basis.

For every term in the color, the product matrix of
\((u s(X))^{-1}p\) is independent of the term and of \(X\).
At a real place the \(X_0\)-parts satisfy
\[
 \bigl((u s(X))^{-1}pF\bigr)_{X_0}
  =\bigl(u^{-1}p_0F\bigr)_{X_0}
       +\bigl(\mathcal S_v^*(\lambda_v+X_v)\bigr)^* Z_i .
\]
Indeed an inverse basis-fixing unipotent with block \(Q\)
adds \(Q^*\) times the rows of dual isotropic coordinates to the \(X_0\)-part.
The \(u^{-1}\) additions at \(p,p_0\) coincide because they
have the same product matrix. The \(s(X)^{-1}\) addition
uses precisely \(Z_i\); the rational lattice difference,
read in the orthonormal basis, has coefficient
\((\mathcal S_v^*\lambda_v)^*\). Indeed, with \(D_v\) the diagonal
matrix of all frame scales and \(D_{v,i}\) its first \(i\) entries
as a diagonal matrix, one has
\[
 Z_i=D_{v,i}^*\Phi_v(\Xi_i)D_v,\qquad
 M\Phi_v(Q_\lambda\Xi_i)D_v
   =M\Phi_v(Q_\lambda)(D_{v,i}^*)^{-1}Z_i
   =(\mathcal S_v^*\lambda_v)^*Z_i.
\]
Here \(D_{v,i}^*=D_{v,i}\) since \(\kappa_v=M\Phi_v(\delta_E)/(2\pi i)\)
is real, and \(\lambda_v=\Phi_v(Q_\lambda)^*\).
The identity follows. This formula treats \(u\) as fixed
before \(X\), which is the order of right translation
on the real basis-fixer fiber.

Since \(u s(X)\) fixes the basis, \(\mathcal V(p,u s(X)y_0)\)
equals \(\mathcal V((u s(X))^{-1}p,y_0)\). Apply the
one-place majorant with base \(u^{-1}p_0\) and parameter
\(\lambda+X\). It has product integral at most
\(\prod_v K_v\), independently of \(u,y_0\). The basis \(\omega\) and the
matrix \(L\) use the common product matrix and \(y_0\).
The unipotent shifts fix the basis pointwise, so the
common leading pairs used in the same cell are preserved.

At each
lattice point the sum of the second boxed cutoffs is at
most \(\mathfrak t\) by basis packing and the right cutoff
partition at a fixed basis. Hence the color sum at
\(u s(X)y_0\) is bounded by
\(\mathfrak t\sum_{\lambda}\varphi(\lambda+X)\) for that
product majorant \(\varphi\). If \(\mathcal D\) is the
measurable fundamental domain of the \(i\) copies of the
definite lattice in the conjugate-transpose coordinates,
write \(\Lambda\) for that full lattice. Conjugate transpose preserves
its real covolume. The parameters attained in the color form a subset
of \(\Lambda\), since they are measured relative to its chosen base term.
Enlarge the nonnegative sum to all of \(\Lambda\). Tonelli gives
\[
 \frac{1}{\operatorname{vol}(\mathcal D)}
  \int_{\mathcal D}\sum_{\lambda\in\Lambda}\varphi(\lambda+X)\,dX
   =\frac{\int\varphi}{\operatorname{vol}(\mathcal D)}
   \le\prod_v K_v.
\]
Use the same domain for each of the at most \(i_\gamma^k\)
colors.

\subsubsection{Summation over the smaller orbits}
\label{sec:orbital-schur-smaller-orbits}

\sectionopening{In this subsection, we use the preceding Schur bounds to sum the smaller orbits required in Section~\ref{sec:marked-pair-weight-sum}.}

For the integration in stages, let \(\mathcal F\) be the
sum over the whole smaller orbit of the second boxed
cutoff times the inverse volume power. This sum is
left invariant under \(U_c^0\), which permutes the terms.
On the fiber through \(y_0\), the first boxed indicator
and the squared Bessel norm are constant. Write its
remaining cutoff as \(a(u)\), so
\(\sum_{\eta\in U_c^0}a(\eta u)=1\).

For each \(X\), right Haar on the nilpotent group gives
\[
 \begin{aligned}
 \int a(u)\mathcal F(uy_0)\,du
 &=\int a(u s(X))\mathcal F(u s(X)y_0)\,du\\
 &=\int a(u)\mathcal F(u s(X)y_0)\,du .
 \end{aligned}
\]
For the second equality insert
\(\sum_\eta a(\eta^{-1}u)=1\) in the first integral,
substitute \(u\mapsto\eta u\) in each term, and sum
\(\sum_\eta a(\eta u s(X))=1\). Tonelli applies as
all weights are nonnegative. Averaging in \(X\) and
using the preceding bound now yields the inequality on
each real fiber. Unipotent disintegration of the first boxed cutoff
times \(\|W^\flat\|^2\) contributes at most \(\nu d_B\).
Indeed \(\int_{U_\infty}\beta_c(uy)\,du=\nu\) by
Section~\ref{sec:marked-class-covolume},
and each cut vector has quotient mass at most the mass
\(d_v\) of the full vector since its multiplier has
modulus at most one.

Combining the colors with \eqref{eq:marked-schur-product} gives
\[
 \sum_{p\in U_c^0\gamma U_{c'}^0}J(p)
 \le \nu d_B\mathfrak t^2\left(\frac{b^{13}}{M^2}\right)^{ikf}
                             (i_\gamma/\sqrt{\mu_\gamma})^k .
\]
Let \(d\mu_p(y)=\beta_c^{\rm box}(y)
\beta_{c'}^{\rm box}(p^{-1}y)\,dy\) and
\(V_p(y)=\prod_v\mathcal V(p_v,y_v)>0\). The pushforward of
\(\mu_p\) by \(y\mapsto p^{-1}y\) is \(\mu_{p^{-1}}\) with the
classes switched, and \(V_{p^{-1}}(p^{-1}y)=V_p(y)^{-1}\).
Define the absolute amplitude
\[
 a(p)=\int\|W^\flat(y)\|\,\|W^\flat(p^{-1}y)\|\,d\mu_p(y).
\]
The integrand contains the Hilbert--Schmidt norm of a rank-one
operator, the product of the two vector norms. On its support
the boxed cutoffs equal the original cutoffs. Insert the factors
\(V_p(y)^{-(k+i)/2}\) and \(V_p(y)^{(k+i)/2}\) in the two
respective factors of the integrand. Cauchy--Schwarz gives
\(a(p)^2\le J(p)J(p^{-1})\) with the original Haar measure.
Write \(\mathcal O=U_c^0\gamma U_{c'}^0\).
Inversion bijects \(\mathcal O\) over \((c,s,c')\) with
\(\mathcal O^{-1}=U_{c'}^0\gamma^{-1}U_c^0\) over
\((c',s^{-1},c)\). Cauchy--Schwarz on finite subsets and
monotone convergence therefore give
\[
 \left(\sum_{p\in\mathcal O}a(p)\right)^2
 \le \left(\sum_{p\in\mathcal O}J(p)\right)
       \left(\sum_{q\in\mathcal O^{-1}}J(q)\right),
\]
where the second \(J\) uses the interchanged classes. In the
definition of the link, \(w'=p^{-1}w\) exchanges the two
quotient coordinates. The identity \(Z_{\gamma^{-1}}=Z_\gamma^*\)
preserves the absolute norms of the leading minors and hence
\(\mu_\gamma\). Exchanging the two bases in each common vector
and rescaling its coefficient pair \((1,\zeta)\) to
\((\zeta^{-1},1)\) gives the same cell rank for \(\gamma^{-1}\).

For \(a\in\Lambda_{c'}\cap A\) and a basis vector \(f\) of
\(F\), the identity \(h(a,\gamma^{-1}f)=h(\gamma a,f)\)
identifies the full pairing lattices and their sublattices generated
by the first \(\min(j,\varrho)\) columns. Thus the inverse orbit
has first index \(i'_\gamma\) and reverse index \(i_\gamma\).
The pairing matrices are constant on each smaller orbit, so
these identities apply to the fixed representatives on both
sides of the inversion. Define
\(g_{\rm min}=N^f\mu_\gamma/(i_\gamma i'_\gamma)\). We obtain
\[
 \left(\sum_p a(p)\right)^2
   \le \left(\nu d_B\mathfrak t^2(b^{13}/M^2)^{ikf}\right)^2
                                                (N^f/g_{\rm min})^k .
\]
With \(n_p,D_j\) as in Section~\ref{sec:marked-dominant-orbits}, \(d_B^2N^{kf}=d_BD_{\rm id}/n_p\) and
\(b^{-8ikf}=\mathfrak t^{-32i}\). The bounds
\(b^{13}/M^2\le(b^4M)^{-1}\), \(d_B\le\mathfrak t^8\) and
\(D_j\ge\mathfrak t^{8-32i}D_{\rm id}\) give
\begin{equation}
       \sum_p a(p)\le
       \nu\mathfrak t^2\left(\frac{D_j}
                  {n_p(M^{2if}g_{\rm min})^k}\right)^{1/2}.
 \label{eq:marked-orbit-schur}
\end{equation}
The last bound on \(D_j\) holds in lower cells because \(i\ge1\),
and in the top cell because
\(D_{\rm main}\ge\mathfrak t^{44}D_{\rm id}\).
For \(\mathfrak g_{\rm norm}(\gamma)\le N^f\mu_\gamma\)
the denominator uses \(G_\gamma=M^{2if}g_{\rm min}\),
so \eqref{eq:marked-orbit-schur} is exactly \(\nu w_\gamma\).
A nonzero cut summand with the original two cutoffs puts
both arguments in their boxes. Nearness is invariant
throughout the smaller orbit, because both basis stabilizers
preserve reduction and length. Hence a non-near smaller
orbit contributes zero. Sum the bound over the original
fixed near representatives of the complementary part.

\subsection{Summing the pair weights}
\label{sec:marked-pair-weight-sum}

\sectionopening{In this subsection, we sum the pair multiplicities with the bounds from
Sections~\ref{sec:marked-dominant-orbits}
and~\ref{sec:marked-schur-orbits}. This completes the estimate for the
corresponding terms of the full Petersson matrices in
Section~\ref{sec:full-orbital-blocks}.}

Use the weights \(w_\gamma\) of Section~\ref{sec:marked-dominant-orbits}.
In each agreeing block \((c,s,c')\), take a family of ordered
pairs from \(\mathcal R_{\Gamma_*}(c,s,c')^2\), using the representatives
fixed in Section~\ref{sec:marked-labelled-frames} for both terms.
Suppose for each such pair of ranks \(j,j'\) that
\((d\eta)^2\le G_\gamma G_{\gamma'}\)
and \(D_j,D_{j'}\le D_0\).
Here \(d=d(L_{12})\) and \(l\) is its corank.
If \(\mathfrak t^{l+6}D_0\eta^{-k}\le B_0\)
for every pair of the family, then
\[
 w_\gamma w_{\gamma'}
 \le \frac{\mathfrak t^4D_0}{n_p d^k\eta^k}
 \le \frac{B_0}{n_p\mathfrak t^{l+2}}d^{-k}.
\]
By Sections~\ref{sec:marked-pair-count} and
\ref{sec:marked-pair-gain} the sum at any \(l\) is at
most \(B_0\). There are at most \(n^3+1\le\mathfrak t\)
coranks, so the total is at most \(\mathfrak t B_0\).

For distinct links, where the codimension \(i=k+1-\dim V_1\) of
Lemma~\ref{lem:discriminant-pair-lattice} satisfies \(i\ge1\), take
\(\eta=(M/2)^{if}\) and \(D_0=D_{\rm main}\).
Using \(l\le i+1\) and
\(\xi=\mathfrak t(2/M)^{kf}\le1\) gives
\[
 \mathfrak t^{l+6}D_0\eta^{-k}
  \le \mathfrak t^7\xi^i D_{\rm main}
  \le \mathfrak t^8(2/b^6)^{kf}D_{\rm id}
  \le \mathfrak t^{-12}D_{\rm id},
\]
where the middle step uses \eqref{eq:marked-main-identity-ratio}.
The rank of the two \(z\)-coordinate forms restricted to a
rational link depends only on the link. It is two at top rank
and one below top rank by Section~\ref{sec:pair-lattice-splitting}.
Thus equal links have the same top-rank status. For equal lower
links \(i=l=0\). Use \(\eta=1\) and
\(D_0=\mathfrak t^{-16}D_{\rm id}\), giving
\(\mathfrak t^{l+6}D_0\eta^{-k}=\mathfrak t^{-10}D_{\rm id}\).
For each ordered pair in either family these bounds give
\[
 w_\gamma w_{\gamma'}\le
 \frac{\mathfrak t^{-10}D_{\rm id}}{n_p\mathfrak t^{l+2}}
 \,d(L_{12})^{-k}.
\]
Apply \eqref{eq:marked-weighted-pairs} once to the union at each
corank \(l\), then sum over the at most \(\mathfrak t\) coranks.
The combined family has total pair weight at most
\begin{equation}\label{eq:marked-other-pair-weights}
                  \mathfrak t^{-9}D_{\rm id}\le b^{-kf}D_{\rm id}.
\end{equation}
For equal top links with meeting norm at least \(z\ge1\),
equation~\eqref{eq:marked-pair-discriminant-gain} gives
\(\eta=z\) and \(l=1\). Take \(D_0=D_{\rm main}\).

\begin{lemma}[Weight of equal top links]\label{lem:marked-equal-link-weight}
Fix \(z\ge1\). In each agreeing block \((c,s,c')\) with
\(Q_{\cap}(s)\ge z\), let \(\mathcal E_z(c,s,c')\) be the ordered pairs
in \(\mathcal R_{\Gamma_*}(c,s,c')^2\) whose terms both have cell rank
\(j=n\) and equal rational links in \((A^\perp/A)^2\).
In every other block put \(\mathcal E_z(c,s,c')=\varnothing\). Then
\begin{equation}\label{eq:marked-equal-pair-weights}
 \sum_{c,s,c'}\,\sum_{(\gamma,\gamma')\in\mathcal E_z(c,s,c')}w_\gamma w_{\gamma'}
 \le\mathfrak t^8 z^{-k}D_{\rm main}=(b^{2f}/z)^kD_{\rm main}.
\end{equation}
\end{lemma}
\begin{proof}
Every pair has \(l=1\) and individual weight at most
\[
 \frac{\mathfrak t^4D_{\rm main}z^{-k}}{n_p d(L_{12})^k}.
\]
Thus \eqref{eq:marked-weighted-pairs} gives
\[
 \sum_{c,s,c'}\,\sum_{(\gamma,\gamma')\in\mathcal E_z(c,s,c')}w_\gamma w_{\gamma'}
 \le \mathfrak t^7D_{\rm main}z^{-k}
 \le \mathfrak t^8D_{\rm main}z^{-k}.\qedhere
\]
\end{proof}

For a stable subset \(S_0\subseteq\Gamma_*\) the cut blocks
satisfy
\begin{equation}\label{eq:needed-near-orbit-size}
 \|A_{S_0}^{\rm cut}(c,s,c')\|_{\rm HS}
 \le \nu\sum_{\gamma\in\mathcal R_{S_0}(c,s,c')} w_\gamma .
\end{equation}
Here \(A_{S_0}^{\rm cut}\) is the cut block of
Section~\ref{sec:marked-full-to-cut}, with the rational
terms restricted to \(S_0\), and \(\mathcal R_{S_0}\) takes
the chosen near representatives in \(S_0\) with cell
ranks over the given agreeing block. Stability means
\(U(F_0)S_0U(F_0)=S_0\). In particular the rational lifts of the
finite unipotents in Section~\ref{sec:marked-orbital-cancellation}
preserve \(S_0\). Fix the same near representatives for every piece.
At a nonagreeing residue take \(\mathcal R_{S_0}=\varnothing\).
Split each cell at
\(\mathfrak g_{\rm norm}=N^f\mu_\gamma\), assigning equality
to the Schur part.

For a unipotent move the product matrix is
multiplied on either side by triangular matrices with diagonal
one, upper triangular on the right and lower triangular on
the left. Hence its leading minors, including the determinant,
are preserved. The ranks are preserved because the changes
of basis preserve the flags and the common leading pairs.
Sections~\ref{sec:marked-dominant-orbits} and
\ref{sec:marked-schur-orbits} give \eqref{eq:needed-near-orbit-size}
on the two stable parts with the same near representatives.
The nonagreeing block vanishes by the finite cancellation
in Section~\ref{sec:marked-orbital-cancellation}.

Now partition the cut sum by links at top rank and group
all lower ranks in one piece. To see stability of a link,
write its two quotient coordinates as the images of
\(w\) and \(\gamma^{-1}w\) in \(A^\perp/A\).
For a move \(\gamma\mapsto u^{-1}\gamma u'\), use
\(u^{-1}w\) and \((u')^{-1}\gamma^{-1}w\). Both unipotents
fix \(A\) as a subspace and act trivially on \(A^\perp/A\),
so the two quotient coordinates are unchanged.

In a fixed block \((c,s,c')\), let \(A_{\rm low}\) be the sum of the
cut amplitudes of ranks \(j<n\), and let \(A_\Lambda\) be the cut
amplitude of rank \(n\) with link \(\Lambda\). By absolute
convergence in Section~\ref{sec:marked-full-to-cut}, for
\(A=A_{\rm low}+\sum_\Lambda A_\Lambda\) we can expand
\[
\begin{aligned}
 \|A\|_{\rm HS}^2-\sum_\Lambda\|A_\Lambda\|_{\rm HS}^2
 &=\|A_{\rm low}\|_{\rm HS}^2
     +2\Re\sum_\Lambda\langle A_{\rm low},A_\Lambda\rangle_{\rm HS}\\
 &\quad+\sum_{\Lambda\ne\Lambda'}
                  \langle A_\Lambda,A_{\Lambda'}\rangle_{\rm HS}\\
 &\le \|A_{\rm low}\|_{\rm HS}^2
      +2\|A_{\rm low}\|_{\rm HS}\sum_\Lambda\|A_\Lambda\|_{\rm HS}\\
 &\quad+\sum_{\Lambda\ne\Lambda'}
           \|A_\Lambda\|_{\rm HS}\|A_{\Lambda'}\|_{\rm HS}.
\end{aligned}
\]
Summing the squares of the top pieces gives the equal-link moment after
normalization. The last sum uses ordered pairs, and the factor \(2\)
accounts for both orders of a lower and a top term. Applying
\eqref{eq:needed-near-orbit-size} to each piece leaves on the right
only pairs of distinct links or equal lower links. For fixed
\(c,c'\), the first term determines
\(s=\operatorname{red}(B_c^{-1}\gamma_f B_{c'})\), and the second
term must have the same residue. Each term belongs to its unique
top link piece or to the single lower piece. Thus every ordered
pair on the right occurs in just one product over its two classes.
Within a top square the first term also determines the common
rational link, so summing the top squares introduces no further
multiplicity.

Combine the block estimates first and then restrict to any
\(\Omega\subseteq G\). By \eqref{eq:marked-other-pair-weights},
\[
 \begin{aligned}
 \mathcal L(\Omega)&=\sum_{s\in\Omega}\sum_{c,c'}
               \sum_\Lambda\|A_\Lambda(c,s,c')\|_{\rm HS}^2,\\
 \nu^{-2}\sum_{s\in\Omega}\sum_{c,c'}
      \|A_{\Gamma_*}^{\rm cut}(c,s,c')\|_{\rm HS}^2
 &\le \nu^{-2}\mathcal L(\Omega)+b^{-kf}D_{\rm id}.
 \end{aligned}
\]
Applying the same size estimate twice within each
top piece and using \eqref{eq:marked-equal-pair-weights}
bounds the normalized equal-link moment by
\((b^{2f}/z)^kD_{\rm main}\) on residues of meeting
norm at least \(z\).
The partition of the blocks is justified by the
absolute convergence in Section~\ref{sec:marked-full-to-cut}.
Thus the pair summation gives the bounds needed in
Section~\ref{sec:full-orbital-blocks}. The reduction to
the cell ranks used in the sum is
\eqref{eq:marked-off-bessel-reduction}.

\section{The common automorphic family and its translated span}
\label{sec:dimension-short-comparison}

\sectionopening{In this section, we assemble the data of Section~\ref{sec:period} and
Sections~\ref{sec:finite-orbital-analysis}--\ref{sec:orbital-cancellation-analysis}
on one family. The moment inequalities of Section~\ref{sec:marked-moments}
then give the dimension bound in Proposition~\ref{prop:marked-upper-estimates}.}

\subsection{Construction of the common automorphic family}
\label{sec:common-filling-assembly}

\sectionopening{In this subsection, we assemble the real and finite data and verify that the resulting
Poincar\'e map has the specified levels at every sufficiently large
scale. The normalization of Section~\ref{sec:period} and the dimension
in Section~\ref{sec:marked-target-dimension} apply to this same
construction.}

Fix \(q_0,d,\epsilon\). For each sufficiently large integer \(x\),
take the primes, weights and skew different of
Section~\ref{sec:datum}. The grading, finite point and
levels in Section~\ref{sec:period} use these choices
and the specified \(u\). The source level is exactly the
preimage of the original principal level under the definite
embedding, by the glued-lattice calculation in
Section~\ref{sec:fixed-vector-common-pairing}. At the chosen split
places reduction of the full integral group is surjective:
lift each prescribed general linear matrix at the chosen
place, use the inverse adjoint at its conjugate and use
the identity at the remaining places. The change of finite isometry has identity reduction at every chosen place.

Use the real signs and weights of Section~\ref{sec:real-marked-filling}.
The different and scales in Section~\ref{sec:datum} give
\[
 |\Phi_v(\delta_E)|\le b,\qquad
 1\le|\kappa_v|\le b^{23},\qquad N\le b^{25m}\le b^{60n}.
\]
Use the adelic point \(x^{\rm pt}\) and the pair \((g_0,t)\)
of Section~\ref{sec:period} in the Poincar\'e map. These
choices precede the selection of forms and residue classes. At each
real place the value space is a determinant twist of the symmetric
power of degree \(a_\tau\le a\) of the standard rank-\(k\)
representation, so its dimension is at most
\((a+k)^k\le b^{2k}\).

In particular \(0<d_B\le b^{2kf}\)
and the source class set is nonempty.
Equation~\eqref{eq:marked-main-identity-ratio} gives both
bounds for \(D_{\rm main}/D_{\rm id}\).
Before imposing the principal level at \(M\), the inverse
image of the glued lattice stabilizer fixes the specified nonzero
covector at each chosen place above \(N\).
The last-row criterion of
Section~\ref{sec:period-finite-coordinates} and transitivity
on nonzero covectors give the index \((N^k-1)^f\ge(N^k/b)^f\). Finally the mass
comparison in Section~\ref{sec:open-equal-link-bound} gives
\[
 C_0\Delta_0\le (1+b^{-9})\nu^2D_{\rm main}\mathcal B,\qquad
 \Delta_0=\frac{v_E^n}{c_G^f}\left(\frac{d_B}{\Lambda_0}\right)^2 .
\]
This is the same constant bound used in the dominant orbit
and in the open equal-link estimates. All these conditions
hold on the constructed Bessel data, with the terminal definite
comparison of Section~\ref{sec:marked-target-dimension}.
Their numerical thresholds depend only on the fixed
\(q_0,d,\epsilon\). The arithmetic comparisons below use the
measured map already specified by \(g_0,t\) and this level.

\subsection{Selecting the translates}\label{sec:marked-span-proof}
\label{sec:translate-selection}

\sectionopening{In this subsection, we convert the two moment inequalities of
Section~\ref{sec:marked-moments} into a lower bound for the dimension
of a sum of translated subspaces. We then pass from the full source to
the ramified subspace of Section~\ref{sec:ramified-source}.}

\begin{lemma}[A dimension bound for translated subspaces]
\label{lem:translated-span-dimension}
Let a finite group \(G\) act unitarily on \(T\) by \(\rho\).
Let \(p\) be the orthogonal projection onto an invariant subspace \(W\),
and let \(B:S\to T\) be linear between Hilbert spaces with
\(1\le d_S=\dim S<\infty\) and \(1\le D=\dim W<\infty\). Write \(P_U\) for orthogonal
projection onto \(U\subseteq S\). Suppose \(0<\eta_1<1/16\) and
\[
 \|pB P_U\|_{\rm HS}^2\ge\dim U-\eta_1 d_S
 \quad(U\subseteq S),\qquad
 D\,\mathbb E_g\|B^*\rho(g)B\|_{\rm HS}^2\le (1+\eta_1)d_S^2.
\]
Then there is \(t\subseteq G\) such that
\[
 |t|d_S\le\tfrac52\left(1+\log\frac1{16\eta_1}\right)D,\qquad
 \dim\sum_{g\in t}p\rho(g)B(S)\ge(1-16\eta_1)D.
\]
\end{lemma}

\begin{proof}
Set \(A_g=\rho(g)BB^*\rho(g)^*\) and
\(C_{\rm av}=\mathbb E_g A_g\). Each \(A_g\) has finite rank.
Cyclicity of the trace gives
\(\operatorname{tr}(A_g A_h)=\|B^*\rho(g^{-1}h)B\|_{\rm HS}^2\).
Since \(h\mapsto g^{-1}h\) preserves the uniform measure,
\[
 \|C_{\rm av}\|_{\rm HS}^2=\mathbb E_{g,h}\operatorname{tr}(A_gA_h)
    =\mathbb E_g\|B^*\rho(g)B\|_{\rm HS}^2 .
\]
For \(F\subseteq W\) of dimension \(\alpha D\), some \(g\) satisfies
\[
 \|P_F\rho(g)B\|_{\rm HS}^2
 \le \operatorname{tr}(P_F C_{\rm av})
 \le \sqrt{\alpha(1+\eta_1)}\,d_S
 \le (1+\alpha+\eta_1)d_S/2.
\]
On the preimage \(U=(p\rho(g)B)^{-1}F\),
\(P_F\rho(g)BP_U=p\rho(g)BP_U=\rho(g)pBP_U\). Restrict to \(U\)
in the near-isometry inequality and compare with the full
Hilbert--Schmidt norm. This gives
\(\dim U\le (1+\alpha+3\eta_1)d_S/2\).
By rank-nullity in \(W/F\), adding the new image adds at least
\((1-\alpha-3\eta_1)d_S/2\) dimensions.

Let \(F_j\) be the span after \(j\) choices, with \(F_0=0\).
Put \(m_j=D-\dim F_j\) and \(a_1=2d_S/(5D)\). Apply the first
moment inequality to \(U=\ker(pB)\). Its left side vanishes, so
\[
 (1-\eta_1)d_S\le\operatorname{rank}(pB)\le D,\qquad
 0<a_1\le \frac{2}{5(1-\eta_1)}<1 .
\]
The dimension increment just proved gives
\[
 m_j>16\eta_1D\quad\Longrightarrow\quad
 m_{j+1}\le(1-a_1)m_j\le e^{-a_1}m_j.
\]
Indeed the increment is at least
\(\tfrac12(m_j/D-3\eta_1)d_S\ge(13/32)(m_j/D)d_S\ge a_1m_j\).
It is positive at every step before termination, so the selected group
element is new. Iteration up to the first
\(m_j\le16\eta_1D\), with \(m_0=D\), gives
\[
 j\le\left\lceil a_1^{-1}\log\frac1{16\eta_1}\right\rceil,\qquad
 jd_S\le\tfrac52\left(1+\log\frac1{16\eta_1}\right)D .\qedhere
\]
\end{proof}

\begin{proposition}[Span of the translates]\label{prop:marked-orbital-span}
For fixed \(q_0,d,\epsilon\) and sufficiently large \(x\), there is
\(t\subseteq G\) such that
\[
 |t|\dim S\le5D\log(1+2b^2),\qquad
 \dim\sum_{g\in t}P_g(S)\ge(1-4/b^2)D .
\]
Here \(P_g=c_B|G|c_{\rm ss}^{-1}p\rho(g)\Pi|_S\), where
\(c_B\) is the positive rational scalar of
Proposition~\ref{prop:shortness-normalization} and \(c_{\rm ss}\ne0\)
is the scalar in \eqref{eq:measured-seesaw-normalized-tensor}.
\end{proposition}

\begin{proof}

Use the Euclidean realization \(C(X)\) of the full
source, where \(C=C_XR_H(\eta^{-1})\), with \(\eta\) as in Section~\ref{sec:period}, is the change of
coordinates in Section~\ref{sec:period}.
The geometric map is \(B:C(X)\to T_{\rm full}\) and the
projection onto the selected types is \(p\). Section~\ref{sec:marked-moments}
gives, for \(U\subseteq C(X)\),
\[
 \|pBP_U\|_{\rm HS}^2\ge\dim U-5b^{-4}d_X,
 \qquad
 D\,\mathbb E_g\|B^*\rho(g)B\|_{\rm HS}^2
          \le(1+9b^{-4})d_X^2 .
\]
The free class coordinates have \(d_X>0\). Applying
the first inequality to \(U=\ker(pB)\) gives
\[
 \dim\ker(pB)\le5b^{-4}d_X,
 \qquad D\ge\operatorname{rank}(pB)\ge(1-5b^{-4})d_X>0 .
\]
The same inequality with \(U=C(S)\) and the codimension estimate of
Section~\ref{sec:ramified-source} give
\[
 \|pBP_{C(S)}\|_{\rm HS}^2\ge (1-b^{-8}-5b^{-4})d_X>0 .
\]
Write \(p=\sum_c p_c\), where \(c\) runs over the selected tuples
at both old rational primes. Orthogonality gives \(p_c B(C(S))\ne0\)
for some \(c\). This component has an allowed type at every chosen
place and a primitive type at one place above each of the two primes.
The commuting normal split Hecke operators away from the old level
preserve the \(c\)-isotypic space, so some simultaneous eigenspace
has nonzero projection of this component. Conversely, for any
nonzero eigenvector \(z\in pT_{\rm full}\), a summand \(p_cz\ne0\)
has the same outside-level Hecke eigenvalues as \(z\), as used in
Section~\ref{sec:main-trace}.
Both error constants are at most \(\eta_1=(8b^2)^{-1}\)
for large \(b\).

The selected translates span all but
\(2D/b^2\) dimensions, and
\[
 |t|d_X\le \tfrac52(1+\log(b^2/2))D
              \le 5D\log(1+2b^2).
\]
Section~\ref{sec:ramified-source} gives
\(\dim(X/S)\le f d_X/(N-1)\le b^{-8}d_X\).
The operator \(R_H(\eta^{-1})\) commutes with the
central scalar averages defining \(S\), so the restricted
source in the Euclidean coordinates is \(C(S)\).
The maps \(p\rho(g)B\) induce the surjection
\[
 \bigoplus_{g\in t} C(X)/C(S)\twoheadrightarrow
 \frac{\sum_{g\in t}p\rho(g)B(C(X))}
      {\sum_{g\in t}p\rho(g)B(C(S))}.
\]
The dimension of the quotient is therefore at most
\[
 |t|\dim(X/S)\le |t|b^{-8}d_X
 \le \tfrac52 b^{-8}(1+\log(b^2/2))D\le 2D/b^2.
\]
With \(v_G=\operatorname{vol}(J)\) and \(\nu\) the
unipotent covolume, the Euclidean coordinates of
Section~\ref{sec:period-petersson-blocks} give
\[
 \begin{aligned}
 B&=(v_G\nu)^{-1/2}\Pi C^{-1},\\
 P_g&=c_B|G|c_{\rm ss}^{-1}p\rho(g)\Pi|_S
     =c_B|G|c_{\rm ss}^{-1}(v_G\nu)^{1/2}p\rho(g)BC|_S.
 \end{aligned}
\]
The scalar is nonzero. Hence
\(\operatorname{im}P_g=p\rho(g)B(C(S))\), giving the
span estimate on the datum of
Section~\ref{sec:common-filling-assembly}.
\end{proof}
For this set \(t\), lift each residue translation integrally at
the specified split places, with the identity at all other finite
places. Its action depends only on its residue by principal-level
invariance. The vectors \(P_g\phi\), \(g\in t\), with \(\phi\) among
the everywhere integral generators of
Lemma~\ref{lem:ramified-short-generators}, span the displayed sum.
A maximal independent subset therefore gives a complex basis of it
before any coefficient place or twist is chosen.
\par\nopagebreak[4]
\part{Theta identities and the upper covolume bound}\label{part:upper}

\section{The Weil action and the measured seesaw}
\label{sec:weil-seesaw}

\sectionopening{In this section, we verify rational invariance and relate the doubled Weil action to the
normalized theta integral. This fixes the scalar and measures in the
comparison with the Poincar\'e map of Section~\ref{sec:period}.}

\subsection{Rational invariance for the Weil action}
\label{sec:marked-rational-theta}

\sectionopening{In this subsection, we verify rational invariance of the
Weil action used in Section~\ref{sec:marked-doubled-action}.
We compute the Levi scalar in a split polarization and
pass to a general Hermitian space through \(B\perp B^-\).}

\begin{lemma}[Composition of the normalized local Weil operators]
\label{lem:weil-vacuum-composition}
Let \(g,h\) belong to the local unitary group of the given Hermitian
space. In the signed real and normalized finite coordinates write
\(T(x)\) for the symmetric Weyl displacement and \(v_0\) for the unit vacuum.
There is a unitary operator \(U_g\) implementing
\(T(x)\mapsto T(gx)\) with
\[
  \langle v_0,U_gv_0\rangle=1/\overline{\Delta_g}.
\]
Here \(\Delta_g=\overline{\det d_g}\) at a real place, where
\(d_g\) is the negative diagonal block of \(g\). At a finite
place \(v\) of \(F_0\), put \(E_v=E\otimes_{F_0}F_{0,v}\) and
\(L=\mathcal O_{E_v}^{\,n}\) in the normalized coordinates. The
norm and character on \(E_v^\times\) are products when \(v\) splits.
Every nonsplit finite place is unramified by
Section~\ref{sec:different-additive-character-scales}.
The finite multiplier is
\(\Delta_g=(\upsilon_v|\cdot|_{E_v}^{-1/2})(\det H_g)\),
where \(H_gL=L\cap g^{-1}L\). Any other choice is \(H_gk\) with
\(k\in\mathrm{GL}(L)\). Both \(\upsilon_v\) and \(|\cdot|_{E_v}\)
are trivial on \(\det k\), since \(\upsilon_v\) is unramified.
Then \(U_gU_h=U_{gh}\), also on the restricted tensor product.
\end{lemma}

\begin{proof}
Write \(e(x)=\langle v_0,T(x)v_0\rangle\), which is real and even
in both models. For the self-dual displacement measure,
\[
 P_0=|v_0\rangle\langle v_0|=\int\overline{e(x)}T(x)\,dx .
\]
Indeed, in Schr\"odinger coordinates \(x=(a,b)\), integration over
the modulation \(b\) is Fourier inversion, and integration over
the translation \(a\) gives \(f\mapsto v_0\langle v_0,f\rangle\).
By the Heisenberg intertwining theorem \cite[no.~10, Theorem~1, p.~157]{WeilOperators}, an isometry \(g\)
has a unitary implementer unique up to a scalar. Its image of
\(P_0\) by conjugation is
\(Q_g=\int e(g^{-1}x)T(x)\,dx\).
It satisfies
\[
 \langle v_0,Q_gv_0\rangle=\int e(x)e(g^{-1}x)\,dx=|\Delta_g|^{-2}.
\]
At a finite place the integral is
\(\operatorname{vol}(L\cap gL)=[L:L\cap g^{-1}L]^{-1}\),
since \(g\) preserves volume and the lattice phases agree.
At a real place \(e(x)=e^{-x^*x/2}\). With \(P\) the positive
coordinate projection, set \(D(a)=P+(I-P)a\). Pseudo-unitarity gives
\((I+a^*a)/2=D(a)^*D(a)\), so the normalized complex Gaussian
integral is \(|\det d_g|^{-2}\) after the substitution \(x=gy\).
Consequently \(U_gv_0=\Delta_gQ_gv_0\) fixes the unique required
phase. For \(k=gh\), covariance gives \(U_gQ_hU_g^{-1}=Q_k\), whence
\[
 U_gU_h=c(g,h)U_k,\qquad
 c(g,h)=\Delta_g\Delta_h\overline{\Delta_k}
                 \langle v_0,Q_kQ_gv_0\rangle,\qquad |c(g,h)|=1 .
\]
Writing \(T(x)T(y)=\beta(x,y)T(x+y)\) and substituting
\(u=-k(x+y),\ v=kx\) in this last double integral gives
\begin{equation}\label{eq:weil-cocycle-integral}
 \langle v_0,Q_kQ_gv_0\rangle
  =\iint e(hx)e(ky)\beta(x,y)e(x+y)\,dx\,dy .
\end{equation}
Indeed \(k^{-1}u=-(x+y),\ g^{-1}v=hx,\ u+v=-ky\) and
\(\beta(u,v)=\beta(x,y)\) by alternation and unitary invariance.

At a finite place write \(c_0=(1-j_{\rm W}^{-1})/2\),
\(\psi_0(z)=\psi_{E,v}(j_{\rm W}\alpha_{\rm W}
\bar\alpha_{\rm W}z)\), and \(e_0(x)=\psi_0(c_0x^*x/2)\).
\(c_0=(1+\sqrt{-d})/2\) is integral and satisfies
\(c_0+\bar c_0=1\). The different identity
of Section~\ref{sec:period-finite-coordinates} makes \(\psi_0\)
self-dual on \(\mathcal O_{E_v}\), and the skew sign of \(j_{\rm W}\)
gives \(\psi_0(\bar z)=\psi_0(z)^{-1}\). In these coordinates
\(\beta(x,y)=\psi_0(y^*x/2)\) and \(e=\mathbf1_L e_0\).
The finite translation fixing the vacuum for integral \(x\) is
\(T_\sharp(x)=e_0(x)T(x)\), as in the finite vacuum formula of
Section~\ref{sec:upper-cm-jet-transport}. Its multiplier is
\(e_0(x)e_0(y)e_0(x+y)^{-1}\beta(x,y)\).
Finite unitary isometries preserve \(x^*x\), so both translations
have the same unitary covariance, including at places above two.
The character \(\psi_0\) is trivial on \(\mathcal O_{E_v}\), so
\(e_0^2=1\) on \(L\). Put
\[
 D=\{(x,y):hx\in L,\ ky\in L,\ x+y\in L\}.
\]
On this compact open lattice \(e_0(hx)=e_0(x)\),
\(e_0(ky)=e_0(y)\) and \(e_0(x+y)=e_0(x+y)^{-1}\). Expanding
the three quadratics in \eqref{eq:weil-cocycle-integral} gives the
integrand \(\psi_0(\bar c_0y^*x)\) on \(D\) and zero elsewhere.

At a split place \(g=(A,A^{-\mathsf t})\). If the elementary divisor
exponents of \(A\) are \(a_i\), the intersection exponents at the
two places are \(\max(0,-a_i)\) and \(\max(0,a_i)\).
Thus \(\Delta_g/|\Delta_g|=\upsilon_w(\det A)^{-1}\) at the chosen
factor \(w\), so \(\Delta_g\Delta_h\overline{\Delta_k}>0\).
Writing \(D=D_+\times D_-\) at the two factors, the phase
\(\psi_0(\bar c_0y^*x)\) is a bicharacter between \(D_+\) and
\(D_-\). Integration over \(D_-\) gives its volume times the
indicator of its annihilator in \(D_+\), which contains an open
subgroup. Hence \(c(g,h)>0\), and therefore \(c(g,h)=1\).
This argument applies at the dyadic split places.

At a nonsplit place let \(q\) be the cardinality of the residue
field of \(F_{0,v}\). The extension is unramified quadratic, and
\(\Delta_g=(-q)^{a_g}\), where
\(a_g=\operatorname{length}_{\mathcal O_{E_v}}L/(L\cap g^{-1}L)\).
Here we used \(\upsilon_v(\varpi)=-1\) on a base uniformizer.
Choose \(m\ge0\) so that \(B=(\varpi^mL)^2\subset D\) and
the phase is constant on cosets of \(B\) in \(D\). The
coordinates of \(x,y\) on \(D\) have valuations bounded below.
Taking \(m\) sufficiently large makes the change in \(\bar c_0y^*x\)
integral after adding a vector of \(B\).
With \(\operatorname{vol}(L)=1\) and residue cardinality \(q^2\)
over \(E_v\), \eqref{eq:weil-cocycle-integral} gives
\[
  q^{4nm}c(g,h)=(-q)^{a_g+a_h+a_k}\mathcal S,\qquad
  \mathcal S=\sum_{(x,y)\in D/B}\psi_0(\bar c_0y^*x).
\]
Simultaneous multiplication by an integral unit of norm \(-1\)
conjugates every summand, so \(c(g,h)\) is real and equals
\(\pm1\). Such a unit exists by norm surjectivity on unramified
units. The Teichm\"uller group of order \(q+1\) and norm one fixes
the summands and acts freely on \((D/B)\setminus\{0\}\), since
\(\xi-1\) is a unit for every nonidentity member \(\xi\).
Consequently \(\mathcal S\equiv1\pmod{(q+1)\mathcal O'}\), where
\(\mathcal O'\) is the ring of integers of the cyclotomic field
generated by the finitely many phase values. Reducing the preceding
equality gives \(c(g,h)\equiv1\pmod{(q+1)\mathcal O'}\).
Since \((q+1)\mathcal O'\cap\mathbb Z=(q+1)\mathbb Z\) and
\(q+1>2\), this forces \(c(g,h)=1\).

At a real place put \(J=2P-I\), \(Q=I-P\), and
\(C_a=(I+a^*a)/2=D(a)^*D(a)\). In these coordinates
\(\beta(x,y)=\exp((y^*Jx-x^*Jy)/2)\).
The integral in \eqref{eq:weil-cocycle-integral}, with measure
\(\pi^{-2n}dx\,dy\) in real Lebesgue coordinates, has matrix
\[
 A=\begin{pmatrix}C_h&P\\ Q&C_k\end{pmatrix},\qquad
 A+A^*=\operatorname{diag}(h^*h,k^*k)+
                      \begin{pmatrix}I&I\\ I&I\end{pmatrix}>0.
\]
Thus its value is \((\det A)^{-1}\) by the complex Gaussian
integral (first for positive Hermitian matrices and then by
holomorphic continuation on matrices with positive Hermitian part).
Write \(a=\left(\begin{smallmatrix}a_{11}&a_{12}\\c_a&d_a
\end{smallmatrix}\right)\) in signed coordinates. Since
\(D(a)=\left(\begin{smallmatrix}I&0\\c_a&d_a\end{smallmatrix}\right)\),
the off-diagonal blocks of \(C_a^{-1}\) are
\(-c_a^*d_a^{-*}\) and \(-d_a^{-1}c_a\).
The Schur complement and \(g=kh^{-1}\) give
\[
\begin{aligned}
 \det A
 &= |\det d_h\det d_k|^2
       \det(I-d_h^{-1}c_hc_k^*d_k^{-*})\\
 &= |\det d_h\det d_k|^2
       \frac{\det(d_hd_k^*-c_hc_k^*)}{\det d_h\det d_k^*}
  =\overline{\det d_h\det d_g}\det d_k .
\end{aligned}
\]
This is precisely \(\Delta_g\Delta_h\overline{\Delta_k}\), also
when either signed summand is zero. Hence \(c(g,h)=1\).
The map \(g\mapsto e(g^{-1}\,\cdot)\) is continuous in \(L^1\).
Since \(\Delta_g\) is continuous, \(U_gv_0=\Delta_g Q_gv_0\)
is norm continuous. Covariance gives continuity on the span of
the Weyl translates of \(v_0\), which is dense by irreducibility.
Unitarity extends this to every vector. The operators preserve
Schwartz space by \cite[nos.~11--12, Proposition~2, pp.~157--159]{WeilOperators}. The finite integral
stabilizers fix \(v_0\) since \(\Delta_g=1\) there, giving the
restricted tensor product representation.\qedhere
\end{proof}

Let \(B\) be the given Hermitian space and \(\omega_B\)
its Weil action before inserting determinant characters.
Write \(J_B\) for its rational Gram matrix, \(X_B\) for its
point matrix, and \(\Gamma_B=X_B^{-1}E^{\dim B}\) in adelic
coordinates before the local rescalings. A rational unitary element
\(\gamma\) acts in these coordinates by \(X_B^{-1}\gamma X_B\).
For \(x,y\in E^{\dim B}\) the symmetric Weyl displacements have
adelic multiplier
\[
 \beta(X_B^{-1}x,X_B^{-1}y)
   =\psi_E\!\left(j_{\rm W} y^*J_B x/2\right)=1 .
\]
Indeed the finite factors of \(\alpha_{\rm W}\bar\alpha_{\rm W}\)
cancel the two normalized coordinate factors, and the signed
real coordinates give the infinite factors of the same rational
trace character. Moreover
\(\beta(z,X_B^{-1}y)^2=\psi_E(j_{\rm W}y^*J_B X_Bz)\).
As \(y\) varies, the covectors \(j_{\rm W}y^*J_B\) fill
\(E^{1\times\dim B}\). The annihilator of \(E\) in
\(\mathbb A_E\) under \(\psi_E\) is \(E\), so triviality
of all these characters is equivalent to \(z\in\Gamma_B\).
In a rational Schr\"odinger polarization the functional
\(\theta_B\) is a nonzero multiple of rational summation.
Indeed rational modulations force its support to be
on the rational lattice by adelic self-duality. At each
support point the modulation differentials span the
real cotangent space and eliminate derivatives of
point masses. Rational translations make all the
point masses equal. Heisenberg covariance implies
that \(\theta_B\omega_B(\gamma)\) has the same rational
Heisenberg invariance for \(\gamma\in U(B)(F_0)\).
Thus
\[
       \theta_B\omega_B(\gamma)=c_B(\gamma)\theta_B,
\]
where \(c_B\) is a character by Lemma~\ref{lem:weil-vacuum-composition}.
We evaluate it using the
multipliers of the action.

First let \(B\) be split, with rational dual isotropic
bases \(k,\ell\) of rank \(p\). In the mixed model
\(f[x]\) is the theta value after displacement by
\(\ell x\) and averaging the \(k\)-modulations on the
additive quotient with probability Haar. In the
Schr\"odinger realization write \(\langle\xi,u\rangle\) for the
character of the modulation parameter \(u\) at \(\xi\). With
\([E^p]=E^p\backslash\mathbb A_E^p\), rational summation gives
\[
 f[x]=C\sum_{\xi\in E^p} f_{\mathrm{Sch}}(x+\xi)
                  \int_{[E^p]}\langle\xi,u\rangle\,du
     =C f_{\mathrm{Sch}}(x).
\]
For fixed \(x\), compact finite support restricts the sum to a
translate of a fractional lattice, on which the archimedean
Schwartz bound is summable. Character orthogonality then gives
the second equality.
In particular
\[
 \theta_B(f)=\sum_{x\in E^p} f[x],\qquad
 \int\overline{f[x]}g[x]\,dx=|C|^2\langle f,g\rangle .
\]
The original theta functional has mean square one on the
Heisenberg translates of a unit vacuum. Orthogonality of
the rational modulation characters and unfolding the
translation variable give \(|C|^2=1\), as computed in
Section~\ref{sec:overlap-fourier-inversion}.

For the Levi \(m(A)\) acting by \(A\) on \(k\) and
\(A^{-*}\) on \(\ell\), covariance gives
\((\omega_B(m(A))f)[x]=c(A) f[A^*x]\). Composing with
the inverse substitution commutes with all Heisenberg
operators, so its effect is scalar. The Jacobian and
unitarity give \(|c(A)|^2=|\det A|_{\mathbb A_E}\).

For the phase, let \(h\) carry the signed orthonormal frame of the
standard hyperbolic model to the adapted frame in \(B\), and test on
\(\omega_B(h)v_0\). By unitarity its matrix coefficient at \(m(A)\) is
\[
 \langle \omega_B(h)v_0,\omega_B(m(A))\omega_B(h)v_0\rangle
   =\langle v_0,\omega_B(h^{-1}m(A)h)v_0\rangle .
\]
This identity holds for an adelic frame change, including one that moves
the original finite lattice. At infinity the change from the
signed orthonormal basis to the hyperbolic coordinates,
up to the ordering of signs, uses normalized sums and
differences. The negative block of the Levi in the
adapted signed basis is \((A^{-*}+A)/2\). Its conjugate
determinant times \(\det A\) is the positive determinant
of \((I+A^*A)/2\). At finite places an integral change
from the hyperbolic Gram to the identity Gram is given
by the matrix
\[
 \begin{pmatrix}I&-cI\\ \varepsilon I&c\bar\varepsilon I\end{pmatrix},
 \qquad c\bar c=-1,\qquad\varepsilon+\bar\varepsilon=1 .
\]

The norm and trace conditions have integral solutions
with \(c\) a unit. This holds even at the places over two:
at an unramified quadratic place use the residue norm
and trace followed by successive lifts, and at a split
place use the product and sum. In our family \(2\) splits
in \(K\) by Section~\ref{sec:datum}. Since \(E=F_0K\),
for every \(v\mid2\) one has
\(E\otimes_{F_0}F_{0,v}=F_{0,v}\times F_{0,v}\), with the
two factors exchanged by conjugation. The choices
\(c=(1,-1)\), \(\varepsilon=(1,0)\) give
\(c\bar c=-1\) and \(\varepsilon+\bar\varepsilon=1\).
Multiplying by the hyperbolic Gram and its adjoint gives the identity
Gram, and the determinant of the change is \(c^p\).
It preserves the elementary divisors of the lattice
intersection. In these coordinates, with
\(L=\widehat{\mathcal O}_E^p\) and
\(\chi_f=\upsilon_f|\cdot|_f^{-1/2}\), the finite
multiplier of \(A\) is
\[
 \Delta_f(A)=\chi_f(\det H),\qquad
 HL=L\cap A^{-1}L.
\]
Since \(\upsilon_f\) is unramified, changing \(H\) by an integral
basis change leaves \(\Delta_f(A)\) unchanged.
Inversion replaces this intersection by \(AHL\),
transpose preserves its elementary divisors and
conjugation conjugates the multiplier. Hence the
finite multiplier of the pair is
\[
 \Delta_f(A^{-*})\Delta_f(A)
     =\overline{\chi_f(\det A)}|\Delta_f(A)|^2 .
\]
Together the real and finite calculations give
\(\Delta\upsilon(\det A)>0\), where \(\Delta\) is the
multiplier of the Levi written in these adapted
coordinates.

The adapted vacuum has slice a nonzero scalar times a
positive real Gaussian and a finite integral-box
indicator. At infinity this follows by computing the
Heisenberg coefficient of the standard Gaussian using
the orthonormal hyperbolic basis: the cross terms cancel
in the product of the paired coordinates. At finite
places use normalized displacement coordinates, dividing
the original coordinates by \(\alpha_{\rm W}\). The
self-dual lattice is \(L\) in these coordinates by the
different identity of Section~\ref{sec:period-finite-coordinates}.
Indeed let \(w\mid p\), \(\mathfrak D_w=\mathfrak D_{E_w/\mathbf Q_p}\), and
\(a_w=(j_{\rm W}\alpha_{\rm W}\bar\alpha_{\rm W})_w\). In the scalar
trace pairing for the usual character of \(\mathbf Q_p\) of conductor
\(\mathbf Z_p\), the annihilator of \(\mathcal O_{E_w}\) is
\(\mathfrak D_w^{-1}\). Fourier inversion gives this ring volume
\(\lvert\mathfrak D_w\rvert_w^{1/2}\) in the self-dual measure.
Scaling the pairing by \(a_w\) changes the volume to
\(\lvert a_w\mathfrak D_w\rvert_w^{1/2}=1\).
Conjugation identifies the dual coordinate with the opposite component
at a split pair. At \(p=d\) one has \(\alpha_{\rm W}=1\) locally and
\(\mathfrak D_w=(\sqrt{-d})\). At \(p\mid q_1\) the conjugate halves
supply the inverse relative different through
\(\alpha_{\rm W}\bar\alpha_{\rm W}\), also at \(p=2\).
The \(k\)- and \(\ell\)-displacements with normalized
parameters in \(L\) fix the vacuum, by isotropy and
the integral change of basis just computed.

In the mixed realization the modulations force support
in the self-dual box and translations force constancy
there. The product of the vacuum slice and its Levi
substitute therefore has positive integral. By the
slice inner-product formula, the phase of \(c(A)\)
is the phase of the vacuum coefficient
\(1/\bar\Delta\) of the Levi in these coordinates.
Thus
\begin{equation}
 (\omega_B(m(A))f)[x]
 =\upsilon(\det A)^{-1}|\det A|_{\mathbb A_E}^{1/2}
                                     f[A^*x].
 \label{eq:marked-hyperbolic-levi}
\end{equation}

For rational \(A\), the scalar is one by the product
formula and the triviality of \(\upsilon\) on \(E^\times\).
The substitution permutes the rational points of
summation. Each element of either opposite unipotent
radical is a commutator with a rational scalar Levi:
conjugation by the Levi with scalar \(2\) multiplies
the additive radical parameter by \(4\) or \(1/4\).
The character \(c_B\) is therefore one on both radicals.

To see that these elements and the Levi generate, write
the image of \(k\) under a rational isometry as \(ka+\ell b\),
so \(a^*b+b^*a=0\). In a complex embedding
\((a+b)^*(a+b)=a^*a+b^*b>0\) since the image has rank
\(p\). The polynomial \(\det(a+Tb)\) is nonzero.
Choose a skew \(t\in E\) where it is nonzero. One
radical changes the upper block to \(A=a+tb\). The
other, with parameter \(-bA^{-1}\), makes the lower
block zero. This parameter is skew by the isotropy
equation after the first operation. The remaining
element stabilizes \(k\) and belongs to its Siegel
parabolic, generated by the Levi and its radical.
This proves \(c_B=1\) for split \(B\).

For general \(B\), the rational diagonal and antidiagonal
split \(B\perp B^{-}\), where \(B^{-}\) has the opposite
form. In direct-sum coordinates \((\gamma,1)\) acts
on product vectors as \(\omega_B(\gamma)\otimes 1\).
To verify the scalar, the real negative minors multiply,
the finite lattice intersections are direct sums, and
the vacuum is a product. The nonzero vacuum coefficient
is therefore the product coefficient. Covariance and
irreducibility identify the two intertwiners.
Rational summation on the product lattice factors. Since each
theta functional is a nonzero multiple of rational summation,
there is a fixed \(C_\oplus\ne0\) with
\[
 \theta_{B\perp B^-}(f\otimes f^-)
       =C_\oplus\,\theta_B(f)\theta_{B^-}(f^-).
\]
Choose \(f^-\) with \(\theta_{B^-}(f^-)\ne0\). Split invariance for
\((\gamma,1)\) cancels the same factor \(C_\oplus\theta_{B^-}(f^-)\)
on the two sides. This proves
\[
       \theta_B\omega_B(\gamma)=\theta_B
                          \quad(\gamma\in U(B)(F_0)).
\]
For the pair \(S,T\) use \(B=\operatorname{Hom}_E(S,T)\).
The determinant characters in the pair action are one
at rational arguments. In particular rational invariance
holds for the normalized theta kernel after all the
automorphic twists used below.

\subsection{The doubled action in the Siegel--Weil formula}
\label{sec:marked-doubled-action}

\sectionopening{In this subsection, we compute the Levi scalar in the doubled Weil action using the
rational normalization of Section~\ref{sec:marked-rational-theta}.
This scalar determines the normalization of the paired lift in
Section~\ref{sec:marked-first-arrow-rallis}.}

Write \(n=\operatorname{rk}_E S\), \(s=\operatorname{rk}_E T\).
Take the mixed model for a rational Lagrangian of
\(S\perp S^{-}\). For a source Siegel Levi \(m(A)\) and a target
isometry \(h\), the argument of the mixed-model function becomes
\(h^{-1}XA\), with determinant \((\det h)^{-n}(\det A)^s\)
on the space of \(s\)-by-\(n\) matrices.
On the isotropic frame in \eqref{eq:marked-hyperbolic-levi},
the action in target coordinates is \(Y\mapsto hYA^*\).
Indeed, \(h\) acts on the output and precomposition by
\(m(A)^{-1}\) acts by \(A^*\) on this frame of source covectors.
Thus the frame block has determinant
\(d=(\det h)^n\overline{\det A}^{\,s}\). Since \(\det h\)
has norm one and \(\upsilon(\bar a)=\upsilon(a)^{-1}\), its scalar is
\[
  \upsilon(d)^{-1}|d|_{\mathbb A_E}^{1/2}
   =\upsilon(\det h)^{-n}\upsilon(\det A)^s
                       |\det A|_{\mathbb A_E}^{s/2}.
\]
The target character \(\upsilon(\det h)^n\) of the doubled pair
cancels the first factor.

The Siegel parabolic preserves the line of evaluation at zero.
Indeed it preserves the modulation subgroup, whose invariant
distributions in this model are the multiples of evaluation at zero.
Its action on this line is a character. Conjugation by the scalar
Levi \(2I\) scales the upper radical parameter by \(4\) or \(1/4\),
so the radical consists of commutators and acts trivially there.
Write \(u(Z)\) for the radical element
\(\ell\mapsto\ell-k j_{\rm W}^{-1}Z\) with \(Z=Z^*\).
Moving the displacement by \(X\) through \(u(Z)\)
adds a \(k\)-displacement absorbed by the average
and gives
\[
 (\omega(u(Z))f)[X]
 =\psi_E\!\left(\tfrac12\operatorname{tr}(Z X^*J_T X)\right)f[X].
\]
Here \(\psi_E=\psi_{F_0}\circ\operatorname{Tr}_{E/F_0}\)
for the rational trace characters. The matrix trace inside the phase is fixed by conjugation.
To identify the unipotent coordinate, convert Yamana's
matrix convention by \(g\mapsto(g^*)^{-1}\).
The ordered basis \((2j_{\rm W}\ell,k)\) has Hermitian
Gram matrix
\(-2j_{\rm W}\left(\begin{smallmatrix}0&I\\-I&0\end{smallmatrix}\right)\).

Removing the scalar gives the skew-Hermitian matrix of
the classical model with the same isometry group. The
upper unipotent \(n(b)\), for \(b=b^*\), sends
\(2j_{\rm W}\ell\) to \(2j_{\rm W}\ell-kb\) and fixes \(k\).
Thus \(Z=b/2\). The mixed variable on \(k\) is unchanged and
\[
 \psi_E\!\left(\tfrac12\operatorname{tr}(Z X^*J_T X)\right)
    =\psi_{F_0}\!\left(\operatorname{tr}(b X^*J_T X/2)\right).
\]
This is the unipotent law in \cite[\S1.3]{YamanaSingular}.
The rational change of basis preserves probability
Haar on the rational additive quotient. The dual Fourier
measures transport with the pairing. In the opposite
polarization the law follows by self-dual Fourier transform,
with its scalar canceled by conjugation.

These are the Schr\"odinger laws of the splitting
character \(\upsilon^s\) with the additive
character. The two opposite parabolics generate over each local
field. For a nonsplit place the block proof of
Section~\ref{sec:marked-rational-theta} works after a
Siegel Levi move and a change of frame. Choose invertible \(P,Q\)
with \(b'=PbQ=\operatorname{diag}(I_j,0)\), and put
\(a'=P^{-*}aQ\). Isotropy gives \(a'_{12}=0\), and joint
injectivity gives invertible \(a'_{22}\). Thus
\[
 \det(a'+Tb')=\det(a'_{11}+TI_j)\det(a'_{22})\ne0
             \quad\text{in }E_v[T].
\]
Choose a skew scalar \(t\) outside its finite set of roots.
In the original frame the required upper radical parameter
is \(N_t=tP^*P\), which is skew, and
\(a'+tb'=P^{-*}(a+N_tb)Q\). Put \(A=a+N_tb\).
Isotropy gives \(A^*b+b^*A=0\), so \(-bA^{-1}\) is skew.
The lower radical with this parameter sends
\(\binom{A}{b}\) to \(\binom{A}{0}\), placing the resulting
isometry in the upper parabolic. At a split place the two
block radicals and their common Levi contain the elementary
transvections and invertible diagonal matrices of
\(\mathrm{GL}_{2n}(F_{0,v})\). The actions
agree on the integral vacua outside a finite set.

At infinity the Fock polynomial space in \cite[pp.~236--237]{YamanaSW}
is relative to the chosen positive and negative orthogonal subspaces.
The corresponding signed decompositions of the doubled source and
of the target can each be carried to the original ones by an isometry:
map signed orthonormal bases to signed orthonormal bases. Write
\(g_0,h_0\) for these real isometries. On the Hom space their action
carries the complex coordinates, the Gaussian and its polynomial
multiples to those of the transported Fock model. Thus the original
tensor is \(\omega(g_0,h_0)\Phi_0\) for a polynomial vector
\(\Phi_0\) in the cited model. The cancellation of
\(\upsilon(\det h_0)^{-n}\) by the doubled target character gives
\((\omega(1,h_0)\Phi)(X)=\Phi(h_0^{-1}X)\). This fixes evaluation
at zero and is absorbed by right translation in the target quotient
integral. By real Iwasawa decomposition one may
take \(g_0\in P_\infty\), the Siegel parabolic, absorbing its
right compact factor into \(\Phi_0\).
Let \(K\) be the maximal compact of the doubled source in
Yamana's model, with the given finite factors, and put
\(K'=g_0Kg_0^{-1}\). Writing \(u=pk'\) gives
\(|a_{K'}(u)|=|a_K(ug_0)|/|a_K(g_0)|\).
With \(s_0=(s-n)/2\), the sections satisfy
\[
 f^{(z),K'}_{\omega(g_0,h_0)\Phi_0}(u)
 = |a_K(g_0)|^{s_0-z} f^{(z),K}_{\Phi_0}(ug_0).
\]
The factor is one at \(z=s_0\). This transports the polynomial
identity also at a holomorphic value obtained by continuation.

In \cite[Theorem~2]{YamanaSW} the division algebra is
\(E\) of index one. The splitting character restricts
to \(\chi_{E/F_0}^{\,s}\). At \(s_0=(s-n)/2\) the normalized
Siegel section is \(g\mapsto(\omega(g,1)f)[0]\), and the theta
integrand is
\[
 \sum_{X\in\operatorname{Hom}_E(E^n,T)}(\omega(g,h)f)[X].
\]
The zero value and the rational sum use the same mixed realization,
so a fixed scalar in this realization multiplies both sides of
the Siegel--Weil identity. In every strict-range
doubling below \(s>2n\), so
\[
         s-\operatorname{Witt}(T)\ge s/2>n .
\]
The ordinary theta integral converges absolutely for every
Schwartz--Bruhat test by \cite[\S51, Proposition~8, p.~75]{WeilSiegel}.
The spherical Eisenstein majorant converges since \(s_0>n/2\).
In Weil's case \({\rm I}_2\), the parameter is \(1/2\), so
his dimension condition is \(s>2n\) \cite[\S\S26--27]{WeilSiegel}.

At an old split place \(v\), put \(k=F_{0,v}\). A full-rank tuple
consists of an injection \(X:k^n\to k^s\) and a surjection
\(Y:k^s\to k^n\), with Gram matrix \(B=YX\) in one split factor.
In \(k^s=k^n\oplus k^n\oplus k^{s-2n}\), take
\[
 X(a)=(a,0,0),\qquad Y(a,b,c)=Ba+b .
\]
These maps realize every \(B\), including \(B=0\), with \(X\)
injective and \(Y\) surjective. For two pairs \((X_i,Y_i)\) with the same
\(B\), match \(X_1a\) with \(X_2a\). Since
\[
  \operatorname{im}X_i\cap\ker Y_i=X_i(\ker B),
\]
extend this matching on the intersections to an isomorphism
\(\ker Y_1\to\ker Y_2\), hence to the sums
\(\operatorname{im}X_i+\ker Y_i\). Complete it by matching lifts
of the same basis of a complement of \(\operatorname{im} B\) in
\(k^n\) under \(Y_1,Y_2\). The resulting \(g\in\mathrm{GL}_s(k)\)
satisfies \(gX_1=X_2\) and \(Y_2g=Y_1\). Thus \(U(T)(k)\) acts
transitively on the full-rank part of every Gram fibre and preserves
the probability Haar measure on the adelic quotient.

For polynomial Fock vectors use \cite[Theorem~2]{YamanaSW}.
For arbitrary real Schwartz vectors use
\cite[\S52, Theorem~5, pp.~76--77]{WeilSiegel}. Both formulas give coefficient
one with probability quotient measure. In
Section~\ref{sec:rank-two-doubling-local-integrals}, the doubled
polynomial Fock datum has definite Hermitian target and the target
quotient is compact. For \(m=n=2\) the holomorphic value is
\(E(g,0)=2I(g)\) with the ordinary theta integral and probability
target measure \cite[Theorem~2, p.~237]{YamanaSW}
\cite[\S1.3, p.~1007, and Proposition~5.8(i), p.~1015]{YamanaSingular}.

\subsection{Translation and scalar normalization in the seesaw}
\label{sec:seesaw-normalization}

\sectionopening{In this subsection, we fix the normalization of the theta operators
used in the seesaw. The source index is the one in
Section~\ref{sec:even-stage-coefficients}, and equivariance carries
the finite group action through the entire composition.}

The lower chain is a composition of the theta
operators \(L_{\xi,j}\) from \(Y_j\) to \(Y_{j+1}\)
with the tame character at the \(j\)-th arrow,
followed by the terminal algebraic coefficient map.
Each \(L_{\xi,j}\) is the original theta integral multiplied
by the source index \(A_{\xi,j}=[K_j^\circ:J_{\xi,j}]\),
where \(J_{\xi,j}\) is the principal subgroup at the
old and twist places. In particular, split reduction gives
\[
 A_{\xi,j}=\bigl|\mathrm{GL}_{r+jt}
                   (\mathcal O_E/(\pi_M Q(\xi)))\bigr|.
\]
Section~\ref{sec:even-stage-coefficients} computes the
transpose with this very index.
Here \(Y_m\) is the lower-chain endpoint, of a different
rank from the dimension-comparison space \(W_m\).
Use the lattice inclusion \(\iota_j\) between the
standard finite coordinates of \(Y_j,Y_{j+1}\).
At the original level places and the additional
twist places, the lattice test has the inclusion
reduction \(\bar\iota_j\). The finite twist factor of
the Weil action on an integral pair is the quotient
of the tame character of its two determinants.

\begin{lemma}[Equivariance of the weighted lower lifts]\label{lem:weighted-lower-equivariance}
For each \(s\in G\) there are integral lifts \(k_j(s)\) with
first-block residue \(\operatorname{diag}(s,I)\) at the old places
and identity at the twist places. They satisfy
\[
 L_{\xi,j}\rho_j(k_j(s))
            =\rho_{j+1}(k_{j+1}(s))L_{\xi,j}.
\]
The successive intertwinings give the action of the original
finite projector on the terminal values.
\end{lemma}

\begin{proof}

For \(s\in G\), choose integral lifts \(k_j(s)\) on
each lower space with reduction \({\rm diag}(s,I)\)
in the first-block coordinates at the original
places and identity at the twist places.
Such lifts exist because these places are split.
At each step
\[
 \overline{k_{j+1}(s)}\,\bar\iota_j
                  =\bar\iota_j\,\overline{k_j(s)}.
\]
The simultaneous translate of the lattice test therefore
preserves its support: integral matrices are carried to
integral matrices, and the inclusion reduction is fixed.

The two reduced determinants agree, since the extra diagonal
block is the identity. Their quotient therefore has trivial
tame character. Away from these places
take the pair to be the identity. The accompanying
real test is unchanged. In the theta kernel, apply the
finite source translation by changing the source integration
variable. The Haar measure is right invariant and the
simultaneous action on the datum has just been shown to
fix it.

Thus
\[
 L_{\xi,j}\rho_j(k_j(s))
            =\rho_{j+1}(k_{j+1}(s))L_{\xi,j}.
\]
The principal kernel \(J_{\xi,j+1}\) in the target variable
fixes the inclusion test and its tame determinant character, as in
Section~\ref{sec:even-stage-coefficients}. Thus the image of
\(L_{\xi,j}\) is fixed by \(J_{\xi,j+1}\). Moreover, for \(s,t\in G\),
\[
 k_j(ts)^{-1}k_j(t)k_j(s)\in J_{\xi,j}.
\]
Together with principal-kernel invariance this gives the action
of \(G\) on the translated coefficients independently of the
integral lifts. The conversion to the algebraic coefficient at
the definite endpoint commutes with this finite action.
Write \(L_\xi\) for the resulting composed lower chain.
Expanding the finite-type projector \(p=\rho(e_0)\)
in the original group algebra gives the identity of
terminal values
\[
 \bigl(L_\xi(|G|p\rho(s)\Pi\phi)\bigr)(h)
  =\sum_{t\in G}|G|e_0(t)
                       (L_\xi\Pi\phi)(h\,k_m(ts)).
\]
Indeed \(k_m(t)k_m(s)\) and \(k_m(ts)\) have the same reduction,
so they act identically on the right \(J_{\xi,m}\)-invariant
endpoint. This uses the same measured \(\Pi\) for every translate.
\end{proof}

Fix \(c_0,C_0\) independently of \(\xi\).
Suppose the unfolded pairing gives, on the original source, \(L_\xi\Pi\phi
=c_0\,\mathcal U_\xi(T_\xi)\phi\), where
\(\mathcal U_\xi\) is the upper chain with its terminal
restriction and is linear in its tensor. Suppose the
tensor is a finite combination of translated
vacua \(T_\xi=C_0\sum_z b_{\xi,z}v_z\), and put \(B_*=M^{f n_m^2}\) for terminal upper
rank \(n_m\). If \(B_*^{2m}b_{\xi,z}\) is integral
off the twist, and the untwisted sum of coefficient
moduli in each complex embedding is at most \(B_*^m\),
then the coefficients \(B_*^{2m}b_{\xi,z}\) give the
integral tensor with untwisted bound \(B_*^{3m}\).

The common scalar multiplying that tensor is
\(c_0 C_0/B_*^{2m}\) when this is nonzero.
If \(c_0 C_0=0\), use the zero integral tensor and scalar
one. In either case the scalar is nonzero and independent
of the twist. Support and invariance under the upper
levels are inherited from the same evaluated tensor,
provided its translated vacua have the prescribed support.
The unfolded pairing is \eqref{eq:seesaw-measured-chain}.
Section~\ref{sec:seesaw-integral-tensor} supplies the tensor
expansion and its coefficient bounds.

\subsection{Rational orbits for consecutive lifts}
\label{sec:adjacent-string-orbits}

\sectionopening{In this subsection, we describe the rational maps that contribute to
the Bessel integral for two consecutive Hermitian spaces. The
character of Section~\ref{sec:common-filling-assembly} determines
these maps up to the Bessel actions, which are used in
Section~\ref{sec:adjacent-unfolding}.}

Use the rational graded coordinates of the two consecutive spaces
\(Y,Y'\), with Hermitian spaces \(V,V'\) and string lengths
\(s,\max(0,s-1)\), respectively. For \(s\ge1\) write \(a_i=s-1-2i\) and
\(a'_b=s-2-2b\) for the string grades, starting at index zero.
The orthogonal definite blocks have grade zero. The Hermitian
forms are conjugate-linear in the first argument. On a
string, \(i^\vee\) denotes the reversed index and the Gram
matrix has entries \(g_i=h(e_i,e_{i^\vee})=1\), except at a
middle index, where \(g_i=-N\). The definite orthonormal
coordinates are their own dual indices. A matrix has
degree at least \(d\) if its entry is zero when the target
grade minus the source grade is less than \(d\).

Let \(c\in F_0^\times\) be the string constant
and \(T:V\to V'\) the degree-\(-1\) map which kills
the definite block and the last string vector. Its
other string entries are \(Te_i=t_i e_i'\), with
\[
 t_i=
 \begin{cases}
  1, -Nc, c & s\text{ odd, in grades }>0, =0, <0,\\
  c, 1       & s\text{ even, in grades }>0, <0.
 \end{cases}
\]
Taking the adjoint in the antidiagonal Gram matrices
gives \(T^\dagger T e_i=t_i t_{s-2-i}(g'_i/g_{i+1})e_{i+1}\)
for \(i<s-1\). This coefficient is \(c\), except
at an odd string's grade zero, where it is \(-Nc\),
and at an even string's grade one, where it is \(-Nc^2\).
In particular all the shift coefficients are nonzero.
Write \(R,R'\) for the Bessel subgroups over \(F_0\).

Every element of \(R\) fixes each nonnegative-grade string basis
vector modulo vectors of strictly higher grade. The same property
holds for \(R'\) on its string. The matrices of these elements
and their inverses have degree at least zero. Indeed the stated
column property preserves the positive flag and, by
orthogonality, the corresponding negative flag.
On the grade-zero quotient the orthogonal complement
of the central string vector, when present, can carry
the definite unitary block. Write \(U_R\subset R\) for
the unipotent group of matrices whose difference from the
identity has degree at least one.

Say that \(w\) agrees with \(T\) on the positive columns
if \(w_{bi}=T_{bi}\) for \(a_i>0\) and \(a'_b<a_i\).
If \(w\) agrees there and \(\gamma\in R(\mathbb A_{F_0})\), then \(w\gamma\)
also agrees: the positive columns are stable under
\(\gamma\), and \(T(\gamma-1)\) has degree at least zero
on them. For \(\delta\in U(V')(\mathbb A_{F_0})\),
\begin{equation}
 \delta w\gamma\text{ agrees with } T
             \quad\Longleftrightarrow\quad \delta\in R'.
 \label{eq:adjacent-agreement-stabilizer}
\end{equation}
For the forward implication put \(\widetilde w=w\gamma\).
In a nonnegative string column \(k\) of \(V'\), use the
positive source pivot \(i\) with \(a_i=a'_k+1\) and \(Te_i=t_i e'_k\).
For every index \(b\) with \(a'_b\le a'_k\), the entry
\(((\delta-1)\widetilde w)_{bi}\) is zero by agreement.

Proceeding in descending order of \(a'_k\) gives
\(0=t_i(\delta-1)_{bk}\):
in column \(i\), all entries in rows of grade at most \(a'_k\)
are already those of \(t_i e'_k\), and entries of
\(\delta-1\) in higher columns vanish in row \(b\) by
induction. These are precisely the defining conditions for
\(R'\). Conversely, if \(\delta\in R'\), then
\((\delta-1)T\) has degree at least zero on the positive
columns, as does \(\delta(\widetilde w-T)\).
This proves \eqref{eq:adjacent-agreement-stabilizer}.

Use the additive character \(\psi_E\) of \(\mathbb A_E/E\)
from the absolute trace and the skew element \(j_{\rm W}\)
of Section~\ref{sec:common-filling-assembly}.
Set
\[
 \beta(x,y)=\psi_E\bigl((j_{\rm W}/2)\operatorname{tr}(y^\dagger x)\bigr).
\]
For target grades \(D=0\) or \(D=a'\), the phase
test at \(w\) requires
\begin{equation}
 \beta(Tu,T)=\beta(wu,w)
 \quad\text{if}\quad
 u\in U_R(\mathbb A_{F_0}),\quad
 \deg(u-1)\ge2,\quad \deg_D(w(u-1))\ge1 .
 \label{eq:adjacent-phase-test}
\end{equation}
Here \(D=0\) assigns grade zero to every target coordinate.
The degree-two condition is automatic for a strictly
grade-raising matrix when \(s\) is odd.

Suppose \(w\) passes \eqref{eq:adjacent-phase-test}. Fix
a source index \(k\) of grade \(d\) and an adelic vector
\(v\) supported in grades at least \(\max(1,d+2)\).
Suppose \(wv\) is supported in target grades
\(D>\max(-1,d)\) and \(we_{k^\vee}\) is supported in
\(D>-a_i\) for every \(i\) in the support of \(v\).
\begin{lemma}[Eichler moment on the supported grades]\label{lem:adjacent-eichler-moment}
Under the stated grade hypotheses and the phase test
\eqref{eq:adjacent-phase-test},
\begin{equation}
 [(T^\dagger T)v]_k=[(w^\dagger w)v]_k .
 \label{eq:adjacent-eichler-moment}
\end{equation}
\end{lemma}

\begin{proof}
Put \(x=\lambda e_{k^\vee}\) and
\(L=v h(x,\cdot)-x h(v,\cdot)\) for arbitrary
\(\lambda\in\mathbb A_E\). Since \(h(x,v)=h(v,v)=0\),
\[
 L^\dagger=-L,\qquad
 L^2=-h(x,x)v h(v,\cdot),\qquad
 L^3=0.
\]
Thus \(u=1+L+L^2/2\) is unitary. The three terms of
\(u-1\) raise the source grade by at least two, and their
images under \(w\) have \(D\)-degree at least one by
the support conditions. For \(A=y^\dagger y\),
\[
 \operatorname{tr}(Au)=\operatorname{tr} A+h(x,Av)-h(v,Ax)
                               -\tfrac12 h(x,x)h(v,Av).
\]
The first and last terms are invariant under conjugation,
and \(h(v,Ax)=\overline{h(x,Av)}\). Since
\(\overline{j_{\rm W}}=-j_{\rm W}\) and \(\psi_E\) is
invariant under conjugation, the phase is
\(\psi_E(j_{\rm W}h(x,Av))\).
Apply \eqref{eq:adjacent-phase-test} for \(y=T,w\).
Since \(h(x,Av)=\bar\lambda g_{k^\vee}(Av)_k\) and
\(j_{\rm W}g_{k^\vee}\ne0\), nondegeneracy of the adelic trace
pairing as \(\lambda\) varies gives
\eqref{eq:adjacent-eichler-moment}.
\end{proof}

Take first \(D=0\) and \(w\) rational. For positive string
indices \(i,l\), use \(v=e_i, k=l^\vee\) in
\eqref{eq:adjacent-eichler-moment}. The resulting identity
is \(h'(Te_l,Te_i)=h'(we_l,we_i)\). On the span of
these columns \(w\) is injective. In fact if \(wv=0\)
for a nonzero vector supported there, let \(i\) be its
least-grade occupied string index and use \(k=i+1\).
The support conditions hold even if \(a_{i+1}\ge0\),
because \(wv=0\). The \(k\)-th coordinate of
\(T^\dagger Tv\) is a nonzero multiple of \(v_i\),
contradicting \eqref{eq:adjacent-eichler-moment}.
The same positive-column space injects under \(T\).
Apply Witt extension \cite[Ch.~7, Theorem~9.1, p.~265]{Scharlau}
to the regular quadratic space \((V',[h'/2])\) over \(E\)
with its CM involution. Symmetrization gives \(h'\), and
division by two identifies Hermitian forms with their quadratic
classes in characteristic zero. Thus the \(E\)-linear isometry
\(we_i\mapsto Te_i\) on the positive-column span extends to
an element of \(U(V')(F_0)\), including when that span is degenerate. In particular the
rational orbit contains a matrix agreeing with \(T\)
on the positive columns. By
\eqref{eq:adjacent-agreement-stabilizer}, the rational
target elements giving agreement have the form
\(R'(F_0)\delta_0\), also after a source change by \(R(F_0)\).

For the second test assume this agreement and take
\(D=a'\). We show that \(w-T\) has degree at least zero,
by induction in increasing target space grade. A possibly
different entry \((b,i)\) has \(a'_b<a_i\le0\).
Use the positive source pivot \(v=e_j\) with
\(Te_j=t_j e'_{b^\vee}\) and \(a_j=1-a'_b\), and
put \(k=i^\vee\), of grade \(-a_i\). Then
\(a_j\ge2-a_i\). The vector \(we_j\) is \(Te_j\)
plus terms of grade at least \(1-a'_b\), by agreement.

The entries of \(we_i\) of grade at most \(-a_j=a'_b-1\)
vanish by induction and by the degree of \(T\).
These are the conditions for
\eqref{eq:adjacent-eichler-moment}, which now gives
\[
 h'(we_i,we_j)=h'(Te_i,Te_j).
\]
The extra terms of \(we_j\) pair only with the
already matched coordinates of \(we_i\) of grade at most
\(a'_b-1\). Their pairing there is zero by the same
degree bound on \(T\). Hence the remaining equation is
\(\overline{(w-T)_{bi}}\,t_j g'_b=0\), and the induction
closes. This includes the definite grade-zero source
columns. In particular all surviving rational terms of
the second test are in the single class of \(T\) modulo
the rational matrices of degree at least zero.

\section{Consecutive theta lifts and their coefficients}
\label{sec:consecutive-theta-lifts}

\sectionopening{In this section, we unfold consecutive theta integrals using the rational orbits of
Section~\ref{sec:adjacent-string-orbits}. The even-step coefficient formula
and the finite support determine the sum evaluated in
Section~\ref{sec:measured-theta-tensor}.}

\subsection{Unfolding the integral for consecutive lifts}
\label{sec:adjacent-unfolding}

\sectionopening{In this subsection, we unfold the adjacent theta integral by the orbit
computation of Section~\ref{sec:adjacent-string-orbits}. The Haar
measures on the original quotients determine the constant in the
resulting expression.}

Use the datum of Section~\ref{sec:adjacent-string-orbits}.
For \(A=\mathbb A_{F_0}\), give the two quotients
\[
   [R]=R(F_0)\backslash R(A),\qquad
   [G']=U(V')(F_0)\backslash U(V')(A)
\]
their inherited Haar measures.
Write \(\rho(w)\) for the Heisenberg translate on
\(\operatorname{Hom}_E(V,V')\), with
\(\rho(w)\rho(z)=\beta(w,z)\rho(w+z)\).
The theta functional \(\theta\) is invariant under rational
Heisenberg translations by reindexing and triviality of the additive
character on rational arguments. For the canonical rational lift,
Weil gives invariance of the theta sum in
\cite[no.~41, (39) and Theorem~6, p.~193]{WeilOperators}.
Section~\ref{sec:marked-rational-theta} verifies the compatibility
of the splitting used here with that rational action.
Use the chosen adelic splitting characters for the pair \(U(V),U(V')\).
For a rational unitary pair its two determinants are rational
norm-one elements, so the full automorphic characters take value
one on them. On each unipotent group the determinant is one
at every place.

For \(D=0,a'\) let \(I_D\) be the space of strictly
positive-degree maps and \(\mathcal C_D\) the rational
maps modulo the rational maps of nonnegative degree.
Reversal of both indices exchanges the degrees \(d,-d\),
so \(I_D\) is isotropic and has as rational dual basis the
opposite coordinates with the inverse Gram entries.
The arrow test is Schwartz--Bruhat in rational
Schr\"odinger coordinates. At infinity the creation
polynomial of the joint harmonic acts on the
Gaussian vacuum.

At the finite places the vacuum is
fixed by a small open lattice of Heisenberg translations
and modulations. Translation invariance makes it locally
constant, and modulation invariance cuts the support to
the compact annihilator lattice. The inclusion test
is a finite sum of translates of that vacuum.
A fixed outer Weil translate preserves this class.
Use the same theta normalization in the two gradings.

For a Schwartz test \(f\) in this realization put
\[
 C_D(f;w)=\int_{I_D(F_0)\backslash I_D(A)}
                      \theta\bigl(\rho(z)\rho(w)f\bigr)\,dz ,
 \qquad \int dz=1 .
\]
For rational \(w\), commutation and invariance give the
integrand \(\beta(z,w)^2\theta(\rho(z)f)\).
Thus \(C_D(f;w)\) depends only on the class of \(w\) in
\(\mathcal C_D\). These are the Fourier coefficients
at the rational dual frequencies of \(z\mapsto
\theta(\rho(z)f)\), with the frequency sign specified
by the displayed integrand.

In particular
\[
                 \theta(f)=\sum_{w\in\mathcal C_D} C_D(f;w).
\]
The sum is absolute. In rational Schr\"odinger
coordinates extending the positive block and its dual basis, the \(I_D\) variable acts by a character on the
first Schr\"odinger coordinate. Orthogonality cuts the
theta sum to one such rational coordinate \(\xi\), leaving
the sum of the test over all rational coordinates in
the remaining Lagrangian. The sum of the coefficient
moduli is bounded by the sum of the modulus of the test
over the full rational Schr\"odinger lattice, times the
fixed modulus of the theta normalization.

The finite test
has compact support, so the relevant rational coordinates
lie in a fixed fractional ideal to a finite power.
Multiply an archimedean Schwartz supremum by
\(\sum_\lambda(1+\|\lambda\|)^{-a}\), with \(a\) greater
than the real dimension. This gives a bound uniform on any
bounded set of archimedean Schwartz tests and bounded
finite tests in a fixed finite-dimensional space of
compactly supported tests.

For \(z_0\in I_D(A)\), translation by \(z_0\) fixes
the averaged theta functional. Hence
\begin{equation}
 C_D(f;w)=\beta(w+z_0,w)C_D(f;w+z_0).
 \label{eq:adjacent-positive-shift}
\end{equation}
For a source rational Bessel element \(\gamma\) and
a rational target element \(\delta\) preserving the
\(D\)-filtration, covariance gives
\begin{equation}
 C_D(\omega(\gamma r,\delta h)\Phi;w)
       =C_D(\omega(r,h)\Phi;\delta^{-1}w\gamma).
 \label{eq:adjacent-rational-ct}
\end{equation}
This is rational theta invariance and change of variable
on the compact quotient of \(I_D\). The change preserves
Haar probability. The same covariance holds for the adelic source
\(u\in U_R(A)\) with \(u-1\) strictly raising the source grade,
and for a strictly \(D\)-grade-raising target unipotent.
Write \(W=\bigoplus_j W_j\) for the graded Hom space.
The form pairs \(W_j\) only with \(W_{-j}\), so
\(I_D=W_{>0}\) and \(I_D^\perp/I_D=W_0\).

Put \(E_D(f)=C_D(f;0)\). In partial Schr\"odinger coordinates this
is the residual theta distribution on the slice at zero.
An \(I_D(A)\)-modulation multiplies a test by an additive
character of the transverse coordinate. At infinity, differentiating
this invariance makes the distribution vanish on multiples of each
real transverse coordinate, so it factors through evaluation at zero.
At the finite places, a test vanishing on the zero slice is covered
by finitely many compact open pieces on each of which a modulation
is a constant different from one. These pieces have zero evaluation.
Invariance under rational Heisenberg translations in \(W_0\)
then makes the residual distribution a multiple of \(E_D\)
by the rational uniqueness in Section~\ref{sec:marked-rational-theta}.
Use the Weil representation of \(U(V)(A)\times U(V')(A)\)
from that section. On the subgroup preserving \(I_D\) and
acting as the identity on \(I_D^\perp/I_D\), covariance gives
these two invariances for \(E_D\circ\omega(g)\). In the residual
translation the added component in \(I_D(A)\) pairs trivially
with \(W_0\). Since \(E_D\ne0\), the scalar in
\(E_D\circ\omega(g)=c(g)E_D\) is a character of this subgroup.

Let \(D_b\) be the target grade of coordinate \(b\) for the chosen
\(D\). Take the rational pair
\[
 d=\bigl(\operatorname{diag}(2^{a_i}),
         \operatorname{diag}(2^{D_b})\bigr).
\]
Opposite grades pair, so the two diagonal matrices are unitary.
On the Hom entry \((b,i)\) the pair acts by \(2^{D_b-a_i}\),
and hence fixes \(W_0\). In each unitary factor let \(N_k\)
consist of the elements whose difference from the identity raises
its own grades by at least \(k\). Both factors' \(N_1\) preserve
\(I_D\) and act trivially on the quotient. For \(g\in N_k\),
\(k\ge1\), put
\[
 h=\exp\bigl((2^k-1)^{-1}\log g\bigr).
\]
The series are finite and \((\log g)^\dagger=-\log g\), so
\(h\) lies in \(N_k\) of the same factor.
Conjugation by \(d\) multiplies the degree-\(k\) part by
\(2^k\). It follows that \([d,h]\equiv g\pmod{N_{k+1}}\).
The character is one on the commutator. Since the filtration
terminates, descending induction gives \(c|_{N_1}=1\) in each
factor. This applies to the two unipotent actions, whose
determinant characters are trivial. Moving \(\rho(w)\) through
their Weil actions gives the stated covariance.

\subsubsection{Phase cancellation and rational orbits}\label{sec:paired-lifts-rational-orbits}

\sectionopening{In this subsection, we use Lemma~\ref{lem:adjacent-eichler-moment} to identify the terms that survive the phase average.}

Let \(W\) be continuous on \([R]\) with
\(W(ur)=\beta(Tu,T)W(r)\) for \(u-1\) of degree at least
two. By \eqref{eq:adjacent-positive-shift} the product
\(F_w(r,h)=W(r)C_D(\omega(r,h)\Phi;w)\) satisfies
\[
 F_w(ur,h)=\frac{\beta(Tu,T)}{\beta(wu,w)}F_w(r,h)
 \qquad\bigl(\deg(u-1)\ge2,\quad\deg_D(w(u-1))\ge1\bigr).
\]
For a class that fails \eqref{eq:adjacent-phase-test},
unfold its rational source orbit. If \(\Gamma_w\) is
its stabilizer, the orbit sum over \([R]\) is the
single integral over \(\Gamma_w\backslash R(A)\), with
quotient measures from counting measure. The sum is
absolutely integrable for each \(h\). Indeed \(R\) is the
semidirect product of its definite unitary Levi and
unipotent radical, so \([R]\) is compact. For each
\(r\) choose a rational left translate into a fixed
compact cover. Equation~\eqref{eq:adjacent-rational-ct}
permutes the classes, and the Schwartz bound applies
to the compact family of translated tests.

Put \(J_D=I_D^\perp\), the space of maps of nonnegative
\(D\)-degree. Thus \(\Gamma_w=\{\gamma\in R(F_0):w\gamma-w\in J_D(F_0)\}\).
Inside the source group let
\[
 H_w=\{u\in U_R:\deg(u-1)\ge2,\quad w(u-1)\in I_D\}.
\]
It is a closed unipotent \(F_0\)-subgroup. Indeed
\(w(uv-1)=w(u-1)v+w(v-1)\), and right composition by
\(v^{\pm1}\) preserves \(I_D\). For \(\gamma\in\Gamma_w\),
\[
 w(\gamma u\gamma^{-1}-1)
   =\{w(u-1)+(w\gamma-w)(u-1)\}\gamma^{-1}\in I_D .
\]
The middle product has \(D\)-degree at least two.
Conjugation also preserves the source degree bound, so \(\Gamma_w\)
normalizes \(H_w\), and \(H_w(A)\cap\Gamma_w=H_w(F_0)\).

Let \(L=\Gamma_w H_w(A)\). In the algebraic normalizer
\(N_R(H_w)(A)\), this is the inverse image of the image of
\(\Gamma_w\) in \((N_R(H_w)/H_w)(A)\). That image is closed
and discrete, since it consists of rational points in an affine
algebraic group. Thus \(L\) is closed. Conjugation by
\(\gamma\in\Gamma_w\) preserves adelic Haar on \(H_w(A)\):
its modulus is the product of the local absolute values of
\(\det(\operatorname{Ad}\gamma|_{\operatorname{Lie}H_w})\in F_0^\times\),
which is one. Hence \(L\) is unimodular, and the invariant
measure on \(\Gamma_w\backslash L\) is the measure on the compact
quotient \(H_w(F_0)\backslash H_w(A)\).
The modular function of \(R(A)\) is trivial on \(L\): for
\(H_w(A)\) the adjoint action is unipotent, and for \(\Gamma_w\)
the product formula applies to the rational adjoint determinant on
\(\operatorname{Lie}R\). Thus the Haar quotient formula applies both
to \(\Gamma_w\backslash L\) and to \(L\backslash R(A)\).

For \(v\in H_w(A)\), the multiplier of \(F_w(vr,h)\) is
\[
 \chi(v)=\frac{\beta(Tv,T)}{\beta(wv,w)} .
\]
If the fibre value is nonzero, two successive translations show
that \(\chi\) is a character, trivial on \(H_w(F_0)\).
A failed phase test gives \(u\in H_w(A)\) with \(\chi(u)\ne1\).
Right translation by \(u\) preserves the fibre measure and multiplies
its integral by \(\chi(u)\), so that integral vanishes. Quotient
integration along \(\Gamma_w\backslash L\) and \(L\backslash R(A)\)
now gives zero for the original orbit, using the absolute integrability
proved above.

For \(D=0\), the first phase test and the rational
orbit calculation show that only the translates of
the classes agreeing with \(T\) on positive columns
contribute. Agreement is a condition on entries with
\(a_i>\max(0,a'_b)\), so it is well defined in both
quotients. For a matrix \(v:V'\to V\), put
\(y_{ib}=v_{ib}\mathbf 1_{a_i>\max(0,a'_b)}\).
The sum of the agreeing terms, for either grading, is
\begin{equation}
 \sum_{w\ {\rm agreeing}} C_D(f;w)=
 \int_{\operatorname{Hom}(V',V)(F_0)\backslash
                  \operatorname{Hom}(V',V)(A)}
 \frac{\theta(\rho(y^\dagger)f)}
      {\psi_E(j_{\rm W}\operatorname{tr}(yT))}\,dv ,
 \qquad \int dv=1 .
 \label{eq:adjacent-common-projection}
\end{equation}
Adjunction reverses the two indices, so \(z_0=y^\dagger\in I_D(A)\).
Translate the variable of the probability integral defining \(C_D\).
Commuting \(\rho(w)\) past \(\rho(z_0)\) gives
\(\beta(w,z_0)^2=\beta(z_0,w)^{-2}\), so
\[
 C_D(\rho(z_0)f;w)
 =\beta(w,z_0)^2 C_D(f;w)
 =\psi_E(j_{\rm W}\operatorname{tr}(yw))C_D(f;w).
\]
Indeed \(\operatorname{tr}(w^\dagger z_0)=
\overline{\operatorname{tr}(yw)}\), \(\overline{j_{\rm W}}=-j_{\rm W}\),
and \(\psi_E(\bar x)=\psi_E(x)\).
After division by the character at \(T\) the frequency is
\(\psi_E(j_{\rm W}\operatorname{tr}(y(w-T)))\).
The coordinates \(y_{ib}\) with \(a_i>\max(0,a'_b)\) vary
independently in \(E\backslash\mathbb A_E\). Their product mean
is one precisely when every corresponding \((w-T)_{bi}=0\), and
zero otherwise. This proves \eqref{eq:adjacent-common-projection}
for both \(D=0\) and \(D=a'\). Indeed, the character identity
and the absolute convergence of the \(I_D\)-Fourier expansion give
\[
 \int\sum_{w\in\mathcal C_D}
         |C_D(\rho(y^\dagger)f;w)|\,dv
      =\sum_{w\in\mathcal C_D}|C_D(f;w)|<\infty
\]
because \(\int dv=1\). Division by the unit character at \(T\)
preserves this majorant. Fubini thus justifies the termwise
projection after the \(I_D\)-average defining each coefficient.

Evaluate this sum again at \(D=a'\). The rational
Bessel group preserves agreement and fixes the class
of \(T\) in \(\mathcal C_{a'}\). On nonzero graded
pieces it fixes the string coordinates by the positive
flag and its dual, and on the grade-zero piece it fixes
the central string coordinate. The definite grade-zero
part is killed by \(T\). Thus \(T(\gamma-1)\) has degree
at least zero. Likewise \((\delta-1)T\) has degree at
least zero for \(\delta\in R'(F_0)\), since \(T\) has
no component in the definite target block.
The second phase test and orbit cancellation leave
only the class of \(T\). By
\eqref{eq:adjacent-agreement-stabilizer}, the translated
agreeing sets for \(D=0\) are indexed by the cosets
\(\delta R'(F_0)\) in \(U(V')(F_0)\). Their integrals
are obtained from the original agreeing set by replacing
\(h\) with \(\delta^{-1}h\) in
\eqref{eq:adjacent-rational-ct}. We obtain
\begin{equation}
 \int_{[R]}W(r)\theta(\omega(r,h)\Phi)\,dr
  =\sum_{\delta\in R'(F_0)\backslash U(V')(F_0)} I(\delta h),
 \qquad
 I(h)=\int_{[R]}W(r)C_{a'}(\omega(r,h)\Phi;T)\,dr .
 \label{eq:adjacent-unfolded-sum}
\end{equation}
The invariance of the class of \(T\) gives left
\(R'(F_0)\)-invariance of \(I\).

\subsubsection{Absolute convergence and unfolding}\label{sec:paired-lifts-absolute}

\sectionopening{In this subsection, we bound the sum over the orbits of Section~\ref{sec:paired-lifts-rational-orbits} before unfolding.}

For integration over \([G']\), choose signed orthonormal
coordinates at each real place, so the Hermitian matrix
satisfies \(J_\tau^*=J_\tau=J_\tau^{-1}\). Measure a rational
vector of the moved finite lattice by \(\max_\tau\|v_\tau\|_2\).
If \(\lambda_1(h)\) is the least such length of a nonzero vector,
put \(H(h)=\max(1,\lambda_1(h)^{-1})\). Choose an integer
\(d_H\ge1\) that clears the Hermitian pairings on the fixed finite
lattice, and set \(\eta=d_H^{-1/2}\).
The moved lattice has fixed covolume, since the determinant of
a unitary isometry has adelic absolute norm one.
If \(H(h)\le 2/\eta\), the first minimum is at least \(\eta/2\).
Minkowski's second theorem
\cite[Ch.~VIII, \S4.3, Theorem~V, p.~218]{Cassels}
bounds the product of the successive minima by a constant times
the covolume, so the last minimum is bounded as well. This gives a bounded rational basis in the adelic
coordinates.

Use a full set of bounded independent
\(\mathbb Z\)-lattice vectors and take an \(E\)-basis
among them. Its finite coordinates are integral in
the chosen basis. The rational Gram matrices of these
bases lie in a compact adelic set and therefore take
finitely many values. Fix a rational basis for each
Gram matrix. The map from the moving basis to that
fixed one is a rational isometry. The change of basis
puts the inverses of the resulting translates in a
compact set of adelic matrices, and thus the translates
themselves in a compact set of the unitary group. Here
inversion on the latter uses the fixed Hermitian form.

If \(H(h)>2/\eta\), take the lattice vectors of length below
\(\eta\). For any two such vectors, their rational Hermitian
pairing multiplied by \(d_H\) is integral at all finite places,
since \(h\) is unitary and both vectors lie in the fixed lattice
in the moved coordinates. At each real place the absolute value is
at most \(d_H\|v_\tau\|_2\|w_\tau\|_2<d_H\eta^2=1\).
The norm of a nonzero algebraic integer has absolute value at least
one. Thus every such pairing is zero and the short vectors span a
totally isotropic space. If \(U_\tau\) denotes the resulting
isotropic space, \(J_\tau U_\tau\) is Euclidean-orthogonal to it
and pairs perfectly with it for the Hermitian form. Dilate the isotropic space
by \(2\) and its opposite by \(1/2\), and fix the joint
orthogonal complement. This is unitary and both it
and its inverse have operator norm at most \(2\).
The short lengths double. Every other length stays at least
\(\eta/2\). The resulting height is at most
\(\max(2/\eta,H(h)/2)\). Repetition gives
\(\delta h=z x\), with \(\delta\) rational,
\(z\) in a fixed compact set and \(x\) purely real, with
\(\|x_\tau\|,\|x_\tau^{-1}\|\le H(h)\).
Indeed there are at most \(\lfloor\log_2 H(h)\rfloor\) dilations,
and each dilation and its inverse has norm at most two.
The compact set and the finite pairing denominator depend only on
the fixed group and lattice of this unfolding.

In the real Schr\"odinger model use the \(L^2\)-norms
of monomials in coordinate multiplication and
differentiation operators. By Heisenberg covariance,
conjugation by the Weil operator sends the linear
operators to linear combinations with coefficients
of the symplectic matrix or its inverse. By
\(L^2\)-unitarity, these norms of the translate of
the test are bounded by the corresponding norms
of the original test times a polynomial in the
matrix and inverse norms.

Sobolev embedding bounds
the weighted suprema above by finitely many such
\(L^2\)-norms. The symplectic action is
\(w\mapsto hwr^{-1}\). Its matrix and inverse are
polynomially bounded on the reduced representatives
with \(r\) in a compact set. At the finite places the
translates of the test over a compact set lie in a
finite set by smoothness of the Weil action. Therefore
\[
          \sum_{w\in\mathcal C_0}|F_w(r,h)|\le C H(h)^B .
\]
Here \(B\) depends on the fixed spaces and test \(\Phi\).
The constant \(C\) also includes \(\sup_{[R]}|W|\), finite because
\([R]\) is compact and \(W\) is continuous. At \(D=0\) every
rational target isometry preserves the filtration. Equation
\eqref{eq:adjacent-rational-ct} therefore permutes all frequencies
under either rational reduction used in this estimate.

For any continuous automorphic \(y\) with
\(\sup_h H(h)^B|y(h)|<\infty\) for all \(B\), the
product majorant is bounded on \([R]\times[G']\).
Each \(I(\delta h)\) is the integral of the sum over its
own agreeing set at \(D=0\), by
\eqref{eq:adjacent-common-projection} and the second orbit
cancellation. The distinct sets are disjoint, whence
\[
 \sum_{\delta}|I(\delta h)y(h)|
 \le |y(h)|\int_{[R]}\sum_{w\in\mathcal C_0}|F_w(r,h)|\,dr
 \le C\operatorname{vol}[R] H(h)^B |y(h)|.
\]
Both quotients have finite volume. Fubini and quotient integration give
\begin{proposition}[Unfolding consecutive lifts]\label{prop:adjacent-paired-unfolding}
For the data let \(W\) be continuous on \([R]\), with
\(W(ur)=\beta(Tu,T)W(r)\) for \(\deg(u-1)\ge2\). Let \(y\) be
continuous on \([G']\) and satisfy
\(\sup_h H(h)^B|y(h)|<\infty\) for every \(B>0\).
The two integrals below converge absolutely and unfolding with the
inherited target measure gives
\begin{equation}
 \int_{[R]} W(r)\int_{[G']}
        \theta(\omega(r,h)\Phi)y(h)\,dh\,dr
   =\int_{R'(F_0)\backslash U(V')(A)} I(h)y(h)\,dh .
 \label{eq:adjacent-paired-unfolding}
\end{equation}
The quotient measure on the right uses counting measure on the rational
subgroup, and the unfolding has scalar one.
\end{proposition}

At the final arrow \(Y_{m-1}\to Y_m\) both strings are empty. All grades
are zero, \(I_D=0\), \(\mathcal C_D=\{0\}\), and \(T=0\). The Bessel
groups are \(R=U(V)\) and \(R'=U(V')\), both definite. Hence
\[
 C_D(f;0)=\theta(f),\qquad
 I(h)=\int_{[R]}W(r)\theta(\omega(r,h)\Phi)\,dr .
\]
Both quotients are compact. Fubini gives
\eqref{eq:adjacent-paired-unfolding} with the original Haar measures.

\subsubsection{Rapid decay of the projected tests}\label{sec:paired-lifts-projected-tests}

\sectionopening{In this subsection, we verify the decay hypothesis of Proposition~\ref{prop:adjacent-paired-unfolding} for the projected test.}

For the target test coming from an integrable discrete
series projector, rapid decay follows as follows.
Let \(y\) be any bounded finite-level form reproduced
by the coefficient projector of the prescribed
integrable real discrete series. For the target tests in the
lower chain, Section~\ref{sec:even-stage-coefficients} verifies
these hypotheses for the conjugate projector. Write the
reproduction identity as
\(y=R(e)y\), where \(e=d_\pi\overline{\langle
v_0,\pi(\cdot)v_0\rangle}\) is the convolution
projector for a unit \(K_\infty\)-finite vector \(v_0\)
of the specified integrable discrete series \(\pi\).

Take the inner product linear in the second variable.
For \(v\) in the Harish-Chandra module of \(\pi\), the map
\[
             v\longmapsto
 d_\pi\int y(hx)\langle\pi(x)v_0,v\rangle\,dx
\]
takes \(v_0\) to \(y\). Denote it by \(T\), and put
\(b_v(x)=\langle\pi(x)v_0,v\rangle\). Write
\(G_\infty=Z_cG_\infty^1\), with \(Z_c\) the compact scalar torus
and \(G_\infty^1\) the product of the special unitary groups. The
product map has finite kernel, and restriction of \(\pi\) to
\(G_\infty^1\) is irreducible since the centre acts by scalars. If
\(\omega_\pi\) is its central character, then
\[
 b_v(zs)=\overline{\omega_\pi(z)}b_v(s)
        \qquad(z\in Z_c,\ s\in G_\infty^1).
\]
Hence the restriction of the nonzero coefficient \(b_{v_0}\) belongs
to \(L^1\cap L^2\) on \(G_\infty^1\). Apply
\cite[arXiv:1610.05483v2, Lemmas~5-3 and~5-5(i), (iii), (v), (vii)]{MuicIntegrable} with
\(p=1\) to this connected semisimple group with finite centre.
After inversion of the coefficient argument, these lemmas give
smoothness in \(L^1\) for the left translates of \(b_v\);
the central factor differentiates by the displayed character.
Since \(b_v(k^{-1}x)=b_{\pi(k)v}(x)\), the change of integration
variable and differentiation give
\(R(k)Tv=T\pi(k)v\) for \(k\in K_\infty\) and
\(dR(X)Tv=T\,d\pi(X)v\).
Each \(Tv\) is bounded, as \(y\) is bounded and the coefficient
is integrable. The quotient has finite volume, so \(Tv\in L^2\).

If \(y\ne0\), the pullback \(\langle Tv,Tw\rangle_{L^2}\)
is a nonzero invariant Hermitian form on the irreducible
Harish-Chandra module of \(\pi\). On each finite-dimensional
\(K_\infty\)-isotypic subspace it is represented, in the Hilbert
inner product of \(\pi\), by a positive semidefinite operator.
Their direct sum intertwines the Lie algebra by infinitesimal
invariance and intertwines \(K_\infty\). Schur's lemma therefore
gives \(\|Tv\|_2^2=c\|v\|_\pi^2\) for one \(c>0\).
Thus \(T\) extends to the Hilbert completion. For a real Lie
algebra element \(X\) and a \(K_\infty\)-finite vector \(v\),
\[
 dR(X)^nTv=T\,d\pi(X)^nv\qquad(n\ge0).
\]
The analytic-vector series for \(v\) converges in norm for small
\(|t|\), and the norm identity gives the same convergence for \(Tv\).
Their sums give \(R(\exp tX)Tv=T\pi(\exp tX)v\).
By density and generation by the identity component and
\(K_\infty\), the extension is equivariant. Its closed image
is a discrete-series summand.

Borel's class-number theorem gives finitely many adelic double
classes at the standard integral compact \(K_f^0\)
\cite[Theorem~5.1, p.~19]{BorelAdele}. For the specified compact open
level \(J\), each such class splits into at most
\([K_f^0:K_f^0\cap J]\) classes at the intersection; hence there
are finitely many arithmetic real quotient components at \(J\). The determinant of each element of its
arithmetic subgroup \(\Gamma\) is an integral unit of absolute value
one at every embedding, hence a root of unity. Thus
\(\Gamma^1=\Gamma\cap G_\infty^1\) has finite index in \(\Gamma\).
Lift to \(\Gamma^1\backslash G_\infty\). The centre acts by a unitary
character, so integration over \(Z_c\) shows that restriction to
\(\Gamma^1\backslash G_\infty^1\) preserves the squared norm up to a
fixed positive scalar. Here the product map has finite central kernel.
The restriction of \(\pi\) is irreducible because \(Z_c\) acts by
scalars, and remains square integrable by the same central integration.

On each component with nonzero image let \(G^1\) be the restriction
of the special unitary group from \(F_0\) to \(\mathbf Q\).
At the finite representative \(h_f\), the compact open
\(J^1=h_fJh_f^{-1}\cap G^1(\mathbb A_f)\) preserves a lattice in
the underlying rational representation, obtained locally by taking
the sum of the translates of an integral lattice under the compact
group and retaining the standard lattice almost everywhere.
Thus \(J^1\) has finite index in the adelic stabilizer of this
lattice in \(G^1(\mathbb A_f)\), and intersection with \(G^1(\mathbf Q)\)
puts \(\Gamma^1\) at finite index in the integral arithmetic group.
For each right differential operator \(D\), the intertwining identity
gives \(D(Tv)=T(d\pi(D)v)\), which is bounded by the integrability
of its coefficient. Together with compact finiteness and the
infinitesimal character this gives Wallach's automorphic-function
hypotheses. His theorem places this tempered summand in the cuspidal
spectrum \cite[Theorem~4.3, p.~233]{WallachConstant}. Each rational parabolic of the
unitary group has the same unipotent radical \(N\) as its intersection
with the derived group, and \(\Gamma\cap N=\Gamma^1\cap N\).
At the finite representative \(h_f\) of a component, put
\(J_N=N(\mathbb A_f)\cap h_fJh_f^{-1}\).
Exponential coordinates and additive approximation give
\(N(\mathbb A_f)=N(F_0)J_N\). The resulting quotient is
\[
 N(F_0)\backslash N(\mathbb A)/J_N
   \simeq (N(F_0)\cap J_N)\backslash N(F_{0,\infty}),
\]
where the intersection uses the finite embedding of \(N(F_0)\).
The quotient integration formula with probability measures
identifies the adelic constant term on this component with
the real one, which vanishes. The restriction of \(T(v_0)\)
on each arithmetic component is a bounded smooth
cuspidal form, hence of moderate growth. The rapid decay theorem
\cite[arXiv:1106.2149v1, Theorems~A and~2.23--2.24]{MillerSchmidDecay} gives \(O_A(\|g\|^{-A})\) on each of the finitely many
Siegel sets for every \(A>0\).
Central translates have the same absolute value.

The finite
lattice is fixed on each real component, so
\(H(h)\ll 1+\|h_\infty\|\) there. Indeed the original
nonzero lattice vectors have a fixed positive minimum
and the inverse real action decreases lengths by at
most \(\|h_\infty\|\). This proves the required decay
in \(H(h)\) for the projected tests.

\subsection{Two successive Bessel integrals}
\label{sec:double-string-unfolding}

\sectionopening{In this subsection, we combine two copies of the unfolding in
Section~\ref{sec:adjacent-unfolding}. The middle Hermitian space then
appears in a single theta integral, and its constant term gives the
paired operator needed in Section~\ref{sec:even-stage-coefficients}.}

Take three lower spaces \(Y,Y',Y''\) with the same \(c,N\), odd initial string length \(s\), and remaining
lengths \(s-1,\max(0,s-2)\). Write \(T,T'\) for the
two string maps and let
\[
 S=(T^\dagger,T')^{\mathsf T}:V'\longrightarrow V^-\oplus V'',
\]
where \(V^-\) has form \(-h\). For an interior pair,
adjunction in the antidiagonal forms gives
\[
 TT^\dagger e'_i
       = t_{s-2-i}t_{i+1}(g'_i/g_{i+1})e'_{i+1}
\]
before the last coordinate. This is the middle shift,
with coefficient \(c\) except at grade one, where it is
\(-Nc^2\). Thus \(TT^\dagger=T'^\dagger T'\), and
\(S^\dagger S=-TT^\dagger+T'^\dagger T'=0\).
The same adjoint computation for even \(s\) gives the
odd target shift, with \(-Nc\) at grade zero. At \(s=1\)
the maps are zero.

Let \(u\in U_{R'}(\mathbb A_{F_0})\) satisfy
\(\deg(u-1)\ge2\). Target covariance replaces \(T\)
by \(u^{-1}T\), and \((u^{-1}-1)T\) has degree at least
one. Equation~\eqref{eq:adjacent-positive-shift} and
\(u^\dagger=u^{-1}\) give
\[
 \begin{aligned}
 I(uh)&=\beta(T,u^{-1}T)I(h),\\
 \beta(T,u^{-1}T)
 &=\psi_E\bigl((j_{\rm W}/2)\operatorname{tr}(TT^\dagger u)\bigr)
   =\beta(T'u,T'),
 \end{aligned}
\]
since \(TT^\dagger=T'^\dagger T'\).
This gives the middle Bessel character on the subgroup
of degree at least two. Section~\ref{sec:two-bessel-constant-terms}
computes the contribution of the degree-one entries.
Continuity follows from the compact \([R]\) integral
and the local Schwartz bounds. We may therefore use
the resulting \(I\) on the next Bessel group.

The adjoint identifies the symplectic spaces
\(\operatorname{Hom}(V,V')\) and
\(\operatorname{Hom}(V',V^-)\). The Hermitian trace
on the latter is the negative conjugate of the trace
pairing on the former, which preserves the alternating
form with \(j_{\rm W}\). Use the corresponding linear
oscillator transport and tensor it with the second datum.
Write \(\mathcal B(f,f')\) for the joined datum. This
intertwines the Heisenberg translations by \(w\) and
\(w^\dagger\).

For the opposite finite basis, multiply
the original basis by the conjugate of
\(\eta\alpha_{\rm W}/\bar\alpha_{\rm W}\), where
\(\eta\in\widehat{\mathcal O}_E^\times\) has norm \(-1\).
Such a unit exists at every finite place since \(E/F_0\)
is unramified or split there. Put
\(\lambda=\eta\alpha_{\rm W}/\bar\alpha_{\rm W}\) and
\(A^*=\bar A^{\mathsf T}\). The pullback on finite oscillator
coordinates is \(A\mapsto\lambda A^*\), with inverse
\(A\mapsto-\lambda A^*\) since \(\lambda\bar\lambda=-1\).
In the split Schr\"odinger coordinates on the conjugate places,
\(\lambda(\alpha_{\rm W}Y)^*=\alpha_{\rm W}\eta Y^*\);
hence the pulled-back test is \(Y\mapsto\Phi(\eta Y^*)\).

For a finite adelic matrix \(X\), use the quadratic phase and
translation of Section~\ref{sec:upper-cm-jet-transport},
\[
 e_0(X)=\psi_0\!\left(c_0\operatorname{tr}(X^*X)/2\right),
 \qquad T_\sharp(X)=e_0(X)\rho_f(\alpha_{\rm W}X),
\]
with the product of the local characters \(\psi_0\) and
\(c_0=(1-j_{\rm W}^{-1})/2\).
Since \(\operatorname{tr}((\eta A^*)^*(\eta A^*))
=-\operatorname{tr}(A^*A)\), one has
\(e_0(\eta A^*)=e_0(A)^{-1}\).
The quadratic phase of a column of two blocks is the
product of their phases. Heisenberg covariance therefore gives
\[
 T_\sharp\!\binom{A}{B}\,\mathcal B(f,f')
   =e_0(A)^2\,\mathcal B\bigl(T_\sharp(\eta A^*)f,T_\sharp(B)f'\bigr).
\]
For integral \(A\), the argument \(c_0\operatorname{tr}(A^*A)\)
is integral, so \(e_0(A)^2=1\) also above \(2\).
For a fractional displacement the displayed phase is retained.
The constant terms below use the unnormalized adelic translations
\(\rho\) of Section~\ref{sec:adjacent-unfolding}.

In rational coordinates the
forward map is ordinary adjunction.
Transporting a rational Schr\"odinger polarization by
this symplectic map gives the same sum over the rational
Lagrangian. The two theta normalizations therefore
satisfy
\[
 \theta\bigl(\mathcal B(f,f')\bigr)
          = z_\theta\,\theta(f)\theta(f'),\qquad z_\theta\ne0 .
\]
The scalar is fixed by the models and is
independent of the two tests.

\subsubsection{The constant terms in the joined block}
\label{sec:two-bessel-constant-terms}

\sectionopening{In this subsection, we group the constant terms before the lower-lift comparison of Section~\ref{sec:seesaw-normalization}.}

Let \(Z_1\) be the zero term along the joined
positive-degree block. The zero term along an isotropic
space \(K\), at a vector in the residual space, is the
theta functional after translating by that vector and
averaging over \(K(F_0)\backslash K(A)\) with probability
measure. The map \((w,w')\mapsto(w^\dagger,w')\) takes
the two positive blocks to the joined block, since
adjunction reverses both grades. It is a rational
bijection and pushes product probability to probability.
Consequently
\begin{equation}
 C_{a'}(f;w) C_{a''}(f';w')
  =z_\theta^{-1}
      \bigl[Z_1\rho((w^\dagger,w'))\mathcal B(f,f')\bigr](0).
 \label{eq:double-constant-term}
\end{equation}

For an element \(n\) of the middle unipotent radical,
the two frequencies in
\eqref{eq:double-constant-term} become \(n^{-1}T,T'n\)
by covariance. Since \(n^\dagger=n^{-1}\), their joined
frequency is \(((n^{-1}T)^\dagger,T'n)^{\mathsf T}=Sn\).
The cross term is trivial because
\(S^\dagger S(n-1)=0\), so
\(\rho(Sn)=\rho(S(n-1))\rho(S)\).
The map \(S(n-1)\) has nonnegative degree. If \(s>1\),
the outer strings have even grades and the middle
string has odd grades. The grade-zero part of Hom
is thus the maps from \(X'\) to the outer grade-zero
space. The two outer middle vectors \(z,z'\) have
norms \(N,-N\) in the joined target. They form the
hyperbolic plane with basis \((z-z')/(2N),z+z'\).
The grade-one middle vector \(e_+\) has
\(Se_+=c(z+z')\): the first coefficient is the adjoint
coefficient \((-Nc)/(-N)\) and the second is \(c\).
It follows that the degree-zero part of \(S(n-1)\)
is \((z+z')b(n)\), with \(b(n)_i=c(n-1)_{+,i}\)
for \(i\in X'\). Its positive part is orthogonal
to the degree-zero part and is absorbed by \(Z_1\).
The product is \(z_\theta^{-1}\) times the value
of \(Z_1\rho(S)\mathcal B(f,f')\) at \((z+z')b(n)\).

The map \(b\) is additive on the middle unipotent,
since the product of two positive-degree differences
has degree at least two. For \(s>1\) it is onto over \(F_0\)
and over \(A\). To prescribe \(b_i\), take \(x=e_+\),
\(v=\sum_i c^{-1}\bar b_i e_i\) in \(X'\), and
\(L=x h'(v,\cdot)-v h'(x,\cdot)\).
Here \(L^\dagger=-L,\ L^3=0\), and \(Le_i=c^{-1}b_i e_+\),
\(L^2e_i=0\) for each definite basis vector \(e_i\).
Hence \(u(t)=1+L+L^2/2\) is unitary, strictly grade-raising
off the identity, and \(b(u(t))=t=(b_i)_i\). Since \(b\)
sends rational elements to rational covectors, it descends
to the left rational quotients. For invariant probability
\(dn\) on the unipotent quotient, right translation by \(u(t)\)
preserves \(dn\) and \(b(nu(t))=b(n)+t\).
Thus \(b_*dn\) is the translation-invariant probability on the
adelic covectors modulo rational covectors. Its average gives a
second zero term \(Z_2\) along \(z+z'\). The two slices
compose to \(\mathcal Z=Z_2 Z_1\), with Haar probability
on the direct sum of the isotropic directions.

For \(s=1\) the two string maps vanish, all three gradings
are zero and there is no plane to remove. The vector
\(z/\pi_N\) has norm \(N/(\pi_N\bar\pi_N)=1\) in \(V^-\).
It is the last rational coordinate of \(X''\), whose other
coordinates are the definite part of \(Y''\). Thus
\(\widehat k=\operatorname{diag}(k,1)\), and both zero-term
operators are the identity in these coordinates.
Section~\ref{sec:marked-neighbor} computes the original lattice
in this basis at the conjugate places above \(N\).
In the interior \(\widehat k=k\).

For the determinant factors, write \(\upsilon_1\) for
the splitting character on norm-one ideles.
On the isotropic block of a parabolic Weil action let
\(A_0\) be its matrix. The scalar of the zero slice,
relative to the action on the residual space, is
\(\upsilon(\det A_0)^{-1}|\det A_0|^{1/2}\).
The modulus is the square root of the Jacobian in
the dual coordinates \(A_0^\dagger x\).
For the phase, use an adapted Gaussian and finite
lattice vacuum.

In the real basis with prescribed signs, the negative
block of the hyperbolic Levi has matrix
\((A_0^{-\dagger}+A_0)/2\). Its conjugate determinant
times \(\det A_0\) is the positive determinant of
\((1+A_0^\dagger A_0)/2\).
At a finite place put \(L=\mathcal O_E^{\dim A_0}\)
in the corresponding local coordinates. The multiplier
for an invertible \(g\) is
\(\chi_f(\det B)\), where \(BL=L\cap g^{-1}L\) and
\(\chi_f=\upsilon_f|\cdot|^{-1/2}\).
For the hyperbolic pair \(A_0^{-\dagger},A_0\) this
has multiplier
\(\overline{\chi_f(\det A_0)}|\Delta_f(A_0)|^2\).

Indeed \(A_0\) carries \(L\cap A_0^{-1}L\) to
\(L\cap A_0 L\), so inversion multiplies \(\Delta_f\)
by \(\chi_f(\det A_0)\). Transposition leaves the
intersection determinant ideal unchanged, by elementary
divisors, and conjugation conjugates \(\Delta_f\).
Thus the hyperbolic multiplier has phase
\(\upsilon(\det A_0)^{-1}\).
The Gaussian overlap and the overlap of the finite
zero slices are positive in the adapted coordinates.
The phase of the scalar is therefore the phase of
this multiplier.

The parabolic radical has trivial
scalar: it acts on the zero slice by a scalar by
Heisenberg covariance and every element of that
radical is a commutator with an isotropic dilation
on its two successive positive-grade quotients.
Equivalently, in the canonical Schr\"odinger
splitting the terms of its phase vanish at the zero
coordinate. This also gives the unipotent covariance
used in Section~\ref{sec:adjacent-unfolding}.

The arrow of parameter \(a\) multiplies the
raw action of \((g,h)\) by
\[
 \frac{a(\det h)}{a(\det g)}\,
 \upsilon_1(\det h)^{\lfloor n/2\rfloor}
 \upsilon_1(\det g)^{-\lfloor n'/2\rfloor},
\]
where \(n,n'\) are the ambient ranks.
The adjoint transport on the first arrow contributes
\(\upsilon_1(\det g)^{n'}\upsilon_1(\det h)^{-n}\).
Indeed, under the transpose of the inverse matrix
of the action on Hom, the opposite sign interchanges
the real positive and negative blocks.

Put \(D(w)=hwg^{-1}\) on the first Hom space. At a real place
\(v\) write \(\Delta_{\epsilon_v}(D_v)=\overline{\det((D_v)_{--})}\).
Jacobi's complementary minor identity and \(|\det D_v|=1\) give
\[
 \Delta_{-\epsilon_v}(D_v^{-1})
   =\overline{\det((D_v^{-1})_{++})}
   =(\det D_v)\Delta_{\epsilon_v}(D_v).
\]
At each finite place \(D_v(L\cap D_v^{-1}L)=L\cap D_vL\).
The determinant has norm one, so \(|\det D_v|_{E_v}=1\), with
both factors included at a split place. The intersection formula gives
\[
 \Delta_f(D_f^{-1})=\upsilon_f(\det D_f)\Delta_f(D_f).
\]
Transposition preserves the principal minors and the finite elementary
divisors. Reordering the matrix indices permutes the corresponding
coordinates in both bases. The original multiplier divided by the
opposed multiplier is therefore \(\upsilon_1(\det D)^{-1}\), with
\(\det D=(\det h)^n(\det g)^{-n'}\).
The exponents on the three determinants in the
joined action are consequently
\[
   (\lceil n'/2\rceil,\
          -\lceil n/2\rceil-\lfloor n''/2\rfloor,\
                     \lfloor n'/2\rfloor).
\]
The factors \(a(\det h)/a(\det g)\) and
\(b(\det k)/b(\det h)\) are retained.

Let \(d_0=(s-1)/2\), let \(p_0\) count the positive
outer grades, and let \(q_0\) be one in the interior
and zero at \(s=1\). On the positive Hom block the
outer determinant occurs once per negative middle
coordinate and the middle inverse determinant once
per positive outer coordinate. Both determinants have
norm one. The first slice changes the exponents by
\((-d_0,p_0,-d_0)\). The plane slice changes them by
\((0,q_0,0)\): the plane is fixed and the column
action is \(h^{-1}\). If \(x,x',x''\) are the final
ranks, then
\[
 n'=x'+2d_0,\qquad n+n''=x+x''+2p_0+2q_0 .
\]
Here \(n''=n+2t\), so each of the pairs \((n,n'')\),
\((x,x'')\) has equal parity. In particular
\[
 -\lceil n/2\rceil-\lfloor n''/2\rfloor+p_0+q_0
 =-(x+x'')/2=-\lceil x/2\rceil-\lfloor x''/2\rfloor .
\]
Subtracting \(d_0\) from the two outer exponents gives
\(\lceil x'/2\rceil,\lfloor x'/2\rfloor\), respectively,
and \(\det\widehat k=\det k\). The joined motions of the definite
groups commute with \(\rho(S)\), since they fix
all the string maps. Thus \(\mathcal Z\) carries
the ambient action to the joined definite action.

\subsubsection{The two quotient volumes}
\label{sec:two-bessel-quotient-volumes}

\sectionopening{In this subsection, we compare the quotient measures after the preceding unfolding with the normalization of Section~\ref{sec:period-quotient-measure}.}

Write \(\kappa_R=\operatorname{vol}[R]/
\operatorname{vol}[U(X)]\) and similarly \(\kappa_{R'}\).
The Levi acts on the unipotent radical with modulus
one, since a determinant of its action is a rational
character of a unitary group. Disintegrating the
original quotient Haar therefore gives the Levi
quotient Haar times this constant and Haar probability
on each unipotent fibre. For odd \(s\), the source
unipotent acts with degree at least two. Its phase in
the first constant term cancels the character of
\(W\) by \eqref{eq:adjacent-positive-shift}.
In the interior the middle mean gives \(Z_2\) by the
surjectivity of \(b\). At \(s=1\) the two unipotent radicals
are trivial and their probability integrals are evaluations.

Let \(\omega^{\rm def}_{a,b}\)
be the joined action on the two definite arrows,
and \(\iota_X,\iota_{X''_0}\) the original definite
group embeddings into the two outer unitary groups.
For the two arrows write
\[
 D_f(r,r')=C_{a'}(\omega_a(r,r')f;T),\qquad
 D_{f'}(r',h)=C_{a''}(\omega_b(r',h)f';T').
\]
The actions include the determinant characters of their respective
arrows. For every continuous \(W\) with the prescribed source
unipotent character, the two integrations just computed give
\begin{equation}
\begin{aligned}
 &\int_{[R']}\!\int_{[R]} W(r)D_f(r,r')
          D_{f'}(r',\iota_{X''_0}(k))\,dr\,dr'\\
 &\quad=\frac{\kappa_R\kappa_{R'}}{z_\theta}
 \int_{[U(X)]}\!\!\int_{[U(X')]}
  W(\iota_X(g))\,
  \theta_{\rm def}\!\left(
     \omega^{\rm def}_{a,b}(g,h,\widehat k)
            \mathcal Z\rho(S)\mathcal B(f,f')\right)\,dh\,dg .
\end{aligned}
 \label{eq:double-unfolded-integral}
\end{equation}
The integrands descend continuously to the compact
Bessel quotients. This justifies the order of the means.
At the endpoint the Bessel subgroup is the full definite unitary group. The coefficient on it is
evaluated at \(\operatorname{diag}(k,1)\) with the
branching normalization of Section~\ref{sec:upper-chain-size}.
The two quotient-mass constants in
\eqref{eq:double-unfolded-integral} are for the
original Haar measures.

\subsection{The coefficients after an even number of steps}
\label{sec:even-stage-coefficients}

\sectionopening{In this subsection, we iterate the paired identity of
Section~\ref{sec:double-string-unfolding}. We keep the source indices
and identify the real coefficients after each pair of theta operators.}

Write \(R_j\) for the Bessel group of \(Y_j\),
\(J_{\xi,j}\) for the principal kernel at the
chosen level and twist places, and
\(A_{\xi,j}=[K_j^\circ:J_{\xi,j}]\) for its index
in the standard compact. Put \(G_j=U(Y_j)\) and write
\(K_{\xi,j}(g,h)=\theta(\omega_{\xi,j}(g,h)\Phi_{\xi,j})\)
for the kernel with the given level test, rank twists and splitting character.
For a source form \(f\) and a test \(z\) of the conjugate target type, the
two operators for the bilinear pairing are
\[
 \begin{aligned}
 (L_{\xi,j}f)(h)&=A_{\xi,j}\int_{[G_j]}K_{\xi,j}(g,h)f(g)\,dg,\\
 (L_{\xi,j}^{\,\mathrm t}z)(g)&=A_{\xi,j}\int_{[G_{j+1}]}
                                  K_{\xi,j}(g,h)z(h)\,dh .
 \end{aligned}
\]
Thus \(\int_{[G_{j+1}]} zL_{\xi,j}f
=\int_{[G_j]}(L_{\xi,j}^{\,\mathrm t}z)f\).
Indeed the compact finite translates give finitely many Schwartz tests.
On reduced representatives the real symplectic matrix for
\(w\mapsto hwg^{-1}\) and its inverse are bounded by powers of the
product of the two heights. The Schwartz seminorm estimate of
Section~\ref{sec:paired-lifts-absolute} gives
\(\lvert K_{\xi,j}(g,h)\rvert\ll H_j(g)^C H_{j+1}(h)^C\)
for the fixed test. The rapid decay of \(f,z\) and finite quotient
volumes make the double integral absolute, so Fubini applies.
In either variable the principal
kernel fixes the lattice-inclusion test and
the determinant character. In the terminal group,
an integral element with block-upper-triangular
reduction at the twist and identity at the chosen
level places has leading blocks that lift to the
preceding unitary groups. The inclusions intertwine those
blocks in the residues. The successive
determinant ratios are the tame characters
on the added diagonal blocks. The initial
block fixes the source vector, which has maximal
level at the twist places. This proves the
terminal finite character.

In the lower grid the added positive blocks
of length \(t=5r\) are centered at \(-w_i\)
and have entries \(-w_i+(t-1)/2-b\) for
\(0\le b<t\). The common half shift at odd
stages leaves all gaps unchanged. The base
tails lie midway between consecutive centers.
Only base entries carry negative sign. Thus
opposite signs have gaps at least
\(\Delta/2-t>500(r+mt)^3\) for large \(x\).
Let \(e_j\) denote the normalized discrete-series coefficient on
\(U(Y_j)_\infty\), with \(R(e)f(g)=\int e(k)f(gk)\,dk\).
Write \(e^\vee(k)=e(k^{-1})=\overline{e(k)}\).
The pair \((Y_j,Y_{j+1})\), \(0\le j<m\), is the horizontal
step in grid \(z=0\) with zero offset, no up blocks and \(j\)
previous down blocks. The ranks are \(r+jt,r+(j+1)t\). Writing \(c\) for the new centre, its positive block has doubled
entries \(2c+t-1-2b\), \(0\le b<t\). The common odd-stage shift cancels
after the opposite determinant exponent is subtracted. Thus the
signed parameter and joint-harmonic calculation of
Section~\ref{sec:four-unitary-real-projections} applies to this step
and, with conjugate types interchanged, to its transpose.

For clarity, take either orientation, with incoming and outgoing
coefficients \((e,e')=(e_j,e_{j+1})\) or
\((\bar e_{j+1},\bar e_j)\). Write \(P\) for its real harmonic and
\(L_P\) for the integral in the incoming variable with the given
finite test. In the real oscillator put
\[
 p_{\rm in}=\int e^\vee(k)\,\omega(k,1)\,dk,\qquad
 p_{\rm out}=\int e'(k)\,\omega(1,k)\,dk .
\]
The local line identity of
Section~\ref{sec:four-unitary-real-projections} gives
\(p_{\rm out}p_{\rm in}P=p_{\rm in}P\) in this orientation.
For a rapid incoming form \(y=R(e)y\), replace \(y(g)\) by
\(\int e(k)y(gk)\,dk\) in the bilinear lift and set \(g'=gk\). The
kernel at \(g'k^{-1}\) contributes \(\omega(k^{-1},1)P\). Inversion
in \(k\) gives \(p_{\rm in}\). Right convolution in the target gives
\[
 L_Py=L_{p_{\rm in}P}y,\qquad
 R(e')L_Py=L_{p_{\rm out}p_{\rm in}P}y=L_Py .
\]
To justify both projections before applying the seesaw, denote the
incoming and outgoing unitary groups by \(G_{\rm in},G_{\rm out}\)
and their ranks by \(a,b'\). Put
\(X(k)=\prod_{\tau}\max(1,\|k_\tau\|)\).
For fixed target point \(h\), source reduction and the coherent
Gaussian estimate in Section~\ref{sec:four-unitary-constant-terms}
give, for the original finite test,
\[
 |K_P(g'k^{-1},hk')|
    \le C_h H(g')^B X(k)^{ab'}X(k')^{ab'} .
\]
Indeed the source and target real moves act on the same Hom
space, whose operator norm is bounded by the product of their
norms, and finite rational reduction leaves a compact set of
finite translates. After \(g'=gk\), an absolute majorant for
the integral over \(g',k,k'\) therefore has integral at most
\[
 C_h\!\int_{[G_{\rm in}]}\! |y(g')|H(g')^B\,dg'\,
 \int_{G_{{\rm in},\infty}}\! |e(k)|X(k)^{ab'}\,dk\,
 \int_{G_{{\rm out},\infty}}\! |e'(k')|X(k')^{ab'}\,dk' .
\]
The first factor is finite by rapid decay, and the coefficient
moments in Section~\ref{sec:four-unitary-real-projections} give
the other two in both orientations. Thus Fubini applies with the
original quotient measure. The finite character and test, together with the
factor \(A_{\xi,j}\), are the same on both sides of each identity. The truncation argument
of Section~\ref{sec:four-unitary-real-projections} now places
\(L_Py\) in the cuspidal coefficient space at its fixed finite
level. Section~\ref{sec:paired-lifts-projected-tests} gives its
rapid decay. Starting at the definite endpoint and applying
this argument successively in the transposed orientation supplies
the rapid test required at every unfolding.

Fix outer points \(g_j=(u_j,y_j)\) and let
\(W_0\) be the original Bessel source coefficient
of \(\phi\) with the functional \(t\).
For \(j\ge0\) define the raw functions
\[
 I_{j+1}(h)=\int_{[R_j]}W_j(r)
 C_{a_{j+1}}\bigl(\omega_{\xi,j}(r g_j,h)
                               \Phi_{\xi,j};T_j\bigr)\,dr,
 \qquad W_{j+1}(r)=I_{j+1}(r g_{j+1}).
\]
The successor covariance of
Section~\ref{sec:double-string-unfolding} gives the required character of
\(W_j\) on \(\deg(u-1)\ge2\) at every step. At an even stage this is the
entire unipotent radical. At an odd stage its degree-one entries contribute
the second isotropic mean of Section~\ref{sec:two-bessel-constant-terms}.
The function \(I_j\) is independent of the
outer \(j\)-th point before its evaluation.
For a rapid automorphic \(z\) of the conjugate target type,
Proposition~\ref{prop:adjacent-paired-unfolding} gives
\[
 \int_{[R_j]}W_j(r)(L_{\xi,j}^{\,\mathrm t}z)(r g_j)\,dr
 = A_{\xi,j}\!
       \int_{R_{j+1}(F_0)\backslash U(Y_{j+1})(A)}
                                      I_{j+1}(h)z(h)\,dh .
\]
This is the unnormalized unfolding for one theta lift with the
source-level factor of the transpose. Put
\(G'=U(Y_{j+1})\) and \(R'=R_{j+1}\). The quotient integral
on its right disintegrates as
\[
 \int_{R'(A)\backslash G'(A)}\!
       \int_{[R']} I_{j+1}(r h)z(r h)\,dr\,d\dot h .
\]
Here \(d\dot h\) is the quotient of the given group Haar
measures, since both groups are unimodular. Absolute convergence
in Proposition~\ref{prop:adjacent-paired-unfolding} applies
also to this disintegration. For \(z=L_{\xi,j+1}^{\,\mathrm t}z'\),
fix a representative \(h=g_{j+1}\) in the outer integral.
The inner source is precisely \(W_{j+1}(r)\), so the next
unfolding contributes \(A_{\xi,j+1}\).

For \(j=2i\), apply \eqref{eq:double-unfolded-integral} to
\[
 f=\omega_{\xi,j}(g_j,g_{j+1})\Phi_{\xi,j},\qquad
 f'=\omega_{\xi,j+1}(g_{j+1},g_{j+2})\Phi_{\xi,j+1} .
\]
Indeed the two original factors have arguments \((r g_j,r'g_{j+1})\)
and \((r'g_{j+1},\iota_{X''_0}(k)g_{j+2})\). The representation
law puts the same middle point \(r'g_{j+1}\) in both. For the raw
function \(I_{j+2}(\iota_{X''_0}(k))\) take \(g_{j+2}=1\).
Join the translated tests and apply the
\(S\)-translation and both slices of
Section~\ref{sec:double-string-unfolding}.
The result is a pure real factor times a joined
finite factor. Indeed in an adapted adelic
Schr\"odinger polarization the slices are
restrictions at zero in the dual coordinates
of the two isotropic directions.

On a pure
tensor this is the product of real and finite
restrictions. Heisenberg translations, Weil
moves and the adjoint transport also factor
by place. Put the fixed theta normalization in
the real factor. Thus the finite factor is
the joined finite zero slice of \(\rho_f(S_f)\)
on the two translated finite vectors, with the
opposite-basis scale. Denote its inverse image
under the finite tensor identification by
\(F_i(\xi,y)\).

This inverse exists since
compactly supported locally constant tests on
a product are finite sums of product tests,
by a finite rectangular partition. The real
factor depends on the \(u_j\) and fixed weights.
Every finite tame character has trivial real
part, so the real factor is independent of
\(\xi,y\).

Write \(V_S,V_T\) for the complex column spaces of the
definite source and target. For the pullback
\(P(Z)\mapsto P(h^{-1}Zg)\) with \(Z:V_S\to V_T\), the
degree-\(D\) polynomials have the Cauchy decomposition
\[
 \operatorname{Sym}^D(V_T^\vee\otimes V_S)
  =\bigoplus_{\substack{\lambda\vdash D\\
                   \ell(\lambda)\le\min(\dim V_S,\dim V_T)}}
                   \mathcal S_\lambda(V_T^\vee)\otimes
                   \mathcal S_\lambda(V_S).
\]
Each summand occurs once \cite[\S4(d), p.~560]{HoweInvariant}.
The two determinant weights of the arrow twist the indicated
factors. Source integration selects the summand of the source
column partition, read in the next rank by padding with zero rows.
For the highest coefficient the polynomial is the product of the minors.

In normalized Bargmann coordinates the Gaussian pairing satisfies
\[
 \int_{\mathbb C^d}\bar z^\alpha z^\nu e^{-\|z\|^2}
            \,\frac{dz}{\pi^d}=\mathbf1_{\alpha=\nu}\,\nu!,
\]
so pairing with a Fock monomial gives its unscaled differential
at the origin. Let \(d_X\) be the source Weyl dimension, \(D\) the
degree of the coupled polynomial, and \(c_{\rm jet}=\pi_M s_{\rm jet}\ne0\)
the scale used in this algebraic lift. Its homogeneous
jet includes \(c_{\rm jet}^D\), so converting that jet back
to the unscaled theta coefficient gives \(c_{\rm jet}^{-D}\).
More generally, write the component paired with the source as a sum
of rectangular Schur coefficients \(Z\mapsto A_i\mathcal S(Z)v_i\),
where \(v_i\) belongs to the source Weyl module and \(A_i\) to the
target dual. These coefficients span the matched compact summand by
the Cauchy decomposition. For real points \(u,u'\), denote the
source and target compact Weyl actions by \(\sigma,\sigma'\).
Substitution in the minors, including the determinant twists,
carries the pair \((A_i,v_i)\) to
\((A_i\sigma'(u'^{-1}),\sigma(u)v_i)\).
At finite points \(g,h\), with finite test \(f\) and finite Weil
action \(\omega_f\), write \(J_cF\) for the Taylor jet at scale \(c\)
from Section~\ref{sec:upper-cm-jet-transport}. The rectangular matrix is
\(\mathcal J_c(F)=\bigl(\mathcal S(\partial_Z)(J_cF)\bigr)(0)\),
with the differential operators applied entrywise. The resulting theta coefficient is
\[
 c_{\rm jet}^{-D}
 A_i\sigma'(u'^{-1})\mathcal J_{c_{\rm jet}}(\omega_f(g,h)f)\sigma(u)v_i .
\]
The incoming source functional at \((u,g)\) is
\(t\sigma(u^{-1})\phi(g)\), where \(\phi\) is the source vector.
Compact integration uses the mass-one identity
\[
  \int\sigma(u)(v_i\otimes t)\sigma(u^{-1})\,du
       =\frac{t(v_i)}{d_X}\operatorname{Id} .
\]
The average commutes with the irreducible source action, and its
trace is \(t(v_i)\). Write
\(K\) for the fixed level on the definite source \(X_j\),
with the prescribed last-row condition modulo \(\pi_N\),
and \(K^\circ\) for its integral stabilizer. For a finite test choose a compact open subgroup \(J\subseteq K\)
fixing the scalar integrand. Its quotient integral for
\(\operatorname{vol}(K^\circ)=1\) is the class sum at \(J\) with
inverse stabilizer orders, multiplied by
\([K^\circ:J]^{-1}=[K^\circ:K]^{-1}[K:J]^{-1}\).
The algebraic lift, normalized by \(\operatorname{vol}(K)=1\),
has the factor \([K:J]^{-1}\) in the same sum. For this upper source
of rank \(n_j\), split reduction gives
\[
 [K^\circ:K]=|\mathrm{GL}_{n_j}(\mathcal O_E/(\pi_M))|\,
                 (N^{n_j}-1)^f .
\]
At each of the \(f\) chosen places over \(N\) the residue group
acts transitively on the nonzero row covectors, and \(K\) fixes the
last covector itself. The opposite component is determined by the
unitary condition.
Thus the scalar from the highest coefficient of the algebraic
lift to the theta integral is
\[
    \frac{c_{\rm harm}}
         {d_X\,c_{\rm jet}^D\,[K^\circ:K]},
\]
where \(c_{\rm harm}\) compares the chosen real harmonic
with the polynomial of minors. All four factors are fixed
at this definite step as \(\xi\) varies.

In particular both
maps have the same finite-datum dependence. For the preceding Schur
expansion the outgoing functional is
\[
   \frac{\sum_i t(v_i) A_i}{d_X c_{\rm jet}^D [K^\circ:K]} .
\]
This contracts the rectangular Schur coefficient with \(t\) and
is independent of the expansion, the finite test and its fixing
subgroup. It includes the cases \(t=0\) and \(t(v_i)=0\).

For the joined real polynomial, the first
block of variables has the transposed
indices of the first arrow and the second
has those of the second arrow. The sign
change of the first form makes the flipped
Fock variables holomorphic in the new basis.
Products of monomials of the two blocks span
all polynomials. Write \(z_{\rm def}\) for the nonzero scalar in the
theta-product identity on the definite triple \((X,X',X'')\). It is
fixed by its rational Schr\"odinger models, by the argument preceding
\eqref{eq:double-constant-term}. On a product polynomial and two
finite tests the residual theta kernel equals \(z_{\rm def}\) times
the product of the definite kernels. Integrating the first and then
the second compact coefficient therefore composes their two
contractions, multiplied by
\(z_{\rm def}\kappa_R\kappa_{R'}/z_\theta\) from
\eqref{eq:double-unfolded-integral}. Extending linearly in the finite
tensor gives the outgoing functional at \(F_i(\xi,y)\).

For a general Schwartz test after the slices,
approximate by Gaussian times polynomials in
the Schwartz topology, by the Hermite expansion.
All the integral evaluations are continuous in
that topology by the compact definite quotients
and absolute theta bounds. For each fixed source functional the possible families of target coefficients,
with the twists, tensors, forms and target points
as indices, form the image of the finite-dimensional
target Weyl dual in a product of copies of \(\mathbb C\).
Choose a complement of the common kernel of these evaluations in
the target dual. Finitely many coordinates give an injective map
on the complement, with a linear inverse on its image. Convergence
in these coordinates therefore gives a limiting functional in the
complement, and convergence in every coordinate gives the identity
for every finite datum. Its dependence on the real points and the
incoming functional is continuous by the same finite inverse and
the compact integral formulas. This also covers a zero-dimensional
image.

At the first pair the incoming functional is the original \(t\).
For the measured period the initial point is
\((g_{0,\infty},c)\) in the finite level trace of
Section~\ref{sec:measured-finite-sum}; the intermediate point is
the common argument of both kernels. At the outgoing place the
real argument can vary with the outgoing finite component held
at \(y_2\). Inductively the even-stage function on
\(U(X_{2i})\) is a coefficient of the first
\(i\) definite pairs with tensors
\(F_0,\ldots,F_{i-1}\). The coefficient
functional depends only on the outer real
points and the fixed datum. Its extension to
\(R_{2i}\) is determined by the unipotent
character. The terminal coefficient is
read at \(\operatorname{diag}(h,1)\), by the
last pair calculation. This statement is
pointwise before any outer integration.
It also applies to the support test for the
finite cuts.

At the terminal space let \(B\) be the
branching map from the Weyl space of
\(X_m\) to that of \(Y_m\). For a terminal
functional \(\ell_u\) and normalized highest
vector \(v_Y\), average
\[
 \int_{U(Y_m)_\infty}
   \sigma_Y(k)(v_Y\otimes\ell_u)
    \sigma_X(\operatorname{diag}(k^{-1},1))\,dk
          =c_{\ell_u} B.
\]
The average intertwines by Haar invariance. The terminal weights
interlace by Lemma~\ref{lem:terminal-branching-norm}. After a common
determinant twist, removing the boxes occupied by the last letter
leaves the smaller diagram, and the removed boxes form a horizontal
strip. Thus the compact branching multiplicity is one. Let \(z\) belong
to the original terminal type and use \(\overline z\) as the test
of the conjugate type in the preceding bilinear transpositions. With the
Hermitian product linear in its second argument,
\[
 \overline{z(k,h)}=\langle z_{\rm alg}(h),\sigma_Y(k)v_Y\rangle .
\]
If \(\psi\) is the upper-chain output, its coefficient at
\((\widehat k,\widehat h)\) is
\(\ell_u(\sigma_X(\widehat k^{-1})\psi(\widehat h))\).
Multiplication by \(\overline{z(k,h)}\) and integration in \(k\)
therefore insert precisely the displayed averaged map. Let \(J_B\)
be the finite class integral of
\[
     \big\langle z_{\rm alg}(h),B\psi(\operatorname{diag}(h,1))\big\rangle
\]
with maximal-compact mass one. At a common fixing level \(J\) this
is its weighted class sum divided by \([K_Y^\circ:J]\), where
\(K_Y^\circ\) is the terminal maximal compact. Here
\(R_m(F_0)=U(Y_m)(F_0)\), so the last quotient
\(R_m(F_0)\backslash U(Y_m)(\mathbb A_{F_0})\) has the ambient
Petersson measure. With probability Haar at infinity the terminal
integral against \(\overline z\) is \(c_{\ell_u}J_B\).
The restriction of Section~\ref{sec:upper-chain-size} is \(B/b_N\),
where \(b_N=\prod_\tau\tau(\bar\pi_N)^{B_\tau}\) and
\(B_\tau=a_\tau-mh_0\) is the exponent denoted by \(d_\tau\)
in that section.
With \(d_Y\) the full terminal Weyl dimension, Schur orthogonality for
the unit vector \(v_Y\) gives
\[
 \langle z,\operatorname{Res}\psi\rangle=(d_Yb_N)^{-1}J_B,
 \qquad c_{\ell_u}J_B=d_Y b_N c_{\ell_u}
                            \langle z,\operatorname{Res}\psi\rangle .
\]
Thus the terminal coefficient uses the same quotient measure and the
same finite representative on both sides.

It separates
into a factor of the real outer points
and the indicated linear expression in the
finite tensors. The two quotient-mass factors
of each pair are included in \(\ell_u\).

\subsection{Finite support on the original Bessel orbits}
\label{sec:marked-seesaw-support}

\sectionopening{In this subsection, we show that the period expression of
Section~\ref{sec:even-stage-coefficients} is supported on a compact
union modulo the original Bessel actions. The lattice inclusions below
give the compact sets used in the class sum of
Section~\ref{sec:measured-finite-sum}.}

For the space \(Y_j\), let \(L_j(x)\) be the lattice in
the chosen basis at the finite point \(x\), and
\(\Lambda_j=L_j(1)\). On its string subspace \(S_j\),
write \(L_j^-=L_j\cap S_j\) and
\(L_j^+=\operatorname{pr}_{S_j}L_j\).
The full lattices \(L_j\) and \(\Lambda_j\) are Hermitian self-dual.
Write \(\#\) for the dual on the nondegenerate string space. A string
vector pairs integrally against \(L_j\) precisely
when it pairs integrally against \(L_j^+\), so
\[
                        (L_j^+)^\#=L_j^-,
                        \qquad (L_j^-)^\#=L_j^+ .
\]
Use \(q_{\rm id}=b_{\rm id}/c_M\), where
\(b_{\rm id}\) is the full modulus idele of the lower
arrow and \(c_M\) is the original level generator.
Then \(q_{\rm id},\bar q_{\rm id}\) are integral and
are units above \(N\). Put \(a_E=\alpha_{\rm W}/\bar c_M\).

We first make a graded representative in a Bessel
orbit. At each local place take a basis \(B\) of
the moved lattice adapted to the decreasing grade
filtration, by extending successive bases of the
saturated intersections. Its degree-zero matrix part
\(D_0\) is invertible and \(BD_0^{-1}-1\) strictly
raises grades. Let \(G=B^*J B\) and \(G_0=D_0^*JD_0\),
where \(J\) is the rational Gram matrix.
The form pairs opposite grades.

Thus \(G-G_0\)
has entries with sum of the two grades negative.
After reversing the grading on one axis, \(G\) is
block triangular, with its diagonal blocks those of
the reversed \(G_0\). Since \(G,G^{-1}\) are integral,
the diagonal blocks and their inverses are integral.
Hence \(G_0\) is unimodular. Write \(G_0^{-1}G=1+P\);
\(P\) is self-adjoint for \(G_0\) and strictly raises
grades. At every finite place \(E/F_0\) is split or
unramified, so choose a local integral
\(\lambda\) with \(\lambda+\bar\lambda=1\).

For \(D=1-\lambda P\), self-adjointness gives
\[
 G_0^{-1} D^*G D
 =(1-\bar\lambda P)(1+P)(1-\lambda P)
 \equiv 1\pmod{P^2}.
\]
Thus a correction doubles the minimum grade increase of the error.
Each correction and its inverse are integral, since \(P\) is integral
and nilpotent. Iterating gives a unitriangular \(H\) with
\(H^*GH=G_0\). For \(d_j=\dim_EY_j\), the isometry \(BHD_0^{-1}\)
takes the graded lattice \(D_0\widehat{\mathcal O}_E^{d_j}\) to
\(L_j\) and acts trivially on the associated graded space. Its inverse
is the required finite Bessel move. Almost everywhere one starts with the
original graded basis and the move is the identity.

For \(0\le j\le m-1\) the string of \(Y_j\) has length
\(s=m-1-j=2\varrho+1-j\). In the basis of normalized string powers
the line ideals for \(0\le i<s\) are
\[
 \Lambda^-_i=p_{s,i}N^{-\mathbf1_{2i\ge s}}
                 \widehat{\mathcal O}_E,\qquad
 \Lambda^+_i=p_{s,i}N^{-\mathbf1_{2i+1\ge s}}
                 \widehat{\mathcal O}_E,\qquad
 p_{s,i}=a_E^{s-i}\bar a_E^{-1-i}.
\]
For \(2i<s\) the coefficient of the \(i\)-th power in the
original string basis is \(c^d\), with \(d=i-\lfloor s/2\rfloor\).
For \(2i\ge s\) it is \(-Nc^d\), with
\(d=i-\lfloor(s-1)/2\rfloor\). These signs make
\(T_{j,f}\) send each surviving power to the power of the same
index, by the three odd-string and two even-string coefficients
of Section~\ref{sec:adjacent-string-orbits}. Put
\(u=(a_E/\bar a_E)^{\lfloor(s+2)/2\rfloor}\).
In the original finite frame the string multiplier is \(u\)
when \(s\) is odd. For even \(s\) it is \(u/a_E\) on the
first half and \(u\bar a_E\) on the second half. Substitution
of the two displayed exponents \(d\) gives
\(p_{s,i}(a_E\bar a_E)^d\) in all cases.

The valuations of \(c\) and \(a_E\bar a_E\) agree,
so the intersection ideals are as displayed.
In odd length the central glued vector has
nonzero component in the last definite line, so its
coefficient must vanish in the string intersection.
With \(u=(a_E/\bar a_E)^{\lfloor(s+2)/2\rfloor}\),
its projection adds \(uz/N\) on the middle line,
giving the extra factor \(N^{-1}\) for the
projection exactly there. For \(s\ge1\) the shift kills the last power. For
\(i<s-1\), \(2i+1\ge s\) is equivalent to \(2i\ge s-1\), and
\(p_{s,i}=a_Ep_{s-1,i}\). Extending a shorter coefficient tuple by
zero proves surjectivity, so
\[
                        T_{j,f}\Lambda_j^+
                              =a_E\Lambda_{j+1}^-.
\]
The quotient \(\Lambda^+/\Lambda^-\) is therefore
zero outside \(N\) and for even length, and one
residue line at each place above \(N\) in odd
length.

\subsubsection{The lattice conditions on a finite slice}\label{sec:seesaw-finite-slices}

\sectionopening{In this subsection, we express the surviving slices of Section~\ref{sec:paired-lifts-rational-orbits} as inclusions of lattices.}

Consider a nonzero finite slice in the step from
\(Y_j\) to \(Y_{j+1}\). In the finite Hom model
write \(L_H\) for the self-dual vacuum lattice,
\(K_I\) for the basis of the positive-degree space
and \(z_T\) for the finite coordinate of \(T_{j,f}\).
In these coordinates use
\(\psi_H(t)=\psi_f(j_{\rm W}\alpha_{\rm W}\bar\alpha_{\rm W}t)\).
The signed Heisenberg operators obey
\[
 \rho_f(x)\rho_f(y)=\psi_H(y^*x)\rho_f(y)\rho_f(x),\qquad
 \rho_f(l)v_{\rm vac}=v_{\rm vac}\quad(l\in L_H).
\]
Let \(\mathfrak m\) be the ideal generated by \(b_{\rm id}\),
\(h=\dim_E\operatorname{Hom}_E(Y_j,Y_{j+1})\), and
\(R=(\mathcal O_E/\mathfrak m)^h\). For the prescribed inclusion
\(a_0\in R\) the test is
\[
 F=\frac1{|R|}\sum_{v\in R}
  \psi_{\mathfrak m}(-a_0\cdot v)\,\rho_f(t_v)v_{\rm vac},
 \qquad
 t_v=\overline{1_{\mathfrak m} b_{\rm id}^{-1}\widetilde v}.
\]
Here \(1_{\mathfrak m}\) selects the modulus places, \(\widetilde v\)
is an integral lift, and \(\psi_{\mathfrak m}\) is the primitive
residue character in the chosen split coordinates. Indeed the sum
of the characters \(\psi_{\mathfrak m}((a-a_0)\cdot v)\), divided
by \(|R|\), is the indicator of \(a=a_0\). Since
\(\bar b_{\rm id}t_v\in L_H\), the commutator of \(\rho_f(b_{\rm id}l)\)
with \(\rho_f(t_v)\) is one for every \(l\in L_H\).
Therefore \(\rho_f(b_{\rm id}l)F=F\).

Let \(g\) be the unitary Hom action of the two finite points, set
\(f_1=\omega_f(g)F\), \(M_1=b_{\rm id}g\) and
\(M_2=\bar b_{\rm id}^{-1}g\). Then \(M_1 M_2^*=1\) and
\(\rho_f(M_1l)f_1=f_1\) by covariance. In the mixed model let
\(\mathcal Z\) denote restriction to zero in the coordinate dual to
the positive-degree space \(I\). With probability Haar on
\(I(F_0)\backslash I(\mathbb A_{F_0})\), character orthogonality
selects exactly this coordinate. On the real and finite tensor the
restrictions factor, so nonvanishing of the adelic constant term
implies \(\mathcal Z(\rho_f(z_T)f_1)\ne0\).
For \(u\in I\) the signed translation coincides with the Heisenberg
translation used to define the zero slice. Indeed its quadratic
normalizing phase is evaluated on \(\operatorname{tr}(u^\dagger u)\),
as in Section~\ref{sec:upper-cm-jet-transport}. Adjunction preserves
the degree, since the Hermitian forms pair opposite graded pieces.
Thus \(u^\dagger u\) strictly raises grades and has trace zero,
including at the places over two. In the mixed model \(\rho_f(u)\)
is modulation by a character that is one at the zero coordinate.
Hence \(\mathcal Z\rho_f(u)=\mathcal Z\). If also \(u\in M_1L_H\),
then \(\rho_f(u)f_1=f_1\) by the lattice covariance above.
Applying \(\mathcal Z\) to
\[
 \rho_f(z_T)\rho_f(u)f_1
   =\psi_H(u^*z_T)\rho_f(u)\rho_f(z_T)f_1
\]
gives \(\psi_H(u^*z_T)=1\).
Choose columns \(C\) so that \(M_2^*K_I C=B_I\)
is an integral basis of the saturation of its
image in \(L_H\).

This is possible over every
local DVR, with the original basis almost everywhere.
Let \(J_I B_I=1\) be an integral left inverse.
For every integral column \(x\), the vector
\(K_I Cx=M_1B_Ix\) is in \(I\cap M_1L_H\). Thus
\(\psi_H(x^* C^*K_I^*z_T)=1\) for all such \(x\), and self-duality
gives \(v=C^*K_I^*z_T\) integral. Taking \(\ell_I=J_I^*v\)
gives \(C^*K_I^*M_2\ell_I=B_I^*J_I^*v=v\). The adjoint uses the
transposed matrix at the conjugate place of a split pair.
Cancelling \(C^*\) yields
\[
                         K_I^*(z_T-M_2\ell_I)=0 .
\]
In the original rational bases the two terms
are \(\alpha_{\rm W}^{-1}T_{j,f}\) and
\(\bar b_{\rm id}^{-1}W\), where \(W\) is the
integral Hom matrix \(\ell_I\) transported to the
actual lattices \(L_j,L_{j+1}\).

Orthogonality to \(I\) means all negative-degree
coordinates of the difference vanish. Thus
\(W-(\bar b_{\rm id}/\alpha_{\rm W})T_{j,f}\)
has degree at least zero. Move the two lattices by Bessel
elements \(r,r'\) to graded representatives. The degree-zero actions
are the identity on each string and preserve the definite blocks.
Since \(T_{j,f}\) kills the definite block and has degree \(-1\)
on the string, \(r'T_{j,f}r^{-1}-T_{j,f}\) has degree at least
zero. Write \(P_d,P'_d\) for the source and target grade projections.
The map \(r'Wr^{-1}\) is integral between the moved lattices. Its
degree-\(-1\) part
\[
 \sum_d P'_{d-1}\,r'Wr^{-1}P_d
     =(\bar b_{\rm id}/\alpha_{\rm W})T_{j,f}
\]
is still integral since these projections preserve the graded
lattices. Factoring through the source string projection and taking
the target string intersection gives
\[
        (\bar q_{\rm id}/a_E)T_{j,f} L_j^+
                                      \subset L_{j+1}^- .
\]
This computation concerns the original finite
factors at these points.

\subsubsection{Exclusion by the ramified character}\label{sec:seesaw-ramified-exclusion}

\sectionopening{In this subsection, we use the source condition of Section~\ref{sec:ramified-source} in the lattice calculation of Section~\ref{sec:seesaw-finite-slices}.}

Let the initial finite point \(y_0\) lie in the chosen level,
so \(L_0(y_0)=\Lambda_0\). If \(I_{j+1}(a,x)\ne0\) at an
archimedean and finite pair, the integral in
Section~\ref{sec:even-stage-coefficients} has a point \(r\) at
which its two factors are both nonzero. The source finite point is
\(r_f y_j\). For \(j=0\), applying \(r_f^{-1}\) recovers the
original graded lattice. For \(j>0\) choose a graded representative
of the source orbit and use the inductive bound.
Keep an arbitrary graded representative of the target orbit.
Duality takes the source bound
\(q_{\rm id}^{j}L_j^-\subseteq\Lambda_j^-\) to
\(\bar q_{\rm id}^{j}\Lambda_j^+\subseteq L_j^+\).
Apply the finite slice inclusion and the standard shift formula,
then dualize on the target, to obtain
\[
 \bar q_{\rm id}^{\,j+1}\Lambda_{j+1}^-\subseteq L_{j+1}^-,
 \qquad q_{\rm id}^{\,j+1} L_{j+1}^+\subseteq\Lambda_{j+1}^+.
\]
A failure of the desired upper bound for
\(L_{j+1}^-\) can occur only at a place above
\(N\) on an odd string. There \(q_{\rm id}\) is a
unit, and the image of a vector witnessing failure generates the
residue line of \(\Lambda^+/\Lambda^-\).
Both inclusions then force \(L^-=L^+=\Lambda^+\)
at that place. At the conjugate place duality
also gives \(L^-=L^+\). In particular the moved
lattice splits orthogonally into its definite
and string parts at the pair of places.

The split option occurs only at an even
stage \(l=j+1<m\).
Let \(w,\bar w\) be its pair above \(N\), choosing
\(w\mid\pi_N\). Let \(z_w\) be the original
norm-one scalar unit supported at that pair,
with primitive reduction at \(w\). Its embedding
\(Z\) acts by that scalar on the definite block
of \(Y_l\) and by one on the string.
The splitting of \(L_l(x_*)\) for the moved
representative \(x_*\) implies
\(x_*^{-1}Zx_*\in J_{\xi,l}\):
it is integral at the split pair and the identity
elsewhere, and the pair is disjoint from the
chosen level and twist places.

Use \(x_*\) at the \(l\)-th outer point
in the even-stage coefficient formula of
Section~\ref{sec:even-stage-coefficients}.
The last lower finite test of that stage,
after the outer move, is fixed by
\(\omega(1,Z)\) by the conjugate level condition.
The exact slice covariance passes this action
to multiplication by \(z_w\) on the last
definite group. Each earlier definite arrow
intertwines scalar translation by \(z_w\) on
its source and target. Indeed their simultaneous
matrix action on Hom is trivial and the finite
determinant characters are unramified, hence
trivial on these units.

In the last arrow the
translation can be moved to the finite test,
where it fixes the evaluated tensor. By linearity
the resulting partial chain \(L\) therefore
satisfies \(L\rho_X(z_w)=L\) on the original
source forms. For \(\phi\) in the ramified
source, the sum of its translates over the
residue-unit action at \(w\) is zero. Thus
\[
       (N-1)L\phi=L\sum_{b=0}^{N-2}\rho_X(z_w)^b\phi=0 .
\]
Thus \(L\phi=0\) in the characteristic-zero space of definite forms.
The same conclusion holds with \(x_*\) replaced by any finite
Levi translate \(h x_*\), since the scalar \(Z\) on the definite
block commutes with that Levi. Every finite Bessel translate is
\(u h x_*\), with \(u\) unipotent. At this even stage the string
has odd length, so \(u-1\) has degree at least two. Successor
covariance multiplies the value at \(h x_*\) by its unipotent
character. Thus the even-stage coefficient vanishes on the entire finite
Bessel orbit at every real point, before outer integration.
If \(x_*=r_*x\) was the graded representative, take the finite
Bessel translate \(r_*^{-1}\). This gives zero at the original
point \((a,x)\), contradicting \(I_l(a,x)\ne0\).

\begin{lemma}[The support on the ramified source]\label{lem:seesaw-ramified-sandwich}
Suppose the source form is ramified and the initial finite point
lies in the chosen level. For \(1\le j\le m-1\), if \(I_j(a,x)\ne0\),
then every finite Bessel translate of \(x\) with graded lattice
\(L_j\) satisfies
\begin{equation}
          \bar q_{\rm id}^{\,j}\Lambda_j^-
                       \subseteq L_j^-\subseteq
               q_{\rm id}^{-j}\Lambda_j^- .
 \label{eq:marked-seesaw-sandwich}
\end{equation}
\end{lemma}

\subsubsection{Compact representatives of the surviving terms}\label{sec:seesaw-support-representatives}

\sectionopening{In this subsection, we use Lemma~\ref{lem:seesaw-ramified-sandwich} to choose representatives for the measured sum.}

To pass from the string lattices to compact
representatives, suppose \(L(x)\) and \(L(y)\)
have the same string intersection. It is
saturated in both lattices. Its local bases
therefore give two primitive bases in the
standard coordinates, with equal Gram matrices.
At a split pair write the bases as
\((A,A')\), \((B,B')\). The equations for one
matrix are \(gA=B\) and \(B'^{\mathsf T}g=A'^{\mathsf T}\).
Over the residue field, first identify the
basis spans with their given bases and extend
the identification on their intersections with
the two kernels to an isomorphism of the kernels.

The two maps \(A'^{\mathsf T}\), \(B'^{\mathsf T}\)
are surjective and agree on these identified
spans by Gram equality. Complete with lifts
of the remaining quotient basis. This gives
an invertible residue solution. Over the DVR
choose \(DA=1\), \(B'^{\mathsf T}R=1\) and put
\[
   g_0=BD+RA'^{\mathsf T}(1-AD).
\]
Lifting the difference of the residue solution
and \(g_0\) to \(H\), replace \(g_0\) by
\(g_0+(1-RB'^{\mathsf T})H(1-AD)\).
The correction is homogeneous for both equations
and gives the chosen invertible reduction.
The conjugate matrix is \(g^{-\mathsf T}\).

At a nonsplit unramified place choose a uniformizer
\(\pi\) of the fixed field; it is also a uniformizer of the
quadratic extension and \(\bar\pi=\pi\). Finite-field Hermitian Witt
extension matches the two primitive residue bases, allowing
degenerate spans. Lift the isometry by successive Gram corrections. At each residue step the
equation \(C+C^*=H\) for Hermitian \(H\) is
solvable using residue trace on the diagonal
and one triangle off it. The trace is
surjective also in residue characteristic two.
If the unitary basis images agree modulo
\(\pi^a\), their discrepancy \(D\) at the
next order satisfies \(A^*D+D^*A=0\).

The prescription \(CA=D\) extends to a
skew-adjoint matrix: prescribe the resulting
skew-Hermitian form on the basis and fill
the complementary block with zero.
Then \(1+\pi^aC\) can be corrected to a
unitary lift starting at order \(a+1\),
so the bases agree one order further.
The limit agrees exactly. The resulting
finite integral isometry \(k\) has
\(y k x^{-1}\) fixing the entire string,
since it fixes its full lattice basis.
It therefore lies in the Bessel Levi.
Thus a common string intersection determines
the Bessel orbit modulo the standard compact.

The bounds \eqref{eq:marked-seesaw-sandwich}
give only finitely many string intersections,
all equal to the standard ones off \(Q\cup\bar Q\), where
\(Q\) is the set of twist places: the quotient of the two bounding
lattices is finite and supported on this union. Start with a graded representative for each occurring
intersection. At a place outside \(Q\cup\bar Q\), its string
intersection agrees with the one at the identity. The local
frame comparison just proved gives a Bessel Levi element on the
left and a standard compact element on the right that take this
component to the identity. At almost every place the component
already belongs to the standard compact, so one may use the
identity Levi element there. At the twist places keep the given
components. The local choices assemble to finite adelic points,
and the resulting representative remains graded everywhere.
Its right standard compact coset meets the same finite Bessel
orbit. This use of \(R_j(\mathbb A_{F_0,f})\) is inside the
integration in stages over the adelic Bessel quotient in
Section~\ref{sec:even-stage-coefficients}. This gives a
compact union at each intermediate stage
in the original basis. At stage zero
use the classes of the original chosen level
modulo \(J_{\xi,0}\).

\section{The measured theta tensor and its integral lattice}
\label{sec:measured-theta-tensor}

\sectionopening{In this section, we evaluate the finite sum from Section~\ref{sec:consecutive-theta-lifts}
with the prescribed Haar measures. The resulting tensor and its comparison
at the additional ramified prime give the local identities used to bound derivatives
in Section~\ref{sec:local-theta-derivatives}.}

\subsection{The measured finite sum}
\label{sec:measured-finite-sum}

\sectionopening{In this subsection, we integrate the paired formula of
Section~\ref{sec:even-stage-coefficients} over the finite classes
furnished by Section~\ref{sec:marked-seesaw-support}. The computation
retains the stabilizer factors and shows how the compact indices
cancel.}

For each intermediate \(U(Y_j)\), use all string
intersections satisfying \eqref{eq:marked-seesaw-sandwich}
that occur for a graded moved lattice. The bounds depend only
on the level and the twist. By
Section~\ref{sec:marked-seesaw-support} this gives a finite
union of standard compact cosets. Its image modulo \(J_{\xi,j}\)
is finite since \(J_{\xi,j}\) has finite index in \(K_j^\circ\).
Choose one class \(xJ_{\xi,j}\) in this image from each
\(R_j(A_f)\)-orbit meeting it. These choices work for every
source form, terminal test and set of outer points. Write \(B_j=R_j(A_f)\) and \(B_j^0=B_j\cap K_j^\circ\).
Let \(\Omega_j\) be the union of the \(B_j^0\)-orbits
of the selected classes. The unitary Levi acts on
the Lie algebra of the unipotent radical with determinant
one: that determinant is a rational character, and
a unitary group has no nontrivial character over \(F_0\).
The Levi and the unipotent radical are unimodular,
so \(R_j\) is unimodular at each place. Take finite
Haar with \(\operatorname{vol}(B_j^0)=1\), rescaling
the real Haar reciprocally to keep the original
Bessel Haar. At an anchor \(xJ_{\xi,j}\) put
\[
 d_j(xJ_{\xi,j})=
 [R_j(A_f)\cap xJ_{\xi,j}x^{-1}:
     R_j(A_f)\cap xJ_{\xi,j}x^{-1}\cap K_j^\circ] .
\]
The Haar measure of the \(r\) with
\(r xJ_{\xi,j}\in\Omega_j\) is this integer on the
represented orbit.

Indeed, with \(Z_c=B_j\cap xJ_{\xi,j}x^{-1}\) at
the anchor \(c=xJ_{\xi,j}\), this set is \(B_j^0 Z_c\).
Its disjoint right translates of \(B_j^0\) are
indexed by \((B_j^0\cap Z_c)\backslash Z_c\).
For a starting point \(b c\) in the same \(B_j\)-orbit
the set is \(B_j^0Z_c b^{-1}\), of the same measure.
Use this anchor's integer on the whole orbit.
If the twist is empty all compact representatives
are integral and \(d_j=1\).
In general \(d_j\) divides a power of the product
of the residue characteristics of the twist places.

Indeed, at each twist place \(w\), let \(K_w(h)\) denote the
principal kernel of depth \(h\) in the integral matrix group.
Choose \(h_w\) with \(x_wK_w(h_w)x_w^{-1}\subseteq K_{j,w}^\circ\).
The level \(J_{\xi,j,w}\) is contained in \(K_w(1)\). The kernel of
\[
 Z_c\longrightarrow
 \prod_{w\mid Q(\xi)}K_w(1)/K_w(h_w),\qquad
 z\longmapsto (x_w^{-1}z_w x_w)_w ,
\]
is contained in \(B_j^0\cap Z_c\), since \(x\) is integral away
from these places. Each successive quotient \(K_w(h)/K_w(h+1)\)
embeds in the additive matrix group over the residue field by
\(1+\varpi_w^hD\mapsto\bar D\). Thus the index \(d_j\) is supported
on the residue characteristics of the twist places.

Choose a nonnegative real cutoff \(\alpha_j\) with
\(\int_{R_j(A_\infty)}\alpha_j(r u)\,dr=1\)
on each real orbit. Such a cutoff comes from the proper
translation action of the closed Bessel subgroup.
The cutoff depends only on the fixed inclusion \(R_j(A_\infty)\subset U(Y_j)(A_\infty)\). Then
\[
        w_j(u,x)=\alpha_j(u)
             \frac{\mathbf1_{\Omega_j}(xJ_{\xi,j})}
                  {d_j(xJ_{\xi,j})}
\]
has full Bessel orbit integral one on the support
of the unfolded integrand. Outside \(\Omega_j\) the
fraction is zero. Quotient integration in stages
can therefore be performed with \(w_j\), using
counting measure for the rational subgroup and
the inherited Haar for the adelic groups.
Indeed, for a nonnegative function \(F\) on
\(R_j(F_0)\backslash G_j(\mathbb A_{F_0})\), integration first
over \([R_j]\) and then over \(G_j(\mathbb A_{F_0})\) gives
\[
 \int w_j(g)\int_{[R_j]}F(rg)\,dr\,dg
 =\int_{R_j(F_0)\backslash G_j(\mathbb A_{F_0})}
       F(h)\left(\int_{R_j(\mathbb A_{F_0})}w_j(b^{-1}h)\,db\right)dh .
\]
To see this, write \(g=b\dot h\) using the quotient integral formula
over \(R_j(\mathbb A_{F_0})\backslash G_j(\mathbb A_{F_0})\).
Right invariance on \([R_j]\) replaces \(F(rb\dot h)\) by
\(F(r\dot h)\) under its integral. The cutoff integral over \(b\)
is constant on the remaining orbit. Tonelli applies to the nonnegative
integrand. Unimodularity and the two cutoff integrals identify the
parenthesis with the indicator of the represented finite orbits.
Those contain the support of \(I_j\). For \(1\le j<m\), with
the earlier outer points fixed and \(z=L_{\xi,j}^{\,\mathrm t}z'\)
the rapid transposed test, apply this identity to \(F=|I_jz|\).
Its right side is \(\int_{R_j(F_0)\backslash U(Y_j)(A)}|I_jz|\),
finite by Proposition~\ref{prop:adjacent-paired-unfolding} at the
preceding arrow. Fubini gives the same identity for \(I_jz\).

Write \(\Xi_j(z)\) for the integral of \(I_j(h)z(h)\)
on \(R_j(F_0)\backslash U(Y_j)(A)\).
Its remaining outer points are understood as
parameters. The product is right \(J_{\xi,j}\)-invariant
when \(z=L_{\xi,j}^{\,\mathrm t}z'\), by the original
level invariance of the last datum in \(I_j\)
and the transpose. Factor the given adelic Haar on \(U(Y_j)\)
as \(du\,dg_f\), with \(dg_f(K_j^\circ)=1\).
Any rescaling of \(dg_f\) is compensated by the reciprocal
rescaling of \(du\). In the integral with \(w_j\) a coset of
\(J_{\xi,j}\) has finite measure \(1/A_{\xi,j}\).
The one-step identity
of Section~\ref{sec:even-stage-coefficients}
has multiplier \(A_{\xi,j}\). These two factors
cancel, leaving exactly
\begin{equation}
 \Xi_j(L_{\xi,j}^{\,\mathrm t}z')=
  \int_{U(Y_j)(A_\infty)}\!\!\alpha_j(u)
       \sum_{c\in\Omega_j}
        \frac{\Xi_{j+1}(z';\,g_j=(u,c))}{d_j(c)}\,du .
 \label{eq:seesaw-cut-peeling}
\end{equation}
For \(1\le j<m\) set \(z=L_{\xi,j}^{\,\mathrm t}z'\).
For fixed earlier outer points the one-step identity at \(g_j\) gives
\[
 A_{\xi,j}\bigl|\Xi_{j+1}(z';g_j)\bigr|
  \le \int_{[R_j]} |I_j(rg_j)z(rg_j)|\,dr .
\]
Integrate this inequality against \(w_j(g_j)\,dg_j\).
On the left, each finite coset contributes
\(A_{\xi,j}^{-1}\). On the right, quotient integration
and the orbit integral of \(w_j\) give
\[
 \int_{U(Y_j)(A_\infty)}\alpha_j(u)
     \sum_{c\in\Omega_j}
       \frac{|\Xi_{j+1}(z';\,g_j=(u,c))|}{d_j(c)}\,du
 \le \int_{R_j(F_0)\backslash U(Y_j)(A)}
                  |I_j(h)z(h)|\,dh <\infty .
\]
The last integral is absolute by
Proposition~\ref{prop:adjacent-paired-unfolding} at the preceding
step, since the transposed test \(z\) is rapid.
For \(j=0\), apply the measured
level formula of Section~\ref{sec:period-measured-pairing} with
\(J_{\xi,0}\subseteq J_M\) to the first transpose.
It gives the sum over \(c\in J_M/J_{\xi,0}\), divided
by \([J_M:J_{\xi,0}]\), of the original periods
at the right translate \(g_0(1,c)=(g_{0,\infty},c)\),
where \(g_{0,\infty}\) is the real point of Section~\ref{sec:period-real-point}.
After unfolding, the factor \(A_{\xi,0}\) leaves
\([K_0^\circ:J_M]\) because indices multiply.
It is independent of \(\xi\).

By the terminal calculation of
Section~\ref{sec:even-stage-coefficients}, the last
period is a product of a function \(a_m(u)\) of the
outer real points and the terminal pairing on the
finite tensor
\(\bigotimes_i F_i(\xi,y)\). At each step of
\eqref{eq:seesaw-cut-peeling}, update the finite
factor by the sum over \(c\in\Omega_j\) with weight
\(d_j(c)^{-1}\). With the earlier real variables
fixed, write \(a(u)b_c\) for the factored integrand
before multiplication by \(\alpha_j(u)d_j(c)^{-1}\).
Set the new real factor to \(\int\alpha_j(u)a(u)\,du\)
when this integral is absolute, and set it to zero
otherwise. This defines a measurable function of the earlier
real points: the integral of \(\alpha_j|a|\) is measurable by
Tonelli, and on its finite locus the signed real and imaginary
parts give the parameter integral. Absolute integrability in
\eqref{eq:seesaw-cut-peeling} at any coset with \(b_c\ne0\)
implies \(\int|\alpha_j a|<\infty\).
In the divergent case every \(b_c\) therefore vanishes
and the same factorization holds. The rule uses
only the fixed real data and works for all twists
and finite tests. Applying it successively proves
the identity successively at the intermediate groups.

Define the evaluated tensor by the weighted sum
\begin{equation}
 T_\xi=\sum_{c_0,\ldots,c_{m-1}}
       \left(\prod_{j=1}^{m-1}d_j(c_j)^{-1}\right)
                    \bigotimes_{i=0}^{m/2-1} F_i(\xi,c),
 \label{eq:seesaw-evaluated-tensor}
\end{equation}
where \(c_0\) runs over \(J_M/J_{\xi,0}\),
\(c_j\in\Omega_j\) for \(j>0\), and the terminal
finite point is one. Each finite tensor uses the
same actual representative as its outer period.
The real cutoffs, integrals, and point \(g_{0,\infty}\) are fixed
before \(\xi\). Define \(c_{\rm meas}\) as the iterated real factor times
the fixed index \([K_0^\circ:J_M]\).
The terminal test \(\overline z\) may have any finite level.
The finite invariance used in \eqref{eq:seesaw-cut-peeling} is that
of \(L_{\xi,j}^{\,\mathrm t}z'\) at \(J_{\xi,j}\) for \(j<m\),
which follows from the inclusion test in its source variable.
Both terminal outputs have the prescribed real type, by the
coefficient identity and the terminal branching map. They are
smooth at the finite places, so their difference is fixed by a
common compact open subgroup. Take \(z\) to be that difference.
The terminal quotient is compact, hence this is an admissible test,
and the equality of pairings gives zero for its squared Petersson
norm. This proves equality of the terminal forms.

\begin{proposition}[Measured seesaw identity]\label{prop:measured-seesaw}
Let \(\phi\in S\) and let \(\xi\) be an allowed twist. With the fixed
real point and original finite levels, the measured lifts satisfy
\begin{equation}
           L_\xi\Pi\phi
                  =c_{\rm meas}\,\mathcal U_\xi(T_\xi)\phi .
 \label{eq:seesaw-measured-chain}
\end{equation}
The scalar \(c_{\rm meas}\) is independent of \(\xi\).
\end{proposition}
For the trivial twist the number of choices of each \(c_j\) is
at most \(M^{f n_m^2}=B_*\), since the classes
lie in \(K_j^\circ/J_{\xi,j}\) and the principal
index there is a product of \(|\mathrm{GL}_{\dim Y_j}(k_w)|\).
The sum of the cut weights is at most \(B_*^m\)
in that case. The expansion of the \(F_i(\xi,c)\) in translated
vacua is given in Section~\ref{sec:seesaw-integral-tensor}.

\subsection{The lattice of an evaluated pair}
\label{sec:marked-neighbor}

\sectionopening{In this subsection, we calculate the lattice obtained when a summand of
Section~\ref{sec:measured-finite-sum} is restricted to the paired
oscillator variables. At the places over \(N\), its two components
satisfy a common residue condition that determines the sum in
Section~\ref{sec:seesaw-integral-tensor}.}

Use the finite bases in a pair from
\(Y_{2i}\) through \(Y_{2i+2}\). Put
\(R_f=\widehat{\mathcal O}_E\) in this
calculation, let \(\pi=\pi_N\) as a
principal idele, and write \(J_N\) for
the idempotent of the places dividing
\(\pi\). Thus \(J_N\bar J_N=0\).
Use \(\eta_-=\eta\), the norm-\(-1\) unit from the opposite finite
basis of Section~\ref{sec:double-string-unfolding}. Its basis scalar
is \(s_-\) below. Put
\[
 \delta=\frac{c_M}{\bar c_M},\qquad
 s_-=\bar\eta_-\frac{\bar\alpha_{\rm W}}{\alpha_{\rm W}},
 \qquad u_j^{\rm str}=(a_E/\bar a_E)^{\lceil(m-j)/2\rceil}.
\]
With \(\beta=\pi_M\) the integral level
generator in valuation, the identities are
\[
 s_-u_{2i}^{\rm str}=\bar\eta_-\delta\,u_{2i+2}^{\rm str},
 \qquad v_w(\bar\beta\delta)=v_w(\beta),
 \qquad \delta=1\text{ above }N .
\]
They follow from \(a_E=\alpha_{\rm W}/\bar c_M\)
and the successive two-step exponents.
In an odd string the lattice contains the
glued vector
\(e/\bar\pi+(u_j^{\rm str}/N)z\)
and \(u_j^{\rm str}z\), where \(e,z\)
are the last definite and central string coordinates. The coefficients of a vector with ordinary
coordinates \(x\) and string coordinate
\(u_j^{\rm str}(\pi^{-1}x_e+d)\)
in this pair of lattice vectors are
\(\bar\pi x_e,d\). Thus membership is tested
by those two integral values and the other
ordinary integral coordinates. The two-by-two
Gram block is
\(\left(\begin{smallmatrix}0&-1\\-1&-N\end{smallmatrix}\right)\).
The even string has unglued definite coordinates.
The lattices in this calculation are over
\(R_f=\prod_w\mathcal O_{E_w}\). At each conjugation orbit a finite
frame identifies the original lattice with a free module over its
integer rings. The outer points of \eqref{eq:seesaw-evaluated-tensor}
act on the two tests before taking their zero slice. The matrices
below are in the fixed original frames. All lattice inclusions can
therefore be checked on the completions, also when the intersection
of a class lattice with its rational space is a nonfree
\(\mathcal O_E\)-module.

Let \(e,e''\) be the last indices of the two
definite blocks in the target of the joined space.
Define
\[
 \begin{split}
 v_0&=e+\bar\eta_-\,\bar J_N e'',\\
 D_N&=J_N(\pi-1)+\bar J_N(\bar\pi^{-1}-1),\\
 \Gamma&=1+D_Nv_0v_0^\dagger,\\
 b_e&=D_N\pi^{-1}\bar J_N,\qquad
 b_{e''}=D_N\pi^{-1}\eta_-J_N .
 \end{split}
\]
Set \(b_j=0\) for the other indices of this
block. The vector \(v_0\) has norm one, and
\(\Gamma\) acts on its line by the norm-one
idele \(1+D_N\). In order \(e,e''\) its
two blocks are
\[
  \begin{pmatrix}\pi&(\pi-1)\eta_-\\0&1\end{pmatrix}
  \quad(w\mid\pi),\qquad
  \begin{pmatrix}\bar\pi^{-1}&0\\
    (\bar\pi^{-1}-1)\bar\eta_-&1\end{pmatrix}
  \quad(w\mid\bar\pi),
\]
and it is the identity elsewhere.
For the column \(x=\Gamma e_j\) with
\(b=b_j\), direct subtraction gives
\[
 \bar\pi x\in R_f^{\dim x},\qquad
 b-\pi^{-1}x_e\in R_f,\qquad
 \bar\eta_-b-\pi^{-1}x_{e''}\in R_f,\qquad
 \delta b=b .
\]
The first and second differences are
\(-[J_N+(1-J_N)\pi^{-1}]\) in the column \(e\)
or \(e''\), respectively, and zero in the
other columns. The bracket is integral.

For an interior pair put \(u=u_{2i}^{\rm str}\).
Lift the column \((x_-,x_+)\) to ordinary
coordinates \((s_-x_-,x_+)\) and an additional
component \(s_-u b(z+z'')\) on the removed
plane. Dividing the first block by its
opposite-basis scale \(s_-\) leaves string
coordinate \(u b\). In the other block it is
\(u_{2i+2}^{\rm str}\bar\eta_- b\).
Both are integral in the required glued lattice
by the two difference congruences. In the
terminal pair the extra target index \(e''\)
is represented by \(z/\pi\) on the negative
source side. Use the diagonal scaling
\(\kappa_{e''}=\delta\), with all other
\(\kappa\)'s one. The corresponding original
source string coordinate is then
\(\pi^{-1}\delta x_{e''}\).
Since \(s_-^{-1}\delta=-\eta_-u\), after
dividing the first block by \(s_-\) its
integrality test is
\[
 -\pi^{-1}x_e-\eta_- \pi^{-1}x_{e''}
 =(b-\pi^{-1}x_e)+\eta_-(\bar\eta_-b-\pi^{-1}x_{e''})
 \in R_f .
\]
The remaining ordinary target entries are
integral. Thus the same lifted procedure
works in both pairs.

Let \(C\) have these lifted columns in the
original finite coordinates. Let \(x_f\) be the finite basis matrix of \(Y_{2i+1}\)
in its graded coordinates, and \(j_0\) the inclusion of its
definite summand. Its definite-coordinate projection is
\[
       R'=j_0^\dagger x_f,\qquad (R')^\dagger=x_f^{-1}j_0 .
\]
Here \(j_0^\dagger\) extracts the definite coordinates, and the
second identity uses the isometry equation for \(x_f\).
The even middle string has unglued definite coordinates, so
\(R'\) is integral.
On matrices the lift \(B:L_{\rm def}\to L_{\rm big}\)
is \(X\mapsto CXR'\), hence integral.

If \(i_f\)
is the orthogonal residual inclusion and \(k_f\)
the isotropic basis in the two slices, its
columns give
\[
       B=i_f g_{\rm nb}+k_f T_{\rm nb},\qquad
       g_{\rm nb}=(\operatorname{diag}\kappa\ \Gamma)\otimes 1 .
\]
Here \(T_{\rm nb}\) is supported above \(N\),
and \(g_{\rm nb}\) is unitary since the diagonal
scaling is one above \(N\) and has norm one.
Write \(I=\operatorname{im}k_f\), and let
\(L_{\rm big}\), \(L_{\rm def}\) denote the
standard finite vacuum lattices.

The quotient
image of \(L_{\rm big}\cap I^\perp\) in
\(I^\perp/I\) is exactly \(g_{\rm nb}L_{\rm def}\).
One inclusion follows from the integral
columns of \(B\). Conversely, for \(x\in L_{\rm def}\) and
\(y\in L_{\rm big}\), the pairing of \(x\) with \(B^\dagger y\)
equals the integral pairing of \(Bx\) with \(y\). Self-duality
of \(L_{\rm def}\) gives \(B^\dagger L_{\rm big}\subseteq L_{\rm def}\).
Since \(i_f^\dagger i_f=1\) and
\(k_f^\dagger i_f=i_f^\dagger k_f=k_f^\dagger k_f=0\),
\[
 B^\dagger(i_f z+k_f v)=g_{\rm nb}^\dagger z=g_{\rm nb}^{-1}z .
\]
Thus every integral lift has \(z\in g_{\rm nb}L_{\rm def}\).

Write \(S_f\) for the joined string shift in the finite
vacuum coordinates, namely \(\alpha_{\rm W}^{-1}\)
times the stacked rational shift in the chosen basis.
It satisfies
\(B^\dagger S_f=0\), because the two string maps vanish on the definite
summand of the middle space, so \(S_f(R')^\dagger=0\). For \(y\) vanishing above
\(N\) one therefore has
\[
                 g_{\rm nb} B^\dagger(S_f+y)=i_f^\dagger y .
\]
On the ambient source and target spaces the
denominators are \(\beta,\bar\beta\),
respectively. At the finite places away
from \(N\), \(\Gamma=1\) and the lifted
columns lie in the corresponding block of
the original inclusion. The additional
diagonal scaling at the terminal index converts
a source bound by \(\beta\) into the target
bound by \(\bar\beta\), by
\(v_w(\bar\beta\delta)=v_w(\beta)\).
More explicitly, let \(d_i\) be \(\beta\) or \(\bar\beta\)
according to the ambient block of coordinate \(i\), and let \(d'_b\)
be the corresponding denominator at residual coordinate \(b\).
At \(w\nmid Nq_\xi\), each nonzero entry
\((B^\dagger)_{bi,w}\) has
\(v_w(d'_b\kappa_b)=v_w(d_i)\) by the block calculation above.
In the split case the entry uses the transpose at \(\bar w\).
Since \(g_{\rm nb}\) acts by \(\kappa_b\) in this coordinate,
\[
 d'_b(g_{\rm nb} B^\dagger y)_{b,w}
 =\sum_{i:\,(B^\dagger)_{bi,w}\ne0}
       (d'_b\kappa_b/d_i)_w\,(B^\dagger)_{bi,w}\,(d_i y_i)_w
       \in\mathcal O_{E_w}
\]
whenever \((d_i y_i)_w\in\mathcal O_{E_w}\) for all \(i\).
At \(w\mid N\) both generators are units and the retained shift
vanishes.

\subsection{Integral coefficients of the measured tensor}
\label{sec:seesaw-integral-tensor}

\sectionopening{In this subsection, we expand the tensor of
Section~\ref{sec:measured-finite-sum} using the lattice of
Section~\ref{sec:marked-neighbor}. The resulting translates have
controlled denominators at coefficient places and an absolute
coefficient bound at complex embeddings.}

Put \(d^{\rm H}_j=\dim(Y_j)\dim(Y_{j+1})\), and let
\(\mathfrak m_\xi=(\beta)\prod_{w\in Q_\xi}\mathfrak p_w\)
be the ideal of the lower residue test. The Fourier
expansion of its one prescribed reduction in
\((\mathcal O_E/\mathfrak m_\xi)^{d^{\rm H}_j}\) gives one
character per element of the dual group, with
coefficient \(\zeta/|\mathcal O_E/\mathfrak m_\xi|^{d^{\rm H}_j}\)
for a root of unity \(\zeta\). Multiplying the finite
vacuum by such a character is its Heisenberg translate
in the conjugate split half. Thus the translation
has denominator \(\bar\beta\) away from \(Q_\xi\)
and its conjugates. At a coefficient prime away from
the twist its coefficient denominator divides
\(|\mathcal O_E/(\beta)|^{d^{\rm H}_j}\). When \(Q_\xi\) is empty
the sum of the absolute coefficients is one under every
complex embedding.

Here is the effect of moving that expansion by an
outer cut point. In the original self-dual finite
lattice \(L\) write \(v_L\) for the vacuum and \(\rho\)
for the signed Heisenberg translations. If \(g\) is
the action on the oscillator and \(H=L\cap g^{-1}L\),
then
\[
 \omega(g)v_L=\frac{\Delta(g)}{[L:H]}
       \sum_{z\in L/H}\rho(gz)v_L ,
 \qquad |\Delta(g)|^2=[L:H],\quad
 \Delta(g^{-1})=\overline{\Delta(g)}.
\]
Indeed the fixed line for \(gL\) is the transport of
the vacuum line.

Projection to that line by averaging
the translations over \(gL\) gives the indicated
sum. The translates by distinct cosets are orthogonal,
so unitarity gives the displayed modulus. The
coefficient at the zero coset is the finite Weil
lattice multiplier divided by \([L:H]\). The signed
translations are the vacuum-normalized ones, so
the other phase factors in multiplying translations
are values of the finite additive character, hence
roots of unity.

For a translated vector \(\rho(x)v_L\),
the shifts become \(g(x+z)\). Outside the twist \(g\)
is an integral isometry, so the denominator of the shift
is unchanged there and the local vacuum multiplier
is one. For the move supported at the twist places, the
lattice multiplier identity of Section~\ref{sec:upper-cm-jet-transport} gives
\[
                  \Delta(g)^2=[L:H]\,\upsilon_f(\det g)^{-1}.
\]
Every prime divisor of \([L:H]\) divides \(q_\xi\). To check the
character at \(t=\det g\), use its ideal \(I\), which satisfies
\(I\bar I=(1)\) and is supported above \(q_\xi\). If
\(I^h=(z)\), then \(y=z/\bar z\) has norm one,
\((y)=I^{2h}\), and \(t^{2h}/y\) is an integral unit idele.
The character on this idele has finite order. Its finite value
on \(y\) is the inverse of the infinity-weight monomial on the
conjugates of \(y\), which are units at every coefficient place
away from \(q_\xi\). Thus \(\upsilon_f(t)\) and \(\Delta(g)/[L:H]\)
are units at every such place. The same argument applies to the
norm-one determinant arguments of the additional pair characters.

The extra
pair characters on the integral components have
finite order on the compact level. At empty
twist all actions are integral and their coefficients
have modulus one in every complex embedding.
The Hilbert norm normalization of the multiplier
above contributes the factor
\[
 \frac{\Delta(g)}{[L:H]}
 =\frac{\Delta(g)}{|\Delta(g)|^2}=\Delta(g^{-1})^{-1}.
\]

The arrow product joins the first shift, by its
adjoint into the negative block, to the second shift
in its original order. In the finite vacuum
coordinates its two blocks are
\[
                  z_-=-\eta_-\,{}^t\bar x_1,\qquad z_+=x_2 .
\]
Since \(\eta_-z_-^*=x_1\), the joined product of the two
normalized translations equals \(e_0(z_-)^{-2}T_\sharp(z)
=e_0(x_1)^2T_\sharp(z)\) on the joined datum, by
Section~\ref{sec:double-string-unfolding}. This phase is a
root of unity. The joined translation \(z\)
therefore has denominator \(\beta\) on the first block
and \(\bar\beta\) on the second. As both are units
above \(N\), the \(N\)-component can be dropped on the
ambient vacuum. The retained and removed parts have
zero Heisenberg pairing since they have disjoint
conjugation-stable supports.

Let \(Z_I\) denote the finite zero slice for the
two isotropic bases, and let \(L_{\rm big}\) be the
ambient vacuum lattice, and let \(S_f\) be the string shift. We explain the finite support test and
its coefficient. If \(Z_I\rho(S_f+z)v_{L_{\rm big}}\ne0\),
its commutator against translations by \(L_{\rm big}\cap I\)
must be trivial. For \(\lambda\in L_{\rm big}\cap I\), the signed
translation fixes the vacuum. Its phase on the zero slice is
\[
 e_0(\lambda)=\psi_0\bigl(c_0h(\lambda,\lambda)/2\bigr)=1
\]
because \(I\) is totally isotropic, including at the places above \(2\).

Put \(A=L_{\rm big}\cap I\) and \(x_0=S_f+z\). The character
condition also holds on all \(R_f\)-multiples of each \(u\in A\).
The annihilator of \(R_f\) under the normalized trace character is
\(R_f\), so \(u^\dagger x_0\in R_f\). At each conjugation orbit
\(A\) is saturated in the free ambient lattice. Extend a basis in
each field factor and extend the integral pairing functional by zero
on the added basis vectors. Hermitian self-duality represents it by
\(\lambda\) in the ambient lattice, with the two field factors paired
at a split orbit. Almost everywhere one may take \(\lambda=x_0\),
which is integral there. The local choices yield
\(\lambda\in L_{\rm big}\) with
\(x_0-\lambda\in I^\perp\), so
\(S_f+z-\lambda=i_f w+k_f v\).
The translation by \(\lambda\) fixes the vacuum,
and the translation by \(k_f v\) acts trivially on
the zero slice.

The result is therefore a phase
times \(\rho(w)Z_I v_{L_{\rm big}}\).
The same argument applies to any translate used
to test whether the zero slice itself vanishes.
The zero slice is nonzero in an adapted split model,
and translates of the finite vacuum span its
Schwartz space. Thus \(Z_Iv_{L_{\rm big}}\ne0\).

For \(a=i_f w+k_f v\in I^\perp\), orthogonality and
isotropy give \(\langle a,a\rangle_{\rm big}=
\langle w,w\rangle_{\rm def}\) in the Hermitian trace pairings
on displacements. Hence the quadratic factors in the
vacuum-normalized translations of \(a\) and \(w\) agree,
including at places above two. The two summands of \(a\)
also have trivial Heisenberg commutator and half-pairing.
Restriction to zero in the coordinates dual to \(I\) thus gives
\[
        Z_I\rho_{\rm big}(i_f w+k_f v)=\rho_{\rm def}(w)Z_I .
\]
For \(w\in g_{\rm nb}L_{\rm def}\),
Section~\ref{sec:marked-neighbor} gives an integral lift \(a\)
fixing the ambient vacuum. Thus the signed translation by \(w\)
fixes its zero slice.
In the Schr\"odinger model the original lattice
invariants are one-dimensional: translation
invariance makes a function constant on additive
lattice cosets and invariance under the dual
characters restricts its support to the integral
coset. Weil transport gives the invariant line
for \(g_{\rm nb}L_{\rm def}\), so the zero slice of
the vacuum is a fixed nonzero multiple of
\(\omega(g_{\rm nb})v_{L_{\rm def}}\).
The residual shift \(w=g_{\rm nb}B^\dagger(S_f+z-\lambda)\)
can be replaced by \(g_{\rm nb}B^\dagger z\) up to
a phase: \(B^\dagger S_f=0\), and
\(g_{\rm nb}B^\dagger\lambda\in g_{\rm nb}L_{\rm def}\)
fixes the slice. Section~\ref{sec:marked-neighbor}
gives the same two denominators for the new shift
and shows that it vanishes above \(N\).

\subsubsection{Translation at the additional level}
\label{sec:seesaw-tensor-translation-N}

\sectionopening{In this subsection, we express the neighbor sum of Section~\ref{sec:marked-neighbor} in translated vacua at \(N\).}

For each \(w\mid\pi\), choose a covector
\(y_w\in\mathbb F_N^{1\times n_{2i+1}}\).
Put \(r_{y,w}=\widetilde y_w/N\) for any integral lift of this row
and \(r_{y,\bar w}=0\). The linked
representative \(n_y\) has last coordinates
\((\eta_-\bar r_y,-\bar r_y)\) in the two blocks.
At \(w\mid\pi\) the matrix \(\Gamma\) has integral
entries, so that factor of the additive quotient is zero.
At \(\bar w\), let \(a,b\) be the two rows in the
distinguished coordinates. Their image modulo the integral lattice is
\[
             (a/\bar\pi,\bar\eta_-a/\bar\pi).
\]
It equals \(n_y\) modulo the lattice precisely
when \(\overline{y_w}\equiv-\pi\bar\eta_-a\pmod{\bar\pi}\)
in every middle coordinate.

The chosen split places have residue
field \(\mathbb F_N\), and \(\pi\bar\eta_-\) is a unit
in the conjugate half. Thus every class occurs once.
In fact \(a=\bar\pi\eta_-\bar r_y\) and
\(b=\bar\eta_-a\) give \(n_y\) exactly. Applying the
inverse arrow product to \(n_y\) gives the shifts
\((-{}^tr_y,-\bar r_y)\). As \(y\) ranges
over all residues this is exactly the original
linked gluing, denoted by \(v_x^{\rm gl}\) when applied
to an upper pair with additional translate \(x\).

The two factors are both supported in the prescribed
last coordinates. The signs may simultaneously be
reversed by \(y\mapsto -y\), which permutes the sum.
All pairings between \(n_y\) and the shift away
from \(N\) vanish. In signed vacuum coordinates the adjoint transport
of a first block \(u\) has phase \(e_0(u)^2\) by
Section~\ref{sec:double-string-unfolding}. For \(n_y\) the first block is supported
on only one split half, so its Hermitian square
is zero and the phase is one.
The entire sum at \(N\) is retained in \(v_x^{\rm gl}\).

Away from \(N\), let \(b_i\) be the number of
classes of \(g_{\rm nb}L_{\rm def}\) modulo its
intersection with \(L_{\rm def}\). On an interior pair
\(\kappa=1\) in every coordinate, so \(b_i=1\). On the terminal
pair only the coordinate \(e''\) is scaled, by \(\delta\) in each
of the \(n_{2i+1}\) middle coordinates. The ideals of \(c_M\)
and \(\beta\) agree, and \((\beta),(\bar\beta)\) are coprime.
Consequently \(\delta R_f/(\delta R_f\cap R_f)\) has order
\(|\mathcal O_E/(\bar\beta)|=M^f\), so \(b_i=M^{f n_{2i+1}}\)
on this pair. Above \(N\) the quotient has the independently
computed \(N^{f n_{2i+1}}\) classes in every pair. These are
already summed in the glued vector.
Put the factor \(b_i\) into the common pair scalar
and assign the \(b_i\) terms coefficients of size
\(b_i^{-1}\), including their phases. To check divisibility away from the twist, set
\(D_M=|\mathcal O_E/(\beta)|=M^f\), using
\(\beta\bar\beta=M\) and \([E:\mathbf Q]=2f\).
The two residue tests have denominators dividing
\(D_M^{d^{\rm H}_{2i}}\) and \(D_M^{d^{\rm H}_{2i+1}}\) there.
The quotient just computed has order \(b_i\mid D_M^{r_i n_{2i+1}}\),
with \(r_i=n_{2i}+n_{2i+2}\). For \(j\le m\) the ranks are
\(\dim X_j=n_j=k+j(t+1)\le n_m\) and
\(\dim Y_j=r+jt=k+m-1+jt\le n_m\). Hence
\[
 d^{\rm H}_{2i}+d^{\rm H}_{2i+1}+r_i n_{2i+1}
 \le n_m^2+n_m^2+2n_m^2.
\]
Consequently the product of the three denominators divides
\(D_M^{4n_m^2}=B_*^4\) away from the twist.

At the empty twist, after factoring out the common pair scalar,
the sum of the absolute coefficients of the glued vectors is at most
one at every complex embedding. Each residue test has absolute
coefficient sum one, and the \(b_i\) weights have absolute sum one.
Opening the linked sum inside a glued vector contributes the separate
factor \(N^{f n_{2i+1}}\) used in Section~\ref{sec:upper-terminal-norm}. The vacuum line scalar, the normalization
of its neighbor sum, and the nonzero normalizing
factors of the two upper vacua depend only on the
original pair. If \(s_{\rm vac}\) is the one-coordinate vacuum normalization, the respective
factors are \(\kappa_j=s_{\rm vac}^{n_j n_{j+1}}\).
More precisely, write \(\mathcal C_i\) for the chosen representatives of
\(g_{\rm nb}L_{\rm def}/(g_{\rm nb}L_{\rm def}\cap L_{\rm def})\).
The unshifted slice coefficient \(c_i\ne0\) is determined by
\[
 Z_I v_{L_{\rm big}}=c_i
       \sum_{\nu\in\mathcal C_i}\rho(\nu)v_{L_{\rm def}} .
\]
Write \(Z_Iv_{L_{\rm big}}=a_i\omega(g_{\rm nb})v_{L_{\rm def}}\).
Multiplication by \(g_{\rm nb}\) identifies
\(L_{\rm def}/(L_{\rm def}\cap g_{\rm nb}^{-1}L_{\rm def})\)
with the quotient represented by \(\mathcal C_i\); hence its order is
\(N^{f n_{2i+1}}b_i\). The finite vacuum expansion gives
\[
 c_i=\frac{a_i\Delta(g_{\rm nb})}{N^{f n_{2i+1}}b_i}.
\]
The linked sum at \(N\) has coefficient one on every residue,
and the \(b_i\) classes away from \(N\) carry weights \(b_i^{-1}\).
The common pair scalar for the two integral vacua is therefore
\[
 \frac{c_i b_i}{\kappa_{2i}\kappa_{2i+1}}
 =\frac{a_i\Delta(g_{\rm nb})}
        {N^{f n_{2i+1}}\kappa_{2i}\kappa_{2i+1}} .
\]

Applying this argument to the tensor in
\eqref{eq:seesaw-evaluated-tensor} gives a fixed nonzero
\(C_0\) and an expansion
\[
       T_\xi=C_0\sum_x b_{\xi,x}v_x^{\rm gl}.
\]
Here \(B_*^{2m}b_{\xi,x}\) is integral at every
coefficient embedding off the twist. Indeed the cut
weights have denominators supported on \(q_\xi\)
and each of the \(m/2\) pairs contributes at most \(B_*^4\).
At the trivial twist the cut sum gives
\[
 \sum_x|\sigma(b_{\varnothing,x})|\le B_*^m
 \qquad(\sigma\in\operatorname{Aut}(\mathbb C)).
\] Terms with coincident translates can be
combined by the triangle inequality.

The expanded tensor is fixed by the original
product of upper levels with the terminal action
one. For \(j<m\) let \(X=X_j\) be the definite block of the
Bessel Levi and let \(\Lambda\) be the original lower lattice.
Write \(L_X\) for the standard definite lattice and \(e\) for its
last basis vector. The Levi element attached to \(g\) in the upper
level fixes the string, so its difference from the identity acts
on \(\Lambda\) by \(A\operatorname{pr}_X\) with \(A=g-I\).
For an even string,
\(\operatorname{pr}_X\Lambda=\Lambda\cap X=L_X\).
For an odd string the two columns \(e/\bar\pi_N+uz/N,uz\)
of Section~\ref{sec:marked-neighbor} give
\[
 \operatorname{pr}_X\Lambda=D_{\bar\pi_N}^{-1}L_X,\qquad
 \Lambda\cap X=D_{\pi_N}L_X,\qquad D_v=\operatorname{diag}(1,\ldots,1,v).
\]
Indeed a linear combination of the two columns with zero string
coordinate has first coefficient divisible by \(N\).
The upper last-row condition gives \(e^t A\equiv0\bmod\pi_N\).
The same condition on \(g^{-1}=g^*\) gives \(Ae\equiv0\bmod\bar\pi_N\).
The last diagonal entry is divisible by \(N\), since
\((\pi_N),(\bar\pi_N)\) are coprime. Thus
\(D_{\pi_N}^{-1}A D_{\bar\pi_N}^{-1}\) is integral in the odd case,
and \(A\) is integral in the even case. In both cases
\(A\operatorname{pr}_X\Lambda\subseteq\Lambda\cap X\), also
for \(g^{-1}\), so the Levi action preserves \(\Lambda\).
At the chosen places above \(M\) both diagonal scalings are units
and \(g\equiv I\bmod\pi_M\), giving identity reduction on \(\Lambda\).

At each intermediate stage the embedded element fixes the string,
preserves \(\Lambda\), and belongs to \(B_j^0=R_j(A_f)\cap K_j^\circ\).
The set \(\Omega_j\) is a union of left \(B_j^0\)-orbits and
\(d_j\) is constant on each orbit by Section~\ref{sec:measured-finite-sum}.
Thus the element permutes the cut classes with their weights.
At the first stage it permutes the chosen level classes. Equivariance
of the two-slice evaluation gives this upper
action on each summand, on both occurrences of
every shared index. Thus it preserves the measured
sum. Before applying the linked sum at \(N\), write each
additional shift as \(x=x^{\rm lev}+x^{\rm away}\), supported above
\(Mq_\xi\) and away from those places, respectively.
Both supports are stable under conjugation and are disjoint.
The signed Heisenberg multiplication law therefore gives
\[
 \rho(x)=\rho(x^{\rm lev})\rho(x^{\rm away}),\qquad
 \rho(x^{\rm away})v_L=v_L ,
\]
where the second identity uses integrality of \(x^{\rm away}\)
and the vacuum normalization of the signed translations.
Thus \(\rho(x)v_L=\rho(x^{\rm lev})v_L\) on each arrow.
Applying the linked sum at \(N\) to these equal tensors
preserves the evaluated tensor.

Apply Section~\ref{sec:seesaw-normalization} to
\eqref{eq:seesaw-measured-chain} with this tensor.
Put \(a_{\rm ss}=c_{\rm meas}C_0/B_*^{2m}\) and define
\[
 (c_{\rm ss},\Psi_\xi(x))=
 \begin{cases}
  (a_{\rm ss},B_*^{2m}b_{\xi,x}) & a_{\rm ss}\ne0,\\
  (1,0) & a_{\rm ss}=0 .
 \end{cases}
\]
In the second case \(c_{\rm meas}=0\), since the two-arrow vacuum
calculation gives \(C_0\ne0\). In both cases \(c_{\rm ss}\ne0\).
Linearity and \eqref{eq:seesaw-measured-chain} give
\begin{equation}
 L_\xi(c_{\rm ss}^{-1}\Pi\phi)
    =\mathcal U_\xi\!\left(\sum_x\Psi_\xi(x)v_x^{\rm gl}\right)\phi
       \qquad(\phi\in S).
 \label{eq:measured-seesaw-normalized-tensor}
\end{equation}
The scalar \(c_{\rm ss}\ne0\) is fixed before the twist.
For each \(\xi\), one tensor gives the identity for every \(\phi\in S\).
For the empty twist
\[
            \sum_x|\sigma\Psi_\varnothing(x)|
                   \le B_*^{3m}=M^{3fm n_m^2}
\]
for every complex embedding \(\sigma\).
The resulting upper lift evaluates the same
source and uses the same terminal branching as
\(L_\xi\Pi\). The equality on the original
ramified source is now in integral tensor form.

\subsection{The comparison at the additional ramified place}
\label{sec:marked-N-lift}

\sectionopening{In this subsection, we compare the two integral realizations at \(N\)
for the tensor of Section~\ref{sec:seesaw-integral-tensor}. The
comparison gives an integral submodule of the ramified source on which
the local theta estimate holds uniformly.}

Use the spaces \(X_j\), their original finite levels, and the
norm equation \(\pi\bar\pi=N\).
These are the upper definite groups and original glued
vacua of Section~\ref{sec:seesaw-integral-tensor}. Fix an \(N\)-adic embedding
\(\jmath\). Put \(P=\Sigma_N\), the set of conjugates of the
places induced by \(\jmath\) and the infinite
embeddings. Put \(S=\{w:w\mid\pi\}\). Each of \(P,S\) contains
one place from every conjugate pair over \(N\), and at each
such place the completion and the residue field are
\(\mathbf Q_N,\mathbf F_N\).

For a tuple \(a\) of matrices at \(P\), let \(\widehat a\) be
its unitary embedding by inverse conjugate transpose at \(\bar P\)
and identity away from \(N\). Put \(\delta_j(P)=\widehat A_j\),
where \((A_j)_w=\operatorname{diag}(I_{n_j-1},N^{\mathbf1_{w\in S}})\)
at \(w\in P\). Write \(\delta_j=\delta_j(P)\) at this embedding. Thus at the chosen half \(P\)
the point takes the standard lattice to \(L_j\), the lattice of
integral vectors whose last coordinate is divisible by \(\pi\).
For right translation \(R_j\) on forms, its expression in these bases is
\(\rho_j(a)=R_j(\delta_j\widehat a\delta_j^{-1})\).
The original local level becomes the stabilizer of
\(e=e_{\rm last}\) modulo \(N\). The conjugate coefficient embedding
uses \(\delta_j(\bar P)\) with the same \(S\) and the same original level.

Let
\[
 B=\prod_{w\in P}\mathbf F_N,\quad
 D_j=\operatorname{diag}(1,\ldots,1,N^{-1}),\quad
 \psi_N(y)=\psi_0(-y/N).
\]
Here \(D_j\) acts at all of \(P\) in the appropriate rank. More
explicitly, let \(e_P\) be the finite idempotent of the places in \(P\)
and let \(\widetilde y\) lift the residue coordinates there. Then
\[
 \psi_{N,P}(y)=
 \psi_f\!\left(-j_{\rm W}\alpha_{\rm W}\bar\alpha_{\rm W}
                       e_P\widetilde y/N\right)=\psi_N(y).
\]
The inverse different identity of Section~\ref{sec:period-finite-coordinates}
makes this independent of the lift and nontrivial on each
\(\mathcal O_{E_w}/N=\mathbb F_N\). Since
\(\bar j_{\rm W}=-j_{\rm W}\) and the trace is invariant under
conjugation, \(\psi_{N,\bar P}(\bar y)=\psi_{N,P}(y)^{-1}\).
All functions of a residue matrix below are extended by
zero off the integral matrices.

\smallskip
\subsubsection{Translation laws and the ramified sum}
\label{sec:ramified-translation-sum}

\sectionopening{In this subsection, we apply the local oscillator action to the summands of Proposition~\ref{prop:measured-seesaw}.}

Let \(v_j\) be a finite Weil vector. In the split
Schr\"odinger model at \(P\), let \(\mathcal B_j\) be the integral
matrix box, and normalize its additive volume to one. The projection
onto its vacuum line at \(N\) is
\[
 v_{P}\otimes v^{(N)}\longmapsto
 \mathbf1_{\mathcal B_j}\otimes
       \left(\int_{\mathcal B_j}v_{P}(T)\,dT\right)v^{(N)}
\]
on pure tensors. Write the projection of \(v_j\) as
\(\mathbf1_{\mathcal B_j}\otimes u_j^{(N)}\), extending the formula
linearly. Inserting \(\Phi\) gives \(\Phi\otimes u_j^{(N)}\).
Write \(L_j(\Phi)\) for the lift of this datum in the
\(L_j,L_{j+1}\) coordinates. A group element supported at \(N\)
therefore acts on the inserted function through its split Weil action.
For \(a,b\) invertible over the local integers, the laws are
\begin{align}
 \rho_{j+1}(b)L_j(\Phi(bT a))\rho_j(a)&=L_j(\Phi(T)),\label{eq:marked-N-law-one}\\
 \rho_{j+1}(D_{j+1})L_j(\mathrm{cov}_e)
 &=\gamma_j L_j(\mathrm{col}_e)\rho_j(D_j).\label{eq:marked-N-law-two}
\end{align}
Here \(\mathrm{cov}_e(T)=\mathbf1_{\operatorname{row}_{\rm last}(\bar T)=e^t}\)
and \(\mathrm{col}_e(T)=\mathbf1_{\bar T e=e}\) on integral
matrices. Put \(n=n_j,\ n'=n_{j+1}\), and write \(\chi_j\)
for the base change of the \((j+1)\)-st lift character.
The substitution \(T\mapsto b^{-1}Ta\) has additive Jacobian
\(J=|\det a|_P^{n'}|\det b|_P^{-n}\). The unitary oscillator acts
on functions by this substitution times
\(\chi_j(\det b/\det a)J^{1/2}\). Changing the source variable in
the lift by right translation gives
\[
 L_j\bigl(\Phi(b^{-1}Ta)\bigr)=
 \chi_j(\det a/\det b)J^{-1/2}\,
 \rho_{j+1}(b)L_j(\Phi)\rho_j(a^{-1}) .
\]

It gives \eqref{eq:marked-N-law-one} since the characters
are unramified there. The equality of functions
\(\mathrm{cov}_e(T)=\mathrm{col}_e(D_{j+1}TD_j^{-1})\)
includes their supports. In blocks the change sends
\(\left(\begin{smallmatrix}A&b\\c&d\end{smallmatrix}\right)\) to
\(\left(\begin{smallmatrix}A&Nb\\c/N&d\end{smallmatrix}\right)\).
Testing the last column and integrality therefore tests
\(c\equiv0,\ d\equiv1\bmod N\) and integrality of \(T\).
The determinant character factor in \eqref{eq:marked-N-law-two}
uses \(a=D_j^{-1}\) and \(b=D_{j+1}^{-1}\) in
\(T\mapsto b^{-1}Ta\). Both determinants have component \(N\) at each
place of \(P\), so their character values cancel and
\[
 \gamma_j=(N^{-|P|})^{-n_{j+1}/2+n_j/2}
          =N^{|P|(n_{j+1}-n_j)/2},\qquad |\gamma_j|_N\le1 .
\] Conjugation by
\(\delta_j\) in these calculations uses the same change on
the datum and the source action. The source and target lattices in the two laws are transported
by \(\delta_j,\delta_{j+1}\), respectively.

A ramified source vector \(u\) on the original level has
zero mean at each scalar group over \(N\). For each
\(w\in P\) write \(u=(1-\rho_j(c_w))u_w\), using a generator
of that local cyclic residue unit group in the scalar action,
supported at \(w\) and its conjugate. Write \(C=\rho_j(c_w)\)
and \(L=N-1\). Since \(C^L=1\) on the level invariants,
\[
 (1-C)\frac1L\sum_{r=0}^{L-1}(L-1-r)C^r
     =1-\frac1L\sum_{r=0}^{L-1}C^r.
\]
The right side fixes \(u\). The scalar translations preserve the
transported integral lattices and \(N\nmid L\), so this formula also
preserves integrality at \(N\). At \(j=0\) the source condition uses \(S\), and at \(w\in P\) it reads through the member of
\(S\) in the same conjugate pair.

Conjugation by \(\delta_j\)
does not change a scalar. Every \(L_j(\Phi)\) for a residue
function carries the difference form of ramification to the
target by \eqref{eq:marked-N-law-one}, taking opposite
scalar changes of variables. Moreover, if the datum is
invariant under scaling \(T\) at one place \(w\), then it
annihilates such a difference at \(w\). This applies to the
indicator that a given matrix column or row reduces to a
vector that vanishes at \(w\). It also applies to the datum
that a specified matrix entry vanishes there.

For a residue column \(y\) in the target let
\(C_y(T)=\mathbf1_{\bar T e=y}\), and for a column in the
source let \(H_y(T)=\psi_N(e^t\bar T y)\), with the integral
support in both cases. Write \(L_j^{\rm d}=L_j(C_e)\).
For the double lift the middle vector \(y\) must be
primitive, meaning nonzero in each residue field. The other
terms in \(\sum_y L_{j+1}(H_y)L_j(C_y)u\) vanish by the
ramified difference property. For primitive \(y\) choose
\(k\) integral invertible with \(\bar k e=y\) by extending
the column to a basis in every residue field and lifting.
The two laws for its insertion in the middle are
\[
 L_j(C_y)u=\rho_{j+1}(k)L_j^{\rm d}u,\qquad
 L_{j+1}(H_y)\rho_{j+1}(k)=L_{j+1}(H_e).
\]
These follow from \eqref{eq:marked-N-law-one}, since
\(\overline{kT}e=\bar k(\bar T e)\) and
\(e^t\overline{T k} e=e^t\bar T y\).
Summing gives
\[
 \sum_y L_{j+1}(H_y)L_j(C_y)u
       = c_{\rm prim}\,L_{j+1}(H_e)L_j^{\rm d}u,
 \qquad c_{\rm prim}=(N^{n_{j+1}}-1)^{|P|}.
\]
Indeed each residue column has \(N^{n_{j+1}}-1\) nonzero
choices at each \(w\in P\), independently.
In particular \(|c_{\rm prim}|_N=1\).

Define the same-level operator
\[
 U_j=\sum_{y:\,y_{\rm last}\in B^\times}
        \psi_N(y_{\rm last})\rho_j(D_j k_y),
\]
where the integral invertible \(k_y\) has both last column
\(y\) and last row \(y^t\) modulo \(N\). For \(y=(v,t)\)
with \(t\) a unit the residue matrix
\(\left(\begin{smallmatrix}1+vv^t/t&v\\v^t&t\end{smallmatrix}\right)\)
has determinant \(t\) by the Schur complement and is such
a lift modulo \(N\). The actual \(U_j\) uses the choice of the lifts. For the actions on mirabolic fixed
vectors only their last columns matter. Conjugating a
mirabolic \(h\) by \(D_j^{-1}\) makes an integral matrix
\(D_j^{-1}hD_j\) whose last row reduces to \(e^t\). Left
multiplication by that matrix permutes the residue \(y\)
with unit last entry, preserving their last entries.
The error after replacing \(k_y\) by the chosen \(k_{\bar
{D_j^{-1}hD_j}y}\) fixes the last column. This proves that
\(U_j\) preserves the original mirabolic fixed vectors.
It commutes with the original level away from \(N\) by
disjoint support. It preserves ramification because it
commutes with the unit scalar actions.

At a ramified source vector one may omit matrices with
nonunit corner from the function \(H_e\), by decomposing
that function according to the corner and using the scalar
difference property at a place where the corner is zero.
Sort the remaining function first by its last column \(y\)
and then, in a separate expression, by its last row \(y^t\).
Using \eqref{eq:marked-N-law-one} in both sorts and
\eqref{eq:marked-N-law-two} for the last-row indicator gives
\[
 \rho_{j+1}(D_{j+1})L_j(H_e)u
     =U_{j+1}L_j^{\rm d}u
     =\gamma_j L_j^{\rm d}U_ju.
\]
These are identities of the same theta lift, so no change of
integral structures occurs when they are used in succession.

\smallskip
\subsubsection{The direct lift at the ramified place}
\label{sec:ramified-direct-lift}

\sectionopening{In this subsection, we identify the gluing of the summed tensors in Section~\ref{sec:ramified-translation-sum}.}

At \(N\) the normalized vacuum projects to the integral
box. In the model at \(P\), the original gluing from \(S\)
is the sum of the tensors
\[
 C_y(T A_j)\otimes H_y(B_{j+2}R),\qquad y\in B^{n_{j+1}},
\]
where \(A_j=\operatorname{diag}(1,\ldots,1,N)\) on \(P\cap S\)
and \(B_{j+2}=\operatorname{diag}(1,\ldots,1,N)\) on
\(P\setminus S\), both with the identity elsewhere.
Write \(e_i\) for the last basis vector at rank \(n_i\), and set
\(t_w=T_w e_j\), \(r_w=e_{j+2}^tR_w\).
Let \(b=(b_v)_{v\in S}\) be the residue column indexing the original
pair of translates in Section~\ref{sec:seesaw-tensor-translation-N},
and let \(\operatorname{red}_w\) denote reduction modulo \(N\) at \(w\).
Write \(\widetilde b_v\) for any integral lift. At \(w\in P\), all
entries outside the indicated last column and last row must be integral.
For the remaining entries the two translated box indicators impose
\[
\begin{array}{c|c|c}
 w & \text{integral entries} & \text{unique original class}\\ \hline
 P\cap S &
 t_w-\widetilde b_w/N,\ r_w &
 b_w=\operatorname{red}_w(Nt_w)\\
 P\setminus S &
 t_w,\ r_w-\overline{\widetilde b_{\bar w}}^{\,t}/N &
 b_{\bar w}=\overline{\operatorname{red}_w(Nr_w)}^{\,t}.
\end{array}
\]
Here conjugation identifies the residue fields at \(w\) and \(\bar w\).
Write \(\psi_{0,w}\) for the restriction of \(\psi_0\) to the \(w\)-factor.
In split coordinates \((a,b^\vee)\), with \(a\) in the position space
at \(w\) and \(b^\vee\) the covector from the conjugate component,
the displacement acts by
\[
 (\mathsf T(a,b^\vee)F)(z)
   =\psi_{0,w}\!\left(b^\vee(z+a/2)\right)F(z+a).
\]
Indeed its multiplier on a second displacement \((a',b'^\vee)\)
is \(\psi_{0,w}((b'^\vee a-b^\vee a')/2)\), the split expression for
\(\beta(x,y)=\psi_0(y^*x/2)\). Each original displacement has only
one nonzero component, so its quadratic phase is one. In the first
row of the table the conjugate translate on \(R_w\) contributes
\(\psi_{0,w}(-r_w\widetilde b_w/N)\). In the second row the first
translate on \(T_w\) contributes
\(\psi_{0,w}(-\overline{\widetilde b_{\bar w}}^{\,t}t_w/N)\).
Both arguments are congruent to \(-r_wt_w\) modulo the local integers.
Since \(\psi_{N,w}(\operatorname{red}_w c)=\psi_{0,w}(-c/N)\) for integral \(c\),
the factor in either row is
\(\psi_{N,w}(\operatorname{red}_w(Nr_wt_w))\).
With \(T'=TA_j\) and \(R'=B_{j+2}R\), the sum is zero unless
\(T',R'\) are integral.
On that support its value is
\[
             \psi_N(e_{j+2}^t\overline{R'}\,\overline{T'}e_j).
\]
Sorting by \(y=\overline{T'}e_j\) recovers the displayed tensor sum
with coefficient one. The semicharacter of each translation is one
by its support on a single isotropic half, and Fourier interchange
preserves the integral box with its self-dual measure of volume one.

For this pair \(j=2i\), let \(\widetilde L_l(\Phi)\) denote
the lift with \(\Phi\) in the standard bases, so
\[
 L_l(\Phi)=R_{l+1}(\delta_{l+1})\,\widetilde L_l(\Phi)\,R_l(\delta_l^{-1}).
\]
Write \(\chi_{l,h}(a)=\chi_l(\det a)|\det a|_P^{-h/2}\)
with the character of arrow \(l\) in Section~\ref{sec:ramified-translation-sum}.
Its split change of variables gives
\[
\begin{aligned}
 \widetilde L_j(C_y(TA_j))
   &=\chi_{j,n_{j+1}}(A_j)\,R_{j+1}(\delta_{j+1}^{-1})L_j(C_y),\\
 \widetilde L_{j+1}(H_y(B_{j+2}R))
   &=\chi_{j+1,n_{j+1}}(B_{j+2})\,
     \rho_{j+2}(D_{j+2}) L_{j+1}(H_y)R_{j+1}(\delta_{j+1}).
\end{aligned}
\]
In the first row take \((a,b)=(A_j,1)\) in \(b^{-1}Ta\);
in the second take \((1,B_{j+2}^{-1})\).
Here \(A_{j+2}B_{j+2}=D_{j+2}^{-1}\), whence
\(\rho_{j+2}(D_{j+2})R_{j+2}(\delta_{j+2})
=R_{j+2}(\widehat B_{j+2}^{-1})\).
On composition the two middle translations cancel. For the
current twist and coefficient embedding the scalar is
\[
 g_{j,\xi,\jmath}=\chi_j(\det A_j)\chi_{j+1}(\det B_{j+2})
           |\det A_j\det B_{j+2}|_P^{-n_{j+1}/2}.
\]
Each \(w\in P\) occurs in exactly one determinant. Through the conjugate embedding
a factor \(N\) on \(P\) reads as \(N^{-1}\) in the norm-one
determinant.

The algebraic character formula therefore
gives an \(N\)-adic factor at most \(N^{2b}\) there,
because the lift weight at that embedding is \(\le2b\)
and the modulus power has \(N\)-adic norm \(\le1\).
Thus \(|\jmath(g_{j,\xi,\jmath})|_N\le N^{2fb}\).
The tame twist has weight zero and is unramified at \(N\).
Its uniformizer values are roots of unity and hence \(N\)-adic units.
For the embedding \(\jmath\circ\mathrm{conj}\), use \(\bar P\)
with its partition by the same \(S\). The argument gives the
same bound at that embedding.
Apply the ramified summation and commute \(U\) through the
direct lifts one at a time.

All intermediate vectors remain ramified. For pair \(i\), the
commutations contribute
\[
 \prod_{j=0}^{2i+1}\gamma_j=N^{|P|(n_{2i+2}-n_0)/2}.
\]
The product for the fixed twist and embedding is therefore
\[
 c_{\xi,\jmath}=
 \prod_{i=0}^{m/2-1}
 \left(g_{2i,\xi,\jmath}\,(N^{n_{2i+1}}-1)^{|P|}
                         N^{|P|(n_{2i+2}-n_0)/2}\right).
\]
The ranks increase, so
\(|\jmath(c_{\xi,\jmath})|_N\le(N^{2fb})^{m/2}=N^{mfb}\).
The scalar is independent of the input vector and the finite
translates. Thus on the ramified source
\[
 \mathcal U_\xi(\Psi)(u)
 =\mathrm{res}_{Y_m}\left(
 c_{\xi,\jmath}\sum_x\Psi_\xi(x)\,
          L_{m-1,x}^{\rm d}\cdots L_{0,x}^{\rm d}
                      U_0^{m/2}u\right).
\]
Here each \(L_{j,x}^{\rm d}\) uses the same translate of the
normalized vacuum as the \(x\)-term of the original chain.
The permitted translates are vacuum-projected above \(N\)
before gluing, as given by admissibility of the chain.

To bound a direct lift, use the standard lattices at the identity. Write \(e_-\) and \(e_+\) for the last basis vectors
in the source and target, respectively. On integral matrices,
finite Fourier inversion on the last column gives
\[
 C_e(T)=N^{-|P|n_{j+1}}
   \sum_{a\in B^{n_{j+1}}}
     \psi_N\!\left(a^t(\bar T e_- -e_+)\right).
\]
Write \(v=\bar T e_- -e_+\). The sum factors over the
places and the row coordinates. Each one-coordinate sum is
\[
 \sum_{a_w\in\mathbb F_N}\psi_{N,w}(a_w v_{w,i})
       =N\mathbf 1_{v_{w,i}=0}.
\]
Thus any single nonzero coordinate kills the sum, including when
its entry in \(B\) is a zero divisor. The inverse character gives
the same identity in the opposite polarization. The denominator
therefore divides the common clearing power
\(N^{[E:\mathbf Q]n_j n_{j+1}}\).

Multiplication by a character of
these residue variables is a Heisenberg translate by
a vector with at most simple poles at the conjugate
places \(\bar P\), by the self-dual pairing. For the
projected vacuum the original translate is either
integral above \(N\), or the projection is zero.
All induced places \(\bar P\) restrict to the same
quadratic prime of \(K\). Choose the quadratic integer
\(\delta=\omega_K-c_0\), where \(\omega_K=(1+\sqrt{-d})/2\)
and \(c_0\in\mathbb Z\) lifts the residue of \(\omega_K\) at
that quadratic prime. The Fourier translations have at most simple
poles above \(N\), all in \(\bar P\). Thus for every surviving
term the total translation \(z\) satisfies
\(\delta z_w\in\mathcal O_{E_w}^{n_jn_{j+1}}\) for every \(w\mid N\).
Write \(\mathfrak n\) for the chosen prime of \(K\).
Then \(v_{\mathfrak n}(\delta)\ge1\), whereas
\(\bar\delta\equiv-\sqrt{-d}\pmod{\mathfrak n}\)
gives \(v_{\bar{\mathfrak n}}(\delta)=0\).
Thus \(\jmath(\bar\delta d)\) is a unit. For the conjugate
coefficient embedding the pole set is \(P\) above
\(\bar{\mathfrak n}\); using \(\bar\delta\) there clears the
poles and gives the same unit bound.
Apply Section~\ref{sec:marked-integral-sections-ordinary-jets}
to these translates of the normalized vacuum.
At its original scale the required factor in the lattice norm is
\((|\jmath(\pi_M)|_N/|\jmath(\bar\delta d)|_N)^{\#\text{boxes}}
\le1\).
After dividing by the Fourier clearing factor, the
direct lattice norm contribution is at most \(N^{2fn_jn_{j+1}}\le N^{fb}\).

The estimate in the standard lattice bounds the kernel at every
pair of finite points. Indeed choose rational translates
of the two points independently so that their finite
coordinates above \(N\) are integral. This is possible
by weak approximation on the two unitary groups,
using their rational Cayley charts. At the integral representatives
the action at \(w\in P\) has the form \(T\mapsto b_w^{-1}Ta_w\)
with \(a_w,b_w\) integral and invertible. It preserves the integral
matrix lattice. Writing \(a_{w,0},b_{w,0}\) for the two reductions,
the transported column condition is
\[
        \bar T a_{w,0}e_-=b_{w,0}e_+ .
\]
The same character expansion applies to this equation. Its linear
forms are integral in the residue coordinates and its constant
phases are roots of unity. Thus the Fourier denominator is still
\(N^{|P|n_{j+1}}\), and the translating vectors have at most simple
poles at \(\bar P\). The unramified character factors above \(N\)
are units. Let \(A\) be the forward action on the displacement space
away from \(N\), extended by the identity above \(N\).
In the finite Weil expansion of Section~\ref{sec:seesaw-integral-tensor},
write \(H=L\cap A^{-1}L\) for the intersection with the vacuum
lattice \(L\). The index \([L:H]\) is prime to \(N\), and the
representatives \(y\in L/H\) can be chosen zero above \(N\).
The splitting multiplier and the pair characters are \(N\)-adic
units by the class-number power calculation in
Section~\ref{sec:marked-gauge-at-M} with coefficient prime \(N\).
Each shift \(A(z+y)\) retains the component \(z_w\) at \(w\mid N\).
Thus \(N^{|P|n_{j+1}}\) times the transformed datum is in the
coefficient valuation ring span of vacua with
\(\delta (A(z+y))_w\in L_w\) for every \(w\mid N\), with the same
\(\delta\) as above. The ordinary jet bound of
Section~\ref{sec:marked-integral-sections-ordinary-jets} applies at
the integral representatives, and rational covariance transports
the bound to the Weyl lattice norms at the original points.
Translating the datum by \((\delta_j,\delta_{j+1})\)
changes the two transported lattices to those at \(g\delta_j\)
and \(h\delta_{j+1}\). In particular, writing \(\|v\|_a\)
for the lattice norm at \(a\), the direct kernel satisfies
\[
 \|K_{j,x}^{\rm d}(g,h)v\|_{h\delta_{j+1}}
          \le N^{fb}\|v\|_{g\delta_j}
                       \qquad\text{for every }g,h .
\]

To sum the products of these kernels, project the
glued tensor onto the original integral vacuum above
\(N\) in each factor. This is the product of the
mass-one additive averages defining the original
vacuum projections. For a linked nonzero residue choose a
place \(w\in P\) where \(y_w\ne0\), identifying residues on conjugate
halves as in the gluing formula. If \(w\in P\cap S\), the first
factor is a translate of the integral box by a last column
\(-\widetilde y_w/N\). Its support is disjoint from the box on
which the vacuum projection integrates. If \(w\in P\setminus S\),
the first factor on that box is the additive character paired with
the nonzero residue. Write \(O_w\) for the integer ring, \(Q_w\) for
its residue cardinality and \(n\) for the number of column entries.
The residue pairing is perfect and the box has additive mass one. Thus
\[
 \int_{O_w^n}\psi_{N,w}(\langle z\bmod N,y_w\rangle)\,dz
 =Q_w^{-n}\sum_{z\in(O_w/NO_w)^n}
                 \psi_{N,w}(\langle z,y_w\rangle)
 =\mathbf1_{y_w=0}.
\]
For \(y_w\ne0\), translation by a residue vector with nontrivial
character value multiplies the sum by that value. This equality
in characteristic zero precedes the coefficient valuation at \(N\).
Only \(y=0\) survives. Both factors have value one on their
original boxes of additive volume one, so the projected tensor is
\[
 \sum_x \Psi_\xi(x)\bigotimes_{j=0}^{m-1} V_j(x),
\]
where \(V_j(x)\) denotes the original translated
vacuum. The original levels preserve the projection
lattices and their probability Haar, so covariance
makes the projection commute with their simultaneous
action. In particular the displayed tensor is
jointly invariant. The direct insertions also commute
with this action. In their conjugated bases the
original levels preserve the integral box, and the
test \(\bar T e=e\) is equivalent to
\(b^{-1}\bar T a e=e\) for the two residue actions.
Hence the tensor with the direct insertions is
jointly invariant under the same levels.

Let \(u=U_0^{m/2}\varphi\) have lattice norm at most one at
\(g\delta_0\) for every \(g\), as required by the
definition of the preimage lattice. For a fixed
terminal point \(g_m\) the product integral of the
\(m\) kernels, summed with \(\Psi_\xi(x)\), is
\[
 \sum_{(g_j)_{j<m}}\sum_x \Psi_\xi(x)
      K_{m-1,x}^{\rm d}(g_{m-1},g_m)\cdots
      K_{0,x}^{\rm d}(g_0,g_1)u(g_0).
\]
Indeed simultaneous level invariance makes the
contracted integrand constant on each product of the
original level classes. Their stabilizers are trivial
at the old odd principal level and the class masses
in the mean integrals are one. Since \(N\) is prime
to the twist modulus, \(|\jmath(\Psi_\xi(x))|_N\le1\).
The \(m\) bounds multiply to \(N^{mfb}\), and the
ultrametric inequality gives that same bound for
the entire sum in the lattice norm at \(g_m\delta_m\).
Multiplication by \(c_{\xi,\jmath}\) gives the bound
\(N^{2mfb}\) before restriction.

For \(h\in U(Y_m)(\mathbb A_{F_0,f})\) put
\(\widetilde h=\operatorname{diag}(h,1)\in U(X_m)(\mathbb A_{F_0,f})\).
Let \(w_\tau\) be the place induced by \(\jmath\circ\tau\), and set
\(\chi_\tau=\mathbf1_{\bar w_\tau\mid\pi}\).
The weight component is read at \(\bar w_\tau\in P\), where
\((\widetilde h\delta_m)^{-1}\) has blocks \(h^{-1}\) and
\(s_\tau=N^{-\chi_\tau}\). The last-coordinate exponent of
the branching matrix \(B_Y\) is
\[
 B_\tau=|\mu(X_m)_\tau|-|\mu(Y_m)_\tau|=a_\tau-mh_0
\]
by Lemma~\ref{lem:terminal-branching-norm}. Since
\(\pi\bar\pi=N\), the restriction scalar satisfies
\[
 \left|\jmath\!\left(\prod_\tau\tau(\bar\pi)^{-B_\tau}\right)\right|_N
      =\left|\prod_\tau s_\tau^{B_\tau}\right|_N .
\]
Indeed \(v_N(\jmath(\tau(\bar\pi)))=v_{\bar w_\tau}(\pi)
=\chi_\tau\), since \(N\) splits completely. Both sides are
\(N^{\sum_\tau\chi_\tau B_\tau}\).
For the conjugate embedding \(\jmath\circ\mathrm{conj}\), form its
own half \(P\) and the corresponding lattice coordinate change
\(\delta_m=\widehat A_m\)
as above. Then \(\chi_\tau\) is replaced by \(1-\chi_\tau\) on both
sides.
For each signed exponent \(B_\tau\), the corresponding scalar
quotient is a unit, so
\[
 u_{\jmath}
 =\jmath\!\left(\prod_\tau\tau(\bar\pi)^{-B_\tau}\right)
                    \prod_\tau s_\tau^{-B_\tau},
 \qquad |u_{\jmath}|_N=1 .
\]
Write \(\rho_X,\rho_Y\) for the Weyl actions in the coefficient
coordinates. The integral branching matrix satisfies
\[
 B_Y\rho_X((\widetilde h\delta_m)^{-1})
 =\left(\prod_\tau s_\tau^{B_\tau}\right)\rho_Y(h^{-1})B_Y .
\]
Thus the restriction in transported coordinates is \(u_\jmath B_Y\).
Since \(B_Y\) has integral entries, the upper lift \(\Phi\)
just bounded satisfies
\[
 \|(\operatorname{res}\Phi)(h)\|_h
       \le \|\Phi(\widetilde h)\|_{\widetilde h\delta_m}
       \le N^{2mfb}
       \qquad\bigl(h\in U(Y_m)(\mathbb A_{F_0,f})\bigr).
\]
Its contribution to the
logarithmic denominator bound is
\[
 e_N\log N=2mfb\log N=O(fm^2bB)=O(faB)=o(faB^3)
\]
by \(\log N=O(mB)\), \(m^2\le4x\), and \(a=xb\), as in
Lemma~\ref{lem:total-local-denominators}.

\subsubsection{Integral clearing on the original source}
\label{sec:ramified-source-clearing}

\sectionopening{In this subsection, we bound the colength required for the direct lift of Section~\ref{sec:ramified-direct-lift}.}

Keep the coefficient embedding \(\jmath\) above \(N\).
Set \(\Delta_N=D_0^{-1}\) on the \(P\) half and use its inverse
conjugate transpose on the other half. Thus the dual point
is \(\delta_0^\vee=\delta_0\Delta_N^{-1}\).
At a place in \(P\cap S\) the original \(\delta_0\)
has last entry \(N\); on \(P\setminus S\) it is the
identity. Every embedding in the coefficient lattice norm
reads the point through \(P\). Therefore
\((\delta_0^\vee)^{-1}\) acts in the Weyl coordinates by integral diagonal matrices. The source representation is polynomial,
with determinant exponent zero, so integrality at \(1\)
implies integrality at \(\delta_0^\vee\).

For a class representative \(c\) of the original source
let \(K_N\) be the product of the integral groups above \(N\).
Choose a rational unitary point \(\gamma_c\) such that
\[
      \gamma_c c\delta_0^\vee\in
                    K_N\cap\Delta_N K_N\Delta_N^{-1}.
\]
Weak approximation here follows from the Cayley transform
\(X\mapsto(1-X)(1+X)^{-1}\) on skew-Hermitian matrices.
The chart is dense in each local unitary group and
the coefficients of \(X\) can be approximated at the finitely
many places.

The required condition is open and nonempty
as a condition on \(\gamma_c\) over the local fields.
Put \(r_c=\gamma_c c\), so \(r_c\delta_0^\vee\in K_N\)
and \(r_c\delta_0\in\Delta_N K_N\).
For \(k\) in the original source level,
\(\delta_0^{-1}k\delta_0\) is integral above \(N\)
and fixes the last column on reduction in \(P\).
At \(P\cap S\) this follows by conjugating the original
condition on the last row, which divides the lower left block
by \(N\) and multiplies the upper right by \(N\).

At \(P\setminus S\) it follows by applying the condition on the last row to \(k^{-1}\) on the other half and taking
the inverse conjugate transpose in the unitary coordinates.
Both coordinates have unit determinant. Thus the unit lattices
at the two points are tested in the same class
bases \(r_c\), with the indicated diagonals. Rational
equivariance transports these tests to arbitrary representatives.

For a Weyl basis index \(I\) let \(w_I\) be the total
degree in the last variable over all embeddings. Integrality
at \(\delta_0^\vee\) gives integral coordinates.
Integrality at \(\delta_0\) requires the \(I\)-th coordinate
to belong to \(N^{w_I}\mathcal O\). Here \(\mathcal O\) is the
valuation ring pulled back along \(\jmath\) to the coefficient field,
with \(v_N(N)=1\). The ratios used in integral column operations
therefore lie in that same field. For a finite tuple of vectors
impose these conditions
one coordinate row at a time. In a nonzero row choose a coefficient
of least valuation, clear the other coefficients by
integral column operations, and multiply the remaining
column by \(N^{w_I}\). A zero row already satisfies
the condition.

The operations preserve the divisibility already imposed on the preceding rows. The determinant contribution is at most
\[
             h_{\rm deg}=\sum_{c,I}w_I\le d_Xfa/k ,
\]
where \(d_X\) is the dimension of the full original source.
Indeed the sum of the degrees in each variable is the
same, by permuting the diagonal in the Weyl representation.
Summing the degrees over all \(k\) variables gives at most
\(fa\) for each \((c,I)\). The source level is free:
a stabilizer at an old odd distinguished prime lies in the
torsion-free principal subgroup proved in
Section~\ref{sec:marked-norm-line}.

The right translation in \(U_0\) at \(y\) is by
\(u_y=\delta_0^\vee k_y\delta_0^{-1}\).
Thus \(u_y\delta_0=\delta_0^\vee k_y\) differs from
\(\delta_0^\vee\) by an integral invertible right factor
above \(N\). For an arbitrary point \(g\) the lattice norm of
the corresponding term \(u(g u_y)\) at \(g\delta_0^\vee\)
is precisely its lattice norm at \(g u_y\delta_0\).
The coefficient \(\psi_N(y_{\rm last})\) is a unit.
It follows that \(U_0\) takes the submodule of integral vectors at
\(\delta_0\) into the one at \(\delta_0^\vee\).
The previous clearing step, applied first to the original
tuple and then to the tuple at each \(U_0\)-iterate, gives
contribution \((m/2+1)h_{\rm deg}\) to obtain integrality at \(\delta_0\) after \(m/2\) iterations. Each clearing acts
on the original tuple by linearity. Since \(a=xb\) and
\(k\ge x\) for sufficiently large \(x\), the contribution is
at most \(mf bd_X\). Intersecting this preimage with
the original ramified source preserves this clearing
bound, since every operation takes integral linear combinations
of the original ramified columns. Explicitly, for a finite adelic element \(d\) put
\[
 \mathcal A(d)=\{u:\|u(g)\|_{gd}\le1
                         \text{ for every finite source point }g\}.
\]
Here the argument of the form is \(g\) and its Weyl lattice is transported
by \(gd\). Put
\[
 \mathcal A'_N=\mathcal A(1)\cap S
     \cap (U_0^{m/2})^{-1}\mathcal A(\delta_0).
\]
For every tuple \((v_1,\ldots,v_s)\) in \(\mathcal A(1)\cap S\)
the operations give a matrix \(H\in M_s(\mathcal O)\) such that
\[
 (v_1,\ldots,v_s)H\in(\mathcal A'_N)^s,\qquad
 \det H\mid N^{(m/2+1)h_{\rm deg}},\qquad
 v_N(\det H)\le mfbd_X .
\]
On a finite coefficient DVR this bounds the length of the
quotient by the intersection with \(\mathcal A'_N\), divided by
the ramification index, for the lattice spanned by any independent
tuple. Indeed the sublattice generated by its transformed tuple has
that normalized index \(v_N(\det H)\).

The preimage is fixed by the level and \(\jmath\)
before choosing the twist or translated arrow data.
The loss in the lattice norm \(2mf b\) has zero saving at \(N\)
since \(N\ge x_1\). The separate branching factor at
infinity occurs in the norm calculation, with the
same branching exponent \(B_\tau=a_\tau-mh_0\).

\section{Integral theta sections and local derivatives}
\label{sec:local-theta-derivatives}

\sectionopening{In this section, we bound the derivatives of the integral theta sections occurring in
Section~\ref{sec:measured-theta-tensor}. The exceptional primes require the
additional Hecke calculation in Section~\ref{sec:exceptional-hecke-analysis}.}

\subsection{Integral sections and ordinary jets}
\label{sec:marked-integral-sections-ordinary-jets}

\sectionopening{In this subsection, we construct an integral CM model and control the
derivatives of translated theta sections in that model. The derivative
estimate applies to the translates in
Section~\ref{sec:seesaw-integral-tensor} and retains their actual
translation denominators.}

Let \(K=\mathbb Q(\sqrt{-d})\), \(d\) a prime congruent to \(7\) modulo \(8\), and \(\beta=(1+\sqrt{-d})/2\). Write \(\Omega_0=\omega_0\) for the period of Section~\ref{sec:upper-cm-jet-transport} and \(\Lambda=\Omega_0\mathcal O_K\) for its CM lattice, \(A\) for its area, and \(t\) for its marked half-period. Modulo \(\Lambda\), the three nonzero points of order two are \(t,\beta t,\bar\beta t\). The scale is fixed by \(\wp''(t)=2\). Put
\[
 x(z)=\wp(z+t)-\wp(t),\quad
 a_1^2=12\wp(t),\quad
 y(z)=\tfrac12(\wp'(z+t)-a_1x(z)).
\]
The equation and differential relations are
\[
 y^2+a_1xy=x^3+x,\qquad
 x'=2y+a_1x,\qquad y'=1+3x^2-a_1y.
\]
For every ideal class and every normalization at a
half-period \(u_i\), the equation in completed-square
coordinates is \(Y^2=X(X^2+3e_i X+1)\).

Its invariant is
\[
 j=256\,\frac{((3e_i)^2-3)^3}{(3e_i)^2-4}.
\]
Complex multiplication makes \(j\) algebraic
\cite[Theorem~12.8]{MilneMF}, so each normalized \(e_i\) and \(a_1\)
is algebraic. The coordinates of the torsion points are
algebraic as well, since the multiplication maps are
defined over \(\overline{\mathbb Q}\).
Here \((x(z),y(z))\) are the Weierstrass coordinates of
\(z+t\) on the group curve. This shifted parametrization sends
\(z=0\) to \((0,0)\) and \(z=t\) to the point at infinity.
Thus \(y(az)\) below evaluates the Weierstrass coordinate at
\(az+t\). The Taylor parameter at \((0,0)\) is \(y\), with
\(y(0)=0\) and \(y'(0)=1\).
At the fixed coefficient place above \(p\), use
\(\iota_p:\mathbb C\simeq\overline{\mathbb Q}_p\hookrightarrow\mathbb C_p\)
of Section~\ref{sec:notation-determinants} and set
\(R=\mathcal O_{\mathbb C_p}\). Apply this embedding coefficientwise
to each Taylor series.

\smallskip
\indent\textit{The discriminant at all primes.}
Let \(e_i\) range over the three half-period values on a CM lattice of any ideal class. If \(u_i\) is the corresponding half-period, set
\[
 s(u_i)=16\frac{(e_j-e_k)^2}{(e_i-e_j)(e_i-e_k)}\qquad(\{i,j,k\}=\{1,2,3\}).
\]
This expression is invariant under homothety. The identity
\(\wp''(u_i)=2(e_i-e_j)(e_i-e_k)\) shows that \(s(u_i)\) is the discriminant after normalizing \(\wp''(u_i)=2\). Taking a common denominator in the sum and multiplying the three expressions gives
\[
 \sum_i s(u_i)=-48,\qquad \prod_i s(u_i)=-2^{12}.
\]
For the quotient by the point \(u_i\), translate \(e_i\) to zero. In the equation \(Y^2=X(X^2+aX+b)\) one has \(a=(e_i-e_j)+(e_i-e_k)\), \(b=(e_i-e_j)(e_i-e_k)\). The quotient by \((0,0)\) has equation \(Y^2=X(X^2-2aX+a^2-4b)\), with dual kernel \((0,0)\). Indeed, the isogeny is given off the kernel by
\((X,Y)\mapsto(X+a+b/X,Y(1-b/X^2))\), as substitution in the equations verifies. The other two old half-periods map to \((0,0)\). Hence the two values of \(s\) at the original and dual kernels are \(16(a^2-4b)/b\) and \(2^8b/(a^2-4b)\), respectively.

Index a family of CM lattices by the ideal classes \(c\), including
the original lattice. On an invertible ideal \(\mathfrak a\),
\((\tfrac12\mathfrak a)/\mathfrak a\) is free of rank one over
\(\mathcal O_K/(2)\simeq\mathbb F_2\times\mathbb F_2\). At each
class use its unique generator \(t\), which is nonzero on both
summands. The idempotents \(\beta,\bar\beta\) give the two other
nonzero points. Homotheties preserve these choices. The quotient by \(\beta t\) replaces an ideal \(\mathfrak a\) by \((1,\beta/2)\mathfrak a\) and carries \(\bar\beta t\) to the nonzero \(\bar\beta\)-half-period of the quotient. In fact \(\beta\) is the projection to the line generated by \(\beta t\) on \(\frac12\mathfrak a/\mathfrak a\), so this is the indicated quotient lattice.

Use \(g\) for the class of \((1,\beta/2)\). With \(s_c,s'_c,s_{0,c}\) the values at \(\beta t,\bar\beta t,t\), the product calculation gives
\[
 s_c s'_{gc}=2^{12},\qquad s_{gc}s_{0,gc}=-s_c.
\]
If \(|s_c|_p>1\), then \(|s'_{gc}|_p\le1\), \(|s_{gc}+s_{0,gc}|_p\le1\), and \(|s_{gc}s_{0,gc}|_p=|s_c|_p>1\). The ultrametric inequality gives \(|s_{gc}|_p=|s_{0,gc}|_p=|s_c|_p^{1/2}>1\). Iterating around the finite cycle contradicts \(|s_c|_p>1\). Interchanging \(\beta,\bar\beta\) proves \(|s'_c|_p\le1\). The sum shows \(|s_{0,c}|_p\le1\), and multiplying the second identity through all classes yields \(\prod_c s_{0,c}=(-1)^{h_K}\). Each \(s_{0,c}\) is a unit. In the chosen Weierstrass model this gives
\(\Delta=a_1^4-64\in R^\times\) and \(a_1\in R\), including at \(p=2\).

\smallskip
\noindent\textit{Theta sections and exact multiplication.}
Set \(H(u,z)=\pi\bar u z/A\). Use the odd sigma section
\(\sigma\) with \(\sigma'(0)=1\) and period law
\[
 \sigma(z+l)=\varepsilon(l)e^{H(l,z)+H(l,l)/2}\sigma(z),\qquad
 \varepsilon(l)=
 \begin{cases}1,&l\in2\Lambda,\\-1,&l\in\Lambda\setminus2\Lambda.\end{cases}
\]
These Appell--Humbert data specify a degree-one line bundle.
Its space of sections has dimension one. With this multiplier,
\(\sigma\) is the reduced theta function of the divisor \([0]\);
the area parameter in \cite[Definition~1.1 and (1.2)--(1.3), p.~190]{BannaiKobayashiTsuji}
is \(A/\pi\) in our notation. For a degree-\(m\) section set
\[
 (T_v^{(m)}G)(z)=e^{-m(H(v,z)+H(v,v)/2)}G(z+v),\qquad
 G[u]_m=(T_u^{(m)}G)(0).
\]
Put \(\theta=T_t^{(1)}\sigma/(T_t^{(1)}\sigma)(0)\). Oddness of
\(\sigma\) and the sign \(\varepsilon(2t)=-1\) in the period law give
\(\theta(-z)=\theta(z)\). Its degree-one semicharacter is
\[
 \chi(l)=\exp(i\operatorname{Im}H(\bar\beta l,l))
 =(-1)^{|l|^2\sqrt d/(2A)}\qquad(l\in\Lambda).
\]
For the sign calculation, the sigma sign at a period \(v\) is \(+1\) if \(v/2\in\Lambda\) and \(-1\) otherwise. Twisting sigma by \(t\) gives the sign at \(al\), for \(a=\beta,\bar\beta\), as the negative of the sigma sign at \(2t+al\). For \(z=t+al/2\) one has \((1-a)z=(1-a)t+|a|^2l/2\notin\Lambda\) since \(|a|^2=(d+1)/4\) is even. Thus \(\chi(al)=1\). The semicharacter rule at \(l=\beta l+\bar\beta l\) gives the displayed formula. Write \(L_\theta\) for the line bundle of \(\theta\).
A section \(G\) of \(L_\theta^{\otimes m}\) satisfies
\(T_l^{(m)}G=\chi(l)^mG\). Both \(\varepsilon\) and \(\chi\)
take values in \(\{1,-1\}\), so \(\sigma^2\) is a section of
\(L_\theta^{\otimes 2}\).

If multiplication by \(a\) is an isogeny into another such lattice of area \(A'\) and \(|a|^2A=NA'\), then \(H'(au,az)=NH(u,z)\) and \(\chi'(al)=\chi(l)^N\). Consequently pullback takes degree \(m\) sections to degree \(mN\), with
\[
 (G\circ[a])[u]_{mN}=G[au]_m .
\]
Also
\(
 (T_v^{(m)}G)[u]_m=e^{\pi i m E(u,v)}G[u+v]_m
\)
where \(E(u,v)=\operatorname{Im}H(u,v)/\pi\). The phase is a root of unity at torsion points since \(E\) is integral on \(\Lambda\). In particular if \(\delta u\in\Lambda\), \(\delta\in\mathcal O_K\setminus\{0\}\), then
\[
 z\longmapsto (T_u^{(1)}\theta)(\bar\delta z)
 \quad\text{has degree }N=\delta\bar\delta .
\]
Indeed \(T_{\bar\delta l}^{(1)}\) commutes with \(T_u^{(1)}\) because \(E(\bar\delta l,u)=E(l,\delta u)\in\mathbb Z\). Its remaining multiplier is \(\chi(\bar\delta l)=\chi(l)^N\).

\smallskip
\noindent\textit{Interpolation on the curve.}
The classical Weierstrass addition and duplication formulas give
\[
 \sigma^2=x\theta^2,\qquad
 \theta(2z)=\theta(z)^4-\sigma(z)^4,\qquad
 \sigma(2z)=x'(z)\theta(z)^4.
\]
For the normalization in these formulas, the addition formula at \((z+t,t)\) proves the first identity after comparing the constant at \(z=0\) of \(\sigma(z)^2/(x(z)\theta(z)^2)\). Indeed this quotient is elliptic with no zero or pole and tends to one there. Twisting twice by \(t\) gives the sigma sign \(-1\) at \(2t\). If \(c=(T_t^{(1)}\sigma)(0)\), translating \(\sigma^2=x\theta^2\) by \(t\) gives \(c^4=x(z)x(z+t)=1\), the last identity being the addition formula at a half-period with \(\wp''(t)=2\). The sigma duplication formula at \(z+t\) now gives the last displayed identity, with the sign at \(2t\) cancelling its minus sign. For the middle identity, apply the sigma addition formula first to \((z+t,v)\) and divide out the common exponentials. This gives \(\theta(z+v)\theta(z-v)=\theta(z)^2(\theta(v)^2-x(z)\sigma(v)^2)\). Set \(v=z\).

For \(w(u)=(\theta[u]_1^2,\sigma[u]_1^2)\), the identities imply \(w(2u)=D(w(u))\) with
\[
 D(U,V)=\big((U^2-V^2)^2,UV(4U^2+a_1^2UV+4V^2)\big).
\]
First, the two coordinates of \(w(u)\) are algebraic
when \(u\) is torsion. The identity \(\sigma^2=x\theta^2\)
gives its projective coordinate as \([1:x(u)]\), or \([0:1]\)
at the unique zero of \(\theta\). Write \(w(u)=\lambda v\)
for a nonzero \(v\in\overline{\mathbb Q}^{\,2}\). Choose
\(a<b\) with \(2^a u\equiv 2^b u\bmod\Lambda\). The section
law at torsion gives \(D^b(w(u))=\zeta D^a(w(u))\) for a
root of unity \(\zeta\). In characteristic zero the two
forms of \(D\) have no common projective zero, since the
curve discriminant is nonzero. For a nonzero coordinate
of \(D^a(v)\) the corresponding coordinate of \(D^b(v)\)
is nonzero as well. The homogeneity of \(D\) then gives
\[
 \lambda^{4^b-4^a}
       =\zeta\,\frac{D^a(v)_i}{D^b(v)_i}\in\overline{\mathbb Q}^{\times}.
\]
This justifies the coefficient-place absolute values on \(w(u)\).

The reduction of these forms has no common projective zero. A zero of the first must have \(U=\pm V\ne0\); the second then vanishes only if \(a_1^2=\pm8\) in the residue field, contrary to \(\Delta\in R^\times\). This also applies in characteristic two. Therefore \(\|D(w)\|_p=\|w\|_p^4\) in the maximum norm, by reducing after dividing by a coordinate of largest norm. The norm is unchanged modulo \(\Lambda\) because sigma and theta have sign multipliers at periods and the remaining torsion phases have norm one. A torsion orbit under doubling takes only finitely many norms, all positive. The identity forces every one to equal one.

The section \(G_3=y\theta^3\) is entire and its twisted torsion values are in \(R\), by the monic quadratic obtained from the curve equation,
\[
 G_3[u]_3^2+a_1\theta[u]_1\sigma[u]_1^2G_3[u]_3
 =\sigma[u]_1^6+\theta[u]_1^4\sigma[u]_1^2 .
\]
For \(u_i=t/2^i\), \(i\ge1\), the reductions of both coordinates of \(w(u_i)\) are nonzero. Indeed both coordinate points map under \(D\) to \((1:0)\), which is fixed, whereas \(D^i(w(u_i))\) reduces to \((0:1)\). Thus \(x(u_i),y(\pm u_i)\in R\). The differences \(x(u_i)-x(u_j)\) are units for \(i<j\), by applying \(D^i\) to the projective reductions and comparing \((0:1)\) with the interior point from \(u_{j-i}\). The normalized duplication identity is
\[
 \sigma[2u_i]_1=x'(u_i)\theta[u_i]_1^4 .
\]
The denominator \(\theta[u_i]_1^4\) and the numerator are units,
including at \(2u_1=t\). Hence \(y(u_i)-y(-u_i)=x'(u_i)\) is
a unit also at \(p=2\).

Let \(G\) be a section of \(L_\theta^{\otimes m}\) whose twisted torsion values are in \(R\). The elliptic function \(G/\theta^m\) has its only possible pole at \(t\), of order at most \(m\). Since \(x,y\) have pole orders two and three there, successive elimination of leading terms gives
\[
 G/\theta^m=P_0(x)+yP_1(x),\qquad
 2\deg P_0\le m,\qquad 3+2\deg P_1\le m
\]
with missing terms omitted. Elimination uses the fact that a function with a sole pole of order one is impossible on an elliptic curve, by the residue theorem.

Evaluating at \(u_i,-u_i\) gives integral values because \(\theta[\pm u_i]_1\) is a unit. Their difference is \(x'(u_i)P_1(x(u_i))\). Unit Vandermonde interpolation implies \(P_1\in R[X]\), then \(P_0\in R[X]\). The curve equation solves for \(x\) in \(R[[y]]\) at the origin, since its partial derivative in \(x\) at \((0,0)\) is a unit. In particular
\[
 \frac{(T_u^{(1)}\theta)(\bar\delta z)}{\theta(z)^N}\in R[[y(z)]]\qquad(\delta u\in\Lambda).
\]
Here the hypothesis on the numerator's torsion values follows from the pullback identity, the torsion translation phase and the norm computation.

For multiplication by any nonzero \(a\in\mathcal O_K\), apply the same interpolation to the pullbacks of \(\theta\) and \(G_3\). It gives
\[
 A_a=\frac{\theta(az)}{\theta(z)^{a\bar a}}\in R[[y(z)]],
 \qquad B_a=\frac{G_3(az)}{\theta(z)^{3a\bar a}}\in R[[y(z)]].
\]
Since \(A_a(0)=1\), division gives \(y(az)=B_a/A_a^3\in R[[y(z)]]\).
The finite Vandermonde systems above also show that whenever the
twisted torsion values of the numerator are algebraic, the resulting
coefficients of \(P_0,P_1\) are algebraic. This applies to these pullbacks
and to the translated numerator above. The differential equations for
\(x,y\), with algebraic initial values, give algebraic Taylor
coefficients by recursion. For \(g=-\log\theta\), the quadratic
coefficient of \(A_\beta\) is
\(\beta(\bar\beta-\beta)g''(0)/2\), so \(g''(0)\) is algebraic.
To compute \(g''\), write the reduced sigma section as
\(\exp(-e_2^*z^2/2)\sigma_{\rm W}(z)\), with
\(\sigma_{\rm W}\) the Weierstrass sigma function
\cite[(1.2), p.~190]{BannaiKobayashiTsuji}. Since
\(-(\log\sigma_{\rm W})''=\wp\) and the translation contributes
a linear term to the logarithm, \(g''(z)=\wp(z+t)+e_2^*\).
Hence \(g''(z)-g''(0)=x(z)\), and \(g(0)=g'(0)=0\) gives
algebraic coefficients for \(g\) and \(\theta\). The translated series
follows from its quotient and rescaling by the nonzero \(\bar\delta\).

For the limits, let \(Z(Y)\) be the compositional inverse of \(y\).
For every \(j\) the coefficients of \(Y^j\) in
\[
 y(aZ(Y)),\qquad \theta(aZ(Y))\exp(sg(Z(Y)))
\]
are polynomials in \(a,s\) over \(\mathbb C_p\). Indeed
\(Z(Y)\in Y\mathbb C_p[[Y]]\) and \(g(Z(Y))\in Y^2\mathbb C_p[[Y]]\).
Writing \(\theta(z)=\sum_m\theta_m z^m\) gives
\[
 [Y^j]\theta(aZ(Y))e^{s g(Z(Y))}
 =\sum_{m+2n\le j}\frac{\theta_m a^m s^n}{n!}
       [Y^j] Z(Y)^m g(Z(Y))^n .
\]
Likewise \([Y^j]y(aZ(Y))\) has degree at most \(j\) in \(a\).
These polynomials are continuous at every coefficient prime.
Since their values lie in the closed valuation ring
\(R\) at \((a,s)=(a,a\bar a)\) for \(0\ne a\in\mathcal O_K\),
they lie in \(R\) at limits of those pairs. At an inert prime
take CM integers tending to \(\zeta\in\mu_{p^2-1}\) in
\(\mathcal O_{K,p}\). Conjugation gives the limit pair
\((\zeta,\zeta^{p+1})\). At a ramified prime take rational
integers tending to \(\zeta\in\mu_{p-1}\subset\mathbf Z_p^\times\).
The limit pair is \((\zeta,\zeta^2)\).
Both specializations use the same normalized even theta section.
For every \(0\ne a\in\mathcal O_K\) the quotient \(A_a\) is even.
Since \(x\) is even and \(y(z)-y(-z)=x'(z)\ne0\), its interpolating
polynomial has \(P_1=0\). Thus \(A_a\in R[x]\) and \(2\deg_x A_a\le a\bar a\).
For odd \(a\bar a\) this gives \(2\deg_x A_a\le a\bar a-1\), both for
\(a=p\) at an inert prime and for \(a=\sqrt{-d}\) at the ramified prime \(p=d\).

In the same way, the even section \(\sigma(az)^2\) has integral quotient \(S_a(x)=\sigma(az)^2/\theta(z)^{2a\bar a}\) with \(\deg S_a\le a\bar a\). At least one coefficient of \(S_a\) is a unit. At \(u=t/a\), its expression before division by \(\theta^{2a\bar a}\) is a sum of monomial sections \(\theta^{2a\bar a-2i}\sigma^{2i}\), each with integral torsion evaluation, while the pullback identity gives the unit value \(\sigma[t]_1^2\) for their linear combination. The quotient \(x(az)=S_a(x)/A_a^2\) is therefore nonzero after reduction and has order at most \(2a\bar a\) in \(y\). Since \(x=y^2+O(y^3)\) over \(R\), the order of the reduction of \(y(az)\) is at most \(a\bar a\). This proves the needed degree bound on the original formal parameter from the torsion interpolation itself. For the inert multiplication-by-\(p\) calculation at the Teichm\"uller-averaged parameter, the reduced series has zero linear term and every allowed term has degree \(1\pmod{p^2-1}\); thus its first nonzero term has degree exactly \(p^2\).

\smallskip
\indent\textit{A basis with denominator prime to the residue characteristic.}
Use the cyclotomic field \(E=K(\zeta)\), with \(\zeta\) of order \(q=q_1\), and its lattice \(\mathfrak a_0=E\cap\alpha_{\rm W}\widehat{\mathcal O}_E\). In this model \(\mathfrak a_0\) is the product of the inverse relative different at the primes over \(\beta\), and the unit ideal elsewhere. It is stable under \(\operatorname{Gal}(E/K)\), it contains \(\mathcal O_E\), and it is self-dual for \(\operatorname{Tr}_{E/K}(x\bar y)\): conjugation interchanges the two sets of primes of the relative different above the split primes dividing \(q\). Write \(\Phi\) for the embeddings over the fixed embedding of \(K\). The powers of \(\zeta\) give Gram sum \(qI\) under these embeddings by the finite geometric sum. For \(p\nmid q\) this already has suitable length \(b=q\).

For \(p\mid q\), let \(j=\lceil\log_2 q\rceil\) and take \(e_i=\bar\beta^j d_i\), where \(d_i\) is trace-dual to \(\zeta^i\), for \(0\le i<[E:K]\). Then \(e_i\in\mathfrak a_0\). Indeed \(d_i\) is in the relative codifferent. The codifferent denominators are supported over \(q\) and are cleared by \(q\), because \(\Phi_q'(\zeta)\) divides \(q\) by differentiating the factorization of \(X^q-1\). For every rational \(\ell\mid q\), \(v_\ell(q)\le j\) and \(\beta\bar\beta\) is divisible by \(\ell\).

Since \(\beta+\bar\beta=1\), the factor \(\bar\beta^j\) clears the codifferent at primes over \(\bar\beta\). Take \(\zeta^{-k}\) for \(0\le k<q\) and all the vectors \(\zeta^k\sigma(e_i)\) with \(0\le k<q\), \(\sigma\in\operatorname{Gal}(E/K)\). Their Gram sum is \(q(1+\sum_i\nu_i)I\), where \(\nu_i=\operatorname{Tr}_{E/K}(e_i\bar e_i)\) is a nonnegative integer by self-duality and reality. Replace each pair \(\zeta^{-i},e_i\) with \(\zeta^{-i}+e_i\). This adds \((\beta^j+\bar\beta^j)I\) by trace duality, since the products of the two square embedding matrices in either order are identity.

The new scalar \(b=q(1+\sum_i\nu_i)+\beta^j+\bar\beta^j\) is positive and prime to \(q\). At any prime over \(\ell\mid q\), one of \(\beta,\bar\beta\) is zero and the other is one, so the last two terms sum to one modulo \(\ell\). In both cases repeat each vector \(b\) times to obtain \(w_i\in\mathfrak a_0\) with Gram sum \(b^2 I\).

\smallskip
\indent\textit{The finite vacuum expansion.}
Let \(L=\mathfrak a_0^n\), embedded in the real Fock coordinates by \(r(\gamma)_{\tau,A}=a_d\tau(\gamma_A)\), where \(a_d=\sqrt{2\pi/\sqrt d}\). This is the self-dual lattice for the Riemann form. Put
\[
 A_{i,A}(z)=\frac{\Omega_0}{b\,a_d}\sum_{\tau\in\Phi}
       \overline{\tau(w_i)}z_{\tau,A},\qquad
 S(z)=\prod_{i,A}\theta(A_{i,A}z).
\]
The area \(A=|\Omega_0|^2\sqrt d/2\) and the basis identity give \(\sum A_{i,A}^*H_C=H_{\rm Fock}\). Also \(A_{i,A}(r(bL))\subset\Lambda\) by trace self-duality. For \(l=b\gamma\in bL\) each period is
\(A_{i,A}(r(l))=\Omega_0\operatorname{Tr}_{E/K}(\gamma_A\overline{w_i})\).
By the frame identity its product of theta signs is
\[
 \begin{aligned}
 \prod_{i,A}\chi(A_{i,A}(r(l)))
 &=(-1)^{\sum_{i,A}|\operatorname{Tr}_{E/K}(\gamma_A\overline{w_i})|^2}\\
 &=(-1)^{b^2\sum_A\operatorname{Tr}_{E/K}(\gamma_A\bar\gamma_A)}
  =(-1)^{\sum_A\operatorname{Tr}_{E/K}(l_A\bar l_A)}.
 \end{aligned}
\]
This is the finite vacuum sign of
Section~\ref{sec:upper-cm-jet-transport}, also when \(b\) is odd.
We have \(S(0)=1\).

The same finite Fourier calculation constructs the product vector and resolves the vacuum. Write \(U_l\) for signed periods and \(T_t\) for ordinary twisted real translations. The groups \(L/bL\) and \(b^{-1}L/L\) are perfectly paired by the translation commutator and both have order \(M_L=b^{2|\Phi|n}\). Let \(\mathcal V_{\rm int}\) be the entire theta image of
\(s_{\rm vac}^{n}v_0^{\otimes n}\), with \(v_0\) the rank-one unit
vacuum of \eqref{eq:upper-CM-vacuum}. On the \(bL\)-periodic space,
\[
 P_tS=\sum_{l\in L/bL}U_lT_tS=d_t \mathcal V_{\rm int}
 \quad(t\in b^{-1}L/L).
\]
Indeed the sum is invariant under \(L\), and the line of theta sections on the principally polarized torus is one-dimensional.

The commutator orthogonality gives
\[
 M_LS=\sum_t T_{-t}P_tS=\sum_t d_tT_{-t}\mathcal V_{\rm int} .
\]
This identifies \(S\) as the entire theta image of a finite-part vector \(F\). In the finite Weil model a rational real translation by \(r(\gamma)\) is induced by the operator at the opposite finite component, by rational invariance of the theta functional. To prove injectivity, let \(f\) have zero entire theta
image. The Gaussian factor is nonzero, so
\(\theta(R(z_\infty)(v_\infty\otimes f))=0\) for every
real displacement \(z_\infty\).

For any finite displacement
\(z_f\), choose a rational \(\gamma\) with \(z_f-\gamma_f\)
in an open translation stabilizer of \(f\).
Rational Heisenberg invariance identifies the evaluation
at \(z_\infty+z_f\), up to a unit phase, with that at
\(z_\infty-\gamma_\infty\). Hence all adelic displacement
evaluations vanish. The polarized norm identity of
Section~\ref{sec:overlap-fourier-inversion} gives
\(\|v_\infty\otimes f\|^2=0\). Since \(v_\infty\ne0\),
\(f=0\), and the identities hold in the finite Weil model.

Evaluating at zero, \(M_L=\sum_t d_t(T_{-t}\mathcal V_{\rm int})(0)\), and substituting \(d_t\mathcal V_{\rm int}=P_tS\) gives
\[
 \mathcal V_{\rm int}=M_L^{-1}\sum_t (T_{-t}\mathcal V_{\rm int})(0)\sum_{l\in L/bL}U_lT_tS .
\]
Each \((T_{-t}\mathcal V_{\rm int})(0)\) is integral. Indeed \eqref{eq:upper-CM-vacuum}
expresses the scaled rank-one vacuum as a sum of finite translates of
the cyclotomic theta product. Tensoring this identity in the \(n\) coordinates
and then translating gives a sum of products of the integral twisted
torsion values, with roots of unity as coefficients.
Since \(p\nmid M_L\), the displayed equality expresses the original
vacuum as an integral combination of translates of \(S\). Each translating
point becomes integral after multiplication by \(b\).

In rank one the integral vacuum is \(s_{\rm vac}v_0\) from
\eqref{eq:upper-CM-vacuum}; in rank \(n\) it is
\((s_{\rm vac}v_0)^{\otimes n}\). The theta evaluations of its
signed finite translates are algebraic integers by the torsion calculation
and the signed translation formula. The expansion of this fixed vector
may use a different \(b\) at each \(p\).

\smallskip
\indent\textit{Ordinary integral jets.}
Let \(f\) lie in the integral span of the signed finite translates of this integral vacuum by \(x\) with \(\delta x\) integral above \(p\), where \(\delta\in\mathcal O_K\setminus\{0\}\). Using the preceding expansion, it suffices to treat finitely many translates of \(F\) with the same condition. In standard finite coordinates a rational vector \(\gamma\)
acts by \(\alpha_{\rm W}^{-1}\gamma\). For each fixed \(x\),
choose an open neighborhood of zero on which \(T(h)F=F\),
\(h\) is integral above \(p\), and the phases of
Lemma~\ref{lem:weil-vacuum-composition} satisfy
\(\beta(x,h)=1\) and \(e_0(x+h)=e_0(x)\).
The two phase conditions are open because the finite additive
character is locally constant. By density choose \(\gamma\) so
that \(h=\alpha_{\rm W}^{-1}\gamma-x\) lies in this neighborhood.
Then \(T(x+h)F=T(x)F\), also with the signed normalization
\(T_\sharp(x)=e_0(x)T(x)\). Above \(p\) we have
\(\delta\alpha_{\rm W}^{-1}\gamma\) integral, so
\(\delta\gamma\in L\otimes\mathbb Z_p\). Only finitely many other
places contribute denominators, so one integer \(m\) prime to \(p\)
makes \(m\delta\gamma\in L\) in every term.

For the translate by \(\gamma\), the factor \((i,A)\) has shift
\(
 t_{i,A}=-(\Omega_0/b)\operatorname{Tr}_{E/K}(\gamma_A\overline{w_i}).
\)
Choose \(k=bm\), so \(k\delta t_{i,A}\in\Lambda\) by trace
self-duality. For each such shift \(s=t_{i,A}\), interpolation with
\(k\delta\) gives
\[
 (T_s^{(1)}\theta)(k\bar\delta Z)
       =\theta(Z)^{k^2 N}\,P_s(y(Z)),\qquad
       P_s\in R[[Y]],\quad N=\delta\bar\delta .
\]
In degrees zero and one this reads
\[
 P_s(0)=(T_s^{(1)}\theta)(0),\qquad
 P_s'(0)=k\bar\delta (T_s^{(1)}\theta)'(0)
\]
because \(\theta(0)=y'(0)=1\) and \(\theta'(0)=0\). This includes shifts
where the translated section vanishes at the origin.
At the argument \(\bar\delta(a_d/\Omega_0)z\), substitute
\(Z=Z_{i,A}=(w_i/(bk))^\dagger z_{\cdot,A}\).
The product \(\prod_{i,A}\theta(Z_{i,A})^{k^2N}\) is independent of
the translating vector because \(k\) was chosen for the entire finite
sum. The remaining factor is an integral sum of the products
\(\prod_{i,A}P_{t_{i,A}}(y(Z_{i,A}))\).

Write \(s_{\rm jet}=p_Fa_d/\omega_0\). In the original normalization the absolute value of \(p_F\) at \(p\) is that of \(\alpha_{\rm W}^{-1}\) at each weight embedding \(a\mapsto\iota_p(\overline{\tau(a)})\).
Indeed Section~\ref{sec:upper-cm-jet-transport} gives \(p_F=s_0/u\),
with \(u\) a unit at every coefficient embedding by the class-cycle
argument, and \(s_0^f=\lambda_f(\alpha_{\rm W})^{-1}\) for
\(\lambda=\upsilon|\cdot|_E^{-1/2}\).
Let \(h>0\) and \(y\in E^\times\) satisfy
\((y)=\mathfrak a_0^h\), and put \(w_\tau^-=w(\iota_p\circ\bar\tau)\),
where \(\bar\tau(a)=\overline{\tau(a)}\). Normalize the valuations
\(v_w\) by \(v_w(p)=1\). The identity
\(s_0^{-fh}=\prod_\tau\overline{\tau(y)}\) in
Section~\ref{sec:upper-cm-jet-transport} gives
\[
 v_p(\iota_p(p_F))=-\frac1f\sum_\tau v_{w_\tau^-}(\alpha_{\rm W})
                   =-v_{w_\tau^-}(\alpha_{\rm W})
                    \qquad(\tau\in\Phi).
\]
All the \(w_\tau^-\) restrict to the same place of \(K\), and
\(\mathfrak a_0\) is invariant over \(K\). For \(p\mid q_1\) write
\(q_1=p^a m'\), \(p\nmid m'\). At either place of \(K\) above \(p\)
the completion of \(K\) is \(\mathbb Q_p\); adjoining
\(\zeta_{m'}\) is unramified. Over this unramified field
\(\Phi_{p^a}(1+X)\) is Eisenstein, and \(\zeta_{p^a}-1\)
generates the ring of integers as an algebra. The different therefore has
normalized valuation
\[
 v_p\bigl(\Phi_{p^a}'(\zeta_{p^a})\bigr)
  =v_p\!\left(\frac{p^a\zeta_{p^a}^{-1}}
                     {\zeta_{p^a}^{p^{a-1}}-1}\right)
  =a-\frac1{p-1}.
\]
Here the last equality follows from \(\Phi_p(1)=p\) and the equal
valuations of \(1-\zeta_p^j\) for \(1\le j<p\). Thus
\[
 v_p(\iota_p(p_F))=
 \begin{cases}
 a-\dfrac1{p-1},&w_\tau^-\mid(\beta),\\[2pt]
 0,&w_\tau^-\mid(\bar\beta).
 \end{cases}
\]
Conjugating the coefficient embedding interchanges the cases.
Since \(w_i\in\mathfrak a_0\), in both cases
\(v_p(\iota_p(p_F\overline{\tau(w_i)}))\ge0\).
This includes \(p=2\). Every choice of the prescribed root
differs by a root of unity.
Before applying the coefficient embedding the fixed jet argument is
\[
 A_{i,A}(s_{\rm jet}z)
       =\frac{p_F}{b}\sum_{\tau\in\Phi}\overline{\tau(w_i)}z_{\tau,A}.
\]
The factors \(a_d\) and \(\Omega_0=\omega_0\) cancel in this identity.
Since \(w_i/(bk)\in\mathfrak a_0\otimes\mathbb Z_p\) and \(p\nmid bk\),
every coefficient of \(\nu=p_F(w_i/(bk))^\dagger\) lies in \(R\).

Let \(g=-\log\theta\), so \(g''=g''(0)+x\) and \(g'(0)=0\).
For \(a=\beta,\bar\beta\), the even interpolation above gives
\(A_a(z)=\theta(az)/\theta(z)^{a\bar a}=P_a(x(z))\), \(P_a\in R[X]\).
Since \(x(z)=z^2+O(z^4)\) by \(\wp''(t)=2\), the quadratic coefficient is
\[
   \frac{a(\bar a-a)}2g''(0)=P_a'(0)\in R .
\]
Subtract the value at \(\beta\) from that at \(\bar\beta\), using
\(\beta+\bar\beta=1\), to get \(\sqrt{-d}\,g''(0)/2\in R\).
Multiplication by \(-2\sqrt{-d}\) gives \(dg''(0)\in R\),
so \(dg''\in R[[y]]\). For \(\mu_i=d\nu_i\), take the ring
generated over \(R\) by integral series in the \(y(\mu_i z)\),
together with \(\theta(\mu_i z)\) and \(g'(\mu_i z)\).

Write \(Y_i=y(\mu_i z)\), \(X_i=x(\mu_i z)\in R[[Y_i]]\).
For \(B\in R[[(Y_i)_i]]\),
\[
 \partial_j B((Y_i)_i)
   =\sum_i(\partial_i B)((Y_i)_i)\,\mu_{i,j}(1+3X_i^2-a_1Y_i).
\]
The other generators satisfy
\[
 \partial_j\theta(\mu_i z)
   =-\mu_{i,j}\theta(\mu_i z)g'(\mu_i z),\qquad
 \partial_j g'(\mu_i z)=\nu_{i,j}(d g'')(\mu_i z).
\]
Thus the ring is preserved by all partial derivatives, and every
element has integral constant term. Applying any polynomial \(Q(\partial)\) with coefficients in \(R\) and evaluating at zero thus gives an integral value for the actual theta vector at scale \(\bar\delta d\,s_{\rm jet}\).

If \(Q\) is homogeneous of degree \(l\), rescaling gives the ordinary bound at \(c\,s_{\rm jet}\),
\[
 \left|Q(\partial)\Theta_f(c\,s_{\rm jet}\,z)\big|_{z=0}\right|_p
 \le \left|\frac{c}{\bar\delta d}\right|_p^l.
\]
For a multi-index \(\mathbf m\), with \(\mathbf m!=\prod_i m_i!\),
this bounds the ordinary derivative by
\[
 \mathbf m!\,[z^{\mathbf m}]\Theta_f(c\,s_{\rm jet}z)
      \in (c/(\bar\delta d))^{|\mathbf m|}R .
\]
For an integral source basis vector \(v\) and a coordinate
covector \(A\) of the integral target basis, put
\(Q_{A,v}(Z)=A\mathcal S_j(Z)v=\sum_{|\alpha|=l}q_\alpha Z^\alpha\).
Let \(F_s,F_t\) embed the source and target Weyl coordinates in
their exterior column tensors. The row of each matrix at a column
filling gives the straightening coefficients of its minor product.
Reading the rows indexed by the standard tableaux gives a left
inverse \(R_tF_t=I\). At \(\tau\), the two diagrams have the same
columns below \(n_j\), and the target has \(c_\tau=m_{\tau,n_j}\)
additional columns of height \(n_j\). Form \(\widehat F_s\) from
\(F_s\) by inserting \(e_1\wedge\cdots\wedge e_{n_j}\) in each
of those slots. Let \(C(Z)\) be the tensor product of the exterior
powers of \(Z\) over the target columns. For a source frame \(X\),
let \(X_\tau\) be the square block of its first \(n_j\) columns and
write \(v_s,v_t\) for the standard minor evaluation vectors.
Cauchy--Binet and the minor expansion defining \(\mathcal S_j\) give
\[
 \begin{aligned}
 F_t v_t(ZX)
  &=\Bigl(\prod_\tau\det(X_\tau)^{c_\tau}\Bigr)
                   C(Z)\widehat F_s v_s(X)\\
  &=\Bigl(\prod_\tau\det(X_\tau)^{c_\tau}\Bigr)
                   F_t\mathcal S_j(Z)v_s(X).
 \end{aligned}
\]
The vectors \(v_s(X)\) with every \(X_\tau\) invertible span the
rational source, since a linear combination of the independent
standard minor polynomials vanishing on this dense open set is zero.
Thus
\[
 C(Z)\widehat F_s=F_t\mathcal S_j(Z),\qquad
 \mathcal S_j(Z)=R_t C(Z)\widehat F_s .
\]
The latter matrix has integral polynomial entries, so
\(q_\alpha\in\mathbf Z\). By the monomial pairing of
Section~\ref{sec:even-stage-coefficients}, the corresponding
entry of the jet matrix is
\[
 j_{Q_{A,v}}(J_{c\,s_{\rm jet}}f)
  =\sum_{|\alpha|=l}q_\alpha\,\alpha!\,
          [z^\alpha]\Theta_f(c\,s_{\rm jet}z)
       \in (c/(\bar\delta d))^l R .
\]
This coordinate computation also applies at \(p=2\) and at every
\(p\mid q_1\) on the same trace lattice.

For a fixed finite sum of the signed translates with algebraic
coefficients, take a positive integer \(D_f\) clearing those coefficients
and \(0\ne\delta\in\mathcal O_K\) with every \(\delta x\) integral at
all finite places. In the finite expansion above, torsion interpolation
gives algebraic coefficients for the translated sections. Their linear
arguments have
coefficients in \(\overline{\mathbf Q}\) after multiplication by
\(s_{\rm jet}\). Thus, for \(c\) an algebraic integer, each
multi-index \(\mathbf m\) satisfies
\[
 D_f\bigl(d\,\mathrm N_{K/\mathbf Q}\delta\bigr)^{|\mathbf m|}
  \,\mathbf m!\,[z^{\mathbf m}]\Theta_f(c\,s_{\rm jet}z)
       \in\overline{\mathbf Z}.
\]
Indeed the local bound for \(D_f f\) applies at every coefficient
embedding, and \(d\,\mathrm N_{K/\mathbf Q}\delta/(\bar\delta d)=\delta\).
The preceding local bound retains the sharper denominator at any single
prime. In particular, in the fixed integral source basis and target
coordinate dual, the polynomial \(Q_{A,v}\) has integral
coefficients. The same integer \(D_f(d\,\mathrm N_{K/\mathbb Q}\delta)^l\)
clears \(j_{Q_{A,v}}(J_{c\,s_{\rm jet}} f)\) for each entry of degree
\(l\), simultaneously at all coefficient embeddings. Applying an
embedding to these cleared coordinates preserves integrality.

To pass from these derivatives to the kernel at finite arguments, let
\(a,b\) be integral invertible above \(p\) and let \(G\) be the induced
isometry \(X\mapsto bXa^{-1}\). In normalized finite coordinates
write \(L_f\) for the standard lattice and \(H=L_f\cap G^{-1}L_f\).
The lattice vacuum formula of Section~\ref{sec:upper-cm-jet-transport}
expresses the image of the translate by \(x\) as a sum of translates
by \(G(x+z)\), \(z\in L_f/H\). The coefficients are
\(\Delta_{G^{-1}}^{-1}\) times translation phases and the pair
characters. Above \(p\), \(H=L_f\) locally and
\(\delta G(x+z)=G(\delta x+\delta z)\) is integral. The determinant
ideal defining \(\Delta_{G^{-1}}\) is a unit at every place over
\(p\), as are \(\det a,\det b\). To evaluate an algebraic finite
character on any of these determinants, take a power making the
ideal principal. Its generator is a unit at every place over \(p\),
so the algebraic weight monomial is a \(p\)-unit. The residual
finite-order factor and all translation phases are roots of unity.
For a norm-one determinant one may take twice the class-group
exponent and a norm-one generator, as in
Section~\ref{sec:upper-cm-jet-transport}. Thus every coefficient is
a \(p\)-unit. The entire span with \(\delta x\) integral above \(p\)
is preserved even when \(a,b\) have denominators away from \(p\).
For arbitrary finite arguments \(g,h\), choose rational isometries
\(\gamma,\eta\) with \(\gamma g,\eta h\) integral invertible above
\(p\). Indeed every prescribed local open set meets the Cayley chart
\(\det(1+u)\ne0\). Multiply a point of the set by a norm-one scalar
arbitrarily near one, avoiding the finitely many roots of
\(\det(1+t u)\) in each field component. Such scalars are given by
\((1+s\zeta)/(1-s\zeta)\) for a nonzero skew scalar \(\zeta\) and
\(s\) in the local fixed field near zero. Weak approximation in the
skew-adjoint matrix space and the rational Cayley map give the
required isometry, also at dyadic places.

Write \(T_g\) for transport of the source coordinates by \(g^{-1}\)
at the coefficient embedding, and similarly on the target. Rational
covariance of the kernel \(K\) gives
\[
 T_{\eta h}K(\gamma g,\eta h)T_{\gamma g}^{-1}
             = T_hK(g,h)T_g^{-1}.
\]
Thus the bound is on the values at the original points. Left rational
translation carries the entire right coset \(gJ\) to \((\gamma g)J\)
at the same source level \(J\), and likewise at the target level.

\subsection{Valuations of the derivative coefficients}
\label{sec:marked-column-gain}

\sectionopening{In this subsection, we use the highest weights to count the columns
entering the derivative estimates of
Section~\ref{sec:marked-integral-sections-ordinary-jets}. This
expresses the local valuation in terms of the prime powers in the
arithmetic sum of Section~\ref{sec:joint}.}

Fix a coefficient prime \(p\) and a twist \(\xi\) with
\(p\nmid q_\xi\). For \(p\ne N\), each displacement in an individual
arrow of the expanded tensor is made integral above \(p\) by
\(\bar\pi_M\). Indeed the two adjoint blocks are returned to their
original order in the individual arrows, the additional translations
are supported above \(Mq_\xi\), and the linked residues in
Section~\ref{sec:seesaw-tensor-translation-N} are supported above
\(N\). At \(p=N\) these residues can have denominator \(N\).
Section~\ref{sec:ramified-direct-lift} treats the resulting sum.

Use \(n_j=k+j(t+1)\), \(t=5r\), \(\Delta=2r^8\),
and put \(w_j(\tau)=\kappa_\tau+(j+\tfrac12)\Delta\).
The weight on \(X_j\) has determinant exponent zero
at \(j=0\) and last entry \(\lfloor n_{j-1}/2\rfloor-w_j\)
for \(j>0\). Its columns of height \(h<n_j\) have counts
\[
 \begin{array}{c|c}
  h&\text{number}\\ \hline
  1&a_\tau=a\mathbf 1_{\Phi_a}(\tau)\\
  k&w_1(\tau)-\lfloor(k+t+1)/2\rfloor\\
  k+i(t+1),\ i\ge1&\Delta-(t+1).
 \end{array}
\]
All other counts vanish. The successive blocks of the
parameter have size \(t+1\), and have centers
\(-w_i(\tau)\). In each such block the parameter and
the highest-weight shift both decrease by one, so
the highest weight is constant there. At adjacent
block boundaries the parameter centers differ by
\(\Delta\), while the shifts differ by \(t+1\).

For \(j\ge1\), put \(l=(t+1)/2\). Subtracting the rank-\(n_j\)
half-sum of positive roots from the parameter gives the entries
on the first two blocks, in zero-based indices,
\[
 \mu_j(b)=-jl+a_\tau\mathbf1_{b=0}\quad(0\le b<k),\qquad
 \mu_j(k+b)=-jl+\lfloor n_1/2\rfloor-w_1(\tau)
                     \quad(0\le b<t+1).
\]
The common shift \(-jl\) cancels in the difference at height \(k\),
giving the table above. It remains in the determinant exponent.
For the diagram below \(n_{j+1}\), its first
\(k\) rows each contain \(\kappa_\tau+O(m\Delta)\)
boxes from the low columns other than height \(1\).
The remaining rows have \(O(m\Delta)\) such boxes.
In particular \(2(m+1)\Delta<r^{16}\) gives the additional term for these rows used in the pointwise size calculation.
For part with index zero its additional term \(b+a/x_1\) absorbs
\(\kappa_\tau\le b\).

For \(j<m\) let \(c_\tau,b_\tau\) count the columns
and the boxes at heights \(1<h<n_{j+1}\), respectively.
The height-\(k\) count is between \(\kappa_\tau\)
and \(\kappa_\tau+2\Delta\), and each subsequent
height has at most \(\Delta\) columns. Thus
\[
 k\kappa_\tau\le b_\tau\le k\kappa_\tau+x^{60},
 \qquad c_\tau\le\kappa_\tau+x^{60},
 \tag{CG1}\label{eq:local-column-gauge}
\]
using \(2k\Delta+n_m^2\Delta\le x^{60}\).

Let \(p<x_1\) be nonsplit in \(K\), \(I_p\) the positive
odd exponents if \(p\) is inert and all positive
exponents if \(p=d\). Put \(v=1\) in the inert case
and \(v=1/2\) when \(p=d\). Write
\(\tau_e=\langle-p^e\rangle\tau\) and
\(\gamma_p=v\sum_{e\in I_p}p^{-e}\).
For \(\mathcal E\subset\Phi\) erase the height-one columns at
\(\mathcal E\), setting \(u_\tau=a\mathbf 1_{\Phi_a\setminus \mathcal E}(\tau)\).
Define \(C_p(\mathcal E)\) as the sum of \(a/p^e\) over
\(\tau,\tau_e\in\Phi_a\setminus \mathcal E\), \(e\in I_p\).
The positive-exponent contribution after erasing is
\[
 D_{\tau,e}=(u_\tau/p^e-u_{\tau_e}-c_{\tau_e})_+
      +(b_\tau/p^e-b_{\tau_e})_+-c_\tau/p^e .
\]
For the full number \(B_j\) of boxes, put
\[
 v\sum_{\tau,\,e\in I_p}D_{\tau,e}+\gamma_p a|\Phi_a\cap \mathcal E|
              -fn_{j+1}^3(1+\log_p\max(B_j,1)).
\]
Denote this value by \(d_{pj}(\mathcal E)\) at inert primes. At \(p=d\)
subtract from the same expression the exponent-zero term
\[
 v\sum_\tau\bigl(\min(u_\tau,u_{\tau_0}+c_{\tau_0})
              +c_\tau+\min(b_\tau,b_{\tau_0})\bigr),
 \qquad \tau_0=\langle-1\rangle\tau .
\]
This specifies the uncharged jet exponent in the
ramified case. The inert exponent is derived by the
null-direction argument in Section~\ref{sec:exceptional-charged-jet}.

Set \(U=a/(pm)\). The saving contribution
\((a/p^e-\kappa_{\tau_e})_+\) for \(\tau\in\Phi_a\)
is at most
\[
 (u_\tau/p^e-u_{\tau_e}-c_{\tau_e})_+ +x^{60}
       +a p^{-e}\mathbf1_{\mathrm{coll}}
       +a p^{-e}\mathbf1_{\tau\in \mathcal E},
\]
where \(\mathrm{coll}\) means that both large parts
survive. When \(\tau\in \mathcal E\) and \(e=1\), the bound
still holds with \(a/x_1\) added to the left, since
\((a/p-\kappa_{\tau_e})_++a/x_1\le a/p\).
When the large part at \(\tau\) survives and its partner does not,
use \eqref{eq:local-column-gauge} and the Lipschitz bound for the positive part.
When both survive the collision term bounds the
left side. At \(\tau\notin\Phi_a\) there is no
saving from a long row. Also \eqref{eq:local-column-gauge} yields
\[
 k(\kappa_\tau/p^e-\kappa_{\tau_e})_+
       \le (b_\tau/p^e-b_{\tau_e})_+ +x^{60}.
\]
For \(p^e\ge x_1\) the saving terms vanish because
\(\kappa_{\tau_e}\ge a/x_1\), so no error is needed
in that range.

Summing contributes at most
\(2x^{60}(1+\log_2 x_1)\le U\) per place.
The negative column term contributes
\(\gamma_p c_\tau\le2(b+x^{60})/p\le U\).
At the ramified prime \(\tau,\tau_0\) cannot both
be \(a\)-places or both be \(b\)-places. Thus each
minimum in the exponent-zero subtraction is at most
\(a/(dm)\), using
\(b_\tau\le ka/x_1+x^{60}\le a/(dm)\) at a non-\(b\)-place.
The same bound applies to \(c_\tau\), so the subtraction
is at most \(2U\) per place. Finally
\(n_{j+1}^3(1+\log_p\max(B_j,1))\le U\).
To see that one cutoff suffices, use \(f,q_1\le x^{12}\),
\(n_{j+1}\le n_m\le x^6\) and \(B_j\le f(a+kb+x^{60})\).
Thus \(\log\max(B_j,1)\ll B+\log x\), with \(B=\log b\).
Uniformly in \(0\le j<m\), \(p\le x_1\) and the erased set,
\[
 \frac{x^{60}(1+\log_2x_1)+
        n_{j+1}^3(1+\log_p\max(B_j,1))}{U}
 \ll \frac{mx_1}{xb}\bigl[x^{60}(1+\log x_1)
                       +x^{18}(1+B+\log x)\bigr]=o(1).
\]
Indeed \(U=xb/(pm)\), \(x_1\le b^{1/10^6}\),
\(m\le2\sqrt x\), and \(B=x^\epsilon+O(1)\).
\begin{lemma}[Local derivative valuation]\label{lem:marked-column-valuation}
Fix \(p<x_1\) nonsplit in \(K\), with \(p\nmid q_\xi\), and \(0\le j<m\).
Let \(\mathcal E\subseteq\Phi\) index the erased height-one columns of the
\(j\)-th map. For fixed \(q_0,d,\epsilon\) and sufficiently large \(x\),
\begin{equation}
 d_{pj}(\mathcal E)\ge g_p/m-C_p(\mathcal E)+v|\Phi_a\cap \mathcal E|a/x_1
                              -faB^2/(mp).
 \label{eq:marked-column-gain}
\end{equation}
\end{lemma}

\begin{proof}
Multiply the positive-exponent comparisons above by \(v\) and
sum over \(\tau,e\). The erased terms \(va/p^e\) sum to
\(\gamma_p a|\Phi_a\cap\mathcal E|\); at \(e=1\) the comparison
retains \(v|\Phi_a\cap\mathcal E|a/x_1\).
Equation~\eqref{eq:local-column-gauge} and the four error bounds
above therefore give
\[
 \frac{g_p}{m}+v|\Phi_a\cap\mathcal E|\frac{a}{x_1}
 \le d_{pj}(\mathcal E)+v C_p(\mathcal E)
        +f(U+U+2U+U).
\]
Now \(v\le1\), \(U=a/(mp)\), and \(B^2\ge5\) for sufficiently
large \(x\).
\end{proof}

Call an inert prime \(p<x_1\) exceptional if \(p\not\equiv\pm1\bmod m\) and
\(o_0=\operatorname{ord}_{q_0}(-p)\) is odd with
\(o_0\le r_p=\lceil\log_p m\rceil\). For such a prime and a coefficient
embedding \(\iota\), put
\[
 \operatorname{ind}_\iota(\tau)
 =\{z\in\mathcal O_E:v_p((\iota\circ\Phi_\tau)(z))>0\}.
\]
Choose one member of each conjugate pair of primes induced by
\(\Phi_a\), and denote this set of representatives by \(\Lambda\).

At \(p=d\) there are no collisions. A partner belongs to the
same fibre modulo \(q_0\) only if \(d^e\equiv-1\bmod q_0\),
which is excluded by Lemma~\ref{lem:auxiliary-residue-class}.
At an inert prime \(p\nmid q_0\), so \(o_0\) is defined and a
collision requires an odd \(e\) divisible by \(o_0\).
For \(p\not\equiv\pm1\bmod m\) and a nonexceptional
prime, such an \(e\) is at least \(r_p+1\), where
\(r_p=\lceil\log_p m\rceil\). For \(p\equiv\pm1\bmod m\)
use \(e\ge1\). The series starting at either lower
bound \(e_0\) is at most \(2fa/p^{e_0}\).
Thus every nonexceptional nonsplit \(p<x_1\) satisfies
\(m C_p(\varnothing)\le L_p=faB^2c_p/p\),
where \(c_p\) is the multiplier in the local error bound.

Now suppose that \(p\) is exceptional. To construct \(\Lambda\),
let \(o\) be the order of \(-p\bmod q_1\).
It is even because \(-p\not\equiv1\bmod m\) and the
reduction modulo \(2^s=m^{22}\) has power-of-two order.
The exceptional condition \(p\not\equiv\pm1\bmod m\) gives
\(v_2(p\pm1)<\log_2m\). Since \(p^o=(-p)^o\equiv1\bmod 2^s\),
\[
 s\le v_2(p^o-1)
    =v_2(p-1)+v_2(p+1)+v_2(o)-1.
\]
The equality follows by factoring differences of squares
after removing the odd part of \(o\). Exactly one of
\(v_2(p-1),v_2(p+1)\) is one, so \(v_2(o)>\log_2m\)
and \(m\mid o\). In
\(\operatorname{Gal}(E/\mathbb Q)
=\operatorname{Gal}(K/\mathbb Q)\times(\mathbb Z/q_1)^\times\),
write \(c_K\) for conjugation on \(K\). The decomposition group
is generated by the arithmetic Frobenius \((c_K,p)\).
Complex conjugation on \(E\) is \(c=(c_K,-1)\). If the inducing
prime were fixed by \(c\), then \(c=(c_K,p)^j\) with odd \(j\).
This would give \((-p)^j=1\bmod q_1\), although its order \(o\)
is even. The two primes are therefore distinct. The identity
\[
 (1,(-p)^j)=(c_K,p)^j c^j
\]
shows that even steps preserve the inducing prime and odd steps
exchange the two primes. Conversely, put
\(w_\tau=\operatorname{ind}_\iota(\tau)\).
Then \(w_\sigma\) belongs to \(\{w_\tau,\bar w_\tau\}\) precisely
when \((1,\sigma\tau^{-1})\) lies in
\(\langle(c_K,p),(c_K,-1)\rangle\). In this subgroup the elements
fixing \(K\) have exponents of the same parity on its two
generators. Their second coordinates are exactly the powers of
\(-p\). Thus
\[
 w_\sigma\in\{w_\tau,\bar w_\tau\}
       \quad\Longleftrightarrow\quad
       \sigma\tau^{-1}\in\langle -p\rangle
                \subseteq(\mathbb Z/q_1)^\times .
\]
Multiplication by \((-p)^{o_0}\) preserves \(\Phi_a\),
so every conjugate pair induced by that set occurs
on both sides. Take one member of each pair to get
\(\Lambda\), with \(\mathcal E\) its inducing set.
Then \(C_p(\mathcal E)=0\). No \(b\)-place induces a prime in either \(\Lambda\) or \(\bar\Lambda\), because \(-1\) is outside the
odd-order subgroup generated by \(-p\bmod q_0\).
On the full set of embeddings inducing \(\lambda\) or \(\bar\lambda\),
the odd step \(\tau\mapsto\tau_1\) exchanges the two halves and
all the weights \(\kappa_\tau\) are \(a/x_1\).
In the column table, only the height-\(k\) multiplicity
\(w_1(\tau)-\lfloor(k+t+1)/2\rfloor\) depends on \(\tau\)
among the heights \(1<h<n_{j+1}\). Thus both the column count
\(c_\tau\) and the box count \(b_\tau\) at these heights agree
on the two halves under \(\tau\mapsto\tau_1\). On the
subsets in \(\Phi_a\), the step \(\tau\mapsto\tau_{o_0}\)
is the bijection exchanging the halves. Thus the sums of the
lift weights \(w_j(\tau)\), the determinant exponents, and the
height-one counts balance separately at each pair.

For \(\tau\in \mathcal E\) use \(c_\tau\le a/x_1+x^{60}\).
If \(p\le m\), then
\(r_p(a/x_1+x^{60})\le a/(pm)\) for large \(x\).
If \(p>m\), then \(r_p=1=o_0\) and \(-p\equiv1\bmod q_0\).
The orbit identity gives \(\mathcal E\subseteq\Phi_a\),
since \(\Phi_a\) is the fibre of \(1\bmod q_0\).
Hence \(c_\tau\le a/x_1+a/(pm)\) for \(\tau\in\mathcal E\).
Thus \(r_p\sum_{\tau\in \mathcal E}c_\tau\) is at most
\(|\mathcal E|a/(pm)+|\Phi_a\cap \mathcal E|a/x_1\).
Also \(|\Lambda|\nu_p\le L_p-fa/(pm)\) for
\(\nu_p=fr_p(n_m+r_p)\), with the sharper bound
\(|\Lambda|\nu_p\le L_p/m-fa/(pm)\) since
\(a/(pm)\) dominates every fixed power of \(x\).
This proves both column comparisons required on
the chosen set of primes \(\Lambda\).

At a nonexceptional small nonsplit prime, use the
uncharged jet bound \(p^{-d_{pj}(\varnothing)}\) from
Section~\ref{sec:exceptional-charged-jet} in the inert case and
Section~\ref{sec:marked-ramified-ordinary-jet} in the ramified case.
Every translate in the evaluated tensor is standard
at that prime, since \(p\nmid MNq_\xi\), and the gluing
is supported at \(N\). For an inert prime take the erased set empty
and all extra exponents zero in \eqref{eq:exceptional-charged-jet}. Inequality~\eqref{eq:marked-column-gain}
gives the full-chain lattice norm
\(p^{m C_p(\varnothing)+L_p-g_p}\le p^{2L_p-g_p}\).
The terminal branching matrix is integral and its
scalar has unit norm at \(p\ne N\).
At a split prime or \(p\ge x_1\), except \(N\),
the saving is zero. Apply the ordinary jet estimate with
\(\delta=\bar\pi_M\), \(c=\pi_M\) and
\(\lvert c/(\bar\delta d)\rvert_p=1\), since \(p\ne d\).
In each case we use the entire integral ramified
source, with denominator contribution zero. In particular every \(p\mid q_1\)
splits in \(K\) by Lemma~\ref{lem:auxiliary-residue-class}, including
\(2\) since \(-d\equiv1\pmod8\). At every such \(p\), the auxiliary
Gram scalar \(b\) in Section~\ref{sec:marked-integral-sections-ordinary-jets}
satisfies \(b\equiv\beta^j+\bar\beta^j\equiv1\pmod p\). Thus its
finite vacuum average has unit denominator. The two valuations of
the period factor computed in Section~\ref{sec:marked-integral-sections-ordinary-jets}
give \(p_F\overline{\tau(w_i)}\) integral at every coefficient embedding.
At \(2\), put \(n_2=v_2(q_1)\). Over the unramified part of the
completion, the polynomial \((1+X)^{2^{n_2-1}}+1\) is Eisenstein
and its derivative at \(\zeta_{2^{n_2}}-1\) has valuation \(n_2-1\).
Thus the two local Weil ideals have valuations
\(-(n_2-1)\) and \(0\). The corresponding period valuations
\(n_2-1\) and \(0\) cancel these possible denominators.
At every such embedding, \(c/(\bar\delta d)=d^{-1}\in\mathbb Z_p^\times\);
hence every ordinary derivative
\(\mathbf m![z^{\mathbf m}]\Theta_f(\pi_Ms_{\rm jet}z)\)
is integral in every degree. The identity matrix therefore clears the
entire upper tuple at \(p\mid q_1\) for every allowed twist prime to
\(p\), with \(v_p(\det I)=e_p=g_p=0\).

\subsection{Twisted derivatives at inert primes}
\label{sec:exceptional-charged-jet}

\sectionopening{In this subsection, we prove the derivative estimate at
each small inert prime. At the exceptional primes
Section~\ref{sec:exceptional-jet-hecke} applies it to Hecke operators.}

Fix an inert prime \(p<x_1\), a twist with \(p\nmid q_\xi\), and the
coefficient embedding \(\iota\) into \(\mathbf C_p\). The prime
divisors of \(MN\) are at least \(x_1\), and every prime divisor of
\(q_1\) splits in \(K\) by Section~\ref{sec:datum}. Thus
\(p\nmid MNq_1q_\xi\). Normalize \(v(p)=1\) and put
\(\mathcal O=\mathcal O_{\mathbf C_p}\). For positive odd \(e\)
write \(\tau_e=\langle-p^e\rangle\tau\) in the chosen CM type
\(\Phi\), and put \(n_j=k+j(t+1)\).
Take any \(\mathcal E\subseteq\Phi\) and nonnegative extra exponents
satisfying the pairing condition of
Section~\ref{sec:exceptional-jet-valuations}. For
\(\mathcal E=\varnothing\) take every extra exponent to be zero.

\indent{\it The ordinary-coordinate exponential.}
The one-variable calculation through \eqref{eq:cm-inert-leading-unit}
holds at every inert prime \(p\). Such a prime is odd since
\(2\) splits in \(K\). We use \(p<x_1\) and \(p\nmid q_\xi\)
when applying the calculation to the translated tensor.
Use the integral section identities and vacuum
resolution of Section~\ref{sec:marked-integral-sections-ordinary-jets}.
In their notation the quotient for \(\theta(pz)\) is
\(R(x)\), \(R(0)=1\), with integral coefficients and
\(2\deg R<p^2\), and \(g=-\log\theta\) satisfies
\(g''=g''(0)+x\), \(g(0)=g'(0)=0\).
The order bound on the reduced multiplication-by-\(p\)
series there also applies after the parameter change below.

Set \(q=p^2\). Average \(y(\zeta z)\) with weight \(\zeta^{-1}\)
and \(g(\zeta z)\) with weight \(\zeta^{-p-1}\) over
\(\zeta^{q-1}=1\), obtaining \(T_0\) and \(g^\flat\), respectively.
The CM pullback formulas extend coefficientwise to the Teichm\"uller
units, for which \(\bar\zeta=\zeta^p\). Applying the formula at
\(\zeta^{-1}\) and substituting \(\zeta z\) gives
\[
 \exp\bigl(\zeta^{-p-1}g(\zeta z)-g(z)\bigr)\in\mathcal O[[y]].
\]
Multiplying over \(\zeta\) and taking the unique constant-one
\((q-1)\)-st root gives \(\exp(g^\flat-g)\in\mathcal O[[y]]\).

Also \(T_0=y+O(y^2)\), so the two parameters give the same integral
series ring. The covariance of \(T_0\) gives
\[
 T_0(pz)=pT_0(z)+b_1T_0(z)^q+
              \sum_{j\ge2}b_jT_0(z)^{1+j(q-1)},\qquad b_1\in\mathcal O^\times .
\]
The order of the reduced parameter under \([p]\) is at most
the degree \(p^2\) of the isogeny. Thus \(b_1\) is a unit. The permitted degrees
above the vanishing linear term start at \(q=p^2\).
Choose \(\varpi^{q-1}=-p/b_1\) in \(\mathbb C_p\) and set
\(Y=T_0(\varpi z)/\varpi\). With \(v(p)=1\) one has
\(v(\varpi)=1/(q-1)\). Multiplying \(\varpi\) by
\(\zeta\in\mu_{q-1}\) leaves \(Y\) and
\(g^\flat(\varpi z)/\varpi^{p+1}\) unchanged by their character
weights \(1,p+1\).
The resulting equation is
\(Y(pz)/p=Y-Y^q+pE(Y)\), with \(E\in U^2\mathcal O[[U]]\).
In each degree \(d\ge2\), division by \(1-p^{d-1}\) determines an
integral coefficient from earlier coefficients. Modulo \(p\), the
equation becomes \(Y-Y^q=z\), whence
\[
 Y\equiv\sum_{u\ge0}z^{q^u}\pmod p .
\]

For clarity, the logarithm calculation giving the nonintegral
exponential is short. Since \(T_0=z+O(z^q)\), integral series in
\(y\) have integral coefficients in degrees below \(q\). For
\(3\le n<q=p^2\), the identity \(n(n-1)g_n=[z^{n-2}]x\)
and \(v_p(n(n-1))\le1\) give \(pg_n\in\mathcal O\).
In degree two the CM pullback gives \(a(\bar a-a)g_2\in\mathcal O\).
For \(a=\beta=(1+\sqrt{-d})/2\) and \(a=\bar\beta\), the difference
of the two multipliers is
\((\beta+\bar\beta)(\bar\beta-\beta)=-\sqrt{-d}\), a unit at the inert
prime \(p\). Hence \(g_2\in\mathcal O\).
Thus \(qg-g(pz)\) lies in \(p\mathcal O[[z]]\) through degree \(q-1\).
Its exponential does also after subtracting the constant, since
\(v(n!)\le(n-1)/(p-1)\) for \(n\ge1\). Comparing coefficients with
\(R(x(z))\), using \(x=z^2+\cdots\) and \(2\deg R<q\), gives
\(R-1\in p\mathcal O[X]\). Its formal logarithm shows
\(qg-g(pz)\in p\mathcal O[[T_0]]\).
Project this equation to degrees congruent to \(p+1\) modulo
\(q-1\). The same covariance in \(T_0\) gives
\[
 qg^\flat(z)-g^\flat(pz)
   =p\bigl(a_0 T_0^{p+1}+T_0^{p+q}H(T_0)\bigr),
 \qquad a_0\in\mathcal O,\quad H\in\mathcal O[[U]] .
\]

\subsubsection{Integral rescaling of the local exponential}
\label{sec:exceptional-rescaled-exponential}

\sectionopening{In this subsection, we normalize the exponential in the local coordinates of Section~\ref{sec:exceptional-charged-jet}.}

Set \(c_\varpi=\varpi^{p+1}/p\) and
\(\Xi=-g^\flat(\varpi z)/c_\varpi\). After rescaling, the leading
scalar \(-p\varpi^{p+1}/(q c_\varpi)\) is a unit, and the
higher terms have the additional factor \(\varpi^{q-1}\in p\mathcal O\).
The coefficient of \(\Xi\) in degree \(d\) is multiplied by
\(1-p^d/q\), a unit congruent to one for \(d\ge p+1\), while all
lower degrees are zero. Thus \(\Xi\) is integral and
\[
 \begin{aligned}
 \theta(\varpi z)&=W(z)\exp(c_\varpi\Xi(z)),
     & W&\in\mathcal O[[z]],\\
 \Xi&\equiv \lambda Y^{p+1}
  \equiv\sum_{e\ \mathrm{odd},\,u\ge0}
           \lambda(z^{p^e}z)^{p^u}\pmod p .
 \end{aligned}
\]
In fact the rescaled equation has the form
\[
 (1-p^{n-2})[z^n]\Xi=
      -a_0[z^n]Y^{p+1}+p R_n,\qquad R_n\in\mathcal O
       \quad(n\ge p+1).
\]
The degrees below \(p+1\) vanish by the Teichm\"uller projection.
Modulo \(p\mathcal O\), multiplying \(Y\) by \(Y^p\) pairs one
even Frobenius index with one odd index. Their minimum \(u\)
and positive odd difference \(e\) give degree \(p^u(p^e+1)\).
Conversely \(u=v_p(n)\) and \(p^e=n/p^u-1\) recover the two
indices. Thus \([z^n]\Xi\equiv\lambda\bmod p\mathcal O\) in
exactly these degrees once \(\lambda\) is shown to be a unit below;
in every other degree the coefficient belongs to \(p\mathcal O\).

We determine the leading coefficient for use in
Section~\ref{sec:twelve-coordinate-inert-expansion}. Put
\(h_j=[z^{p-1}]x(z)^j\). Since the projection defining
\(g^\flat\) retains degree \(p+1\) and \(g''=g''(0)+x\),
\begin{equation}\label{eq:cm-inert-leading-unit}
 [z^{p+1}]\Xi=-p[z^{p+1}]g=-\frac{h_1}{p+1}.
\end{equation}
The identities \((x')^2=4x^3+a_1^2x^2+4x\) and
\(x''=6x^2+a_1^2x+2\) follow from the normalized Weierstrass
equation. For \(j\ge0\) they give
\[
 (4j+6)h_{j+2}+(j+1)a_1^2h_{j+1}+(4j+2)h_j
   =[z^{p-1}](x^jx')'=p[z^p](x^jx')\in p\mathcal O .
\]
Indeed \(x\) and \(x'=2y+a_1x\) have integral coefficients in
degrees below \(p^2\) by the preceding parameter calculation.
Put \(k_p=(p-1)/2\). Then \(h_0=0\) and \(h_{k_p}=1\) since
\(x=z^2+O(z^3)\). For \(0\le j\le k_p-2\),
\(4j+6=2(2j+3)\) is a unit. If \(h_1\) vanished in the
residue field, the recurrence would force \(h_{k_p}\) to vanish.
For \(p=3\) we already have \(h_1=h_{k_p}=1\).
Thus \([z^{p+1}]\Xi\) is a unit for every inert \(p\), and
comparison with \(Y^{p+1}\) shows that \(\lambda\) is a unit.

For the theta series of standard level-one data on the matrix space,
the standard-vacuum resolution writes the series at normalized scale
\(a_d/\Omega_0\) as a product of
\(\theta(u_i^*z)^{m_i}\) times an integral series in \(y(u_i^*z)\),
with \(u_i\) integral above \(p\) and \(\sum_i m_i u_i u_i^*=1\).
The resolution uses the powers of \(\zeta_{q_1}\), whose Fourier
orthogonality gives the scalar basis matrix \(q_1 I\). Repeat each
vector \(q_1\) times and use the basis scale \(q_1^{-1}\).
Since \(p\nmid q_1\), the signed finite Fourier inversion from the
basis product to the integral vacuum has unit denominator at \(p\).
Write \(v_i\) for the repeated vectors, so that their Gram sum
is \(q_1^2 I\). Take the rational representatives of the finite
translations as in
Section~\ref{sec:marked-integral-sections-ordinary-jets}. They
are integral above \(p\), so one may take \(\delta=1\) there and
clear their torsion shifts by a common integer \(k_0\) prime to \(p\),
divisible by \(q_1\). The integral section formula then gives
\[
 u_i=v_i/(q_1k_0),\qquad m_i=k_0^2,\qquad
 \sum_i m_i u_i u_i^*=I
\]
in the embedding coordinates, separately for each matrix entry.
All terms have this same theta prefactor.

For each odd \(e\), the Frobenius congruence in the field
gives \(\overline{\tau(u)}^{p^e}\equiv\tau_e(u)\pmod p\)
in the coefficient embedding for \(u\) integral above \(p\).
On \(K\) the odd Frobenius cancels conjugation, and on the
cyclotomic part their composite acts by \(-p^e\). Since
\(p\nmid d q_1\), each completion of \(E\) above \(p\) is
unramified, with maximal ideal generated by \(p\). Thus the
residue-field identity gives the displayed congruence in
\(\mathcal O/p\mathcal O\). Reduction of the basis identity gives
\[
 \sum_i m_i(u_i^*z)^{p^e}(u_i^*z)
 \equiv \Gamma_e(z):=\sum_{\tau,b,a}
              z_{\tau,b,a}^{p^e}z_{\tau_e,b,a}\pmod p .
\]
The integers \(m_i\) satisfy \(m_i^{p^u}\equiv m_i\pmod p\).
Raising the preceding congruence to \(p^u\) therefore gives the
coupling \(\Gamma_e^{p^u}\) with coefficient \(\lambda\).
The difference of the two integral exponents is \(pR(z)\), with
\(R(0)=0\). For \(a\ge1\), Legendre's formula gives
\[
 v\!\left((pc_\varpi)^a/a!\right)
       =s_p(a)/(p-1),
\]
where \(s_p(a)\) is the sum of the base-\(p\) digits.
Thus \(\exp(pc_\varpi R)\) is an integral unit as a formal
series: only finitely many powers contribute in each degree.
In the linear forms of Section~\ref{sec:marked-integral-sections-ordinary-jets},
the factors \(a_d\) and \(\Omega_0\) cancel at the analytic scale
\(\pi_Ms_{\rm jet}\). The remaining scalar is
\(\alpha=\iota(c\,p_F)\), with \(c=\pi_M\).
All primes of \(q_1\) split in \(K\), so here \(p\nmid q_1\)
and \(p_F\) is a unit at \(\iota\). Also \(p\nmid M\) and
\(c\bar c=M\), so \(c\) and \(\alpha\) are units.
Homogeneity permits the choice
\(l_{eu}=\lambda\alpha^{(p^e+1)p^u}\) for every \(e,u\). Thus the actual Taylor series
\(G\), in the ordinary matrix coordinates of the arrow, satisfies
\begin{equation}
 G(\varpi z)=W_0(z)\exp\left(
   c_\varpi\sum_{e\ {\rm odd}\ge1,\,u\ge0}
                       l_{eu}\Gamma_e(z)^{p^u}\right),
 \qquad W_0\in\mathcal O[[z]],\quad l_{eu}\in\mathcal O .
 \label{eq:exceptional-jet-normal-form}
\end{equation}
For each tuple \(x\) in \eqref{eq:measured-seesaw-normalized-tensor},
let \(t(y)\) be the tuple of linked shifts indexed by the residues
\(y\) of Section~\ref{sec:seesaw-tensor-translation-N}. Write \(V_j(u)\)
for the original integral vacuum of the \(j\)-th arrow translated by \(u\).
That section gives the tensor identity
\begin{equation}\label{eq:linked-translate-tensor}
 \sum_x\Psi_\xi(x)v_x^{\rm gl}
   =\sum_{x,y}\Psi_\xi(x)
       \bigotimes_{j=0}^{m-1}V_j\bigl(x_j+t_j(y)\bigr).
\end{equation}
The supports of \(x\) and \(t(y)\) are disjoint, so their mutual
Heisenberg phase is one. Equal tuples on the right have coefficients
equal to sums of the original \(\Psi_\xi(x)\), integral above
\(p\nmid q_\xi\). Since \(p\nmid MNq_\xi\), the local factor of every
summand above \(p\) is the same vacuum. In particular
\eqref{eq:exceptional-jet-normal-form} applies to each arrow.
Applying the linear chain operator and terminal restriction to this
tensor identity gives equality on the full source space. Thus it
also applies to the individual summands of
\(Y_\Lambda^m u+p^{j_p}v\), \(u,v\in M_\iota\), in
Section~\ref{sec:exceptional-jet-hecke}. The identification with
\(L_\xi c_{\rm ss}^{-1}\Pi\) on the ramified subspace then follows
from \eqref{eq:measured-seesaw-normalized-tensor}.

\subsubsection{Valuations of the jets}
\label{sec:exceptional-jet-valuations}

\sectionopening{In this subsection, we bound the Schur coefficients of the expansion in Section~\ref{sec:exceptional-rescaled-exponential}.}

Let \(n=n_{j+1}\), \(\kappa_{\tau,h}\) be the number of columns of
height \(h<n\), and \(B_j\) the total number of their boxes.
Erase height-one columns at \(\mathcal E\), writing \(B_j'\) for the remaining
number of boxes. Let \(t_\tau,c_\tau,b_\tau\) count the remaining
height-one columns, the other columns, and the boxes of those other
columns, respectively. For each relevant \(\tau,e\) define
\[
 \ell_{\tau,e}=\min(t_\tau/p^e,t_{\tau_e}+c_{\tau_e})
     +c_\tau/p^e+\min(b_\tau/p^e,b_{\tau_e}).
\]
Consider directions \(y_{\tau,a,k}\in\mathcal O\) and a fixed
target basis \(I'\). Give the row index \(b\) of a matrix entry in a column of height \(h\) at \(\tau\)
an extra exponent \(w_{\tau,h,b}\in\mathbf Z_{\ge0}\). Require that one factor
has positive extra exponent when the paired entries have the same row index and their columns have heights \((1,h)\)
at \((\tau,\tau_e)\) if \(\tau\in \mathcal E\) or if \(h=1\) and
\(\tau_e\in \mathcal E\). Let \(C(I')\) sum the extra exponents over its boxes.

Use only heights with \(\kappa_{\tau,h}>0\) in the original
diagram, including its erased height-one columns. Variables at
absent heights have no derivatives in the minor product and may
be set to zero. For each occurring \((\tau,h)\) change the first
\(h\) directions by a matrix \(U_{\tau,h}\) over \(\mathcal O\).
For the positive odd \(e\) with \(p^e\le B_j\), the changes can make all
Frobenius pairings \(\sum_a\beta_a^{p^e}\beta_a'\)
divisible by \(p\) except the pairs with a positive extra exponent. The contribution in determinant
valuation, with each height weighted by its column count, is at
most \(\sum_{\tau,e:p^e\le B_j}\ell_{\tau,e}\). For fixed \((\tau,e)\)
use the following constraints.

\begin{enumerate}
\item For an unerased height one at the source, either multiply that direction
by a \(p^e\)-th root of \(p\), or impose at every nonexempt height
on the partner side one linear congruence against the original
source height-one direction. Choose the cheaper option.
\item For source height above one against partner height one,
impose the linear congruence
\(\sum_a\sqrt[p^e]{y_{\tau_e,a,0}}\,\beta_a=0\pmod{p^{1/p^e}}\)
on the source. Taking \(p^e\)-th powers is additive modulo \(p\).
Use the original partner direction even if its place is erased,
since this pair has source height above one.
\item For heights above one on both sides, make all directions on
the cheaper side divisible by \(p^{1/p^e}\) on the source or
by \(p\) on the partner.
\end{enumerate}

Each linear constraint contributes at most the valuation of its modulus.
Choose a coefficient of least valuation and clear the other components
by integral subtraction before multiplying the pivot direction by the
modulus. Collect the forms in the original directions and impose them
successively on the changed directions. All earlier constraints persist
under integral changes of the directions. A changed height-one direction is a
scalar multiple of its original direction, so the constraints give the
required pairings after all changes. The three contributions are the three
summands of \(\ell_{\tau,e}\).

In \eqref{eq:exceptional-jet-normal-form} substitute
\[
 z_{\tau,b,a}=\sum_{\substack{h<n,\,\kappa_{\tau,h}>0\\ k<h}}
          p^{w_{\tau,h,b}}\beta_{\tau,h,k,a}s_{\tau,h,k,b}.
\]
For fixed row \(b\), heights \(h,h'\) and directions \(k,k'\),
the corresponding pairing modulo \(p\) has coefficient
\[
 p^{p^e w_{\tau,h,b}+w_{\tau_e,h',b}}
   \sum_a\beta_{\tau,h,k,a}^{p^e}\beta_{\tau_e,h',k',a}.
\]
A positive extra exponent makes this zero; the change of directions
imposes the congruence for every other pair. Each term is divisible
by \(p\) even if its two variables coincide. For \(p^e>B_j\)
the coupling has degree above \(B_j\). In the quotient by monomials
of total degree greater than \(B_j\), the substituted exponent before
multiplication by \(c_\varpi\) therefore belongs to \(p\mathcal O[[s]]\).
The same factorial estimate for \((pc_\varpi)^a/a!\) proves
integrality of every coefficient through degree \(B_j\).
There are at most \(fn^3\) variables \(s\). For a multi-index
\(\alpha\) of total degree \(B_j\), the factorial estimate is
\[
 v(\alpha!)\ge \frac{B_j}{p-1}
       -fn^3\bigl(1+\log_p\max(B_j,1)\bigr).
\]
Apply the product of column determinants in the differential variables
\(\partial_s\) indexed by \(I'\). On a single determinant the chain
rule replaces its derivative matrix by the original matrix times
\(U_{\tau,h}^{\mathsf T}\), with row factors
\(\varpi p^{w_{\tau,h,b}}\). Thus the derivative is the original
minor-product jet \(J_y\) times
\(\varpi^{B_j}p^{C(I')}\prod_{\tau,h}\det(U_{\tau,h})^{\kappa_{\tau,h}}\).

Using the factorial estimate and the determinant contribution gives
\[
 \begin{aligned}
 v(J_y)+C(I')&\ge
 \gamma_p B_j-\sum_{\tau,e:p^e\le B_j}\ell_{\tau,e}
   -fn^3\bigl(1+\log_p\max(B_j,1)\bigr),\\
 \gamma_p&=\frac1{p-1}-\frac1{p^2-1}.
 \end{aligned}
\]
If \(B_j=0\) use the integral constant term. Otherwise enlarge the
nonnegative deduction sum to all positive odd \(e\). Since
\(\gamma_p=\sum_{e\ {\mathrm{odd}}\ge1}p^{-e}\), the resulting bound is
\[
 v(J_y)+C(I')\ge d_{pj}(\mathcal E),
\]
Indeed \(B'_j=\sum_\tau(t_\tau+b_\tau)\) and
\(B_j-B'_j=a|\Phi_a\cap\mathcal E|\). For the contribution
\(D_{\tau,e}\) of Section~\ref{sec:marked-column-gain},
\[
 \frac{t_\tau+b_\tau}{p^e}-\ell_{\tau,e}
 =(t_\tau/p^e-t_{\tau_e}-c_{\tau_e})_+
 +(b_\tau/p^e-b_{\tau_e})_+-c_\tau/p^e=D_{\tau,e}.
\]
Summation over odd \(e\ge1\) gives
\(\gamma_p B_j-\sum_{\tau,e}\ell_{\tau,e}
 =\gamma_pa|\Phi_a\cap\mathcal E|+\sum_{\tau,e}D_{\tau,e}\).
The series converge absolutely.

The same bound applies to each Schur jet entry \(J_{I',I}\).
Indeed the product of the flag minors of \(z_\tau y_\tau\)
expands as \(\sum_I A_I(y)S_{I',I}(z)\), where \(S\) is the
matrix of the Weyl map and the \(A_I\) are the source flag-minor
basis multiplied by \(\prod_\tau\det(Y_\tau)^{c_{\tau,j}}\).
Here \(Y_\tau\) is the square block of the first \(n_j\) direction
columns of the \(n_j\)-row matrix \(y_\tau\), and \(c_{\tau,j}\)
is the number of target columns of height \(n_j\). Each such target
column contributes one source top exterior power by the minor
identity of Section~\ref{sec:marked-integral-sections-ordinary-jets}.
For the diagram the only possible
heights other than one are \(k+h(t+1)\), so there are no heights
strictly between \(n_j\) and \(n_{j+1}\).
The \(A_I\) are independent after reduction.

Order the variables by the direction column and then by increasing
basis index, putting the earlier variable first in lexicographic
order. Every flag minor has leading coefficient one on its diagonal
monomial. In their product the variables of direction \(d\) give the
multiset of entries in tableau row \(d\); its weak order recovers that
row. Distinct tableaux therefore give distinct leading monomials in
every residue characteristic. The common determinant factor is nonzero
after reduction, so multiplication by it preserves independence.
An integral change of source basis preserves this independence.
The target index \(I'\) is kept in the standard tableau basis.

If all \(J_{I',I}\) vanish there is nothing to bound. Otherwise choose
\(I_0\) maximizing \(|J_{I',I_0}|_p\). The polynomial
\[
  J_{I',I_0}^{-1}\sum_I A_I(y)J_{I',I}
\]
has integral coefficients and nonzero reduction modulo the
maximal ideal by this independence.
The residue field of \(\mathbf C_p\) is infinite, so some integral
\(y\) gives it a unit value. Then
\( |J_y|_p=\max_I|J_{I',I}|_p\). The bound for every \(J_y\) gives
\begin{equation}
           v(J_{I',I})+C(I')\ge d_{pj}(\mathcal E) .
 \label{eq:exceptional-charged-jet}
\end{equation}

\subsubsection{The resulting Hecke operators}
\label{sec:exceptional-jet-hecke}

\sectionopening{In this subsection, we insert the valuations of Section~\ref{sec:exceptional-jet-valuations} into the integral Hecke action.}

Fix an exceptional \(p\) as in Section~\ref{sec:marked-column-gain} and
its coefficient embedding \(\iota\). The conjugate pairs above \(p\)
are split over \(F_0\), and \(p\) is a uniformizer on both sides.
For the construction below fix any allowed twist \(\xi\) with \(p\nmid q_\xi\).
Among the primes induced by \(\Phi_a\), choose one member of each
conjugate pair that occurs. The choice depends on \((p,\iota)\) and
the original datum. Write \(\Lambda\) for this set and put
\(\mathcal E=\{\tau:\operatorname{ind}_{\iota}(\tau)\in\Lambda\}\).
Odd partners induce conjugate primes, so
\(C_p(\mathcal E)=0\).

At \(\lambda\in\Lambda\) each allowed target neighbor has the form
\(k_1\operatorname{diag}(pI_i,I_{n_{j+1}-i})k_2\), with
\(k_1,k_2\in\mathrm{GL}_{n_{j+1}}(\mathcal O_\lambda)\) and \(0\le i\le r_p\).
Its lattice index is \(q_\lambda^i\). Write \(J_\lambda\) for the set
of \(i\) positions with elementary divisor \(p\). Form the product of
the chosen neighbors, written \(k_1Dk_2\) with \(D\) diagonal at
all places above \(p\), using the inverse conjugate transpose on the
opposite halves. Let \(K_-\) and \(K_+\) be the products of the
integral source and target groups above \(p\).

For arbitrary finite arguments \(g,h\), weak approximation on the
rational Cayley charts gives rational source and target isometries
\(\gamma,\eta\) with
\[
 (\gamma g)_p\in K_-,\qquad
 k'=(\eta h k_1D)_p\in K_+\cap D_p^{-1}K_+D_p .
\]
Both conditions specify nonempty open sets. At the transformed
arguments \((\gamma g,\eta h k_1Dk_2)\) both acting points are
integral above \(p\). The translated vacuum then has integral
displacements there. Away from \(p\) the vacuum expansion of
Section~\ref{sec:seesaw-integral-tensor} has \(p\)-unit indices
and character values. Thus the datum is still an integral
combination of the allowed vacuum translates, and the source lattice is standard. For the target,
\[
 (\eta h)_p=D_p^{-1}(D_p k' D_p^{-1})k_{1,p}^{-1}
                         \in D_p^{-1}K_+ .
\]
Thus its inverse acts on the transported value by \(D_p\) followed by
an integral invertible matrix. By rational covariance, the charged
row estimate for the transformed jet applies at the original
arguments \(g,h\). The weight embedding \(\iota\circ\bar\Phi_\tau\)
induces the place conjugate to \(\operatorname{ind}_\iota(\tau)\).
The diagonal entries above \(\lambda,\bar\lambda\) at a position in
\(J_\lambda\) are \(p,p^{-1}\), respectively. A box read at
the conjugate place therefore contributes \(\lvert p^{-1}\rvert_p=p\)
when \(\operatorname{ind}_\iota(\tau)=\lambda\), and
\(\lvert p\rvert_p=p^{-1}\) when \(\operatorname{ind}_\iota(\tau)=\bar\lambda\).
Thus its contribution to the base-\(p\) logarithm of the lattice norm multiplier is
\[
 \mathbf1_{b\in J_\lambda}
 \left(\mathbf1_{\operatorname{ind}_\iota(\tau)=\lambda}
       -\mathbf1_{\operatorname{ind}_\iota(\tau)=\bar\lambda}\right).
\] At each pair take
\[
 w_{\tau,h,b}=
 \begin{cases}
  \mathbf1_{h=1}\mathbf1_{b\notin J_\lambda},
       &\operatorname{ind}_\iota(\tau)=\lambda,\\
  \mathbf1_{b\in J_\lambda},
       &\operatorname{ind}_\iota(\tau)=\bar\lambda .
 \end{cases}
\]
Set \(w_{\tau,h,b}=0\) at every other embedding.
An odd partner has the same \(J_\lambda\) on the other side.
For an exempt pair the assigned exponent is one on the inducing
side if the entry is outside \(J_\lambda\), and one on the
conjugate side if it is in \(J_\lambda\). This verifies the
condition on the extra exponents for every neighbor. The determinant weights cancel,
since their sums on the two sets of embeddings inducing \(\lambda\)
and \(\bar\lambda\) agree. The sum of the box exponent and the extra exponent
is zero on the conjugate side, is one in each height-one box on
the inducing side, and contributes at most \(r_p\) per other column there
by strict increase down each column. If
\[
 A_\lambda=a\,|\{\tau\in\Phi_a:
                     \operatorname{ind}_\iota\tau=\lambda\}|,\qquad
 L_j(\mathcal E)=\sum_{\tau\in \mathcal E}\sum_{h\ne1,\,h<n_{j+1}}\kappa_{\tau,h},
\]
then \eqref{eq:exceptional-charged-jet} bounds the neighbor arrow
by the factor in the lattice norm
\[
           p^{\sum_{\lambda\in\Lambda}A_\lambda
                         +r_pL_j(\mathcal E)-d_{pj}(\mathcal E)}.
\]

Fix a uniformizer \(\varpi_\lambda\). Normalize Haar measure on \(\mathrm{GL}_d(E_\lambda)\) by
\(\operatorname{vol}(\mathrm{GL}_d(\mathcal O_\lambda))=1\).
For \(d=n_j\), let \(T_{j,\lambda,i}\) denote the double-coset operator
of \(\operatorname{diag}(\varpi_\lambda I_i,I_{d-i})\), \(0\le i\le d\).
Write \(\chi_{j,\lambda}\) for the local base change of the norm-one
character of the lift \(X_j\to X_{j+1}\), obtained by composition with
\(u\mapsto u/\bar u\). With \(q_\lambda=z_\lambda^2=p^{2h_\lambda}\)
and \(2h_\lambda\le f\), put
\[
 T_{j,\lambda}(U)=\sum_i(-U)^i q_\lambda^{\binom i2}T_{j,\lambda,i}.
\]
On source forms spherical at \(\lambda\), the identity is
\[
 \begin{aligned}
 T_{j+1,\lambda}(U)\Theta_{x,j}
  &= R_{j,\lambda}(U)\Theta_{x,j}
             T_{j,\lambda}(z_\lambda^{t+1}U),\\
 R_{j,\lambda}(U)&=\prod_{k=0}^t
     (1-q_\lambda^k z_\lambda^{n_j}\chi_{j,\lambda}(\varpi_\lambda)U).
 \end{aligned}
\]
In the split type-II correspondence for irreducible representations, the retained
Langlands factors are the contragredients of the smaller factors
\cite[Corollary~6.3(2)]{MinguezHowe}. We verify the operator identity
on the local Schwartz function, using the matrix calculation of
Section~\ref{sec:main-trace}. The pair \(\lambda,\bar\lambda\) avoids
the original level and the twist. The translating coordinates of the
weak-chain tensor are integral there by
Section~\ref{sec:seesaw-integral-tensor}, and the rank twists are
unramified by Section~\ref{sec:splitting-level-characters}. Thus the
local factor in self-dual coordinates is the indicator of
\(M_{n',n}(\mathcal O_\lambda)\), where \(n=n_j\) and \(n'=n+t+1\).
Put \(q=q_\lambda\), \(z=z_\lambda\), and \(c=\chi_{j,\lambda}(\varpi_\lambda)\).
For the additive matrix measure the map \(Z\mapsto h^{-1}Zg\)
has modulus \(J=|\det h|_\lambda^{-n}|\det g|_\lambda^{n'}\), by
its action on the \(n\) columns and \(n'\) rows. The unitary
Schr\"odinger action on functions uses \(J^{1/2}\) because the
substitution in the integral of the squared norm contributes \(J^{-1}\).
Thus the target coset contributes \(q^{ni/2}\), and the inverse source
coset contributes \(q^{n'i/2}\). Both ranks are even, since \(k=r+1-m\)
and \(t+1=5r+1\) are even. Write \(\nu\) for the character of this arrow, so
\(\nu^{\rm bc}(a)=\nu(a/\bar a)\). In either rank \(d=n,n'\),
\[
 (\nu\upsilon_1^{d/2})^{\rm bc}\upsilon^{-d}=\nu^{\rm bc},
 \qquad \upsilon_1^{\rm bc}(a)=\upsilon(a/\bar a)=\upsilon(a)^2
\]
by Section~\ref{sec:splitting-level-characters}. The two character factors
are consequently both \(c^i\). Indeed the local action, with source matrix
\(g\) and target matrix \(h\), is
\[
 (\omega(g,h)\phi)(Z)=
 \chi_{j,\lambda}(\det h/\det g)
 |\det h|_\lambda^{-n/2}|\det g|_\lambda^{n'/2}\,
 \phi(h^{-1}Zg).
\]
For the positive target coset and the inverse source coset this gives
\(\alpha^i=(z^nc)^i\) and \(\beta^i=(z^{n'}c)^i\), respectively.
In \(\lambda\)-coordinates the positive source coset at
\(\bar\lambda\) is the inverse coset at \(\lambda\), by inverse transpose.

At an integral matrix \(Z\), let \(s\) be its rank modulo
\(\varpi_\lambda\). The target cosets contributing to the indicator
are indexed by the \(i\)-planes in \(\ker\bar Z^{\mathsf t}\),
and the inverse source cosets by the \(i\)-planes in \(\ker\bar Z\).
Put \(K_d=\mathrm{GL}_d(\mathcal O_\lambda)\), \(a_i=\operatorname{diag}(\varpi_\lambda I_i,I_{d-i})\)
and \(D_i^{(d)}=K_da_iK_d\). Choose \(\gamma_V\in D_i^{(d)}\) with
\(\ker\bar\gamma_V=V\) for each \(i\)-plane \(V\) of the residue space.
The coset \(K_d\gamma_V\) is determined by \(\gamma_V^{-1}\mathcal O_\lambda^d\), whose image
under multiplication by \(\varpi_\lambda\) and reduction is \(V\).
Thus \(D_i^{(d)}=\bigsqcup_V\gamma_V^{\mathsf t}K_d\) and
\((D_i^{(d)})^{-1}=\bigsqcup_V\gamma_V^{-1}K_d\), with each coset of measure one.
If \(Z\gamma_V^{-1}\) is integral, so is \(Z\), since \(\gamma_V\)
is integral. At integral \(Z\), the further condition is
\(V\subseteq\ker\bar Z\). Transposition gives the target criterion
for \(\gamma_V^{-\mathsf t}Z\), including vanishing outside the
integral matrix box.
The finite \(q\)-binomial theorem thus gives the polynomials
\(P_{n'-s}(\alpha U)\) and \(P_{n-s}(\beta U)\), where
\(P_b(V)=\prod_{k=0}^{b-1}(1-q^k V)\).
Since \(q^{t+1}\alpha/\beta=z^{t+1}\),
\[
 P_{n'-s}(\alpha U)=P_{t+1}(\alpha U)\,
                         P_{n-s}(\beta z^{t+1}U).
\]
For the theta integral write \(B(\phi,v)\) for the lift of a source
form \(v\) with datum \(\phi\), and \(\rho_d\) for right translation.
Translation in the two variables gives
\[
 B(\omega(g,1)\phi,v)=B(\phi,\rho_n(g^{-1})v),\qquad
 B(\omega(1,h)\phi,v)=\rho_{n'}(h)B(\phi,v).
\]
The target is spherical since \(K_{n'}\) fixes \(\phi\). Inversion carries
Haar measure on \((D_i^{(n)})^{-1}\) to Haar measure on \(D_i^{(n)}\),
since \(\mathrm{GL}_n(E_\lambda)\) is unimodular. Integrating the first
identity therefore gives the positive source operator at \(\lambda\) on \(v\).
Bilinearity preserves the scalar coefficients in the kernel identity and
gives the asserted formula for \(\Theta_{x,j}\).
The unramified character value \(c\) is a unit at \(\iota\). In fact, let
\(a_\lambda\) be the idele supported at \(\lambda\) with component
\(\varpi_\lambda\). Evaluation of the weight character on
\(a_\lambda/\bar a_\lambda\) gives
\[
 |\iota(c)|_p
 =p^{\sum_{\operatorname{ind}_\iota(\tau)=\bar\lambda}w(\tau)
      -\sum_{\operatorname{ind}_\iota(\tau)=\lambda}w(\tau)}=1 ,
\]
where \(w\) is the weight of this arrow character. To obtain the
first equality, raise the idele ideal to a principal power and use the
infinity character on a generator as in
Section~\ref{sec:splitting-level-characters}. Finite unit characters
have absolute value one. The local uniformizer has \(v_p=1\) since
\(E_\lambda/\mathbf Q_p\) is unramified. The twist has weight zero,
and the two weight sums agree by Section~\ref{sec:marked-column-gain}. For the fixed twist, \(R_{j,\lambda}\) is the same for every translate
in the support of the chain datum. Its coefficients are \(\iota\)-integral
and its constant term is one. For the original initial segment
\(C=C_{x,j+1}=\Theta_{x,j}\cdots\Theta_{x,0}\), iteration gives
\[
 \left(\prod_{a=0}^{j}R_{a,\lambda}
       (z_\lambda^{(j-a)(t+1)}U)\right)^{-1}
       T_{j+1,\lambda}(U)C
       =C T_{0,\lambda}(z_\lambda^{(j+1)(t+1)}U).
\]
The substitution in each factor is by an integer multiple
of \(U\), so the product still has integral coefficients and
constant term one. Its formal reciprocal is integral by the
coefficient recurrence. Since \(r_p\le m\le k=n_0\),
the coefficient of \(U^{r_p}\) on the right is
\((-1)^{r_p}p^{u_j}C T_{0,\lambda,r_p}\), with
\[
 u_j=h_\lambda r_p\bigl(r_p-1+(j+1)(t+1)\bigr)
      \le f r_p(n_m+r_p)=\nu_p .
\]
Here \(2h_\lambda\le f\) and \(j<m\). Multiplying the coefficient
identity by \((-1)^{r_p}p^{\nu_p-u_j}\) gives integral coefficients
\(c_{j,\lambda,i}\) satisfying
\[
 \sum_{i\le r_p}c_{j,\lambda,i}
       T_{j+1,\lambda,i}C_{x,j+1}
    =p^{\nu_p} C_{x,j+1}T_{0,\lambda,r_p},
 \qquad \nu_p=fr_p(n_m+r_p).
\]
Put \(H_\lambda=p^{A_\lambda}T_{0,\lambda,r_p}\) and
\(Y_\Lambda=\prod_{\lambda\in\Lambda}H_\lambda\).
Move the outgoing Hecke sums into the data by Weil equivariance,
using the operators \(p^{A_\lambda-\nu_p}\sum_i
c_{j,\lambda,i}T_{j+1,\lambda,i}\) on the target side of each datum.
At each distinguished pair every supported datum has the local vacuum
factor. A double-coset operator preserves its \(K_{n_{j+1}}\)-invariance
and commutes with the source action and with the level actions at the
other places. The resulting sums preserve the original level equivariance.
The operations at different pairs commute. The coefficients
\(c_{j,\lambda,i}\) may depend on \(\xi\) but are common to all the
supported translates for that twist. The source operator \(Y_\Lambda\)
uses only the fixed arithmetic datum and \((p,\iota,\Lambda)\). The resulting arrow applied
to \(C_{x,j}\) is \(C_{x,j+1}Y_\Lambda\), so the modified full
chain gives exactly the original chain on \(Y_\Lambda^m\).
The factors \(p^{A_\lambda-\nu_p}\) cancel the preceding
\(\sum A_\lambda\) in the exponent in the lattice norm. Taking the maximum over
the integral Hecke coefficients gives the per-arrow exponent
\[
              |\Lambda|\nu_p+r_pL_j(\mathcal E)-d_{pj}(\mathcal E).
\]

The column comparison in
Section~\ref{sec:marked-column-gain} gives the following
two bounds with \(L_p=faB^2/p\) and \(e_p=2L_p\). At these
primes \(p\not\equiv\pm1\bmod m\), so \(c_p=1\) in the definition of \(L_p\).
\[
 \begin{aligned}
 d_{pj}(\mathcal E)&\ge g_p/m+|\Phi_a\cap \mathcal E|a/x_1-L_p/m,\\
 |\Lambda|\nu_p+r_pL_j(\mathcal E)&\le L_p/m+|\Phi_a\cap \mathcal E|a/x_1 .
 \end{aligned}
\]
The first uses zero collisions on the unerased large parts. Together
they give an exponent at most \((2L_p-g_p)/m\). The composite
therefore has lattice norm at most \(p^{2L_p-g_p}\) on
\(Y_\Lambda^m M_\iota\), where \(M_\iota\) is the original
integral source. With no erased columns, the same comparison
gives \(p^{m C_p+L_p-g_p}\) on \(M_\iota\), where
\(C_p\le2fa/p\) by the sum of \(a/p^e\) over the colliding
large parts. Thus on \(p^{j_p}M_\iota\), with
\(j_p=\lceil2mfa/p\rceil\), the exponent is at most \(L_p-g_p\).
Ultrametricity gives \(p^{2L_p-g_p}\) on the sum. The terminal
restriction uses the original integral Weyl matrix and the original
scalar, which is a unit here because \(p\ne N\).

Proposition~\ref{prop:source-hecke-colength} gives
\(H_\lambda M_\iota\subseteq M_\iota\) for the scaled operator
\(H_\lambda=p^{A_\lambda}T_{0,\lambda,r_p}\) and
\[
 \operatorname{length}_p\!\left(M_\iota/
        (H_\lambda M_\iota+p^{j_p}M_\iota)\right)
 \le d_X\left(\frac{12A_\lambda}{mp}+\frac{j_p}{m^2f}\right).
\]
Compute normalized lengths over the integers \(\mathcal O_F\)
of a finite extension \(F/\mathbb Q_p\) containing the entries of
all \(H_\lambda\) and of the operators \(R_H(z_w)\), \(w\mid\pi_N\),
defining the ramified source, in the original class and Weyl bases.
These entries use finite sums of algebraic characters and matrix entries in
the local coefficient fields. Divide each length by \(e(F/\mathbb Q_p)\).
Smith normal form identifies it with the sum of the invariant-factor
valuations for \(v(p)=1\), which is unchanged after scalar extension.
On the resulting finite quotient \(R=M_{\iota,F}/p^{j_p}M_{\iota,F}\)
the cokernel lengths of two endomorphisms bound the length for their
composition. For endomorphisms \(A,B\), the exact
sequence
\[
 0\longrightarrow AR/ABR\longrightarrow R/ABR
                  \longrightarrow R/AR\longrightarrow0
\]
has first term a quotient of \(R/BR\), via \(A\).
After normalization the length inequality is unchanged.
Apply it successively to all factors in the \(m\) copies
of the product over \(\Lambda\). The normalized length of
\(M_\iota/(Y_\Lambda^m M_\iota+p^{j_p}M_\iota)\) is at most
\[
 d_Xm\sum_{\lambda\in\Lambda}
    \left(\frac{12A_\lambda}{mp}+\frac{j_p}{m^2f}\right)
 \le d_X(14fa/p+1/m)\le d_X L_p .
\]
Here \(\sum A_\lambda\le af\), \(|\Lambda|\le f\), and the last
inequality is the parameter inequality
\(14fa+p/m\le faB^2\).

Put \(J=Y_\Lambda^mM_\iota+p^{j_p}M_\iota\).
Over this same field \(F\), by Section~\ref{sec:ramified-source},
the characteristic-zero source \(S_F\) is the common kernel of
\[
 \sum_{j=0}^{N-2}R_H(z_w)^j=(N-1)A_w
                         \qquad(w\mid\pi_N).
\]
Each \(z_w\) has identity components above \(p\), since \(p\ne N\).
Thus \(R_H(z_w)\) and its inverse preserve the transported integral
source lattice, and the displayed sums have integral matrices. Their
common kernel on \(M_{\iota,F}\) is \(M_{\iota,F}\cap S_F\), even
when \(p\mid N-1\), since the kernel is taken before reduction.
Put \(J_F=Y_\Lambda^mM_{\iota,F}+p^{j_p}M_{\iota,F}\).
The inclusion
\[
 (M_{\iota,F}\cap S_F)/(J_F\cap S_F)
                   \hookrightarrow M_{\iota,F}/J_F
\]
bounds the length of its source. The valuation ring \(R_\iota\) is flat
over \(\mathcal O_F\), since it is torsion-free over a DVR.
Scalar extension therefore preserves these kernels and intersections.
The normalized colength on the original integral ramified source is
bounded by the same \(d_X L_p\), also when \(p\mid N-1\).
The tuple clearing is given in Section~\ref{sec:short-images-proved}. The resulting module is fixed after the chosen arithmetic datum
and \((p,\iota)\) and before any twist or translate. For every allowed \(\xi\) with
\(p\nmid q_\xi\), the local datum at the pair is spherical and
the preceding estimates apply on this same module to all evaluated
translates. The twist can change the unramified character values
in the transfer identity; their weight-zero factors are roots of
unity, hence \(p\)-adic units even if their orders are divisible by
\(p\). The operators \(H_\lambda\) defining the module act on the
original source independently of these values.
The set of representatives and the module may depend on the
coefficient embedding. The loss and denominator bounds depend on \(p\) alone.

\section{Hecke operators and colength at exceptional primes}
\label{sec:exceptional-hecke-analysis}

\sectionopening{In this section, we construct the source Hecke operators that supply
the divisibilities in Section~\ref{sec:exceptional-charged-jet}.
The local choices of Section~\ref{sec:marked-column-gain} ensure that
the colength bounds hold on the same source lattice for all permitted twists.}

Keep the exceptional prime \(p\), its coefficient embedding
and the set of representatives \(\Lambda\) of Section~\ref{sec:exceptional-charged-jet}.
The two members above each chosen place of the real subfield are
distinct by the decomposition-group calculation in
Section~\ref{sec:marked-column-gain}. The exceptional prime \(p\) is
inert in \(K\), whereas \(d\) ramifies and every divisor of \(q_1\)
splits in \(K\). Thus \(p\nmid d q_1\), and
\(E=K(\zeta_{q_1})\) is unramified above \(p\). Since \(p<x_1\)
and every prime divisor of \(MN\) is at least \(x_1\), the local
factor of \(J_H\) at the pair of each \(\lambda\in\Lambda\) is
the hyperspecial group \(\mathrm{GL}_k(\mathcal O_{E_\lambda})\). Let \(X_0\) be the rank-\(k\) Hermitian
space and \(X\) the full complex space of algebraic automorphic
forms at its original level \(J_H\) of Section~\ref{sec:ramified-source}.
Set \(R_\iota=\{z\in\overline{\mathbf Q}_p:v_p(z)\ge0\}\) with
\(v_p(p)=1\). In the given coefficient embedding let \(M=M_\iota\)
be the \(R_\iota\)-lattice of integral forms. We may calculate the
matrices over \(\mathcal O=\mathcal O_{\mathbf C_p}\). Their maximal minors
below determine the colength over \(R_\iota\).
Equivalently, fix a finite extension \(F/\mathbf Q_p\) containing the
images of \(E\) at the coefficient embeddings in the source, the
entries of the finitely many original class transitions, and
\(w=\sqrt p\). The matrices of all Levi elements on the active
factors have entries in \(\mathcal O_F\), since the active
embeddings read the same unramified completion. For every integer
\(l\ge1\), the tests by these matrices modulo \(w^l\) depend only
on the Levi elements modulo \(p^{\lceil l/2\rceil}\). Hence each
intersected kernel at a fixed depth is an intersection of finitely
many kernels. The occurrence-matrix operators in
Section~\ref{sec:exceptional-single-average} and their images of the
specified tensor basis vectors define the trace functional over this
same DVR. Flat scalar extension preserves the kernels. Smith form with
\(v(p)=1\) preserves the valuations of the elementary divisors.
Put \(d_X=\dim_{\mathbf C}X\), \(r=r_p=\lceil\log_p m\rceil\), and
\(A_\lambda=a\,\#\{\tau\in\Phi_a:\operatorname{ind}_\iota(\tau)=\lambda\}\).
In the unnormalized double-coset convention the operator of
Section~\ref{sec:exceptional-jet-hecke} is
\(Y_\lambda=H_\lambda=p^{A_\lambda}T_{0,\lambda,r}\).

\indent{\it The Frobenius orbit of embeddings.}
Fix \(\lambda\in\Lambda\) and choose \(\tau\in\Phi_a\) inducing it. Put \(u=\langle-p\rangle\),
\(e=\operatorname{ord}_{q_0}(-p)\), and
write \(N_{\rm orb}=\operatorname{ord}(u^e)\) for the orbit length in this calculation. The infinite places inducing
\(\lambda\) or \(\bar\lambda\) are the \(u\)-orbit of
\(\tau\), by decomposition at the inert prime of the
imaginary quadratic field. The defining \(q_0\)-condition
of \(\Phi_a\) is preserved by \(u^j\) exactly when
\(e\mid j\). Thus
\[
 t(s)=u^{es}\tau,\qquad s\in\mathbb Z/N_{\rm orb}\mathbb Z
\]
lists precisely the original \(\Phi_a\)-places above
\(\lambda,\bar\lambda\), each once. The integer \(e\) is odd,
and the step \(u^e=\langle-p^e\rangle\) exchanges the two
primes. With \(t(0)=\tau\), the even places induce \(\lambda\).
Hence \(N_{\rm orb}\) is even and \(A_\lambda=aN_{\rm orb}/2\).

Write \(W=(\mathcal O_E)_\lambda=E\cap\mathcal O_{E_\lambda}\)
for the localization and \(\widehat W=\mathcal O_{E_\lambda}\)
for its completion. As \(p\) is unramified, \(pW=\lambda W\);
the common residue field is \(W/pW=\widehat W/p\widehat W\).
Put \(P=p^e\). At an even embedding use the embedding of \(t(s)\),
and at an odd embedding compose that embedding with complex
conjugation on the field. These embeddings \(\iota_s\)
all induce \(\lambda\). On the field the step
\(\langle-p^e\rangle\) followed by conjugation is the
\(e\)-th power of arithmetic Frobenius at \(p\).
Indeed in \(E=K\mathbb Q(\zeta_{q_1})\) the step fixes
\(K\) and takes \(\zeta_{q_1}\) to \(\zeta_{q_1}^{-p^e}\).

Conjugation changes the cyclotomic sign and acts
nontrivially on \(K\), where \(p\) is inert and \(e\)
is odd. More explicitly, with
\(\sigma_p=(c_K,p)\) and \(c=(c_K,-1)\), the alternating
embeddings are
\(\iota_s=\iota_0\circ u^{es}c^s=\iota_0\circ\sigma_p^{es}\).
Arithmetic Frobenius satisfies the residue congruence modulo
\(\lambda W=pW\). Elements of \(\mathcal O_E\setminus\lambda\)
have unit images under every \(\iota_s\), so the congruence
extends to \(W\) and holds modulo \(p\mathcal O\).
Consequently
\[
 \iota_{s+1}(x)\equiv\iota_s(x)^P\pmod{p\mathcal O},
 \qquad q=|W/\lambda W|=p^{eN_{\rm orb}}=P^{N_{\rm orb}} .
\]
Indeed arithmetic Frobenius is \(u\) times conjugation,
and conjugation has order two and is independent of
\(\operatorname{Gal}(E/K)\). The residue degree is
\(\operatorname{lcm}(\operatorname{ord}(u),2)=eN_{\rm orb}\), since
\(N_{\rm orb}\) is even.

To check the exceptional order condition,
write \(m=2^H\). The order \(eN_{\rm orb}\) satisfies
\(u^{eN_{\rm orb}}\equiv1\bmod 2^{22H+v_2(q_0)}\).
It is even, because an odd exponent would give
\(v_2(u^{eN_{\rm orb}}-1)=v_2(u-1)<H\).
For an even exponent the factorization of the difference
of powers gives
\[
 v_2(u^{eN_{\rm orb}}-1)=v_2(u-1)+v_2(u+1)+v_2(eN_{\rm orb})-1.
\]
One of the first two terms is \(1\), and both are less
than \(H\) since \(p\not\equiv\pm1\bmod m\).
It follows that \(v_2(eN_{\rm orb})\ge H\), so \(m\mid eN_{\rm orb}\). This construction keeps \(\lambda\) and its parameters
fixed.

\indent{\it Coordinates for the exact operator.}
Represent each class \(x\) at the original level by a point
\(a_x\) integral above \(p\). For any representative \(c_x\)
prescribe \(\gamma_p\in K_{X_0,p}(c_x)_p^{-1}\), with
\(K_{X_0,p}\) the product of the integral groups above \(p\).
Since \(p\nmid M\), one may simultaneously prescribe the old
principal subgroups at the primes dividing \(M\).
Weak approximation on the rational Cayley chart of
Section~\ref{sec:ramified-source-clearing} gives \(\gamma\)
in these open sets. Then \(a_x=\gamma c_x\) represents the
same class.
The original level has trivial class stabilizers. Evaluation
at these representatives therefore identifies \(M\) with
the module of all integral coordinate vectors at the classes.
Indeed left rational translation in the lattice norm cancels the
rational transformation of the value, and right level
translation preserves the integral lattice norm.

Since \(p\) generates the maximal ideal of \(\widehat W\),
the uniformizers \(\varpi_\lambda\) and \(p\) give the same
hyperspecial double coset and the same unramified character value.
Put \(d=\operatorname{diag}(pI_r,I_{k-r})\) at \(\lambda\).
Put \(K_\lambda=\mathrm{GL}_k(\widehat W)\),
\(\widehat{\mathcal P}=K_\lambda\cap dK_\lambda d^{-1}\), and
\[
 \mathcal P=\widehat{\mathcal P}\cap\mathrm{GL}_k(W)
   =\{g\in\mathrm{GL}_k(W):d^{-1}gd\in\mathrm{GL}_k(W)\}.
\]
The incoming edges at \(x\) are indexed by
\(K_\lambda/\widehat{\mathcal P}\). Reduction identifies this
quotient with the \((k-r)\)-planes over the residue field. It
gives the same flag space for \(\mathrm{GL}_k(W)/\mathcal P\):
the upper right block vanishes modulo \(p\), and every invertible
residue matrix lifts to a matrix over \(W\) with unit determinant.
Choose representatives \(\beta_b\in\mathrm{GL}_k(W)\) in this way.
With \(\operatorname{vol}(K_\lambda)=1\),
\[
 K_\lambda d K_\lambda=\bigsqcup_b\beta_b dK_\lambda,\qquad
 \operatorname{vol}(\beta_b dK_\lambda)=1,\qquad
 \operatorname{vol}(\widehat{\mathcal P})
          =[K_\lambda:\widehat{\mathcal P}]^{-1}.
\]
Thus every edge has coefficient one in the unnormalized Hecke sum.
At \(\lambda\) replace the component of \(a_x\) by \(\beta_b d\)
and make the dual change
at \(\bar\lambda\). This gives the same Hecke coset as
right multiplication by \(a_{x,\lambda}^{-1}\beta_b d\).
Write the resulting point as \(\gamma_e a_{y(e)}v_e\),
with \(\gamma_e\) rational unitary and \(v_e\) of original
level. At \(\lambda\) this identity gives
\[
        \gamma_e^{-1}=g_e d^{-1}\beta_b^{-1},
        \qquad g_e\in\mathrm{GL}_k(W).
\]
At \(\lambda\), \(g_e=\gamma_e^{-1}\beta_b d\) has entries
in \(E\), and equals \(a_{y(e),\lambda}v_{e,\lambda}\) in
\(\mathrm{GL}_k(\widehat W)\). Thus \(g_e\in\mathrm{GL}_k(W)\).
At each \(w\mid p\) outside \(\{\lambda,\bar\lambda\}\), the
modified representative still has component \(a_{x,w}\), so
\[
 \gamma_{e,w}^{-1}=a_{y(e),w}\,v_{e,w}\,a_{x,w}^{-1}
        \in\mathrm{GL}_k(\mathcal O_{E_w}).
\]
Here both class representatives are integral above \(p\) and the
original level is hyperspecial there.

For use on the dual edges, the map
\(e\mapsto(y(e),(g_e^{\mathsf T})^{-1}\mathcal P)\)
is a bijection. If two edges have the same image,
\(g_e^{\mathsf T}(g_{e'}^{\mathsf T})^{-1}\in\mathcal P\).
Taking inverse transpose in the conjugation test defining
\(\mathcal P\) shows that
\(d g_e^{-1}g_{e'}d^{-1}\) is integral with integral inverse.
Let \(x,x'\) be the two termini and put \(\delta=\gamma_e\gamma_{e'}^{-1}\).
At \(\lambda\), \(\delta=\beta_b(dg_e^{-1}g_{e'}d^{-1})\beta_{b'}^{-1}\)
is integral with integral inverse. The dual component has the same
property. Write \(J\) for the original level. Away from the pair,
the two factorizations at their common origin give
\(a_x^{-1}\delta a_{x'}=v_e^{-1}v_{e'}\) locally.
At the pair the representatives are integral and \(J\) is hyperspecial.
Hence \(a_x^{-1}\delta a_{x'}\in J\), so \(x=x'\). Now
\(\delta\) is in the trivial rational stabilizer of \(a_xJ\), so \(\gamma_e=\gamma_{e'}\).
The displayed equation for \(\gamma_e^{-1}\) then gives
\(\beta_b^{-1}\beta_{b'}\in\mathcal P\). This proves
injectivity, and equal finite cardinalities give bijectivity.

At an embedding \(\tau\) the one-row weight has degree
\(a_\tau=a\mathbf1_{\Phi_a}(\tau)\) and determinant exponent zero.
Every column has height one, so its evaluation vector in the
integral tableau coordinates is \(e(z)=(z^\alpha)_{|\alpha|=a_\tau}\).
For the symmetric-power matrix
\(C(A)_{\beta,\alpha}=[z^\beta](A^{\mathsf T}z)^\alpha\),
coefficient comparison gives
\[
       e(Az)=C(A^{\mathsf T})^{\mathsf T}e(z).
\]
This is an identity over \(\mathbb Z[A]\), also for \(a_\tau=0\).
Use \(\rho_s(g)\)
for the symmetric power of \(\iota_s(g)\) at even \(s\)
and of \(\iota_s(g^{-\mathsf T})\) at odd \(s\), and let
\(\rho_{\rm let}(g)=\bigotimes_s\rho_s(g)\).
On even embeddings, \(p\) times the transpose of the
conjugated rational matrix is
\(p\iota_s(\gamma_e^{-1})\), with middle factor
\(p d^{-1}\). On odd embeddings the transpose instead is
\(\iota_s(\gamma_e^{\mathsf T})\), with middle factor \(d\).
Outside this orbit, every positive-degree factor lies away from
\(\lambda,\bar\lambda\), where \(\gamma_e\) and its inverse are
integral. A factor of degree zero has trivial action even above the pair.
Therefore the \(p^{A_\lambda}\)-scaled source matrix for
this edge is the tensor product of an integral invertible
transition on the other Weyl factors and the transpose of
\[
   \rho_{\rm let}(g_e)\,
       \operatorname{diag}_{\mu}(p^{\operatorname{dep}(\mu)})\,
       \rho_{\rm let}(\beta_b^{-1}).
\]
Number the coordinates from \(0\) to \(k-1\). Then
\(\operatorname{dep}(\mu)\) counts exponents in coordinates
\(r,\ldots,k-1\) on even embeddings and \(0,\ldots,r-1\) on odd
embeddings. This follows by taking symmetric powers of the two
diagonal matrices. In particular \(Y_\lambda\) preserves \(M\).
The formula is the edge sum on the original class lattice.

\indent{\it Size conditions for the local calculation.}
Set \(\theta=(m^2f)^{-1}\).
The polynomial estimates below use
\[
  0<2r\le k,\quad 0<\theta\le1,\quad a^4\le q,\quad
             (N_{\rm orb}k\log q/\theta)^5 P^{3r}\le a .
\]
These hold uniformly for the primes in \(\Lambda\) for all
sufficiently large \(x\). Indeed \(m\asymp\sqrt x\),
\(k\asymp x^5\), \(f=O(x^{11}\log x)\), \(N_{\rm orb}\le f\),
\(r=r_p=O(\log x)\), and \(e\le r\).
In particular \(2r\le k\) and \(r^2m\le k\).
Since \(a=xb\) and \(\log b=x^\epsilon+O(1)\),
\[
            a^4\le b^5\le\exp(m)\le p^{eN_{\rm orb}}=q .
\]
The exceptional prime \(p\) is inert in \(K\), whereas all
primes below \(13\) divide \(q_0\) and split in \(K\) by
Lemma~\ref{lem:auxiliary-residue-class}. Thus \(p\ge13\).
Since \(p<x_1\), we have \(\log q\le eN_{\rm orb}\log x_1\).
Also \(P^{3r}\le x_1^{3r^2}\le\sqrt b\), because
\((\log x)^2\log x_1=o(\log b)\).
Finally \((eN_{\rm orb}^2 k\log x_1/\theta)^5\) is bounded by a
fixed power of \(x\), uniformly in \(p,\lambda\).
Its product with \(\sqrt b\) is at most \(a\) for large
\(x\). This gives the second large-degree inequality.

Put
\[
      C_{\rm div}=N_{\rm orb}a\left(P^{-2r}+\frac{r^2}{kP}\right).
\]
Let \(w^2=p\) on extending scalars if necessary. Let
\(b_\mu\) be the diagonal symmetric-power pairing, with
integral complementary weights \(b'_\mu\) satisfying
\(b_\mu b'_\mu=p^{h_0}\).
For the exact common exponent put \(h_0=N_{\rm orb}k^2eN_{\rm orb}\).
If \(M_s=a!/\prod_i\mu_{s,i}!\), take
\(b_\mu=p^{h_0}/\prod_sM_s\) and \(b'_{\mu}=\prod_s M_s\).
Since \(a<q=p^{eN_{\rm orb}}\), the factorial valuation formula
bounds \(v_p(M_s)\) by \(keN_{\rm orb}\). This proves integrality.
The multinomial theorem proves that substitution by a
matrix and by its transpose are adjoint for the weights
\(1/M_s\), and hence for their scaled tensor product.
Finally \(k^2eN_{\rm orb}\,P^{2r}\le a\) follows from the displayed
large-degree bound and \(eN_{\rm orb}\le\log q\), so
\(h_0\le N_{\rm orb}a/P^{2r}\le C_{\rm div}\).
Write \(s_{\mathcal P}=[\mathrm{GL}_k(W):\mathcal P]\)
and \(d_{\rm let}\) for the number of monomials of degree \(a\)
in each of the \(N_{\rm orb}\) embedding factors. Thus
\[
 d_{\rm let}=\binom{a+k-1}{k-1}^{N_{\rm orb}}.
\]

Write \(W_0=\operatorname{diag}_\mu(w^{\operatorname{dep}(\mu)})\).
A saturated halfway restriction is a vector \(W_0F\)
with \(F\) over the coefficient fraction field such that
\(W_0\rho_{\rm let}(g)F\) is integral in every parahoric
direction. To see that representatives suffice, let \(u\in\mathcal P\).
In the linear substitution at each embedding, an entry with shallow
row and deep column is divisible by \(p=w^2\). For even embeddings
this is the upper-right block of \(u\), and for odd embeddings it is
the lower-left block of \(u^{-\mathsf T}\). Conjugating by the
diagonal scaling of deep coordinates by \(w\) divides those entries
by \(w\) and preserves integrality. Applying this also to \(u^{-1}\)
and taking symmetric and tensor powers gives
\[
 W_0\rho_{\rm let}(u)W_0^{-1}\in\mathrm{GL}_{d_{\rm let}}(\mathcal O).
\]
Choose \(h_{\rm clr}\ge0\) such that \(F_{\rm int}=w^{h_{\rm clr}}F\) is integral.
The condition on a saturated restriction is precisely
\[
 W_0\rho_{\rm let}(g)F_{\rm int}\in
     w^{h_{\rm clr}}\mathcal O^{d_{\rm let}}
 \qquad(g\in\mathrm{GL}_k(W)).
\]

\subsection{The depth count for the endpoint restrictions}
\label{sec:exceptional-source-depth-count}

\sectionopening{In this subsection, we bound the valuations of the elementary divisors
of the two restrictions of Section~\ref{sec:exceptional-hecke-analysis}.
The proof counts coefficients satisfying the divisibility conditions
in coordinates adapted to the parahoric.}

Work over the residue coefficient field \(K_0\) of
characteristic \(p\) for the cycle of embeddings above. Write
\(k=2r+k_1\). A variable \(X_{b,s}\), \(0\le b<k\), is deep
when \(s\) is even and \(b\ge r\), or \(s\) is odd and \(b<r\).
Let \(D_j\) be the space of degree-\(a\) forms in each embedding
whose translates in all parahoric directions have no
coefficient of deep degree less than \(j\).
The condition then holds after every \(\mathrm{GL}_k(W)\)
translate. In fact the parahoric preserves the deep-degree
filtration modulo \(p\) in both parities. The action on the reduced coefficients factors through
\(\mathrm{GL}_k(W/pW)\).

Treat the complementary variables as coefficients in the field
\[
 L_0=K_0\bigl(X_{b,s}:
     (s\text{ even},\, b<r)\ \text{or}\ (s\text{ odd},\, b\ge r)\bigr).
\]
For each deep variable \(v=(b,s)\) put
\[
 \kappa_v=\begin{cases}P,&b<2r,\\P^{2r},&b\ge2r.\end{cases}
\]
Give a monomial in the deep variables weight
\(\sum_v\mu_v/\kappa_v\).
We show that any nonzero form in \(D_j\) has a nonzero
coefficient of weight at least \(j\).

Let \(B\) run through matrices over \(\mathbb F_q\) supported
on entries \((b,b+r)\), \(b<r\), and
\((b,t)\), \(b<r,\ 2r\le t<k\). There are \(r(k_1+1)\)
independent entries. Such matrices satisfy \(B^2=0\).
The translation by \((1+B)^{-1}\) fixes every complementary variable. It takes a deep variable \(X_v\) to \(X_v-z_v(B)\),
where \(z_{t,s}\) at even \(s\) is
\(\sum_{b<r}X_{b,s}\iota_s(B_{bt})\), for \(t\ge r\),
and \(z_{b,s}\) at odd \(s\) is
\(-\sum_{t\ge r}X_{t,s}\iota_s(B_{bt})\), for \(b<r\).
For \(F\in D_j\) the deep degree of every term in the translate
is at least \(j\), so \(F(X-z(B))\in(X_v:v\text{ deep})^j\).
Using \(-B\) gives
\(F\in(X_v-z_v(B):v\text{ deep})^j\) at each point \(z(B)\).

We claim that all the coordinate functions
\(\iota_i(B_{bt})\), for any fixed embedding \(i\), are
\(L_0\)-linear combinations of the functions \(z_{v}(B)^{P^u}\)
with \(s+u=i\) and \(P^u<\kappa_v\).
For a column \(t\ge2r\), choose the \(r\) distinct integers
\(0\le u<2r\) for which \(s=i-u\) is even.
The \(r\) values \(z_{t,i-u}(B)^{P^u}\) have coefficient
matrix
\[
                 \bigl(X_{b,i-u}^{P^u}\bigr)_{u,b}.
\]
This uses \(\iota_s(y)^{P^u}=\iota_{s+u}(y)\) on the
residue field. The determinant is nonzero. Distinct
permutations use different monomials, because within any
one column the variables are distinguished either by
their embedding or by the distinct exponent \(P^u\).

Thus these values solve for every entry of the column. For
\(b<r\) the remaining entries satisfy
\[
 \begin{aligned}
 z_{b+r,i}(B)&=X_{b,i}\iota_i(B_{b,b+r}) &&(i\text{ even}),\\
 -z_{b,i}(B)&=X_{b+r,i}\iota_i(B_{b,b+r})
       +\sum_{2r\le t<k}X_{t,i}\iota_i(B_{bt}) &&(i\text{ odd}).
 \end{aligned}
\]
The displayed coefficients of \(\iota_i(B_{b,b+r})\) are nonzero
in \(L_0\). Subtracting the known sum in the odd case proves the
claim, also when \(k=2r\) and the sum is empty. When \(r=1\),
the interpolation matrix for each \(t\ge2r\) has a single
nonzero entry \(X_{0,i-u}^{P^u}\).

It follows first that the points \(z(B)\) are distinct.
It also follows, by taking the \(P\)-power or the
\(P^{2r}\)-power of the coordinate expression at its
original embedding, that for every deep variable \(v\)
\[
             z_v(B)^{\kappa_v}=g_v(z(B)),
\]
where \(g_v\) is a polynomial of weight strictly less
than one. Indeed the power moves the coefficient
functions to the corresponding later embedding, where the
claim expresses them using the allowed lower powers.
There are \(N_{\rm orb}r\) deep coordinates with \(b<2r\) and
\(N_{\rm orb}k_1/2\) others, the latter all in even embeddings.
Since \(q=P^{N_{\rm orb}}\), the number of points equals
\[
 q^{r(k_1+1)}
   =P^{N_{\rm orb}r}(P^{2r})^{N_{\rm orb}k_1/2}
   =\prod_v\kappa_v .
\]

For completeness, the weighted interpolation argument
for these points is as follows. Let
\(\psi_v=X_v^{\kappa_v}-g_v\), of weight one.
The ideals of the points are pairwise comaximal, as are
their \(j\)-th powers. The map to truncated polynomials
of degree less than \(j\) at all these points is therefore
surjective. Division using \(X_v^{\kappa_v}=\psi_v+g_v\)
expresses each polynomial in the span of
\(X^\alpha\psi^\gamma\), \(0\le\alpha_v<\kappa_v\).
In each replacement the term with \(g_v\) strictly lowers
weight and the other term removes \(X_v^{\kappa_v}\)
from the monomial preceding \(\psi^\gamma\).

At each step the pair
\((\operatorname{wt}(\alpha)+|\gamma|,\operatorname{wt}(\alpha))\)
decreases in lexicographic order on pairs of elements of
\(P^{-2r}\mathbb Z_{\ge0}\). Thus division terminates without
increasing total weight. For \(j\ge1\), all products with
\(|\gamma|\ge j\) vanish in the jets, since each \(\psi_v\)
vanishes at every point. Let \(d\) be the number of deep
coordinates. Translation to each point identifies its jet algebra
with \(L_0[X_v]/(X_v)^j\), where the ideal is generated by all
monomials of total degree \(j\). The remaining products span the
target and number
\[
       \Bigl(\prod_v\kappa_v\Bigr)\binom{d+j-1}{d}
           =\dim_{L_0}\bigoplus_B L_0[X_v]/(X_v-z_v(B):v)^j .
\]
Their images are therefore independent. If \(F\) had weight less than
\(j\), division would use only products with
\(|\gamma|<j\). Its zero image would force \(F=0\).
The change of complementary variables to coefficients is
injective, so it follows that a nonzero original
coefficient of \(F\) has weight at least \(j\).

Projection onto the monomial coordinates of weight
at least \(j\) is consequently injective on \(D_j\).
A monomial of weight \(u\) is counted at no more than
\(u\) positive integer thresholds. In the uniform
average over monomials of degree \(a\) in each embedding,
each variable has mean exponent \(a/k\). There are
\(N_{\rm orb}r\) deep variables with \(b<2r\), and at most \(N_{\rm orb}k\)
other deep variables. Thus the mean weight is at most
\[
    \frac ak\left(\frac{N_{\rm orb}r}{P}+\frac{N_{\rm orb}k}{P^{2r}}\right)
       \le N_{\rm orb}a\left(P^{-2r}+\frac{r^2}{kP}\right)=C_{\rm div}.
\]
Summing the injective dimension bounds gives
\(\sum_{j\ge1}\dim D_j\le C_{\rm div}d_{\rm let}\),
as required for both endpoint restrictions.

\subsection{The rank-two defect estimate}
\label{sec:exceptional-root-defect}

\sectionopening{In this subsection, we make the root calculation needed for the source
Hecke operator of Section~\ref{sec:exceptional-hecke-analysis}. The
two-dimensional calculation bounds the coefficients on the saturated
restrictions of Section~\ref{sec:exceptional-source-depth-count}.}

Fix two distinct indices \(b,b'\) and an embedding \(i\).
Set \((X_s,Y_s)=(X_{b,s},X_{b',s})\) for even \(s\) and
\((X_s,Y_s)=(X_{b',s},X_{b,s})\) for odd \(s\).
Thus \(Y_s\) is the positive variable in the defect.
Put \(P_s=P^{(s-i)\bmod N_{\rm orb}}\). For an exponent \(\alpha\)
define
\[
 D=q\sum_s P_s\alpha(X_s),\qquad
 E=q\alpha(Y_i)+\sum_{s\ne i}P_s\alpha(Y_s).
\]
The defect of \(\alpha\) at \(i\) is \((E-D/q)/(q-1)\).

First work over the integral coefficient ring with
\(w^2=p\). Let \(c\) give weights on the variables with
\(1/2\le c(Y_s)-c(X_s)\le1\).
For a coefficient \(a_\alpha\), give the term
\(a_\alpha X^\alpha\) total weight
\(\operatorname{ord}_w(a_\alpha)+c\cdot\alpha\), where for \(a\ne0\)
\[
 \operatorname{ord}_w(a)=\max\{j\in\mathbf Z_{\ge0}:a\in w^j\mathcal O\}
                         =\lfloor2v(a)\rfloor ,
 \qquad \operatorname{ord}_w(0)=+\infty .
\]
Writing \(o=\operatorname{ord}_w\), we have
\[
 o(ab)\ge o(a)+o(b),\qquad o(a+b)\ge\min\{o(a),o(b)\}.
\]
Thus lower bounds on total weight pass to products and sums.
Over the residue field every nonzero coefficient has order zero.
Suppose all terms of all translates of \(F\) by the
root \(\mathrm{SL}_2(W)\) have weight at least
\(l+c\cdot\alpha\), where the coefficient \(a_\alpha\) of \(X^\alpha\)
is nonzero and \(l=\operatorname{ord}_w(a_\alpha)\).
Thus \(a_\alpha\ne0\bmod w^{l+1}\).
Over the residue field take \(l=0\) and omit reduction.
We prove \(Eq\le D\).

Put \(\epsilon_s=-1\) if \(s\) has the parity of \(i\),
and \(\epsilon_s=1\) otherwise. In a new variable \(T\)
make the substitution \(\mathcal A_T\) given by
\[
 \begin{aligned}
 X_i&\longmapsto T^qX_i,&Y_i&\longmapsto X_i,\\
 X_s&\longmapsto T^{qP_s}X_s-\epsilon_s T^{P_s}Y_s,&
 Y_s&\longmapsto Y_s-\epsilon_s X_s &&(s\ne i).
 \end{aligned}
\]
All other variables are fixed. Let \(\bar\alpha\) be the
exponent obtained by merging \(Y_i\) into \(X_i\), and
let \(\Phi(T)\) be the coefficient at \(\bar\alpha\)
of \(\mathcal A_TF\), reduced modulo \(w^{l+1}\).

For the expansion at infinity put
\(\rho=(c(Y_i)-c(X_i))/q\). Give \(T\) weight \(-\rho\),
\(X_s\) weight \(c(X_s)+qP_s\rho\), and every other
variable its original weight. The main term of each
substituted variable has its original weight. At \(i\)
this uses \(c(X_i)+q\rho=c(Y_i)\). At \(s\ne i\)
the cross terms have strictly greater weight, since
\[
 P_s\rho\le 1/P<1/2,\qquad qP_s\rho\ge P/2>1 .
\]
The variable weight of \(T^D X^{\bar\alpha}\) is \(c\cdot\alpha\).
Any contribution using a cross term has strictly greater total
weight than \(l+c\cdot\alpha\), so its integral coefficient order
is at least \(l+1\). For a contribution using only main terms
the exponents away from \(X_i,Y_i\) are unchanged, and
\[
 \alpha(X_i)=D/q-\sum_{s\ne i}P_s\bar\alpha(X_s),\qquad
 \alpha(Y_i)=\bar\alpha(X_i)-\alpha(X_i).
\]
Hence \([T^D]\Phi=a_\alpha\bmod w^{l+1}\) is nonzero.
For degree \(n>D\) the variable weight drops by \((n-D)\rho>0\);
the coefficient order is again at least \(l+1\). Hence \(\deg\Phi=D\).

For each finite point take \(r_0\in W\), put
\(t_0=\iota_i(r_0)\), and shift \(T\mapsto T+t_0\).
The root matrix is
\(\left(\begin{smallmatrix}0&1\\-1&-r_0\end{smallmatrix}\right)\)
if \(i\) is even, with \(+r_0\) in the last entry if \(i\)
is odd. Denote it by \(g_{r_0}\), embedded by the identity on
the other coordinates. It has determinant one. As \(r_0\) runs through
residue representatives, its second column \((1,\pm r_0)^{\mathsf T}\)
represents the \(q\) affine directions in \(\mathbf P^1(W/pW)\).
The remaining direction is the second column of the identity. Put
\(\delta_s=-1\) at even \(s\) and \(\delta_s=1\) at odd \(s\).
Reading its two columns at even embeddings and those of its
inverse transpose in the reversed order at odd embeddings gives
\[
 \rho_{\rm let}(g_{r_0})X_s=\delta_sY_s,\qquad
 \rho_{\rm let}(g_{r_0})Y_s=
   -\delta_s\bigl(X_s+\epsilon_s\iota_s(r_0)Y_s\bigr).
\]

Compose that root substitution with \(\mathcal A_{T+t_0}\).
The result on each of \(X_s,Y_s\) is the signed image
under \(\mathcal A_T\) of the swapped variable,
up to errors on the image of \(Y_s\) with coefficients
in \(p\mathcal O[T]\).
Indeed for \(s\ne i\) the change in \(\mathcal A_T(X_s)\)
modulo \(p\) is
\[
 t_0^{P_s}(X_s-\epsilon_sY_s)
       =-\epsilon_s\iota_s(r_0)\mathcal A_T(Y_s)
                         \pmod p .
\]
Here the shifts of \(T^{qP_s}\) and \(T^{P_s}\) have
the same constant modulo \(p\), by the Frobenius cycle.
At \(i\) the change is \(t_0^qX_i=t_0X_i\) modulo \(p\),
and \(-\epsilon_i=1\). The errors consist only of
\(p\)-multiples of polynomials times \(X_s,Y_s\), or
of \(X_i\) after merging. This follows by expanding the
two shifted \(p\)-powers. Over the residue field the
errors vanish.

For this finite expansion give \(T\) weight \(+\rho\),
\(X_s\) weight \(c(X_s)\), \(Y_s\) weight
\(c(Y_s)-P_s\rho\) for \(s\ne i\), and \(Y_i\) weight
\(c(Y_i)-q\rho\). The signed image of the swapped
variable has weight at least the weight of the source vector. For \(s\ne i\) the required two inequalities
are \(c(Y_s)-P_s\rho\ge c(X_s)\) and
\(c(X_s)+qP_s\rho\ge c(Y_s)\); at \(i\) they
hold with the stated \(\rho\). Each error still has
weight at least the source weight, because a factor
\(p=w^2\) contributes two, while the difference
lost in either variable is at most \(1+1/P\).
Applying this to the inverse root translate of \(F\)
shows that all terms of \(\mathcal A_{T+t_0}F\)
have weight at least \(l+c\cdot\alpha\) in this
finite expansion. Now \(T^E X^{\bar\alpha}\) has weight
\(c\cdot\alpha\). Every lower \(T\)-degree at that
exponent has zero coefficient modulo \(w^{l+1}\).
Thus \(\Phi\) vanishes at \(t_0\) to order at least \(E\).

Choose one \(r_0\in W\) lifting each of the \(q\) classes of \(W/pW\).
Distinct lifts differ by a unit of \(W\), whose image under \(\iota_i\)
is a unit even modulo \(w^{l+1}\). Therefore the ideals generated by
the distinct \(T-t_0\) are pairwise comaximal in the quotient
polynomial ring, as are their \(E\)-th powers. Their product divides
\(\Phi\). Since this product is monic of degree \(Eq\),
multiplication by it raises the degree of every nonzero polynomial
by \(Eq\), even over \(\mathcal O/(w^{l+1})\). Thus \(Eq\le D\). This proves
that the defect of the extreme term is nonpositive.
Taking \(l=0\) over the residue field with the height
weights proves the least-height fact used in the flag
count.

For the ordinary depth estimate take an ordered root
\(b<r\le b'\), in particular the boundary root \((r-1,r)\).
In both parities \(X_s\) is shallow and \(Y_s\) is deep, so
\(\operatorname{dep}(Y_s)-\operatorname{dep}(X_s)=1\).
Fix \(i\) and change these weights by
subtracting one third of the defect functional at \(i\).
The increments of the changed weights are still in
\([1/2,1]\): the increments of the defect are
\(1+2/(q-1)\) at \(i\) and \(2P_s/(q-1)\) elsewhere,
all in \([0,3/2]\).

Suppose every root translate initially has order plus
ordinary depth at least \(v_0\). If the changed bound
failed in even one translate, choose an extreme term of
minimal order plus changed weight over all translates.
Such a minimum exists, because the exponents lie in a
finite set and nonnegative integral orders with weight
below \(v_0\) are bounded. The estimate just proved says that the extreme term has nonpositive
defect. Its ordinary depth plus order would therefore
be below \(v_0\), a contradiction.
Thus
\[
    \operatorname{ord}_w(a_\beta)+\operatorname{dep}(\beta)
             \ \ge\ v_0+\operatorname{defect}_i(\beta)/3 .
\]

In particular let \(v=w^{-h_{\rm clr}}W_0F\) be a saturated
halfway restriction with \(F\) integral and
\(w^{-h_{\rm clr}}W_0\rho_{\rm let}(g)F\) integral in every
direction \(g\). For a block-diagonal Levi element \(g_L\) and
any root element \(g_2\), the product \(g_2g_L\) is an allowed direction.
The depth hypothesis therefore holds for every root translate of
\(F_L=\rho_{\rm let}(g_L)F\) with \(v_0=h_{\rm clr}\).
Since \(g_L\) preserves the block degrees, its action commutes with
\(W_0\). For every nonzero coefficient this gives
\[
 \operatorname{ord}_w\bigl((\rho_{\rm let}(g_L)v)_\beta\bigr)
  = \operatorname{ord}_w\bigl((F_L)_\beta\bigr)
       +\operatorname{dep}(\beta)-h_{\rm clr}
  \ge \operatorname{defect}_i(\beta)/3 .
\]
Here the equality uses \(\lfloor 2v_p(z)+d\rfloor
=\lfloor 2v_p(z)\rfloor+d\) for \(d\in\mathbf Z\).
If \(\operatorname{defect}_i(\beta)\ge3u\) for
\(u\in\mathbf Z_{\ge0}\), this coefficient belongs to \(w^u\mathcal O\).
Writing \(e_\beta\) for the monomial vector and using the integral
pairing weights \(b_\beta\), we also have
\[
 \langle\rho_{\rm let}(g_L^{-1})e_\beta,v\rangle
   = b_\beta (\rho_{\rm let}(g_L^{-\mathsf T})v)_\beta
   \in w^u\mathcal O ,
\]
since \(g_L^{-\mathsf T}\) belongs to the same Levi.

\subsection{Counting forms with bounded local defect}
\label{sec:exceptional-defect-packing}

\sectionopening{In this subsection, we use the two-variable estimate of
Section~\ref{sec:exceptional-root-defect} to bound a dimension in the
residue representation. We count the permitted monomials after
imposing the same conditions under all elements of the Levi.}

Write \(b_0=r-1\) and take the defect bound \(3t\) at the
root \((b_0,r)\), where \(t\) is an integer with \(t\ge2(k+1)\). Work in the residue
coefficient field of characteristic \(p\). Let \(S\) be
the degree-\(a\) forms in the embedding factors for which every translate
by the reduction of the two-block Levi has all defects along this root
at most \(3t\). The reduction is
\(\mathrm{GL}_r(\mathbb F_q)\times\mathrm{GL}_{k-r}(\mathbb F_q)\)
acting through the prescribed Frobenius embeddings into \(K_0\).
The other simple roots lie inside the Levi.

The rank-two estimate of Section~\ref{sec:exceptional-root-defect}
gives the following fact at the roots inside one Levi block. Among the nonzero coefficients
of all Levi translates of a nonzero polynomial, choose
one of least height, with height of \(X_{b,s}\) given
by \(b\) on even embeddings and \(-b\) on odd embeddings.
The defect at that coefficient is nonpositive at each
embedding for every simple root inside the block.
The required root calculation applies without a
homogeneity requirement on the polynomial, since the
root substitutions are linear and preserve the degree
at each embedding. 

A least-height coefficient occurs at a representative
of a Levi flag. Write \(g=u x\), with \(u\) lower triangular
modulo \(p\) and \(x\) chosen for the quotient by those matrices
acting on the left. On even
embeddings an off-diagonal term of \(u\) strictly increases
the variable index, and on odd embeddings inverse transpose
strictly decreases it. Each increases the height.
Thus at the least-height coefficient \(u\) contributes
only its diagonal action, with a nonzero factor.
The coefficient at \(x\) is therefore already nonzero.

For a positive integer \(M=q^j\), take a family
\((F_u)_u\) in \(S\) indexed by the vectors \(u\) with
\(0\le u_{b,s}<M\). In each embedding put
\[
  \gamma_{u,b,s}=u_{b,s}
        +\mathbf1_{b=b_0}\,M\left(
             k-\left\lfloor\frac{\sum_c u_{c,s}}M\right\rfloor
                              \right).
\]
Then \(X^{\gamma_u}\) keeps the original residue \(u\)
modulo \(M\), while the degree of \(\gamma_u\) in each
embedding belongs to \([Mk,M(k+1))\).
Pack the family as
\[
               F_{\rm pack}(X)=\sum_u X^{\gamma_u}F_u(X^M).
\]
This linear map is injective by distinct exponent residues.
The residue entries of every translated matrix lie in the
images of \(\mathbb F_q\), so the substitutions commute
with \(X\mapsto X^M\). Consequently any nonzero coefficient
at \(\alpha\) in a Levi translate of the packed polynomial
comes from a term with exponent
\(\alpha=M\beta+\eta\), where \(\beta\) occurs in a translate
of some \(F_u\) and \(\eta\) occurs in the translate of its
monomial \(X^{\gamma_u}\). The degrees of \(\eta/M\) are in
\([k,k+1)\), because linear substitutions preserve them.

The defects along \((b_0,r)\) of \(\beta\) are at most \(3t\).
In the defect of \(\eta/M\), the term with factor \(q\) is
bounded by \((k+1)q/(q-1)\); the singled-out term with factor \(-1\)
is nonpositive. Including the other weighted differences gives
\[
 \frac{k+1}{q-1}\left(q+\sum_{l=1}^{N_{\rm orb}-1}P^l\right)
   =(k+1)\left(1+\frac1{P-1}\right)\le2(k+1).
\] Thus each of these defects of \(\alpha/M\) is at most
\(3t+2(k+1)\le4t\). Its degree at each embedding is in
\([a+k,a+k+1]\).

We move the boundary-root bound into the zero-defect cone
by a common integral shift. Put \(\mu_b=P^{2(b_0-b)}\)
for \(b\le b_0\) and \(\mu_b=0\) otherwise.
Put \(\nu_{b,s}=\mu_b\) on even embeddings and \(P\mu_b\)
on odd embeddings. Put \(\varepsilon_s=(-1)^s\). The cyclic steps
\(\sigma_{b,s}=
\varepsilon_s(P\nu_{b,s+1}-\nu_{b,s})\) are \((P^2-1)\mu_b\) on even
embeddings and zero on odd embeddings. It is integral and
nonnegative, with degree at most \(P^{2r}-1\) at every
embedding. Along a simple root, its defects are
\(P^2\mu_{b+1}-\mu_b\) on even embeddings and
\(P(\mu_{b+1}-\mu_b)\) on odd embeddings. Thus they are
at most \(-1\) at the boundary root and at most zero
at all other simple roots.

Write \(\mathcal C_0\) for the cone in which all simple-root
defects are nonpositive, and \(\mathcal D\) for the nonnegative
degree slab with the bounds below. For a nonzero packed
family, choose a least-height
coefficient \(\alpha\) at a flag representative \(x\).
The fact about roots inside the Levi and the last
estimates show
\[
  \alpha/M+4t\sigma\in\mathcal C_0\cap\mathcal D,\qquad
  D_s=a+k+1+4t\sum_b\sigma_{b,s}\le a+4tP^{2r}.
\]
Since \(\alpha\) and \(4t\sigma\) are integral arrays, the point lies in
\(M^{-1}\mathbb Z^{kN_{\rm orb}}\).

Therefore the linear map that lists the coefficients
at \(M(y-4t\sigma)\) in the flag translates, for
\(y\in\mathcal C_0\cap\mathcal D\cap M^{-1}\mathbb Z^{kN_{\rm orb}}\),
is injective on the packed families. Points giving
negative coordinates can be skipped.
It follows that \(M^{kN_{\rm orb}}\dim S\) is bounded by the
number of Levi flags times the number of these lattice
points. Divide by \(M^{kN_{\rm orb}}\) and let \(M=q^j\to\infty\).
The bounded polytope has boundary of measure zero,
and its normalized lattice counts converge to its volume
by the ordinary Riemann sums. This is the required
lattice-point bound for the volume calculation.

\subsection{The volume factor in the residue-defect estimate}
\label{sec:exceptional-pieri-volume}

\sectionopening{In this subsection, we compute the volume ratio in the dimension estimate of
Section~\ref{sec:exceptional-defect-packing}. Ordering the
exponents gives interlacing inequalities. Their volume
normalization follows by determinant expansion.}

Use the same parity sign \(\varepsilon_s=(-1)^s\).
For permutations \(\sigma_s\in S_k\), set
\[
 (T_\sigma \nu)_{b,s}
       =\varepsilon_s(P\nu_{b,s+1}-\nu_{\sigma_s(b),s}).
\]
Write \(q=P^{N_{\rm orb}}\). The sum of the entries at each embedding
is independent of the permutations. Let \(\mathcal D\)
be any slab of nonnegative vectors whose degree in the
embedding \(s\) is in \([D_s-1,D_s]\).
In coordinates \(T_1\nu\), the zero-defect cone
\(\mathcal C_0\) is given by the additional inequalities
\[
 \begin{cases}
  P\nu_{b+1,s+1}\le \nu_{b,s}, &\varepsilon_s=1,\\
  \nu_{b+1,s}\le P\nu_{b,s+1}, &\varepsilon_s=-1,
 \end{cases}
 \qquad (b<k-1).
\]
Indeed put \(\beta=T_1\nu\) and
\(d_s=\nu_{b,s}-\nu_{b+1,s}\). The signed weighted
difference of the two relevant beta coordinates is
\(d_s-Pd_{s+1}\). Summing with weights
\(P^{(s-i)\bmod N_{\rm orb}}/(q-1)\) gives \(-d_i\), because the
cycle closes with \(P^{N_{\rm orb}}=q\). The remaining term of
the defect is \(\beta_{b+1,i}\) if \(i\) is even and
\(\beta_{b,i}\) otherwise. Thus zero defect is the
displayed pair of inequalities.
We prove
\[
 \operatorname{vol}(\mathcal C_0\cap\mathcal D)
       =\frac{\operatorname{vol}(\mathcal D)}
             {\prod_{b=1}^k(1+q+\cdots+q^{b-1})}.
\]

First take \(\nu\) with distinct entries in each embedding,
avoiding also equality in adjacent comparisons. These
restrictions discard only finitely many hyperplanes.
For each \(s\) let \(M_s\) be the \(k\)-square matrix of
indicators of the inequalities
\(\varepsilon_s(P\nu_{b,s+1}-\nu_{c,s})>0\).
The product \(\prod_s\det M_s\) is unchanged if the entries
in any one embedding are independently permuted, since this
permutes the rows of one determinant and the columns of
another with the same sign.

Sort the entries of \(\nu\)
decreasingly in each embedding. For two decreasing sequences
\(x,y\), the determinant of \((1_{x_b>y_c})_{b,c}\)
is the indicator of interlacing
\(x_0>y_0\ge x_1>y_1\ge\cdots\).
Indeed each row consists of zeros followed by ones,
and the first positions of the ones are nondecreasing.
The determinant vanishes if two rows coincide or there
is a zero row. In the remaining case the first positions
are \(0,\ldots,k-1\) and the matrix is triangular with
diagonal ones. For the reversed inequality transpose
and interchange the two sequences. Thus the product of
determinants indicates the cone inequalities and
positivity on the sorted arrays.

Conversely the cone inequalities and positivity force
the entries of \(\nu\) to be decreasing in each embedding,
by the same interlacing inequalities. Therefore the
number of tuples of permutations of the embeddings sending
\(\nu\) to a point satisfying these conditions equals that
product of determinants. Expanding the determinants
gives the signed sum over the tuples \(\sigma\) of the
indicators of \(T_\sigma\nu>0\), with sign
\(\prod_s\operatorname{sgn}\sigma_s\). Multiply by the
degree-slab condition, common to all terms, and integrate.

To compute the resulting volumes, eliminate the embedding
vectors cyclically in
\(P\nu_{s+1}=\Pi_{\sigma_s}\nu_s+\varepsilon_s\beta_s\),
where \(\Pi_{\sigma_s}\) is the corresponding permutation
matrix. Eliminating \(N_{\rm orb}-1\) vectors gives diagonal factors
\(P\) on each eliminated coordinate. The final block is
\(P(1-q^{-1}\Pi_\sigma)\), with \(\Pi_\sigma\) the product
round the cycle. Thus
\[
 |\det T_\sigma|
       =P^{kN_{\rm orb}}|\det(1-q^{-1}\Pi_\sigma)|.
\]
Each determinant is nonzero since permutation eigenvalues
have absolute value one and \(q>1\).
For a prescribed product \(\pi\) there are \((k!)^{N_{\rm orb}-1}\)
tuples yielding it.

Each permutation of the embedding coordinates has Jacobian of
absolute value one. Integrating the indicator identity gives
\[
 \frac{(k!)^{N_{\rm orb}}}{|\det T_1|}
       \operatorname{vol}(\mathcal C_0\cap\mathcal D)
 =\operatorname{vol}(\mathcal D)
       \sum_{\sigma\in(S_k)^{N_{\rm orb}}}
       \frac{\prod_s\operatorname{sgn}(\sigma_s)}{|\det T_\sigma|}.
\]
Grouping by the product permutation yields
\[
 \operatorname{vol}(\mathcal C_0\cap\mathcal D)
   =\operatorname{vol}(\mathcal D)\frac{(1-q^{-1})^k}{k!}
       \sum_{\pi\in S_k}
       \frac{\operatorname{sgn}\pi}{|\det(1-q^{-1}\Pi_\pi)|}.
\]
The same determinant calculation holds for each choice of the
embedding signs. We can therefore evaluate the factor, which depends
only on \(q,k\), with one positive embedding, \(P=q\), and the
simplex \(\sum_b\beta_b\le1\), \(\beta_b\ge0\), of positive volume.
In this case the
step map is scalar multiplication by \(q-1\), and the
cone is the cone of \(\beta_b\ge q\beta_{b+1}\),
\(\beta_{k-1}\ge0\). Put
\(d_b=1+q+\cdots+q^b\) and \(\beta_k=0\).
The triangular map
\(\beta_b\mapsto d_b(\beta_b-q\beta_{b+1})\)
maps this cone to the nonnegative orthant, preserves the
sum of coordinates and has determinant \(\prod_b d_b\).
This proves the claimed factor for every \(N_{\rm orb}\).

Here is the consequence of the corresponding lattice-point
bound. Section~\ref{sec:exceptional-defect-packing} shows that the bounded-defect space \(S\) has dimension
at most the number \(F_L\) of flags of the two-block Levi
times \(\operatorname{vol}(\mathcal C_0\cap\mathcal D)\),
with \(1\le D_s\le a+E\) and \(E=4tP^{2r}\).
The parahoric index times \(F_L\) is the number of all
flags of \(\mathbb F_q^k\), namely
\(\prod_{b=1}^k(1+q+\cdots+q^{b-1})\).
A complete flag is uniquely its \((k-r)\)-plane with the
flags on that plane and on the quotient. Thus the volume
factor cancels this index and
\[
       [\mathrm{GL}_k(W):\mathcal P]\dim S
                  \le \operatorname{vol}(\mathcal D).
\]
The nonnegative vectors in one embedding with degree in
\([D_s-1,D_s]\) have volume
\((D_s^k-(D_s-1)^k)/k!\le D_s^{k-1}/(k-1)!\).
If \(8N_{\rm orb}ktP^{2r}\le\theta a\), then
\(N_{\rm orb}(k-1)E/a\le\theta/2\). Therefore
\[
 \operatorname{vol}(\mathcal D)
 \le \exp\!\left(\frac{N_{\rm orb}(k-1)E}{a}\right)
                  \left(\frac{a^{k-1}}{(k-1)!}\right)^{N_{\rm orb}}
 \le (1+\theta)\binom{a+k-1}{k-1}^{N_{\rm orb}}.
\]
The last step uses \(0<\theta\le1\) and
\(\exp(\theta/2)\le1+\theta\).

\subsection{An integral average in one matrix variable}
\label{sec:exceptional-single-average}

\sectionopening{In this subsection, we construct one integral functional on polynomials
in a generic matrix. The invariance and vanishing properties are
chosen for the Levi in Section~\ref{sec:exceptional-hecke-analysis},
while the congruence to evaluation at the identity will control the
pairing in Section~\ref{sec:exceptional-edge-average}.}

Let \(W\) be the local ring of the distinguished prime, with
residue field of order \(q=p^{eN_{\rm orb}}\), and fix one embedding into the coefficient ring \(\mathcal O\).
Let \(L_W\subseteq\mathrm{GL}_k(W)\) be the two-block Levi.
For integral weights \(\omega_i\) on the indices, call a
monomial \(X^\alpha\) in the \(k^2\) matrix entries bad if
its row weight \(\sum_{i,j}\omega_i\alpha_{ij}\) exceeds
a given threshold \(t_0\). Suppose \(D+k\le q\) and \(D<q\).
There is an integral functional \(I\) on polynomials in
\(X\) which is invariant under conjugation by \(L_W\),
kills bad monomials, and satisfies
\[
    I(F)\equiv p^{k^2eN_{\rm orb}}F(1)\pmod{\mathfrak b}
\]
whenever \(\deg F\le D\) and every Levi conjugate of \(F\)
has its bad coefficients in the ideal \(\mathfrak b\).

Here is the construction. Work in degree \(d\le D\), on
\(U_d=(\mathcal O^k)^{\otimes d}\).
Let \(\mathcal S_d\) be the Schur algebra with basis the
occurrence matrices of pairs of tensor indices with the
same sequence of block parameters. For such an occurrence
matrix \(\alpha\), the corresponding basis element is the
coefficient of \(X^\alpha\) in
\((X_{\rm block}^{\mathsf T})^{\otimes d}\), where off-block
entries of \(X_{\rm block}\) are set to zero. The basis
elements have disjoint supports in the matrix entries.
For independent block matrices \(X,Y\), put
\(A(X)=(X^{\mathsf T})^{\otimes d}\). Comparing coefficients
of \(X^\alpha Y^\beta\) in \(A(X)A(Y)=A(YX)\) expresses each product
of two occurrence operators as an integral combination of these operators.
For an exponent involving an off-block variable the occurrence operator is zero.

The tensors of matrices in \(L_W\) span \(\mathcal S_d\)
over \(\mathcal O\). To verify this, first work over \(W\).
If a homogeneous polynomial of degree \(d\) over
\(\mathbb F_q\) in the block entries vanishes on all
invertible block matrices, its product with the block
determinant vanishes at every point of the affine space.
This product is homogeneous of degree \(d+k\le q\).
A homogeneous polynomial of degree at most \(q\) that
vanishes on all of this affine space is zero.

Indeed
reduction by the relations \(X_i^q=X_i\) gives a
polynomial of degree less than \(q\) in each variable,
and such polynomials are tested injectively by evaluation
on the grid. In degree \(q\) only the pure powers
\(X_i^q\) have to be changed, and their reductions of
degree one cannot cancel the other homogeneous terms.
The block determinant is a nonzero polynomial. Thus the
evaluations at invertible block matrices span the dual of
the space of degree-\(d\) polynomials over \(\mathbb F_q\).
Lift the matrices to \(L_W\), apply Nakayama to the finite
free coefficient space over \(W\), and extend along the
embedding.

For \(|\alpha|=d\), put
\(u_\alpha=p^{k^2eN_{\rm orb}}/\binom{d}{\alpha}\). These weights are
integral. By the factorial valuation formula each of the
\(eN_{\rm orb}\) relevant floors contributes less than \(k^2\) to
\(v_p\binom d\alpha\). Use \(u_\alpha\) to pair the
coefficients of polynomials in the generic matrix
variables. Let \(F^\sharp\in\mathcal S_d\) be the operator
whose coefficient at the occurrence matrix \(\alpha\) is
\(u_\alpha[X^\alpha]F\).
For a matrix \(A\) set \(\ell_A(X)=\sum_{i,j}A_{ij}X_{ij}\), and
let \(A_{\rm block}\) be its block-diagonal part. Writing \(E_\alpha\)
for the occurrence operators and using
\(u_\alpha\binom d\alpha=p^{k^2eN_{\rm orb}}\) gives
\[
 (\ell_A^d)^\sharp=p^{k^2eN_{\rm orb}}\sum_\alpha A^\alpha E_\alpha
                 =p^{k^2eN_{\rm orb}}(A_{\rm block}^{\mathsf T})^{\otimes d}.
\]
Since \(g\in L_W\) is block diagonal,
\(\ell_A(g^{-1}Xg)=\ell_{g^{-\mathsf T}Ag^{\mathsf T}}(X)\), and
taking the block-diagonal part commutes with this change of \(A\).
Thus the identity
\[
           \bigl(F(g^{-1}Xg)\bigr)^\sharp
            =g^{\otimes d}F^\sharp(g^{-1})^{\otimes d}
\]
holds on the \(d\)-th powers of linear forms. These span in
characteristic zero by polarization. Both sides are integral on
integral polynomials, so equality over the fraction field proves
the stated integral identity.
For a tensor index \(j=(j_1,\ldots,j_d)\) let its weight
be \(\sum_l\omega_{j_l}\). Let \(U_+\) be the submodule
generated under \(\mathcal S_d\) by all tensors of
weight \(>t_0\). It is finitely generated and torsion-free,
hence free over the valuation ring. It is invariant
under both actions in the last conjugation identity.

Define \(I_d(F)\) to be the trace of \(F^\sharp\) on
\(U_d\) minus its trace on \(U_+\). This is integral and
conjugation invariant. If \(F=X^\alpha\) is bad, its
operator is zero on tensor indices of weight at most
\(t_0\), because the first indices of \(X^\alpha\) are the
column indices of \(F^\sharp\). On the higher-weight basis tensors its image lies in
\(U_+\), since \(F^\sharp\in\mathcal S_d\) preserves \(U_+\).
Thus the two traces agree by the cyclic trace
identity for a map factoring through \(U_+\).
Thus \(I_d\) kills bad monomials.

Suppose now all the bad coefficients of every conjugate
of \(F\) are in \(\mathfrak b\). For a high-weight basis
tensor \(v\), write \(E_\alpha\) for the occurrence operator of
\(\alpha\). Then
\[
 F^\sharp v=\sum_{\alpha:\,E_\alpha v\ne0}
       [X^\alpha]F\,u_\alpha E_\alpha v \in \mathfrak b U_+.
\]
Indeed the input indices in each nonzero term have the weight of
\(v\), hence \(\alpha\) is bad. Also \(E_\alpha\in\mathcal S_d\),
so \(E_\alpha v\in U_+\) before reduction.
For a translate \(g^{\otimes d}v\) the same assertion
follows by the conjugation identity. Such translates
generate \(U_+\) by the integral Schur spanning just
proved. Hence \(F^\sharp U_+\subseteq\mathfrak b U_+\),
so its trace there belongs to \(\mathfrak b\).
Finally the trace on \(U_d\) of \((X^\alpha)^\sharp\)
is zero unless \(\alpha\) is diagonal. For a diagonal
\(\alpha\) there are \(\binom d\alpha\) identical diagonal
coefficients, so that trace is \(p^{k^2eN_{\rm orb}}\).
Thus the trace on \(U_d\) is \(p^{k^2eN_{\rm orb}}\) times
evaluation at 1 in degree \(d\). Sum \(I_d\) on the homogeneous
components for \(0\le d\le D\) and set the functional equal to zero in
higher degrees. Conjugation preserves degree.

Conjugation invariance also kills every monomial with
different row and column weights. Indeed conjugation by
the diagonal matrix with entries \((1+p)^{\omega_i}\)
scales that monomial by a nontrivial integer power of
\(1+p\) when the weights differ. The coefficient ring
has characteristic zero and no zero divisors.

\subsection{The simultaneous average on the coefficient space}
\label{sec:exceptional-edge-average}

\sectionopening{In this subsection, we tensor the functionals of
Section~\ref{sec:exceptional-single-average} across the embeddings in
the local representation. The estimates of
Sections~\ref{sec:exceptional-root-defect} and
\ref{sec:exceptional-defect-packing} yield both a pairing congruence
and a dimension bound for its range after reduction.}

Keep \(N_{\rm orb}=2h\), \(P=p^e\), \(q=P^{N_{\rm orb}}\), and \(r\) as above.
Put \(c_0=k^2eN_{\rm orb}\), \(h_0=N_{\rm orb}c_0\), \(c_1=5h_0\), and
\(t=3kc_1\). Write \(\rho_b=1_{b=r}-1_{b=r-1}\).
For a fixed copy \(i\), give a variable in row \(b\) and
column \(c\) of the generic matrix in copy \(j\) weight
\[
 \omega_i(j,b,c)=P^{-h}
 \begin{cases}
  \rho_b, & j=i,\\
  q^{-1}P^{(j-i)\bmod N_{\rm orb}}(\rho_b-\rho_c), & j\ne i.
 \end{cases}
\]
Denote the weight of a monomial in all copies by
\(\omega_i(\nu)\). Write \(\operatorname{def}_i\) for the
defect at an embedding at the boundary root, in the convention
of the rank-two calculation.

Take one generic \(k\)-square matrix \(X_j\) for each copy.
At embedding \(s\) use the product of the entrywise
\(P^d\)-th powers of \(X_{s-d}\), in decreasing order of
\(d\), with \(0\le d<h\) for a shallow row and \(0\le d\le h\)
for a deep row. In odd embeddings take the cofactor matrix
\(\operatorname{adj}(M)^{\mathsf T}\) of the corresponding product
\(M\) before extracting the row. Assemble these rows into a matrix
at each embedding. Let \(\mathcal B\) be the tensor product of the
degree-\(a\) symmetric powers of these matrices. Its evaluation at \(X_j=1\) for all \(j\)
is \(1\). For a change of basis \(g\), replace the factor
in embedding \(s\) at position \(d\) by
\[
 g_s\bigl(g_{s-d}^{-1}X_{s-d}g_{s-d}\bigr)^{[P^d]}g_s^{-1},
              \qquad g_s=\iota_s(g).
\]
Brackets denote entrywise power. Call the resulting
substitution \(\mathcal B_g\). Modulo \(p\) these
conjugations cancel, because the \(P^d\)-power map
is a ring map there and takes \(g_{s-d}\) to \(g_s\).
The difference from the original factor has all
coefficients divisible by \(p=w^2\).

We give the weight estimate on \(\mathcal B_g\).
Fix \(|\epsilon|\le1/2\) and give the generic variables
weights \(\epsilon\omega_i\), adding the \(w\)-order to
the weight of each monomial coefficient. Every monomial
of any single entry of the factor at \(d\) has degree
\(P^d\) and weight between \(-1\) and \(1\), since
\(|\epsilon\omega_i(j,b,c)|\le P^{-h}\) and \(d\le h\).
Adding \(w\)-order two therefore makes every error term
weigh at least as much as the original entry monomial.
In particular each entry has weight at least
\(\epsilon P^d\omega_i(s-d,b,c)\).

For the matrix product put
\[
 A_{s,d}(b)=\frac{\epsilon\rho_b}{P^h}
                       \frac{P^{d+(s-d-i)\bmod N_{\rm orb}}}q,
                 \qquad u=\frac1{2P^h}.
\]
At \(d=0\), \(|A_{s,0}(b)|\le u\). For each \(d\le h\)
one has
\[
 \left|A_{s,d+1}(b)-A_{s,d}(c)
                   -\epsilon P^d\omega_i(s-d,b,c)\right|
                   \le u .
\]
Write \([v]\) for the residue of \(v\) in
\(\{0,\ldots,N_{\rm orb}-1\}\) and put \(j=[s-d-i]\).
For \(j>0\) one has \([j-1]=j-1\), whereas for \(j=0\)
one has \(P^{1+[j-1]}=q\). Substitution gives
\[
 A_{s,d+1}(b)-A_{s,d}(c)-\epsilon P^d\omega_i(s-d,b,c)
     =-\mathbf1_{j=0}\,\epsilon\rho_c P^{d-h}/q .
\]
Its absolute value is at most \(u\), including at \(d=h\).
Along a path in a product of \(m\le h+1\) factors,
the potentials telescope, so all entries in row
\(b\) weigh at least \(A_{s,m}(b)-(h+2)u\).
In an odd embedding a cofactor omits its own row and
takes one entry in each other row. The potentials
sum to zero over the rows. Thus in both parities the
row bound is
\(\varepsilon_s A_{s,m}(b)-k(h+2)u\), with
\(\varepsilon_s=1\) for even embeddings and \(-1\)
otherwise.

For \(m=h+1\) on deep rows and \(h\) on shallow rows,
the term \(\varepsilon_s A_{s,m}(b)\) is exactly
\(\epsilon(1-1/q)\) times the coefficient of the
variable \((b,s)\) in the defect at \(i+h\).
For shallow rows the power simplifies to
\(P^{(s-i-h)\bmod N_{\rm orb}}\) after removing \(P^h\).
For deep rows it is the same power except at
\(s=i+h\), when it is \(q\), giving the additional
deep term of the defect. The signs are those of the
two root variables, and the other potentials vanish.
Taking the symmetric powers and then the tensor
product adds these weights according to the row
exponents \(\beta\). Hence, for every monomial \(\nu\)
in an entry of row \(\beta\) of \(\mathcal B_g\), that
coefficient is divisible by \(w\) to the nonnegative
ceiling of
\[
 \epsilon\bigl((1-1/q)\operatorname{def}_{i+h}(\beta)
                     -\omega_i(\nu)\bigr)-N_{\rm orb}a k(h+2)u.
\]
The row estimate uses exactly \(a\) factors in each
embedding, one row factor for each variable counted
with its exponent.

The size hypotheses give \(N_{\rm orb}a k(h+2)+a+k<P^h\).
Indeed \(P^h\ge a^2\), \(h+2\le N_{\rm orb}k\), and
\(a\ge(N_{\rm orb}k)^5\ge(N_{\rm orb}k)^2+2\), with \(k<a\).
Thus the loss \(N_{\rm orb}a k(h+2)u\) is at most \(1/2\).
The sum of the cyclic geometric weights of the
defect is \(1/(P-1)\), and each embedding has degree
\(a\), so
\(|\operatorname{def}_{i+h}(\beta)|/q\le2a/q\le1\).
Taking \(\epsilon=1/2,-1/2\) shows that if the
coefficient survives modulo \(w^{c_1}\), then
\[
 |\operatorname{def}_{i+h}(\beta)-\omega_i(\nu)|
                       \le 2c_1+2\le t.
\]

The degree in each generic copy is at most
\(D=N_{\rm orb}a k(h+2)P^h\). This follows by adding the
degrees of the factors in each product, using at
most \(k\) factors of product entries for a cofactor
and \(a\) row factors in each embedding. The same bound
holds in any basis. Our strict size inequality gives
\(D+k<q\), so the one-copy average applies.
At each copy \(j\) use its embedding and row weights
\(\rho_b\), with bad threshold \(2tP^h\).
Let \(I_j\) be the resulting invariant functional,
and let \(I_j^0\) be \(p^{c_0}\) times evaluation at
the identity. Both are invariant under the Levi.

Apply some such invariant functionals in all
copies but \(j\), leaving a polynomial in \(X_j\).
If a coefficient of it survives modulo \(w^{c_1}\),
some original term survives and its other-copy
monomials are not killed by the chosen functionals.
Each has equal row and column weights by invariance
under diagonal conjugation. Consequently for the
surviving term \(\omega_j\) is just its row weight
in copy \(j\), divided by \(P^h\).
The last weight bound has two consequences.
A defect above \(3t\) at \(j+h\) forces the row weight in this copy above \(2tP^h\), and row weight above \(2tP^h\) forces defect above \(t\).

For \(g\) in the Levi, write \(C_g\) for the substitution
\(X_s\mapsto g_s^{-1}X_sg_s\) in every copy, and put
\(H_s=g_s\) for even \(s\), \(H_s=g_s^{-\mathsf T}\) for odd \(s\).
In each product defining \(\mathcal B_g\) the adjacent factors
\(g_s^{-1}g_s\) cancel. At an odd embedding we use the polynomial identity
\[
 \operatorname{adj}(g_s A g_s^{-1})^{\mathsf T}
     =g_s^{-\mathsf T}\operatorname{adj}(A)^{\mathsf T}g_s^{\mathsf T}.
\]
The two products of lengths \(h,h+1\) contribute to different
depth blocks. Projection onto either block commutes with \(H_s\).
Taking the symmetric powers and their tensor product therefore gives
the identity over the integral polynomial ring
\[
 C_g\mathcal B=\rho_{\rm let}(g)^{-1}\mathcal B_g\rho_{\rm let}(g).
\]
Fix \(j\) and let \(l_s\) be any of the invariant functionals above
for \(s\ne j\). For an integral coefficient vector \(v\), put
\[
 F_g(X_j)=\left(\bigotimes_{s\ne j}l_s\right)
                (\mathcal B_g)\,\rho_{\rm let}(g)v.
\]
Applying these functionals to the polynomial identity leaves
\[
 F_1(g_j^{-1}X_jg_j)=\rho_{\rm let}(g)^{-1}F_g(X_j).
\]

Write \(\mathcal I=\bigotimes_s I_s\) and \(V=\mathcal I(\mathcal B)\),
applying the functional entrywise. Invariance in every copy and the
integral basis identity give, for every \(g\in L_W\),
\[
 \rho_{\rm let}(g)Vv=\mathcal I(\mathcal B_g)\rho_{\rm let}(g)v.
\]
For a row of defect \(>3t\) at \(i\), contract first all copies
except \(j=i-h\). Modulo \(w^{c_1}\) every remaining monomial
in that row is bad in copy \(j\) by the first weight implication.
The functional \(I_j\) kills those monomials, and its integrality
keeps all other coefficients in \(w^{c_1}\mathcal O\). Thus the row
of every Levi translate of \(Vv\) vanishes modulo \(w^{c_1}\).

For the pairing, replace \(I_j^0\) by \(I_j\) one
copy at a time. Fix the other copies and a saturated
halfway vector \(v'\). Pair the uncontracted
polynomial with \(v'\) using the diagonal weights.
Its degree is at most \(D\). Write this polynomial as
\(\Psi(X_j)=\sum_\beta b_\beta v'_\beta F_{1,\beta}(X_j)\).
The adjunction of \(\rho_{\rm let}(g^{-1})\) with
\(\rho_{\rm let}(g^{-\mathsf T})\) for the diagonal pairing gives
\[
 [X_j^\alpha]\Psi(g_j^{-1}X_jg_j)
 =\sum_\beta [X_j^\alpha]F_{g,\beta}(X_j)\,
       b_\beta\bigl(\rho_{\rm let}(g^{-\mathsf T})v'\bigr)_\beta .
\]
For a bad \(\alpha\), any first factor surviving modulo \(w^{c_1}\)
forces the defect of \(\beta\) to exceed \(t=3kc_1\).
Since \(g^{-\mathsf T}\) belongs to the same Levi, the rank-two estimate
of Section~\ref{sec:exceptional-root-defect} makes the second factor
divisible by \(w^{c_1}\). Thus all the bad
coefficients of every conjugate of the pairing
polynomial are in \((w^{c_1})\). The congruence
property of \(I_j\) shows that replacing \(I_j^0\)
does not change the pairing modulo \(w^{c_1}\).
With all identity evaluations the contraction
is \(p^{N_{\rm orb}c_0}1=p^{h_0}1\). This proves the desired
pairing congruence for \(V\).

Finally let \(C\) consist of the integral coefficient vectors
whose \(L_W\)-translates have coefficients of defect \(>3t\)
divisible by \(w^{c_1}\). The rank-two estimate puts the
saturated halfway restrictions in \(C\), and the full
contraction takes values there. Suppose an integral \(v\)
satisfies \(w^{c_1-1}v\in C\). Lift any residue Levi element
\(\bar g\) to \(g\in L_W\) by lifting its two invertible blocks
with unit determinants. For each coordinate \(\alpha\) of
defect \(>3t\) at some embedding,
\[
 w^{c_1-1}(\rho_{\rm let}(g)v)_\alpha\in w^{c_1}\mathcal O
   \quad\Longrightarrow\quad
 (\rho_{\rm let}(g)v)_\alpha\in w\mathcal O .
\]
Reduction commutes with the symmetric powers and inverse
transposes of the Levi blocks. Since \(v_p(w)>0\), all
these coordinates of \(\bar g\bar v\) vanish in the residue
field, so \(\bar v\in S\).
The flag count uses the inequalities:
\(t\ge2(k+1)\) and
\[
 8N_{\rm orb}ktP^{2r}
    =120N_{\rm orb}^2 k^4 eN_{\rm orb}\,P^{2r}
    \le (N_{\rm orb}k\log q)^5P^{3r}\le \theta a .
\]
For the first inequality use \(eN_{\rm orb}\le\log q\),
\(N_{\rm orb},k\ge2\), and \(\log q\ge40\log2>27\), which
follows from \(q\ge a^4\) and \(a\ge(N_{\rm orb}k)^5\).
The final bound follows from the large-degree
hypothesis and \(0<\theta\le1\).
The flag estimate therefore gives
\(s_{\mathcal P}\dim S\le(1+\theta)d_{\rm let}\).
This supplies both the congruence and the residue
dimension used in the edge Smith argument.

\subsection{The colength of the source Hecke operator}

\sectionopening{In this subsection, we combine the two restriction estimates of
Section~\ref{sec:exceptional-source-depth-count} with the congruence
of Section~\ref{sec:exceptional-edge-average}. Smith normal form gives
an explicit bound on the length of the quotient by the image of the
Hecke operator and the prescribed power of \(p\).}

The depth count of Section~\ref{sec:exceptional-source-depth-count}
bounds the sum of the dimensions of the residue spaces
above all positive depths by \(C_{\rm div}d_{\rm let}\).
Use \(C,V,S\) from Section~\ref{sec:exceptional-edge-average}.
In particular \(s_{\mathcal P}\dim S\le(1+\theta)d_{\rm let}\),
\(C\) contains the saturated halfway restrictions, and
\(w^{c_1-1}u\in C\) implies \(\bar u\in S\). The pairing
of \(V\) with a saturated halfway restriction agrees with
that of \(p^{h_0}1\) modulo \(w^{c_1}\).

Here is the linear algebra that uses those estimates.
On an edge the two restrictions are the origin value
acted on by \(W_0\rho_{\rm let}(g_e^{\mathsf T})\) and the
terminus value acted on by
\(W_0\rho_{\rm let}(\beta_b^{-1})\), respectively.
At either endpoint every class is incident to an edge, by the forward
coset description and the dual-edge bijection. On a fixed edge the
restriction at that class is invertible over the coefficient field,
since \(W_0\) and each change of coordinates are invertible there.
Both restrictions therefore have column rank \(d_X\).
Factor them as \(R_1=Q_1D_1\), \(R_2=Q_2D_2\), with \(D_i\) square
and the saturated bases \(Q_i\) admitting integral left inverses \(L_i\).
Write \(\mathcal E\) for the integral module of edge values. Since
\(L_iQ_i=1\), an element \(Q_i c\) of the fraction-field image is
in \(\mathcal E\) exactly when \(c=L_iQ_i c\) is integral. Thus
\(Q_i\mathcal O^{d_X}\) is the saturation of \(R_i\mathcal O^{d_X}\) in \(\mathcal E\).
Put \(g_{1,e}=g_e^{\mathsf T}\), \(g_{2,e}=\beta_b^{-1}\). For either
endpoint, a vector \(Q_i c=R_iD_i^{-1}c\) of the saturation has,
at the edge \(e\), active component \(W_0\rho_{\rm let}(g_{i,e})F\),
where \(F\) is the corresponding component of \(D_i^{-1}c\) at that class.
Fix the same coordinate in the remaining tensor factors. Integrality
of the saturated column in \(\mathcal E\) gives
\(W_0\rho_{\rm let}(g_{i,e'})F\in\mathcal O^{d_{\rm let}}\) for
every edge \(e'\) incident to this class, with this same \(F\).
As \(e'\) runs through those edges, the \(g_{i,e'}\) represent
\(\mathcal P\backslash\mathrm{GL}_k(W)\), since the cosets
\(g_{i,e'}^{-1}\mathcal P\) exhaust \(\mathrm{GL}_k(W)/\mathcal P\).
For \(g\in\mathrm{GL}_k(W)\), write \(gg_{i,e}=u g_{i,e'}\) with
\(u\in\mathcal P\). Then
\[
 W_0\rho_{\rm let}(gg_{i,e})F
 =\bigl(W_0\rho_{\rm let}(u)W_0^{-1}\bigr)
      W_0\rho_{\rm let}(g_{i,e'})F
      \in\mathcal O^{d_{\rm let}}.
\]
At the fixed edge \(e\), multiply
\(\rho_{\rm let}(g_{i,e})F\) by a sufficiently large power of
\(w\) to clear its finitely many denominators. The displayed
identity already gives the integral tested values of the
unscaled vector in every direction \(g\). Thus the integral
edge values of each \(Q_i\) are saturated halfway restrictions
for the common local test in Section~\ref{sec:exceptional-edge-average}.

Let \(B=\operatorname{diag}(b_\mu)\) at each edge and \(B_X\) be the
same diagonal matrix repeated at the classes. Write \(U=\bigoplus_e U_e\)
for the transitions on the remaining tensor factors, acting as the
identity on the active factors, and set \(A_0=Q_2^{\mathsf T}U B\).
On the tensor product of the active factors the identity
\[
 (W_0\rho_{\rm let}(\beta_b^{-1}))^{\mathsf T}B
      W_0\rho_{\rm let}(g_e^{\mathsf T})
 =\bigl(\rho_{\rm let}(g_e)W_0^2\rho_{\rm let}(\beta_b^{-1})\bigr)^{\mathsf T}B
\]
follows from \(B\rho_{\rm let}(g_e^{\mathsf T})=\rho_{\rm let}(g_e)^{\mathsf T}B\).
Here \(W_0^2=\operatorname{diag}(p^{\operatorname{dep}(\mu)})\).
Tensoring with \(U_e\) and summing incoming edges gives, in the
original class coordinates,
\[
 D_2^{\mathsf T}A_0 Q_1D_1=R_2^{\mathsf T}UBR_1=Y_\lambda B_X.
\]
Its image is therefore contained in \(Y_\lambda M\) after scalar extension.
To bound \(v(\det D_i)\), write \(t_j\) for the Smith factors of
one restriction, and \(v_j\) for the corresponding input basis. All
parahoric cosets occur at that endpoint by the two edge descriptions
above. Thus \(t_j\) divides every coefficient of
\(W_0\rho_{\rm let}(g)v_j\) in each active factor at every class,
for every \(g\in\mathrm{GL}_k(W)\), by the parahoric invariance of
the weighted lattice proved above. The remaining tensor factors
change by integral invertible matrices.

If \(2v(t_j)>\delta\) for an integer \(\delta\ge0\), every
coefficient of depth at most \(\delta\) in the unweighted translates
of \(v_j\) vanishes modulo the maximal ideal. Thus \(\bar v_j\)
belongs to the direct sum of \(d_X/d_{\rm let}\) copies of
\(D_{\delta+1}\). These reduced basis columns are independent.
Since \(v(p)=1\) and \(w^2=p\), we obtain
\[
 \begin{aligned}
 2v(\det D_i)
 &=\sum_j 2v(t_j)
 \le\sum_{\delta\ge0}\#\{j:2v(t_j)>\delta\}\\
 &\le\frac{d_X}{d_{\rm let}}\sum_{\delta\ge0}\dim D_{\delta+1}
 \le C_{\rm div}d_X .
 \end{aligned}
\]
Applying this at both endpoints gives
\(v(\det D_1)+v(\det D_2)\le C_{\rm div}d_X\).

Let \(B'=\operatorname{diag}(b'_\mu)\) and let \(\mathcal V\) act
by \(V\) in every active edge factor. Since \(BB'=p^{h_0}1\) and
\(L_2Q_2=1\), we have \(A_0 B'U^{-1}L_2^{\mathsf T}=p^{h_0}1\).
The matrix \(B_0=\mathcal V B'U^{-1}L_2^{\mathsf T}\) has edge
values in \(C\). The matrix \(U\) acts on the remaining tensor
factors and commutes with \(B,\mathcal V,B'\). The pairing congruence, tested against
the active edge values of the columns of \(Q_2\), gives
\[
 \begin{aligned}
 A_0B_0
 &=Q_2^{\mathsf T}B\mathcal V B'L_2^{\mathsf T}\\
 &\equiv p^{h_0}Q_2^{\mathsf T}BB'L_2^{\mathsf T}
  =p^{2h_0}1\pmod{w^{c_1}}.
 \end{aligned}
\]
The values of \(Q_1\) also belong to \(C\).

Take Smith form for the image of
\(B_0\) in the complement of \(Q_1\) given by \(L_1\).
The projection \(1-Q_1L_1\) preserves \(C\), since integral
linear combinations of the columns of \(Q_1\) belong to \(C\).
Keep those nonzero divisors \(t_j\) dividing
\(w^{c_1-1}\), and write \(u_j\) for their primitive columns.
Then \(w^{c_1-1}u_j\) lies in \(C\) on each edge, so
each \(u_j\) reduces into \(S\) at every edge.
The columns of \(Q_1\) also reduce into \(S\).
These two lists together are independent modulo the maximal
ideal, since the primitive columns are basis vectors in the
complement. There are \(s_{\mathcal P}d_X/d_{\rm let}\)
copies of the residue space at the edges, including the multiplicities
from the remaining tensor factors. Thus
\[
 d_X+\#\{u_j\}\le
 \frac{s_{\mathcal P}d_X}{d_{\rm let}}\dim S\le(1+\theta)d_X.
\] The other Smith terms are all divisible by
\(w^{c_1-1}\).

Let \(G\) be the sum of the image of \(A_0Q_1\) and the
span of the kept \(A_0u_j\). The near-inverse congruence
gives \(p^{2h_0}\mathcal O^{d_X}\subset
G+w^{c_1-1}\mathcal O^{d_X}\).
Since \(c_1=5h_0\) and \(h_0\ge2\),
\(w^{c_1-1}=p^{2h_0} w^{h_0-1}\).
The image of \(p^{2h_0}\mathcal O^{d_X}\) in \(\mathcal O^{d_X}/G\)
is finitely generated and equals its multiple by \(w^{h_0-1}\).
Nakayama's lemma therefore gives \(p^{2h_0}\mathcal O^{d_X}\subseteq G\).
For finitely generated lattices \(L'\subseteq L''\) of the same
rank over the coefficient valuation ring, successive elimination
at an entry of least valuation diagonalizes the inclusion. Thus
\(\operatorname{length}_p(L''/L')=\sum_i v(t_i)\), where
\(t_i\) are the elementary divisors and \(v(p)=1\).
This normalized length is additive along inclusions by multiplicativity
of the determinants of the basis changes. Write \(L=\mathcal O^{d_X}\),
\(z=p^{j_p}\), and \(\Phi=A_0Q_1\). The truncated sums of
images used here are finitely generated and contain \(zL\).

The containment \(p^{2h_0}L\subset G\)
bounds the normalized length of \(L/(G+zL)\) by
\(2h_0d_X\). The quotient between \(\Phi L+zL\) and
\(G+zL\) is generated by the kept columns, each killed
by \(z\), so its normalized length is at most
\(j_p\#\{u_j\}\). Replacing \(\Phi L\) by \(\Phi D_1L\)
adds at most \(v(\det D_1)\), since the intervening
quotient is a quotient of \(L/D_1L\).

Put \(T=\Phi D_1L+zL\). Since \(D_2^{\mathsf T}(zL)\subseteq zL\),
\(D_2^{\mathsf T}T\subseteq D_2^{\mathsf T}\Phi D_1L+zL\).
The injection \(D_2^{\mathsf T}\) identifies \(L/T\) with
\(D_2^{\mathsf T}L/D_2^{\mathsf T}T\).
The inclusion \(zD_2^{\mathsf T}L\subseteq D_2^{\mathsf T}T\)
gives the full rank needed for the normalized lengths.
Additivity along
\(D_2^{\mathsf T}T\subseteq D_2^{\mathsf T}L\subseteq L\) gives
\[
 \operatorname{length}_p\bigl(L/(D_2^{\mathsf T}\Phi D_1L+zL)\bigr)
 \le v(\det D_2)+\operatorname{length}_p(L/T).
\]
The sum of these endpoint valuations is at most
\(C_{\rm div}d_X\) and \(h_0\le C_{\rm div}\).
This proves
\[
  \operatorname{length}_p\!\bigl(M/(Y_\lambda M+p^{j_p}M)\bigr)
           \le (3C_{\rm div}+\theta j_p)d_X .
\]
In a class basis of \(M\), let
\(H=[Y_\lambda\mid p^{j_p}I_{d_X}]\). Its maximal minors
generate the zeroth Fitting ideal of the quotient. Smith form
therefore gives its normalized length, with \(v_p(0)=+\infty\), as
\[
 \min_{\#J=d_X}v_p(\det H_J)
   =\sum_{i=1}^{d_X}\min\{v_p(s_i(Y_\lambda)),j_p\},
\]
where the \(s_i\) include the zero Smith factors. Indeed, if \(P Y_\lambda Q=\operatorname{diag}(s_i)\) is a Smith form,
left multiplication by \(P\) and right multiplication by
\(\operatorname{diag}(Q,P^{-1})\) preserve the ideal of maximal minors
and keep the second block equal to \(p^{j_p}I\). In these coordinates
each nonzero maximal minor selects one entry from each row.
Selecting an entry of lesser valuation in every row realizes the sum
above, including rows with \(s_i=0\). In the original class basis the
minimum ranges over all column sets and can depend on the coefficient
place. The valuations of its minors are unchanged after any extension
of the coefficient field with
\(v_p(p)=1\). Hence the bound applies to the original local module.

Finally \(p^r\ge m\), \(r\ge1\), and \(P\ge p\) give
\(P^{2r}\ge mp\). The size bound \(r^2m\le k\) gives
\(r^2/(kP)\le1/(mp)\). Hence
\[
 3 C_{\rm div}\le \frac{6N_{\rm orb}a}{mp}=\frac{12 A_\lambda}{mp}.
\]
\begin{proposition}[Colength of the source Hecke operator]\label{prop:source-hecke-colength}
For each \(\lambda\in\Lambda\), \(Y_\lambda\) preserves the original
source lattice \(M\). With \(j_p=\lceil 2mfa/p\rceil\), its normalized
colength modulo \(p^{j_p}\) satisfies
\[
 \operatorname{length}_p\bigl(M/(Y_\lambda M+p^{j_p}M)\bigr)
 \le d_X\left(\frac{12A_\lambda}{mp}+\frac{j_p}{m^2 f}\right).
\]
\end{proposition}

\section{Norm and index estimates for the translated images}
\label{sec:integral-poincare-bounds}

\sectionopening{In this section, we combine the derivative estimates of
Sections~\ref{sec:local-theta-derivatives} and~\ref{sec:exceptional-hecke-analysis}
with the calculations at the level primes. The common scalar of
Proposition~\ref{prop:shortness-normalization} gives the upper estimates
on the family of Section~\ref{sec:common-filling-assembly}.}

\subsection{Derivatives at the ramified prime}
\label{sec:marked-ramified-ordinary-jet}

\sectionopening{In this subsection, we specialize the integral CM calculation of
Section~\ref{sec:marked-integral-sections-ordinary-jets} to \(p=d\).
The ramified uniformizer changes the exponent in the estimate of
Section~\ref{sec:marked-column-gain}, including the contribution at
exponent zero.}

Put \(p=d\), \(\pi=\sqrt{-d}\), and normalize \(v(p)=1\).
Here \(v(\pi)=1/2\), \(p<x_1\), and \(p\nmid MNq_1\).
Fix any allowed twist \(\xi\) with \(p\nmid q_\xi\).
By Section~\ref{sec:seesaw-integral-tensor}, the translations in the
evaluated tensor can be supported on the primes dividing \(MNq_\xi\).
The signed action of an integral displacement fixes the local vacuum.
Thus each arrow has the standard level-one datum at \(p\).
At the actual scale \(\pi_M s_{\rm jet}\) the algebraic multiplier in
the normalized Taylor variables is \(\iota_p(\pi_Mp_F)\).
Section~\ref{sec:marked-integral-sections-ordinary-jets} gives
\(v_p(\iota_p(p_F))=0\), since \(p\nmid q_1\) and \(\alpha_{\rm W}\)
is integral invertible above \(p\). Also \(\pi_M\bar\pi_M=M\) with
\(p\nmid M\), so this multiplier is a unit at either coefficient
embedding.

On the CM curve with its chosen half-period, put \(g=-\log\theta\). Average
\(y(az)\) and \(g(az)\) over \(a^{p-1}=1\), with weights
\(a^{-1}\) and \(a^{-2}\), respectively, and denote the results
by \(T\) and \(g^\flat\). At the ramified place these roots lie
in \(\mathbb Q_p\), so their CM conjugates are themselves.
The integral CM pullbacks give \(y(az)\in y\mathcal O[[y]]\) for
each such \(a\), with linear term \(a y\). Thus
\(T=(p-1)^{-1}\sum_a a^{-1}y(az)=y+O(y^2)\) over \(\mathcal O\).
Since \(a\bar a=a^2\), the corresponding pullback quotient is
\[
 Q_a(z)=\exp\bigl(a^2g(z)-g(az)\bigr)\in1+y\mathcal O[[y]].
\]
Substituting \(y(az)\) for \(y\) in the quotient at \(a^{-1}\) gives
\[
 \prod_{a^{p-1}=1}Q_{a^{-1}}(az)
    =\exp\bigl((p-1)(g^\flat(z)-g(z))\bigr)\in1+y\mathcal O[[y]].
\]
In each positive degree the equation for its constant-one
\((p-1)\)-st root has leading coefficient \(p-1\in\mathcal O^\times\).
It follows recursively that \(\exp(g^\flat-g)\in\mathcal O[[y]]
=\mathcal O[[T]]\). The torsion calculation in
Section~\ref{sec:marked-integral-sections-ordinary-jets} bounds the
order of the reduced parameter under \([\pi]\) by \(\pi\bar\pi=p\).
The linear coefficient vanishes on reduction and all allowed degrees
are \(1\pmod{p-1}\). The reduced order is therefore \(p\), giving
\[
 T(\pi z)=\pi T(z)+b_1T(z)^p+
             \sum_{j\ge2}b_jT(z)^{1+j(p-1)},\qquad b_1\in\mathcal O^\times.
\]
Choose \(\varpi^{p-1}=-\pi/b_1\) in \(\mathbb C_p\) and write
\(Y=T(\varpi z)/\varpi\). Thus \(v(\varpi)=1/(2(p-1))\).
If the chosen root is multiplied by \(\zeta\in\mu_{p-1}\), then
\(T(\zeta z)=\zeta T(z)\) and
\(g^\flat(\zeta z)=\zeta^2g^\flat(z)=\zeta^{p+1}g^\flat(z)\).
Both \(Y\) and \(g^\flat(\varpi z)/\varpi^{p+1}\) are therefore
independent of the root. The equation
\(Y(\pi z)/\pi=Y-Y^p+\pi E(Y)\), with \(E\in U^2\mathcal O[[U]]\),
gives
\[
 Y\in\mathcal O[[z]],\qquad Y\equiv\sum_{j\ge0}z^{p^j}\pmod\pi.
\]
Indeed the recursion divides in degree \(n>1\) only by the unit
\(1-\pi^{n-1}\), and the reduction is the unique zero-constant
solution of \(Y-Y^p=z\).

The pullback \(\theta(\pi z)\) has degree \(\pi\bar\pi=p\).
The elliptic quotient by \(\theta(z)^p\) is even and has its
only possible pole at the theta zero, of order at most \(p\).
The integral interpolation of
Section~\ref{sec:marked-integral-sections-ordinary-jets} expresses it
as \(P_0(x)+yP_1(x)\). Comparing \(z,-z\) and using
\(y(z)-y(-z)=x'(z)\ne0\) gives \(P_1=0\). Since the pole
of \(x\) has order two and \(p\) is odd, we obtain
\[
 \exp(pg-g(\pi z))=R(x),\qquad R\in\mathcal O[X],\qquad
 R(0)=1,\qquad 2\deg R\le p-1 .
\]
For \(2\le n<p\), \(\pi[z^n]g\) is integral. In
degree two use the two CM pullbacks at \(\beta,\bar\beta\)
and subtract their quadratic coefficients. For \(n>2\)
use \(g''=g''(0)+x\), since the factors \(n(n-1)\) are
units and \(x\in\mathcal O[[T]]\) with \(T=z+O(z^p)\).
It follows that \(pg-g(\pi z)\) belongs to \(\pi\mathcal O[[z]]\)
through degree \(p-1\), since \(p=-\pi^2\).

Exponentiation
preserves the congruence through those degrees because
\(k/2-v(k!)\ge1/2\) for \(k\ge1\). Triangular comparison
using \(x=z^2+\cdots\) gives \(R-1\in\pi\mathcal O[X]\).
Taking the formal logarithm gives \(pg-g(\pi z)\in\pi\mathcal O[[T]]\).
The component of degree \(2\) modulo \(p-1\) thus satisfies
\[
 pg^\flat(z)-g^\flat(\pi z)
 =\pi(aT^2+T^{p+1}H),\qquad a\in\mathcal O,\quad H\in\mathcal O[[T]].
\]
With \(c_\varpi=\varpi^{p+1}/p\) and
\(\Xi=-g^\flat(\varpi z)/c_\varpi\), the scalar
\(-\pi\varpi^2/(pc_\varpi)\) is a unit. The preceding identity
gives, coefficientwise, \((1-\pi^n/p)\Xi_n\equiv
\lambda[z^n]Y^2\pmod\pi\) for some \(\lambda\in\mathcal O\).

Indeed \(-\pi\varpi^2/(pc_\varpi)=-\pi/\varpi^{p-1}=b_1\);
the higher terms contain an additional factor \(\varpi^{p-1}\in\pi\mathcal O\).
Since \(p=-\pi^2\), the multiplier is \(2\) for \(n=2\)
and \(1+\pi^{n-2}\equiv1\pmod\pi\) for \(n>2\).
The constant and linear terms vanish. Hence
\[
 \Xi\in\mathcal O[[z]],\qquad
 \Xi\equiv\lambda(Y^2-z^2/2)\pmod\pi .
\]
In particular the coupling indices now include \(e=0\).
The coefficients of \((z^{p^e}z)^{p^j}\) are \(2\lambda\)
for \(e>0\), \(\lambda\) for \(e=0,j>0\), and \(\lambda/2\)
for \(e=j=0\).

Use the same resolution of the integral vacuum as for the
exceptional jet. The basis entries are integral above \(p\),
since \(p\nmid q_1\). On \(K\), conjugation is identity
modulo \(\pi\); on the cyclotomic generator, the \(p^e\)-power
is the corresponding Frobenius. Therefore
\[
 \overline{\tau(w)}^{p^e}\equiv \tau_e(w)\pmod\pi,
 \qquad\tau_e=\langle-p^e\rangle\tau.
\]
In fact \(p=d\nmid q_1\), so the completion of the cyclotomic
factor at \(p\) is unramified, and adjoining \(\pi^2=-p\)
gives ramification index two and uniformizer \(\pi\). Thus equality
of residues for these elements of \(E\) is congruence modulo
\(\iota_p(\pi)\mathcal O_{\mathbf C_p}\) after embedding, including
at \(e=0\). Denominators prime to \(p\) preserve it, and either
conjugate of \(\pi\) generates the same ideal. The basis identity
turns the single-curve terms into
\(\Gamma_e=\sum_{\tau,b,a}z_{\tau,b,a}^{p^e}z_{\tau_e,b,a}\).
The congruence error has zero constant term and coefficients
divisible by \(\pi c_\varpi\). Since \(v(\pi c_\varpi)=1/(p-1)\),
Legendre's formula gives
\[
 v\bigl((\pi c_\varpi)^n/n!\bigr)=s_p(n)/(p-1)\ge0 ,
\]
where \(s_p(n)\) is the sum of the base-\(p\) digits. Thus its
formal exponential is integral. Scaling by the actual integral
factor \(\pi_Mp_F\) preserves this form.
We obtain the same matrix normal form as in the exceptional
case, now for all \(e\ge0\) and with
\(v(\varpi)=1/(2(p-1))\).

Since \(\lambda\in\mathcal O\),
a direction that kills the coupling modulo \(\pi\) also kills its
integral multiple \(\lambda\), including \(\lambda=0\).
Use the null-direction proof with no extra exponents or discarded columns and
with each congruence modulus \(p\) there replaced by \(\pi\).
The coefficient ring modulo \(\pi\) has characteristic \(p\),
so the same linear congruences, using \(p^e\)-th roots
of \(\pi\) where needed, work also for \(e=0\). Each
determinant contribution is multiplied by \(v(\pi)=1/2\). In detail,
with \(u_\tau=a\mathbf1_{\Phi_a}(\tau)\) and \(c_\tau,b_\tau\)
as in the column counts, write
\[
 E_{\tau,e}
 =\min(u_\tau/p^e,u_{\tau_e}+c_{\tau_e})
       +c_\tau/p^e+\min(b_\tau/p^e,b_{\tau_e}).
\]
The factorial bound gives the lower bound for every minor jet,
and therefore every Schur entry,
\[
 \frac{B_j}{2(p-1)}
 -\frac12\sum_{\tau,e\ge0}E_{\tau,e}
 -fn_{j+1}^3(1+\log_p\max(B_j,1)).
\]
If \(B_j=0\), the constant term is integral.

Here \(B_j=\sum_\tau(u_\tau+b_\tau)\) since no columns
were erased. Subtracting \(v(\varpi^{B_j})=B_j/(2(p-1))\) from
the factorial contribution \(B_j/(p-1)\) leaves
\(\tfrac12\sum_{e\ge1}B_j/p^e\). For each \(e\ge1\),
\[
 (u_\tau+b_\tau)/p^e-E_{\tau,e}
   =(u_\tau/p^e-u_{\tau_e}-c_{\tau_e})_+
     +(b_\tau/p^e-b_{\tau_e})_+-c_\tau/p^e
   =D_{\tau,e}.
\]
The \(e=0\) coupling gives solely \(-\tfrac12\sum_\tau E_{\tau,0}\),
with partner \(\langle-1\rangle\tau\). This is the required
\(d_{pj}(\varnothing)\), including the exponent-zero subtraction.
The argument applies to every arrow with the same original tensor.
Transport of bases and composition of the finite class sums preserve
the lattice norm estimates. Since \(p\ne N\), the final branching
matrix is integral and its extra scalar is a unit. Thus the
restricted chain has valuation at least \(g_p-L_p\) at \(p=d\),
since \eqref{eq:marked-column-gain} gives
\(\sum_{j<m}d_{pj}(\varnothing)\ge g_p-L_p\) when
\(C_p(\varnothing)=0\), as proved for \(p=d\) in
Section~\ref{sec:marked-column-gain}.

\subsection{Integrality at the two original level primes}
\label{sec:marked-gauge-at-M}

\sectionopening{In this subsection, we treat the two primes dividing \(M\). The tensor
in Section~\ref{sec:seesaw-integral-tensor} can have denominators
there, so we use the actual translation parameter in
Section~\ref{sec:marked-integral-sections-ordinary-jets} to check the
unit bound.}

Let \(\ell\mid M\), so \(\ell=p_0\) or \(p_1\). These primes
split in \(K\) and are distinct from \(d,N\). For a twist
prime to \(\ell\), the admissible translations of the actual
tensor satisfy \(\bar\pi_M x\) integral above \(\ell\).
This is the datum of standard level \(\delta=\bar\pi_M\);
the argument allows translations with denominators at those places.
The genuine arrow scale in the notation of the ordinary jet bound is
\(c=\pi_M\), so for each homogeneous Schur entry of degree \(B_j\)
the ordinary multiplier is
\[
 \left|\frac{c}{\bar\delta d}\right|_\ell^{B_j}
 =|d^{-1}|_\ell^{B_j}=1.
\]
For the source and target arguments separately, weak approximation
in the Cayley chart provides rational isometries whose left translates
of those arguments are integral at every place above \(\ell\).

Write \(g_1,h_1\) for these integral translates. The action on
displacements is \(A(Z)=h_1 Zg_1^{-1}\), as in
Section~\ref{sec:marked-integral-sections-ordinary-jets}.
The pullback on functions uses \(A^{-1}\).
Let \(L\) be the vacuum lattice and \(H=L\cap A^{-1}L\).
At every split place \(w\) above \(\ell\), \(AL_w=L_w\).
Thus \(H_w=L_w\) and \(\ell\nmid[L:H]\). The vacuum expansion of
Section~\ref{sec:seesaw-integral-tensor} expresses the translate of
\(T(x)v_L\) as a sum of translates \(T(A(x+y))v_L\), \(y\in L/H\).
Their coefficients are \(\Delta_A/[L:H]\) times the character of the
pair and a root of unity. Write \(a_H=\det H_A\) with \(H_A L=H\).
Put \(\vartheta=\upsilon_E|\cdot|_{\mathbb A_E}^{-1/2}\), where
\(\upsilon_E\) is the full splitting character.
Lemma~\ref{lem:weil-vacuum-composition} gives \(\Delta_A=\vartheta_f(a_H)\).
The character \(\vartheta\) is unramified at the finite places. At infinity,
\(\vartheta_\infty((z_v)_v)=\prod_v\overline{z_v}^{-1}\).
Put \(h_E=|\operatorname{Cl}(E)|\) and write \(a_H^{h_E}=\alpha u\)
with \(\alpha\in E^\times\) and \(u\in\widehat{\mathcal O}_E^\times\).
Since \(\vartheta\) is trivial on \(E^\times\),
\[
  \Delta_A^{h_E}=\vartheta_f(a_H)^{h_E}=\prod_v\overline{\Phi_v(\alpha)} .
\]
The generator \(\alpha\) is a unit at every place above \(\ell\)
since \(a_H\) is. Thus
\(v_\kappa(\Delta_A)=0\) at every \(\kappa\mid\ell\). For the determinant
\(t\) of either finite argument, Section~\ref{sec:marked-level-characters}
gives \(t^{h_E}=\alpha_1 u_1\) in the norm-one torus by the
surjection from \(\operatorname{Cl}(E)\) onto its finite quotient.
Here \(\alpha_1\in T_1(F_0)\) and \(u_1\in T_{1,f}^{0}\).
The finite character on \(u_1\) is a root of unity, and its value on
\(\alpha_1\) is the inverse infinite weight monomial on a unit
at every place above \(\ell\). Thus the full pair character is a
\(\ell\)-unit, even at a place of nontrivial tame character.
At each \(w\mid\ell\), \(E\)-linearity gives
\(\delta A(x+y)=A(\delta x+\delta y)\in L_w\).
Thus the span of the allowed translated vacua over the coefficient
valuation ring is preserved, including under the actions at the other
finite places.

The ordinary jet estimate applies at the two integral translates.
Rational equivariance transports it to the two original arguments in
their respective Weyl lattices. The original principal levels give free class sums. In the
normalized tensor, the zero coefficients are integral. In the
case \(a_{\rm ss}\ne0\) the coefficients
\(\Psi_\xi(x)=B_*^{2m}b_{\xi,x}\) include the point multiplicities.
Section~\ref{sec:seesaw-integral-tensor} clears the two residue-test
denominators and the neighbor denominator at each pair. Since \(\ell\nmid q_\xi\),
these coefficients are integral also at the two old primes.

The final branching has integral entries
and unit scalar at \(\ell\ne N\), since \(\pi_N\bar\pi_N=N\).
Consequently the full integral source has lattice norm at most one after restriction.
The saving at \(\ell\mid M\) is zero because its positive partner set is
empty for a split prime. The distinct prime \(N\) uses the original
preimage clearing and has saving zero from \(N\ge x_1\).

\subsection{Norm and index bounds for the Poincar\'e images}
\label{sec:short-images-proved}

\sectionopening{In this subsection, we combine the integral estimates with the complex
norm bound of Section~\ref{sec:upper-chain-size}. The normalized
tensor of Section~\ref{sec:seesaw-integral-tensor} gives the
upper realization of \(L_\xi(c_{\rm ss}^{-1}\Pi)\) in
\eqref{eq:measured-seesaw-normalized-tensor}. Here the finite
tensor in \eqref{eq:seesaw-evaluated-tensor} is the argument of
\(\mathcal U_\xi\). In Section~\ref{sec:marked-chain-composition}
the symbol \(T_\xi\) denotes the weighted map \(L_\xi\), with
\(\Theta_\xi=L_\xi|_V\). We apply the scalar of
Proposition~\ref{prop:shortness-normalization}.}

\begin{lemma}[Class sum of the upper tensor]\label{lem:upper-tensor-class-sum}
The contracted upper integral with the jointly invariant tensor of
Section~\ref{sec:seesaw-integral-tensor} is the sum over the product
of the original upper class sets, with every class given mass one.
The tensor can then be expanded at those representatives.
\end{lemma}

\begin{proof}
The tensor is jointly invariant under the upper levels with the last
action one. At an intermediate
class variable the two adjacent kernels transform by
inverse coefficient actions, and the contracted product
with the invariant tensor is unchanged by the level.
The source factor is already invariant in its original
level. Thus the contracted integrand can be summed over
the Cartesian product of the original upper classes.

Each quotient measure at a definite stage is multiplied
by that stage's level index in the definition of the
algebraic lift. Its mass on a class is the reciprocal
of the class stabilizer order. All the intermediate
levels are free by reduction at an old odd distinguished place,
so every class has mass one in this calculation.
One can therefore expand the finite tensor after choosing
these class representatives. Individual translated
vacuum terms need not have the separate level invariance
of the combined tensor. Their terms are evaluations at
the same representatives in this one finite sum.

Apply \eqref{eq:linked-translate-tensor} at these representatives.
Each linked shift \(t_j(y)\) is supported above \(N\) and
\(Nt_j(y)\) is integral there. The shifts \(x_j\) vanish above \(N\),
so their mutual Heisenberg phases are one. Collecting equal tuples
preserves coefficient integrality off the twist; the translation
denominators and supports away from \(N\) are unchanged.
Linearity of the upper integral and terminal restriction gives
this expansion of the upper operator on the full original source.
Equation~\eqref{eq:measured-seesaw-normalized-tensor} identifies
its value with \(L_\xi(c_{\rm ss}^{-1}\Pi\phi)\) for \(\phi\in S\).
At \(p=N\) the estimate on \(S\) uses
Section~\ref{sec:ramified-direct-lift}. For the exceptional Hecke operators,
the local Hecke operators act on the spherical place
and commute with the original level actions on these
data, so the modified tensor has the same joint invariance
and the same class-sum formula. A one-arrow lattice norm bound
in the transported coordinates therefore composes
inside each summand in the class sum. The ultrametric
inequality then preserves the resulting bound under
summing the classes and the integral tensor coefficients.
\end{proof}

\begin{proposition}[Norm and index on the original translates]\leavevmode\label{prop:translated-norm-and-index}
Fix \(q_0,d,\epsilon\) as in Section~\ref{sec:datum}. For every sufficiently large integer \(x\) there
is a spanning set \(\mathcal S\) of everywhere integral forms in \(S\)
such that, for every complex automorphism \(\sigma\) and \(g\in G\),
\[
 \|\sigma\Theta_1P_g\phi\|^2\le e^{A-m\mathcal R+L}
       \qquad(\phi\in\mathcal S).
\]
For each \(t\subseteq G\) and each independent tuple \(y\) of images
\(P_g\phi\) with \(g\in t\) and \(\phi\in S\) everywhere integral,
\[
 \operatorname{Ind}(y)\le2|t|\dim(S)\,faB^3 .
\]
At each coefficient place over \(p\), its clearing matrix works for
every allowed twist of level prime to \(p\).
\end{proposition}

\begin{proof}

Fix \(p\) and a coefficient place \(\kappa=(p,\iota)\). Write
\(S_\iota\) for the ramified subspace in these coefficient coordinates
and \(M_\iota\) for the full integral source lattice.
When \(p\ne N\) and either \(p\) is split or \(p\ge x_1\), the
ordinary derivative estimate gives lattice norm at most one on
\(M_\iota\cap S_\iota\). At a nonexceptional inert prime \(p<x_1\)
or at \(p=d\) use the same module, with the bounds in
Sections~\ref{sec:marked-column-gain} and~\ref{sec:marked-ramified-ordinary-jet},
respectively. In the exceptional
case use
\[
 J_\iota\cap S_\iota,\qquad
 J_\iota=Y_\Lambda^m M_\iota+p^{j_p}M_\iota,
\]
as in Sections~\ref{sec:exceptional-charged-jet}
and~\ref{sec:exceptional-hecke-analysis}. For
\(\varphi=Y_\Lambda^m u+p^{j_p}v\) in the intersection,
apply the upper tensor operator of Lemma~\ref{lem:upper-tensor-class-sum}
on the full integral source \(M_\iota\). Its values on
\(Y_\Lambda^m u\) and \(p^{j_p}v\) have lattice norms at most
\(p^{2L_p-g_p}\) and \(p^{L_p-g_p}\), respectively.
Linearity and the ultrametric inequality bound its value on
\(\varphi\) by \(p^{2L_p-g_p}\), also when \(p\mid N-1\).
It is this sum \(\varphi\in S_\iota\) to which the measured
seesaw identity \eqref{eq:measured-seesaw-normalized-tensor} applies.
At \(N\), use the preimage under \(U_0^{m/2}\)
of the integral lattice at \(\delta_0\), intersected with
the original integral ramified source. These modules
are fixed after the embedding and before choosing any twist
or translate. At an exceptional prime the natural map gives
\[
 (M_\iota\cap S_\iota)/(J_\iota\cap S_\iota)
                      \hookrightarrow M_\iota/J_\iota .
\]

For a tuple \(\phi_1,\ldots,\phi_s\) of integral preimages,
work over the ring of integers \(O\) of a finite extension of
\(\mathbf Q_p\) over which the lattice, the subspace, the Hecke matrices
and the tuple are defined. At an exceptional prime choose a Smith
basis for the inclusion in Section~\ref{sec:exceptional-jet-hecke}.
After extension to \(R_\iota=\iota(\mathcal O_\kappa)\) this writes
\(M_\iota=\bigoplus_i R_\iota e_i\) and
\(J_\iota=\bigoplus_i d_i R_\iota e_i\), with all \(d_i\ne0\)
because \(J_\iota\) contains \(p^{j_p}M_\iota\).
In a nonzero coordinate of the tuple choose a pivot of least
valuation, subtract its multiples from the other columns, and
multiply the pivot column by \(d_i\). The ratios used to clear
the other columns belong to \(R_\iota\). In a zero coordinate
no operation is needed. These operations preserve previous
divisibilities and membership in \(S_\iota\).
Writing \(\Phi:R_\iota^s\to M_\iota\cap S_\iota\) for the map
with the given tuple as columns, the resulting matrix \(C\)
therefore satisfies
\[
  \Phi C(R_\iota^s)\subseteq J_\iota\cap S_\iota,\qquad
  v_\iota(\det C)\le\sum_i v_\iota(d_i)
           =\frac{\operatorname{length}_O(M_{\iota,O}/J_{\iota,O})}
                  {v_O(p)}\le d_XL_p.
\]
Here \(M_{\iota,O},J_{\iota,O}\) are the two lattices over \(O\)
and \(v_O\) assigns value one to a uniformizer of \(O\).
The bound is unchanged on finite extension, since the length
and \(v_O(p)\) are multiplied by the same ramification degree.
All operations pull back to \(\mathcal O_\kappa\).
At \(p=N\), Section~\ref{sec:ramified-source-clearing} provides
the clearing directly, with determinant valuation at most
\((m/2+1)d_Xfa/k\le mfbd_X\).
In the other cases the whole integral module suffices. The
determinant valuation in each case is at most \(d_X e_p\). For every allowed
twist prime to \(p\) the chain lattice norm after restriction
is at most \(p^{e_p-g_p}\) on this module.

Let \(z=c_{\rm ss}^{-1}\Pi\varphi\), using the common
scalar of Section~\ref{sec:seesaw-integral-tensor}.
Then \eqref{eq:seesaw-measured-chain} and the integral tensor expansion
identify \(L_\xi z\) with the restricted upper chain at \(\varphi\).
Set
\[
       U_g\varphi=|G|\,\rho(e_0)\rho(g)z .
\]
By the lower-chain equivariance of
Section~\ref{sec:seesaw-normalization},
\(\Theta_\xi U_g\varphi\) is the sum of the values of
\(L_\xi z\) at the finite translates. More precisely,
\[
 (\Theta_\xi U_g\varphi)(h)
  =\sum_{s\in G}|G|e_0(s)(L_\xi z)(h k_m(sg)),
\]
because \(U_g\varphi\in V\) and \(\Theta_\xi=L_\xi|_V\).
The coefficients are those of \(|G|e_0\). They are algebraic integers:
for each irreducible character the coefficient of the
central idempotent multiplied by \(|G|\) is its degree
times the character on the inverse class.

The indicated
terminal translations are integral at the places over
\(p\). The transported lattice norm for a point and its right
translate by an integral invertible matrix is the same
on a fixed value. Thus the same bound \(p^{e_p-g_p}\)
holds after the projected translate at every terminal
point, with no change of the local source module.
This uses the original \(\Theta_\xi\) and its prescribed
level, so its unit lattices are the ones tested in the index bound.

Lemma~\ref{lem:ramified-short-generators},
Proposition~\ref{prop:upper-terminal-norm} and
Lemma~\ref{lem:integral-finite-type-projection} give
\[
 \|\sigma\Theta_1 U_g\varphi\|^2
 \le \exp\!\left(
     A-m\mathcal R+2\sum_{p\in F}g_p\log p+2E_0+x^{99}\right).
\]
Here \(K_{\rm in}=E_0\) for the integral generators,
the terminal restriction contributes another \(E_0\), and the passage
from pointwise to class norm costs at most
\(\log h_{\rm term}+4\log|G|\le x^{99}\).
The exponents in the lattice norm are
\(2L_p-g_p=e_p-g_p\) at nonsplit \(p<x_1\) other than \(d\),
at most \(L_p-g_p\le e_p-g_p\) at \(d\), and \(2mfb=e_N\)
at \(N\), where \(g_N=0\). They are zero elsewhere. Using
\(\sum_p e_p\log p\le E_0\), the rational rescaling in
Proposition~\ref{prop:shortness-normalization} now gives
the short length and index bounds for \(P_g=c_B U_g\)
on the original family. For the given ordered tuple \(y\), let \(P\)
be a permutation matrix grouping its columns by their chosen translates.
If \(C_g\) clears the group for \(g\), take
\[
 U_\kappa=P\operatorname{diag}(C_g)P^{-1},\qquad
 a_p=|t|d_Xe_p .
\]
In the original column order this matrix applies each clearing within
its own group. The determinant bound \(v_\iota(\det U_\kappa)\le a_p\)
holds with the same \(a_p\) at every embedding above \(p\).
The scalar \(c_B\) leaves the cleared images integral since
\(e_p-g_p-v_p(c_B)\le0\), and
\[
 |t|d_X\sum_p e_p\log p\le2|t|\dim(S)E_0.
\]
At the infinite places the same scalar gives logarithmic squared norm
at most \(A-m\mathcal R+4E_0+2b+x^{99}\le A-m\mathcal R+L\)
by Proposition~\ref{prop:shortness-normalization}.

For the set \(t\) of Proposition~\ref{prop:marked-orbital-span},
choose a basis \(y_i=P_{g_i}\phi_i\) of the sum of the images from
\(g_i\in t\) and \(\phi_i\in\mathcal S\). Such a basis exists because
\(\mathcal S\) spans \(S\) and \(c_B\ne0\).
Make this choice before the coefficient place or twist. For each
\(\sigma\in\operatorname{Aut}(\mathbb C)\), use the same \(n\)
selected columns. Hadamard's inequality in the positive untwisted
terminal class pairing gives
\[
 \det\operatorname{Gram}(\sigma\Theta_1 y)
 \le\prod_{i=1}^{n}\|\sigma\Theta_1P_{g_i}\phi_i\|^2
 \le \exp\bigl(n(A-m\mathcal R+L)\bigr).
\]
Combining this inequality with the index bound gives
\[
 \operatorname{Cov}(y)\le n(A-m\mathcal R+L),\qquad
 \operatorname{Ind}(y)\le
  10D\log(1+2b^2)E_0=D I_0 .\qedhere
\]
\end{proof}

\section{The mixed Gram determinant and the product formula}
\label{sec:gram-product-formula}

\sectionopening{In this section, we use the product formula to pass from valuations of
mixed pairings to a global determinant inequality. We use the
upper covolume bounds of Section~\ref{sec:integral-poincare-bounds}, the
critical-value bounds of Section~\ref{sec:common-critical-twist}, and the
lattice completion of Section~\ref{sec:chain-determinants-completion}.}

Use the normalized absolute value \(|\ell|_\kappa=\ell^{-1}\)
at a coefficient place \(\kappa=(\ell,\iota)\), where
\(\iota:\mathbb C\simeq\overline{\mathbb Q}_\ell\) is a field
isomorphism. Write
\(\bar\kappa=(\ell,\iota\circ\mathrm{conj})\). Say that a
family \(f_\kappa\) has height at most \(j\) if there is a
single \(H\ne0\) and bounds \(\alpha\in\mathbb R,\ e_\ell\ge0\)
such that
\[
 |\sigma H|\le e^\alpha,\qquad
 |H f_\kappa|_\kappa\le\ell^{e_\ell},\qquad
 \alpha+\sum_{\ell\in U} e_\ell\log\ell\le j
\]
for every complex automorphism, every coefficient place and
every finite set \(U\) of primes.

Height upper bounds add under
multiplication by multiplying the multipliers and adding the bounds.
Suppose that \(z\ne0\) is a constant family of height at most \(j\)
and \(|\sigma z|\le e^u\) for every \(\sigma\).
Choose its multiplier \(H\) and put \(h=Hz\). Since
\(|\sigma h|\le e^{\alpha+u}\) for every \(\sigma\), \(h\) is algebraic.
Indeed, a transcendental number can be carried to a transcendental
number of arbitrarily large modulus by extending a bijection of
transcendence bases to the algebraic closures. Likewise a fixed cleared
coordinate that is bounded at every coefficient embedding above a
single \(\ell\) is algebraic. If it were transcendental, its images
could include \(\ell^{-m}\theta\) for every \(m\ge0\), with
\(\theta\in\overline{\mathbb Q}_\ell\) transcendental over \(\mathbb Q\).
Both algebraically closed fields have cardinality continuum;
the algebraic closure of \(\mathbb Q(\mathcal B)\) has cardinality
\(\max(\aleph_0,|\mathcal B|)\) for a transcendence basis \(\mathcal B\).
Thus both transcendence degrees are continuum. Thus any embedding of a number
field in either field extends to a map of the required type. A finite
list of such fixed algebraic coordinates lies in a single number field.

Write \(K=\mathbb Q(h)\), \(d_h=[K:\mathbb Q]\) and
\(q_h=N_{K/\mathbb Q}(h)\). With \(|\ell|_\lambda=\ell^{-1}\)
at every \(\lambda\mid\ell\), the bounds give
\[
 \begin{aligned}
 \log|q_h|
   &=\sum_{\tau:K\hookrightarrow\mathbb C}\log|\tau h|
       \le d_h(\alpha+u),\\
 \log|q_h|_\ell
   &=\sum_{\lambda\mid\ell}[K_\lambda:\mathbb Q_\ell]\log|h|_\lambda
       \le d_h e_\ell\log\ell.
 \end{aligned}
\]
The first sum counts both members of each pair of complex embeddings.
For every \(\ell\), the local degrees in the second sum satisfy
\(\sum_{\lambda\mid\ell}[K_\lambda:\mathbb Q_\ell]=d_h\).
Take the finite set
\[
 U=\{\ell:\ |h|_\lambda\ne1
              \text{ for some finite place }\lambda\mid\ell\text{ of }K\}.
\]
The rational product formula
\cite[Theorem~7.13]{MilneANT} applied to \(q_h\) gives
\[
 0=\log|q_h|+\sum_{\ell\in U}\log|q_h|_\ell
 \le d_h\Bigl(u+\alpha+\sum_{\ell\in U} e_\ell\log\ell\Bigr)
 \le d_h(u+j).
\]

Let \(w=(y_1,\ldots,y_n,z_1,\ldots,z_s)\) be a basis of the
target space \(V\) of the Poincar\'e map. The comparison lift
\(\Theta_\xi\) has domain \(V\) and is injective. Use the fixed
Hermitian metric on \(V\) and the class sum on the target of each
lift. For a twist \(\xi\), including the trivial twist \(1\), put
\[
 G_\xi=\det\operatorname{Gram}(\Theta_\xi w),\qquad
 Z_\xi=\det\nolimits_V(\Theta_\xi^*\Theta_\xi).
\]
Choose an orthonormal basis \(e\) of \(V\), write \(w=eP\) and
set \(M_\xi=\operatorname{Gram}(\Theta_\xi e)\). Then
\[
 G_\xi=\det(P^*M_\xi P)
      =|\det P|^2 Z_\xi
      =\det\operatorname{Gram}(w)\,Z_\xi.
\]
\begin{proposition}[Mixed Gram determinant inequality]\label{prop:mixed-gram-inequality}
For the basis \(w=(y,z)\) and the injective comparison lifts above,
suppose that the following bounds hold.
\begin{itemize}
\item For every \(\sigma\in\operatorname{Aut}(\mathbb C)\), the Gram determinant of the coordinatewise
  conjugates \(\sigma\Theta_1 y_i\) is at most \(e^c\), and
  each \(\sigma\Theta_1 z_i\) has squared norm at most
  \(e^{b_{\rm short}}\). The positive definite Hermitian matrix
  of the trivial lift's coordinate metric is rational.
\item Fix \(a_\ell\ge0\), independent of the place above \(\ell\),
  with \(\sum_{\ell\in U}a_\ell\log\ell\le\delta\) for every
  finite set of primes \(U\). At each \(\kappa\) above \(\ell\)
  take a clearing matrix \(U_\kappa\) for the \(y_i\), with
  \(\kappa\)-integral entries and
  \(|\det U_\kappa|_\kappa\ge\ell^{-a_\ell}\).
  Its cleared vectors lift
  integrally at \(\kappa\) for every allowed twist at \(\ell\).
\item For each \(\kappa\) choose an allowed twist \(\xi_\kappa\)
  and write \(\gamma=\gamma(\xi_\kappa)\ne0\).
  The vectors \(\gamma\Theta_{\xi_\kappa}z_i\) are integral
  at both \(\kappa\) and \(\bar\kappa\). The pairings of
  the \(\bar\kappa\)-integral lattice against the
  \(\kappa\)-integral lattice have absolute value at most \(1\)
  at \(\kappa\).
\item The families \((\gamma\bar\gamma)^{-1}\) and
  \(Z_1/Z_{\xi_\kappa}\) have heights at most \(j_\gamma,j_Z\).
\end{itemize}
Then
\[
                  0\le c+2\delta+s(b_{\rm short}+j_\gamma)+j_Z .
\]
\end{proposition}

\begin{proof}
\smallskip
\indent\textit{The terminal mixed Gram form.}
The pairing in the determinant is at the original lift
level, principal at the old chosen places and parabolic at
the twist places. The dual delta used for the period at the end of the path
has the smaller principal level at the twist places.
The pointwise Weyl lattice norm checked by those deltas is the
lattice norm of the same lifted vector. Thus the check applies
to its values on the class representatives at the lift level.

The terminal archimedean group is compact, so each class
stabilizer at the lift level is finite by discreteness of the
rational points in the adelic group. At any one old odd
completely split place it embeds in
\(1+p_iM_{n_{\rm term}}(\mathbf Z_{p_i})\) after conjugation by the
class representative. A finite-order element of this
group is one. Indeed, for \(g=I+p_i^aB\) with \(a\ge1\) and some
entry of \(B\) a unit, the minimum valuation of the entries of
\(g^m-I\) is \(a+v_{p_i}(m)\) for \(m\ge1\), by the binomial expansion
first for \(p_i\nmid m\) and then for successive \(p_i\)-th powers.
Thus every class stabilizer is trivial, and the lift-level Gram
pairing is the plain class sum.
Write \(K_{\rm par}\) for this level and \(K_{\rm pr}=K(\pi_MQ(\xi))\)
for the level of \eqref{eq:original-dual-delta-readout}. Over the
class of \(g\), the principal-level representatives are \(gk\)
for \(k\in K_{\rm par}/K_{\rm pr}\). Both vectors acquire the
same unitary level character at \(k\), so the principal-level
class sum multiplies the parabolic one by
\([K_{\rm par}:K_{\rm pr}]=I_\xi/I_{\rm par}(\xi)\),
as computed in Section~\ref{sec:marked-chain-composition}.
We check the lattice at each representative \(g\) of the parabolic
level, using its pointwise dual-delta value.
Write \(J_{\rm c}\) for the signed complement in
Section~\ref{sec:completion-terminal-dual-coordinates}, and
\(G_\vee\) for the Gram matrix of the dual tensor basis. Its
isometry on the two Weyl spaces gives the primal Gram matrix
\[
 J=J_{\rm c}^{\mathsf T}G_\vee J_{\rm c}.
\]
Indeed both sides are the restriction of the exterior tensor
inner product in the primal integral basis. This is the metric
in the upper estimate of Section~\ref{sec:short-images-proved}.
The principal-level delta reads the bilinear form with matrix
\(G_\vee J_{\rm c}\). The integral coordinate extraction in
Section~\ref{sec:completion-terminal-dual-coordinates} therefore
recovers the components in this same primal basis.
At a representative \(g\), choose both transported Weyl
coordinates at their own coefficient embeddings.

Write
\(T\) for the transport by \(g^{-1}\) at the embedding \(\kappa\), and \(T^\#\) for the equivalent transport on
the conjugate coefficient vector. With \(a\) the local
matrices of \(g\) read at \(\kappa\), unitarity reads
\(T=\rho(a^{-1})\), \(T^\#=\rho(a^{\mathsf t})\).
The integral rational Weyl Gram matrix \(J\) satisfies
\(\rho(a^{\mathsf t})^{\mathsf t}J=J\rho(a)\), by the
transpose rule on the exterior tensor construction.
Hence \((T^\#)^{\mathsf t}JT=J\).
Consequently the mixed term at \(\kappa\) is
\[
 \iota\bigl(\langle v_-,v_+\rangle_J\bigr)
    =\bigl(T^\#\iota(\overline{v_-})\bigr)^{\mathsf t}
                                      J\bigl(T\iota(v_+)\bigr).
\]
The two transported vectors have integral coordinates when \(v_-\)
has lattice norm at most one at \(\bar\kappa\) and \(v_+\) has
lattice norm at most one at \(\kappa\). Since \(J\) is integral,
their scalar pairing is integral as well.
Summing over the free classes is integral. For the trivial
twist the coordinate metric is a direct sum of the same
rational Gram matrices, which is the metric used in the
archimedean exterior-product inequality.

Use the convention that the first variable of the Hermitian
form is conjugate-linear. Write the basis vectors as columns.
At a coefficient place \(\kappa=(\ell,\iota)\), a clearing matrix
\(U_\kappa\) acts by column combinations.
The entries of \(\iota(U_\kappa)\) lie in a common finite extension
of \(\mathbb Q_\ell\). Regard the matrix in its ring of integers
\(\mathcal O\), with integer-valued valuation \(v_{\mathcal O}\).
Smith normal form gives
\[
 \operatorname{Fitt}_0(\mathcal O^n/U_\kappa\mathcal O^n)
      =(\det U_\kappa),\qquad
 -\log_\ell|\det U_\kappa|_\kappa
      =\frac{\operatorname{length}_{\mathcal O}
                 (\mathcal O^n/U_\kappa\mathcal O^n)}
             {v_{\mathcal O}(\ell)} .
\]
This normalized length is unchanged by finite extension of the coefficient
field. Under scalar extension both the module length and the valuation
of \(\ell\) acquire the same relative ramification factor. For a fixed \(\kappa\), use the same allowed twist
\(\xi=\xi_\kappa\) at \(\kappa\) and \(\bar\kappa\); both lie over
\(\ell\). Put \(M=\operatorname{Gram}(\Theta_\xi w)\) and
\[
 C_+=\operatorname{diag}(U_\kappa,\gamma I_s),\qquad
 C_-=\operatorname{diag}(U_{\bar\kappa},\gamma I_s).
\]
By the clearing hypotheses for this tuple and \(\ell\nmid q_\xi\),
the columns of \((\Theta_\xi w)C_+\) are integral at \(\kappa\),
and those of \((\Theta_\xi w)C_-\) are integral at \(\bar\kappa\).
Their mixed pairing matrix \(C_-^*MC_+\) therefore has
integral entries at \(\kappa\).
Its determinant is
\(\overline{\det U_{\bar\kappa}}\,(\gamma\bar\gamma)^s
G_\xi\det U_\kappa\). Since \(|\bar t|_\kappa=|t|_{\bar\kappa}\),
\[
 \left|(\gamma\bar\gamma)^sG_\xi\right|_\kappa
 \le |\det U_{\bar\kappa}|_{\bar\kappa}^{-1}
     |\det U_\kappa|_\kappa^{-1}\le\ell^{2a_\ell}.
\]
Thus \(((\gamma\bar\gamma)^sG_{\xi_\kappa})_\kappa\)
has height at most \(2\delta\), with multiplier \(1\).
Now
\[
 G_1=((\gamma\bar\gamma)^sG_{\xi_\kappa})\,
      (\gamma\bar\gamma)^{-s}\,\frac{Z_1}{Z_{\xi_\kappa}},
\]
so the constant \(G_1\) has height at most
\(2\delta+s j_\gamma+j_Z\).

Let \(B\) be the coordinate matrix with columns \(\Theta_1 w\)
and \(Q\) the rational matrix of the coordinate metric. Write
\(\bar\sigma=\mathrm{conj}\circ\sigma\circ\mathrm{conj}\). Entrywise
application of \(\sigma\) gives
\[
 \sigma G_1=\det\bigl((\bar\sigma B)^*Q(\sigma B)\bigr).
\]
Write \(H_\tau(v)=\det\operatorname{Gram}(\tau\Theta_1v)\).
Gram--Schmidt and the norm bounds give
\[
 H_\tau(w)\le H_\tau(y)\prod_{i=1}^s\|\tau\Theta_1z_i\|^2
                         \le e^{c+sb_{\rm short}} .
\]
Cauchy--Schwarz in the \((n+s)\)-th exterior power gives
\[
 |\sigma G_1|^2\le H_\sigma(w)H_{\bar\sigma}(w)
                         \le e^{2(c+sb_{\rm short})}.
\]
Since \(G_1\ne0\), the
preceding product-formula argument yields
\[
 0\le c+2\delta+s(b_{\rm short}+j_\gamma)+j_Z.
\]
The same twist is needed on the two sides of each mixed
pairing. The individual clearing matrices may depend on the
two coefficient places. The volume and denominator heights
are separate arithmetic estimates on the allowed twists.
\end{proof}

Proposition~\ref{prop:marked-completion-vectors} gives the spanning
set and the clearing-scalar height. Proposition~\ref{prop:chain-determinant-height}
gives the determinant-ratio height. Section~\ref{sec:marked-determinant-assembly}
substitutes them into this inequality to prove
Proposition~\ref{prop:marked-determinant}.

\part{Selmer bounds and Galois congruences}\label{part:congruences}

\section{Length bounds after removing auxiliary Selmer conditions}
\label{sec:selmer-relations}

\sectionopening{In this section, we estimate the change in the length of the cocycle
module when auxiliary local relations are discarded. The local bounds
below are used in the polynomial estimate of
Section~\ref{sec:selmer-polynomials}.}

Let \(\mathcal O\) be the ring of integers in a finite extension
of \(\mathbb Q_\ell\). Normalize \(v(\ell)=1\) and put
\(\mathcal B=\mathcal O/\ell^M\mathcal O\) for \(M\ge1\).
In this section lengths are measured by
\(\lambda(X)=\operatorname{length}_{\mathcal O}(X)/v_{\rm ord}(\ell)\),
where \(v_{\rm ord}\) gives a uniformizer valuation one.

\begin{lemma}[Length lost by imposing local relations]\label{lem:selmer-relations}
Let \(\rho\) be a continuous integral Galois representation
of rank \(n\). At each place in a finite list let \(\phi_j\)
be an arithmetic Frobenius lift or its inverse, let \(\gamma_j\)
belong to inertia, and suppose \(\rho\bmod\ell^M\) is
unramified there. All these places have residue characteristic
different from \(\ell\) and residue cardinality \(1\bmod\ell^N\).
Let \(s_j,y\in\mathcal O\), and let \(\psi\) be a continuous
scalar character into \(\mathcal B^\times\). Suppose
\[
 \psi(\phi_j)=1+y,\quad \psi(\gamma_j)=1+s_jy,\qquad
 s_j\ne0\ \Longrightarrow\ \rho(\phi_j)\equiv I\bmod\ell^M .
\]
For every element of inertia at the listed places, suppose
its scalar \(\alpha\) under \(\psi\) satisfies
\[
  \sum_{a=0}^{\ell^N-1}\alpha^a=0\qquad\text{in }\mathcal B .
\]
Put \(R=\psi(\rho\bmod\ell^M)\). Let \(C\) be an
\(\mathcal O\)-submodule of the continuous \(R\)-cocycles
with values in \(\mathcal B^n\). Write \(C_0\) for its submodule
defined by \(z(\gamma_j)=s_jz(\phi_j)\) for all \(j\). Then
\[
 \lambda(C/C_0)\le
     v\!\left(\prod_j\det((1+y)\rho(\phi_j)-I)\right).
\]
\end{lemma}

\begin{proof}
For \(z\in C\), write \(A_g(u)=z(g)+R(g)u\) for the affine
action on \(\mathcal B^n\). If \(g\) belongs to a listed inertia
group, \(R(g)=\alpha I\), so
\[
 A_g^{\ell^N}(u)=
   \alpha^{\ell^N}u+\left(\sum_{a=0}^{\ell^N-1}\alpha^a\right)z(g)=u.
\]
Indeed \(\alpha^{\ell^N}-1=(\alpha-1)\sum_{a<\ell^N}\alpha^a\).
The finite continuous affine image of wild inertia is trivial:
it is a finite group of order a power of the residue characteristic
and all its elements have order dividing \(\ell^N\).
Let \(q_j\) be the residue cardinality and \(f_j\) an arithmetic
Frobenius lift. The tame inertia relation gives
\(A_{f_j}A_{\gamma_j}A_{f_j}^{-1}=A_{\gamma_j}^{q_j}\)
\cite[\S5.1, Lemma~5.1]{BCT}. Since
\(q_j\equiv1\pmod{\ell^N}\) and \(A_{\gamma_j}^{\ell^N}=1\),
both \(A_{f_j}\) and its inverse commute with \(A_{\gamma_j}\).

Put \(F_j=(1+y)\rho(\phi_j)-I\in M_n(\mathcal O)\) and let
\(\bar F_j\) be its reduction in \(M_n(\mathcal B)\). Evaluating the commuting
affine maps at zero and applying the condition on
\(\rho(\phi_j)\) gives, on \(\mathcal B^n\),
\[
 \bar F_j z(\gamma_j)=s_jy z(\phi_j),\qquad s_j\bar F_j=s_jyI.
\]
These identities give the injection
\[
 C/C_0\lhook\joinrel\longrightarrow\bigoplus_j\ker(\bar F_j:\mathcal B^n\to\mathcal B^n),
 \qquad [z]\longmapsto(z(\gamma_j)-s_jz(\phi_j))_j.
\]
Let \((u_{j,i})_i\) be the Smith diagonal of \(F_j\) over
\(\mathcal O\), with \(v(0)=+\infty\).
The kernel of multiplication by \(u\) on \(\mathcal B\) has
ordinary length \(\min\{Mv_{\rm ord}(\ell),v_{\rm ord}(u)\}\).
Thus
\[
 \lambda(C/C_0)\le\sum_j\lambda(\ker\bar F_j)
  =\sum_{j,i}\min\{M,v(u_{j,i})\}
  \le v\!\left(\prod_j\det F_j\right).
\]
The middle sum is finite also for singular \(F_j\).
\end{proof}

For the auxiliary twist constructed from homomorphisms
\(\chi_w\) into \(\mathbb Z/\ell^N\mathbb Z\), use variables
\(x_w\in y\mathcal O\) and put
\(\psi(g)=\prod_w(1+x_w)^{\chi_w(g)}\) in \(\mathcal B\).
The scalar condition in Lemma~\ref{lem:selmer-relations}
follows if \(\sum_{a<\ell^N}\alpha^a\) is divisible by
\(\ell^M\) for every \(\alpha\in1+y\mathcal O\).
In particular each \((1+x_w)^{\ell^N}=1\) in \(\mathcal B\).
This makes the character well defined. With a first
coordinate \(x_0=y\), the other coordinates are \(x_w=s_wy\).
At each subsequent place, the character values needed for
the lemma are the additive values of the \(\chi_w\) in
\(\mathbf Z/\ell^N\mathbf Z\), namely \(1\) in coordinate \(0\) and \(0\) elsewhere
at \(\phi_j\), and \(1\) in coordinate \(j\) and \(0\)
elsewhere at \(\gamma_j\). These statements concern the
chosen pairs at both conjugate places, with inverse lifts
where used. They give the two evaluations of \(\psi\) in the
lemma. The polynomial whose valuation bounds the length
loss is consequently \(\prod_j\det((1+Y)\rho(\phi_j)-I)\).

For a common lower bound \(v(y)\ge a_0>0\), the choice of
\(N\) ensuring the scalar condition can be made uniformly.
For \(z\in y\mathcal O\) the polynomial identity
\[
 \sum_{a=0}^{\ell^N-1}(1+z)^a
      =\sum_{j=1}^{\ell^N}\binom{\ell^N}{j}z^{j-1}
\]
shows that the right side has valuation at least
\(N-\sup_{j\ge1}\{\log_\ell j-(j-1)a_0\}\).
Here \(j\binom{\ell^N}{j}=\ell^N\binom{\ell^N-1}{j-1}\),
and the supremum is finite. At the first place and its
conjugate use the untwisted specialization \(y=0\), with
slope zero. The scalar condition then reduces to \(N\ge M\).
In this specialization the kernel matrices only depend on
the reductions of the first Frobenius matrices. Thus their
length contribution is bounded by \(\sum_j v(\det(\rho(g_j)-I))\)
for any chosen \(g_j\) with
\(\rho(g_j)\equiv\rho(\phi_j)\bmod\ell^M\), with the sum
over the two first-place parameters. This uses the kernel bound
before replacing the lifts in the determinant.

\subsection{The denominator bound for odd cocycles}
\label{sec:odd-cocycle-budget}

\sectionopening{In this subsection, we bound the odd cocycles in the
paired Tate-twist estimate of Section~\ref{sec:triangle-character-cost}.
We reduce their estimation to the annihilator of
Section~\ref{sec:odd-positive-annihilator} and the detection argument
of Section~\ref{sec:odd-large-detection}.}

Set \(N_c=dq_1\), \(t=5r\), and \(y=N_c+t+3\).
The field \(E\) of Section~\ref{sec:datum} embeds in
\(\mathbb Q(\zeta_{N_c})\). For a coefficient DVR
\(\mathcal O\) as above, use ordinary length in this
subsection. Write \(Z^-_{e,\nu}\) for the module of continuous
cocycles \(G_E\to(\mathcal O/\ell^\nu)(-e)\) that vanish on
inertia outside \(\ell\). For \(e>0\) require
\[
 z(cgc^{-1})=-z(g),\qquad
 z|_{D_v\cap\ker\chi_{\rm cyc}}=0
\]
for every decomposition group \(D_v\) at a prolongation
of a place above \(\ell\), using the full \(\ell\)-adic
cyclotomic character. For \(e=0\) require
vanishing on inertia above \(\ell\) and impose no sign.
Lengths are taken on these modules of cocycles themselves.

For \(e=0\) the action is trivial. Since \(E\) has no real
places, global reciprocity \cite[Ch.~V, Theorems~5.3 and~5.5]{MilneCFT} gives
\[
 Z^-_{0,\nu}=\operatorname{Hom}(\mathrm{Cl}(E),\mathcal O/\ell^\nu).
\]
Write \(h_E=|\mathrm{Cl}(E)|\).
If the class group has cyclic factors of orders \(d_i\), then
\[
 \operatorname{length}_{\mathcal O}Z^-_{0,\nu}
   =\sum_i\min\{\nu v_{\rm ord}(\ell),v_{\rm ord}(d_i)\}
   \le v_{\rm ord}(h_E).
\]

\smallskip
\indent\textit{Detection at small coefficient primes.}
Suppose that \(\ell\le y^2\). There is one list of at most
\(y^{20}\) elements of \(G_E\) whose evaluations detect every
cocycle unramified outside \(\ell\) with values in
\((\mathcal O/\ell^\nu)(-e)\), for all \(\nu,e,\mathcal O\).
Put \(m_0=N_c\ell^2\), \(K'=\mathbb Q(\zeta_{m_0})\) and
\(s'=[K':\mathbb Q]\le m_0\). Use representatives of
\(G_E/G_{K'}\), together with lifts of generators of the
maximal abelian quotient of \(G_{K'}\) of exponent dividing
\(m_0\) unramified outside \(\ell\). Denote this quotient by
\(D_{m_0}\) and put \(U_w=\mathcal O_{K'_w}^{\times}\).
The ray class sequence and global reciprocity give
\[
 \prod_{w\mid\ell}U_w/U_w^{m_0}\longrightarrow D_{m_0}
 \longrightarrow \mathrm{Cl}(K')/\mathrm{Cl}(K')^{m_0}
 \longrightarrow 1 .
\]
Indeed, take the inverse limit over moduli supported above \(\ell\)
and then the quotient by \(m_0\)-th powers. The kernel before taking
powers is the image of the local units, and \(K'\) has no real places
\cite[Ch.~V, Theorems~1.7, 5.3 and~5.5]{MilneCFT}. Thus
\[
 |D_{m_0}|\le h_{K'}\prod_{w\mid\ell}
            |K'_w{}^\times/(K'_w{}^\times)^{m_0}|.
\]
The valuation factor, the cyclic group of roots of unity, and the
free part of the principal units give the bound
\(m_0^{2+[K'_w:\mathbb Q_\ell]}\) at \(w\). The free rank is
\([K'_w:\mathbb Q_\ell]\), since the local logarithm identifies
sufficiently deep principal units with a full additive lattice,
also for \(\ell=2\).

Their product is at most \(m_0^{3m_0}\).
In a number field \(L'\) of degree \(s'\), each class has an
integral ideal of norm at most \(\sqrt{|D_{L'}|}\) by
Minkowski's bound \cite[Theorem~4.3]{MilneANT}. There are at most
\(h^{s'}\) ideals of norm \(h\), since over each rational
prime \(p\) there are at most \(s'\) prime ideals and the
number of exponent assignments is at most
\((1+v_p(h))^{s'}\). Hence
\(h_{L'}\le (2|D_{L'}|)^{s'+1}\).
Using \(|D_{K'}|\le m_0^{s'}\) bounds the size of the
class-field quotient by \(2^{m_0^3+5m_0^2+m_0}\).
Thus it has a generating set of size at most
\(m_0^3+5m_0^2+m_0\). With the coset representatives the
total is at most \(m_0^4\le y^{20}\), since \(m_0\ge4\)
and \(m_0\le y^5\).

If a cocycle \(z\) vanishes at those elements, induct on the
power of \(\ell\) annihilating its values. The induction
hypothesis applied to \(m_0z\) gives \(m_0z=0\).
On these values \(G_{K'}\) acts trivially, since it fixes
the \(m_0\)-th roots of unity. The restriction of \(z\)
is now a homomorphism through the indicated class-field
quotient. It vanishes on the generators, hence on \(G_{K'}\),
and the coset representatives give vanishing on all of
\(G_E\) by the cocycle identity. This includes \(\ell=2\).

\smallskip
\indent\textit{Assembly of the bound.}
For the positive twists, take \(\mathcal O\) to be the
integers of a finite extension of \(\mathbb Q_\ell\)
containing a primitive root of unity of order
\(\exp((\mathbb Z/N_c)^\times)\). For \(0<e<t\),
Section~\ref{sec:odd-positive-annihilator} provides an
annihilator \(Y_e\le 2^{4y^6}\) independent of
\(\ell,\nu,\mathcal O\) under this condition.
For \(\ell\le y^2\) the same set of at most \(y^{20}\) evaluations
works for all exponents and depths. For \(\ell>y^2\),
Section~\ref{sec:odd-large-detection} gives at most \(y^{20}\)
evaluations for each fixed \(e,\nu,\mathcal O\).
In either case the bound on the number of evaluations is uniform.
Evaluation therefore embeds \(Z^-_{e,\nu}\) into
at most \(y^{20}\) copies of
\((\mathcal O/\ell^\nu)[Y_e]\), of length at most
\(v_{\rm ord}(Y_e)\) each. Together with the zero twist this
gives one integer
\[
 X_{\rm odd}=h_E\left(\prod_{0<e<t}Y_e\right)^{y^{20}},
 \qquad
 \operatorname{length}_{\mathcal O} Z^-_{e,\nu}
      \le v_{\rm ord}(X_{\rm odd}),\qquad
 \log X_{\rm odd}\le\log h_E+4y^{27}\log 2 .
\]
For an arbitrary coefficient ring \(\mathcal O\), adjoin the
specified roots of unity and let \(\mathcal O'\) be the resulting
ring of integers. In an \(\mathcal O\)-basis of \(\mathcal O'\),
the cocycle, inertia, cyclotomic-kernel and parity conditions hold
coordinatewise because the coefficient action is scalar. Thus
\[
 Z^-_{e,\nu}(\mathcal O')=
 Z^-_{e,\nu}(\mathcal O)\otimes_{\mathcal O}\mathcal O'.
\]
Both the length and the valuation of \(X_{\rm odd}\) multiply
by the ramification index. Hence the normalized bound
\(v_\ell(X_{\rm odd})\) per exponent holds over \(\mathcal O\).
The ideal-class estimate above, applied to \(E\), and
\(\operatorname{rd}(E)\le N_c\) give
\[
 \log h_E=O(f^2\log N_c),\qquad
 N_c=O_d(x^{11}\log x),\quad t=O(x^5).
\]
Consequently, for every finite set of primes \(U\),
\[
 \sum_{\ell\in U}v_\ell(X_{\rm odd})\log\ell
                   \le\log X_{\rm odd}\le x^{2000}
\]
for all sufficiently large \(x\). Summing the odd contributions of
at most \(t\) ordered exponents therefore contributes at most
\(t x^{2000}\) to the logarithmic denominator bound.

\smallskip
\indent\textit{Application to a paired cocycle test.}
For \(e>0\) let \(S\) be a module of cocycles unramified
outside \(\ell\) satisfying the local vanishing above \(\ell\).
Suppose evaluations at \(g_0\) and a finite set of points
detect all unramified-outside-\(\ell\) cocycles of this type.
Put \(d_e=\chi_{\rm cyc}(g_0)^{-e}-1\).
Assume \(d_e z(g_0)=0\) for \(z\in S\), and that at each test
point \(a\) the paired Frobenius conditions give
\(2(z+z\circ c)(a)=0\), where
\((z\circ c)(g)=z(cgc^{-1})\).
Under these hypotheses
\(z\mapsto z\circ c-z\) sends \(S\) to \(Z^-_{e,\nu}\).
The local conditions are invariant under conjugation by
\(c\).

In the kernel \(z=z\circ c\), hence \(4z(a)=0\) at
all test points. The cocycle \(4d_e z\) therefore vanishes
at every test point and at \(g_0\), so it is zero.
If \(4d_e\mid C\), in particular for the prescribed multiple of \(4d_e^2\),
the kernel is contained in \(S[C]\), and consequently
\[
 \operatorname{length}_{\mathcal O}(S/S[C])
                  \le \operatorname{length}_{\mathcal O}Z^-_{e,\nu}.
\]
The paired evaluations and separation scalar act on the same
cocycle module over \(\mathcal O\). The coefficient descent above
gives the bound for this ring, also at \(\ell=2\).

\subsection{A uniform annihilator for the positive twists}
\label{sec:odd-positive-annihilator}

\sectionopening{In this subsection, we construct a common integer that controls the
positive exponents appearing in Section~\ref{sec:odd-cocycle-budget}.
The distribution identity for the Bernoulli polynomial also supplies
the auxiliary units needed in
Section~\ref{sec:odd-jacobi-derivatives}.}

Fix \(e>0\), abbreviate \(N_c\) to \(N\), and put
\(u=N+e+3\). We construct an integer \(0<X_e\le 2^{u^6}\).
For each prime \(\ell\) it supplies a positive integer
\(a\equiv1\pmod N\), prime to \(\ell\), with
\((a^e-1)(a^{e+1}-1)\mid X_e\).
For a character \(\psi\) of \((\mathbf Z/N)^\times\)
with \(\psi(-1)=-(-1)^e\) and values in \(\mathcal O\), set
\[
 S_m(a,\psi)=
 \sum_{b\in(\mathbf Z/Nm)^\times}
  b^e\psi(b)
  \left(\left\lfloor a\left\langle\frac b{Nm}\right\rangle\right\rfloor
       -\left\lfloor a\left\langle\frac{-b}{Nm}\right\rangle\right\rfloor\right)
       \pmod m ,
\]
where angle brackets denote fractional part.
For \(m=\ell^\nu\), this scalar in \(\mathcal O/m\mathcal O\)
divides \((a^{e+1}-1)X_e\) for every positive
\(a\equiv1\pmod N\) prime to \(\ell\).

\smallskip
\indent\textit{Construction of the integer.}
For the Bernoulli polynomials put
\[
 B(b)=(e+2)!N^{e+1}B_{e+1}(b/N),\qquad
 P(X)=(e+1)!N^{e+2}(B_{e+2}(X/N)-B_{e+2}(0)).
\]
Both have integer coefficients. In fact \((j+1)!B_j\)
is integral by the Bernoulli recursion. In the coefficient
of \(X^k\) in \(P\), \(j=e+2-k\). For \(k\ge2\),
\((j+1)!\mid(e+1)!\), and for \(k=1\) the
binomial coefficient contributes \(e+2\).
The identities \(P'=B\) and
\(P(b+N)-P(b)=N(e+2)!b^{e+1}\) now give
\[
 (e+2)!\sum_{i=0}^{m-1}(b+Ni)^{e+1}
                  \equiv m B(b)\pmod{m^2}.
\]
Indeed, summing the difference gives
\[
 N(e+2)!\sum_{i=0}^{m-1}(b+Ni)^{e+1}
   =P(b+Nm)-P(b)=Nm B(b)+(Nm)^2 R,\qquad R\in\mathbf Z .
\]
Division in \(\mathbf Z\) by \(N\) gives the congruence also when
\(\ell\mid N\).

Extend \(\psi\) by zero on the nonunits modulo \(N\). In
particular \(T_m\) includes multiples of \(\ell\) when
\(\ell\nmid N\). Use the sums
\[
 T_m=\sum_{0\le b<Nm}\left\lfloor\frac{ab}{Nm}\right\rfloor
                   b^e\psi(b),\qquad
 U_\psi=\sum_{0\le b<N}B(b)\psi(b).
\]
Multiplication by \(a\) permutes residues modulo \(Nm\)
and preserves \(\psi\). Modulo \((Nm)^2\), expand
\[
 \bigl(ab-Nm\lfloor ab/(Nm)\rfloor\bigr)^{e+1}.
\]
Multiply by \((e+2)!\), sum, group by the residues modulo
\(N\), and cancel \(m\) in the congruence modulo \(m^2\).
This gives
\begin{equation}\label{eq:odd-voronoi-congruence}
 (a^{e+1}-1)U_\psi
      \equiv(e+2)!(e+1)Na^e T_m\pmod m .
\end{equation}

On \(G=(\mathbf Z/N)^\times\) define the integer matrix
\[
 A_{s,t}=2B(t/s)+(-1)^e\mathbf1_{t=s}+\mathbf1_{t=-s},
\]
with arguments read in \([0,N)\). Work first in the cyclotomic number field of character values,
using a complex embedding to test nonvanishing. The vector
\((\psi(s))_{s\in G}\) has eigenvalue \(2U_\psi+(-1)^e+\psi(-1)\).
Grouping the Dirichlet series by residue and applying the
Hurwitz special value gives
\(U_\psi=-(e+2)!(e+1)N L(-e,\psi)\).
For \(N=1\) the sole character has \(L(s,\psi)=\zeta(s,1)\);
here \(B_{e+1}(0)=B_{e+1}(1)\) since \(e+1\ge2\).
At primitive conductor one the zeta reflection formula
\cite[\href{https://dlmf.nist.gov/25.4.E2}{(25.4.2)}]{DLMF}
gives vanishing at \(-e\) exactly when \(e\) is even.
At primitive conductor greater than one the parity factor
in the Dirichlet functional equation gives vanishing at
\(-e\) exactly when \(\psi(-1)=(-1)^e\). In both cases
the Euler product at \(1+e\) converges to a nonzero value.
For a character induced by \(\psi^*\), every omitted factor
\(1-\psi^*(p)p^e\) is nonzero. We have used
\cite[\href{https://dlmf.nist.gov/25.11.E14}{(25.11.14)}
and \href{https://dlmf.nist.gov/25.15.E3}{(25.15.3)}--\href{https://dlmf.nist.gov/25.15.E5}{(25.15.5)}]{DLMF}.
The eigenvalues are therefore \(2(-1)^e\) on the vanishing
parity and \(2U_\psi\ne0\) on the other parity, so
\(\det A\ne0\).

Take
\[
 X_e=|\det A|\prod_{a\in\{1+N,1+2N\}}(a^e-1)(a^{e+1}-1).
\]
The two choices for \(a\) are relatively prime.
The Bernoulli Fourier series
\cite[\href{https://dlmf.nist.gov/24.8.E3}{(24.8.3)}]{DLMF} bounds \(|B_{e+1}(x)|\) by \(17(e+1)!\) on \([0,1]\).
With \(d'=\varphi(N)\) and
\(M'=(e+2)!N^{e+1}(e+1)!\), this gives
\(|\det A|\le d'!(64M')^{d'}\).
Use \(d'\le u\), \(M'\le u^{3u}\) and bound the four
factors by \((2u)^u\) each. Then
\(\log_2 X_e\le3u^3+5u^2+10u\le u^6\) for \(u\ge5\).

For the odd parity the eigenvalue is \(2U_\psi\).
In the identity coordinate the character vector has entry \(1\).
Apply \(\operatorname{adj}(A)A=\det(A)I\) there and embed its
cyclotomic integer coefficients in \(\mathcal O\). This shows \(X_e=2U_\psi q_\psi\) with
\(q_\psi\in\mathcal O\). Multiply
\eqref{eq:odd-voronoi-congruence} by \(2q_\psi\) to
obtain
\[
 (a^{e+1}-1)X_e\equiv 2wT_m\pmod m,\qquad
 w=q_\psi(e+2)!(e+1)Na^e.
\]
The same \(w\) works at every power of \(\ell\).
Write \(T_{\ell^j}=V_j+\ell^e\psi(\ell)T_{\ell^{j-1}}\),
where \(V_j\) is the sum on units modulo \(N\ell^j\).
Subtract the scaled congruence at level \(j-1\), using
\(e\ge1\), to get
\[
 (1-\ell^e\psi(\ell))(a^{e+1}-1)X_e
                                  \equiv 2wV_j\pmod{\ell^j}.
\]
For \(j=1\) the level-zero congruence is automatic.
Replacing the representative \(b\) by that of \(-b\)
changes \(b^e\psi(b)\) to its negative modulo \(\ell^j\),
so \(S_{\ell^j}(a,\psi)=2V_j\) in the quotient.
The scalar \(1-\ell^e\psi(\ell)\) is a unit there. This
proves the asserted divisibility, also for \(\ell=2\).

\subsubsection{The Jacobi quotient}\label{sec:odd-cocycles-jacobi-quotient}

\sectionopening{In this subsection, we compute the split-prime Jacobi quotient for the odd cocycles of Section~\ref{sec:odd-cocycle-budget}.}

Let \(m=\ell^\nu\) with \(\nu>0\), let
\(p\equiv1\pmod{Nm}\) be prime, and fix a place \(B\)
above \(p\). Write \(\theta_j\) for the multiplicative characters of
\(\mathbf F_p^\times\), indexed by \(j\in\mathbf Z/(p-1)\). Fix a nontrivial additive
character \(\xi\) and put \(G_j=\sum_{b\ne0}\theta_j(b)\xi(b)\).
Normalize \(v_B(p)=1\). Inertia fixes the multiplicative
characters and acts on \(\xi\) by scaling its argument.
The resulting scalars for \(G_j\) multiply in \(j\).
Thus \(G_{j+k}/(G_jG_k)\) belongs locally to the unramified
subfield and has integral order with \(v_B(p)=1\).
It follows that \(j\mapsto v_B(G_j)\bmod\mathbf Z\) is additive.
For \(j\ne0\), reduction of the additive character to
one gives positive order, and
\(G_j\sum_{b\ne0}\theta_{-j}(b)\xi(-b)=p\)
gives order less than one. The product identity follows
by summing first over the ratio of the two nonzero
arguments: the inner additive sum is \(p-1\) at ratio
one and \(-1\) otherwise. For \(j=0\), \(G_0=-1\).
Hence \(j\mapsto v_B(G_j)\bmod\mathbf Z\) identifies the
character group with the \((p-1)\)-torsion of \(\mathbf Q/\mathbf Z\).

If \(x\) corresponds to \(j\), put
\[
 J(x)=\frac{G_j^a/G_{aj}}{G_{-j}^a/G_{-aj}} .
\]
A Galois change of additive character cancels in each
ratio \(G_j^a/G_{aj}\). A change of multiplicative character
sends \(j\) to its cyclotomic multiple. Thus if \(x\) has
denominator dividing \(Nm\), then \(J(x)\in
\mathbf Q(\zeta_{Nm})^\times\), and its conjugate under
\(b\in(\mathbf Z/Nm)^\times\) is \(J(bx)\).
The Gauss product shows that these numbers are units off
\(p\). By the least nonnegative orders above,
\[
 v_B(J(x))=\lfloor a\langle x\rangle\rfloor-
                  \lfloor a\langle-x\rangle\rfloor .
\]

\subsubsection{Annihilation on character components}\label{sec:odd-cocycles-character-annihilation}

\sectionopening{In this subsection, we apply reciprocity to the Jacobi quotient of Section~\ref{sec:odd-cocycles-jacobi-quotient}.}

Write \(C=\mathcal O/\ell^\nu\mathcal O\),
\(K'=\mathbf Q(\zeta_{Nm})\) and
\(\Theta_\psi=\chi_{\rm cyc}^e\psi\).
For this annihilation argument take \(a\in\{1+N,1+2N\}\)
prime to \(\ell\). Such a choice exists since
\(2(1+N)-(1+2N)=1\); at \(\ell=2\) one may take \(1+2N\).
The product defining \(X_e\) contains
\((a^e-1)(a^{e+1}-1)\) as an integer factor.
Let \(z_\psi:G_{\mathbf Q}\to C\) be a cocycle for
\(\Theta_\psi^{-1}\), with \(z_\psi(c)=0\) at complex
conjugation \(c\). Assume it vanishes on inertia off
\(\ell\) lying in \(G_{\mathbf Q(\zeta_N)}\), and on
the elements of decomposition groups over \(\ell\)
that fix \(\zeta_N\) and have \(\chi_{\rm cyc}=1\).
Suppose \(\Theta_\psi(c)=-1\). On \(G_{K'}\) the character
\(\Theta_\psi^{-1}\) is trivial in \(C^\times\), so
\(z_\psi|_{G_{K'}}\) is an additive homomorphism.
At any place over \(\ell\), choose inertia \(\delta\)
with \(\chi_{\rm cyc}(\delta)=a\). This is possible by
cyclotomic total ramification over \(\mathbf Q_\ell\).
It fixes \(\zeta_N\), since inertia fixes the prime-to-\(\ell\)
part and \(a\equiv1\pmod N\) handles the \(\ell\)-part.
For \(g\in G_{K'}\) in that decomposition group,
\([\delta^{-1},g]\) fixes \(\zeta_N\) and has trivial
cyclotomic character. Its cocycle value is
\((a^e-1)z_\psi(g)\) since \(g\) acts trivially on \(C\).
It follows that \(2(a^e-1)z_\psi|_{G_{K'}}\)
vanishes on the decomposition groups above \(\ell\),
the inertia groups off \(\ell\), and any real
decomposition groups.

Apply global reciprocity over \(K'\) to this homomorphism
and the principal idele \(J(1/(Nm))\)
\cite[\href{https://www.jmilne.org/math/CourseNotes/CFTc.pdf\#page=187}{Ch.~V, Proposition~5.2}
and \href{https://www.jmilne.org/math/CourseNotes/CFTc.pdf\#page=188}{Theorem~5.3(a)}]{MilneCFT}.
At the chosen place over \(p\) write \(\phi\) for an
arithmetic Frobenius and choose lifts \(s_b\) of
\(b\in(\mathbf Z/Nm)^\times\). The other places above
\(p\) are \(s_b^{-1}B\). Their Frobenius values are
\[
 z_\psi(s_b^{-1}\phi s_b)=b^e\psi(b)z_\psi(\phi)
\]
because \(\phi\) acts trivially on \(C\). Their orders on
the principal idele are \(v_B(J(b/(Nm)))\).
The sum of local Artin images is zero, giving
\[
 S_m(a,\psi)\,2(a^e-1)z_\psi(\phi)=0 .
\]
Write \((a^{e+1}-1)X_e=S_m k\) in \(C\) and
\(X_e=(a^e-1)(a^{e+1}-1)d_0\) in \(\mathbf Z\).
Multiplication by \(kd_0\) gives \(2X_e^2 z_\psi(\phi)=0\).

The cocycle \(2X_e^2z_\psi\) is zero on \(G_{K'}\)
by Chebotarev's density theorem
\cite[\href{https://www.jmilne.org/math/CourseNotes/ANT.pdf\#page=149}{Theorem~8.31, p.~147}]{MilneANT}.
Indeed take a finite normal quotient
carrying its values and the cyclotomic quotient modulo
\(Nm\). On the kernel of the latter quotient, conjugacy
changes a cocycle value by a unit, and each conjugacy
class in the kernel is represented by a Frobenius above
a rational prime split in \(K'\). So \(2X_e^2z_\psi\)
factors through the abelian group \(\operatorname{Gal}(K'/\mathbf Q)\).
For a cocycle on that group with value zero at \(c\),
commuting an element with \(c\), which acts by \(-1\),
shows that twice every value is zero. We obtain
\(4X_e^2 z_\psi=0\).

\subsubsection{Integral recombination of character components}\label{sec:odd-cocycles-integral-recombination}

\sectionopening{In this subsection, we pass from Section~\ref{sec:odd-cocycles-character-annihilation} to the integral lattice.}

For \(z\in Z^-_{e,\nu}\), let \(s_v\) be lifts of
\(v\in(\mathbf Z/N)^\times\), with \(s_1=1\). Form the
character components
\[
 z_\psi(g)=\sum_v \Theta_\psi(s_v)
                  z(s_v[c,g]s_v^{-1}).
\]
Here \([c,g]=cgc^{-1}g^{-1}\). Each conjugated commutator lies in
\(G_{\mathbf Q(\zeta_N)}\subseteq G_E\) and acts
trivially on \(C\). For \(h\in G_E\) acting trivially,
the cocycle rule gives
\(z(khk^{-1})=\chi_{\rm cyc}(k)^{-e}z(h)\) for \(k\in G_E\).
It follows that each summand as a function of \(s_v\) is
unchanged by left multiplication by \(G_{\mathbf Q(\zeta_N)}\).
Using \([c,gg']=[c,g]g[c,g']g^{-1}\) and reindexing
\(s_vg\) proves the cocycle rule for \(\Theta_\psi^{-1}\).
The value at \(c\) is zero. At an element \(g\) in
either prescribed local subgroup, the two conjugates
\(s_vcgc^{-1}s_v^{-1}\) and \(s_vgs_v^{-1}\) are in
the corresponding prescribed local subgroups of \(G_E\).
Their values under \(z\) vanish, and their difference
is the commutator value. Thus the local conditions above
hold for every component.

Assume \(\mathcal O\) contains the roots of unity of
order the exponent of \((\mathbf Z/N)^\times\).
Character orthogonality holds as an identity in
\(\mathcal O\) before reduction. For \(g\in G_E\) it gives
\[
 \sum_\psi z_\psi(g)=\varphi(N)(z(cgc^{-1})-z(g))
                             =-2\varphi(N)z(g).
\]
If the odd components vanish, the even components are
invariant under conjugation by \(c\) by their cocycle rule.
Comparing the sum at \(g\) and \(cgc^{-1}\) gives
\(4\varphi(N)z=0\). Apply this to \(4X_e^2 z\), whose
odd components vanish by the previous argument. Hence
\begin{equation}\label{eq:odd-positive-annihilator}
              Y_e Z^-_{e,\nu}=0,\qquad
              Y_e=16\varphi(N)X_e^2\le 2^{4u^6}.
\end{equation}
The bound is uniform in \(\ell\), \(\nu\), and the coefficient
ring. Its integral character argument includes the primes
dividing \(2\varphi(N)\).

\subsection{Detection at large coefficient primes}
\label{sec:odd-large-detection}

\sectionopening{In this subsection, we use the integer of
Section~\ref{sec:odd-positive-annihilator} to choose a torsion order
at a large coefficient prime. Reciprocity and Jacobi derivatives will
then detect the nonzero odd component from
Section~\ref{sec:odd-cocycle-budget}.}

Fix \(0<e<t\), a character \(\psi\) with
\(\psi(-1)=(-1)^{e+1}\) as in
Section~\ref{sec:odd-positive-annihilator}, a coefficient
ring \(O=\mathcal O\), and the component depth \(\nu\).
Let \(C_{\psi,\nu}\) consist of the continuous cocycles
\(G_{\mathbf Q}\to(O/\ell^\nu O)((\chi_{\rm cyc}^e\psi)^{-1})\)
vanishing at \(c\). On restriction to \(G_{\mathbf Q(\mu_{N_c})}\), they
vanish on inertia away from \(\ell\) and on the cyclotomic kernel in
each decomposition group above \(\ell\), as in
Section~\ref{sec:odd-cocycles-character-annihilation}.
Put \(h=v_\ell(X_e)+1\). We prove
\[
 \exists g_1,\ldots,g_h\in G_{\mathbf Q},\qquad
 \{z\in C_{\psi,\nu}: z(g_i)=0\ (1\le i\le h)\}=0.
\]
It suffices to detect \(C_{\psi,\nu}[\ell]\): in a nonzero
evaluation kernel the last nonzero \(\ell\)-multiple of any
nonzero element is killed by \(\ell\).
The inequality \(\ell>y^2\) gives
\(e+1<\ell-1\) and \(\ell\nmid N_c\).
For the detection argument choose a positive \(a\equiv1\bmod N_c\)
with primitive residue modulo \(\ell\), by the Chinese remainder
theorem. For this \(a\), Section~\ref{sec:odd-positive-annihilator}
gives
\[
 (a^{e+1}-1)X_e\in S_{\ell^h}(a,\psi)\,(O/\ell^h O).
\]
Indeed the Voronoi congruence there gives this membership after
multiplication by the unit \(1-\ell^e\psi(\ell)\).
Then \(2,\varphi(N_c\ell),a^e-1,a^{e+1}-1\) are units
in \(O\). At a place above \(\ell\) choose
\(\gamma\) in inertia, fixing \(\mu_{N_c}\), with
\(\chi_{\rm cyc}(\gamma)=a\). The existence follows by
restricting to the compositum of the unramified
prime-to-\(\ell\) cyclotomic extension and the
totally ramified \(\ell\)-power cyclotomic extension.
The value of the dual character
\(\Upsilon=\chi_{\rm cyc}\Theta_\psi=\chi_{\rm cyc}^{e+1}\psi\) at \(\gamma\) is
\(a^{e+1}\), with \(\Upsilon(\gamma)-1\) a unit.

Here and below, a cocycle with finite coefficient depth
has open zero set, which is a subgroup. If \(w(g_0)=0\)
and \((\chi_w(g_0)-1)w\ne0\), the zero set cannot contain
the subgroup fixing any cyclotomic modulus. Indeed, vanishing on that subgroup would give
\(w(g_0g)=w(gg_0)\), since \((g_0g)^{-1}gg_0\) belongs to it.
The cocycle identity would then give
\((\chi_w(g_0)-1)w(g)=(\chi_w(g)-1)w(g_0)=0\) for every \(g\).
For a nonzero odd component take \(g_0=c\), so its character
difference is \(-2\). For a nonzero normalized dual cocycle
take \(g_0=\gamma\). These differences are units in the
large-prime range.

Suppose no \(h\) evaluations detect the odd component.
For any set \(T\) of at most \(h\) primes
\(r\equiv1\bmod N_c\ell\), use their Frobenius elements
as tests and take the last nonzero \(\ell\)-multiple of a
missed cocycle. This gives a nonzero component \(z\), killed
by \(\ell\), trivial at all test Frobenii. At these split
primes the action on the values of \(z\) is trivial.
It vanishes on inertia, so it vanishes on each tested
decomposition group. Conjugating a tested Frobenius multiplies
its value by the character at the conjugating element, since
the Frobenius action on the coefficients is trivial.
The same applies to inertia, so \(z\) kills every decomposition
group above each \(r\in T\).

For sets \(T\) of primes split modulo \(N_c\ell\), call
\(T\) tight at depth \(h\) when every dual cocycle in
\(O/\ell^h\) killed by \(\ell\) and unramified over
\(\mathbb Q(\mu_{N_c})\) outside \(\ell,T\) is unramified
there outside \(\ell\). The empty set is tight.
If \(T\) is tight and \(z\) is nonzero on a decomposition group
at a new split \(r\) and trivial above \(T\), the
modulo-\(\ell\) reciprocity argument in Section~\ref{sec:odd-tight-reciprocity} implies that
\(T\cup\{r\}\) is tight.

Choose a finite set \(Q'\) of primes
\(r\equiv1\bmod aN_c\ell^h\) detecting every pair \((u,z)\)
in which \(z\) is a nonzero \(\ell\)-killed odd component and \(u\) is a
dual cocycle at depth \(h\), unramified outside \(\ell\)
over \(\mathbb Q(\mu_{N_c})\), and nonzero on the subgroup
fixing \(\mu_{aN_c\ell^h}\). Both cocycle spaces are finite.
Over the cyclotomic field killing the coefficient action,
each restriction is a homomorphism to a finite additive
coefficient module. Its kernel fixes an extension of degree
at most the order of that module, unramified outside \(\ell\).
There are finitely many such extensions
\cite[\href{https://www.jmilne.org/math/CourseNotes/ANT.pdf\#page=153}{Theorem~8.42, p.~151}]{MilneANT}.
Each restriction and the values on finitely many coset
representatives then have only finitely many possibilities.

On the subgroup fixing
\(\mu_{aN_c\ell^h}\) the two zeros are proper subgroups.
Their union is proper. If neither contains the other,
a product of elements exclusive to the two subgroups
lies in neither. Choose such an element and pass to a finite
Galois quotient through which the two cocycles, their coefficient
actions and the cyclotomic quotient factor. Chebotarev gives
infinitely many unramified rational primes in its conjugacy class.
Every such prime splits completely in
\(\mathbf Q(\mu_{aN_c\ell^h})\). Conjugating a prime of the finite
extension above it, and then prolonging that prime to
\(\overline{\mathbf Q}\), realizes the prescribed element as the
image of arithmetic Frobenius. Both cocycle values are nonzero
at this same prolongation. Choosing a distinct prime for each
pair gives \(Q'\). The preceding abelian
calculation ensures that every nonzero normalized dual
cocycle here qualifies for this pair test.

Choose a nonzero odd component killed by \(\ell\) and
apply its cyclotomic nonvanishing to the modulus
\(N_c\ell^h\prod_{r\in Q'}r\). Chebotarev gives a prime \(p\)
split at this modulus, a distinguished place \(B\mid p\), and
a nonzero component value at its Frobenius.
This uses the same finite-quotient argument as above
with one cocycle. In particular \(p\notin Q'\), and
reciprocity with \(T=\varnothing\) makes \(\{p\}\) tight.
Section~\ref{sec:odd-jacobi-derivatives} gives \(\kappa_S\) for
\(S\subseteq Q'\) at depth \(h\), with this \(p,B,\gamma\).
It uses \(r\equiv1\bmod aN_c\ell^h\) for \(r\in Q'\) and
\(p\equiv1\bmod N_c\ell^h\prod_{r\in Q'}r\).

Prove by descending induction for \(i=h,\ldots,0\) that
\(\ell^i\kappa_S\) is unramified outside \(\ell\) over
\(\mathbb Q(\mu_{N_c})\) whenever \(|S|\le i\) and
\(\{p\}\cup S\) is tight. At \(i=h\) it is zero.
For \(i<h\), if \(u=\ell^{i+1}\kappa_S\ne0\), then it
is unramified by induction and normalized at \(\gamma\).
Since \(|\{p\}\cup S|\le h\), choose \(z\ne0\), killed
by \(\ell\), trivial above \(\{p\}\cup S\).
Use the pair test to choose \(r\in Q'\), a prolongation
\(A\mid r\) and a lift \(\phi\) of arithmetic Frobenius at
\(A\) with \(u(\phi),z(\phi)\ne0\). It lies outside \(S\)
and \(z\) makes \(\{p\}\cup S\cup\{r\}\) tight.
Since \(|S\cup\{r\}|\le i+1\), the induction hypothesis at
\((i+1,S\cup\{r\})\) gives vanishing on inertia above \(r\).
Apply the Frobenius--inertia identity of
Section~\ref{sec:odd-frobenius-inertia} at the same prolongation to obtain
\[
 u(\phi)=\ell^{i+1}\kappa_S(\phi)
   =\ell^{i+1}\kappa_{S\cup\{r\}}(\tau_r)=0,
 \qquad \tau_r\in I_A,
\]
a contradiction. Thus \(u=0\) as a cocycle and
\(\ell^i\kappa_S\) is killed by \(\ell\). By construction it
is unramified outside \(\ell,\{p\}\cup S\), so tightness
removes the additional primes. This proves the induction.

At \(i=0\) the cocycle \(\kappa_\varnothing\) kills
inertia at \(B\). The base derivative identity gives
\((\Upsilon(\gamma)-1) S_{\ell^h}(a,\psi)=0\) in
\(O/\ell^h\). Cancel the unit and use the Voronoi
divisibility of Section~\ref{sec:odd-positive-annihilator} to get \(\ell^h\mid(a^{e+1}-1)X_e\)
in \(O\). The factor \(a^{e+1}-1\) is a unit.
This contradicts \(h=v_\ell(X_e)+1\), in the
normalized valuation \(v_\ell(\ell)=1\).
This proves existence of the \(h\) component tests.

For each \(\psi\) with \(\psi(-1)=(-1)^{e+1}\), take its
\(h\) tests \(g\).
Evaluate the original cocycle at all \(s_u[c,g]s_u^{-1}\)
in the character sums of
Section~\ref{sec:odd-positive-annihilator}. At most
\(\varphi(N_c)^2(v_\ell(X_e)+1)\le y^{20}\)
points are used, since \(v_\ell(X_e)\le y^6\).
If \(z\) vanishes there, every component of total
character \(\chi_{\rm cyc}^e\psi\) of sign \(-1\) vanishes. Then \(4\varphi(N_c)z=0\) and
\(4\varphi(N_c)<\ell\), so \(z=0\).

Thus these evaluations detect the original module for
\(\ell>y^2\).

\subsubsection{Local reciprocity at an auxiliary prime}
\label{sec:odd-tight-reciprocity}

\sectionopening{In this subsection, we apply local reciprocity at a split auxiliary
prime as part of the detection argument of
Section~\ref{sec:odd-large-detection}. The local calculation
prescribes how an unramified cocycle is evaluated against the
auxiliary units.}

Let \(K=\mathbf Q(\mu_{N_c\ell})\) and suppose
a split prime \(r\equiv1\bmod N_c\ell\) lies outside \(T\).
Let \(x\) be an \(\ell\)-killed odd component, trivial on
the decomposition groups above \(T\) and nonzero at one
Frobenius over \(r\). On the decomposition groups over
\(r\) the action on its values is trivial, and \(x\) is
unramified there. Conjugation therefore shows
\(2(a^e-1)x(\phi)\ne0\) for an arithmetic Frobenius
\(\phi\) at any \(B\mid r\).

Take an inertia element
\(\tau\) at \(B\). Let \(\kappa\) be a dual cocycle for
\(\Upsilon=\chi_{\mathrm{cyc}}^{e+1}\psi\), killed by
\(\ell\) and unramified over \(\mathbf Q(\mu_{N_c})\) off
\(\ell,T\cup\{r\}\).
Write \(z=\kappa(\tau)\).
Take an \(\mathbf F_\ell\)-basis \(b_i\) of \(Oz\)
inside the \(\ell\)-torsion of the coefficient module
and extend its coordinate maps \(L_i\) to that whole
\(\mathbf F_\ell\)-space.
On \(G_K\) the cocycle action on these coefficients
is trivial.

With \(\zeta_\ell\) fixed, Hilbert's Theorem~90
\cite[\href{https://www.jmilne.org/math/CourseNotes/CFTc.pdf\#page=76}{Ch.~II, Proposition~1.22}]{MilneCFT}
applies over the finite Galois extension of \(K\) cut out by
\(L_i\circ\kappa|_{G_K}\). It gives \(r_i\) in that
extension with
\[
 g(r_i)/r_i=\zeta_\ell^{L_i(\kappa(g))}\qquad(g\in G_K).
\]
Hence \(y_i=r_i^\ell\in K^\times\).
Apply Artin reciprocity over \(K\) to \(y_i\) and
\(2(a^e-1)x|_{G_K}\). This unramified homomorphism is
zero on the decomposition groups above \(\ell\) by the
cyclotomic commutator calculation in
Section~\ref{sec:odd-positive-annihilator}, and on those
above \(T\) by hypothesis. Elsewhere off \(r\) the
root \(r_i\) is inertia-fixed, so \(y_i\) has normalized
order divisible by \(\ell\).

Such orders contribute zero
since the homomorphism is killed by \(\ell\). Thus only
the places above \(r\) contribute.
For a cyclotomic lift \(s_u\), with
\(u\in(\mathbf Z/N_c\ell)^\times\), put
\(m_i(u)=v_B(s_u y_i)\) with \(v_B(r)=1\).
Fix \(\pi^\ell=r\) and write
\(\tau(\pi)/\pi=\zeta_\ell^j\).
Since \(r\) splits in \(K\), \(m_i(u)\in\mathbb Z\).
The relations \(r_i^\ell=y_i\) and \(\pi^\ell=r\) give
\[
 v_B\bigl(\pi^{-m_i(u)}s_u(r_i)\bigr)
     =-m_i(u)/\ell+v_B(s_u y_i)/\ell=0.
\]
The \(\ell\)-th power of this local unit lies in \(K\), which
\(\tau\) fixes. Thus its quotient under \(\tau\) belongs to
\(\mu_\ell\) and reduces to \(1\) at \(B\) by inertia.
Reduction is injective on \(\mu_\ell\) since \(r\ne\ell\).
The equality of exponents is therefore
\[
 j m_i(u)=uL_i(\kappa(s_u^{-1}\tau s_u))\quad
                    \text{in }\mathbf F_\ell .
\]
Since \(\tau\) fixes \(\mu_{N_c\ell}\), its dual
coefficient action is trivial.

The cocycle identity
gives \(\kappa(s_u^{-1}\tau s_u)=\Upsilon(s_u)^{-1}z\).
Together these identities read in the chosen basis as
\[
 (u^e\psi(u))\,\sum_i j m_i(u)b_i=z.
\]
The weighted Artin identity gives
\[
 \left(\sum_u m_i(u)u^e\psi(u)\right)
                         2(a^e-1)x(\phi)=0 .
\]
Each scalar in parentheses thus lies in the maximal
ideal \(\mathfrak m_O\). Summing the basis identity
over \(u\) yields
\(\varphi(N_c\ell)z\in\mathfrak m_O Oz\).
Since \(\varphi(N_c\ell)\) is a unit in this range,
the cyclic module \(Oz\) vanishes by Nakayama.
This holds at all places over \(r\), so \(\kappa\) is
unramified off \(\ell,T\). If \(T\) is tight, it follows
that \(T\cup\{r\}\) is tight.

\subsubsection{The Jacobi derivative cocycles}
\label{sec:odd-jacobi-derivatives}

\sectionopening{In this subsection, we differentiate the Jacobi family built from the
distribution identity of Section~\ref{sec:odd-positive-annihilator}.
The resulting cocycles satisfy the relations needed for the detection
calculation of Section~\ref{sec:odd-large-detection}.}

We use Kolyvagin's derivatives of Gauss sums in the form discussed
by Rubin~\cite{RubinGauss}. Use \(N=N_c\), \(m=\ell^h\), \(L_1=Nm\),
and the primes
\(r\in Q'\), \(p\) and place \(B\) from the descent. In
the construction of \(J\) in Section~\ref{sec:odd-positive-annihilator}
use the chosen multiplier \(a\). All the arguments used
below have denominators dividing \(p-1\).
Since \((a,r)=1\), the distribution quotient
\(\prod_{b\bmod r}J(x+b/r)/J(rx)\) is a root of
unity. Its valuations disappear by the floor
distribution identity at \(B\), its conjugate
versions at all primes above \(p\), and
the units away from \(p\).

The equality of absolute values
\(|G_j|=|G_{-j}|\) at each embedding makes every absolute
value of the quotient one, so Kronecker's theorem applies.
The floor distribution follows by writing
each floor as \(a\langle x\rangle-\langle ax\rangle\)
and using
\(\sum_{b\bmod r}\langle x+b/r\rangle
  =(r-1)/2+\langle rx\rangle\) twice for
\(x\) and \(-x\). If \(r\mid p-1\) is prime and
\(r b=0\) in \(\mathbf Q/\mathbf Z\), the Gauss
sums at the shifted indices for \(x+b\) and
\(x\) are congruent at each place above \(r\).

The extra character takes \(r\)-th root
values, which reduce to \(1\), and the
Gauss sums are units there. Therefore
\(J(x+b)/J(x)\) is a principal unit there.
For the derivative system, let \(m=\ell^h\),
choose \(y(x)^m=J(x)\), and define \(F_g(x)\)
by
\[
 \zeta_m^{F_g(x)}
   =\frac{\gamma(y(x))}{(g\gamma g^{-1})(y(x))}.
\]
On the \(M_S=L_1\prod_{r\in S}r\)-torsion
this quotient has order dividing \(m\), since
the action of the commutator on \(J(x)\)
is trivial. Use the sum
\[
 \mathcal D_S F
   =\sum_{b\in(\mathbf Z/L_1)^\times}
      b^e\psi(b)
      \left[\prod_{r\in S}
       \sum_{u\in(\mathbf Z/r)^\times}
         \log_r(u)\,T_{u/r}\right]
             F\left(\frac{b}{L_1}\right),
 \qquad (T_bF)(x)=F(x+b).
\]
For each \(r\), \(\log_r\) is a fixed surjective homomorphism
to \(\mathbf Z/m\), available since \(m\mid r-1\).

The exponent \(F_g(x)\) is independent of the
choice of \(m\)-th root \(y(x)\): the two conjugate
Galois operators act alike on \(\mu_m\).
Write \(b_{g,S}\in(\mathbf Z/M_S)^\times\) for the exponent
of \(g\) on \(\mu_{M_S}\). For \(b=b_{g,S}\) and \(M_Sx=0\),
\(g(y(x))\) and \(y(bx)\) differ by an \(m\)-th root of
unity. Cancelling that root gives, for \(g'\in G_{\mathbf Q}\),
\[
 F_{gg'}(bx)-F_g(bx)=bF_{g'}(x).
\]
The quotient on the left is the \(g\)-translate of
\(\gamma y(x)/(g'\gamma g'^{-1})y(x)\).

Let \(T\) be a set of the allowed auxiliary primes and write
\(\mathcal D_T^{\mathbf d}\) for the same sum with weights
\(\log_s(u)+d_s\), \(\mathbf d\in(\mathbf Z/m)^T\).
For an auxiliary prime \(q\notin T\), set
\(\boldsymbol\lambda_T(q)=(\log_s(q))_{s\in T}\).
Since \(q\equiv1\bmod L_1\), the base indices are fixed
and substituting \(u'=qu\bmod s\) gives
\[
 \mathcal D_T^{\mathbf d}\bigl(F(q\,\cdot)\bigr)
       =\mathcal D_T^{\mathbf d-\boldsymbol\lambda_T(q)}F .
\]
The quotient \(\prod_{b\bmod q}y(x+b/q)/y(qx)\) is a
root of unity by the Jacobi distribution identity. The two
Galois actions defining \(F_g\) agree on roots of unity, so
\(\sum_{b\bmod q}F_g(x+b/q)=F_g(qx)\).
For \(S=T\sqcup\{q\}\), changing only \(d_q\) by
\(h_0\in\mathbf Z/m\) therefore gives
\[
 \mathcal D_S^{\mathbf d+h_0\mathbf e_q}F_g-
 \mathcal D_S^{\mathbf d}F_g
 =h_0\bigl(\mathcal D_T^{\mathbf d|_T-\boldsymbol\lambda_T(q)}F_g
                 -\mathcal D_T^{\mathbf d|_T}F_g\bigr),
\]
where \(\mathbf e_q\) is the coordinate vector at \(q\).
Starting with \(T=\varnothing\), induction removes all shifts.
For any \(b\) prime to \(M_S\), also substitute \(a'=ba\bmod L_1\)
in the base sum. The original base weight is
\((b^e\psi(b))^{-1}(a')^e\psi(a')\), and the same shift calculation
at the auxiliary primes gives
\[
 \mathcal D_S\bigl(F_g(b\,\cdot)\bigr)
       =(b^e\psi(b))^{-1}\mathcal D_S F_g .
\]
With \(b=b_{g,S}\), applying this identity to the equation for
\(F_{gg'}-F_g\) proves
\[
 \kappa_S(gg')-\kappa_S(g)
      =b^{e+1}\psi(b)\kappa_S(g'),\qquad
 \kappa_S(g)=\mathcal D_S F_g .
\]
Thus \(\kappa_S\) has coefficients in \(O/mO\) with character
\(\chi_{\mathrm{cyc}}^{e+1}\psi\). It is continuous because each
sum uses finitely many roots.
The cocycle at \(\gamma\) is zero termwise.
Away from \(\ell,p,S\), inertia fixing
\(\mu_N\) fixes \(\mu_{M_S}\). For \(g\) in
this inertia, \(g y(x)/y(x)\) is a principal
unit in \(\mu_m\), since \(y(x)\) is a unit
and its \(m\)-th power is fixed. The same
holds for \(g\gamma y(x)/\gamma y(x)\), since
\((\gamma y(x))^m=J(b_{\gamma,S}x)\). Here
\(b_{\gamma,S}\equiv a\bmod L_1\): \(\gamma\) has
\(\ell\)-adic cyclotomic value \(a\), fixes \(\mu_N\), and
\(a\equiv1\bmod N\). For each \(r\in S\) one has
\(b_{\gamma,S}\equiv1\bmod r\), since \(\gamma\) belongs
to inertia above \(\ell\ne r\).
Injectivity of reduction on \(\mu_m\) gives
\(F_g(x)=0\).

For the base derivative \(S=\varnothing\), take
\(x\in L_1^{-1}\mathbf Z/\mathbf Z\). Then
\(J(x)\in\mathbf Q(\mu_{L_1})\) and
\(b_{\gamma,\varnothing}x=ax\).
At the chosen prolongation \(B\mid p\), take a tame inertia element
\(\tau\) that multiplies the \(m\)-th root \(\varpi_m\)
of \(p\) by \(\zeta_m\). It exists by total ramification
of the uniformizer extension over the unramified field
containing \(\mu_m\). It fixes \(\mu_{L_1}\). Write
\[
 s_x=v_B(J(x))=
          \lfloor a\langle x\rangle\rfloor-\lfloor a\langle-x\rangle\rfloor .
\]
Then \(y(x)/\varpi_m^{s_x}\) is a unit whose \(m\)-th
power lies in \(\mathbf Q(\mu_{L_1})\). The quotient
of this unit under \(\tau\) is in \(\mu_m\) and
reduces to one. Therefore
\(\tau y(x)/y(x)=\zeta_m^{s_x}\), and likewise
\(\tau\gamma y(x)/\gamma y(x)=\zeta_m^{s_{ax}}\),
since \(\gamma(y(x))^m=J(ax)\).
It follows that \(F_\tau(x)=a s_x-s_{ax}\).
Summing over \(x=b/L_1\) and reindexing
\(b\) by \(ab\), one obtains
\[
 a^e\psi(a)\,\mathcal D_\varnothing F_\tau
  =\bigl(a^{e+1}\psi(a)-1\bigr)
             S_m(a,\psi).
\]
This proves the base Stickelberger identity
on the same \(B\).
\subsubsection{Comparison of Frobenius and inertia}\label{sec:odd-frobenius-inertia}

\sectionopening{In this subsection, we compute the Frobenius evaluation of the Jacobi
cocycles of Section~\ref{sec:odd-jacobi-derivatives} in terms of
inertia. This gives the precise local relation used with the
reciprocity calculation of Section~\ref{sec:odd-tight-reciprocity}.}

Retain the Jacobi family and \(m=\ell^h\).
Fix \(r\in Q'\setminus S\), any prolongation \(A\mid r\)
and any lift \(\phi\) of arithmetic Frobenius at \(A\).
Let \(x\) be an \(M_S\)-torsion point. Choose a unit
\(v\bmod r\) with \(\log_r(v)=1\) and choose
\(\tau\) in the inertia group at \(A\) acting
by \(v^{-1}\) on the \(r\)-torsion roots of
unity. It acts trivially on \(\mu_{M_S}\)
because \(r\nmid M_S\). Use additive notation
for \(\overline{\mathbf Q}^{\,\times}\) in the following
calculation. Write \(Y_b=y(x+b/r)\), \(Y_*=y(rx)\),
and \(\lambda(u)\in\mathbf Z\) for any integral
lift of \(\log_r(u)\).

For \(\delta\in G_{\mathbf Q}\) set
\[
 W_\delta=\sum_{b\bmod r}Y_b-(\delta\phi\delta^{-1})Y_0
       +\sum_{u\in(\mathbf Z/r)^\times}
            \lambda(u)\bigl((\delta\tau\delta^{-1})Y_{vu}-Y_u\bigr).
\]
At the conjugate place \(\delta A\), all \(Y_b\)
are units. Its inertia acts trivially on
their reductions and its Frobenius acts on
the residue field by the \(r\)-th power.
With \(Z_b=Y_b-Y_0\), reduction of \(W_\delta\)
modulo principal units there gives
\[
 \sum_{u\in(\mathbf Z/r)^\times}
       \bigl(1+\lambda(v^{-1}u)-\lambda(u)\bigr) Z_u .
\]
Each coefficient is divisible by \(m\);
the Jacobi congruence says \(m Z_u\) is a
principal unit at every place above \(r\).

Thus \(W_\delta\) is a principal unit at
\(\delta A\). This uses the Gauss-sum
congruence at the conjugate place as well
as at \(A\).
Set \(R_u=\tau Y_{vu}-Y_u\) and \(R_*=\phi Y_0-Y_*\).
The action of \(\phi\) on \(\mu_{M_S}\) is the
\(r\)-th power, as reduction is injective there.
Thus \(R_u,R_*\) are \(m\)-torsion by Jacobi
equivariance. The distribution quotient
gives that \(m(\sum_bY_b-Y_*)\) is torsion,
hence \(\sum_bY_b-Y_*\) is torsion.

It follows that \(W_1\) is a root of unity,
so \(\gamma W_1\) is a power of \(W_1\)
and is still a principal unit at \(A\).
Also \(\gamma W_{\gamma^{-1}}\) is a
principal unit at \(A\). The actions of
\(\gamma\) and any of its conjugates agree
on all roots of unity. Expand the difference
of the two terms, using \(\phi Y_0=Y_*+R_*\)
and \(\tau Y_{vu}=Y_u+R_u\), to obtain
\[
 \gamma W_1-\gamma W_{\gamma^{-1}}
  =\sum_u\lambda(u)
        \bigl(\gamma Y_u-(\tau\gamma\tau^{-1})Y_u\bigr)
       -\bigl(\gamma Y_*-(\phi\gamma\phi^{-1})Y_*\bigr).
\]
The \(R_u\) and \(R_*\) terms cancel by the
agreement on roots of unity.

The right side
is \(m\)-torsion and a principal unit.
Since \(r\ne\ell\), reduction on \(\mu_m\)
is injective. The exponent therefore gives
\[
 \sum_{u\in(\mathbf Z/r)^\times}
       \log_r(u) F_\tau(x+u/r)=F_\phi(rx).
\]
Here \(r\notin S\) is prime to \(M_S\) and
\(r\equiv1\bmod L_1\), so \(r^e\psi(r)=1\).
The scaling identity just proved integrates the local equation as
\[
 \mathcal D_{S\cup\{r\}}F_\tau
  =\mathcal D_S(F_\phi(r\,\cdot))
  =\mathcal D_SF_\phi .
\]
This gives \(\kappa_{S\cup\{r\}}(\tau)=\kappa_S(\phi)\) for
the given \(\phi\) and \(\tau\in I_A\). In particular the identity
holds for the prolongation and Frobenius lift selected by simultaneous
detection in Section~\ref{sec:odd-large-detection}.

\section{Polynomial bounds for the Selmer lengths}
\label{sec:selmer-polynomials}

\sectionopening{In this section, we represent the finite cocycle conditions by a matrix
over a polynomial ring in the twist variables. Its determinantal
ideals control the lengths from Section~\ref{sec:selmer-relations}
after specializing the twist.}

Retain the normalized valuation and length from
Section~\ref{sec:selmer-relations}. Fix a number field \(E\),
a finite set \(\Sigma\) of finite places, an integral
continuous representation \(\rho\) of \(G_E\), and finitely
many continuous additive characters
\(\chi_w:G_E\to\mathbb Z/\ell^N\mathbb Z\), where \(N\ge1\).
For a finite extension of coefficient rings \(\mathcal O\subset\mathcal O'\)
we use \(\lambda(C)=\operatorname{length}_{\mathcal O'}(C)/
\operatorname{ord}_{\mathcal O'}(\ell)\), so \(v(\ell)=1\)
in every extension.
For each \(g\in G_E\) set
\[
 R_X(g)=\rho(g)\prod_w(1+X_w)^{[\chi_w(g)]},
\]
where the brackets denote the representative in
\([0,\ell^N-1]\). Take a fixed collection of homogeneous linear
relations with coefficients in \(\mathcal O\), each involving
only finitely many evaluations.
Over the ring of integers \(\mathcal O'\) of a finite
extension of \(\operatorname{Frac}(\mathcal O)\), denote by
\(C(x)\) the continuous functions into
\((\mathcal O'/\ell^M)^n\) satisfying the cocycle equations
for \(R_x\), the chosen linear relations, and vanishing
on every inertia group at the finite places outside \(\Sigma\).
At every integral \(x\) the equation is
\(z(gh)=z(g)+R_x(g)z(h)\). It is a homogeneous linear
equation on the values of \(z\) even when exponent reduction
modulo \(\ell^N\) introduces a carry. At a character
specialization these are the usual cocycle equations.

\begin{lemma}[Polynomial minors for Selmer lengths]\label{lem:sharp-selmer}
There is a finite set of integral polynomials \(P_j\) such that
\[
            \lambda(C(x))=\min_j v(P_j(x))
\]
over every such coefficient ring and for each integral \(x\).
If \(j_0\) attains the minimum at \(x_0\), then
\[
 \lambda(C(x_0))=v(P_{j_0}(x_0)),\qquad
 \lambda(C(x))\le v(P_{j_0}(x))\quad\text{for every integral }x.
\]
\end{lemma}

\begin{proof}
Let \(H\) be the open normal subgroup in the kernels of
all \(\chi_w\) and of \(\rho\bmod\ell^M\), and let \(E_H\)
be its fixed field. Every \(z\in C(x)\) is additive on \(H\),
since \(R_x(h)=I\bmod\ell^M\) there. The continuous
homomorphisms from \(H\) to \(\mathbb Z/\ell^M\) unramified
outside the places above \(\Sigma\) form a finite set.
Indeed their kernels cut out cyclic extensions of \(E_H\)
of degree at most \(\ell^M\) with this ramification bound,
and those form a finite set by Hermite's finiteness theorem
\cite[Theorem 8.42]{MilneANT}. Each such field admits finitely
many homomorphisms. Let \(H_1\) be the core in \(G_E\) of the intersection
of their kernels in \(H\). The dual into
\(\mathbf Z/\ell^M\mathbf Z\) separates points of the finite
abelian group \((\mathcal O'/\ell^M)^n\). Applying its
characters to each \(z|_H\) gives
\[
 z|_{H_1}=0,\qquad
 z(gh)=z(g)+R_x(g)z(h)=z(g)\quad(h\in H_1).
\]
Thus \(z\) factors through \(G_E/H_1\).

The finite quotient \(G_E/H_1\) depends only on the original
data and \(M,N\), and works for every \(\mathcal O'\) and \(x\).
In the original cocycle, inertia and linear equations replace each
value \(z(g)\) by the coordinate indexed by \(gH_1\). Conversely,
any solution on the quotient pulls back to a continuous solution
on \(G_E\). The rows of these equations generate a submodule of
a finite free \(\mathcal O[X]\)-module. Since \(\mathcal O[X]\)
is Noetherian, finitely many of the rows span it.
Let \(A(X)\) be their matrix, with \(q\) columns. For every
specialization the kernel of \(A(x)\) acting on
\((\mathcal O'/\ell^M)^q\) is exactly \(C(x)\), since all the
equations are equivalent to this finite selection after
evaluation. Put
\[
 A^+(X)=\begin{pmatrix} A(X)\\ \ell^M I_q\end{pmatrix},
 \qquad
 C(x)=\ker\!\left(A^+(x)\bmod\ell^M\right).
\]
The bottom \(q\) rows give the constant maximal minor
\(\ell^{Mq}\), so \(A^+(x)\) has full column rank at every
specialization. For \(q=0\) the empty minor is \(1\), and the
length is zero. For \(q>0\) its Smith invariants \(u_1,\ldots,u_q\) divide \(\ell^M\), since its
row module contains \(\ell^M(\mathcal O')^q\). Multiplication by
\(u_i\) on \(\mathcal O'/\ell^M\) has kernel
\((\ell^M/u_i)\mathcal O'/\ell^M\mathcal O'\) of normalized
length \(v(u_i)\). The maximal minors \(P_j(x)\) generate the
ideal \((\prod_i u_i)\) by Smith reduction. Hence
\[
 \lambda(C(x))=\sum_{i=1}^q v(u_i)=\min_j v(P_j(x)).\qedhere
\]
\end{proof}

\subsection{The line through the untwisted point}

\sectionopening{In this subsection, we restrict the polynomial matrix of
Section~\ref{sec:selmer-polynomials} to a line through the untwisted
specialization. At zero the distinguished pair has already been
released. The determinant product bounds the cost of releasing
the relations at the other auxiliary pairs along this line.}

Let \(w_0\) be the distinguished place, and fix integral slopes
\(s_w\in\mathcal O\) for \(w\ne w_0\). Define the line by
\[
 x_{w_0}(Y)=Y,\qquad x_w(Y)=s_wY\quad(w\ne w_0).
\]
The first coordinate has coefficient one, while the untwisted
local relation at \(w_0\) used slope zero. At all other pairs assume
the Frobenius and inertia values, norm congruences, and the conditions
on \(\rho\bmod\ell^M\) of Lemma~\ref{lem:selmer-relations}.
Write \(C_0(x)\) for the module obtained by imposing
the auxiliary relations at these pairs, so the first
place is excluded from this set of relations. Suppose
\(\lambda(C_0(0))\le\beta\). This is the untwisted
length bound after accounting for the first place and its
conjugate.

Let \(J_{\rm aux}\) index the auxiliary pairs where the
relations are imposed. For each \(w\in J_{\rm aux}\), put
\(\phi_{w,1}=\phi_w\), \(\phi_{w,2}=(\phi_w^c)^{-1}\). Both
evaluations of Lemma~\ref{lem:selmer-relations} enter the product
\[
 D_{\rm aux}(Y)=
   \prod_{w\in J_{\rm aux}}\prod_{i=1}^2
                     \det((1+Y)\rho(\phi_{w,i})-I).
\] Choose by Lemma~\ref{lem:sharp-selmer}
a polynomial \(G\) for \(C_0\) that is sharp at \(0\).
For every \(y\) satisfying the scalar condition of
Lemma~\ref{lem:selmer-relations}, the exact sequence for
\(C_0(x)\subset C(x)\) gives
\[
   \lambda(C(x))\le v(G(x))+v(D_{\rm aux}(y)).
\]
The degree of \(D_{\rm aux}\) is \(D_0=2n|J_{\rm aux}|\) and
its leading coefficient \(\prod_{w,i}\det\rho(\phi_{w,i})\)
is a unit. The bound on \(|J_{\rm aux}|\) from
Section~\ref{sec:auxiliary-descent} is independent of \(N\).
Fix a positive rational radius \(a_0\) and require
\(N\min(1,a_0)\ge M\) before selecting the auxiliary places at
level \(\ell^N\). For \(1\le j\le\ell^N\) put \(k=v_\ell(j)\).
As \(k\le N\) and \(j-1\ge k\), for \(v(z)\ge a_0\) we have
\[
 v\!\left(\binom{\ell^N}{j}z^{j-1}\right)
 \ge N-k+(j-1)a_0\ge N\min(1,a_0).
\]
This gives the scalar condition of Section~\ref{sec:selmer-relations}
on the whole line in every coefficient extension.

For \(q(Y)=\sum_k q_k Y^k\) put
\(v_{a_0}(q)=\min_k(v(q_k)+ka_0)\). With \(g(Y)=G(x(Y))\), the
constant coefficient of \(g\) and the leading coefficient of
\(D_{\rm aux}\) give
\[
 v_{a_0}(g)\le v(G(0))\le\beta,\qquad
 v_{a_0}(D_{\rm aux})\le D_0 a_0 .
\]
Choose \(y_0\) of valuation \(a_0\). For \(q=g,D_{\rm aux}\), divide
\(q(y_0U)\) by a coefficient of least valuation to obtain a
nonzero reduction. Lift a nonzero residue avoiding the zeros of both
reductions to a unit \(u\) in a finite coefficient extension.
At \(y=y_0u\) both evaluations attain their Gauss valuations.
Choose a minor \(P\) for \(C\) attaining its minimum at \(x=x(y)\).
It satisfies
\[
 v(x_w)\ge a_0,\qquad v(P(x))\le\beta+D_0a_0,
 \qquad \lambda(C(x'))\le v(P(x'))\quad\text{for every }x'.
\]

\subsection{Transfer to divisible coefficients}

\sectionopening{In this subsection, we apply the finite-length estimate of
Section~\ref{sec:selmer-polynomials} to divisible coefficient modules
at torsion specializations. We check that the same inertia and linear
conditions identify the finite modules with their torsion subgroups.}

At a tuple of \(\ell^N\)-power torsion points
\(x'_w=\zeta_w-1\) with \(\zeta_w^{\ell^N}=1\), the
specialized twist is a scalar character. Let \(C_\infty\)
be the module of cocycles for this twist with coefficients in
\((\operatorname{Frac}(\mathcal O')/\mathcal O')^n\),
with the same inertia and linear conditions. Write \(\varpi'\)
for a uniformizer and \(e'=\operatorname{ord}_{\mathcal O'}(\ell)\).
The equivariant coefficient isomorphism
\[
 \mathcal O'/\ell^M\mathcal O'\longrightarrow
 (\operatorname{Frac}(\mathcal O')/\mathcal O')[\ell^M],
 \qquad b\longmapsto b/\ell^M\pmod{\mathcal O'}
\]
preserves all homogeneous relations, hence identifies \(C(x')\)
with \(C_\infty[\ell^M]\). Every continuous cocycle has finite
image since \(G_E\) is compact and its target is discrete.
If a cocycle has annihilator \((\varpi'^b)\) with \(b>e'M\),
multiplication by \(\varpi'^{b-e'M}\) gives a cyclic submodule
of \(C(x')\) of normalized length \(e'M/e'=M\). Hence, for the
full set of maximal minors in Lemma~\ref{lem:sharp-selmer},
\[
 \min_j v(P_j(x'))<M\quad\Longrightarrow\quad
 C_\infty=C_\infty[\ell^M],\qquad
 \lambda(C_\infty)=\lambda(C(x'))=\min_j v(P_j(x')).
\]
In particular, the minor \(P\) chosen on the test line gives
\(\lambda(C_\infty)\le v(P(x'))\) whenever \(v(P(x'))<M\).
This holds over every finite extension containing the coefficients
and the torsion tuple.

\section{Simultaneous evaluation on torsion}
\label{sec:torsion-specialization}

\sectionopening{In this section, we choose torsion specializations at which finitely
many polynomials satisfy their valuation bounds simultaneously. The
choice will apply to the Selmer minors of
Section~\ref{sec:selmer-polynomials} and the group-ring evaluations of
Section~\ref{sec:joint-twist}.}

Let \(v\) be the additive valuation on \(\overline{\mathbb Q}_\ell\),
normalized by \(v(\ell)=1\), and write \(\mathcal O\) for its
valuation ring. The following evaluation lemma gives a torsion
exponent before the auxiliary primes or their polynomials are
chosen.

\begin{lemma}[Simultaneous torsion specialization]\label{lem:common-torsion-specialization}

Fix a variable bound \(s_0\in\mathbb Z_{\ge0}\), an integer
\(m\ge1\), a positive rational radius \(a\), a tolerance
\(\epsilon>0\), and a finite list of real bounds \(B_i\ge0\).
Use the convention \(v(0)=+\infty\). Set
\[
 \rho=a/m,\qquad
 B=m\sum_i B_i,\qquad
 \delta=\frac{\min(1,\epsilon)}{s_0+2},
 \qquad \rho_h=\varphi(\ell^h)^{-1}\quad(h\ge1).
\]
Let \(H_1\) be the least positive integer with \(\rho_{H_1}\le\rho\).
Any \(H\ge H_1\) satisfying
\[
 B+1 < \delta\sum_{h=H_1}^{H}\frac1{h+1}
\]
works. For example one can take
\[
 H=\left\lceil(H_1+1)\exp((B+1)/\delta)\right\rceil
\]
since
\[
 \sum_{h=H_1}^H\frac1{h+1}
 \ge\log\frac{H+2}{H_1+1}>\frac{B+1}{\delta}
\]
for this choice.
For every family of integral polynomials \(P_i\) in the same
\(s\le s_0\) variables and points \(x_i\) such that
\[
 v(x_{i,w})\ge a,\qquad v(P_i(x_i))\le B_i,
\]
there are roots of unity \(\zeta_w^{\ell^H}=1\) with
\(v(\zeta_w-1)\le \rho\) and
\[
     v\big(P_i((\zeta_w^j-1)_w)\big)
                       \le v(P_i(x_i))+\epsilon
                  \qquad(1\le j\le m,\ \text{all }i).
\]

\end{lemma}

\begin{proof}

\indent\textit{One coordinate.}
For a nonzero integral polynomial \(g\), write
\[
 v_t(g)=\min_b\{v(g_b)+bt\}.
\]
This valuation is multiplicative on products in additive
notation. To see this for rational \(t\), scale the variable
by an element of valuation \(t\), then divide each scaled
polynomial by a coefficient of minimum valuation. The reductions
are nonzero polynomials over a field and have nonzero product.
By factorization, with leading coefficient \(g_{\rm top}\),
\[
      v_t(g)=v(g_{\rm top})+\sum_\alpha\min(t,v(\alpha)),
\]
where roots are counted with multiplicity. In particular
\(v_t(g)\le v(g(y))\) for \(v(y)\ge t\), and \(v_t(g)\)
is increasing in \(t\).

Suppose \(v_\rho(g)\le B+1\). The polynomial is integral, so \(v_0(g)\ge0\). Subtracting
the root formula at radius zero from that at radius \(\rho\)
gives a sum of the nonnegative terms
\(\min(\rho,v(\alpha))-\min(0,v(\alpha))\).
Since the radii \(\rho_h\le\rho\) are distinct, the multisets
\(R_h\) of roots with \(v(\alpha)=\rho_h\) obey
\[
 \sum_{h=H_1}^{H}|R_h|\rho_h\le v_\rho(g).
\]
Thus there is an \(h\) with
\(|R_h|\rho_h(h+1)\le\delta\), by the harmonic-sum inequality.
For each \(\alpha\in R_h\), the polynomial
\[
            F_\alpha(Y)=(1+\alpha+Y)^{\ell^h}-1
\]
is monic and has linear coefficient of valuation \(h\),
because \(v(\alpha)>0\).

Its roots include
\(Y=\zeta-1-\alpha\) for the primitive \(\ell^h\)-th roots
\(\zeta\). All roots of \(F_\alpha\) are of positive valuation.
Applying the Gauss valuation at radius \(1\),
\[
 \sum_{\zeta\ {\rm primitive}}
               \min\{1,v(\zeta-1-\alpha)\}
            \le v_1(F_\alpha)\le h+1 .
\]
There are \(1/\rho_h\) such roots of unity. Sum over \(\alpha\in R_h\)
and average to find one \(\zeta\) with the sum of the truncated
valuations at most \(|R_h|\rho_h(h+1)\le\delta\).

Since \(\delta<1\), each truncation is inactive. All primitive
roots of this order have \(v(\zeta-1)=\rho_h\), by taking
the norm of \(1-\zeta\) or using the cyclotomic polynomial
at \(1\). For \(y=\zeta-1\), the roots of \(g\) with valuation
different from \(\rho_h\) have
\(v(y-\alpha)=\min\{\rho_h,v(\alpha)\}\). The other roots give
\[
       v(g(y))-v_{\rho_h}(g)
           =\sum_{\alpha\in R_h}
                    \big(v(y-\alpha)-\rho_h\big)\le\delta.
\]

Now let there be finitely many polynomials \(g_i\), with their
old arguments \(u_i\) of valuation at least \(\rho\)
and \(\sum_i v(g_i(u_i))\le B+1\). Apply this argument to
\(g=\prod_i g_i\).
Since \(v_{\rho_h}(g_i)\le v_\rho(g_i)\le v(g_i(u_i))\) and
\(v(g_i(y))\ge v_{\rho_h}(g_i)\), multiplicativity gives
\[
 \sum_i \max\{0,v(g_i(y))-v(g_i(u_i))\}
 \le \sum_i \bigl(v(g_i(y))-v_{\rho_h}(g_i)\bigr)
 =v(g(y))-v_{\rho_h}(g)\le\delta .
\]

\indent\textit{All coordinates and powers.}
For each pair \((i,j)\), \(1\le j\le m\), replace \(P_i(X)\)
by the integral polynomial
\[
            Q_{i,j}(Y)=P_i(((1+Y_w)^j-1)_w).
\]
For each old coordinate \(x=x_{i,w}\) the monic polynomial
\((1+Y)^j-(1+x)\) has \(j\) integral roots whose product
up to sign is \(x\). For \(x\ne0\), the sum of their valuations
is \(v(x)\), so one root has valuation at least \(a/j\ge\rho\).
For \(x=0\) choose the root \(0\). Use these roots as the test point of \(Q_{i,j}\).
Its evaluation is \(P_i(x_i)\); the sum of the initial
valuations over \((i,j)\) is at most \(B\).

Process the coordinates one at a time and freeze the others at
their current values for each test point. On the coordinate
currently being replaced the old valuations are at least \(\rho\).
The sum of the evaluation valuations after \(k\) replacements
is at most \(B+k\delta\le B+1\). Hence the one-coordinate
argument applies with the same \(H\) at every step, including
when the previous replacements have valuations below \(\rho\).
They remain integral, so the resulting one-variable polynomials
have integral coefficients. Each is nonzero since its evaluation
at the current test point has finite valuation. The common
final points are \((\zeta_w-1)_w\). Summing the positive changes bounds
the increase at each \((i,j)\) by \(s\delta\le\epsilon\). This proves the assertion.
\end{proof}

\section{A common twist for the Selmer and Euler-factor conditions}
\label{sec:joint-twist}

\sectionopening{In this section, we choose one twist for the cocycle conditions of
Section~\ref{sec:selmer-polynomials} and the Euler-factor relations of
the norm systems. The auxiliary primes and torsion specialization are
arranged so that both estimates apply to this same character.}

Fix a coefficient place \(\kappa=(\ell,\iota)\), normalize
\(v(\ell)=1\), and take a finite family of the Selmer data of
Section~\ref{sec:selmer-polynomials}. Use a common finite coefficient
DVR for their lattices. In this section \(N\) denotes the exponent of
the auxiliary level \(\ell^N\).
Fix also finitely many norm systems, each consisting of a family
\((\Theta_Q)_Q\) of integral group-ring elements over
\(\overline{\mathbf Q}_\ell\), indexed by the allowed auxiliary sets.
For each system let \(\phi_w\) denote a lift of geometric Frobenius at
the new place \(w\), in the reciprocity convention of
Section~\ref{sec:first-auxiliary}. For a character \(\psi\)
factoring through level \(Q\), require that its evaluation at \(\Theta_{Q\cup\{w\}}\)
is its evaluation at \(\Theta_Q\) multiplied by
\[
 E_w(\vartheta(\phi_w)\psi(\phi_w)),\qquad
 E_w(A)=\prod_{e\in\mathcal E}
              (1-Aq_w^{-e})(1-A^{-1}q_w^{-e}).
\]
Here \(q_w\) is the residue cardinality, \(\vartheta:G_E\to\mathcal O^\times\) is the continuous
character of that system, unramified at the allowed auxiliary places,
and \(\mathcal E\) is its fixed finite list of integer exponents.
Thus \(\vartheta(\phi_w)\) is independent of the lift of Frobenius
at each such place. In the application
to the signed systems of Section~\ref{sec:critical-signed-systems},
\(\psi\) is the inverse \(\ell\)-adic avatar of the auxiliary
norm-one character. Both signs are included among the systems before
the choice of \(\psi\).

\subsection{A specialization of a norm system}

\sectionopening{In this subsection, we express each group-ring evaluation from
Section~\ref{sec:joint-twist} as a polynomial in the chosen character
coordinates. The norm relation determines the divisibility along the
line used for the simultaneous specialization.}

Use one list \(J\) of auxiliary primes with distinguished
prime \(w_0\), equipped with characters \(\chi_w\) of order
dividing \(\ell^N\) as in Section~\ref{sec:selmer-relations}.
Write \(\chi_w:G_E\to\mathbf Z/\ell^N\mathbf Z\) additively.
It comes from a character of the norm-one classes at level \(w\),
is zero on the decomposition groups above
\(\ell\), and takes the value \(1\) at the chosen
tame generator at \(w\). The other characters are zero there.
At the Frobenius and its conjugate inverse for
every subsequent prime, the first character has value \(1\)
and the other characters are zero. All \(q_w\) are
\(1\bmod\ell^N\). For a norm-one class represented by
\(u/\bar u\), read \(\chi_w(\operatorname{Art}u)\)
using the induced class character. Put \(h=\ell^N\) and take
\([c]_h\in\{0,\ldots,h-1\}\) for \(c\in\mathbf Z/h\mathbf Z\).
The assignment on class elements
\[
 [u/\bar u]\longmapsto
    \prod_{w\in J}(1+X_w)^{[\chi_w(\operatorname{Art}u)]_h}
\]
is multiplicative in
\(\mathcal O[X_w:w\in J]/((1+X_w)^h-1:w\in J)\).
Define \(P_J\) in the polynomial ring by summing these displayed
representatives with the coefficients of \(\Theta_J\). Each variable
has degree less than \(h\). For \(X_w=\zeta_w-1\) with \(\zeta_w^h=1\)
the value is the evaluation at \(\prod_w\zeta_w^{\chi_w}\).

Restrict \(P_J\) to the line with \(X_{w_0}=Y\) and
all other variables zero, and call the resulting polynomial
\(P_J^0(Y)\). Let \(R_0(Y)\) denote the polynomial for the
same one-variable evaluation of \(\Theta_{\{w_0\}}\).
Put \(j_0=|\mathcal E|(|J|-1)\). With \(A_w=\vartheta(\phi_w)\),
the product
\[
 R_1(Y)=\prod_{w\ne w_0}\prod_{e\in\mathcal E}
       \bigl(1-A_w q_w^{-e}(1+Y)\bigr)
       \bigl(1+Y-A_w^{-1}q_w^{-e}\bigr)
\]
has degree \(2j_0\) and unit leading coefficient.
The norm relation iterated away from \(w_0\) gives
\[
 (1+Y)^{j_0}P_J^0(Y)\equiv R_0(Y)R_1(Y)
                     \pmod{(1+Y)^{\ell^N}-1}.
\]
Indeed at \(Y=\zeta-1\) the character is \(\zeta^{\chi_{w_0}}\).
Its value on the Frobenius at \(w\ne w_0\) is
\(\zeta\), so each pair of factors of \(R_1\) equals
\(\zeta\) times the Euler-factor pair. This gives equality
at every root of the modulus, which is monic and has simple
roots. Polynomial division gives the congruence over the
valuation ring.

For the first-prime Euler factor, fix \(g_0\) with the
following finite bound
\[
 \beta=v(1(\Theta_\varnothing))
                 +2|\mathcal E|v(\vartheta(g_0)-1).
\]
The second term is zero for an empty exponent list.
Choose \(M'>\beta+1\), and require \(N\ge M'\) and
\(v(\vartheta(\phi_{w_0})-\vartheta(g_0))\ge M'\).
If \(\mathcal E\) is nonempty, each of the two factors
for each \(e\) at the first prime and trivial character
has valuation \(v(\vartheta(g_0)-1)\). Indeed \(q_{w_0}^{-e}\)
is \(1\bmod\ell^{M'}\), and \(\vartheta(\phi_{w_0})\)
and its inverse agree with \(\vartheta(g_0)\) and its
inverse, respectively, to that precision.
For a nonempty list the integrality of the augmentation and
\(M'>\beta\) imply \(v(\vartheta(g_0)-1)<M'\).
For an empty list the Euler product is one. In both cases the norm
relation gives \(v(R_0(0))=\beta\).

Let \(0<a_0\le1\) be rational. The Gauss valuation
of \(R_0R_1\) at radius \(a_0\) is at most
\(\beta+2j_0a_0\), by the constant coefficient of \(R_0\)
and the leading coefficient of \(R_1\).
For \(v(y)=a_0\), the binomial expansion shows
\[
 v((1+y)^{\ell^N}-1)\ge Na_0 .
\]
Indeed, for \(1\le j\le\ell^N\), put \(k=v_\ell(j)\le N\).
The identity \(j\binom{\ell^N}{j}=\ell^N\binom{\ell^N-1}{j-1}\)
gives
\[
 v\!\binom{\ell^N}{j}+ja_0\ge N-k+ja_0
      \ge Na_0,
 \qquad (N-k)(1-a_0)+(j-k)a_0\ge0.
\]
The coefficient bound just proved gives
\(v_{a_0}((1+Y)^{\ell^N}-1)\ge Na_0\).
The quotient in the cyclic polynomial congruence is integral, so
the difference of its two sides has Gauss valuation at least
\(Na_0\). Since \(v_{a_0}(1+Y)=0\), whenever
\(\beta+2j_0a_0<Na_0\) this gives
\[
 v_{a_0}(P_J^0)=v_{a_0}(R_0R_1)\le\beta+2j_0a_0 .
\]
For finitely many such systems \(i\) on \(J\), choose \(c\) of
valuation \(a_0\) and divide each \(P_{i,J}^0(cZ)\) by a coefficient
of minimum valuation. All these reductions are nonzero. Their product
has finitely many roots in \(\overline{\mathbf F}_\ell\). Lift a nonzero
residue outside this set to a unit \(u\), passing to a finite
unramified extension of the coefficient field containing \(c\) and the
polynomials if needed. At \(y=cu\) we obtain simultaneously
\[
 v(P_{i,J}^0(y))=v_{a_0}(P_{i,J}^0)
      \le\beta_i+2|\mathcal E_i|(|J|-1)a_0 .
\]
The other coordinates are zero on this line. The degree of the
product of reductions need only be known after \(J\) is chosen.

\subsection{The distinguished auxiliary prime}
\label{sec:first-auxiliary}

\sectionopening{In this subsection, we use the class-field conditions of
Section~\ref{sec:joint-twist} to choose the first auxiliary prime. Its
Frobenius lies in the specified neighborhood, and its local character
has the prescribed restrictions at the other places.}

Let \(E\) be CM and let \(g_0\in G_E\) lie in the normal
subgroup generated by inertia at a finite set of places
\(R_{\rm in}\) away from \(\ell\). Fix a neighborhood
\(\mathcal W\) of \(g_0\), a finite set of forbidden
primes, and a finite set \(S_{\rm sp}\) of places at which
the new character must be trivial on decomposition groups.
Require the conjugate closure of \(S_{\rm sp}\) to be
disjoint from \(R_{\rm in}\).
For each \(N\ge1\) there is a prime of prime residue
cardinality \(1\bmod\ell^N\), outside the forbidden set
with its conjugate and distinct from its conjugate,
equipped with a character of the kind used in
Section~\ref{sec:joint-twist}. Its geometric Frobenius lift
can be chosen in \(\mathcal W\) with character value zero.

To prove this, choose a finite set \(T\), stable under
conjugation, containing \(S_{\rm sp}\) and the places
above \(\ell\), and avoiding \(R_{\rm in}\) and its
conjugate. Require the primes in \(T\) to generate the
ideal class group.
Such a set exists by Chebotarev in the Hilbert class
field, which gives a prime in every ideal class outside
any finite set \cite[Ch.~V, Theorems~3.23 and~5.5]{MilneCFT}. Let \(\mathbb A_{E,f}^\times(T)\) denote the finite ideles
that are units outside \(T\). Since \(T\) generates the class group,
the ideal of any finite idele differs from an ideal supported on
\(T\) by a principal ideal. Thus
\[
 \mathbb A_{E,f}^\times=E^\times\mathbb A_{E,f}^\times(T),
 \qquad E^\times\cap\mathbb A_{E,f}^\times(T)
                      =\mathcal O_{E,T}^\times .
\]

Put \(h=\ell^N\) and take the field obtained by adjoining
\(\mu_h\) and the \(h\)-th roots of a finite generating
set of the \(T\)-units. The generating set exists by
the unit theorem, since valuation at \(T\) has kernel
the ordinary unit group and image in \(\mathbb Z^T\).
This Kummer extension is Galois over \(E\) and unramified
at each place of \(R_{\rm in}\): the \(T\)-units are units
there and \(h\) is prime to the residue characteristic.

Their \(h\)-th roots and \(\mu_h\) can be obtained in
an unramified extension by lifting the distinct roots
over a residue field extension. Choose a finite Galois
extension of \(E\) such that the coset of the kernel
containing \(g_0\) lies in \(\mathcal W\), and enlarge it
by the Kummer extension. Prescribe the conjugacy class of \(g_0^{-1}\)
for arithmetic Frobenius in this compositum. Chebotarev gives
\(w_0\) and a place above it with this arithmetic element in the
quotient \cite[Ch.~V, Theorem~3.23]{MilneCFT}. Its inverse lifts to
geometric Frobenius in \(\mathcal W\). We may avoid \(T\),
the forbidden primes and their conjugates, and the
primes ramified over \(\mathbb Q\). We may also take
\(w_0\) of residue degree one over \(\mathbb Q\),
since prime ideals of higher absolute residue degree
have density zero. Then \(w_0\ne\bar w_0\),
as an unramified prime fixed by the nontrivial
conjugation of \(E/E^+\) has relative residue degree two.

The Frobenius of \(w_0\) fixes the Kummer extension, so
\(q_{w_0}=1+hk\) for an integer \(k\).
For each generating \(T\)-unit \(u\), its root \(u^{1/h}\)
reduces to a nonzero element fixed by the residue
Frobenius. Hence the reduction of \(u\) is an \(h\)-th
power in the residue field of \(w_0\). The discrete
logarithm of residue units modulo \(h\) kills these
\(h\)-th powers and therefore all \(T\)-units.

The residue logarithm on the \(w_0\)-component of a
\(T\)-idele kills \(\mathcal O_{E,T}^\times\), so the displayed
quotient identifies it with a character \(\theta\) of the finite idele
classes. Its kernel contains an open subgroup of finite ideles, giving
continuity. Extend it trivially at the complex places. By reciprocity
\cite[Ch.~I, Theorem~1.1 and Ch.~V, Theorems~5.3 and~5.5]{MilneCFT},
\(\theta(u/\bar u)\) defines an anticyclotomic Galois
character, trivial on the decomposition groups at
\(S_{\rm sp}\) and above \(\ell\), and unramified away
from \(w_0,\bar w_0\).
Indeed a local idele supported at a place of \(T\) and its
conjugate have unit component \(1\) at \(w_0\), even for a
local uniformizer. Thus the character vanishes on the full
decomposition groups there. If \(u\) is a unit idele supported at
\(w_0\), the \(w_0\)-component of \(u/\bar u\) is \(u_{w_0}\), since
\(w_0\ne\bar w_0\). Its residue logarithm maps onto \(\mathbb Z/h\),
including for \(h=2^N\). Normalize the logarithm to read \(1\) on the
chosen tame generator \(\gamma_{w_0}\). It then reads \(1\) on
\((c\gamma_{w_0}c^{-1})^{-1}\) as well.

Multiplying the Frobenius lift by a power of this tame
generator makes the character value of Frobenius zero.
This correction does not alter the chosen finite Frobenius,
as the compositum used for Chebotarev is unramified
at \(w_0\). The corrected lift still lies in \(\mathcal W\).

\subsection{Selmer length descent and the uniform bound}
\label{sec:auxiliary-descent}

\sectionopening{In this subsection, we iterate the local relations of
Section~\ref{sec:selmer-relations} until the remaining module can be
treated by Section~\ref{sec:terminal-duality}. The argument bounds the
number of added primes independently of the larger torsion order.}

Fix the coefficient ring \(\mathcal O\), its uniformizer
\(\varpi\), the finite data indexed by \(i\) of ranks \(n_i\),
and \(M\). Set \(e=\operatorname{ord}_{\mathcal O}(\ell^M)\).
If necessary, extend \(\mathcal O\) by a finite unramified
extension so that its residue field has more than \(2en_i+1\)
elements for every \(i\). This preserves \(e\) and the ranks.
Use this ring throughout the descent, before computing the initial
length bound or choosing the auxiliary exponent \(N\). In this subsection lengths are unnormalized.
For an auxiliary list with prescribed congruences \(J\) and \(1\le a\le e\), let
\(C_{i,a}(J)\) be the untwisted cocycles modulo \(\varpi^a\)
with the original linear conditions, allowed ramification
at \(J\) and its conjugates, and the given congruences at each
auxiliary pair. Let \(C_{i,a}(\varnothing)\) have just the
original ramification and linear conditions. The lengths
are finite by the finite-presentation argument.

Apply Section~\ref{sec:first-auxiliary} at exponent \(N>M\)
to the inertial reference \(g_0\) in the prescribed common
neighborhood \(\mathcal W\), with \(S_{\rm sp}=\varnothing\).
The neighborhood and its precisions are specified below
before increasing \(N\).
Forbid every place above the rational primes dividing the
norms of the original bad places. Start with slope zero at both members of the distinguished
pair. For each datum \(i\), depth \(a\), and eligible new pair,
consider the following condition.

The two Frobenius values
of every member of \(C_{i,a}(J)\), at that prime and its
conjugate inverse, can be realized as the two tame values
of a cocycle with the old congruences and ramification also
allowed at the new prime pair. If this condition holds for every eligible pair for \(i\)
at depth \(a\), and some cocycle in \(C_{i,a}(J)\) has a unit
coordinate at some group element, Section~\ref{sec:terminal-duality}
bounds the length by \(G_i+L_i\), where \(G_i,L_i\)
are the datum's global and original-place terms, in
unnormalized length.

An eligible prime \(P\) has prime
norm, is distinct from its conjugate, and both \(P\) and \(\bar P\)
are outside the forbidden places and the previously chosen pairs. The prime \(P\) does not lie over
the rational norm of \(w_0\). Its character is trivial
on the decomposition groups above \(\ell\) and at both members
of each old pair other than \(\{w_0,\bar w_0\}\), and vanishes on its
own Frobenius lift.

On this lift, the old character at
\(w_0\) reads \(1\) and the other old characters read \(0\).
The new prime has the same tame-generator and norm
conditions as the original list. For the selected datum
\(i\), the representation modulo \(\ell^M\) is trivial on
both Frobenius parameters. This makes every nonzero slope
admissible in that datum.

Here is the local calculation used to choose a slope.
At an eligible new prime pair let \(C\) be the enlarged
group with the old congruences, and for this datum and depth write
\(A,B:C\to(\mathcal O/\varpi^a)^{2n_i}\) for the two tame
and the two Frobenius evaluations respectively. Put
\[
 C'=\sum_{u\in\mathcal O}\ker_C(uA-B).
\]
For \(c\in\ker(uA-B)\) one has \(Bc=A(uc)\), so
\(B(C')\subseteq A(C')\). Let \(L'\) be the inverse
image of \(A(C')\) in \(\mathcal O^{2n_i}\), with reduction
\(\pi:L'\twoheadrightarrow A(C')\).

It is free of rank \(2n_i\) since it contains
\(\varpi^a\mathcal O^{2n_i}\).
Lift \(\pi\) through \(A:C'\twoheadrightarrow A(C')\)
to \(\sigma:L'\to C'\), then lift \(B\sigma\) through \(\pi\)
to an endomorphism \(U'\) of \(L'\).
This gives \((uA-B)\sigma=\pi(uI-U')\).
Exclude the roots of the characteristic polynomial
of \(U'\) modulo \(\varpi\).
There are at most \(2n_i\) excluded residues. Outside
these residues \(uI-U'\) is invertible, so the image
of \(uA-B\) on \(C'\) equals \(A(C')\).

Since its kernel
on \(C\) is contained in \(C'\), the length identity for
the two surjections onto \(A(C')\) gives
\[
 \operatorname{length}\ker_C(uA-B)
   =\operatorname{length}(C'\cap\ker A)
   \le \operatorname{length}\bigl(\ker_C A\cap B^{-1}(A(C))\bigr).
\]
The representations are unramified at the new prime
pair, and the additive inertia restrictions with these
coefficients factor through the pro-\(\ell\) tame quotient.
Thus \(\ker_C A=C_{i,a}(J)\), since the chosen generators
detect vanishing on that quotient. For a unit \(u\) the
new congruence has slope \(u^{-1}\), and the inequality bounds
its length by that of the old cocycles satisfying the
Frobenius-value property at the new pair.

If the Frobenius-value property fails at any datum
and depth, choose an eligible prime where the condition fails.
For this datum exclude the union of the bad residues over
\(1\le a\le e\) and zero. The residue-field bound supplies
one unit \(u\) for the chosen datum at all depths. Give it
slope \(u^{-1}\), and give the other data slope zero. For the selected datum the upper bound is the length of the old
submodule \(C_{i,a}(J)\cap B^{-1}(A(C))\). At the failing depth
this is a proper submodule, hence the length drops by at least one.
For the other data, slope zero recovers unramifiedness at the new
pair, so their lengths do not increase. The sum of lengths over all
\(i,a\) decreases strictly.
For the initial singleton with slope zero this sum is
at most
\[
 B_0=\sum_i\sum_{a=1}^{e}
                  \operatorname{length}C_{i,a}(\varnothing).
\]
The same inertia argument justifies this bound at the
first prime. It follows that the process terminates
with at most \(1+B_0\) primes, independently of \(N\).

At termination the Frobenius-value property holds at every
depth for each eligible prime. To apply the primitive bound
to the group at depth \(e\), consider any depth \(a\) whose
group has no member with a unit value coordinate. Every
value then lies in \(\varpi\mathcal O/\varpi^a\).
The isomorphism
\(T/\varpi^{a-1}T\xrightarrow{\,\varpi\,}\varpi T/\varpi^aT\)
therefore identifies the group at depth \(a-1\) with that at depth
\(a\). Indeed the same integral homogeneous linear equations define
the cocycle, local, and auxiliary conditions at both depths, and
\(\varpi t=0\bmod\varpi^a\) exactly when
\(t=0\bmod\varpi^{a-1}\). Descend until there is a unit
value coordinate or the depth is zero. Thus the terminal
bound gives length at most \(G_i+L_i\) for \(C_{i,e}(J)\).
The first-place bound of Section~\ref{sec:selmer-relations}
then bounds releasing the distinguished pair at zero.
After dividing by \(\operatorname{ord}_{\mathcal O}(\ell)\)
this gives the bound \(\beta_i\) for the untwisted length
used in Section~\ref{sec:selmer-polynomials}. The argument
uses the hypotheses on the representation of Section~\ref{sec:terminal-duality}.

\subsection{Duality for the final local representation}
\label{sec:terminal-duality}

\sectionopening{In this subsection, we analyze the residual classes at the end of the
descent in Section~\ref{sec:auxiliary-descent}. The Cartier duality
calculation gives the necessary relation between the two
representations at the chosen torsion depth.}

Consider an untwisted specified space at a torsion depth \(a\)
with \(V=T/\varpi^aT\), \(T=\mathcal O^n\), and
\(a\le\operatorname{ord}_{\mathcal O}(\ell^M)\).
Put \(\Gamma=G_E\) and write \(g^c=cgc^{-1}\) for the
specified complex conjugation. For the Cartier dual put
\(V^*=\operatorname{Hom}_{\mathbf Z}(V,\mu_{\ell^N})\),
with scalar action \((t b)(v)=b(tv)\) and the contragredient
Galois action twisted by the cyclotomic character. Here
\(\ell^M V=0\) and \(N>M\). Suppose that on the original lattice
\(T\) there is an element \(s\in\Gamma\) satisfying
\[
 \rho(s)=\rho((s^c)^{-1})=\xi I,\qquad
 \chi_{\rm cyc}(s)=1,\qquad \xi^2-1\in\mathcal O^\times.
\]
Suppose \(s\) lies in inertia at an original place
prime to \(\ell\) and the matrices of \(\rho\) span the
full matrix algebra over \(\mathcal O\). The cocycles \(z\) satisfy
\(z(s)=0\). At termination assume the Frobenius-value property
of Section~\ref{sec:auxiliary-descent} at this depth. We prove
its length bound for a group containing a cocycle \(z\) with
a unit coordinate at some group element.

Let \(S'\) contain the original ramification set of the representation,
including the places above \(\ell\), and the later auxiliary
primes and their conjugates. The first auxiliary pair is outside
\(S'\). For \(v\in S'\), let \(A_v\) be the restrictions to
the local group of cocycles satisfying the original
local relations and the old congruences, with neither a ramification
constraint nor the normalization at \(s\).

Assume the original local conditions contain the coboundaries.
The congruence conditions also contain them, since inertia acts trivially at
the auxiliary primes and Frobenius acts trivially when the
slope is nonzero. Let \(\mathcal L_v\) be the image of \(A_v\)
in local cohomology. At each finite \(v\notin S'\), impose the
unramified condition. Here \(V\) is unramified and its order is prime
to the residue characteristic. The unramified subgroups for \(V\) and
\(V^*\) are exact annihilators under local Tate duality
\cite[Ch.~I, Theorem~2.6]{MilneADT}. In particular, the dual condition
is unramified at the first auxiliary pair and at each eligible new
pair outside \(S'\).
The local length estimate needed here is
\[
  \sum_{v\in S'}\operatorname{length}A_v
       \le (|S'|+f)\operatorname{length}V+L_i .
\]
The old auxiliary places contribute at most
\(\operatorname{length}V\) each because the Frobenius value
determines the tame value by the congruence and determines the
rest of the cocycle with that inertia value. Thus the
remaining local bounds concern the excess \(L_i\) at the original
places away from \(\ell\) and the crystalline rank bound at
the places above \(\ell\).

The cocycle module \(C\) injects into \(H^1_{\mathcal L}(E,V)\).
At the first pair the zero slope forces unramifiedness,
and if a coboundary vanishes at \(s\) its defining vector
is killed by the unit \(\xi-1\). Write
\(D^\vee_{\mathcal L}=H^1_{\mathcal L^\perp}(E,V^*)\)
for the ordinary Selmer group using the local Tate
orthogonal conditions. Since all infinite places of \(E\) are complex, the
archimedean factors below apply also at \(\ell=2\). The Greenberg--Wiles formula
\cite[Theorem 2.19]{DDT}, in unnormalized length, gives
\[
 \operatorname{length}H^1_{\mathcal L}(E,V)
 = \operatorname{length}D^\vee_{\mathcal L}
   +\sum_{v\in S'}\bigl(\operatorname{length}A_v-\operatorname{length}V\bigr)
   - f\operatorname{length}V .
\]
Indeed both global invariants vanish by the action of \(s\).
Since \(A_v\) contains the coboundaries, the local ratio
\(|\mathcal L_v|/|H^0(E_v,V)|\) is \(|A_v|/|V|\).
Unramified places give ratio one and complex places give
\(1/|V|\). Taking logarithms to the residue field cardinality
gives the displayed identity. In particular
\(\operatorname{length}C\le\operatorname{length}D^\vee_{\mathcal L}+L_i\).

Put \(U=\ker(\rho\bmod\ell^M)\cap\ker(\chi_{\rm cyc}\bmod\ell^N)\).
For a dual class \([b]\), write \(y=b|_U\).
This is additive and equivariant under conjugation, and is
independent of the representative since \(U\) acts trivially.
Restriction is injective on dual classes. Indeed its kernel
inflates from \(H^1(\Gamma/U,V^*)\), where the image of \(s\)
is central and acts by \(\xi^{-1}\). Comparing a cocycle
at \(sg\) and \(gs\) shows that such a cocycle is a coboundary
because \(s-1\) acts invertibly on \(V^*\).

Let \(H\) be generated by \([g,s]\) and \([g,s^c]\),
\(g\in\Gamma\), with \([g,s]=gsg^{-1}s^{-1}\).
The identity \([ag,s]=a[g,s]a^{-1}[a,s]\), also for \(s^c\),
gives normality, and conjugation by \(c\) interchanges the two sets
of generators. Both scalar actions commute with the whole image
on \(T/\ell^MT\); the cyclotomic character kills commutators.
Thus \(H\subset U\cap U^c\).
The identity \(z([g,s])=(1-\xi)z(g)\) and equivariance on
\(U\) show that the \(\mathcal O\)-span of \(z(H)\) is \(V\):
its reduction contains a nonzero vector, is stable under
all matrices, and Nakayama's lemma applies.

Identify \(\mu_{\ell^N}\) additively using one fixed primitive
root of unity. For \(t\in\mathcal O\) and \(u,u'\in U\cap U^c\)
put
\[
 F_t(u,u')=y(u')(t z(u))-y(u^c)(t z(u'^c)).
\]
For \(h\in H\), apply the prime construction of
Section~\ref{sec:first-auxiliary} to the target
\(\gamma_{w_0}h\). The commutators defining \(H\) lie
in the normal subgroup generated by inertia at the
original place of \(s\) and its conjugate. Together
with the first prime these places form \(R_{\rm in}\).

Let \(J\) index the old auxiliary pairs, with distinguished
representative \(w_0\). Take
\[
 S_{\rm sp}=\bigcup_{w\in J\setminus\{w_0\}}\{w,\bar w\}.
\]
The distinguished pair \(\{w_0,\bar w_0\}\) belongs to \(R_{\rm in}\),
so the splitting conditions apply only to the other pairs.
Exclude the old prime pairs, the original forbidden
primes and all places above the rational norm of \(w_0\).
The target and its conjugate are in \(U\). Every old
character kills \(h\), since \(H\) is generated by
commutators. At \(\gamma_{w_0}\) only the first character
is nonzero and reads \(1\). Prescribe these finitely many
residual and character evaluations, as well as the locally constant values of \(y,z\) and their conjugates, in the
Frobenius neighborhood.

The new character vanishes on its own Frobenius,
on every \(G_{E_v}\) for \(v\mid\ell\), and on \(G_{E_w}\)
and \(G_{E_{\bar w}}\) for \(w\in J\setminus\{w_0\}\).
It is unramified at the original level places away from \(\ell\).
Anticyclotomicity
preserves the character values on the inverse conjugate
parameter. Hence this produces an eligible pair with
\(\phi,\phi^c\in U\). The cocycles \(y,z\) vanish at the
first tame pair because both local classes there are
unramified and inertia acts trivially. Thus
\(F_t(h,h)=F_t(\phi,\phi)\).

Choose \(z'\) by the terminal property, with
\(z'(\gamma_P)=z(\phi)\) and
\(z'((\gamma_P^c)^{-1})=z((\phi^c)^{-1})\).
Apply global reciprocity to \(b\cup(tz')\). Its local
invariants vanish on \(S'\) by orthogonality, and
outside \(S'\) and the new pair by unramifiedness.
At the new pair the coefficient actions are trivial and
\(b\) is unramified. The local invariant reads the Frobenius
value of \(b\) on the tame value of \(tz'\), with respect to
the generator determined by the fixed root of unity.

For this formula, subtract from \(z'\) the Kummer cocycle
of a uniformizer times the normalized tame value of \(z'\).
The difference is unramified and its cup product with \(b\)
has zero invariant. To check the sign, put \(d_0=\ell^N\), choose a local
uniformizer \(\pi\) and \(r^{d_0}=\pi\), and write
\(\kappa_\pi(g)=g(r)/r\).
For the resulting unramified character \(\chi\) with values in
\(\mathbf Z/d_0\mathbf Z\), take integer lifts
\(\widetilde\chi(g)\). The coboundary of
\(g\mapsto r^{\widetilde\chi(g)}\) at \((g,h)\) is
\[
 \pi^{(\widetilde\chi(g)+\widetilde\chi(h)-\widetilde\chi(gh))/d_0}
 \kappa_\pi(g)^{\widetilde\chi(h)} .
\]
Thus graded commutativity identifies
\(\chi\cup\kappa_\pi\) with \(\pi\cup\delta(\chi/d_0)\)
in the Brauer group, where \(\chi/d_0\) takes values in
\(\mathbf Q/\mathbf Z\).
The local invariant is \(\chi(\phi)/d_0\) for arithmetic
Frobenius \(\phi\)
\cite[Ch.~III, Proposition~3.6 and Lemma~3.7]{MilneCFT}.
For the chosen geometric Frobenius and its geometric conjugate
the invariant has the opposite sign at both places.
The inverse conjugate of a normalized tame generator is
normalized by the same root of unity at the conjugate place.
The two chosen generators therefore differ from the
normalized generators by the same unit modulo \(\ell^N\).
Multiplying the reciprocity sum by this common unit gives
\[
 0=y(\phi)(t z(\phi))+
      y(\phi^c)(t z((\phi^c)^{-1}))=F_t(\phi,\phi).
\]
Here \(z((\phi^c)^{-1})=-z(\phi^c)\).
Hence \(F_t(h,h)=0\) on \(H\).

Additivity makes \(F_t\) antisymmetric on \(H\times H\).
Conjugating the first argument by \(s\) replaces \(t\) by
\(\xi t\), and conjugating the second replaces \(t\) by
\(\xi^{-1}t\), by the two scalar actions at \(s,s^c\).
Indeed
\[
 F_{\xi t}(h,k)=F_t(shs^{-1},k)
      =-F_t(k,shs^{-1})=-F_{\xi^{-1}t}(k,h)
      =F_{\xi^{-1}t}(h,k).
\]
Replacing \(t\) by \(\xi t\) gives \(F_{(\xi^2-1)t}=0\).
Because \(\xi^2-1\) is a unit, \(F_t\) vanishes for every \(t\).
It follows that \(z(h)\mapsto y(h^c)\) extends to a unique
\(\mathcal O\)-linear map \(V\to(V^*)^c\).

In fact, if \(\sum_h c_hz(h)=0\), then for every \(k\in H\)
and \(t\in\mathcal O\), the identities \(F_{tc_h}(h,k)=0\) give
\[
 \Big(\sum_h c_h y(h^c)\Big)(t z(k^c))
       = y(k)\!\Big(t\sum_h c_h z(h)\Big)=0 .
\]
The elements \(t z(k^c)\) generate \(V\) additively since
\(H^c=H\) and \(z(H)\) spans over \(\mathcal O\).
Equivariance follows from the conjugation formulas, and the
map depends \(\mathcal O\)-linearly on \([b]\) by uniqueness.
If the map vanishes, \(y|_H=0\); for \(u\in U\) the commutator
\([s,u]\) then gives \(y(u)=y(sus^{-1})=\xi^{-1}y(u)\).
Thus \(y=0\). We have proved the injection
\[
 D^\vee_{\mathcal L}\hookrightarrow
       \operatorname{Hom}_{\mathcal O[\Gamma]}(V,(V^*)^c).
\]
The bound \(G_i\) for the length on the right gives the
terminal bound \(G_i+L_i\) for \(C\). This proof uses the
stated local length and hypotheses on the representation. The first-place determinant losses
are paid separately at the untwisted specialization.

\section{Galois representations, finite types, and local lengths}
\label{sec:galois-types-and-traces}

\sectionopening{In this section, we attach Galois representations to the original Hecke eigenspaces and
compute the normalizations at the finite places. The resulting representations
supply the data for Section~\ref{sec:crystalline-auxiliary-selection}.}

\subsection{The original Hecke system and the normalized trace}
\label{sec:main-trace}

\sectionopening{In this subsection, we isolate the admissible finite type in the
original Hecke eigenspace using
Section~\ref{sec:marked-target-dimension}. We then identify the
normalization of the Frobenius trace needed for the attachment in
Section~\ref{sec:cofinite-spherical-galois}.}

Let \(z\ne0\) have Hecke system \(b=(b_{w,i})\) in the
final projected space. On the principal level the finite-type idempotent is a sum of the central idempotents attached
to admissible tuples of finite parameters. Since the sum fixes
\(z\), some term acts nontrivially on \(z\). The corresponding
component keeps \(b\): the unramified Hecke operators act
at places away from the residue group and commute
with its action. The joint trace comparison of
Section~\ref{sec:marked-target-dimension} gives a nonzero
component on the definite group with the same prescribed finite type and outside-level Hecke system.

The terminal rank-\(r\) space \(W_m\) is definite with standard
form and finite basis. Its level is principal at the primes
dividing \(\pi_M\) and the prescribed maximal integral compact
elsewhere, as in Section~\ref{sec:marked-target-dimension}.
The last-covector condition at \(\pi_N\) is imposed on the
rank-\(k\) source \(X\) of Section~\ref{sec:ramified-source}.
Let \(\sigma\) be the compact
coefficient representation of the weight,
\(v\) its unit highest vector and \(d_\sigma=\dim\sigma\).
On a form reproduced by the real projector set
\[
            Lz(g_f)=d_\sigma\int_{U(W_m)_\infty}
                             z(k,g_f)\sigma(k)v\,dk ,
                  \qquad \int dk=1 .
\]
The compact projector kernel in these coordinates is
\(d_\sigma\langle v,\sigma(k)v\rangle\), with the
inner product linear in its second variable. Indeed
\(\langle v,\sigma(k^{-1})Lz(g_f)\rangle=z(k,g_f)\)
by the projector equation and change of variable
in Haar. Thus \(L\) is injective. For a rational unitary
element \(\gamma\), the change \(k\mapsto\gamma_\infty k\)
uses left rational invariance of \(z\) and gives
\(Lz(\gamma_f g_f)=\sigma(\gamma_\infty)Lz(g_f)\).
This is the coefficient transformation of the algebraic form.
The integral commutes with every finite right translation,
so it preserves the original level, finite type projector
and outside-level Hecke system.

\smallskip
\noindent\textit{The split polynomial.}
Let \(w,\bar w\) be a split pair off the level and
twist of a descent with source rank \(r_0\) and target rank
\(n=r_0+t_0\). The local Schr\"odinger test
there is the matrix-lattice vacuum
\(\Phi_w=\mathbf 1_{\operatorname{Mat}_{n\times r_0}
(\mathcal O_w)}\). The other local factors are unchanged.
For an \(i\)-space \(U\) in a residue vector space of
dimension \(k\), choose an integral \(\gamma_U\) with
elementary divisors \(\varpi_w\) in \(i\) positions and
\(1\) elsewhere, whose reduction has kernel \(U\).
These kernels parameterize the Hecke representatives.
The target representatives are read by transpose
and those at the opposite source place by inverse.
The split Schr\"odinger rules give the two terms
\[
      \nu_w^i\Phi_w(\gamma_U^{-\mathsf T}X),\qquad
      (\nu'_w)^i\Phi_w(X\gamma_U^{-1}).
\]
The target determinant modulus factor is \(q^{r_0i/2}\),
the source one is \(q^{ni/2}\), and the splitting
characters at the uniformizer give the other factors
of \(\nu_w^i,(\nu'_w)^i\) specified below.

The two terms require \(X\) integral. For integral
\(X\) the tests are \(U\subseteq\ker\bar X^{\mathsf T}\)
and \(U\subseteq\ker\bar X\), respectively, where here
\(\bar X\) denotes reduction. Indeed \(\gamma_U^{-1}\)
has only a simple pole and the image of
\(\varpi_w\gamma_U^{-1}\) modulo \(\varpi_w\) is \(U\).
Let \(s\) be the residue rank of \(X\). For
\(P_d(Z)=\prod_{j=0}^{d-1}(1-q^jZ)\), the subspace
identity
\[
          P_d(Z)=\sum_{i=0}^{d}(-Z)^i q^{\binom i2}
                                             \genfrac{[}{]}{0pt}{}{d}{i}_q
\]
gives target and source series \(P_{n-s}(\nu_wY)\)
and \(P_{r_0-s}(\nu'_wY)\).

\subsubsection{Hecke multiplicities and the trace}
\label{sec:normalized-trace-multiplicities}

\sectionopening{In this subsection, we compute the multiplicities used to normalize the original trace in Section~\ref{sec:main-trace}.}

The identity follows by counting \(i\)-spaces
contained in and transverse to a fixed hyperplane,
which gives the recurrence for multiplying by
\(1-q^{d-1}Z\). Because \(n=r_0+t_0\),
\[
 P_{n-s}(\nu_wY)
       =P_{t_0}(\nu_wY) P_{r_0-s}(q^{t_0}\nu_wY).
\]
Thus on the actual test the target Hecke polynomial
action is the source one at \(q^{t_0}\nu_wY/\nu'_w\)
times \(P_{t_0}(\nu_wY)\). Moving the source action
through the bilinear theta pairing inverts the coset under
the inherited Haar, changing its component at \(\bar w\)
to \(w\). Both actions preserve the original data
at the other local factors. This proves the split
Hecke comparison when the lift is nonzero.

Apply the split polynomial identity to a nonzero lifted
Hecke eigenform. Fix a descent
pair with ranks \(r_0\) and \(n=r_0+t_0\), with \(n\) even,
and let \(\mu\) be its unitary character before the tame
twist \(\vartheta\). At an unramified split place \(w\) of norm
\(q\), use the reciprocal Hecke polynomial
\[
 P_{b,w}(Y)=\sum_{i=0}^{r_0}(-Y)^i q^{\binom i2} b_{w,i}.
\]
Here \(b_{w,0}=1\). For the lift let \(H_w(X)\)
be the monic degree-\(n\) polynomial in the Galois
normalization.

The split Hecke comparison needed here is
\[
 Y^nH_w(Y^{-1})=
     \prod_{a=0}^{t_0-1}(1-\nu_wq^aY)\,
         P_{b,w}(q^{t_0}\nu_wY/\nu'_w).
\]
The nonzero scalars \(\nu_w,\nu'_w\) are the rank-\(r_0\)
and rank-\(n\) Frobenius values of the splitting characters.
For a norm-one character \(\chi\) and a rank \(n'\),
this value is
\[
 q^{n'/2}\chi(\varpi_w/\bar\varpi_w)\upsilon(\varpi_w)^{-n'},
\]
where \(\upsilon\) is the splitting character of the Weil
data. In the even rank the character for the lift is
\(\chi=\mu\vartheta(\upsilon|_{U(1)})^{n/2}\). By
Section~\ref{sec:splitting-level-characters},
\(\upsilon(\varpi_w/\bar\varpi_w)=\upsilon(\varpi_w)^2\), so
\(\nu'_w=q^{n/2}(\mu\vartheta)(\varpi_w/\bar\varpi_w)\),
where \(\varpi_w\) denotes the uniformizer idele.

Comparing linear terms gives the exact formula
\[
 -\nu_w^{-1}[X^{n-1}]H_w
    =\frac{q^{t_0}}{\nu'_w} b_{w,1}
             +\sum_{a=0}^{t_0-1}q^a .
\]
In the geometric Frobenius convention,
\(\chi_{\rm cyc}(\mathrm{Fr}_w)^{-1}=q\).
The scalar \(q^{t_0}/\nu'_w\) is the product of the
characters at \(\mathrm{Fr}_w\) obtained from
\(\vartheta^{-1}\) and from the algebraic character
\[
 u\longmapsto \mu^{-1}(u/\bar u)|u|^{(r_0-t_0)/2}
\]
by reciprocity in the same convention. This follows
from \(|\varpi_w|=q^{-1}\) and \(n=r_0+t_0\).
Write \(\theta\) for the avatar of the latter character
and \(\vartheta^{-1}_{\ell}\) for the former avatar.

Let \(a\) be the source splitting character of the original
pair. The normalizing avatar \(\alpha_{\rm nor}\) is that of
the algebraic Hecke character
\[
     x\longmapsto a(x/\bar x)^{-1}\upsilon(x)^{r_0}|x|^{r_0/2}.
\]
On the chosen CM type, the last two factors have infinity
type \(z^{r_0}\), by the infinity type of \(\upsilon\).
The first factor has integral powers of \(z,\bar z\) since
\(a\) is a character of the norm-one torus. Its finite unit
action is trivial at the original and auxiliary level places, since the
source character is a prescribed power of the norm-one
restriction of \(\upsilon\). At the selected geometric
Frobenius the avatar is \(\nu_w^{-1}\). Multiplying the
Hecke character by its conjugate gives \(|x|^{r_0}\), hence
\(\alpha_{\rm nor}\alpha_{\rm nor}^{\,c_E}=\chi_{\rm cyc}^{r_0}\).

If \(R\) is the representation attached to \(b\) and
\(R_{\rm lift}\) to \(H_w\), put
\(\widetilde R_{\rm lift}=\alpha_{\rm nor}R_{\rm lift}\).
After transporting coefficients by \(\iota\), the
normalized trace equality is
\[
 \operatorname{tr}\widetilde R_{\rm lift}(g)
    =\vartheta^{-1}_{\ell}(g)\operatorname{tr}((\theta R)(g))
                    +\sum_{a=0}^{t_0-1}\chi_{\rm cyc}(g)^{-a}.
\]
Use the cofinite Frobenius comparison of
Section~\ref{sec:cofinite-spherical-galois} for both systems.
In every finite quotient of \(G_E\), the conjugacy class of
the image of \(g\) has positive Chebotarev density
\cite[Ch.~V, Theorem~3.23]{MilneCFT}.
The primes of \(E\) of absolute residue degree greater than one
have density zero, and every remaining unramified degree-one
prime is split over \(F_0\). Using both members of split pairs
and omitting the two systems' finite excluded sets, we obtain
the linear coefficient identity on dense conjugacy classes.
Both sides are continuous class functions on \(G_E\), giving the
displayed identity for every \(g\).
By semisimplicity and the independence of irreducible trace
characters, the main summand of \(\widetilde R_{\rm lift}\) at
twist \(\vartheta\) is \(\vartheta_\ell^{-1}\otimes\theta R\).
The representation \(\theta R\) and its stable lattice can be fixed
before choosing \(\vartheta\); at \(\vartheta=\xi_\kappa^j\)
the scalar twist of this lattice is exactly \(\psi^j\), where
\(\psi\) is the inverse avatar of \(\xi_\kappa\).

\smallskip
\noindent\textit{Removing auxiliary character factors.}
Suppose two continuous semisimple representations \(R_i\)
in the same coefficient field have characteristic polynomial
at the prescribed Frobenius elements equal, in the limit, to
a common monic polynomial times the linear factors of a list
of characters \(U_i\). Assume the two character lists are
disjoint as Galois characters. On taking first coefficients
and limits, the virtual characters \(\operatorname{tr}R_i-
\operatorname{tr}U_i\) agree on \(G_E\).
Characters of pairwise nonisomorphic irreducible representations
over an algebraically closed field of characteristic zero are
linearly independent.

Indeed the algebra generated by the
group on their direct sum is finite-dimensional and its
faithful module is semisimple, so its simple factors are the
full matrix algebras on the distinct simple modules. Their
trace functions are independent on these factors.
For any constituent in \(U_i\), comparison with the
other list shows that its multiplicity in \(R_i\)
is at least its multiplicity in \(U_i\).
Thus both \(R_i\) contain \(U_i\) and have the same
semisimple remainder \(R\). The summands are continuous
as restrictions to invariant subspaces.

This proves
the cancellation at the level of representations.
The resulting characteristic polynomials satisfy
\(\operatorname{charpol}(R_i(g))=
\operatorname{charpol}(R(g))\operatorname{charpol}(U_i(g))\).
At the specified Frobenius elements, division by the
monic character factor recovers the common attachment
polynomial in the limit. Indeed coefficients of the monic
quotient are obtained by successively subtracting products
of higher coefficients, so this division is continuous.

If both larger attachments are unramified, crystalline,
or potentially crystalline at a place, these properties pass
to the remainder as a direct summand. The jumps at a
coefficient embedding are obtained by subtracting the
multiset of jumps of the corresponding character list.
At an unramified place, if the two larger annihilating
polynomials are \(P(X)Q_i(X)\) and \(Q_0,Q_1\) are
coprime there, \(P\) annihilates the remainder by the
B\'ezout identity applied to its Frobenius matrix.
These reductions require the two larger attachments and
the local polynomial comparisons on the form.

\subsection{Galois attachment and spherical compatibility}
\label{sec:cofinite-spherical-galois}

\sectionopening{In this subsection, we attach a Galois representation to the spherical
eigensystem of Section~\ref{sec:main-trace}. The contragredient
compact type fixes the normalization, which is used again for the
Frobenius tests of Section~\ref{sec:triangle-unit-input}.}

Use an eigensystem on a definite group of rank \(n\),
with the chosen compact coefficient representation
\(\sigma_{\rm alg}\). Its scalar coefficient is
\(\langle\sigma_{\rm alg}(k)v,Lz(g_f)\rangle\).
At the chosen CM embedding the coefficient action on Weyl
coordinates is \(k\mapsto S_\mu(\bar k)\), with the determinant
power included in \(S_\mu\); here \(\bar k=(k^{-1})^{\mathsf T}\)
on the definite unitary group. Right translation acts in the
conjugate-linear first variable, so the compact automorphic
constituent has type \(\sigma_{\rm alg}^{\vee}\simeq S_\mu\).
Its Harish--Chandra entries are
\(\mu_q+(n-1-2q)/2\), \(0\le q<n\); at the conjugate embedding
their multiset is negated. This fixes the sign of the raw jumps
in Section~\ref{sec:coefficient-hodge-tate-intervals}. One may choose an irreducible automorphic
constituent containing a component of the eigenvector:
the space is compact modulo rational points, and the
outside-level eigenvalues and the specified finite type
are preserved by taking components. By the compact
coefficient realization in Section~\ref{sec:main-trace}
this is the same finite-level Hecke system.

Apply \cite[Corollaire~5.3, p.~43 of the author's version]{LabesseCM}.
The corollary specifies local compatibility at each split place and at
each unramified place. In the paired split coordinates the transfer
of \(\pi_v\) is \(\pi_v\boxtimes\pi_v^\vee\), including at a ramified
split component. At a spherical inert place the transfer uses the
unramified Satake parameter
\cite[\S2.4, Theorem~2.2, pp.~1040--1041]{GuerberoffUnitary}. Section~\ref{sec:datum} gives \([F_0:\mathbb Q]=\varphi(q_0)3L_0Q\),
a multiple of four, so the degree condition for the definite group
also holds in even rank. The extension \(E/F_0\) is unramified at
finite places. The definite groups used here have the
standard positive Hermitian matrix. At a nonsplit
finite place choose a unit \(u\) of norm \(-1\), using
the unit norm surjectivity of
Section~\ref{sec:splitting-level-characters}. The
vectors \(e_{2j-1}+u e_{2j}\), \(1\le j\le\lfloor n/2\rfloor\),
span a totally isotropic subspace of maximal dimension.
Thus the group is quasi-split there. At a split place it
is \(\mathrm{GL}_n\).
The real groups are compact, so their representations
are discrete series and Labesse's condition \((*)\) holds.

The rank-\(n\) transfer is an isobaric sum of discrete
conjugate-self-dual representations over \(E\).
It agrees with the local base change at every split
finite place and every inert place where the representation
is spherical. Its infinitesimal character is the base change
of the compact weight. By M{\oe}glin--Waldspurger
\cite[theorem of the introduction, pp.~605--606]{MWResidual}, each block of rank \(ab\) is a
length-\(b\) Speh representation from a unitary cuspidal
\(\sigma\) of rank \(a\), with exponents
\(u_j=(b-1)/2-j\), \(0\le j<b\).
Conjugate self-duality and uniqueness of the cuspidal
support give \(\sigma^c\simeq\sigma^\vee\).
In both complex coordinates at infinity, the
infinitesimal weights of the full base change lie in
\((n-1)/2+\mathbb Z\) and are pairwise distinct.
Subtracting \(u_j\) shows that the weights of \(\sigma\)
are regular and lie in \((n-b)/2+\mathbb Z\).

Use the unitary character \(\upsilon\) of
Section~\ref{sec:splitting-level-characters}. It is everywhere
unramified at finite places, satisfies \(\upsilon^c=\upsilon^{-1}\),
and has infinity type \(z/|z|\) at the chosen CM embeddings. Set
\[
 \psi=\begin{cases}
 1,& n-a-b+1\text{ even},\\
 \upsilon,& n-a-b+1\text{ odd}.
 \end{cases}
\]
Then \(\sigma\otimes\psi\)
is regular algebraic and conjugate self-dual, since
its infinitesimal exponents lie in \((a-1)/2+\mathbb Z\).

For its polarizing character on \(F_0\) use \(1\) if
\(a\) is even and the quadratic character of \(E/F_0\)
if \(a\) is odd. Its pullback by the norm is trivial,
and its sign at each real place is \((-1)^a\).
The sign in the Galois multiplier is then
\((-1)^{1-a}(-1)^a=-1\).
Hence the hypotheses of Barnet-Lamb--Gee--Geraghty--Taylor
\cite[Theorem~2.1.1]{BLGGT} apply for every coefficient
prime \(\ell\), including \(2\).

For each \(j\) the character
\[
 \alpha_j=\psi^{-1}|\cdot|^{(a-n)/2+u_j}
\]
is algebraic. Put \(d_j=a-n+b-1-2j\) and let \(e_\psi\)
be \(0\) for \(\psi=1\) and \(1\) for \(\psi=\upsilon\).
Since \(d_j\equiv n-a-b+1\equiv e_\psi\pmod 2\), its
infinity character on the chosen CM type is
\[
 z^{(d_j-e_\psi)/2}\bar z^{(d_j+e_\psi)/2}.
\]
The two integral exponents exchange at the conjugate embedding.
Let \(\rho_{\ell,\iota}\) be the direct sum over all blocks and \(j\)
of
\[
           r_{\ell,\iota}(\sigma\otimes\psi)
                         \otimes r_{\ell,\iota}(\alpha_j).
\]
Here \(r_{\ell,\iota}(\alpha_j)\) is the class field
avatar in the same geometric Frobenius convention.
The twists total \((1-a)/2+(a-n)/2+u_j=(1-n)/2+u_j\).
Thus this direct sum, on the chosen split coordinate,
has the local Langlands normalization of the base
change twisted by \(|\det|^{(1-n)/2}\). This is
the direct sum computation in the proof of
\cite[proof of Theorem~2.3, pp.~1041--1042]{GuerberoffUnitary}. The summand representations
here come from the cuspidal attachments of
\cite[Theorem~2.1.1]{BLGGT} for every \(\ell\), including when
\(\ell=2\) or a place above \(\ell\) is nonsplit.
At a nonsplit spherical place
\(w\nmid\ell\), Labesse gives the unramified base change.
The factors \(\psi_w\) and \(\psi_w^{-1}\) cancel in each
shifted block, so the direct sum has trivial
Weil--Deligne inertia and zero monodromy. Thus its Galois representation
is unramified at \(w\).

At an unramified split place \(w\) let \(\lambda_i\)
be the normalized Satake parameters. With Haar volume
one for the maximal compact and right translation by
the positive diagonal coset with \(j\) uniformizers,
the eigenvalue is
\(t_{w,j}=q^{j(n-j)/2}e_j(\lambda_1,\ldots,\lambda_n)\).
The Hall--Littlewood polynomial for \((1^j)\) is \(e_j\)
by \cite[Ch.~III, (2.8)]{MacdonaldHall}.
With \(q^{-s_i}=\lambda_i\), the Hecke algebra map and its spherical
specialization in \cite[Ch.~V, (2.7) and (3.3)]{MacdonaldHall} give
\[
 q^{-\binom{j}{2}} e_j\bigl(q^{(n-1)/2}\lambda_1,\ldots,
 q^{(n-1)/2}\lambda_n\bigr)
 = q^{j(n-j)/2}e_j(\lambda_1,\ldots,\lambda_n).
\]
\begin{proposition}[Spherical Frobenius polynomial]\label{prop:marked-frobenius-polynomial}
For an unramified split place away from \(\ell\) and the finite
excluded set of the attached system, transport the Hecke and
Satake coefficients to \(\overline{\mathbb Q}_\ell\) by \(\iota\).
The characteristic polynomial of geometric Frobenius in the
chosen split coordinate is
\begin{equation}\label{eq:actual-spherical-galois}
\begin{aligned}
 \det(XI_n-\rho_{\ell,\iota}(\operatorname{Frob}_w^{\rm geom}))
 &=\sum_{j=0}^n(-1)^j q^{\binom j2}t_{w,j}X^{n-j}\\
 &=\prod_{i=1}^n(X-q^{(n-1)/2}\lambda_i).
\end{aligned}
\end{equation}
\end{proposition}
\begin{proof}
On the unramified quotient geometric reciprocity sends the class of
\(\varpi_w\) to geometric Frobenius.
The local comparison above twists the unitary parameter by
\(|\det|^{(1-n)/2}\), whose value on one \(\varpi_w\) is
\(q^{(n-1)/2}\). The Frobenius eigenvalues are therefore
\(\alpha_i=\iota(q^{(n-1)/2}\lambda_i)\), counted with multiplicity.
For \(0\le j\le n\), the minuscule coset formula gives
\[
 e_j(\alpha_1,\ldots,\alpha_n)
   =\iota\!\left(q^{j(n-1)/2}e_j(\lambda_1,\ldots,\lambda_n)\right)
   =q^{\binom j2}\iota(t_{w,j}).
\]
Expansion of \(\prod_i(X-\alpha_i)\) gives
\eqref{eq:actual-spherical-galois} with the indicated coefficient transport.
\end{proof}
All coefficients on the left are computed at \(w\); the inverse
Satake tuple belongs to \(\bar w\). Each parity-corrected cuspidal
attachment lies in a regular compatible system by
\cite[Theorem~2.1.1(3) and \S5.1, final paragraph]{BLGGT}.
After a finite coefficient extension its integral models at the
coefficient places are given by \cite[Lemma~5.3.1(3)]{BLGGT}.
The compatible systems of the polarized characters
are regular because their rank is one. For each fixed \(\ell\),
the direct sum with the specified character twists is a continuous
representation into \(\mathrm{GL}_n(\overline{\mathbb Q}_\ell)\).

\begin{lemma}\label{lem:continuous-attachment-stable-lattice}
A continuous representation of a profinite group in
\(\mathrm{GL}_n(\overline{\mathbf Q}_\ell)\) is defined over a
finite extension \(K_\ell/\mathbf Q_\ell\) and preserves an
\(O_{K_\ell}\)-lattice of rank \(n\).
\end{lemma}
\begin{proof}
Exhaust \(\overline{\mathbf Q}_\ell\) by an increasing sequence of finite
extensions, using finiteness of the extensions of each fixed degree.
The inverse images of their matrix groups are closed subgroups
covering the given group. By Baire one has interior and hence finite
index. Adjoin the matrix entries of finitely many coset representatives
to its coefficient field to obtain \(K_\ell\). For an initial
\(O_{K_\ell}\)-lattice \(L_0\), compactness gives \(a\ge0\) with
\[
 L_0\subseteq\sum_g\rho(g)L_0\subseteq\ell^{-a}L_0 .
\]
The middle module is invariant and finite free of rank \(n\).
\end{proof}
Apply the lemma to each attachment. For finitely many systems take
the compositum of their coefficient fields and extend the lattices
by scalars. Thus the coefficients
in the \(\ell\)-adic embedding are integral and the
polynomial annihilates the actual Frobenius operator.
These claims hold at every such split place outside
a finite set including the coefficient and level places.

Fix \(\ell\) and the original level. For the family in
Section~\ref{sec:triangle-unit-input}, fix the representations and
the union of their finite excluded sets before choosing the depth
and the auxiliary primes. At the resulting raised level,
Section~\ref{sec:integral-raised-traces} uses the finite union for its
own systems before applying Chebotarev at that level.
An existing semisimple attachment
specified by convergence to the same Frobenius
polynomials has the same trace by continuity on
the dense selected Frobenius classes. Independence
of irreducible characters as in
Section~\ref{sec:main-trace} identifies the two
semisimple representations. In particular
\eqref{eq:actual-spherical-galois} applies to the
same systems used in the residual comparison.
It supplies the spherical part of
Section~\ref{sec:triangle-unit-input}.

\subsection{Inertial parameters of the finite types}
\label{sec:inertial-labels}

\sectionopening{In this subsection, we read the inertia eigenvalues of the attachment
in Section~\ref{sec:cofinite-spherical-galois} from the original
finite types. These local data verify the spectral conditions for the
length estimates of Section~\ref{sec:main-row-budgets}.}

For a split place of residue characteristic different
from \(\ell\), the attachment of
Section~\ref{sec:cofinite-spherical-galois} has the full
Frobenius-semisimplified Weil--Deligne comparison in the
chosen coordinate. Labesse supplies local base
change at split places, and \cite[Theorem~2.1.1(2)]{BLGGT}
gives the comparison for each parity-corrected cuspidal block.
In rank at least two the monodromy assertion is
\cite[Theorem~1.1, p.~2312]{CaraianiMonodromy}. The regular algebraic cuspidal
representation \(\sigma\otimes\psi\) is cohomological at the complex
infinite places by \cite[Lemme~3.14, p.~114]{ClozelMotifs}.
The auxiliary quadratic character in that lemma is trivial at every
complex place since \(\mathbf C^\times\) is connected. Conjugate self-duality
therefore gives temperedness at all finite places by
\cite[Theorem~1.2, p.~2312]{CaraianiMonodromy}. Rank-one blocks are
unitary characters. Removing the unitary parity twist preserves
temperedness. Thus the factors of a length-\(b\) Speh block,
ordered by \(u_0>\cdots>u_{b-1}\), form a standard module
on each finite general linear group. By the Langlands quotient
description, its parameter is the sum of the shifted tempered
parameters. Each discrete automorphic block is unitary, so its
irreducible local component is unitarizable. Bernstein's theorem
\cite[Theorem~2.5(i), p.~346]{TadicUnitary} gives irreducibility of the
normalized induction of these local components across the discrete
blocks. By \cite[Corollary~A.7(i), p.~375]{TadicUnitary}, its
Langlands datum is the union of the block data with multiplicities.
Its Weil--Deligne parameter is therefore the direct sum of the block
parameters, including their monodromy operators. The algebraic
twists in the direct sum cancel the parity corrections
at this very place. This proves the required comparison
on the same attachment. Frobenius semisimplification preserves
the Weil--Deligne inertia action and the monodromy operator.
When the latter vanishes, the inertia action on the original
\(\ell\)-adic representation at \(w\nmid\ell\) equals its
finite Weil--Deligne inertia action.

We describe the level-zero local calculation for the
finite type of that attachment. Let \(k=\mathbb F_q\)
be the residue field, \(K_w=\mathrm{GL}_n(\mathcal O_w)\),
and \(K_{w,1}\) its principal kernel. For an irreducible
factor of the squarefree type polynomial of degree \(a\), the
finite type is the genuine cuspidal character
\(\tau_{\theta}=(-1)^{a-1}R_T^{\theta}\) of \(\mathrm{GL}_a(k)\),
where \(T(k)=\mathbb F_{q^a}^{\times}\) and \(\theta\) is indexed by a root
as in Section~\ref{sec:marked-inertia-characters}. Its Frobenius
orbit has size \(a\), since the root has degree \(a\) and the
character parametrization commutes with Frobenius. This is the
primitivity hypothesis in \cite[\S2.4]{EmertonGeeHerzig}. Full
multiplicative order is required only at the separately chosen
primitive type polynomials. The corresponding depth-zero supercuspidal has on
inertia the \(a\) Frobenius conjugates of this character,
with zero monodromy
\cite[Proposition~2.4.1(i), p.~10 of the authors' version]{EmertonGeeHerzig}.
This is a characteristic-zero statement.

For the full parameter, the \(K_{w,1}\)-invariants of a
parabolic induction are the finite Harish--Chandra
induction of the inducing \(K_{w,1}\)-invariants \cite[Proposition~6.18]{ShottonBM}.
Indeed the maximal compact and the parabolic give
Iwasawa decomposition, and reducing the intersection
with the parabolic gives the finite one. Invariants
are exact in characteristic zero by compact averaging.
Use the finite irreducible constituent \(\sigma_{\mathcal P}\) of
\cite[Definitions~6.6--6.7 and Proposition~6.9]{ShottonBM}.
A tame Frobenius orbit of multiplicity one contributes its finite
cuspidal representation. On a determinant block of size \(d\), untwist
its determinant character. The trivial representation occurs in the
finite induction from the Borel; by Frobenius reciprocity it is absent
from any induction having a Steinberg block of size greater than one.
Thus it is the constituent labelled by \((1^d)\). Finite parabolic
induction combines the constituents on the distinct cuspidal lines.

For an irreducible local representation containing this type, take
its Langlands standard module. A supercuspidal factor with nonzero
principal-kernel invariants contributes precisely the finite cuspidal
representation of its depth-zero type, by the Cartan decomposition
in Section~\ref{sec:raised-tame-frobenius}. The finite cuspidal
support then identifies their multiset
\cite[\href{https://maths.dur.ac.uk/users/jack.g.shotton/papers/gln.pdf\#page=36}{\S6.3}]{ShottonBM}.
If \(\mathcal P'\) is the
family of segment lengths in this module, its principal-kernel
invariants are the finite induction \(\pi_{\mathcal P'}\) of the
corresponding generalized Steinberg blocks. Exactness implies
\(\operatorname{Hom}(\sigma_{\mathcal P},\pi_{\mathcal P'})\ne0\).
The finite multiplicity formula gives
\(\mathcal P'\preceq\mathcal P\) on every line
\cite[Corollaries~6.10 and~6.19, Theorem~6.20]{ShottonBM}.
All the partitions in \(\mathcal P\) have only singleton parts, so
the monodromy is zero. This is the forward implication of
\cite[Theorem~3.7(1) and Corollary~6.21(1)]{ShottonBM} for the
specified type. For the smooth complex calculation one may use
an auxiliary coefficient prime different from the residue
characteristic and transport its Weil--Deligne parameter to
\(\mathbf C\). The local--global comparison at the prescribed
coefficient prime is applied to this same complex parameter.
In particular the squarefree type gives one depth-zero factor per
irreducible polynomial, and a determinant block of size \(d\)
gives \(d\) rank-one factors on its tame line.
There is no nontrivial segment on a line occurring once.

On a repeated line, twist out the common tame character.
A segment of length \(a>1\) on the remaining rank-one unramified
line is an unramified twist of Steinberg. Put
\(H_a=\mathrm{GL}_a(k)\) and let \(B_a\) denote its upper Borel.
Reducing the flag functions at the principal kernel and taking
invariants by compact averaging gives
\[
 \operatorname{St}_{H_a}=\mathbf C[H_a/B_a]\Big/
       \sum_{B_a\subsetneq P\subseteq H_a}\mathbf C[H_a/P].
\]
Each standard parabolic summand is pulled back along its flag
projection. For \(a>1\), the summand \(P=H_a\) contains the
constants, which are the only \(H_a\)-invariants in the numerator.
Exactness gives \((\operatorname{St}_{H_a})^{H_a}=0\). By Frobenius
reciprocity the finite Harish--Chandra induction from
such segments cannot contain the trivial representation
of the repeated block.

Fix an old place \(w\in S_M\). Let \(P\subset\mathrm{GL}_{r_0+t_0}(k_w)\)
stabilize the image of the original rank-\(r_0\) block, and let
\(\tau_w\) be its selected lower residue type. The pointwise stabilizer
of the original columns contains the unipotent radical and the whole
complementary \(\mathrm{GL}_{t_0}(k_w)\); its action there is trivial by
the original level. On the first block the projector uses right
translation multiplied by \(\chi_{M,w}^{j}\circ\det\), as in
Section~\ref{sec:top-coefficient-clearing}. Thus ordinary right
translation there has type \(\tau_w\otimes(\chi_{M,w}^{-j}\circ\det)\).
The raised level has the same local conditions at \(w\). On every
upper constituent at either level, finite Frobenius reciprocity gives
\[
 \operatorname{Hom}_{\mathrm{GL}_{r_0+t_0}(k_w)}\!\left(
  \operatorname{Ind}_{P}^{\mathrm{GL}_{r_0+t_0}(k_w)}
   \bigl((\tau_w\otimes(\chi_{M,w}^{-j}\circ\det))\boxtimes\mathbf1\bigr),
  \pi_w^{K_{w,1}}\right)\ne0 .
\]
For \(r_0=r\), the lower type excludes the roots of order dividing
\(q_2\), so its rescaled roots exclude \(1\). For \(r_0=1\), the
unscaled polynomial is \(X-1\) and \(\chi_{M,w}^{-j}\) has order
\(q_2/\gcd(q_2,j)>1\), since \(1\le j\le m<q_2\). In both cases
the two blocks have disjoint cuspidal supports. At a twist place of
exact ramification the two determinant characters are distinct. In finite parabolic
Mackey theory, a double coset moving a positive-rank subblock between
distinct cuspidal lines would require the corresponding two Jacquet
modules to share a cuspidal factor. Their supports are disjoint by the parabolic product formula
\cite[\S2, Theorem~2, and \S8, Theorem~14]{Green}.
Thus only the identity double coset can contribute when the Levi blocks
are grouped by cuspidal line in the same order on both sides.
The multiplicity of the prescribed type in the finite induction
is the product of its multiplicities on the distinct cuspidal
lines \cite[Proposition~6.9 and Corollary~6.10]{ShottonBM}.
The zero multiplicity on the repeated line therefore remains
zero after adjoining the other lines. Thus each segment in the determinant
block has length one and its monodromy is zero.
The same proof applies to a pair of distinct tame
determinant characters on two blocks. In particular
at a split auxiliary prime with ramification
the normalized characters are \(\vartheta_{\ell}^{-1}\)
on the original rank block and \(1\) on the added block.
For each inertia element \(g\) the polynomial
\((X-\vartheta_{\ell}^{-1}(g))(X-1)\) annihilates this type.
Away from \(\ell\) this is actual Galois inertia, since monodromy is
zero. At the coefficient place the polynomial calculation concerns the
linear inertia on potentially crystalline periods, as proved below.

\subsubsection{The inertia characters}
\label{sec:marked-inertia-characters}

\sectionopening{In this subsection, we compute the inertia characters of the finite type on the generator used in the trace comparison of Section~\ref{sec:main-trace}.}

To specify the action of the generator, fix the faithful
multiplicative character \(\jmath:\bar k^\times\to
\mathbb C^\times\) from the parameter convention.
If \(\jmath(\xi)=\exp(2\pi i b/(q^a-1))\), the Green
character of the degree-\(a\) field is
\(\jmath|_{\mathbb F_{q^a}^\times}^{\,b}\).
For each \(a\), let \(F_a/E_w\) be unramified of degree \(a\).
On inertia set
\[
 \omega_a(s)=\overline{\operatorname{Art}_{F_a}^{-1}(s^{\rm ab})}
                      \in\mathbb F_{q^a}^{\times},
\]
using geometric reciprocity and abelianization in \(W_{F_a}\).
For \(a\mid D\) the norm relation gives
\(\omega_a=\omega_D^{(q^D-1)/(q^a-1)}\).
Take \(D\) divisible by all degrees in the finitely many
systems at this place. The tame quotient is procyclic, and
\(\omega_D\) maps it onto \(\mathbb F_{q^D}^{\times}\).
Choose one tame generator \(s\) with
\(\jmath(\omega_D(s))=\exp(2\pi i/(q^D-1))\).
For every such system the Green character on \(s\) then equals
\(\jmath(\omega_a(s))^b=\jmath(\xi)\).

Its Frobenius conjugates are the \(q\)-th powers.
Fix any \(w\in S_M\) above \(p_\nu\ne\ell\), and put \(q=p_\nu\). Let \(c_w\in k_w[X]\) be the
original monic type polynomial of a rank-\(n\) main summand. Write
\(u_{j,w}\in k_w^\times\) for the parameter corresponding under the
degree-one Green convention above to \((\chi_{M,w}^{\,j})^{-1}\),
where \(\chi_{M,w}\) is the original character of order \(q_2\)
and \(q_2^2\mid q-1\). At \(w\) the character \(\theta\) on the lower summand
\(\theta R\) of Section~\ref{sec:main-trace} has tame character
\(\chi_{M,w}^{-j}\), by the choice of \(\mu\). The character
\(\alpha_{\rm nor}\) multiplying the upper attachment has trivial
smooth finite-unit character there. On that attachment the first
parahoric block already carries \(\chi_{M,w}^{-j}\circ\det\), so its
inertia roots receive the same single multiplier \(u_{j,w}\) as those
of \(\theta R\) on the lower side. The absolute norm factor and
the auxiliary twist are trivial on inertia there. Pullback from
\(k_w^\times\) to the degree-\(a\) extension uses the norm exponent
\((q^a-1)/(q-1)\). Thus the same multiplier \(u_{j,w}\) acts on
each root, and the normalized polynomial is
\[
 C_w(X)=u_{j,w}^{\,n}c_w(X/u_{j,w}).
\]
The specified tame generator has eigenvalues \(\iota\jmath(\xi)\)
for the roots \(\xi\) of \(C_w\). In rank \(r\) the raw parameter
is squarefree with nonzero roots and excludes all \(q_2\)-th roots
of unity in \(k_w\). Since \(u_{j,w}^{-1}\) is one of these
excluded values, none of the scaled roots equals one. In rank one
\(c_w=X-1\), so \(C_w=X-u_{j,w}\). Its root has order
\(q_2/\gcd(q_2,j)>1\), since \(1\le j\le m<q_2\).

A root of \(C_w\) has degree \(h\le n\) over the residue
field and order dividing \(q^h-1\). It follows that the
inertia eigenvalues at the tame generator have orders
dividing \(q^{n!}-1\). The local comparison gives zero monodromy
and a sum of tame characters. Hence wild inertia acts trivially,
and the image of tame inertia is finite cyclic, generated by
the image of the chosen generator. Its order divides
\(q^{n!}-1\), which is prime to the residue characteristic.

For any eigenvalue \(\alpha\ne1\) of order dividing \(q^h-1\),
the identity \(\sum_{a=0}^{q^h-2}\alpha^a=0\) gives
\[
            v_\ell(\alpha-1)\le v_\ell(q^h-1):
\]
subtract the geometric sum at \(1\) and divide
by \(\alpha-1\), with integral quotient.
Since \(q^h-1\) divides the order of
\(\mathrm{GL}_{rt}(\mathbb F_q)\), multiplying the
eigenvalue differences proves the determinant bound
\(n v_\ell(A_{\rm aut})\).
At the conjugate place use the inverse conjugate of
the chosen element. In the polarization trace relation
the cyclotomic factor is trivial on this inertia.

Applying
it to every power of the chosen element gives equal power
traces at the two parameters. Newton's identities give equal
characteristic polynomials and hence the same estimate.
In this way the determinant bounds concern the specified
tame generator at each chosen place and its inverse
conjugate. Choosing one distinguished place on a rational prime
different from \(\ell\) supplies the common reference
element in the bound.

\subsubsection{The matrix span at the coefficient place}
\label{sec:marked-inertia-matrix-span}

\sectionopening{In this subsection, we use the inertia characters of Section~\ref{sec:marked-inertia-characters} to obtain the matrix span for the trace comparison.}

For the full matrix span, choose \(p_i^{\rm lev}\ne\ell\).
The selected tuple at \(p_i^{\rm lev}\), as defined in
Section~\ref{sec:marked-residue-projection}, has a primitive
type polynomial at a distinguished place over \(p_i^{\rm lev}\).
Its original root generates
\(\mathbb F_{q^r}^\times\). Multiplication by any element of
order dividing \(q_2\) preserves this property: the exponent
of the multiplier in a generator of this group is a multiple
of \((q^r-1)/q_2\). Every prime dividing \(q^r-1\)
divides this quotient, since \(q_2^2\mid q-1\).
Thus the new exponent remains prime to \(q^r-1\).
The normalized roots form the \(q\)-Frobenius orbit
of a generator \(\alpha_0\).

This orbit stays distinct after every \(\ell\)-power.
Indeed \(r\) is a positive power of \(3\) with \(r\ge3\).
Put \(Q'=q^{r/3}\). The factor
\((q^r-1)/(q^{r/3}-1)=Q'^2+Q'+1\) is congruent
to \(3\bmod 9\), since \(9\mid q-1\), and exceeds \(3\).
It has at least two distinct prime divisors.
It cannot divide a power of \(\ell\), so no
\(\alpha_0^{\ell^a}\) belongs to \(\mathbb F_{q^{r/3}}\).
Every proper subfield over \(\mathbb F_q\) lies in that
field. Under the faithful character map into the
coefficient field the orbit consequently has pairwise
unit differences: a root of unity reduces to \(1\)
only when its order is an \(\ell\)-power, by the
separability of roots of order prime to \(\ell\).
The same argument applies if the tame generator raises
all eigenvalues to an exponent prime to their orders.

At such a place, the orbit is the full spectrum of the
inertia generator. To see this from the annihilating polynomial,
one eigenvalue occurs, all its \(q\)-powers occur by tame
Frobenius conjugation, and there are \(r\) eigenvalues counted
with multiplicity. Raising the generator to the exponent
\((q^r-1)/(q-1)\) gives a scalar whose root-of-unity order
is \(q-1\), since every normalized root is a Frobenius
conjugate of a generator.

For the rank-one main pair, its own normalized root
already gives a scalar with order \(q_2/\gcd(q_2,j)\).
Let \(\xi\) be the resulting scalar and \(n_{\rm up}\) the upper
rank. For \(a\ge0\) and \(1\le h\le4n_{\rm up}+2\), the
difference \(\xi^{h\ell^a}-1\) is a unit. Otherwise the order of
\(\xi\) would divide \(h\ell^{a+b}\) for some \(b\), giving
\(q_2\mid jh\ell^{a+b}\) in either case.
Take \(p_*\in\{2,5\}\) different from \(\ell\). Then
\(p_*^s\mid q_2=(q_0d)^{s+6}\). The main rank is \(r\) or \(1\),
and \(n_{\rm up}=r+t\) or \(2z\le2r\), respectively. Thus
\[
 jh\le(4n_{\rm up}+2)m\le(24r+2)m=O(x^{11/2})
       <x^{11}\le m^{22}=2^s\le p_*^s .
\]
This contradicts \(p_*^s\mid jh\). In particular \(\xi^2-1\) is a
unit. The scalar comes from the original tame inertia.

For lift \(j\) and real place \(\tau\), the two conjugate jumps at
index \(u\) are \((t-1)/2\pm(\lambda_u(\tau)+w_j(\tau))\) at
rank \(r\), and \(z-1\pm c_\tau\) at rank one, by
\eqref{eq:separated-hodge-tate-intervals} and
\eqref{eq:hodge-tate-crystalline-defect}. The two sums are
\(t-1\) and \(2z-2\), respectively, both nonnegative.
Each pair thus contains at most one negative jump.

\smallskip
\noindent\textit{At the coefficient place.}
Apply \cite[Theorems~1.1 and~3.1, pp.~170 and~172]{BLGGTLocalII}
to each parity-corrected regular algebraic conjugate-self-dual
cuspidal block of rank at least two. Theorem~3.1 allows every
coefficient prime \(\ell\), including \(2\). At \(w\mid\ell\)
its conclusion identifies
the semisimplified Weil representations and bounds the Galois
monodromy partition by the automorphic partition on each
unramified twist class. For the fixed compact type the
automorphic monodromy on the whole parameter is zero by
\cite[Theorem~3.7(1) and Corollary~6.21(1)]{ShottonBM}.
In the standard-module description above that operator is
the direct sum of the monodromy operators of the shifted
cuspidal blocks. Every one of these operators is therefore zero,
with the same conclusion for repeated summands.
On a class of multiplicity \(d\), the resulting partition
bound on the Galois side is \(\lambda\preceq(1^d)\).
Its largest part is at most one, so its monodromy also vanishes.
The Frobenius-semisimplified parameters consequently agree.
For a rank-one block let \(\delta\) be its local avatar at
\(F=E_w\), \(w\mid\ell\), and fix a uniformizer \(\varpi\).
At the chosen coefficient embedding extend the coefficient field
finitely to \(L\) containing all images of \(F\) and the values of
the smooth units character \(\varepsilon_w\). In geometric
reciprocity the defining formula for the avatar on units is
\[
 \delta(\operatorname{Art}_F(u))=\varepsilon_w(u)
       \prod_{\tau:F\hookrightarrow L}\tau(u)^{b_\tau},
       \qquad u\in\mathcal O_F^\times .
\]
Here \(b_\tau\) are the exponents of the algebraic factor.
Extend \(\varepsilon_w\) by \(\varepsilon_w(\varpi)=1\) and view it
as a finite Galois character via the same Artin map.
The Lubin--Tate character \(\chi_\varpi\) satisfies
\(\chi_\varpi(\operatorname{Art}_F(u))=u\) on units and
\(\chi_\varpi(\operatorname{Art}_F(\varpi))=1\). Hence
\[
 \delta=\varepsilon_w\,\delta_{\rm ur}
                 \prod_{\tau}(\tau\circ\chi_\varpi)^{b_\tau},
 \qquad \delta_{\rm ur}|_{I_F}=1 .
\]
The unramified factor is crystalline of weight zero, and
\(\tau\circ\chi_\varpi\) is crystalline with filtration jump
\(-1\) at \(\tau\) and \(0\) at the other embeddings
\cite[Hypothesis~6.0.1, Definition~6.2.5 and Example~6.2.6(1),(3),
pp.~1091--1092, 1095--1096]{KedlayaPottharstXiao}. In particular the jump at
\(\tau\) for \(\delta\) is \(-b_\tau\); for \(F=\mathbf Q_\ell\)
and \(\varpi=\ell\), the Lubin--Tate character is \(\chi_{\rm cyc}\).
Over the field cut out by \(\varepsilon_w\), a basis \(e\) of the
finite-character line gives the period \(1\otimes e\).
Its descent action at \(s\in I_F\) is multiplication by
\(\varepsilon_w(s)\). The remaining factor is crystalline over
\(F\) and has trivial linear inertia on its periods. Thus the
inertial character on the potentially crystalline periods of
\(\delta\) is exactly its smooth units character under reciprocity.
With \(f_w=[k_w:\mathbf F_\ell]\), the eigenvalue of
\(\varphi^{f_w}\) on the crystalline factor
\(\delta\varepsilon_w^{-1}\) is
\[
 \delta(\operatorname{Art}_F(\varpi))
                \prod_{\tau:F\hookrightarrow L}\tau(\varpi)^{-b_\tau}
\]
by \cite[Example~6.2.6(1),(3)]{KedlayaPottharstXiao}.
Since \(\varepsilon_w(\varpi)=1\), this is the value of the smooth
Weil--Deligne character at the geometric Artin image of \(\varpi\).
Period spaces commute with the finite extension of coefficients,
since tensoring with a finite-dimensional coefficient space
commutes with the defining Galois invariants. Thus
these attachments are potentially crystalline at
the coefficient places, with the specified
tame inertial types. When the finite unitary representation
is spherical, Labesse's local base change \(\Pi_w\) is spherical.
Put \(K_s=\mathrm{GL}_s(\mathcal O_{E_w})\). By the Speh
description in Section~\ref{sec:cofinite-spherical-galois},
\(\Pi_w\) is a subquotient of the normalized parabolic induction
of the \(\sigma_w|\det|^{u_j}\) over all blocks and shifts.
Taking \(K_n\)-invariants is exact by compact averaging.
By Iwasawa decomposition, evaluation at \(1\) identifies the
invariants of the induction with
\(\bigotimes_{\rm blocks}\bigotimes_{j=0}^{b-1}\sigma_w^{K_a}\),
since the modulus and norm twists are trivial on the integral Levi.
As \(\Pi_w^{K_n}\ne0\), each \(\sigma_w^{K_a}\ne0\).
The unramified parity character preserves these fixed vectors. By
\cite[Theorem~2.1.1(4)]{BLGGT} each parity-corrected attachment
is crystalline there. The algebraic characters \(\alpha_j\) in
the direct sum have trivial smooth units character. Taking
\(\varepsilon_w=1\) in the preceding factorization gives crystallinity
of their avatars over \(E_w\). With arithmetic reciprocity, the Lubin--Tate character on units
is inversion, as in the proof of \cite[Proposition~B.4(i), pp.~35--36]{ConradLocalCharacters}.
When the smooth units character is trivial, the avatar therefore
restricts in this convention to the inverse algebraic character
described above. The same
criterion applies to the normalizing characters off the level.
Tensor products and direct sums preserve crystallinity.
In particular every place above \(2\) is spherical when
\(\ell=2\), since the original level and the auxiliary twist
places are odd and the finite level there is hyperspecial.

At \(w\mid\ell\) the Frobenius comparison concerns the
Frobenius-semisimplified Weil--Deligne parameter. On each shifted
block the two parity characters cancel, giving the local parameter
of \(\sigma_w|\det|^{(1-n)/2+u_j}\) after transport of coefficients.
For the inverse auxiliary avatar
\(\psi_{\rm aux}=\prod_{u\in Q}\zeta_u^{\chi_u}\) in
Section~\ref{sec:joint-twist}, each additive character \(\chi_u\)
is zero on the entire decomposition group \(G_{E_w}\) by
Sections~\ref{sec:selmer-relations} and~\ref{sec:first-auxiliary}.
Thus \(\psi_{\rm aux}^{h}|_{G_{E_w}}=1\) for every integer \(h\).
Every power preserves the local representation, its period inertia
and its labeled Hodge--Tate numbers, also when \(\ell\mid M\).

If \(\ell=p_i\), complete splitting of the original level prime
gives \(E_w=\mathbb Q_\ell\) at every \(w\mid\ell\).
Let \(B_w(s)\) denote the linear inertia action on the
potentially crystalline periods. Each eigenvalue of \(B_w(s)\)
from a degree-\(a\) tame parameter has order dividing
\(\ell^a-1\), with \(a\le r\). Polarization gives
\[
 B_{\bar w}(c_E s c_E^{-1})\sim (B_w(s)^{-1})^{\mathsf T},
\]
since the Tate power in the polarization has trivial linear period
inertia. Thus \(c_Es^{-1}c_E^{-1}\) at \(\bar w\) has the same
eigenvalue multiset as \(s\) at \(w\). The common exponent
\[
              e_r=\mathop{\mathrm{lcm}}_{1\le a\le r}
                                    (\ell^a-1)
\]
kills inertia on these periods. Put \(m_r=\operatorname{lcm}(1,\ldots,r)\).
Since \(\ell=p_i>r\), one has \(\ell\nmid m_r e_r\).
The unramified extension of degree \(m_r\) contains \(\mu_{e_r}\),
since \(e_r\mid\ell^{m_r}-1\). Adjoining an \(e_r\)-th root of a
uniformizer then gives a totally ramified extension of degree \(e_r\).
The composite is Galois of degree \(m_r e_r\), prime to \(\ell\).
Since \(m_r\mid r!\), it is contained in the chart field of
Section~\ref{sec:tame-cup-original-periods}, obtained by taking
the unramified part of degree \(r!\) and the same root.

Here the Swan conductor is that of the Weil--Deligne representation.
Its wild inertia acts trivially on the potentially crystalline periods,
since the finite inertia exponent \(e_r\) is prime to \(\ell=p_i\).
The tame inertia of the extension maps into
the \(e_r\)-th powers of the original tame quotient.
For a de Rham representation the Weil--Deligne
representation is obtained from its potentially
semistable periods. Its restriction here has trivial
inertia and zero monodromy, so the representation is
crystalline over this extension. The same extension
works for the finitely many types under consideration,
including the degree-one normalizing characters.

Fix the individual tuple of old types for the descent pair before
varying the upper system. For each inertia element the upper
spectrum, including its determinant block, consists of at most
\(r+t\) eigenvalues. An eigenvalue of order dividing
\(\ell^a-1\) lies in the unramified field of degree \(a\le r\).
The product of the distinct monic minimal polynomials over
\(\mathbf Q_\ell\) of these eigenvalues is in \(\mathbf Z_\ell[X]\)
and has degree at most \((r+t)r\). Finite inertia acts
semisimply in characteristic zero, so this polynomial annihilates
the linear period action throughout the fixed type, and hence
on subrepresentations and direct sums of arbitrary multiplicity.
Inversion of the eigenvalues gives the same degree bound on Tate
duals since the Tate character has trivial linear period inertia.
At the fixed comparison \(\kappa=(\ell,\iota)\), every
\(\mathbf Q_\ell\)-embedding of the local coefficient field
preserves this polynomial, since its coefficients lie in
\(\mathbf Z_\ell\). For another comparison above \(\ell\)
take the minimal polynomials of the resulting eigenvalues anew.
Their common degree bound \((r+t)r\) depends only on the ranks.
This gives the inertial polynomials of
Section~\ref{sec:relative-crystalline-colength} for each
coefficient place with the same bound.

\smallskip
\noindent\textit{The paired jumps.}
The raw rank-\(n\) attachment satisfies
\(R^c\simeq R^\vee\chi_{\rm cyc}^{\,1-n}\):
in each Speh string \(j\) pairs with \(b-1-j\),
and \(\alpha_j^c\alpha_{b-1-j}=|\cdot|^{a-n}\).
Combine this with the polarization of the
rank-\(a\) cuspidal attachment. At the given embedding let \(v\)
be the infinitesimal exponent of \(\sigma\) in the contragredient
compact type, and let \(b_\psi\) be the exponent of the parity
character there. It is \(0\) for \(\psi=1\) and
\(\pm\tfrac12\) for \(\psi=\upsilon\).
By \cite[Theorem~2.1.1(3)]{BLGGT} and the rank-one
computation above, the algebraic factors
\(r_{\ell,\iota}(\sigma\otimes\psi)\) and
\(r_{\ell,\iota}(\alpha_j)\) have integral jumps
\[
 \frac{a-1}{2}-(v+b_\psi),\qquad
 -\left(\frac{a-n}{2}+u_j-b_\psi\right).
\]
Their sum is \((n-1)/2-(v+u_j)\). The parentheses run through
the infinitesimal exponents of the rank-\(n\) base
change.

In the two conjugate embeddings the raw
jumps match \(u\) with \(n-1-u\).
The normalizing character and its conjugate have
product \(\chi_{\rm cyc}^{n-t_0}\) in the rank-\(n\)
attachment of a descent with rank difference \(t_0\).
Thus the normalized jumps pair as \(u,t_0-1-u\).
The main summand keeps this pairing after subtracting
the extra character string, which is invariant
under the same involution. Hence at most one
jump per pair is negative. The rank-one chart has
the same property, with jumps \(z-1\pm c_w\).

\subsection{Local length bounds for a main summand}
\label{sec:main-row-budgets}

\sectionopening{In this subsection, we use the inertial eigenvalues of
Section~\ref{sec:inertial-labels} and the stable coefficient lattices
to bound the lengths needed in Section~\ref{sec:terminal-duality}. The
estimates are normalized by the chosen valuation at \(\ell\).}

Fix the coefficient place \(\kappa=(\ell,\iota)\) and the
original main summand of rank \(n\le r\). We describe the
precise data on the summand that suffice for the previous
count. Work over a finite \(\ell\)-adic coefficient field,
enlarged to contain the eigenvalues and the images
of the designated embeddings of the local fields.
The representation is unramified outside the places above
\(\ell p_0p_1\), is de Rham at the places above \(\ell\),
and has trace relation
\[
 \operatorname{tr}\rho(g^c)
    =\chi_{\rm cyc}(g)^{-d_{\rm pol}}\operatorname{tr}\rho(g^{-1}),
 \qquad 1\le d_{\rm pol}+1\le t,\quad t=5r .
\]
Assume the prescribed tame type gives a cyclic Frobenius orbit of
\(n\) distinct inertia eigenvalues whose pairwise differences
are units at a distinguished place prime to \(\ell\).

There the norm is \(p_\nu\) for \(\nu\in\{0,1\}\).
Assume also a scalar inertia element \(s\) away from \(\ell\)
whose scalar \(\xi\) has \(\xi^2-1\) a unit. At the original
places above \(p_0p_1\) away from \(\ell\) there is a bound
on the exponent of the original inertia image prime to the
residue characteristic. At each such place a specified
inertial element has determinant of \(\rho(g)-I\) of
valuation at most \(n u_\ell\), where
\[
 A_{\rm aut}=\prod_{w\in S_M}
                      |\mathrm{GL}_{rt}(k_w)|,\qquad
 u_\ell=v_\ell(A_{\rm aut}).
\]
There are \(f=[F_0:\mathbb Q]\) chosen places over each
\(p_i\), so \(|S_M|=2f\). Local ramification can occur on both
members of a pair, giving at most \(4f\) original places.
At the common reference \(g_0\) and its inverse conjugate
require the same determinant bounds. These are local
hypotheses on the original summands.

Take a stable lattice, for example the span of the orbit of
a coefficient lattice, which is finitely generated because the
image is compact and the orbit is bounded. The integral
spectral idempotents of the tame generator split it
into rank-one summands, by Lagrange interpolation with unit
denominators. Frobenius permutes those lines with unit
multipliers in integral bases. The idempotents and the
powers of Frobenius give all the matrix units over the
valuation ring.

For \(g_1\) the \(n\)-th power of this arithmetic Frobenius,
\(\rho(g_1)=\alpha I\), since on each line its multiplier is the
product of the same \(n\) units around the cycle, and
\(\chi_{\rm cyc}(g_1)=p_\nu^n\).
The trace relation and the full matrix span show that
\[
 \rho(g)=\alpha_g I\quad\Longrightarrow\quad
 \rho((g^c)^{-1})=\chi_{\rm cyc}(g)^{d_{\rm pol}}\alpha_g I .
\]
Indeed for every \(x\), apply
\(\operatorname{tr}\rho((u^c)^{-1})
     =\chi_{\rm cyc}(u)^{d_{\rm pol}}\operatorname{tr}\rho(u)\)
to \(u=gx^{-1}\), then use the formula at \(x^{-1}\).
Pairing the difference of the claimed matrices with
all \(\rho(x^c)\) gives trace zero. The trace pairing
on the full matrix algebra is nondegenerate.
At \(s\), the cyclotomic value is one by inertiality
away from \(\ell\), giving the two evaluations in
Section~\ref{sec:terminal-duality}.

For \(V=T/\varpi^aT\), an equivariant map
\(V\to(V^*)^c\) is determined by its value on the
first coordinate vector, whose orbit spans \(V\).
The scalars of \(g_1\) on source and target are
\(\alpha\) and \(\alpha p_\nu^{n(d_{\rm pol}+1)}\), respectively.
Hence the evaluation lies in the torsion killed by
\(1-p_\nu^{n(d_{\rm pol}+1)}\). For \(0\ne u\in\mathcal O\),
the scalar action on the Cartier dual gives
\[
 V^*[u]=\operatorname{Hom}_{\mathbb Z}(V/uV,\mu_{\ell^N}),\qquad
 \operatorname{length}_{\mathcal O}V^*[u]
 =\operatorname{length}_{\mathcal O}(V/uV)
 =n\min\{a,\operatorname{ord}_{\mathcal O}(u)\}.
\]
The middle equality follows from equality of cardinalities:
\(V/uV\) has exponent dividing \(\ell^N\), as in
Section~\ref{sec:terminal-duality}.
Dividing by \(\operatorname{ord}_{\mathcal O}(\ell)\)
gives a bound of \(nv_\ell(u)\). For \(1\le n(d_{\rm pol}+1)\le rt\),
\(p_\nu^{n(d_{\rm pol}+1)}-1\) divides
\(|\mathrm{GL}_{rt}(\mathbb F_{p_\nu})|\).
This bounds the normalized global term by \(n u_\ell\).

At an original place away from \(\ell\), write \(b_v\)
for the given exponent of the inertia image.
On cocycles modulo \(\varpi^a\), the affine action on \(V\)
shows that every inertia element raised to \(b_v\) acts
by translation. Thus its order divides \(b_v\ell^j\)
for some common \(j\), prime to the residue characteristic.
The affine action kills wild inertia. Evaluation at a tame
generator has image of length at most \(\operatorname{length} V\).
On its kernel the cocycle vanishes on inertia, and its
Frobenius value is inertia-invariant and determines the
cocycle. It therefore lies in the kernel of
\(\rho(g)-I\) on \(V\) for the chosen inertia element.
By Smith normal form the excess normalized length is
at most \(v_\ell\det(\rho(g)-I)\le n u_\ell\).
Summing over \((S_M\cup\bar S_M)\setminus\{v:v\mid\ell\}\)
therefore costs at most \(4fn u_\ell\) in excess of the
baseline local lengths.

At \(v\mid\ell\), impose the full annihilator of the integral
cocycles whose rational classes lie in \(H^1_f\).
For the de Rham representation \(W\) over the coefficient field
\(K\), put \(D=D_{\rm dR}(W)\). Since \(K\) contains the images
of all embeddings \(\tau:E_v\hookrightarrow K\), the decomposition
\(E_v\otimes_{\mathbf Q_\ell}K=\prod_\tau K\) gives
\[
 d_v=\dim_K(D/\mathrm{Fil}^0D)
     =\sum_\tau\dim_K(D_\tau/\mathrm{Fil}^0D_\tau).
\]
Divide the \(\mathbf Q_\ell\)-dimension formula of Bloch--Kato
\cite[Corollary~3.8.4, p.~355]{BK} by \([K:\mathbf Q_\ell]\). This gives
\[
 \dim_K H^1_f(E_v,W)=d_v+\dim_K H^0(E_v,W),\qquad
 \dim_K B^1(E_v,W)=n-\dim_K H^0(E_v,W).
\]
Thus the preimage of \(H^1_f\) in the rational cocycles has
dimension \(n+d_v\). To see the
lattice assertion, finitely many coordinate evaluations
embed the rational cocycle space in a finite product of
coefficient fields. The integral submodule is then finite
free.

Each rational cocycle can be cleared into it by
compactness of its image. The quotient of the finite-support
evaluation module by the annihilator is a submodule of
the dual of this lattice, free of rank at most the same rank.
Allowed values modulo \(\varpi^a\) define homomorphisms
from that quotient to \(\mathcal O/\varpi^a\), so the
local length is at most \(a\) times the rank.
Thus the sum over the places above \(\ell\) satisfies
the required crystalline bound if the total number of
negative de Rham jumps is at most \(fn\), as follows
from a bound of \(n\) at each complex conjugate pair of
embeddings.

Adding the global term, the original-place excess and
the two first-prime determinant contributions give a normalized
bound of at most
\[
                  (4f+3)n u_\ell\le(4f+3)r u_\ell.
\]
This bound uses the original lattice and all its original local data. For a finite family of such summands a
common torsion cutoff is \(M=(4f+3)r u_\ell+2\).
After putting the covectors into one finite coefficient field,
an unramified enlargement increases the residue cardinality
without increasing the ramification index. It can therefore
meet the residue-field requirement in the length descent
for all the covectors at once.

The extended local equations still annihilate lattice-valued
crystalline cocycles after twisting by a character constructed in
Section~\ref{sec:joint-twist}. At each \(v\mid\ell\) its local
restriction is trivial: every auxiliary character is trivial on
the full decomposition group there by Section~\ref{sec:first-auxiliary},
and so is every product of their powers. To verify the
claim after a finite coefficient extension
\(\mathcal O\subseteq\mathcal O'\), apply an
\(\mathcal O\)-linear functional \(\mathcal O'\to\mathcal O\)
to the values of a cocycle. Applying its rational extension
to the crystalline period-ring equation shows that the
result satisfies the original local condition. Evaluating
a finite-support relation commutes with the functional,
and these functionals separate the elements of the finite free
module \(\mathcal O'\). This proves the required inclusion
of crystalline cocycles in the module tested by the
extended equations. Framing at the scalar inertia element
is preserved because the twist is unramified at its place.

\subsection{The order of the simultaneous choices}
\label{sec:simultaneous-choice-order}

\sectionopening{In this subsection, we specify which choices in
Sections~\ref{sec:auxiliary-descent} and \ref{sec:selmer-relations}
are made before the torsion order grows. This allows the finite Selmer
family and the norm systems of Section~\ref{sec:joint-twist} to use
one specialization.}

\begin{proposition}[One auxiliary twist for the Selmer and norm tests]\label{prop:joint-norm-selmer-twist}
Fix a coefficient place \(\kappa=(\ell,\iota)\), the finite family
of original integral Selmer data satisfying the hypotheses of
Sections~\ref{sec:terminal-duality} and~\ref{sec:main-row-budgets},
the finite family of norm systems of Section~\ref{sec:joint-twist},
\(m\), and \(0<\eta<1\). Use the same element \(g_0\) of inertia at an original place
away from \(\ell\) for the first-prime tests in every system. Specifically, require the
two first-prime determinant bounds of
Section~\ref{sec:main-row-budgets} at \(g_0,(g_0^c)^{-1}\), and
\(\vartheta_i(g_0)\ne1\) for every norm system whose exponent list
is nonempty. Assume all norm-system augmentations are nonzero.
One can choose a common bound \(\ell^H\) for the torsion orders before
the auxiliary level \(\ell^N\) and the data \(J,(\chi_w),(P_i)\) at that level. There is then
one tuple \((\zeta_w)_{w\in J}\) with
\[
 \zeta_w^{\ell^H}=1,\qquad \psi=\prod_{w\in J}\zeta_w^{\chi_w},\qquad
 v(P_i((\zeta_w^j-1)_w))\le\beta_i+\eta
                  \quad(\text{all }i,\ 1\le j\le m).
\]
The bounds \(\beta_i\) are independent of \(N\).
Each \(\psi^j\) is ramified on \(J\cup\bar J\). For a Selmer test
of rank \(n_i\), let \(O'\) be the ring of integers in any finite
extension of the original coefficient field containing the tuple
\((\zeta_w)_{w\in J}\). Put
\(T_{i,j}=(L_i\otimes_{O_0}O')\otimes\psi^j\), where \(O_0\)
is the original coefficient DVR, and put \(K'=\operatorname{Frac}(O')\).
Let \(\mathcal S_{i,j}\) consist of the continuous cocycles with
coefficients in \(T_{i,j}\otimes_{O'}(K'/O')\) vanishing at
the scalar inertia element of this summand and on inertia
outside the places above \(\ell p_0p_1\) and \(J\cup\bar J\).
At each \(v\mid\ell\) impose all \(O'\)-linear relations of
finite support in the cocycle evaluations that annihilate the integral
cocycles of \(T_{i,j}\) whose rational class lies in
\(H^1_f(E_v,T_{i,j}\otimes_{O'}K')\).
For the uniformizer valuation \(v_{O'}\), the valuation
\(v_\ell(\ell)=1\), and a tame generator \(\gamma_w\), one has
\[
 \frac{\operatorname{length}_{O'}\mathcal S_{i,j}}{v_{O'}(\ell)}
       \le (4f+3)n_i u_\ell+\eta,\qquad
 2\sum_{w\in J}v_\ell\bigl(\psi^j(\gamma_w)-1\bigr)\le\eta
       \qquad(1\le j\le m).
\]
\end{proposition}

\begin{proof}

Write \(\beta_i^{\rm Sel}=(4f+3)n_i u_\ell\) for the bounds
after releasing the distinguished pair. Adjoin zero to
each finite list of bounds before taking its maximum. Fix an
integer \(M>1+\max(\{0\}\cup\{\beta_i^{\rm Sel}\}_i)\)
for the Selmer rows. Make the finite unramified coefficient
extension of Section~\ref{sec:auxiliary-descent} at this \(M\).
Over the resulting fixed DVR the descent bounds \(|J|\) by
\[
 1+\sum_i\sum_{a=1}^{M v_O(\ell)}
        \operatorname{length}_O C_{i,a}(\varnothing)
\]
over the fixed coefficient DVR \(O\) of Section~\ref{sec:auxiliary-descent}.
This integer is independent of \(N\). For these \(J\) the congruence
Selmer lengths at zero satisfy the bounds after the first
place has been released as in Section~\ref{sec:selmer-relations}.
Independently,
for each norm system with character \(\vartheta_j\)
choose an integer \(M'_j>\beta_j^{\rm norm}+1\), using
its augmentation and first Euler-factor bound.

In the descent require the first Frobenius in
\(\mathcal W\), the intersection at \(g_0\) of
\[
 \begin{split}
 \rho_i(g)&\equiv\rho_i(g_0)\pmod{\ell^M},\\
 \rho_i((g^c)^{-1})&\equiv
                \rho_i((g_0^c)^{-1})\pmod{\ell^M},\\
 \vartheta_j(g)&\equiv\vartheta_j(g_0)\pmod{\ell^{M'_j}},
 \end{split}
\]
over all original main data and signed norm systems.
This is an open neighborhood by continuity. At
\(q_{w_0}\equiv1\bmod\ell^N\) and
\(N\ge\max(\{0\}\cup\{M'_j\}_j)\), the last conditions
preserve each first norm bound.
The \(M'_j\) are independent of the Selmer depth \(M\).
All these conditions refer to the original main summands.
The later trace systems at the raised level retain their
prescribed type and contain the twisted main systems. Write
\(\beta_i\) for the combined list of the Selmer and norm
bounds.

Take \(s_0\ge1\) bounding \(|J|\), the degrees of the
shifted Frobenius products in the Selmer tests, and
\(2|\mathcal E|(|J|-1)\) for each norm system. For a prescribed
number \(m\) of powers and \(0<\eta<1\), choose a positive rational
\(a_0\) with \(s_0a_0\le\eta/2\). Apply Lemma~\ref{lem:common-torsion-specialization}
with \(s_0,m,a_0\), tolerance
\(\eta/2\), and the list \(\beta_i+1\). This fixes \(H\)
before the primes or their polynomials. Choose an integer
\(N>\max\{H,M\}\) with \(Na_0\ge M\) and, for each norm system \(j\),
\[
 N\ge M'_j,\qquad Na_0>\beta_j^{\rm norm}+s_0a_0 .
\]
Since \(a_0<1\), the first scalar inequality gives
\(N\min(1,a_0)\ge M\) on every Selmer test line in
Section~\ref{sec:selmer-polynomials}. All these lower bounds
are fixed before the primes are chosen. Then choose \(J\),
the characters, and all the polynomials. Each test has a
point of depth at least \(a_0\) where its valuation is at most
\(\beta_i+s_0a_0\). One tuple \((\zeta_w)_w\) now suffices
for all tests and all powers \(1\le j\le m\).

In particular
\[
 \psi=\prod_{w\in J}\zeta_w^{\chi_w},\qquad
 v(P_i((\zeta_w^j-1)_w))\le\beta_i+\eta
                    \quad (1\le j\le m).
\]
Let \(O'\) be any coefficient DVR in the statement and let \(R\)
be the ring of integers of the compositum of \(\operatorname{Frac}(O')\)
with the field used to choose the slopes. Over \(R\),
Lemma~\ref{lem:sharp-selmer} bounds the normalized length by the
valuation of the same minor. Since \(\beta_i^{\rm Sel}+\eta<M\),
Section~\ref{sec:selmer-polynomials} gives this bound for divisible
coefficients.

At \(v\mid\ell\) the twist is trivial. Write \(T\) for the local
\(O'\)-lattice and \(T_R=T\otimes_{O'}R\).
Expand an integral crystalline cocycle in an \(O'\)-basis of \(R\).
The coordinate projections applied to its period primitive give
primitives for its coordinates, and scalar extension preserves each
primitive. Thus
\[
 Z_f(T_R)=Z_f(T)\otimes_{O'}R .
\]
On any finite set of evaluations the same basis expansion gives
the corresponding equality for the modules of relations annihilating
these cocycles. Indeed a relation over \(R\) vanishes on the
extended cocycles from \(T\) exactly when each coefficient relation
vanishes on \(Z_f(T)\). These cocycles span \(Z_f(T_R)\) over \(R\).
The same argument modulo \(\ell^a\), for every \(a\ge1\), applies
to the images of the cocycles in the finite coefficient modules,
using the reduced basis of \(R/\ell^a R\) over \(O'/\ell^a O'\).

The twist values and original matrices lie in \(O'\), so the global
cocycle equations, inertia conditions and normalization at the scalar
element also hold coordinatewise. Since
\((\operatorname{Frac}(O')/O')\otimes_{O'}R=\operatorname{Frac}(R)/R\),
this gives \(\mathcal S_{i,j}(R)=\mathcal S_{i,j}(O')\otimes_{O'}R\).
The coordinate argument applies to every finite extension of
coefficient DVRs \(O_0\subseteq O\) over which the selected twisted
lattice is obtained by scalar extension from \(O_0\). It gives
\begin{equation}\label{eq:specialized-selmer-coefficient-extension}
 \frac{\operatorname{length}_O\mathcal S_{i,j}(O)}{v_O(\ell)}
 =\frac{\operatorname{length}_{O_0}\mathcal S_{i,j}(O_0)}{v_{O_0}(\ell)},
\end{equation}
since the ramification index multiplies both the length and the
uniformizer valuation of \(\ell\). Taking \(O_0=O',\ O=R\) gives the
same normalized bound over \(O'\).
The norm tests are the values of \(\psi^j\) on their systems. Each \(v(\zeta_w-1)\le a_0/m\).
If \(\zeta_w\) has order \(\ell^h\) and \(j=\ell^u j'\) with
\(\ell\nmid j'\), then
\[
 \frac{j}{\varphi(\ell^h)}\le a_0<1
 \quad\Longrightarrow\quad u<h,\qquad
 v(\zeta_w^j-1)=\frac{\ell^u}{\varphi(\ell^h)}\le a_0 .
\]
Thus \(\psi^j\) is ramified at each \(w\in J\). At \(\bar w\) the
values on the chosen generators have the same valuation, by
anticyclotomicity and invariance under prime-to-\(\ell\) powers. Hence
\[
 \sum_{w\in J\cup\bar J}v\bigl(\psi^j(\gamma_w)-1\bigr)
                         \le 2|J|a_0\le\eta .\qedhere
\]
\end{proof}

\section{Crystalline bounds and auxiliary prime selection}
\label{sec:crystalline-auxiliary-selection}

\sectionopening{In this section, we bound the relative cocycle lengths for the lattices of
Section~\ref{sec:galois-types-and-traces} and select auxiliary primes
through paired Kummer evaluations. These choices enter the trace determinant
in Section~\ref{sec:trace-determinant-specialization}.}

\subsection{Crystalline colength for the relative cocycles}
\label{sec:relative-crystalline-colength}

\sectionopening{In this subsection, we bound the local indices in the relative cocycle
modules from the original potentially crystalline lattices and their
Tate duals. Appendices~\ref{sec:tame-cup-appendix}
and~\ref{sec:crystalline-colength-appendix} give the local bounds used
in Section~\ref{sec:main-row-budgets}.}

Let \(O=\mathcal O\) be the common coefficient DVR.
For a main or character lattice \(L\), write
\(L^\dagger=\operatorname{Hom}_O(L,O)(1)\).
Fix the lift number \(j\), its descent pair, the upper weight and the
individual tuple of old finite types selected by \(Q\) in
Section~\ref{sec:top-coefficient-clearing} before varying the upper
Hecke system \(x\). For a rank-one pair this fixes \(z\) and the
character of weights \(-c_\tau\) in \eqref{eq:hodge-tate-crystalline-defect}. The auxiliary pairs at
positive depth satisfy \(\ell^\nu\mid\mathrm Nw-1\), so no such
place lies over \(\ell\). Raising the level at these pairs preserves
each old type and preserves the hyperspecial level at
\(v\mid\ell\nmid M\). Indeed at an old split place let \(P\)
be the stabilizer in the residue group of the image of the
original block inclusion. The subgroup fixing its columns pointwise
contains the unipotent radical of \(P\) and its complementary
general linear block. Its prescribed character is trivial on the
unipotent radical. On the remaining block the projector \(Q\)
selects the chosen finite type with the determinant twist in
Section~\ref{sec:top-coefficient-clearing}. Thus each upper constituent
in this subspace contains the indicated Levi type on its
unipotent invariants. By finite Frobenius reciprocity it contains
the corresponding finite parabolic induction: the cuspidal lines
of the two blocks are disjoint, so this induction is irreducible
\cite[Proposition~6.9]{ShottonBM}. This verifies the compact-type
hypothesis of Section~\ref{sec:inertial-labels} on every upper
constituent at the raised level. At \(v\mid\ell\) with \(\ell\nmid M\), the original level
contains the integral stabilizer of the standard unimodular Hermitian
lattice. Raising the level away from \(\ell\) preserves this stabilizer.
The quadratic algebra of local integers is \'etale. After this \'etale
base change the unitary group of the lattice is \(\mathrm{GL}_n\),
since the inverse adjoint determines the second factor from the first.
Smoothness and connected reductive geometric fibres descend, so the
stabilizer is hyperspecial also at residue characteristic two.
The spherical comparison of Section~\ref{sec:inertial-labels}
therefore gives crystallinity over \(E_v\).
For \(\ell\mid M\) the annihilating polynomials for inertia on the
potentially crystalline periods at each old place are common to the
whole family for this tuple. For each \(v\mid\ell\), define \(\mathcal Q_v\) by the jump support
of the full upper list in Section~\ref{sec:coefficient-hodge-tate-intervals}
and the local conditions in Section~\ref{sec:inertial-labels}.
Off \(M\) the local class requires crystallinity over \(E_v\).
For \(v\mid\ell\) with \(\ell\mid M\) it requires potential
crystallinity and, for every \(s\in I_v\),
\[
 B_Y(s)^{e_r}=1,\qquad q_{v,s}\bigl(B_Y(s)\bigr)=0
\]
on the linear period action. Here
\(e_r=\operatorname{lcm}_{1\le a\le r}(\ell^a-1)\)
and \(q_{v,s}\in\mathbf Z_\ell[X]\) is the common monic
polynomial for the upper type of degree at most \((r+t)r\)
from Section~\ref{sec:marked-inertia-matrix-span}.
Since \(e_r\) is prime to \(\ell\), the first equation
forces trivial wild inertia. The jumps at each embedding
must lie in the specified support, with arbitrary multiplicities.
As checked below, this class contains the restrictions of \(L\)
and of every \(U_x\) for this pair and is closed under equivariant
inclusions into finite direct sums of its members.
The local class, including the jumps and inertial
polynomials, is fixed before the cocycles are chosen.
At depth \(\nu=0\) the coefficient quotient is zero.
Take \(\nu\ge1\) for the finite-level argument below.

Write \(Z_f(A)\) for the integral cocycles of a local
lattice \(A\) that become coboundaries over \(B_{\max}\).
This agrees with the usual crystalline condition:
\(\varphi(B_{\max})\subseteq B_{\rm cris}\subseteq B_{\max}\)
\cite[\S13, before Theorem~13.1]{BergerPadic}. Applying \(\varphi\otimes1\)
to a \(B_{\max}\)-primitive gives a \(B_{\rm cris}\)-primitive
for the same cocycle, and the opposite inclusion uses the
ring inclusion.

Define the cup annihilator of \(\mathcal Q_v\) in
\(Z_f(L^\dagger)\) using all members \(Y\in\mathcal Q_v\),
all depths \(\nu\), and all \(O\)-linear forms
\(Y\otimes L^\dagger\to O(1)\) equivariant modulo
\(\ell^\nu\). For each form compose with any
\(\mathbf Z_\ell\)-linear functional \(O\to\mathbf Z_\ell\).
The test is that the cup with each \(Z_f(Y)\), reduced
at \(a=\ell^\nu\), vanish in \(H^2(E_v,\mu_a)\).
By Kummer and Hilbert~90 this is the \(a\)-torsion
subgroup \(\operatorname{Br}(E_v)[a]\).
Let \(\Delta_v(L)\) be the \(O\)-length of the quotient by this
annihilator. Put \(J_{\rm pair}=\operatorname{Hom}_O(Y\otimes_O L^\dagger,O(1))\).
Since \(Y\otimes_O L^\dagger\) is finite free, every equivariant pairing
modulo \(a\) lifts to an \(O\)-linear form \(q\in J_{\rm pair}\).
Equivariance gives the integral cocycle \(z_q(g)=(gq-q)/a\)
with \(B_{\max}\)-primitive \(q/a\).
Replacing \(q\) by \(q+ab\) adds \((g-1)b\) to the cocycle.
The map \((j,c)\mapsto(j+cb,c)\) identifies the extension for \(q+ab\)
with that for \(q\) and satisfies
\[
 a(j+cb)+cq=aj+c(q+ab).
\]
Thus the integral evaluation is independent of the lift under
this identification.

\smallskip
\indent\textit{The common integral presentation.}
A relative cocycle \(z:G_E\to L/aL\) obeys all the
\(O\)-linear relations among \(d(g)=(u_x(g)-1)_x\).
Let \(W\) be the finite product of the matrix lattices
of the \(u_x\), acted on by left multiplication, and
let \(U\subset W\) be generated by the \(d(g)\).
The map \(d(g)\mapsto z(g)\) defines \(U\to L/aL\).
It is equivariant since
\[
                  h\,d(g)=d(hg)-d(h).
\]
Choose a Smith basis of \(W\) in which \(U\) is generated
by \(\varpi^{s_i}e_i\) for a subset of its indices, with
\(\varpi\) a uniformizer of \(O\). Choose \(b_1\ge1\) with
\((b_1-1)v_{\rm ord}(a)\ge\max_i s_i\) if \(U\ne0\).
Then \(a^{b_1}W\cap U\subseteq aU\). The map kills this intersection
and extends by zero on \(a^{b_1}W\) to \(Y=U+a^{b_1}W\).
Both summands are stable, so the extended map is equivariant.
The columns embed \(Y\) into copies of the upper lattices, and
\(a^{b_1}W\subseteq Y\subseteq W\) gives
\(Y[1/\ell]=W[1/\ell]\) at each completion. Thus \(Y\) has the
required jumps and local period conditions, including when \(\ell=2\).
The identity tuple \(I=(1)_x\in W\subset Y[1/\ell]\) satisfies
\(d(g)=gI-I\), so the integral cocycle \(d\) lies in
\(Z_f(Y)\). Therefore the original
relative cocycle at \(v\) is the image of a crystalline
cocycle in this same class under an equivariant map
\(t:Y\to L/aL\).

\smallskip
\indent\textit{Finite Cartier duality.}
Write \(A_\nu,A_\nu^\dagger\) for the reductions of
\(Z_f(L),Z_f(L^\dagger)\) at \(v\), and \(C_\nu\subseteq
A_\nu^\dagger\) for the reduction of the cup annihilator.
The trace against a generator of the inverse different
of \(\operatorname{Frac}(O)/\mathbf Q_\ell\) gives a
perfect \(\mathbf Z_\ell\)-pairing between \(L,L^\dagger\)
with values in \(\mathbf Z_\ell(1)\). Its reduction gives
Cartier duality between \(L/aL\) and \(L^\dagger/aL^\dagger\).
Write bars for the images of cocycles in cohomology.

The image of any relative cocycle as above annihilates
\(\bar C_\nu\). Lift \(t\) to an \(O\)-linear map \(q:Y\to L\)
by lifting the images of an \(O\)-basis of \(Y\). Write
\(b_q(y,\lambda)=\lambda(qy)\). For the usual action
on \(\operatorname{Hom}_O(Y,L)\) one has
\[
   (gq-q)(y)=g q(g^{-1}y)-q(y)\in aL,\qquad
   (g b_q-b_q)(y,\lambda)=\lambda((gq-q)y).
\]
Thus \(b_q\bmod a\) is equivariant, and the divided integral
cocycle has rational primitive \(b_q/a\). Replacing \(q\)
by \(q+au\) changes that cocycle by the integral coboundary
of \(b_u\). Use the inverse-different trace functional of the
Cartier pairing among the allowed functionals.
Cup functoriality gives the claimed orthogonality
in the Brauer group, hence in \(H^2(E_v,\mu_a)\) by the
Kummer sequence and Hilbert 90. The same reasoning with
the reduction of \(L\) itself applies to \(A_\nu\).

The two finite conditions \(\bar A_\nu,\bar A_\nu^\dagger\)
are exact annihilators. For the cardinality calculation let
\(V=L[1/\ell]\), \(h_0=\dim_{\mathbf Q_\ell}H^0(E_v,V)\),
\(d_f=\dim_{\mathbf Q_\ell}H^1_f(E_v,V)\), and let \(Q_L\)
be the torsion of \(H^1(E_v,L)\).
The integral finite condition is the preimage of the rational
one \cite[(3.7.3)]{BK}, hence is saturated and contains \(Q_L\).
Its intersection with \(aH^1(E_v,L)\) is \(aH^1_f(E_v,L)\), so
\(\bar A_\nu\simeq H^1_f(E_v,L)/aH^1_f(E_v,L)\).
The cohomology sequence for multiplication by \(a=\ell^\nu\) gives
\[
 0\longrightarrow H^0(E_v,L)/aH^0(E_v,L)
 \longrightarrow H^0(E_v,L/aL)
 \longrightarrow Q_L[a]\longrightarrow0.
\]
The free \(\mathbf Z_\ell\)-ranks of \(H^1_f(E_v,L)\) and
\(H^0(E_v,L)\) are \(d_f\) and \(h_0\), respectively.
Since \(|Q_L/aQ_L|=|Q_L[a]|\), it follows that
\[
 |\bar A_\nu|=\ell^{\nu d_f}|Q_L[a]|,\qquad
 |H^0(E_v,L/aL)|=\ell^{\nu h_0}|Q_L[a]|.
\]
Use the same identities on the dual. The dimension formula for de Rham
representations \cite[Corollary~3.8.4, p.~355]{BK} and the complementary
filtrations give
\[
 d_f+d_f^\dagger-h_0-h_0^\dagger
       =[E_v:\mathbf Q_\ell]\operatorname{rank}_{\mathbf Z_\ell}L.
\]
Rational orthogonality \cite[Proposition~3.8]{BK} gives integral
orthogonality by the injection
\[
 H^2(E_v,\mathbf Z_\ell(1))\longrightarrow
                              H^2(E_v,\mathbf Q_\ell(1)).
\]
The local invariant identifies this map with
\(\mathbf Z_\ell\hookrightarrow\mathbf Q_\ell\) by passing to the
inverse limit of the invariants on \(H^2(E_v,\mu_{\ell^j})\) and
inverting \(\ell\). The finite \(H^1(E_v,\mu_{\ell^j})\) satisfy
the Mittag--Leffler condition, so this limit computes the integral
group.
For \(M=L/aL\), finite local duality
\cite[Ch.~I, Corollary~2.3]{MilneADT} gives
\[
 |H^2(E_v,M)|=|H^0(E_v,L^\dagger/aL^\dagger)|.
\]
The local Euler characteristic \cite[Ch.~I, Theorem~2.8]{MilneADT} gives
\[
 |H^1(E_v,M)|=|H^0(E_v,M)|\,|H^2(E_v,M)|\,
               |M|^{[E_v:\mathbf Q_\ell]} .
\]
Substituting the cardinalities above yields
\[
 \bar A_\nu\subseteq(\bar A_\nu^\dagger)^\perp,\qquad
 |\bar A_\nu|\,|\bar A_\nu^\dagger|=|H^1(E_v,L/aL)|.
\]
The finite Tate pairing is perfect, so the inclusion is an equality.

If \(S_\nu\) denotes the local image of the relative
cocycles, coboundaries modulo \(a\) lift to integral
coboundaries. The kernel of passage to cohomology on
\(S_\nu\) is therefore contained in \(A_\nu\), and
\[
 \operatorname{length}_O\bigl(S_\nu/(S_\nu\cap A_\nu)\bigr)
 \le \operatorname{length}_O(\bar C_\nu^\perp/\bar A_\nu)
 = \operatorname{length}_O(\bar A_\nu^\dagger/\bar C_\nu)
 \le \Delta_v(L).
\]
The two middle quotients have equal cardinality by the perfect
Tate pairing. They are \(O\)-modules because the inverse-different
trace pairing satisfies \(\langle bx,y\rangle=\langle x,by\rangle\)
for \(b\in O\). For any finite \(O\)-module \(D\), writing
\(f_O=[O/\mathfrak m_O:\mathbf F_\ell]\), one has
\[
 \log_\ell |D|=\operatorname{length}_{\mathbf Z_\ell} D
             = f_O\operatorname{length}_O D .
\]
This gives the equality of \(O\)-lengths. The last quotient is
an image of the defining defect quotient.

Every member of \(A_\nu\) obeys all the integral linear
relations of \(Z_f(L)\). Thus imposing those linear
relations on \(S_\nu\) contributes at most \(\Delta_v(L)\).
At each place \(v\) this argument uses the restrictions of the
actual integral lattices \(Y,L,L^\dagger\) to \(G_{E_v}\).
If \(v\ne\bar v\), repeat it on \(G_{E_{\bar v}}\) with the
same integral coefficient maps. The polarization compares the
rational spectra there. The local cup pairings at the two places
are formed separately.

\begin{lemma}[Crystalline colength over the coefficient prime]\label{lem:relative-crystalline-sum}
If \(c_{\rm crys}(L,\nu)\) denotes the global colength for the
crystalline relations at the places above \(\ell\), then
\begin{equation}\label{eq:relative-crystalline-colength}
                   c_{\rm crys}(L,\nu)\le\sum_{v\mid\ell}\Delta_v(L).
\end{equation}
\end{lemma}
\begin{proof}
Write \(Z\) for the relative cocycle module at depth \(\nu\) and
\(Z_{\rm lin}\) for its submodule satisfying the crystalline relations at
all \(v\mid\ell\). With \(S_{\nu,v}\) the local restriction of \(Z\)
and \(A_{\nu,v}\) the reduction of the integral crystalline cocycles, put
\(Z_{\rm red}=\bigcap_{v\mid\ell}\operatorname{res}_v^{-1}(A_{\nu,v})\).
For finitely supported \(F_g\in\operatorname{End}_O(L)\), write
\(R_F(z)=\sum_{g\in G_{E_v}}F_g z(g)\). The relations are
those with \(R_F(Z_f(L))=0\), reduced modulo \(a=\ell^\nu\).
They all vanish on \(A_{\nu,v}\), so
\(Z_{\rm red}\subseteq Z_{\rm lin}\). Restriction induces an injection
\[
 Z/Z_{\rm red}\ \hookrightarrow\
  \bigoplus_{v\mid\ell} S_{\nu,v}/(S_{\nu,v}\cap A_{\nu,v}).
\]
The preceding finite-duality estimate bounds each summand's
\(O\)-length by \(\Delta_v(L)\), giving the assertion for
\(c_{\rm crys}(L,\nu)=\operatorname{length}_O(Z/Z_{\rm lin})\).
\end{proof}

\smallskip
\indent\textit{Summing the defects.}
On the descent pair the upper lattice for its
own parameter is \(T\oplus\bigoplus_{0\le e<t'}O(-e)\), with
\(t'\le t\). Here \(r'=\operatorname{rank}_O T\) is \(r\) in
\eqref{eq:separated-hodge-tate-intervals} and \(1\) in
\eqref{eq:hodge-tate-crystalline-defect}.
At a coefficient embedding the jumps of this direct
sum are the multiset sum of the original main jumps
and \(0,\ldots,t'-1\). In the Tate dual the jumps \(u\)
are replaced by \(-u-1\) and the inertia operator is
the transpose of the inverse period operator.
For this \(\ell\), take the coefficient fraction field to be a
finite Galois extension of \(\mathbb Q_\ell\) containing the matrix
entries, local type eigenvalues, all embedded completions \(E_v\)
and their residue-field Teichm\"uller roots. Take the normal closure
of their finite compositum. This gives one field for all places
above \(\ell\).
For a coefficient embedding extended by
\(g\in G_{\mathbf Q_\ell}\), read the jumps at the
embedding \(g^{-1}e_v\) of \(E\), where \(e_v\) is the
chosen local embedding, transported to the complex side
by the fixed comparison. This choice depends only on
the coefficient embedding, because \(g^{-1}e_v(E)\)
lies in the coefficient field. The upper class uses
the full original upper list at this embedding.

Use the local class \(\mathcal Q_v\) defined above.
An
equivariant inclusion into a sum of these lattices
at each embedding gives a jump multiset contained, with multiplicity,
in that of the ambient direct sum. Its jump support therefore lies in the
prescribed union, and the inertia operators satisfy the same polynomials. Indeed one may read
the crystalline periods after a common finite
crystallizing extension, using the exact crystalline
period functor there. In particular the direct
summands of the upper lattice for its own parameter belong
to that same class.

\subsubsection{Hodge--Tate intervals at the coefficient places}\label{sec:coefficient-hodge-tate-intervals}

\sectionopening{In this subsection, we determine the intervals for the lattices used in Appendix~\ref{sec:crystalline-colength-appendix}.}

Write \(h=r^8\), \(G=xb/x_1\), \(Q_*=kt\), \(W_*=4t\), and
\(\eta=h^2+1\). For the fixed lift \(1\le j\le m\) put
\[
 w_j(\tau)=\kappa_\tau+j\Delta+h,\qquad
 \delta_j=j\Delta+h-2t,\qquad S_j(\tau)=\kappa_\tau+\delta_j .
\]
Let \(\varepsilon_\sigma\) be \(1\) if \(\sigma\) restricts to the
chosen embedding of \(K\), and \(-1\) otherwise, and let \(\tau\)
be its infinite place. Write \(J_\sigma\) for the multiset of
main jumps. For the rank-\(r\) pair, \(t'=t\) and
\[
 J_\sigma=\left\{\frac{t-1}{2}-\varepsilon_\sigma
                 \bigl(\lambda_q(\tau)+w_j(\tau)\bigr)
                   :0\le q<r\right\}.
 \tag{H1}\label{eq:separated-hodge-tate-intervals}
\]
Indeed, if \(n=r+t\), the upper highest-weight entry \(\mu_q\)
gives the doubled parameter \(P_q=2\mu_q+n-1-2q\).

The raw jump is \((n-1-\varepsilon_\sigma P_q)/2\).
The normalizing character subtracts \(r-r_-\) on the chosen CM
type and \(r_-\) on its conjugate, where \(r_-\) counts the
negative lower signs. Thus the normalized formula replaces
\(n-1,P_q\) by \(t-1,P_q-2r_-+r\). By the parameter
correspondence of the descent pair, the latter shifted
entries consist of \(2(\lambda_q+w_j)\) and \(t-1-2a\) for
\(0\le a<t\). In the order of this parameter list the second
block reads \(a\) on the chosen CM type and \(t-1-a\) on its
conjugate. Reindexing the latter by \(a\mapsto t-1-a\) gives the
same Tate string at both embeddings. The fixed global summand
\(O(-a)\) has jump \(a\) at every embedding, and its Tate dual
\(O(a+1)\) has jump \(-a-1\).
Equivalently one gets \eqref{eq:separated-hodge-tate-intervals} directly by lowering the original
rank-\(r\) jumps by \((r-t)/2+\varepsilon_\sigma w_j\).

For a rank-one pair the parameter \(z\) and its character
of weight \(-c_\tau\) satisfy
\[
 m<z\le r,\qquad w_j+r-t\le c_\tau\le w_j,\qquad
 t'=2z-1,\qquad J_\sigma=\{z-1-\varepsilon_\sigma c_\tau\}.
 \tag{H2}\label{eq:hodge-tate-crystalline-defect}
\]
The last equality is the jump of the character obtained by
normalizing the trivial rank-one system. In both cases
\(\varepsilon_{\bar\sigma}=-\varepsilon_\sigma\), with \(\tau\)
unchanged. Corresponding main jumps sum to \(t'-1\), so their
total negative count over the two embeddings is at most \(r'\).

For \(0\le q<k-1\) we have \(\lambda_q=k/2-1-q\) by
Section~\ref{sec:real-marked-filling}. The corresponding jump in
\eqref{eq:separated-hodge-tate-intervals} therefore satisfies
\[
 |\varepsilon_\sigma u+w_j|
 =\left|\varepsilon_\sigma\frac{t-1}{2}-\lambda_q\right|
 \le \frac{t+k-3}{2}<t .
\]
The rank-one formula \eqref{eq:hodge-tate-crystalline-defect} gives
\(
 |\varepsilon_\sigma u+w_j|
 =|\varepsilon_\sigma(z-1)+w_j-c_\tau|
 \le (r-1)+(t-r)<t .
\)
For the remaining \(m\) indices of \eqref{eq:separated-hodge-tate-intervals},
\[
 \lambda_{k-2+i}+w_j=(j-i)\Delta+h=(2(j-i)+1)h
 \qquad(1\le i\le m).
\]
Since \(1\le |2(j-i)+1|\le2m-1\), either CM sign gives
\[
 \frac h2\le h-\frac{t-1}{2}\le |u|
 \le (2m-1)h+\frac{t-1}{2}<2mh\le h^2 .
\]
For every jump with \(|\varepsilon_\sigma u+w_j|<t\), the reverse
triangle inequality gives \(|u|>w_j-t>2h^2\), since
\(w_j\ge G\ge99h^3\).

At a coefficient embedding \(\gamma\), use \(\sigma\) as above
and put \(\mathcal U_\gamma=\operatorname{supp}(J_\sigma)
 \cup\{0,\ldots,t'-1\}\). For the Tate dual of an individual
summand \(L\), let \(\mathcal U'_{\gamma,L}\) be its actual
jump set. These sets lie in opposite intervals:
\[
\begin{aligned}
 \mathcal U_\gamma&\subset[-\eta,\eta]\ \cup\
        -\varepsilon_\sigma[S_j(\tau),S_j(\tau)+W_*],\\
 \{-u-1:u\in\mathcal U_\gamma\}&\subset[-\eta,\eta]\ \cup\
         \varepsilon_\sigma[S_j(\tau),S_j(\tau)+W_*].
\end{aligned}
\]
The \(1\) in \(\eta\) absorbs the dual shift for low jumps, and
the far-interval half-width \(2t\) absorbs it for high jumps.
For the main summand, every negative member of
\(\mathcal U_\gamma\cup\mathcal U'_{\gamma,L}\) is at most
\(-h/2\). Its Tate dual has at most
\(\lfloor r'm/k\rfloor\) jumps in \([-\eta,\eta]\), counted
with multiplicities: at full rank there are at most \(m\)
tail indices, and at rank one there are none. For an
individual Tate character of jump \(a\), the dual set is
\(\{-a-1\}\subset[-\eta,0]\).

\smallskip
\indent\textit{Residue-field parameters and sparsity.}
Write \(E_v\) for the embedded local field at \(v\mid\ell\),
with residue field \(k_v\), degree \(d_v\), ramification index
\(e_v^{\rm ram}\) and residue degree \(f_v\).
Then \(e_v^{\rm ram}f_v=d_v=[G_{\mathbf Q_\ell}:G_{E_v}]\)
and \(d_v\le[E:\mathbf Q]\le h\), since
\([E:\mathbf Q]\le 2q_1=O_{q_0}(x^{11}\log x)\) and \(h=r^8\ge x^{40}\).
Choose a primitive \((\ell^{f_v}-1)\)-st root \(\zeta_v\)
in the maximal unramified subfield of \(E_v\). A coefficient
embedding \(\gamma\) has parameter \(\gamma(\zeta_v)\).
The chart is indexed by all \(i\in k_v^\times\), identified with
their Teichm\"uller lifts. The parameters occurring at embeddings
are precisely \(\{\zeta_v^{\ell^h}:0\le h<f_v\}\): the
restrictions to \(\mathbf Q_\ell(\zeta_v)\) are its Frobenius
powers, each of which extends to an embedding of \(E_v\).
Witt Frobenius sends \(i\) to \(i^\ell\), so the level-\(a\) test
at \(i\) reads \(\varphi^{-a}i\).
Call a parameter thick if it occurs at an embedding whose read
infinite place lies in \(T\), where \(\kappa_\tau=b\).

If \(\ell\nmid q_1d\), the completion is unramified and
\(\zeta_v\) generates it. Two embeddings with the same parameter
agree on the completion, hence read the same \(\sigma\).
On the orbit of \(\zeta_v\) the parameter therefore determines
both thickness and sign.
If \(\ell\nmid q_1d\), then \(\ell>7\) since
\(2\cdot3\cdot5\cdot7\mid18480\mid q_0\mid q_1\).
Every prime dividing \(q_1d\) divides \(q_0d\), so is at most
\(Q_*\) for large \(x\).
If \(\ell\mid q_1\), then \(\ell\mid q_0\). By
Lemma~\ref{lem:auxiliary-residue-class}, \(K\) splits at \(\ell\),
so all embeddings at a fixed completion have the same
restriction to \(K\) and therefore the same sign. At \(\ell=d\), the
completion of the cyclotomic field is unramified, and adjoining
\(\sqrt{-d}\) is totally ramified quadratic. Each embedding of the
maximal unramified subfield consequently extends to two embeddings
of \(E_v\) with opposite CM signs. If two thick embeddings at this
completion had opposite signs, the ratio of their cyclotomic
coordinates on the chosen CM type would be \(-d^a\bmod q_0\)
for some \(a\), since Frobenius acts by \(d\) and conjugation
negates the coordinate. Both thick coordinates are \(-1\bmod q_0\),
so \(d^a\equiv-1\bmod q_0\).
This is impossible by
\(-1\notin\langle d\bmod q_0\rangle\) in
Lemma~\ref{lem:auxiliary-residue-class}.
At \(\ell\mid q_1\), set \(a=v_\ell(q_1)\) and
\(q'=q_1/\ell^a\). At either place of \(K\) the completion is
\(\mathbf Q_\ell\). The roots of order \(q'\) generate the unramified
extension of degree \(\operatorname{ord}_{q'}(\ell)\), and
\(\Phi_{\ell^a}(1+X)\) remains Eisenstein over it. At \(\ell=d\),
the roots of order \(q_1\) generate an unramified extension, over
which \(X^2+d\) is Eisenstein. Thus
\[
 (e_v^{\rm ram},f_v)=
 \begin{cases}
   (\varphi(\ell^a),\operatorname{ord}_{q'}(\ell)),&\ell\mid q_1,\\
   (2,\operatorname{ord}_{q_1}(d)),&\ell=d .
 \end{cases}
\]
These degrees are independent of the coefficient embedding.
In the ramified cases assign all parameters the common sign of
the thick embeddings, with either choice if there are none.

Suppose \(\ell\le Q_*\) and three thick parameters are
\(\omega,\omega^{\ell^a},\omega^{\ell^{a+b}}\), where
\(1\le a,b\le\lfloor\log_\ell h^3\rfloor\).
Represent them by coefficient embeddings extended to
\(g_0,g_1,g_2\in G_{\mathbf Q_\ell}\).
Put \(D=\ell^{f_v}-1\) and \(q'=(q_1)_{\ell'}\).
The image under \(e_v\) of the prime-to-\(\ell\) part of
\(\zeta_{q_1}\) has order \(q'\) and is a power of \(\zeta_v\),
so \(q'\mid D\). Write \(g_i(\zeta_v)=\zeta_v^{u_i}\) with
\(u_i\in(\mathbf Z/D\mathbf Z)^\times\). The three parameters give
\(u_1=\ell^a u_0\), \(u_2=\ell^b u_1\) modulo \(D\). The corresponding
embedding is \(g_i^{-1}e_v\), so its exponent on the root of
order \(q'\) is multiplied by \(u_i^{-1}\). If \(\tau_i\) denotes
its infinite place on the chosen CM type and \(\varepsilon_i\)
its sign on \(K\), complex conjugation multiplies the exponent
by \(-1\). Consequently, modulo \(q'\),
\[
 \tau_1/\tau_2=\varepsilon_1\varepsilon_2\,u_2/u_1
             =\varepsilon_1\varepsilon_2\ell^b,\qquad
 \tau_0/\tau_2=\varepsilon_0\varepsilon_2\,u_2/u_0
             =\varepsilon_0\varepsilon_2\ell^{a+b}.
\] Since
\[
 0<b<a+b\le 2\lfloor\log_\ell h^3\rfloor\le24s<L_0 ,
\]
where \(2\log_\ell h^3\le48\log_2 r\le24s\) for large \(x\),
by \(r<3x^5\) and \(s\ge11\log_2x\).
Also \(L_0=27\cdot3^{\lfloor\log_3s\rfloor+1}>27s\).
Thus the sparsity proved in Section~\ref{sec:datum} excludes
this triple. At \(\ell=2\), two thick parameters at any distance
\(0<i<L_0\) would already give \(4^i\equiv1\pmod{3L_0}\):
their signed ratio is \(\pm2^i\), while both places are
\(-1\bmod 3L_0\) and \(3L_0\mid q'\). Since
\(v_3(4^i-1)=1+v_3(i)\), this would force \(L_0\mid i\).

For a parameter take the union of the jump sets over its
coefficient embeddings and adjoin \(0\). At an index outside
the orbit of \(\zeta_v\), the two original sides consequently
have jump sets \(\{0\}\) after adjoining zero. This index is
thin. Assign it the thin center
\(G+\delta_j\), with either sign in the unramified case and
the common sign chosen above in the ramified case.
When \(\ell>Q_*\)
the resulting intervals are the displayed ones with the
sign and center read on the parameter. When \(\ell\le Q_*\)
the low radius can be enlarged to \(2G\), since for a thin
embedding \(\eta\le G\) and \(G+\delta_j+W_*\le2G\) put both
of its intervals in \([-2G,2G]\). In particular at \(\ell=d\) a
thick residue parameter also has an embedding of \(E_v\) with the
opposite CM sign. That embedding is thin and all its jumps
belong to \([-2G,2G]\). The other intervals still have the sign
and center read on the thick parameters. For
dual sets the signs of the separated intervals are reversed.

\smallskip
\indent\textit{The Kummer coordinate and the scale inequalities.}
Choose a uniformizer \(\pi_v\) of \(E_v\) and compatible roots
\(\underline{\pi_v}\). In Fontaine's ring
\(A_{\inf}=W(\mathcal O_{\mathbf C_\ell^\flat})\) use
\(Z=[\underline{\pi_v}]\), and embed \(W(k_v)[[Z]]\) by
evaluation, which is injective modulo \(\ell\) by the
valuation of the first nonzero term. Write \(\theta\) for
Fontaine's map \cite[Lemma~3.10 and Remark~3.11]{BMSIntegral}. Lift the unit
\(-\ell/\pi_v^{e_v^{\rm ram}}\) to a unit \(A_v\) of
\(W(k_v)[[Z]]\) by the uniformizer expansion. Then
\(Z^{e_v^{\rm ram}}A_v+\ell\) generates \(\ker\theta\):
it lies in that ideal, and its reduction modulo \(\ell\)
has the same valuation as the reduction of the generator
\([\underline{\ell}]-\ell\).

Their quotient is therefore
a unit. If \(e_v^{\rm ram}=1\), choose \(\pi_v=-\ell\)
and \(A_v=1\). The group fixing \(W(k_v)\) pointwise and
fixing \(\pi_v\) is \(G_{E_v}\), since \(E_v=W(k_v)[1/\ell](\pi_v)\).
If \(\underline\epsilon\) is a compatible
system of primitive \(\ell^n\)-th roots of unity, for every
integer \(a\ge0\) some element of \(G_{E_v}\) sends \(Z\)
to \([\underline\epsilon]^aZ\). Indeed
\(X^{\ell^n}-\pi_v\) is Eisenstein over \(E_v\), so the
closed sets imposing this action at level \(n\) are
nonempty and nested. Compactness produces the result.
The residue action of the arithmetic Frobenius is
\(y\mapsto y^{\ell^{f_v}}\).

The size conditions for these intervals are
\[
\begin{gathered}
 Q_*>0,\quad
 99\bigl(Q_*^2+(W_*+1)Q_*+Q_*+W_*+m\bigr)\le h,\\
 99h^3\le G,\qquad 99(Gh^3+G)\le b,\qquad
 2\log(4b)<m,\\
 0\le\delta_j<\eta,\qquad
 99((\eta+W_*)Q_*+\eta+h)\le h^3 .
\end{gathered}
\]
Indeed \(Q_*\le5r^2\), \(W_*=20r\), \(m<r\), \(h=r^8\),
and \(\delta_j\le(2m+1)h\). For the comparisons with \(G,b\)
use \(b\ge x_1^{10^6}\), \(x_1\ge x^{10^6}\),
\(G/b=x/x_1\), and \((\log b)^{100}\le x\).
At a prescribed level prime \(\ell\mid M\), the congruence
\(\ell\equiv1\bmod q_1d\) makes \(E_v=\mathbf Q_\ell\),
so all coefficient embeddings read the same \(\sigma\).
The tame interval threshold holds because \(\ell\ge b^{10}\):
\[
 S_j(\tau)+6(W_*+\eta)+r^2((r+t)r+1)^2+2<\ell .
\]

For a local cup exponent \(q_v\), meaning that
\(\ell^{q_v}Z_f(L^\dagger)\) lies in the cup annihilator,
the normalized defect is at most
\[
 \frac{\Delta_v(L)}{v_{\rm ord}(\ell)}
 \le q_v\,\operatorname{rank}_O Z_f(L^\dagger)
 \le 2q_v [E_v:\mathbf Q_\ell]\operatorname{rank}_O L.
\]
Put \(K_O=\operatorname{Frac}(O)\) and
\(V=L^\dagger\otimes_O K_O\). The rational cocycles \(Z_f(V)\)
are the inverse image of \(H^1_f(E_v,V)\) under passage to
cohomology. Their subspace of coboundaries has dimension
\(\dim_{K_O}V-\dim_{K_O}H^0(E_v,V)\).
The representation \(V\) is de Rham since it is potentially
crystalline. The Bloch--Kato formula
\cite[Corollary~3.8.4]{BK} on the underlying
\(\mathbf Q_\ell\)-representation, with dimensions divided by
\([K_O:\mathbf Q_\ell]\), gives
\[
 \dim_{K_O}Z_f(V)
 =\dim_{K_O}V+
   \dim_{K_O}\bigl(D_{\rm dR}(V)/\operatorname{Fil}^0\bigr)
 \le (1+[E_v:\mathbf Q_\ell])\operatorname{rank}_O L .
\]
Here \(D_{\rm dR}(V)\) has \(\mathbf Q_\ell\)-dimension
\([E_v:\mathbf Q_\ell]\dim_{\mathbf Q_\ell}V\).
The common kernel of all evaluations on \(G_{E_v}\) is zero.
While a finite collection has nonzero joint kernel, evaluation at
a point where a cocycle in that kernel is nonzero strictly lowers
its dimension. The resulting finite evaluations embed
\(Z_f(L^\dagger)\) into a finite direct sum of copies of
\(L^\dagger\), so it is finite free over \(O\).
Every rational continuous cocycle has compact image, and a power of
\(\ell\) makes it integral. Hence \(Z_f(L^\dagger)\) has rank
\(\dim_{K_O} Z_f(V)\).

For the fixed \(\mathcal Q_v\), apply the same definition of
\(\Delta_v\) to the sum of the character lattices. Then
\begin{equation}\label{eq:fixed-character-cup-annihilator}
 \Delta_v\!\left(\bigoplus_{a=0}^{t'-1}O(-a)\right)
 =\sum_{a=0}^{t'-1}\Delta_v(O(-a)).
\end{equation}
Indeed each allowed pairing restricts to one on each dual summand
at the same depth, and the cup of a cocycle tuple is the sum of the
component cups with the same partner and coefficient functional.
Conversely, every test on one summand extends by zero to the others.
Thus the cup annihilator is the direct sum of the individual
annihilators over the original coefficient DVR.

The estimates needed at \(v\) are normalized
defects at most
\([E_v:\mathbf Q_\ell]r'F_\ell\) for \(T\) and
\([E_v:\mathbf Q_\ell] F_\ell/t\) for each character,
where
\[
 F_\ell=10^6\frac{m}{k}
     \left(bv_\ell(M)+v_\ell\!\left((2b)!\!
        \prod_{\substack{b-30h_2^3<j\le b+30h_2^3\\j>0}}
                          j^{h_2^3}\right)\right),
 \qquad h_2=r^8 .
\]
\subsubsection{Uniform bounds for the local defects}\label{sec:coefficient-local-defects}

\sectionopening{At the primes dividing the original level we use Appendix~\ref{sec:tame-cup-appendix}. At the other coefficient primes we use Appendix~\ref{sec:crystalline-colength-appendix}, which applies the lifts of Appendix~\ref{sec:crystalline-lifting-appendix} to the intervals of Section~\ref{sec:coefficient-hodge-tate-intervals}.}

Put
\[
 \lambda_\ell=v_\ell\!\left((2b)!
   \prod_{\substack{b-30h^3<n\le b+30h^3\\ n>0}}n^{h^3}\right),
 \qquad K_\ell=\left\lfloor\frac{2\lambda_\ell}{Q_*}\right\rfloor .
\]
Then \(\ell(\lambda_\ell+1)>2b\), since each multiple of
\(\ell\) up to \(2b\) contributes at least one to the factorial
valuation. If \(\lambda_\ell<h^3\), no integer in
\((b-30h^3,b+30h^3]\) is divisible by \(\ell\). In this case set
\(N_\ell=6(2(\lambda_\ell+1)+\eta+W_*)\). The scale inequalities
give \(\delta_j+N_\ell<30h^3\) and
\[
 \ell>\frac{2b}{h^3+1}>G+\delta_j+N_\ell.
\]
For either center \(S=\kappa_\tau+\delta_j\), no nonzero multiple
of \(\ell\) is within \(N_\ell\) of \(0\) or of either signed
center \(\pm S\). For \(S=G+\delta_j\) all such integers have
absolute value smaller than \(\ell\). For \(S=b+\delta_j\) take
the absolute value in the separated intervals to reduce to the preceding
exclusion near \(b\).

For fixed \(q_0,d\), the prime divisors of \(q_1d\) are those
of \(q_0d\), since \(q_1/q_0=3L_0\,2^s\) and \(6\mid q_0\).
As \(Q_*\to\infty\), every \(\ell>Q_*\) is therefore unramified
in \(E\) for sufficiently large \(x\). Moreover
\[
 \ell\le Q_*\quad\Longrightarrow\quad
 \lambda_\ell\ge\lfloor2b/\ell\rfloor
       \ge\lfloor2b/Q_*\rfloor\ge h^3 ,
\]
since \(b/(Q_*h^3)\to\infty\).
For these small primes the factorial bound also implies
\[
 Q_*(K_\ell+3)>2\lambda_\ell+2>4b/Q_*.
\]
Since \(b/Q_*^2\) grows faster than any power of \(r\), this gives
\(K_\ell\ge 2b/Q_*^2\) for large \(x\). In particular
\(2(b+\delta_j)<h^2 K_\ell\).
For the normalization of the local defect, write \(\rho\) for the
rank of the summand being tested and \(c_*=1\) on the main
summand, \(c_*=t\) on each Tate character. Put \(\mu=\lfloor\rho m/(kc_*)\rfloor\), and write
\(d_v=[E_v:\mathbf Q_\ell]\) and \(e_O=v_{\rm ord}(\ell)\).
For \(\ell\nmid M\), Proposition~\ref{prop:crystalline-off-level-colength}
applies to the actual dual of \(L\) and the whole upper jump set at the
fixed lift. It gives
\[
 \frac{\Delta_v(L)}{e_Od_v}
 \le 10^5(\rho K_\ell m+\mu\lambda_\ell)
 \le 3\cdot10^5\frac{\rho m\lambda_\ell}{kc_*}
 \le \rho F_\ell/c_* .
\]
The last two comparisons use \(Q_*=kt\), \(c_*\le t\), whence
\(0\le K_\ell\le 2\lambda_\ell/(kt)\le2\lambda_\ell/(kc_*)\), and
\(\mu\le\rho m/(kc_*)\).
For \(\ell\mid M\), Proposition~\ref{prop:tame-cup-exponent}
with \(\eta=h^2+1\), \(W_*=4t\) and
\(k_0=r^2((r+t)r+1)^2\) gives the sufficient exponent
\[
 q_{\rm tame}=3r^2((r+t)r+1)^2+150(h^2+4t+2).
\]
For \(i=0,1\), Proposition~\ref{prop:marked-level-primes} gives
\(p_i\ge b^{10}>\max(2b,b+30h^3)\) for sufficiently large \(x\).
Since \(M=p_0p_1\) with distinct primes,
\[
 v_{p_i}(M)=1,\qquad \lambda_{p_i}=0,\qquad
 F_{p_i}=10^6\frac{m}{k}b,\qquad k=r+1-m .
\]
At each \(v\mid p_i\), complete splitting gives \(d_v=1\) and
\(\operatorname{rank}_O Z_f(L^\dagger)\le2\rho\). The inclusion of
\(p_i^{q_{\rm tame}}Z_f(L^\dagger)\) in the annihilator gives
\[
  \frac{\Delta_v(L)}{e_Od_v}\le 2\rho q_{\rm tame}
                              \le \rho F_{p_i}/c_* .
\]
The last inequality follows from \(c_*\le t\), \(m\ge1\), and
\(ktq_{\rm tame}=o(b)\).

With these estimates, \eqref{eq:relative-crystalline-colength}
and \(\sum_{v\mid\ell}[E_v:\mathbf Q_\ell]=2f\) give the
following normalized bound at every depth \(\nu\):
\[
 \frac{c_{\rm crys}(T,\nu)+\sum_{0\le j<t'}c_{\rm crys}(O(-j),\nu)}
      {r'e_O}
 \le 2f\Bigl(1+\frac{t'}{r't}\Bigr)F_\ell\le 4fF_\ell .
\]
This is the crystalline contribution per unit rank in
\eqref{eq:display-rb1}.
A uniform cup exponent bounds the index by the rank of
the cocycle lattice. In the small-loss range for the main summand,
Appendix~\ref{sec:crystalline-colength-appendix} bounds the colength
of a common sublattice directly. Both bounds hold at every finite
pairing depth \(\nu\).

For sufficiently large \(x\), \(b+30h^3<2b\).
Thus \(F_\ell=0\) for \(\ell>2b\) with \(\ell\nmid M\).
The same primes have zero defect by the local bound.
For a finite set \(I\) of coefficient primes, put
\(B=\log b\) and
\(P_h=\prod_{b-30h_2^3<j\le b+30h_2^3,\ j>0}j\).
By nonnegativity and factorization over the finite support we have
\[
 \sum_{\ell\in I}4fF_\ell\log\ell
 \le\sum_{\ell\ {\rm prime}}4fF_\ell\log\ell
 = 4\cdot10^6 f\frac{m}{k}
       \bigl(b\log M+\log((2b)!)+h_2^3\log P_h\bigr).
\]
The half-open interval in \(P_h\) contains at most \(60h_2^3\)
positive integers. Since \(h_2^6=r^{48}<(3x^5)^{48}\) and
\(b\ge x_1^{10^6}\ge x^{10^{12}}\), for fixed \(q_0,d,\epsilon\)
and all sufficiently large \(x\) we have \(120h_2^6\le b\). Hence
\[
 h_2^3\log P_h
 \le 60h_2^6\log(b+30h_2^3)\le 120h_2^6 B\le bB .
\]
Both primes dividing \(M\) satisfy \(p_i\ge b^{10}>2b\), so
\(\lambda_{p_i}=0\). Appendix~\ref{sec:crystalline-colength-appendix}
treats the primes \(\ell\nmid M\) in three disjoint cases,
\(\ell\le Q_*\), \(\ell>Q_*\) with \(\lambda_\ell\ge h_2^3\),
and \(\ell>Q_*\) with \(\lambda_\ell<h_2^3\).
The triangle in Section~\ref{sec:triangle-unit-input} uses
\(Y_{\rm tr}^2=h_2^2>Q_*\). Throughout \(Q_*<\ell\le h_2^2\),
\(\lambda_\ell\ge\lfloor2b/\ell\rfloor\ge h_2^3\), since
\(b\) dominates \(h_2^5\). Thus these primes use the two-start
colength estimate and the small-prime triangle estimate. The
contribution of \(\lambda_\ell\) is counted once per rational
coefficient prime in the displayed sum.
Also \(M\le b^{22}\) and \(\log((2b)!)\le 4bB\).
Dividing the right side by \(fxbB^3/(mk)\) gives at
most \(108\cdot10^6 m^2/(xB^2)\), which tends to zero
because \(m^2\le4x\) and \(B\to\infty\).

\subsection{Cyclotomic-kernel evaluations for the auxiliary primes}
\label{sec:triangle-unit-input}

\sectionopening{In this subsection, we refine the finite detection argument of
Section~\ref{sec:odd-cocycle-budget} to choose elements in the
cyclotomic kernel. The normalization of Section~\ref{sec:main-trace}
lets these elements test the representations in the prime selection of
Section~\ref{sec:triangle-selection}.}

The finite set of evaluation elements in Section~\ref{sec:odd-cocycle-budget}
constructed by adjoining \(\mu_{N_c\ell^2}\) works for every
\(\ell\). The bound \(y^{20}\) on its size was used only for
small \(\ell\). To put its points in the cyclotomic kernel,
let \(P\) be that finite set and take a cyclotomic generator
\(g_{\rm gen}\), meaning that its cyclotomic value generates the
image of \(G_E\) topologically. For \(p\in P\) choose \(q_p\)
in the closure of \(\langle g_{\rm gen}\rangle\) with \(\chi_{\rm cyc}(q_p)=\chi_{\rm cyc}(p)\).
Such a choice exists by continuity and compactness.
Set \(s_p=pq_p^{-1}\), so \(s_p\in\ker\chi_{\rm cyc}\).
The zero set of a continuous cocycle is a closed subgroup.
If it contains \(g_{\rm gen}\) and every \(s_p\), it contains all
\(p=s_pq_p\), so the cocycle is zero. This procedure detects
all unramified-outside-\(\ell\) cocycles in the stated
cyclotomic characters, without a parity condition. The
chosen points are independent of \(e,\nu,\mathcal O\).

Put \(Y_{\rm tr}=r^8=\Delta/2\), so \(N_c+t+3\le Y_{\rm tr}\) for
large \(x\). The image of \(\chi_{\rm cyc}\) on \(G_E\)
is procyclic, since it is a closed subgroup of
\(\mathbb Z_\ell^\times\) for odd \(\ell\), and of
\(1+4\mathbb Z_2\) at \(\ell=2\) since \(i\in E\).
Choose \(\gamma_\ell\in G_E\) whose cyclotomic value
\(\xi_\ell\) generates that image topologically. If
\(w\) is any value in the image, \(\xi_\ell^a-1\)
divides \(w^a-1\), by taking limits of powers of \(\xi_\ell\)
in the closed ideal generated by \(\xi_\ell^a-1\).
For \(\ell\le Y_{\rm tr}^2\) take \(w=p_i^{\rm lev}\) from
arithmetic Frobenius at a distinguished place not above \(\ell\).
For \(0<a<t\), the integer \(w^a-1\) divides \(A_{\rm aut}\)
since it is a factor of the order of the general
linear group of rank \(rt\). For \(\ell>Y_{\rm tr}^2\), inertia at
\(\ell\) fixes \(E\subseteq\mathbb Q(\zeta_{N_c})\) and
its cyclotomic reduction fills \(\mathbb F_\ell^\times\).
Choose \(w\) reducing to a generator. Then \(w^a-1\) is
a unit for \(0<a<t\). Thus, with
\[
 h_\ell=\begin{cases}v_\ell(A_{\rm aut}),&\ell\le Y_{\rm tr}^2,\\
                 0,&\ell>Y_{\rm tr}^2,\end{cases}
\]
every \(\xi_\ell^a-1\), \(0<a<t\), divides \(\ell^{h_\ell}\).

For the triangle take \(g_{\rm gen}=g_{\rm tr}=\gamma_\ell^{-1}\),
whose cyclotomic value also generates the image, and use the characters \(\chi_e=\chi_{\rm cyc}^{-e}\). Then \(\chi_e(g_{\rm tr})=\xi_\ell^e\).
For the actual rank increment \(t'\le t=5r\) of the fixed pair, define
\[
 \Delta_e=\prod_{\substack{0\le j<t'\\j\ne e}}
              (\xi_\ell^e-\xi_\ell^j),\qquad C=\ell^{Y_{\rm tr}h_\ell}.
\]
The scalar \(4\Delta_e\) divides \(C\) for \(0\le e<t'\), and
\(4(\xi_\ell^e-1)^2\) divides \(C\) for \(0<e<t'\).
Indeed, for \(j\ne e\) the difference \(\xi_\ell^e-\xi_\ell^j\)
is a unit times \(\xi_\ell^{|j-e|}-1\), with valuation
at most \(h_\ell\). At
\(\ell=2\), \(4\mid \xi_\ell-1\) implies \(v_2(4)\le h_\ell\);
at odd primes \(4\) is a unit. The inequalities
\(t-1+1\le Y_{\rm tr}\) and \(2+1\le Y_{\rm tr}\) prove the claims.

Let \(R\) be the raw lower attachment of rank \(n\) and
\(T=\vartheta_\ell^{-1}\otimes\theta R\) the main summand in the
normalized trace of Section~\ref{sec:main-trace}. At rank one the raw
attachment is trivial and \(\theta\) supplies the inverse old tame
character. Sections~\ref{sec:marked-inertia-characters}
and~\ref{sec:marked-inertia-matrix-span} supply the off-\(\ell\)
inertia element \(\gamma\) such that \(T(\gamma)=\xi I\),
\(\chi_{\rm cyc}(\gamma)=1\), and
\(\xi^{a\ell^\nu}-1\) is a unit for \(1\le a\le n\).
Here the auxiliary character is unramified at the original inertia used for \(\gamma\). Fix any \(g\in G_E\)
and depth \(\nu\). The eigenvalues of \(T(g)\) are integral,
since the image is compact and preserves a lattice over
its coefficient field. For \(\lambda\) in this multiset,
at most one \(i\in\{0,\ldots,n\}\) satisfies
\(\overline{\xi^{i\ell^\nu}\lambda}=1\), where bars denote
residue classes. Indeed the \(n+1\) powers of
\(\overline{\xi^{\ell^\nu}}\) are distinct. Choose \(i\)
outside all the forbidden indices, and set \(\delta=\gamma^i\).
Then
\[
 \chi_{\rm cyc}(\delta)=1,\qquad
 \det(1-T(\delta^{\ell^\nu}g))
                   =\prod_\lambda(1-\xi^{i\ell^\nu}\lambda)
                         \in\mathcal O^\times .
\]
The unit condition persists on an open neighborhood of
\(\delta^{\ell^\nu}g\) by continuity. This is the determinant
condition when selecting new Frobenius points for the
triangle.

An original-level system \(x\) with normalized representation \(N_x\)
is a weak lift when \(\det(1-N_x(g))=0\) for every \(g\in G_E\).
At an off-level split spherical place with residue norm
\(q\), let \(N_H=r_0+t'\) be the rank of a weak lift and
\(t'>0\) the chain length for this descent pair. Write \(H\) for its monic Hecke
polynomial. Suppose the normalized weak-lift representation
has eigenvalue \(1\) at geometric Frobenius \(\phi\), and
the unnormalized representation at \(\phi\) is killed by \(H\)
in the coefficient embedding. If \(R^{\rm nor}=\chi R\), the unit
eigenvalue of \(R^{\rm nor}(\phi)\) gives the root
\[
 X_0=\chi(\phi)^{-1},\qquad |X_0|=q^{r_0/2}
\]
in the complex realization, by the normalization of
Section~\ref{sec:main-trace}.
At this split place the original local unitary-group representation is
spherical and unitarizable. Apply the classification before the common
arithmetic shift.

Tadi\'c's spherical unitary classification
\cite[Theorem~3.2, p.~54]{TadicSpherical} expresses the Satake
roots in strings. In arithmetic normalization a string
of length \(h\) has terms
\[
 a q^{(N_H-1)/2+\theta-(h-1)/2+j},
 \qquad 0\le j<h,\quad |a|=1,\quad -\tfrac12<\theta<\tfrac12 .
\]
The generators in the classification are unitary determinant
characters and complementary pairs of their norm twists.
The norm of a uniformizer is \(q^{-1}\); reversing each
string writes its powers of \(q\) in increasing order.
The shift \((N_H-1)/2\) is the arithmetic normalization
of Section~\ref{sec:cofinite-spherical-galois}.
If \(X_0\) is at position \(u\) in this string, with the
positions starting at zero, subtracting the common shift gives
\[
 \theta-(h-1)/2+u=(1-t')/2,qquad
 2\theta=h-t'-2u\in\mathbb Z .
\]
The strict bound \(|\theta|<1/2\) forces \(\theta=0\)
and \(h=t'+2u\). The \(h-1-u=t'-1+u\) steps above
\(X_0\) include the distinct roots \(X_0q^j\) for
\(0\le j<t'\).

Put \(B_{t'}(X)=\prod_{0\le j<t'}(X-q^j)\). In the coefficient embedding the exact chain gives
\[
 H(\chi^{-1}X)=B_{t'}(X)U(X),\qquad \chi\in\mathcal O^\times,
 \qquad \deg H=N_H,
\]
where \(\chi\) is the normalizing character at the chosen
\(\phi\). Let \(C(X)\) be its normalized Galois characteristic
polynomial. Assume the coefficientwise comparison
\(C\equiv\chi^{N_H}H(\chi^{-1}X)\) in the residue field,
with the rescaled polynomial integral, and a residual
factorization \(C\equiv S B_{t'}\) with \(S(1)\) a unit.
At a depth-\(\nu\) triangle prime, \(q\equiv1\bmod\ell^\nu\).
For \(\ell=2\) one also has \(q\equiv1\bmod4\), since the rational
prime splits completely in \(E\supset\mathbb Q(i)\). Thus for
\(0\le i<j<t'\),
\[
 v_\ell(q^j-q^i)=v_\ell(q-1)+v_\ell(j-i)\ge\nu .
\]
Indeed an exponent prime to \(\ell\) leaves the valuation of
\(q-1\) unchanged, and raising to the \(\ell\)-th power increases
it by one by the binomial formula, using \(4\mid q-1\) at \(2\).
Thus the roots are distinct in characteristic zero and all reduce to
\(1\) modulo \(\ell^\nu\); in particular
\(\overline B_{t'}=(X-1)^{t'}\) in the residue field. Monic division in \(\mathcal O[X]\), followed
by cancellation of \(\overline B_{t'}\) in the residue polynomial
ring, gives
\[
 \chi^{N_H}U\in\mathcal O[X],\qquad
                    \chi^{N_H}U\equiv S .
\]
Thus \(U(1)\) is a unit and \(U/U(1)\) is an integral
polynomial of value one at \(1\). For finitely many
residual-near weak lifts at the same prime, multiply
these polynomials to obtain a single such polynomial
divisible by their normalized quotients. 

For an original system \(x\) put \(C_x(g;X)=\det(X-N_x(g))\).
The product uses the weak lifts satisfying
\[
 C_x(g;1)=0,\qquad
 C_x(g;X)\equiv C_{x_0}(g;X)\pmod{\mathfrak m_{\mathcal O}}
 \qquad (g\in G_E),
\]
where \(x_0\) is the distinguished normalized lift and the
congruence is coefficientwise. Before choosing the auxiliary primes,
take a basis of the original-level space for each of the finitely many
cofinite spherical systems. By definition each basis vector has the
specified spherical eigenvalues away from finitely many split places.
Exclude the union of these exceptions and the exceptions in the Galois
comparisons of Section~\ref{sec:cofinite-spherical-galois}.
At each remaining split place every spherical operator acts by its
specified scalar on the full original space for each system. In
particular, at each subsequently chosen auxiliary prime, the root chain
used in Section~\ref{sec:top-coefficient-clearing} applies to every
vector of each original system in the product.
The trace identity of Section~\ref{sec:main-trace},
applied to every power of \(g\), and Newton's identities
in the characteristic-zero fraction field give the distinguished polynomial as the product of
that of \(T(g)\) and
\(\prod_{0\le j<t'}(X-\chi_{\rm cyc}(g)^{-j})\).
In the determinant-unit neighborhood the first factor
has value a unit at \(1\). Thus the residual quotient
condition holds for all the near systems there. The
weak-lift condition gives eigenvalue \(1\) in each
normalized representation itself. The local argument
applies simultaneously to all these systems.

Eligible Frobenius elements of prime norm are dense
in this neighborhood, also after any further finite
exclusion. For any point \(h\) of the given open set, choose an
open normal subgroup \(V\triangleleft G_E\) with \(hV\) in that set.
Let \(L/\mathbb Q\) be the normal closure of the fixed field of \(V\)
and put \(h_L=h|_L\in\operatorname{Gal}(L/E)\).
Chebotarev over \(\mathbb Q\) gives a positive density of
unramified rational primes with arithmetic Frobenius conjugate to
\(h_L^{-1}\) in \(\operatorname{Gal}(L/\mathbb Q)\)
\cite[Theorem~8.31]{MilneANT}.
Choose a prime of \(L\) over such a prime whose arithmetic Frobenius
is exactly \(h_L^{-1}\), and prolong it to \(\overline{\mathbb Q}\).
Since \(E/\mathbb Q\) is Galois and \(h_L|_E=1\), the rational prime
splits completely in \(E\). Its geometric Frobenius at the chosen
prolongation belongs to \(G_E\) and restricts to \(h_L\), hence belongs
to \(hV\). Excluding the rational primes below the earlier places and
the finite exceptional sets gives distinct primes of \(F_0\). This proves abundance at
\(\delta^{\ell^\nu}g\) for
Section~\ref{sec:triangle-selection}, on the
stated scalar inertia and distinguished lift data.
The finite family and its off-level comparisons
are fixed before choosing the depth \(\nu\).

\subsection{The fresh Kummer maps}
\label{sec:fresh-kummer-maps}

\sectionopening{In this subsection, we construct the Kummer maps for adding a prime to
the places allowed to ramify. Their local restrictions are the ones
prescribed in the inductive selection of
Section~\ref{sec:triangle-selection}.}

Fix \(\nu>0\) and \(e>0\). A \(\nu\)-split place is a
degree-one prime \(w\nmid\ell\) of \(E\) with a chosen prolongation
and geometric Frobenius there fixing \(\mu_{\ell^\nu}\).
Take \(k_0\) earlier such places of type \(e\) with distinct
restrictions to \(F_0=E^+\). Let \(S_0\) be the set of their primes in
\(F_0\) together with the primes over \(\ell\). Put
\(F_\nu=E(\mu_{\ell^\nu})\), \(U=G_{F_\nu}\), and
\(\kappa=\chi_{\rm cyc}\chi_{E/F_0}\) on \(G_{F_0}\),
where the second factor has values in \(\{\pm1\}\).
Then \(\kappa|_{G_E}=\chi_{\rm cyc}\) and
\(\kappa(c)=1\) at every real conjugation. It is
unramified off \(\ell\) since \(E/F_0\) is finite
unramified. Let
\[
 M_\nu=(\mathbb Z/\ell^\nu\mathbb Z)(\kappa^{-e}),
 \qquad A=M_\nu[\ell],
 \qquad M_\nu^\vee=\operatorname{Hom}(M_\nu,\mu_{\ell^\nu}).
\]
The earlier primes of \(F_0\) split in \(F_\nu\),
so both coefficient modules are trivial there.

Use primal local conditions of zero at real places and
at the earlier primes, vanishing on
\(D_v\cap\ker\kappa\) at \(v\mid\ell\), and unramified
elsewhere. These give conditions on cohomology classes
as coboundaries obey the required vanishings.
Let \(P\) be a finite set of auxiliary primes split in
\(F_\nu\), outside \(S_0\), where the primal condition is
relaxed to the full local group. Write \(Z_\nu(P)\) for the
resulting group of actual cocycles. The same conditions
for \(A\) give \(Z_A(P)=Z_\nu(P)[\ell]\).

To choose \(P\), take all residual dual classes with
coefficients in \(A^\vee=\operatorname{Hom}(A,\mu_\ell)\)
unramified outside \(S_0\) and infinity. These classes
form a finite group. Their restrictions to \(U\) are
homomorphisms. For each nonzero restriction, take an
element of \(U\) detecting it. For a cocycle representative
\(b\), the map \(g\mapsto(b(g),g|_{F_\nu})\) is a homomorphism
to \(A^\vee\rtimes\operatorname{Gal}(F_\nu/F_0)\), since the
action on \(A^\vee\) factors through this finite Galois group.
Intersect the kernels for the finitely many representatives with \(U\).
In the resulting finite quotient of \(G_{F_0}\), each detecting
element has trivial image in \(\operatorname{Gal}(F_\nu/F_0)\).
Chebotarev applied to its inverse conjugacy class gives a prime
splitting completely in \(F_\nu\). Choose its prolongation so that
geometric Frobenius has the detecting image. The restriction of
\(b\) to \(U\) is independent of the representative, since \(U\)
acts trivially. Exclude \(S_0\) and the earlier choices at each step.

At finite \(v\notin S_0\), the coefficients are unramified
of order prime to the residue characteristic. For either \(A\)
or \(M_\nu\), the unramified conditions for the module and its
Cartier dual are exact annihilators
\cite[Ch.~I, Theorem~2.6]{MilneADT}. At \(P\) the relaxed
primal condition is the full group, with zero annihilator.
Thus the residual dual classes are unramified outside \(S_0\)
and infinity and vanish on the chosen Frobenius elements at \(P\).
By the choice of \(P\) their restrictions to \(U\) vanish.
Inflation--restriction gives the inclusion in
\(H^1(\operatorname{Gal}(F_\nu/F_0),A^\vee)\).
The Galois group is generated by at most two elements.
For odd \(\ell\) it is a subgroup of the product of
the order-two group of \(E/F_0\) and the cyclic units
modulo \(\ell^\nu\). For \(\ell=2\), use \(E=F_0(i)\):
at \(\nu=1\) the group has order two and for \(\nu\ge2\)
it is a subgroup of the units modulo \(2^\nu\).
Cocycles are determined by their values on generators,
so the residual dual group has dimension at most two.

Apply Greenberg--Wiles to \(A\) with the relaxed conditions
\cite[Theorem~2.19]{DDT}, which applies to finite modules
of arbitrary order. At \(P\) the local factor is
\(\ell\) by the unramified and tame directions in local
\(H^1\). Each earlier prime contributes \(\ell^{-1}\).
At a real place \(v\), the local condition is zero and
\(\kappa(c)=1\), so \(|\mathcal L_v|/|H^0(F_{0,v},A)|=\ell^{-1}\),
also for \(\ell=2\).
At \(v\mid\ell\), an allowed cocycle factors through the
procyclic \(\kappa(D_v)\) and has at most \(\ell\) possible
values on a lift of a generator. After division by
\(\ell/|H^0(F_{0,v},A)|\) coboundaries the local factor is
at most one. Elsewhere the unramified factor is one. At \(\ell=2\), the character of
\(E=F_0(i)\) is the sign modulo \(4\) of
\(\chi_{\rm cyc}\), so \(\kappa\) takes values in
\(1+4\mathbb Z_2\) and its local images are procyclic.
The global coboundary factor \(\ell/|H^0(F_0,A)|\)
cancels the numerator of the global invariant ratio.
With the residual dual bound this gives
\[
  |Z_\nu(P)[\ell]|\le\ell^{|P|-k_0-f+3}
                         \le \ell^{|P|-k_0-1},
\]
since \(f\ge4\).

At full depth let \(B\subset M_\nu^P\) be the image
of tame evaluation at \(P\) on \(Z_\nu(P)\). Let
\(B^\vee\subset(M_\nu^\vee)^P\) be the Frobenius evaluations
of the global dual classes with the orthogonal conditions
for the unrelaxed primal group. Put \(d_\nu=\ell^\nu\).
At every \(v\in P\) both coefficient modules are trivial.
Use a single primitive \(d_\nu\)-th root of unity \(\zeta_{d_\nu}\) for
the Kummer coordinates, transported to each completion by its chosen
prolongation. For a uniformizer \(\varpi_v\) of \(F_{0,v}\), choose
\(\tau_v\) in tame inertia with
\[
 \tau_v(\varpi_v^{1/d_\nu})/\varpi_v^{1/d_\nu}=\zeta_{d_\nu}.
\] A primal class with
tame value \(t\) differs from \(t\kappa_{\varpi_v}\) by an
unramified class. Two unramified classes have zero cup product,
since the unramified quotient has \(\ell\)-cohomological dimension one.
For the unramified dual coordinate \(\chi\), the calculation of
Section~\ref{sec:terminal-duality} gives
\[
 \operatorname{inv}_v(\chi\cup t\kappa_{\varpi_v})
      =-t\,\chi(\phi_v)/d_\nu\quad\text{in }\mathbf Q/\mathbf Z ,
\]
where \(\phi_v\) is geometric Frobenius. Indeed the local
Artin map in \cite[Ch.~III, Proposition~3.6 and Lemma~3.7]{MilneCFT}
sends \(\varpi_v\) to arithmetic Frobenius on the unramified
quotient. Reversing the cup factors changes the sign at every \(v\).
In the chosen root-of-unity coordinates the product of the
local pairings is the perfect pairing on the two copies of the
coefficient modules indexed by \(P\). Global cup reciprocity
gives \(B^\vee\subset B^\perp\).
The kernels of tame and Frobenius evaluation are the
unrelaxed primal and the strict dual groups, respectively.
The ratio of the two Greenberg--Wiles identities for
\(M_\nu\) with and without relaxation at \(P\) yields
\[
                     |B|\,|B^\vee|=\ell^{\nu|P|}.
\]
Indeed each relaxed prime changes the local numerator of the
Greenberg--Wiles formula from the order of the unramified line
to the order of the full local group, a factor \(\ell^\nu\).
The global invariant terms cancel. Since
\(|B^\perp|=\ell^{\nu|P|}/|B|\), this proves \(B^\vee=B^\perp\).

Every character of \(A^P/(B\cap A^P)\) extends to
\(M_\nu^P/B\) by divisibility of \(\mathbb Q/\mathbb Z\),
and therefore arises from one of the dual classes.
The coefficient map here is restriction
\[
 M_\nu^\vee\longrightarrow A^\vee,\qquad
 b\longmapsto b|_{\ell^{\nu-1}M_\nu}.
\]
It is Galois-equivariant. In coordinates from the chosen roots,
\(b(\ell^{\nu-1}x)=b(x)^{\ell^{\nu-1}}\in\mu_\ell\),
so it reduces the exponent of \(b\) modulo \(\ell\).
On \(U\) coboundaries in \(M_\nu^\vee\) vanish. The
resulting homomorphism on \(U\) determines these values on
\(A^P\) since the Frobenius elements of \(P\) belong to \(U\).
Moreover \(|B[\ell]|\le |Z_\nu(P)[\ell]|\), by surjectivity
on the quotients by multiplication by \(\ell\) and the
equality of the kernel and cokernel sizes of multiplication
on a finite abelian group. Hence there are at least
\(\ell^{k_0+1}\) different reductions of the restricted
dual classes. Some restriction is outside the span of
the \(k_0\) earlier maps modulo \(\ell\).

The restricted cocycle on \(U\) is a homomorphism \(\psi\)
to \(\mathbb Z/\ell^\nu\mathbb Z\) in the chosen root-of-unity
coordinates. Conjugation acts through
\(\chi_{\rm cyc}\kappa^e\). This gives oddness under \(c\)
and exponent \(e+1\) under conjugation by \(\gamma_\ell\).
The global dual class is unramified at every finite place
outside \(S_0\). At a real place its local condition is unrestricted,
the orthogonal complement of the zero primal condition.
The exploratory set \(P\) can be discarded.

For the reciprocity law, let \(z\) be an even cocycle over
\(E\) allowed to ramify at the later split pairs and above
\(\ell\), vanishing on \(D_v\cap\ker\chi_{\rm cyc}\) for every
\(v\mid\ell\), with earlier tame and type-\(e\) Frobenius tests
zero. Extend it to \(G_{F_0}\) by
\(\widetilde z(gc)=z(g)\) for \(g\in G_E\).
This is a cocycle by evenness and \(\kappa(c)=1\).
Use \(2\widetilde z\) as the primal cocycle. At every
real conjugation \(2\widetilde z\) vanishes by squaring.

For \(g\in D_v\cap\ker\kappa\) at \(v\mid\ell\),
\(2\widetilde z(g)=\widetilde z(g^2)=0\), since
\(g^2\in D_v\cap G_E\) has cyclotomic value one.
At an earlier prime of type \(e\), unramifiedness and the
zero Frobenius value give zero on its decomposition group.
Thus the terms at \(S_0\) and infinity vanish by the orthogonal
local conditions. At earlier primes of other types and all
remaining finite primes outside the new pairs, both classes
are unramified. Their cup inflates from the procyclic residue
Galois group of \(\ell\)-cohomological dimension one, so its
local term is zero. At each new prime of \(F_0\) split in \(E\),
take \(\tau_i^+\) with
\(\tau_i^+(\varpi_i^{1/d_\nu})/\varpi_i^{1/d_\nu}=\zeta_{d_\nu}\)
for a uniformizer of \(F_{0,v_i}\), using the same primitive root.
Represent the local class of \(z\) by a Kummer element
\(b_i\in F_{0,v_i}^\times\) in these coordinates. The
\(d_\nu\)-th roots of units lie in the maximal unramified
extension, since \(v_i\nmid\ell\). Thus
\(z(\tau_i^+)=v_i(b_i)\bmod d_\nu\).
Local reciprocity pairs the class with the unramified
dual value as \(\pm\psi(\phi_i)v_i(b_i)\), with the same sign
at every prime from the geometric Frobenius and cup-product
conventions. This gives the term \(\pm2\psi(\phi_i)z(\tau_i^+)\)
for \(2\widetilde z\), and proves the reciprocity relation. Inertia acts trivially
on the coefficients, so each arbitrarily prescribed tame generator
has value an \(A_\nu\)-multiple of \(z(\tau_i^+)\), with
\(A_\nu=\mathbf Z/\ell^\nu\mathbf Z\).
For \(O/\ell^\nu\)-valued cocycles apply the relation
after all linear functionals \(O/\ell^\nu\to
\mathbb Z/\ell^\nu\mathbb Z\), which separate values since
\(O/\ell^\nu\) is free over \(\mathbb Z/\ell^\nu\mathbb Z\). This proves the
fresh-map assertion at every positive depth.

\subsection{The ordered prime selection from the Kummer relations}
\label{sec:triangle-selection}

\sectionopening{In this subsection, we combine the cyclotomic-kernel evaluations of
Section~\ref{sec:triangle-unit-input} with the Kummer maps of
Section~\ref{sec:fresh-kummer-maps}. The order of the choices ensures
that each newly selected prime satisfies the conditions imposed by the
earlier ones.}

Keep \(g_{\rm tr}\), \(C=\ell^{Y_{\rm tr}h_\ell}\), and a finite
set of evaluation elements in the cyclotomic kernel \(P_0\) as in
Section~\ref{sec:triangle-unit-input}. For \(\ell\le Y_{\rm tr}^2\),
the finite detection bound of Section~\ref{sec:odd-cocycle-budget} gives
\[
 m_0=N_c\ell^2\le Y_{\rm tr}^5,\qquad
 |P_0|\le m_0^4\le Y_{\rm tr}^{20}.
\]
For the actual rank increment \(t'\le t\), order the
pairs \((e,s)\), \(0<e<t'\), \(s\in P_0\), lexicographically
with the exponent first. The uniform bound uses \(t\).
Fix this list before choosing the depth. For each \(\nu>0\), put
\(A_\nu=\mathbb Z/\ell^\nu\mathbb Z\) and
\[
 G_\nu=G_{E(\mu_{\ell^\nu})}
      =\{g\in G_E:\chi_{\rm cyc}(g)\equiv1\pmod{\ell^\nu}\}.
\]

At each step Section~\ref{sec:fresh-kummer-maps}
gives a Kummer map \(\psi_j\) on \(G_\nu\) with values in
\(A_\nu\) for the earlier split pairs of type \(e\), and
a finite excluded set of primes of \(F_0\). The map is a continuous homomorphism, is odd
under complex conjugation, and conjugation by \(\gamma_\ell\)
multiplies its values by \(\xi_\ell^{e+1}\). Its reduction
modulo \(\ell\) is outside the span of the previous maps of
type \(e\).

The required reciprocity law is as follows.
Let an even cocycle \(z\) of type \(e\) vanish on the full
cyclotomic kernel in the decomposition groups above \(\ell\).
Suppose it is unramified at the previous pairs and outside
\(\ell\) and a finite family \(\{w_i,\bar w_i\}\) of further
\(\nu\)-split pairs with pairwise distinct restrictions to \(F_0\),
all outside the excluded set. If it vanishes
on the earlier Frobenius elements of type \(e\), then
\[
             \sum_i 2\psi_j(\phi_i)z(\tau_i^+)=0.
\]
Here \(\phi_i\) is geometric Frobenius and \(\tau_i^+\)
is the generator normalized over \(F_0\) in the chosen decomposition
group as in Section~\ref{sec:fresh-kummer-maps}.

The primes over \(\ell\) and the previously specified
local tests belong to the excluded set of the new map.
The set of usable Frobenius points must be abundant
near a translate of every prescribed point. At each depth
\(\nu\) and \(g\in G_E\) there must be
\(\delta\in\ker\chi_{\rm cyc}\) such that every open
neighborhood of \(\delta^{\ell^\nu}g\) contains usable
Frobenius elements of new places of prime norm outside any
given finite excluded set. At those places use the
determinant-unit and Hecke conditions of
Section~\ref{sec:triangle-unit-input}.

Under this Frobenius condition, consider the index \(j=(e,s)\).
The earlier maps of type \(e\) and the new \(\psi_j\)
are independent modulo \(\ell\). Their joint evaluation
on \(G_\nu\) has image all of the appropriate power
of \(A_\nu\): modulo \(\ell\) this follows by the absence
of an annihilating linear combination, and Nakayama's
lemma gives the assertion over \(A_\nu\).
Choose \(a\in G_\nu\) on which the new map reads one
and the earlier maps of type \(e\) read zero.
For each \(k\in\{0,\ldots,t'-1\}\setminus\{e\}\), replace
\(a\) successively by
\[
 (\gamma_\ell a\gamma_\ell^{-1})^{\xi_\ell^{-1}}
                                  a^{-\xi_\ell^k}.
\]
The exponents are integer representatives modulo \(\ell^\nu\).
Every map of type \(d\) is multiplied by
\(\xi_\ell^{d}-\xi_\ell^k\). Thus the resulting \(a_*\in G_\nu\)
has every earlier value zero and
\(\psi_j(a_*)=\Delta_e\bmod\ell^\nu\).

Apply the auxiliary-prime theorem at
\(g=s csc^{-1} a_*(c a_*c^{-1})^{-1}\).
All factors and \(\delta^{\ell^\nu}\) lie in \(G_\nu\).
An odd Kummer map reads \(2\psi(a_*)\) at
\(\delta^{\ell^\nu}g\), since the first two factors
cancel and \(\ell^\nu\) kills its value on \(\delta\).
An even cocycle of type \(e\), restricted to \(G_\nu\),
is additive and reads \(2z(s)\) there.
The finitely many Kummer evaluations persist on an open
neighborhood.

For the even evaluations take the finitely
many \(A_\nu\)-valued cocycles unramified outside \(\ell\)
using the finite evaluation set. Their values persist on one
open neighborhood. Applying the additive \(A_\nu\)-linear
duals of \(O/\ell^\nu\), which separate its elements,
gives the same equalities for the \(O/\ell^\nu\)-valued
cocycles, on that neighborhood for every coefficient ring.
Set \(h=\delta^{\ell^\nu}g\).
Intersect the open kernels in \(G_\nu\) of the earlier and current Kummer
maps and the restrictions of these finitely many cocycles. Each restriction
is additive on \(G_\nu\). Take an open normal subgroup of \(G_E\)
inside this intersection, shrinking it so that its coset at \(h\)
also preserves the residual determinant of \(1-T(h)\).
The Frobenius construction of Section~\ref{sec:triangle-unit-input}
now gives one prime and prolongation satisfying all these conditions.
Exclude \(\ell\), the rational primes under the original level,
the twist support, the exceptional places of the finitely many
original systems, the excluded sets of the Kummer maps chosen so
far, and every earlier pair. These are finitely many rational
primes at this step. Even when the same Frobenius class is
prescribed again, it still has positive density after these
exclusions. Use the conjugate prolongation at \(\bar w_j\). Writing
\(\phi'_j=(c\phi_jc^{-1})^{-1}\), oddness and evenness on \(G_\nu\) give
\[
 \psi_i(\phi'_j)=\psi_i(\phi_j),\qquad
 z(\phi'_j)=-z(\phi_j)
\]
for the Kummer maps and even cocycles just tested. The geometric
Frobenius at the conjugate prolongation is \(c\phi_jc^{-1}\).
The paired test uses its inverse.
Its prime norm implies it is split over \(F_0\), and
\(\ell^\nu\mid\mathrm Nw_j-1\).
It has
\[
 \psi_i(\phi_j)=0\ (i<j),\qquad
 \psi_j(\phi_j)=2\Delta_e,\qquad
 z(\phi_j)=2z(s)
\]
for the stated unramified even cocycles.

For an even cocycle ramified only on the new pairs and
above \(\ell\), first impose vanishing at all earlier tame
generators. Evenness gives vanishing at their conjugates,
so the cocycle is unramified at those earlier pairs.
Impose the earlier Frobenius tests of type \(e\).
The reciprocity law for \(\psi_j\) now applies to the
current and later pairs. The later coefficients vanish by
construction. The remaining term is
\(2\psi_j(\phi_j)z(\tau_j^+)=4\Delta_e z(\tau_j^+)\).
Write \(C=u_e\,4\Delta_e\) in \(\mathbb Z_\ell\). Multiplication
by \(u_e\) gives \(Cz(\tau_j^+)=0\) at every depth, including
those for which \(\Delta_e\equiv0\pmod{\ell^\nu}\).
The value at the prescribed \(\tau_j\) is a multiple of the
value at \(\tau_j^+\) in the cyclic pro-\(\ell\) tame quotient
modulo \(\ell^\nu\), so \(Cz(\tau_j)=0\) as well. This gives the ordered tame test and the reads of
the original evaluation map at the new Frobenius elements. The
total loss exponent counts \((t'-1)|P_0|\) ordered indices,
\(|P_0|\) evaluation points and one reference point. Thus
\[
  (t'|P_0|+1)Y_{\rm tr}h_\ell
   \le (t|P_0|+1)Y_{\rm tr}h_\ell\le Y_{\rm tr}^{22}h_\ell .
\]
When \(\ell>Y_{\rm tr}^2\) this holds because \(h_\ell=0\).
There is no bound on the norms of the new primes in
this argument.

\subsection{The paired tame and Frobenius tests}
\label{sec:triangle-character-cost}

\sectionopening{In this subsection, the successive primes of Section~\ref{sec:triangle-selection} give two
local evaluations for each relative cocycle. We combine the tame and
Frobenius evaluations with the original local conditions to bound the
length of their common kernel.}

At depth \(\nu\), let \(Z\) be an \(O\)-submodule of the
relative character cocycles for \((O/\ell^\nu)(-e)\) of
Section~\ref{sec:relative-crystalline-colength}. Put
\(\chi_e=\chi_{\rm cyc}^{-e}\). With \(\gamma_\pi\) the prescribed
scalar inertia element and \(g_{\rm tr}\) the reference element of
Section~\ref{sec:triangle-unit-input}, assume
\[
 z\in Z,\quad z(\gamma_\pi)=0
 \quad\Longrightarrow\quad
 (\chi_e(g)-1)z(g_{\rm tr})=0\quad(g\in G_E).
\]
The congruence modules used below satisfy this implication by the
normalizer calculation.
Suppose the old prescribed tame tests and their conjugate
inverses vanish and cover all original and auxiliary
ramification away from \(\ell\) except the new pairs
\(w_j,\bar w_j\). Outside these places the matrices
\(u_x(g)-1\) vanish for inertia \(g\), so evaluation
at \(g\) is a relative linear relation. At an old place
away from \(\ell\), the character cocycle on inertia is
a homomorphism to an \(\ell\)-primary finite group.
It factors through pro-\(\ell\) tame inertia. The topological tame generator and its conjugate inverse
therefore give vanishing throughout those groups.
Thus \(Z\) is unramified outside \(\ell\) and the new pairs.

Let \(\tau_j\) generate pro-\(\ell\) tame inertia at
\(w_j\), and let the chosen \(\phi_j\) fix
\(\mu_{\ell^\nu}\). Set \(\tau'_j=(c\tau_jc^{-1})^{-1}\)
and \(\phi'_j=(c\phi_jc^{-1})^{-1}\). Thus \(\phi_j\) is
geometric Frobenius and \(\phi'_j\) is arithmetic Frobenius at
the conjugate prolongation. For \(q_j=\mathrm Nw_j\),
\[
 \chi_{\rm cyc}^{-1}(\phi_j)=q_j,\qquad
 \chi_{\rm cyc}^{-1}(\phi'_j)=q_j^{-1},\qquad
 q_j\equiv q_j^{-1}\equiv1\pmod{\ell^\nu}.
\]
Let \(Z_0\subseteq Z\) require vanishing on the local
sets above \(\ell\) used for \(Z^-_{e,\nu}\) in
Section~\ref{sec:odd-cocycle-budget}. This concerns
inertia for \(e=0\), and the full cyclotomic kernel
in each decomposition group for \(e>0\).
The crystalline linear relations imply this vanishing.
For positive \(e\), the cyclotomic image is infinite
and the period \(t_{\rm cyc}^e\otimes v_{-e}\) generates the de Rham
module in \(\operatorname{Fil}^0\), where \(v_{-e}\) is a basis of the
Tate character and \(t_{\rm cyc}\) is the cyclotomic period. Bloch--Kato's formula
gives \(H^1_f(E_v,\operatorname{Frac}(O)(-e))=0\).
Thus an integral crystalline cocycle is a rational
coboundary and vanishes on that kernel. For \(e=0\) the formula
gives dimension one over \(\operatorname{Frac}(O)\)
\cite[Corollary~3.8.4, p.~355]{BK}.
Let \(f_v\) be the residue degree of \(E_v/\mathbf Q_\ell\).
The equation \(\varphi^{f_v}u-u=1\) has a solution in
\(W(\overline{\mathbf F}_\ell)\subset B_{\rm cris}\).
Solve successively modulo the powers of \(\ell\) using
surjectivity of \(x\mapsto x^{\ell^{f_v}}-x\) over
\(\overline{\mathbf F}_\ell\). Since arithmetic Frobenius acts on
this Witt ring by \(\varphi^{f_v}\), the unramified cocycle
sending Frobenius to \(1\) is the coboundary of \(u\) in
\(B_{\rm cris}\). Hence the finite condition is the unramified line.
Each asserted evaluation is
therefore an integral crystalline relation. It remains
a relation modulo \(\ell^\nu\), including for cocycles
that only obey the relations. Conjugation transfers
vanishing from the designated decomposition group to
its conjugates since the character at an element of the
local set is one. Consequently
\[
         \operatorname{length}_O(Z/Z_0)
                               \le c_{\rm crys}(O(-e),\nu).
\]

For \(z\in Z_0\), paired agreement at either
\((\tau_j,\tau'_j)\) or \((\phi_j,\phi'_j)\) is
equivalent to vanishing of \(z+z\circ c\) at the
first member. In fact the character acts trivially on
both pairs and hence inversion negates the value.
The sum \(z+z\circ c\) is an even cocycle with the same
ramification and local vanishing. Suppose the new
places have the ordered positive types and satisfy the conditions of
Section~\ref{sec:triangle-selection}, with the separation
scalar \(C=\ell^{Y_{\rm tr}h_\ell}\).
Its even-cocycle test at an index of type \(e\) says
that vanishing on earlier tame generators and earlier
Frobenius elements of type \(e\) forces the value at
the current tame generator to be killed by \(C\).
For the paired tests this gives agreement up to \(C\)
at \(\tau_j,\tau'_j\), using all earlier tame pairs
and only earlier Frobenius pairs of type \(e\).

Now let \(S\subseteq Z_0\) require the two tame values
at every new pair to vanish individually and the
Frobenius values to agree for type \(e\). Killing the
tame generators removes all new ramification as above.
For the reference test use the normalizing element for the covector congruences.
Write \(d(g)\) for simultaneous augmentation coordinates
in a finite generating list \(g_i\) for the tuples
\((u_x(g)-1)_x\). Thus
\[
            u_x(g)-1=\sum_i d(g)_i(u_x(g_i)-1)
                              \qquad\text{for every }x.
\]
For a group-algebra covector \(v\), its value is
\(\sum_i v_i(g_i-1)\). The congruence ideal contains the
common action kernel of the \(u_x\) and the left
multiples of the normalizer \(N_0\) by the augmentation
ideal.

Here the normalizer is
\[
 E_\pi=c_\pi(\gamma_\pi-1),\qquad
 N_0=E_\pi+(1-E_\pi)(g_{\rm tr}-1),
\]
where \(c_\pi=(y_\pi-1)^{-1}\) for the scalar
\(y_\pi\) of \(\gamma_\pi\). Its inverse is integral
since the reduction of the scalar differs from one.
The covector \(v_0=c_\pi d(\gamma_\pi)+(1-E_\pi)d(g_{\rm tr})\)
has value \(N_0\) modulo the common action kernel.
Hence \((g-1)v_0\) satisfies the congruences for every \(g\).
Evaluation at \(\chi_e\) kills \(E_\pi\), since
\(\gamma_\pi\) is in inertia off \(\ell\).
The character cocycle given by a functional annihilating
the vectors satisfying the congruences therefore satisfies
\[
    (\chi_e(g)-1)\bigl(c_\pi z(\gamma_\pi)
                                      +z(g_{\rm tr})\bigr)=0 .
\]
This gives the required implication for the coordinate quotients
in the Fitting estimate, as verified there.

For \(z\in S\) the value on \(\gamma_\pi\) is zero
by unramifiedness, giving
\((\chi_e(g_{\rm tr})-1)z(g_{\rm tr})=0\).
For \(e>0\), use the triangle's detecting points
besides \(g_{\rm tr}\), its prescription that twice the
even value at each point is an appropriate Frobenius
value of type \(e\), and its divisibility
\(4(\chi_{\rm cyc}(g_{\rm tr})^{-e}-1)^2\mid C\).
The paired-cocycle argument in
Section~\ref{sec:odd-cocycle-budget} gives
\(\operatorname{length}_O(S/S[C])\le
\operatorname{length}_O Z^-_{e,\nu}\).
For \(e=0\) there is no index of type zero, and
\(S\) is already everywhere unramified, hence a
submodule of \(Z^-_{0,\nu}\). Thus in both cases
\begin{equation}\label{eq:triangle-character-cost}
 \operatorname{length}_O(Z/Z_0)+\operatorname{length}_O(S/S[C])
       \le c_{\rm crys}(O(-e),\nu)+v_{\rm ord}(X_{\rm odd}).
\end{equation}
Together with the tests for unramified cocycles,
evaluations at both tame generators of each new pair
detect all of \(Z\): their vanishing removes the new
ramification, and the remaining tests give zero.
The reference test uses the original congruence ideal.
The prime triangle uses the scalar inertia, distinguished
lift and residual comparison of the same representation.

\section{Measured norm tests and the integral trace algebra}
\label{sec:norm-tests-and-trace-algebra}

\sectionopening{In this section, we normalize the two theta norm tests for the critical values and construct
the integral trace algebra on the raised levels. The norm tests are used in
Section~\ref{sec:common-critical-twist}, and the algebra is used in
Section~\ref{sec:hecke-trace-comparison}.}

\subsection{Norm tests from a definite line}
\label{sec:marked-norm-line}

\sectionopening{In this subsection, for \(z\ge2\), the theta lift from the definite line gives a measured
norm involving the critical value at \(z\). We compute the finite
factors and the archimedean normalization here, and establish
integrality for the same test in
Section~\ref{sec:marked-both-norm-jets}.}

Let \(p=(\nu,z)\) be a critical pair on the prescribed old level
with \(z\ge2\). Its tame residue characters at the old places
may be trivial. The infinity weights of \(\nu\) are
\(-c_v\) with \(c_v\ge z\). Use the standard definite
line and the standard definite group of rank \(a_0=2z\)
over the same \(E/F_0\). In this section \(S_M\)
denotes the chosen halves of the old split level places.
Let \(Q\) be a set of the auxiliary split places
allowed in the tests, disjoint from the old places
and their conjugates. Write
\[
 q_Q=\prod_{w\in Q}q_w,\qquad
 S_Q=(S_M\cup Q)\cup(\bar S_M\cup\bar Q),\qquad
 C_G=\prod_{w\in S_M}|\mathrm{GL}_{rt}(k_w)| .
\]
Let \(\tau\) be a finite-order character of the adelic norm-one
quotient, trivial at infinity, unramified away from \(Q\) and tame at
\(Q\). Its ramification set may be smaller than the chosen modulus.
For a norm-one character \(\chi\) write \(\chi^{\rm bc}(a)
=\chi(a/\bar a)\). We use the finite \(L\)-function of this
base change and omit the places in \(S_Q\).

Use the vacuum multiplied by \(s_{\rm vac}^{a_0}\).
At each old place restrict the column to the reduction
of the specified basis vector. At a place of \(Q\)
prescribe the distinguished coordinate to be one.
Use the same arrow scale \(\pi_M=\bar\alpha\) and
\(s_{\rm jet}=\gamma=p_Fa_d/\omega_0\) from \eqref{eq:upper-CM-jet-transport}. Use the jet of weight
\(\bigotimes_v\operatorname{Sym}^{c_v-z}(\mathbf C^{a_0})\)
at scale \(\pi_Ms_{\rm jet}\), with the fixed character
\(\nu\).

Denote the resulting kernel with trivial auxiliary target
character by \(K_Q(g,h)\). Here \(g\) is on the definite line
and \(h\) on the rank-\(a_0\) space, both with the original measures.
For line classes \(g_i\) take the old principal-level
classes times integral lifts of inverse unit residues
at \(Q\). Set \(u_i=\det(1,g_i)\). The values of
\[
                   \Phi_\tau(h)=
       \sum_i \tau(u_i)^{-1}K_Q(g_i,h)
\]
are the class lift of the line character, since summing
over the inverse unit residues gives the test on which
the distinguished coordinate is nonzero and weighted by
the tame character. The sum on the line has no stabilizer
denominators at the old principal level. Its quotient
stabilizer is trivial by reduction at an old odd split
prime in rank one.

On the target the level \(J_Q\) is the old vector fixer
and the inverse image of the parabolic preserving the complementary
hyperplane at each new place. For \(k\in J_Q\), write \(t_w(k)\in k_w^\times\)
for its action on this one-dimensional quotient and \(\psi_{w,\tau}\)
for the residue character of \(\nu\tau\). In the integral column,
the inverse pullback satisfies \((k^{-1}X)_1=t_w(k)^{-1}X_1\) modulo
the maximal ideal. Therefore
\[
 \Phi_\tau(hk)
   =\prod_{w\in Q}\psi_{w,\tau}(t_w(k))^{-1}\Phi_\tau(h).
\]
The old vector fixers preserve the prescribed column.
If the auxiliary target character \(a\) is used, the lift
is multiplied by \(a(\det h)\), so the displayed transformation
law acquires the factor \(a(\det k)\).
All factors are unitary, and the squared Weyl norm is invariant
under \(J_Q\) for either choice. The class pairing
is the sum of the transported Weyl pairings divided by the finite
class stabilizer orders. Every such order divides \(C_G\).
Indeed the stabilizer is finite, since the real group is definite and its rational integral points are discrete.
After conjugation into an old compact group at
an odd split place with \(E_w=\mathbf Q_{p_i}\), the
kernel of reduction is torsion-free. In that kernel
a nonidentity \(1+p_i^sA\), for \(s\ge1\) and \(A\)
nonzero modulo \(p_i\), has the same depth \(s\) after
raising to a power prime to \(p_i\) and depth \(s+1\)
after raising to \(p_i\), by the binomial expansion.
Reduction therefore injects the finite stabilizer
into the finite linear group. The rank is at most \(rt\).

\begin{proposition}[Measured norm of the line lift]
\leavevmode\label{prop:marked-line-norm}
For the line lift \(\Phi_\tau\) with the original class pairing
at \(J_Q\),
\begin{equation}
       \|\Phi_\tau\|_{J_Q}^2
           =c_p\,q_Q^{a_0-1}
                 L^{S_Q}(z,((\nu\tau)^{\rm bc})^{-1}),
       \qquad c_p>0 .
 \label{eq:marked-rank-one-norm}
\end{equation}
The scalar \(c_p\) is independent of \(Q\) and \(\tau\).
\end{proposition}
\begin{proof}
Put \(u=\upsilon|_{U(1)}\), \(b=(\nu\tau)u^z\).
Use the column test \(f\) after the sum over inverse unit
residues at \(Q\), and put \(f_g=\omega_b(g,1)f\), where
the source action includes \(b(\det g)^{-1}\).
Let \(\eta=(1,\eta_f)\), \(\eta_f\bar\eta_f=-1\), be
the opposite-form coordinate of
Section~\ref{sec:marked-graph-doubling}. In the two original
blocks the doubled test is
\[
     f_g(X)\,\overline{f_{g'}(\eta^{-1}Y)} .
\]
Here the phase of the opposite theta functional is chosen
to give the conjugate evaluation. Its modulus is one by
Section~\ref{sec:marked-rational-theta}.

If the auxiliary character on the common target in
the first factor is \(a\), its opposed character in
the rank-one negative factor is \(a^{-1}u\), by the
intrinsic determinant multiplier of that factor.
Their product is \(u\). Thus this is the target action
of Section~\ref{sec:marked-doubled-action} with Siegel
rank one, with both source translations already in the
test. The partial Fourier transform into the split
polarization uses the original dual additive measures.
The doubled line has a one-dimensional Lagrangian. In
\cite[Theorem~2]{YamanaSW} the Siegel index is therefore \(n=1\),
the Hermitian dimension is \(a_0=2z\), and the Witt index is zero.
Thus \(a_0>n\), \(s_0=(2z-1)/2\), and evaluation at the identity
of the double has coefficient one for target quotient probability Haar.

Every rational isotropic line of the doubled source
is a graph of an isometry \(\delta\), since both
projections are injective by definiteness. Evaluation
at zero in the transverse polarization evaluates the two
blocks at \((M\widehat\delta,\eta M)\) in the original
point coordinates of Section~\ref{sec:marked-graph-doubling}.
Thus the second argument returns to \(M\). Rational
invariance and unitarity give
\[
 \langle f_{g'},\omega_b(\delta,1)f_g\rangle
      =\langle f,\omega_b(g'^{-1}\delta g,1)f\rangle .
\]
The source character here is \(b(\det(g'^{-1}\delta g))^{-1}\),
since \(b(\det\delta)=1\).
Use the original target volume to undo the
probability normalization. The quotient formula
unfolds the rational graph sum to the ambient line
group, leaving one original source quotient volume.
The absolute local sum at a split place outside \(S_Q\)
is bounded by \(\sum_{j\in\mathbb Z}q_w^{-z|j|}\).
The product converges for \(z>1\), and the remaining
local factors are finite. Together with the compact
real groups this justifies the absolute unfolding
and factorization.

Here are the local quantities, with each standard
finite compact of mass one. At a split place outside \(S_Q\), a
line element with valuation pair \((j,-j)\) has unitary
Schr\"odinger coefficient of absolute value
\[
 q_w^{-a_0j/2}q_w^{-a_0\max(0,-j)}=q_w^{-a_0|j|/2}.
\]
The two factors are the square root of the Jacobian and the volume
of the overlap of the integral column boxes. The split Weil character factor,
after multiplication by the source character, is
\(A_w^{-j}\) for \(A_w=(\nu\tau)^{\rm bc}(\varpi_w)\):
indeed the determinant factor from the splitting is
\(\upsilon(\varpi_w)^{2zj}\) and the source factor
contains its inverse.

The sum of the overlaps is
\[
 \sum_{j\in\mathbb Z}A_w^{-j}q_w^{-z|j|}
 =
 \frac{1-q_w^{-2z}}{(1-A_wq_w^{-z})(1-A_w^{-1}q_w^{-z})}.
\]
At a nonsplit finite place \(v\) of \(F_0\), write \(w\) for
the place of \(E\) above it and \(q_v=|k_v|\). The norm-one
line group is its integral compact, the test is fixed, and its
local integral is one. Since \(E_w/F_{0,v}\) is unramified
quadratic, \(q_w=q_v^2\) and a uniformizer of \(F_{0,v}\) is also
one of \(E_w\). The base change of \(\nu\tau\) is trivial on this
uniformizer. Hence its inverse has Euler factor
\((1-q_w^{-z})^{-1}\), whose product with
\(1-q_v^{-2z}\) is one.
Let \(S_Q^+\) be the finite places of \(F_0\) below \(S_Q\),
and let \(I_v\) denote the local integral. Since \(a_0=2z\)
and \(S_M\cup Q\) contains one place above each \(v\in S_Q^+\),
\[
 \prod_{\substack{v\nmid\infty\\v\notin S_Q^+}} I_v
 =\frac{L^{S_Q}(z,((\nu\tau)^{\rm bc})^{-1})}
        {\zeta_{F_0}^{S_Q^+}(2z)}
 =\frac{L^{S_Q}(z,((\nu\tau)^{\rm bc})^{-1})}{\zeta_{F_0}(2z)}
       \prod_{w\in S_M\cup Q}(1-q_w^{-a_0})^{-1}.
\]

At a distinguished place a nonunit moves the primitive
column support off itself. The old column has
squared mass \(q_w^{-a_0}\), and among source
units its support is fixed exactly by reduction
one. On this subgroup every allowed tame character is trivial,
so the old matrix coefficient integral is \(q_w^{-a_0}/(q_w-1)\)
for each old residue character, including the trivial one. At a new place the sum over
the inverse unit residues replaces the point condition
by the character-weighted condition \(X_1\ne0\) on the
distinguished coordinate \(X_1\).

Every source unit
preserves this support. The residue multiplier cancels
the source character, so the coefficient integral is
the mass \((q_w-1)/q_w\). The same statement holds with
a trivial new-place character: the residue sum still
gives the same support with coefficient one. Thus in
that case the new quotient level and volume factor
are the same as for an exact-support place. The passage to the
target class pairing multiplies by \(q_w^{a_0}-1\)
at an old place and \((q_w^{a_0}-1)/(q_w-1)\)
at a new place. These are the coset counts for
the two target stabilizers. Turning the line Haar
lift into the old principal-level class sum multiplies
its norm square by the square of the old source index,
whose old local factor is \(q_w-1\).

At infinity
the line character cancels the rotation of the
selected highest monomial \(X_1^{c_v-z}\) and its
squared norm is \((c_v-z)!\) in the Fock pairing.
For a highest covector of norm one, Schur orthogonality on the compact
group gives the squared scalar coefficient integral as the squared Weyl
norm divided by \(\dim\bigl(\bigotimes_{v\mid\infty}
\operatorname{Sym}^{c_v-z}(\mathbb C^{a_0})\bigr)\). Thus passing to the Weyl norm multiplies by this dimension.
At an old place the two class-sum source indices, the source
coefficient integral, the target index and the removed Euler numerator give
\[
 (q_w-1)^2\frac{q_w^{-a_0}}{q_w-1}
           \frac{q_w^{a_0}-1}{1-q_w^{-a_0}}=q_w-1 .
\]
At a new place the corresponding factor is
\[
       \frac{q_w-1}{q_w}
       \frac{q_w^{a_0}-1}{q_w-1}
       (1-q_w^{-a_0})^{-1}=q_w^{a_0-1}.
\]
In particular, with \(V_j\) the volume of the standard definite
rank-\(j\) quotient for these Haar measures, put
\[
 t_v=c_v-z,\qquad
 d_\infty=\prod_{v\mid\infty}\binom{t_v+a_0-1}{a_0-1}.
\]
For the vacuum and jet scales specified above the constant in
\eqref{eq:marked-rank-one-norm} is
\[
 c_p=\frac{V_1V_{a_0}}{\zeta_{F_0}(a_0)}\,
  d_\infty\prod_{v\mid\infty}t_v!\,
  \prod_{w\in S_M}(q_w-1)\,
  \left|s_{\rm vac}^{a_0}
     (\pi_Ms_{\rm jet})^{\sum_{v\mid\infty}t_v}\right|^2 .
\]
All factors are positive and independent of \(Q,\tau\).
\end{proof}

For the integral coefficient test put
\[
   d_p=M^{2fa_0}(Md_K)^{\sum_v(c_v-z)},\qquad
   F_{p,Q}(\tau)=C_Gd_p^2 q_Q^{1-a_0}\|\Phi_\tau\|_{J_Q}^2 .
\]
Section~\ref{sec:marked-both-norm-jets} proves the unit
bound on \(d_pK_Q\) in the transported coordinates at
every embedding away from \(q_Q\), including the old level primes
and \(2\). To apply it in the pairing, let \(A\) be the tuple of
original coefficient embeddings of \(h\), \(G_{\rm W}\)
the integral Weyl Gram matrix, and \(S\) the weight
representation. The conjugate value of \(h^{-1}\)
is \(A^t\), since the standard definite form is unitary.
The pairing identity
\[
                 S(A^t)^tG_{\rm W}S(A^{-1})=G_{\rm W}
\]
follows on tensor powers of the coordinate pairing
and hence on the Weyl space.

For the Gram form
linear in its second entry write
\[
 \begin{split}
  F_{p,Q}(\tau)&=\sum_{i,j}C_{ij}(Q)\tau(u_j u_i^{-1}),\\
  C_{ij}(Q)&=\frac{C_Gd_p^2}{q_Q^{a_0-1}}
       \sum_{[h]}\frac{\operatorname{Gram}(K_Q(g_j,h),K_Q(g_i,h))}
                              {|\operatorname{Aut}(h,J_Q)|}.
 \end{split}
\]
Indeed \(\overline{\tau(u_j)^{-1}}\tau(u_i)^{-1}=\tau(u_j u_i^{-1})\).
Use the given embedding on the second vector and its composite with
complex conjugation on the first. The displayed Gram identity then
evaluates both cleared vectors in integral coordinates at the same
class point. The integer \(q_Q\) is a unit at each prime being tested. Thus the stabilizer bound
and the unit bound give integral \(C_{ij}(Q)\) at
every embedding away from \(q_Q\). Use the same representatives
for all \(\tau\), since the unitary level character cancels in the
squared norm of the complete lift.
Equation~\eqref{eq:marked-rank-one-norm} proves
\[
 \frac{F_{p,Q}(\tau)}{F_{p,\varnothing}(1)}
 =
 \frac{L^{S_Q}(z,((\nu\tau)^{\rm bc})^{-1})}
      {L^{S_\varnothing}(z,(\nu^{\rm bc})^{-1})}.
\]
The denominator is positive by absolute convergence
and pairing the two conjugate Euler factors at each
split place. Using \(\tau\) and \(\tau^{-1}\) gives the
two sign tests with the same clearing constant.

\subsection{The rank-two determinant norm}
\label{sec:marked-norm-determinant}

\sectionopening{In this subsection, when \(z=1\), we instead evaluate a determinant period on the rank-two
group. Its ratio to the untwisted period will enter the same integral
comparison as the norm of Section~\ref{sec:marked-norm-line}.}

Use the same finite bases and the sets \(S_M,Q,S_Q\), and
\(q_Q,C_G\) as in the preceding subsection. Let \(p=(\nu,1)\)
be a ramified critical pair, with infinity weights \(-c_v\),
where \(c_v\ge1\). Put \(k_v=c_v-1\) and
\[
 \alpha=\nu\tau\,\upsilon|_{U(1)},\qquad
 P(X)=\prod_{v\mid\infty}\det(X_v)^{k_v}.
\]
Use the original oscillator on the pair of standard definite
rank-two spaces. Write \(H=U(I_2)\), and \(K_2\) for the
stabilizer of its standard finite lattice. Give each finite
compact factor of \(K_2\) and the real compact group mass one.
Write \(v_2\) for the volume of \([H]=H(F_0)\backslash
H(\mathbb A_{F_0})\) for this Haar normalization.

The finite test is \(\chi_{\mathfrak m}(\det X)\) in the
original split Schr\"odinger coordinates at the modulus
\(\mathfrak m=(\pi_M)^2\prod_{w\in Q}\mathfrak p_w\).
The residue character
\(\chi_{\mathfrak m}\) is the finite unit character of
\(\alpha\) in these coordinates, extended by zero on
nonunits. At every other place the finite test is the
original vacuum. Denote by \(f_Q\) this finite test with the
real Fock vector \(P\), and by \(\theta_{f_Q}(g,h)\) its
kernel with the two characters \(\alpha,\alpha\).
The determinant period is
\begin{equation}\label{eq:determinant-marked-period}
 \Lambda_Q(\tau)
   =\sum_{h\in H(F_0)\backslash H(\mathbb A_f)/K_2}
       \frac{\theta_{f_Q}(1,h)}{|\operatorname{Aut}(h,K_2)|}.
\end{equation}

\indent{\it The local actions and the quotient integral.}
At a split place let \(X\) be the \(2\)-by-\(2\)
Schr\"odinger variable. Identify both groups there with
\(\mathrm{GL}_2(F_{0,v})\) by the same chosen member of the
split pair. The action is
\[
 (\omega(g,h)f)(X)=
 B(\det h)B(\det g)^{-1}
 |\det g|_v|\det h|_v^{-1}f(h^{-1}Xg),
 \qquad B=\alpha^{\rm bc}\upsilon^{-2}=(\nu\tau)^{\rm bc}.
\]
Here \(B\) and \(\upsilon\) are read on the chosen split
factor. The powers of the absolute values are the
square roots of the two matrix Jacobians. The original
split self-dual coordinate measure gives the integral
matrix box mass one. The splitting
factor on each side is \(\upsilon^{-2}\), and
\(\upsilon(t/\bar t)=\upsilon(t)^2\) on this split
coordinate.

On integral units \(B\) is precisely the
finite character. In particular the test is
fixed by both unit compacts: multiplying the argument
of \(\chi_{\mathfrak m}\circ\det\) on the left by
\(h^{-1}\) and on the right by \(g\) cancels the two
unit character factors. This holds also on the
noninvertible residue matrices, on which both sides
vanish. Away from the chosen places, the original
unit compacts fix the vacuum and the characters are
unramified. At infinity \(\alpha\) has weight \(-c_v+1=-k_v\). The bare Fock substitution on
\(P(X_v)=\det(X_v)^{k_v}\) contributes
\((\det h_v)^{k_v}(\det g_v)^{-k_v}\). The factors
\(\alpha_v(\det h_v)\alpha_v(\det g_v)^{-1}\) cancel these powers.
The oscillator vector and the two theta variables are
therefore invariant under their respective real compact
groups and under \(K_2\) at finite places.

The kernel is left rationally invariant in each adelic
variable. In the finite evaluation in
\eqref{eq:determinant-marked-period}, this includes the
real factor of a rational move: substitution in
\(P\) multiplies it by the inverse power of the rational
determinant, which cancels the real multiplier of that
move. The definite quotients have compact sets of
representatives. The finite translates of the fixed
test over each such set are a finite family, by its
open stabilizer. On the real compact groups its
Schwartz seminorms are uniformly bounded. Accordingly
the rational theta sum converges absolutely and
uniformly on these sets of representatives, and the
kernel is continuous and bounded. Haar disintegration of a class of
representative \(h\) gives the mass
\(|\operatorname{Aut}(h,K_2)|^{-1}\), since the real
compact and \(K_2\) have mass one. Thus
\(\Lambda_Q(\tau)=F_Q(1)\), where
\[
 F_Q(g)=\int_{[H]}\theta_{f_Q}(g,h)\,dh .
\]

\indent{\it Constancy on the original group.}
Fix a chosen split place and regard the integral over
\(h\) as a linear functional in the local test, with
all other factors and the source point fixed. On
\(C_c^\infty(\mathrm{GL}_2(F_{0,v}))\) this functional
satisfies
\[
 T(f(l^{-1}X))=\lambda(\det l)\,T(f),
 \qquad \lambda(t)=B(t)^{-1}|t|_v .
\]
Indeed right translation by \(l\) of the target point
preserves the quotient Haar, and its oscillator action
has scalar \(B(\det l)|\det l|_v^{-1}\).
For a small compact open \(L\) on which \(\lambda\circ
\det\) is trivial, replace each test by its finite
expansion in indicators of left cosets \(lL\).

This
shows that on every such test the functional is a
constant times \(\int f(X)\lambda(\det X)\,d^\times X\).
More explicitly, for any finite set of tests one may
choose one \(L\) so that all their compact supports
are disjoint finite unions of its cosets. The constant
for these tests is \(T(\mathbf 1_L)/\operatorname{vol}(L)\),
by the displayed covariance. Both the test
and its conjugate by any fixed matrix satisfy this
condition since their supports are compact subsets
of the open invertible locus. Simultaneous translation
of source and target by the same matrix sends the
test to its conjugate and has scalar one.

Conjugation
preserves the weighted integral. Absorb the target
translation into its Haar integral to obtain right
invariance of \(F_Q\) by the full source split group
at this place. The argument applies at every place
of \(S_M\cup Q\), including a place of a larger
modulus at which \(\tau\) is unramified.

The group \(\operatorname{SU}(I_2)\) is simply connected
and absolutely almost simple. At an old split place
it is \(\operatorname{SL}_2\), which is noncompact.
Strong approximation \cite[Theorem~7.12]{PlatonovRapinchuk} away from
this place implies
\(F_Q(sg)=F_Q(g)\) for every \(s\in
\operatorname{SU}(I_2)(\mathbb A_{F_0})\). Indeed,
approximate \(s\) by rational \(s_0\) away from the
place so that \(g^{-1}s_0^{-1}sg\) is integral at all
other finite places and tends to the identity at the
real places. Apply rational left invariance, right
invariance by the full group at the omitted place,
and continuity. The determinant map is onto the
adelic norm-one group, by \(u\mapsto\operatorname{diag}(u,1)\).
The lift therefore factors through the determinant.

At a split place outside \(S_Q\) both tests are the same vacuum. The minuscule source and target actions on
it agree when the target acts at \(w\) and the source
acts at \(\bar w\). Namely choose integral elementary
matrices \(\gamma_U\) with determinant valuation one
and reduction kernel \(U\), for the residue lines
\(U\subset k_w^2\). Use \(\gamma_U^{\mathsf T}\) on
the target coordinate at \(w\), and \(\gamma_U^{-1}\)
on the source coordinate at \(w\), the latter being
the positive elementary class at \(\bar w\).

Take
the corresponding right coset representatives. The
target and source translates are respectively
\[
       q A\,\mathbf 1_{\mathrm{int}}(\gamma_U^{-\mathsf T}X),
       \qquad
       q A\,\mathbf 1_{\mathrm{int}}(X\gamma_U^{-1}),
       \qquad A=B(\varpi_w).
\]
They require integral \(X\), and then require \(U\)
to lie respectively in the kernel of \(\bar X^{\mathsf T}\)
and \(\bar X\). Indeed the image of
\(\varpi_w\gamma_U^{-1}\) on reduction is \(U\).
The two sums are therefore equal because these
kernels have the same dimension.

The group on the integrated variable can be moved
through the integral, including the inverse transpose
if that coordinate is used. Thus the corresponding
source Hecke operator acts on \(F_Q\) by \(q+1\).
On a determinant function it acts by \(q+1\) times
translation by the determinant of the diagonal with
one uniformizer (or its inverse). This proves
invariance by that translation and hence by the whole
split group. At a nonsplit finite place every norm-one
element is integral and \(\operatorname{diag}(u,1)\)
belongs to the standard compact. At infinity the
real compact invariance was verified above. This proves
that \(F_Q\) is constant, so
\[
             \|F_Q\|_{[H]}^2=v_2|\Lambda_Q(\tau)|^2 .
\]

\subsubsection{Doubling and the local integrals}
\label{sec:rank-two-doubling-local-integrals}

\sectionopening{In this subsection, we apply the doubled action of Section~\ref{sec:marked-doubled-action} to the determinant norm.}

For this equal-rank calculation put
\(u=\upsilon|_{U(1)}\) and write \(\theta_c\) for the
kernel at the same vector \(f_Q\) with pair characters
\((u,u)\). With \(\chi=(\nu\tau)\circ\det\) one has
\[
       \theta_{f_Q}(a,h)=\chi(a)^{-1}\chi(h)\theta_c(a,h).
\]
Here \(h\) is originally in the target slot. To double that slot as a
source, take the rational adjoint \(X\mapsto X^*\) from the original
Hom space to the Hom space with reversed factors and opposite second
form. The opposite sign preserves the Heisenberg pairing. In the
opposite basis of Section~\ref{sec:marked-equal-rank-vanishing} its
unitary intertwiner \(J\) takes the two characters \((u,u)\) to
\((u^{-1}u^2,u^{-1}u^2)=(u,u)\). The rational theta functional
changes by a unit scalar, since the adjoint identifies the rational
lattices and the normalized vacuum translates have mean square one.
The same \(J\) acts on both factors. Thus its theta phase cancels,
and for \(v_h=\omega_c(1,h)f_Q\) the identity-graph pairing is
\(\langle Jv_h,Jv_{h'}\rangle=\langle v_h,v_{h'}\rangle\).
The passage to the opposite form leaves the rational unitary group,
its integral compact and its quotient Haar unchanged.

On negating a rank-two factor, conjugation inverts
its character \(b\), and the intrinsic Weil multiplier
on the conjugate Hom action contributes \(u^2\) on
each determinant in the respective numerator or
denominator. Hence the opposed character is \(b^{-1}u^2\).
Both centered characters \(u\) are unchanged.
On doubling the integrated group, the pair characters
are \(u\) for the doubled source and \(u^2\) for the
common rank-two target, in the convention that the
action has the target determinant in the numerator.

Since \(\upsilon(\bar y)=\upsilon(y)^{-1}\) for every idele \(y\),
we have \(u^{\mathrm{bc}}(y)=\upsilon(y/\bar y)=\upsilon(y)^2\).
Thus both splitting characters are trivial by the identities
\[
 u^{\mathrm{bc}}\upsilon^{-2}=1,\qquad
 (u^2)^{\mathrm{bc}}\upsilon^{-4}=1 .
\]
In the doubled model the source-character factor
contributes \((u^{\mathrm{bc}})^{-1}=\upsilon^{-2}\)
to the untwisted Siegel Levi factor of
Section~\ref{sec:marked-doubled-action}, so the centered
section has the trivial splitting character. The
doubled theta kernel on the pullback is the product of
the two centered theta kernels with the first factor
conjugated. Both inputs use the above \(f_Q\) before
conjugation. On the pullback the two real factors are
\(\overline P\) and \(P\) in the corresponding Fock models.
In common graph coordinates the finite tensor is
\(\overline{f_{Q,f}}\otimes f_{Q,f}\), where \(f_{Q,f}\) is the
finite determinant test at
\(\mathfrak m=(\pi_M)^2\prod_{w\in Q}\mathfrak p_w\).
At the diagonal this gives \(|f_{Q,f}(X)|^2\), including when
the determinant character is nonquadratic.
In taking the partial Fourier transform, the passage from
midpoint and difference to the two arguments is
\[
       (Z,U)\longmapsto (Z+U/2,Z-U/2).
\]
Its determinant in each scalar coordinate is \(-1\). At \(v\mid2\)
the place splits in \(E/F_0\) by Section~\ref{sec:marked-rational-theta}.
For the two vacuum indicators on \(L=M_2(\mathcal O_v)\), the
conditions in these coordinates are
\[
 U\in L,\qquad Z\in U/2+L.
\]
Thus Fourier transformation in \(Z\) with negative sign gives
\(\mathbf1_L(U)\mathbf1_{L^\perp}(\Xi)
  \psi_v(-\langle U/2,\Xi\rangle)\), since \(L\) has mass one.
In particular the phase is one at the graph evaluation at zero.
More generally, after the original translations, evaluation at
\(U=\Xi=0\) integrates \(\overline{v_h(Z)}v_{h'}(Z)\) in the
original measure. This explains the graph pairing also when
\(2\) is not a local unit.

At a real place \(v\), use the oscillator complex structure given by
the two definite summands of the doubled source and the definite
target. The corresponding maximal compact of \(U(2,2)\) is
\(U(2)\times U(2)\). In the opposite and original Fock factors,
respectively, \(\overline{P_v}\otimes P_v\) is a polynomial of total
degree \(4k_v\) in creation operators applied to the tensor vacuum.
The partial Fourier transform \(\mathcal F\) to the mixed realization
intertwines the Weyl displacements. Differentiating gives
\(\mathcal F a^*(u)\mathcal F^{-1}\) as the creation operator for the
transported complex structure, and carries the tensor vacuum to its
vacuum line. Choose a real isometry taking the ordered positive and
negative subspaces in Yamana's doubled model to these two subspaces.
Write it as \(p_\infty k\) by Iwasawa decomposition, with
\(p_\infty\) in the Siegel parabolic and \(k\) in the standard
maximal compact \(K_\infty\). Since \(k\) preserves each subspace,
\(p_\infty\) effects the same change. The target is definite.
Thus \(\Phi_0=\omega_c(p_\infty^{-1})\Phi\) is polynomial in
Yamana's Fock model, where \(\Phi\) is the transformed tensor and
\(\omega_c\) denotes the centered doubled action.
With \(K'_\infty=p_\infty K_\infty p_\infty^{-1}\) and the same
finite compacts, Section~\ref{sec:marked-doubled-action} gives
\[
 f_{\Phi,K'}^{(s)}(g)
   =|a_K(p_\infty)|^{-s}f_{\Phi_0,K}^{(s)}(gp_\infty),
 \qquad
 \Theta(g,h;\Phi)=\Theta(gp_\infty,h;\Phi_0).
\]
Both rank parameters in Yamana are \(2\), so the height factor is one at \(s_0=0\).
Both evaluations use the same centered Weil action and its trivial
Siegel splitting character. The graph variables evaluate the section
of the original transformed tensor \cite[pp.~236--237]{YamanaSW}.

Write \(E_c\) for the Eisenstein series of this section
and \(\iota(h,h')\) for the graph embedding. The norm
integrand carries \(\chi(h)^{-1}\chi(h')\) outside the
Eisenstein series, because
\(\overline{\theta_{f_Q}(a,h)}\theta_{f_Q}(a,h')\)
has that factor. The Hermitian and Siegel ranks are
both two and the Hermitian space is definite. Thus
\cite[Theorem~2, p.~237]{YamanaSW} at \(s=0\), with the probability
target measure of \cite[\S1.3, p.~1007, and Proposition~5.8(i), p.~1015]{YamanaSingular}, gives
\[
 \|F_Q\|_{[H]}^2
    =\frac{v_2}{2}
        \int_{[H]\times[H]}
          \chi(h)^{-1}\chi(h')E_c(\iota(h,h'),0)
                                      \,dh\,dh' .
\]
Indeed the squared norm integrates the product of the two theta
kernels in the common variable \(a\). Yamana's identity reads
\(E_c=2\int_{[H]}\theta_c^{\rm dbl}\,da/v_2\).
Thus integration in \(a\) contributes \(v_2/2\). The two remaining
quotient integrations use the original measures.

All rational maximal isotropic subspaces are graphs of
rational isometries by definiteness of the two factors.
The graph zero term is the centered coefficient
\(\langle f_Q,\omega_c(1,g)f_Q\rangle\) at \(s=0\):
the self-dual partial Fourier transform at zero gives
the inner product on the common polarization.
For \(\Re s>1\), extend the coefficient by the
parabolic height to the power \(s\), using the compact
stabilizing the original doubled self-dual lattice.

Its Iwasawa decomposition is the integral Witt
completion of an isotropic direct summand. The change
to the classical skew coordinates in
Section~\ref{sec:marked-doubled-action} can be nonintegral,
so the compact is transported with the lattice.
The choice of extension leaves the value at zero
unchanged: two holomorphic sections agreeing at zero
differ by \(s\) times a holomorphic section, and in
rank two with trivial inducing character
\cite[Theorem~1(2), p.~236]{YamanaSW} gives no pole at zero
for the series of the latter section. The diagonal copy of \(H\)
stabilizes the identity graph. Its action on this Lagrangian has
determinant of absolute value one over \(E\), since
\(\det h\,\overline{\det h}=1\). The centered inducing character is
trivial, so left diagonal translation preserves the flat section
for every \(s\). Thus on \(\Re s>1\) the rational graph sum unfolds
with relative variable \(g=h^{-1}\gamma h'\), \(\gamma\in H(F_0)\).

The external factor equals \(\chi(g)\) since
\(\chi(\gamma)=1\). The relative variable has ambient Haar, and
the remaining quotient variable has volume \(v_2\). Writing \(Z(s)\)
for the unfolded integral including \(\chi(g)\), continuation and
the constancy already proved give
\[
 v_2|\Lambda_Q(\tau)|^2=\|F_Q\|_{[H]}^2
       =\frac{v_2^2}{2}Z(0).
\]
For evaluation under the integrals, choose finite
open compacts in the two factors fixing the pullback
of the flat section and the external characters.
The definite finite class sets at these levels are
finite. At a fixed pair of representatives the
real restriction is in a fixed finite sum of compact
matrix coefficient spaces. Indeed the embedded real
groups lie in the doubled maximal compact and the
flat height is right invariant under that compact.

In the Fock model the polynomial tensor lies in the sum of its finitely
many homogeneous degrees. The doubled real maximal compact preserves
each of these finite-dimensional spaces; the flat extension has the same
compact restriction. Thus the section is finite under its right action.
Finitely many point evaluations recover the coefficients
in this fixed space, and every such evaluation is
holomorphic at zero by \cite[Theorem~1(2), p.~236]{YamanaSW}.
Integration against the external characters on
the compact quotients is a fixed linear functional
on the finite sum of these spaces. It therefore
commutes with evaluation at zero.

Here are the local integrals of this section. Invert
the ambient variable, preserving its unimodular Haar.
The height is unchanged by inversion since the unitary
matrix and its inverse have equal absolute lattice
intersection index. At the split places the centered operator is
\(\omega_c(1,g^{-1})\) and acts by \(X\mapsto gX\)
with its intrinsic half-density. The external
character becomes \(\chi(g)^{-1}\). In particular
on elementary divisors \(j,k\) it is \(A^{-j-k}\),
where \(A=(\nu\tau)^{\mathrm{bc}}(\varpi_w)\).
In fact
\(\chi(g)^{-1}\omega_c(1,g^{-1})=\omega_{\alpha,\alpha}(1,g^{-1})\).
In particular, if \(g\in K_{2,v}\) at a place of the
determinant test, this coefficient is the squared norm of that test. At a finite place the doubled compact preserves the doubled
self-dual lattice. A graph intersects it in
\[
 \{(gx,x):x\in\Lambda\cap g^{-1}\Lambda\},
 \qquad \Lambda=\mathcal O_{E,v}^2.
\]
The inverse index of \(\Lambda\cap g^{-1}\Lambda\)
in \(\Lambda\) is the absolute determinant of the
graph coordinate in the Levi, hence the parabolic
height.

This follows by taking an integral basis of
the graph intersection and completing it to a
hyperbolic basis of the doubled lattice. A primitive
isotropic summand in the split or unramified
self-dual lattice admits such a completion: choose
vectors paired dually with its basis and subtract
multiples of the isotropic basis to cancel their Gram
matrix. In the unramified case the diagonal equations
use surjectivity of the integral trace, also at
residue characteristic two.

At a split pair outside \(S_Q\) with norm \(q\) and
\(A=B(\varpi_w)\), the Cartan representative with
elementary divisors \(j\ge k\) has height
\(q^{-|j|-|k|}\), using both members of the pair.
The overlap of the two integral matrix boxes, with
the square-root Jacobian in the local action, gives
the coefficient
\[
              A^{-j-k}q^{-|j|-|k|}.
\]
For example in a coordinate with valuation \(j\)
the matrix Jacobian factor is \(q^{-j}\),
and the overlap condition in the two entries of that row
has probability \(q^{-2\max(0,-j)}\).
The double coset has Haar mass \(1\) if \(j=k\),
and mass \(q^{j-k}(1+q^{-1})\) otherwise, by counting
free lines over \(\mathcal O_v/\varpi_v^{j-k}\).
With \(x=q^{-1-s}\), \(t=A^{-1}x\), \(u=Ax\), summing the weighted Cartan double cosets gives
\[
 \begin{aligned}
 &\frac1{(1-t)(1-qt)}+\frac1{(1-u)(1-qu)}-1
       +\frac{(1+q^{-1})q^2tu}{(1-qt)(1-qu)}\\
 &\hspace{12mm}=
 \frac{(1-x^2)(1-qx^2)}
      {(1-Ax)(1-A^{-1}x)(1-qAx)(1-qA^{-1}x)}.
 \end{aligned}
\]
For \(j\ge k\ge0\), set \(j=k+d\). The local sum on this cone is
\[
  \sum_{k\ge0}t^{2k}\Bigl(1+(1+q^{-1})\sum_{d\ge1}(qt)^d\Bigr)
   =\frac{1}{(1-t)(1-qt)}.
\]
The cone \(0\ge j\ge k\) gives the same expression in \(u\), with
the origin counted a second time. The sector \(j>0>k\) gives the
last term.

At an unramified nonsplit place let \(q\) be the norm
in \(F_0\), so the norm in \(E\) is \(q^2\).
An integral hyperbolic basis of the standard local
plane exists because the unit norm is onto. It may be
constructed by first taking a primitive isotropic
vector and then making its integral dual vector
isotropic using the trace. The special unitary group in the hyperbolic basis is
\(\mathrm{SL}_2(F_{0,v})\) under conjugation by a suitable
diagonal matrix. Indeed determinant one forces the
diagonal entries to be fixed by conjugation and the
off-diagonal ones to be skew. Choose a skew integral
unit \(\delta\), for example \(z-\bar z\) with the residue
of \(z\) outside the fixed residue field. Conjugating by
\(\operatorname{diag}(\delta,1)\) gives the asserted
isomorphism and identifies the integral compact with
\(\mathrm{SL}_2(\mathcal O_v)\), including above two.

Let \(K_v\) be the unitary integral compact, let
\(K_v^1=K_v\cap\operatorname{SU}(2)\), and put
\(a_j=\operatorname{diag}(\varpi^j,\varpi^{-j})\), \(j\ge0\).
The determinant of \(K_v\cap a_j K_v a_j^{-1}\)
is all norm-one units. In fact the integral isometries
\(\operatorname{diag}(z,\bar z^{-1})\) with
\(z\in\mathcal O_{E_w}^{\times}\) commute with \(a_j\)
and have determinant \(z/\bar z\). A Hilbert 90
preimage of a norm-one unit can be made a unit by
dividing by a power of the fixed-field uniformizer.
In particular every unitary double coset has special
unitary Cartan representative \(a_j\) and
\[
 [K_v:K_v\cap a_jK_va_j^{-1}]
  =[K_v^1:K_v^1\cap a_jK_v^1a_j^{-1}]
  =q^{2j}(1+q^{-1})\quad(j>0).
\]
The last equality counts free lines modulo
\(\varpi^{2j}\). The identity coset has mass one.

The graph height at \(a_j\) is \(q^{-2j}\) by the
intersection lattice formula. In the split-plane
polarization the centered matrix coefficient is the
overlap of the integral box in \(E_w^2\) with its
dilation by the rank-one Levi. In the convention
\(X\mapsto\varpi^j X\), the old box is the smaller
box. The intrinsic half-density is
\(|\varpi^j|_{E_w}^{2/2}=q^{-2j}\) and the centered
splitting character is trivial. The external character
is one on \(a_j\). Thus the vacuum coefficient is
\(q^{-2j}\) in this same model.
The local sum is
therefore
\[
     1+(1+q^{-1})\sum_{j\ge1}q^{-2sj}
              =\frac{1+q^{-1-2s}}{1-q^{-2s}} .
\]

At a chosen split place, overlap of the two
determinant tests requires \(X,gX\in
\mathrm{GL}_2(\mathcal O_v)\) and thus \(g\in
\mathrm{GL}_2(\mathcal O_v)\). Conversely for every such
\(g\) the finite determinant phase cancels by the local
action checked above. The height is one, the ambient
compact has mass one, and the integral of the squared
test is exactly
\[
 \frac{|\mathrm{GL}_2(k_w)|}{q_w^4}
          =(1-q_w^{-1})(1-q_w^{-2}).
\]
In the original point tests, put \(e=1\) or \(2\),
\(R_e=\mathcal O_v/\varpi_v^e\) and
\(\delta_y(X)=\mathbf1_{y+\varpi_v^e M_2(\mathcal O_v)}(X)\).
For \(g\in K_{2,v}\) and \(y,y'\in\mathrm{GL}_2(R_e)\),
\[
 \int \delta_y(X)\delta_{y'}(gX)\,dX
       =q^{-4e}\,\mathbf1_{\{y'=\bar g y\}} .
\]
Here \(q=q_w\) and \(\bar g\) denotes reduction modulo \(\varpi_v^e\).
Weighting these two tests by the conjugate and ordinary determinant
characters gives \(B(\det g)\) on every surviving pair. The external
factor is its inverse. Thus the mass is exactly
\(q^{-4e}|\mathrm{GL}_2(R_e)|=q^{-4}|\mathrm{GL}_2(k_w)|\),
also at an inserted split place with trivial unit character. At infinity the parabolic
height is one on the two compact groups, and the
local integral is \(\|\det(X_v)^{k_v}\|^2\).
In the normalization of Bargmann's monomial basis
\cite[(1.6)]{BargmannTransform}, polar coordinates give
\[
 \frac1\pi\int_{\mathbb C}z^a\bar z^b e^{-|z|^2}\,dA(z)
      =\delta_{ab}\,a! \qquad(a,b\in\mathbb Z_{\ge0}).
\]
Applying this identity in the four matrix entries gives
\[
 \|\det(X)^{k}\|^2
  =\sum_{j=0}^k{\binom{k}{j}}^2(j!)^2((k-j)!)^2
  =k!(k+1)! .
\]

Let \(\beta=(\nu\tau)^{\rm bc}\) and let \(\xi\) be the
quadratic character of \(E/F_0\). The local formulas
give the absolute Euler product for \(\Re s>1\).
For example the nonidentity contributions at
places away from \(S_Q\) are \(O(q^{-\Re s})\) uniformly
on closed half-planes \(\Re s\ge1+\delta\).
The remaining finite local factors are constant in \(s\). The product away from \(S_Q\) is
\[
 \frac{L^{S_Q}(s,\beta^{-1})L^{S_Q}(s+1,\beta^{-1})}
 {L_{F_0}^{S_Q}(2s+1,\xi)\,\zeta_{F_0}^{S_Q}(2s+2)} .
\]
In a superscript on an \(F_0\)-function, omit the
places below \(S_Q\). At a spherical split place the parameter
in \cite[\S11.6, pp.~823--824]{GanQiuTakeda} is our normalized induction
parameter \(s\). The standard factor there is evaluated at
\(s+\tfrac12\): the two Satake exponents of a determinant
character of \(\mathrm{GL}_2\) are \(\pm\tfrac12\), giving the
shifts \(s,s+1\) in the numerator. The two places above \(v\)
supply the factors with \(A\) and \(A^{-1}\).
The split denominator is \(\zeta_v(2s+1)\zeta_v(2s+2)\),
as in the Cartan sum above.

At an inert \(v\), put \(q=\mathrm{N}v\).
The character \(\beta\) is unramified and trivial on the fixed-field
uniformizer, the residue norm upstairs is \(q^2\), and
\(\xi_v(\varpi_v)=-1\). Hence the local Euler quotient is
\[
 \frac{(1+q^{-2s-1})(1-q^{-2s-2})}
      {(1-q^{-2s})(1-q^{-2s-2})}
   =\frac{1+q^{-2s-1}}{1-q^{-2s}},
 \qquad \Re s>1 .
\]
The individual nonsplit integrals
\((1+q^{-1-2s})/(1-q^{-2s})\) have poles at zero. The Euler product
was taken on \(\Re s>1\). Meromorphic continuation
of the Eisenstein series and of the Hecke functions
extends the equality of their integrals from
\(\Re s>1\).

The continued Eisenstein pullback may
be evaluated under the compact quotient integrals
at the holomorphic Siegel--Weil point. At zero the
local masses at the omitted places insert exactly the omitted reciprocal
denominator factors. The positive residues at one in
\(\zeta_E(s)=\zeta_{F_0}(s)L_{F_0}(s,\xi)\) show that
\(L_{F_0}(1,\xi)>0\). Thus the denominator
\(L_{F_0}(1,\xi)\zeta_{F_0}(2)\) is positive, giving
\begin{equation}\label{eq:determinant-marked-ratio-product}
 |\Lambda_Q(\tau)|^2
    =D_p L^{S_Q}(0,\beta^{-1})L^{S_Q}(1,\beta^{-1}),
 \qquad D_p>0 ,
\end{equation}
where \(D_p\) depends only on the original field,
Haar and real test at \(p\). With the stated probability
real measures and unscaled Bargmann test it is
\(v_2\prod_{v\mid\infty} k_v!(k_v+1)!/
(2L_{F_0}(1,\xi)\zeta_{F_0}(2))\).
Using the fixed normalization of the theta vector
throughout multiplies this constant by a fixed
positive factor. The same normalization is used
with and without \(Q\).

For a twist of exact support \(Q\), the conductor
of \(\beta^{-1}\) contains both members of every
pair with exponent one. The old factors are
fixed, since \(p\) is ramified at each old place
and the auxiliary support is disjoint. The archimedean
gamma factors are \(\Gamma_{\mathbb C}(s+c_v)\) and
are finite and nonzero at zero and one. With \(\mathfrak f_\beta\) the conductor, set
\(D_\beta=|\operatorname{disc}(E)|\,
  \mathrm N_{E/\mathbb Q}\mathfrak f_\beta\).
The completed function is
\[
 D_\beta^{s/2}\prod_{v\mid\infty}\Gamma_{\mathbb C}(s+c_v)
                      L(s,\beta^{-1}).
\]
Tate's global functional equation
\cite[\S2.5, Theorem~4.4.1 and \S4.5]{TateThesis} exchanges this function
with the completion for \(\beta\) at \(1-s\), up to a factor of
modulus one. On the real axis the two finite \(L\)-values are complex
conjugates. Exact support gives
\(D_\beta/D_{\nu^{\rm bc}}=\prod_{w\in Q}q_w^2=q_Q^2\);
the archimedean factors are the same for every twist. Since
\(\Gamma_{\mathbf C}(1+c_v)/\Gamma_{\mathbf C}(c_v)=c_v/(2\pi)\),
the absolute values at the two Mellin points satisfy
\[
  |L(0,\beta^{-1})|=D'_p q_Q |L(1,\beta^{-1})|,\qquad
  D'_p=D_{\nu^{\rm bc}}^{1/2}\prod_{v\mid\infty}\frac{c_v}{2\pi}>0.
\]
In this calculation the omitted factors of the full
Hecke function at the exact ramified support are all
one.

For real \(s>1\), applying \(3+4\cos\theta+\cos(2\theta)\ge0\)
to prime powers outside the conductors gives
\[
  \zeta_E^S(s)^3 |L^S(s,\beta^{-1})|^4 |L^S(s,\beta^{-2})|\ge1 ,
\]
where \(S\) contains all ramified primes. The omitted local factors
are finite and nonzero at one. The first factor has pole order three,
and \(L^S(s,\beta^{-2})\) is holomorphic there since its infinity
weights are nonzero. A zero of \(L^S(s,\beta^{-1})\) at one would
force the left-hand side to tend to zero. Hence \(L(1,\beta^{-1})\ne0\). These
facts give \(\Lambda_{\varnothing}(1)\ne0\) and
\begin{equation}\label{eq:determinant-marked-ratio}
 \frac{|\Lambda_Q(\tau)|^2}{|\Lambda_{\varnothing}(1)|^2}
   =q_Q\left|
       \frac{L^{S_Q}(1,((\nu\tau)^{\rm bc})^{-1})}
            {L^{S_\varnothing}(1,(\nu^{\rm bc})^{-1})}
       \right|^2 .
\end{equation}
Keep the character \(\tau\) fixed and insert a split place
\(w\) at which \(\beta\) is unramified. Put \(q=q_w\),
\(A=\beta(\varpi_w)\) and \(f=\mathbf1_{M_2(\mathcal O_w)}\).
The new local test is \(\mathbf1_{\mathrm{GL}_2(\mathcal O_w)}\).
For each line \(U\subset k_w^2\) choose \(h_U\) of determinant
\(\varpi_w\) with \(h_U\mathcal O_w^2\) the inverse image of \(U\).
Counting the lines containing the image of a matrix of rank
zero, one or two gives, on \(M_2(E_w)\),
\[
 \mathbf1_{\mathrm{GL}_2(\mathcal O_w)}(X)
 =f(X)-\sum_U f(h_U^{-1}X)+qf(\varpi_w^{-1}X).
\]
All other local tests and both characters are fixed. Let \(T\)
be the quotient integral defining the period with variable test
at \(w\). Right invariance of target Haar and the split oscillator
action give
\[
 T(f(h^{-1}\,\cdot))=|\det h|_w B(\det h)^{-1}T(f).
\]
The \(q+1\) line terms each contribute \(q^{-1}A^{-1}T(f)\);
the last test, with \(h=\varpi_w I_2\), contributes
\(q^{-2}A^{-2}T(f)\) before its coefficient \(q\).
Consequently
\[
 \Lambda_{Q\cup\{w\}}(\tau)=
 (1-A^{-1})(1-q^{-1}A^{-1})\Lambda_Q(\tau).
\]
Since \(|A|=1\), the squared period is multiplied by
\begin{equation}\label{eq:rank-two-euler-deletion}
      (1-A)(1-A^{-1})(1-q_w^{-1}A)(1-q_w^{-1}A^{-1}).
\end{equation}
This identity is a multiplication identity and remains
valid when the multiplier or the period is zero.

\subsubsection{Integral coefficients of the rank-two norm}
\label{sec:rank-two-integral-coefficients}

\sectionopening{In this subsection, we express the norm of Section~\ref{sec:rank-two-doubling-local-integrals} in the original finite coordinates.}

Let \(h_c\) be representatives in
\eqref{eq:determinant-marked-period}. For
\(y\in\mathrm{GL}_2(\mathcal O_E/\mathfrak m)\), let
\(\widetilde y\) be the split integral lift in the
original coordinates at the chosen places, with
identity components elsewhere, and set
\(h_{c,y}=h_c\widetilde y\), \(u_{c,y}=\det h_{c,y}\).
Let \(K^{\rm pt}_{Q,1}(h)\) be the kernel at the same
modulus with the point indicator of the identity and
with \(\tau=1\). Expand the determinant character
as the sum of the point indicators at invertible
matrices \(y\), weighted by
\(\chi_{\mathfrak m}(\det y)\). Write \(\delta_y\) for the indicator of \(X\equiv y\) at this
modulus, and \(\alpha_1=\nu\upsilon|_{U(1)}\). The split action gives
\[
 \omega_{\alpha,\alpha}(1,\widetilde y)\delta_1
       =\chi_{\mathfrak m}(\det y)\delta_y,\qquad
 \theta^{\alpha,\alpha}_{\delta_1}(1,h_{c,y})
       =\tau(u_{c,y})\,\theta^{\alpha_1,\alpha_1}_{\delta_1}(1,h_{c,y}).
\]
In the second equality the modulus of \(\delta_1\) still includes
every prime in \(Q\), even where the character under evaluation is
unramified. Thus
\[
 \Lambda_Q(\tau)=
    \sum_{c,y}\frac{K^{\rm pt}_{Q,1}(h_{c,y})\tau(u_{c,y})}
                    {|\operatorname{Aut}(h_c,K_2)|}.
\]
The stabilizer order here divides \(C_G\), by reduction
at an old odd completely split place as in the
preceding subsection.

Use the fixed nonzero integral-jet scale
\[
 s_p=(M^2)^{8f}s_{\rm vac}^{4}
           (M^2d_Ks_{\rm jet})^{2\sum_v k_v}
\]
for the original rank-two point test at modulus
\(\mathfrak m\). Section~\ref{sec:marked-both-norm-jets} gives
\[
 \iota(s_p K^{\rm pt}_{Q,1}(h_{c,y}))
       \in\mathcal O_{\overline{\mathbb Q}_\ell}
     \qquad(\ell\nmid q_Q)
\]
at every coefficient embedding, with the same \(s_p\)
for all \(Q,c,y\).

Set
\[
 a_{c,y}(Q)=
   \frac{C_Gs_p K^{\rm pt}_{Q,1}(h_{c,y})}
        {|\operatorname{Aut}(h_c,K_2)|}.
\]
Its coefficients are integral away from \(q_Q\) by this jet bound. With \(i,j\) ranging over the pairs
\((c,y)\), the squared norm system is
\[
 F_{p,Q}(\tau)=|C_Gs_p\Lambda_Q(\tau)|^2
      =\sum_{i,j}a_i(Q)\overline{a_j(Q)}\,\tau(u_i u_j^{-1}).
\]
The conjugate coefficients are integral by applying
the bound at the conjugate embedding. Formula
\eqref{eq:determinant-marked-ratio} gives its
exact-support ratio, with the same scale \(s_p\)
on numerator and denominator.

\subsection{The integral trace algebra on the auxiliary levels}
\label{sec:integral-raised-traces}

\sectionopening{In this subsection, the attachment of Section~\ref{sec:cofinite-spherical-galois}
identifies traces with Hecke data at a cofinite set of primes. We use
the integral Hecke action there to place traces at all Galois elements
in one coefficient algebra.}

Fix a coefficient embedding above \(\ell\) and any of the
finitely raised levels under consideration. The old type and twist
characters are fixed. Raising changes only the invariance at the
auxiliary split places away from \(\ell\). Let
\(V_{\rm orig}\subseteq V_{\rm raised}\) be these two spaces.
Their Weyl lattices are defined by the same inequalities on the
values in transported coordinates. Consequently
\[
 L_{\rm orig}=V_{\rm orig}\cap L_{\rm raised},\qquad
 V_{\rm orig}\cap aL_{\rm raised}=aL_{\rm orig}\quad(a\ne0),
\]
by dividing the vector by \(a\) in \(V_{\rm orig}\).
On the ambient integral space all local translations
supported away from \(\ell\) preserve this lattice norm. In fact
the lattice norm of a value at a point is read in the transported
coordinates of that point at the places above \(\ell\).
Right translating by a point with identity components above
\(\ell\) changes neither those coordinates nor the list of
points where the lattice norm is bounded. The finite local sums
that preserve the raised invariants therefore act on the
raised integral lattice.

For a system \(x\) of spherical eigenvalues at almost all split
places, identify systems that agree off a finite set. Write
\(W_x\) for the space of forms with these eigenvalues at almost
all split places. For each level the corresponding space is
\(V\cap W_x\). Thus
\[
 V_{\rm orig}\cap W_x\subseteq V_{\rm raised}\cap W_x .
\]
The set \(Y=\{x:V_{\rm raised}\cap W_x\ne0\}\) is finite by
finite dimensionality and orthogonality of the distinct systems.
Auxiliary translations preserve \(W_x\), since they commute
with spherical Hecke operators away from their support.
Define \(\mathbb T\) to
be the ring of functions \(Y\to\overline{\mathbb Q}_\ell\)
whose diagonal operators preserve the raised integral lattice.
It contains the scalar valuation ring. It is closed in the
product of the coefficient topologies: for each vector in
the integral lattice and each test in the lattice coordinates, the
result of the diagonal operator is a continuous linear
function of the finitely many scalars. The condition that
this result be in the valuation ring is closed, and the
intersection of all these conditions is closed. This
argument uses the original lattice norm at every class point,
so it is valid with stabilizers as well.

For this fixed raised level apply
Section~\ref{sec:cofinite-spherical-galois} to the finite set \(Y\).
Choose a basis of each \(V_{\rm raised}\cap W_x\).
Remove the selected triangle primes, the exceptional places of
these finitely many basis vectors, and the excluded places of the
Galois comparisons for \(Y\) from the split places away from \(\ell\).
On the remaining set, an original eigenvector and its image in the
raised space have the same spherical Hecke polynomials, since
the inclusion commutes with the local translations there.
Density and semisimplicity therefore identify their Galois attachments.
For newly occurring systems, the comparison at a selected triangle
prime uses the split local parameter of
Section~\ref{sec:raised-tame-frobenius}.
On the common good set the minuscule
spherical sums act by the stipulated scalar on each outside-level
isotypic space. As finite sums of the integral translations above,
their diagonal scalar functions belong to \(\mathbb T\). The coefficients of the normalized
Hecke polynomial are integral multiples of these
operators, with powers of \(q_w\) as in
Section~\ref{sec:cofinite-spherical-galois}.
Fix \(g\in G_E\). The cofinite Frobenius comparison
for each member of \(Y\), and Chebotarev,
approximate the coefficients of
\(\det(X-R_x(g))\) simultaneously, for all \(x\in Y\),
by these Hecke coefficients.

Indeed for any required
precision, continuity gives an open normal subgroup
on whose coset at \(g\) the coefficients agree to that
precision. The functions are conjugacy invariant.
Chebotarev's theorem \cite[Theorem~8.31]{MilneANT} gives
positive density for the conjugacy class in the finite quotient.
The prime ideals of \(E\) of residue degree at least two over
\(\mathbb Q\) have Dirichlet density zero, since their sum of
\((\mathrm N w)^{-s}\) for \(s>1\) is bounded by
\([E:\mathbb Q]\sum_{p}p^{-2}\). A degree-one prime outside the
finite ramified set must split over \(F_0\), by multiplicativity
of residue degrees. Thus the required class occurs in the common
good set.
Closedness now gives each coefficient function in
\(\mathbb T\).

For an integral character \(a\), putting \(N_x=aR_x\)
multiplies the coefficient of \(X^{n-i}\) by \(a(g)^i\);
hence these twisted characteristic coefficients also
belong to \(\mathbb T\). In particular \(x\mapsto
\operatorname{tr}(N_x(g))\) belongs to \(\mathbb T\)
for every \(g\), including products involving the selected
tame generators and Frobenius lifts. The trace functions
used below are integral linear combinations of these.
For instance expansion of \(p\) in monomials writes
\[
 \operatorname{tr}\bigl(p(N_x(\phi))(N_x(\tau)-1)N_x(g)\bigr)
 =\sum_i p_i\bigl(\operatorname{tr}(N_x(\phi^i\tau g))
                     -\operatorname{tr}(N_x(\phi^i g))\bigr).
\]
All of these diagonal operators commute with the local
translations preserving the decomposition into
outside-level eigensystems. Thus the raised local comparison
uses the same integral endomorphism ring and original
invariant lattice throughout.

\section{Hecke operators and the common trace coefficient}
\label{sec:hecke-trace-comparison}

\sectionopening{In this section, we identify the tame and Frobenius operators in the algebra of
Section~\ref{sec:integral-raised-traces}. Spherical compatibility gives a
common Hecke coefficient to which the Fitting and trace arguments of
Section~\ref{sec:trace-determinant-specialization} apply.}

\subsection{Tame and Frobenius operators at a raised prime}
\label{sec:raised-tame-frobenius}

\sectionopening{In this subsection, we represent tame inertia and Frobenius by operators on the
\(K_1\)-invariants at the same split place. The resulting trace
identities refine the integral trace comparison of
Section~\ref{sec:integral-raised-traces} by retaining the local finite
type.}

Fix a coefficient prime \(\ell\), a chosen split place \(w\nmid\ell\)
off the original level, and a congruence depth \(\nu\). Write \(k=k_w\),
\(q=|k|\), \(H=\mathrm{GL}_n(k)\), \(K=\mathrm{GL}_n(\mathcal O_w)\),
and \(K_1=\ker(K\to H)\). At \(w\) the original level is \(K\)
and the raised level is \(K_1\). The other local factors stay fixed.
We assume \(\ell^\nu\mid q-1\) and fix the prescribed geometric
Frobenius lift \(\phi\) and tame generator \(\tau\) in the local
Weil group. In its abelianization choose \(\varpi\in E_w^\times\)
and \(u\in\mathcal O_w^\times\) corresponding to their images
under geometric local reciprocity. The Frobenius degree of \(\phi\)
gives \(v_w(\varpi)=1\). Let \(\lambda:k^\times\to A=\mathbb Z/\ell^\nu
\mathbb Z\) be the chosen homomorphism with \(\lambda(\bar u)=1\).
Put \(d=\mathrm{diag}(\bar u,1,\ldots,1)\).

Write \(\rho\) for the \(H\)-action on the raised space and
\(\mathcal A\) for the integral endomorphisms of its lattice.
The raised lattice is the relevant invariant space in the ambient
integral lattice on which all translates at \(w\) act. Section~\ref{sec:integral-raised-traces} gives the
integral diagonal trace action on this same lattice.
The outside-level eigensystems decompose the raised space over
characteristic zero, and the local operators at \(w\) preserve each
summand.

On every raised system \(x\) use its semisimple attachment
\(R_x^{\rm up}\) of rank \(n\) and the full split-place comparison from
Section~\ref{sec:inertial-labels}. At a split level place this
comparison follows from the same split base change on each
constituent of the eigenvector. The cofinite comparison of
Section~\ref{sec:cofinite-spherical-galois} identifies the
semisimple attachments of all constituents with system \(x\).
Choose the normalization \(N_x=aR_x^{\rm up}\), where \(a(\tau)=1\).
For each Galois element \(g\), write
\[
 F_x=N_x(\phi),\qquad T_{x,\tau}=N_x(\tau),\qquad M_x=N_x(g),
\]
and let \(C_p,D_p\in\mathcal A\) be the diagonal operators
whose entries for an integral polynomial \(p\) are
\[
 C_p(x)=\operatorname{tr}(p(F_x)M_x),\qquad
 D_p(x)=\operatorname{tr}(p(F_x)(T_{x,\tau}-1)M_x).
\]
All scalars are evaluated at the same embedding.

Let \(\bar P\subset H\) be the stabilizer of \(ke_1\),
\(J\subset K\) its inverse image, and \(I_1\subset K\)
the inverse image of the upper unitriangular group.
For a character \(\psi\) of \(A\), put \(\chi=\psi\circ\lambda\).
On the \(H\)-type \(\chi\circ\det\), the inverse determinant
twist has a \(K\)-fixed vector. Its parameter is spherical.
Twisting back gives zero monodromy, and \(\tau\) acts on
the entire Galois representation by \(\psi(1)\), the eigenvalue
of \(\rho(d)\) on the same type.

Consider now the \(J\)-type \(j\mapsto\chi(\bar j_{11})\)
with \(\chi\ne1\). In an irreducible local constituent \(\pi\),
finite Frobenius reciprocity gives a common constituent of
\(\pi^{K_1}\) and \(\operatorname{Ind}_{\bar P}^H(\chi\boxtimes1)\).
The latter embeds into the Borel induction from
\((\chi,1,\ldots,1)\) by transitivity of induction.
Write \(\pi\) as a subquotient of an induction from its
supercuspidal support. Iwasawa decomposition and exactness of
compact invariants give the \(K_1\)-invariants of that induction
as finite induction from the corresponding principal-kernel
invariants of its factors.

For a supercuspidal factor of rank \(b\), set
\(K_b=\mathrm{GL}_b(\mathcal O_w)\).
If its \(K_{b,1}\)-invariants are nonzero, it has depth zero
at the vertex of the standard lattice.
By \cite[Propositions~6.6 and~6.8, pp.~114--115, 117--118]{MoyPrasadJacquet},
it is compactly induced from a representation of a maximal parahoric
normalizer containing an inflated finite cuspidal representation.
In \(\mathrm{GL}_b(E_w)\) every such vertex is a lattice class;
choose a basis for the lattice to take its parahoric to \(K_b\).
The normalizer \(E_w^\times K_b\) acts on the same finite cuspidal
space, since the centre acts by scalars. To calculate the invariants
of this compact induction, take Cartan representatives \(s=\mathrm{diag}(\varpi^{v_1},\ldots,
\varpi^{v_b})\), with \(v_1\le\cdots\le v_b\), modulo
\(E_w^\times K_b\) on the left and \(K_b\) on the right.
For any strict cut \(v_i<v_{i+1}\), the reduction of
\(K_b\cap sK_{b,1}s^{-1}\) contains the upper rectangular
unipotent radical at that cut, since conjugation by \(s^{-1}\)
multiplies the entry in row \(a\le i\), column \(b'>i\), by
\(\varpi^{v_{b'}-v_a}\). The intersection is the same for
\(s k_0\), \(k_0\in K_b\), because \(K_{b,1}\) is normal in
\(K_b\). Cuspidality forces the value of every
\(K_{b,1}\)-invariant vector at those cosets to vanish.
The central Cartan class gives the finite cuspidal factor.

\subsubsection{The restriction to the tame subgroup}
\label{sec:raised-prime-tame-restriction}

\sectionopening{In this subsection, we calculate the tame subspace using the local Hecke data of Section~\ref{sec:cofinite-spherical-galois}.}

Now use permutation representatives in finite Mackey theory
between the Borel and the block parabolic. For each representative,
the intersection in each conjugated Levi block contains the
unipotent radical of a Borel of that block. A finite cuspidal of
rank greater than one has no such invariants, because that radical
contains the radical of a proper parabolic. Thus every block has
rank one. Comparing the characters on the diagonal torus gives
the tuple \((\chi,1,\ldots,1)\), up to permutation.
Unramified twists preserve restriction to \(\mathcal O_w^\times\).
The entry with character \(\chi\ne1\) therefore forms a segment
of length one on its own cuspidal line. The classification
by segments and the product over distinct lines give
\[
 \pi=\operatorname{Ind}_{P}^{\mathrm{GL}_n(E_w)}
                    (\chi_0\boxtimes\pi_1)
             \quad\text{with normalized induction},
\]
where \(P\) is the stabilizer of the first coordinate line,
\(\chi_0|_{\mathcal O_w^\times}=\chi\) via reduction, and
\(\pi_1\) has all-unramified rank-one support
\cite[Theorem~6.1 and Proposition~8.6, pp.~188, 197]{ZelevinskyInducedII}.

Apply finite Mackey theory to \(\bar P\backslash H/\bar P\).
For \(n\ge2\) the two representatives are \(1\) and \(s=(12)\).
Choose \(b\in k^\times\) with \(\chi(b)\ne1\). The element
\(h=\operatorname{diag}(1,b,1,\ldots,1)\) belongs to
\(\bar P\cap s\bar P s^{-1}\) and acts trivially by the prescribed
character \(j\mapsto\chi(j_{11})\). On the conjugated inducing
representation it acts by \(\chi(b)\), since
\(s^{-1}hs=\operatorname{diag}(b,1,\ldots,1)\) acts on the
first block alone. The intertwining space for this coset is
therefore zero. The identity coset gives
\(\pi_1^{\mathrm{GL}_{n-1}(\mathcal O_w)}\).
In particular \(\pi_1\) is spherical and evaluation at 1
is injective on the required \(J\)-type.
For \(n=1\) there is only \(\chi_0\) and the evaluation is immediate.

Set \(t=\mathrm{diag}(\varpi,1,\ldots,1)\) and let \(U\)
be the unnormalized operator of \(I_1tI_1\).
For \(b\in k^{n-1}\), put
\(u_b=\left(\begin{smallmatrix}1&\widetilde b\\0&I_{n-1}\end{smallmatrix}\right)\).
Both \(I_1\cap tI_1t^{-1}\) and \(J\cap tJt^{-1}\) impose
precisely the extra divisibility by \(\varpi\) on the off-diagonal
entries of the first row. Thus the cosets are \(u_btI_1\) and \(u_btJ\), respectively. The \(J\)-character agrees with
its \(t\)-conjugate on the intersection because the \((1,1)\)
entry is unchanged. Left multiplication by \(j\in J\) permutes
the cosets and contributes its character to each term, so \(U\)
preserves the type. In the induced model every term has value
\[
 f(u_bt)=q^{-(n-1)/2}\chi_0(\varpi)f(1).
\]
There are \(q^{n-1}\) terms. Evaluation at 1 therefore shows
that the \(U\)-scalar is \(q^{(n-1)/2}\chi_0(\varpi)\).

Both blocks have zero monodromy, and the complementary block
has trivial inertia. On the unique \(\psi(1)\)-inertia line the
Frobenius scalar after the twist \(|\det|^{(1-n)/2}\) is
\(q^{(n-1)/2}\chi_0(\varpi)\). Frobenius semisimplification
preserves this scalar on the one-dimensional line. Hence
\(T_{x,\tau}-1\) acts only on this line, where \(F_x\) has the scalar
of \(\widetilde U=a(\phi)U\). Write
\(\alpha_\phi=a(\phi)q^{(n-1)/2}\chi_0(\varpi)\).
Then \(p(F_x)(T_{x,\tau}-1)=p(\alpha_\phi)(T_{x,\tau}-1)\),
since the image of \(T_{x,\tau}-1\) lies on that line.
Right multiplication by \(M_x\) and taking the trace yield
\(D_p=p(\widetilde U)D_1\) on the given \(J\)-type for every \(g\).
This and the determinant-type calculation apply on the
corresponding types of every raised eigensystem, by taking
irreducible components.

For an integer-valued weight \(f\) on \(A\), define
\[
 \begin{aligned}
 L_J(f)&=\sum_{k\in\bar P}f(\lambda(k_{11}))\rho(k),&
 L_H(f)&=\sum_{h\in H}f(\lambda(\det h))\rho(h),\\
 B(f)&=p(\widetilde U)L_J(f).
 \end{aligned}
\]
Left multiplication by \(I_1\) preserves the line weight, so
the image of \(L_J(f)\) is \(I_1\)-fixed. The translates in
the coset sum \(U\) act on the ambient integral lattice, and
their sum on \(I_1\)-invariants preserves those invariants.
As \(K_1\subset I_1\), \(B(f)\) is an integral operator on
the raised lattice made from translates at \(w\).
Let \(\mathcal E(f)=B(f)D_1-L_J(f)D_p\). Fourier expansion
in \(A\) and the evaluations on the two types give
\[
 L_H(f)D_p=(\rho(d)-1)L_H(f)C_p,\qquad
 \mathcal E(f(\cdot-n!)-f)=0.
\]
The weight \(\psi^{-1}\) has image in the \(\psi\)-type, and
its coefficient in \(f(\cdot-n!)-f\) is multiplied by
\(\psi(n!)-1\). The equality is
tested in characteristic zero between the integral operators.

Let \(\omega:A\to\mathbb Z\) give the signed representatives.
Taking \(f(s)=\omega(s)(\omega(s)+n!)\) gives
\(f(s-n!)-f(s)\equiv-2n!\,\omega(s)\pmod{\ell^\nu}\), so
\(2n!\mathcal E(\omega)\in\ell^\nu\mathcal A\).
Let \(\sigma\) be a reduction homomorphism on \(\mathcal A\)
which kills every \(X\) with \(2n!X\in\ell^\nu\mathcal A\).
For example, use reduction modulo \(\ell^\nu/(2n!)\) when
that quotient belongs to the coefficient valuation ring.
Write \(\mathcal N_w=\sum_{h\in H}\rho(h)\).
The determinant identity with the constant weight gives
\(\mathcal N_wD_1=0\) exactly.

Choose representatives \(r_s\) with \(r_s\bar P\) running
over the lines, and set
\[
 C=\sum_s\rho(r_s),\quad C^\vee=\sum_s\rho(r_s^{-1}),\quad
 \partial=-CB(\omega)C^\vee.
\]
Write \(hr_s=r_{hs}k_{h,s}\) with \(k_{h,s}\in\bar P\).
In the group algebra the coefficient of \(h\) in
\(CL_J(\omega)C^\vee\) is
\(\sum_s\omega(\lambda((k_{h,s})_{11}))\).
The product \(\prod_s(k_{h,s})_{11}\) is a homomorphism in
\(h\) by the coset cocycle identity and is independent of the
representatives because their changes telescope.
On \(\operatorname{diag}(b,1,\ldots,1)\) the product is \(b\).

For each line other than \(ke_1\) choose a representative
vector with its first nonzero coordinate among coordinates
2 through \(n\) equal to 1. The diagonal element preserves
this normalization, so the only contribution is that of
\(ke_1\). A transvection has residue-characteristic order,
prime to \(|k^\times|\), and therefore contributes 1.
These elements generate \(H\), so the product is the
determinant. It follows that
\[
 CL_J(\omega)C^\vee\equiv L_H(\omega),\qquad
 (\rho(d)-1)L_H(\omega)\equiv-\mathcal N_w
                          \pmod{\ell^\nu\mathcal A}.
\]
The second congruence uses
\(\omega(s-1)-\omega(s)\equiv-1\pmod{\ell^\nu}\).

Since the diagonal trace operators commute with the local action,
the two congruences and \(\sigma(\mathcal E(\omega))=0\) give
\[
 \begin{aligned}
 \sigma(\partial D_1)
 &=-\sigma\bigl(CL_J(\omega)C^\vee D_p\bigr)
   =-\sigma\bigl(L_H(\omega)D_p\bigr)\\
 &=-\sigma\bigl((\rho(d)-1)L_H(\omega)C_p\bigr)
   =\sigma(\mathcal N_w C_p).
 \end{aligned}
\]
If \(p(X)=(X-q^b)H_0(X)\), with \(b\ge0\) and \(H_0\) integral,
then \(q^b\equiv1\pmod{\ell^\nu}\) and
\(\sigma(\partial D_1)=\sigma(\mathcal N_w C_{(X-1)H_0})\).
The operator \(\partial\) is independent of \(g\) and uses
the original unnormalized norm on the full raised lattice.
At every selected place the argument uses the same depth
and reduction, with the hypotheses stated at the start of
the subsection.

\subsubsection{Integral clearing from the tame action}
\label{sec:raised-prime-integral-clearing}

\sectionopening{In this subsection, we use the tame restriction of Section~\ref{sec:raised-prime-tame-restriction} to pass to integral vectors.}

For use with clearing, suppose an integral endomorphism
\(\Psi\) and scalar \(c\) satisfy
\[
 \sum_{k\in\bar P}\rho(k^{-1})\Psi\rho(k)=c\,1.
\]
Let \(\mathcal N_{\mathrm{all}}\) be the product of the selected
full norms, including \(\mathcal N_w\), and suppose the
spherical annihilation on the supported original systems gives
\(\delta\mathcal N_{\mathrm{all}}B(1)=0\), with \(1\) the
constant weight. For \(h\in H\) set
\(A_h=\sum_s\omega(\lambda((k_{h,s})_{11}))\rho(r_s^{-1})\).
The equations
\[
 B(f)\rho(k)=B(f(\cdot-\lambda(k_{11}))),\qquad
 r_{hs}^{-1}h=k_{h,s}r_s^{-1}
\]
and the signed weight congruence give
\[
 \partial(\rho(h)-1)\equiv C B(1) A_h\pmod{\ell^\nu\mathcal A}.
\]
Here the difference of the weights at \(k_{h,s}\) is
congruent to the negative of its logarithm.

The double-coset
calculation for the trivial \(J\)-character gives
\(\rho(k)B(1)=B(1)\). Applying the identity for \(\Psi\)
to \(B(1)\) yields \(c CB(1)=\mathcal N_w\Psi B(1)\).
Set \(Z=\Psi B(1)A_h\). Then
\[
 \sigma\bigl(c\partial(\rho(h)-1)\bigr)=\sigma(\mathcal N_w Z).
\]
Before applying \(C\) to that identity for \(\Psi\), multiply
on the left by \(\delta\mathcal N_{\mathrm{all}}\). Norm
invariance gives
\[
 c\delta\mathcal N_{\mathrm{all}}B(1)
       =|\bar P|\delta\mathcal N_{\mathrm{all}}\Psi B(1).
\]
The left side vanishes. Integral endomorphisms have no integer
torsion, so \(\delta\mathcal N_{\mathrm{all}} Z=0\).

\subsection{Spherical root chains on the original level}
\label{sec:marked-spherical-root-chain}

\sectionopening{In this subsection, we factor the normalized spherical polynomial into successive roots in
the pro-unipotent Iwahori space. The construction preserves the
original spherical vector and keeps the local normalization used in
Section~\ref{sec:raised-tame-frobenius}.}

In this subsection \(r\) is the local rank and \(t\) is the length
of the specified chain of roots. Let \(K_0=\operatorname{GL}_r(\mathcal O_w)\) with
residue field of size \(q\), and use the fixed ordered
coordinates and uniformizer \(\varpi\) at the
auxiliary place. Let \(I_1\subset K_0\) be the
pro-unipotent Iwahori. Work in the given smooth
unitary representation with positive invariant
pairing \(B\), conjugate linear in the first
argument. Let \(x\ne0\) be the original
\(K_0\)-fixed eigenvector, with eigenvalues \(t_j\)
for the unnormalized spherical coset sums of
\(\operatorname{diag}(\varpi I_j,I_{r-j})\).

The polynomial is
\[
  P(X)=\sum_{j=0}^r a_jX^{r-j},
  \qquad a_j=(-1)^j q^{j(j-1)/2}t_j,\quad t_0=1 .
\]
Put \(d=\operatorname{diag}(\varpi,1,\ldots,1)\).
The operator \(U\) in the derivative congruence is
the unnormalized \(I_1\)-double-coset sum of
\(d\). On \(I_1\)-fixed vectors, it is
\(q^{r-1} M_{I_1}\rho(d)M_{I_1}\), where \(M_{I_1}\)
is probability Haar averaging.
Indeed \(I_1\cap d I_1d^{-1}\) imposes the
vanishing modulo \(\varpi\) of the off-diagonal
entries in the first row, giving index \(q^{r-1}\).
The adjoint is
\(U^\dagger=q^{r-1}M_{I_1}\rho(d^{-1})M_{I_1}\).

For \(m,j\ge0\), let \(T_{m,j}\) be the sum of
\(\rho(g)x\) over the lattices \(L=g\mathcal O_w^r\)
between \(\varpi\mathcal O_w^r\) and
\(\mathcal O_w^r\), of colength \(j\), containing
the first \(m\) coordinate vectors. A representative
\(g\) gives a well-defined vector because \(x\) is
\(K_0\)-fixed. In particular \(T_{0,j}=t_jx\).
Put \(T_j=T_{1,j}\).
These vectors are fixed by \(I_1\), which
permutes the relevant lattice sets.

The interval of lattice sums gives a forward
recursion. For each residue covector \(\xi\) with first
entry \(1\), take \(n_\xi\in I_1\) the identity
with first row lifting \(\xi\).
For a lattice in \(T_{0,j+1}\), its occurrence in
\(n_\xi^{-1}dT_j\) means that its reduction lies in
the hyperplane \(\xi=0\). To check the condition,
\(n_\xi L\) must have first residue coordinate zero;
then \(d^{-1}n_\xi L\) is integral, contains
\(\varpi\mathcal O_w^r\) and contains \(e_1\)
since \(n_\xi L\supset\varpi\mathcal O_w^r\).
The colength has dropped by one.
There are no such covectors if \(L\) contains \(e_1\).
Otherwise the requirement that \(\xi\) vanish
on a subspace of dimension \(r-j-1\) and take
value \(1\) on \(e_1\) leaves \(j\) free
coordinates, hence \(q^j\) covectors. Applying
\(M_{I_1}\), which kills the action of each
\(n_\xi^{-1}\) on the left, gives
\begin{equation*}\tag{SC1}\label{eq:display-sc1}
             U T_j=q^j(t_{j+1}x-T_{j+1}) .
\end{equation*}
In this notation \(T_0=x\) and \(T_r=0\).

For the inverse count use the representation
\(\rho'(g)=\rho((g^{\mathsf T})^{-1})\).
Under duality the lattice \(L\) becomes
\(L^\vee=(g^{\mathsf T})^{-1}\mathcal O_w^r\)
for the standard integral bilinear pairing.
For \(L\) containing the first \(m\) vectors,
the first \(m\) dual coordinates are integral.
With
\(C_m=\operatorname{diag}(I_m,\varpi^{-1}I_{r-m})\),
the lattice \(C_m^{-1}L^\vee\) is between
\(\varpi\mathcal O_w^r\) and \(\mathcal O_w^r\),
contains those same coordinate vectors,
and its colength is \(r-m-j\).
This is a bijection, by applying duality
and the inverse diagonal scaling.
It proves
\[
 T'_{m,j}=\rho(C_m)T_{m,r-m-j},
\]
where primes refer to \(\rho'\).
In the forward count for \(\rho'\) test by
the linear map \(\Psi=M_{I_1}\rho(C_1^{-1})\).

For every first-row unipotent \(n_\xi\) of
that count, \(C_1^{-1}n_\xi^{\mathsf T}C_1\)
is in \(I_1\): its off-diagonal first-column
entries are multiplied by \(\varpi\).
Hence \(\Psi\rho'(n_\xi^{-1})=\Psi\),
as required in the same hyperplane count.
Use \(d=C_0^{-1}C_1\). For \(1\le i<r\), take the dual
colength \(j=r-1-i\). The left side of the dual count,
after applying \(\Psi\), is
\(q^{r-1}\Psi\rho'(d)T'_{1,j}=U^\dagger T_i\), since
\(d\) commutes with \(C_1\) and \(T_i\) is \(I_1\)-fixed.
The two lattice sums on the right give
\[
 \Psi T'_{0,j+1}=t_i M_{I_1}\rho(d^{-1})x,\qquad
 \Psi T'_{1,j+1}=T_{i-1}.
\]
The multiplicity remains \(q^j=q^{r-1-i}\), so
\begin{equation*}\tag{SC2}\label{eq:display-sc2}
 U^\dagger T_i
   =q^{r-1-i}\bigl(q^{-(r-1)}t_i U^\dagger x-T_{i-1}\bigr).
\end{equation*}
At the end, \(T_{r-1}\) is the single lattice
of \(C_1^{-1}\), and \(t_r x=\rho(\varpi I_r)x\).
Thus \(t_r\ne0\),
\[
 U^\dagger x=q^{r-1}t_r^{-1}T_{r-1},
 \qquad U T_{r-1}=q^{r-1}t_rx,\qquad
 U U^\dagger x=q^{2r-2}x .
\]
Write \(P_i(X)=\sum_{j=0}^i a_j X^{i-j}\).
Then \(P_r=P\) and \(P_{i+1}=XP_i+a_{i+1}\).

\begin{lemma}[The spherical root relation]\label{lem:marked-spherical-annihilation}
For the spherical vector \(x\), the Hecke operator \(U\), and
its polynomials \(P_i\), one has
\[
 P_i(U)x=(-1)^i q^{\binom i2}T_i\quad(0\le i\le r),\qquad P(U)x=0 .
\]
\end{lemma}
\begin{proof}
For \(i=0\) both sides are \(x\). Using \eqref{eq:display-sc1}
and the definition of \(a_{i+1}\), the induction step is
\[
\begin{aligned}
 P_{i+1}(U)x
 &=(-1)^i q^{\binom i2+i}(t_{i+1}x-T_{i+1})+a_{i+1}x\\
 &=(-1)^{i+1}q^{\binom{i+1}{2}}T_{i+1}.
\end{aligned}
\]
At \(i+1=r\), the only lattice of colength \(r\) between
\(\varpi\mathcal O_w^r\) and \(\mathcal O_w^r\) is
\(\varpi\mathcal O_w^r\), which omits the first basis vector.
Thus \(T_r=0\) and \(P(U)x=0\).
\end{proof}

For \(1\le i\le r\), \eqref{eq:display-sc2} and the endpoint
identity give
\begin{equation*}\tag{SC3}\label{eq:display-sc3}
 U^\dagger P_i(U)x
      =a_i q^{-i} U^\dagger x
        +q^{r-2}P_{i-1}(U)x .
\end{equation*}
For \(i=r\) both terms on the right cancel,
using the formula for \(U^\dagger x\).
For any scalar \(c\), set \(P_{r+1}=XP_r\) and
\(D_i(X)=(P_i(X)-P_i(c))/(X-c)\) for \(0\le i\le r+1\).
Then \(D_{i+1}=P_i+cD_i\) for \(0\le i\le r\).
By \eqref{eq:display-sc3}, induction gives
\[
 U^\dagger D_{i+1}(U)x
    =q^{-i}P_i(qc)U^\dagger x
         +q^{r-2}D_i(U)x,\qquad 0\le i\le r .
\]
The case \(i=0\) is \(D_1=1\), \(D_0=0\).
In the step use
\(q^{-i}P_i(qc)=a_i q^{-i}
          +c q^{-(i-1)}P_{i-1}(qc)\).

If \(P(c)=0\), put \(z=D_r(U)x\).
Since the constant term of \(P\) is nonzero,
\(c\ne0\). The divided difference identity
and \(P(U)x=0\) give \(Uz=cz\) and
\[
 cU^\dagger z=q^{-r}P(qc)U^\dagger x+q^{r-2}z .
\]
Pair this last identity first with \(U^\dagger x\)
in the first variable. From \(UU^\dagger x=q^{2r-2}x\)
one obtains
\[
        c(q^r-1) B(x,z)=P(qc) B(x,x).
\]
Pairing with \(z\) in the second variable instead
and substituting this equality gives
\[
 (|c|^2-q^{r-2}) B(z,z)=
       \frac{|P(qc)|^2}{q^r(q^r-1)}B(x,x).
\]
Since \(x\ne0\), positive definiteness gives \(B(x,x)>0\).
If \(P(qc)\ne0\), the first identity gives \(z\ne0\), and the
second gives \(|c|^2>q^{r-2}\) since \(B(z,z)>0\).

If \(P(qc)=0\), the first identity gives
\(B(x,z)=0\), and the operator identity
gives \(U^\dagger z=q^{r-2}c^{-1}z\).
Therefore \(B(U^j x,z)=0\) for every \(j\),
by moving the adjoint. Since \(z\) is a
polynomial in \(U\) applied to \(x\), it is
orthogonal to itself, and \(z=0\).
Thus when \(c,qc\) are both roots,
\((P/(X-c))(U)x=0\).
Suppose now \(P(X/a)\) is
\(\prod_{i=0}^{t-1}(X-q^i)\,R(X)\),
with \(a\ne0\). Let \(H\) be any polynomial divisible by \(R\).
At the auxiliary prime indexed by \(j\) we use the normalized
annihilating polynomial \(H_j\) constructed from the complementary
factors in Section~\ref{sec:triangle-unit-input}.
If \(t=1\), \((aU-1)R(aU)=P(U)\) already annihilates \(x\)
and hence its cyclic space. For \(t>1\), the minimal
polynomial of \(U\) on the cyclic
space \(\mathcal C=\mathbb C[U]x\)
divides \(P/(X-q^i/a)\) for every \(i<t-1\),
by the calculation for adjacent roots.

For \(t>1\) these roots are distinct, so taking the greatest
common divisor cancels one copy of each indicated factor from \(P\).
In either case \((aU-q^{t-1})R(aU)\) kills
\(\mathcal C\), and so does
\(p(aU)\) for \(p(X)=(X-q^{t-1})H(X)\).
The space \(\mathcal C\) is spanned by
\(T_0,\ldots,T_{r-1}\) and is stable under
\(U^\dagger\) by \eqref{eq:display-sc2}.
The adjoint of \(p(aU)\) applied to \(x\)
belongs to \(\mathcal C\) and is orthogonal
to it, hence is zero. For every \(I_1\)-fixed
\(y\), therefore, \(B(x,p(aU)y)=0\).
Apply this to \(y=L_J(1)v\), with \(a=a(\phi_w)\), on a system
in the support of the diagonal operator \(\delta\) for which the
displayed factorization holds. Since \(L_J(1)v\) is \(I_1\)-fixed,
\(B(1)v=p(aU)L_J(1)v\) is orthogonal to every original spherical
vector in this system. Over \(\mathbf C\) the product of full norms
\(\mathcal N_{\mathrm{all}}\) is \(\prod_w|H_w|\) times the
orthogonal projection onto the original spherical subspace at all
selected places. It follows, system by system, that
\(\delta\mathcal N_{\mathrm{all}}B(1)=0\) on the raised lattice,
as required in Section~\ref{sec:raised-prime-integral-clearing}.

The same calculation also supplies the
full root chain from its initial modulus.
The adjoint of the \(j\)-th spherical double-coset
sum is the \((r-j)\)-th such sum times
\(\rho(\varpi^{-1}I_r)\), by inverting the
diagonal representative and permuting coordinates.
Thus \(\overline{t_j}=t_r^{-1}t_{r-j}\), since
the original vector is a joint eigenvector and
\(t_r\) is the unitary central eigenvalue.
Substitution into \(P\) gives the reciprocal
identity, for \(z\ne0\),
\[
 q^{r(r-1)/2}\overline{P(z)}
   =(-\bar z)^r\bar t_r\,P(q^{r-1}/\bar z).
\]
For example, reindex the right-side term of
order \(r-j\) by \(j\) and use
\(\binom{r}{2}+\binom{j}{2}
=\binom{r-j}{2}+(r-1)j\).
The roots are therefore stable under
\(z\mapsto q^{r-1}/\bar z\).
If \(P(c)=0\) and \(|c|^2q^k=q^r\), this
reflection interchanges \(cq^i\) and
\(cq^{k-1-i}\).
For \(0<i\le(k-1)/2\), the previous root
\(cq^{i-1}\) has squared modulus \(<q^{r-2}\).
The endpoint bound forces \(cq^i\) to be a root.
Induction covers the lower half and reflection
covers the upper half. 
\subsection{The common coefficient of the Hecke polynomial}
\label{sec:marked-ribet-coefficient}

\sectionopening{In this subsection, we construct one polynomial over the integral Hecke algebra of
Section~\ref{sec:integral-raised-traces} on the entire finite set of
systems. The coefficient of interest is isolated by specializing the
common relations on the selected lifts.}

Fix the descent pair and a nonzero original lift generator with
system \(\pi\). Identify systems by their spherical eigenvalues at
all but finitely many split places, as in
Section~\ref{sec:cofinite-spherical-galois}. Let \(\mathcal X_0\)
be the systems of the original nonzero lift generators for the
chosen twist. At a chosen auxiliary level let \(Y\) be the full
set of systems occurring there, as in
Section~\ref{sec:integral-raised-traces}, and put
\(\mathcal X=Y\cap\mathcal X_0\). The form defining \(\pi\) remains
in its system after raising level, since both realizations agree
at all split places outside the old and auxiliary levels. Thus
\(\pi\in\mathcal X\). On this intersection use the block
decomposition and twisted lattice of the original lift. Write \(r'\) for
the main rank and \(t'\) for the rank difference. Use
the \(\ell\)-adic embedding in the congruence construction and put
\[
 G=C_G,\qquad v=v_\ell(G),\qquad
 H=Y_{\rm tr},\qquad h=h_\ell,\qquad C=\ell^{Hh},\qquad 0\le h\le v.
\]
Let \(J\) be the ordered set of new primes and
\(\mathrm{ty}\colon J\to\{1,\ldots,t'-1\}\) the monotone type.

The coefficient of \(X_J=\prod_{j\in J}X_j\) in the common
Ribet polynomial at \(\pi\) is \(d\). The bound for the
actual passage to the Hecke algebra is
\[
 v_\ell(G^2d)
 \le \mathrm{sel}+2t'\mathrm{ram}+r'\mathrm{loc}_\ell+\mathrm{err}_\ell ,
 \qquad \mathrm{loc}_\ell=4fF_\ell ,
 \tag{RB1}\label{eq:display-rb1}
\]
where \(f\) is the number of infinite places of \(E\), and
\[
 \mathrm{err}_\ell
  =t\,v_\ell(X_{\rm odd})
   +8 ft\,v+2v+\bigl((4t+1)H^{22}+H\bigr)v.
\]
For the selected twist \(\eta\), write \(\psi_\eta\) for
the \(\ell\)-adic Galois avatar of \(\eta^{-1}\).
If \(\eta=\xi_\kappa^j\) and \(\psi\) is the inverse avatar
of \(\xi_\kappa\), then \(\psi_\eta=\psi^j\). Here
\(\mathrm{ram}=\sum_{w\in Q}v_\ell(\psi_\eta(\sigma_w)-1)\),
where \(Q\) contains one place of each pair in the twist support
and \(\sigma_w\) is the tame generator there. For the main
lattice \(L_0\) below, \(\mathrm{sel}\) bounds the normalized length
of the divisible-coefficient Selmer module with \(z(\gamma_\pi)=0\),
that is,
\(\operatorname{length}_{O_0}(\mathcal S_0)\le
\mathrm{sel}\,v_{O_0}(\ell)\).
Fix \(Q,\eta\) and the lattice carrying this Selmer bound before
choosing the congruence depth \(\nu\) and the distinct auxiliary set
\(J\) for the trace comparison. In the application,
Proposition~\ref{prop:joint-norm-selmer-twist} gives the bound for
each power of the same twist, independently of \(\nu\) and \(J\).

For the relative traces, \(c\) may be any common multiple of the
ambient compact stabilizer orders. For the original type projector
and the congruence depth of Section~\ref{sec:top-coefficient-clearing},
fix \(c=G=C_G\) before selecting the raised level. Reduction at
one old split place and block inclusion show that every ambient
stabilizer order divides \(C_G\). The product of the block inclusions
at \(w\in S_M\) also gives
\(\prod_{w\in S_M}|\mathrm{GL}_{r'}(k_w)|\mid C_G\), so \(cQ\)
preserves the original lattice. With \(v_\ell(\ell)=1\),
\eqref{eq:display-rb1} bounds \(v_\ell(c^2d)=2v+v_\ell(d)\).

For this raised family take the finite compositum of the fields of
definition of the representations \(N_y\), \(y\in Y\), the prescribed
lattices, characters, congruence polynomials and block identifications.
Adjoin the torsion values and all conjugates of the completions of
Section~\ref{sec:relative-crystalline-colength}. Use the valuation ring
\(O\) of its normal closure over \(\mathbb Q_\ell\) for the entire family. On the distinguished
lift \(\pi\) take \(L_\pi=L_0\otimes_{O_0}O\), where \(L_0\)
is the main-summand lattice over \(O_0\) carrying the selected
Selmer bound. Let \(\mathcal S_0\) be the divisible-coefficient
Selmer cocycle module for \(L_0\), normalized by \(z(\gamma_\pi)=0\). In
\eqref{eq:display-rb1} we use
\[
 \mathrm{sel}=\frac{\operatorname{length}_{O_0}\mathcal S_0}
                         {v_{O_0}(\ell)},\qquad
 \mathrm{ram}=\sum_{w\in Q}v_\ell(\psi_\eta(\sigma_w)-1),
\]
with the same representatives and tame generators. Replacing the tame generator raises its \(\ell\)-power root of unity
to an \(\ell\)-adic unit power and leaves the valuation unchanged.
Equation~\eqref{eq:specialized-selmer-coefficient-extension}
bounds the normalized length for \(L_\pi\) by \(\mathrm{sel}\).

For \(x\in\mathcal X\), let \(R_x^{\rm main}\) be the rank-\(r'\)
summand before the auxiliary twist. For the source attachment \(R\)
of Section~\ref{sec:normalized-trace-multiplicities} this summand is
\(\theta R\). On each other lift take its specified twisted lattice.
In the semisimple block decomposition add the standard Tate lattices
in the prescribed order. Thus the integral model is
\[
 u_x=T_x\oplus\bigoplus_{e=0}^{t'-1}\xi_e,\qquad
 T_x=\psi_\eta R_x^{\rm main},\quad \xi_e=\chi_{\rm cyc}^{-e}.
\]
For each \(y\in Y\setminus\mathcal X\) choose any
\(G_E\)-stable \(O\)-lattice in \(N_y\) and let \(u_y\) be
its representation. Such a lattice is the \(O\)-span of
the translates of an initial lattice. Compactness of the
continuous image bounds all translates in one scalar
multiple of the initial lattice, so this span is a full
finite free \(O\)-module and is stable. All these models
have the specified rank \(n=r'+t'\).
If \(O_*\subseteq O\) is the valuation ring
used for the original systems in the prime selection, restriction
of the normalized valuation gives
\[
  \ell^\nu O\cap O_*=\ell^\nu O_*,
  \qquad \mathfrak m_O\cap O_*=\mathfrak m_{O_*}.
\]
Hence a ramified enlargement of coefficients after raising the
level preserves the original depth-\(\nu\) congruences and
the residual determinant-unit test. The newly occurring systems
enter the trace comparison at the chosen level.

The element \(\gamma_x\) acts on \(T_x\) by \(y_x\) with
\(y_x-1\in O^\times\), and trivially on all \(\xi_e\).
The \(T_\pi(g)\) span \(\operatorname{Mat}_{r'}(O)\).
The original local conditions in
Section~\ref{sec:relative-crystalline-colength} hold on
every \(u_y\) by inclusion in \(N_y\). The auxiliary primes
are away from \(\ell\). At places over \(\ell\) the original
type projector and the original determinant conditions
are kept on the raised space, so all \(u_y\) belong to the
same class \(\mathcal Q_v\), by the inertial and de Rham
comparisons there. Since \(\psi_\eta\) is trivial on the full
decomposition group at each \(v\mid\ell\) by
Section~\ref{sec:first-auxiliary}, the integral matrices on the chosen
main lattice satisfy \(T_\pi(g)=R_\pi^{\rm main}(g)\) for every
\(g\in G_{E_v}\). Globally \(T_\pi=\psi_\eta R_\pi^{\rm main}\)
on \(L_\pi\) carries the selected Selmer bound.

Use the doubled pairs
\((\tau_j,\phi_j)\) and
\((\overline{\tau_j}^{-1},\overline{\phi_j}^{-1})\), with
the same variable \(X_j\) for the two pairs. The raw
rank-\(n\) attachments have polarization exponent \(1-n\).
The normalizing avatar of Section~\ref{sec:main-trace}
has product with its conjugate \(\chi_{\rm cyc}^{r'}\).
Consequently
\[
                 N_y^{c_E}\simeq N_y^\vee\chi_{\rm cyc}^{\,1-t'}
                           \qquad (y\in Y).
\]

Put \(\chi=\chi_{\rm cyc}^{-1}\) and
\(g^\sharp=(c_E g c_E^{-1})^{-1}\). On each \(y\in Y\),
polarization and cyclicity of trace give
\[
 \operatorname{tr}N_y((b^\sharp)^i g)
  =\chi(g)^{t'-1}\chi(b^\sharp)^{i(t'-1)}
                       \operatorname{tr}N_y(b^i g^\sharp),
\]
since the conjugate inverse of \((b^\sharp)^i g\) is
\(g^\sharp b^i\). At \(j\), let \(H_j\) be the product of the
normalized quotients \(U/U(1)\) constructed in
Section~\ref{sec:triangle-unit-input}. Thus \(H_j\in O[X]\)
and \(H_j(1)=1\), also for the empty product.
Write \(P=(X-1)H_j=\sum_i P_iX^i\).
With \(u=(u_y)_{y\in Y}\), use
\(R_j(g)=\operatorname{tr}\bigl(u((\tau_j-1-X_jP(\phi_j))g)\bigr)\)
for the tuple of trace relations on \(Y\),
and put
\[
 R_j^\sharp(g)=\operatorname{tr}
    u\bigl((\tau_j^\sharp-1-X_j P(\phi_j^\sharp))g\bigr).
\]
Thus the tame and Frobenius elements are replaced by their conjugate
inverses, with the trace tested at the same \(g\).
Tame inertia gives \(\chi(\tau_j^\sharp)=1\). The
auxiliary Frobenius condition gives
\(\chi(\phi_j^\sharp)\equiv1\bmod\ell^\nu\).
Therefore \(b_i=(\chi(\phi_j^\sharp)^{i(t'-1)}-1)/
\ell^\nu\) belongs to \(O\), and
\[
 R_j^\sharp(g)_y=\chi(g)^{t'-1}
   \left(R_j(g^\sharp)_y-X_j\sum_i P_i b_i\ell^\nu
                           \operatorname{tr}N_y(\phi_j^i g^\sharp)\right).
\]
The terms in parentheses belong to the original trace
ideal by its ordinary congruences and its \(\ell^\nu\)-multiples
of traces. The coefficients are common to \(Y\).
The same trace identities hold after changing to the
integral realizations \(u_y\).

Since the auxiliary Frobenius fixes the \(\ell^\nu\)-th roots of
unity, its value under each cyclotomic power \(\xi_e\) is
\(1\bmod\ell^\nu O\), also at the conjugate inverse. By integrality
and \(H_j(1)=1\), this gives
\(H_j(\xi_e(\phi))-1\in\ell^\nu O\) on both sides.
Take \(a=g_{\rm tr}=\gamma_\ell^{-1}\) from
Section~\ref{sec:triangle-unit-input}. Its values are
\(\xi_e(a)=\xi_\ell^e\). The cyclotomic image has finite
index in \(\mathbb Z_\ell^\times\) since \(E/\mathbb Q\) is
finite, so its topological generator has infinite order.
The values are pairwise distinct. By Section~\ref{sec:triangle-unit-input},
each \(\xi_\ell^b-1\) for \(0<b<t\) divides \(\ell^h\).
Each difference below is a unit multiple of such a factor.
Since \(t'-1\le H\), its product satisfies
\[
 \prod_{e'\ne e}\bigl(\xi_e(a)-\xi_{e'}(a)\bigr)\mid C
 \qquad(0\le e<t').
\]
Use the tame elements \(\delta_w\) at the places of \(S_M\)
away from \(\ell\) and at the twist support \(Q\), together with their
conjugate inverses, as in Section~\ref{sec:coefficient-old-tame-relations}.
The cyclotomic character is unramified at these places, so
\(\xi_e(\delta_w)=1\) for every \(e\). The finite types give
\((u_y(\delta_w)-1)o_w(u_y(\delta_w))=0\) for every \(y\in Y\).
At an old place over \(\ell\), Section~\ref{sec:relative-crystalline-colength}
uses the inertial polynomials on the potentially crystalline periods.
Let \(s'\) consist of the test set \(P_0\) of
Section~\ref{sec:triangle-selection}, the reference \(a\), and
the two tame generators at each new prime pair. The triangle
gives, with \(n_{\rm test}\) counting the tests and \(a\),
\[
 (|J|+n_{\rm test})Hh\le H^{22}h,\qquad
 n_{\rm det}:=|s'|\le 2(|J|+n_{\rm test}).
\]

\subsubsection{The main-block rows}
\label{sec:coefficient-main-block-covectors}

\sectionopening{In this subsection, we use the trace relations above and the congruence of Section~\ref{sec:triangle-unit-input} to select the main-block rows for Section~\ref{sec:trace-determinant-specialization}.}

Put \(\Lambda=O[X_j:j\in J]\) and \(A=\Lambda[\mathrm{Gal}(\overline E/E)]\).
Use one finite list \(g_i\), \(i\in I\), generating the
augmentation module
\[
 \mathcal M=\operatorname{span}_O\{(u_y(g)-1)_{y\in Y}:g\in G_E\}
                   \subseteq\prod_{y\in Y}\operatorname{Mat}_n(O).
\]
This module is finitely generated over the DVR \(O\).
Include the old tame elements and \(s'\) in the list.
Let \(z(g)\in O^I\) express \((u_y(g)-1)_{y\in Y}\)
in these generators. Extend each \(u_y\) to \(A\) by \(\Lambda\)-linearity, and put
\(I_{\rm aug}=\ker(A\to\Lambda)\) for the augmentation ideal.
Write \((\tau_{j,+},\phi_{j,+})=(\tau_j,\phi_j)\) and
\((\tau_{j,-},\phi_{j,-})=(\tau_j^\sharp,\phi_j^\sharp)\).
For \(\epsilon\in\{+,-\}\), put
\[
 \begin{aligned}
 b_{j,\epsilon}&=\tau_{j,\epsilon}-1
       -X_jH_j(\phi_{j,\epsilon})(\phi_{j,\epsilon}-1),\\
 \varepsilon_\pi&=(y_\pi-1)^{-1}(\gamma_\pi-1),&
 \rho_0&=\varepsilon_\pi+(1-\varepsilon_\pi)(a-1).
 \end{aligned}
\]
Define the left ideal \(\mathfrak b\subset A\) and its module of rows by
\[
 \begin{aligned}
 \mathfrak b&=\Bigl(\bigcap_{y\in Y}\ker u_y\Bigr)+\ell^\nu I_{\rm aug}
              +\sum_{j,\epsilon} A b_{j,\epsilon}+I_{\rm aug}\rho_0,\\
 \mathcal N&=\Bigl\{r\in A^I:\ \sum_i r_i(g_i-1)\in\mathfrak b\Bigr\}.
 \end{aligned}
\]
Let \(\psi_0\) be evaluation at \(X_j=0\) followed by the
\(T_\pi\) action.

The spanning hypothesis makes \(\psi_0\)
onto \(\operatorname{Mat}_{r'}(O)\).
The zero-variable image of the row
\(z(\tau)-X_jH_j(\phi)z(\phi)\) is \(z(\tau)\) acting
by scalar entries. The row
\[
 (\gamma_\pi-1)\bigl((y_\pi-1)^{-1}z(\gamma_\pi)
                                  +(1-\varepsilon_\pi)z(a)\bigr)
\]
satisfies the congruences and has image \(z(\gamma_\pi)\) with scalar
entries, since \(\psi_0(\varepsilon_\pi)=1\). These statements
about values are read modulo the common kernel.
The same kernel contains the values of the cocycle vectors
\(z(gg')-z(g)-g z(g')\), the coordinate rows
\(z(g_i)-\mathbf e_i\), and the rows of all linear
relations among \((u_y(g)-1)_{y\in Y}\).

Let \(L\subset O^{I\times r'}\) be the span of the individual
scalar rows of \(\psi_0(\mathcal N)\). For
\(w\in(K/O)^{I\times r'}\) annihilating \(L\), \(K=\operatorname{Frac}(O)\),
put \(c_w(g)=\sum_i z_i(g)w_i\). The displayed relations
give the cocycle rule for \(T_\pi\), the simultaneous
linear relations, and \(c_w(\tau)=c_w(\gamma_\pi)=0\).
Also \(c_w(g_i)=w_i\), so this map is injective.
Put \(a_\nu=\ell^\nu\). The rows \(a_\nu\mathbf e_i\) belong to
\(\mathcal N\), so \(a_\nu w_i=0\) for every \(i\).
Over the fixed DVR \(O\), multiplication by \(a_\nu\) induces
\[
 (K/O)[a_\nu]\xrightarrow{\sim}O/a_\nu O,\qquad
           x+O\longmapsto a_\nu x+a_\nu O .
\]
Write \(\widetilde c_w\) for the corresponding cocycle in
\(L_\pi/a_\nu L_\pi\). By the simultaneous linear relations its
evaluation on the augmentation tuple factors through
\(\mathcal M/a_\nu\mathcal M\). In a Smith basis of
\(W_{\rm mat}=\prod_{y\in Y}\operatorname{Mat}_n(O)\), write
\(\mathcal M=\bigoplus_i \varpi^{s_i}Oe_i\). For
\(n_0\ge v_O(a_\nu)+\max(\{0\}\cup\{s_i\})\) one has
\(\mathcal M\cap\varpi^{n_0}W_{\rm mat}\subseteq a_\nu\mathcal M\).
Continuity of the tuple of original representations therefore gives
continuity of \(\widetilde c_w\) and \(c_w\).

Sections~\ref{sec:triangle-selection} and~\ref{sec:raised-tame-frobenius}
choose the triangle places \(J\cup\bar J\) off \(\ell\), the
original level, and \(Q\cup\bar Q\) for the chosen twist.
Thus \(T_\pi\) is unramified at both members of every pair.
Its cocycle on each inertia group is a continuous additive
homomorphism into the finite \(\ell\)-primary torsion module.
If \(p\ne\ell\) is the residue characteristic, wild inertia is
pro-\(p\) and maps to zero, as do the prime-to-\(\ell\) factors of tame inertia.
The two equations \(c_w(\tau_{j,\epsilon})=0\) test generators of
the remaining pro-\(\ell\) quotient at the respective places.
Hence the inertia restrictions at these places vanish. Off the
prescribed support all the tuple differences vanish on inertia.
Let \(D_\nu\) be the image of \(w\mapsto\widetilde c_w\) in the
relative cocycles for this tuple. With \(A_{\nu,v}\) the reduction of
the integral crystalline cocycles, put
\(D_{\nu,\mathrm{red}}=D_\nu\cap
\bigcap_{v\mid\ell}\operatorname{res}_v^{-1}(A_{\nu,v})\).
Restriction induces
\[
 D_\nu/D_{\nu,\mathrm{red}}\hookrightarrow
 \bigoplus_{v\mid\ell} S_{\nu,v}/(S_{\nu,v}\cap A_{\nu,v}).
\]
The local bounds in the proof of \eqref{eq:relative-crystalline-colength}
give \(\operatorname{length}_O(D_\nu/D_{\nu,\mathrm{red}})
\le\sum_{v\mid\ell}\Delta_v(T_\pi)\).
For \(z\in D_{\nu,\mathrm{red}}\), choose at each \(v\mid\ell\) an
integral crystalline cocycle \(c_v\) lifting \(\operatorname{res}_v z\).
Under the inverse coefficient identification
\[
 L_\pi/a_\nu L_\pi\hookrightarrow L_\pi[1/\ell]/L_\pi,\qquad
            y+a_\nu L_\pi\longmapsto a_\nu^{-1}y+L_\pi,
\]
its local restriction is \(a_\nu^{-1}c_v\bmod L_\pi\) and satisfies
the homogeneous crystalline relations. The inertia tests and
\(z(\gamma_\pi)=0\) are preserved. This gives an injection of
\(D_{\nu,\mathrm{red}}\) into \(\mathcal S_O\) on actual cocycles.
Since evaluation at the \(g_i\) recovers every \(w_i\), the annihilator
is isomorphic to \(D_\nu\).

At the fixed depth \(\nu\) and over the coefficient ring \(O\)
of this raised family, the ordinary \(O\)-length of the
annihilator is at most
\[
 b_T=\operatorname{length}_O\mathcal S_O
               +\sum_{v\mid\ell}\Delta_v(T_\pi),
\]
where \(\mathcal S_O\) is the Selmer module with divisible
coefficients and \(z(\gamma_\pi)=0\).
For each depth and corresponding raised family, the coefficient
extension bound and the main-block estimate of
Section~\ref{sec:coefficient-local-defects} give
\[
 \frac{b_T}{v_O(\ell)}\le\mathrm{sel}+2fr'F_\ell ,
\]
using \(\sum_{v\mid\ell}[E_v:\mathbf Q_\ell]=2f\).
The injection is applied at each such depth. The common twist and
its Selmer bound were chosen at the earlier cutoff \(M\).
Smith form identifies this annihilator
with \(\operatorname{Hom}_O(O^{I\times r'}/L,K/O)\).
A free summand in the quotient would give a \(K/O\)
summand of infinite length in the annihilator. For a
summand \(O/(\varpi^a)\) the two lengths are both \(a\).
It follows that \(L\) has full rank and colength
at most \(b_T\).

Choose an \(O\)-basis of \(L\). Each basis vector is an \(O\)-linear
combination of scalar rows of the matrices \(\psi_0(r)\), \(r\in\mathcal N\).
To put it in row \(c\) of a block row, multiply each contributing \(r\)
on the left by a lift of the corresponding matrix unit \(E_{c,d}\),
where \(d\) indexes the chosen scalar row. Sum with the given coefficients.
The lifts exist by surjectivity of \(\psi_0\), and
\(\mathcal N\) is a left \(A\)-submodule. Group the basis vectors into
\(|I|\) blocks of \(r'\) rows. This gives \(B_0\in M_I(A)\) with
rows in \(\mathcal N\) and
\[
 v_O\bigl(\det\psi_0(B_0)\bigr)
       =\operatorname{length}_O(O^{I\times r'}/L)\le b_T .
 \tag{RB2}\label{eq:display-rb2}
\]
Thus the whole matrix, rather than a multiple of its
determinant by a rank-dependent denominator, realizes
the cocycle length bound.

\smallskip\noindent
\textit{The character matrices and the top coefficient.}
Order the generator coordinates as the disjoint blocks
\(I_{\rm old}\sqcup I_{\rm aux}\sqcup s'\). The first
block consists of the old tame tests and the last
of the test set. Adjoin a finite list \(I_{\rm aux}\)
of generators as needed to span the simultaneous
augmentation tuples. The coordinate relation
\(z(g_i)-\mathbf e_i\) remains valid if a generator
is repeated. Let \(\eta_e:A\to\Lambda\) evaluate each
group element by \(\xi_e\), with the polynomial
variables fixed. We need rows satisfying the congruences \(U_i\in\mathcal N\),
\(i\in I_{\rm aux}\), whose \(I_{\rm aux}\)-parts are
the respective unit vectors under every \(\eta_e\).

First fix \(e\). Let \(P_e\) be the image of
\(\mathcal N\) under \(\eta_e\) restricted to the
\(I_{\rm aux}\) coordinates. It contains
\(\ell^\nu\Lambda^{I_{\rm aux}}\). Put
\(B_e=O^{I_{\rm aux}}/(P_e\cap O^{I_{\rm aux}})\).
For a functional \(b:B_e\to O/\ell^\nu O\), extend it
by zero on the two other constant coordinate blocks.
The cocycle, simultaneous linear, coordinate and normalizer
rows have \(\eta_e\)-values in \(O^I\). Their auxiliary parts
therefore belong to \(P_e\cap O^{I_{\rm aux}}\), so \(b\) kills
these rows. In particular \(c(g)=b(z(g))\) satisfies the
cocycle and simultaneous linear relations. The Smith form
argument for the simultaneous tuple module gives continuity.
The coordinate relations give zero at all old tame and test points,
including both tame generators of each new pair.
The vectors \((g-1)((y_\pi-1)^{-1}z(\gamma_\pi)
+(1-\varepsilon_\pi)z(a))\) give
\[
   (\xi_e(g)-1)((y_\pi-1)^{-1}c(\gamma_\pi)+c(a))=0 .
\]
These are precisely the conditions for the injectivity criterion
in Section~\ref{sec:triangle-character-cost}, so
\(c=0\). Evaluating at the coordinate generators
gives \(b=0\). Since \(B_e\) is a finite-length
\(O\)-module killed by \(\ell^\nu\), finite DVR
duality gives \(B_e=0\). Thus each unit vector of
\(I_{\rm aux}\) has a lift \(U_{e,i}\in\mathcal N\)
under this evaluation.

For simultaneous lifting, start with \(U_i=0\).
At the next character \(e\) replace \(U_i\) by
\[
 U_i+\sum_{j\in I_{\rm aux}}
               (\mathbf1_{j=i}-(U_i)_j)\,U_{e,j}.
\]
The coefficients multiply on the left in \(A\), so the
result still satisfies the congruences. Its evaluation on \(I_{\rm aux}\) is the
unit vector at \(e\). At every previously processed
character all coefficients of the correction evaluate
to zero. Induction gives \(U_i\) as required. This
construction is integral.

For each \(e\) eliminate the auxiliary block of a
row \(w\) satisfying the congruences by subtracting
\(\sum_{i\in I_{\rm aux}} w_i U_i\) before evaluating.
This operation makes the auxiliary part zero at \(e\).
It defines a test-coordinate image in
\(\Lambda^{s'}\), the same image used in the next
subsection. For \(J_e=\{j:\mathrm{ty}(j)=e\}\), its
Fitting calculation gives
determinant sums \(d_e(X)\) and polynomials \(\Gamma_e,\zeta_e\)
such that
\[
 d_e+\zeta_e=\Gamma_e\prod_{j\in J_e}X_j,\qquad
 \zeta_e\prod_{\mathrm{ty}(j)<e}X_j\in(X_j^2:j\in J),\qquad
 v_O(\Gamma_e(0))\le b_e+(|J_e|+n_{\rm det})v_O(C).
 \tag{RB3}\label{eq:display-rb3}
\]
Write each \(d_e\) as a weighted sum of test-block
determinants. Each is formed from rows that lift
to the congruence module after elimination at \(e\).
Choose one such matrix \(B_e^{\rm test}\) for each \(e\)
and denote the eliminated lifts by \(w_{e,l}\), \(l\in s'\).
In the common test block use
\(\sum_e q_e^{\rm sep} w_{e,l}\), with the separation element
\[
 q_e^{\rm sep}=\frac{C}{\prod_{e'\ne e}(\xi_e(a)-\xi_{e'}(a))}
                    \prod_{e'\ne e}(a-\xi_{e'}(a)).
\]

The values \(\xi_e(a)\) are distinct in \(\operatorname{Frac}(O)\).
Their difference product divides \(C\) in \(O\) by
Section~\ref{sec:triangle-unit-input}. Form the displayed quotient
in \(O\), before reducing at depth \(\nu\). Thus \(q_e^{\rm sep}\in A\) and
\(\eta_{e'}(q_e^{\rm sep})=C\mathbf1_{e=e'}\), also when two characters
have the same reduction.
Use the common \(U_i\) in the auxiliary block.
For \(w\in I_{\rm old}\), use the row supported in
that coordinate with entry \(o_w(\delta_w)\).
The congruences follow since \((\delta_w-1)o_w(\delta_w)\)
annihilates every system. In the order
\(I_{\rm old},I_{\rm aux},s'\), evaluation at \(\eta_e\) gives
\[
 \begin{pmatrix}
  \operatorname{diag}_{w\in I_{\rm old}}(o_w(1)) & 0 & 0\\
  * & I_{I_{\rm aux}} & *\\
  * & 0 & C B_e^{\rm test}
 \end{pmatrix}.
\]
Elimination with the common \(U_i\) gives the zero in the lower
middle block. Expansion along the old coordinates gives determinant
\(c_0\det B_e^{\rm test}\), where
\(c_0=(\prod_w o_w(1))C^{n_{\rm det}}\).
In particular the separation factor is used only on
the test block. This establishes the simultaneous
specialization for each tuple of determinant summands.

For a fixed tuple of test-block determinants let \(M\) be the
corresponding simultaneous character matrix. In the integral group
algebra \(A\) put
\[
 k_{\rm char}=\prod_{x\in\mathcal X}^{\rm ord}
                 \frac{y_x-\gamma_x}{y_x-1},\qquad
 B=(1-k_{\rm char})B_0+k_{\rm char}M.
\]
Each denominator is a unit. The factor indexed by \(x\) acts by
zero on \(T_x\), and every factor acts by one on each \(\xi_e\).
The multiplications are on the left of the entries, so the rows
of \(B\) remain in the left \(A\)-module \(\mathcal N\).
The main and character blocks of its evaluations are \(T_x(B_0)\)
and \(\eta_e(M)\), respectively. Use this same \(B\) for the trace
determinant on the full family \(Y\). Sum over the tuples
of chosen test-block determinants with the product
of their weights from the individual character sums. The common trace
determinant construction applied to these congruence matrices gives a single multiplier \(f\) for all
\(\det(1-u_x(g))\), with specialization on \(\mathcal X\)
\[
 f_x(X)=C\,c_0^{t'}D_x(X)\prod_e d_e(X),
\]
where \(D_x\) is the block determinant. For the fixed \(M\) the trace calculation of
Section~\ref{sec:marked-ribet-trace-determinant} uses
\(q_0^{\rm sep} k_{\rm char}\). On each lift this element acts by \(C\) on
\(\xi_0\) and by zero on the main block and all other characters.
Since \(\xi_0(\rho_0)=0\), the trace of this element times the
Ribet polynomial evaluated at \(\rho_0\) is \(C\) times its
constant term. This gives the
displayed block product.

Replace the factors \(d_e\) in increasing order by
\(d_e+\zeta_e\). At each replacement the difference
is a multiple of
\(\zeta_e\prod_{\mathrm{ty}(j)<e}X_j\). It is zero
modulo \((X_j^2:j\in J)\). Hence on \(\mathcal X\) all the
coefficients of proper divisors of \(X_J\) vanish and
\[
 d=Cc_0^{t'} D_\pi(0)\prod_e\Gamma_e(0).
\]
Write \(L_{0,e}\) and \(b_{A,e}\) for \(L_0\) and \(b_A\)
in the Fitting calculation for \(\xi_e\) below, and put
\[
 m_T=\operatorname{length}_O(O^{I\times r'}/L),\qquad
 m_e=\operatorname{length}_O(O^{s'}/L_{0,e}).
\]
Since \(\Gamma_e(0)=\varpi^{b_{A,e}}C^{|J_e|}d_T(0)\) with
\(v_O(d_T(0))=m_e\) in that calculation, the product gives
\[
\begin{aligned}
 v_O(d)={}&m_T+\sum_e(b_{A,e}+m_e)+t'\sum_w v_O(o_w(1))\\
          &+(1+t'n_{\rm det}+|J|)v_O(C).
\end{aligned}
\]
Here \(\sum_e |J_e|=|J|\), including when some \(J_e\) is empty.
The estimates \(m_T\le b_T\) and
\(b_{A,e}+m_e\le b_e+n_{\rm det}v_O(C)\) now give
\[
 v_O(d)\le b_T+\sum_e b_e+
 t'\left(\sum_w v_O(o_w(1))+2n_{\rm det} v_O(C)\right)
 +(|J|+1)v_O(C).
 \tag{RB4}\label{eq:display-rb4}
\]

\subsubsection{The old tame relations in the coefficient comparison}\label{sec:coefficient-old-tame-relations}

\sectionopening{In this subsection, we use the tame operators of Section~\ref{sec:raised-tame-frobenius} to bound the remaining coefficient contributions.}

Fix the coefficient prime \(\ell\), the chosen embedding, and the
prescribed twist at this embedding. Recall \(C_G=A_{\rm aut}\) and
\(v=v_\ell(C_G)\). Use all the places above \(M\) prime to
\(\ell\), together with the twist support \(Q\). At an old place
above \(p_i\), the normalized upper representations have finite
inertia and every inertia eigenvalue has order dividing
\(p_i^a-1\) for some \(1\le a\le r\). The conjugate-inverse
element has the same bound, by the local types on the
conjugate place. The integer \(p_i^a-1\) divides the order of
the corresponding general linear group in \(C_G\).
Finite order in characteristic zero implies semisimplicity.
Thus, on every upper representation and for either choice of
inertia element, \(X^{C_G}-1\) annihilates its action.
Take \(o_w(X)=\sum_{j=0}^{C_G-1}X^j\), so \(o_w(1)=C_G\).

The old local projector is imposed on the whole raised
space. At a twist place the prescribed two determinant
types give the quadratic annihilator of
Section~\ref{sec:inertial-labels} on every raised system.
Choose the tame generator \(\sigma_w\) used in the ramification
bound. For \(g=\sigma_w\) or its conjugate inverse the
quadratic relation permits \(o_g(X)=\psi_\eta(g)-X\). The
anticyclotomic identity gives
\(\psi_\eta((c_Egc_E^{-1})^{-1})=\psi_\eta(g)\). Hence the two copies
give the same valuation \(v_\ell(\psi_\eta(\sigma_w)-1)\).
These annihilations pass to all the fixed integral models
by their injections into the representations over the
coefficient field. There are at most \(4f\) old places above
the two rational level primes and at most twice that many
old inertia elements in the list. Thus on this same family
\[
 \sum_g v_\ell(o_g(1))\le 2\mathrm{ram}+8f v .
\]
Places above \(\ell\) omitted here are handled by the
crystalline colength bound. The auxiliary triangle uses
separate places outside both the old and the twist level.

Divide \eqref{eq:display-rb4} by \(v_O(\ell)\). In normalized units,
\[
 \frac{b_T}{v_O(\ell)}=\mathrm{loss}_T+\mathrm{sel}_T,
 \qquad
 \frac{b_e}{v_O(\ell)}=\mathrm{loss}_e+\mathrm{odd}_e .
\]
The local defect estimate bounds the normalized sum of all losses
by
\[
 \sum_{v\mid\ell}[E_v:\mathbf Q_\ell](r'+t'/t)F_\ell
     =2f(r'+t'/t)F_\ell\le4fr'F_\ell
     =r'\mathrm{loc}_\ell .
\]
Here \(r'\ge1\) and \(t'\le t\). The normalized
\(\mathrm{sel}_T\) is at most \(\mathrm{sel}\);
the normalized \(\mathrm{odd}_e\) at most
\(v_\ell(X_{\rm odd})\); and the normalized
old tame sum at most \(2\mathrm{ram}+8fv\).
Finally the two triangle inequalities give
\[
 (2t'n_{\rm det}+|J|+1)Hh
 \le (4t+1)(|J|+n_{\rm test})Hh+Hh
 \le \bigl((4t+1)H^{22}+H\bigr)v .
\]
Adding \(2v\) to \eqref{eq:display-rb4} proves \eqref{eq:display-rb1}.
The next two subsections establish the character Fitting
estimate and construct \(f\) on \(Y\). The annihilation
of \(f\) by the detecting determinants modulo the trace
ideal is proved on \(Y\), and the coefficient bound is
computed by the block evaluation at \(\pi\).
The quantitative bound uses the local defect bounds
of Section~\ref{sec:relative-crystalline-colength}
and the Selmer bound on the prescribed representation.

\section{Fitting ideals and trace determinant specialization}
\label{sec:trace-determinant-specialization}

\sectionopening{In this section, we specialize the common coefficient of
Section~\ref{sec:hecke-trace-comparison} on the chosen characters.
The resulting Fitting inclusion and trace determinant give the
finite coefficient used with the norm tests in
Section~\ref{sec:common-critical-twist}.}

\subsection{The Fitting element for the ordered characters}

\sectionopening{In this subsection, the specialization on a character summand of
Section~\ref{sec:marked-ribet-coefficient} reduces to an inclusion
between Fitting ideals. We derive it by eliminating the unit relations
and applying the bounds on the paired evaluations of
Section~\ref{sec:triangle-character-cost}.}

Fix the character \(\xi_e\) in the common Ribet
construction and write \(I=J_e\), \(s=s'\),
\(R=O[X_j:j\in J]\), \(p=\ell^\nu\) and \(F=R^s\).
Write \(\mathfrak q=(X_j^2:j\in J)\) and
\(m_{<e}=\prod_{\mathrm{ty}(j)<e}X_j\), with empty product \(1\).
Let \(M_e\subset F\) be the image obtained by evaluating the congruence
covectors at \(\xi_e\), eliminating the auxiliary coordinates with the
unit covectors, and retaining the test coordinates. Put
\[
                  N=M_e+(\mathfrak q:m_{<e})F .
\]
Write \(\Phi:O^{I_{\rm coord}}\to F\) for the
inclusion of constant coordinate covectors followed
by this evaluated elimination and projection to the
test block. Here \(I_{\rm coord}\) is the full list
of the common augmentation generators.

Write \(y(g)=\Phi(z(g))\).
The cocycle vectors show that the forms on \(O^{I_{\rm coord}}\)
modulo \(p\) which annihilate \(\Phi^{-1}N\)
give the character cocycles by \(z\).
The coordinate relations \(z(g_i)-\mathbf e_i\)
show that \(z\) recovers each such form from its
cocycle. On the test positions \(g\in s\) one has
\(y(g)\equiv\mathbf e_g\pmod N\).
The two congruences are \(y(\tau_j^\pm)\equiv
 X_j y(\phi_j^\pm)\pmod N\), because
\(H_j(\xi_e(\phi_j^\pm))\equiv1\pmod p\) and
\(pF\subset N\). Indeed the \(p\)-multiples of the pure test rows
are congruence covectors and have zero auxiliary coordinates.
For an earlier \(i<j\) with \(j\in I\),
either \(i\in I\) or \(\mathrm{ty}(i)<e\); in the latter
case \(X_iF\subset N\) by the colon ideal.

For each \(R\)-submodule \(L\) containing \(N\), put
\[
 A_L=O^{I_{\rm coord}}/\Phi^{-1}L,\qquad
 Z_L=\operatorname{Hom}_O(A_L,O/(p)).
\]
The finite \(O\)-module \(A_L\) is killed by \(p\).
A functional \(\lambda\in Z_L\) gives the cocycle
\(c_\lambda(g)=\lambda([z(g)])\). The relations at \(g_i\) give
\(c_\lambda(g_i)=\lambda([\mathbf e_i])\), so this identifies \(Z_L\)
with its image in the character cocycles. We use the following determinant
operation throughout. For \(w\in F\), if \(aw\in L\) for
\(a\in R\), then \(a\,\operatorname{Fitt}_0(F/(L+Rw))
\subseteq\operatorname{Fitt}_0(F/L)\).
Expand a determinant of covectors in \(L+Rw\).
Terms with two or more covectors proportional to \(w\)
are zero, and in each remaining term the factor \(a\)
can be placed in the \(w\)-covector if it is present. The
Fitting ideal here is generated by all determinants of
\(|s|\) covectors, so this proves the inclusion.

For each old tame generator \(\delta=g_i\),
\(\Phi(z(\delta)-\mathbf e_i)\in N\) and
\(\Phi(\mathbf e_i)=0\). Thus \(c_\lambda(\delta)=0\)
for every \(\lambda\in Z_L\). The simultaneous linear relations
make \(c_\lambda\) factor through \(\mathcal M/p\mathcal M\).
Continuity follows from the Smith basis of
Section~\ref{sec:relative-crystalline-colength}.
The row \((g-1)v_0\) lies in \(\mathcal N\). With
\(a=g_{\rm tr}\), its evaluation at \(\xi_e\) is the constant row
\[
 m_g=(\xi_e(g)-1)\bigl(c_\pi z(\gamma_\pi)+z(a)\bigr),
 \qquad \xi_e(\gamma_\pi)=1 .
\]
Eliminating the auxiliary coordinates and retaining the test
coordinates gives \(\Phi(m_g)\in N\subseteq L\). Applying
\(\lambda\) gives
\((\xi_e(g)-1)(c_\pi c_\lambda(\gamma_\pi)+c_\lambda(a))=0\).
In particular the required reference relation holds for all
\(\lambda\in Z_L\) and we may take \(Z=Z_N\) in
Section~\ref{sec:triangle-character-cost}.
Let \(A\subset G\) be the union of the inertia groups at
\(v\mid\ell\) if \(e=0\), and of
\(D_v\cap\ker\chi_{\rm cyc}\) at \(v\mid\ell\) if \(e>0\).
Write \(Z^0=Z_0\) for the submodule of \(Z_N\) that vanishes on \(A\). The \(O\)-submodule
of \(A_N\) spanned by the \(z(a)\), \(a\in A\),
has length \(b_A=\operatorname{length}(Z_N/Z^0)\),
by finite DVR duality. In fact \(O/(p)\) identifies
with the \(p\)-torsion of \(K/O\), and duality with
\(K/O\) is exact on finite \(O\)-modules.
Choose cyclic direct summands of that submodule, and
lift their generators to constant vectors \(w_i\) in the
preimage of the submodule in \(O^{I_{\rm coord}}\).
If the summand has length \(a_i\),
\(\varpi^{a_i}\Phi(w_i)\in N\), where \(\varpi\) is
a uniformizer. Set \(K_1=N+\sum_iR\Phi(w_i)\).
Applying the determinant operation to the successive
augmentations gives
\[
 \varpi^{b_A}\operatorname{Fitt}_0(F/K_1)
     \subseteq \operatorname{Fitt}_0(F/N).
\]
The generators \(y(a)\), \(a\in A\), belong to \(K_1\).
Consequently \(Z_{L}\subseteq Z^0\) whenever \(L\supset K_1\).

For \(j\in I\) set \(v_j=y(\phi_j^+)-y(\phi_j^-)\).
Starting at \(K_1\), add \(Rv_j\) in increasing order of
\(j\). Immediately before adding \(v_j\), call the
current module \(L\). Every \(c\in Z_L\) vanishes
on the differences of the earlier tame generators:
for \(i\in I\) this follows from \(v_i\in L\) and the
two congruences, and for \(i\notin I\) from \(X_iF\subset N\).
It also vanishes on the earlier Frobenius differences
of type \(e\), by \(v_i\in L\). The triangle condition
therefore gives \(C(c(\tau_j^+)-c(\tau_j^-))=0\).
The finite dual \(Z_L\) separates classes of \(A_L\).
Applied to the constant vector \(C(z(\tau_j^+)-
z(\tau_j^-))\), this gives
\[
 C\bigl(y(\tau_j^+)-y(\tau_j^-)\bigr)\in L,
 \qquad CX_jv_j\in L.
\]
After all augmentations, let \(T\) be the resulting
module. The determinant operation gives
\[
 C^{|I|}\left(\prod_{j\in I}X_j\right)\operatorname{Fitt}_0(F/T)
   \subseteq \operatorname{Fitt}_0(F/K_1).
\]

Let \(L_0\subset O^s\) be the image of \(T\) at \(X=0\).
It contains \(pO^s\). The evaluation map \(v\mapsto\Phi(v)(0)\)
surjects onto \(O^s/L_0\), since the images of \(z(g)\), \(g\in s\),
are the test basis modulo \(L_0\). Its kernel modulo \(L_0\) is
\(\Phi^{-1}\bigl(T+(X_j:j\in J)F\bigr)\), because
\(L_0\) is the image of \(T\) at zero. Thus
\(Z_1=Z_{T+(X_j:j\in J)F}\) is the dual of \(O^s/L_0\). Every member of \(Z_1\) is
in \(Z^0\), vanishes at both tame generators for each
\(j\) by the two congruences specialized at zero, and vanishes
on the Frobenius differences for \(j\in I\) by \(v_j\in T\).

Thus \(Z_1\subseteq S\) for the module \(S\) of
Section~\ref{sec:triangle-character-cost}, and
\(Z_1/Z_1[C]\) injects into \(S/S[C]\).
Our choice \(Z=Z_N\) gives \(b_A=\operatorname{length}_O(Z/Z_0)\). Equation
\eqref{eq:triangle-character-cost} therefore gives
\[
 b_A+\operatorname{length}_O(Z_1/Z_1[C])
             \le c_{\rm crys}(O(-e),\nu)+v_{\rm ord}(X_{\rm odd})=b_e .
\]
Put \(e_O=v_O(\ell)\). The modulus \(p=\ell^\nu\) generates
\((\varpi^{e_O\nu})\) in \(O\). Since \(Z_1\subset (O/(p))^s\)
via the test basis, the contribution of its \(C\)-torsion is
\[
 \operatorname{length}_O Z_1[C]\le
 |s|\min\{v_O(C),e_O\nu\}\le |s|v_O(C).
\]
Choose a basis of \(L_0\) and lift its elements to \(T\).
Their determinant \(d_T\in\operatorname{Fitt}_0(F/T)\)
satisfies
\[
 v_O(d_T(0))
   =\operatorname{length}(O^s/L_0)
   =\operatorname{length}(Z_1)
   \le b_e-b_A+|s|v_O(C).
\]
Take \(\Gamma_e=\varpi^{b_A}C^{|I|}d_T\).
The two determinant inclusions prove
\[
 \Gamma_e\prod_{j\in I}X_j\in\operatorname{Fitt}_0(F/N),
 \qquad
 v_O(\Gamma_e(0))\le b_e+(|I|+|s|)v_O(C).
\]
In a determinant with rows in \(N=M_e+(\mathfrak q:m_{<e})F\),
every expansion term containing a row in \((\mathfrak q:m_{<e})F\)
belongs to \((\mathfrak q:m_{<e})\). Hence
\[
 \operatorname{Fitt}_0(F/N)\subseteq
     \operatorname{Fitt}_0(F/M_e)+(\mathfrak q:m_{<e}).
\]
Writing \(\Gamma_e\prod_{j\in I}X_j=d_e+\zeta_e\) accordingly gives
a polynomial combination \(d_e\) of determinants from \(M_e\) and
\(m_{<e}\zeta_e\in\mathfrak q\), as in \eqref{eq:display-rb3}.

\subsection{The trace determinant under the Frobenius relations}
\label{sec:marked-ribet-trace-determinant}

\sectionopening{In this subsection, we use the trace pairing on the algebra of
Section~\ref{sec:marked-ribet-coefficient} to obtain a determinant
over its coefficient ring. Its annihilation property holds on the full
family of representations and permits specialization at the chosen
system.}

Use the coefficient DVR \(O\) of Section~\ref{sec:marked-ribet-coefficient}
and the full raised family \(Y\) of Section~\ref{sec:integral-raised-traces}.
The diagonal copy of \(O\) preserves the raised lattice. Put
\[
 \mathbb T_O=\mathbb T\cap O^Y,\qquad
 S=\mathbb T_O[X_j:j\in J]\hookrightarrow
 L=\prod_{y\in Y}O[X_j:j\in J].
\]
The map records each coefficient at each system. For every \(g\in G_E\)
the characteristic coefficients of \((u_y(g))_{y\in Y}\) lie in
\(\mathbb T\) by the Hecke comparison and in \(O^Y\) by the stable
integral models. The local polynomials \(H_j\) have coefficients in \(O\).
Projection to any \(y\in Y\) fixes the diagonal \(O[X_j:j\in J]\).
In this trace calculation write \(G=G_E\), and extend the tuple
of representations to \(u:S[G]\to\operatorname{Mat}_n(L)\).
Define the scalar ideal and its trace annihilator by
\[
 \mathcal J=\bigl(\ell^\nu\operatorname{tr}u(g):g\in G\bigr)_S
     +\bigl(R_j(g):j\in J,\ g\in G\bigr)_S,
\]
\[
 \mathfrak a=\{w\in S[G]:\operatorname{tr}(u(wv))
                      \in\mathcal J\text{ for every }v\in S[G]\}.
\]
Explicitly the local relations are
\[
 R_j(g)=\operatorname{tr}u\bigl((\tau_j-1-X_jP_j(\phi_j))g\bigr),
 \qquad P_j(U)=(U-1)H_j(U) \quad (g\in G,\ j\in J).
\]
The relations
\(R_j^\sharp(g)\) also lie in \(\mathcal J\) by the computation
in Section~\ref{sec:marked-ribet-coefficient}.
The trace condition defining \(\mathfrak a\subset S[G_E]\)
is two-sided by cyclicity of the trace,
and it contains \(\ker u\).
Write \(y_i=g_i-1\), and let
\(\mathcal I=\sum_i S[G]y_i\) be the left ideal they generate.

Let \(m=|I_{\rm coord}|\), and let \(B\) be an \(m\)-by-\(m\)
matrix over \(S[G]\) with rows given by the congruence covectors.
View \(y=(g_i-1)_i\) as a column. The trace relations give
\[
 By\equiv B_1y\,\rho_0\pmod{\mathfrak a}
\]
for a matrix \(B_1\) over \(S[G]\).
Use the coefficient inclusion
\(O[X_j:j\in J]\hookrightarrow S\) to view the rows in \(S[G]^m\).
To test membership in \(\mathfrak a\) it suffices by \(S\)-linearity
to multiply on the right by each \(h\in G\). For \(p=\ell^\nu\) we have
\[
 \operatorname{tr}u\bigl(p(g-1)h\bigr)
    =p\operatorname{tr}u(gh)-p\operatorname{tr}u(h)\in\mathcal J .
\]
The tests on the paired generators are \(R_j(h),R_j^\sharp(h)\);
the common kernel maps to zero. For any
\(v=\sum_g c_g(g-1)\) in the augmentation ideal, with finite support,
the simultaneous coordinates give
\[
 v\rho_0\equiv
  \sum_i\Bigl(\sum_g c_g z_i(g)\Bigr)y_i\rho_0\pmod{\ker u}.
\]
Applied row by row, this constructs the required \(B_1\).

The characteristic coefficients of \(u(g)\) lie in \(S\).
The same holds for \(u(w)\), \(w\in S[G]\), and for
all flattened block matrices over this image, also
with polynomial coefficients. For square matrices \(A_1,\ldots,A_h\)
of the same size over a commutative ring, the primitive-word form
of Amitsur's identity
\cite{AmitsurCharpoly}\cite[p.~301, (6)]{ReutenauerSchutzenberger}
gives, in a central formal variable \(t\),
\[
 \det(1-t\textstyle\sum_i A_i)=
 \prod_{[v]\ {\rm primitive}}
       \det(1-t^{|v|}A_v),
\]
where the product is over primitive cyclic words and
\(A_v\) is the ordered product of the matrices of a
representative word.

Only words of length at most \(k\) contribute to the coefficient
of \(t^k\). This coefficient is an integral polynomial in the
characteristic coefficients of the word matrices. Applying it to scalar
multiples of \(u(g)\) proves the assertion for \(u(w)\).
For block matrices, split into individual blocks
\(E_{ij}\otimes u(w_{ij})\).
Every product is zero or \(E_{ab}\otimes u(w)\).
If \(a\ne b\), the matrix \(1-tE_{ab}\otimes u(w)\)
is block upper unitriangular when \(a\) precedes \(b\).
Its determinant is one over the coefficient ring.
If \(a=b\), its only contributing diagonal block is
\(u(w)\). This proves coefficient integrality
in \(S\) and in \(S[T]\).

Put
\[
 B_T=B-TB_1,\qquad
 P(T)=\det(u(B_T))\in S[T].
\]
Flatten \(u(B_T)\) to an \(mn\)-by-\(mn\) matrix \(M(T)\);
the variable \(T\) is central. With \(N_B=mn\ge1\), write
\(\det(zI-M(T))=z^{N_B}+\sum_{k=0}^{N_B-1}c_k(T)z^k\)
and \(c_{N_B}(T)=1\). Then
\[
 P(T)=(-1)^{N_B}c_0(T),\qquad
 \operatorname{adj}M(T)=(-1)^{N_B-1}
            \sum_{k=1}^{N_B}c_k(T)M(T)^{k-1}.
\]
Lifting the right-hand polynomial in the block entries gives a
matrix \(B_T'\) over \(S[G][T]\) such that
\[
 B_T'B_T\equiv P(T)\mathbf 1\pmod{\ker u}
\]
coefficientwise. View \(S[G]/\mathfrak a\) as an
\(S[G][T]\)-module by left multiplication from
\(S[G]\) and right multiplication by \(\rho_0\)
for \(T\). The two actions commute by associativity.
The covector relation says that \(B_T\) annihilates \(y\)
in this module, so the adjugate relation gives
\[
 y_iP(\rho_0)\in\mathfrak a,\qquad
 \mathcal I P(\rho_0)\subset\mathfrak a .
\]
Taking \(T=0\) in the definition yields
\(P(0)=\det(u(B))\).

Use the ordered product \(k_{\rm char}\) of
Section~\ref{sec:marked-ribet-coefficient}. On each lift it
acts by zero on the main summand and by one on the characters.
The separation element \(q_0^{\rm sep}\) acts on \(\xi_e\)
by \(C\mathbf1_{e=0}\). Thus \(\rho_1=q_0^{\rm sep}k_{\rm char}\)
acts on \(\xi_0\) by \(C\) and vanishes on all other summands
of each lift. Now
put \(f_B=\operatorname{tr}(u(\rho_1 P(\rho_0)))\).
The trace belongs to \(S\), by the characteristic
coefficient result applied to the group-algebra
element inside the trace. For \(w\in\mathcal I\)
put \(b=\det(u(w))\in S\). The adjugate of \(u(w)\)
is the image of an element \(w'\in S[G]\), again
by Cayley--Hamilton with coefficients in \(S\).

Thus \(u(w'w)=b\mathbf1\), and
\[
 b f_B
   =\operatorname{tr}\bigl(u(\rho_1\,w'wP(\rho_0))\bigr)
   \in\mathcal J ,
\]
since \(w'wP(\rho_0)\in\mathfrak a\).
At a lift system \(x\), one has \(\xi_0(\rho_0)=0\),
and \(u_x(\rho_1)\) is zero off that character. Write
\(\mu_0=\xi_0(\rho_1)\), the value of the character projector
multiplied by the first factor of the common multiplier. Hence
\[
 f_B(x)
  =\mu_0P(0)(x)
  =\mu_0\det(T_x(B))\prod_{e=0}^{t'-1}\det(\xi_e(B)).
\]
For each \(g\in G\) take
\(w_g=-\sum_i z(g)_i(g_i-1)\in\mathcal I\).
Then \(u(w_g)=1-u(g)\) on the full family, so that
\(\det(1-u(g))f_B\in\mathcal J\) with the indicated sign for every \(g\). Combining
the matrices satisfying the congruences and summing their integral weights
gives the common \(f\) used above. The argument
does not introduce factorial denominators or require
changing the coefficient embedding or the system set.

\subsection{The highest squarefree coefficient and the original test space}
\label{sec:top-coefficient-clearing}

\sectionopening{In this subsection, the polynomial of Section~\ref{sec:marked-ribet-coefficient} must be
evaluated on the space at the original level. We use the spherical
root construction of Section~\ref{sec:marked-spherical-root-chain} to
track the top squarefree coefficient under this restriction.}

\smallskip
\indent\textit{The type projection and the order of the congruence.}
Write \(r'\) for the source rank of this descent pair and
\(G_c=\prod_{w\in S_M}\mathrm{GL}_{r'}(k_w)\).
For the irreducible finite types \(\lambda\) with the
prescribed parameters the coefficient of \(g\in G_c\) in the
projector is
\[
 |G_c|^{-1}\sum_\lambda \dim(\lambda)\,
                              \operatorname{tr}\lambda(g^{-1}).
\]
Each trace is a sum of roots of unity. The residue action of
\(g\) is right translation by an integral lift of the first
block of \(g\), multiplied by \(\prod_{w\in S_M}\psi_w(\det g_w)\).
The scalar is a root of unity. Let \(k\) be the finite lift and
\(\rho_{\rm W}\) the Weyl representation in the transported coordinates.
For each finite point \(x\) the value of its right translate satisfies
\[
 \rho_{\rm W}(x_\ell^{-1})y(xk)
  =\rho_{\rm W}(k_\ell)\rho_{\rm W}((xk)_\ell^{-1})y(xk).
\]
The first matrix on the right preserves the integral Weyl lattice:
\(k_\ell\) and its inverse are integral at every place over \(\ell\),
including the selected split places when \(\ell\mid M\).
Thus every summand preserves the lattice before the division by
\(|G_c|\) in the characteristic-zero idempotent.

The lift is supported at the specified old places and their conjugates.
The block embeddings in the groups defining \(C_G\) show
\(|G_c|\mid C_G\). Hence \(C_G Qy\) is integral for integral
\(y\), also when \(\ell\) is an old prime. This conclusion
holds in the original level with the type condition imposed:
one takes the given level and character space at the old and
twist places, its \(Q\)-fixed subspace, and the spherical
invariants outside those places. Raising auxiliary level
changes only the last invariance condition.

Fix the chosen original system \(x_0\). There are finitely
many systems at the original level in the fixed type space.
For each original system \(x_k\) outside the weak lifts,
the definition in Section~\ref{sec:triangle-unit-input} gives
\(g_k\in G_E\) with \(\lambda_k=\det(1-N_{x_k}(g_k))\ne0\). These values are
integral: the powers of the continuous Galois matrix are
bounded by compactness, so its eigenvalues have absolute
value at most one. Put \(n_0=\prod_k\lambda_k\), with the
empty product equal to one.

The triangle construction
first gives its index set \(I\), with a choice of primes
and root chains at every positive depth for the weak
original systems congruent to \(x_0\). Put
\[
 n_1=C_G^2n_0,\qquad
 n_2=2n!n_1^{|I|+1}.
\]
At the fixed embedding, \(n_2\) is a nonzero element of the
original coefficient valuation ring. Choose an integer
\(\nu\ge\max\{1,\lceil v_\ell(n_2)\rceil\}\), where
\(v_\ell(\ell)=1\). This choice remains valid on extending that
valuation to the raised coefficient field.
Only now select the actual depth-\(\nu\) primes and the
common raised trace algebra. At each raised system use
\(\beta_k(x)=\det(1-N_x(g_k))\) with the same \(g_k\).
The scalar \(\eta=\ell^\nu/(2n!)\) satisfies
\[
 v_\ell(\eta)=\nu-v_\ell(2n!)\ge (|I|+1)v_\ell(n_1).
\]
In the coefficient DVR \(O\) of any resulting raised family this is
\[
 v_O(\eta/n_1^{|I|+1})
 =v_O(\ell)\bigl(\nu-v_\ell(2n!)-(|I|+1)v_\ell(n_1)\bigr)\ge0 .
\]
Hence \(\eta\in n_1^{|I|+1}O\) at the depth fixed from the
original systems and their index set.

Use the full raised trace algebra \(\mathbb T\) acting on its
integral lattice \(L\). The raised places are indexed by \(I\).
Write \(\mathcal N_j,\partial_j\) for the local norm and
derivative from the local comparison at \(j\), and
\(\mathcal N=\prod_{j\in I}\mathcal N_j\). The groups at
different places commute, as do their local operator algebras.
The original-level space is the fixed space for the product
\(H=\prod_{j\in I}G_j\), where \(G_j\) is the full finite residue group at \(j\).
Its integral lattice is \(L^H\): at either level the integral test
requires \(\rho_{\rm W}(x_\ell^{-1})y(x)\) to lie in the same Weyl
lattice for every finite point \(x\). Passing to the original level
adds precisely the \(H\)-invariance condition to this pointwise test.
Use \(c=C_G\) for the common multiple of the ambient compact stabilizer orders.

For the local polynomial \(H_j\) of
Section~\ref{sec:marked-ribet-coefficient}, put
\(P_j(X)=(X-1)H_j(X)\). At a raised system \(x\) and \(g\in G_E\),
\[
 \begin{aligned}
 t_{j,g}(x)&=\operatorname{tr}\bigl((N_x(\tau_j)-1)N_x(g)\bigr),\\
 s_{j,g}(x)&=\operatorname{tr}\bigl(P_j(N_x(\phi_j))N_x(g)\bigr).
 \end{aligned}
\]
Under \(J=I,\ X_j=\alpha_j\), extend the scalar ideal
from \(S\) to \(\mathcal J_{\mathbb T}
=\mathcal J\,\mathbb T[\alpha_j:j\in I]\) by the coefficient
inclusion \(\mathbb T_O\hookrightarrow\mathbb T\). Its generators are the traces \(\ell^\nu\operatorname{tr}N(g)\)
and \(t_{j,g}-\alpha_js_{j,g}\).
In this calculation reductions are modulo
\(\eta=\ell^\nu/(2n!)\) in the coefficient valuation ring.
On the full raised lattice,
Section~\ref{sec:raised-prime-tame-restriction} gives
\[
 \mathcal N_jt_{j,g}=0,\qquad
 \partial_j t_{j,g}\equiv\mathcal N_js_{j,g}\pmod\eta.
\]
The local operator in the second relation uses
\((X-q_j^{t'-1})H_j(X)\). Since \(q_j\equiv1\bmod\ell^\nu\),
its trace reduces to the one with \(P_j\) modulo \(\eta\).
Let \(A\in\mathbb T\) be the indicator
of the residual congruence class of the chosen lift.
It belongs to the integral diagonal algebra. To construct it,
for each system outside this class choose an integral
characteristic-polynomial coefficient separating it modulo
the maximal ideal from the chosen lift. Subtract the value
at the outside system and divide by the unit value at the
chosen lift.

The product \(T\) of these factors is zero
outside the class and is a unit on the class. The polynomial
\[
         1-\prod_{x\,\text{in the class}}(1-X/T(x))
\]
evaluated at \(T\) gives the indicator. All scalars here
are in the coefficient valuation ring. Put
\[
 \delta_0=\prod_k(\lambda_k-\beta_k),\qquad
 n_0=\prod_k\lambda_k,\qquad n_1=c^2 n_0,\qquad
 \delta'=A\delta_0 .
\]
On every lift system the trivial character \(\xi_0\) gives
\(\beta_k=0\) for all \(k\) and \(\delta_0=n_0\).
At the nonweak original system \(x_k\), the \(k\)-th factor in
\(\delta_0\) is zero. Thus an original system in the support
of \(\delta'\) is weak and residually congruent to \(x_0\).
For each \(j\in I\) the complementary spherical factor for
that system divides \(H_j\) by the product construction in
Section~\ref{sec:triangle-unit-input}. The depth-\(\nu\) prime
choice gives the spherical root chain on every such system.

We construct the relative trace operators and verify the spherical
annihilation on these systems. Section~\ref{sec:raised-prime-integral-clearing}
then yields \(Z_{j,h}\) for \(h\in G_j\) with
\[
 c\partial_j(\rho(h)-1)\equiv\mathcal N_j Z_{j,h}\pmod\eta,
 \qquad \delta'\mathcal N Z_{j,h}=0 .
\]
For the relative trace on endomorphisms, see~\cite[\S13, pp.~507--508]{HigmanOperators}.
To construct the local and full relative traces, take \(K\)
to be the ambient standard maximal finite compact.
Reduction at the original auxiliary primes gives
\(\pi:K\to H\), where \(H=\prod_{j\in I}G_j\) uses
the full rank \(n\) at each place. Put \(K^1=\ker\pi\).
On \(\Gamma\backslash G_f/K^1\), with \(\Gamma\)
the rational left group, let \(D\) be \(H\) or a
line stabilizer in one factor. The stabilizer at
\(x=\Gamma gK^1\) is
\[
            D_x=D\cap \pi(K\cap g^{-1}\Gamma g).
\]
Indeed a lift in \(K\) of a residue element fixing \(x\)
differs from a rational stabilizer by an element of
\(K^1\). Its order therefore divides the order of the
ambient compact rational stabilizer. At an old split place over the
odd prime \(p_i\) its integral action injects under reduction into
\(\mathrm{GL}_n(\mathbb F_{p_i})\): the principal congruence kernel
in \(\mathrm{GL}_n(\mathbb Z_{p_i})\) is torsion-free.
Block inclusion for \(n=r'+t'\le rt\) and the definition of \(C_G\)
in Section~\ref{sec:marked-norm-line} give the asserted divisibility.

On each \(D\)-orbit define an integer-valued weight
\(m_D\) equal to \(c/|D_x|\) at one point and zero
elsewhere. Each orbit point occurs \(|D_x|\) times
under translation by \(D\), so
\(\sum_{d\in D}m_D(xd^{-1})=c\).
Multiplication by the pullback of this weight commutes
with every right \(K^1\)-translation. The old and twist
compact factors, including conjugate values, lie in
\(K^1\) since their places are off the auxiliary support.
The multiplier therefore respects the original projection
\(Q\), both character equations and the auxiliary
principal levels. It preserves \(L\) because it is an
integer scalar at each point with specified lattice coordinates. For this
endomorphism \(\Psi_D\) the orbit identity gives
\[
             \sum_{d\in D}\rho(d^{-1})\Psi_D\rho(d)=c\,1_L .
\]
In particular \(\mathcal N\Psi_H v=cv\) on the
original-level invariants. Use \(D\) equal to the
line stabilizer in the local derivative calculation.

For the spherical annihilation at \(j\), let \(q_j\) be the residue
cardinality and \(a=a(\phi_j)\). Here \(t'=n-r'>0\). The actual
operator is \(B_j(1)=p_j(aU)L_J(1)\), where
\(p_j(X)=(X-q_j^{\,t'-1})H_j(X)\). Its input \(L_J(1)y\) is
\(I_1\)-fixed. The congruence \(q_j\equiv1\bmod\ell^\nu\)
gives \(p_j(X)\equiv (X-1)H_j(X)\bmod\ell^\nu\), the polynomial
appearing in \(s_{j,g}\). Over the characteristic-zero space, if \(W_x\) is a cofinite
Hecke eigenspace in the raised space, the action of \(H\) preserves
\(W_x\) and
\[
             \mathcal N|_{W_x}=|H|\operatorname{pr}_{W_x^H}.
\]
Here \(W_x^H\) is the original-level space of that system. In
particular, the norm kills the systems which have no original-level
vectors. The integral operator is the sum of the \(H\)-translations.
For an original spherical vector \(x\) on a system supported by
\(\delta'\), the norm is self-adjoint and acts on \(x\) by \(|H|\).
Therefore
\[
 \langle x,\mathcal N B_j(1)y\rangle
             =|H|\langle x,B_j(1)y\rangle=0
\]
by the adjoint orthogonality in
Section~\ref{sec:marked-spherical-root-chain}.
The image of \(S_j=\mathcal N B_j(1)\) lies in the
full original-level invariants. For every original-level
\(w\), the vector \(\delta'w\) is a sum of supported
vectors, and hence \(\langle\delta'w,S_j y\rangle=0\).
The local operators defining \(S_j\) commute with
the diagonal trace operator \(\delta'\). Applying the orthogonality to \(w=S_jy\) and to \(S_j(\delta'y)\)
gives \(\langle\delta'S_jy,S_j\delta'y\rangle=\|\delta'S_jy\|^2=0\).
This proves \(\delta'\mathcal N B_j(1)=0\) on the
characteristic-zero space and on the integral lattice.
Together with the local Higman operator it gives the
two stated properties of \(Z_{j,h}\).

Put \(\Delta=c^2\delta'\) and
\[
 \widetilde{\mathcal N}=\prod_{j\in I}(\mathcal N_j+\alpha_j\partial_j),
 \qquad \widetilde{\mathcal N}_S=[\alpha_S]\widetilde{\mathcal N},\qquad
 \alpha_S=\prod_{j\in S}\alpha_j .
\]
In the polynomial operator ring modulo \(\eta\) and
\(\alpha_j^2=0\), the two local trace equations give
\(\widetilde{\mathcal N}(t_{j,g}-\alpha_j s_{j,g})=0\). Coefficients in
\(\mathbb T\) commute with the local operators. The remaining generators
\(\ell^\nu\operatorname{tr}N(g)\) are zero modulo \(\eta\) since
\(\ell^\nu=(2n!)\eta\). Hence
\([\alpha_I](\widetilde{\mathcal N} B)\equiv0\pmod\eta\)
for every \(B\in\mathcal J_{\mathbb T}\). The quotient by \(\alpha_j^2\) does
not change any squarefree coefficient.

For \(i\ge |S|\) there is an integral \(Y_{S,i}\) such that
\[
 \Delta^i\widetilde{\mathcal N}_S\equiv Y_{S,i}\mathcal N\pmod\eta,
 \qquad Y_{\varnothing,i}=\Delta^i .
\]
Indeed the empty case is \(\widetilde{\mathcal N}_\varnothing=\mathcal N\).
For \(S\ne\varnothing\) put \(X=c\delta'\Delta^{i-1}\widetilde{\mathcal N}_S\).
For a place \(j\notin S\), \(X(\rho(h)-1)=0\) for \(h\in G_j\),
since \(\mathcal N_j\) is a factor. For \(j\in S\), apply the
local derivative congruence at \(j\) and commute \(\delta'\)
to the right of \(\widetilde{\mathcal N}_{S\setminus\{j\}}\).
By induction, \(X(\rho(h)-1)\) reduces to
\(Y_{S\setminus\{j\},i-1}\mathcal N\delta' Z_{j,h}=0\).
Thus \(X\rho(h)\equiv X\) for all \(h\in H\). Multiplying
the full Higman identity on the left by \(X\) gives
\(cX\equiv X\Psi_H\mathcal N\). Since \(cX=\Delta^i\widetilde{\mathcal N}_S\),
take \(Y_{S,i}=X\Psi_H\) with this order of composition.
The Hecke operators \(\Delta D_S\) in the calculation below commute
with \(\mathcal N\) and act before this multiplier. When \(I=\varnothing\)
only \(Y_{\varnothing,i}\) occurs.

Let \(D(\alpha)=f((\alpha_j)_{j\in I})\), where \(f\) is
the common trace multiplier of Section~\ref{sec:marked-ribet-trace-determinant}
on the same raised space, with distinguished lift \(\pi=x_0\).
Under \(J=I\) and \(X_j=\alpha_j\), the relations \(R_j(g)\)
become \(t_{j,g}-\alpha_js_{j,g}\). For the augmentation element \(w_{g_k}\in\mathcal I\)
constructed there, \(u(w_{g_k})=1-u(g_k)\), so
\(\det u(w_{g_k})=\beta_k\). The trace determinant gives
\(\beta_kD\in\mathcal J_{\mathbb T}\). On a lift system \(x\), the successive
replacements \(d_e\mapsto d_e+\zeta_e\) give
\[
 f_x\equiv Cc_0^{t'}D_x(0)\prod_e\Gamma_e(0)X_J
                         \pmod{(X_j^2:j\in J)} .
\]
Indeed the \(e\)-th difference contains
\(\zeta_e\prod_{\mathrm{ty}(j)<e}X_j\), which vanishes in the
squarefree quotient by \eqref{eq:display-rb3}. Thus
\(D(x)\equiv [\alpha_I]D(x)\,\alpha_I\) in that quotient and
\([\alpha_I]D(x_0)=d\). Write \(D_S=[\alpha_S]D\), and put
\(i=|I|\), \(b=i+1\). 
The trace relation \(\beta_kD\in\mathcal J_{\mathbb T}\) gives
\((\delta_0-n_0)D\in\mathcal J_{\mathbb T}\) by expanding the product.
The first-order lower coefficient vanishing means \(D_S=0\)
on the lift systems when \(S\ne I\).
The original-level choice above gives
\(2n!\,n_1^b\mid\ell^\nu\) in the coefficient valuation ring
before the raised prime selection.

Let \(l,v\) be original-level vectors with \(d_0l\in L\) for
an integral scalar \(d_0\).
Assume \(\delta_0v=n_0v\), \(D_Sv=0\) for \(S\ne I\), and
\(A\delta_0l=n_0 Av\).
In the raised lattice \(L\) of the prescribed old type, take
\[
 G=[\alpha_I](\widetilde{\mathcal N}D)
    =\sum_{S\subseteq I}\widetilde{\mathcal N}_{I\setminus S}D_S
    \in\operatorname{End}_O(L).
\]
Indeed only complementary subsets contribute to the squarefree
monomial \(\alpha_I\). Each \(D_S\in\mathbb T\) acts by a scalar
on every cofinite Hecke eigenspace, which the auxiliary local
operators preserve by Section~\ref{sec:integral-raised-traces}.
Thus \(D_S\) commutes with those operators. Reduction on
\(L/\eta L\) preserves the coefficient identity.
Since \(\Delta^b-n_1^b A^b\) is divisible by
\(\delta_0-n_0\), the vanishing modulo the congruence ideal and the factorization
of \(\widetilde{\mathcal N}_S\) give the congruence below. In each
summand \(\Delta\) commutes with the local norm coefficients, so put
\[
 \Delta^b\widetilde{\mathcal N}_{I\setminus S}D_S
   =(\Delta^i\widetilde{\mathcal N}_{I\setminus S})\Delta D_S
\]
before substituting the factorization of the parenthesized operator. Thus
\[
 n_1^b A^b G\equiv
 \Delta^b G\equiv
 \sum_{S\subseteq I}Y_{I\setminus S,i}\mathcal N\Delta D_S
                         \pmod{\eta}.
\]
Apply this congruence to \(u=\Psi_H(d_0 l)\in L\). Since \(d_0 l\) is fixed by \(H\),
\(\mathcal N u=c d_0 l\).
The term for \(S\) on the right is
\[
 n_1 c d_0\,Y_{I\setminus S,i} A D_S v .
\]
Only \(S=I\) survives. Since \(Y_{\varnothing,i}=\Delta^i\),
the sum is \(n_1^b c d_0 A^b D_Iv\).
There is therefore a vector \(z\in L\) with
\[
 n_1^b(A^bGu-cd_0A^bD_Iv)=\eta z.
\]
Since \(n_1\ne0\), \(\eta/n_1^b\in O\) and \(A^bGu\in L\),
division in the ambient vector space gives
\(cd_0 A^bD_Iv=A^bGu-(\eta/n_1^b)z\in L_{\rm raised}\).
The functions \(A,D_I\) act by scalars on each original summand
\(V_{\rm orig}\cap W_x\). Hence this vector belongs to
\(V_{\rm orig}\cap L_{\rm raised}=L_{\rm orig}\) by
Section~\ref{sec:integral-raised-traces}.
Thus the exact clearing function before the type denominator
is
\[
                         s=A^bD_I=AD_I .
\]
At the chosen lift system \(x_0\) it has value \(D_I(x_0)\).

For the fixed descent pair and twist, let \(Z\) be the span
at the original level of a collection of good lift eigenvectors.
Each generator belongs to one old finite type and one outside-level
Hecke system. Several generators may belong to the same joint
summand. The distinguished system \(x_0\) is the system of one
of these generators. Taking all the generators gives the full good
lift space.
The old central idempotents split the projected source into its type
summands. The commuting spherical operators give the isotypic
decomposition there, and the finite-type and spherical intertwining
relations of Sections~\ref{sec:upper-good-projector} and
\ref{sec:galois-types-and-traces} give the stated supports on the
lifts. Choose a good ambient finite-type idempotent \(e_0\)
whose image contains the types of the selected generators.
Write \(K_w(X_w)\) for the
original spherical Hecke polynomial at the lift root, for
each split \(w\) outside the descent level and its conjugate.
The false-lift subspaces are
\[
            F_w=\operatorname{ran}e_0\cap Z^\perp
                                 \cap\ker K_w(X_w).
\]
Let \(Q\) be the self-adjoint projection to the fixed lift type
used for \(x_0\), with \(\operatorname{ran}Q\subseteq
\operatorname{ran}e_0\).

On a lift generator of finite type tuple \(\lambda\), the central
idempotent \(Q\) acts by one or zero according as \(\lambda\) occurs
in its defining sum or does not. Likewise each orthogonal isotypic
projection \(E_x\) acts on the generator by one or zero. Thus both
preserve \(Z\), and self-adjointness gives preservation of \(Z^\perp\).
Moreover \(Q\) commutes with each \(E_x\), since the old residue
actions commute with the spherical Hecke actions away from the level.
For an original integral test \(y\perp F_w\) for all \(w\),
take \(l=Qy\) and \(v=Q\pi_Zy\), where \(\pi_Z\) is orthogonal
projection. Choose one split place \(w\nmid\ell\) off the level,
the exceptional places of bases of the original summands, and the
excluded places of their Galois comparisons, as in
Section~\ref{sec:integral-raised-traces}. Write \(q=\mathrm Nw\), let \(P_x\) be the
monic degree-\(n\) spherical polynomial of system \(x\) there, and
put \(\nu_w=\alpha_{\rm nor}(\mathrm{Fr}_w)^{-1}\) as in
Section~\ref{sec:normalized-trace-multiplicities}. On a weak system,
\(\det(1-N_x(\mathrm{Fr}_w))=0\) and hence \(P_x(\nu_w)=0\).
Since \(|\nu_w|^2=q^{r'}\), the spherical reciprocity identity
of Section~\ref{sec:marked-spherical-root-chain} gives
\[
 P_x(q^{n-1}/\overline{\nu_w})=
 P_x(q^{t'-1}\nu_w)=0 .
\]
Consequently the reciprocal Hecke operator \(K_w(X_w)\) acts
by \(X_w^nP_x(X_w^{-1})=0\) on that system at the original
spherical level, with \(X_w=(q^{t'-1}\nu_w)^{-1}\). The projections \(E_x,Q,\pi_Z\) therefore commute. For a weak
original system \(x\), the vector \(u_x=E_xQ(1-\pi_Z)y\)
belongs to \(F_w\) at such a place \(w\), since
\(\operatorname{ran}Q\subseteq\operatorname{ran}e_0\).
As \(E_xQ(1-\pi_Z)\) is an orthogonal projection, we obtain
\(\|u_x\|^2=\langle y,u_x\rangle=0\) by the test hypothesis.
Hence \(l=v\) on the weak original systems. Moreover \(v\) is
supported on lift systems since \(Z\) is spanned by the lift
generators and \(Q\) preserves each of their isotypic components.
There \(\delta_0=n_0\) and \(D_S=0\) for \(S\ne I\); on the
nonweak original systems \(\delta_0=0\). This proves all three
required relations for \(l,v\).

The integral realization of \(Q\) at its base level uses
\(cQy\in L\) for integral \(y\). Taking \(d_0=c\) gives the following bound.

\begin{proposition}[Clearing on the original test space]
\label{prop:clearing-on-original-tests}
The operator \(c^2sQ\pi_Z\) is integral on the original integral
tests perpendicular to \(F_w\) for all \(w\). Its eigenvalue on
\(Z\cap W_{x_0}\cap\operatorname{ran}Q\), the joint summand
containing the selected lift generator, satisfies
\[
 v_\ell\bigl(c^2s(x_0)\bigr)
 \le \mathrm{sel}+2t'\mathrm{ram}
                  +r'\mathrm{loc}_\ell+\mathrm{err}_\ell .
\]
\end{proposition}
\begin{proof}
Take \(l=Qy\), \(v=Q\pi_Zy\) and \(d_0=c\) in the original-level
congruence above. The type denominator gives \(d_0l\in L\); the
orthogonality to every \(F_w\) gives \(A\delta_0 l=n_0Av\).
The lift supports give \(\delta_0v=n_0v\) and \(D_Sv=0\) for
\(S\ne I\). Hence \(cd_0A^bD_Iv=c^2sQ\pi_Zy\) belongs to
\(L_{\rm orig}\). On the indicated joint summand \(Q\) and
\(\pi_Z\) act as the identity and \(A(x_0)=1\). Its eigenvalue is
\(c^2D_I(x_0)=c^2d\); \eqref{eq:display-rb1} gives the valuation
bound with \(r'\) the source rank and \(t'=n-r'\).
\end{proof}
For \(Z=0\) the determinant is the empty product. Otherwise take a basis
\(b_1,\ldots,b_d\) from the lift generators, each in a single finite-type summand
and a single outside-level Hecke eigenspace. For each \(b_i\) use its clearing
operator \(D_i=c^2s_iQ_i|_Z\) from the proposition. Write \(\mathrm{sel}_i\)
for the normalized Selmer length in this application. The ranks \(r',t'\)
and the terms \(\mathrm{ram},\mathrm{loc}_\ell,\mathrm{err}_\ell\) are fixed
for this lift image and twist. Put
\[
 \mathrm{sel}_Z=\max_{1\le i\le d}\mathrm{sel}_i,\qquad
 B_{\rm tr}=\mathrm{sel}_Z+2t'\mathrm{ram}
                      +r'\mathrm{loc}_\ell+\mathrm{err}_\ell .
\]
Write \(D_i b_k=d_{ik}b_k\). Indeed the type projectors act on a fixed
type by zero or one, and the trace coefficients act by scalars on
each Hecke system. The selected diagonal coefficient satisfies
\(v_\ell(d_{kk})\le B_{\rm tr}\).

For each \(k\), choose a nonzero entry \(a_k\) of least valuation in
the column \((d_{ik})_i\), using the fixed coefficient embedding.
The reduction of \(\sum_i (d_{ik}/a_k)Y_i\) is a nonzero linear form.
Over \(\overline{\mathbb F}_\ell\) the product of these forms is
nonzero at some tuple. Its entries lie in a finite field
\(\mathbb F_{\ell^a}\). Lift zero entries to zero and nonzero entries
to their Teichm\"uller lifts, roots of unity of order dividing
\(\ell^a-1\). Transport these coefficients \(u_i\) back by the chosen
embedding. The individual clearing operators and the column minima
are taken at that same embedding.
The operator \(D_Z=\sum_i u_iD_i\) still clears \(\pi_Z y\) for every
original integral \(y\) orthogonal to all \(F_w\). Its eigenvalue on
\(b_k\) has valuation \(v_\ell(a_k)\le v_\ell(d_{kk})\), so
\[
                   v_\ell(\det\nolimits_Z D_Z)\le d B_{\rm tr}.
\]
\nopagebreak[4]
The \hyperref[sec:simultaneous-original-clearing]{application in
Section~\ref*{sec:compression-trivial-norm}} uses this combination
on the common original test space. Section~\ref{sec:critical-integral-pairings}
multiplies \(D_Z\) by the stabilizer multiple \(C_G\) when clearing
the class pairing.
\par\nopagebreak[4]
\part{Period comparison and the lower determinant bound}\label{part:lower}

\section{Period vanishing and the integral comparison lattice}
\label{sec:comparison-vanishing-lattice}

\sectionopening{In this section, we use the finite support of Section~\ref{sec:adjacent-unfolding} to
determine the terms that survive in the comparison periods. The integral
first factors and their common scalar give the lattice used in
Section~\ref{sec:comparison-operator-first-lifts}.}

\subsection{Vanishing of the comparison periods}
\label{sec:marked-comparison-period-vanishing}

\sectionopening{In this subsection, the finite support of the data in Section~\ref{sec:adjacent-unfolding}
forces certain kernels to vanish after integration. We establish the
vanishing with the original local characters and measures before using
it in Section~\ref{sec:original-comparison-tests}.}

Fix \(1\le j\le m\) and an allowed twist \(\vartheta\).
Use the pair \(W=W_j,V=I_n\), the coordinate inclusion
\(\iota_0:E^r\hookrightarrow E^n\) into the first \(r\) rows,
and the
split sets \(S_M,Q\) of chosen and twist places. The
places in \(S_M\) are the divisors of \(\pi_M\), and every
residue character on the exact support \(Q\) is nontrivial. These sets
and their conjugates are disjoint in the construction.
Let \(S=S_M\cup Q\). At \(w\in S\) write \(k_w\)
for its residue field and \(\chi_w\) for the exact residue
character in the descent pair. Fix the false lift \(\Xi\):
it lies in the image of the finite-type idempotent, is perpendicular to the
image of the projected lift at this level, and is killed by the Hecke polynomial at a split place off \(S\cup\bar S\).

At a chosen place and a nonzero residue covector \(\lambda\) of
length \(r=\dim W\), put
\[
 H_\lambda=\{1+u\lambda:\det(1+u\lambda)\ne0\},\qquad
 \mathcal A_\lambda
    =\sum_{s\in H_\lambda}\chi_w(\det s)^{-1}s .
\]
In a basis with \(\lambda=(1,0,\ldots,0)\), \(H_\lambda\)
is the product of the lower unipotent radical for the
\((1,r-1)\) parabolic and the first \(\mathrm{GL}_1(k_w)\).
The operator \(|H_\lambda|^{-1}\mathcal A_\lambda\)
first projects to unipotent invariants and then to type
\(\chi_w\) on the first \(\mathrm{GL}_1\).

For a squarefree Green parameter \(c_w\), the finite
Harish--Chandra support is the collection of cuspidal
types of the irreducible factors of \(c_w\). Any
\(\mathrm{GL}_1\) of type \(\chi_w\) in those unipotent
invariants would give, by finite parabolic induction and
Frobenius reciprocity, the factor of degree one named by
\(\chi_w\). This contradicts the condition on the selected parameters:
\(\chi_w^{q_2}=1\) and no \(q_2\)-torsion root of
\(k_w^\times\) occurs in the polynomial. The
\(q_2\)-power condition is for the original \(\chi_w\),
since the auxiliary twist is unramified at \(w\) and hence is
trivial on \(\mathcal O_{E_w}^{\times}\).
Thus the weighted sum kills every selected finite type at \(w\)
and hence kills \(\Xi\). This is the parabolic product
in Green's classification \cite[\S2, pp.~410--411, and \S8, Theorem~14]{Green}.

Let \(\mathcal D\Xi\) be the normalized descent with the original
block test. First consider it in the Petersson space at the old
principal level. Choose \(w\mid\pi_M\).
In the split Schr\"odinger coordinates its test is supported on
\(X\equiv\binom{I_r}{0}\pmod w\). The top minor is a unit, so
\(X:E_w^r\to E_w^n\) is injective throughout the support.
Changing the source or target point at \(w\) multiplies \(X\)
on the right or left by an invertible matrix and preserves this
injectivity, including at every representative of the definite
target classes.

Fix such a representative \(c\), its real polynomial \(P\), and
write \(H=\operatorname{ht}_{W}(g)\).
The source-variable estimate of
Section~\ref{sec:four-unitary-constant-terms} starts from
\(|L_{g,c}(P)|\ll_c H^{A_0}s(P)\) for the theta evaluation and
a Schwartz seminorm \(s\). For \(H>2\) take a shortest
isotropic source vector \(v\). In the moved finite coordinates
the modulation detected by a definite target basis is nonzero,
since \(Xv_w\ne0\). Thus the isotropic average annihilates the
evaluation for every real test. Its real displacements are
\(O_c(H^{-1})\). Iterating the difference from the average
\(A_0+B\) times as in that section gives, for every integer \(B\ge0\),
\begin{equation}
    |L_{g,c}(P)|\ll_{c,B} H^{-B}s_B(P).
    \label{eq:original-block-descent-decay}
\end{equation}
On bounded heights the initial estimate suffices. The class
sum is finite, so the descent is rapid and square integrable.
Its real type is the one computed in
Section~\ref{sec:bottom-descent-real-type}, using the transposed
coefficient projection of Section~\ref{sec:four-unitary-real-projections}.
At the old places the source principal level fixes the inclusion
test. At a twist place the character of the invertible top block
cancels the inverse determinant character on the integral source
group, and away from these places the vacuum is fixed.

The pairing of \(\mathcal D\Xi\) against
any source form \(y\), up to the fixed nonzero normalizing
scalar, is the pairing of \(\Theta y\) with \(\Xi\).
Write \(e_{\rm src}\) and \(e_{\rm tar}=e_0\) for the actions of
the same central idempotent of
\(\prod_{w\in S_M}\mathrm{GL}_r(k_w)\) on the source and target.
For a tuple \(s\) lift \(s_w\) on the source and
\(\operatorname{diag}(s_w,I_{n-r})\) on the target. Their simultaneous
action fixes the reduction of the block test, since
\(\operatorname{diag}(s_w,I_{n-r})\iota_0=\iota_0s_w\).
The simultaneous Weil multiplier on this test is
\(\prod_{w\in S_M}\chi_w(\det s_w)^{-1}\), as in the split calculation
of Section~\ref{sec:critical-hecke-systems}. Write
\(\rho_{\rm tar}(s)\) for the target translation multiplied by
\(\prod_{w\in S_M}\chi_w(\det s_w)\). The twist is trivial on the
unit groups at the old places. At each \(w\in Q\), the character
of the invertible block cancels the inverse character of the
integral source translation. Thus the simultaneous change of
variables gives
\(\rho_{\rm tar}(s)\Theta=\Theta\rho_{\rm src}(s)\).
For an integral source lift \(k\) of \(s\), apply the same change
of variables against compactly supported smooth source tests and
take adjoints. The two residue actions are unitary, hence
\[
   (\mathcal D\Xi)(gk)
     =\mathcal D\bigl(\rho_{\rm tar}(s)\Xi\bigr)(g).
\]
The equality of these continuous functions follows from the
equality of the pairings against the smooth tests.
Summing with the central idempotent coefficients gives
\(e_{\rm tar}\Theta=\Theta e_{\rm src}\). Adjointness and
self-adjointness of the two projectors give
\[
 e_{\rm src}\mathcal D\Xi=\mathcal D e_{\rm tar}\Xi=\mathcal D\Xi
 \qquad\text{in } L^2([U(W)]) .
\]
Each summand retained by \(e_{\rm src}\) has a primitive residue type
at a split old place, by Section~\ref{sec:marked-residue-projection}.
To check cuspidality, take a proper rational parabolic with unipotent
radical \(N\). At this place use Iwasawa decomposition in computing
the constant term at any translate. Right invariance of the integral on
\([N]\) allows insertion of the probability average over a compact
conjugate of \(N(F_{0,w})\cap K_w\) on the local vector, with the
parabolic in standard position. Its reduction is the unipotent radical
of a proper residue parabolic and its average on the cuspidal residue
type is zero. Hence the constant term vanishes for every selected
summand. Together with the real-line and level conditions above this
puts \(\mathcal D\Xi\) in the original projected source. Taking
\(y=\mathcal D\Xi\) in the adjunction, the pairing
\(\langle\Theta y,\Xi\rangle\) vanishes by \(\Xi\perp Z\).
Thus \(\|\mathcal D\Xi\|_2^2=0\). The descent is continuous
at every source point. Indeed, over a compact set of source points and
compact representatives of the definite target quotient, the finite
translates of its test form a finite set. Their rational support lies
in a fixed fractional lattice. The real Schwartz seminorms of the
translated test are uniformly bounded there, so a summable inverse
power on that lattice dominates the kernel. Dominated convergence
gives continuity after pairing with \(\Xi\). The quotient Haar measure
has full support, hence \(\mathcal D\Xi(g)=0\) for every \(g\). At the
original block finite datum the kernel at the harmonic highest vector
is this normalized descent up to conjugation and a nonzero scalar.

\subsubsection{The real type in the bottom descent}\label{sec:bottom-descent-real-type}

\sectionopening{In this subsection, we determine the type for the descent used in Section~\ref{sec:original-comparison-tests}.}

At a place \(\tau\), order the entries of the source
Harish--Chandra parameter as
\(\varsigma_0>\cdots>\varsigma_{r-1}\).

Write \(a=k+j-1\), \(b_-=m-j\), \(n=r+t\), and
\(w=w_j(\tau)\). Then \(a+b_-=r\). The matrices for
the positive source sector use the first \(a\) source
coordinates in decreasing order. Those for the
negative sector use the last \(b_-\) in increasing order.
The two families of column counts are
\[
 \begin{split}
 k^+(h)&=\bigl(\varsigma_{h-1}
       -\max(\varsigma_h,(t-1)/2-w)-1\bigr)_+,\\
 k^-(h)&=\bigl(\min(\varsigma_{r-h-1},-w-(t-1)/2)
                       -\varsigma_{r-h}-1\bigr)_+ .
 \end{split}
\]
Here one uses the first formula only through height \(a\)
and the second only through height \(b_-\).
In these formulas use \(\varsigma_r=-\infty\) and
\(\varsigma_{-1}=+\infty\) at a missing endpoint. Write \(\mu^+_\tau(h)=k^+(h)\) and \(\mu^-_\tau(h)=k^-(h)\),
extended by zero at height zero and beyond the corresponding sign block.
The inequality
\(\varsigma_i+w\ge \Delta/2\) for \(i<a\), and
\(\varsigma_i+w\le-\Delta/2\) for \(i\ge a\),
makes all entries beyond these heights zero.

In the compact
target the column multiplicity at height \(h\) is
\(k^+(h)+k^-(n-h)\), and the determinant exponent is
\(-\sum_h k^-(h)\).
The differences of consecutive twice Harish--Chandra
parameters are twice these column multiplicities plus two.
At position \(a\), the determinant exponent cancels the
whole negative contribution. It follows that after
subtracting \(2u-1\), where
\(u=b_-+w-\lfloor r/2\rfloor\), the compact twice
parameter is the doubled source list with the inserted block
\[
 -2w+t-1,\ -2w+t-3,\ \ldots,\ -2w-t+1
\]
between the two sectors. This is the regular discrete-series
theta correspondence with one inserted compact block,
after the known determinant shifts of the Weil action.
In Harish--Chandra units the retained entries are
\(\eta_i=\varsigma_i+w\) and the inserted entries are
\(\gamma_b=(t-1)/2-b\), \(0\le b<t=n-r\). Thus
\[
 \min_i|\eta_i|\ge\Delta/2\ge(t+1)/2,\qquad
 \min_{i,b}|\eta_i-\gamma_b|
       \ge(\Delta-t+1)/2>500(r+n)^3
\]
for sufficiently large \(x\). The positive source sector has \(a\)
positive entries and the negative sector has \(b_-\) negative
entries. They all give positive target signature, and the inserted
string contributes \(t\) further positive directions. Hence the
target is \(U(a+b_-+t,0)=U(n,0)\).
Use Li's regular discrete-series correspondence
\cite{LiThetaCohomology} in the formulation of
\cite[Theorem~2.7(b), p.~641]{PaulFirstOccurrence}.
Here \(r\) is odd and \(n=r+t=6r\) is even, so the genuine
Harish--Chandra lattice in this normalization is
\(\mathbb Z+(r+n+1)/2=\mathbb Z\). The entries \(\eta_i\) are
integers since both \(\varsigma_i\) and \(w\) are integral, and the
inserted entries \(\gamma_b\) are integral because \(t\) is odd.
The four counts are \((a,0,0,b_-)\), the rank difference is \(t\),
and the first signature case is \((a+b_-+t,0)\).
It gives a nonzero theta lift with compact parameter
\[
 (\eta_0,\ldots,\eta_{a-1},
      \gamma_0,\ldots,\gamma_{t-1},
      \eta_a,\ldots,\eta_{r-1}).
\]
The last string is empty if \(b_-=0\). Restoring the determinant
shifts gives the compact representation computed above.
For this regular range see also the proof of
\cite[Lemma~8.2]{IchinoRealTheta}.

The polynomial is the product of powers of leading
minors in the two sectors. For the positive sector use the
first target rows in the prescribed order and the positive source
columns just specified. For the negative sector use the last target
rows in reverse order and the negative source columns in their
prescribed order. The two sets of target rows are disjoint since
\(a+b_-=r<n\).

Each leading minor is a nonzero polynomial. A raising operator
annihilates it or produces a determinant with two equal rows or
columns. A mixed Laplacian differentiates one positive and one
negative variable with the same target row. Their row sets are
disjoint, so all these contractions vanish. The product is therefore
jointly harmonic and highest for the two compact orders. Its target
weight is the list of positive row lengths followed, after the zero
coordinates, by the negatives of the reversed negative row lengths.
The splitting character cancels the vacuum determinant shift.

It is exactly the compact representation with column
multiplicities and determinant exponent described above.
The polynomial is the image of its highest vector
by the Cartan product, since evaluating both highest
vectors on the ordered coordinate arrays gives these same
leading minors. Thus projection to that compact target
highest line is nonzero. Write \(\tau\) for this compact type and
\(M_\tau\) for its multiplicity space. On \(\tau\widehat\otimes M_\tau\)
the compact group generates
\(\mathcal B(\tau)\otimes 1\). Its real dual partner generates the
commutant \(1\otimes\mathcal B(M_\tau)\)
by \cite[Theorem~6.1, p.~550]{HoweDuality}. Projection onto the chosen highest
line in \(\tau\) therefore gives one irreducible source representation.
This proves the occurrence used in the compact bump argument.

The local oscillator correspondence
expresses the lift with any polynomial Fock vector at that
same finite datum as a real convolution of the harmonic lift. The normalized descent therefore
vanishes there for every polynomial Fock vector.
For this assertion let \(K=U(n)\) at one real place, with
probability Haar measure, and let \(v\) be the unit highest vector
in the compact type \(\tau\) of the harmonic. Set
\(c(k)=\langle v,\tau(k)v\rangle\), \(d_\tau=\dim\tau\).
Schur orthogonality gives the matrix unit
\[
 p=\int_K d_\tau\overline{c(k)}\,\omega(1,k)\,dk
\]
on the oscillator. The conjugate scalar coefficient \(y\) of
\(\Xi\) satisfies \(R(d_\tau c)y=y\). In the bilinear descent,
changing \(h\) to \(hk^{-1}\) replaces \(c(k)\) by
\(c(k^{-1})=\overline{c(k)}\); hence the descent of \(Q\)
equals the descent of \(pQ\). Take the product over the real places.

The projected joint harmonic is the original highest vector.
Its orbit under the smaller real group has dense linear span
in the matched irreducible theta representation.
A polynomial vector projected to it is compact-type finite.
For its finitely many compact types the image of compactly
supported smooth bump averages of the projected harmonic is
dense and has finite dimension after the compact-type projection.
It is therefore the entire finite-type subspace.

The compact-type
projection can be absorbed in the smooth bump by averaging
over the real maximal compact. This expresses the projected
polynomial exactly as one compact bump average. Moving that
average to the source variable gives the claimed convolution.
The finite vector stays fixed and the compact support justifies
the interchange on the definite target quotient. Trivializing the second character of the
kernel multiplies the descent by a character of \(\det g\),
so leaves its vanishing intact.

Here is the extension from the block datum to every
finite test off the original vacuum places.
Fix a real polynomial and a source point \(g\).
Call the corresponding descent pairing \(B(\Psi,\Xi)\),
linear in \(\Psi\) and conjugate-linear in \(\Xi\).
At \(w\in S\), write a matrix as \(z=(z',z'')^{\mathsf t}\)
in the blocks, with \(z'\) an \(r\)-by-\(r\) block.
Say that the block ratio is allowed at \(w\) if
\(z'\) is invertible over the local field and
\(z''z'^{-1}\) is integral, with zero reduction at
\(w\in S_M\). The block function \(\Psi_{\rm bl}\) on
integral matrices is the indicator of the inclusion
modulo \(w\in S_M\), multiplied at \(w\in Q\) by
\(\chi_w(\det z')\), extended by zero at singular
reductions. It has the original vacuum at the other
places. The conjugate place is prescribed by the split
unitary identification.

Suppose first that \(z\) has allowed block ratio at every
place of \(S\). Take the finite source element \(A\)
whose component at each \(w\in S\) is \(z'^{-1}\)
and which is the identity off \(S\cup\bar S\).
Then \(\Psi_0(x)=\Psi_{\rm bl}(xA)\) takes the value one
at \(z\). A source translation shows
\(B(\Psi_0,\Xi)\) is the original block descent at \(gA\)
times an invertible Weil scalar. Thus it is zero.

The nonzero support of \(\Psi_0\) is one orbit of the
target level, with its character action, in the
\(S\)-coordinates. Indeed after right multiplication by
\(A\), a matrix in the support extends to an integral
invertible block matrix with that matrix as its first
columns. At \(S_M\) these columns are the coordinate
columns modulo the residue place, and at \(Q\) the completed
matrix is block lower triangular with the indicated first
block. These are exactly the two level conditions. On
this orbit all equivariant functions are determined by
their value at the original inclusion multiplied by
\(A^{-1}\).

Let \(C=\prod_{w\in S}(z_w+\varpi_w^{a_w}M_{n,r}(\mathcal O_w))\).
First suppose that the ratio is allowed at every point of \(C\).
Then every possible first block \(z'+\varpi_w^{a_w}B\) is invertible.
In particular every entry of \(\varpi_w^{a_w}z'^{-1}\) lies in the
maximal ideal. Otherwise, for an entry \(u\) of modulus at least one,
take \(B\) to have the single transposed entry \(-u^{-1}\).
The rank-one determinant formula then gives determinant zero.
It follows that \(x'z'^{-1}\equiv I\bmod\varpi_w\) and
\(x''z'^{-1}\equiv z''z'^{-1}\bmod\varpi_w\) for every \(x\in C\).
Thus \(C\) is contained in the orbit supporting \(\Psi_0\). Average
\(\mathbf1_C\) over this same compact target level, including
the conjugate of its level character. The translates have
finite span and factor through a finite quotient.
Equivariance of the Weil action and covariance of \(\Xi\)
give the same value of \(B(\cdot,\Xi)\) before and after
averaging. The averaged function is supported on this orbit.
The average is therefore a scalar multiple of \(\Psi_0\).
It follows that \(B(\mathbf1_C,\Xi)=0\).

Otherwise choose a point \(z\in C\) where the ratio fails at a
place \(w\), and use that point as the center of the same ball. There
exists an integral covector \(l=(\lambda,\lambda')\) with
\(lz=0\), \(\bar\lambda\ne0\), and \(\bar\lambda'=0\)
when \(w\in Q\). If \(z'\) is singular take an integral
multiple of a nonzero covector of its left kernel, normalized
to have a unit entry. Otherwise use a covector of
\((z''z'^{-1},-I)\) whose first part has a coefficient of
modulus at least one at \(S_M\), or strictly greater than
one at \(Q\).

Divide by an entry of largest modulus in
that first part. This gives exactly the asserted integrality
and residue conditions. For each \(s=1+\bar u\bar\lambda\)
in \(H_{\bar\lambda}\), lift \(\bar u\) to an integral
column \(u\) and take the target shear
\[
                 h_s=1+\iota_0(u)l .
\]
It is integral invertible since its determinant reduces
to \(\det s\), and fixes \(z\). Extend it by the identity
at the other split places with the corresponding
conjugate components.

Since \(h_s\) fixes \(z\) and preserves
\(M_{n,r}(\mathcal O_w)\), it fixes each such additive ball. Its level
decomposition is upper block triangular at a chosen
place, and at a twist place the reduction of the
off-diagonal top-right block is zero. By that decomposition,
its action on \(\Xi\) is \(\chi_w(\det s)^{-1}\) times the
action of \(s\) at a chosen place, and just
\(\chi_w(\det s)^{-1}\) at a twist place.

The split Weil action of \(h_s\) on the target variable
preserves \(\mathbf1_C\). Its determinant is a unit, and the
rank-twisted source splitting character is unramified on this
unit group. Thus the determinant multiplier as well as the
modulus factor equals one. Changing the target integration
variable gives \(B(\mathbf1_C,h_s\Xi)=B(\mathbf1_C,\Xi)\).
At a chosen place the sum of these actions kills \(\Xi\)
by \(\mathcal A_{\bar\lambda}\). At a twist place it also
kills \(\Xi\), since \(\det:H_{\bar\lambda}\to k_w^\times\)
is onto with equal fibers, and \(\chi_w\ne1\) on its
exact support. Adding the pairings therefore gives
\(|H_{\bar\lambda}|B(\mathbf1_C,\Xi)=0\).

Every compactly supported locally constant test on the
finite split coordinates is a finite linear combination
of these cosets, each of which is covered by one of the two cases.
This proves the bottom vanishing for all finite
Schwartz--Bruhat factors with the prescribed vacuum off
\(S\cup\bar S\). Fix one of these finite tests \(f_f\), let
\(u\) be the real Schr\"odinger vector, and let \(C\) be a compact set
of source points. The same lattice majorant justifies
\[
 \sup_{g\in C}|B_g(u\otimes f_f,\Xi)|
 \le C_{C,f_f}\,\operatorname{vol}([U(V)])\|\Xi\|_\infty p_j(u)
\]
for a fixed real Schwartz seminorm \(p_j\). In fact the translates on
the two compact sets are bounded by a higher seminorm of \(u\);
summing the resulting inverse power costs a constant independent of
\(u\). For any real Schwartz datum, the Taylor polynomials times the
Gaussian in Section~\ref{sec:four-unitary-real-projections} converge to
its Fock vector in every weighted absolute moment. The coherent-state
bound of Section~\ref{sec:four-unitary-constant-terms} bounds \(p_j\)
of the inverse Bargmann transform by such a moment. Passing to this
limit proves vanishing for that finite test. Every Schwartz--Bruhat
datum with the prescribed off-level vacuum is a finite sum of such
tensors: on the remaining finite support it is constant on the
cosets of a common compact open additive subgroup. This establishes
the all-datum bottom vanishing used in the tower. 

\subsubsection{Extension through the hyperbolic planes}\label{sec:bottom-descent-planes}

\sectionopening{In this subsection, we extend the type of Section~\ref{sec:bottom-descent-real-type} through the original tower.}

Keep the original pair, level, support \(S\) and false lift
\(\Xi\) of the bottom-descent argument above.
Write \(V=I_n=I_{d_1}\oplus I_{d_2}\) for the prescribed gluing
blocks, with \(2r+2d_1+3<n\), and let \(\mathbb H^l\) be the
standard hyperbolic summand of rank \(2l\).
In the rank-\(r\) chain \(n=r+t\), \(t=5r\), and
\(d_1=r+9m+1\). The displayed inequality is then
\(4r+18m+5<6r\), which follows from the original bound
\(100m+9<k\) and \(r+1=k+m\).

Let \(Y\) denote the conjugate of the matrix coefficient
of \(\Xi\) at its original highest vector.
The chosen split place of the false-lift test is \(w\)
with norm \(q\), off \(S\cup\bar S\).
Conjugation of this scalar coefficient commutes with right
translation and conjugates the spectral argument of the
Hecke polynomial. Thus \(Y\) is spherical at \(w\) and
\(\mathcal P_{n,w}(X)Y=0\), where
\[
 |X|=q^{1-n+r/2},\qquad
 \mathcal P_{a,w}(X)
     =\sum_i(-X)^iq^{i(i-1)/2}T_{w,i}.
 \tag{F1}\label{eq:display-f1}
\]
The \(T_{w,i}\) are the original right-coset sums,
using compact mass one.

With \(t'=n-r\), the Frobenius value of the original source
character in Section~\ref{sec:normalized-trace-multiplicities} is
\[
 \nu_w=\alpha_{\rm nor}(\mathrm{Fr}_w)^{-1}
 =q^{r/2}\chi_W(\varpi_w/\bar\varpi_w)\upsilon(\varpi_w)^{-r}.
\]
The auxiliary twist enters \(\nu'_w\). The last factor of the
complementary string gives \(X_w=(q^{t'-1}\nu_w)^{-1}\), the
parameter in Proposition~\ref{prop:clearing-on-original-tests}.
The right-coset sums have real coefficients. On the conjugate
matrix coefficient their polynomial in \eqref{eq:display-f1}
therefore vanishes at \(X=\overline{X_w}\), of the stated modulus.

First suppose that the descents to \(W\oplus\mathbb H^{l-1}\)
vanish for all data vacuum off \(S\cup\bar S\), for the
source character \(a\upsilon^{l-1}\), where \(\upsilon\) denotes
the norm-one restriction and \(a=\chi_W=\upsilon^{r_-}\) for the
negative index \(r_-\) of \(W\). Set
\(W_l=W\oplus\mathbb H^l\), \(r_l=r+2l\).
Split off the last plane, of Gram matrix
\(J_{\mathbb H}=\left(\begin{smallmatrix}0&1\\1&0\end{smallmatrix}\right)\).
For its finite adelic coordinates take the fixed integral trace-one
element \(e=c_{\mathrm{tr}}\) and a unit \(c\) of norm \(-1\)
by Section~\ref{sec:splitting-level-characters}. The matrix
\[
 B_f=\begin{pmatrix}1&-c\\ e&c\bar e\end{pmatrix},
 \qquad B_f^*J_{\mathbb H} B_f=I_2,\qquad \det B_f=c
\]
identifies the hyperbolic lattice with the original finite unitary
lattice, including at dyadic places. At infinity use
\(2^{-1/2}\left(\begin{smallmatrix}1&-1\\1&1\end{smallmatrix}\right)\),
which is orthogonal and carries \(J_{\mathbb H}\) to
\(\operatorname{diag}(1,-1)\). Removing successive planes only
reorders these blocks and leaves the summand \(W\) fixed.

Write \(g=n_0(u,g_0)k\) in Iwasawa form. At a split place the
integral compact moves an isotropic line to the specified line by
extending its primitive vector inside the kernel of its primitive
dual covector to a basis. At a nonsplit finite place the extension
is unramified. For a primitive isotropic vector take an integral
dual vector with pairing one, then subtract a multiple of the
isotropic vector to make the dual vector isotropic, using the
surjective integral trace. These two vectors split off a unimodular
plane. Its orthogonal complement admits an integral orthonormal
basis: lift a vector of nonzero norm over the quadratic residue
field, normalize by surjectivity of the unit norm, and repeat.
A nondegenerate residue Hermitian form has such a vector even in
characteristic two, by nondegeneracy of the residue trace pairing.
These local choices are integral at almost all places and give the
adelic decomposition. Let \(\theta_0\) denote the isotropic zero
slice of the theta functional. The compact
factor carries the vacuum subspace at every place off
\(S\cup\bar S\) to itself, and the radical \(n_0\) acts as the
identity on the zero slice.

The slice of the Levi action
is the Weil action of the smaller pair at \(g_0\), times
\[
 \bigl((1)_{n}(\det u)\bigr)^{-1}
                            |\det u|_{\mathbb A_E}^{n/2}.
 \tag{F2}\label{eq:display-f2}
\]
Here \((1)_n\) is the rank-\(n\) space character associated
to the trivial target automorphic character. Write \(b_d(z)=
b(z/\bar z)\upsilon_E(z)^{-d}\) for the splitting character
of a norm-one character \(b\) at rank \(d\). Since
\(\upsilon_E(\bar z)=\upsilon_E(z)^{-1}\), its restriction
\(\upsilon\) satisfies \(\upsilon(z/\bar z)=\upsilon_E(z)^2\). Hence
\[
 (a\upsilon^l)_{r+2l}(z)
    =(a\upsilon^{l-1})_{r+2l-2}(z)=a_r(z).
\]
Let \(s_0\) denote evaluation at zero in the mixed coordinate
of the removed plane and let \(\omega_j\) denote the Weil action
for \((W\oplus\mathbb H^j,V)\) with characters
\((a\upsilon^j,1)\). For \(g=n_0(u,g_0)k\), the radical has
zero Heisenberg displacement and trivial quadratic phase at this
coordinate. The Levi substitutes \(x\mapsto h^{-1}xu\) and
contributes \eqref{eq:display-f2}. Consequently
\[
 s_0\bigl(\omega_l(g,h)\phi\bigr)
  =(1)_n(\det u)^{-1}|\det u|_{\mathbb A_E}^{n/2}
      \,\omega_{l-1}(g_0,h)\,
           s_0\bigl(\omega_l(k,1)\phi\bigr).
\]
At each finite place off \(S\cup\bar S\), the integral compact
preserves the vacuum line. In the integral plane frame the mixed
vacuum is the integral-box indicator times the smaller vacuum,
up to the fixed nonzero normalization. Evaluation at zero therefore
preserves the prescribed vacuum projection. The scalar is independent
of \(h\), so the integral of this term against \(Y\) vanishes by
the induction hypothesis on the same datum.

The remainder gives a bounded descent on \(U(W_l)\).
Here are the analytic details. Write
\(\Lambda=\theta-\theta_0\) and denote by \(k_0\) the
rational isotropic column of the plane. The zero slice is
the probability mean of the theta functional along the
adelic isotropic translations
\[
 Z=\{z k_0^*J:z\in E^{n\times1}\}
                     \subset\operatorname{Hom}_E(W_l,V).
\]
Each mixed-model character indexed by a rational vector is trivial
on \(Z(E)\), since the global additive character is trivial on
rational arguments. On \(Z(\mathbb A)/Z(E)\) the zero index gives
the constant character, and every nonzero index gives a nontrivial
character by the nondegenerate trace pairing. Their probability
means are therefore one and zero, respectively. Thus \(\Lambda\)
is invariant under the rational translations \(Z(E)\) and has
probability mean zero on every datum.

\subsubsection{The mixed-model height bound}\label{sec:bottom-descent-height}

\sectionopening{In this subsection, we bound the mixed theta sum used in the vanishing argument after Section~\ref{sec:bottom-descent-planes}.}

Let \(\delta\) be the first Euclidean minimum of the
moved integral source lattice, with its original positive
Euclidean archimedean metric.

Put \(H=\max(1,\delta^{-1})\).
If \(\delta<1/4\), the \(E\)-span of the vectors of length
\(<4\delta\) is totally isotropic. Indeed, the pairings of
such vectors are integral at every finite place and have
absolute value \(<1\) at every embedding, so vanish by
the product formula. At each infinite place \(\tau\), with sign
matrix \(J_\tau\), the completion \(I_\tau\) of this isotropic subspace
and \(J_\tau I_\tau\) are Euclidean orthogonal and dual for the signed
pairing. At every \(\tau\) simultaneously, in the Euclidean splitting
into \(I_\tau,J_\tau I_\tau\), and their common orthogonal complement,
multiply the first factor by two, the second by \(1/2\), and the
third by one.

This is a
signed unitary map of norm at most two. It doubles the
lengths of the short vectors and decreases the others
by at most a factor of two. The new first minimum
is therefore at least \(2\delta\). Iteration of the inverse
moves on the adelic point puts its minimum above \(1/4\),
using a purely real product of norm \(O(H)\).
For the resulting point, the lattices in the Euclidean coordinates
have fixed covolume, since \(g\) is unitary, and first minimum
at least \(1/4\). Minkowski's second theorem bounds their
last minimum. Choose a full set of independent lattice
vectors of bounded length. They span a sublattice of bounded
index by Hadamard's inequality. Its possible enlargements,
in the coordinates of these vectors, lie between
\(\mathbb Z^N\) and \(\frac1{D}\mathbb Z^N\), where \(N\) is the
real dimension and \(D\) is a common multiple of the possible
indices. There are finitely many such enlargements. Choosing
a basis in each gives a uniformly bounded \(\mathbb Z\)-basis
of the moved lattice, with uniformly bounded inverse basis
matrix since its covolume is fixed.

In the corresponding rational basis, the matrices for
multiplication by a fixed integral basis of \(\mathcal O_E\)
have bounded integer entries. The signed Hermitian pairings
of the basis vectors lie in the fixed fractional ideal of
pairing values and are bounded at every embedding, so they
also have only finitely many possible values.
Two bases with the same multiplication matrices and
pairings are carried into one another by an \(E\)-linear
unitary map defined over \(F_0\). Fix one rational basis
and one finite lattice representative for each realized
set of these values. After the corresponding rational
left translation, the finite point lies in that representative
times the compact stabilizer of the original finite lattice.
The real point and its inverse are bounded because they
carry the fixed real basis to the bounded moved basis.
This gives a compact set of representatives for all points
with first minimum at least \(1/4\).
The definite target has a positive uniform lower bound
on its first minimum, by integrality of the signed
pairing and the product formula. The same argument
gives a compact set of representatives there.

The theta distribution on a rational Schr\"odinger
polarization is a rational lattice sum. On a fixed finite
Schwartz space it is bounded by a finite-order Schwartz
seminorm at infinity, by summing a decaying power on
the rational lattice. Compact sets of finite translates
of such a space are contained in another fixed finite
Schwartz space. At infinity compact sets of Weil
translates act boundedly on each fixed seminorm, allowing a
higher fixed seminorm on the source vector.

For a general real
move the conjugation formula on the position and
differentiation operators is its linear symplectic action.
Together with unitarity on \(L^2\), this bounds a
fixed Hermite Sobolev norm of the image by a polynomial
in its matrix norm times a fixed Hermite Sobolev norm
of the source vector. The real symplectic action of the pair is
the linear Hom action, whose norm is polynomial in
the two matrix norms. Rational invariance of theta
and the preceding compact reduction consequently give
\[
  |\theta(\omega(g,h)\phi)|\le C H^A p(\phi)
\]
with one fixed Schwartz seminorm \(p\) on the finite
Schwartz spaces in use, uniformly in the definite
target \(h\).

In the short-translation argument each difference is
conjugated past \(\omega(g,h)\), after choosing the rational
representative with the fixed short isotropic vector.
The transformed finite shifts fix the test because its
finite stabilizer contains the fixed ideal used for the
coefficient representatives. The real shifts have size
\(O(H^{-1})\) and are in one isotropic family, so commute.
The fundamental theorem of calculus along a real
Heisenberg one-parameter group bounds a difference by
its length times the seminorm of the corresponding
position or derivative operator. Repeating for \(A\)
or more differences gives a factor \(H^{-A}\) in a
higher fixed seminorm of the original polynomial Fock test.
The following averaging identity applies this estimate directly to
the full theta functional.

If \(H>2\), choose a shortest rational integral source
vector of length at most \(2/H\). Its norm pairing is
integral at finite places and of absolute value less than
one at the archimedean places, so it is isotropic by the
product formula. Witt's theorem moves it to \(k_0\) by a
rational translate \(\gamma\) of the source. Rational
invariance of the full kernel allows the descent at
\(\gamma g\) in place of \(g\); the zero-slice integral at
that point is again zero. At the definite target, a rational
basis in its moved lattice has uniformly bounded lengths:
Euclidean covolume is fixed and the first minimum is
bounded below, so Minkowski's second theorem bounds the
last minimum. With this basis \(d_j\), the maps
\(d_j k_0^*J\) span \(Z\), are integral in the moved Hom
lattice, and have real lengths \(O(H^{-1})\).

Choose a fixed compact open additive ideal in the finite
coefficient coordinates which makes the transformed finite
shifts act trivially on the fixed Schwartz vector. Such an
ideal works uniformly because the moved finite rank-one maps
are integral. Additive approximation gives representatives for
the adelic quotient with finite coefficients in this ideal,
and the archimedean coefficients can then be taken in a fixed
compact fundamental domain for the corresponding rational
coefficient lattice.

Give \([Z]=Z(E)\backslash Z(\mathbb A)\) probability Haar measure.
For \(j\ge1\) and the moved datum \(f\), rational periodicity and
isotropy give
\[
 (\theta-\theta_0)(f)
  =\int_{[Z]^j}\theta\!\left(\prod_{i=1}^j
                     (1-\varrho(z_i))f\right)\prod_{i=1}^j dz_i.
\]
Indeed each nonempty term of the expanded product averages to
\(\theta_0(f)\), since translation preserves Haar measure on
\([Z]\); their signed coefficients sum to \(-1\). The rational basis
\(d_j k_0^*J\) identifies this probability measure with that on the
coefficient quotient, independently of the chosen basis. After conjugation
each finite shift fixes the test and each real shift has length
\(O(H^{-1})\). Integrating its Heisenberg generator gives one factor
\(H^{-1}\) per difference in the Schwartz seminorm bound.

For \(H>2\), at the rational representative just chosen the
difference estimate bounds the remainder by \(C_j H^{A-j}p_j'(\phi)\),
uniformly in the definite target point, for a higher real Schwartz
seminorm and the fixed finite datum. The zero slice contributes zero
after pairing with \(Y\), by the induction hypothesis. At \(H\le2\)
bound the full theta kernel by its initial estimate. Thus for every
integer \(B\ge0\), choosing \(j\ge A+B\) gives the descent bound
\[
 |F_l(g)|\le C_{l,B,f_f}\operatorname{vol}([U(V)])\|Y\|_\infty
                         H(g)^{-B}p_{l,B}(\phi_\infty).
\]
Here the seminorm is on the original real test and the finite factor
is fixed. In particular the boundedness required below holds for
each polynomial datum at this level.

The local Hecke transfer now kills this bounded function.
For a polynomial Fock datum with the required off-level
vacua write \(F_l\) for its descent to \(W_l\),
\(k=n-r_l\), \(\alpha=\mathrm{frob}_{r_l}(a\upsilon^l,\bar w)\),
and \(\beta=\mathrm{frob}_n(1,\bar w)\).
The vacuum at \(w,\bar w\) makes \(F_l\) spherical and
gives the exact adjoint transfer
\[
 0=\prod_{i=0}^{k-1}(1-q^i\alpha X)
          \mathcal P_{r_l,w}(X')F_l,\qquad
                   X'=q^k\alpha\beta^{-1}X .
 \tag{F3}\label{eq:display-f3}
\]
For a target double coset \(D=K d K\), write
\(T_DY(h)=\int_D Y(hs)\,ds\), with \(\operatorname{vol}K=1\).
For the kernel \(H_g(h)=\theta(\omega_l(g,h)\phi)\), right
invariance of the target quotient measure gives the bilinear identity
\[
 \int_{[U(V)]} H_g(h)(T_DY)(h)\,dh
 =\int_{[U(V)]}\!\left(\int_{D^{-1}}H_g(hs)\,ds\right)Y(h)\,dh .
\]
The double coset and the target quotient are compact. Inversion
preserves Haar measure on the unimodular local group. The inverse
coset may therefore be decomposed afresh into right \(K\)-cosets
on the spherical kernel. In the \(\bar w\)-coordinate the split
inverse-transpose relation gives the positive target coset; on the
source it reads the inverse coset there as the positive one at \(w\). The determinant factors on the vacuum have
moduli \(|\alpha|=q^{r_l/2}\), \(|\beta|=q^{n/2}\).
There is no conjugation of \(X\) at this step, since the
kernel integration with \(Y\) is linear.
To check the local rectangular vacuum identity in \eqref{eq:display-f3},
use the chosen coordinate at \(v=\bar w\), with residue
cardinality \(q\). Write \(L=\mathcal O_v^a\). For each
\(U\subseteq L/\varpi_vL\) take an integral matrix \(\gamma_U\)
with elementary divisors \(1,\varpi_v\) and reduction kernel \(U\).
Its transpose and inverse determine the image lattices
\[
 \gamma_U^{\mathsf t}L/\varpi_vL=U^\perp,\qquad
 \varpi_v\gamma_U^{-1}L/\varpi_vL=U .
\]
The first equality is the row-space identity over the residue
field, and the second follows by reducing
\(\gamma_U(\varpi_v\gamma_U^{-1}L)=\varpi_v L\) and comparing
dimensions. Right cosets are determined by their image lattices.
Use \(\gamma_U^{\mathsf t}\) on the large group. On the
small group the inverse double coset uses \(\gamma_U^{-1}\).
Their actions on the vacuum matrix \(Z\in M_{n,r_l}(E_v)\) are
the substitutions \(\gamma_U^{-\mathsf t} Z\) and \(Z\gamma_U^{-1}\),
respectively. Taking \(Y'=Z^{\mathsf t}\), \(a=n\) on the large
side and \(Y'=Z\), \(a=r_l\) on the small side, both tests have
the form \(Y'\gamma_U^{-1}\). In a basis adapted to \(U\),
\[
 Y'\gamma_U^{-1}\text{ integral}
 \quad\Longleftrightarrow\quad
 Y'\text{ integral and }U\subseteq\ker\overline{Y'} .
\]
Indeed the columns in the directions of \(U\) must be divisible by
\(\varpi_w\); the remaining columns must be integral. Thus both
operators vanish off the integral matrices. On an integral matrix
of residue rank \(a_0\), the contributing subspaces lie in kernels
of dimensions \(n-a_0\) and \(r_l-a_0\), respectively. Both base
change characters at \(\bar w\) are unramified, since the original
source character and \(\upsilon\) are unramified there, and the target
character is trivial. A subspace of dimension \(i\) therefore
contributes determinant factor \(\alpha^i\) or \(\beta^i\).
For a residue space of dimension \(d\), the finite binomial identity is
\[
 \sum_{U\subset k_w^d}(-s)^{\dim U}q^{\binom{\dim U}{2}}
                  =\prod_{j=0}^{d-1}(1-q^js).
\]
With \(\beta X'=q^k\alpha X\), the two sums are consequently related by
\[
 \prod_{j=0}^{n-a_0-1}(1-q^j\alpha X)
 =\prod_{j=0}^{k-1}(1-q^j\alpha X)
    \prod_{j=0}^{r_l-a_0-1}(1-q^j\beta X').
\]
This gives \eqref{eq:display-f3} for all \(X\) in the original Haar
normalization.

The absolute values in \eqref{eq:display-f1} and
\eqref{eq:display-f3} give
\[
 q^k|\alpha X|=q^{1-l}\le1,\qquad
 q^{r_l+1}|X'|
      =q^{(r_l+r+4-n)/2}\le1 .
\]
Since \(i<k\), we have \(|q^i\alpha X|\le q^{-1}<1\), so each
factor in that prefactor is nonzero. The second inequality follows from \(2r+2d_1+3<n\)
and \(l\le d_1\). For a bounded spherical function \(F\)
on the source, the absolute sum of the coefficients of
the Hecke right-coset operator \(\mathcal P_{r_l,w}(X')\)
is
\[
 \sum_{U\subset k_w^{r_l}}q^{\binom{\dim U}{2}}|X'|^{\dim U}
  =\prod_{i<r_l}(1+q^i|X'|)
  \le\exp(q^{r_l}|X'|)\le e^{1/2}<2.
\]
The term from \(U=0\) has weight one and acts as the
identity.

Isolating it proves
\(\|F_l\|_\infty\le(e^{1/2}-1)\|F_l\|_\infty\), so
\(F_l=0\). This uses no spectral decomposition of the
noncompact source. With each finite factor fixed, the compact-set
Schwartz seminorm bound and polynomial approximation proved at the
bottom apply also at rank \(r_l\). They extend this pointwise
vanishing to every Schwartz datum with the prescribed off-level
vacuum. Induction starting
with the bottom-descent argument above gives
vanishing at \(l=d_1\). Every conjugate matrix coefficient
of \(\Xi\) has the same vanishing: translates of the highest
vector span the irreducible Weyl space, and a real
translation of the coefficient corresponds to the inverse
real translation of the datum, leaving its finite
projection fixed. 

Next we pass from these descents to the prescribed block
period. Write \(T\) for the complement of \(S\cup\bar S\).
For a continuous function \(Y_0\) on \([U(V)]\) fixed by
the integral group at \(T\), suppose first that its
descents to the split space \(\mathbb H^{d_1}\) with
character \(\upsilon^{d_1}\) vanish for all data vacuum at
\(T\). Fix \(h_1\in U(I_{d_1})(\mathbb A)\) and set
\[
 g_0=\operatorname{diag}(h_1,1),\qquad
 I(h)=\int_{[U(I_{d_2})]}
                 Y_0(\operatorname{diag}(h_1,h_2)h)\,dh_2 .
\]
Let \(o_1\) be the original orthonormal first-block basis.
For \(u\in\operatorname{Mat}_{n\times d_1}(\mathbb A_E)\) use
the original absolute-trace character \(\psi_E\) to define
\[
  x(g)=g^{-1}o_1 h_1,\qquad
  W(x)=|\mathrm{vac}(\alpha_0 x)|^2,\qquad
  \psi_u(x)=\psi_E\!\left(\delta\sum_{i,j}
                             u^{\circ}_{ij}x_{ij}\right).
\]
Take \(j_{\rm W}=-1/\sqrt{-d}\) and the finite idele
\(\alpha_{\rm W}\) of Section~\ref{sec:period-finite-coordinates}. Put
\(\alpha_f=(\alpha_0)_f=\alpha_{\rm W}\) and
\((\alpha_0)_\tau=(2\pi\operatorname{Im}\tau(j_{\rm W}))^{-1/2}\).
The vacuum support in each finite coordinate is
\(\alpha_f\widehat{\mathcal O}_E\), so \(W(x)>0\) precisely when
\(x_f\) is integral. The idele \(\delta=j_{\rm W}\alpha_0\bar\alpha_0\) generates the
inverse different at finite places, and \(u^\circ\) is obtained
by setting the finite components of \(u\) at \(T\) to zero.

Parametrize the Siegel center by the skew block \(C-C^*\),
where \(C\in\operatorname{Mat}_{d_1}(\mathbb A_E)\).
The additive matrix quotient by rational matrices maps onto the
central quotient, carrying probability Haar to probability Haar.
Write \(G_z=z^*Jz\) for the Gram matrix of a rational term.
The central multiplier in the mixed model is
\[
 \psi_E\left(-\tfrac12j_{\rm W}\operatorname{tr}
            (G_z(C-C^*))\right)
       =\psi_E\bigl(-j_{\rm W}\operatorname{tr}(G_z C)\bigr).
\]
Here \(\psi_E\) is the character through the absolute trace,
\(\overline{j_{\rm W}}=-j_{\rm W}\), and
\(\operatorname{tr}(G_z C^*)=\overline{\operatorname{tr}(G_z C)}\).
Multiplying the zero descent by
\(\psi_E(j_{\rm W}\operatorname{tr}(I_{d_1}C))\) and averaging
therefore retains exactly the bases with \(G_z=I_{d_1}\):
\[
 \int \psi_E\bigl(j_{\rm W}\operatorname{tr}((I_{d_1}-G_z)C)\bigr)\,dC
      =\mathbf1_{G_z=I_{d_1}} .
\]
For \(d_1>0\), these matrices are injective. For \(d_1=0\)
there is just the empty matrix.
They are one rational orbit by Witt's theorem, with
stabilizer the original rational \(U(I_{d_2})\).
Fix the Siegel Levi point and take compact sets of representatives
for the definite quotient and for the additive unipotent quotient.
The finite supports of the translated mixed slices impose
\(z\in\mathfrak a^{n d_1}\) for a fixed fractional ideal
\(\mathfrak a\) of \(E\). Their real Schwartz seminorms are
uniformly bounded on these compact sets, and the unipotent
character has modulus one. Thus the sum of absolute values is at most
\[
 C_N\sum_{z\in\mathfrak a^{n d_1}}(1+\|z\|)^{-N}<\infty,
 \qquad N>n d_1[E:\mathbb Q].
\]
The last inequality follows by counting lattice points in dyadic
annuli. The bounded factor \(Y_0\) preserves this majorant.
We may therefore take the coefficient term by term and unfold
the surviving orbit to the quotient by its rational stabilizer
with the original ambient Haar measure.

Read the coefficient at the Siegel Levi element with
matrix \(\alpha_0 h_1\), and use the Heisenberg translate
of the original vacuum by the isotropic vector whose
adapted coordinates are \(\alpha_0\overline{u^\circ}\).
It is vacuum at \(T\). In the integral hyperbolic frame
used above, the finite vacuum on the position coordinates
is the indicator of \(\alpha_{\rm W}\widehat{\mathcal O}_E^{nd_1}\)
times its fixed nonzero normalization. Its opposite lattice
is dual by the trace identity of
Section~\ref{sec:period-finite-coordinates}, including above \(2\).
At an infinite place the orthogonal plane rotation gives
a positive radial Gaussian. On the frame orbit \(x^*x=I_{d_1}\)
its value and the local factor of \(W(x)\) are fixed
positive constants, since \(\|x\|^2=d_1\) there.

For the pure momentum translation the Schr\"odinger
multiplier is the full character of the displacement paired
against the position. Here its argument is
\[
 j_{\rm W}\operatorname{tr}
    \bigl((\alpha_0\overline{u^\circ})^*(\alpha_0 x)\bigr)
       =\delta\sum_{i,j}u^\circ_{ij}x_{ij}.
\]
Thus the translated slice on this orbit is a common
nonzero scalar times \(W(x)\psi_u(x)\).
The target action on the position frame is \(x\mapsto g^{-1}x\)
because the pure hyperbolic source character satisfies
\((\upsilon^{d_1})_{2d_1}=1\).
The remaining Levi scalar is the nonzero factor
\eqref{eq:display-f2} for \(\alpha_0 h_1\), independent of the frame.
For \(h_2\in U(I_{d_2})(\mathbb A)\),
\(\operatorname{diag}(h_1,h_2)g_0^{-1}g
=\operatorname{diag}(1,h_2)g\).
The original Haar disintegration therefore gives the inner
period \(I(g_0^{-1}g)\). Canceling the common scalars, the
vanishing coefficient gives
\[
 \int_{U(I_{d_2})(\mathbb A)\backslash U(V)(\mathbb A)}
          W(x(g))\psi_u(x(g))I(g_0^{-1}g)\,d\dot g=0
                           \qquad\text{for every }u .
 \tag{F4}\label{eq:display-f4}
\]
The quotient measure is the quotient of the original
group Haar measures, and the inner integral defining
\(I\) uses the original adelic quotient Haar measure.

The positive weight in \eqref{eq:display-f4} detects the period at the
chosen basis. On its support the finite basis is integral
and orthonormal. Integral Witt extension at each finite
place extends this basis to an integral unitary matrix.
At a split place this is completion of a direct summand.
At an inert place the basis splits off a unimodular orthogonal
summand. The integral orthonormalization in
Section~\ref{sec:period-finite-coordinates} identifies its complement
with the standard one, including in residue characteristic two.

The finite nonsplit extensions are unramified by
Section~\ref{sec:different-additive-character-scales}. The residue
and integral trace maps are surjective also in residue characteristic
two, as in Section~\ref{sec:splitting-level-characters}. At infinity
orthonormal frames extend by compact unitary transformations.

More precisely, put \(G=U(V)\), \(H_2=U(I_{d_2})\), and
\(K=G_\infty\prod_{v<\infty}G(\mathcal O_v)\) for the standard
lattice. When \(W(x(g))>0\), choose \(k\in K\) with
\(x(g)=k^{-1}o_1\). Then \(gk^{-1}o_1=o_1h_1\), whence
\[
 g\in H_2(\mathbb A)g_0K.
\]
The support of the integral in \eqref{eq:display-f4} is therefore
contained in the compact image of \(g_0K\) in
\(H_2(\mathbb A)\backslash G(\mathbb A)\). The inner period is
unchanged by left multiplication of \(g\) by \(H_2(\mathbb A)\),
by right invariance of its quotient measure. Two integral
representatives with the same characters \(\psi_u\)
have identical basis coordinates off \(T\), by nondegeneracy
of the additive pairing. Multiplying on the right by
the \(T\)-component of their relative integral matrix
then identifies their full bases and leaves \(Y_0\) fixed.
Thus \(\overline I\) is a continuous function of the array
of these characters on the compact character image.

On that compact image the span of the characters separates points,
contains \(1=\psi_0\), and satisfies
\(\psi_u\psi_v=\psi_{u+v}\) and \(\overline{\psi_u}=\psi_{-u}\).
By Stone--Weierstrass there are finite character sums \(p_j\)
converging uniformly to \(\overline I\). Equation~\eqref{eq:display-f4}
gives \(\int W(x)p_jI=0\), and
\[
 \left|\int W(x)(p_j-\overline I)I\right|
 \le\|p_j-\overline I\|_\infty\int W(x)|I|\longrightarrow0
\]
by the unfolded absolute convergence. Thus \(\int W(x)|I|^2=0\).
Continuity and full support of quotient Haar measure
give \(W(x)|I|^2=0\) everywhere. At \(g=g_0\) the basis
is \(o_1\), which has strictly positive vacuum weight.
Hence \(I(1)=0\), the desired block period.
This argument uses the chosen basis and complement
throughout. 

\subsubsection{Vanishing for the gluing adjoint}\label{sec:gluing-block-period-vanishing}

\sectionopening{In this subsection, we apply the vanishing calculation of Section~\ref{sec:bottom-descent-height} to the original false-lift coefficient.}

Allow a finite right translate \(\gamma\) of \(\Xi\) supported
on \(S\), as in the prescribed divided lattice.
The proved vanishing for \(W\oplus\mathbb H^{d_1}\) remains
valid after this translate because the transported datum
is unchanged at \(T\). Let \(\phi_2\) be the actual
\(W,I_{d_2}\) datum of the slice, with target character
\(b\), and fix its source point \(g\).

Choose a \(W,I_{d_1}\) datum \(\phi_1\) vacuum at \(T\)
whose kernel, with target character \(b^{-1}\), is
nonzero at \((g,h_1)\). Choose a rational Schr\"odinger polarization with position
space \(D\) over \(F_0\), and take the original finite vacuum.
Its finite Weil translate at \((g_f,h_{1,f})\) is a nonzero locally
constant Schwartz function \(f_f\) on \(D(\mathbb A_{F_0,f})\).
Additive strong approximation gives \(x_0\in D(F_0)\) with
\(f_f(x_0)\ne0\). In rational coordinates the compact support
of \(f_f\) lies in the profinite completion of a fractional
lattice. The rational points in that support therefore form a
discrete set in \(D(F_{0,\infty})\). Take a smooth real bump
supported in a neighborhood meeting this set only in \(x_0\),
with value one there, and pull it back by the inverse real Weil
action at \((g_\infty,h_{1,\infty})\). The rational theta sum
then has exactly one nonzero term, with value \(f_f(x_0)\) up to
the nonzero normalization of the distribution. This choice
keeps the finite vacuum at \(T\).

For a conjugate scalar coefficient \(Y\) of \(\Xi\), set
\[
 Y_0(h)=Y(h\gamma)\,
     \theta^{a,1}_{W,V}(\phi_1\boxtimes\phi_2;g,h),
 \qquad a=\upsilon^{m-j},\quad 1=b^{-1}b.
\]
This is continuous and left rationally invariant. If \(k\) belongs
to the original integral group supported on \(T\), its support is
disjoint from that of \(\gamma\). Hence \(k\gamma=\gamma k\) and
\(Y(hk\gamma)=Y(h\gamma)\). Both finite tests are vacuum at \(T\)
and \(a\) is unramified there. Thus \(\omega(1,k)\) fixes the tensor
at \(T\), also after the commuting source action at any \(g\).
Consequently \(Y_0(hk)=Y_0(h)\), as required in
\eqref{eq:display-f4}.
The split-space descents of \(Y_0\) vanish, since the
product with the \(\mathbb H^{d_1},V\) kernel is the
\(W\oplus\mathbb H^{d_1},V\) kernel with the box datum
and source character \(a\upsilon^{d_1}\).
On the block \(I_{d_1}\oplus I_{d_2}\), the \(W,V\)
kernel factors as the product of the two chosen kernels.
Its first factor is the fixed nonzero value at \((g,h_1)\).
Cancel it from the block period and conjugate to obtain
zero for the pairing of the original translated
matrix coefficient with the conjugate \(W,I_{d_2}\)
kernel. To apply this to the vector appearing in the gluing
map, take the first-block class representative with real component
one. Let \(b_t\) be the negative tableau polynomials, and write
\(b(U X)=S_-(U)b(X)\) for the column of these polynomials.
For a real second-block point \(u\) and finite point \(h_f\), put
\[
 \theta_t(u)=\theta_{b_t}(g;(u,h_f)),\qquad
 v(u)=(\overline{\theta_t(u)})_t .
\]
In the target action the character of weight \(m-j\) cancels the
vacuum determinant of the \(m-j\) negative source directions.
The negative matrix variables undergo the substitution
\(X\mapsto u^{-1}X\). The entries of the tableau substitution matrix are
polynomials with real coefficients in its matrix argument. Thus with the conjugate Weyl action \(\rho_-(u)=S_-(\bar u)\) this gives
\[
            v(u)=\rho_-(u^{-1})v(1).
\]
These identities hold place by place, also for the constant
polynomial when \(j=m\).

For the fixed first factor let \(C(v)=\operatorname{Cartan}(y_1,\bar v)\).
In the coefficient action of
Section~\ref{sec:cofinite-spherical-galois}, each complemented
minor of height \(n-h\) contributes a determinant and an
inverse-transpose minor of height \(h\). At \(\tau\) the exponent
\(-\sum_h\mu^-_\tau(h)\) cancels the determinants. The first
block of \(\operatorname{diag}(1,\bar u)\) is the identity,
and the second inverse-transpose block is
\((\bar u)^{-\mathsf T}=u\). Thus
\[
 \rho_V(\operatorname{diag}(1,u))C(v)=C(\rho_-(u)v).
\]
Write \(\Xi_f=\Xi(\operatorname{diag}(h_{1,f},h_f)\gamma)\)
and let \(e_t\) be the tableau coordinate vectors. The two
conjugate-linear maps, \(C\) and the first variable of the pairing,
give
\[
 \langle C(v(1)),\Xi_f\rangle
   =\sum_t\overline{\theta_t(u)}\,
     \langle\rho_V(\operatorname{diag}(1,u))C(e_t),\Xi_f\rangle .
\]
Indeed \(\rho_V(\operatorname{diag}(1,u))C(v(u))=C(v(1))\).
Each term on the right is the conjugate theta kernel times an
original scalar coefficient at the prescribed finite translate.
Its measured block period vanishes by the preceding calculation.

Write \(\varpi=\pi_M\). The finite function is fixed by
\(J(\varpi^2)\).
On the kernel vector, a right integral action
precomposes the residue function by the inverse
of the matrix reduction at that level. This reduction
is the identity, and the determinant factor in the
kernel transformation is one for these integral
translates, since the source splitting character is
unramified there. Conjugation by a divided lattice point with
\(\varpi\mathcal O^n\subset L_i\mathcal O^n\subset\mathcal O^n\)
loses at most one factor of \(\varpi\) in the conjugated
difference from the identity.

Hence it is the identity
modulo \(\varpi\) at the chosen places. At a twist
place the divided lattice point is the identity and the
second block fixes the original first-block covectors,
with level character one. This verifies the exact
right invariance in the original finite level.
Normalize the maximal finite compact and the compact real
group to have mass one, or use the same common positive
constant on both sides.

Each class at \(J(\varpi^2)\) then
has measure \( [K_{\rm std}:J(\varpi^2)]^{-1}/|\Gamma_c|\).
The stabilizers \(\Gamma_c\) are trivial by the original
principal level trace argument in Section~\ref{sec:ramified-source},
which applies in rank at most \(r+t\). Thus the vanishing
measured integral is the plain sum over the second-block
classes of the given gluing map times one positive constant.
For a second-block class \(c\), let \(h_c\) be its finite
representative, and put \(\kappa(g,h_c)\) equal to the vector of
\emph{conjugates} of the second-block theta kernels in the negative
tableau basis. Here the two characters are
\(\upsilon^{m-j},\upsilon^{m/2-j}\) and the residue test is the
inclusion modulo \(\varpi^2\). Write \(\eta_{i,g}\) for the tuple
supported on choice \(i\) with value \(\overline{\kappa(g,h_c)}\)
at \((i,c)\). We use the inner product linear in the second variable.

\begin{lemma}[Vanishing against the gluing adjoint]
\leavevmode\label{lem:marked-gluing-period-vanishing}
For the gluing map \(G_\vartheta\) of
Section~\ref{sec:original-comparison-tests}, each false lift \(\Xi\)
at this twist satisfies, for every \(i,g\),
\[
 \overline{K^\circ_{(G_\vartheta^*\Xi)_i}(g)}
    =\langle\eta_{i,g},G_\vartheta^*\Xi\rangle=0 .
\]
Here \(K^\circ\) is the component pairing of
Section~\ref{sec:v2-measured-gluing-adjoint}.
\end{lemma}

For a fixed class \(c\), let \(b_a\) be the standard tableau
polynomials in the negative coordinates and put
\(l_{c,a}=\theta_{b_a}(g,(1,h_c))\) with the given finite test.
At each infinite place write \(s_1,\ldots,s_{m-j}\) for the
negative source indices in the prescribed reverse order. A tableau
column \(1\le i_1<\cdots<i_h\le d_2\) then gives the minor
\[
    \det\bigl(Z_{i_\alpha,s_\nu}\bigr)_{1\le\alpha,\nu\le h}.
\]
Here \(h\le m-j\) by Section~\ref{sec:bottom-descent-real-type}.
The substitution uses distinct source indices and retains all
\(d_2\) target indices, so it is injective on the polynomial
ring in the original matrix entries. Thus the standard tableau products
form a basis of the negative test polynomials. When \(j=m\), the
diagram is empty and its basis is \(1\). The real positive definite
Gram matrix \(G_-\) is invertible. By the tableau expansion in
Section~\ref{sec:original-comparison-tests},
\(\sum_a(G_-\overline{\zeta_c})_a b_a\) is the Weyl
polynomial of \(\overline{\zeta_c}\), and these polynomials
span exactly the specified tests as \(\zeta_c\) varies.
Thus \(\kappa(g,h_c)=\overline{l_c}\), and
\[
 \overline{K^\circ_{\zeta}(g)}
   =\sum_c\overline{l_c}^{\mathsf t}G_-\zeta_c
   =\sum_c\langle l_c,\zeta_c\rangle_{G_-}.
\]
This is the asserted first equality at
\(\zeta=(G_\vartheta^*\Xi)_i\). Write
\(y_i(h)=\overline{c''\Upsilon_{\vartheta,i}(h)}\) for the first
factor constructed below. By adjunction the inner product is
\[
 M^{l_1}\!
 \sum_{[h]\bmod J_\Lambda}
 \sum_c \left\langle
  \operatorname{Cartan}\bigl(y_i(h),\overline{\kappa(g,h_c)}\bigr),
  \Xi\bigl(\operatorname{diag}(h,h_c)L_i\bigr)
 \right\rangle_{\rm W}.
\]
The stabilizer weights are one at these free levels. For each
fixed \(h\), the inner sum is the vanishing block period just proved.

\subsection{The lattice for the comparison operator}
\label{sec:marked-comparison-lattice}

\sectionopening{In this subsection, we decompose the target lattice into the two Hermitian summands used
in Section~\ref{sec:original-comparison-tests}. The calculation
identifies the subgroup and quotient of the original level needed for
the tensor contraction.}

On the \(j\)-th pair \(W_j,V_j\) the target splits
in its prescribed coordinate order as \(V''\oplus V^\circ\), of
ranks
\[
 d_+=r+9m+1,\qquad d_-=t-9m-1 .
\]
The first \(r\) coordinate indices of \(V''\) label the underlying
basis of \(W_j\). The same index set is used on the source of the
exchange to form
\[
 W'=W_j\perp I_j^-\perp I_{3m-j}\perp W^\flat,
 \qquad W^\flat=I_{6m+1}.
\]
The signs in this decomposition belong to \(W'\). All coordinate
directions in the first target factor \(V''=I_{d_+}\) are positive.

The second factor uses
\(V^\circ\) with its original level \(J(\pi_M^2)\).
Let \(I_e,I_z\) denote the inclusions of \(r\) columns
in the fixed coordinates of \(V'',V^\circ\),
respectively. Let \(a_0=3m\) be the number of rows of \(\Lambda\), the rank
of the auxiliary Hermitian space of signature \((2m,m)\). The residue choices are matrices
\(\Lambda\in M_{a_0,d_+}(\mathcal O_E/(\pi_M^2))\)
and a level-\(\pi_M^2\) class \(h^\flat\) of \(U(W^\flat)\)
with its chosen short rows. For each choice use the level \(J_\Lambda\)
on \(V''\) that is \(1\) modulo \(\pi_M\),
fixes \(\Lambda\) modulo \(\pi_M^2\), and has the
identity columns off \(W_j\) modulo \(\pi_M^2Q\).
Thus \(J_\Lambda\) contains the indicated principal
level and is contained in the original free
chosen level.

Here is the matrix of the lattice point in this
description. It suffices to construct it at each
place dividing \(\pi_M\), in its split matrix realization, and extend by inverse conjugate transpose.
Put \(\pi=\pi_M\) locally and lift \(\Lambda\) to an
integral matrix. Let \(D\) be a basis matrix for the
row lattice generated by the rows of
\(\binom{\Lambda}{\pi I_{d_+}}\).
This row lattice contains \(\pi\mathcal O^{d_+}\)
and lies in \(\mathcal O^{d_+}\), so both \(D\) and
\(B_+=\pi D^{-1}\) are integral. The equality of
the row lattices gives, for all matrices \(X,Y\)
of compatible size, the equivalence modulo \(\pi^2\)
\[
 D X=D Y
 \quad\Longleftrightarrow\quad
 \binom{\Lambda}{\pi I_{d_+}}X
      =\binom{\Lambda}{\pi I_{d_+}}Y .
\]
Write \(J_z=I_zI_z^{\mathsf t}\).

Set
\[
 A=\begin{pmatrix}D&-D I_e I_z^{\mathsf t}\\
                 0&\pi J_z+1-J_z\end{pmatrix},
 \qquad
 C=\begin{pmatrix}B_+&I_e I_z^{\mathsf t}\\
                 0&J_z+\pi(1-J_z)\end{pmatrix}.
\]
The product is \(CA=\pi I\). Permute the rows of \(A\) and the source columns of \(C\)
inversely, so that the bottom \(z\)-rows of \(A\)
are the prescribed inclusion rows of the original
pair. Take \(L_i=C\) at this place, and
use the identity away from the chosen pairs.
In particular \(L_i\) is chosen using the old
modulus only and is independent of the twist.

To verify the exact block congruence, write \(x_+,x_-\)
for the blocks of \(x\) and put \(x_z=I_z^{\mathsf t}x_-\).
The condition \((\pi L_i^{-1})x=\pi I_{\rm incl}\)
modulo \(\pi^2\) says
\[
       x_z\equiv I\pmod\pi,\qquad
       (1-J_z)x_-=0,\qquad D(x_+-I_e x_z)=0
                                      \pmod{\pi^2}.
\]
The first two conditions determine the unique invertible matrix
\(Y=x_z^{-1}\) modulo \(\pi^2\) with \(Y\equiv I\)
modulo \(\pi\) for which \(x_-Y=I_z\).
Then the last condition is \(D(x_+Y-I_e)=0\).
By the row-lattice equivalence it is exactly
\[
       x_+Y\equiv I_e\pmod\pi,\qquad
       \Lambda x_+Y=\Lambda I_e\pmod{\pi^2}.
\]
The implications reverse by multiplication by \(Y^{-1}\).

Thus the original inclusion indicator equals the sum
over \(Y\equiv I\bmod\pi\) of the two block indicators.
At every residue point this sum has zero or one term.
The character of a member of this congruence group
is trivial, since the original determinant twist factors
through reduction modulo \(\pi\). When the class periods
are assembled, the right changes by these matrices
give equal integrals on the original projected source,
which is fixed by this old congruence group. Consequently
the unnormalized sum in that period contributes its
actual order, once.

For the denominator estimate, a change of columns in
\(V''\) puts the kernel of \(\Lambda\bmod\pi\) on
all but at most \(a_0\) coordinates. In these coordinates
the corresponding columns of \(D\) are zero modulo
\(\pi\), by the same equality of row lattices. Thus
\(L_i^{-1}=\pi^{-1}A\) is integral on every
first-block basis vector except at most \(a_0\);
on the second block all \(d_-\) vectors may be divided
once.

Cofactor matrices at the conjugate places are
given by the unitary inverse transpose convention, and
the integral inverse there follows from the displayed
integral \(C\). For a height-\(s\) column of the first
or second Weyl diagram, the exterior minor contains
at most \(a_0\) or \(d_-\) divided variables,
respectively. The maximal heights in the two diagrams are at most
\(a=k+j-1<d_+\) and \(b_-=m-j<d_-\), so no full-height columns
have been omitted. This gives the common clearing exponent
\[
 l_1=\sum_\tau\left(a_0\sum_{s<d_+}\mu_\tau^+(s)
                  +d_-\sum_{s<d_-}\mu_\tau^-(s)\right).
\]
Since \(\pi_M\bar\pi_M=M\), the rational integer
\(M^{l_1}\) clears these denominators at every coefficient embedding.
To check the Cartan map on Weyl vectors, take the positive exterior
columns in \(V''\) and the dual negative columns in \((V^\circ)^*\),
included in the dual of the whole space by extension by zero on \(V''\).
Order the positive basis first and the negative basis last, in the
orders of Section~\ref{sec:bottom-descent-real-type}. For each negative
column use the integral exterior isomorphism
\[
 \bigwedge^s V^*\simeq\bigwedge^{n-s}V\otimes(\det V)^{-1},
 \qquad n=d_++d_- .
\]
With the volume signs, the highest tensor of the two blocks now has
the leading wedge in every column of the target diagram. The entire
tensor product of the two block Weyl spaces lies in the target Weyl
span: it is spanned in characteristic zero by the orbit of this tensor
under the two independent block linear groups. Hence the normalized
Cartan projection is the identity on these included tensors. The signed
complement maps preserve the integral exterior coordinates, and the
target tableau embedding has an integral left inverse by
Section~\ref{sec:integral-terminal-branching}. After transport by
\(L_i^{-1}\) the tensor remains in the same target Weyl span by
equivariance. On a positive column the transport is exterior power of
\(L_i^{-1}\), and on a dual negative column it is exterior power of
\(L_i^{\mathsf t}\). This last matrix is integral at the chosen split
place, and at the conjugate place its restriction to an \(s\)-fold
wedge has denominator at most \(\bar\pi^s\), with \(s\le d_-\).
Together with the positive-column bound this proves the asserted
\(M^{l_1}\)-clearing in the actual target Weyl coordinates. For a unitary matrix \(h\) in the definite standard bases,
\(h^{-\mathsf t}=\bar h\). Thus in the lattice norm at \(h\) the
second Cartan factor acts by
\((h^{-1})^{-\mathsf t}=h^{\mathsf t}=\bar h^{-1}\).
Its input norm is therefore computed on the conjugate vector at
\(\iota\circ\mathrm{conj}\) when the first and output norms are
computed at \(\iota\). This also applies at \(\ell\mid M\).

\subsection{A common integral scale for the first factors}
\label{sec:marked-first-factor-bounds}

\sectionopening{In this subsection, from the lattice decomposition of
Section~\ref{sec:marked-comparison-lattice}, we bound denominators and
archimedean norms of the first factors on one integral scale. These
bounds are uniform over all the finite classes in the subsequent
comparison.}

Write \(d_K\) for the positive quadratic parameter \(d\) of the
datum. In this subsection the letter \(d\) denotes a polynomial degree.
Fix the original \(j\) and the gluing choices of
Section~\ref{sec:marked-comparison-lattice}.
The space \(V''\) has rank \(d_+=r+9m+1\).
Use the same modulus \(\mathfrak m_\vartheta=\varpi^2Q\)
and the original split residue coordinates.
The first factor is a vector \(\overline{c''\Upsilon_{\vartheta,i}(h)}\)
of weight \(\mu^+\) at \(h\in U(V'')(\mathbb A_f)\),
where \(\Upsilon_{\vartheta,i}\) is the left class series
on the original levels. The scalar \(c''\) is
independent of \(\vartheta,i,h\).

We use the following normalized definite kernel.
Let \(K(P,\Psi;g,h)\) be the kernel on the standard definite
pair \((I_l,I_n)\) of positive ranks \(l,n\), at finite group
points \(g,h\) with trivial
archimedean components, the prescribed algebraic infinity types and
tame characters, and a polynomial \(P\) homogeneous of degree \(d\). Put
\[
 s_{l,n,d}
   =(M^2)^{2f ln}\,s_{\rm vac}^{ln}
                  (M^2d_Ks_{\rm jet})^d .
 \tag{I1}\label{eq:display-i1}
\]
Here \(s_{\rm vac}\) is the rank-one vacuum scale,
\(s_{\rm jet}\) the jet scale, and \(d_K\) the
positive integer of the original quadratic field.
All factors are nonzero. Let \(q_Q=\operatorname{Nm}Q\).

Let \(\mathfrak m\) be the split modulus of the individual kernel.
If \(q_QM^2\in\mathfrak m\), the residue
values of \(\Psi\) are integral away from \(q_Q\), and
\(\ell\nmid q_Q\), then at integral points above \(\ell\)
one has
\[
           |s_{l,n,d}K(P,\Psi;g,h)|_\iota\le1
 \tag{I2}\label{eq:display-i2}
\]
whenever \(P\) has integral coefficients at \(\iota\).
Indeed the residue datum applied to the integral vacuum
is a finite Fourier expansion in translates of that
vacuum, with translation coordinates integral after
multiplication by \(q_QM^2\).

The denominators of Fourier inversion divide
\(\lvert\mathcal O_E/\mathfrak m\rvert^{ln}\),
which divides \(\operatorname{Nm}((q_QM^2))^{ln}\) since
\((q_QM^2)\subseteq\mathfrak m\). As \(\ell\nmid q_Q\),
their \(\ell\)-parts are cleared by
\((M^2)^{2fln}=\operatorname{Nm}((M^2))^{ln}\).
At \(\iota\) the same unit \(q_Q\) can be removed from the
translation denominator. The resulting scaled datum belongs
to the local module of standard data at level \(M^2\). The two integral
group points preserve it by the finite Weil action
and the prescribed norm-one character weights.
The integral finite transport of Section~\ref{sec:seesaw-integral-tensor}
and the ordinary derivatives of Section~\ref{sec:critical-integral-pairings}
give integrality at scale \(M^2d_Ks_{\rm jet}\) on that module. Applying the homogeneous polynomial
gives precisely the scale \eqref{eq:display-i1}, by the homogeneous jet formula in
Section~\ref{sec:marked-integral-sections-ordinary-jets}.

We also use the conjugate size part of the same jet
estimate. For the prescribed characters suppose that \(l,n\le\Delta\),
\(d\le M\), and that the split modulus \(\mathfrak m\) contains
\((M^2)\). A bound \(A_0\) for all
conjugates of the residue values and a bound \(B_0\)
for the conjugate polynomial on the polydisc
of radius \(M^4\) give
\[
   |\sigma(s_{l,n,d}K(P,\Psi;g,h))|
               \le e^{x_1}A_0B_0 .
 \tag{I3}\label{eq:display-i3}
\]
Each normalized finite Fourier coefficient is an average of the
residue values times roots of unity. Its conjugate absolute value is
at most \(A_0\). The total is therefore at most
\[
 |\mathcal O_E/\mathfrak m|^{ln} A_0
 \le (M^2)^{2fln}A_0,
\]
since \((M^2)\subseteq\mathfrak m\) and \([E:\mathbb Q]=2f\).
At the unit ideal the finite group has one element and the constant
datum gives the vacuum itself.

\subsubsection{Bounds for the normalized derivatives}
\label{sec:first-factor-derivative-bounds}

\sectionopening{In this subsection, we apply the derivative estimates of Part~\ref{part:upper} to the first factors of the original lift.}

The normalized derivative estimate applies to each
translated vacuum and to every finite point of the
definite groups. The usual homogeneous rescaling
uses
\[
 \frac{(M^2d_K)^2}{\Omega}(2d+2)\le (M^4)^2
\]
in its coefficient ball, where \(\Omega\) is the original
positive period parameter. At radius \(R\) in the squared
norm of that ball, homogeneity bounds the polynomial
by \(B_0 e^{dR}\), because \(\sqrt R\le e^R\).
The jet exponential factor is
\[
 \exp\bigl(C_0(q_1ln+1)^2+fln\log(2d+2)\bigr)
\]
in the original numerical notation. Multiplying by
the Fourier cardinality and the denominator-clearing
factor gives at most \(e^{x_1}\) by the cutoff inequality.

Here is the polynomial to which to apply \eqref{eq:display-i2}--\eqref{eq:display-i3}.
Write \(w_{\mathrm{off}}=a/x_1\) for the weight at the
non-\(b\)-places. The parameter bounds give \(w_{\mathrm{off}}\le b\).
The total determinant count in the specified space is bounded by \(q_*\), where
\[
 q_*=f(b+w_{\mathrm{off}}+\Delta+m\Delta).
\]
Let \(Y\) be the square matrix of short covectors used for the
\(W^\flat=I_{6m+1}\) class of \(h^\flat\) in choice \(i\).
The definite lattice \(L=E^{6m+1}\cap h^\flat
\widehat{\mathcal O}_E^{\,6m+1}\) is Hermitian self-dual. Its Euclidean
first minimum is at least one, since the squared norm
\(\operatorname{Tr}_{F_0/\mathbf Q}\langle v,v\rangle\) of a nonzero
\(v\in L\) is a positive integer. Apply the transference bound in
Section~\ref{sec:peak-hom-lattice} and take the columns of \(Y^{\mathsf t}\)
to be the resulting \(E\)-basis in \(L\). Then
\[
 (h^\flat)^{-1}Y^{\mathsf t}\in
       M_{6m+1}(\widehat{\mathcal O}_E),\qquad
 \max_{\tau,a,b}|\tau(Y_{ab})|
      \le 6d_Kq_1\bigl(f(6m+1)\bigr)^2\le M^4 .
\]
This one choice gives the finite integrality and the bound at every
complex embedding used below.
The product of stacked determinants consists of constant
covectors, some given by \(Y\), over variable covectors.

Separate
the \(W_j\) positive variables from the ones on
\[
 W^{\mathrm{add}}=I_{3m-j}\oplus W^\flat\oplus I_{3m}.
\]
Let \(Z\) denote the positive variables of \(W_j\), let \(b_s(Z)\)
be the standard tableau products of minors in those variables, and
let \(G_{\rm W}\) be the Gram matrix of the corresponding integral
Weyl basis. For the stacked determinant product \(P\), write
\[
 P=\sum_s d_s(Y)b_s(Z),\qquad d(Y)=G_{\rm W}\zeta(Y),
\]
where \(\zeta(Y)\) gives its Weyl coordinates over the
coefficient field.
Both the coefficients of this expansion and the entries of
\(\zeta(Y)\) are polynomials in the added variables. The number of
variable columns remaining in a determinant attached
to height \(h\) is \(d_+-h\). For \(1\le h<k-1\) the initial parameter entries satisfy
\(\varsigma_{h-1}-\varsigma_h=1\), hence the formula for
\(k^+(h)\) in Section~\ref{sec:bottom-descent-real-type}
gives \(\mu^+_\tau(h)=0\). Also \(\mu^+_\tau(0)=0\) by its definition.
Thus \(d_+-h\le d_+-(k-1)=10m+1\) at every contributing height,
using \(r=k+m-1\). Their common homogeneous degree is
\[
 d=\sum_{\tau,h}(d_+-h)\mu^+_\tau(h)
           \le(10m+1)q_*\le M .
\]
The constant rows have entries \(0,1\), or \(\tau(Y_{ab})\)
for one of the chosen complex embeddings \(\tau:E\hookrightarrow\mathbb C\).

Write \(R_Y=\mathbb Z[\tau(Y_{ab}):\tau,a,b]\) and denote the added
variables by \(X_{\mathrm{add}}\). Expanding first along the fixed
unit rows and then along the remaining \(Z\)-columns expresses \(P\) as a
sum of products of
column minors in \(Z\) with coefficients in \(R_Y[X_{\mathrm{add}}]\).
In the transposed diagram convention of
\cite[Theorem~II.2.16, p.~232]{ABWSchur}, the exterior factors are our columns. The Schur map alternates their
entries and multiplies along our rows. Its standard tableaux give a
basis over \(\mathbb Z\), so straightening these products gives
\(d_s(Y)\in R_Y[X_{\mathrm{add}}]\). We control the denominator of
\(\zeta(Y)=G_{\rm W}^{-1}d(Y)\) by the complementary Gram form below.
\subsubsection{Inverse Gram matrices in the chosen coordinates}
\label{sec:first-factor-inverse-gram}

\sectionopening{In this subsection, we pass from the derivative bounds of Section~\ref{sec:first-factor-derivative-bounds} to the integral Gram scale.}

For the complementary shape let \(P^\vee\) be the product over
\(\tau,h\) of the leading minors of size \(d_+-h\), each to the
power \(\mu^+_\tau(h)\). Put
\[
        H=P^\vee(\partial)P^\vee(0)
           =\sum_\alpha\alpha!\,\bigl([X^\alpha]P^\vee(X)\bigr)^2 .
\]
Complementation in exterior powers, with reversed
order of the covectors and application of the determinant tensor,
is an integral isometry between the Gram spaces
for the complementary shapes. On one exterior column of rank
\(d_+\), put \(J_h e_I=\operatorname{sgn}(I,I^c)e_{I^c}\). The
cofactor identity is
\[
 J_h\bigwedge^h g=(\det g)\bigwedge^{d_+-h}(g^{-\mathsf t})J_h .
\]
Reorder the complementary columns by decreasing height. The leading
column wedges go to the final coordinate wedges, which are the images
of the leading complementary wedges under the reversing permutation,
up to sign. Thus the signed tensor isometry carries the rational span
of the highest orbit onto its complementary span. Integral tableau
straightening identifies both Weyl coordinate lattices
with the intersections of these rational subspaces and their
integral tensor lattices. Indeed the rows of the integral filling
matrix indexed by the semistandard tableaux form an identity
submatrix, giving an integral left inverse for each embedding.
Write \(U\) for its unimodular coordinate matrix and
\(G^\vee\) for the complementary Gram matrix. The Fock Gram
identity in Section~\ref{sec:original-comparison-tests} gives
\(H(G^\vee)^{-1}\) as an integral differential pairing. Thus
\[
 G_{\rm W}=U^{\mathsf t}G^\vee U,\qquad
 HG_{\rm W}^{-1}=
 U^{-1}\bigl(H(G^\vee)^{-1}\bigr)U^{-\mathsf t}
 \ \in\operatorname{Mat}(\mathbb Z).
\] Thus \(H\zeta(Y)=(HG_{\rm W}^{-1})d(Y)\)
is integral in the original Weyl basis at every embedding
where the entries of \(Y\) are integral.
The complementary shape has \(d\) boxes. For \(d>0\) the
product of factorials of the minor sizes is at most \(d!\), so
\(\|P^\vee\|_1\le d!\) and \(0<H\le(d!)^3\le d^{3d}\).
For \(d=0\), \(H=1\).
For fixed \(j\), all choices of the short rows have the same \(d,H\).

For the norm estimate, specialize the added variables
to a point of the \(M^4\)-polydisc at the chosen conjugate embedding. The short-covector entries of \(Y\)
also have absolute values at most \(M^4\) there.
Set \(a=k+j-1\). For \(h<a\), the contributing heights
satisfy \(a-h\le j\le m\). The unit rows specified in
Section~\ref{sec:gluing-real-matched-factors} thus include all
the positive \(W_j\) columns after the first \(h\). Expand
along those unit rows. For \(h=a\), the short rows are
supported on \(W^\flat\), so the \(a\) positive \(W_j\)
columns are already variable columns with zero short-row entries.
In the resulting determinant let \(C\) be the matrix of the
specialized columns remaining beside the \(h\) columns in \(Z\).
Laplace expansion in these \(h\) columns gives coefficients
\(\pm\det C_{I^c}\), where the \(h\)-element row sets \(I\)
are restricted to the variable rows. Cauchy--Binet on all row
sets and Hadamard's inequality give
\[
 \sum_I|\det C_{I^c}|^2\le\det(C^*C)
      \le\prod_s\|C_s\|_2^2 .
\]
The number of specialized columns is at most
\(\dim W^{\mathrm{add}}+(k+j-1)-h\le12m+1\).
Each squared column length is at most
\((r+11m+1)M^8\le M^9\).
Let \(c\) be the tensor product of these exterior coefficient
vectors after specialization, and let \(F\) be the integral filling
matrix above. Straightening gives \(d(Y)=F^{\mathsf t}c\).
Since \(G_{\rm W}=F^{\mathsf t}F\), the exterior tensor
\[
 F\zeta(Y)=F(F^{\mathsf t}F)^{-1}F^{\mathsf t}c
\]
is the orthogonal projection of \(c\). Its norm is at most the
product of the exterior coefficient norms at every complex embedding.

Put \(V=H\zeta(Y)\) and apply \(\sigma\) to the coefficients
to form \(V^\sigma\). For \(\|z\|_\infty\le M^4\),
\[
 \big\|V^\sigma(z)\big\|_{\rm Gram}\le B:=H M^{9(12m+1)q_*/2}.
 \tag{I4}\label{eq:display-i4}
\]
Extend \(\sigma\) to \(\overline{\mathbf Q}\) and take a dual
row \(u\) of Gram norm one with algebraic entries. Put
\(Q=(\sigma^{-1}u)V\). With the other arguments fixed, let
\(K_{\mathrm{alg}}(P)=s_{l,n,d}K(P,\Psi;g,h)\) be the
algebraic normalized kernel of \eqref{eq:display-i3}. Then
\[
 Q^\sigma(z)=uV^\sigma(z),\qquad
 \sigma(K_{\mathrm{alg}}(Q))
       =u\bigl(\sigma(K_{\mathrm{alg}}(V_s))\bigr)_s .
\]
Apply \eqref{eq:display-i3} with \(B_0=B\) to \(Q\).
Algebraic rows are dense in the dual Gram unit sphere, so taking
their supremum gives the same norm bound for the vector of
normalized kernel values.

To expand the first factor, let
\[
 c_0=s_{m,3m,0},\qquad
 c_1=s_{j,d_+,0},\qquad
 c_2=s_{\dim W^{\mathrm{add}},d_+,d}.
\]
The scalar theta kernel \(K_0\) on \((I_m,I_{3m})\) has the unit
ideal as its finite modulus and \(1\) as both polynomial and finite
datum. Its characters are \(\chi_M^j\vartheta\) and \(1\). The second
kernel \(K_1\) is on \((I_j,V'')\) with polynomial
one, trivial first character, and the prescribed inclusion indicator
evaluated at the norm-minus-one residue unit times its argument. Its second
character is the inverse of the first-factor character
times \(\upsilon^{d_+}\). The negative-source identity
in the Weil model gives
\(\overline{K_1}=\epsilon K_{I_j^-}\), where
\(|\epsilon|=1\) and \(K_{I_j^-}\) is the negative-source kernel
in \(\Upsilon_{\vartheta,i}\).

Indeed the conjugating map of Section~\ref{sec:v2-first-exchange-phase}
takes a residue datum \(D(X)\) to \(\overline{D(\eta^{-1}X)}\).
For \(D(X)=\mathbf1_{\rm incl}(\eta X)\) this is the original
inclusion indicator. On the rectangular Hom space the determinant
of \(X\mapsto hXg^{-1}\) is \((\det h)^j(\det g)^{-d_+}\).
Thus the characters \((1,b_1)\) become
\((\upsilon^j,b_1^{-1}\upsilon^{d_+})\), with \(b_1\) the second
character of \(K_1\). The phase \(\epsilon\) comparing the two
normalized theta functionals depends only on this fixed pair.
The third kernel \(K_2\) is the linear map on polynomials
for \((W^{\mathrm{add}},V'')\) with the first-factor
character on both sides. Write its residue columns as
\((z_0,z_\flat,z_3)\) on \(I_{3m-j}\oplus W^\flat\oplus I_{3m}\),
and put \(R_0=\mathcal O_E/(\varpi^2)\). Let
\(\phi:\mathcal O_E/\mathfrak m_\vartheta\to R_0\) be reduction.
Writing \(J_{0,\flat}\) for the coordinate inclusion of the first
two blocks in \(V''\), the residue function is
\[
 \mathbf1_{(z_0,z_\flat)=J_{0,\flat}\bmod\mathfrak m_\vartheta}
 \sum_{A\in\mathrm{GL}_{3m}(R_0)}
       \psi\bigl(\operatorname{tr}(\Lambda\,\phi(z_3)A)\bigr).
\]
This is the factor on the added columns in
\eqref{eq:gluing-determinant-test-polynomial}. Its values are integral and every
conjugate absolute value is at most
\[
 \mathcal N=|\mathrm{GL}_{3m}(\mathcal O_E/(\varpi^2))|.
\]
The double class sum is over \(I_{3m}\) and \(I_m\)
at their integral levels. Both class weights are the
reciprocals of the actual stabilizer orders. Put \(n_A=C_G\):
reduction at an old odd split prime embeds each stabilizer into a
finite general linear group of rank at most \(rt\), so its order divides \(n_A\).
Expanding that sum and the coordinate minor product
gives the identity
\[
 \begin{split}
 c''&=n_A^2c_2 H c_0\overline{c_1}\epsilon\ne0,\\
 \overline{c''\Upsilon_{\vartheta,i}(h)}
   &=n_A^2\sum_{w,v}^{\mathrm{cl}}
      \overline{(c_0K_0(v,w))\,
         \overline{c_1K_1(h)}\,
         c_2K_2(H\zeta(Y);w,h)} .
 \end{split}
 \tag{I5}\label{eq:display-i5}
\]
In \eqref{eq:display-i5} the arguments of \(K_2\) on the first group
are \(\operatorname{diag}(1,h^\flat,w)\).

For \eqref{eq:display-i2}, integral representatives can be chosen
by rational translation at all places over \(\ell\).
When moving the representative of the short covectors by
\(\gamma_1\), transport the matrix \(Y\) to
\(Y\gamma_1^{\mathsf t}\).
This matrix is integral at the conjugate embedding
because
\[
       \gamma_1Y^{\mathsf t}
              =(\gamma_1 h^\flat)((h^\flat)^{-1}Y^{\mathsf t})
\]
is a product of two integral matrices there.
The last factor is integral by the original choice
of short covectors from that moved lattice.
Rational covariance of the kernel brings the transformed
value back to the lattice norm at the original \(h\).
For source and target rational elements \(\Gamma,\gamma'\)
and character exponents \(e_\tau,e'_\tau\), the exact
complex identity is
\[
 \begin{aligned}
 K_P(\Gamma g,\gamma'h)
   &=\left(\prod_\tau
       \tau(\det\gamma')^{-e_\tau}
       \tau(\det\Gamma)^{e'_\tau}\right) K_{P^{A,B_1}}(g,h),\\
 P^{A,B_1}(X)&=P\bigl((A_\tau^{\mathsf t}X_\tau
                                      B_{1,\tau}^{\mathsf t})_\tau\bigr),\\
 A_\tau&=\tau(\gamma'^{-1}),\qquad B_{1,\tau}=\tau(\Gamma).
 \end{aligned}
\]

Rational invariance of theta and the compact Fock action
at the inverse real pair give the identity before any
coefficient embedding. In the stacked case use
\(B'=\operatorname{diag}(1,B_1)\) on the full column set. Each stacked determinant changes by
\(\det(A)\det(B')\). Put \(e^{(1)}_j=j+\lfloor(t+k-1)/2\rfloor\). There are
\(m_\tau=w_j(\tau)-e^{(1)}_j\) stacked determinants at \(\tau\),
the first row length of \(\mu^+_\tau\). The first-factor character has
weight \(-m_\tau\) at \(\tau\) on both sides. Thus the determinant
multiplier \(\tau(\det\gamma')^{-m_\tau}\tau(\det\Gamma)^{m_\tau}\)
cancels the character multiplier in the displayed transformation. The remaining
Schur action on conjugate coordinates is exactly
the lattice norm transformation at \(\gamma'h\).
Thus no polynomial coefficient integrality is
asserted before moving the rational representatives.

In the conjugated product in \eqref{eq:display-i5}, use the embedding
\(\iota\circ\mathrm{conj}\) for \(K_0,K_2\) and
\(\iota\) for \(K_1\).
The scalar kernels have the indicated integral
values after moving their representatives, with
constant polynomials. Both exponents of \(K_0\)
are zero. In \(K_1\) the source point stays at one, and its
trivial first character makes the multiplier on the determinant
of the rational target change equal to one.
Thus the rational moves preserve each scalar value. The vector kernel is integral
in the original lattice norm by the covariance just proved.
The factors in the two class denominators divide
\(n_A^2\). Therefore \eqref{eq:display-i5} has lattice norm at most one for
every \(\iota\) above every \(\ell\nmid q_Q\).

For fixed \(\vartheta,i,h\), choose one such \(\ell\) where \(h\)
is integral at every place above \(\ell\). The change to the transported
lattice coordinates and its inverse are then integral, so each fixed coordinate of
\(\overline{c''\Upsilon_{\vartheta,i}(h)}\) is bounded by one for
all field isomorphisms \(\iota:\mathbb C\simeq\overline{\mathbb Q}_\ell\).
A transcendental coordinate could be sent to an element of negative
valuation by extension of transcendence bases. Thus these cleared
coordinates are algebraic. In applying \eqref{eq:display-i3} to
\eqref{eq:display-i5}, the two barred factors are read at
\(\sigma\circ\mathrm{conj}\) and the middle factor at \(\sigma\).
Both are among the automorphisms in the conjugate estimate.

Suppose the auxiliary conductor is empty, \(Q=(1)\), as for
the trivial twist. The three moduli in
\eqref{eq:display-i5} are \((1),(\pi_M^2),(\pi_M^2)\), each containing
\((M^2)=(\pi_M^2\bar\pi_M^2)\). Since \(1\le j\le m\le r\), the
ranks \(m,3m,j,d_+\) and \(\dim W^{\mathrm{add}}=12m-j+1\) are
positive and at most \(\Delta=2r^8\). The degrees are \(0,0,d\)
with \(d\le M\). Both scalar kernels have residue and polynomial
bounds one. Their conjugate norms are at most \(e^{x_1}\) by
\eqref{eq:display-i3}. For the vector kernel the residue bound is
\(\mathcal N\) and the polynomial bound of \eqref{eq:display-i4}
is \(B\), giving \(e^{x_1}\mathcal N B\).
For every choice of short rows, trace matrix and finite group
points, the triangle inequality and the bound
\(\lvert\Gamma\rvert^{-1}\le1\) on each stabilizer weight
therefore bound the squared conjugate norm of \eqref{eq:display-i5} by
\[
 \bigl(n_A^2 h_{3m}h_m e^{3x_1}\mathcal N B\bigr)^2,
\]
where \(h_a\) denotes the number of classes of \(I_a\)
at its integral level. The integral-level counting argument in the proof
of Lemma~\ref{lem:terminal-class-count}, at \(n=m,3m\), gives
\(\log h_n=O((2fn)^3\log x)\le x^{53}\) for sufficiently large \(x\). Also
\(\mathcal N\le (M^2)^{2f(3m)^2}\) and
\(\log B\le 3d\log M+
            \frac92(12m+1)q_*\log M\).
The logarithm of the displayed square is thus at most
\[
 4\log n_A+4x^{53}+6x_1+
 8f(3m)^2\log M+
            \bigl(6d+9(12m+1)q_*\bigr)\log M
        \le fmb(\log b)^2 .
\]
Indeed \(q_*=O(fb)\), \(d=O(fmb)\), and \(\log M\le22\log b\).
The definition of \(n_A\) gives \(\log n_A\le f(rt)^2\log M\).
The ranks and \(f\) grow at most polynomially in \(x\), while
\(x_1\le b^{1/10^6}\) by Section~\ref{sec:marked-level-prime-selection}.
The displayed left side is therefore \(o(fmb(\log b)^2)\)
for fixed \(q_0,d_K,\epsilon\) as \(x\) tends to infinity.
Thus the scale in \eqref{eq:display-i5} gives both first-factor bounds
for the gluing. Its nonvanishing precedes
both the three-exchange identity and the proof of injectivity of the
common period used to establish \(c'_\vartheta\ne0\).

\section{The comparison operator and the first theta lifts}
\label{sec:comparison-operator-first-lifts}

\sectionopening{In this section, we define the contraction on the lattice of
Section~\ref{sec:comparison-vanishing-lattice} and compute the Gram forms
and paired periods that enter the exchange of
Section~\ref{sec:normalized-period-exchange}.}

\subsection{The comparison operator on the original lattice}
\label{sec:original-comparison-tests}

\sectionopening{In this subsection, we contract the target tensor against the first factors of
Section~\ref{sec:marked-first-factor-bounds}. The operator takes
values in the measured direct sum over the original classes, and the
vanishing of Section~\ref{sec:marked-comparison-period-vanishing}
identifies the period terms it retains.}

Keep the actual twist \(\vartheta\), original descent
level and quotient measures fixed. The tuple space
\(\mathcal H\) is the direct sum over the gluing choices
and the \(J(\pi_M^2)\)-classes in \(V^\circ\) of the
second-factor Weyl spaces with their Gram forms.
This direct sum has its plain class metric, since the
second-factor level is free. Put \(n_A=C_G\).
For the gluing classes the bottom rows \(\pi_M I\)
of the stabilizer matrix force the reduction of every
element of \(J_\Lambda\) to be the identity at each
old place. A definite class stabilizer is finite.

At
one old odd split place it injects into a conjugate
of the torsion-free principal congruence subgroup,
so it is trivial. The same argument applies to the
second-factor \(J(\pi_M^2)\)-classes.
Consequently the reciprocal stabilizer weights in the
original gluing operator \(G_\vartheta\) are all one,
and the operator is
\[
 G_\vartheta f
   =M^{l_1}\sum_{i,c}\sum_{[h]\bmod J_\Lambda}
       \delta_{(h,c_f)L_i}
        \bigl(\phi''_{\vartheta,i}(h)\cdot f(i,c)\bigr).
\]
The dot is the integral Cartan product in the original
block order.

The notation \(\delta_z(v)\) means the
class form obtained by K-equivariant averaging of the
specified value \(v\) at \(z\); its values sum over the
class stabilizer without division. In the target pairing
against an equivariant value this multiplicity cancels
the reciprocal target stabilizer order, so
\(\langle\delta_z(v),a\rangle=\langle v,a(z)\rangle\).
The first-factor
coefficient is \(\overline{c''\Upsilon_{\vartheta,i}(h)}\)
from Section~\ref{sec:marked-first-factor-bounds}.
Its lattice norm is at most one. For the second factor use the
finite family of integral vectors constructed below. Section~\ref{sec:marked-comparison-lattice}
shows that the power \(M^{l_1}\) clears the divided
lattice in the Cartan product; the point-form evaluation
sums unitary quotient transports with unit level
characters, so preserves lattice norm. This verifies the
integral gluing for these definite coordinates.

To construct this family, work first with the fixed
column shape \(\mu^-\). Let \(b_i(X)\) denote the standard
minor products and \(F\) the integral matrix of the
embedding of their dual coordinates into the tensor product
of the exterior powers for the columns. Explicitly, the
row of \(F\) at a tableau \(T\) is the list of
coefficients of its product of minors in the \(b_i\).

The matrices here use the original space enumerations.
Writing \(G\) for the Gram matrix in the coordinate
space, Cauchy--Binet in every column gives
\[
 K(X,Y)=\sum_T b_T(X)b_T(Y)
       =\sum_{i,k}G_{ik}b_i(X)b_k(Y),\qquad G=F^{\mathsf t}F.
\]
Here \(K(X,Y)\) is the product of the prescribed
leading minors of \(Y^{\mathsf t}X\).
In the sum over \(T\), choose an increasing list of indices
independently in each exterior column, including repeated columns.
The indices \(i,k\) enumerate semistandard tableaux. Straightening
each exterior filling in this basis gives the integral rows of \(F\).
For a standard filling \(T_i\), \(F_{T_i,k}=\delta_{ik}\).
The monomial coefficients \(D_{i\alpha}=[X^\alpha]b_i(X)\) instead give
the differential pairing
\[
 H_{ik}=\left.b_i(\partial)b_k(X)\right|_{X=0}
       =\sum_\alpha D_{i\alpha}\,\alpha!\,D_{k\alpha},\qquad
 h_-=\left.b_{\rm high}(\partial)b_{\rm high}(X)\right|_{X=0}.
\]
The last integer is positive because the highest minor
product takes the value one on the standard array, and
the Fock pairing of a polynomial with real coefficients
with itself is a positive weighted sum of its squared
coefficients.

The Gram identity holds for any of these column shapes.
Write \(b\) for the column of basis polynomials and \(b(gX)=S(g)b(X)\)
for the change of basis under an invertible matrix \(g\) at each
embedding. The identity \(K(gX,Y)=K(X,g^{\mathsf t}Y)\) gives
\[
 S(g)^{\mathsf t}G=GS(g^{\mathsf t}).
\]
The differential pairing moves \(g\) from one argument to
\(g^{\mathsf t}\) on the other, so
\(S(g)H=HS(g^{\mathsf t})^{\mathsf t}\).
It follows that \(HG\) commutes with the Schur action and is scalar
by irreducibility. At the standard array \(I\), put \(v=b(I)\).
Then \(v^{\mathsf t}Gv=K(I,I)=1\) and \(b(X)^{\mathsf t}Gv=K(X,I)\)
is the highest minor polynomial. If \(HG=cI\), its differential
square is \(v^{\mathsf t}GHGv=c\). Thus \(HG=h_-I\) for the
displayed shape, with the same identity and its own highest
minor square for any complementary shape.

\subsubsection{The comparison in the Fock model}
\label{sec:original-comparison-fock}

\sectionopening{In this subsection, we evaluate the comparison operator on the monomials in the metric of Section~\ref{sec:first-factor-inverse-gram}.}

For the Hermitian Fock pairing the adjoint of substitution by a
complex matrix \(a\) is substitution by \(a^*\). For the bilinear
differential pairing used to define \(H\) it is substitution by
\(a^{\mathsf t}\). Since the basis polynomials have real integral
coefficients, the two pairings have the same Gram matrix in this basis.
The preceding computation gives \(HG=GH=h_-I\). Thus for a vector
\(e_T\) in the orthonormal basis of the full exterior tensor product
the vector \(\eta_T=HF^{\mathsf t}e_T\) has integral
coordinates and
\[
   F\eta_T=h_-\operatorname{pr}_{\operatorname{im}F}(e_T).
\]
It follows that
\[
 \sum_T\langle\eta_T,\xi\rangle_G\eta_T=h_-^2\xi,
 \qquad \|\eta_T\|_G^2\le h_-^2.
 \tag{A}\label{eq:original-comparison-test-covector}
\]
For a given prime \(\ell\nmid\operatorname{Nm}Q\),
move each original second-factor representative \(h_c\)
by a rational unitary element \(\delta_c\), so that
\(\delta_c h_c\) is integral above \(\ell\).

For an adelic point \(g\) and the integral lattice stabilizer
\(K_v\), require \(\delta_v\in K_vg_v^{-1}\) at the places \(v\)
of \(F_0\) above \(\ell\). At any other finite set of places one
can also require \(\delta_v\in g_vJ_vg_v^{-1}\) for prescribed
open level subgroups \(J_v\), which keeps \(\delta_vg_v\) in \(g_vJ_v\).
Each of these open sets contains a unitary point \(u_v\) with
\(\det(1+u_v)\ne0\): multiplying a point by a norm-one scalar
arbitrarily close to one avoids the finitely many excluded
eigenvalues. The Cayley parameters
\(A_v=(1-u_v)(1+u_v)^{-1}\) are skew-Hermitian.
Approximate them in their \(F_0\)-vector space by one rational
parameter \(A\) and set \(\delta=(1-A)(1+A)^{-1}\).
Continuity gives all the prescribed open conditions. The inverses
are in the local fields, including at the places above two.

Only finitely many places above \(\ell\)
are involved. Transport \eqref{eq:original-comparison-test-covector} at coordinate \((i,c)\) by
\[
 \eta_{c,T}=\overline{\rho(\delta_c^{-1})\overline{\eta_T}},
 \qquad
 \rho(\delta_c)\overline{\eta_{c,T}}=\eta_T .
\]
Here \(\rho\) is the Weyl action of weight
\(\mu^-\) and determinant exponent zero. The transport
is Gram-unitary, since the rational point is unitary at
every embedding and the Gram matrix is real.
For \(x\in E\), the CM involution gives
\[
 \sigma(\overline{\tau(x)})=(\sigma\tau)(\bar x)
                     =\overline{(\sigma\tau)(x)} .
\]
Thus the rational unitary equation for \(\delta_c\) holds at
\(\sigma\tau\) also when \(\sigma\) acts nontrivially on \(K\).
The rational Weyl Gram matrix and its exterior tensor realization
then make the transported action Gram-unitary at every \(\sigma\).

Therefore each tuple \(f_a\), supported at the single
coordinate \((i,c)\) with value \(\eta_{c,T}\), has squared
norm at most \(h_-^2\) at every \(\sigma\), and these tuples
have sum of rank-one operators \(h_-^2I\) in \(\mathcal H\).
Rational covariance of the lattice norm and the displayed integral
coordinate identity show that the value of \(f_a\) at
\(h_c\) has lattice norm at most one at every embedding above
\(\ell\).

In this calculation the embedding applied to the unbarred coordinate is the conjugate of the embedding applied to
\(\rho(\delta_c)\overline{\eta_{c,T}}\).
The same family therefore works at both coefficient embeddings.
Write \(v_a=G_\vartheta f_a\).
All \(v_a\) are integral at every embedding above \(\ell\)
when \(\ell\) is prime to the twist modulus.
Define the comparison operator by
\[
 U_\vartheta=n_A h_-^2G_\vartheta G_\vartheta^*
   =n_A\sum_a v_a\langle v_a,\,\cdot\,\rangle .
\]
Indeed the class coordinates of \(\mathcal H\) have the unweighted
direct-sum Gram metric. Thus for every target vector \(z\),
\[
 \sum_a v_a\langle v_a,z\rangle
 =G_\vartheta\!\left(\sum_a f_a
                   \langle f_a,G_\vartheta^*z\rangle_{\mathcal H}\right)
 =h_-^2G_\vartheta G_\vartheta^*z
\]
by \eqref{eq:original-comparison-test-covector} on each coordinate.

At the chosen embedding use the conjugate coefficient embedding in
the first Gram factor. The two transported coordinates
cancel by the transpose identity for the Weyl Gram matrix.
The rational class weights have denominators dividing
\(n_A\). For integral \(z\), the coefficients
\(n_A\langle v_a,z\rangle\) are therefore integral.
This proves preservation of the original integral
target lattice by \(U_\vartheta\).

For the fixed twist let \(y\) belong to the original projected
automorphic space on \(W_j\), and set \(Z=\operatorname{im}u_{j,\vartheta}\)
on that space.
For \(\xi\in F_w\), Lemma~\ref{lem:marked-gluing-period-vanishing}
applies at that same \(w\) and twist. In its notation the tuple
\(\eta_{i,g}\) has entries \(l_c=\overline{\kappa(g,h_c)}\) in
the \(i\)-th summand and zero entries elsewhere. With \(G_-\)
the negative tableau Gram matrix, the lemma gives
\[
  \overline{K^\circ_{(G_\vartheta^*\xi)_i}(g)}
    =\sum_c\overline{l_c}^{\mathsf t}G_-(G_\vartheta^*\xi)_{i,c}
    =\langle\eta_{i,g},G_\vartheta^*\xi\rangle=0 .
\]

For \(z=u_{j,\vartheta}y\) and \(\xi\in F_w\), insert
\(\zeta=G_\vartheta^*\xi\in\mathcal H\) in the unscaled
gluing formula~\eqref{eq:unscaled-gluing-pairing}. Each factor
\(K^\circ_{\zeta_i}(g)\) vanishes by the individual pairing above
and its complex conjugate. Hence
\(\langle G_\vartheta^*\xi,G_\vartheta^*u_{j,\vartheta}y\rangle_{\mathcal H}=0\).
Multiplying by \(n_Ah_-^2\) and taking the adjoint gives
\begin{equation}
 \langle \xi,U_\vartheta z\rangle=0
 \qquad(\xi\in F_w,\ z=u_{j,\vartheta}y,\ \text{every }w).
 \label{eq:block-period-compression-orthogonality}
\end{equation}
Thus \(U_\vartheta\) sends
integral vectors in \(Z\) to integral tests for the same
false-lift subspaces. In particular a clearing operator
\(\Lambda\) with the bound above makes
\(\Lambda\pi_Z U|_Z\) preserve the original integral lift
lattice. This proves the required compatibility of the
clearing and comparison operators on that lattice.

\indent{\it Nonvanishing of the real first-arrow constant.}
At one real place let \(P\ne0\) be the prescribed joint
harmonic polynomial. In the source action the variables
are grouped by target index. In each such group, the
source signs give the two Fock signs up to a common
reversal. The polynomial is termwise harmonic:
\(\partial_\alpha\partial_\beta P=0\) when the two
variables have opposite signs in a single such group. Use the
source discrete-series parameter \(\lambda\) so that its
lowest compact highest weight is the weight of \(P\)
for the oscillator action, including the determinant
twist. In the real factor of the doubled integral
\(\overline{e_j(g)}/e_j(1)\) is then the conjugate of the
normalized coefficient \(c_\lambda(g)\) of the unit
highest vector of this discrete series. This convention
is the contragredient on the scalar form side.

If \(z\) is a positive definite element of the source
unitary group, its action on the matrix space is
positive definite, since on each group of variables with
a common target index it is
a copy of \(z\), its transpose, or its inverse
transpose in the indicated sign conventions. Its
blocks on the positive and negative Fock variables
are positive definite. Its determinant twist is one,
since \(\det z=1\).

In the Gaussian Bargmann kernel
for this action the two pure quadratic terms involve
only opposite-sign pairs with the same target index. The
mixed bilinear term has positive definite coefficient
given, up to transposition with the variables, by the
inverse of the complex-linear part of the action.
Indeed that linear part is the direct sum of the
positive sign blocks, with complex conjugation on
the negative variables. The scalar of the kernel
is the positive vacuum coefficient. It is the
reciprocal determinant of the negative principal
block in the unitary normalization.

\subsubsection{The Gaussian pairing}
\label{sec:original-comparison-gaussian}

\sectionopening{In this subsection, we evaluate the corresponding Gaussian integral in the Fock model of Section~\ref{sec:original-comparison-fock}.}

Pair the Gaussian kernel against \(P\) in both
Bargmann variables by differentiation at zero.
Every contribution from a pure quadratic term
vanishes by termwise harmonicity. The remaining
mixed kernel pairs \(P\) with its linear substitution
by the positive definite coefficient matrix \(T\).
For \((S_aP)(X)=P(aX)\) one has \(S_aS_b=S_{ba}\) and
\(S_a^*=S_{a^*}\) in the Hermitian Fock metric. With the
positive definite square root \(D=T^{1/2}=D^*\) this gives
\[
             \langle P,S_TP\rangle_F=\|S_DP\|_F^2>0 ,
\]
since \(D\) is invertible and \(P\ne0\).
Consequently
\[
                   \langle P,\omega(z)P\rangle>0 .
\]
In particular, one can compute this coefficient with both oscillator
signs present. Write the positive element as \(z=kak^{-1}\), with
\(k\) compact and \(a\) a product of real two-dimensional boosts.
Write \(P_k=\omega(k)^{-1}P=\sum_\alpha p_\alpha Z^\alpha\).
Compact substitution preserves the mixed derivative conditions
within each group of variables with a common target index.
For a boost parameter \(t_l\), the two participating coordinates
have opposite signs. Its Bargmann kernel, with the Gaussian in the
measure, is
\[
 \frac{1}{\cosh t_l}\exp\!\left(
   \varepsilon_l\tanh t_l(Z_+Z_- -\bar Y_+\bar Y_-)
  +\frac{Z_+\bar Y_++Z_-\bar Y_-}{\cosh t_l}\right).
\]
Differentiation in \(t_l\) gives the action of
\(\varepsilon_l(Z_+Z_- -\partial_+\partial_-)\), and at zero the
kernel is the identity. Passing to the inverse transpose on a dual
factor may change \(\varepsilon_l\in\{1,-1\}\). In either case the
two pure quadratic terms vanish in the pairing with \(P_k\).
Index all boost pairs in the Hom space, with repetitions over target
coordinates. Orthogonality of the Bargmann monomials gives
\[
 \langle P,\omega(z)P\rangle
 =\sum_\alpha \alpha!\,|p_\alpha|^2
     \prod_l(\cosh t_l)^{-1-\alpha_{l,+}-\alpha_{l,-}}>0 .
\]
Coordinates fixed by \(a\) retain their usual monomial norm.
The determinant twist is one at \(a\), and the compact scalar at \(k\)
cancels that at \(k^{-1}\).

To deduce the assertion about the projector, let
\(K_\infty\) be the maximal compact of this source
group, with probability Haar, and \(d_K\) the dimension
of its lowest compact type. Insert Sun's positive-factor
formula for this coefficient, as in the real coefficient
decay calculation, and average the integrand over
\(K_\infty\times K_\infty\). At a positive factor \(z\)
the result is
\[
 \begin{split}
 &\int_{K_\infty^2}\overline{c_\lambda(k_1^{-1}zk_2)}
                   \langle P,\omega(k_1^{-1}zk_2)P\rangle
                                                  \,dk_1\,dk_2\\
 &\qquad =\frac{1}{d_K}\int_{K_\infty}
 m_{-\lambda-\rho}(u^{-1}zu)^{-1}
                   \langle P,\omega(u^{-1}zu)P\rangle\,du.
 \end{split}
\]
In fact Sun's formula is an integral of rank-one
compact projectors times the displayed weight of the
positive factor. For a fixed \(u\), substitute
\(k_1\mapsto uk_1,\ k_2\mapsto uk_2\) and apply compact
Schur orthogonality on each side of the oscillator
matrix coefficient.

The two factors \(d_K^{-1}\)
cancel one of the factors from the compact integral
formula, leaving the displayed \(d_K^{-1}\).
A preceding compact factor in the polar decomposition
is absorbed by Haar in \(k_1\).
For \(Z=u^{-1}zu\) we have \(\det Z=1\) and \(Z^{-1}=JZJ\).
Jacobi's complementary-minor identity gives \(D_I(Z)=D_{I^c}(Z)\).
Since \(\{i:-a_i\ge t\}=\{i:a_i\ge1-t\}^c\) for integral \(a\),
the definition in Section~\ref{sec:real-coefficient-decay} gives
\(m_{-a}(Z)=m_a(Z)>0\). Thus either sign of \(a=\lambda+\rho\)
gives the displayed weight, with \(\rho\) indexed in decreasing
\(\lambda\) order. The right side is strictly positive. Absolute integrability follows from
integrability of the original sufficiently separated
discrete-series coefficient and the unitary norm
bound for \(\omega\). Integrating over the real group
proves nonvanishing. The tensor product over the real
places gives the nonzero constant in the original
doubled first-arrow formula.

\subsection{Gram forms for the first lifts}
\label{sec:marked-first-arrow-rallis}

\sectionopening{In this subsection, we compute the Gram form of the initial theta lift with the finite and
archimedean data fixed. Along with the paired period of
Section~\ref{sec:marked-square-analytic-pairing}, this identity
determines the normalization for composing the lifts.}

Use the split doubled source and the diagonal target action of
\cite[\S11.3, pp.~820--821]{GanQiuTakeda}. The tensor of a Weil datum with the
conjugate datum on the opposite source gives the product of the two
theta kernels on the diagonal target. The strict-range formula
of Section~\ref{sec:marked-doubled-action} applies to each pair of
separately translated Schwartz--Bruhat tests at the ranks below.
Unfolding the graph terms gives the Gram form with the specified
vectors and quotient measures.

Put \(j=i+1\), \(\mathfrak m=\pi_M Q(\xi)\) and let \(S\)
be the union of its support and the conjugate support.
Write \(r\) for the original source rank and \(d=r+t\)
for the rank at the first arrow, so \(d>2r\).
Let \(K_0\) be the original standard compact in the finite
source group. The finite bases identify this group and
\(K_0\) for every space in the calculation. Give every local
factor of \(K_0\) mass one. Denote the source principal
level at the old places by \(J_M\) and at the whole
modulus by \(J_{\mathfrak m}\).

The first arrow \(a_{i,0}\) uses the character of index \(j=i+1\).
At target rank \(d\) its rank twist is
\[
       b_j=\mu_j^0\xi^j\upsilon_1^{\lfloor d/2\rfloor},\qquad
       \mu_j^0=\chi_0(\upsilon_1^{-\Delta}\chi_M)^j,\qquad
       \upsilon_1=\upsilon|_{U(1)}.
\]
For an allowed \(m\)-fold twist the local component of \(\xi^j\)
is nontrivial at each new place for every \(1\le j\le m\).
Thus \(Q(\xi^j)=Q(\xi)\), including the empty support when
\(\xi=1\). For the single lift on \(W_j\) at \(\xi^j\), the rank twist
of its defining character is the same expression. This
is the character in the source factor of each Weil
matrix coefficient. In column zero the doubled entries are the stem \(2\lambda_a\)
at every offset. Passing to the next row negates all signature signs,
so \(J\) changes to \(-J\) while the ordered parameter stays fixed.
The equations defining \(U(J)\) and \(U(-J)\) coincide. Compact
roots are determined by equality of two signature signs, which is
preserved by simultaneous negation. Thus the compact Cartan, positive
roots, lowest compact highest line and its normalized coefficient
agree under this identification. The group equation
\(g^*Jg=J\) is the same for the two forms over \(F_0\), so their
diagonal rational subgroups also agree. The finite change of frame is
scalar and acts trivially by conjugation on the group and principal
level. The identity on the resulting adelic quotient preserves its
normalized Haar measure and the right translations, giving the isometry
\(j_i\) of the triangle base copies. The finite target character is unramified
off \(S\) on both source spaces.

\subsubsection{The masses of the first lifts}\label{sec:first-lift-sharp-masses}

\sectionopening{In this subsection, we compute the local measures required for Proposition~\ref{prop:first-lift-gram}.}

At a chosen member of a chosen split pair use the
original Schr\"odinger matrix \(X\in M_{d,r}(E_w)\).
For \(g\) in the source integral compact, the action on
the residue test is
\[
       F(X)\longmapsto b_{j,w}(\det g)^{-1} F(Xg).
\]
Here the finite character is read on the chosen split
coordinate. The absolute-value factor on such \(g\)
is one and the unramified Weil splitting has trivial
unit factor. The opposite coordinate is already
accounted for by the dual basis in this model.
On the triangle the test at every distinguished place is
the point \(\binom{I_r}{0}\). On the single lift it is
this point at an old place, and at a new place it is
\[
       F_w\left(\binom UV\right)=\chi_w(\det U),
       \qquad U\in M_r(k_w),\quad V\in M_{t,r}(k_w).
\]
The value is zero for noninvertible \(U\).
Here \(\chi_w\) is the unit character of \(b_{j,w}\).
At the other places both data are the original vacuum.

Every matrix in either prescribed support has an integral
left inverse. If \(X\) and \(Xg\) belong to their
respective supported integral boxes, applying these left
inverses shows \(g,g^{-1}\) are integral. For the point
support the overlap then occurs exactly when \(g\equiv
I_r\). On the new block support the transformation
\(\chi_w(\det(Ug))=\chi_w(\det U)\chi_w(\det g)\)
cancels the phase. The stabilizer of the coefficient of
the block test is therefore the entire integral source
group at that place.

At an old place the stabilizer is
the principal subgroup. Normalizing the residue
measure by the original integral matrix box of mass
one, put
\[
 B_\xi=|M_{d,r}(\mathcal O_E/\mathfrak m)|,\qquad
 w_0=[K_0:J_{\mathfrak m}]
       =|\mathrm{GL}_r(\mathcal O_E/\mathfrak m)|.
\]
The matrix-coefficient integrals on the two products are exactly
\[
 c_{\triangle}=\frac{B_\xi^{-1}}{w_0},\qquad
 c_{\rm single}
    =B_\xi^{-1}
       \frac{\prod_{w\mid Q(\xi)}
                   |\mathrm{GL}_r(k_w)|q_w^{rt}}
            {\prod_{w\mid\pi_M}|\mathrm{GL}_r(k_w)|}.
\]
Indeed the coefficient is the squared norm of the
local test on its stabilizer and zero elsewhere.
The squared norm of a block test is its support
proportion, and that support has \(rt\) unconstrained
entries and an invertible \(r\)-block. There is no
extra Haar factor at a conjugate distinguished place.

For the single lift let \(J_T\subset K_T^0\) be the target
level in the standard finite compact, and let \(D_\infty\) be the
dimension of the target Weyl representation. Give \(K_T^0\) and the
compact real target group mass one. For two lifted Weyl-valued forms
\(F_1,F_2\) let \(f_1,f_2\) be their scalar coefficients at the
prescribed unit highest vector. With \(c\) ranging over the finite
target classes at \(J_T\), quotient disintegration and compact Schur
orthogonality give
\[
 \langle f_1,f_2\rangle_{\rm Pet}
 =\frac{1}{[K_T^0:J_T]D_\infty}
   \sum_c\frac{\langle F_1(c),F_2(c)\rangle_{\rm W}}
                    {|\operatorname{Aut}(c)|}.
\]
Indeed each class has measure
\(\operatorname{vol}(J_T)/|\operatorname{Aut}(c)|\), the unitary
level character cancels in the product, and the real integral
contributes \(1/D_\infty\). Hence the class pairing multiplies the
Petersson pairing by \([K_T^0:J_T]D_\infty\). Only the finite index
in this multiplier varies with the twist. Write
\[
 F_w=\prod_{h=t+1}^{t+r}(1-q_w^{-h}).
\]
At an old place the index of the target vector-basis stabilizer is the number of full-rank \(d\)-by-\(r\)
matrices, namely \(q_w^{dr}F_w\).

At a new place write the upper blocks of a residue matrix as
\((A\ B)\), in block sizes \((r,t)\). If it preserves the support
of the test, then \(A\) is invertible by testing \((U,V)=(I_r,0)\).
If \(Bb\ne0\) for some \(b\in k_w^t\), put \(a=-A^{-1}Bb\) and choose
\(V\in M_{t,r}(k_w)\) with \(Va=b\). Then \((A+BV)a=0\), contradicting
support preservation at \((I_r,V)\). Thus \(B=0\).
Conversely, an invertible lower block triangular matrix sends
\(U\) to \(AU\), so it acts on the test through the determinant
character of \(A\), with inversion for the inverse action.
The resulting parabolic preserves the span of the first \(r\) coordinate
covectors under the dual action and the kernel of the top-block
projection under the column action.
There are \(q_w^{dr}F_w\) ordered independent covector
frames and \(|\mathrm{GL}_r(k_w)|\) ordered bases of their
span. The target index is therefore
\(q_w^{dr}F_w/|\mathrm{GL}_r(k_w)|\). The target index times
\(c_{\rm single}\) consequently contributes
\(F_w/|\mathrm{GL}_r(k_w)|\) at each old place and
\(q_w^{rt}F_w\) at each new one. The first group of
factors is fixed and the second is precisely
\[
               c_Q(\xi)=\prod_{w\mid Q(\xi)}q_w^{rt}F_w .
\]
For the triangle Petersson pairing there is no such
target conversion at the first arrow.

\subsubsection{Doubling and the common series}\label{sec:first-lift-doubling-series}

\sectionopening{In this subsection, we use the masses of Section~\ref{sec:first-lift-sharp-masses} to compare the two Gram forms.}

Let \(A\) be the complement of \(S\) and put
\(G_A=U_r(\mathbb A_{F_0,A})\), with the finite
restricted product interpreted in the original finite
bases. Let \(K_A\) be its standard integral compact.
Represent \(K_A\) itself by \(1\) in \(G_A/K_A\).
Extend the chosen representatives \(\gamma\) by the identity at \(S\) and write
\[
 u_j(\gamma)
     =\bigl(b_j(\det\gamma)\Delta(\gamma)^d\bigr)^{-1},
 \qquad
 R_j=\sum_{\gamma\in G_A/K_A}u_j(\gamma)R(\gamma).
\]
The series acts on \(K_A\)-fixed vectors. Its coefficient is constant
on every double coset. Indeed, for the standard lattice \(\Lambda\)
and \(k,k'\in K_A\),
\[
 \Lambda\cap(k\gamma k')^{-1}\Lambda
          =(k')^{-1}(\Lambda\cap\gamma^{-1}\Lambda).
\]
The determinant ideals of the two intersection lattices coincide.
Thus their multipliers \(\Delta\) coincide by the unramified formula
of Section~\ref{sec:reproduction-vacuum-coefficients}.
Also \(b_j\) is trivial on determinants of \(K_A\).
This applies at every conjugation orbit in \(A\). Both source levels
contain \(K_A\), so \(R(\gamma k)y=R(\gamma)y\) for \(k\in K_A\).
On each of the \(d\) target covectors the source acts by
\(\gamma^{-\mathsf t}\). In the finite unitary bases this is
\(\bar\gamma\). Conjugation preserves the absolute value and
conjugates the unitary character in the intersection multiplier, so
\(\Delta(\gamma^{-\mathsf t})=\overline{\Delta(\gamma)}\) on each
conjugation orbit. The vacuum coefficient
\(1/\overline{\Delta(\gamma^{-\mathsf t})}\) thus contributes
\(\Delta(\gamma)^{-1}\) per covector.
In particular its absolute value is
\[
 |u_j(\gamma)|=\left(\prod_{w\in A}
     [\gamma\Lambda_w:\gamma\Lambda_w\cap\Lambda_w]\right)^{-d/2},
\]
where the product includes both members of split pairs.
This follows by multiplying the integral lattice
overlap by its unitary square-root factor. The
complex multiplier also agrees under the finite lattice isometry
\(\jmath_f:\Lambda\longrightarrow\Lambda'\). Indeed, with
\(\gamma'=\jmath_f\gamma\jmath_f^{-1}\) and
\(B\Lambda=\Lambda\cap\gamma^{-1}\Lambda\), one has
\[
  (\jmath_f B\jmath_f^{-1})\Lambda'=\Lambda'\cap(\gamma')^{-1}\Lambda',
  \qquad \det(\jmath_f B\jmath_f^{-1})=\det B,\qquad
  \det\gamma'=\det\gamma .
\]
Thus both the splitting character in \(\Delta\) and the factor
\(b_j(\det\gamma)\) have the same values on the two groups.

Absolute convergence and injectivity can both be
seen by a lattice count in the present ranks
\(t=5r\), \(d=6r\). At a split pair a lattice at the
chosen place determines its dual at the other.
If it is not contained in \(\Lambda_w\), choose
\(x\) in it with \(\varpi_wx\) primitive in \(\Lambda_w\).
Enlarge \(\Lambda_w\) by \(x\) and take its dual
at the opposite place. The new intersection with the
target lattice maps onto the quotient of the enlarged lattice by
\(\Lambda_w\), generated by the class of \(x\) and of order \(q_w\).
Its kernel is the old intersection at the chosen place.
The dual lattice shrinks at the opposite place, but contains
the old intersection there by integral pairing with \(x\).
Since the target lattices are fixed, the index at the chosen
place is divided by \(q_w\) and the index at the opposite place
is unchanged. Thus the product of the two indices is divided by
\(q_w\). If containment is strict in the other
direction, start on the opposite member. Each enlargement is
specified by the line of \(\varpi_w x\) in the residue
space of dimension \(r\). Thus there are at most
\(2(1+q_w+\cdots+q_w^{r-1})\) moves at one split pair.

At a nonsplit unramified place the two lattices are
self-dual. Choose \(\varpi_w=\varpi_v\) in the fixed field. If the
lattices differ, take \(x\) as above in
\(\Lambda'_w\) and put \(y=\varpi_wx\). The lattice
\(\Lambda_0=\{z\in\Lambda_w:h(y,z)\in\varpi_w\mathcal O_w\}\)
has index \(q_w\). The lattice
\(\Lambda_0+\mathcal O_w x\) is integral because
\(h(x,\Lambda_0)\subseteq\mathcal O_w\) by definition and
\(h(x,x)\in\mathcal O_w\) since \(x\in\Lambda'_w\).
As \(\varpi_w x=y\in\Lambda_0\) and \(x\notin\Lambda_w\),
it contains \(\Lambda_0\) with index \(q_w\). Its covolume
equals that of \(\Lambda_w\), so integrality and self-duality of
\(\Lambda_w\) imply that it too is self-dual.
Put \(P=\Lambda_w\cap\Lambda'_w\). Then \(P\subseteq\Lambda_0\)
since \(h(x,P)\subseteq\mathcal O_w\), and \(y\in P\).
The new intersection maps onto
\((\Lambda_0+\mathcal O_w x)/\Lambda_0\), since it contains \(x\),
and the kernel is \(P\). Hence it equals \(P+\mathcal O_w x\)
and the intersection index is divided by \(q_w\).
For a fixed residue line
of \(y\), two choices of \(x\) differ by \(z\in\Lambda_w\)
after fixing that line's representative. Modulo \(\varpi_v\) one has
\[
 \varpi_v\bigl(h(x+z,x+z)-h(x,x)\bigr)
    \equiv\operatorname{Tr}_{k_w/k_v}h(\bar y,\bar z).
\]
The residue pairing with \(\bar y\ne0\) identifies
\(\Lambda_w/\Lambda_0\) additively with \(k_w\). Since
\(k_w/k_v\) is finite separable of degree two, its trace is
surjective, also in characteristic two. Thus at most
\(|k_v|=\sqrt{q_w}\) classes give an integral square for this line.
Each new self-dual lattice has an integral orthonormal basis.
Its reduction is a nondegenerate Hermitian space over \(k_w/k_v\).
If all residue norms were zero, expanding the norm of \(x+ay\)
would give \(\operatorname{Tr}_{k_w/k_v}(a h(x,y))=0\)
for all \(a,x,y\), contrary to nondegeneracy of the residue
trace pairing. A lift of a vector with nonzero norm has
unit norm and can be scaled to norm one by surjectivity
of the norm on units
(Section~\ref{sec:splitting-level-characters}).
Its orthogonal complement in the lattice is self-dual.
Induction gives the basis, including in residue characteristic two.

For each such fixed move choose an isometry taking
the new standard lattice back to the original one.
Repeat on the moved \(\Lambda'\). The product of
indices in \(|u_j|\) decreases by a factor at least
\(q_w\) at each step. Reversing the resulting word
recovers the original lattice coset. Its coefficient
is bounded by the product of the step weights
\(q_w^{-d/2}\).
The sum of those step weights over all choices of a
single step is
\[
 \eta \ll [E:\mathbb Q]\sum_{n\ge2} n^{-2r}
             \ll [E:\mathbb Q]\,2^{-2r}.
\]
Indeed either count above is \(O(q_w^r)\), and there
are at most \([E:\mathbb Q]\) places of a given prime
power norm.

Here \([E:\mathbb Q]\ll x^{12}\) and
\(r\ge x^5\), so \(\eta<1/3\) for all sufficiently
large \(x\). Choose one word for each nonidentity coset. The geometric series gives
absolute convergence in operator norm and
\[
 \|R_j-I\|\le\sum_{\gamma\ne1}|u_j(\gamma)|
                \le\frac{\eta}{1-\eta}<\frac12 .
\]
The estimate holds on the whole \(K_A\)-fixed space, including the
full triangle base and the two projected spaces. On each of these
spaces the identity coset acts by the identity, so the corresponding
operator is injective. More precisely, let \(\eta_0\) be the same
sum of one-step weights allowing every finite place. It has the
same bound \(O([E:\mathbb Q]2^{-2r})\). Each \(b_j\) is unitary,
so its value on a coset contributes a factor of absolute value
one. Omitting the old and new places only decreases the sum.
Thus \(\|R_j-I\|,\|R_j^W-I\|\le\eta_0/(1-\eta_0)\) with
one bound for all the allowed supports and all \(1\le j\le m\).
Let \(\ell(\gamma)\) be the least number of the preceding neighboring
lattice changes joining \(\Lambda\) to \(\gamma\Lambda\). The
standard compact preserves adjacency and fixes \(\Lambda\), so
\(\ell\) is constant on double cosets. Write \(R_{j,<L}\) for the
absolutely convergent subsum with \(\ell(\gamma)<L\).
Every omitted coset has chosen reducing word of length at least \(L\);
hence on either realization
\[
 \|R_j-R_{j,<L}\|\le\sum_{h\ge L}\eta_0^h
      =\frac{\eta_0^L}{1-\eta_0}\qquad(L\ge1).
\]

Equivalently, index the finite conjugation orbits in \(A\) by places
\(v\) of \(F_0\). Let \(b_v\) be the sum of the absolute values of
the nonidentity local coset coefficients, and let \(\eta_v\) be the
sum of the one-step weights there. The same word count gives
\(b_v\le\eta_v/(1-\eta_v)\) and \(\sum_v\eta_v\le\eta_0<1/3\).
Multiplicativity of the coefficient and monotone convergence give
\[
 \sum_{\gamma K_A}|u_j(\gamma)|
     =\prod_{v\in A}(1+b_v)<\infty .
\]
At a split place the two dual lattices enter one local factor.
In particular the places over \(N\) occur in this product:
\((q_\xi,MN)=1\), the glued rank-\(r\) lattice is unimodular
there, and the last-covector condition is on the rank-\(k\) source
of Section~\ref{sec:ramified-source}.

The doubled action is identified in Section~\ref{sec:marked-doubled-action}.
The doubled source has dimension \(2r\) and a Lagrangian of dimension
\(r\). Thus the Siegel index in \cite[Theorem~2]{YamanaSW} is \(n=r\).
For the target \(T\) of dimension \(d=6r\) and Witt index \(\iota(T)\),
\[
 d-\iota(T)\ge d/2=3r>r,\qquad s_0=(d-r)/2 .
\]
Weil's convergent identity applies to the full Schwartz--Bruhat
space with coefficient one and probability target measure
\cite[\S52, Theorem~5, pp.~76--77]{WeilSiegel}, by the verification in
Section~\ref{sec:marked-doubled-action}. The doubled datum is the original
tensor times the conjugate tensor on the opposite
source. In a rational Schr\"odinger polarization the
partial Fourier transform in the difference variable
gives the product of the theta kernels on the diagonal
group, and its graph zero term gives the original
unitary matrix coefficient. The finite opposite-basis
scale cancels against its inverse in the conjugate tensor.
The graph change of variables preserves the probability
measure on the additive quotient, by the calculation
in the \hyperref[sec:marked-graph-doubling]{principal-level doubling calculation}.

With source character \(b_j\), the opposite source carries
\(b_j^-=b_j^{-1}\upsilon_1^d\). Indeed its conjugate character
is \(b_j^{-1}\), and the intrinsic multiplier on the opposite
Hom space supplies \(\upsilon_1^d\), as in
Section~\ref{sec:marked-graph-doubling}. Write
\(f_g=\omega_{b_j}(g,1)f\), and take the antiunitary \(C\) of that
calculation, normalized by \(\theta_-(Cw)=\overline{\theta_+(w)}\).
The doubled test is \(f_g\otimes C f_{g'}\). The source multipliers
are \(b_j(\det g)^{-1}\) and \(b_j^-(\det g')^{-1}\).
The graph parameter \(Z\mapsto X_T^{-1}J_T^{-1}ZJ_S X_S\)
preserves probability Haar between the rational additive quotients.
The graph indexed by \(\delta\in U_r(F_0)\) therefore gives
\(\langle f,\omega_{b_j}(g'^{-1}\delta g,1)f\rangle\).

For clarity, the other rational Lagrangians give no
term on either of the present tests. If a Lagrangian of
the doubled source is not a graph, its intersection
with the first source has a nonzero vector. Indeed
if its second projection is injective it is an
isomorphism, and the isotropy equation then makes
the first projection an isometry. A vector in that
intersection yields a Heisenberg translation in the
first tensor factor alone under which the doubled
zero-term functional is invariant.

In the split
Schr\"odinger coordinates at an old place this
translation is multiplication by
\(\psi_w(\operatorname{tr}(Z Xc))\), for arbitrary
\(Z\), where \(c\ne0\) is the transported first-source
vector at \(w\). Invariance forces the functional to
be supported on \(Xc=0\). To see this directly,
partition a compact open support in \(X\) into finitely
many compact open pieces on each of which some one
multiplier is a constant different from one. Test
invariance on each piece. The support has injective
\(X\), so all pieces are of this kind. This is the
usual open-graph unfolding with a full-rank
Schr\"odinger test.

In the convergent range the absolute sum of the
Siegel--Weil sections is bounded by a Schwartz
seminorm of the transformed datum times the
convergent spherical Eisenstein series at the
identity. This follows by writing each section in
Iwasawa coordinates and bounding its restriction to
the maximal compact by that seminorm. At reduced
representatives of the two source quotients the
seminorm of the tensor of the Weil translates grows
polynomially in the reduction height. At fixed source levels the
finite class representatives range over finite sets. Their remaining
compact factors have finite orbits on the locally constant finite
data. Thus the seminorm estimates apply to both original residue
tests for all such representatives.

The source
forms in the prescribed real coefficient types are
rapidly decreasing by the coefficient-projector
argument in Section~\ref{sec:adjacent-unfolding}.
That argument applies to every bounded finite-level
form of these types. For a fixed-level \(L^2\) vector reproduced by
the projector, truncate its values and project again. The projected
truncations converge in \(L^2\) and lie in the same finite-dimensional
cusp space of the given compact type and infinitesimal character,
as in Section~\ref{sec:marked-deep-separation}. This gives a rapid
representative also for the full projected space.
For a finite right translate \(t\),
choose an integer \(c\ge1\) such that the fixed finite lattice
satisfies \(c\widehat\Lambda\subset t\widehat\Lambda
\subset c^{-1}\widehat\Lambda\). Multiplication by \(g_f\)
and intersection with the rational vectors give the corresponding
inclusions for the moved lattices. Hence
\(c^{-1}H(g)\le H(gt)\le cH(g)\). Such a translate preserves rapid
decay, commutes with the real projector, and changes its finite
level \(J\) to \(tJt^{-1}\). In particular both source forms may
be used even when only the second is fixed by \(J_{\mathfrak m}\).

Before unfolding, the absolute product of the original theta
kernels is integrable in all three quotient variables. Indeed,
with each source test paired with its own conjugate, the
pointwise convergent Siegel--Weil formula above gives
\[
 \int_{[U_d]}|\theta_f(g,h)|^2\,dh\ll_f H(g)^{2B}.
\]
Here the diagonal doubled section has the preceding polynomial
bound and the target quotient volume is fixed. Cauchy--Schwarz
in \(h\) bounds the integral of the absolute mixed product at
\(g,g'\) by \(C H(g)^B H(g')^B\), increasing \(B\) for the two
fixed data. Integration against the rapid source forms is finite.
The absolute graph sum satisfies the same product majorant by the
Iwasawa bound. Thus Fubini applies first to the theta product and
then to the graph sum, as in Section~\ref{sec:marked-graph-doubling}.
Unfolding the rational graph sum by the quotient
Haar formula leaves the integral over the original
adelic source of the Weil coefficient times
\(\langle x,R(g)y\rangle\), multiplied by the original
target quotient volume.

At infinity a rank-one discrete-series coefficient
projector gives
\[
 \langle x,R(g_\infty)z\rangle
     =\frac{\overline{e_j(g_\infty)}}{e_j(1)}
                                        \langle x,z\rangle .
\]
This is the compression of the real right regular
action to the unit highest coefficient line, by
Schur orthogonality. The source form in the second variable is fixed by \(J_{\mathfrak m}\) for the triangle
and by \(J_M\) for the single lift. Fubini therefore
gives the local masses computed above and the same series over the remaining
places. The real factor is
\[
 \operatorname{vol}([U_d])
  \int_{U_r(F_0\otimes\mathbb R)}
   \langle v_j,\omega(g_\infty,1)v_j\rangle
      \frac{\overline{e_j(g_\infty)}}{e_j(1)}\,dg_\infty .
\]

For the two pairs use their respective joint harmonic
and quotient volume. These constants are independent
of \(\xi\), since the real parameters and vectors do
not include the tame twist. At a real place let \(d_\lambda\)
be the formal degree for the given Haar measure. With
\(c_\lambda\) in the oscillator convention above, put
\(p_\lambda=\omega(d_\lambda\overline{c_\lambda})\).
Schur orthogonality makes this an orthogonal projection, and
the compact average in Section~\ref{sec:original-comparison-gaussian}
gives strict positivity for the indicated harmonic:
\[
 \int\langle P,\omega(h)P\rangle\overline{c_\lambda(h)}\,dh
      =d_\lambda^{-1}\|p_\lambda P\|^2>0.
\]
The same calculation applies to the first true
return: the rank \(r+1\) source and rank \(r+t\)
target have the common entries and inserted block
computed in the signed parameter comparison.
Their joint harmonic has only one Fock sign in
each group of variables with the same target index,
so both pure quadratic terms vanish in the
Bargmann pairing as above.

For \(x,y\) in the original full triangle base, \(y\) is fixed
by \(J_{\mathfrak m}\), as required in the unfolding. This gives the
following two identities in the original Petersson metrics.

\begin{proposition}[Gram forms of the first lifts]
On\label{prop:first-lift-gram} the full triangle base, with \(\mathcal V_i\) the convergent
spherical series,
\[
  w_0\langle a_{i,0}j_i x,a_{i,0}j_i y\rangle
           =C_i B_\xi^{-1}\langle x,\mathcal V_i y\rangle .
\]
For \(j=i+1\), the restriction of \(\mathcal V_i\) to the projected source is \(R_j\). On the projected source of the single lift,
\[
    \langle u_{j,\xi^j}x,u_{j,\xi^j}y\rangle_{\rm class}
          =A_j c_Q(\xi)\langle x,R_j^W y\rangle ,
\]
where \(R_j^W\) uses the same coefficients and cosets
as \(R_j\), acting on \(W_j\). On each right-hand side use the
Petersson pairing on the corresponding source. The two restrictions
satisfy \(\det R_j=\det R_j^W\).
\end{proposition}

For either lift \(U\), write \(K\) for its spherical operator and
\(c\) for the scalar multiplying its quadratic form. The Haar, Gaussian
and finite factors give \(c>0\). The preceding series estimate gives
\[
 \|U y\|^2=c\operatorname{Re}\langle y,Ky\rangle
    \ge c(1-\|K-I\|)\|y\|^2\ge\tfrac{c}{2}\|y\|^2
\]
in the respective source and target metrics.

\subsubsection{The two spherical actions}\label{sec:first-lift-spherical-actions}

\sectionopening{In this subsection, we compare the actions of the coefficient series of Proposition~\ref{prop:first-lift-gram} on the projected spaces.}

For a double coset \(K_A\gamma K_A\), its right cosets
are represented by \(k\gamma\) with \(k\) ranging over
\(K_A/(K_A\cap\gamma K_A\gamma^{-1})\). Because
\(u_j\) is bi-invariant, the portion of the series
on this double coset is exactly \(u_j(\gamma)\)
times its unnormalized Hecke operator. These
operators commute with the finite residue action on
the projected subspace.

Take finite subsums \(h_n\) over complete double cosets exhausting
the series. They belong to the common finite spherical Hecke algebra,
including its factors at nonsplit places. Their actions converge in
operator norm on the two projected spaces by absolute convergence
and unitarity. For each positive integer \(k\),
\eqref{eq:marked-full-spherical-trace} applied to \(h_n^k\) gives,
on passing to the limit,
\[
       \operatorname{tr}(R_j^k)=\operatorname{tr}((R_j^W)^k).
\]
The same trace identity at the unit gives equality of dimensions.
Newton's identities now give \(\det R_j=\det R_j^W\) with the
original Petersson metrics on both spaces.

For the first arrow \(a=a_{i,0}\), take \(u\in K_0\) and
write \(u_{\mathfrak m}\) for its reduction at the chosen split
members of \(\mathfrak m=\pi_MQ(\xi)\). Lift
\(\operatorname{diag}(u_{\mathfrak m},I_t)\) to an integral target
element \(h\), using inverse transpose at the conjugate members
and the identity away from these pairs. Writing \(h_{\mathfrak m}\)
for its chosen reduction and \(X_0=\binom{I_r}{0}\), one has
\[
                 h_{\mathfrak m}^{-1}X_0u_{\mathfrak m}=X_0 .
\]
Thus the simultaneous action preserves the point indicator.
The tame determinant factors cancel on the matched reductions,
and the unramified factors fix the vacuum on integral elements off
the modulus. More generally, writing \(\lambda\) for the unit scalar
of the simultaneous action, covariance and the change of source
variable \(g=g'u\) give \(R(h)a=\lambda aR(u)\). The integral
elements normalize the respective principal levels. With their
original Petersson pairings,
\[
 \langle aR(u)x,aR(u)y\rangle
   =\langle R(h)ax,R(h)ay\rangle=\langle ax,ay\rangle .
\]
Pairing with every vector of the full triangle base gives
\(R(u)^{-1}a^*aR(u)=a^*a\). In particular the return preserves
the \(J_M\)-fixed subspace of the \(J_{\mathfrak m}\)-fixed base.
Identify the base copy by the isometry \(j_i\). For every \(y\) in
this entire original base, Proposition~\ref{prop:first-lift-gram}
and the uniform estimate \(\|\mathcal V_i-I\|<1/2\) give
\[
 \langle y,a^*ay\rangle=\|ay\|^2
 =c\,\operatorname{Re}\langle y,\mathcal V_i y\rangle
 \ge \tfrac{c}{2}\|y\|^2,\qquad
 c=\frac{C_i}{w_0 B_\xi}>0 .
\]
Thus the return is injective on the original base level.
The first adjoint descent gives this same return
on the next isometric base copy.

\subsection{The period pairing for the four unitary groups}
\label{sec:marked-square-analytic-pairing}

\sectionopening{In this subsection, we write the composition of the two theta lifts as a paired period on
the direct sum of the source and target spaces. This places both
choices of intermediate space in the normalization fixed in
Section~\ref{sec:marked-first-arrow-rallis}.}

Let \(U_1\) and \(U_2\) be the source and target Hermitian
spaces and \(A\) and \(A'\) the two middle spaces. All quotient
integrals below use the measures of the theta operators.
Write \([B]=U(B)(E^+)\backslash U(B)(\mathbb A_{E^+})\).
Set \(W=U_1\perp U_2\), with the source block first.
Write \(X=(X_1,X_2)\) and \(Y=(Y_1^{\rm t},Y_2^{\rm t})^{\rm t}\)
for the two specified point arrays. Write
\[
 T_X=D_{2,X_2}L_{1,X_1},\qquad
 T'_Y=L'_{2,Y_2}D'_{1,Y_1}.
\]
Each arrow uses the indicator of the indicated block at the
modulus and the original vacuum away from it. The two lower
real polynomials are \(P'_1,P'_2\), and the two upper real
polynomials are \(P_1,P_2\).

Write \((a_i,b_i)\) for the upper splitting pairs and
\((b'_i,a'_i)\) for the lower splitting pairs. Put
\(N=\dim W=\dim A+\dim A'\) and
\(\gamma=\nu^{-1}\upsilon^{N/2}\), where \(\upsilon\)
is restricted to the norm-one ideles. Each pair of opposite
rank corrections has exponents \(\lfloor p/2\rfloor\) and
\(\lceil q/2\rceil\), with \(p+q=N\) even, in one order or
the other. Their sum is \(N/2\), and the product of the
two factors before rank correction is \(\nu^{-1}\).
Applying this calculation at the two kept ranks and at the
two middle ranks gives
\[
 a_1a_2=a'_1a'_2=b_1b'_1=b_2b'_2=\gamma.
\]
In particular \((\alpha,\alpha')=(b_1,b'_1)\) satisfies the
character hypothesis of Proposition~\ref{prop:normalized-period-exchange}.
Apply that proposition on \(W\) with \(V=A\) and \(U=A'\). Its
rank, signature and weight hypotheses for the chosen squares are
checked in Section~\ref{sec:four-unitary-marked-application}. For the fixed lower real tensor
\(v_-=P'_1P'_2\) (with its standard Fock Gaussian), the
comparison supplies a real Fock vector \(v_+\) of finite
moments, before choosing the admissible modulus or character.
For matched finite arrays \((\Phi,\Psi)\) it gives
\begin{equation}
 \int_{[A]}K_{W,A}^{\gamma,b_1}(\Phi,v_+;g,h)\,dh
 =
 c_\xi\int_{[A']}K_{A',W}^{b'_1,\gamma}(\Psi,v_-;h',g)\,dh'
 \qquad(g\in[W]).
 \label{eq:marked-square-paired-period}
\end{equation}
Here \(c_\xi\ne0\). If \(\dim A-\dim A'=2t\), the comparison
asserts
\[
 |c_\xi|^2 =
 C\prod_{w\mid\mathfrak m}q_w^{\,t^2}
       \prod_{k=1}^t\big|L^{S_{\mathfrak m}}(k,\nu)\big|^2,
 \qquad C>0.
\]
The constant \(C\) and \(v_+\) may depend on the ordered
real spaces, the fixed lower datum and archimedean weight.
The finite characters and arrays are varied after this choice.
Equation~\eqref{eq:marked-square-paired-period} is a period
comparison for two kernels, with the individual original
quotient measures of \(A\) and \(A'\).

Set \(\nu_0=b_1/b_2=b'_2/b'_1\).
Our convention is
\(K^{a,b}(g,h)=a(\det h)b(\det g)^{-1}K^{1,1}(g,h)\).
Multiplying the upper kernel with a pure real tensor by
\(\nu_0(\det g_2)\) changes the factor
\(b_1(\det g_2)^{-1}\) to \(b_2(\det g_2)^{-1}\).
Multiplying the lower kernel by the same factor changes
\(b'_1(\det g_2)\) to \(b'_2(\det g_2)\).
Use the product normalization of the rational theta functional on
the two Hom blocks, as in Section~\ref{sec:gluing-real-matched-factors}.
In a product rational polarization the rational points form the product
of the two rational lattices. The two Heisenberg actions commute and
their tensor product is the block action. Thus the normalized theta
functional on a pure tensor is the product of its evaluations on the
two factors. The Gaussian,
the real pure tensor and the finite point indicator are products
in those blocks. Thus the two corrections have the same sign.

Take the source vector \(f\) and the target vector \(z\) in their original real
types, and set \(f_1=f\), \(f_2=\overline z\). Pairings of
automorphic forms are linear in the second argument.
For a real Fock vector \(v\) set
\[
 H_{\xi,X,f,z}(v)=
 \int_{[U_1]\times[U_2]\times[A]}
  f_1(g_1)f_2(g_2)\nu_0(\det g_2)
  K_{W,A}^{\gamma,b_1}
       (\mathbf{1}_X,v;g_1\oplus g_2,h)\,dg_1\,dg_2\,dh .
\]
Only primitive \(X\) are needed. Indeed, the support of a matched
\(\Phi\) is formed by the last rows of \(k^{-1}\), and the
support of its partner \(\Psi\) by the first columns of \(k\),
for invertible \(k\). Thus \(X\) has a right inverse and \(Y\)
has a left inverse at every factor of the nonempty modulus.

We prove a bound uniform in the integrated middle
variable. Write \(\rho=\operatorname{ht}_W(g)\) and
\(H=\operatorname{ht}_A(h)\) for the reciprocal minima
as in Section~\ref{sec:adjacent-unfolding}. In the real
Schr\"odinger model write \(u\) for the image of \(v\).
With the finite test fixed, reduction in both variables
and Weil covariance of position and differentiation
give the bound
\[
             |L_{g,h}(u)|\le C(\rho H)^A s_j(u)
\]
for the evaluation \(L_{g,h}\) of the kernel and a fixed
Schwartz seminorm \(s_j\). Indeed the reduced finite
translates range over a finite set and the real changes
of variable have matrix and inverse norms bounded by
powers of \(\rho H\). Take \(A\) to be a positive integer.

\subsubsection{Uniform bounds for the two kernels}
\label{sec:four-unitary-constant-terms}

\sectionopening{In this subsection, we bound the kernels in \eqref{eq:marked-square-paired-period} by averaging over a rational isotropic subspace.}

When \(H>2\), take a shortest vector \(d\) of the moved
middle lattice. Its integral Hermitian norm has absolute
value below one at every infinite place, so it is isotropic.
In the kept lattice take an \(E\)-basis of vectors \(c_j\)
with lengths \(O(\rho^a)\), for a fixed \(a\ge1\).
Indeed the real lattice has fixed covolume and first
minimum at least \(\rho^{-1}\), so Minkowski's second
theorem bounds its last minimum by \(O(\rho^{D-1})\)
in real dimension \(D\). Among \(D\) independent
\(\mathbb Z\)-lattice vectors one can take an \(E\)-basis.

For the upper kernel average the theta functional over
the adelic quotient of the rational isotropic space
spanned by \(dc_j^*J_W\), with probability Haar.
In the moved finite coordinates these displacements
are integral. As in Section~\ref{sec:marked-mixed-absolute-summability},
the modulations at a split place read
\(d_h^*\widetilde X C_g\), where \(C_g\) is the moved
basis matrix and \(\widetilde X\) is a local matrix in
the support with residue \(X\). This covector is nonzero
over the local field: right invertibility of the residue
implies right invertibility of \(\widetilde X\).

For the lower kernel take the short isotropic vector \(c\)
in the middle source and a basis \(d_j\) in the kept target.
The modulations read \(D_g^*\widetilde Y c_{h'}\),
nonzero by left invertibility of \(\widetilde Y\)
over the local field. In either case, partition
the compact finite support into pieces on each of which
one modulation is a constant different from one.
Invariance of the averaged functional under this
modulation makes it vanish on the piece for every real
vector \(u\).

Rational Heisenberg invariance and additive strong
approximation allow computing the average with bounded
real coefficients of the rational generators and with
finite coefficients in a fixed ideal of
\(\widehat{\mathcal O}_E\). The ideal can be chosen so
that every such finite displacement fixes the finite
test, by integrality of the generators in the moved
coordinates. This yields an average \(A_h\) of real
translations and modulations of size \(O(\eta)\),
where \(\eta=\rho^a/H\). The vanishing gives
\(L_{g,h}(A_hu)=0\) for every \(u\). For \(\eta\le1\),
differentiating the real Heisenberg law along the segment
from zero to each displacement gives
\(s_j((1-A_h)u)\le C_j\eta s_{j'}(u)\) with a higher
fixed seminorm \(s_{j'}\). Now
\[
 |L_{g,h}(u)|=|L_{g,h}((1-A_h)^A u)|
    \le C(\rho H)^A\eta^A s_{j''}(u)
    =C\rho^{A(a+1)}s_{j''}(u).
\]
For \(\eta>1\) use \(H<\rho^a\) in the initial bound.
For \(H\le2\) the initial bound suffices.

In the Fock coordinates the coherent-state
formula is
\[
 v(w)=\pi^{-d}\int_{\mathbb C^d}v(y)
 e^{-(|w|^2+|y|^2)/2+\bar y\cdot w}\,dy.
\]
The inverse Bargmann transform of the exponential
is a translated and modulated normalized Gaussian,
up to a unit scalar, by the Gaussian integral in the
Bargmann transform. By Leibniz's rule and translation
of the polynomial weights, its \(s_{j''}\)-seminorm
is \(O((1+|y|)^N)\) for a fixed \(N\).

Integration in that seminorm gives
\(s_{j''}(u)\le C p_N(v)\). Thus
\[
 \left|K_{W,A}^{\gamma,b_1}(\mathbf{1}_X,v;g,h)\right|
 \le C_X p_N(v)\operatorname{ht}_{W}(g)^B
 \qquad\text{uniformly in }h.
\]
More generally, for every integer \(D\ge0\), the target height can
be given the factor \(H^{-D}\). For \(H>2\) and
\(\eta=\rho^a/H\le1\), use \(A+D\) differences. The preceding
Schwartz estimate becomes
\[
  C_D(\rho H)^A\eta^{A+D}s_{j_D}(u)
       =C_D\rho^{A+a(A+D)}H^{-D}s_{j_D}(u).
\]
If \(H\le2\) or \(\eta>1\), then \(H\le2\rho^a\), and the initial
bound gives the same inequality after increasing \(C_D\).
The lower kernel uses the integral left inverse of \(Y\) in the
transposed model, in place of the integral right inverse of \(X\).
The proof applies termwise to each finite sum of such residue tests.
In particular the estimate is uniform in \(h'\) for
\(K_{A',W}^{b'_1,\gamma}(\mathbf{1}_Y,v_-;h',g)\).
Here \(p_N\) is the weighted absolute Fock \(L^1\) moment in
the prescribed normalization. In the original block coordinates
the moved integral lattice at \(g_1\oplus g_2\) is the direct sum of
the two moved integral lattices. The length of a vector in the
orthogonal sum dominates the length of each component, and each
summand embeds isometrically. Thus its first minimum is the lesser
of the two first minima and
\[
 \operatorname{ht}_W(g_1\oplus g_2)
    =\max(\operatorname{ht}_{U_1}(g_1),\operatorname{ht}_{U_2}(g_2)).
\]
The source tests are rapidly decreasing, and
the quotient measures are finite. These inequalities give
absolute integrability and, for each fixed
\(j=(\xi,X,f,z)\), the continuity bound
\[
            |H_j(v)-H_j(v')|\le C_j p_{N_j}(v-v').
\]
The constants \(C_j,N_j\) can vary with \(j\). In particular these
bounds do not change the normalization of the quotient measures.

With lower real datum \(v_-\), the block product formula and
Fubini give the integral of the two transposed lower kernels
against \(f_1,f_2\), paired bilinearly on \([A']\).
By transposition of the second lower arrow this is exactly
\(\langle z,T'_Y f\rangle_{[U_2]}\).
The identical upper calculation with datum \(P_1P_2\) gives
\(\langle z,T_X f\rangle_{[U_2]}\). These transpositions use
the product quotient measure and introduce no indices.

For an arbitrary upper pure polynomial tensor \(Q_1Q_2\),
the real edge convolution identity proved below expresses the two lifts on
\([A]\) as \(C_{\alpha_1}F_1\), \(C_{\alpha_2}F_2\), where
\(F_j\) are the lifts with \(P_j\), and
\[
 C_{\alpha}F(h)=\int_{U(A)_\infty}\alpha(k)F(hk)\,dk .
\]
The two compactly supported real smoothing functions depend on
the real edge and \(Q_1,Q_2\), and precede the finite arrays and forms.
The distinguished real types of \(F_1\) and \(F_2\) are conjugate.
Primitivity concerns the combined matrix \(X=(X_1,X_2)\) or
\(Y=(Y_1^{\mathsf t},Y_2^{\mathsf t})^{\mathsf t}\).
The Gaussian estimate below justifies the individual coefficient
projections and proves square integrability of the two lifts.

Consider one edge of source and target ranks \(a,b'\),
with the joint harmonic \(P\) and a rapidly
decreasing source test. At each infinite place use the
lowest compact type coefficient projectors \(e_S,e_T\)
for the source and target lines, with the inverse in the
source convolution as prescribed by the invariant
pairing. Both kernels are integrable with weight
\(\|k\|^{ab'}\). Indeed the opposite-sign gaps in the
types exceed \(500d_*^3\), where \(d_*\) bounds
the individual ranks, by the half-gap comparison of
the base and inserted blocks. For the whole grid take \(d_*=r+1+2mt\), \(t=5r\). Indeed the node
\((z,a_1,a_2)\) has rank \(r+z+(t+z)a_1+(t-z)a_2\), with
\(0\le a_1,a_2\le m\) and \(z=0,1\). This includes the two nodes
at \((a_1,a_2)=(m,m)\). The gap before the inserted blocks is at
least \(\Delta/2\), and the subtractions are bounded by \(5r+2+2m\).
Since \(m=O(\sqrt x)\), \(r\ge x^5\), and \(\Delta=2r^8\),
\[
 500(2d_*)^3+5r+2+2m<\Delta/2
\]
for sufficiently large \(x\). Reversing or conjugating the parameters
preserves the gaps. Section~\ref{sec:real-coefficient-decay}
gives at least the exponent \(ab'+8d_*^2+1\).

The Cartan density of \(U(p,q)\), \(p+q\le d_*\),
in coordinates \(t_1\ge\cdots\ge t_{\min(p,q)}\ge0\)
is a product of powers of \(\sinh\alpha(t)\).
The restricted roots have at most \(2d_*^2\) total
multiplicity and \(|\alpha(t)|\le2t_1\), so the
density is \(O(e^{4d_*^2t_1})\). Since
\(\|k\|=e^{t_1}\) on the Cartan part, the stated
moment is integrable. The compact cases are immediate.

Fix the finite test \(f_f\) in a rational Schr\"odinger polarization.
Its compact support at the finite places restricts the rational
positions to a single fractional lattice \(\Lambda\) in the real
position space. Let \(G\) be the fixed Schr\"odinger Gaussian and let
\(u\) be the inverse Bargmann transform of the weighted Fock vector
\(v\). A coherent vector at \(z\) gives, up to a modulation, the
translate \(G(\,\cdot-t_z)\). Consequently
\[
 \sum_{q\in\Lambda}|f_f(q)u(q)|
 \le \|f_f\|_\infty p_0(v)
          \sup_t\sum_{q\in\Lambda}|G(q-t)| .
\]
The Gaussian series converges locally uniformly in \(t\) and is
\(\Lambda\)-periodic, so its supremum is finite. Here \(u(q)\)
is the pointwise coherent integral. For \(p_0(v)<\infty\) it
gives a continuous position representative by dominated convergence
against \(|v(z)|\). For a Hilbert Fock vector it agrees with the
inverse Bargmann transform by the coherent resolution of the identity.
Tonelli applied to the Gaussian majorant proves the displayed
inequality also for the projected vectors used below.

Write \(d=ab'\) at one real place, let \(v_0\) be the Fock vacuum,
and set \(p_0(v)=\pi^{-d}\int|v(z)|\,dz\), including the Gaussian
in \(v\). For the real linear symplectic transformation \(\sigma\)
of the Hom space, write its positive singular values in pairs
\(e^{t_l},e^{-t_l}\) with \(t_l\ge0\), \(1\le l\le d\).
The two orthogonal symplectic factors in its singular-value
decomposition preserve \(p_0\) and preserve the vacuum line.
In the resulting coordinates the Gaussian kernel of
Section~\ref{sec:original-comparison-gaussian} gives the transformed
vacuum, up to phase, as
\[
 \prod_{l=1}^d (\cosh t_l)^{-1/2}
  \exp\!\left(-\frac{|z_l|^2}{2}
            +\frac{\tanh t_l}{2}z_l^2\right).
\]
Integration in \(\Re z_l,\Im z_l\) contributes a factor
\((1-\tanh^2t_l)^{-1/2}=\cosh t_l\) relative to the unsqueezed
Gaussian. Thus
\[
 \frac{p_0(\omega(\sigma)v_0)}{p_0(v_0)}
       =\prod_l(\cosh t_l)^{1/2}\le\|\sigma\|^{d/2}.
\]
For a coherent vector \(\rho(y)v_0\), covariance gives
\(\omega(\sigma)\rho(y)v_0=\rho(\sigma y)\omega(\sigma)v_0\).
Heisenberg displacements preserve \(p_0\). Applying the coherent
integral first to the input and then the fixed-lattice theta bound
to its image therefore gives, for the theta functional \(L\),
\[
 |L(\omega(\sigma)v)|\le C_L p_0(\omega(\sigma)v)
       \le C'_L\|\sigma\|^{d/2}p_0(v).
\]
The phases of the unitary splitting characters have modulus one.
Reduce the source point as in Section~\ref{sec:adjacent-unfolding}.

For a fixed target
\(h\) the resulting kernel at \((g,hk)\) is bounded
by
\[
             C_h\operatorname{ht}_S(g)^B
                      \prod_{v\mid\infty}\|k_v\|^{ab'} .
\]
This follows from
\(\|\sigma_v\|\le\|(h k)_v\|\|g_v^{-1}\|\) and the
finitely many finite translates in the reduction.
Integrating the rapid source test gives the same
bound in \(k\) for the lift. For either edge write
\(\mathcal L_f(v;h)=\int_{[S]}f(g)K(v;g,h)\,dg\) for its
bilinear source pairing after the determinant correction above.
For \(k\in U(S)_\infty\) and \(k'\in U(T)_\infty\), use
\(R(k)f(g)=f(gk)\). Covariance and the source substitution \(g'=gk\)
give the two translation laws
\[
 \mathcal L_f(\omega(k,1)v;h)=\mathcal L_{R(k^{-1})f}(v;h),
 \qquad
 \mathcal L_f(\omega(1,k')v;h)=\mathcal L_f(v;hk').
\]
Put \(e_S^\vee(k)=e_S(k^{-1})\). In the scalar realization of
the source forms the unit highest line is cut out by
\(R(e_S^\vee)\), with the contragredient coefficient for the
bilinear pairing, as in Section~\ref{sec:first-lift-doubling-series}.
Thus \(R(e_S^\vee)f=f\) on every copy of that line. The analogous
weighted bound in the source real translation gives Fubini against
\(e_S\); the target bound gives Fubini against \(e_T\). Integrating
the two laws yields
\[
 \mathcal L_f(p_Sv;h)=\mathcal L_f(v;h),\qquad
 \mathcal L_f(p_Tv;h)=\bigl(R(e_T)\mathcal L_f(v;\cdot)\bigr)(h).
\]
The calculation holds for every multiplicity vector at the original
level, including the conjugate test on the second edge. All the
translations used in this projection have finite component one.

\subsubsection{The real line projections}
\label{sec:four-unitary-real-projections}

\sectionopening{In this subsection, we derive the real convolution identity, using the uniform bounds of Section~\ref{sec:four-unitary-constant-terms} to justify integration.}

On the real oscillator, let \(p_S,p_T\) be the line
projections given by these coefficients, with the
dual source representation used in the pairing.
We first verify the parameter rule and \(p_S P\ne0\)
for the following block insertions, including the column edges of
Section~\ref{sec:four-unitary-marked-application}. At fixed
offset \(o\) in grid \(z\in\{0,1\}\), a node with coordinates
\((u,v)\) has rank \(N=r+z+u(t+z)+v(t-z)\).
The triangle position \((i,n)\) uses offset \(o=i\) and grid
coordinates \((0,n)\).
Its unlowered doubled entries are the stem \(2\lambda_a+z\)
and the block entries \(2c+d-1-2b\), \(0\le b<d\),
with \(d=t+z\) for an up block and \(d=t-z\) for a down block.
The lowering \(\epsilon_N\) is \(1\) at odd \(u+v\)
in grid \(0\), and \(0\) otherwise. For either kind of
edge the opposite rank \(N'\) has parity \(\epsilon_N\).

For a down insertion at \(c\), the opposite determinant
exponent is \(c+\lfloor N'/2\rfloor\). Its shifted entry
is therefore \((u-\epsilon_N)-2(c+\lfloor N'/2\rfloor)+N'
=u-2c\) for an unlowered entry \(u\). For an up insertion
at \(c\), the conjugate parameter and opposite exponent
\(-c+\lceil N'/2\rceil\) give
\((-u+\epsilon_N)-2(-c+\lceil N'/2\rceil)+N'=-u+2c\).
Thus the matched shifted entries coincide. The inserted
entries form the length \(d\) compact block with a single
signature sign. The signature of a retained entry contains
one factor \(-1\) per absent down centre above it and
per absent up centre below it. Inserting \(c\) changes this
factor precisely for a negative shifted entry. With coordinates
\((o,u,v)\) for the offset and the two block counts, the insertions are
\[
\begin{array}{c|c|c}
 \text{step}&\text{centre}&\text{length}\\ \hline
 (0,u,v)\longrightarrow(0,u+1,v)& w_{u+1}&t+z\\
 (o,u,v)\longrightarrow(o,u,v+1)&-w_{o+v+1}&t-z\\
 (o+1,0,v)\longrightarrow(o,0,v+1)&-w_{o+1}&t-z .
\end{array}
\]
In a rectangular grid \(o=0\); in a triangle \(u=0\).
For the last insertion give the smaller space the opposite form
to put both sign rules over \((-1)^o\).
At stem index \(p\), let \(a=(p+1-k)_+\). The first \(a\)
down centres are above this entry, and \(\min((a-o)_+,v)\)
of them are present. Its sign is
\((-1)^{o+a+\min((a-o)_+,v)}\).
A present up block has no missing centre contributing to its sign.
For a present down block the \(o\) missing down centres above it
cancel the offset sign. Thus the inserted blocks are positive in
the arithmetic signatures. Opposite signs at a node require a
missing centre strictly between the two entries. In doubled
coordinates each entry is at least \(\Delta-t-2\) from that
centre, by the half-gap and the block widths in
Section~\ref{sec:marked-deep-separation}. This gives the undoubled
gap used in Section~\ref{sec:four-unitary-constant-terms}.

For the prescribed opposite characters and parameters
use the substitution
\[
 (l,l',s,s';\gamma,\gamma')
 \longmapsto (-l,-l',-s,s';N-\gamma,N'-\gamma'),
\]
where \(s,s'\) are the two sign vectors for this insertion.
The shifted entries change sign and the smaller signature
returns to its original one. The retained rule is again
equal signs at positive entries and opposite signs at
negative entries. On reversing an arrow the two commuting
actions still act on the same harmonic polynomial. Interchange of
the two factors in the bilinear theta pairing takes both coefficient
representations to their contragredients. The two kernels of the
reversed arrow are consequently \(\overline{e_T},\overline{e_S}\)
in this order.

To apply Li's regular discrete-series correspondence \cite{LiThetaCohomology}
in the form of \cite[Theorem~2.7(b), p.~641]{PaulFirstOccurrence}, remove the
determinant twists and use the tensor realization
\(\operatorname{Hom}(S,T)=S^\vee\otimes T\):
on the source dual the compact weights change sign and their order
in each compact factor reverses. More explicitly, if \(s_a,s'_b\)
are the two signatures, the vacuum weights after removing the scalar
twists are \(-s_a\sum_b s'_b/2\) and \(s'_b\sum_a s_a/2\)
in the Hom model. Dualizing the first factor gives exactly the two
vacuum weights of \cite[Lemma~2.5]{PaulFirstOccurrence}.
The doubled stem entries have parity \(z\), since \(k\) is even
and the \(\lambda_a\) are integral. Every block has the same
parity, since \(t\) is odd and its length is \(t\pm z\).
Shifting by \(2c\) and negating preserve this parity. Thus for
rank difference \(d=t\pm z\) the genuine Harish--Chandra entries
lie in \(\mathbf Z+(d-1)/2\), as required on both groups. Let \(n_{\epsilon,\eta}\) count the retained entries
with signature \(\epsilon\) and shifted sign \(\eta\) on the smaller side.
They contribute
\[
     (n_{+,+}+n_{-,-},\ n_{-,+}+n_{+,-})
\]
to the larger signature. The remaining entries are
\((d-1)/2,(d-3)/2,\ldots,(1-d)/2\), all in one compact factor.
In fact \(|\xi_a|\ge\Delta-1\) at every retained index.
For a tail entry of excess \(h\), its difference from the doubled
down centre \(-2w_j\) is \((2j-2h-1)\Delta+z\), by the stem
formula in Section~\ref{sec:marked-deep-separation}. An up centre
is positive and has distance greater than \(\Delta\) from the tail.
The doubled head has modulus at most \(k+1\), whereas
\(w_j\ge\Delta^2\). For an earlier block the distinct centres
differ by at least \(\Delta\), and its doubled half-width is
at most \(t\). This gives \(|\xi_a|\ge 2\Delta-t\) on all
retained blocks. In particular \(|\xi_a|/2\ge(d+1)/2\).
Use the whole-grid rank bound
\(N,N'\le d_*=r+1+2mt\) of
Section~\ref{sec:four-unitary-constant-terms}. Its inequality
\(500(2d_*)^3+5r+2+2m<\Delta/2\) implies the required undoubled
separation \((\Delta-t-2)/2>500(N+N')^3\).
The bound \(|\xi_a|/2\ge(d+1)/2\) and the two choices of the
signature of the inserted block are the hypotheses of
\cite[Theorem~2.7(b)]{PaulFirstOccurrence}. The stronger separation
is used for the coefficient moment in
Section~\ref{sec:four-unitary-constant-terms}. Simultaneous reversal of the two Hermitian
forms preserves the correspondence. Reversing just one form
conjugates the oscillator and takes the contragredient on both
representation factors
\cite[Proposition~2.9 and Corollary~2.10]{PaulFirstOccurrence}.

Write \(\xi_a\) for the doubled shifted entry at a retained index
and \(\xi_b\) also for the doubled entries of the inserted string.
Let \(C_a\) and \(R_b^{\rm wt}\) be the relative lowest compact
weights of the dual smaller factor and the larger factor. The
lowest-type formula \(\lambda+\rho-2\rho_K\) and the vacuum
weights above give
\[
 \begin{aligned}
 2C_a&=-\xi_a+s_a
       \left(\sum_{a'}s_{a'}\operatorname{sgn}(\xi_a-\xi_{a'})
                                    +\sum_{b'}s'_{b'}\right),\\
 2R_b^{\rm wt}&=\xi_b-s'_b
       \left(\sum_{b'}s'_{b'}\operatorname{sgn}(\xi_b-\xi_{b'})
                                    +\sum_{a'}s_{a'}\right).
 \end{aligned}
\]
Here \(a'\) ranges over the smaller parameter and \(b'\) over the
larger parameter. Set \(\eta_a=\operatorname{sgn}\xi_a\), so
\(s'_a=s_a\eta_a\) at a retained index. For another retained index
\(a'\), its contribution to \(2(R_a^{\rm wt}+C_a)\), divided by
\(s_as_{a'}\), is
\[
 (1-\eta_a\eta_{a'})\operatorname{sgn}(\xi_a-\xi_{a'})
                            +\eta_{a'}-\eta_a=0 .
\]
For an inserted index \(b\), the contribution divided by \(s_as'_b\)
is \(1-\eta_a\operatorname{sgn}(\xi_a-\xi_b)=0\).
At an inserted index \(b\), the retained terms in \(2R_b^{\rm wt}\)
cancel because \(\operatorname{sgn}(\xi_b-\xi_a)=-\eta_a\).
The inserted entries have one compact sign and spacing two, so
\[
 s'_b\sum_{b'\ {\rm inserted}}s'_{b'}\operatorname{sgn}(\xi_b-\xi_{b'})
 =\xi_b .
\]
Thus \(R_a^{\rm wt}=-C_a\) at retained indices and \(R_b^{\rm wt}=0\)
on the inserted string. Also
\[
 |2R_a^{\rm wt}-\xi_a|\le N+N'-1<|\xi_a|,
 \qquad \operatorname{sgn}R_a^{\rm wt}=\eta_a .
\]
In either compact order let \(L_i^*\) be the decreasing doubled
parameter entries and \(\mu_i\) the relative weights. For \(i<j\)
of the same signature, let \(c\) be the number of indices of that
signature in \((i,j]\). Their common vacuum correction cancels.
The successive entry gaps are at least two, so
\[
 \mu_i-\mu_j=(L_i^*-L_j^*)/2+(j-i)-2c\ge 2(j-i-c)\ge0 .
\]
In each sign quadrant take the nested products of minors on the
matched positions of degree \(|R_a^{\rm wt}|\ge h\), \(h\ge1\).
These are initial segments in the compact order for positive weights
and in the dual order for negative weights. In the signed Fock
coordinates every compact raising operator kills a minor or replaces
an index by one already present in it. At the matrix with ones on the
matched pairs and zeros elsewhere the minors have determinants
\(\pm1\). Each source or target index occurs with one oscillator sign,
so the mixed derivatives with opposite signs and a common source
or target index kill \(P\). The two lowest compact types occur
together with multiplicity one in the joint harmonics by
\cite[Theorem~2.2(1) and Lemma~2.5]{PaulFirstOccurrence}.

The positive Bargmann pairing and compact Schur averages of
Section~\ref{sec:original-comparison-gaussian} now apply
to either group as source. On \(\operatorname{Hom}(S,T)\), a
positive source matrix \(g\) acts by \(X\mapsto Xg^{-1}\)
and a positive target matrix \(h\) acts by \(X\mapsto hX\).
Their matrices on the tensor factors are \(g^{-\mathsf T}\) on
the dual source and \(h\) on the target, both Hermitian positive.
At either positive factor the determinant twist is one, the pure quadratic terms
vanish by the mixed derivative calculation, and the
mixed term gives the squared Fock norm after an invertible
substitution, times the positive vacuum coefficient.
The two compact Schur averages against the source matrix
coefficient consequently have strictly positive integral
at each positive factor. Integrate in the polar decomposition
and multiply by the positive formal degree.
The coefficient moment verified above justifies this step
and gives \(\langle P,p_S P\rangle>0\).

Write \(c_S(g)=\langle v_S,\pi_S(g)v_S\rangle\) for the coefficient
in the source oscillator convention. The convolution calculation of
Section~\ref{sec:real-reproducing-coefficient} gives
\(p_S\omega(g)p_S=c_S(g)p_S\), and similarly for the target.
For any unit \(u\in\operatorname{im}p_S\), its source translates
therefore have the same Gram matrix as the translates of \(v_S\).
Their closed span is a copy of \(V_S\) under
\(\pi_S(g)v_S\mapsto\omega(g)u\). In particular \(p_SP\) belongs
to the full source discrete-series isotypic component.

Let \(Q_S\) project onto this component of the oscillator Hilbert
space \(\mathcal H_\omega\).
This is a central projection in the source von Neumann algebra. On
\(Q_S\mathcal H_\omega=V_S\widehat\otimes M\), the source algebra
is \(\mathcal B(V_S)\overline\otimes1\). Howe's commutant theorem
\cite[Theorem~6.1, p.~550]{HoweDuality}, restricted by this central
projection, identifies the target algebra with
\(1\overline\otimes\mathcal B(M)\). Thus \(M\) is irreducible and
admissible. Scalar twists leave these algebras unchanged.

Here admissibility follows from Harish--Chandra's admissibility
theorem for irreducible unitary representations of real reductive
groups, applied to the target and its metaplectic cover.
The map \(Q_S\) intertwines both derived actions on Fock polynomials
and has dense polynomial image in \(V_S\widehat\otimes M\).
Projection onto any fixed \(K_S\times K_T\)-type preserves the
polynomials, since the compact action preserves degree. The
image fills this finite-dimensional isotypic space by density and
admissibility. Consequently \((V_S)_{K_S}\boxtimes M_{K_T}\) is a
quotient of the polynomial oscillator module. Its target factor is
the one identified by the preceding discrete-series correspondence,
using the smaller group first and reciprocity for the reversed
orientation.

On this component \(p_S=P_{v_S}\otimes 1\), where \(v_S\) spans
the source compact highest line. Since \(p_S\) commutes with the
target action, \(0\ne p_SP\) has the target compact highest weight
computed above. For a noncompact real factor the corresponding
lowest compact type occurs once by
\cite[Theorem~8.5, p.~42]{AS} on the derived group.
Adjoining the compact central torus with its fixed character
preserves this multiplicity. For a definite factor the highest line
is one-dimensional by highest-weight theory. Thus the highest line
of the target type is one-dimensional and \(p_Tp_SP=p_SP\).

The nonzero vector \(p_SP\) is cyclic in the target
representation. Smearing its translates by smooth
compactly supported functions has dense linear span.
Project this span onto the finite sum of target compact
types of \(p_S Q\), for a polynomial \(Q\). This sum is
finite dimensional by admissibility. The projected span
is therefore the whole sum, and its projection is still
obtained by smooth compactly supported convolution, by
averaging the convolution functions over the target
compact group with the compact-type character projectors.
Thus one such function sends \(p_SP\) to \(p_SQ\).

Use the weighted Fock estimate of Section~\ref{sec:four-unitary-constant-terms}
to evaluate this identity by theta summation. On a compact neighborhood
of any real group element the Gaussian quadratic forms of the moved
vacuum have uniformly positive real part in the decaying exponent.
Consequently the moved vacuum is continuous in \(p_0\) by dominated
convergence. For each fixed displacement the same holds for its
translate, by covariance and continuity of Heisenberg translations
and modulations in \(L^1\). Their \(p_0\)-norms are bounded on this
neighborhood independently of the displacement. Integrating against
the coherent coefficient of \(v\) gives strong \(p_0\)-continuity,
with \(p_0(\omega(g)v)\le C\|g\|^{ab'}p_0(v)\) as above. In the Hilbert norm the action is
strongly continuous and unitary. The coefficient moment proved
in Section~\ref{sec:four-unitary-constant-terms} and the compact
support of the smoothing function therefore give convergence
of the source projection and smoothing integrals in
\(L^1\cap L^2\) with norm \(p_0+\|\cdot\|_2\).
Thus \(p_SQ=\omega(\alpha)p_SP\), for the chosen smoothing
function \(\alpha\), holds in this space: the \(L^2\)
identity identifies the two \(L^1\) representatives.

At fixed finite datum, the coherent Gaussian lattice
bound in Section~\ref{sec:four-unitary-constant-terms}
makes theta evaluation continuous for \(p_0\). At the
translated source points its bound is multiplied by at most
the fixed power of the source height displayed there,
which is integrable against the rapidly decreasing source
test. Hence we may apply the pointwise theta lift to the
smoothing identity and interchange the convergent integrals.
Taking the product at the infinite places gives the asserted
replacement of the real polynomial lifts. The same weighted
coefficient bound applies to the target projector \(p_T\).
Evaluating \(p_Tp_SP=p_SP\) therefore gives
\(R(e_T)F=F\) as an absolutely convergent pointwise integral
for the lift \(F\) of \(P\), at the same finite datum and
with the same quotient Haar measures.

Fix a compact open subgroup \(J\) of the finite target group
that fixes the chosen finite Weil vector with its splitting
character. Smoothness gives such a \(J\). Then \(F\) and all
\(F_j^{\rm tr}=F\mathbf 1_{\{|F|\le j\}}\) are right
\(J\)-invariant. Put \(F_j=R(e_T)F_j^{\rm tr}\).
The quotient has finite volume, so
\(F_j^{\rm tr}\in L^2\), and convolution gives
\[
  \|F_j\|_\infty\le j\|e_T\|_1,\qquad R(e_T)F_j=F_j .
\]
The last identity follows from \(e_T*e_T=e_T\) by Fubini.
The convolution gives a continuous representative on each real
component, since translations of \(e_T\) are continuous in \(L^1\).
Real convolution commutes with the finite right action of \(J\), so
the same \(J\) fixes every \(F_j\).

Apply Section~\ref{sec:paired-lifts-projected-tests} with \(y=F_j\).
On each arithmetic component its coefficient integral intertwines
the Harish--Chandra module of the target discrete series with
right translation and sends the distinguished vector to \(F_j\).
The matrix coefficients and their derivatives are integrable, so
the image consists of smooth, compact-finite automorphic forms
of the given infinitesimal character. Each is bounded by
\(\|F_j\|_\infty\) times the \(L^1\)-norm of its coefficient,
including the formal-degree factor, and hence is square-integrable.
Restricting to the arithmetic quotients of the derived group
as in Section~\ref{sec:paired-lifts-projected-tests}, Wallach's
theorem gives cuspidality of this tempered module
\cite[Theorem~4.3, p.~233]{WallachConstant}. The adelic constant
terms are the corresponding real integrals by the same unipotent
approximation. All \(F_j\) have the same
finite level, archimedean compact type and infinitesimal character.
As in Section~\ref{sec:paired-lifts-projected-tests}, there are finitely
many arithmetic quotient components. On each, the compact centre
acts by the fixed character, so restriction to the derived group
determines the form. The restricted compact types form a finite set,
and the kernel of the derived infinitesimal character has codimension
one in the centre of its enveloping algebra. Intersecting the
arithmetic subgroup of each component with the special unitary
group gives an arithmetic subgroup of the derived group. Inversion
identifies its left quotient with the right quotient in
\cite[Ch.~I, \S2, Theorem~1]{HarishChandraAutomorphic}, intertwining
right translations with left translations. Boundedness gives the
required moderate growth. Finiteness on these components therefore
places all \(F_j\) in a single finite-dimensional cusp space.

For each fixed \(h\), all truncations have the integrable majorant
\[
  C_h|e_T(k)|\prod_{v\mid\infty}\|k_v\|^{ab'}
\]
by the preceding coefficient moment bound.
Hence \(R(e_T)F_j^{\rm tr}(h)\to R(e_T)F(h)=F(h)\).
Take finitely many evaluations that are injective
on the finite-dimensional space. Convergence at
these points gives convergence of the coefficients
in a basis, so \(F\) lies in that same space.
In particular the two individual upper lifts and
the two transposed lower lifts are square-integrable.
They are rapidly decreasing by Section~\ref{sec:paired-lifts-projected-tests}.

For fixed real elements \(s,t\) of the source and target groups,
write \(\omega_S(s)=\omega(s,1)\), \(\omega_T(t)=\omega(1,t)\).
The central source projection \(Q_S\) commutes with the source
compact group. Hence \(Q_SP\) has the compact highest weight of
\(P\), that of the lowest compact type of \(V_S\). The same
multiplicity result gives one occurrence in \(V_S\), with highest
line \(\mathbb C v_S\).
Thus \(Q_SP\in\mathbb C v_S\widehat\otimes M\) and
\(Q_SP=p_SP\), for the full multiplicity space \(M\).
Since \(p_S=p_SQ_S\) and \(Q_S\) commutes with \(\omega_S(s)\),
the source compression and the two translation laws give
\[
 p_S\omega_S(s)P=c_S(s)p_SP,\qquad
 \mathcal L_f(\omega_T(t)\omega_S(s)P)=c_S(s)R(t)F .
\]
With \(e_T^t(k)=e_T(t^{-1}kt)\), substitution \(k=tlt^{-1}\)
in the convolution gives
\[
 R(e_T^t)R(t)=R(t)R(e_T).
\]
Right translation by \(t\) preserves the finite level and commutes
with the left unipotent integrations defining the constant terms.
The lattice norms at \(h\) and \(ht\) differ by bounded factors
depending on \(t,t^{-1}\), so \(H(ht)\asymp_t H(h)\) and
\(R(t)F\) is again rapidly decreasing.

Write \(e_A\) for the normalized
matrix coefficient of right translation on the real line of
\(F_1\). With the projector kernel \(e\) of this type, its
argument convention is \(e_A(g)=e(g^{-1})/\overline{e(1)}\).
Then
\[
 \int_{[A]} F_2(u)F_1(ug)\,du
      = e_A(g)\int_{[A]}F_2(u)F_1(u)\,du
       \qquad(g\in U(A)_\infty).
\]
Indeed \(R(e)F_1=F_1\), and the conjugate distinguished type
of \(F_2\) gives \(R(e)\overline{F_2}=\overline{F_2}\).
The integral on the left is
\(\langle\overline{F_2},R(g)F_1\rangle\).
Inserting the self-adjoint projector on both sides and using
\(R(e)R(g)R(e)=e_A(g)R(e)\) from
Section~\ref{sec:real-reproducing-coefficient} proves the identity
on every multiplicity space.

The change \(u=hk_2\) preserves the Haar measure of \([A]\).
Consequently, for real edge data \(Q_1,Q_2\),
\begin{align*}
 H_j(Q_1Q_2)
 &= r(Q_1,Q_2)\langle z,T_X f\rangle_{[U_2]},\\
 r(Q_1,Q_2)
 &=\int_{U(A)_\infty^2}\alpha_1(k_1)\alpha_2(k_2)
        e_A(k_2^{-1}k_1)\,dk_1\,dk_2 .
\end{align*}
Fubini here also follows directly from Cauchy--Schwarz,
unitarity of the real translations on \([A]\), and compact
support of the smoothing functions. The integrals are
bilinear, so neither smoothing function is conjugated.
Every polynomial in both blocks is a sum of pure tensors,
and the same conclusion holds with a scalar \(r(Q)\).

To approximate \(v_+\), use the signed Fock coordinates
\(\zeta(w)\) and Gaussian \(G(w)=e^{-|w|^2/2}\).
For each \(m\), the degree-\(m\) Taylor kernel
\[
 G(w)G(y)\sum_{k=0}^{m}
             \frac{(\overline{\zeta(y)}\cdot\zeta(w))^k}{k!}
\]
integrated against \(v_+(y)\,dy/\pi^d\) gives a polynomial
Fock vector \(v_m(w)\), where \(d\) is the complex dimension.
The full exponential is the Fock reproducing kernel, by
the Gaussian integral for orthogonal monomials.
The absolute value of each partial kernel and of its limit
is at most \(e^{-(|w|-|y|)^2/2}\).

For every \(N\) there is
a constant \(C_{d,N}\) such that
\[
 \int_{\mathbb C^d}(1+|w|)^N e^{-(|w|-|y|)^2/2}\,dw
 \le C_{d,N}(1+|y|)^{N+2d-1}.
\]
Indeed, use polar coordinates and bound the polynomial in
the radial variable \(r\) by a constant times
\((1+|y|)^{N+2d-1}(1+|r-|y||)^{N+2d-1}\).
For \(d=0\) the Fock space consists of constants.
Thus dominated convergence in \((y,w)\) proves
\(p_N(v_m-v_+)\to 0\) for every \(N\), using the finite
moments of \(v_+\). Set
\(B_j=\langle z,T_X f\rangle_{[U_2]}\).

If some \(B_{j_0}\ne0\), the identities for polynomials give
\[
 H_j(Q)B_{j_0}=H_{j_0}(Q)B_j .
\]
Continuity for each \(j\) gives the same formula at \(v_+\);
put \(\kappa=H_{j_0}(v_+)/B_{j_0}\in\mathbf C\). If all the \(B_j\)
vanish, all the \(H_j(v_+)\) vanish, and we may take
\(\kappa=0\). In both cases
\[
       H_j(v_+)=\kappa\langle z,T_X f\rangle_{[U_2]}
       \qquad\text{for every }j.
\]
This gives one \(\kappa\) for the entire family.

Integrate \eqref{eq:marked-square-paired-period} at
\(g=g_1\oplus g_2\) against
\(f_1(g_1)f_2(g_2)\nu_0(\det g_2)\).
Expanding the two finite point arrays and applying the preceding
identities yields
\[
  \kappa\sum_X\Phi(X)\langle z,T_X f\rangle_{[U_2]}
  =c_\xi\sum_Y\Psi(Y)\langle z,T'_Y f\rangle_{[U_2]} .
\]
The target pairing is nondegenerate. Thus the two composites
on the original four Hermitian spaces satisfy the matched operator identity
\[
          \sum_Y\Psi(Y)T'_Y
       =\frac{\kappa}{c_\xi}\sum_X\Phi(X)T_X .
\]
This argument permits \(\kappa=0\). Its needed nonvanishing
for the forward triangle follows from injective first returns
and the weighted chain argument.

\subsubsection{Application to the unitary groups}
\label{sec:four-unitary-marked-application}

\sectionopening{In this subsection, we apply the two pairings of Section~\ref{sec:four-unitary-real-projections} with the original levels and measures.}

For the forward square the nodes \(U_1,U_2,A,A'\) are
\[
 (i,n+1),\qquad(i+1,n+1),\qquad(i,n+2),\qquad(i+1,n),
\]
respectively. These are triangle positions, so each has first
grid coordinate zero. For grid parameter \(z_0\in\{0,1\}\),
the rank in column \(j\) is \(d_j=r+z_0+j(t_0-z_0)\).
Hence
\(\dim W=2d_{n+1}=d_n+d_{n+2}\) and the exchange
defect is \(t=t_0-z_0\).
The signs of the added blocks are positive. In the stem,
the head stays on the same side of every inserted centre at
each infinite place. The relative positions of the tail entries
and centres are independent of \(\tau\), since the common
\(\kappa_\tau\)-shift cancels. Thus all negative counts here are
constant across the infinite places. On a fixed base coordinate
with excess \(a\), the sign in column \(j\)
at offset \(i\) is
\[
            (-1)^{i+a+\min((a-i)_+,j)} .
\]
Moving \((i,n+1)\) to \((i,n+2)\) and \((i+1,n+1)\)
to \((i+1,n)\) flips both signs exactly when \(a-i\ge n+2\).
In that case the two old signs are opposite. Thus the number
of negative signs agrees between the two diagonals, as the
period comparison requires. That number is at most
\(2(r+z_0)<d_{n+1}\), using \(t_0\ge r+6\).
The inverse quotient character in the comparison has weight
\((n+1)\) times the original gap at each real place. The
inequality \(d_{n+1}\le(n+1)\Delta\) gives
\(\dim W\le 2(n+1)\Delta+2q\), with \(q\) the number
of negative signs in \(W\).

For the reversed square in the \(z_0=0\) triangle,
the source and target are \((i+1,n)\) and \((i,n+2)\).
The upper and lower middle nodes in the ordered comparison
are then \((i+1,n+1)\) and \((i,n+1)\). Their ranks are
equal, so \(t=0\) in the measured period comparison.
Interchanging the two diagonals preserves their common negative
index \(q\) and their common total rank \(2d_{n+1}\).
The inverse quotient character has weight \((n+1)\Delta\), so
\[
       2q+2\le 2d_{n+1}\le2(n+1)\Delta+2q .
\]
Apply the same paired calculation with the two reversed Fourier
edge data, with the two row blocks in their prescribed order. The modulus of the
period scalar is constant on this family.

To identify the adjoint consequence, let \(A_0,A_1\) be
the two reversed lifts with the Fourier data of the original
test and its extended test, \(B_0,B_1\) the corresponding
transposed descents, and \(L_+,L_-\) the upper and lower
ordinary lifts. Let \(S_n,S_{n+1}\) be the two original
central reversal translations. Write the two reflection identities,
the reversed square and the central intertwining as
\begin{align*}
 \langle B_0x',y\rangle&=d_0\langle x',A_0 S_n^{-1}y\rangle,\\
 L_+ A_0&=g_0 A_1 L_-,\\
 L_-S_n^{-1}&=\eta S_{n+1}^{-1}L_-,\qquad |\eta|=1,\\
 \langle x,A_1 S_{n+1}^{-1}u\rangle&=d_1\langle B_1x,u\rangle .
\end{align*}
For \(x'=L_+^*x\), the second-linear convention gives
\(\langle x',v\rangle=\langle x,L_+v\rangle\).
It follows directly that
\[
     \langle B_0x',y\rangle
       =d_0g_0\eta d_1\langle B_1x,L_-y\rangle .
\]
In particular none of the four scalars is conjugated in this
passage. The factors in \(d_0,d_1\) come from the two matrix
reversals \eqref{eq:display-mr}. For matched arrays the scalar
of the reversed square is \(\kappa/c_\nu\), with \(\kappa\) the real
projection scalar and \(c_\nu\) the equal-rank period scalar.
For the Fourier data in \(A_0,A_1\),
Section~\ref{sec:period-finite-determinant-sum} gives
\[
 |g_0|=\frac{|\kappa|}{|c_\nu|\,|\mathcal G|\,|R|^{d_n t_0}},
 \qquad R=\mathcal O_E/(\pi_MQ(\xi)),\qquad |\mathcal G|^2=|R|^{t_0^2}.
\]
Here the triangle has \(z_0=0\), so its column increment is \(t_0\).

\section{Real restrictions and self-dual Fourier tests}
\label{sec:local-period-transforms}

\sectionopening{In this section, we compute the local restrictions of the paired periods of
Section~\ref{sec:marked-square-analytic-pairing}. The real slice and the
finite self-dual transform give the normalization for the exchange in
Section~\ref{sec:normalized-period-exchange}.}

\subsection{The real restriction of the larger period}
\label{sec:marked-real-slice-period}

\sectionopening{In this subsection, the real integrals in the paired periods of
Section~\ref{sec:marked-square-analytic-pairing} involve the original
Fock vectors. We evaluate the larger and smaller real slices with
their quotient measures to determine the factor in the exchange.}

Let \(W_0=H\perp W\) have rank \(2j\), let \(B_0=U^-\)
have rank \(b\), and let \(X=H^b\perp W_0\). Use the original
two isotropic bases of \(X\) and the original mixed
Schr\"odinger--Fock models for its pairs with
\(B_0\perp W_0\) and \(B_0\). Fix the lower real Schwartz
vectors \(v\) for \((W_0,B_0)\) and \(v_1\) for
\((W_0,W_0)\). Fix the common infinity weight
\(\mu=(\mu_\tau)_\tau\). It is positive on the admissible
family, so the exponents \(2\mu_\tau\) below are nonnegative
integers. All matrices in the real formula are families at
the infinite places. Let \(J=\operatorname{diag}(J_1,J_2)\)
be the matrix of \(B_0\perp W_0\), and set \(x_0=(1;0)\).

\medskip
\noindent\emph{The two real vectors.}
Near \(x_0\), choose a smooth unitary \(k(y)\) carrying
\(\operatorname{im}x_0\) onto \(\operatorname{im}y\).
For the graph of \(T:\operatorname{im}x_0\to W_0\), the
derivative at zero of the graph chart on the action of
\[
 \exp\begin{pmatrix}0&-J_1^{-1}T^*J_2\\T&0\end{pmatrix}
\]
is \(T\mapsto T\). The real inverse function theorem produces
the section on a neighborhood. Multiplying its Lie-algebra
coordinates by a cutoff and exponentiating extends \(k\)
smoothly while keeping it unitary. If \(b>0\), take at each real
place a ball about \(x_{0,\tau}\) with closure in this chart such
that on the closure
\[
 \|(x_{0,\tau}^\dagger x_{0,\tau})^{-1}
      (x_\tau^\dagger x_\tau-x_{0,\tau}^\dagger x_{0,\tau})\|<1/2 .
\]
Such balls exist by continuity and nondegeneracy of the fixed first
block. Take a smooth nonnegative \(\sigma\), supported in their product,
with \(\sigma(x_0)=1\). For an empty first block use its unit function.
These choices use only the prescribed real bases.

For \(x^\dagger=x^*J\) choose a smooth matrix function
\(\iota\) equal to inversion on a neighborhood of
\(x_0^\dagger x_0\) and of all the Gram matrices from the
support of \(\sigma\). This is obtained by multiplying
inversion by a smooth cutoff supported where the determinant
is nonzero. Set
\[
 P(x)=x\,\iota(x^\dagger x)x^\dagger x_0,\qquad
 x^\flat=k(P(x))^{-1}x.
\]
The large real vector \(\Phi\) and small real vector \(\Lambda\)
have the following slices in the initial polarization:
\begin{align}
 \Phi[x]
  &=\kappa\,\sigma(P(x))\sigma(x^\flat)
       \prod_\tau\overline{\det (x^\flat_\tau)_1}^{\,2\mu_\tau}
       \ \omega_{j-\mu}(1,k(P(x)))(v\boxtimes v_1),        \label{eq:marked-real-slice-big}\\
 \Lambda[B]
  &=\kappa\,\sigma(x_0B)
       \prod_\tau|\det B_\tau|^{2(j+\mu_\tau)}\,v .         \label{eq:marked-real-slice-small}
\end{align}
In the first formula \((x^\flat_\tau)_1\) is the first
\(b\)-block. The second-group determinant exponent in
\(\omega_{j-\mu}\) is \(j-\mu_\tau\) at each place.

These sections are Schwartz in the mixed model. In the first
formula \(x=k(P(x))x^\flat\), so the support is contained
in the compact image of the two bump supports. All
scalar factors are smooth, including at rank-deficient matrices,
since no inverse is taken outside the invertible Gram
neighborhood. On that compact support the Weil operators
form a smooth compact family on Schwartz functions.
For every derivative and Schwartz seminorm this follows
locally from the linear-change, quadratic-multiplier
and Fourier formulas for the metaplectic action.

The small section is supported on a compact set because
\(B\mapsto x_0B\) is a closed linear embedding.
The mixed-model isomorphism therefore gives Schwartz
vectors with these slices. Equivalently, in signed Fock
coordinates the condition on the real source can be stated as
finiteness of all weighted absolute moments. To see that
this condition gives the Schwartz estimates, use the
coherent-state reproducing integral. The modulus of its
kernel in real coordinates \(x,y\) is
\(e^{-|x-y|^2/2}\), and differentiation contributes a
polynomial in \(x,y\). For any prescribed weighted
derivative, bound the resulting Gaussian times the
weight by a fixed polynomial in \(|y|\), uniformly in
\(x\), and use finiteness of the corresponding moment.

In the opposite source polarization the zero slice is the
integral of the slices in the initial polarization, with
the Fourier measure of the chosen bases. To compute its
normalization for \(\Lambda\), write
\(t_\Lambda(B)=\sigma(x_0B)\prod_\tau|\det B_\tau|^{2(j+\mu_\tau)}\).
Choose \(\kappa\) so that the real integral of \(\kappa t_\Lambda\)
is the real integral of the scalar vacuum slice. The integral
of \(t_\Lambda\) is positive since it is continuous, nonnegative,
and equal to one at the identity.

The Gaussian vacuum integral is nonzero. With any finite
complementary test, the adelic opposite zero slice factors
into this real integral and the same finite slice integral
as for the vacuum. The linear change of slice variables
carries the self-dual Haar to the same positive multiple of the
product real and finite Haar in these two evaluations.
Equality of the real integrals therefore gives equality for
every finite test. The vacuum has the normalized opposite
zero slice equal to the smaller
vacuum. Consequently the opposite zero slice of \(\Lambda\)
is \(v\). This normalization uses the actual mixed model and
places no condition on the finite data later used in the period.

\medskip
\noindent\emph{A single translate for every test.}
Fix a real point \(h\in U(B_0\perp W_0)_\infty\) and put
\(x_h=h^{-1}x_0\), \(y=P(x_h)\), and \(\rho(h)=\sigma(y)\).
Here \(P\) uses the true inverse, since the Gram matrix of
\(x_h\) equals that of \(x_0\). The function \(\rho\) is
nonnegative and equals one at the identity. It depends
only on the plane \(\operatorname{im}(h^{-1}x_0)\).

If \(\rho(h)\ne0\), \(y\) is a basis spanning the plane
of \(x_h\). Since \(k(y)x_0\) and \(x_h\) have the same
Gram matrix and span, there is \(u_1\in U(B_0)_\infty\)
such that \(k(y)x_0u_1=x_h\). The isometry \(hk(y)\)
preserves the first summand and its orthogonal complement.
Write
\[
 hk(y)=\operatorname{diag}(u_1^{-1},u_2),
 \qquad p=\operatorname{diag}(u_1,u_2^{-1}).
\]
Regard \(p\) as an adelic point with finite component \(1\).
For every invertible \(B\) in the chosen adelic basis,
the relevant real basis is
\[
 x=(ph)^{-1}x_0B=x_hu_1^{-1}B,\qquad
 k(y)^{-1}x=x_0B,\qquad phk(y)=1 .
\]
If \(\sigma(x_0B)\ne0\), the Gram matrix of \(x\) is that of
\(x_0B\) and the inversion in \(P(x)\) is genuine.

The
orthogonal projection is independent of the choice of basis
of the same plane, so \(P(x)=y\).
If \(\sigma(x_0B)=0\), a nonzero value in
\eqref{eq:marked-real-slice-big} would force
\(k(P(x))^{-1}x\) to lie in the bump. Its Gram matrix is
the Gram matrix of \(x\), so again the projection is genuine
and \(P(x)=y\). This makes the second bump
\(\sigma(x_0B)=0\), a contradiction.
If \(\rho(h)=0\), choose \(p=1\). The same argument
shows that the large slice at \((ph)^{-1}x_0B\) vanishes,
since otherwise the first bump would equal \(\rho(h)\).
\begin{lemma}[Real slice at a common translate]\label{lem:common-real-slice}
For each \(h\in U(B_0\perp W_0)_\infty\), there is
\(p=p(h)\in(U(B_0)\times U(W_0))_\infty\), with finite component \(1\),
such that for every invertible \(B\) in the chosen adelic basis,
\begin{equation}
 \omega_{j-\mu}(1,ph)\Phi[(ph)^{-1}x_0B]
    =\rho(h)\kappa\,\sigma(x_0B)
       \prod_\tau\overline{\det B_\tau}^{2\mu_\tau}
       \ (v\boxtimes v_1).                                  \label{eq:marked-real-slice-cancel}
\end{equation}
\end{lemma}
\medskip
\noindent\emph{The product and denominator.}
Let \(\nu\) be an admissible split tame character of weight
\(\mu\), with residue character \(\chi\). Its tame support
consists of the split pairs represented by a nonempty set \(S\);
\(\chi_w\ne1\) for every \(w\in S\).
Put \(\gamma_0=\nu^{-1}\upsilon^j\). For \((W_0,W_0)\) at the
finite places use the original determinant datum
\(\phi_\chi=\mathbf1_{\mathrm{GL}_{2j}(O_S)}
\chi^{-1}(\det)\) at \(S\), with the original vacuum off
the support. Let \(F\) be any finite Schwartz vector
for \((W_0,B_0)\). The real vectors are evaluated against
the theta distribution on the finite vector
\(F\boxtimes\phi_\chi\) and its small member \(F\),
respectively. In the large functional the Heisenberg
translate is \((c,0)\), the group point is
\((g,y_0\oplus l)\), and the splitting characters are
\((\gamma_0,\gamma_0\beta_0)\).

On the small functional
these data read \(c,(g,y_0)\) and \((\gamma_0,\beta_0)\),
for any remaining splitting character \(\beta_0\).
The block product of the rational lattice sums, together
with the direct-sum Heisenberg action, factors the theta
functional of \(v\boxtimes v_1\) into that small functional
times the equal-rank kernel of \(v_1,\phi_\chi\) at \((g,l)\)
with characters \((\gamma_0,\gamma_0)\).
This formula is for the original coordinates.
The second Heisenberg move is zero and the determinant
twist is the product of the twists of the two blocks.

The equal-rank kernel has absolutely convergent period
on \([W_0]\) by the primitive support of \(\phi_\chi\)
and the target-bounded theta estimate. Its period
is independent of \(g\), as follows. For a source element
\(s\) supported at the chosen split places, simultaneously
translate the target by \(s\). The split Weil action on the
square finite matrix reads
\(\phi(X)\mapsto\phi(s^{-1}Xs)\) with scalar one, since
the two ranks and characters are equal.

Set \(\lambda(a)=\nu^{\rm bc}(a)|a|^j\) on the split
local factors. If \(I_g(\phi)\) is the integrated kernel,
the split target action and right invariance give
\[
 I_g\bigl(\phi(l^{-1}\,\cdot)\bigr)=\lambda(\det l) I_g(\phi).
\]
Thus \(I_g\) on compactly supported locally constant
functions on the invertible matrices depends through
\(\int\lambda(\det X)\phi(X)\,d^\times X\).
Indeed, choose a common compact open subgroup \(K\)
on whose right cosets the tests are constant and on which
\(\lambda\circ\det=1\).

The equality on each coset
\(lK\) follows from the target-action formula, giving
the assertion by finite linearity. This weighted integral
is invariant under conjugation by \(s\). The determinant
period is therefore right invariant at the chosen split
places, and off them it is invariant under the original
standard compact level away from both selected split halves.

Strong approximation \cite[Theorem~7.12]{PlatonovRapinchuk}, applied to
\(\operatorname{SU}(W_0)\) now gives invariance under left
translation by \(\operatorname{SU}(W_0)(\mathbb A_{E^+})\).
Since \(W_0=H\perp W\), its rank \(2j\) is at least two.
The special unitary group is simply connected and absolutely almost
simple, and at every selected split place it is the noncompact group
\(\mathrm{SL}_{2j}\).
Explicitly, for \(s'\in\operatorname{SU}(W_0)(\mathbb A_{E^+})\),
approximate \(s'\) outside the selected places by a rational
point \(s_0\) so that \(g^{-1}s_0^{-1}s'g\) lies in the
original level at the other finite places and tends to
the identity at infinity. At the selected places the full
right factor is allowed. Automorphy and continuity give
the claimed left invariance.

The period therefore equals its normalized radical mean
on the one-plane parabolic of \(W_0\).
For this equal-rank datum the one-plane unfolding
has lower source and target ranks equal, with empty deleted
block. At its source Levi point \((a_0,g_0)\), the
unfolded integrand uses the idele variable \(a\) solely in
the combination \(aa_0\): the source and target Levi
powers combine to \(|aa_0|^{h+1}\), multiplied by the
mixed slice factor \(|aa_0|^{-h/2}\) in lower rank \(h\).

The character and the slice variable are likewise read at
\(aa_0\). Right translation in the absolutely convergent
idele-class integral proves that the mean is independent
of \(a_0\). Hence the period at \(\operatorname{levi}(a_0,1)\)
equals its value \(I_0\) at \(1\). The determinant of the
Levi point is \(a_0/\bar a_0\). Local Hilbert 90, with unit
representatives away from finitely many places, gives every
norm-one idele this way. As two points of equal
determinant differ by the special unitary group, the
equal-rank period is the constant \(I_0\).
This uses the equal-rank one-plane unfolding on the same
finite determinant datum.

Finally let \(z_B\) be the idele with component \(\det B_\tau\)
at infinity and \(1\) at every finite place. With the
normalized complex absolute value,
\[
 |z_B|^j=\prod_\tau|\det B_\tau|^{2j},\qquad
 \nu^{\rm bc}(z_B)
       =\prod_\tau(\det B_\tau/|\det B_\tau|)^{2\mu_\tau}.
\]
Consequently
\[
 \prod_\tau\overline{\det B_\tau}^{2\mu_\tau}
     =\frac{\prod_\tau|\det B_\tau|^{2(j+\mu_\tau)}}
                 {\nu^{\rm bc}(z_B)|z_B|^j}.
\]
Use \eqref{eq:marked-real-slice-small} and
\eqref{eq:marked-real-slice-cancel}, and integrate the
equal-rank factor over its original \([W_0]\)-measure.
The result is the small functional of \(\Lambda[B]\)
multiplied by
\[
              \frac{\rho(h)I_0}{\nu^{\rm bc}(z_B)|z_B|^j}.
\]
This holds simultaneously for every \(B,\beta_0,F,c,g,y_0\),
with the fixed choice of \(p(h)\). The real vectors, bump,
normalizing constant \(\kappa\), and correction \(p(h)\)
are chosen before the character and modulus.

\subsection{Extension of the self-dual Fourier tests}
\label{sec:marked-primitive-fourier-transfer}

\sectionopening{In this subsection, to pass from the real computation of
Section~\ref{sec:marked-real-slice-period} to an adelic identity, the
finite data must match under self-dual Fourier transform. We establish
the required Fourier identity on the full-rank locus and then extend
the period formula.}

Let \(D\) be the \(b\)-plane basis of the large source
\(X\), and \(A\) the additional plane inside its
complement. The \(D\)-basis embeds from the complement of
\(A\). All conjugate-transpose and Fourier coordinates
are those of these bases. Let \(S\) be the nonempty set
of selected split places of \(E\), and let \(P\) be the
finite places outside \(S\cup\bar S\). Use the integral compact level at \(P\), and give every
unipotent rational quotient its probability Haar measure.

In this subsection put \(2j=\dim W_0\),
\(N=\dim X=2b+2j\) and \(a=\dim(B_0\perp W_0)=b+2j\), so \(N=a+b\).
Write \((K_D,I_D,L_D)\) for the ordered Witt basis with \(K_D\)
spanning the given isotropic subspace. In these coordinates the
Gram pairing with \(K_D\) is the row block \(s=L_D^\vee\) of
the inverse basis matrix. Thus \(sK_D=0\), \(sI_D=0\) and
\(sL_D=I_b\). At the selected places put
\[
 \Psi_0(Y)=\chi^{-1}(\det(sY))
       \quad(Y\text{ integral},\ \det(sY)\text{ a unit}),
\]
with zero value elsewhere.

\smallskip\noindent\emph{The common-section space and the equivariant difference.}
Let \(F_S=\prod_{w\in S}E_w\) and let \(O_S\) be its integral
subring. Put \(\lambda(u)=\nu^{\rm bc}(u)|u|^{N/2}\)
on \(F_S^\times\), embedding the selected components as ideles
with value one elsewhere. Indeed the split target change
\(X\mapsto l^{-1}X\) acts on each of the \(N\) columns and
contributes \(|\det l|^{-N/2}\) to the Weil action.
For a full-column-rank matrix \(Y\), define a
section as follows. In each rank write
\(G=\mathrm{GL}(F_S)\) and \(K=\mathrm{GL}(O_S)\).
Choose a fixed set of representatives simultaneously for
\(G/K\) and \(K\backslash G\). Indeed each compact double
coset has equally many left and right cosets by unimodularity,
with each coset of measure one when \(K\) has measure one.

Each left and right coset in \(KgK\) meet, since
\(k_1gk_2\in k_1gK\cap Kgk_2\).
Pair the two finite sets and select an element in each
intersection. Define
\[
         S_Y(\Psi)=\sum_r \lambda(\det r)\Psi(Yr).
\]
On \(X^{\mathsf t}\) use the same construction with
\(\mathrm{GL}_a\). This defines
\[
          S_X^{\mathrm{cov}}(\Phi)=
                 \sum_r\lambda(\det r)\Phi(r^{\mathsf t}X).
\]
Here the tests are compactly supported in the open full-rank
locus. For fixed full-rank \(Y\), choose a left inverse \(A\).
If \(Yr\) is in that compact support, then \(r=A(Yr)\) is
bounded, and every limit of these \(r\)'s is invertible since
the corresponding limit of \(Yr\) has full rank. Thus \(r\)
ranges over a compact subset of \(\mathrm{GL}_b(F_S)\), which
meets only finitely many cosets \(r\mathrm{GL}_b(O_S)\). The covector sum follows by transpose. In the column space \(W\subset E_w^N\) of a
full-rank matrix, the intersection \(W\cap O_w^N\) is
a saturated lattice. Its basis gives the integral
reference for comparing the fractional lattice generated
by the matrix columns.

If \(Y=Y_0t\) with \(Y_0\) a primitive integral basis of this
saturated module at every place, the coset \(Kt\) is independent
of the primitive basis. The set on which \(Kt=Kr\) is given by
primitivity of \(Yr^{-1}\). This is a clopen condition since its entries
are integral and its reduction has full column rank at each place.
The disjoint open sets for \(Kr\in K\backslash G\) cover the
full-rank locus. A compact support therefore meets only finitely
many of them. The same argument applies to \(X^{\mathsf t}\).
The integral determinant covariance is
\(\Psi(Yu)=\lambda(\det u)^{-1}\Psi(Y)\) for
\(u\in \mathrm{GL}_b(O_S)\), and
\(\Phi(vX)=\lambda(\det v)^{-1}\Phi(X)\) for
\(v\in \mathrm{GL}_a(O_S)\). Integral matched functions
have these covariances by changing the complementary
block of an integral completion. Both are preserved by
all split source translations. Each weighted summand is therefore
independent of the choice of representative in \(rK\).

For \(k\in\mathrm{GL}_N(F_S)\), take \(Y\) to be its first
\(b\) columns and \(X\) the last \(a\) covectors of \(k^{-1}\).
Take pairs of Schwartz--Bruhat functions on the full-rank
loci with those integral covariances and satisfying
\[
                S_X^{\rm cov}(\Phi)=\lambda(\det k)S_Y(\Psi)
                \quad\text{for every }k.
\]
For an integral matched pair, choose bases \(a_0X,Yb_0\) of
the saturated covector and column modules. The corresponding covector
and column sections each have only one potentially nonzero
term, up to the integral change of bases. Its weight cancels
the determinant covariance on the integral group. The primitive
map \(a_0X:O_S^N\to O_S^a\) is onto at every factor, hence has an
integral right inverse \(Z\). Its kernel has basis \(Yb_0\), so
\(k_0=(Yb_0,Z)\in\mathrm{GL}_N(O_S)\). In the fixed block order,
\[
 k^{-1}k_0=\begin{pmatrix}b_0&*\\0&a_0^{-1}\end{pmatrix},
 \qquad \det k_0=\det k\,\det b_0/\det a_0 .
\]
The matching law on \(k_0\) gives the displayed section identity. The section space is stable under
\[
 (\Phi,\Psi)\longmapsto
 \bigl(\lambda(\det g)\Phi(\,\cdot\, g),
                    \Psi(g^{-1}\,\cdot\,)\bigr).
\]
In fact its identity at \(k\) is the old one at \(g^{-1}k\),
multiplied by \(\lambda(\det g)\).

\subsubsection{Factorization through the Fourier test}
\label{sec:selfdual-fourier-factorization}

\sectionopening{In this subsection, we use the real restriction of Lemma~\ref{lem:common-real-slice} to compare the large period and its Fourier test.}

The large period \(\mathcal L(\Phi)\) factors linearly through
\(\Psi\) on this space. To prove this, split the compact
support of \(\Phi\) into the finitely many clopen coset pieces
\(\Phi_r\) for which \(X^{\mathsf t}r^{-1}\) is primitive
integral. If \(X^{\mathsf t}=Y_0t\) for a primitive basis
\(Y_0\), this is exactly \(t\in Kr\). Thus \(K\backslash G\)
indexes the partition. At a primitive covector \(X\), the term
\(\Phi_r(r^{\mathsf t}X)\) equals \(\Phi(r^{\mathsf t}X)\),
and the same representatives enumerate \(G/K\) in the section sum.

Translate each piece to primitive integral covectors,
with datum \(\Phi_r(r^{\mathsf t}X)\). By changing the target
variable at the selected split places and using the determinant
multiplier of the Weil action,
\[
 \mathcal L(\Phi_r)
        =\lambda(\det r)\mathcal L\bigl(\Phi_r(r^{\mathsf t}\,\cdot\,)\bigr).
\]
Consequently \(\mathcal L(\Phi)\) is the large period of the
function \(S_X^{\rm cov}(\Phi)\) on primitive integral covectors,
extended by zero. At each such \(X\), choose an integral
completion \(k_X\) with that last block.

The section
identity reads this function directly from \(\Psi\) as
\(\lambda(\det k_X)S_{Y(k_X)}(\Psi)\). Thus it depends
linearly on \(\Psi\). For convergence, apply the ordinary theta
estimate at the end of Section~\ref{sec:marked-mixed-absolute-summability}
with the full rational source basis. A compactly supported covector
test of full row rank satisfies its surjectivity hypothesis at a
selected split place. The transpose of a full-column-rank test
satisfies it at the conjugate place for the small integral. Cutting
into clopen pieces and multiplying the matrix argument by a fixed
invertible matrix on either side preserve compact support and full
rank, also for a nonintegral matrix. For each resulting finite test
and compact set of source points, the absolute value of the whole
theta kernel has a common constant majorant on the finite-volume
target quotient. Hence the changes of target variable use absolutely
convergent integrals and the original Haar measure. Fix \(\nu\) and
take \(c_{\mathrm{top},\nu}\) to be the explicit scalar in
\eqref{eq:period-open-orbit-scalar}. For the small period
\(\mathcal L_{\rm sm}\) on the original column test define
\[
 \mathcal D(\Psi)=\mathcal L(\Phi)
                    -c_{\mathrm{top},\nu}\mathcal L_{\rm sm}(\Psi).
\]
The factorization makes this independent of the choice of matching
covector test \(\Phi\). It is linear on the image of the section
space in the column functions, including when \(c_{\mathrm{top},\nu}=0\).

Let \(d\) be the finite idele equal to \(\det g_S\) at the selected
places and to one elsewhere. On the large side the covectors
transform by \(X\mapsto Xg_S^{-1}\) with section multiplier
\(\lambda(d)^{-1}\). Since \(a+b=N\), the combined factor is
\[
 \bigl((\nu^{\rm bc})^{-1}\upsilon_E^b\bigr)(d)^{-1}
       |d|^{a/2}\lambda(d)^{-1}
       =\bigl(\upsilon_E(d)|d|^{1/2}\bigr)^{-b}.
\]
On the small side the transposed columns at the conjugate
places read \(\bar d^{-1}\). Their rank-\(b\) character is
\(\upsilon_E^b\). As \(\upsilon_E(d\bar d)=1\), the factor is
\[
 \upsilon_E(\bar d^{-1})^{-b}|\bar d^{-1}|^{b/2}
       =\bigl(\upsilon_E(d)|d|^{1/2}\bigr)^{-b}.
\]
Both calculations use the original quotient measures in
\(\mathcal L\) and \(\mathcal L_{\rm sm}\). This gives the stated
covariance of \(\mathcal D\).

The self-dual tests first force only those central coefficients
of \(f=\mathcal D(\Psi_0)\) whose index is primitive at the
selected split places to vanish. The finite Fourier combination
\(B\) will then satisfy \(\mathcal D B=0\), and the translates
of its dilate span the locally constant matched column data.
For rank \(b=0\) there is one column datum up to scalar and
one central coefficient with empty determinant equal to one.
The self-dual test already gives its vanishing on the full
parabolic, which in this case is the whole group. The line
shifts and the Gauss sum are unnecessary in that case.

In these coordinates the covariance is
\[
 \mathcal D(\Psi\circ g_S)(xg_S)
   =(\upsilon_E(\det g_S)|\det g_S|^{1/2})^{-b}
                 \mathcal D(\Psi)(x).
\]
For every such pair the compact-set bound above gives continuity
of both periods by dominated convergence. Rational automorphy follows
from that of the theta kernels. At \(P\) the fixed vacua and the
unramified determinant characters give the indicated right compact
invariance. These properties consequently hold for \(\mathcal D\).

\medskip
\noindent\emph{Interchange of the two line averages.}
For \(T\in M_b(E)\), set \(n_D(M)=1+K_D(M-M^*)s\). With
probability Haar measure on the additive quotient the coefficient is
\[
 F_T(g)=\int_{M_b(E)\backslash M_b(\mathbb A_E)}
    \psi_E(\operatorname{tr}(TM))\, f(n_D(M)g)\,dM .
\]
Since \(sK_D=0\), the map \(M\mapsto n_D(M)\) is additive.
Translation in the integral gives
\[
 F_T(n_D(Z)g)=\psi_E(-\operatorname{tr}(TZ))F_T(g).
\]
For Hermitian \(Z\) we have \(n_D(Z)=1\), so a nonzero
coefficient requires \(T^*=-T\) by nondegeneracy of the trace pairing.
Let \(n_+(a)\) and \(n_-(b')\) be the shifts using the
two opposite isotropic lines of \(A\), where \(a,b'\)
are covectors of length \(b\). Define
\[
 \phi_g^+(a)=\int_{E^b\backslash\mathbb A_E^b}
                 F_T(n_-(b')n_+(a)g)\,db',
\]
and use the opposite order for \(\phi_g^-\). Take the dual
isotropic basis of \(A\) with cross pairing one. In the two
orders of multiplication the upper-right blocks of the products
are \(-a^*b'\) and \(-(b')^*a\), respectively. Hence
\[
 n_+(a)n_-(b')=n_D(-a^*b')\,n_-(b')n_+(a).
\]
Moving \(n_+(a)\) to the left in the opposite mean introduces
\(n_D(a^*b')\). Its central multiplier therefore gives
\[
 \phi^+_{n_+(a)g}(u)=\phi_g^+(u+a),\qquad
 \phi^-_{n_+(a)g}(b')
       =\psi_E(-\operatorname{tr}(T a^*b'))\phi^-_g(b').
\]
For a rational covector \(a_0\), move \(n_+(a_0)\) across the
opposite shift and use rational automorphy. The value
\(\phi_g^+(a+a_0)\) is then the Fourier coefficient at
\(-Ta_0^*\) of \(b'\mapsto F_T(n_-(b')n_+(a)g)\), with
the displayed Fourier sign. If \(\det T\ne0\), these are
all the characters of the covector quotient. For fixed \(a\) its
continuous function \(b'\mapsto F_T(n_-(b')n_+(a)g)\) is therefore
determined everywhere by \(\phi_g^+(a+a_0)\), by Fourier uniqueness.
Commuting the two shifts with the stated central character and
integrating over the other compact covector quotient recovers
\(\phi_g^-\). Reversing the lines recovers \(\phi_g^+\), and these
operations respect scalar multiplication. The line and central
averages integrate continuous functions against probability Haar
on compact sets of representatives.

\medskip
\noindent\emph{Shears on the actual standard datum.}
At a point fixing the \(D\)-basis away from \(P\), the
partial mean at a shifted point, as a linear evaluation
on the equivariant data, is invariant under
\[
                  Y\longmapsto Y+u w sY
                  \qquad(w\in M_{1,b}(E_S)),
\]
where \(u\) is the moved extra-line vector at \(S\).
Indeed these are the split components of the opposite
root shifts used in the mean. They have determinant
one in the split covariance, and right translation
of the covector quotient preserves its Haar measure.
The covector here is still \(s\): it is the adjoint covector of
the fixed basis, and the added extra-line vector is
orthogonal to it. Consequently \(su=0\).
For a general allowed rank-one split shear the change
in the local datum is \(Y\mapsto Y+u rY\), with \(ru=0\).
The other split half has the conjugate unitary coordinates.

Such a shear also fixes the value of this linear evaluation
on \(\Psi_0\). Average the two functions over the restricted
maps \(Y\mapsto Y+u r t sY\) with \(t\) in the compact
additive group \(M_{\dim X,b}(O_S)\). If for some covector
\(\rho\) the point \(Y+u\rho\) lies in \(\operatorname{supp}
\Psi_0\), then \(sY\) is integral invertible and
\(t_0=(Y+u\rho)(sY)^{-1}\) is integral. Since \(ru=su=0\),
\[
 Y+urY+urtsY=Y+ur(t+t_0)sY .
\]
Translation by \(t_0\) preserves Haar on the compact group,
so the two averages agree. If no such \(\rho\) exists the
averages vanish. The averages are finite linear combinations
of translates because the datum is smooth for the compact
group. The restricted invariance therefore proves the needed
invariance under the full split shear with the exact datum.

For either extra line this gives partial-mean invariance under
the adelic line radical, between the points under consideration.
First, its centre fixes the partial mean for rational parameters
by automorphy and for parameters supported over the selected split
pairs by the shear argument. A central parameter at these pairs is
specified by its coordinates at \(S\), with the conjugate coordinates
determined by the skew-Hermitian condition. The rational and selected
parameters together are dense in the full adelic centre by additive
strong approximation.
Averaging over this centre now writes the partial mean as the
coefficient of the combined isotropic basis with index
\(\operatorname{diag}(0,T)\).

The mixed blocks of its centre
give the opposite line shift already integrated in the mean.
A rational parameter in the line radical acts on the combined
basis by an upper block unipotent and preserves this bordered
index by congruence. For the noncentral column parameter \(c\),
take support on \(\bar S\). Write \(i\) for the middle-block
inclusion, \(i^\vee\) for its coordinate projection, and \(J_0\)
for the middle Hermitian matrix in the Witt basis. At \(w\in S\),
\(c_w=0\), so the linear term with image in the middle block
and the quadratic central term vanish. The action on this half is
\[
  1+u_w(c^*)_w J_{0,w}i_w^\vee,
  \qquad i_w^\vee u_w=0 .
\]
It is the shear just treated on the selected datum. The opposite
component is determined by unitarity.
Let \(N_1\) be the extra-line radical, and let \(H_S\) be its
subgroup with column parameter supported on \(\bar S\) and centre
supported over the selected pairs.
For \(h\in H_S\), multiplying by a selected central element
still gives a single shear on the selected datum at the original
point. Rational central translations preserve the partial mean
by automorphy. Density in the centre and continuity then show
that the bordered coefficient at \(hg\) equals the partial mean
at \(g\), for every original point \(g\) under consideration.

The bordered coefficient is invariant under \(N_1(F_0)\) by
the rational congruence above. The product \(N_1(F_0)H_S\) is
dense in \(N_1(\mathbb A_{F_0})\). For the column parameter use
additive strong approximation away from \(\bar S\ne\varnothing\).
The product of a rational column lift and a column lift supported
on \(\bar S\) differs from the lift of their sum by a central
term supported over the selected pairs. The rational and selected
central parameters then approximate the centre as above.
Continuity gives equality of the bordered coefficient at \(ng\)
with the partial mean at \(g\) for every \(n\in N_1(\mathbb A_{F_0})\).
At a point \(ng\) under consideration the partial mean equals its
bordered coefficient, which gives the stated invariance.

\medskip
\subsubsection{Integral extension of the Fourier test}
\label{sec:selfdual-fourier-integral-extension}

\sectionopening{In this subsection, we extend the test of Section~\ref{sec:selfdual-fourier-factorization} in the local self-dual coordinates.}

Suppose \(F_T(g)\ne0\) at a point fixing the \(D\)-basis
away from \(P\), with \(\det T\ne0\). The line interchange
gives a nonzero shifted mean. At the places in \(P\), act
on the right by the integral group so that the
action on the combined isotropic basis is upper block
triangular for the flag given by the extra line and the
\(D\)-plane. This is the integral flag decomposition:
intersect the two isotropic spaces of the transported flag
with the standard lattice, take primitive bases in succession,
and complete to an integral Witt basis.

For a primitive isotropic basis matrix \(X\), unimodularity
and primitivity give integral dual columns \(Y\) with
\(X^*JY=1\). Choose integral \(H\) with
\(H+H^*=Y^*JY\), using the upper triangular entries off the
diagonal and surjectivity of the integral trace on the diagonal.
Every nonsplit finite place of \(E/F_0\) is unramified by
Section~\ref{sec:different-additive-character-scales}.
For an unramified extension this surjectivity follows from the
surjective residue trace: a lift with unit trace can be divided by
that trace. The split trace is the sum of the two components.
Now \(Y-XH\) is isotropic and still dual to \(X\), so the
hyperbolic span splits off with unimodular complement.
This construction includes residue characteristic two. The compact action is trivial at the
other places. The complementary factor of the upper-triangular
action is in the extra-line radical, which may be discarded
by the invariance just proved. This leaves a point which at
\(P\) agrees with a block action on the chosen bases with
scalars \(t,t'\) on the two lines and matrix \(\gamma\) on
the \(D\)-basis, where \((t')^*t=1\). Call it a straight
point. It still fixes the \(D\)-basis off \(P\).

Put \(\delta=j_{\rm W}\alpha_0\bar\alpha_0\), the additive scale, and set \(C=\delta^{-1}\gamma^*T\gamma\) at
\(P\). At every straight point where the shifted mean is
nonzero, \(C\) is integral at \(P\).
Indeed some shifted value of the coefficient is nonzero at a
point with the same \(D\)-basis action. The compact level
allows all integral central translations at \(P\) there.
At \(w\in P\), take the integral matrix parameter
\(M=aE_{ji}\), \(a\in\mathcal O_w\), supported at that place.
Write \(z(M)\) for the central translation with skew block
\(M-M^*\). In the defining coefficient the matrix parameter
\(U\) has weight \(\psi_E(\operatorname{tr}(TU))\).
This weight factors through \(U-U^*\): for \(H=H^*\),
\(\operatorname{Tr}_{E/E^+}\operatorname{tr}(TH)=0\).
At the coefficient point \(g\), we have \(gK_D=K_D\gamma\)
on \(P\); unitarity gives \(sg^{-1}=\gamma^*s\).
Thus \(g n_D(M)g^{-1}=n_D(\gamma M\gamma^*)\) for \(M\)
supported there. Changing the integration variable to
\(U+\gamma M\gamma^*\) gives
\[
 F_T(g n_D(M))=\psi_{E,w}(-\delta_w aC_{ij,w})F_T(g).
\]
Here \(\operatorname{tr}(MC)=a C_{ij,w}\), so the multiplier
is \(\psi_{\rm sd}(-a C_{ij,w})\),
where \(\psi_{\rm sd}(u)=\psi_E(\delta u)\) on these finite
coordinates. For the absolute trace character the trace dual of
\(\mathcal O_w\) is \(\mathfrak D_{E_w/\mathbb Q_p}^{-1}\),
and \(\delta_w\mathcal O_w=\mathfrak D_{E_w/\mathbb Q_p}^{-1}\)
by Section~\ref{sec:period-finite-coordinates}. Thus the annihilator
of \(\mathcal O_w\) for this scaled pairing is \(\mathcal O_w\).
Since the compact translation fixes the nonzero coefficient,
\(C_{ij,w}\) is integral. The central parameter introduces no
factor of two: for any integral skew block \(Z\) one can take
\(M=c_{\rm tr} Z\), since \(c_{\rm tr}+\bar c_{\rm tr}=1\)
and \(M-M^*=Z\). This works in residue characteristic two as well. Choose a straight point with
nonzero shifted mean minimizing the positive integer
\[
                 I(C)=\prod_{w\in P}q_w^{v_w(\det C)}.
\]
The factors are one almost everywhere. For a fixed straight \(g\),
the two functions \(\phi_g^+\) and \(\phi_g^-\) are either both zero
or both nonzero, by Fourier uniqueness in the line interchange.
Straight points preserve the two opposite lines at \(P\), so exchanging
them preserves this set of points and the \(D\)-block \(\gamma\).
Thus a single minimizing point works for both functions.

For a covector \(a\) with \(\phi_g^+(a)\ne0\), one has
\begin{equation}
              \gamma^{-1}a^*t\in O_P^b.                    \label{eq:primitive-transfer-lattice-plus}
\end{equation}
At \(w\in P\) the bottom block row of the shifted basis matrix is
\((-a^*t,\gamma)_w\). Its image on \(\mathcal O_w^{b+1}\) is
\[
 \Lambda'_w=\gamma_w\mathcal O_w^b+(a^*t)_w\mathcal O_w .
\]
Choose a basis \(\gamma'_w\) of \(\Lambda'_w\). The surjection onto
this free lattice splits. A basis of its rank-one kernel, followed by
lifts of the columns of \(\gamma'_w\), gives an integral invertible
column change that turns the bottom row into \((0,\gamma'_w)\).
At almost every place \(a,t,\gamma,\gamma^{-1}\) are integral,
so \(\Lambda'_w=\mathcal O_w^b\). All the column changes are
integral and together give an element of the product compact.
Extend the resulting integral change on the combined isotropic plane
by its inverse adjoint on the opposite plane and by the identity on the
complement. This is a right compact move. The remaining upper
off-diagonal block lies in the extra-line radical, under which the
shifted mean is invariant. For the inverse column change \(M\) we have
\(-a^*t=\gamma'M_{21}\) and \(\gamma=\gamma'M_{22}\).
Thus \(\gamma_w\mathcal O_w^b\subseteq\Lambda'_w\), and the inclusion
is strict at some \(w\) if \eqref{eq:primitive-transfer-lattice-plus} fails.
The new Gram is integral since the new mean is nonzero,
and
\[
 C=M_{22}^*\,(\delta^{-1}(\gamma')^*T\gamma')\,M_{22}.
\]
Write \(C'=\delta^{-1}(\gamma')^*T\gamma'\).
Since \(P\) is stable under conjugation and \(q_w=q_{\bar w}\),
\[
 I(C)=I(C')\prod_{w\in P}q_w^{2v_w(\det M_{22})}>I(C').
\]
All the valuations on the right are nonnegative and at least
one is positive, contradicting minimality.
The opposite argument gives
\(\gamma^{-1}(b')^*t'\in O_P^b\) whenever
\(\phi_g^-(b')\ne0\).

Take any column \(x\) with \(Cx\in O_P^b\), and set
\(a=(\gamma x(t')^*)^*\) on \(P\) and zero at the
other places. On the support of \(\phi_g^-\) the
multiplier of the shift \(n_+(a)\) is one. Indeed the
scaled argument of its additive character is
\(-(\gamma^{-1}(b')^*t')^*Cx\), which is integral.
Off \(P\) the argument is zero. Thus \(a\) is a period of
\(\phi_g^+\) by the interchange. Apply
\eqref{eq:primitive-transfer-lattice-plus} to a nonzero
point of \(\phi_g^+\) and its translate by \(a\), and
subtract. Since \((t')^*t=1\), this gives \(x\in O_P^b\).
Use the columns of \(C^{-1}\) as \(x\). The minimizing
Gram is therefore unimodular at every place in \(P\).

A nonzero shifted mean gives a nonzero coefficient at a
point with this same \(D\)-basis action on \(P\) and identity
on that basis elsewhere, because both line shifts fix \(D\).
Suppose now that \(\delta^{-1}T\) is primitive at \(S\).
Since \(T^*=-T\) and \(\bar\delta=-\delta\), the matrix
\(D_0=\delta^{-1}T\) and its inverse are Hermitian. At a split
pair, conjugate transposition therefore transports integrality
of both matrices from \(w\in S\) to \(\bar w\).
Write \(T=j_{\rm W}\beta\) and \(q=\alpha_0 q_{\rm Lev}\)
for the parabolic Levi matrix of the point. Its finite action on
the \(D\)-basis is \(\gamma\) at \(P\) and the identity on the
selected split pairs. Hence
\[
 Q=(q_{\rm Lev}^*\beta q_{\rm Lev})_f
   =(\delta^{-1}q^*Tq)_f,\qquad
 Q|_P=C,\quad Q|_{S\cup\bar S}=D_0|_{S\cup\bar S}.
\]
This Hermitian matrix and its inverse are integral everywhere
finite. The integral unitary basis construction in
Section~\ref{sec:marked-mixed-absolute-summability} gives
\(k\in\mathrm{GL}_b(\widehat{\mathcal O}_E)\) with \(k^*Qk=I_b\). If \(\beta\) is outside the
real congruence class of \(B_0=U^-\), its coefficient
vanishes. For the chosen real slice this follows from its
support and the empty small orbit in
Section~\ref{sec:marked-mixed-absolute-summability}.
For the arbitrary real vectors in
Section~\ref{sec:marked-equal-rank-vanishing}, the compact
averaging argument there gives this vanishing.

Otherwise Landherr's theorem \cite[pp.~246--247]{LandherrHermitian}
gives \(d\in\mathrm{GL}_b(E)\) with \(\beta=d^*J_{B_0}d\):
the forms have the same real signatures and are isometric at
each finite place. Let \(m_D(d)\) be the rational Levi element
acting by \(d\) on \(K_D\). In the coefficient integral put
\(M=dM'd^*\). Rational automorphy and cyclicity of trace give
\[
 n_D(M)m_D(d)=m_D(d)n_D(M'),\qquad
 F_{j_{\rm W}J_{B_0}}(m_D(d)g)=F_T(g).
\]
The change is an automorphism of the rational additive quotient,
so it preserves probability Haar. On \(M_b(E)\) its
\(E\)-linear determinant is \((\det d\,\overline{\det d})^b\),
whose absolute values over the places of \(E\) have product one.
The moved Levi matrix is \(dq=\alpha_0 X_{B_0}Y\) with
\(Y=X_{B_0}^{-1}d q_{\rm Lev}\). From
\((X_{B_0}^*J_{B_0}X_{B_0})_f=I_b\) we get
\[
 (Y_f k)^*(Y_f k)=k^*Qk=I_b .
\]
This is precisely the finite basis condition of the self-dual
test at index \(j_{\rm W}J_{B_0}\), so the coefficient vanishes. The
contradiction proves that \(F_T(1,g_f)=0\) when
\(\delta^{-1}T\) is primitive at \(S\) and \(g_f\) is
trivial there and at the conjugate split places.

\medskip
\noindent\emph{Finite detection of the other coefficients.}
At \(S\) use the opposite ordering \((K,L)=(L_D,K_D)\)
in this calculation. For the corresponding coordinate rows of
the inverse basis matrix we have \(s=K^\vee=L_D^\vee\),
\(K^\vee L=L^\vee K=0\) and \(K^\vee K=L^\vee L=I_b\).
In particular \(LzK^\vee=K_DzL_D^\vee\) is the original central
block on the chosen split half.
For \(z\in M_b(E_S)\) put
\[
           \Psi_z(Y)=\Psi_0(Y+LzK^\vee Y).
\]
Extend \(z\) by zero outside \(S\). The component of \(n_D(z)\)
on the chosen half is \(1+LzK^\vee\), and its determinant is one.
Split covariance gives
\(\mathcal D(\Psi_z)(x)=\mathcal D(\Psi_0)(x n_D(-z))\).
At the indicated points \(g=(1,g_f)\) with \(g_f=1\) on
\(S\cup\bar S\), this inverse shear commutes with \(g\).
Apply the central translation formula with \(Z=-z\), and put
\(\psi_S(u)=\prod_{w\in S}\psi_{E,w}(\delta_w u_w)\). Then
\[
 \Lambda_{T,g}(\Psi_z)
   =\psi_S(\operatorname{tr}((\delta^{-1}T)_S z))
                  \Lambda_{T,g}(\Psi_0),
\]
where \(\Lambda_{T,g}\) is the central coefficient of
\(\mathcal D\) at \(g\) as a linear form of the datum. Since
\(\delta_w\mathcal O_w=\mathfrak D_{E,w}^{-1}\), the conductor
of \(u\mapsto\psi_{E,w}(\delta_wu)\) is \(\mathcal O_w\).
Thus \(r\mapsto\psi_{E,w}(\delta_w\varpi_w^{-1}\widetilde r)\)
is a well-defined nontrivial additive character of each \(k_w\).
Let \(k_S\) be the product of residue fields, and let
\(\widetilde r\) be integral lifts of their matrices
with \(\widetilde1=1\). Put
\[
 B(Y)=\sum_{r\in M_b(k_S)}
          \chi^{-1}(\det r)\Psi_{\varpi^{-1}\widetilde r}(Y).
\]
The character is extended by zero on nonunits. This detector
uses the residue field even if the level at \(w\) is
\(\mathfrak p_w^2\). Indeed, on the deeper ring a determinant
character inflated from \(k_w\) is constant on
\(r_0+\varpi_w M_b(k_w)\). Summation over this fiber
against a primitive additive character of the deeper ring gives
zero unless the Fourier index lies in \(\varpi_w M_b(\mathcal
O_w/\mathfrak p_w^2)\), by additive orthogonality. The shifts
\(\varpi^{-1}\widetilde r\) in \(B\) pair directly with
indices reduced modulo \(\varpi\). The full depth of a matched
datum enters the subsequent compact average through \(c\).

If \(\Lambda_{T,g}(\Psi_0)\ne0\), invariance under
integral shears forces \((\delta^{-1}T)_S\) integral.
Write \(A=(\delta^{-1}T)_S\). Since \(\psi_S\) is trivial
on \(O_S\), \(\psi_S(\operatorname{tr}(A\varpi^{-1}\widetilde r))\)
depends only on \(r\). By the primitive test, \(\bar A_w\)
is singular for some \(w\in S\). The tame residue character
\(\chi_w\) is nontrivial at each selected place. Choose a unit
\(a\) with \(\chi_w(a)\ne1\), a nonzero \(v\in\ker\bar A_w\)
and a residue-field linear form \(\ell\) with \(\ell(v)=1\).
Take \(u_w=1+(a-1)v\ell\) and \(u_{w'}=1\) for \(w'\ne w\).
Then \(\bar A u=\bar A\) and \(\det u_w=a\).
Reindex the Gauss sum by \(r\mapsto ur\).
The additive term is unchanged and the determinant character
changes by a nontrivial scalar. Thus the sum is zero.
This also handles a singular rational \(T\), and all other
terms with \(\Lambda_{T,g}(\Psi_0)=0\) are already zero.

The skew coefficients exhaust the characters of the central
quotient in the original convention. Fourier uniqueness
gives \(\mathcal D B(g)=0\) at these points.
Their rational translates are dense by weak approximation for
the unitary group at the specified real and split finite places,
with all other finite components unrestricted.
Weak approximation here follows from its rational Cayley chart
and weak approximation for the vector space of skew matrices.
Continuity and automorphy give \(\mathcal D B=0\).

Translate by \(Y\mapsto-\varpi^{-1}Y\) on the split source.
Put \(p=L^\vee\) and \(U=-\varpi^{-1}K^\vee Y\). The
argument of \(\Psi_0\) in the \(r\)-summand is
\[
 X_r=\varpi^{-1}(L\widetilde r U-Y),\qquad K^\vee X_r=U .
\]
A nonzero term requires \(U\in\mathrm{GL}_b(O_S)\),
\(r\in\mathrm{GL}_b(k_S)\), and \(X_r\) integral. It follows
that \(Y\) is integral and
\(\overline{pY}=r\bar U\), so at most one \(r\) contributes.
Conversely, write \(Y=LV+\varpi Z\) with \(pZ=0\).
If \(V=pY\) and \(K^\vee Z\) are integral invertible matrices
and \(Z\) is integral, then \(U=-K^\vee Z\) is invertible,
\(r=\bar V\bar U^{-1}\) gives \(X_r\) integral, and its
weighted value is \(\chi^{-1}(\det V)\). Thus, writing
\(\Omega\) for this set of \(Y\), we have the exact formula
\[
 G(Y):=B(-\varpi^{-1}Y)
     =\mathbf1_\Omega(Y)\chi^{-1}(\det(pY)).
\]
On integral matrices membership in \(\Omega\) depends only on
\(Y\bmod\varpi^2\): first test \((1-Lp)Y\equiv0\bmod\varpi\),
then test invertibility of \(pY\) and of
\(K^\vee(1-Lp)Y/\varpi\) after reduction. At
\(x_0=L-\varpi K\) the term \(r=1\) has argument \(K\), so
\(G(x_0)=1\). Apply the compact averaging identity
\eqref{eq:period-original-compact-average} with \(F\) multiplied
by the indicator of \(\Omega\), and with this \(x_0,p\).
Its positive scalar gives all the integral matched data from \(G\).

\section{The normalized period exchange}
\label{sec:normalized-period-exchange}

\sectionopening{In this section, we compare the periods using the local factors of
Section~\ref{sec:local-period-transforms}. The real identity and its
extension to the finite tests supply the comparison used in
Section~\ref{sec:period-gluing-operator}.}

\subsection{The normalized comparison of the two periods}
\label{sec:marked-partner-period}

\sectionopening{In this subsection, we compare the paired periods of
Section~\ref{sec:marked-square-analytic-pairing} by removing
hyperbolic planes. The real factor from
Section~\ref{sec:marked-real-slice-period} and the finite Fourier
transform of Section~\ref{sec:marked-primitive-fourier-transfer} fix
the final scalar.}

Let \(W,V,U\) have ranks \(N,a,b\), their original standard
unimodular models at the finite places, and respective negative
indices \(q(W),q(V),q(U)\) independent of the infinite place.
Assume \(N=a+b\), \(a-b=2t\) for \(t\in\mathbb Z\),
\(q(W)=q(V)+q(U)\), and
\[
 2q(W)+2\le N\le2\mu_\tau+2q(W)\qquad(\tau\mid\infty).
\]
Fix the ordered splitting into \(b,a\) coordinates and a real
Schwartz vector \(y\) for \((U,W)\), equivalently a Fock vector with
all weighted absolute moments by Section~\ref{sec:marked-real-slice-period}.
For the admissible split tame character \(\nu\), put
\[
 \gamma=\nu^{-1}\upsilon^{N/2},\qquad
 M_\nu(j)=\prod_{w\in S}q_w^{\,j^2}
                 \prod_{k=1}^{j}|L^S(k,\nu)|^2
 \quad(j\ge0),\qquad M_\nu(-j)=M_\nu(j)^{-1}.
\]
The positive-index factors are nonzero by the argument
preceding \eqref{eq:plane-removal-ratio}, so the reciprocal
defines the negative-index values. Here \(\upsilon_E\)
is the fixed splitting Hecke character and \(\upsilon\)
is its restriction to the norm-one ideles.
Let \(S\ne\varnothing\) consist of one chosen member of each split
pair in the modulus. Write \(\nu^{\rm bc}(z)=\nu(z/\bar z)\).
Admissibility requires \(\nu^{\rm bc}\) to be unramified at every
finite place outside \(S\cup\bar S\). At each \(w\in S\) its
restriction to \(\mathcal O_w^\times\) is the inflation of a
nontrivial character of \(k_w^\times\). Its restriction at
\(\bar w\) is determined by conjugation and inversion. The product
\(L^S\) omits exactly these pairs from the base-change \(L\)-series. The characters have weight
\(\mu\) at infinity, which is fixed. The Haar measures and
vacua throughout are the original ones for the models.

We give the assembly from the Fourier and constant-term
identities, specifying the versions needed to normalize the
comparison. For the top calculation use
\[
 W_0=H\perp W,\qquad X=H^b\perp W_0,\qquad
 V_{\mathrm{top}}=U^-\perp W_0 .
\]
The small side is \((X,U^-)\), with \(U^-\) the opposite
space. Write \(j_0=\dim W_0/2\). We use the real slice
of Section~\ref{sec:marked-real-slice-period} with \(j=j_0\).
Given the two fixed real vectors, one for the small side and
one \(z_0\) for \((W_0,W_0)\), the slice identity on the
original self-dual finite lattices has coefficient
\[
  \frac{\rho(h)I_{W_0,\nu}(z_0)}
                 {\nu^{\rm bc}(z_B)|z_B|^{j_0}} .
\]
Here \(B\) is the invertible small slice matrix, \(z_B\) is
\(\det B\) at infinity and \(1\) at finite places, and
\(I_{W_0,\nu}(z_0)\) is the equal-rank \([W_0]\)-period
for the finite determinant datum \(\chi^{-1}(\det)\).

The real factor and the real coordinate corrections are chosen
before \(\nu,B\) and the finite data. The self-dual Fourier
unfolding and the primitive-coefficient vanishing and spanning
on its test space give the identity on standard partners
with scalar
\[
 c_{\mathrm{top},\nu}
    =R_0\kappa_0\varepsilon\,\nu^{\rm bc}(d_0)
                I_{W_0,\nu}(z_0),\qquad
 \kappa_0=(v_{\mathrm{big}}/v_{\mathrm{sm}})|d_0|^{j_0}.
 \tag{1}\label{eq:period-open-orbit-scalar}
\]
The nonzero vacuum constants \(v_{\mathrm{big}},v_{\mathrm{sm}}\)
and the coordinate-change idele
\(d_0=\det(\alpha_0 X_{U^-})\) belong to the original splittings, with \(X_{U^-}\) the small adelic basis.

At the selected places let \(B_D=(L_D,K_D,I_D)\) be the
integral basis matrix in this order, with \(s=L_D^\vee\).
For the stacked covectors \(C_s(Z)=\binom{s}{Z}\),
\[
 C_s(Z)B_D=\begin{pmatrix}I_b&0\\ ZL_D&Z(K_D,I_D)\end{pmatrix}.
\]
Consequently, on the integral unit support,
\[
 \chi^{-1}(\det C_s(Z))
  =\varepsilon\,\chi^{-1}\bigl(\det[Z(K_D,I_D)]\bigr),
 \qquad \varepsilon=\chi(\det B_D),\qquad |\varepsilon|=1 .
\]
The column partner is \(\chi^{-1}(\det(sY))\) on the support
specified in Section~\ref{sec:marked-primitive-fourier-transfer}.
The positive integral \(R_0\) is the integral of
\(\rho(h)\mathbf 1_{H_fK_f^{\rm std}}(u)\) on the quotient
by the block stabilizer. In particular the compact set in this
integral is the full fixed standard level. The exponent here
is \(j_0\) because the complement in the mixed model is \(W_0\).

\subsubsection{Absolute summability of the mixed terms}\label{sec:marked-mixed-absolute-summability}

\sectionopening{In this subsection, we justify the unfolding needed for Proposition~\ref{prop:normalized-period-exchange}.}

In this mixed model write \(K\) for the isotropic columns spanning
\(D\), denoted by \(L\) in Section~\ref{sec:selfdual-fourier-integral-extension},
and \(K'\) for the opposite columns. With \(n_0=\dim W_0\) and \(I\)
the \(n_0\) middle columns, \(A=(K,I)\) is a rational basis of
\(K^\perp\). Let \(A_f^{\rm rat}\) give the finite coordinates of these
rational columns in the standard model. The adapted frame
\(A_f=(K_f,I_f)\) uses the chosen adelic basis of \(W_0\) in its middle
block, so \(A_f=A_f^{\rm rat}\operatorname{diag}(1_b,x_{W_0,f})\)
with \(x_{W_0,f}\) invertible. For a source point \(g\) in the
parabolic fixing the span of \(K\), one has
\[
 A_g=g_f^{-1}A_f^{\rm rat},\qquad A_g C_g=A_f .
\]
Indeed the columns of \(g_f A_f\) and \(A_f^{\rm rat}\) are bases of
the same space with the same initial isotropic subspace. Their
change of basis \(C_g\) is invertible and upper triangular by blocks.

For a target of rank \(d\), the condition used below is that
\(Y A_g\) have a right inverse over the local field for every matrix
\(Y\) in the support at one fixed split place. For the large test,
\(d=b+n_0\), and the integral basis matrix at a chosen \(w\in S\) is
\(P=(K'_f,K_f,I_f)_w\). The covector \(s_w\) from
Section~\ref{sec:marked-primitive-fourier-transfer} is the first row
block of \(P^{-1}\). Thus
\[
 \binom{s_w}{Y}P=
 \begin{pmatrix}1_b&0\\ YK'_{f,w}&YA_{f,w}\end{pmatrix}.
\]
On the support of the determinant test both square basis changes are
integrally invertible. Hence \(YA_{f,w}\) is invertible, and
\(C_{g,w}(YA_{f,w})^{-1}\) is a right inverse to \(YA_{g,w}\).

For the small test, \(d=b\). Read the transpose of its column
variable at \(\bar w\), identifying \(E_{\bar w}\) with \(E_w\) by
conjugation. In the standard finite Hermitian form,
\(s_w^{\mathsf t}\) is then \(K_{f,\bar w}\), since \(s\) is the
Hermitian dual of \(K\). Thus the character is supported where
\(\det(s_wY^{\mathsf t})=\det(YK_{f,\bar w})\) is a unit. With
\(J_K=\binom{1_b}{0}\) the required right inverse is
\[
       C_{g,\bar w}J_K(YK_{f,\bar w})^{-1}.
\]
If the target has rank zero its theta sum has one term. These
calculations use a single selected split pair. The vacuum at all
other places enters through its integral lattice and open stabilizer.
For a change of finite coordinates both the support and the frame are
transported in the product \(YA_g\). Full row rank of \(Y\) alone
suffices for the ordinary period with a full source basis. In a
mixed slice it is the restriction of \(Y\) to the moved
\(K^\perp\) that must be surjective.

For a fixed source point let \(q_x\) denote the theta functional
of the mixed slice at the rational matrix \(x\). Let \(h\) be the
target point and \(H=H(h)\) its reciprocal lattice minimum cut off
below at one. There are \(A_0,M\) and \(C_g\) such that, for an
arbitrary real Schwartz vector \(u\) and the fixed finite datum
\(f_f\),
\[
 \sum_x |q_x(\omega(g,h)(u\otimes f_f))|
       \le C_g H^{A_0}p_M(u).
\]
Here \(p_M\) is a fixed finite-order Schwartz seminorm.

Indeed the
absolute sum is invariant under rational target translation,
which permutes the rational matrices \(x\). Write \(\Lambda_h\) for the
rational vectors integral in the moved finite target lattice and
\(\|v\|_h=\max_\tau\|v\|_{h,\tau}\) in its signed orthonormal real bases.
The finite Hermitian form is integral on \(\Lambda_h\), since the original
finite lattice is integral and \(h_f\) is unitary. Thus vectors of length
below one pair to zero: every conjugate of their integral pairing has
modulus below one. At each infinite place double their complex span,
halve its image under the sign matrix, and fix the joint Euclidean
orthogonal complement, as in Section~\ref{sec:paired-lifts-absolute}.
This simultaneous unitary move sends the first minimum \(\lambda\) to
at least \(\min(2\lambda,1/2)\). After at most \(\lfloor\log_2 H\rfloor\)
moves the height is at most two. The compact reduction of
Section~\ref{sec:bottom-descent-height} therefore gives \(\gamma h=z a\)
with \(\gamma\) rational, \(z\) in a fixed adelic compact set and
\(a\) purely real, with \(\|a_\tau^{\pm1}\|\le H\) at every \(\tau\).
In any rational
polarization of the smaller mixed block, the absolute sum of the
slice theta evaluations is bounded by the sum of the absolute
Schwartz values over the rational coordinates of both blocks.

The finite coordinates in the support range over a compact set.
In a fixed rational basis it is contained in
\(m^{-1}\widehat{\mathcal O}_E^{\,r}\) for some positive integer \(m\),
so the rational coordinates in the sum lie in \(m^{-1}\mathcal O_E^{\,r}\).
Summing \((1+\|x\|+\|y\|)^{-M}\) for sufficiently large \(M\)
gives the claimed bound before the real move. The position and
differentiation operators transform by the real symplectic
matrix under Weil conjugation. Their finite products, Sobolev's
inequality and unitarity give the stated polynomial factor for
the real move, increasing the fixed order \(M\) if needed. The
finite translates in the compact reduction are taken from a
finite set by smoothness.

For \(H>2\) take a shortest nonzero \(d\in\Lambda_h\). The rational
lattice is discrete in the product of its archimedean completions,
so the minimum is attained and \(\|d\|_h=H^{-1}\). Its Hermitian norm
\(\langle d,d\rangle\) belongs to \(\mathcal O_{E^+}\) and
\(\lvert\tau\langle d,d\rangle\rvert\le \|d\|_{h,\tau}^2\le H^{-2}<1\)
for every \(\tau\). Hence \(d\) is isotropic. Multiply the columns \(c_j\) of the fixed rational basis
of \(K^\perp\) by a single positive integer to put them in the
moved source lattice.

Their lengths are bounded in terms of
\(g\). Let \(Z\) be the rational isotropic span of the maps
\[
                         d\,c_j^*J_X .
\]
These maps have integral finite coordinates and real lengths
\(O_g(H^{-1})\) in the moved Hom space. Use the opposite
frame \(K'\), and write \(J_{\rm tar}\) for the rational target
form and \(\Theta_0\) for the theta functional
of the remaining block. For \(y\in Z(E)\) one has \(yK=0\).
In rational coordinates the mixed Heisenberg formula gives
\[
 q_x(\rho(y)f)=
 \psi_E\!\left(j_{\rm W}\operatorname{tr}((yK')^*J_{\rm tar}x)\right)
 \Theta_0\bigl(\rho_0(yI) f[x]\bigr)=q_x(f).
\]
Indeed the argument of the additive character lies in \(E\),
and the rational translation \(yI\) preserves \(\Theta_0\).
As \(d^*J_{\rm tar}d=0\), the Heisenberg cocycle on
\(Z(\mathbb A_E)\) is trivial. Consequently, for any finite
sum \(L=\sum_x\epsilon_xq_x\), \(|\epsilon_x|=1\), the formula
\[
 L_0(f)=\int_{Z(E)\backslash Z(\mathbb A_E)}
                    L(\rho(z)f)\,dz
\]
defines a \(Z(\mathbb A_E)\)-invariant functional using the
probability Haar measure.

The mean annihilates the translated datum. In the moved finite
coordinates, the covector \(d_h^*Y A_g\) is nonzero at a selected
split place whenever the datum at \(Y\) is nonzero. For a right inverse \(R\) of \(YA_g\) we have
\[
                  (d_h^*YA_g)R=d_h^*\ne0 .
\]
In the opposite split component the corresponding
Heisenberg shift multiplies the Schr\"odinger datum by the
additive character paired with this covector. For every support
point one can choose a parameter for which this character is
different from one. It is constant on a small open neighborhood.
A finite partition of the compact support into such neighborhoods,
and invariance of \(L_0\) under the corresponding shifts, show
that \(L_0\) vanishes on each piece, with any real part.

For clarity the estimate for the difference from this mean is
as follows. Choose a positive integer \(m_0\) such that all
finite Heisenberg shifts in \(m_0\widehat{\mathcal O}_E\) in
the chosen integral coordinates fix the fixed finite datum.
For each coefficient of the displayed generators, reduce
modulo \(E\) to an adele with finite part in
\(m_0\widehat{\mathcal O}_E\) and bounded real part. This uses
additive strong approximation followed by a fundamental
parallelepiped of \(m_0\mathcal O_E\).

These measurable representatives fix the finite test after conjugation
and give real displacements \(s\) of size \(O_g(H^{-1})\).
The finite rank condition gives the exact vanishing of the mean.
The displacement bound uses the shortest vector at the real places
and the bounded real representatives for the scalar coefficients.
Let \(A_h\) be their probability average on the real test and set
\(L_{g,h}(u)=L(\omega(g,h)(u\otimes f_f))\). Weil covariance gives
\[
 L_{g,h}(A_hu)=L_0(\omega(g,h)(u\otimes f_f))=0 .
\]
This holds for every real Schwartz vector, so
\(L_{g,h}(u)=L_{g,h}((1-A_h)^n u)\) for every \(n\ge0\).

Use real self-dual coordinates and the seminorms
\[
 p_j(u)=\max_{|\alpha|+|\beta|\le j}
             \sup_t |t^\alpha\partial^\beta u(t)|.
\]
For \(s=(a,b)\), the generator of the Weyl displacement
\(\rho(\lambda s)\) is \(D_s=a\cdot\nabla+2\pi i\,b\cdot t\)
in coordinates with translation by \(+a\) and modulation by \(+b\).
Leibniz's rule gives \(p_j(D_su)\le C_j\|s\|p_{j+1}(u)\).
Translation of the polynomial weights and differentiation of the
modulation bound \(p_j(\rho(\lambda s)u)\le C_{g,j}p_j(u)\)
for \(0\le\lambda\le1\) on this bounded set of displacements. Integrating
\(\partial_\lambda\rho(\lambda s)u=\rho(\lambda s)D_su\), then
averaging and iterating, gives
\[
 p_j((1-A_h)u)\le C_{g,j}H^{-1}p_{j+1}(u),\qquad
 p_M((1-A_h)^nu)\le C_{g,M,n}H^{-n}p_{M+n}(u).
\]
The initial estimate for \(L_{g,h}\) is independent of the finite set
of \(x\)'s and their phases.

Taking phases so that the finite
sum is the sum of the absolute values, then increasing the
finite set, gives
\[
       \sum_x |q_x(\omega(g,h)(u\otimes f_f))|
                 \le C_{g,n}H^{A_0-n}p_{M+n}(u).
\]
At \(H\le2\) use the first bound. In particular this is bounded
on the target quotient, which has finite volume. The radical
quotient is compact. Over a compact
set of representatives its source points stay in the
parabolic, so the same rank proof applies.

The constants
can be taken uniformly there. The bounded representative
set has compact finite projection, so a single integer clears
the coordinates of the fixed rational basis throughout, and
its real basis coordinates are uniformly bounded. The initial
Schwartz bound is uniform on that compact set as well.
The central parameter itself only multiplies each mixed term
by a unit phase. Thus the absolute sum is integrable also
in the radical variables.

Tonelli and Fubini now justify
both the extraction of the central coefficient term by term
and the subsequent unfolding of the surviving rational orbit.
For the ordinary period, take a full rational source basis
instead of the basis of the isotropic perpendicular and use the
whole theta functional in place of \(q_x\). Suppose that, at one
split place and in the orientation for the integrated group,
every matrix in the support of \(f_f\) has full row rank over
the local field. Multiplication by the matrix of the moved full
source basis preserves this rank at every source point. The mean
then vanishes by the same finite modulation calculation. The
initial bound, the denominator clearing for the fixed source
basis, and the real bounds for its columns are uniform for \(g\)
in a compact subset \(C\) of the source group. Taking sufficiently
many differences gives a Schwartz seminorm \(p_{N_B}\) with
\begin{equation}\label{eq:partner-exchange-ordinary-bound}
  \sup_{g\in C}\big|\theta(\omega(g,h)(u\otimes f_f))\big|
        \le C_{C,f_f,B}\,H(h)^{-B} p_{N_B}(u)
              \qquad (B\ge0).
\end{equation}
For a rank-zero integrated group the quotient consists of one
point and continuity follows from the ordinary lattice sum.
For positive rank the quotient has finite volume by
\cite[Theorem~5.6(i), p.~21]{BorelAdele}. A rational character of the
unitary group factors through its determinant torus, and inversion
by the nontrivial element of \(\operatorname{Gal}(E/E^+)\) on the
character lattice of that norm-one torus leaves only the zero
character invariant. Thus \(B=0\) gives an integrable bound uniform
on \(C\), and dominated convergence gives continuity of the
ordinary period. By the coherent-state formula of
Section~\ref{sec:four-unitary-constant-terms}, the same bound
can be written with a sufficiently high absolute Fock moment
in place of \(p_{N_B}\). These constants are for the fixed
finite datum.

Here is the scalar computation on the original self-dual tests.
Put \(B_0=U^-\) and \(F_0=E^+\). For its rational Hermitian matrix
\(\beta=J_{B_0}\) the skew central Fourier index in
Section~\ref{sec:selfdual-fourier-integral-extension} is
\(T=j_{\rm W}\beta\). Let the parabolic Levi matrix be
\(q=\alpha_0 X_{B_0}Y_0\), where \(Y_{0,f}k\) is unitary
for an integral finite invertible matrix \(k\).
For a rational frame \(x\) with target Gram matrix \(G_x=x^*Jx\),
use the central variable \(C-C^*\) of
Section~\ref{sec:bottom-descent-height}. Its action contributes
\(\psi_E(-j_{\rm W}\operatorname{tr}(G_x C))\).
The Fourier test contributes \(\psi_E(\operatorname{tr}(TC))\).
With probability Haar on the additive matrix quotient,
\[
 \int_{\operatorname{Mat}_b(E)\backslash\operatorname{Mat}_b(\mathbb A_E)}
       \psi_E\!\left(\operatorname{tr}((T-j_{\rm W}G_x)C)\right)\,dC
     =\mathbf1_{G_x=\beta}.
\]
Over \(F_0=E^+\), Witt extension gives one orbit of ordered
isometric embeddings on each side, with representatives \((1;0)\) and \(1\), and pointwise stabilizers
\(U(W_0)(F_0)\) and \(1\). Indeed fixing the ordered frame fixes
its span pointwise, while its orthogonal complement can move by
an arbitrary isometry.

Let \(c\in\operatorname{Mat}_{2j_0\times b}(\mathbb A_E)\) be the
noncentral radical coordinate and \(J_{\rm std}\) the Gram matrix
of \(W_0\) in its adelic basis. In the adapted isotropic frame,
the radical sends a middle vector \(v\) to \(v\) plus the first
isotropic frame applied to \(c^*J_{\rm std}v\).
On a mixed slice with frame value \(x'\) in target coordinates,
the induced Heisenberg translation is therefore \(x'c^*J_{\rm std}\).
For the two representatives the frame values are
\(\binom{X_{B_0}^{-1}}{0}\) and \(X_{B_0}^{-1}\).
Thus with \(d_c=X_{B_0}^{-1}c^*J_{\rm std}\) the translations are
\(\binom{d_c}{0}\) and \(d_c\), respectively. Write \(z\) for the
unreduced central parameter. In the ordered Witt basis the
upper-right block of the radical element is
\[
  -\frac{z+z^*}{2j_{\rm W}}-\tfrac12c^*J_{\rm std}c .
\]
The mixed slice at a frame of Gram matrix \(G_x\) has multiplier
\(\psi_E(\tfrac12\operatorname{tr}(zG_x))\): in the Heisenberg
phase the quadratic \(c\)-term has zero absolute trace since
\(\operatorname{tr}(c^*J_{\rm std}cG_x)\) is fixed by conjugation
and \(j_{\rm W}\) changes sign. Replacing \(z\) by \(z-2j_{\rm W}C\)
adds \(C-C^*\) to the upper-right block and multiplies the slice by
\(\psi_E(-j_{\rm W}\operatorname{tr}(G_xC))\). Both frames have
\(G_x=\beta\), so the multipliers agree before integration.
The equal-rank block in the large integrand has zero translation.
This gives the large unfolded quotient by \(U(W_0)(F_0)\) and the
small unfolded integral over \(U(B_0)(\mathbb A_{F_0})\), with the
original Haar normalizations.

For the slices of Lemma~\ref{lem:common-real-slice}, a nonzero
large slice in \eqref{eq:marked-real-slice-big} requires
\(x^\flat\in\operatorname{supp}\sigma\). Unitarity of \(k(P(x))\)
gives \((x^\flat)^*Jx^\flat=x^*Jx\). For \(B=x_0^*Jx_0\) and \(G=x^*Jx\), the Gram bound
\(\|B^{-1}(G-B)\|<1/2\) keeps \(B+u(G-B)\) invertible for
\(0\le u\le1\). Hence \(G\) has the signature of \(J_{B_0}\)
at every real place.
The source frame changes act by invertible congruence and preserve
rank and signature. Thus an index of different signature or deficient
rank contributes zero. A nondegenerate index with a different signature
has no small orbit. At a singular index every small square frame
is singular since \(\det(B^*J_{B_0}B)=\det(J_{B_0})|\det B|^2\);
its bump \(\sigma(x_0B)\) vanishes.

At a selected split pair write the source columns in the integral
order \((K'_f,K_f,I_f)\) of the mixed model. Put
\(A=YK'_f\), \(B_1=YK_f\), and \(C=YI_f\) on the chosen half.
The stacked determinant of the standard test is, up to the fixed
basis determinant, \(\det[\,B_1\ C\,]\). Its dependence on \(A\)
is therefore exactly the indicator of the integral matrix box.
The partial Fourier transform in \(A\) reads the frame on the
conjugate half, since \(K'_f\) is dual to \(K_f\). In the normalized
coordinates its entrywise calculation is
\[
 \int_{O_w}\psi_{E,w}(\delta_w a t)\,da=\mathbf1_{O_w}(t),
 \qquad \delta=j_{\rm W}\alpha_0\bar\alpha_0,\quad
                 \operatorname{vol}(O_w)=1 .
\]
It holds also with the inverse additive character. Thus both halves
of the frame must be integral. For the small transposed test at
\(\bar w\) the determinant depends on \(YK_{f,\bar w}\), since
\(s_w^{\mathsf t}=K_{f,\bar w}\) in the paired coordinates. The
other line block has the same integral-indicator transform. The
calculation uses the full integral box also at a place of exponent
two; only the multiplicative character is read in the residue field.

A nonzero finite slice forces the large group point to
lie in \(H_fK_f^{\rm std}\), where \(H_f\) is the block
group for \(B_0\perp W_0\). Indeed multiply the integral
basis by \(k\) to obtain an integral orthonormal basis.
It extends to a unitary integral basis by splitting off its
unimodular orthogonal complement. At split places extend
the chosen basis and its dual by a complementary basis.
At an unramified nonsplit place, reduction of the
complement is a nondegenerate Hermitian space over the
quadratic residue extension. A vector of nonzero norm
exists by the residue trace pairing: if all norms vanished,
polarization after multiplying one argument by every residue scalar
would make every pairing zero. Lift this vector and scale its unit
norm to one. The residue-field norm is onto, and for \(j\ge1\) the
norm of \(1+\varpi^jz\) is congruent to
\(1+\varpi^j\operatorname{Tr}(z)\bmod\varpi^{j+1}\) for a
fixed-field uniformizer. The residue trace is onto also in
characteristic two. Successive lifting and completeness give the
required unit norm. Splitting off each resulting norm-one line
diagonalizes the unimodular complement. This proves the claim in all the finite models, which are split or unramified.

On a supported finite coset choose its representative in
\(K_f^{\rm std}\). Right translation by this compact fixes each
large mixed term. Indeed, at a chosen split place let \(A\) be the
integral target matrix. In the determinant test
\[
 \binom{s_w}{A^{-1}Y}
       =\operatorname{diag}(1_b,A^{-1})\binom{s_w}{Y}.
\]
The matrix preserves the integral box. The determinant
character \(\chi_w^{-1}\) produces the factor \(\chi_w(\det A)\)
on the prime-power test, using the inflated residue character. The target character \(\nu^{-1}\) supplies the
inverse factor, and the unramified splitting and modulus factors
are one on units. The opposite split component is determined by
unitarity; away from these pairs the integral vacuum is fixed.
Thus the finite Weil datum itself is fixed before the mixed
functional is applied, so on this coset its formula is the one at
the identity. Let \(u\in U(B_0)(\mathbb A_{F_0})\) be its small
group variable, \(B=u^{-1}Y_0\) and \(X_f=B_f\).
Both data vanish unless \(X_f\) is integral. On this support
\(X_fk\) is integral unitary, so its adjoint is its integral inverse.
Since \(k\) is integral invertible, \(x=\det X_f\) is a unit at
every finite place.
At \(w\in S\), the large residual determinant in the adapted
row order is
\[
 \det\begin{pmatrix} X_{f,w}&Z_1\\0&Z_2 \end{pmatrix}
       =\det X_{f,w}\det Z_2 .
\]
The large slice has the integral vacuum on \(Z_1\) and the
original determinant datum on the \(W_0\)-block \(Z_2\), with
scalar \(\varepsilon_w\chi_w^{-1}(x_w)\), where
\(\varepsilon_w=\chi_w(\det B_{D,w})\). On the small side the
residual variable has the integral vacuum with scalar
\(\chi_w^{-1}(\overline{x_{\bar w}})\). Here conjugation identifies
\(E_{\bar w}\) with \(E_w\). The quotient of the two factors
depending on the frame is
\[
  \frac{\chi_w^{-1}(x_w)}{\chi_w^{-1}(\overline{x_{\bar w}})}
  =\chi_w(x_w/\overline{x_{\bar w}})^{-1}.
\]
Its product over \(S\) is \(\nu^{\rm bc}(x)^{-1}\), since away
from these split pairs \(\nu^{\rm bc}\) is trivial on units. The
product of the fixed phases \(\varepsilon_w\) is \(\varepsilon\)
in \eqref{eq:period-open-orbit-scalar}.
Put \(n=\dim X\), \(d=\dim V_{\rm top}\), so that
\(n=d+b\) and \(d-b=2j_0\), and put
\(\gamma_n=\nu^{-1}\upsilon^{n/2}\).
The standard large pair is \((\gamma_n,\gamma_n)\).
The signed adjoint of the pair \((\gamma_n,\gamma_n)\)
in the reverse order has characters
\((\gamma_n^{-1}\upsilon^n,\gamma_n^{-1}\upsilon^b)\)
by the adjoint calculation below. Dividing its kernel by
\(\gamma_n(\det g)\) in the source variable gives the
standard small pair \((\nu\upsilon^{n/2},\upsilon^b)\).
Since \(K^{a_1,a_2}(g,h)\) has character multiplier
\(a_1(\det h)a_2(\det g)^{-1}\), the ratio of the large
and small target twists at \(u\) is \(\nu(\det u)^{-2}\).
For the source Levi use
\(\beta_k(z)=\beta(z/\bar z)\upsilon_E(z)^{-k}\)
as in \eqref{eq:display-f2}. Then
\[
 (\gamma_n)_d=(\nu^{\rm bc})^{-1}\upsilon_E^b,\qquad
 (\upsilon^b)_b=\upsilon_E^b ,
\]
because \(\upsilon(z/\bar z)=\upsilon_E(z)^2\). At
\(\det q\) the quotient of the two Levi factors is therefore
\((\gamma_n)_d^{-1}(\upsilon^b)_b|\cdot|^{(d-b)/2}
=\nu^{\rm bc}|\cdot|^{j_0}\).
Since
\(z_B=\det(Y_0)(\det u)^{-1}(\det X_f)^{-1}\)
as ideles, and \(\nu^{\rm bc}(\det u)=\nu(\det u)^2\),
the resulting phase and positive factor are exactly
\[
 \frac{\nu^{\rm bc}(\det q)}
      {\nu^{\rm bc}(z_B)\nu^{\rm bc}(\det X_f)\nu(\det u)^2}
        =\nu^{\rm bc}(d_0),\qquad
 \frac{|\det q|^{j_0}}{|z_B|^{j_0}}=|d_0|^{j_0}.
\]
The remaining block stabilizer quotient integral of the
bump gives the factor \(R_0\) of \eqref{eq:period-open-orbit-scalar}.

\subsubsection{The measure on the open orbit and the signed comparison}

\sectionopening{In this subsection, the first stage of the comparison in
Section~\ref{sec:marked-partner-period} unfolds over the full-rank
orbit. We keep track of the quotient measure and the source involution
to identify its signed scalar.}

Let \(G=U(B_0\perp W_0)(\mathbb A_{E^+})\),
\(L=U(B_0)(\mathbb A_{E^+})\),
\(H_1=U(W_0)(\mathbb A_{E^+})\), and
\(H=L\times H_1\) inside \(G\). This is the stabilizer of the
nondegenerate subspace \(B_0\).
Let \(\Gamma=\{1\}\times U(W_0)(F_0)\) be the stabilizer of
its ordered rational frame in the rational group.
The groups \(G\) and \(H\) are unimodular, and \(H\) is closed.
Give the homogeneous space \(H\backslash G\) the invariant quotient
measure determined by \(dG\) and \(dL\,dH_1\)
via \cite[Theorem~2.51 and (2.52)]{FollandAbstract}. The same Haar disintegration, first by
\(\Gamma\) and then by \(H\), gives for an integrable
unfolded term \(F\)
\[
 \int_{\Gamma\backslash G} F
  =\int_{H\backslash G}\int_L\int_{\Gamma\backslash H_1}
          F((l,h_1)g)\,dh_1\,dl\,d\bar g.
\]
This follows from the quotient integral formula first
for compactly supported lifts to \(G\), then for
integrable functions by applying the positive formula
to their absolute values. Right translation by each
component of \(p(h)\) in the two inner integrals preserves
the original Haar, including the measured \([W_0]\). For \(B_0=0\) the frame is
empty, \(H_1=G\), and \(L\) is the trivial group. If \(dL\) has
mass \(v_0\) on that point, the quotient integral formula gives
\(H\backslash G\) mass \(1/v_0\), since \(dH_1=dG\).
Fix the source point and the matched self-dual test. Denote by
\(A_{\rm sm}\) its small Fourier coefficient, unfolded as the
integral over \(L\) of the term at the identity basis.
Applying \eqref{eq:marked-real-slice-cancel} inside the two inner
integrals gives
\[
 \int_L\int_{\Gamma\backslash H_1} F((l,h_1)g)\,dh_1\,dl
 =\kappa_0\varepsilon\nu^{\rm bc}(d_0)I_{W_0,\nu}(z_0)
      A_{\rm sm}\,\rho(g_\infty)\mathbf1_{H_fK_f^{\rm std}}(g_f).
\]
The \(H_1\)-integral contributes the equal-rank period, and the
\(L\)-integral is precisely \(A_{\rm sm}\) with its original
Haar measure. The only dependence on the outer coset is in the
last two factors.

This weight descends to \(H\backslash G\). The finite
indicator is left \(H_f\)-invariant. At infinity
\(\rho(h)\) depends only on the span of \(h^{-1}x_0\).
On its support \(hk(P(h^{-1}x_0))\in H_\infty\), so
the corresponding cosets lie in the compact image of
\(k(\operatorname{supp}\sigma)^{-1}\).
At finite places the supported cosets are the image
of \(K_f^{\rm std}\), also compact. The weight is
positive on a neighborhood of the identity coset.
The translate \(p(h)\) of Lemma~\ref{lem:common-real-slice}
is smooth where \(\rho(h)>0\) by its block formula and is taken to
be the identity on the complementary closed set. It is therefore
measurable. The quotient integral of this continuous compactly
supported weight is the positive finite \(R_0\) in
\eqref{eq:period-open-orbit-scalar}. When \(b=0\), the basis is empty,
the determinant of its empty block is one, \(H=G\), and the outer
quotient is a single point. The same formula integrates its unit
weight with the quotient measure induced by the fixed Haar measures.
Consequently integrating the last display gives the large
coefficient \(c_{\rm top,\nu} A_{\rm sm}\), with the same
\(R_0\) even if \(I_0=0\).

For completeness, at each plane removal the old target
need not initially exhibit the next hyperbolic
summand. If its rank and negative count are those of
the new target plus one plane, form \(H\perp V\) with
its standard chosen basis and permute the basis so that
its signs agree entry by entry with the old target.
At infinity the forms are thus identified.

At finite
places the forms have the same standard
unimodular model, and Landherr's theorem
\cite[pp.~246--247]{LandherrHermitian} gives a rational isometry \(C\). If \(X_{\rm old}\)
and \(X_{\rm ad}\) are their adelic bases, write
\(X_{\rm old}=C^{-1}X_{\rm ad}u\) with \(u\) adelic
unitary. Their rational groups in chosen coordinates
are conjugate by \(u\). Their theta distributions have
the same rational vectors after the corresponding
unitary change of point.

The real representation and
the finite vacuum and test, in chosen coordinates,
are unchanged because they depend on the common
signs and on the given covector ordering.
Transport the theta functional along the isometry. Its rational
summation agrees up to a scalar with that in the adapted
polarization. The unit-vacuum mean-square normalization gives
this scalar modulus one, as in the calculation following
\eqref{eq:display-mr}. With transported functionals the target
argument changes by \(h\mapsto uh\) and the character contributes
\(\alpha(\det u)^{-1}\), where \(\alpha\) is the character of
the kept source.
Conjugation preserves the counting measure on rational
points and left translation preserves ambient Haar,
so this map preserves the induced quotient measure.

Thus the same finite covector datum occurs, and only a unit
phase changes in the scalar. In the block-reindexed
adapted model the plane basis and the new target
inclusion have their prescribed integral and
orthonormal forms. We compute the two constant terms in these
coordinates. Through the entire source
sequence use the fixed plane blocks and fixed
complement in the original adapted basis.
Successive zero slices are evaluations in the product
mixed model, so they give its combined zero slice.
Thus the same coordinates supply the adapted vectors at
each removal. At the selected finite places let \(\pi_i\) be
the middle projection at the \(i\)-th removal and \(j_i\) its
integral basis inclusion, so \(\pi_i j_i=1\). Given a lower
covector parameter \(p\), use successively \(p\pi_i\) on the
preceding spaces. For two consecutive removals the upper
parameter is \(p\pi_{i+1}\pi_i\), and restriction along
\(j_i j_{i+1}\) recovers \(p\). Each projection and inclusion
is defined before reduction to any chosen residue depth. Thus
the induction applies to every lower standard parameter, with
the same real vectors and scalar. If \(p\) has deficient rank
in a selected residue field, both standard determinant tests
vanish there. For the empty parameter the same identities use
the unique zero-row matrix.

For a single plane removal write \(e\) for the new defect
\(h-N_1/2\), using the new target rank \(h\) and source rank
\(N_1\). The old defect is \(e+1\). Under the quotient radical
mean the standard partner covector \(p\) is read on the old spaces
as \(p\,\mathrm{pr}_{W_1}\), where \(\mathrm{pr}_{W_1}\)
is the projection in the dual basis. The two
constant-term formulas on this datum have the form
\[
       \mathrm{CT}\,\Theta_{\mathrm{old}}^+
                  =\gamma_{\mathrm{pl},\nu}\Theta_{\mathrm{new}}^+,
       \qquad
       \mathrm{CT}\,\Theta_{\mathrm{old}}^-
                  =\kappa_A\Theta_{\mathrm{new}}^- .
 \tag{2}\label{eq:plane-removal-constant-terms}
\]
The lower small vector is the zero slice of the old one.
The large vector is the cone slice at \(s(e)=e+1\) for
\(e\ge0\), or its dual cone slice at \(s(e)=-e\) for \(e<0\).
The small vacuum slice scalar \(\kappa_A\ne0\) depends on
the splitting and the fixed measures. These operations
choose the new real vectors before the finite data and twist.

\subsubsection{The two radical terms after removing one hyperbolic plane}
\label{sec:marked-plane-removal-proof}

\sectionopening{We compute the two unipotent integrals on the standard
finite data after removing one hyperbolic plane. The large constant term
gives the equal-rank Levi invariance used in
Section~\ref{sec:marked-real-slice-period}. The two constant terms together
give \eqref{eq:plane-removal-constant-terms}.}

Use the adapted decompositions \(\widetilde W=H\perp W_1\)
and \(\widetilde V=H\perp V\) on the original tower.
Put \(N_1=\dim W_1\), \(h=\dim V\), and \(e=h-N_1/2\).
Retain the given middle and isotropic bases. Let \(p\)
be the integral covector parameter in the lower basis. The two
upper data use \(p\operatorname{pr}_{W_1}\), where the
projection is the middle dual covector of the source decomposition.
The upper determinant test on covectors and its column partner
are the functions \(\chi^{-1}(\det[p\operatorname{pr}_{W_1};Y])\)
and \(\chi^{-1}(\det(p\operatorname{pr}_{W_1}Z))\),
respectively, restricted to integral matrices with unit
determinant in the indicated formula.

\medskip
\noindent\emph{Finite restriction of the upper space.}
Write \((K,I,L)\) and \((K^\vee;I^\vee;L^\vee)\)
for the column and dual covector bases of the source, and use
the same notation with subscript \(V\) for the target.
At the selected split places these are integral inverse
bases. Extend the target basis to the stacked test by
the identity on the fixed covectors of \(p\). The resulting
square change of basis \(G\) takes the first and last
line coordinates to \(K_V,L_V\) on the moving block
and acts on the middle by the identity on the fixed
covectors and by \(I_V\) on the moving covectors.

Let
\(R=(K^\vee;I^\vee;L^\vee)\); both \(G\) and \(R\)
have integral inverses. In the two partial Fourier models,
the matrix inside the determinant is \(G M_\pm R\), with
\[
 M_+=
 \begin{pmatrix}
   a & z_1 &z_3\\
   0 &[p;Y]&z_2\\
   0 & 0 &u
 \end{pmatrix},\qquad
 M_-=
 \begin{pmatrix}
   u & z_1 &z_3\\
   0 &[p;Y]&z_2\\
   0 & 0 &a
 \end{pmatrix}.
\]
The column \(z_2\) is zero on the fixed covectors. The entry \(z_3\)
is the coordinate of \(K_VL^\vee\).
These block expressions are in the chosen order of the covectors:
multiplication by \(R\) reads the three dual source
coordinates, and multiplication by \(G\) restores the
target coordinates on the three moving blocks.
On the fixed covectors the diagonal line blocks vanish and
the middle \(p\) is read by precisely
\(\operatorname{pr}_{W_1}=I^\vee\).

Thus restriction to the line leaves the old
fixed covectors on the same middle coordinates. Integrality
of the determinant test is equivalent to integrality
of \(M_\pm\), and
\(\det(GM_\pm R)=\det(RG)\,a\,u\,\det[p;Y]\).
The basis phase is
\(\varepsilon=\chi^{-1}(\det(RG))\), of absolute value
one. The two diagonal line entries merely switch places.

The slice at \(\alpha_0 a\) in the self-dual
measure Fourier transforms \(z_1,z_2,z_3,u\). The frequency is
\(\bar a\) at \(u\) in the scaled coordinates. Each of \(z_1,z_2,z_3\)
has zero frequency and enters the test only through its integral
lattice indicator. The self-dual Haar gives each of these lattices
mass one at the chosen split places. Let \(\varpi=(\varpi_w)_{w\in S}\)
on the chosen split side. The \(u\)-integral at \(w\in S\)
is
\[
 \int_{O_w^\times}\chi_w(u)^{-1}\psi_w(\bar a_w u)\,du
 =
 \begin{cases}
  g_w\,\chi_w(\varpi_w\bar a_w),&
                   \varpi_w\bar a_w\in O_w^\times,\\
  0,&\text{otherwise},
 \end{cases}
 \qquad |g_w|=q_w^{-1/2}.
\]
The additive character has conductor \(O_w\) in this
coordinate. To see the conductor in the formula, write the
units as the fibers \(u_0+\varpi_w O_w\) of reduction to
\(k_w^\times\). The character \(\chi_w^{-1}\) is constant
on each fiber. If \(v_w(\bar a_w)<-1\), its additive integral
there is zero. If \(v_w(\bar a_w)\ge0\), the integral is
\(q_w^{-1}\sum_{u_0\in k_w^\times}\chi_w(u_0)^{-1}=0\).
At \(v_w(\bar a_w)=-1\) the field Gauss sum has absolute
value \(q_w^{1/2}\), and each fiber has mass \(q_w^{-1}\).
More explicitly, if \(\psi_{0,w}(t)=\psi_w(\widetilde t/\varpi_w)\)
on the residue field, then
\[
 g_w=q_w^{-1}\sum_{t\in k_w^\times}\chi_w(t)^{-1}\psi_{0,w}(t).
\]
For \(v=\varpi_w\bar a_w\in O_w^\times\), replacing \(\psi_{0,w}(t)\)
by \(\psi_{0,w}(vt)\) multiplies the sum by \(\chi_w(v)\), by the substitution \(t\mapsto v^{-1}t\).
The calculation applies at every modulus exponent \(d_w\ge1\),
since the determinant character always factors through \(k_w^\times\).
Together with the diagonal \(a\) in \(M_\pm\) it gives
the standard lower covector datum on the same \([p;Y]\),
with scalar \(l_\pm\mathfrak g\,F_\varpi(a_f)\).
Here \(l_\pm\) is the nonzero vacuum slice constant in
the indicated Fourier model,
\[
 \mathfrak g=\varepsilon\prod_{w\in S}g_w,\qquad
 F_\varpi(a_f)=\chi^{-1}(a_S)\chi((\bar a\,\varpi)_S)
\]
on \(a_S,(\bar a\,\varpi)_S\in O_S^\times\) and
\(a\) integral outside \(S\cup\bar S\), and zero otherwise.
At a selected split pair, write \(x_1,x_2\) for the two
coordinates, identifying the second with \(E_w\) by conjugation.
The scaled alternating trace pairing is
\(\psi_w(x_1y_2-x_2y_1)\), since the scale
\(j_{\rm W}\bar\alpha_0\alpha_0\) changes sign under conjugation.
Both additive measures give \(O_w\) mass one. Put
\(g(\theta)=q_w^{-1}\sum_{t\in k_w^\times}\theta(t)\psi_{0,w}(t)\).
The transform of the local factor
\(\chi_w^{-1}(x_1)\chi_w(\varpi_w x_2)\), with its two unit
conditions, has the same support and the same values:
the \(x_1\)-integral contributes
\(g(\chi_w^{-1})\chi_w(\varpi_w y_2)\), and substitution
\(u=\varpi_w x_2\) in the \(x_2\)-integral contributes
\(q_w g(\chi_w)\chi_w^{-1}(-y_1)\).
Indeed \(q_w g(\chi_w)g(\chi_w^{-1})=\chi_w(-1)\), by
complex conjugation and \(|g(\chi_w)|^2=q_w^{-1}\).
This also works with the inverse Fourier character, since the
two minus signs cancel in the original test. Away from the
chosen pairs the integral is the original vacuum slice. Thus
the two line tests have the same \(F_\varpi\) with their
fixed vacuum constants \(l_\pm\). In particular both line data have zero
value at the zero line coordinate.

\medskip
\noindent\emph{The large radical mean.}
At a source Levi point \((u,g)\), take \(u=1\), or allow
arbitrary \(u\) if \(h=N_1\). Both original kernels can be
integrated in either order over the compact source radical
quotient and the finite-volume target quotient. Indeed the
large determinant test has right-invertible moving rows and the
small column test gives a left-invertible \(Z\) at a conjugate
split place, since \(pZ\) is invertible there. The ordinary whole-kernel estimate at the end of
Section~\ref{sec:marked-mixed-absolute-summability}, or the primitive
row and column estimates of Section~\ref{sec:four-unitary-constant-terms},
is uniform as
the radical point ranges over compact representatives. A target
of rank zero contributes one value.

Now fix the target point. On the compact radical representatives
the translated archimedean tests form a bounded set in Schwartz
space, and their finite parts range over a finite set. In rational
coordinates for the mixed slice and its residual polarization,
the rational indices belong to a fixed fractional lattice.
The sum of \((1+\|x\|+\|y\|)^{-M}\) on this lattice is finite
for sufficiently large \(M\). Thus the mixed expansion at this
target point can be integrated term by term over the radical.
The resulting vector terms are absolutely summable at the fixed
point and can be grouped by their rational lines. The central
average of the source line retains exactly
the rational vectors of zero Hermitian norm in
\(\widetilde V\). Its other radical coordinates produce
the line zero-term functional for each retained nonzero
vector. For an isotropic vector \(x=\xi k_V\), \(\xi\in E^\times\),
write \(c\) for the noncentral radical coordinate in the
adapted Witt basis. The central action on this slice is trivial.
The action on the residual theta test is the Heisenberg
translation by \(x c^* J_{W_1}\), where \(J_{W_1}\) is the
Gram matrix of the retained source space in a rational basis.
The map \(c\mapsto x c^*J_{W_1}\) is a rational additive
isomorphism onto \(Z=\operatorname{Hom}_E(W_1,Ek_V)\)
and sends the normalized Haar of the compact rational quotient
to normalized Haar. Take the self-dual mixed polarization in
which \(Z\) acts by modulation. Its opposite block is
\(Z^\vee\), paired with \(Z\) by the additive character of the
Weil representation. For each rational index \(\lambda\) in
that block the mean contributes
\begin{equation}\label{eq:radical-target-line-zero-term}
 \int_{Z(F_0)\backslash Z(\mathbb A_{F_0})}
        \psi\bigl(\langle z,\lambda\rangle\bigr)\,dz
           =\mathbf 1_{\lambda=0}.
\end{equation}
The absolute bound above permits this termwise integration.
It leaves the residual theta on the target-line zero slice,
with scalar one. The conjugate polarization replaces the
character by its inverse and gives the same mean. For the zero vector the source slice at a
chosen \(w\in S\) restricts the moving covector matrix \(Y\) to
\(YK=0\), with \(K\) the first isotropic basis column. The fixed rows
also kill \(K\), since \(I^\vee K=0\). Thus \([p I^\vee;Y]K=0\),
contrary to invertibility on the support and \(K^\vee K=1\).
This source slice vanishes before taking a target-line term. Moving
the target point preserves its vanishing by equivariance of the zero slice.

\subsubsection{Unfolding the isotropic line orbit}
\label{sec:radical-isotropic-line-orbit}

\sectionopening{In this subsection, we group the nonzero radical terms before the equal-rank vanishing in Section~\ref{sec:marked-equal-rank-vanishing}.}

Witt extension gives one orbit of rational isotropic lines.
Group the nonzero vectors by this orbit, then by their
coordinate \(\xi\in E^\times\) on the line.
Let \(P\subset U(\widetilde V)\) stabilize \(E k_V\), and let
\(N\) be its unipotent radical. In the ordered Witt decomposition,
\[
 m(b,h_0)=\operatorname{diag}(b,h_0,\bar b^{-1}),\qquad
 P=N\rtimes\bigl(\operatorname{Res}_{E/F_0}\mathrm{GL}_1\times U(V)\bigr).
\]
The stabilizer of \(k_V\) is \(N\rtimes U(V)\); the rational
scalars moving \(k_V\) on its line form \(E^\times\).
In coordinates for \(N\) with a column \(c\) in \(V\) and a
central coordinate \(t\) on an \(F_0\)-line, conjugation by
\(m(b,h_0)\) acts by
\((c,t)\mapsto(h_0 c\bar b,b\bar b t)\). Its absolute
Jacobian is \(|b|_{\mathbb A_E}^{h+1}\), since
\(|\det h_0|_{\mathbb A_E}=1\). Hence for
\(\tilde h=n m(b,h_0)k\) the Iwasawa density is
\(|b|_{\mathbb A_E}^{-(h+1)}\,dn\,d^\times b\,dh_0\,dk\)
up to its fixed positive Haar scalar. Here \(k\) ranges over the
original maximal compact of \(U(\widetilde V)(\mathbb A_{F_0})\);
the finite Iwasawa decomposition in these lattices is proved in
Section~\ref{sec:bottom-descent-planes}. Since
\(m(b,h_0)^{-1}k_V=b^{-1}k_V\), the slice coordinate is
\(a=b^{-1}\), and its density is \(|a|_{\mathbb A_E}^{h+1}\).

For the positive line term at \(\xi k_V\), write \(S_{\xi k_V}\)
for the source slice and \(T_0\) for the subsequent target zero
slice. For \(n\in N(\mathbb A_{F_0})\) the mixed model gives
\[
 S_{\xi k_V}\omega(1,n)f=\omega_{\rm low}(1,n)S_{\xi k_V}f,
 \qquad T_0\omega_{\rm low}(1,n)=T_0.
\]
In the first equality \(n^{-1}k_V=k_V\) and \(\det n=1\).
For the second, in the target line model the displacement from the
off-diagonal column is linear in the slice coordinate, and the
central phase is quadratic in it. At coordinate zero the displacement
and phase exponent vanish, and \(n\) induces the identity on
\(k_V^\perp/Ek_V\). Any remaining scalar on this zero slice is a
character of \(N\) invariant under \(m(2,1)\). Conjugation multiplies
the centre by \(4\), so the character is trivial on the centre.
On the vector quotient conjugation multiplies by \(2\), so the
character there is trivial as well. Thus the full target radical fixes the term
before summing over \(\xi\). Its quotient by \(N(F_0)\) has a fixed
positive finite volume: the column and central coordinates have
compact additive adelic quotients, with product Haar since the
group law acts by successive translations in these coordinates.
When \(h=0\), this central line remains and the column is absent.

Let \(\ell_\pm(z;h_0,k)\) be the lower kernel at the
line slice at the actual parameter \(z\),
target point \(h_0\) and old target compact point \(k\);
the retained source point is \(g\).
Put \(x=au\), where \(u=1\) unless \(h=N_1\).
In the positive model, first slicing on the source gives
\(|u|^{(h+2)/2}\) from the old target rank. Slicing on the
target then gives \(|a|^{-N_1/2}\) from the lower source rank.
In the opposite order these factors are \(|a|^{-(N_1+2)/2}\)
and \(|u|^{h/2}\). The line arguments for a rational index
\(\xi\in E^\times\) are \(\xi x\) and \(\xi\bar x^{-1}\),
respectively. To compute the unitary factors, write
\(\gamma_j=\nu^{-1}\upsilon^{j/2}\) and
\(\eta_k(z)=\eta(z/\bar z)\upsilon_E(z)^{-k}\).
A source-line Levi of parameter \(u\) contributes
\(\eta_{\dim(\mathrm{target})}(u)^{-1}\), and a target-line Levi
of parameter \(b=a^{-1}\) contributes
\(\eta_{\dim(\mathrm{source})}(a^{-1})\), with the character of
the opposite space in each case. Removing a plane lowers its
space character by one power of \(\upsilon\).
For the source first these two factors are
\((\gamma_{N_1+2})_{h+2}(u)^{-1}\) and
\((\gamma_{N_1})_{N_1}(a^{-1})\); for the target first they are
\((\gamma_{N_1+2})_{N_1+2}(a^{-1})\) and
\((\gamma_{N_1})_{h}(u)^{-1}\). In both orders their product is
\[
  \nu^{\rm bc}(au)\upsilon_E(u)^{h-N_1}=\nu^{\rm bc}(x),
\]
since \(u=1\) or \(h=N_1\). Including the target
orbit density gives the two scalar factors
\[
 \begin{aligned}
  \nu^{\rm bc}(x)|a|^{h+1}|u|^{(h+2)/2}|a|^{-N_1/2}
     &=\nu^{\rm bc}(x)|x|^{e+1},\\
  \nu^{\rm bc}(x)|a|^{h+1}|a|^{-(N_1+2)/2}|u|^{h/2}
     &=\nu^{\rm bc}(x)|x|^{e}.
 \end{aligned}
\]
In the second expression \(x'=\bar x^{-1}\) satisfies
\(\nu^{\rm bc}(x')=\nu^{\rm bc}(x)\) and \(|x'|=|x|^{-1}\),
so the exponent in the actual negative line variable is \(-e\).
These coordinates give the following integral for the
contribution of one line modulo its rational stabilizer, with
the exponent \(s=e+1\) supplied by the densities. We verify its
absolute convergence before unfolding the sum over the orbit. Put
\[
 Z(s)=\int \nu^{\rm bc}(au)|au|^s
       \sum_{\xi\in E^\times}\ell_+(\xi au;h_0,k)
                 \,d^\times a\,dh_0\,dk,
\]
where \(a\in E^\times\backslash\mathbb A_E^\times\),
\(h_0\in U(V)(F_0)\backslash U(V)(\mathbb A_{F_0})\), and
\(dk\) is Haar on the original maximal compact.
The scalar \(c_{\rm unf}>0\) is the Iwasawa Haar scalar times
the volume of \(N(F_0)\backslash N(\mathbb A_{F_0})\). The
quotients use counting measure on their rational subgroups.
In particular the full \(E^\times\)-sum unfolds with this
quotient measure, including all integral units and roots of unity.
For \(V=0\) the \(h_0\)-integral has the prescribed mass of the
trivial group. The Fourier
formula in the two opposite line coordinates gives
\[
 |a|\sum_{\xi\ne0}\ell_+(\xi a;h_0,k)
   =\sum_{\xi\ne0}\ell_-(\xi\bar a^{-1};h_0,k).
\]
To identify the two polarizations, write
\(A=EK,\ A'=EL,\ B=EK_V,\ B'=EL_V\) in the ordered rational
Witt bases. The common zero slice removes the isotropic subspace
\[
 D=\operatorname{Hom}(W_1\oplus A',B)
       \oplus\operatorname{Hom}(A',V)
       \ \subset\operatorname{Hom}(\widetilde W,\widetilde V).
\]
Extend each map by zero on the other displayed summands. The trace
pairing of the matrix units pairs opposite Witt indices. Thus
\[
 D^\perp/D=\operatorname{Hom}(W_1,V)\ \perp\
               (E K_VK^\vee\oplus E L_VL^\vee).
\]
The last two maps have zero norms and pair to one for the Hermitian
trace form. In the first order of slicing the translating map is
\(aK_VK^\vee\), and in the reverse order it is \(aL_VL^\vee\).
The remaining Weyl change is therefore Fourier transformation in this
single adelic line with its original self-dual measure. It fixes the
\(\operatorname{Hom}(W_1,V)\) coordinates. This description also
applies when that summand is zero. For fixed \(h_0,k,g\), the remaining theta sum
commutes with this transform: its rational finite coordinates range
in a fractional lattice, and Schwartz decay in the rational and
line variables bounds the absolute sum of the line integrals.
Write \(P_+=K_VK^\vee\), \(P_-=L_VL^\vee\), and abbreviate the fixed
lower and compact arguments of \(\ell_\pm\). Let \(F\) be the
common complementary zero slice and \(\Theta\) its residual theta
functional. The commutator of \(\mathcal R(P_+x)\) and
\(\mathcal R(P_-y)\) is
\(\langle x,y\rangle=\psi_E(j_{\rm W}\bar y x)\), since \(P_+,P_-\)
pair to one. In the first polarization these operators translate
by \(x\) and multiply by \(\langle\,\cdot\,,y\rangle\),
respectively. In the second their roles are reversed and the
commutator is inverted. Thus the restrictions of the displaced
theta value to the two isotropic lines have the expansions
\[
 \begin{aligned}
 \Theta(\mathcal R(P_-y)F)&=\sum_{\xi\in E}
                             \ell_+(\xi)\langle\xi,y\rangle,\\
 \Theta(\mathcal R(P_+x)F)&=\sum_{\eta\in E}
                             \ell_-(\eta)\langle x,\eta\rangle^{-1}.
 \end{aligned}
\]
Both sums at the identity equal \(\Theta(F)\). More generally, put
\([E]=E\backslash\mathbb A_E\) with probability Haar. The opposite
slice is the mean after modulation and translation in the first
polarization. For every \(y\in\mathbb A_E\) it gives
\[
 \begin{aligned}
 \ell_-(y)
 &=\int_{[E]}\Theta(\mathcal R(P_+x)\mathcal R(P_-y)F)\,dx\\
 &=\int_{[E]}\sum_{\xi\in E}
           \ell_+(\xi+x)\langle\xi+x,y\rangle\,dx
   =\int_{\mathbb A_E}\ell_+(z)\langle z,y\rangle\,dz .
 \end{aligned}
\]
The last Haar measure has rational covolume one. In this Fourier
convention the transform of \(z\mapsto f(az)\) is
\(y\mapsto |a|^{-1}\widehat f(\bar a^{-1}y)\).
The line function and its transform are Schwartz--Bruhat:
derivatives and polynomial weights in that line are bounded
after summing the remaining rational lattice by taking sufficiently
many further powers of Schwartz decay in the lattice variables.
Their periodizations converge uniformly on compact representatives.
Adelic Poisson summation \cite[Lemma~4.2.4]{TateThesis} therefore
gives the displayed equality before the substitution \(z=\alpha_0a\)
in the finite line test.
Both zero-frequency contributions are zero by the selected finite
unit support of the respective slices.
Put \(a'=\overline{au}^{-1}\). Then
\(\nu^{\rm bc}(a')=\nu^{\rm bc}(au)\) and
\(|a'|=|au|^{-1}\). Multiplicative Haar on the idele class group
is invariant under translation, conjugation and inversion.
Substituting the Poisson formula shows that \(Z(s)\) uses the
negative line sum with exponent \(1-s\).
In the equal-rank case, \(a\mapsto au\) shows directly
that the constant term is independent of \(u\).

For \(s>1\) unfold the sum over \(E^\times\) in the positive-sign formula for \(Z(s)\) or in the negative-sign formula for
\(Z(1-s)\). From the finite factor just obtained, the
idele translation is \(K=\alpha_0 c_0^{-1}\), where
\(c_0\) is \(\bar\varpi_w\) at \(\bar w\) and one
elsewhere. Thus \(|K|=|\alpha_0|Q_S\), with
\(Q_S=\prod_{w\in S}q_w\).
In the original line slices substitute \(z=Ka\). Since
\((c_0)_\infty=1\), this gives
\[
    (z/\alpha_0)_\infty=a_\infty,\qquad
    F_\varpi((z/\alpha_0)_f)=F_1(a_f),
\]
where \(F_1\) uses the same formula with \(\varpi=1\).
Multiplicative Haar is invariant, and the character and modulus
supply the factor \(\nu^{\rm bc}(K)|\alpha_0|^s Q_S^s\).

Normalize the local multiplicative measures by
\(\operatorname{vol}_{d^\times a_w}(\mathcal O_{E_w}^\times)=1\)
and take their restricted product on the finite ideles.
The original finite measure is \(v_\pm\) times this product, where
\(v_\pm>0\) is its mass on \(\widehat{\mathcal O}_E^\times\), fixed before \(S\).
At \(w\in S\cup\bar S\) the translated function is supported on units
and \(\nu^{\rm bc}F_1=1\) there, so the normalized local integral is one.
Elsewhere the integral is
\(\sum_{n\ge0}\nu_w^{\rm bc}(\varpi_w)^n q_w^{-ns}\), since each
multiplicative shell has mass one. The full finite factor is therefore
\(v_\pm L^S(s,\nu^{\rm bc})\), abbreviated to \(v_\pm L^S(s,\nu)\) below.
In absolute value the sum is bounded
by the partial Dedekind zeta series at \(s>1\).
The real integral is the cone vector \(C_\pm(s)\)
formed by integrating the real line slice against
\(\nu^{\rm bc}(a_\infty)|a_\infty|^s\) and the original
target compact Haar. For a finite point \(k\) of the original
target compact, let \(k_S\) be its component on the chosen split side.
In the upper covector model set
\(\phi_p(Y)=\chi^{-1}(\det[p\operatorname{pr}_{W_1};Y])\)
on the integral unit support and zero elsewhere. Multiplication by
\(\operatorname{diag}(I,k_S^{-1})\) preserves that support and gives
\[
  \phi_p(k_S^{-1}Y)=\chi(\det k_S)\phi_p(Y),\qquad
  (\nu^{-1}\upsilon^{\dim(\widetilde W)/2})(\det k)
                 =\chi(\det k_S)^{-1}.
\]
The second equality uses the unramified splitting character and the
prescribed tame character on units. The two factors cancel. The
unselected compact factors preserve the integral vacuum. Thus the
finite compact fixes the standard vector in the twisted Weil action,
and its images in the two slice formulas have the same invariance.

For either polarization let \(v_\pm(k_\infty,a_\infty)\)
be the real line slice before insertion of the Mellin weight.
The common zero slice is restriction in a fixed mixed model. After
the possible one-line Fourier transform it is still Schwartz jointly
in the line and the remaining coordinates. As \(k_\infty\) varies
in the maximal compact these functions form a bounded Schwartz family.
Thus for every remaining Schwartz seminorm \(p_j\) and integer \(M\),
\[
 \sup_{k_\infty}p_j(v_\pm(k_\infty,a_\infty))
       \ll_{j,M}\prod_{\tau\mid\infty}(1+|a_\tau|)^{-M},
\]
where the absolute values are ordinary complex moduli. For
\(0<\sigma_0\le\Re s\le\sigma_1\), the weighted majorant at
each place is \(r^{2\sigma_0}\) for \(r\le1\) and
\(r^{2\sigma_1-M}\) for \(r>1\), with radial measure \(dr/r\).
Choosing \(M>2\sigma_1\) gives absolute convergence in every
seminorm, also after any fixed power of \(\log|a_\infty|\).
Schwartz completeness gives the cone vector.
For either order of slicing, the retained covector test has
full row rank on its support at a selected split place, by the
invertibility of \([p;Y]\). The fixed lower source point preserves
this rank by invertible right multiplication. Use the full lower
source basis in the kernel estimate of
Section~\ref{sec:four-unitary-constant-terms}.
Equation~\eqref{eq:partner-exchange-ordinary-bound} thus bounds the
absolute value of the lower theta kernel uniformly on the target
quotient by a fixed Schwartz seminorm. For rank zero use evaluation
on the one-point quotient. At \(z=\alpha_0a\) the finite slice
is \(l_\pm\mathfrak g F_\varpi(a_f)\) times the retained covector
test. Together with the real seminorm estimate this gives, for
each \(M>0\),
\[
 \sup_{h_0,k}|\ell_\pm(\alpha_0a;h_0,k)|
 \ll_{g,M}|F_\varpi(a_f)|
                 \prod_{\tau\mid\infty}(1+|a_\tau|)^{-M}.
\]
Choose an integral ideal \(\mathfrak d\)
clearing the finite line-coordinate denominators in both slices.
Such an ideal exists because their finite tests are compactly
supported and have only finitely many translates under the compact
variable \(k_f\). In the normalized finite multiplicative measure,
the integral of \(|z_f|^\sigma\) on
\(\mathfrak d^{-1}\widehat{\mathcal O}_E\cap\mathbb A_{E,f}^\times\)
is \(\mathrm N(\mathfrak d)^\sigma\zeta_E(\sigma)\), by summing the
valuation shells. For \(\sigma>1\) and \(M>2\sigma\) this gives
\[
 \int_{\mathbb A_E^\times}\!\!\int_{[U(V)]\times K_{\widetilde V}}
      |z|^\sigma|\ell_\pm(z;h_0,k)|\,dh_0\,dk\,d^\times z
 \ll \mathrm N(\mathfrak d)^\sigma\zeta_E(\sigma)
     \prod_{\tau\mid\infty}\int_0^\infty
           r^{2\sigma}(1+r)^{-M}\,\frac{dr}{r}<\infty .
\]
Here \(K_{\widetilde V}\) denotes the original maximal compact of
the target, and the constant includes \(\operatorname{vol}([U(V)])\)
and the fixed measure factors. The bound justifies
termwise unfolding for each sign in its own half-plane
\(\operatorname{Re}s>1\).
For the Mellin integral on the original quotient, write
\[
 D(a,h_0,k)=\nu^{\rm bc}(au)
       \sum_{\xi\in E^\times}\ell_+(\xi au;h_0,k),
 \qquad r=|au|.
\]
The positive unfolding bounds the integral of
\(|D(a,h_0,k)|r^{\sigma_+}\) for \(\sigma_+>1\).
For \(\sigma_-<0\), apply Poisson summation to the sum
before taking its absolute value. The zero terms vanish as above.
With \(a_-=\overline{au}^{-1}\), the Poisson identity and
unitarity of \(\nu^{\rm bc}\) give
\[
 \begin{aligned}
 |D(a,h_0,k)|r^{\sigma_-}
    &=|a_-|^{1-\sigma_-}
      \left|\sum_{\xi\in E^\times}\ell_-(\xi a_-;h_0,k)\right|\\
    &\le |a_-|^{1-\sigma_-}
      \sum_{\xi\in E^\times}|\ell_-(\xi a_-;h_0,k)|.
 \end{aligned}
\]
The negative unfolding bounds the integral because
\(1-\sigma_->1\) and \(d^\times a_-=d^\times a\).

Consequently \(r^\sigma\le r^{\sigma_-}+r^{\sigma_+}\) for
\(\sigma_-\le \sigma\le \sigma_+\) gives an integrable majorant
for the grouped integral. Widening the interval absorbs every
power of \(\log r\), giving holomorphy of the Mellin integral.

To finish the unfolding, put \(G=U(\widetilde V)\) and
\([G]=G(E^+)\backslash G(\mathbb A_{E^+})\). Let \(P_L\)
stabilize the chosen rational isotropic line \(L\), and write
\(B_L(\tilde h)\) for the sum of the radical terms indexed by the
nonzero rational vectors on \(L\). Rational covariance permutes these
vectors, so \(B_L\) is left \(P_L(E^+)\)-invariant. Each line
term is invariant under the adelic target radical by the zero-slice
calculation above. In the Iwasawa
coordinates above, \(|B_L|\) times the orbit density is
\(|D(a,h_0,k)|\,r^{e+1}\), with the fixed measure scalar.
The radical quotient is compact. The preceding strip bound
therefore gives
\[
 \int_{P_L(E^+)\backslash G(\mathbb A_{E^+})}|B_L|
  = \int_{[G]}\sum_{\delta\in P_L(E^+)\backslash G(E^+)}
                   |B_L(\delta\tilde h)|\,d\tilde h <\infty .
\]
The equality is the quotient integration formula for nonnegative
functions. The zero-vector term vanishes, and Witt extension gives
one orbit on the remaining isotropic lines. We can now integrate
the grouped line expansion and unfold it to
\(c_{\rm unf}Z(e+1)\), including at \(e+1=0,1\).

In both cases the value required is \(Z(e+1)\). If \(e\ge0\),
the positive cone has exponent \(s(e)=e+1\). If \(e<0\),
the negative cone has exponent \(s(e)=1-(e+1)=-e\).
The boundary values at \(e=0,-1\) are
\[
 Z(1)=\lim_{t\downarrow1}Z(t),\qquad
 Z(0)=\lim_{t\downarrow1}Z(1-t),
\]
respectively. The real cone factors are continuous at one by
their Schwartz bounds, and the nontrivial Hecke \(L\)-function
is regular there.
The exact scalar multiplying the lower measured
period in the respective model is
\[
 \gamma_{\pm,\nu}
  =c_{\rm unf}\,l_\pm\,v_\pm\,
       (\mathfrak g\sqrt{Q_S})\,\nu^{\rm bc}(K)\,
       |\alpha_0|^{s(e)}Q_S^{\,s(e)-1/2}L^S(s(e),\nu).
\]
In particular \(|\mathfrak g\sqrt{Q_S}|=1\).
This gives the first identity of \eqref{eq:plane-removal-constant-terms} and its stated
constant with all measure-dependent quantities
\(c_{\rm unf},l_\pm,v_\pm\) fixed before the twist.

\medskip
\subsubsection{The small radical mean and its zero term}\label{sec:radical-local-line-terms}

\sectionopening{In this subsection, we compute the second radical integral of Section~\ref{sec:marked-plane-removal-proof} and show that only its zero term remains.}

For the small side read the transposed column datum
at \(\bar S\). In the original paired bases the
middle projection read at \(S\) is the transpose
of the middle inclusion read at \(\bar S\).
Consequently the test in the old basis only reads
the lower complementary block of the transpose
inside \(\det(pZ)\). All other factors are the
integrality conditions of the original lattice.
Its source-line slice at \(\alpha_0x\) is zero unless
\(x_f\) is integral, and otherwise equals the standard
lower column datum, times the real slice at \(x_\infty\)
and the fixed vacuum slice constant \(\kappa_A\ne0\).
Indeed the extra sliced coordinates have exactly the
constant integral-indicator test of the vacuum in
this projection basis. In particular the slice at zero
uses exactly the same lower finite datum at \(\bar S\).

Take the source radical mean of this small kernel at
\(\operatorname{levi}(1,g)\).
The central average again retains the isotropic rational
vectors of the fixed small space \(U\).
At a nonzero isotropic vector, apply
\eqref{eq:radical-target-line-zero-term} with the transported
line in the small target. This first takes the source slice and then
the target zero term. The surviving rational indices have zero
component in the opposite block. Hence the resulting theta
functional is invariant under all Heisenberg translations along
the transported target line.

The transported line vector at any selected place
is nonzero. In the finite Schr\"odinger polarization
at \(\bar S\) this invariance is invariance under all
multiplications by the additive characters paired
against \(x^*Z^{\rm t}\), where \(x\) is that line
coordinate and \(Z\) the column variable.
On the support, \(pZ\) is invertible, so \(Z\) is
left-invertible; hence \(x^*Z^{\rm t}\ne0\).
Nondegeneracy of the additive pairing makes one of
these multipliers different from one at each support
point.

It is constant on a sufficiently small compact
open neighborhood. Compactness gives a finite
disjoint refinement of this cover, and on each part
the invariant functional must be zero. This proves
vanishing of the whole nonzero isotropic term, with
no size estimate and no change of finite datum.
Anisotropic terms have already vanished in the
central mean.

Only the zero source slice remains. At that slice the
retained group action is the old one restricted to
\(W_1\) and \(U\). The source splitting character drops
by the prescribed one power of \(\upsilon\) when the
rank drops by two, which is the lower rank twist.
The finite part is exactly the same standard column
datum and the small real vector is the zero slice.
Thus the small radical mean is \(\kappa_A\) times
the lower measured period, as asserted in \eqref{eq:plane-removal-constant-terms}.
For \(U=0\) this argument simply retains the unique
zero vector on the zero-group quotient.

In particular,
\[
 |\gamma_{\mathrm{pl},\nu}|=C_{\mathrm{pl}}L_\nu^\sharp(s(e)),
 \qquad L_\nu^\sharp(s)=|L^S(s,\nu)|\prod_{w\in S}q_w^{s-1/2},
 \tag{3}\label{eq:plane-period-local-norm}
\]
where \(C_{\mathrm{pl}}>0\) includes the original line-orbit,
slice and unit-volume factors, all fixed before the twist.

For the real section needed below, let \(x_B\) be the
target isotropic vector. Give the lower vector \(v\)
a target-line mixed-model extension \(v_B\) by tensoring
with a Schwartz function of value one at zero in the
removed coordinates. At each real place \(x_{B,\tau}\ne0\), and
Witt extension gives the product orbit
\[
 \prod_{\tau\in\Phi}
 \{x_\tau\in\widetilde V_\tau\setminus\{0\}:h_\tau(x_\tau,x_\tau)=0\}.
\]
Its orbit maps are submersions. Take smooth local sections
\(u_x x_B=x\) with \(u_{x_B}=1\) on this product near \(x_B\).
At every real place the nonzero isotropic cone is a smooth
real hypersurface near \(x_{B,\tau}\): the differential of
\(x\mapsto \langle x,x\rangle\) there is the nonzero real form
\(v\mapsto2\Re\langle x_{B,\tau},v\rangle\). A submanifold chart
gives a smooth retraction onto the cone, so the local section extends
to an ambient neighborhood by composition with this retraction.
Choose a smooth nonnegative \(\beta\) with \(\beta(x_B)=1\)
supported in a product of balls whose closures lie in that
neighborhood and whose radii are less than \(\|x_{B,\tau}\|/4\).
Multiplication by \(\beta\) and extension by zero gives a global
Schwartz slice.
Prescribe on the source-line mixed model the slice
\[
            \beta(x)\omega_{N_1/2-\mu}(1,u_x)v_B .
\]
On the support of the cone integral at \((\kappa,a_\infty)\)
we have \(x=\kappa_\infty^{-1}x_Ba_\infty\).
Then \(\kappa_\infty u_x\) sends \(x_B\) to
\(x_Ba_\infty\).

Let \(h\in U(V)_\infty\) be its induced isometry on
\((E_\infty x_B)^\perp/(E_\infty x_B)\).
The diagonal blocks in the adapted Witt basis are
\(a_\infty,h,\bar a_\infty^{-1}\).
The earlier target radical calculation applies to the real zero
slice: the displacement and central phase exponent vanish
there. Thus only the three Levi blocks act on that slice.
At \(\tau\) put \(u_\tau=a_\tau/|a_\tau|\), where
\(|a_\tau|_{\mathbb C}=|a_\tau|^2\) is the normalized local absolute
value. The isotropic block on
\(\operatorname{Hom}(W_{1,\tau},\mathbb C x_{B,\tau})\) is
\(a_\tau I_{N_1}\). Equation~\eqref{eq:marked-hyperbolic-levi} gives
its intrinsic multiplier \(u_\tau^{-N_1}|a_\tau|_{\mathbb C}^{N_1/2}\).
The determinant of the three Levi blocks divided by \(\det h_\tau\)
is \(a_\tau/\bar a_\tau=u_\tau^2\). The source slice reduces the
target exponent from \(N_1/2+1-\mu_\tau\) to \(N_1/2-\mu_\tau\).
Relative to the action of \(h_\tau\), the zero-slice multiplier is
therefore
\[
 u_\tau^{-N_1}|a_\tau|_{\mathbb C}^{N_1/2}
       (u_\tau^2)^{N_1/2-\mu_\tau}
   =|a_\tau|_{\mathbb C}^{N_1/2}u_\tau^{-2\mu_\tau}.
\]
The phase over all \(\tau\) cancels \(\nu^{\rm bc}(a_\infty)\), and
\(|a_\infty|=\prod_\tau |a_\tau|_{\mathbb C}\).
The translation by \(h\) disappears from the lower
period by right invariance on its original quotient.

The cone period of the constructed real vector is
therefore the prescribed lower period times
\[
 B(s)=\int_{\kappa,a_\infty}\beta(\kappa_\infty^{-1}x_Ba_\infty)
            |a_\infty|^{s+N_1/2}\,d\kappa\,d^\times a_\infty .
\]
On the support, Euclidean unitarity of \(\kappa_\infty\) gives
\(3/4<\|\kappa_\tau^{-1}x_{B,\tau}a_\tau\|/\|x_{B,\tau}\|
      =|a_\tau|<5/4\) at every real place, where \(|a_\tau|\)
is the ordinary complex modulus. At \((\kappa,a_\infty)=(1,1)\)
the bump is one, so the integrand is positive on an open set
of positive measure for every real \(s\). Differentiation in \(s\)
on the fixed compact annuli proves that \(B\) is entire. The same
bump and normalization apply at every modulus and finite character.
At a fixed positive exponent we may divide the real
vector by \(B(s)\) to have exactly the desired period.
For the continued plus identity at a negative exponent
the same section works without altering its finite test.

 The Euler product proves nonvanishing for real \(s>1\).
Taking logarithms of Euler products and using
\(3+4\cos\theta+\cos 2\theta\ge0\) prime power by prime power gives
\[
 1\le \zeta_E^S(s)^3 |L^S(s,\nu^{\rm bc})|^4
                  |L^S(s,(\nu^{\rm bc})^2)| .
\]
All three products omit \(S\cup\bar S\). The pole of
\(\zeta_E^S\) at \(1\) is simple. A zero of
\(L^S(s,\nu^{\rm bc})\) at \(1\) would therefore make the right
side tend to zero as \(s\downarrow1\). The square character is nontrivial
because the weights are positive, so its \(L\)-function is
regular at \(1\) by the Hecke continuation theorem.

Since \(M_\nu(j+1)/M_\nu(j)=L_\nu^\sharp(j+1)^2\)
for \(j\ge0\), and \(M_\nu(-j)=M_\nu(j)^{-1}\), \eqref{eq:plane-period-local-norm} gives
\[
  |\gamma_{\mathrm{pl},\nu}|^2 M_\nu(e)
          =C_{\mathrm{pl}}^2 M_\nu(e+1).
\]
Let \(\epsilon_{\mathrm{ad},\nu}\) be the unit multiplier of
the preceding target transport, with the convention that the upper
identity in the adapted coordinates reads
\(\Theta_{\mathrm{old}}^+
 =\epsilon_{\mathrm{ad},\nu}c_{\mathrm{old},\nu}\Theta_{\mathrm{old}}^-\).
It is the product of the functional and character multipliers of the
chosen isometry and is defined even when \(c_{\mathrm{old},\nu}=0\).
For an already adapted model take \(\epsilon_{\mathrm{ad},\nu}=1\).
Taking the two means in \eqref{eq:plane-removal-constant-terms} gives
\[
 c_{\mathrm{new},\nu}
  =\epsilon_{\mathrm{ad},\nu}c_{\mathrm{old},\nu}
        \kappa_A/\gamma_{\mathrm{pl},\nu}.
\]
The positive normalizing factor
changes from \(D_{\mathrm{old}}\) to
\(D_{\mathrm{old}}C_{\mathrm{pl}}^2/|\kappa_A|^2\), so
\[
 \frac{D_{\mathrm{new}}|c_{\mathrm{new},\nu}|^2}{M_\nu(e)}
   =\frac{D_{\mathrm{old}}|c_{\mathrm{old},\nu}|^2}{M_\nu(e+1)} .
 \tag{4}\label{eq:plane-removal-ratio}
\]
The recurrence is valid for every integral defect allowed by
the rank and signature conditions. If \(e=-j<0\), then
\[
 \frac{M_\nu(e+1)}{M_\nu(e)}
   =\frac{M_\nu(j)}{M_\nu(j-1)}
   =L_\nu^\sharp(j)^2.
\]
Thus the negative case uses the positive Mellin exponent
\(j=-e\), including the later application with \(e=-k\) at
\((J,0,J)\). Both positive Mellin exponents are at least one,
and \(\nu\) is nontrivial. Thus the factor
\(\gamma_{\rm pl,\nu}\) of \eqref{eq:plane-period-local-norm}
is nonzero. The two scalars \(c_{\rm old,\nu}\) and
\(c_{\rm new,\nu}\) can vanish until the choice of the equal-rank
vector below.

At the top the factor is \(D_{\mathrm{top}}=
(R_0|\kappa_0|)^{-2}\) and the value of this invariant
is
\[
          R_{W_0,\nu}(z_0)
           =|I_{W_0,\nu}(z_0)|^2/M_\nu(j_0).
 \tag{5}\label{eq:period-invariant-value}
\]
The small target \(U^-\) stays fixed. After \(i\) removals,
\(0\le i\le b+1\), the source and large target ranks are
\(N_i=N+2(b+1-i)\) and \(a_i=a+2(b+1-i)\).
The large target has \(q(V)+b+1-i\) negative signs, since
\(q(U^-)=b-q(U)\). Thus its defect is
\(e_i=a_i-N_i/2=j_0-i\), from \(e_0=j_0\) to \(e_{b+1}=t\).
Writing \(C_i,\kappa_{A,i}\) for the constants of step \(i\),
the telescoped identity is
\[
 |c_{\mathrm{final},\nu}|^2
   =|c_{\mathrm{top},\nu}|^2
      \frac{M_\nu(t)}{M_\nu(j_0)}
      \prod_{i=1}^{b+1}\frac{|\kappa_{A,i}|^2}{C_i^2}.
\]
It applies also when \(b=0\), with one removal, and when
\(a=0\), with final defect \(-N/2\).
Successive source zero slices compose, since they are
evaluations of separate blocks in the product mixed model.
Let \(y^{\rm ad}\) be the image of the prescribed vector \(y\) for
\((U,W)\) under the signed real adjoint, so its pair is \((W,U^-)\).
This coordinate change preserves each coordinate modulus up to
permutation and preserves the Gaussian and Lebesgue measure. The
weighted Fock moments of \(y^{\rm ad}\) are finite, with the
coordinatewise moment indices permuted. In the one-plane
mixed model on \((W_0,U^-)\), extend
\(y^{\rm ad}\) by a Schwartz tensor of value one at zero to obtain
\(v\). Use this \(v\) in \eqref{eq:marked-real-slice-small}, whose
opposite zero slice after the other \(b\) planes is \(v\) by its
integral normalization. The combined zero slice is consequently
\(y^{\rm ad}\). On the large side take the successive cone slices of
\eqref{eq:marked-real-slice-big} at their prescribed exponents
\(s(e)\). Their real integrals use only the fixed infinity type
\(\mu\), the splittings and the two initial real vectors. Thus, for a
fixed \(z_0\), they give one large vector before the modulus and the
finite matched arrays are chosen. After \(b+1\) removals
\eqref{eq:plane-removal-ratio} holds at \((W,V,U^-)\) with defect
\(t\) and small vector \(y^{\rm ad}\); adjoint transport returns it
to the prescribed \(y\) below.

\medskip
\noindent\emph{From standard data to the actual finite data.}
In the product \(F_S=\prod_{w\in S}E_w\), put
\(O_S=\prod_{w\in S}\mathcal O_w\). Let \(p\) be an
integral block of size \(b\) and use the two standard
data \(\chi^{-1}(\det[p;X])\) and
\(\chi^{-1}(\det(pY))\), each supported on integral
arguments with unit determinant in the indicated formula.
For \(X\) the last covectors of \(k^{-1}\) and \(Y\) the
first columns of \(k\), one has
\[
                   \det[p;X]\det k=\det(pY).
\]
Thus these are matched, with factor \(\chi(\det k)\).

For use with the equivariant difference of the two theta
periods, we prove finite generation.
A locally constant matched column function \(\Psi\) is
supported on primitive integral columns. It satisfies
\(\Psi(Yu)=\chi(\det u)^{-1}\Psi(Y)\) for
\(u\in\mathrm{GL}_b(O_S)\): replace a completion \(k\)
by \(k\,\mathrm{diag}(u,1)\).
Choose \(x_0\) primitive and an integral \(p\) with
\(px_0=1\). Let \(F(Y)=\chi^{-1}(\det(pY))\) on its
standard support. Take a common depth \(c\ge1\) of local
constancy of \(\Psi\), and the product uniformizer
\(\varpi=(\varpi_w)_{w\in S}\).

Set
\[
 d=x_0p+\varpi^{-c}(1-x_0p),\qquad K=\mathrm{GL}_N(O_S).
\]
If \(F(dkY)\ne0\), put \(Z=dkY\) and \(u=pZ\).
Then \(u\) is invertible and
\[
 Y=k^{-1}x_0u+\varpi^c k^{-1}(Z-x_0u),\qquad
 \Psi(Y)=\Psi(k^{-1}x_0)F(dkY).
\]
Conversely any nonzero value of \(F(dkY)\) forces \(Y\)
primitive integral. Thus
\begin{equation}\label{eq:period-original-compact-average}
 \int_K\Psi(k^{-1}x_0)F(dkY)\,dk=c_F\Psi(Y)
\end{equation}
everywhere, where \(c_F=\int_K|F(dkx_0)|\,dk>0\).
For a datum \(m(Z)F(Z)\) with nonnegative locally constant
multiplicity \(m\), the same calculation replaces \(c_F\)
by \(\int_K m(dkx_0)|F(dkx_0)|\,dk\), which is positive when
\(m(x_0)>0\).
To see constancy of \(c_F\) on primitive \(Y\), use the
transitivity of \(K\) on primitive columns. Positivity
follows on a neighborhood of \(k=1\).

All functions of \(k\) in the average are locally constant.
The translated \(F\)'s themselves range over finitely
many functions on \(K\), by smoothness of the translation
action on locally constant compactly supported functions.
Consequently this is a finite linear combination of translates
of \(F\), or of the function with multiplicity in the preceding
formula. For the detector \(G=B(-\varpi^{-1}\,\cdot)\) the average
in fact factors through \(\mathrm{GL}_N(O_S/\varpi^{c+2}O_S)\).
If \(k_1\equiv k_2\bmod\varpi^{c+2}\) and \(Y\) is integral,
then \(dk_1Y-dk_2Y\in\varpi^2 M_{N,b}(O_S)\). Both evaluations
vanish on nonintegral \(Y\), since \(d^{-1}\) is integral.
The factor \(\Psi(k^{-1}x_0)\) is constant on the same cosets.
If this finite quotient has order \(h\), the scalar multiplying
\(\Psi(Y)\) is \(h^{-1}\sum_{\bar k}|G(dkx_0)|\ge h^{-1}\):
the identity coset contributes \(|G(x_0)|=1\). The split covariance of the difference \(\mathcal D\)
takes the form
\[
 \mathcal D(\Psi\circ g_\Sigma)(xg_\Sigma)
    =\lambda_b(\det g_\Sigma)^{-1}\mathcal D(\Psi)(x),
 \qquad \lambda_b(u)=(\upsilon_E(u)|u|^{1/2})^b.
\]
Its kernel is therefore stable under the translates.
The common-section calculation in
Section~\ref{sec:marked-primitive-fourier-transfer} identifies
the domain with the original integral matching and its
split translates. This proves the comparison for every
locally constant matched pair.

For a modulus at \(S\), the reduction
\(O_S\to\mathcal O_E/\mathfrak m\) is surjective and detects
units by the Chinese remainder theorem. Each invertible
residue matrix lifts to an invertible integral matrix by
lifting its entries and using the determinant test. A
primitive integral covector or column block extends to such a
matrix by completion over each local ring. The residue
matching law and its support thus lift exactly, by pullback
and extension by zero, to the same law and support on
\(O_S\). The unit character is the same product of nontrivial
residue characters on the selected support \(S\).
More explicitly, write \(\mathfrak m=\prod_{w\in S}\mathfrak p_w^{d_w}\)
with arbitrary \(d_w\ge1\) and \(S\cap\bar S=\varnothing\). The
character on \((\mathcal O_E/\mathfrak m)^\times\) is the product of
the tame characters on \(k_w^\times\), pulled back by reduction.
The two matched functions may depend on all residue coordinates
modulo \(\mathfrak p_w^{d_w}\). In particular the comparison on
columns determines every original covector value over the finite
ring. For a primitive row matrix \(X\), take a completion \(k\)
with \(X\) the last rows of \(k^{-1}\), and let \(Y\) be the
first columns of \(k\). Another completion is
\(k\left(\begin{smallmatrix}u&v\\0&1_a\end{smallmatrix}\right)\)
with \(u\) invertible. Its determinant factor \(\chi(\det u)\)
cancels the factor \(\chi(\det u)^{-1}\) in \(\Psi(Yu)\).
Thus \(\Phi(X)=\chi(\det k)\Psi(Y)\) is independent of the
completion. Both functions vanish off the primitive locus.
Their pullbacks are locally constant, so the finite spanning
formula applies for every \(d_w\ge1\), including the old exponent
\(b_0\) in Section~\ref{sec:deep-transport-arrays}.
Outside \(S\cup\bar S\) the prescribed vacua are identical.
In particular the change from the finite ring to the local
integral model introduces no level indices or measure ratios.

\medskip
\subsubsection{The equal-rank normalizer}\label{sec:marked-equal-rank-vanishing}

\sectionopening{In this subsection, we use the invariant period ratio \eqref{eq:plane-removal-ratio} to normalize the comparison in Section~\ref{sec:gluing-three-exchanges}.}

Write \(N=2k+2q(W)\), so that \(\mu_\tau\ge k\ge1\).
First take a positive definite space \(J\) of rank \(2k\) and
\(Y=H\perp J\). For each \(\nu\) the equal-rank period
\(I_{Y,\nu}(z)\) with determinant datum is constant in its
unitary variable by the special-unitary invariance and the one-plane
Levi calculation of Section~\ref{sec:marked-real-slice-period}.

Fix a split tame character of the indicated weight.
Suppose that every equal-rank period of \((Y,Y)\) with
the finite determinant datum is zero.
In the large space \(X=H^{2k}\perp Y\), with target
\(J^-\perp Y\), consider an arbitrary real Schwartz
vector. First use the self-dual Fourier test with rational Hermitian
index \(J_{J^-}\) and its prescribed adelic basis. The frame
\(\binom{1}{0}\) has orthogonal complement \(Y\), and its displacement
\(\binom{d_c}{0}\) has zero \(Y\)-component. On the supported finite
cosets the upper triangular determinant in
Section~\ref{sec:marked-mixed-absolute-summability} leaves precisely
the original square determinant datum on \(Y\), up to the scalar
from the frame block. A rational congruence transports the index
and the coefficient point together as in
Section~\ref{sec:selfdual-fourier-integral-extension}.

Fix an outer point of the unfolded integral and denote the moved
real mixed slice by \(w\).
Use the weighted Fock coordinates of
Section~\ref{sec:four-unitary-constant-terms}. For a block of complex
dimension \(d_i\), write \(du_i=\pi^{-d_i}d\operatorname{Leb}(u_i)\).
The Gaussian is included in the values of \(w\), and the unit
coherent vector is
\(v_u(z)=\exp(-(|z|^2+|u|^2)/2+\bar u\cdot z)\).
With the inner product linear in the second argument, the
reproducing formula in the two blocks gives
\[
 w(u_1,u_2)=\langle v_{u_1}\boxtimes v_{u_2},w\rangle,
 \qquad
 w=\int w(u_1,u_2)\,v_{u_1}\boxtimes v_{u_2}\,du_1\,du_2 .
\]
In the real Schr\"odinger model each \(v_u\) is a translate and
modulation of the fixed Gaussian, with parameters linear in the real
coordinates of \(u\). Translation of the weights and differentiation
of the modulation give \(p_j(v_u)\ll_j(1+\|u\|)^j\) for the
seminorms of Section~\ref{sec:marked-mixed-absolute-summability}.
With the first finite vector and test point fixed, the finite-order
theta bound therefore gives polynomial growth of \(T(v_{u_1})\).
All weighted absolute moments of the Schwartz slice \(w\) are finite,
so the reproducing integral and the contraction below converge in
every Schwartz seminorm. In particular, the contraction
\[
 w_2=\int w(u_1,u_2)T(v_{u_1})v_{u_2}\,d u_1\,d u_2
\]
is a real Fock vector with all moments finite.
Write \(T_{Y,h}(v_{u_2})\) for the second theta evaluation at
\(h\in[Y]\) with its finite determinant datum. At any chosen
\(w\in S\ne\varnothing\) its support consists of invertible square
matrices, so it has full row rank in the orientation for integration
over \([Y]\). The ordinary period bound
\eqref{eq:partner-exchange-ordinary-bound}, with the source point
fixed and height exponent zero, bounds \(T_{Y,h}(v_{u_2})\)
uniformly in \(h\) by a fixed Schwartz seminorm of \(v_{u_2}\). Translation and modulation of the Gaussian bound this
seminorm by a polynomial in \(\|u_2\|\). Thus for some \(C,M\), with
the first test and point fixed,
\[
 \begin{aligned}
 &\int_{[Y]}\int
   |w(u_1,u_2)T(v_{u_1})T_{Y,h}(v_{u_2})|\,du_1\,du_2\,dh\\
 &\quad\le C\,\operatorname{vol}([Y])
       \int |w(u_1,u_2)|(1+\|u_1\|+\|u_2\|)^M\,du_1\,du_2
 <\infty .
 \end{aligned}
\]
The tensor kernel product and Fubini identify the inner period
with the assumed zero period of \(w_2\). For an arbitrary upper
vector the original unfolded term is absolutely integrable by
Section~\ref{sec:marked-mixed-absolute-summability}, so its
disintegration over the open orbit integrates this zero inner
period.

For the primitive Fourier indices of
Section~\ref{sec:marked-primitive-fourier-transfer}, the Gram matrix is
invertible over \(E\), since their scaled reductions are invertible at
the chosen places. At the points used in the finite detection argument,
the coefficient of the detector at any remaining \(T\) is the coefficient
of \(\Psi_0\) times the matrix Gauss sum. If \(\Psi_0\) has a nonzero
coefficient, the scaled index is integral at every selected place.
At a singular reduction the change \(r\mapsto ur\) in that local
sum, with \(\bar A u=\bar A\) and \(\chi_w(\det u)\ne1\), gives zero.
This applies in particular to every singular rational index. Hence
only represented primitive indices require the orbit calculation.
For such an index outside the negative real class, choose a real
place at which its nondegenerate span has \(s_\tau\ge1\) positive
directions. Put this span in standard signed position by a unitary
change of coordinates. In the smaller mixed Fock space the source
\(Y=H\perp J\) has \((p_Y,q_Y)=(2k+1,1)\).
The target has \(b+q_Y\) negative coordinates, where \(b=2k\).
The span uses \(b-s_\tau\) of them, so the compact group on the
\(q_Y+s_\tau\) unused negative coordinates fixes the entire frame.
The real determinant exponent after the source slice is
\(m_\tau=(\dim Y)/2-\mu_\tau=k+1-\mu_\tau\).

Fix a point of the signed polynomial Fock coordinates, and let
\(Z_-\) be its \(q_Y\) negative-source evaluation columns in the
unused negative coordinates. Choose a unit column \(r\) there with
\(r^{\mathsf t}\overline{Z_-}=0\). These are at most \(q_Y\)
linear conditions in dimension \(q_Y+s_\tau>q_Y\).
The circle \(k_\eta=I+(\bar\eta-1)\bar r r^{\mathsf t}\), \(|\eta|=1\),
lies in the same fixed compact stabilizer. Its determinant is
\(\bar\eta\), and the negative-block vacuum multiplier and the twist
give
\[
 (\det k_\eta)^{-p_Y+m_\tau}
       =\eta^{p_Y-m_\tau}=\eta^{k+\mu_\tau}.
\] On the signed polynomial coordinates with negative target sign,
substitute \(k_\eta^{\mathsf t}\) on the negative-source columns
and \(k_\eta^{-1}\) on the positive-source columns.
The first substitution fixes \(Z_-\), since
\(\bar r^{\mathsf t}Z_-=0\). The second multiplies the line
component along \(\bar r\) by \(\eta\) and fixes its complement. If a polynomial vector
is invariant under the stabilizer's compact
subgroup, its evaluation at this point is
therefore the mean on this circle of a
polynomial in \(\eta\) all of whose exponents
are at least \(p_Y-m_\tau>0\).
This mean is zero. As the point was arbitrary,
there is no invariant polynomial vector.
The compact group preserves each homogeneous polynomial space.
Its probability average is a bounded projection on the Fock Hilbert
space, so polynomial density makes this projection zero. For a real
Schwartz vector the compact orbit map is continuous in the Schwartz
topology. Its Haar integral there has the same Hilbert image, hence is
zero in Schwartz space. The continuous mixed theta functional may
therefore be applied to this average.

Unfold the Fourier term over the rational basis
stabilizer as in the self-dual calculation.
Disintegrate through the adelic stabilizer and
average over the compact subgroup just used,
conjugated back to the given basis.
Right multiplication by the compact subgroup
preserves the stabilizer Haar, since the modular
character is trivial on compact groups.
The mixed-model slice commutes with this
stabilizer action at the chosen vector, with
the lower source rank twist. Thus its compact
average vanishes by the Fock calculation.
The coefficient therefore vanishes off the negative
class for an arbitrary upper real vector. If no
rational basis realizes the index, it is an empty sum.

Apply Section~\ref{sec:selfdual-fourier-integral-extension} to
the large period and the zero functional on the small side.
The two cases above give vanishing for every self-dual test.
The primitive coefficient transfer then gives vanishing of the
large period at every point for every upper real vector and
standard covector datum.
To descend this vanishing, take any proposed lower
period \(I\) and the common real bump section on
the cone, whose scalar is \(B(s)>0\) for every real \(s\).
For this lift of the given lower vector, with its fixed covectors
and source point, write \(Z(s)\) for the grouped Mellin integral of
Section~\ref{sec:radical-isotropic-line-orbit}, absorbing the positive
unfolding constant. The positive-sign unfolding gives
\[
       Z(s)=d_\nu|\alpha_0|^s Q_S^{s-1/2}
                     L^S(s,\nu)B(s) I
       \qquad (\Re s>1),\qquad d_\nu\ne0 .
\]
The positive and Poisson-dual estimates there give an integrable
majorant for the grouped integral on each closed vertical strip.
Thus \(Z(s)\) is entire. Its value at \(e+1\) is the constant
term of the upper period on this lifted datum.

The right side is holomorphic on \(\Re s>-k\) by
Tate's continuation theorem \cite[Theorem~4.4.1 and \S4.5]{TateThesis}.
The nonzero angular weights exclude a norm character. Holomorphic
continuation gives the identity throughout this half-plane. At each needed integer \(s=e+1>-k\) this
Hecke factor is nonzero. For \(s>1\) the Euler product converges
absolutely. At \(s=1\) use the nonvanishing argument with
\(3+4\cos\theta+\cos(2\theta)\) above. For \(s\le0\) the
infinity factors are \(\Gamma_{\mathbb C}(s+\mu_\tau)\)
and \(\Gamma_{\mathbb C}(1-s+\mu_\tau)\) at the two
reflected points. Their arguments are positive,
and the Hecke functional equation reduces the assertion
to nonvanishing at \(1-s\ge1\) for the dual character.

The Euler factors removed at \(S\cup\bar S\)
are already the ramified factors in that equation.
By the upper vanishing the left side is zero at
\(s=e+1\), so \(I=0\).
For these removals, use the rational hyperbolic splitting
of \(J^-\perp J\) given by its diagonal and antidiagonal.
In the original finite orthonormal bases of the two copies
the Gram matrix is \(\operatorname{diag}(I,I)\).
At each finite place choose integral \(u_-\) and \(c\) with
\(u_-\bar u_-=-1\) and \(c+\bar c=1\).
The unit norm and integral trace are surjective in the split
and unramified quadratic algebras, also at residue characteristic two.
The two isotropic frames in the original coordinates have columns
\[
       (u_-u,u),\qquad (-u_-cv,(1-c)v),
\]
for the standard basis vectors \(u,v\) of the two coordinate
spaces. Their cross pairing is \(\bar u^{\mathsf t}v\).
The block matrix with these columns has determinant \(u_-^{2k}\),
an integral unit. It therefore preserves the finite lattice
and its additive measure. At every split place fix
\(u_-=(1,-1)\) and \(c=(1,0)\) in the product algebra,
independently of the selected set \(S\).
Then \(u_-^{2k}=(1,1)\), so the finite determinant
character of this basis change has value one on \(S\).
Using these rational and finite bases, apply the preceding
descent at the successive \(2k+1\) planes from
\((H^{2k}\perp Y,J^-\perp Y)\) to \((J,0)\).
After removal \(1\le i\le2k+1\) the ranks are \(6k+2-2i\) and
\(4k+2-2i\). Thus the positive-model continuation uses
\(s_i=e_i+1=k+2-i\), with \(e_i=k,k-1,\ldots,-k\).
The smallest argument is \(1-k>-k\). For \(s_i\le0\)
the two gamma arguments satisfy \(s_i+\mu_\tau\ge1\)
and \(1-s_i+\mu_\tau\ge1+\mu_\tau\).
At the last removal the positive-model exponent is \(1-k\) and
the dual exponent is \(k\). The conductor of \(\nu^{\rm bc}\) has
exponent one at both members of each pair in the chosen support and
is trivial elsewhere, so its norm is \(Q_S^2\), independently of
the depths of the residue tests. Taking absolute values in the
completed functional equation gives
\[
 Q_S^{1/2-k}|L^S(1-k,\nu)|
 =|\operatorname{disc}(E)|^{k-1/2}
  \prod_{\tau\mid\infty}
   \frac{\Gamma_{\mathbb C}(k+\mu_\tau)}
        {\Gamma_{\mathbb C}(1-k+\mu_\tau)}
     Q_S^{k-1/2}|L^S(k,\nu)|.
\]
Here the dual character has the conjugate finite \(L\)-value at the
real argument \(k\). Both gamma factors are regular, including at
\(k=1\), where the dual value is taken at one.
In the final test take the real scalar vector \(1\), source point
\(g=1\), and fixed covector matrix \(p=I_{2k}\) at the selected
places. The moving matrix has zero rows, so the stacked determinant
is \(\det p=1\). Every finite factor equals one and the rational
theta sum on the zero Hom space evaluates its scalar at the sole
point. The period therefore has absolute value
\(v_0=\operatorname{vol}([U(0)])>0\), also after transport by a
unit phase. This contradiction proves
the required nonvanishing of some \(I_{Y,\nu}(z)\).

If the family with the fixed infinity type is nonempty, choose one
reference character \(\nu_*\) in it, at any admissible modulus.
Fix a single real vector \(z\) on \((Y,Y)\) with
\(R_{Y,\nu_*}(z)>0\), as just proved. The standard exchange
constructed by \eqref{eq:period-open-orbit-scalar}--\eqref{eq:plane-removal-ratio}, allowing its scalar to vanish,
applies at \((J,J,0)\) and \((J,0,J)\), with defects \(k\)
and \(-k\), respectively. In the reverse exchange the upper target
is \(J^-\perp H\perp J\). Use the rational and finite hyperbolic bases above.
The \(2k+1\) removals end at the zero target also at dyadic
places. The split basis choice above applies to the same
residue-character test in the reverse comparison. The last step has \(e=-k\) and uses
\(L_\nu^\sharp(k)\) in \eqref{eq:plane-period-local-norm}.
Use the real image under the first map as the argument of the
second, before varying \(\nu\).
In the forward exchange \(b=0\) and \(N=2k\). The empty Gram
matrix has determinant one. The two frame shifts and the central
Fourier group have dimension zero, so the sole Fourier coefficient
is the period itself and the real frame is \(B\in\mathrm{GL}_0\).
The quotient by the full block stabilizer is a point of mass
\(v_0^{-1}\), cancelling the inner \(U(0)\) mass in the
open-orbit disintegration. For \(a=0\) the top Fourier index has
size \(b=N\), and the target first has rank zero at the final removal.
Write
\[
 D(X)=\mathbf1_{\mathrm{GL}_N(O_S)}(X)\chi^{-1}(\det X)
\]
on integral square matrices, extending the character by zero on
nonunits. The forward law \(\Phi(k^{-1})=\chi(\det k)\Psi(0)\)
gives \((\Phi,\Psi)=(D,1)\) when \(\Psi(0)=1\).
In the reverse direction take \(p=I_N\). The law
\(\Phi(0)=\chi(\det k)\Psi(k)\) gives \((\Phi,\Psi)=(1,D)\)
when \(\Phi(0)=1\). In each case every integral matched
pair at the prescribed depth is a scalar multiple, by this law
and support on the primitive locus. Both tests have the original
vacuum outside the chosen pairs. The adjoint change of variables
below uses an isometry of the original quotients and multiplication
by a unit scalar. It applies to the preliminary period identity
at either sign of the defect, including when the scalar of that
identity is zero. With \((\alpha,\alpha')=(1,\gamma)\), write
\(K_\nu=K_{J,J}^{\gamma,\gamma}\) for the kernel on the fixed
real image and the square determinant test \(D\). In the forward
exchange, the source character gives
\[
 K_{J,J}^{\gamma,1}(g,h)=\gamma(\det g)K_\nu(g,h).
\]
The zero-source integral is \(\gamma(\det g)v_0\). Cancelling
this character uses no nonvanishing assertion about \(c_+\).
In the reverse exchange the small kernel is \(K_\nu(h,g)\);
the rank-zero target kernel with characters \((\gamma,1)\)
is the fixed scalar coordinate \(\zeta_0\).
Let \(B,C>0\) denote the final normalization factors
of \eqref{eq:plane-removal-ratio} in these two exchanges.
With the positive zero-group quotient volume \(v_0\),
the identities read
\begin{align*}
 \int_{[J]}K_\nu(g,h)\,dh&=c_+v_0,
 & B|c_+|^2&=M_\nu(k)R_{Y,\nu}(z),\\
 c_-\int_{[J]}K_\nu(h,g)\,dh&=v_0\zeta_0,
 & C|c_-|^2&=M_\nu(-k)R_{Y,\nu}(z).
\end{align*}
The product \([J]^2\) is compact because \(J\) is definite.
On compact sets of representatives the finite vectors range over a
finite set and the real Schwartz seminorms are uniformly bounded.
The theta kernel is therefore bounded and continuous on the product.
Set \(v_J=\operatorname{vol}([J])\). Integrating the first
identity in \(g\), then using the second identity with the order of
integration reversed, gives
\[
 c_-c_+v_0v_J
   =c_-\int_{[J]^2}K_\nu(g,h)\,dg\,dh
   =v_0\zeta_0 v_J .
\]
Thus \(c_+c_-=\zeta_0\). The scalar \(\zeta_0\) is the value of the
real vector produced by the reverse exchange on the zero Hom space;
the point-quotient measure has been accounted for by \(v_0\).
Multiplication, using \(M_\nu(k)M_\nu(-k)=1\), gives
\[
                     R_{Y,\nu}(z)^2=BC|\zeta_0|^2.
\]
At \(\nu_*\) the right side is strictly positive. Since
\(R_{Y,\nu}(z)=|I_{Y,\nu}(z)|^2/M_\nu(k+1)\ge0\) for every
admissible character, the equality of squares gives
\(R_{Y,\nu}(z)=R_{Y,\nu_*}(z)>0\) for all moduli and characters
in the family.

An equal-rank step from rank \(2j\) to rank \(2j+2\)
uses \eqref{eq:plane-removal-constant-terms} on the standard equal-rank determinant datum.
The larger period is constant, hence equals its normalized
constant term. For every given lower vector a real
isotropic-line bump supplies the same upper cone section
before \(\nu\), normalized by its positive real integral
at \(s=j+1\). Equations \eqref{eq:plane-period-local-norm}--\eqref{eq:plane-removal-ratio} say that the ratio at
the upper vector is a fixed positive multiple of the
lower ratio.

Iterate through \(q(W)\) planes. Identify the
resulting space with \(W_0\) by the rational
Hermitian isometry and changes of adelic basis in the
model. The rational group, adelic compact
normalization and quotient Haar are transported, and the
constant period changes at most by a phase. The rational
isometry follows from Landherr's theorem
\cite[pp.~246--247]{LandherrHermitian}, since the signatures and the
finite local forms agree. This gives a single \(z_0\) for
which \eqref{eq:period-invariant-value} is a fixed positive \(R_0'\).
If the admissible character family is empty the assertion
has no instances and a vacuum can be chosen.

\medskip
\noindent\emph{Return to the original small quotient.}
The standard small comparison uses \((W,U^-)\).
Write \(X_U,X_W\) for the chosen adelic bases. Fix a finite integral
unit \(c_f\) with \(c_f\bar c_f=-1\). Such units exist at every
finite place since \(E/F_0\) is split or unramified there and the
norm on units is onto. For the opposite form \(-J_U\) take
\[
 X_{U^-}=dX_U,\qquad d_\infty=1,\qquad
 d_f=\bar c_f\bar\alpha_0/\alpha_0 .
\]
The finite Gram matrix of this basis is the identity because
\(d_f\bar d_f=-1\). For a rational map \(A:U\to W\) put
\[
 T=X_W^{-1}AX_U,\qquad
 T'=X_{U^-}^{-1}(J_U^{-1}A^*J_W)X_W .
\]
At infinity \(T'\) is the signed adjoint of \(T\). At finite places
the two original Gram matrices in the chosen bases are identities,
so \(T'_f=d_f^{-1}T_f^*\). In the scaled finite coordinates this is
\[
 T_f/\alpha_0=c_f(T'_f/\alpha_0)^* .
\]
The rational adjoint is bijective, with inverse
\(B\mapsto J_W^{-1}B^*J_U\).
The opposite Hermitian sign cancels the sign from
adjunction in the Heisenberg form. The adjoint preserves
the finite vacuum lattice and the real Gaussian, so the
two normalized theta distributions differ by a phase.
Indeed, in a rational Schr\"odinger polarization, rational
modulation and translation characterize the distribution
as a multiple of the sum over the rational points.
The normalized mean square on vacuum translates fixes the
absolute value of that multiple.

For characters \((a_1,a_2)\) on \((U,W)\), the flipped
characters on \((W,U^-)\) are
\[
       (a_2^{-1}\upsilon^N,\ a_1^{-1}\upsilon^b).
\]
With the ranks interchanged, applying this character change a
second time gives
\[
 ((a_1^{-1}\upsilon^b)^{-1}\upsilon^b,\,
  (a_2^{-1}\upsilon^N)^{-1}\upsilon^N)=(a_1,a_2).
\]
The real complementary-minor identity and the finite
elementary divisors give the multiplier ratio
\[
 \frac{\upsilon(\det h)^N}{\upsilon(\det g)^b},
 \qquad (h,g)\in U(U)\times U(W).
\]
The absolute determinant powers cancel because the
unitary determinants have norm one. For any finite scaled matrix
\(Z\) of size \(b\)-by-\(N\), let
\(Z_S^\dagger=(\overline{Z_{\bar w}}^{\mathsf t})_{w\in S}\), where
conjugation identifies \(E_{\bar w}\) with \(E_w\). The right
covariance of the matched column datum gives the pointwise identity
\[
 \Psi(Z_S^\dagger)=\chi((c_f)_S)^b\,\Psi((c_f)_S Z_S^\dagger).
\]
It holds on the full matrix space since the unit scalar preserves
the support. Thus the flipped test reads the transpose at the
conjugate places and its phase is independent of the matched functions.

The forms \(J_U,-J_U\) have the same rational unitary matrices.
The basis change \(X_{U^-}=dX_U\) is central, so conjugating those
matrices into the two adelic models gives the same arithmetic
subgroup. The identity on adelic unitary matrices preserves the
chosen compact group and its Haar normalization, hence also the
quotient measure. Thus the original small integral is
related by a unit scalar and no volume ratio. More
precisely, starting with \((\gamma,\gamma)\) on
\((U,W)\), the standard small period obeys
\[
 \gamma(\det g)\Theta_{\rm small}(g)
   =\varepsilon_{\rm ad}\int_{[U]}
        K_{U,W}^{\gamma,\gamma}(\Psi,y;h,g)\,dh,
      \qquad |\varepsilon_{\rm ad}|=1 .
\]
On the standard large side, multiplication by
\(\gamma(\det g)\) changes the characters from
\((\gamma,\gamma)\) to \((\gamma,1)\).

Absorb \(\varepsilon_{\rm ad}\) into the scalar in \eqref{eq:plane-removal-ratio}.
Since \eqref{eq:period-invariant-value} is \(R'_0>0\), the resulting scalar \(c_\nu\)
is nonzero, with
\[
                   |c_\nu|^2=(R'_0/D_{\rm final})M_\nu(t).
\]
Finally, the product of the two rank-specific space
characters equals \((\nu^{\rm bc})^{-1}\) precisely when
\(\alpha\alpha'=\gamma\). This follows by multiplying
the two powers of the splitting character and using
the surjectivity of \(x\mapsto x/\bar x\) on the norm-one
ideles, by local Hilbert 90, choosing units at the unramified
finite places.

Replacing \((1,\gamma)\) by
\((\alpha,\alpha')\) on the two middle groups multiplies
both kernels by \(\alpha(\det g)^{-1}\).
\begin{proposition}[Normalized period exchange]\label{prop:normalized-period-exchange}
Suppose \(t\ge0\). For the fixed spaces \(W,V,U\), Haar measures,
infinity weight \(\mu\) and real Schwartz vector \(y\), one can choose
a real Schwartz vector \(x\) for \((W,V)\) and
\(R'_0/D_{\rm final}>0\) before the split tame modulus and
character \(\nu\).
For every such modulus \(\mathfrak m\) with nonempty chosen split
support \(S\), and every admissible \(\nu\) whose tame character
on \(k_w^\times\) is nontrivial for every \(w\in S\), there is
a scalar \(c_\nu\ne0\) such that,
for all \(\alpha,\alpha'\) with \(\alpha\alpha'=\gamma\), all
matched arrays \((\Phi,\Psi)\) over \(\mathcal O_E/\mathfrak m\)
for its induced residue character, and every adelic \(g\), both
periods are absolutely convergent and
\[
\begin{aligned}
 \int_{[V]}K_{W,V}^{\gamma,\alpha}(\Phi,x;g,h)\,dh
   &=c_\nu\int_{[U]}K_{U,W}^{\alpha',\gamma}(\Psi,y;h,g)\,dh,\\
 |c_\nu|^2&=(R'_0/D_{\rm final})M_\nu(t).
\end{aligned}
\]
\end{proposition}

For the two periods in the proposition, the chosen tests have
full row rank on \(\Phi\)'s support and full column rank on
\(\Psi\)'s support at every \(w\in S\). Apply
\eqref{eq:partner-exchange-ordinary-bound} to the covector test
and to the column test in the conjugate transpose model. At
rank zero the integral is over the one-point quotient.
The stronger mixed-slice estimate is used for the top central
Fourier coefficients, with their adapted frames. For the two
one-plane constant terms the ordinary estimate precedes the
compact radical average, and on the large side the sum over
rational lines is unfolded with the grouped estimate in
Section~\ref{sec:radical-isotropic-line-orbit}.

\subsection{The finite extension of the period identity}
\label{sec:marked-square-finite-extension}

\sectionopening{In this subsection, after the analytic exchange of
Proposition~\ref{prop:normalized-period-exchange}, finite Fourier inversion
gives an identity on the entire matrix space. We determine the
resulting next datum, including the inverse determinant character and
its scalar.}

Let \(R=\mathcal O_E/(\pi_MQ(\xi))\), a product of finite
fields at the chosen split members, and let \(\chi\)
be the character of \(R^\times\) of the quotient
character. Extend multiplicative characters by zero
on nonunits. A paired residue datum on \(a\)-by-\((b+a)\)
matrices \(X\) and \((b+a)\)-by-\(b\) matrices \(Y\) is
supported, respectively, on the last \(a\) covectors of
\(k^{-1}\) and the first \(b\) columns of \(k\) for
\(k\in\mathrm{GL}_{a+b}(R)\), and satisfies
\[
                     \Phi(X)=\chi(\det k)\Psi(Y)
\]
on such pairs. The order of the covectors in this definition is
\(b,a\). The matched point data from the
measured square is the following. The scalar tested
composites \(P(X),P'(Y)\) on the prescribed four Hermitian spaces
satisfy
\[
 \begin{split}
 P(gX)&=\chi(\det g)P(X),\qquad
 P'(Yh)=\chi(\det h)P'(Y),\\
 \sum_Y\Psi(Y)P'(Y)&=c\sum_X\Phi(X)P(X)
                 \qquad\text{for paired data }(\Phi,\Psi).
 \end{split}
\]
Here the test functional in the top covector is arbitrary;
the scalar \(c\) is independent of it and of the
residue arrays. The first line is the finite
determinant covariance. It holds on all residue arrays:
the fixed-point tests define individual lifts against the rapid source
forms. The coefficient projector argument of
Section~\ref{sec:four-unitary-real-projections} applies to each fixed
finite datum by the Gaussian theta bound and gives square integrability
of the resulting lifts on the middle group. Changing the middle
basis in that pairing preserves its quotient measure and gives the
displayed determinant multiplier. The second line uses the measured
identity only on the full-rank supports of paired data.

For the triangle at the original level take \(i,n\ge0\) and
\(i+n+2\le m\). In column \(n\) the rank is \(r+nt\). The two kept spaces at \((i,n+1)\) and
\((i+1,n+1)\) have rank \(r+(n+1)t\); the two
integrated spaces at \((i+1,n)\) and \((i,n+2)\)
differ by rank \(2t\). The negative counts on the
two diagonal pairs agree. For a stem entry with
excess parameter \(a_*\), its sign in covector \(o\) and
column \(j\) is
\[
                 (-1)^{o+a_*+\min(j,(a_*-o)_+)}.
\]
If \(a_*\le i+n+1\), the two column changes leave
both signs unchanged.

If \(a_*\ge i+n+2\), both
signs change and they were opposite. All added
entries are positive. The common negative count is
at most \(2r<r+(n+1)t\). Moreover the positive
weight of the inverse quotient is
\((n+1)\Delta\), greater than \(r+(n+1)t\).
These are the rank, sign and weight inequalities for
the exchange. In the reversed square interchange the
two roles of the diagonal pairs. The two integrated
ranks now agree, so the rank defect is zero.

Write \(u=\upsilon|_{U(1)}\) and
\(\mu_j=\mu_0(u^{-\Delta}\chi_M)^j\). The horizontal
and left-pair determinant twists at ranks \(p,q\)
in the square respectively have factors
\[
  \mu_{i+n+2}\xi^{i+n+2}u^{\lfloor p/2\rfloor},
  \qquad
  \mu_{i+1}^{-1}\xi^{-(i+1)}u^{\lceil q/2\rceil}.
\]
The second is the inverse of the first-column
character at rank \(q\), multiplied by \(u^q\).
Since \(p+q=2(r+(n+1)t)\), its product with the
first is the rank twist of
\((\mu_{i+n+2}/\mu_{i+1})\xi^{n+1}\).
In the split Weil characters the two rank
corrections therefore cancel, leaving the base
change of that quotient. The paired finite
determinant character is its inverse.
The two lift indices satisfy \(1\le i+1<i+n+2\le m\),
and their difference satisfies \(1\le n+1\le m-1\).
The order constraint for powers up to \(m\) gives nontriviality
at every chosen old or new place. The finite character used
above is therefore the character on \(R^\times\) of the square.

\subsubsection{Fourier transform in the complementary blocks}
\label{sec:period-finite-fourier-blocks}

\sectionopening{In this subsection, we calculate the finite transform for the exchange identity in Proposition~\ref{prop:normalized-period-exchange}.}

Write \(X=(X_1,X_2)\) and \(Y=\binom{Y_1}{Y_2}\)
with square blocks \(X_2,Y_1\). For a fixed
\(w\in M_{a,b}(R)\) use
\[
 \Phi(X)=\chi^{-1}(\det X_2)\,
             \mathbf1_{\{X_2^{-1}X_1=w\}},
 \qquad
 \Psi(Y)=\chi^{-1}(\det Y_1)\,
             \mathbf1_{\{-Y_2Y_1^{-1}=w\}},
\]
defined to be zero if their square block is not
invertible. These functions are paired. For supported
pairs \(XY=0\). For a given \(k\), \(Y_1\) is invertible
if and only if \(X_2\) is invertible, by complementary
inverse blocks. When they are invertible, block
elimination gives
\[
 X_2^{-1} X_1=-Y_2Y_1^{-1},\qquad
             \det k=\det Y_1/\det X_2 .
\]
Conversely, a supported \(X\) is completed to the last
covectors of an inverse by adjoining a complementary covector
basis, since \(X_2\) is invertible. A supported \(Y\)
is completed to a column basis in the same way.
The verification works over \(R\) factor by factor.

Using the covariance on each supported basis, the two
weighted sums reduce to
\begin{equation}\label{eq:triangle-square-block-test}
 |\mathrm{GL}_b(R)|P'\binom{I_b}{-w}
     =c\,|\mathrm{GL}_a(R)|P(w,I_a).
\end{equation}
For example \(Y=\binom{I_b}{-w}Y_1\); its character
cancels \(\chi^{-1}(\det Y_1)\) and every \(Y_1\)
contributes the same scalar. The left block is similar.

Now take the successive ranks \(N,N+t,N+2t\) at
the square. In the transposed-lift square \(Y\) has
\(N\) columns, \(X\) has \(N+2t\) covectors, and the
fixed outer blocks adjoined to each datum are the
original inclusions. The abstract sizes in the next
calculation are \(n=N+t\), \(a=N\), and \(b=t\). Write
\(p\) for a datum of size \(n\)-by-\(a\) and define
\begin{equation}\label{eq:triangle-exact-next-datum}
 p^+\begin{pmatrix}\eta&s\\0&I_b\end{pmatrix}
       =\frac{p(-\eta)}{|R|^{nb}},
 \qquad s\in M_{n,b}(R),
\end{equation}
with value zero on all other matrices. In the order
\(n,a,b\), the two full arguments are
\[
 Y=\begin{pmatrix}-\eta\\I_a\\0\end{pmatrix},\qquad
 X=\begin{pmatrix}I_n&\eta&s\\0&0&I_b\end{pmatrix}.
\]
To obtain them in \eqref{eq:triangle-square-block-test},
write the basis in the order \(a,(n,b)\) and use
\(w=\binom{\eta}{0}\). Reindex the ambient module from
\((a,n,b)\) to \((n,a,b)\) in both the domain and the range of
the completion \(k\). Its matrix is conjugated by the block
permutation, so its determinant and the matching character are
unchanged.

Multiplying the covector argument by
\(\begin{psmallmatrix}I_n&s\\0&I_b\end{psmallmatrix}\)
introduces the displayed \(s\). This shear has
determinant one, so the covariance preserves the
tested value. Average over \(s\) and then sum with
weight \(p(-\eta)\) over \(\eta\). Write
\(\iota_d=\binom{I_d}{0_{b\times d}}\) for \(d=a,n\).
With \(v=-\eta\) this gives
\[
 |\mathrm{GL}_a(R)|\sum_{v\in M_{n,a}(R)}
        p(v)P'\binom{v}{\iota_a}
 =c\,|\mathrm{GL}_{n+b}(R)|\sum_{Z\in M_{n+b,a+b}(R)}
        p^+(Z)P\bigl((\iota_n,Z)\bigr).
\]
The sum over \(s\) has \(|R|^{nb}\) terms, giving precisely the
denominator in \eqref{eq:triangle-exact-next-datum}. In column
notation the weights are \(w_n,w_{n+2}\).

For this forward square write \(k_{i,n}\) for the archimedean
projection scalar of Section~\ref{sec:four-unitary-real-projections},
and \(c_{i,n}(\xi)=c_\xi\) for its nonzero measured partner scalar.
The former is independent of \(\xi\). The paired-data scalar is
\(c=k_{i,n}/c_{i,n}(\xi)\), as in \eqref{eq:chain-gram-first-comparison}.
The two intermediate ranks differ by \(2t\). In
Proposition~\ref{prop:normalized-period-exchange} the real constant
is fixed, the old primes occur in both moduli, and each new place
contributes \(q_w^{t^2}\) to \(M_\nu(t)\). Thus
\[
 \frac{|c_\xi|^2}{|c_1|^2}
     =q_\xi^{t^2}Q_{i+1,i+n+2}(\xi)^2,
\]
where \(Q_{i+1,i+n+2}\) is the product over the critical points
\(1,\ldots,t\) of the absolute ratios for the quotient character
at \(\xi^{n+1}\). When the original source \(W\ne0\),
\eqref{eq:chain-gram-propagation} gives \(k_{i,n}\ne0\).
The reciprocal conjugate
\(\beta'_{i,n}(\xi)=\overline{c_{i,n}(\xi)/k_{i,n}}\)
then has this same squared-modulus ratio.

\indent{\it Fourier extension for the reversed square.}
Write \(\mathcal F p(Z)=\sum_X p(X)\psi(\operatorname{tr}(ZX))\)
for the unnormalized matrix Fourier transform with
the original primitive additive character of \(R\).
For block \(Z=\begin{psmallmatrix}A&B\\C&D\end{psmallmatrix}\)
compatible with \eqref{eq:triangle-exact-next-datum},
character orthogonality gives
\[
  \mathcal F(p^+)(Z)=
       \mathbf1_{\{C=0\}}\mathcal Fp(-A)\,
                                     \psi(\operatorname{tr}D).
\]
Indeed the trace in the Fourier phase is
\(\operatorname{tr}(A\eta)+\operatorname{tr}(Cs)
+\operatorname{tr}(D)\). The mean over \(s\) is zero
unless \(C=0\), by primitivity and the matrix trace
pairing. The change of variable \(\eta\mapsto-\eta\)
accounts for the sign.

\subsubsection{The finite determinant Gauss sum}
\label{sec:period-finite-determinant-sum}

\sectionopening{In this subsection, we evaluate the determinant factor in Section~\ref{sec:period-finite-fourier-blocks}.}

We use the following exact determinant Gauss sum over a
residue field \(k\) of size \(q\). If a function \(f\)
on \(M_b(k)\) transforms by the nontrivial character
\(\chi(\det g)\) under left multiplication, it is zero
on singular matrices. A singular matrix has a left
stabilizer with arbitrary determinant, obtained by
scaling a line complementary to its image. Thus
\(f(D)=\chi(\det D)f(I_b)\) on all matrices with
the zero extension convention.

The additive Fourier
transform of \(\chi\circ\det\) transforms by the inverse
determinant character on the opposite side. It is
also zero on singular matrices, using their right
stabilizer. Consequently it equals
\(\mathcal G\,\chi^{-1}\circ\det\), where
\(\mathcal G\) is its value at \(I_b\). For the convention
\(\mathcal F f(Y)=\sum_D f(D)\psi(\operatorname{tr}(YD))\),
average the defining sum for \(\mathcal G\) over left multiplication
by upper unitriangular matrices \(u\). This preserves \(\chi(\det D)\).
The trace is
\(\operatorname{tr}(uD)=\operatorname{tr}D+\sum_{i<j}u_{ij}D_{ji}\).
The mean vanishes unless every \(D_{ji}\) with \(i<j\) is zero.
For upper triangular \(D\) the strictly upper entries give
\(q^{b(b-1)/2}\) choices and the nonzero diagonal entries sum
independently. Thus for \(b\ge1\)
\[
 \mathcal G=q^{b(b-1)/2}
       \left(\sum_{s\in k^\times}\chi(s)\psi(s)\right)^b ,
\]
the full-rank case of
\cite[arXiv:1105.4513v1, (2.12)--(2.13) and Theorem~2.1]{LiHuMatrixGauss}.
For \(b=0\) the matrix sum has its single empty term, of value one.
Parseval for the unnormalized Fourier transform on
the \(b^2\)-dimensional matrix space gives
\[
                   |\mathcal G|^2=q^{b^2}.
\]
The original modulus is squarefree, and the quotient character is nontrivial in each field factor,
because its old and new tame powers have full support.
The additive character is primitive in each factor.
Taking products proves the same formula with \(|R|\)
in place of \(q\).

In the reversed square the abstract block sizes are
\(n=N\), \(a=N+t\), \(b=t\). Thus the two
linear-group orders in \eqref{eq:triangle-square-block-test}
are equal after the shear extension, since \(a=n+b\).
Put \(p_*=\mathcal F(\bar p)\). In the sum weighted by
\(\mathcal F(\overline{p^+})\), the last formula leaves
only matrices with lower-left block zero. At fixed
\(\eta,s\), the tested value as a function of the
lower-right block \(D\) transforms by
\(\chi(\det g)\) under \(D\mapsto gD\), by applying
\(\operatorname{diag}(I_n,g)\) to the whole covector argument.

The sum over \(D\), with its Fourier phase, is therefore
\(\mathcal G\) times the value at \(D=I_b\).
The extra factor \(|R|^{nb}\) converts the unrestricted
sum over \(s\) into the mean in \(p_*^+\). Hence
\[
 \sum_X \mathcal F(\overline{p^+})(X)P(\iota,X)
     =\mathcal G|R|^{nb}\sum_X p_*^+(X)P(\iota,X),
\]
where \(\iota\) is the fixed inclusion adjoined to the
covector argument. The squared modulus of the coefficient
here is \(|R|^{b(b+2n)}\). Combining it with the
paired-data law for \(p_*\) gives the reversed
commutation coefficient with the additional inverse
factor \((\mathcal G|R|^{nb})^{-1}\). The remaining
paired-data scalar for that square has twist-independent
modulus because its two middle ranks agree, and the
critical product is empty. Thus the reversed
commutation scalar has modulus a fixed constant times
a fixed half-integral power of \(q_\xi\).

\indent{\it Matrix reversal in the original two bases.}
Let \(S,T\) have ranks \(u,v\), and let
\(R=\mathcal O_E/\mathfrak m\) have the original selected
split support. Denote by \(K_{S,T}^{a,b}(D,P;g,h)\)
the theta kernel with residue datum \(D\) and real
Fock polynomial \(P\). Write
\(\psi_0(\bar x)=\psi_f(j_{\rm W}\alpha_f\bar\alpha_f
x/t_{\mathfrak m})\), where \(t_{\mathfrak m}\) is the
chosen generator on this support. Let
\(s_{\mathfrak m}=\overline{t_{\mathfrak m}\alpha_f}/
(t_{\mathfrak m}\alpha_f)\) at the finite places.
With the unnormalized matrix Fourier transform
\[
 \widehat{\bar D}(Z)=
       \sum_{X\in\operatorname{Mat}_{v\times u}(R)}
           \overline{D(X)}\,\psi_0(\operatorname{tr}(ZX)),
\]
the exact kernel identity is
\[
 \overline{K_{S,T}^{a,b}(D,P;g,h)}
   =\tau_{\mathfrak m}
      K_{T,S}^{b,a}(\widehat{\bar D},\operatorname{flip}(\bar P);
                                                  h,g s_{\mathfrak m}),
 \qquad
 |\tau_{\mathfrak m}|^2=|R|^{-uv}.
 \tag{MR}\label{eq:display-mr}
\]
For source signs \(s_j\) and target signs \(t_i\), the adjoint
on the opposite source is \(Y=-J_SX^*J_T\). On each coordinate
\(Y_{ji}=-s_jt_i\overline{X_{ij}}\). The signed Fock coordinate
in this position is either the coordinate or its conjugate, with the
same choice in the reversed pair. Thus conjugation followed by
adjunction sends a polynomial to
\[
   (\operatorname{flip}(\bar P))(Z)
             =\bar P\bigl((-s_jt_i Z_{ji})_{i,j}\bigr),
\]
where the bar on \(P\) conjugates its coefficients. For a minor on
row and column sets \(I,J\) of size \(h\), the substitution contributes
\((-1)^h\prod_{i\in I}t_i\prod_{j\in J}s_j\) times the minor
on \(J,I\) in the transposed order.

First conjugate with the source replaced by \(S^-\).
Negating the Hermitian form negates the symplectic form
on \(\operatorname{Hom}(S,T)\), so conjugation intertwines
the Heisenberg actions with the two opposite central
characters. In the chosen coordinates the real sign of
every variable changes, so conjugation replaces its
holomorphic polynomial coefficients by \(\bar P\) on the
opposite signs. At the finite places the coordinate change
is multiplication by \(u_0^{-1}\), for the fixed unit
\(\bar u_0u_0=-1\).

Its preservation of the integral boxes
changes \(D(X)\) to \(\overline{D(u_0^{-1}X)}\).
The two determinant weights become
\(a^{-1}\upsilon^u,b^{-1}\upsilon^v\).
These powers are the change from each split Weil
normalization to its complex conjugate: the underlying
splitting character of \(E^\times\backslash\mathbb A_E^\times\)
is conjugated to its inverse, while the fixed standard
splitting on a rank \(h\) space contributes \(\upsilon^h\).
The rational matrices coincide under passage to \(S^-\)
when the coordinate change is undone. Equivariance for
rational Heisenberg translations identifies the transported
functional with a scalar multiple of the normalized theta
functional. To determine the modulus, write \(J\) for the
unitary or antiunitary transport. It carries the rational
translation lattice onto the new one and preserves probability
Haar on the quotient. In the self-dual Schr\"odinger measures,
\[
 \int |\theta'(\rho'(x)Jv_0)|^2\,dx
    =\|Jv_0\|^2=1 .
\]
Indeed, integration in the modulation variable gives the
sum of squared Schr\"odinger values by Parseval and integration
in the displacement variable unfolds this to their squared
norm. The same identity holds for the original unit vacuum.
Consequently the scalar has modulus one. The adjoint
transport below has the same property before the finite dilation.

Next apply the Hermitian adjoint on the opposite-source
pair, which returns to the original \(S\) as the second
space. In rational coordinates it is
\(X\mapsto J_{S^-}^{-1}X^*J_T\), a bijection of the rational
matrix sets. The signed real transpose and the original
finite basis scalar are the entrywise transports in the
line calculation. On the original source group the basis
scalar moves the finite argument by
\(\bar\alpha_f/\alpha_f\). Swapping again inverts the
determinant weights and inserts the same respective rank
powers, so the two weights are now \(b,a\).

To verify the finite factor, use a point indicator at
\(X_0\) after the opposite-source transport. Its pre-dilation
vacuum displacement has box \(X-X_0/t_{\mathfrak m}\)
integral at the selected places. The split source
dilation by \(t_{\mathfrak m}^{-1}I_u\) gives the box
\((X-X_0)/t_{\mathfrak m}\) integral. In rank \(u\)
the determinant norm of the dilation is \(|R|^u\).
The split Weil factor is its \(v/2\)-power, giving
squared modulus \(|R|^{uv}\) for the multiplier from
the displaced vacuum to the point datum.

The Hermitian transpose sends the displacement to the
conjugate places, where the Heisenberg relation turns it
into multiplication on the original integral box by
\(\psi_0(\operatorname{tr}(Z(-\bar u_0)X_0))\).
This is the entrywise trace pairing, so the computation
for one column extends to the full matrix without a second
measure factor on the conjugate split side. The normalized
theta transport under adjoint is a unit phase. The source
dilation moves to the group argument and combines with
the basis scalar, giving precisely \(s_{\mathfrak m}\).
Solving for the original kernel inverts the displacement
multiplier, giving the stated value of \(|\tau_{\mathfrak m}|^2\).
Locally the annihilator of \(t_w\mathcal O_w\) for
\(\psi_{E,w}(\delta_w xy)\) is \(t_w^{-1}\mathcal O_w\).
It follows that the pairing on \(\mathcal O_w/t_w\mathcal O_w\)
given by \(\psi_{E,w}(\delta_w xy/t_w)\) is perfect. For
rectangular matrices, testing one matrix unit at a time gives the
same conclusion for the trace pairing on the two opposite shapes.
At an unselected place the self-dual integral box has volume one
and transforms to itself.

Finally sum the point identities with coefficients
\(\overline{D(u_0^{-1}X_0)}\). The substitution
\(X_0=u_0X\) cancels the shift in the additive character
because \(-\bar u_0u_0=1\). This gives exactly
\(\widehat{\bar D}\) under the unchanged primitive character
\(\psi_0\), proving \eqref{eq:display-mr}. On applying the
conjugated transform twice, orthogonality for the primitive trace
pairing gives
\[
 \sum_Z\overline{\widehat{\bar D}(Z)}
          \psi_0(\operatorname{tr}(XZ))
 =\sum_{X'}D(X')\sum_Z\psi_0(\operatorname{tr}((X-X')Z))
 =|R|^{uv}D(X).
\]
The two real signed transpositions also compose to the identity,
since \((-s_jt_i)^2=1\) at each entry. Put \(z=s_{\mathfrak m}\)
and write \(\tau_{S,T},\tau_{T,S}\) for the scalars in the two
orders, with the corresponding characters. Applying
\eqref{eq:display-mr} twice gives the original kernel as
\(\overline{\tau_{S,T}}\tau_{T,S}|R|^{uv}\) times that kernel
at \((gz,hz)\). The simultaneous action on \(\operatorname{Hom}(S,T)\)
is the identity, and the determinant characters contribute
\[
 K_{S,T}^{a,b}(D,P;gz,hz)
       =a(z^v)b(z^u)^{-1}K_{S,T}^{a,b}(D,P;g,h).
\]
Here \(z\bar z=1\), so this factor has absolute value one.
The same right translations preserve quotient Haar and every
principal level by centrality, including at the two halves of a split
place where \(z\) need not be integral.
Write \(J(P)=\operatorname{flip}(\bar P)\). On edge \(j=0,1\)
let \(P_j,P'_j\) be the distinguished original and reversed
harmonics, and \(p_j,p'_j\) the respective products of the
commuting source and target line projectors. In each orientation
Section~\ref{sec:four-unitary-real-projections} identifies the
source isotypic component with \(V_S\widehat\otimes V_T\).
The two projectors there select the tensor product of the two
highest lines. Thus \(p_jP_j\) and \(p'_jP'_j\) span their
respective lines. In a pairing with the two prescribed form types
one may insert both projectors on the oscillator datum.
Indeed moving each coefficient integral to its form fixes that
form, with the inverse coefficient in the bilinear source pairing.
The weighted coefficient moments and the two rapid tests justify
both interchanges as in Section~\ref{sec:four-unitary-constant-terms}.

Define the real-place scalars by
\[
 p'_0 J(P_0)=\overline{k_0}\,p'_0P'_0,\qquad
 p_1J^{-1}(P'_1)=\overline{k_1}\,p_1P_1 .
\]
On the first edge apply \eqref{eq:display-mr} to \(P_0\) and
compress \(J(P_0)\) in the reversed kernel, which is in the
linear leg. The tested factor is \(\tau_{\mathfrak m}\bar k_0\).
On the second edge apply \eqref{eq:display-mr} to \(J^{-1}(P'_1)\),
so that the reversed kernel has polynomial \(P'_1\).
Solve for that kernel using \(\tau_{\mathfrak m}\ne0\).
The original kernel is conjugated, so its compression gives the
tested factor \(\tau_{\mathfrak m}^{-1}k_1\), with the ranks
of the second edge. Both identities hold also for a zero
compression. The two scalars \(k_j\) depend only on the real data. In
\eqref{eq:chain-terminal-gram} the two transforms use
\(R=R_0(\xi)\). With \(s_0,s_1\) these tested factors,
\(\gamma\) the reversed-square factor, and \(\eta\) the unitary
translation phase in that equation, their product satisfies
\[
   |s_0\gamma\eta s_1|^2=C_{i,n}|R_0(\xi)|^{t^2},
   \qquad C_{i,n}\ge0 ,
\]
where \(C_{i,n}\) is independent of \(\xi\).

\section{The leading coefficient and the four-block test}
\label{sec:four-block-coefficient}

\sectionopening{In this section, we estimate the leading Taylor coefficient on the prescribed Weyl lattice
and construct the four-block test. The reference norm determines the
normalization of the critical ratios in Section~\ref{sec:common-critical-twist}.}

\subsection{A lower bound for a coefficient in twelve coordinates}
\label{sec:marked-twelve-sharp-coefficient}

\sectionopening{In this subsection, the four-block test of Section~\ref{sec:v2-four-block-main-numerator}
requires the leading Taylor coefficient from its first twelve
coordinates. We compare this coefficient with its Weyl norm in the
prescribed lattice and normalization.}

Write \(V_0=s_{\rm vac}\), \(\Lambda_{12}=\mathfrak a_0^{12}\), and
\(\delta=\bar\pi_M\). Use the original \(12\)-coordinate
Weil vacuum \(v\). The Taylor map \(\mathcal T\) includes
the scale \(\pi_Ms_{\rm jet}\).
For nonnegative integers \(\lambda_\tau\) satisfying
\(\kappa_\tau\le\lambda_\tau\le\kappa_\tau+x^{42}\), put
\[
 n_\#=\lceil6b/k\rceil,\quad
 B_\tau=\prod_{j=1}^{n_\#}((\lambda_\tau-n_\#)_++j),
 \qquad P_p=\sum_\tau S_p(\kappa_\tau,\tau).
\]
We construct one nonzero \(t_0\), depending on the
datum and the \(\lambda_\tau\), with
\[
 \log|\sigma t_0|\le fmb/k\qquad\text{for every }\sigma.
\]
At each embedding above a prime \(p\), we will find
an exponent \(\alpha\) satisfying
\[
 \sum_c\alpha_{\tau,c}=\lambda_\tau
\]
and a residue cell \(x_0\bmod\pi_M\) with indicator \(I_{x_0}\)
such that
\begin{equation}
 v_p\!\left(\frac{V_0^{12}}{t_0}
      \prod_\tau\lambda_\tau!\,
       [Z^\alpha]\mathcal T I_{x_0}\right)\le P_p.
 \label{eq:marked-rank-twelve-sharp}
\end{equation}

\begin{lemma}[Field of definition of the coefficient test]\label{lem:normalized-elliptic-jet-base-change}
For every residue cell \(x_0\) and multi-index \(\alpha\), the coefficient
\(c_{\alpha,x_0}=V_0^{12}[Z^\alpha]\mathcal T I_{x_0}\) is algebraic.
For any embedding \(\iota:\overline{\mathbb Q}\hookrightarrow\mathbf C_p\),
its image is the coefficient obtained from the corresponding CM sections
using the base-changed invariant differential, torsion translations and
finite character sums.
\end{lemma}

\begin{proof}
Work over a number field defining the CM elliptic curve, its invariant
differential, its CM action and the distinguished point \(t=-t\ne0\).
Since \(t=-t\), the bundle \(\mathcal O([t])\) is symmetric. Its
semicharacter has order at most two, and its Hermitian form satisfies \(H(az,w)=H(z,\bar a w)\)
for \(a\in\mathcal O_K\). The section has value one at the origin.
Consequently its reduced theta series and its normalized torsion
translates have algebraic coefficients in the logarithm coordinate
\(z\) of this differential
\cite[Prop.~2.1 and Thm.~2.9 in the authors' version]{BannaiKobayashiTheta}.
More precisely, let \(F\) be the field of definition above, enlarged
to contain \(K\). Since \(0\notin[t]\), the local equation at zero can
be taken to be \(1\); the normalized value there is \(1\).
If \(N_u u\in\Lambda_C\), put \(n=2N_u\) and write
\(\pi_C:\mathbb C\to C(\mathbb C)\) for uniformization.
Retain the division point \(Q_u=\pi_C(u/(2n))\in C[4N_u^2]\)
and its double \(P_u=2Q_u\). The symmetric bundle \(\mathcal O([t])\)
and these points determine the rational translation function
\(R_u=f^M_{[t],P_u}\in F(C[4N_u^2])(C)\), with
\[
 R_u(\pi_C(v))=\frac{(U_u^M\vartheta_C)(nv)}{\vartheta_C(nv)},\qquad
 (U_u^M\vartheta_C)(z)
   =e^{-\pi(z\bar u+|u|^2/2)/A_C}\vartheta_C(z+u).
\]
The canonical rational function and the displayed analytic formula
are those of \cite[Cor.~2.4 and Lemma~2.7 in the authors' version]{BannaiKobayashiTheta}.
The trivialization at \(Q_u\) and its inverse cancel their scalar
factors in that construction, and its pullback description
\cite[Prop.~2.8 in the authors' version]{BannaiKobayashiTheta}
identifies the function from the divisor and division points on
\(C^\iota\) with \(\iota R_u\).
The period law for the semicharacter \(\chi_C\) gives
\[
 U_{u+\lambda}^M\vartheta_C
 =\chi_C(\lambda)
  \exp\!\left(\frac{\pi(u\bar\lambda-\lambda\bar u)}{2A_C}\right)
       U_u^M\vartheta_C,\qquad \lambda\in\Lambda_C .
\]
The multiplier for a different representative of the torsion
point is a root of unity, since
\(\operatorname{Im}(u\bar\lambda)/A_C\in N_u^{-1}\mathbb Z\)
and \(\chi_C(\lambda)\in\{1,-1\}\).

Write \(y_{\rm W}\) for the affine coordinate at \(t=(0,0)\)
on the Weierstrass model with group origin at infinity. Its
shifted pullback \(y(P)=y_{\rm W}(P+t)\) satisfies
\(y(0)=0\) and \(y'(0)=1\) in the invariant logarithm coordinate
by Section~\ref{sec:marked-integral-sections-ordinary-jets}.
Take \(g_u=y\) if \(\pi_C(u)=t\), when the translated divisor has
a simple zero at the origin, and take \(g_u=1\) otherwise.
Theorem~2.9 of the same reference gives
\(g_u^{-1}U_u^M\vartheta_C\in F(C[4N_u^2])[[z]]^\times\),
where \(g_u\) is pulled back by \(\pi_C\).
Multiplication by \(g_u\) recovers the translated series, including
its zero constant term when \(\pi_C(u)=t\).

The forms in Lemma~\ref{lem:upper-cm-transport} satisfy
\[
 A_{k,i}(\pi_Ms_{\rm jet}Z)=
 \frac{\pi_Mp_F}{q_1}\sum_\tau\overline{\tau(\zeta_{q_1}^k)}Z_{\tau,i},
 \qquad p_F=s_0/u\in\overline{\mathbb Q}^{\times}.
\]
Here \(s_0\) has a power in \(K\), and \(u\) is obtained by taking a root
of the product of algebraic isogeny multipliers and discriminants in the
same lemma. The finite translates of its rank-one seed are products of the
preceding torsion series in these forms, up to roots of unity. This
holds at every finite translate by rational approximation modulo the
open stabilizer. Tensoring \eqref{eq:upper-CM-vacuum} twelve times
therefore gives the asserted algebraicity for every translate of
\(V_0^{12}v\).

On the residue group \(D\) at the prescribed places the identity
\[
 {\bf1}_{\bar x=x_0}=|D|^{-1}
         \sum_{\chi\in\widehat D}\chi(-x_0)\chi(\bar x)
\]
expresses the cell indicator times the vacuum as a finite combination
of its translates by the dual lattice. The perfect pairing is given in
Section~\ref{sec:twelve-coordinate-finite-pairing}. All coefficients
and substitutions in this expression lie in a number field. In the
chart at \(t\), let \(x(Y)\in YF[[Y]]\) solve
\(x^3+x-a_1xY-Y^2=0\). It is unique since the derivative in \(x\)
at zero is one. The invariant differential gives
\[
 z=\log_C(y)=\int_0^y\frac{dY}{1+3x(Y)^2-a_1Y}=y+O(y^2).
\]
The inverse \(y(z)\) reads the point \(\pi_C(z)+t\).
Accordingly, substitute \(y(z/n)\) in \(R_u\circ T_{-t}\),
where \(T_{-t}(P)=P-t\), to obtain \(R_u(\pi_C(z/n))\).
This and multiplication by \(\vartheta_C(z)\) recover the translated
series. In each degree the substitutions commute with \(\iota\),
including when \(p\mid n\), since they are over a characteristic-zero
field. Coefficient extraction then gives \(\iota(c_{\alpha,x_0})\).
\end{proof}

\smallskip
\noindent\textit{The one-variable tests.}
Use the integral parameter \(y=z+O(z^2)\)
and \(R=\mathcal O_{\mathbf C_p}\).
The formal group is integral: the translated identity
is the section \((0,0)\) of the good Weierstrass model,
and the group law transported by this section gives a
series in \(R[[y_1,y_2]]\). Differentiation preserves
\(R[[y]]\) since \(y'=1+3x^2-a_1y\).
Let \((\pi,q)\) be \((p,p)\) if \(p\) splits,
\((p,p^2)\) if it is inert, and \((\sqrt{-d},p)\)
if \(p=d\). The ordinary integral-section calculation
gives \(y(az)\in R[[y]]\) for \(a\in\mathcal O_K\) on the formal
group at its identity. For \(a\ne0\) let \(i(a)\)
be the order of the reduced parameter of \([a]\). The CM isogeny
extends over the good model, with degree \(a\bar a\) on the special
fibre: its kernel is contained in the finite kernel of \([a\bar a]\),
and the degree is the degree of the pullback of a degree-one line
bundle, constant in this smooth proper family. Its local
multiplicity at the identity is \(i(a)\), so \(i(a)\le a\bar a\).
For split \(p\), take a principal power
\(\mathfrak p^h=(a)\) of the chosen split prime.
Then \(a\bar a=p^h\), \(\bar a\) is a unit at the chosen embedding, and \(i(\bar a)=1\). Composition of reduced
parameter series gives \(i(p)^h=i(a)i(\bar a)\le p^h\).
Thus \(i(\pi)\le q\) in every case.

Average \(a^{-1}y(az)\) over \(a^{q-1}=1\). These
Teichm\"uller values are limits of CM integers at the chosen embedding, so the result \(T=y+O(y^2)\) is integral.
Root-of-unity covariance gives
\[
 T(\pi z)=\pi T(z)+b_1T(z)^q+
             \sum_{j\ge2}b_jT(z)^{1+j(q-1)},
 \qquad b_1\in R^\times ,
\]
where the reduced order bound gives the unit assertion.
This also works at \(p=2\), when \(q-1=1\).
For \(\varpi^{q-1}=-\pi/b_1\) put
\(\rho=v_p(\varpi)=v_p(\pi)/(q-1)\).
The coefficients of \(Y=T(\varpi z)/\varpi\) are
integral. Its equation is \(Y(\pi z)/\pi=Y-Y^q+\pi E(Y)\)
with \(E\in X^2R[[X]]\), and the recursion divides
in degree \(n>1\) by the unit \(1-\pi^{n-1}\).
Thus \(Y\equiv\sum_{j\ge0}z^{q^j}\pmod\pi\).

Put \(c_n=\max_{H\in R[[y]]}|n![z^n]H|_p\).
For \(H=\sum_j h_jy^j\), only \(0\le j\le n\) contribute.
Since \(|h_j|_p\le1\), the ultrametric inequality gives
\(c_n=\max_{0\le j\le n}|n![z^n]y^j|_p\). A monomial attains this
maximum. Integral differentiation
gives \(c_n\le1\), and the integral group law gives
\(c_{m+n}\le c_mc_n\) by taking the \(z^m w^n\)
coefficient of \(H(z+w)\) and multiplying by \(m!n!\).
Since \(T(z)=\varpi Y(z/\varpi)\) belongs to \(R[[y]]\),
\([z^n]T=\varpi^{1-n}[z^n]Y\) and \([z^{q^j}]Y\) is a unit.
Taking \(H=T\) and using Legendre's formula gives
\[
 c_{q^j}\ge p^{-\gamma_p(q^j-1)},\qquad
 \gamma_p=\frac1{p-1}-\rho .
\]
The case \(n=0\) gives \(c_0=1\). For \(n\ge1\), write
\(q^j=an+b\), \(0\le b<n\). Submultiplicativity gives
\(c_{q^j}\le c_n^a\). Let \(j\) tend to infinity to get
\(c_n\ge p^{-\gamma_p n}\).
Consequently in each degree \(n\) some \(H\in R[[y]]\) has
\(v_p(n![z^n]H)\le\gamma_p n\).

Put \(g=-\log\theta\). The normalized pullback
\(\exp(a\bar a\,g(z)-g(az))\in R[[y]]\)
extends to limits of \((a,a\bar a)\).
For \(\delta'\in\mathbf Z_p\) approximate \(\delta'\)
by ordinary integers, obtaining the same assertion
with \(\bar a a=(\delta')^2\).
There are four elements of \(\mathbf Z_p\)
whose squares sum to \(-1\), by four squares
modulo \(p^u\) and compactness. Adjoin \(1\) and \(0\)
to this list to give \(\delta_0,\ldots,\delta_5\),
with \(\delta_0=1\), and \(\sum_c\delta_c^2=0\).
It follows that
\(\prod_{c=0}^5\theta(\delta_c z)\in R[[y]]^\times\)
by multiplying the six pullbacks.

\smallskip
\noindent\textit{The weighted Gram identity and coefficient calculation.}
Use vectors \(w_i\) in the trace-dual lattice
\(\mathfrak a_0=E\cap\alpha_{\rm W}\widehat O_E\),
whose positive integral multiplicities \(c_i\) satisfy
\[
 \sum_i c_i\tau(w_i)\overline{\sigma(w_i)}
       =M\delta_{\tau,\sigma},\qquad c_i\ge6fb.
\]
For this family take the \(q_1\) powers of \(\zeta_{q_1}\)
and the trace-dual mixed family constructed for ordinary
jets, before squaring either scalar. Their scalar
constants are \(q_1\) and \(h\) respectively, with
\((h,q_1)=1\). They can both be taken at most \(b\).
In the mixed family the dual vectors are
\(\bar\beta^{\lceil\log_2q_1\rceil}b_i/\Phi'_{q_1}(\zeta)\),
with \(b_i\) the coefficients of \(\Phi_{q_1}(X)/(X-\zeta)\).
Write \(u=q_1/q_0\). The polynomials satisfy
\(\Phi_{q_1}(X)=\Phi_{q_0}(X^u)\): if the \(u\)-th
power of a root of unity has order \(q_0\), its order
is \(q_1\), since every prime dividing \(u\) divides \(q_0\).
The coefficient bound for \(\Phi_{q_0}\) and polynomial
division at \(|\zeta|=1\) give \(|b_i|\le q_1 2^{q_0}\)
at every embedding. Differentiation gives
\[
 |\Phi'_{q_1}(\zeta)|
     =u|\Phi'_{q_0}(\zeta^u)|\ge c(q_0)u ,
\]
where \(c(q_0)>0\) is the minimum at the primitive
\(q_0\)-th roots. These are bounds at fixed \(q_0\).
Thus \(h\) is polynomial in \(q_1,d^{\lceil\log_2q_1\rceil}\),
whereas \(b\) exceeds every fixed power of \(x\).
Put \(m_0=\lceil6fb\rceil\).
The scale \(M\ge b^{20}\) ensures
\(M-m_0(q_1+h)\ge q_1h\).
Write the left side as \(u q_1+v h\) with \(u,v\ge0\):
choose \(0\le v<q_1\) in the required residue class first.
The two multiplicities are \(m_0+u,m_0+v\).

Write \(\mu_i(Z)=p_F\sum_\tau
\overline{\tau(w_i)}Z_\tau\). At every coefficient embedding these forms
have integral coefficients by the period-factor
valuation. At least one \(\mu_i\) has only unit coefficients.
At weight embeddings outside \(\beta\)-primes use the
vector \(1\). At the primes over \(\beta\) use the
unchanged vector
\(\zeta\bar\beta^{\lceil\log_2q_1\rceil}/\Phi'_{q_1}(\zeta)\);
the denominator cancels the period factor in valuation.
Either case applies to all weight embeddings together,
since their restrictions to \(K\) are equal.
At an inert prime the \(\mu_i\) span the coordinate forms
integrally: the columns for the powers \(1,\zeta,\ldots,\zeta^{f-1}\)
give a Vandermonde determinant which is a unit since \(p\nmid q_1\)
and the \(f\) distinct \(q_1\)-th roots \(\tau(\zeta)\) remain distinct
in the residue field.

Put
\(\Theta(Z)=\prod_i\theta(\mu_i Z)^{c_i}\) and let
\(\mathcal A\) be the integral power series in
\(y(\mu_i Z)\). For a degree vector \(\nu\) write
\(q_N(\nu)\le u\) if some coefficient of combined degree
\(\nu\) in a product of \(N\) independent copies \(H\Theta\),
\(H\in\mathcal A\), has valuation at most \(u\).
Products on distinct copies add bounds. An integral linear
substitution of copies within each coordinate preserves the
combined degree, and each new coefficient is a finite
integral combination of the old ones.

For the six-copy general bound,
\(U(Z)=\prod_{c=0}^5\Theta(\delta_c Z)\) is a unit of
\(\mathcal A\) by the one-variable identity. In the
first copy take \(H=H_1(\mu_{i_0}Z)/U(Z)\) with
\(\mu_{i_0}\) the unit form, and in the other five take
\(H=1\). Substitute \(Z_c=\delta_c Z\). The product
becomes \(H_1(\mu_{i_0}Z)\), whose coefficient is obtained
from the one-variable test and multinomial formula.
Thus
\begin{equation}
 q_6(\nu)+v_p(\nu!)\le\gamma_p|\nu|,
 \qquad \nu!=\prod_\tau\nu_\tau!.
 \label{eq:marked-twelve-six}
\end{equation}
At inert \(p\), the spanning property also gives
\(q_6(\nu)\le0\): multiply the integral series with
linear terms \(Z_\tau\), to powers \(\nu_\tau\), by
\(\Theta\). The lowest term in total degree is the
desired monomial with unit coefficient.

\subsubsection{The inert expansion and valuation}
\label{sec:twelve-coordinate-inert-expansion}

\sectionopening{In this subsection, we apply the original expansion of Section~\ref{sec:exceptional-charged-jet} to the coefficient in twelve coordinates.}

The inert normal form has
\(\Theta(\varpi Z)=W(Z)\exp(K\Xi(Z))\),
\(K=\varpi^{p+1}/p\), with \(W\) an integral unit
and \(\Xi\) reducing modulo \(p\) to
\[
 \sum_{e\text{ odd},\,j\ge0}\lambda_{e,j}
       \left(\sum_\tau Z_\tau^{p^e} Z_{-p^e\tau}\right)^{p^j},
 \qquad \lambda_{e,0}\in R^\times.
\]
Since \(2\) splits in \(K\), an inert \(p\) is odd.
Equation~\eqref{eq:cm-inert-leading-unit} gives a unit for the
\(z^{p+1}\) coefficient. In \(Y^{p+1}\bmod p\) each
\(z^{p^e}z\), for \(e\) positive and odd, occurs once: the two
Frobenius indices are \(0\) and \(e\). For integral \(w_i\) at
this inert prime the odd Frobenius congruence in
Section~\ref{sec:exceptional-rescaled-exponential} and the weighted
Gram identity give
\[
 \sum_i c_i\,\overline{\tau(w_i)}^{p^e}\overline{\sigma(w_i)}
       \equiv \sum_i c_i((-p^e\tau)(w_i))\overline{\sigma(w_i)}
       =M\mathbf1_{\sigma=-p^e\tau}\pmod p .
\]
Each scaled form supplies the additional factor
\(p_F^{p^e+1}\).
Both \(M\) and the period factor are units at an
inert prime, so the transferred coefficient is a unit.

Let selected starts \(\tau\) have odd exponents \(e_\tau\)
and ends \(u_\tau=-p^{e_\tau}\tau\). For each selected \(\tau\),
assume that for every selected \(\sigma\) and every positive odd
\(e\le e_\tau\),
\[
 -p^e\sigma=u_\tau \quad\Longrightarrow\quad
 (\sigma,e)=(\tau,e_\tau).
\]
Use this coupling
\(j_\tau\) times, with degrees
\(\nu'=\sum_\tau p^{e_\tau}j_\tau[\tau]\) and
\(\nu''=\sum_\tau j_\tau[u_\tau]\); put \(J_c=\sum j_\tau\).
Apply the substitution using
\[
 G(Z',Z'')=\Theta(\varpi(Z'+Z''))
                 \prod_{c=1}^5\Theta(\varpi\delta_c Z').
\]
The logarithmic derivative in \(Z''_u\), divided by \(K\), is
\[
 (\partial_u\Xi)(Z'+Z'')
     +K^{-1}(W^{-1}\partial_u W)(Z'+Z'').
\]
Here \(v_p(K)=1/(p-1)-1<0\). If \(L\) is integral and
\(\partial=\partial_{Z''_u}\), then
\[
 K^{-1}\partial(GL)
   =G\bigl(L\,\partial_u\Xi+
      K^{-1}(\partial L+L W^{-1}\partial_u W)\bigr),
\]
where \(W\) and \(\Xi\) are evaluated at \(Z'+Z''\).
Let \(\mathfrak m_R=\{a\in R:v_p(a)>0\}\). In the last
recursion \(K^{-1}\in\mathfrak m_R\) since
\(v_p(K^{-1})=(p-2)/(p-1)\). Reducing in \(R/\mathfrak m_R\) and
iterating gives \(K^{-J_c}\partial_{Z''}^{\nu''}G=GH\), where
\[
 H\bmod\mathfrak m_R
     =\prod_\tau \bigl((\partial_{u_\tau}\Xi)(Z'+Z'')
                         \bmod\mathfrak m_R\bigr)^{j_\tau}.
\] The recursion uses ordinary iterated derivatives over characteristic
zero, including for \(j_\tau\ge p\). The two sets of variables
\(Z',Z''\) are independent even if a start equals an end. In the
residue field the derivative of a loop monomial \(X_u^{p^e+1}\)
is \(X_u^{p^e}\), with coefficient one.

Set
\(Z''=0\) and all \(Z'_\tau\) outside the starts to zero.
The unique lowest term of each relevant derivative
is \(\bar\lambda_{e_\tau,0}(Z'_\tau)^{p^{e_\tau}}\).
In the residue field the terms with \(j>0\) vanish on
differentiation, as does the derivative of every factor
\((Z'_\sigma)^{p^e}\). The surviving summands have
\(-p^e\sigma=u_\tau\); the selection condition leaves just
\((\tau,e_\tau)\) through degree \(p^{e_\tau}\).
Also \(G(Z',0)=U(\varpi Z')\) is integral with
constant term \(1\), by \(T(\varpi z)=\varpi Y(z)\).
After restriction to the starts, every other term of a derivative
has strictly greater total degree than \(p^{e_\tau}\). The product
of its prescribed lowest monomials has degree \(|\nu'|\). Choosing
any higher term, including one from \(U(\varpi Z')-1\), gives
larger total degree. Thus the coefficient of \(Z'^{\nu'}\) in
\((GH)(Z',0)\bmod\mathfrak m_R\) is
\(\prod_\tau\bar\lambda_{e_\tau,0}^{j_\tau}\ne0\).
Evaluation of the ordinary derivatives at \(Z''=0\)
multiplies the corresponding coefficient by \(\nu''!\), including
when \(p\mid\nu''!\). Thus the coefficient of \(Z'^{\,\nu'}Z''^{\,\nu''}\)
in \(G\) has valuation exactly
\(J_c(1/(p-1)-1)-v_p(\nu''!)\).
Undo the scale \(\varpi\) and use the integral substitution
back to the six independent copies. With
\(\nu=\nu'+\nu''\) this gives
\[
 q_6(\nu)\le-J_c,\qquad
 q_6(\nu)+v_p(\nu'!)+v_p(\nu''!)\le\gamma_p|\nu|-J_c.
\]
Indeed undoing the common dilation subtracts
\(\rho|\nu|=|\nu|/(p^2-1)\) from the valuation. For the first bound
use \(|\nu|\ge(p+1)J_c\). For the second use
\(|\nu''|=J_c\) and \(v_p(\nu'!)\le |\nu'|/(p-1)\).
The exact coefficient here is that of the substituted
series. On the original independent variables this
proves existence at the indicated valuation bound.

We now show
\begin{equation}
 q_{12}(\lambda)+\sum_\tau v_p(\lambda_\tau!)
              \le P_p+2\sum_\tau v_p(B_\tau).
 \label{eq:marked-twelve-coeff}
\end{equation}
At a split prime \(q=p\) and \(v_p(\pi)=1\), giving
\(\gamma_p=0\). Also \(1-\chi_K(p^e)=0\) for every \(e\), so
\(P_p=0\). This applies at \(p=2\) because \(d\equiv7\pmod8\).
Use \eqref{eq:marked-twelve-six} and the constant coefficients of
six further copies.
For inert \(p>k\) use every start with \(e_\tau=1\),
\(i_\tau=\lfloor(\lambda_\tau-n_\#)_+/p\rfloor\), and
\(j_\tau=\min(i_\tau,\lambda_{-p\tau})\).
Here \(1\) is the least positive odd exponent. Multiplication by
\(-p\) permutes the places, so each end has a unique start.
Since \(i_\tau\le n_\#\), the outgoing plus incoming
degree at each place is at most
\((\lambda_\tau-n_\#)_++\min(n_\#,\lambda_\tau)=\lambda_\tau\).
The remaining six copies give valuation at most zero on
unused degree. Put \(m_\tau=(\lambda_\tau-n_\#)_+\).
Legendre's formula gives
\[
 v_p(\lambda_\tau!)=i_\tau+v_p(i_\tau!)
                  +v_p(\lambda_\tau!/m_\tau!)
             \le i_\tau+2v_p(B_\tau).
\]
Indeed \(i_\tau!\mid n_\#!\mid B_\tau\), as \(B_\tau/n_\#!\)
is a binomial coefficient, and \(\lambda_\tau!/m_\tau!\mid B_\tau\)
also when \(\lambda_\tau<n_\#\). Also \(i_\tau-j_\tau\le
(\kappa_\tau/p-\kappa_{-p\tau})_+\) because
\(n_\#\ge x^{42}\). Summing proves \eqref{eq:marked-twelve-coeff}.

For fixed \(d\) and sufficiently large \(x\) we have \(d\le k\),
so the remaining nonsplit primes all satisfy \(p\le k\). Select the pairs with
\(2e<L_0\) whose two places have weight \(b\), with
\(e\) among the allowed positive exponents.
There are none at \(d\): both endpoints of a weight-\(b\) pair
would be \(-1\bmod q_0\), forcing \(-d^e\equiv1\bmod q_0\),
whereas \(-1\notin\langle d\bmod q_0\rangle\) by
Lemma~\ref{lem:auxiliary-residue-class}.
At inert \(p\) the partners are \((-p)^e\tau\).
Put \(D(\tau)=(-p)\tau\). Since \(p\nmid q_1\), this is a
permutation. Section~\ref{sec:datum} excludes three weight-\(b\)
places \(\tau,D^i\tau,D^j\tau\) with \(0<i<j<L_0\),
even if the places coincide. Two outgoing edges of different
exponents \(e<e'\) at the same start give such a triple with
\((i,j)=(e,e')\). For two edges ending at the same place,
equal exponents give the same start by bijectivity. If \(e<e'\),
the orbit from the start of the longer edge meets the other
start and the common end at \((e'-e,e')\). Finally two
successive edges give \((e,e+e')\), with \(e+e'<L_0\);
this includes a loop followed by itself.

If \(-p^{e'}\sigma=u_\tau\) for a selected start and an odd
\(0<e'\le e_\tau\), this too is a selected edge, since
\(2e'<L_0\) and both endpoints have weight \(b\).
Uniqueness at the common end gives
\((\sigma,e')=(\tau,e_\tau)\). In particular a shorter
return from \(\tau\) would give the prohibited triple
\((\tau,u_\tau,u_\tau)\) at exponents \((e',e_\tau)\).
Take \(j_\tau=\lfloor\lambda_\tau/p^e\rfloor\) on
these pairs; their ends have degree at least \(b\ge j_\tau\).
Use the second coupling bound and
\eqref{eq:marked-twelve-six} on the remaining degrees.
At a place write \(l_1+l_2+l_3=\lambda_\tau\) for the
starting, ending and residual degrees. The factorial loss is
\[
 v_p(\lambda_\tau!)-\sum_{i=1}^3v_p(l_i!)
  =\sum_{h\ge1}\left(\left\lfloor\lambda_\tau/p^h\right\rfloor
       -\sum_{i=1}^3\left\lfloor l_i/p^h\right\rfloor\right)
  \le 3(\log b+2).
\]
Indeed each summand is at most two and vanishes for
\(p^h>\lambda_\tau\); here \(\lambda_\tau\le 2b\) and \(p\ge3\).
At the ramified prime \(d\) the selected set is empty and
the six-copy general bound applies to the entire residual degree.

For each \(\tau\),
\[
 R_\tau=\sum_{e\text{ allowed}}\min(\kappa_\tau/p^e,\kappa_{-p^e\tau})
\]
has its selected term at most \(j_\tau+1\);
all other terms with \(2e<L_0\) are at most \(a/x_1\),
and the tail is at most \(2b/(kp)\).
For the tail, \(p^{\lceil L_0/2\rceil-1}\ge2^s\ge k\).
For inert \(p\) the allowed exponents are positive and odd,
and for \(p=d\) they are all positive. In both cases
\(v_p(\pi)\sum_{e\ \text{ allowed}}p^{-e}=\gamma_p\), so subtracting
the positive parts defining \(P_p\) gives
\[
 \gamma_p\sum_\tau\kappa_\tau-P_p=v_p(\pi)\sum_\tau R_\tau .
\]
Writing \(D_\tau=v_p(\lambda_\tau!)-\sum_i v_p(l_i!)\)
for the three-factor loss at \(\tau\), the two six-copy bounds give
\[
 q_{12}(\lambda)+\sum_\tau v_p(\lambda_\tau!)-P_p
 \le \sum_\tau\bigl(\gamma_p(\lambda_\tau-\kappa_\tau)
            +v_p(\pi)R_\tau+D_\tau\bigr)-J_c .
\]
The selected contributions to \(\sum_\tau v_p(\pi)R_\tau\)
are at most \(J_c+f\), as \(v_p(\pi)\le1\). The remaining
contributions and the factorial errors are bounded place by place by
\[
 3(\log b+2)+x^{42}+1+L_0a/x_1+2b/(kp)
       \le3b/(kp)\le v_p(B_\tau).
\]
Here the scales give
\(3(\log b+2)+x^{42}+1\le a/x_1\),
\((L_0+1)a/x_1\le b/k^2\), and \(b/(kp)\ge1\).
The last inequality follows by counting at least
\(\lfloor n_\#/p\rfloor\) multiples of \(p\) in \(B_\tau\).
This proves \eqref{eq:marked-twelve-coeff} also here.

\smallskip
\noindent\textit{The exact geometric section and finite cells.}
At a coefficient furnished by \eqref{eq:marked-twelve-coeff},
expand each \(H\) in monomials in \(y(\mu_i Z_c)\).
These series have zero constant term, so only the finitely many
monomials of parameter degree at most \(|\lambda|\) contribute.
Their scalar coefficients belong to \(R\). By the ultrametric
inequality, some product of unweighted monomials has a coefficient
of the required multi-index with no greater valuation. Thus the
ensuing geometric section requires only the selected integer
exponents.
Writing their exponents \(u_{i,c}\), the corresponding
curve sections are
\[
       (y\theta^3)^{u_{i,c}}\theta^{c_i-3u_{i,c}} .
\]
The powers are nonnegative since \(3|\lambda|\le6fb\le c_i\).
The sections have degree \(c_i\) and integral twisted
torsion values. Evaluate them in the Fock coordinates
at the forms \(\Omega_0 w_i^\dagger/(\bar\delta\,a_d)\)
on the corresponding coordinates. Their weighted pullback
Hermitian forms sum to the original Fock form by the
weighted Gram identity and \(|\bar\delta|^2=M\).

In the real coordinates \(r(\lambda)_{\tau,c}=a_d\tau(\lambda_c)\),
trace self-duality gives
\[
 A_{i,c}(r(\bar\delta\lambda))
  =\Omega_0\operatorname{Tr}_{E/K}(\lambda_c\bar w_i)
  \in\Omega_0\mathcal O_K\qquad(\lambda\in\Lambda_{12}).
\]
Write \(H_C(l,z)=\pi\bar lz/A_C\) and let \(H_{\rm Fock}\)
denote the sum of the pullbacks of \(c_iH_C\). The sign
\(\chi_C(l)=\exp(i\,\operatorname{Im}H_C(\bar\beta l,l))\)
computed for the chosen two-torsion point gives
\[
 \prod_{i,c}\chi_C(A_{i,c}L)^{c_i}
   =\exp\!\bigl(i\,\operatorname{Im}H_{\rm Fock}
                       (\bar\beta L,L)\bigr).
\]
This is exactly the prescribed signed period law on
\(\bar\delta\Lambda_{12}=\pi_M\Lambda_{12}\).
To realize this section in the finite Weil model, put
\(L=\Lambda_{12}\). Self-duality and
\(E_{\rm alt}(\bar\delta x,y)=E_{\rm alt}(x,\delta y)\)
for the alternating polarization give
\((\bar\delta L)^\#=\delta^{-1}L\). The index
\([L:\bar\delta L]=M^{12f}\) is the product of the polarization
elementary divisors on \(\bar\delta L\).
Indeed the polarization kernel is
\((\bar\delta L)^\#/\bar\delta L\), of order
\[
 [\delta^{-1}L:L][L:\bar\delta L]=M^{24f}.
\]
For the positive polarization in characteristic zero, the dimension
of the space of sections is the square root of this order
\cite[\S1, (II)--(IV), p.~289]{MumfordTheta}.
For \(u\in\delta^{-1}\widehat L/\widehat L\), the translates
\(T(u)v\) have this multiplier. The signed translations by \(L\) give them distinct
commutator characters by the dual-lattice equality. Their entire
images are nonzero by injectivity of the finite theta realization in
Section~\ref{sec:marked-integral-sections-ordinary-jets}, so they form
a basis. The section is therefore the image of a finite-part vector
\(F_*\). At \(\pi_Ms_{\rm jet}\) its forms are
exactly \(\mu_i\) because \(\bar\delta=\pi_M\).
A finite translate of \(F_*\) can be read after rational
approximation modulo its open stabilizer. Each factor
then evaluates, up to a root of unity, at the twisted value at
\(-A_{i,c}r(\gamma)\) for the rational displacement \(\gamma\).
Choose an integer \(b_0\ge1\) with \(b_0\gamma\in\Lambda_{12}\).
Since \(M=\delta\bar\delta\),
\[
 M b_0 A_{i,c}r(\gamma)
  =\Omega_0\delta\operatorname{Tr}_{E/K}(b_0\gamma_c\bar w_i)
  \in\Omega_0\mathcal O_K .
\]
The twisted torsion values of both \(\theta\) and \(y\theta^3\)
are integral by Section~\ref{sec:marked-integral-sections-ordinary-jets}.
Thus each theta displacement evaluation of \(F_*\) at real
argument zero and arbitrary finite displacement is integral.

\subsubsection{The finite coefficient pairing}
\label{sec:twelve-coordinate-finite-pairing}

\sectionopening{In this subsection, we pair the resulting coefficient of Section~\ref{sec:twelve-coordinate-inert-expansion} against the finite test.}

Put \(L=\Lambda_{12}\), \(A=L/\bar\delta L\) and
\(B=\delta^{-1}L/L\). The Weil commutator pairs \(A\) and \(B\)
perfectly by \((\bar\delta L)^\#=\delta^{-1}L\) and \(L^\#=L\).
Both groups have order \(\mathrm N_{E/\mathbf Q}(\pi_M)^{12}=M^{12f}\).
Their rational representatives give all classes of the completed
lattices, since completion preserves finite quotients.
Write \(y=\alpha_{\rm W}^{-1}a\) and \(x=\alpha_{\rm W}^{-1}u\). With the
translations of Section~\ref{sec:upper-cm-jet-transport}, put
\[
 S(a)=T_\sharp(y)=e_0(y)\rho_f(a),\qquad T(u)=\rho_f(u).
\]
For \(a,a'\in L\) the normalized coordinates \(y,y'\) are integral.
The multiplier \(\psi_0(\bar c_0(y')^*y)\) is therefore one, including
at places above \(2\). Thus \(S(a)S(a')=S(a+a')\).
Rational Heisenberg invariance identifies the signed period law
on \(\bar\delta L\) with \(S(a_0)F_*=F_*\) for \(a_0\in\bar\delta L\).
The ordinary displacement satisfies
\[
 T(u)S(a)T(-u)=\psi_0(y^*x)S(a).
\]
This is the commutator pairing just described. In particular
\(\bar\delta L\) fixes \(T(-u)F_*\), so the sum over \(A\) is
independent of the lifts of its elements. For each \(u\), the sum
\[
 d_u v=\sum_{a\in A} S(a)T(-u)F_*
\]
satisfies the signed period law on \(L\), whose theta space
is the vacuum line. Conjugation by \(T(u)\) multiplies \(S(a)\)
by the commutator character at \(u\). Its sum over \(B\) is zero
for \(a\ne0\) and \(|B|\) for \(a=0\). Thus
\[
 M^{12f}F_*=\sum_{u\in B}d_u\,T(u)v .
\]
If the vacuum vanished on the whole finite slice, rational Heisenberg
translation would give \(v(r_\infty,z_f)=0\) for every rational
\(r\in E^{12}\) and every finite argument \(z_f\): translation by
\(-r\) takes the real argument to zero and changes the value by a
nonzero phase. Since \(E^{12}\) is dense in its archimedean vector
space, continuity would give \(v=0\). Thus the finite slice has a
nonzero value.

Fix such a point \(z_f^*\) and put
\(C_v=V_0^{12}v(0,\alpha_{\rm W}z_f^*)\ne0\).
Local constancy and finite adelic additive approximation allow the
finite translate at this point to be represented by a rational one.
Its normalized value is therefore a degree-zero translated CM
coefficient, fixed before the coefficient embedding.
At this point \(d_u v(0,\alpha_{\rm W}z_f^*)\) is a sum of
finite evaluations of \(F_*\) multiplied by roots of unity.
It is therefore integral also at the primes dividing \(M\).
In the normalized finite coordinates \(x=\alpha_{\rm W}^{-1}u\)
we have \(\delta x\) integral. The nonintegral components of
\(x\) are at the conjugates of the places of \(\pi_M\).
Subtracting the other components acts on the vacuum by
a root-of-unity phase. Take \(Z\) in the chosen Schr\"odinger
lattice. The commutator above makes \(T(u)v\) an eigenvector
of the integral position translation by \(Z\), with eigenvalue
\(\psi_0(-Z^*x)=\psi_0(x^*Z)\). Thus on that lattice its
additive multiplier relative to the vacuum is \(\psi_0(x^*Z)\),
up to the fixed phase from the discarded components.
It factors through reduction modulo \(\pi_M\) since
\(\pi_M\bar x=\overline{\delta x}\) is integral. The perfect
pairing of \(A\) and \(B\) gives all the residue characters.
For the original residue cells write
\(T(u)v=\sum_x\omega_u(x)I_x\), where each nonzero \(\omega_u(x)\)
is a root of unity. Set \(v_*=v(0,\alpha_{\rm W}z_f^*)\) and
\(a_x=\sum_{u\in B}(d_u v_*)\omega_u(x)\). Each product
\(d_u v_*\) is a sum of the integral algebraic torsion values above
with root-of-unity coefficients, so \(v_p(a_x)\ge0\). Multiplying
the Fourier identity by \(v_*\) and applying the scaled coefficient
functional gives
\[
 M^{12f}C_v[Z^\alpha]\mathcal T F_*
     =\sum_x a_x\,V_0^{12}[Z^\alpha]\mathcal T I_x .
\]
The jet of \(F_*\) on the left is the chosen coefficient of the CM
section in the algebraic forms \(\mu_i\). It is nonzero. Each
scaled cell coefficient on the right is algebraic by
Lemma~\ref{lem:normalized-elliptic-jet-base-change}. Applying the
coefficient embedding to this identity and taking the minimum
valuation on the right gives a cell \(x_0\) with
\[
 v_p(V_0^{12}[Z^\alpha]\mathcal T I_{x_0})
 \le v_p(M^{12f}C_v)+v_p([Z^\alpha]\mathcal T F_*).
\]
Set \(t_0=M^{12f}C_v\prod_\tau B_\tau^2\).

It is fixed before the coefficient prime and embedding, and nonzero.
Put \(C_{\rm int}=2+\log M_K\) and take
\(C_P=\max\{C_K,C_{\rm int}+2d+6\}\), with the constants of
Lemma~\ref{lem:upper-gaussian-kernel}. The rank-one conjugate
vacuum bound there, applied on the twelve coordinates, gives
\[
 \log|\sigma t_0|
 \le12f\log M+12C_Pq_1^2+
             2f n_\#\log(2b+n_\#)\le fmb/k .
\]
Indeed \(n_\#\le 7b/k\) eventually and
\(\log(2b+n_\#)\le B+\log 3=o(m)\).
Also \(\log M=O(B)\), \(q_1\le x^{12}\), and \(b\) exceeds every
fixed power of \(x\). These bounds hold at one cutoff for the
fixed \(q_0,d,\epsilon\), independently of the coefficient prime.
Now \eqref{eq:marked-twelve-coeff} proves
\eqref{eq:marked-rank-twelve-sharp}. The polynomial
degree, section, and finite Fourier steps all use the
same \(\Lambda_{12},\alpha_{\rm W},\pi_M\) as the numerator test.

\smallskip
\subsubsection{CM transport of the untwisted numerator}\label{sec:rank-one-cm-transport}

\sectionopening{In this subsection, we apply Lemma~\ref{lem:upper-cm-transport} to the coefficient of the original numerator test.}

Fix \(\sigma\in\operatorname{Aut}(\mathbb C/K)\).
Instantiate the exact translated jet identity \eqref{eq:upper-CM-jet-transport}
and the vacuum identity \eqref{eq:upper-CM-vacuum} on the standard Hom
space of rank \(2z\). Their proof depends on the rank
only via the coordinate product, so the same construction
uses the normalized scale \(t_\sigma\) with
\(|t_\sigma|^2=\Omega^{-1}\).
Keep the exact identities \eqref{eq:upper-finite-vacuum-action} and the finite multiplier
law, rather than their subsequent bounds on coefficient sizes.
Write the resulting scalar translate as \(b_\sigma=u_\sigma a_\sigma\)
where \(u_\sigma\bar u_\sigma=1\) and \(a_\sigma\) is an
integral unit equal to \(1\) on the old chosen places.
This is possible because \(b_\sigma\bar b_\sigma\) is a unit.
In every split pair choose one member \(w\), using the prescribed
member at each old pair. Set \(a_{\sigma,w}=1\) and
\(a_{\sigma,\bar w}=\overline{(b_\sigma\bar b_\sigma)_w}\).
At a conjugation-fixed place \(b_\sigma\) itself is a unit,
so use it as \(a_\sigma\).

In the untwisted column datum the finite coefficient function
is the original chosen-column indicator and has values \(0,1\).
After removing the scalar action of \(u_\sigma\), the exact
semilinear operator transports each signed translate by
\(x\) to the one by \(a_\sigma x\), and multiplies the vacuum
by one fixed nonzero scalar \(c_\sigma\).
Let \(D\) be the additive residue matrix space at the chosen
places, and let \(Y_0\in D\) be the prescribed inclusion. On the
integral box its indicator has the expansion
\[
 \mathbf1_{Y=Y_0}=|D|^{-1}\sum_{\chi\in\widehat D}
                         \chi(-Y_0)\chi(Y).
\]
The modulations representing \(\chi(Y)\) are supported at the
opposite places. Since \(a_\sigma=1\) on the chosen side,
\eqref{eq:upper-finite-vacuum-action} changes their characters to
\(\sigma\chi\). Semilinearity changes their coefficients to
\((\sigma\chi)(-Y_0)\). Permuting the characters recovers the same
indicator times \(c_\sigma\). At the other places there is only
the vacuum.

Write \(b_X=\nu_X\upsilon_1^z\), \(b_Y=\nu_Y\upsilon_1^z\),
where \(\upsilon_1=\upsilon|_{U(1)}\), and let \(g'\) denote
the inverse source action on the \(2z\) Hom coordinates.
For a finite norm-one element \(g\), put \(D_1=\Delta_1(g^{-1})\)
for the one-coordinate multiplier. The squared-multiplier identity
in Section~\ref{sec:upper-cm-jet-transport} gives
\(D_1^2=|D_1|^2\upsilon_1(g)\), where \(|D_1|^2\) is a rational
lattice index. The \(2z\) identical coordinates give
\(\Delta(g')=D_1^{2z}=|D_1|^{2z}\upsilon_1(g)^z\), whence
\[
 d_\sigma(g')=\frac{\sigma(\Delta(g'))}{\Delta(g')}
     =\left(\frac{\sigma(\upsilon_1(g))}{\upsilon_1(g)}\right)^z .
\]
After writing \(\nu_Y(g)=\sigma(\nu_X(g))\) on the finite
ideles, the source factor in a corresponding kernel term is
\[
 \sigma\bigl(b_X(g)^{-1}\bigr)d_\sigma(g')
   =b_Y(g)^{-1}
   =b_Y(u_\sigma^{-1})\,b_Y(gu_\sigma^{-1})^{-1}.
\]
The last inverse character belongs to the transported kernel
at \((gu_\sigma^{-1},h)\). Multiplication by \(u_\sigma^{-1}\)
permutes the classes of the norm-one torus at the same old
level. Summing over this permutation leaves
\(b_Y(u_\sigma^{-1})\) outside the class sum.

Let \(X\) be the original untwisted datum and \(Y\) its transport.
Write \(X^{\rm raw},Y^{\rm raw}\) for their unscaled class sums,
\(J_X\) for their common total degree, and \(h'\) for the target
Hom action at \(h\).
Write \(p_\sigma\) for the permutation of the weight embeddings
given by \(\sigma\Phi_\tau=\Phi_{p_\sigma(\tau)}\).
Relabel \(Z\) by \((R_\sigma Z)_{\tau,j}=Z_{p_\sigma^{-1}(\tau),j}\),
and put \(i^\sigma_\tau=i_{p_\sigma^{-1}(\tau)}\).
The Schur coordinates are products of polynomials with integer
coefficients in these variables. If \(P_{X,i}\) and
\(P_{Y,i^\sigma}\) are the corresponding coordinate polynomials,
then \(P_{Y,i^\sigma}(R_\sigma Z)=P_{X,i}(Z)\).
Their total degree \(J_X\) gives the scalar \(t_\sigma^{J_X}\).
With \(A_X=s_{\rm jet}^{J_X}X^{\rm raw}\), the exact transport
at the unchanged target point \(h\) is
\[
 \sigma(A_X(h)_i)=t_\sigma^{J_X}c_\sigma
           b_Y(u_\sigma^{-1})\,
        \frac{\sigma(\Delta(h'))}{\Delta(h')}\,
                           Y^{\rm raw}(h)_{i^\sigma}.
\]
The scalar \(c_\sigma\) and the Taylor scale were fixed before
the character and degree. This proves the transport identity
needed for the mixed pairing in the untwisted reference norm.
For the twisted numerator the two point tests are constructed
separately at \(\kappa\) and \(\bar\kappa\). The twelve-coordinate
cell is evaluated at the identity, and the other block carries
the finite twist in \eqref{eq:display-nb2}. This keeps the same
\(\eta\), multiplier \(t_0\), and main line at both embeddings.

\subsection{The four-block test and the untwisted reference norm}
\label{sec:v2-four-block-main-numerator}

\sectionopening{In this subsection, we insert the coefficient of
Section~\ref{sec:marked-twelve-sharp-coefficient} in the original
four-block theta test. The resulting numerator is measured against the
untwisted norm on the same space.}

Fix a main pair \(p=(\nu,z)\), its twist \(\eta\) of
exact support \(Q\), and the original lift \(X\) from \(U(1)\)
to \(U(2z)\). Let \(C\) be the old \(U(1)\) class group,
\(h_1=|C|\), and \(V_0=s_{\rm vac}\).
Use its original ordered blocks of sizes \(12,4,4,2z-20\).
Write \(\lambda_\tau=c_\tau-z\) for the part sizes and put
\(J=(\pi_M s_{\rm jet})^{\sum_\tau\lambda_\tau}\).
The normalized main lift is
\[
 \phi_\eta=V_0^{2z}J\sum_{g\in C}K_X(g),
\]
where \(K_X(g)\) is the raw kernel form. The class
stabilizers in this source sum are trivial.

\smallskip
\subsubsection{The derivative estimate}\label{sec:four-block-derivative}

\sectionopening{In this subsection, we use the coefficient bound of Section~\ref{sec:marked-twelve-sharp-coefficient} for the point test.}

The required scalar \(t_0\ne0\) depends on the part sizes
and the fixed datum, with \(\log|\sigma t_0|\le fmb/k\) for
every complex automorphism. At each coefficient embedding
\(\kappa=(\ell,\iota)\), the sharp derivative bound gives an
exponent \(\alpha\) on the first twelve coordinates with
\(\sum_i\alpha_{\tau,i}=\lambda_\tau\), and a cell \(x_0\)
in the original residue space modulo \(\pi_M\), such that
\[
 \left| \frac{V_0^{12}}{t_0}
        \left(\prod_\tau\lambda_\tau!\right)
        [Z^\alpha]\mathcal T I_{x_0}\right|_\kappa
               \ge \ell^{-P_\ell},\qquad
 P_\ell=\sum_\tau S_\ell(\kappa_\tau,\tau).
 \tag{NB1}\label{eq:display-nb1}
\]
Here \(I_{x_0}\) is the point indicator at the chosen places
and the original vacuum elsewhere, and \(\mathcal T\) reads
the Taylor series at scale \(\pi_Ms_{\rm jet}\).
By Lemma~\ref{lem:normalized-elliptic-jet-base-change}, this
coefficient is the image under \(\iota\) of the normalized
algebraic jet. The degree-zero assertion for a translated
vacuum gives \(C_v\in\overline{\mathbb Q}^{\times}\) and
hence \(t_0\in\overline{\mathbb Q}^{\times}\).
For a main pair the defining inequalities of
Section~\ref{sec:critical-hecke-systems} give
\(w_j(\tau)+r-t\le c_\tau\le w_j(\tau)\) and \(m<z\le r\).
Since \(w_j(\tau)=\kappa_\tau+(2j+1)r^8\) and \(t=5r\), we have
\[
 \kappa_\tau\le w_j(\tau)-t\le c_\tau-z
            \le w_j(\tau)\le\kappa_\tau+x^{42}.
\]
The last inequality follows for all \(1\le j\le m\) from
\((2m+1)r^8=O(x^{81/2})\) by
Section~\ref{sec:marked-level-prime-selection}.
Thus Section~\ref{sec:marked-twelve-sharp-coefficient}
applies at the actual degree \(\lambda_\tau=c_\tau-z\);
it fixes \(t_0\) before the coefficient embedding and gives
\eqref{eq:display-nb1} on the original residue space.

\smallskip
\subsubsection{The point test on the original target}\label{sec:four-block-point-test}

\sectionopening{In this subsection, we insert the derivative of Section~\ref{sec:four-block-derivative} in the test on four blocks.}

Use the empty rank-four pair \(p_{\rm aux}=(\chi_M\upsilon_1^{-2},2)\)
in the second block at twist \(\eta\), and two untwisted empty
partitions in the third and fourth blocks. Write \(\upsilon_1\) for the restriction to the norm-one torus.
For an even rank \(2z'\) and an old weight-zero character \(\psi\),
the empty partition is the pair \((\psi\upsilon_1^{-z'},z')\), with
rank-twisted character \(\psi\). These even ranks lie
between \(4\) and \(t\), since \(m<z\le r\). The character
on the first block is the quotient of the target
character of \(X\) by the other three characters.

Set \(T_\alpha=\prod_\tau\lambda_\tau!/\alpha!\), with
\(\alpha!=\prod_{\tau,i}\alpha_{\tau,i}!\).
Let \(P_i\) be the dual Schur basis of products of minors
for this weight. Its column minors have size one.
In the original Weyl coordinates take \(v_i=[Z^\alpha]P_i\),
so \(v\) is integral. The Fock polynomial \(\mathcal P(v)\)
in these coordinates satisfies
\[
 \langle\mathcal P(v),P_i\rangle_{\rm Fock}
       =\left(\prod_\tau\lambda_\tau!\right)v_i,\qquad
 \langle Z^\alpha,P_i\rangle_{\rm Fock}=\alpha![Z^\alpha]P_i .
\]
The first identity is the Weyl to dual Schur pairing
with the hook product of the one-row partition. The second follows by
monomial orthogonality. The pairing with all \(P_i\)
therefore gives \(\mathcal P(v)=T_\alpha Z^\alpha\). The contragredient
substitution fixes the coefficient functional when the
first twelve columns are identity columns. Equivalently
the transpose substitution fixes its monomial in the
Fock pairing. This gives invariance of \(v\) under the
last three embedded groups. There is no determinant
factor in this Weyl weight.

At the old chosen places write \(Y=(x_0,e_2,e_3,e_1)^{\mathsf t}\),
with \(e_i\) the included unit vector in the indicated
empty block. Write \(e_{\rm old}\) for the original included
vector of \(X\) and \(e_{\rm second}\) for the included
vector of the second block. Choose the reduction of the
integral point lift \(A\) by
\[
 A e_{\rm old}=Y \quad(w\in S_M),\qquad
 e_{\rm old}^{\mathsf t}A^{-1}=e_{\rm second}^{\mathsf t}
                      \quad(w\in Q).
\]
At an old place \(Y\) is nonzero by its last component
and extends to a basis. At a new place the prescribed
nonzero covector extends to an invertible matrix of
covectors. Split reduction gives the local integral lifts.
They define \(A\) by taking the identity at all remaining finite
places of \(F_0\).

Put \(x=\operatorname{diag}(1,h_2,h_3,h')A\) and let
\(J_Q,\chi\) be the target level and its character.
Write the covariance of a form as
\(f(\gamma yk)=\chi(k)\rho(\gamma)f(y)\), where
\(\Gamma\) is the rational group and \(\rho\) its
coefficient action. Define \(F(h',h_2,h_3)=F_x(v)\) by
\[
 F_x(v)(y)=\sum_{\substack{\gamma\in\Gamma\\
                                  x^{-1}\gamma y\in J_Q}}
        \chi(x^{-1}\gamma y)\rho(\gamma^{-1})v .
\]
The sum is finite by definiteness and discreteness.
It has the stated covariance by changing variables
and is supported on the class of \(x\).
For \(k=x^{-1}\gamma y\), transporting the summand to
the lattice norm at \(y\) acts on \(v\) by \(k^{-1}x^{-1}\)
and a root of unity. The inverse block diagonal fixes \(v\).
Both \(A^{-1}\) and \(k^{-1}\) preserve the local Weyl lattice
at every coefficient place above \(\ell\). Thus each summand has
lattice norm at most one there, as does their sum.

At \(x\) the sum is over \(S_x=\Gamma\cap xJ_Qx^{-1}\).
For \(\gamma\in S_x\), put
\(T_\gamma=\chi(x^{-1}\gamma x)\rho(\gamma^{-1})\).
It is unitary in the original definite pairing and
fixes \(f(x)\) by covariance. Since the class weight
at \(x\) is \(|S_x|^{-1}\), the pairing gives
\[
 \langle F_x(v),f\rangle_{J_Q}
   =\frac{1}{|S_x|}\sum_{\gamma\in S_x}
                       \langle T_\gamma v,f(x)\rangle
   =\langle v,f(x)\rangle .
\]
Move \(A\) from the kernel argument into the finite datum
by the split Weil action. The determinant of \(A\) occurs in the target variable
with the character \(\chi_W=1\) of the definite source line.
The local lattice action thus introduces no determinant scalar.
The datum now acts on \(A^{-1}\) times the residue column.
On \(\pi_M^2Q\) the two displayed residue conditions give
the prescribed column indicator at \(S_M\) and the
tame character of the second-block covector at \(Q\).

This is exactly the product of the four
block functions, pulled back to that modulus. The
archimedean polynomial is \(T_\alpha Z^\alpha\) on the
first block and \(1\) on the other three. The rational
theta sum and the Weil character both factor across
these orthogonal blocks. If \(K_i\) are the raw empty
kernels and \(K_0\) is the selected first kernel, this
proves
\[
 \langle F(h',h_2,h_3),K_X(g)\rangle
   = T_\alpha K_0(g)K_2(g,h_2)K_3(g,h_3)K_1(g,h').
 \tag{NB2}\label{eq:display-nb2}
\]
In this formula the kernel is in the linear leg of the
pairing. The same \(g\) serves every block.

For an empty block \(p_i\) of rank \(a_i\), let
\(\phi_i=V_0^{a_i}\phi_i^{\rm raw}\). The denominator
\(d_i=d_{p_i}=M^{2fa_i}\) is the denominator of
Section~\ref{sec:marked-norm-line} at degree zero.
The local jet degree vanishes at each infinite place because the
polynomial of this block is \(1\).
Put
\[
 A_i=C_Gd_i^2\|\phi_i\|^2/Q_i^{a_i-1},
 \quad Q_2=\mathrm NQ\text{ for the twisted block},
 \quad Q_1=Q_3=1.
\]
Here \(Q_i\) in this display denotes only the norm of the
auxiliary modulus on that block. The class averaging
operator on \(F\) weights its value at \(h\) by
\[
 \frac{C_Gd_i^2}{|\operatorname{Aut}(h)|\,Q_i^{a_i-1}}\,
          \langle 1,\phi_i(h)\rangle_{\rm Weyl}.
\]
These weights are integral at all the embeddings above
\(\ell\), since \(\ell\nmid\mathrm NQ\), \(d_i\) clears
the normalized empty lift and \(C_G\) clears the class
weights.

The weighted point form is independent of the empty-block
representatives. For a right level change \(h_i\mapsto h_i k_i\)
the new point is \(xk\), with \(k=A^{-1}\operatorname{diag}(k_i)A\)
and the other diagonal blocks unchanged.
At the old places \(\operatorname{diag}(k_i)\) fixes \(Y\)
in the residue space. At the new places the defining covector is
on the second block, so \(k\) belongs to \(J_Q\) and its character
is the character of the changing block. The point form acquires
\(\chi(k)^{-1}\) by its definition and the averaging weight
acquires \(\chi(k)\). Left rational changes act on the point
vector through the block action fixing \(v\). The stabilizer
weights therefore apply to the actual classes even when several
choices give the same target class.

In the pairing against \(K_X(g)\) the weights
are conjugated. They replace \(K_i(g,h)\) in \eqref{eq:display-nb2} by
\(A_i/(h_1 V_0^{a_i})\). To see this, a normalized
central translation by \(u^{-1}\) sends \(K_i(g)\)
to \(K_i(gu)\): the simultaneous scalar cancels on
the Hom variable, leaving the inverse source twist.
This translation is unitary and fixes
\(\phi_i^{\rm raw}=\sum_g K_i(g)\), as right translation
permutes the source classes. It follows that
\(\langle\phi_i^{\rm raw},K_i(g)\rangle\) is independent
of \(g\); summing over \(g\) gives its value
\(\|\phi_i^{\rm raw}\|^2/h_1\).
Combining this with the conjugated scalar
\(\overline{V_0^{a_i}}\) in the weights gives precisely
\(A_i/(h_1V_0^{a_i})\).

\subsubsection{The last-block class sum}\label{sec:four-block-class-sum}

\sectionopening{In this subsection, we evaluate the last block of the test in Section~\ref{sec:four-block-point-test}.}

Apply the measured period argument in
Section~\ref{sec:marked-comparison-period-vanishing}
to the original rank-one descent datum \(O\) of the
four-block lift \(X\). The source has rank \(1\), the
target has rank \(n=2z\), and the two blocks for this
application are the retained first \(12+4+4=20\)
coordinates and the last \(n-20\) coordinates. Thus
the dimension condition is \(2+2\cdot20+3=45<n\), since
\(z>m\) and \(m\ge22\) for all sufficiently large \(x\).
At infinity the line parameter is zero at every
place, the determinant exponent is zero, and the
target space of length \(\lambda_\tau\) is matched
by the homogeneous Fock polynomials of that degree.
Indeed the inverse source scalar acts on that homogeneous
degree by its \(-\lambda_\tau\)-th power, and the inverse
character of source weight \(-\lambda_\tau\) cancels it.
The positive Fock variables give the one-row partition
\((\lambda_\tau,0,\ldots,0)\), with no negative-sector
columns or mixed Laplacians. In each degree this is the irreducible
symmetric power of the compact target, including its trivial
representation in degree zero. The source circle has different
characters in different degrees, giving the required correspondence.

Fix a split \(w\) with \(w,\bar w\) off the original
chosen and twist places. For this test a false lift
\(y\) lies in
\[
 \operatorname{im}(e_{X-1})\cap
       \ker\mathcal P_{n,w}(X_*)\cap\ker D_O ,
 \qquad
 X_*=\bigl(q_w^{n-2}\operatorname{frob}_1(1,w)\bigr)^{-1},
\]
where \(D_O\) is the normalized descent for \(O\), and
\(\operatorname{frob}_1(1,w)\) is the rank-one
unramified value of the trivial \(\chi_W\)
with the original splitting character included.
If \(\Upsilon\) is the original Hecke splitting
character of \(E\), then
\(\operatorname{frob}_1(1,w)=q_w^{1/2}\Upsilon(\varpi_w)^{-1}\).
In particular \(|X_*|=q_w^{3/2-n}\).
The conjugate matrix coefficient of \(y\) is killed
at \(\overline{X_*}\), as in \eqref{eq:display-f1}. The coefficient
operator here is the one of this rank-one datum,
with \(\chi_W=1\), even when \(\eta\) is
nontrivial in the other character.

The idempotent \(e_{X-1}\) is the product of the means
for the normalized actions of the old groups
\(\mathrm{GL}_1(k_v)=k_v^\times\).
Its image consists of the invariant vectors for these actions.
The raw stabilizer shear includes the factor \(\chi_v(s)^{-1}\)
in the calculation of Section~\ref{sec:bottom-descent-real-type}.
If \(\lambda\) is a nonzero residue covector at a chosen
place \(v\), the rank-one set \(1+u\lambda\) with
nonzero determinant is all of \(k_v^\times\).
The action of its hyperplane sum on \(y\) is therefore
\(\left(\sum_{s\in k_v^\times}\chi_v(s)^{-1}\right)y=0\).
Here \(\chi_v\) is the exact old main-pair tame character,
which is nontrivial, while \(\eta\) has trivial old tame
components.

The same weighted sum in the product of
old local groups is zero if any one of its covectors is
nonzero. At the new places the exact-support
hypothesis gives nontrivial residue characters.
The bottom descent vanishing for the original block is
the condition \(D_O y=0\). These are the three conditions for
the bottom extension to arbitrary Schwartz data vacuum
off the original support. They replace the rank-\(r\)
argument with selected Green parameters and the orthogonal
complement of the projected image.

For the following tower argument, the source at step
\(l\), \(1\le l\le20\), has rank \(1+2l\).
With \(k_l=n-(1+2l)\), the estimates of \eqref{eq:display-f3} read
\[
 q_w^{k_l}|\alpha X|=q_w^{1-l}\le1,\qquad
 q_w^{\,1+2l+1}|X'|=q_w^{(6+2l-n)/2}\le1,
\]
where \(X=\overline{X_*}\). The latter follows from
\(n\ge46\). Hence the same hyperbolic-plane induction
gives vanishing of all descents through twenty planes,
and the argument from \eqref{eq:display-f4} gives the measured
second-block periods with retained dimension \(20\).
The bottom and tower arguments use the fixed \(O\),
fixed target point, original finite vacua and original
right-coset operators throughout.

For the application to the test, take the fixed
first-block point \(\operatorname{diag}(1,h_2,h_3)\)
of the retained twenty coordinates. Take the additional
right translate to be \(1\). The integral point lift \(A\)
in \(F_x(v)\) then gives the period point
\(\operatorname{diag}(1,h_2,h_3,h')A\), with finite
translation supported at the old and auxiliary split places. The
last-block kernel \(K_1(g,h')\) is the degree-zero kernel,
whose polynomial is \(1\), whose target dimension is
\(n-20\), and whose rank-twisted character is the
untwisted old factor \(\psi_1\).
Write \(b_X=\nu\eta\upsilon_1^z\) for the rank-twisted character
of the main pair. For its original Schwartz datum the kernel identity
\[
 K^{1,1}(g,h)=b_X(\det g)K^{1,b_X}(g,h)
\]
multiplies the descent at fixed \(g\) by a nonzero scalar, conjugated
when the kernel is in the conjugate leg. Thus \(D_Oy=0\) gives the
bottom vanishing with second character \(1\) on the original datum.
The old residue sums and the exact-support condition above extend this
vanishing to the Schwartz tests used in the tower.
In the block period take the last character \(b=\psi_1\) and the
auxiliary retained character \(b^{-1}\). Their product is \(1\)
for every allowed \(\psi_1\). The auxiliary retained kernel is chosen
only for this vanishing argument.

The finite datum of \(K_1(g,h')\) has the original vacuum at
the new places and elsewhere off the old support. For each fixed
source point \(g\) and \(h_{\rm ret}=\operatorname{diag}(1,h_2,h_3)\),
Section~\ref{sec:gluing-block-period-vanishing} constructs an auxiliary
twenty-coordinate kernel nonzero at \((g,h_{\rm ret})\). This scalar is
independent of \(h'\), so it cancels from the period for each false
lift with the prescribed last-block datum, before summing over
\(g,h_2,h_3\). The datum of \(K_1\) can be read on the nominal
modulus \(\pi_M^2\), with the prescribed column
indicator depending only on reduction modulo \(\pi_M\).
In the measured period use the conjugate kernel vector
\(\overline{K_1(g,h')}\) and the conjugate-linear
coefficient map \(c\mapsto\bar c\,v\).
This map is equivariant for the real last-block group:
the empty weight is trivial, and the retained coefficient
vector \(v\) is fixed by that group. Since the class
pairing is conjugate-linear in its first argument, the
resulting scalar integrand is
\[
 \overline{K_1(g,h')}\,
       \langle F(h',h_2,h_3),y\rangle .
\]
The point-form reproducing identity evaluates the second
form at the finite point
\[
          \operatorname{diag}(1,h_2,h_3,h')A .
\]

Use the actual last-block vector-fixer level \(J_1\) of
the untwisted degree-zero lift. It is integral and fixes
\(e_1\) at the old chosen places, so fixes the
last-block datum. At an old place the embedded element
fixes the prescribed column \((x_0,e_2,e_3,e_1)^{\mathsf t}\).
Consequently \(A^{-1}\operatorname{diag}(1,k)A\)
fixes the originally included column. At a new place
this block embedding fixes the prescribed twist covector
because that covector lies in the second block of the
four-block splitting. Conjugating by \(A\) therefore
fixes the originally included covector there and has level
character \(1\). Off those places the element is
integral. This is the right invariance needed for the
period-to-class-sum passage. Rational left invariance
follows from the same kernel and form covariance.

The measured period is zero. Compact integration at
\(J_1\) multiplies its class sum by the common nonzero
factor \([K_{\rm std}:J_1]^{-1}\) with maximal compact
mass one. Retaining the actual stabilizer weights gives
\[
 \sum_{[h']_{J_1}}
       \frac{\overline{K_1(g,h')}}
                  {|\operatorname{Aut}_{J_1}(h')|}
              \langle F(h',h_2,h_3),y\rangle=0 .
\]
This is for every original \(g,h_2,h_3\), before the
outer three sums. The last-block sum has its vector-fixer
level, so there is no need to assert freeness of this
class space. Expanding the degree-zero lift in its source
class sum now gives the required orthogonality of
the averaged test to \(y\).
The two remaining scalar averages preserve this
orthogonality, and all three averages preserve integrality.

Denote the test after the three averages by \(\Xi\).
Substituting the replacement identities into \eqref{eq:display-nb2} gives
\[
 h_1^3V_0^{2z-12}\langle\Xi,K_X(g)\rangle
        =T_\alpha K_0(g)A_1A_2A_3 .
 \tag{NB3}\label{eq:display-nb3}
\]

\smallskip
\subsubsection{Selection of the class character}\label{sec:four-block-class-character}

\sectionopening{In this subsection, we use the class sum of Section~\ref{sec:four-block-class-sum} to choose characters for the numerator estimate.}

Before selecting the two untwisted empty-block characters,
put
\[
 \chi_0=\frac{\nu\eta\upsilon_1^z}
   {(\chi_M\upsilon_1^{-2})\eta\upsilon_1^2}
       =\nu\upsilon_1^z/\chi_M .
\]
These are the rank-twisted characters of the main and second
blocks, respectively. The common finite twist cancels in the
quotient. At \(g=1\) every character has value one, so the
twelve-coordinate jet is exactly the untwisted cell coefficient
of Lemma~\ref{lem:normalized-elliptic-jet-base-change}. Let \(s(g)\)
be the unscaled first kernel for \(\chi_0,x_0,\alpha\).
Choose the actual source representatives \(g_a\) of \(C\),
and let \(z_0=g_{1_C}^{-1}\). The representatives
\(g_a z_0\) have the identity point in the identity class.
For \(\psi\in\widehat C\) put
\[
 B_\psi=h_1^{-1}\sum_{a\in C}
                  \psi(\det(g_a z_0))s(g_a z_0).
\]
On these representatives the class of \(\det(g_a z_0)\) is
\(a\), and \(g_{1_C}z_0=1\). Since \(h_1=|C|\), character
orthogonality gives
\[
   \sum_{\psi\in\widehat C}\psi(\det(g_a z_0))
       =h_1\mathbf1_{\{a=1_C\}},\qquad
   \sum_{\psi\in\widehat C} B_\psi=s(1).
\]
This identity uses only the values of \(s\) on the chosen
representatives. Applying the ultrametric inequality to the last
sum gives a \(\psi\) with \(|B_\psi|_\kappa\ge |s(1)|_\kappa\),
also for \(\ell\mid h_1\). The three class weights forming \(\Xi\)
are integral by the calculation preceding
Section~\ref{sec:four-block-class-sum}.

Let \(\mathcal C_0\) be the finite group of allowed old tame
weight-zero characters. Each \(\psi\in\widehat C\) belongs
to this group by inflation. Indeed, on a norm-one idele take
\(\psi\) of the class of its finite component. This character
is trivial on rational points, on the archimedean torus and on
the old principal level. Its finite unit restrictions thus factor
through the residue groups at \(S_M\) and are trivial elsewhere.
Fix \(L_{\rm geom}\ne0\) by choosing a root
\[
 L_{\rm geom}^{|\mathcal C_0|}
       =\prod_{\chi\in\mathcal C_0} A_1(\chi)A_3(\chi).
\]
These are untwisted degree-zero values on the last and third
blocks with the same old places. Proposition~\ref{prop:marked-line-norm}
at \(Q=\varnothing\) gives a positive norm for every such character,
including the trivial one: the old vector-fixer computation gives
\(q_w-1>0\) and the omitted Euler product is absolutely convergent
and nonzero at \(z'\ge2\). Section~\ref{sec:marked-both-norm-jets}
gives the denominator bound with \(q_Q=1\).
In particular \(L_{\rm geom}\ne0\). As \(\chi\) runs over
\(\mathcal C_0\), so does \(\psi\chi^{-1}\), and hence
\[
 \prod_{\chi\in\mathcal C_0}
      |A_1(\chi)A_3(\psi\chi^{-1})|_\kappa
      =|L_{\rm geom}|_\kappa^{|\mathcal C_0|}.
\]
At least one factor is at least the geometric mean. Thus
one factorization \(\psi=\psi_1\psi_3\) has
\(|A_1(\psi_1)A_3(\psi_3)|_\kappa\ge|L_{\rm geom}|_\kappa\).
Use these two characters on the corresponding blocks.
The first-block target-leg character becomes \(\chi_0/\psi\).
Its multiplier on the source variable is
\[
 ((\chi_0/\psi)(\det g))^{-1}
                 =\psi(\det g)\chi_0(\det g)^{-1}.
\]
Hence \(K_0(g)=\psi(\det g)s(g)\).

Write \(L_{\rm twist}=A_2(1)\), and let
\(\rho_\eta=A_2(\eta)/A_2(1)>0\). For the empty rank-four
pair \(p_{\rm aux}=(\chi_M\upsilon_1^{-2},2)\), the exponents
\(c_v-z\) vanish, so the denominator \(d_2=M^{8f}\) is
independent of \(\eta\). The norm identity
\eqref{eq:marked-rank-one-norm} gives
\[
 \frac{\|\phi_2(\eta)\|^2}{\|\phi_2(1)\|^2}
       =q_Q^3 R_{p_{\rm aux},+}(\eta),\qquad
 \rho_\eta=R_{p_{\rm aux},+}(\eta).
\]
On the source, rational equivariance and the trivial
source weight give left invariance of the full kernel \(K_X\).
The old principal level fixes its old column. At a new
place, multiplication of the residue covector value by the
unit character cancels the inverse character in the source
Weil action. At all remaining places the units fix the vacuum.
Thus \(K_X\) descends to the old class group \(C\), and
\(\phi_\eta\) can be evaluated at \(g_a z_0\). Sum \eqref{eq:display-nb3} there and take its linear-leg
normalization. The result is
\[
 \frac{h_1^2}{\rho_\eta}\langle\Xi,\phi_\eta\rangle
   =V_0^{12}J T_\alpha L_{\rm twist}
                              A_1A_3 B_\psi .
 \tag{NB4}\label{eq:four-block-fourier-average}
\]
Here the third power of \(h_1\) in \eqref{eq:display-nb3} loses one
power in passing to the Fourier average.

No complex
conjugate of \(V_0\) or \(J\) occurs in \eqref{eq:four-block-fourier-average}.
The coefficient identity at \(g=1\) is
\((J/\alpha!)s(1)=[Z^\alpha]\mathcal T I_{x_0}\),
since differentiation multiplies the Taylor coefficient
by \(\alpha!\).
In particular \(JT_\alpha s(1)=
(\prod_\tau\lambda_\tau!)[Z^\alpha]\mathcal T I_{x_0}\).
Using \(v_\kappa(B_\psi)\le v_\kappa(s(1))\) and
\(v_\kappa(A_1A_3)\le v_\kappa(L_{\rm geom})\) in
\eqref{eq:four-block-fourier-average}, then \eqref{eq:display-nb1}, gives
\[
 v_\kappa\langle\Xi,\phi_\eta\rangle
 \le v_\kappa(t_0L_{\rm twist}L_{\rm geom}\rho_\eta)
       +P_\ell-2v_\ell(h_1).
\]
This bound also applies for \(\ell\mid h_1\).
For \(\ell\in\{2,p_0,p_1\}\), \(P_\ell=0\) because \(\ell\) splits
in \(K\), so every factor \(1-\chi_K(\ell^e)\) in \(S_\ell\) vanishes.

\begin{proposition}[The four-block numerator test]\label{prop:four-block-numerator}
For each coefficient embedding \(\kappa\) above \(\ell\) and
each allowed twist \(\eta\) of exact support \(Q\) with
\(\ell\nmid q_Q\), there is a test
\(\Xi_{\kappa,\eta}\) integral at every coefficient place above \(\ell\)
and orthogonal, for every split off-level \(w\), to
\(\operatorname{im}e_{X-1}\cap\ker\mathcal P_{n,w}(X_*)\cap\ker D_O\),
where \(D_O\) is descent to the full line form space. It satisfies
\[
 |\langle\Xi_{\kappa,\eta},\phi_\eta\rangle|_\kappa
     \ge |c_0\rho_\eta|_\kappa\ell^{-P_\ell},
 \qquad c_0=t_0L_{\rm twist}L_{\rm geom}\ne0 .
\]
The constant \(c_0\) is fixed for the original main pair and
\(\log|\sigma c_0|\le J_*\) at every complex automorphism \(\sigma\).
\end{proposition}
At the untwisted levels let \(h_{\rm s},h_{\rm t}\) be the numbers
of classes of the source line and the target of rank \(a_i\le 2r\).
Both levels contain the principal subgroup modulo \(M\).
Apply the lattice count in the proof of
Lemma~\ref{lem:terminal-class-count} in these ranks. Its integral
class bound is \(O(u^3\log x)\) on the logarithmic scale for
\(u=2fn\), \(n\in\{1,a_i\}\); refinement modulo \(M\)
adds at most \(2fn^2\log M\). Thus
\(\log h_{\rm s}+\log h_{\rm t}=O(x^{52})\) uniformly in the block.
For each of \(\sigma\) and \(c\sigma c\), apply the translated-vacuum
bound in the proof of Lemma~\ref{lem:upper-gaussian-kernel}
in \(a_i\) coordinates with polynomial one. The transported
\(V_0^{a_i}\)-scaled vacuum has scalar of modulus at most
\(\exp(C_{\rm int}q_1^2a_i)\). The ordinary Gaussian estimate
is \(\exp((2d+6)(q_1a_i+1)^2)\).
The untwisted column test has at most \(M^{fa_i}\) finite Fourier
terms. Including the degree-zero integral gives
\[
 U_i=M^{fa_i}\exp\!\big(C_P(q_1a_i+1)^2+fa_i\log4\big)
\]
for the absolute value of either conjugated kernel. In the Gram sum
use \(\sigma\) on the first kernel and \(c\sigma c\) on the kernel
whose value is conjugated before applying \(\sigma\). The stabilizer
weights are at most one and the weight-zero Weyl line has norm one.
Consequently, for every old character in this sum,
\[
 |\sigma A_i(\chi)|\le C_Gd_i^2 h_{\rm t}h_{\rm s}^2 U_i^2.
\]
Here \(\log C_G\le f(rt)^2\log M\),
\(\log d_i=2fa_i\log M\), and
\(\log U_i=O((q_1a_i+1)^2+fa_i(1+\log M))\).
The logarithm of the displayed bound is at most \(x^{53}\) eventually.
It follows that
\[
 \log|\sigma L_{\rm geom}|\le 2x^{53},\qquad
 \log|\sigma c_0|\le fmb/k+3x^{53}\le J_*.
\]
Indeed \(x^{53}=o(fmb/k)\) and
\(J_*/(fmb/k)=xB^3/m^2\ge B^3/4\). The cutoff depends only on
the fixed \(q_0,d,\epsilon\).

For a fixed twist \(\eta\) of exact support \(Q\) with \(\ell\nmid q_Q\),
the construction applies to the same \(\eta\)
at \(\kappa\) and \(\bar\kappa\); only the local exponent,
cell, Fourier character and test may change. The
fixed multiplier and the main line stay the same.

\smallskip
\subsubsection{Untwisted norms at both conjugate embeddings}\label{sec:four-block-reference-norm}

\sectionopening{In this subsection, we control the untwisted norms in the finite comparison of Proposition~\ref{prop:four-block-numerator}.}

Let \(m<z\le r\) be the half-rank of the fixed main pair.
Take \(a=2z\) and its original level at \(Q=\varnothing\).
In \eqref{eq:display-nb5}--\eqref{eq:display-nb6}
take \(\sigma\in\operatorname{Aut}(\mathbb C/K)\). Both \(\sigma\)
and \(\sigma^\#=c\sigma c\) then preserve the chosen CM type, the
embeddings extending the fixed embedding of \(K\).
For its untwisted datum \(X\) write \(J_X=\sum_\tau\lambda_{X,\tau}\) and
\(A_X=s_{\rm jet}^{J_X}X^{\rm raw}\), with \(X^{\rm raw}\)
the source class sum. The CM transport of
Section~\ref{sec:rank-one-cm-transport}, evaluated at the same
target point \(h\), gives
\[
 \sigma(A_X(h)_i)=
 t_\sigma^{J_X}c_\sigma b_Y(u_\sigma^{-1})
               d_\sigma(h')\,Y^{\rm raw}(h)_{i^\sigma}.
 \tag{NB5}\label{eq:display-nb5}
\]
Here \(b_Y\) is the rank-twisted character of the conjugate datum.
For \(\sigma\) fixing \(K\), let \(p_\sigma\) be the permutation
defined by \(\sigma\Phi_\tau=\Phi_{p_\sigma(\tau)}\) on the CM
embeddings used for the Fock coordinates.
Then \(\lambda_{Y,p_\sigma(\tau)}=\lambda_{X,\tau}\)
and \(i^\sigma_{p_\sigma(\tau)}=i_\tau\),
and \(d_\sigma(h')=\sigma(\Delta(h'))/\Delta(h')\) for the target Hom action
\(h'\). The scale satisfies \(|t_\sigma|^2=\Omega^{-1}\).
The vacuum scalar \(c_\sigma\ne0\) depends only on \(\sigma\),
the rank and the original level.

\(\sigma^\#=c\sigma c\) gives the same
datum \(Y\). Its permutation agrees with that of \(\sigma\)
since \(E/\mathbb Q\) is abelian. To check the finite
character, write any finite norm-one idele as \(y/\bar y\)
by local Hilbert 90, taking unit preimages at almost all
places. Choose a finite ray modulus so that its ray units
map into the kernel of the finite character. A power
\(y^h\) is a principal idele \(a\) times a ray unit.
Automorphy expresses the corresponding \(h\)-th power
of the finite character value as a product of integer
powers of the embeddings of \(a/\bar a\).
All conjugates of these numbers have absolute value
one. Thus for each finite character value \(v\),
\(\sigma^\#(v)=\overline{\sigma(v^{-1})}=\sigma(v)\).
Permuting its infinity type and conjugating its finite
part therefore gives one automorphic character, on the
same level, for both complex automorphisms.

In each Gram term, \(D_{h'}=\Delta(h')\overline{\Delta(h')}\)
is a positive integral lattice index. Hence
\[
 d_\sigma(h')\overline{d_{\sigma^\#}(h')}=\frac{\sigma(D_{h'})}{D_{h'}}=1.
\]
The class weights and Gram matrix are
rational and the place permutation reindexes the basis.
Using \eqref{eq:display-nb5} in the two legs gives the equality
\[
 |\sigma\langle A_X,A_X\rangle|
       =|c_\sigma c_{\sigma^\#}|\,
              \Omega^{-J_X}\|Y^{\rm raw}\|^2 .
 \tag{NB6}\label{eq:display-nb6}
\]
The two character values in \eqref{eq:display-nb5} have modulus one.

Apply the proof of Proposition~\ref{prop:marked-line-norm} at the
untwisted level to the raw kernel. For the transported datum \(Y\), let
\(\lambda_{Y,\tau}=c_{Y,\tau}-z\) denote its local degree and put
\(W_Y=\bigotimes_\tau\operatorname{Sym}^{\lambda_{Y,\tau}}(\mathbb C^a)\).
Use a highest covector of norm one for the compact invariant Hermitian
form. The local monomial has Bargmann
square norm
\[
 \pi^{-a}\int_{\mathbb C^a}|Z_1|^{2\lambda}
                      e^{-\|Z\|^2}\,dZ=\lambda!.
\]
For a unit highest covector, compact Schur orthogonality gives the
scalar coefficient's square integral as the squared Weyl norm divided
by \(\dim W_Y\). Write \(v_j=\operatorname{vol}([U(j)])\) in the
original measures of Section~\ref{sec:marked-norm-line} and
\(S^0=S_M\cup\bar S_M\). Then
\[
 \begin{aligned}
 C_{a,S_M}
  &=\frac{v_a v_1}{\zeta_{F_0}(2z)}\prod_{w\in S_M}(q_w-1),\\
 \|Y^{\rm raw}\|^2
  &= C_{a,S_M}\dim W_Y\left(\prod_\tau\lambda_{Y,\tau}!\right)
                 L^{S^0}(z,(\nu_Y^{\rm bc})^{-1}).
 \end{aligned}
\]
Indeed, for \(q=\mathrm N v\) at a split pair \(w,\bar w\) off
\(S^0\), put \(A_w=\nu_Y^{\rm bc}(\varpi_w)\). The two Euler factors of the base change multiply to
\(|1-A_wq^{-z}|^{-2}\), with \(|A_w|=1\). At an inert place
the base change is trivial on the uniformizer of \(F_{0,v}\), and
its factor \((1-q^{-2z})^{-1}\) cancels the local zeta denominator.
The finite extension \(E/F_0\) is unramified at every finite place by
the inertia calculation in Section~\ref{sec:different-additive-character-scales}. At an old
split pair, even for a ramified tame character, the character is
trivial on the principal-unit stabilizer. The local mass, the
target index, the square of the source index and the removed zeta
factor contribute together
\[
 \frac{q^{-a}}{q-1}(q^a-1)(q-1)^2(1-q^{-a})^{-1}=q-1.
\]
These factors are positive. For \(|A|=1\), \(0<t<1\),
the power series for \(\log(1-At)\) gives
\(\bigl|\log|1-At|\bigr|\le-\log(1-t)\). Summing over places yields
\[
 \left|\log L^{S^0}(z,(\nu_Y^{\rm bc})^{-1})\right|
   \le\log\zeta_E(z)\le [E:\mathbb Q]\log\zeta(2)<4f .
\]
The vacuum and jet factors satisfy
\[
 \|\phi_X\|^2=|V_0^a|^2 M^{J_X}\langle A_X,A_X\rangle,
 \qquad |\pi_M|^2=M\in\mathbb Q_{>0}.
\]
In particular the common term at \(\sigma\) is
\[
 C_{\sigma,a,S_M}
 =\log\!\left(\left|\sigma\!\left(|V_0^a|^2\right)\right|\,
            |c_\sigma c_{\sigma^\#}|\, C_{a,S_M}\right).
\]
Equation~\eqref{eq:display-nb6} and the permutation of the local
degrees give
\[
 \log|\sigma\|\phi_X\|^2|
  =C_{\sigma,a,S_M}+r_X+
     \sum_\tau \bigl(\log(\lambda_{X,\tau}!)
                     +\lambda_{X,\tau}\log(M/\Omega)\bigr),
\]
where \(r_X=\log\dim W_Y+\log L^{S^0}(z,(\nu_Y^{\rm bc})^{-1})\).
Thus \(r_X\ge-4f\), and \(r_X\le4f\) for the empty partition,
when \(\dim W_Y=1\).

For the empty reference of this same rank, take
\(q_{\rm ref}=(\nu_{\rm ref},z)\) with
\(\nu_{\rm ref}=\upsilon_1^{-z}\), where \(\upsilon_1=\upsilon|_{U(1)}\).
Its old unit characters are trivial, its infinity weights are
\(-z\), and its rank-twisted character is \(\nu_{\rm ref}\upsilon_1^z=1\).
The source retains its old principal level and the target its
vector fixer \(J_\varnothing\). Since \(J_{\rm ref}=0\) and the
coefficient space has dimension one, the measured identity gives
\[
 \|\phi_{\rm ref}\|^2=|V_0^a|^2 C_{a,S_M}
    L^{S^0}(z,(\nu_{\rm ref}^{\rm bc})^{-1})>0 .
\]
The positive Euler product for the main datum also gives
\(\|\phi_X\|^2>0\). By Section~\ref{sec:marked-both-norm-jets}
at \(Q=\varnothing\), the values of
\(d_{\rm ref}\phi_{\rm ref}\), with \(d_{\rm ref}=M^{4fz}\), are
integral at every coefficient embedding above each \(\ell\).
In the class pairing, use the two embeddings \(\kappa\) and
\(\bar\kappa\) for the conjugate factors. The Weyl Gram matrix
is integral and every stabilizer order divides \(C_G\), by
Section~\ref{sec:marked-norm-line}. Hence
\[
 I_{\rm ref}=C_Gd_{\rm ref}^2\|\phi_{\rm ref}\|^2\ne0,
 \qquad v_\kappa(I_{\rm ref})\ge0
\]
for every \(\kappa\). The fixed multiplier
\(H=I_{\rm ref}/\|\phi_X\|^2\) therefore clears the
untwisted norm at every prime. Subtract the two
conjugate logarithms. The common constant cancels
and the main part sizes are at least \(\kappa_\tau\).
Monotonicity of \(J(n)=\log(n!)+n\log(M/\Omega)\)
gives
\[
 \log|\sigma H|\le \log(C_Gd_{\rm ref}^2)+8f
                     -\sum_\tau J(\kappa_\tau)
             \le J_*-\sum_\tau J(\kappa_\tau).
\]
Indeed \(\log(C_Gd_{\rm ref}^2)\le
f((rt)^2+8z)\log M\). Since \(z\le r\), \(\log M\le22\log b\)
and \(f,r,t,k,m\) grow at most polynomially in \(x\), this
upper bound and \(8f\) are \(o(J_*)\) for fixed \(q_0,d,\epsilon\),
where \(J_*=fxb(\log b)^3/(mk)\).
If \(\sigma\) acts nontrivially on \(K\), then
\(c\sigma\in\operatorname{Aut}(\mathbb C/K)\) and
\(|\sigma H|=|(c\sigma)H|\). Applying the bound to \(c\sigma\)
on the same old level proves it for every complex automorphism.
This is the fixed multiplier estimate
for the untwisted term of the reciprocal height.

\section{The common twist and bounds on critical values}
\label{sec:common-critical-twist}

\sectionopening{The two measured norm tests of
Section~\ref{sec:norm-tests-and-trace-algebra} give the signed systems
for Proposition~\ref{prop:joint-norm-selmer-twist}. That proposition
selects one twist per coefficient place and bounds the Selmer lengths
and signed critical ratios for its prescribed powers. We then apply
Proposition~\ref{prop:clearing-on-original-tests} to those twisted
summands to bound the projected pairings.}

Fix \(q_0,d\) and the exponent of Section~\ref{sec:datum}, and take \(x\)
sufficiently large for these parameters. Fix the resulting admissible
datum. For a coefficient place \(\kappa=(\ell,\iota)\)
put \(\bar\kappa=(\ell,\iota\circ\mathrm{conj})\) and normalize
\(v_\kappa(\ell)=1\). Write \(S_M\) for the old chosen places and
\(A_{\rm aut}=\prod_{w\in S_M}|\mathrm{GL}_{rt}(k_w)|\). For a
ramified critical pair \(p=(\nu,z)\), sign \(\epsilon=\pm1\), and
twist \(\eta\) of exact support \(Q\), put \(q_Q=\prod_{w\in Q}\mathrm Nw\),
taking one place of each conjugate pair. Define
\[
 R_{p,\epsilon}(\eta)
 =\frac{|L^{S_M\cup Q}(z,((\nu\eta^\epsilon)^{\rm bc})^{-1})|}
 {|L^{S_M}(z,(\nu^{\rm bc})^{-1})|}.
\]
The superscript omits the conjugates of the indicated places as well.
Let \(\beta\) be either inverse base-change character in the ratio.
Omit a conjugation-stable finite set \(S'\) containing the ramified
places. For real \(s>1\), the Euler product converges absolutely
and has the factors at a split pair and an inert place, respectively,
\[
 |1-\beta(\varpi_w)q^{-s}|^{-2},
 \qquad (1-q^{-2s})^{-1},
\]
where \(q\) is the residue norm below the pair or inert place.
The first identity uses \(\beta^c=\beta^{-1}\) and unitarity.
At the inert place, base change is trivial on the uniformizer
from \(F_0\). This applies to the two partially omitted products
appearing in the ratio as well.

At \(z=1\) the critical pair has \(c_v\ge1\) at each infinite
place by Section~\ref{sec:notation-comparison}. In its chosen
complex coordinate the inverse base-change character restricts
to the unit circle by \(u\mapsto u^{2c_v}\). The finite-order twist
leaves this restriction unchanged. Both \(\beta\) and \(\beta^2\)
are therefore nontrivial on the ideles of module one.
Hecke continuation gives holomorphy of their \(L\)-functions at one
\cite[Theorem~4.4.1 and \S4.5]{TateThesis}. For real \(s>1\),
\[
 \zeta_E^{S'}(s)^3|L^{S'}(s,\beta)|^4|L^{S'}(s,\beta^2)|\ge1
\]
by the logarithmic Euler-product argument of
Section~\ref{sec:rank-two-doubling-local-integrals}. A zero at one would
make the left side tend to zero, since \(\zeta_E\) has a simple pole.
Removing finitely many Euler factors preserves nonvanishing at every
\(z\ge1\). The common twist chosen below controls the valuations of
these ratios.

\smallskip
\subsection{Signed systems for the critical pairs}
\label{sec:critical-signed-systems}

\sectionopening{In this subsection, we form the systems for the common torsion choice of Proposition~\ref{prop:joint-norm-selmer-twist}.}

Write each of the two coefficient formulas as
\(F_{p,Q}(\eta)=\sum_i c_i(Q)\eta(u_i)\).
Fix any \(\ell\nmid q_Q\). The coefficient formulas give
integrality at every coefficient embedding. If some \(c_i(Q)\)
were transcendental, choose a transcendental
\(u\in\mathbb Z_\ell^\times\). A field isomorphism
\(\mathbb C\simeq\overline{\mathbb Q}_\ell\) extending the chosen
embedding of the original algebraic data can send \(c_i(Q)\)
to \(u/\ell\), by extension of transcendence bases as in
Section~\ref{sec:gram-product-formula}. The negative valuation
contradicts integrality. For an algebraic coefficient the same
integrality at every place away from \(q_Q\) then gives
\(c_i(Q)\in\overline{\mathbb Z}[1/q_Q]\), where
\(\overline{\mathbb Z}\) is the ring of algebraic integers.
Use the signed array
\[
 \Theta_{p,\epsilon,Q}
    =\sum_i\iota(c_i(Q))[u_i^{-\epsilon}]_Q.
\]
For the inverse avatar \(\psi=\eta^{-1}_\ell\), evaluate a norm-one
class \([u]_Q\) on an Artin preimage of \(u=y/\bar y\).
In the construction of the auxiliary characters, every
\(\chi_w(\operatorname{Art}y)\) is \(\theta_w(y/\bar y)\)
for an idele class character \(\theta_w\). Changing \(y\) by a
fixed-field idele therefore leaves the value unchanged, also when
\(\ell=2\). Evaluation at the resulting inverse avatar gives
\(\iota(\eta(u)^{-1})\).
Therefore \(\psi(\Theta_{p,\epsilon,Q})=\iota(F_{p,Q}(\eta^\epsilon))\).
The sign is in the array, so each sign is evaluated with the same
\(\psi\).

The unit-valued character in the Euler norm relation is
the avatar \(\vartheta_\epsilon\) of \(B_\epsilon=(\nu^{-\epsilon})^{\rm bc}\).
For a finite idele \(x\) of \(E\) its formula in geometric reciprocity is
\[
 \vartheta_\epsilon(\operatorname{Art}_E x)
  =\iota(B_{\epsilon,f}(x))
      \prod_{\tau\in\Phi}
       \left(\frac{(\iota\tau)(x_\ell)}
                   {(\iota\bar\tau)(x_\ell)}\right)^{\epsilon c_\tau},
\]
where each embedding reads its induced component above \(\ell\).
On principal ideles the product is the transported infinity component,
so it cancels the finite character. The finite values of \(\nu\) lie
in one number field by the norm-one class calculation of
Section~\ref{sec:marked-level-characters}. The displayed map is
continuous and factors through the compact quotient of the finite
ideles by the closure of \(E^\times\), by global reciprocity.
Its valuation has compact image in \(\mathbf R\), hence is zero.
The algebraic correction in the displayed formula uses only
components at \(\ell\). Away from \(\ell\) and the old level the finite
character is trivial on local units, so the avatar is unramified
there as required in Section~\ref{sec:joint-twist}. At a split pair outside \(S_M\cup\bar S_M\cup Q\cup\bar Q\)
and of residue characteristic different from \(\ell\), put
\(A_w=\iota((\nu\eta^\epsilon)^{\rm bc}(\varpi_w))\).
In the geometric Artin convention we have
\[
 (\vartheta_\epsilon\psi)(\operatorname{Fr}_{w,\rm geom})
       =A_w^{-\epsilon},\qquad
 (\vartheta_\epsilon\psi)(\operatorname{Fr}_{w,\rm ar})
       =A_w^\epsilon .
\]
Indeed \(\nu^{-\epsilon}\eta^{-1}=(\nu\eta^\epsilon)^{-\epsilon}\),
and the two Frobenius elements are inverses. Each
Euler-factor pair is invariant under \(A_w\mapsto A_w^{-1}\). On an
unramified enlargement of the nominal modulus the rank-one
system deletes the pair at exponent \(z\), and the rank-two
system deletes the pairs at exponents \(0,1\). The rank-one
normalizing power of \(q_Q\) cancels the local mass at
each new place even when the residue character is trivial.
In rank two the local multiplier at \(0,1\) is
\eqref{eq:rank-two-euler-deletion}, with the same
fixed jet scalar for the two levels. These are multiplication
identities even if a deleted factor vanishes.

For clarity, let \(\mathcal C_Q\) be the finite norm-one
class group with maximal integral units off \(Q\) and
principal units at \(Q\). The old level is part of the
fixed norm test and of \(\nu\).
Write \(P:\mathcal C_{Q\cup\{w\}}\to\mathcal C_Q\) for
the quotient and let \(b_w\) be the class of \(y_w/\bar y_w\)
where \(y_w\) is the idele that is a uniformizer at \(w\)
and one elsewhere. Put \(a_w=\iota(\nu^{\rm bc}(y_w))\) and
\(\mathcal E_p=\{z\}\) for \(z\ge2\), \(\mathcal E_p=\{0,1\}\)
for \(z=1\). Then
\[
 P_*\Theta_{p,\epsilon,Q\cup\{w\}}=
 \Theta_{p,\epsilon,Q}
 \prod_{e\in\mathcal E_p}
    (1-q_w^{-e}a_w^{-\epsilon}[b_w])
    (1-q_w^{-e}a_w^\epsilon[b_w]^{-1}).
\]
For every \(j\ge1\), evaluation on \(\psi^j\) sends
\(a_w^{-\epsilon}[b_w]\) to
\[
 a_w^{-\epsilon}\iota(\eta(b_w))^{-j}
 =\iota\bigl((\nu\eta^{j\epsilon})^{\rm bc}(y_w)\bigr)^{-\epsilon}
 =(\vartheta_\epsilon\psi^j)(\operatorname{Fr}_{w,\rm geom}).
\]
The reciprocal group unit gives the inverse value. Thus only the
auxiliary character is raised to \(j\). For \(z=1\), if this
Frobenius value is one, the factor at exponent zero makes
\(F_{p,Q\cup\{w\}}(\eta^{j\epsilon})=0\) at this enlarged modulus.
The identity still holds for that character.
The levels for the two class groups differ only at the split
place \(w\). Hence \(\ker P\) is the image of the integral norm-one
units there. If a power of a character on
\(\mathcal C_{Q\cup\{w\}}\) is trivial on these units, that power
descends to \(\mathcal C_Q\), where the identity applies to it.
The characters of \(\mathcal C_Q\) have trivial residue part at
\(w\) and separate its group algebra in characteristic zero. Indeed, take a finite extension \(K_\ell\)
containing the displayed coefficients and all character values of
\(\mathcal C_Q\). If \(D=\sum_c d_c[c]\) is the difference of the two sides,
Fourier inversion in \(K_\ell\) gives
\[
 |\mathcal C_Q|d_c=\sum_{\chi\in\widehat{\mathcal C_Q}}
                         \chi(c)^{-1}\chi(D)=0.
\]
Thus each \(d_c=0\) in \(K_\ell\). The identity holds in its valuation
ring coefficient by coefficient, including when \(\ell\mid|\mathcal C_Q|\),
since \(q_w,a_w\) are units and the arrays are integral.
It supplies the norm relation of Section~\ref{sec:joint-twist}.

For \(Q=J\) and the auxiliary characters \(\chi_w\) of that section,
choose finite idele representatives \(y_i\) with
\(u_i=y_i/\bar y_i\) in \(\mathcal C_Q\). With \(N\) the auxiliary exponent,
put \(n_{iw}=[-\epsilon\chi_w(\operatorname{Art}y_i)]_{\ell^N}\),
the representative between \(0\) and \(\ell^N-1\). The integral polynomial
\[
 P_{p,\epsilon,Q}(X)=\sum_i\iota(c_i(Q))
                              \prod_{w\in Q}(1+X_w)^{n_{iw}}
\]
has value \(\iota(F_{p,Q}(\eta^{j\epsilon}))\) at
\(X_w=\zeta_w^j-1\) whenever \(\zeta_w^{\ell^N}=1\) for each \(w\) and
\(\prod_{w\in Q}\zeta_w^{\chi_w}\) is the inverse avatar of \(\eta\).
The signs therefore use one torsion tuple.

For exact support the two norm formulas give
\[
 \frac{F_{p,Q}(\eta^\epsilon)}{F_{p,\varnothing}(1)}
 =\begin{cases}
     R_{p,\epsilon}(\eta),&z\ge2,\\
     q_QR_{p,\epsilon}(\eta)^2,&z=1.
   \end{cases}
\]
At \(Q=\varnothing\), the line test prescribes the old column modulo
\(\pi_M\) and uses its vector fixer \(J_\varnothing\). The rank-two
period uses the determinant character modulo \((\pi_M)^2\) and the full
standard compact \(K_2\) in \eqref{eq:determinant-marked-period}.
The augmentations for both signs are the images under \(\iota\) of
\[
 F_{p,\varnothing}(1)=
 \begin{cases}
     C_Gd_p^2\|\Phi_1\|_{J_\varnothing}^2,&z\ge2,\\
     |C_Gs_p\Lambda_\varnothing(1)|^2,&z=1 .
 \end{cases}
\]
Proposition~\ref{prop:marked-line-norm} gives positivity in the first
case. In the second, \eqref{eq:determinant-marked-ratio-product}
and the functional equation there give \(\Lambda_\varnothing(1)\ne0\).
With \(q_\varnothing=1\), the integral coefficient formulas for
both tests apply at every coefficient place. Thus
\(a_p=F_{p,\varnothing}(1)\) is a positive algebraic integer.
It is fixed before the auxiliary places and its valuation at a
coefficient place may be positive.
For a twist of exact support, put
\(D=F_{p,Q}(\eta^\epsilon)/F_{p,\varnothing}(1)\).
Both evaluations are algebraic by their finite character sums.
If \(z\ge2\), the critical ratio is \(D\). If \(z=1\), it is the
positive real root of \(X^2-D/q_Q\), hence is algebraic as well.
At every coefficient embedding with \(\ell\nmid q_Q\) this gives
\begin{equation}\label{eq:critical-norm-ratio-valuation}
 v_\kappa(R_{p,\epsilon}(\eta))=
 \begin{cases}
    v_\kappa(D),&z\ge2,\\
    \tfrac12 v_\kappa(D),&z=1 .
 \end{cases}
\end{equation}

\smallskip
\subsection{Original Hecke systems and their integral models}
\label{sec:critical-hecke-systems}

\sectionopening{In this subsection, we define the lift systems corresponding to the arrays of Section~\ref{sec:critical-signed-systems}.}

Fix the coefficient embedding above \(\ell\). Let \(\mathcal A\) be the
set of pairs \((j,b)\), where \(1\le j\le m\) and \(b\) is an
outside-level Hecke system in the projected space on \(W_m\).
Let \(\mathcal B\) be the set of \((j,p)\), where \(p=(\nu,z)\)
is main for \(j\), so that
\[
 m<z\le r,\qquad \nu_{\rm tame}=\chi_M^j,\qquad
 \mathrm{wt}_v(\nu)=-c_v,\qquad w_j(v)+r-t\le c_v\le w_j(v).
\]
The spaces of principal-level forms are finite dimensional.
For a fixed infinity type and fixed finite level the characters in
\(\mathcal B\) differ by characters of a finite norm-one class
quotient. The bounds on \(c_v,z\) and \(j\) therefore make
\(\mathcal B\) finite.

For each lower system \(b\) in \(\mathcal A\), choose once a nonzero
form of system \(b\) in the image of a central idempotent for a
selected finite type on \(W_m\). This is
possible because the sum of the idempotents for the selected types fixes the projected space, and each idempotent commutes with the off-level Hecke
operators. The residue action on its algebraic form is the one
obtained by reducing at the prescribed level, since the principal
kernel fixes the form and the algebraic conversion commutes with
right translation.

Choose the attachment \(R_b\) of a constituent of this form on
\(W_m\). Its local type parameters are those of the chosen constituent.
For \((j,b)\in\mathcal A\) set
\[
 T_{j,b}=\theta_j R_b,\qquad
 \theta_j=\operatorname{av}_\ell
       \bigl(u\longmapsto\mu_j^{-1}(u/\bar u)|u|^{(r-t)/2}\bigr),
\]
where \(\mu_j\) is the original descent character of lift \(j\).
The local--global comparison gives the scaled lower parameters of
\(T_{j,b}\). On every nonzero lift, Section~\ref{sec:main-trace}
identifies the normalized upper trace with these parameters and the
string \(\chi_{\rm cyc}^{-a}\), \(0\le a<t\).
The normalization changes the raw polarization exponent \(1-r\)
to \(1-t\), since \(\theta_j\theta_j^c=\chi_{\rm cyc}^{r-t}\).
Its jumps are \eqref{eq:separated-hodge-tate-intervals}.

For \((j,p)\in\mathcal B\), use the constant system on the
original norm-one line. The raw attachment is the trivial
character, with jump \(0\) and parameter \(X-1\) at every distinguished
place. Put
\[
 t_p=2z-1,\qquad
 T_{j,p}=\operatorname{av}_\ell
      \bigl(u\longmapsto\nu^{-1}(u/\bar u)|u|^{\,1-z}\bigr).
\]
This is the same normalization with lower rank \(1\) and upper
rank \(2z\). Its polarization exponent is \(1-t_p\) and its
single jump at \(\sigma\) is \(z-1-\varepsilon_\sigma c_v\).
In particular
\[
 |\varepsilon_\sigma(z-1-\varepsilon_\sigma c_v)+w_j(v)|
     \le t .
\]
Indeed \(0\le w_j(v)-c_v\le t-r\) and \(z-1<r\).
As \(w_j(v)>2(r^8)^2+t\), there are no small jumps, and the two
conjugate jumps add to \(2z-2\). For every allowed twist the
lift of the constant line form is nonzero by
Proposition~\ref{prop:marked-line-norm}: its squared norm is a positive
scalar times an absolutely convergent, nonvanishing Euler product
at \(z>m\ge1\). Section~\ref{sec:main-trace} therefore gives the
string of length \(t_p\) for its upper system.

The unscaled parameter \(X-1\) is the trivial character of each
\(\mathrm{GL}_1(k_w)\). It also fixes the upper lift in the
normalized residue action. For a tuple \(g\) of residue scalars
choose the indicated lower integral lift and the upper integral
lift of its first block. Their simultaneous Weil action on the
column multiplies it by
\(\prod_{w\in S_M}\nu_w(\det g_w)^{-1}\).
In the split matrix coordinates, if \(i_0\) is the
original block inclusion, the two lifts act on the support by
\(X\mapsto k'Xk^{-1}\), up to the fixed transpose convention.

Its reduction fixes \(i_0\), since \(\operatorname{diag}(g_w,1)i_0
=i_0g_w\). The integral matrix lattice is preserved and both
determinants are units. The two splitting-character factors have
opposite rank powers of the splitting character and leave
\(\nu_w(\det g_w)^{-1}\) at \(w\).
The same identity holds after the auxiliary twist, whose
tame part is trivial on \(S_M\). Changing the source variable
in the theta integral shows that the upper translation multiplied
by \(\prod_{w\in S_M}\nu_w(\det g_w)\) acts on the lift as
the lower translation acts on the source. The constant vector is
fixed. Hence the mean over the \(\mathrm{GL}_1(k_w)\) fixes
the line, which is the required \(X-1\) idempotent. Its
normalized inertial root instead names \((\chi_{M,w}^j)^{-1}\)
in the original parameter convention. It is nontrivial since
\(1\le j\le m<q_2\).

At rank \(r\), the same residue equivariance on the
vector column gives the corresponding assertion for every
component in a selected finite type of a lower eigenform on \(W_j\).
Decompose the projected space on \(W_j\) into simultaneous
eigenspaces for the off-level split spherical Hecke operators
and the finite central type projectors. These operators commute
and are normal in the original Petersson inner product.
Fix the old residue tuple and let \(b\) occur on \(W_j\).
For each other spherical character \(b'\) in the union of the two
finite spectra choose a test \(h_{b'}\) away from the old level
with \(b(h_{b'})\ne b'(h_{b'})\). In the common spherical algebra put
\[
 h_b=\mathop{\prod}_{b'\ne b}
       \frac{h_{b'}-b'(h_{b'})\,1}{b(h_{b'})-b'(h_{b'})},
\]
with multiplication by convolution and empty product \(1\).
The spherical operators commute with their adjoints, so the spectrum
is semisimple and \(h_b\) projects onto the entire \(b\)-eigenspace.
The trace identity of Section~\ref{sec:marked-target-dimension}
at the fixed old tuple gives the same positive rank on \(W_j\)
and \(W_m\), including all multiplicities.

The hypotheses of that trace argument
hold for both source spaces, since the same two finite places with primitive types
and the prescribed common regular real parameter occur there.
For every nonzero \(y\) in such a piece and every admissible
\(\xi\), Proposition~\ref{prop:first-lift-gram} gives
\[
 \|u_{j,\xi^j}y\|_{\rm class}^2
   \ge \tfrac12 A_j c_Q(\xi)\|y\|^2>0.
\]
This uses the first-lift estimate on the original single-lift space.
Residue equivariance gives the specified upper type.
For each lower constituent on \(W_j\) with system \(b\),
choose a constituent on the definite group \(W_m\) realizing its
tuple and system by this trace comparison. Its rank-\(r\) base
change is cuspidal. Indeed at an old split place with primitive type
the local component is supercuspidal, and the split base change of
\cite[Corollaire~5.3]{LabesseCM} has this component in the chosen
coordinate. Two or more global discrete blocks would give a
subquotient of induction from a proper local parabolic. A single
length-\(d'\) Speh block with \(d'>1\) would likewise give a
quotient of induction from \(d'\) blocks of rank \(r/d'<r\)
by \cite[theorem of the introduction, pp.~605--606]{MWResidual}. Supercuspidality excludes both
possibilities. The general discrete-block construction of
Section~\ref{sec:cofinite-spherical-galois} also applies at the
larger target ranks.

To compare two such definite constituents, choose one prime \(\ell'\nmid M\) and one
coefficient isomorphism there. Their raw semisimple attachments
have equal traces at cofinite split Frobenius classes.

Prime ideals of \(E\) of absolute degree at least two have
Dirichlet density zero, since for \(s>1\) their sum of inverse
norms to the power \(s\) is at most
\([E:\mathbf Q]\sum_p p^{-2s}\). At a degree-one prime away from
ramification the underlying rational prime splits completely in
the Galois field \(E\). Thus the indicated traces agree on a dense
set of Frobenius classes by Chebotarev. Brauer--Nesbitt identifies
the two semisimple representations.
At every \(w\in S_M\) their Frobenius-semisimplified
Weil--Deligne parameters agree by Section~\ref{sec:inertial-labels}.
Both use the twist by \(|\det|^{(1-r)/2}\) of the local Langlands
parameter and geometric reciprocity. Removing that common twist
identifies their smooth local components.
For the common tame generator of
Section~\ref{sec:marked-inertia-characters}, the eigenvalues on
the raw inertia representation are the images under the chosen
coefficient isomorphism of
\(\{\jmath(\xi):c_w(\xi)=0\}\). Since \(\jmath\) is injective and
\(c_w\) is monic and squarefree, this multiset recovers \(c_w\).
Consequently the same system \(b\) gives one original tuple, also
when the working coefficient prime \(\ell\) divides \(M\), since
the comparison of these characteristic-zero types used
\(\ell'\nmid M\).
Return to the prescribed coefficient embedding \(\kappa\).
For each nonzero lifted component of system \(b\), original tuple \(c\)
and allowed auxiliary twist \(\eta\), put
\(\psi_\eta=\operatorname{av}_\ell(\eta^{-1})\).
At every split Frobenius outside the two finite exceptional sets,
the split polynomial in Section~\ref{sec:normalized-trace-multiplicities}
uses the eigenvalues \(b_{w,i}\) of the chosen constituent on \(W_m\).
Taking its linear coefficient and using density as above gives
\[
 \operatorname{tr}\widetilde R_{\rm lift}(g)
 =\psi_\eta(g)\operatorname{tr}T_{j,b}(g)
     +\sum_{a=0}^{t-1}\chi_{\rm cyc}(g)^{-a}
        \qquad(g\in G_E).
\]
Semisimplicity identifies this normalized attachment with
\(\psi_\eta\otimes T_{j,b}\oplus
  \bigoplus_{a=0}^{t-1}\overline{\mathbf Q}_\ell(-a)\)
over \(\overline{\mathbf Q}_\ell\). The finite twist has jump zero
and is unramified at the old places. Thus its main summand has the
jumps in \eqref{eq:separated-hodge-tate-intervals} and the same
old inertial parameters there, read on potentially crystalline periods
at the places over \(\ell\). For \(\eta=\xi_\kappa^j\) the scalar is
\(\psi_\eta=\psi^j\). Residue equivariance places the component in
its own upper type projector for \(c\). The lifts of these components
span the original projected lift image.
Therefore the finite list supplies a main summand for every
generator used in clearing on the original lift image.

\smallskip
\indent\textit{Integral models and the full matrix span.}
Take \(K_\ell/\mathbf Q_\ell\) to be the normal closure of the
finite compositum for the original matrices, prescribed lattices,
scalar values and local type eigenvalues.
Apply Lemma~\ref{lem:continuous-attachment-stable-lattice} and extend
the resulting stable lattices by scalars to this compositum.

Let \(L\) be any such invariant lattice over a ring \(O\) containing
the relevant inertia eigenvalues. At rank \(r\) choose a primitive
place above an old prime distinct from \(\ell\). The preceding
primitive-parameter calculation gives \(r\) mutually incongruent
eigenvalues for the specified tame generator \(s\).
Their Lagrange idempotents belong to the \(O\)-span of the
group image and split \(L\) into free lines of rank one.

An arithmetic Frobenius \(\phi\) at this place permutes the
lines in one cycle and has unit entries between the lines
in integral bases, since its inverse preserves the lattice.
Thus the products of the idempotents and powers of \(\phi\)
give all matrix units up to units. The span is an algebra
and is \(\operatorname{End}_O(L)\). The operator \(\phi^r\)
acts on every line by the same product of the \(r\) unit
multipliers, with cyclotomic value \(p_i^r\).

At rank one the full span follows from the identity and
the same Frobenius is itself scalar, with cyclotomic value \(p_i\).
After reduction over any extension of the residue field, the
matrix units still span. Consequently every lattice in either
rank has absolutely irreducible reduction. A scalar unit-valued
twist preserves the span, since each matrix is multiplied by
a unit. In particular, after passing to any common finite coefficient
DVR \(O'\), every stable lattice in this twisted representation is
homothetic to the selected one \(L'\). Scale the other lattice
\(L''\) so that \(L''\subseteq L'\) but
\(L''\not\subseteq\varpi' L'\). Its image in \(L'/\varpi' L'\)
is a nonzero invariant subspace, hence is the whole reduction.
Nakayama's lemma gives \(L''=L'\). This scalar identification
preserves the cocycle equations and the normalized local conditions.

The old type calculation bounds the inertial determinants at
every original place prime to \(\ell\), with the same chosen tame
generator at a distinguished place and its conjugate inverse. In each
summand one may choose the scalar inertia element supplied by
that calculation, with scalar \(\zeta\) satisfying
\(\zeta^{h\ell^a}-1\in O^\times\) for \(a\ge0\),
\(1\le h\le4\operatorname{rank}(L)+2\).
This element is chosen for the summand and need not be the
common reference tame generator \(g_0\).
The twists are unramified at the scalar inertia place, so on
the twisted lattice the same element acts by \(\zeta\).
The unit \(\zeta-1\) gives the normalization
\(z(\gamma)=0\). Adding the coboundary of
\(-(\zeta-1)^{-1}z(\gamma)\) imposes this condition uniquely.

Use one common \(g_0\) and torsion cutoff for
\(\mathcal A\amalg\mathcal B\).
The joint choice applies every power of the selected character to
\emph{every} original summand in this union. The ramified critical
pairs and their signs give the separate norm tests. The character
is trivial on the full decomposition groups above \(\ell\), so
the selected local lattice is the scalar extension of the original
one there.

For a selected twisted lattice \(L_0\) over \(O_0\), the
normalized divisible-coefficient Selmer length after any finite
coefficient extension is bounded by its value over \(O_0\),
by \eqref{eq:specialized-selmer-coefficient-extension}. After choosing
the torsion tuple and the raised systems, this applies to the common
normal closure of Section~\ref{sec:marked-ribet-coefficient}.

\subsection{The common specialization and its valuations}
\label{sec:critical-common-specialization}

\sectionopening{In this subsection, we use Proposition~\ref{prop:joint-norm-selmer-twist} for the signed systems of Section~\ref{sec:critical-signed-systems}.}

Let \(\mathcal P_{\rm ram}\) be the set of ramified critical pairs.
The norm tests include each pair with both signs.
For the datum fixed at \(x\), their number \(2|\mathcal P_{\rm ram}|\) is at most
\[
 2t(4b)^f|\mathcal C_{\varnothing}|\prod_{w\in S_M}(q_w-1).
\]
Indeed \(1\le z\le t\), each integer \(c_v\) lies in \([z,4b]\),
and each old residue character is on \(k_w^\times\). Two characters
with the same infinity weights and residue characters differ by
a character of the finite norm-one class group
\(\mathcal C_{\varnothing}\) of Section~\ref{sec:critical-signed-systems}.
Thus the family is fixed before the auxiliary level is chosen.
In the common specialization there is one chosen Selmer minor for each
member of \(\mathcal A\amalg\mathcal B\) and one group-algebra
polynomial for each signed norm test. Counting these chosen polynomials
at each power \(1\le h\le m\) gives
\[
 m\bigl(|\mathcal A|+|\mathcal B|+2|\mathcal P_{\rm ram}|\bigr)
\]
evaluations with multiplicity.

Choose \(p_i\in\{p_0,p_1\}\) distinct from \(\ell\), possible
because the two old rational primes are distinct. For the common
reference choose one \(w_*\in S_M\) over \(p_i\) and its designated
tame generator \(g_0\). Write \(q=\mathrm Nw_*\), and use the
common multiple \(D\) of the degrees of the original type factors at
\(w_*\) in Section~\ref{sec:marked-inertia-characters}. The element
\(\omega_D(g_0)\) generates \(\mathbb F_{q^D}^{\times}\), so its norm
\(\omega_1(g_0)\) generates \(\mathbb F_q^{\times}\). The old-type
determinant bound applies to every original summand at both \(g_0\)
and \((g_0^c)^{-1}\) by that same section. A ramified critical
pair has nontrivial residue character at every place in \(S_M\).
Its inverse has the same property. For each signed
system its unit-valued character sends \(g_0\)
to a nontrivial root of unity \(\zeta\) of order \(h\mid q_{w_*}-1\).
The identity
\[
 h=(1-\zeta)\sum_{i=0}^{h-1}(h-1-i)\zeta^i
\]
gives \(v_\kappa(\zeta-1)\le v_\ell(q_{w_*}-1)\).
Thus the first-prime bound is the augmentation valuation
plus at most \(4v_\ell(q_{w_*}-1)\), as there are at most two
Euler exponents.

Fix the norm systems, the coefficient ring for the
original Selmer lattices and avatar congruences, the
Selmer truncation, and the neighborhood of \(g_0\)
before choosing the auxiliary exponent. The uniform
bound on the number of auxiliary places permits this
exponent to be chosen before the actual auxiliary primes.
After selecting that finite prime set, the polynomials
for specialization have integral coefficients in
\(\overline{\mathbb Q}_\ell\). They and the finitely many
test points lie in a finite extension \(L\) of the fraction field
\(K_0\) of the fixed Selmer coefficient DVR. The bound \(\ell^H\)
in Proposition~\ref{prop:joint-norm-selmer-twist} precedes the
auxiliary level \(\ell^{N_{\rm tor}}\), where \(N_{\rm tor}\)
denotes the exponent called \(N\) in that proposition. Thus all possible torsion coordinates
already belong to \(K_0(\mu_{\ell^H})\). After choosing \(Q=J\),
one may evaluate all the Selmer minors and all signed norm polynomials
in \(L(\mu_{\ell^H})\), with the same normalized valuation.
The witness points for the different polynomials may differ.
Lemma~\ref{lem:common-torsion-specialization} replaces them by one
tuple for all rows, both signs and all powers \(1\le j\le m\).
Its tolerance is \(2^{-\ell-1}\); the separate witness-point errors
are at most \(s_0a_0\le2^{-\ell-1}\), giving \(2^{-\ell}\) in total.
For example, at \(\ell=2\) take \(a_0=1/(8s_0)\) and an integer
\[
 N_{\rm tor}>\max\{H,8s_0M,\max_i M'_i,\max_i(8s_0\beta_i^{\rm norm}+s_0)\},
\]
with the last two maxima over the norm tests. Each selected coordinate
of order \(2^h\) has \(2^{1-h}\le a_0/m\). Consequently
\(v_2(\zeta_w^j-1)=2^{1-h+v_2(j)}\le a_0\) for every
\(1\le j\le m\).
For the triangle and turn factors of Section~\ref{sec:critical-exchange-height},
the auxiliary powers are respectively \(j-i\) and \(j+1\).
Their old tame powers are \(j-i\) and \(i+j+2\), between \(1\)
and \(2m<q_2\), and their weight bounds include the required critical
points. The same list contains \(p_{\rm aux}=(\chi_M\upsilon_1^{-2},2)\)
from Section~\ref{sec:four-block-point-test}. Its old unit character
at each \(w\in S_M\) is the character of exact order \(q_2\) from
Section~\ref{sec:marked-level-characters}. Its infinity weights are
\(-2\), and its new twisting power is the same \(j\) as for the main pair.

The forbidden set consists of the places over \(\ell,p_0,p_1\)
and the rational level prime \(N\) of
Proposition~\ref{prop:marked-level-primes}. This prime is fixed
before the exponent of the auxiliary level is chosen.
The resulting degree-one support \(Q\)
is disjoint from its conjugate and avoids this set. The
norm-one quotient law identifies the selected product of
characters with the inverse avatar of one torus character
\(\xi_\kappa\). At \(w\in Q\), if \(\zeta_w\) has order \(\ell^h\), then
\[
 \frac{1}{\varphi(\ell^h)}=v_\ell(\zeta_w-1)\le\frac{a_0}{m}.
\]
Since \(a_0<1\), we have \(\ell^h>m\). For \(1\le j\le m\),
\(\zeta_w^j\ne1\) and
\(v_\ell(\zeta_w^j-1)=\ell^{v_\ell(j)}v_\ell(\zeta_w-1)\le a_0\).
The other tame characters are trivial on inertia at \(w\).
Thus each power has exact support \(Q\) on the chosen side and
\(\bar Q\) on the conjugate side.

At the rational prime \(d\) the cyclotomic projection of the
decomposition group is generated by \(d\bmod q_1\). It omits \(-1\)
by Lemma~\ref{lem:auxiliary-residue-class} after reduction modulo
\(q_0\), so every place of \(F_0\) over \(d\) splits in \(E\).
The residue degree over \(\mathbb Q\) there is at least two since
\(d\bmod q_0\) has order two. The prime \(N\) splits completely by
Proposition~\ref{prop:marked-level-primes}. The places over \(d,N\)
therefore lie outside \(S_M\cup\bar S_M\cup Q\cup\bar Q\).
The same separation holds at each rational prime \(p\mid q_0\).
Indeed \(p\) splits in \(K\), so every place above it splits in \(E/F_0\).
The absolute residue degree in \(E\) is even at \(2\) since \(\zeta_3\in E\).
At odd \(p\mid q_0\) it is a multiple of \(\operatorname{ord}_{2^s}(p)>1\),
since \(\zeta_{2^s}\in E\) and \(1<p<2^s\) for all sufficiently large \(x\).
Thus no such place belongs to the degree-one old or auxiliary sets.
The Euler products retain the two split factors at every place over
\(q_0dN\), with the actual unramified character values on uniformizers. Their
local conductor exponents are zero in both characters. Consequently
only the places in \(Q\cup\bar Q\) enter the conductor quotient
\(q_Q^2\) in \eqref{eq:determinant-marked-ratio}. This calculation
uses the complex Hecke characters also when \(\ell=d\) or \(\ell=N\).

Write \(u_\ell=v_\ell(A_{\rm aut})\). Proposition~\ref{prop:joint-norm-selmer-twist}
gives the following bounds in the trace comparison for each
selected power \(\eta=\xi_\kappa^j\), with \(\psi\) the inverse avatar
of \(\xi_\kappa\):
\[
 \begin{aligned}
 \mathrm{sel}&\le (4f+3)r u_\ell+2^{-\ell}=e_\ell,\\
 \mathrm{ram}&=\sum_{w\in Q}v_\ell\bigl(\psi^j(\gamma_w)-1\bigr)
                   \le |Q|a_0\le 2^{-\ell-1}.
 \end{aligned}
\]
The first line bounds the normalized Selmer lengths with divisible
coefficients, since each main summand has rank at most \(r\).
The second sum uses one member of each auxiliary pair.
The norm-one quotient \(y\mapsto y/\bar y\) gives
\(\psi(cgc^{-1})=\psi(g)^{-1}\). The conjugate generators have
equal valuations since \(v_\ell(u^{-1}-1)=v_\ell(u-1)\) for units.
The character is trivial on each \(G_{E_v}\) over \(\ell\), so the
local crystalline lattices there are obtained by scalar extension.
At \(Q\cup\bar Q\) the residue characteristic is different from
\(\ell\), and the contribution is that of the tame character.
Changing a tame generator raises its \(\ell\)-power root to a
prime-to-\(\ell\) power and preserves the valuation. The
coefficient-extension inequality above keeps the Selmer bound in
normalized units. The norm tests independently give
\[
 v_\kappa F_{p,Q}((\xi_\kappa^j)^\epsilon)
  \le v_\kappa F_{p,\varnothing}(1)
          +4v_\ell(q_{w_*}-1)+2^{-\ell}.
\]
Since \(a_p\ne0\), the right side is finite, and the numerator
is nonzero at the selected specialization.
Since \(\ell\nmid q_Q\), the exact-support formula
\eqref{eq:critical-norm-ratio-valuation} gives
\[
 \begin{aligned}
 (1+\mathbf1_{z=1})v_\kappa R_{p,\epsilon}(\xi_\kappa^j)
  &=v_\kappa F_{p,Q}((\xi_\kappa^j)^\epsilon)
       -v_\kappa F_{p,\varnothing}(1)\\
  &\le 4v_\ell(q_{w_*}-1)+2^{-\ell}.
 \end{aligned}
\]
Consequently
\[
 v_\kappa R_{p,\epsilon}(\xi_\kappa^j)\le e_\ell .
 \tag{RT1}\label{eq:display-rt1}
\]
More precisely, put \(d_z=1+\mathbf1_{z=1}\) and
\(b_\ell=4v_\ell(q_{w_*}-1)+2^{-\ell}\). At \(\eta=\xi_\kappa^j\),
integrality of \(F_{p,Q}(\eta^\epsilon)\) above \(\ell\) gives
\[
 -\frac{v_\kappa(a_p)}{d_z}
 \le v_\kappa R_{p,\epsilon}(\eta)
 \le \frac{b_\ell}{d_z}.
\]
The same bounds hold at \(\bar\kappa\) for this \(\eta\), since
\(R_{p,\epsilon}(\eta)\) and \(a_p\) are real.
For \(z=1\), including at \(\ell=2\), this gives
\[
 \max_{\lambda\in\{\kappa,\bar\kappa\}}
 \left|R_{p,\epsilon}(\eta)^{-1}\right|_\lambda
 \le\ell^{b_\ell/2},\qquad |x|_\lambda=\ell^{-v_\lambda(x)} .
\]
Here \(q_{w_*}-1\mid A_{\rm aut}\), so
\(4v_\ell(q_{w_*}-1)\le (4f+3)r\,v_\ell(A_{\rm aut})\). In particular,
at the same \(\eta=\xi_\kappa^j\) the second block of
Section~\ref{sec:four-block-class-character} satisfies
\[
 \frac{A_2(\eta)}{A_2(1)}
 =\frac{F_{p_{\rm aux},Q}(\eta)}{F_{p_{\rm aux},\varnothing}(1)}
 =R_{p_{\rm aux},+}(\eta).
\]

\smallskip
\subsection{Integral clearing and projected pairings}
\label{sec:critical-integral-pairings}

\sectionopening{In this subsection, we apply Proposition~\ref{prop:clearing-on-original-tests} to the models of Section~\ref{sec:critical-hecke-systems}.}

The first-order comparison on the selected twisted summands
gives the denominator-clearing operator on the projected image and the rank-one
projected pairing estimate below. With its local defect \(F_\ell\), odd-class integer \(X_{\rm odd}\),
original rank-one level class number \(h_1\), and \(Y_{\rm tr}=r^8\)
from Section~\ref{sec:triangle-unit-input}, put
\[
 \begin{aligned}
 C_\ell={}&4fF_\ell+(1+2t)e_\ell+t v_\ell(X_{\rm odd})\\
   &+\bigl(8ft+3+(4t+1)Y_{\rm tr}^{22}+Y_{\rm tr}\bigr)v_\ell(A_{\rm aut}),\\
 B_\ell={}&C_\ell+v_\ell(h_1),\qquad J_\ell=B_\ell\log\ell .
 \end{aligned}
\]
On the lift image \(Z\), the operator \(C_GD_Z\) clears the
orthogonal projections of integral tests perpendicular to all \(F_w\).
Its determinant valuation is at most \(r'\dim(Z)C_\ell\), where \(r'\)
is the source rank of this lift. This follows by summing the
eigenvalue bound below. For the rank-\(r\) source, \(\dim Z\le\dim V\).
In rank one let \(\Pi\) project to the normalized
line \(\mathbb Cu\) for \(\eta=\xi_\kappa^j\).
The two tests \(a,a'\) must be integral at every coefficient embedding over \(\ell\) and perpendicular to the false
lifts defined by the kernel of the full descent.

The following pairing and central-average calculation
gives the line estimate
\[
 v_\kappa\langle\Pi a,\Pi a'\rangle\ge-B_\ell .
 \tag{RT2}\label{eq:display-rt2}
\]
\smallskip
\indent\textit{Integral pairings after the first-order comparison.}
Let \(D_Z\) be the denominator-clearing operator on the projected image obtained
from the first-order coefficients. The original class
pairing is the sum of the transported Weyl pairings with
denominators the stabilizer orders. The transport and
the integral Weyl evaluations use the same original
bases. The left test is integral also at the conjugate embedding used for complex conjugation. Each stabilizer
order divides \(C_G\) by the reduction argument in
Section~\ref{sec:marked-norm-line}. Thus
\[
 v_\kappa\langle a,C_G D_Z\pi_Z a'\rangle\ge0
\]
for the original tests, with the right test perpendicular
to the false lifts. The integer multiplication does not
affect the integral image property. The valuation per
eigenvalue is at most
\[
 \begin{split}
 &\mathrm{sel}_Z+2t'\mathrm{ram}+4fr'F_\ell+\mathrm{err}_\ell+u_\ell\\
 &\qquad\le r'\bigl(4fF_\ell+(1+2t)e_\ell+
                                     \mathrm{err}_\ell+u_\ell\bigr).
 \end{split}
\]
Here \(\mathrm{sel}_Z,\mathrm{ram}\le e_\ell\) since the selected twist
bounds each \(\mathrm{sel}_i\), and \(1\le r'\le r\), \(t'\le t\). This is the allowance \(r'C_\ell\)
in the determinant on the projected image. The \(C_G\) for this pairing
is in addition to the \(C_G^2\) in the first-order coefficient.

\smallskip
\indent\textit{Rank-one averaging on the original line.}
Let \(T=U(1)\), let \(K_1\) be the original chosen line level,
and let \(s_c\) represent the \(h_1\) classes of
\(T(F_0)\backslash T(\mathbf A_{F_0,f})/K_1\).
Take the representatives integral at \(\ell\).
To do this, write the inverse local component of each proposed
representative as \(a_v/\bar a_v\) by Hilbert 90, using a
split pair of factors at a split place. Approximate the
\(a_v\) at the places above \(\ell\) by a single \(a\in E^\times\).
Replacing the representative by its product with
\(a/\bar a\) puts the local components in the open integral compact.

For the main pair \((\nu,z)\) and its twist \(\eta\), let
\(\beta=(\nu\eta)\upsilon_1^z\) be the full target-leg
character. For a norm-one scalar \(x\) let \(T(x)\) act by
right translation by \(x I_{2z}\), multiplied by \(\beta(x)\).
This operator is unitary, commutes with the type projector
and spherical operators, and has adjoint \(T(x^{-1})\).
For the original line kernel \(K(g)\) one has
\[
 T(x^{-1})K(g)=K(gx).
\]
Indeed the simultaneous scalars \((x^{-1},x^{-1}I_{2z})\)
fix the Hom matrix variable. In the Weil action
their multiplier is \(\beta(x^{-1})^{-1}\), which
cancels that of \(T(x^{-1})\). The line stabilizers are trivial.
With \(s=V_0^{2z}J\ne0\) as in Section~\ref{sec:v2-four-block-main-numerator},
the raw and normalized lifts satisfy
\[
 u_{\rm raw}=\sum_c K(s_c),\qquad u=\phi_\eta=s\,u_{\rm raw}.
\]
Translation permutes the classes, so \(T(x)u_{\rm raw}=u_{\rm raw}\)
and \(T(x)u=u\).

Here \(K(g)=K_X(g)\) is the raw kernel of the four-block
test, with the original inclusion column at the old places, the
character on the distinguished unit coordinate at \(Q\), and the
vacuum elsewhere. In particular the normalized descent to the full
line form space is the same \(D_O\) as in
Proposition~\ref{prop:four-block-numerator}. Indeed contraction
with the polynomial of the one-row Weyl type and summation over
the original \(J_Q\)-classes give
\[
 (\mathcal D y)(g)=(D_Oy)(g)
    =\sum_{[h]_{J_Q}}
       \frac{\langle K_X(g,h),y(h)\rangle_{\rm Weyl}}
            {|\operatorname{Aut}(h,J_Q)|}
    =\langle K(g),y\rangle_{J_Q}.
\]
Set \(S=\sum_c T(s_c)\). Commutativity of the torus gives
\[
 \mathcal D(Sy)(g)
    =\left\langle\sum_c T(s_c^{-1})K(g),y\right\rangle
    =\langle u_{\rm raw},y\rangle
    =\bar s^{-1}\langle u,y\rangle .
\]
The equality of the kernel sum uses that \(gs_c\) runs
over precisely the same line classes.
For Proposition~\ref{prop:clearing-on-original-tests} take
\(Z=\mathbb C u\) and \(e_0=Q=e_{X-1}\).
Indeed \(u\) is a nonzero multiple of the lift of the
original constant line form. Its Hecke system \(x_0\)
and its \(X-1\) type are those of
Section~\ref{sec:critical-hecke-systems}, and \(\pi_Z=\Pi\).
At every split off-level place \(w\), the descent identity and
the commutation of \(S\) with \(Q\) and \(K_w(X_w)\) give
\[
 S\bigl(\operatorname{ran}Q\cap u^\perp\cap\ker K_w(X_w)\bigr)
   \subseteq\operatorname{ran}Q\cap\ker\mathcal D\cap\ker K_w(X_w).
\]
The original test \(a'\) is perpendicular to the right-hand side.
Hence \(S^*a'\perp F_w\) for this line at every such \(w\).
Also \(S^*a'\) is integral at every embedding above \(\ell\):
each scalar translation by \(s_c^{-1}\) preserves the
integral local Weyl lattice norm. For a finite norm-one idele \(x\)
integral above \(\ell\), its class in the original quotient is killed
by \(h_1\). Write \(x^{h_1}=\gamma k\) with \(\gamma\in T(F_0)\)
and \(k\in K_1\). Then \(\gamma\) is an \(\ell\)-unit and automorphy gives
\[
 \beta_f(x)^{h_1}=\beta_\infty(\gamma)^{-1}\beta_f(k).
\]
The first factor is a monomial in conjugates of an \(\ell\)-unit.
The second is a root of unity, since \(\beta_f\) is smooth on
the compact \(K_1\). Thus \(v_\kappa(\beta_f(x))=0\) at every
coefficient embedding above \(\ell\).
Since \(Su=h_1u\), orthogonal projection gives
\(\Pi S^*a'=h_1\Pi a'\). Let \(s_{\rm tr}=AD_I\) be the Hecke
coefficient of Section~\ref{sec:top-coefficient-clearing} for the
system \(x_0\) of \(u\). Here \(Qu=u\), \(c=C_G=A_{\rm aut}\), and
the eigenvalue of \(C_Gc^2s_{\rm tr}Q\) on \(\mathbf C u\) is
\[
 d_u=C_Gc^2s_{\rm tr}(x_0)\ne0,\qquad
 v_\kappa(d_u)\le
 u_\ell+\mathrm{sel}+2t'\mathrm{ram}+4fF_\ell+\mathrm{err}_\ell
 \le C_\ell .
\]
Proposition~\ref{prop:clearing-on-original-tests}
applies to the original integral test \(S^*a'\) just verified.
It gives
\(c^2s_{\rm tr}Q\Pi S^*a'\in L_{\rm orig}\)
and the valuation bound with \(r'=1\) before the extra
\(u_\ell=v_\ell(C_G)\).
In each class the two transports satisfy
\((T^\#)^{\mathsf t}G_{\rm W}T=G_{\rm W}\) for the integral
Weyl Gram matrix, as in Section~\ref{sec:gram-product-formula}.
Pairing this integral image at \(\kappa\)
with \(a\), integral at \(\bar\kappa\), and clearing the class
stabilizers by \(C_G\) gives
\[
 d_u h_1\langle\Pi a,\Pi a'\rangle
   =\langle a,C_Gc^2s_{\rm tr}Q\Pi S^*a'\rangle
   \in\mathcal O_\kappa .
\]
The scalar acts in the second, linear variable. Taking valuations gives
\eqref{eq:display-rt2}.

For the class-number factor, put \(T(\widehat{\mathcal O}_E)
=\prod_{v<\infty}T(\mathcal O_{E_v})\), with the analogous
interpretation at split places. The norm-one ideal map from
classes at this maximal compact into \(\mathrm{Cl}(E)\)
has kernel of order at most \(2^f\). For a principal ideal
\((a)\) in the kernel, \(a\bar a\) is a unit of \(F_0\).
Changing \(a\) by a unit changes this unit by a norm, and
if the resulting unit class is trivial there is a norm-one
generator.

Thus the kernel injects into
\(\mathcal O_{F_0}^\times/N_{E/F_0}(\mathcal O_E^\times)\).
The norm subgroup contains \((\mathcal O_{F_0}^\times)^2\).
The totally real field \(F_0\) has degree \(f\), so its unit group
has rank \(f-1\) and torsion \(\{\pm1\}\). Its quotient by squares
has size \(2^f\). The index from this maximal
compact to \(K_1\) is at most \((MN)^{2f}\), since \(K_1\)
contains the principal level \(MN\), and reduction modulo
\(MN\mathcal O_E\) has at most that many values. Hence
\[
 h_1\le 2^f h_E(MN)^{2f}.
\]
By the Minkowski bound \cite[Thm.~4.3]{MilneANT}, each class has
an integral representative of norm at most
\(D=(4/\pi)^f(2f)!\sqrt{|\Delta_E|}/(2f)^{2f}\).
Thus
\[
 h_E\le\sum_{\mathrm N\mathfrak a\le D}1
 \le D^2\zeta_E(2)\le 2^{2f}|\Delta_E|\zeta(2)^{2f}.
\]
All logarithms on the right-hand side of the bound for
\(h_1\) are polynomially bounded in \(x\) by the parameter
bounds in Section~\ref{sec:datum}.

For every finite set \(U\) of rational coefficient primes, factorization
gives \(\sum_{\ell\in U}v_\ell(H)\log\ell\le\log H\) for each
positive integer \(H\). Also \(C_G=A_{\rm aut}\) and
\(\sum_{\ell\in U}2^{-\ell}\log\ell\le\sum_{n\ge1}n2^{-n}=2\).
The explicit defect bound of Section~\ref{sec:coefficient-local-defects}
and the expressions for \(C_\ell,B_\ell\) give
\[
\begin{aligned}
 \sum_{\ell\in U}J_\ell
 &\le 108\cdot10^6\frac{fm}{k}\,b\log b
       +(1+2t)\bigl((4f+3)r\log A_{\rm aut}+2\bigr)\\
 &\quad +t\log X_{\rm odd}+\log h_1\\
 &\quad +\bigl(8ft+3+(4t+1)Y_{\rm tr}^{22}+Y_{\rm tr}\bigr)
          \log A_{\rm aut}.
\end{aligned}
\]
The first term divided by \(J_*=fxb(\log b)^3/(mk)\) is at most
\(108\cdot10^6m^2/(x(\log b)^2)\to0\), since \(m^2\le4x\).
All other terms are polynomially bounded in \(x\), whereas \(J_*\)
and \(x_1\) each exceed every fixed power of \(x\) for fixed
\(q_0,d,\epsilon\). Each allowance bounds the local term at every
coefficient embedding over its rational prime.

For every finite set \(U\) of coefficient primes,
\[
 \sum_{\ell\in U}e_\ell\log\ell\le
  (4f+3)r\log A_{\rm aut}+2\le x_1,\qquad
 \sum_{\ell\in U}J_\ell\le J_*=\frac{fa(\log b)^3}{mk}.
 \tag{RT3}\label{eq:display-rt3}
\]
Here \(\log A_{\rm aut}\le (rt)^2 f\log M\), since the original
set contains one split half above each place over the two level
primes. The allowances are independent of the coefficient embedding,
and the twist \(\xi_\kappa\) may be chosen at each \(\kappa\). The first-order comparison is applied at
\(\kappa\) on that chosen twist and does not use the twist
selected at \(\bar\kappa\).

\smallskip
\subsection{Reciprocal heights of the main critical values}
\label{sec:critical-main-value-height}

\sectionopening{In this subsection, we use the local estimates of Section~\ref{sec:critical-integral-pairings} and the numerator test of Proposition~\ref{prop:four-block-numerator}.}

For each ramified pair \(p=(\nu,z)\), each sign \(\epsilon\), and
\(\eta=\xi_\kappa^j\), equations \eqref{eq:display-rt1} and
\eqref{eq:display-rt3} give
\[
 v_\kappa(R_{p,\epsilon}(\eta)^{-1})\ge-e_\ell,\qquad
 \operatorname{Height}\bigl((R_{p,\epsilon}(\xi_\kappa^j)^{-1})_\kappa\bigr)
 \le x_1
\]
with multiplier \(1\). This includes \(z=1\), where the valuation of
the squared-period ratio is \(2v_\kappa(R_{p,\epsilon}(\eta))\)
because \(q_Q\) is an \(\ell\)-adic unit.

For the stronger estimate fix a main pair for \(j\), so \(m<z\le r\),
and take the positive sign. At both conjugate embeddings
use \(\eta=\xi_\kappa^j\) and the normalized lift
\(u=\phi_\eta\) of Section~\ref{sec:v2-four-block-main-numerator}.
Summing the inverse unit residues in the class formula of
Section~\ref{sec:marked-norm-line} gives this same vector
\(\Phi_\eta\), with vacuum and jet scalar
\(V_0^{2z}(\pi_Ms_{\rm jet})^{\sum_\tau(c_\tau-z)}\).
Use \(u_1=\phi_1\) at \(Q=\varnothing\) and
\(\rho=R_{p_{\rm aux},+}(\eta)>0\).
The numerator estimate gives a fixed \(c_0\ne0\) with
\(\log|\sigma c_0|\le J_*\) for every complex automorphism
and, for any coprime modulus and any coefficient embedding, an integral
test satisfying
\[
 v_\kappa\langle a,u\rangle\le v_\kappa(c_0\rho)+P_\ell,
 \qquad P_\ell=\sum_\tau S_\ell(\kappa_\tau,\tau).
 \tag{RT4}\label{eq:display-rt4}
\]
Here \(S_\ell\) is the local prime sum of
Section~\ref{sec:short-normalization}. The test is integral at every coefficient place above this
\(\ell\) and has the false-lift orthogonality of \eqref{eq:display-rt2}. This is the four-block
numerator test on the original rank-one lift.

Apply \eqref{eq:display-rt4} at \(\bar\kappa\) to the same \(\eta,u\) and
call the resulting test \(a'\). Using that \(\rho\) is real
and conjugating this inequality gives, for
\(I=\langle a,u\rangle\langle u,a'\rangle\ne0\),
\[
 v_\kappa(I)\le v_\kappa(c_0\bar c_0)+2v_\kappa(\rho)+2P_\ell.
\]
Apply \eqref{eq:display-rt2} to these two tests at \(\kappa\).
Proposition~\ref{prop:marked-line-norm} uses \(J_Q\) for \(u\)
and the old vector fixer \(J_\varnothing\) for \(u_1\). At the
trivial twist \(Q=\varnothing\) and \(q_\varnothing=1\). The vacuum and jet
scales agree in both norms. Thus
\[
 \frac{\|u\|^2}{\|u_1\|^2}=q_Q^{2z-1}R_{p,+}(\eta),\qquad
 \langle\Pi a,\Pi a'\rangle=\frac{I}{\|u\|^2}.
\]
Both norms are nonzero by the positive Euler products. It follows that
\[
 R_{p,+}(\eta)^{-1}
 =q_Q^{2z-1}\|u_1\|^2
       \frac{\langle\Pi a,\Pi a'\rangle}{I}.
 \tag{RT5}\label{eq:display-rt5}
\]
For algebraic \(u\), one has \(v_\kappa(\bar u)=v_{\bar\kappa}(u)\).
The numerator test at \(\bar\kappa\) uses the same \(\eta\), whose
modulus satisfies \(\ell\nmid q_Q\). The projection estimate
\eqref{eq:display-rt2} is at \(\kappa\), with the twist selected there.

For the untwisted term use the integral algebraic norm
\(I_{\rm ref}\ne0\) and the fixed empty reference pair of rank
\(2z\) from Section~\ref{sec:four-block-reference-norm}. Its multiplier
\(H=I_{\rm ref}/\|u_1\|^2\) satisfies
\[
 H\|u_1\|^2\text{ integral at every coefficient embedding},\qquad
 \log|\sigma H|\le J_*-\sum_\tau J(\kappa_\tau).
 \tag{RT6}\label{eq:display-rt6}
\]
Put \(h=Hc_0\bar c_0\). Since the bound for \(c_0\) applies at
\(\sigma\) and \(\sigma\circ\mathrm{conj}\),
\[
 \log|\sigma h|\le 3J_*-\sum_\tau J(\kappa_\tau).
\]
Using \(H\|u_1\|^2=I_{\rm ref}\), \(v_\kappa(q_Q)=0\) and the bound
on \(v_\kappa(I)\) in \eqref{eq:display-rt5} gives
\[
 \begin{aligned}
 v_\kappa(hR_{p,+}(\eta)^{-1})
 &\ge v_\kappa(I_{\rm ref})+v_\kappa\langle\Pi a,\Pi a'\rangle
             -2v_\kappa(\rho)-2P_\ell\\
 &\ge -B_\ell-2e_\ell-2P_\ell
 \end{aligned}
\]
by \eqref{eq:display-rt1}, \eqref{eq:display-rt2} and
\eqref{eq:display-rt6}. Also
\[
 2\sum_{\ell\in U}P_\ell\log\ell\le\sum_\tau S(\kappa_\tau,\tau).
\]
Indeed \(2S_\ell(\kappa_\tau,\tau)\log\ell\) is the sum of the terms
of \(S(\kappa_\tau,\tau)\) at \(n=\ell^e\). Terms with \(\ell^e>a\)
vanish, since \(\kappa_\tau/\ell^e<1\) whereas every weight is at least
one. All terms are nonnegative.
By \eqref{eq:display-rt3}, the fixed multiplier has total allowance at most
\[
 3J_*-\sum_\tau J(\kappa_\tau)+J_*+2x_1+
       \sum_\tau S(\kappa_\tau,\tau)
 \le 6J_*-\sum_\tau\bigl(J(\kappa_\tau)-S(\kappa_\tau,\tau)\bigr),
\]
since \(x_1\le J_*\).

\begin{proposition}[Reciprocal height of the main critical ratio]
\label{prop:main-critical-height}
Let \((\xi_\kappa)_\kappa\) be the chosen family. For every main pair
\(p\) and its lift index \(j\),
\[
 \operatorname{Height}\bigl((R_{p,+}(\xi_\kappa^j)^{-1})_\kappa\bigr)
      \le 6J_* - d_{\rm dir},
 \qquad
 d_{\rm dir}=\sum_\tau\bigl(J(\kappa_\tau)
                                      -S(\kappa_\tau,\tau)\bigr).
\]
\end{proposition}

\smallskip
\subsection{Height of the exchange product}
\label{sec:critical-exchange-height}

\sectionopening{In this subsection, we use the same auxiliary twists as in Proposition~\ref{prop:main-critical-height} for the exchange factors.}

Use the untwisted lift character \(\mu_j^0=\chi_0(\upsilon_1^{-\Delta}\chi_M)^j\).
Its weight at \(\tau\) is \(-w_j(\tau)\), where
\(w_j(\tau)=\kappa_{j,\tau}=\kappa_j(\tau)=\kappa_\tau+(j+\tfrac12)\Delta\).
The following table lists the characters and endpoints in reverse order
of the measured exchanges of Section~\ref{sec:gluing-three-exchanges}.
\[
 \begin{array}{c|c|c}
 \text{character}&\text{upper point}&\text{sign}\\ \hline
 \mu_j^0\upsilon^{j+(t-11m-1)/2}&7(m/2)+k/2&+\\
 \upsilon^{-m}\chi_M^j&m&+\\
 \upsilon^{-m/2}\chi_M^{-j}&m/2&-
 \end{array}
\]
In this table \(\upsilon=\upsilon_1=\upsilon_E|_{T_1(\mathbb A_{F_0})}\).
Lowering the point keeps the old character and tame support.
For the interval with upper endpoint
\(z_*=(7m+k)/2=(r+6m+1)/2\), put
\(n_j=j+(t-11m-1)/2\). For large \(x\) one has \(r\ge6m+1\);
hence \(z_*\le r\) and \(0\le n_j\le t-r\), since \(t=5r\)
and \(1\le j\le m\). At the \(k/2\) integers
\(7m/2<u\le z_*\) one therefore has \(m<u\le r\), old tame
character \(\chi_M^j\), and
\[
 w_j(v)+r-t\le c_v=w_j(v)-n_j\le w_j(v),\qquad
 c_v\ge r\ge u.
\]
Here \(w_j(v)>2(r^8)^2+t\), as in
Section~\ref{sec:critical-hecke-systems}. Thus
Proposition~\ref{prop:main-critical-height} applies to these points.
For the remaining \(7m/2+m+m/2=5m\) pairs the original tame powers
are \(\pm j\) and their weights are at least their respective upper
points. Equation~\eqref{eq:display-rt1} bounds each reciprocal
height by \(x_1\), using multiplier \(1\) and \eqref{eq:display-rt3}. Write \((\chi_i,z_i,\epsilon_i)\) for the three rows of the table,
with \(\epsilon_i\in\{1,-1\}\). The three exchange identities of
Section~\ref{sec:gluing-three-exchanges} give
\[
 \frac{H_j(\eta)}{H_j(1)}
  =q_Q^{\sum_i z_i^2}
    \left(\prod_{i=1}^3\ \prod_{1\le u\le z_i}
                   R_{(\chi_i,u),\epsilon_i}(\eta)\right)^2.
\]
Write \(M_\lambda^S(z)\) for the factor in
Section~\ref{sec:marked-partner-period} with split support \(S\).
The character on row \(i\) is \((\chi_i\eta^{\epsilon_i})^{-1}\),
whose base change is \(((\chi_i\eta^{\epsilon_i})^{\rm bc})^{-1}\).
Thus the quotient for this row is
\[
 \frac{M_{(\chi_i\eta^{\epsilon_i})^{-1}}^{S_M\cup Q}(z_i)}
      {M_{\chi_i^{-1}}^{S_M}(z_i)}
 =q_Q^{z_i^2}\prod_{u=1}^{z_i} R_{(\chi_i,u),\epsilon_i}(\eta)^2 .
\]
At a new pair the numerator omits both Euler factors and the
untwisted denominator retains both. The conductor counts the chosen
half once, with exponent \(\sum_{u=1}^{z_i}(2u-1)=z_i^2\)
by \eqref{eq:plane-period-local-norm}. The old conductor factors cancel.
For \(p_{i,u}=(\chi_i,u)\), set
\[
 U_{i,u}=\frac{F_{p_{i,u},Q}(\eta^{\epsilon_i})}
                  {F_{p_{i,u},\varnothing}(1)} .
\]
For exact support \(Q\), Section~\ref{sec:critical-signed-systems}
gives \(U_{i,u}\in\overline{\mathbf Q}^{\times}\),
\(R_{p_{i,1},\epsilon_i}(\eta)^2=U_{i,1}/q_Q\),
and \(R_{p_{i,u},\epsilon_i}(\eta)=U_{i,u}\) for \(u\ge2\).
Each of the three positive endpoints includes \(u=1\), so
\[
 \frac{H_j(\eta)}{H_j(1)}
 =q_Q^{z_*^2+m^2+(m/2)^2-3}
       \prod_{i=1}^3\left(U_{i,1}\prod_{u=2}^{z_i}U_{i,u}^2\right)
 \in\overline{\mathbf Q}^{\times}.
\]

For each \(\kappa\) the chosen \(Q\) avoids \(\ell\), so the reciprocal
conductor power has valuation zero there and requires only multiplier \(1\).
Put \(A_j=\{u\in\mathbf Z:7m/2<u\le z_*\}\). For \(u\in A_j\)
let \(X_{j,u}\) be the original untwisted rank-\(2u\) datum for the third
measured exchange, the first row of the table. Use its numerator constant \(c_{0,j,u}\) from
Proposition~\ref{prop:four-block-numerator} and put
\[
 h_{j,u}=\frac{C_G M^{8fu}\|\phi_{(\upsilon_1^{-u},u)}\|^2}
                  {\|\phi_{X_{j,u}}\|^2},
 \qquad
 h_j=\prod_{u\in A_j}(h_{j,u}c_{0,j,u}\bar c_{0,j,u})^2.
\]
Both norms use the same old level. The reference in the numerator
has row degree zero and denominator \(M^{4fu}\), as in
Section~\ref{sec:four-block-reference-norm}. Thus \(h_j\) is fixed
before \(\kappa\) and \(Q\). At \(\eta=\xi_\kappa^j\), the \(k/2\)
main factors have sign \(+\). The remaining \(5m\) factors of either
sign use multiplier \(1\). Using \eqref{eq:display-rt1}
and \eqref{eq:display-rt5}--\eqref{eq:display-rt6} before summing gives
\[
 v_\kappa\!\left(h_j\frac{H_j(1)}{H_j(\eta)}\right)
 \ge -k(B_\ell+2e_\ell+2P_\ell)-10m e_\ell .
\]
Since \(\log|\sigma h_j|\le k(3J_*-\sum_\tau J(\kappa_\tau))\),
the finite-prime sums above give
\[
 \operatorname{Height}\left(
   \left(\frac{H_j(\xi_\kappa^j)}{H_j(1)}\right)^{-1}\right)
      \le 10mx_1+k(6J_*-d_{\rm dir}).
\]
For \(1\le i<j\le m\) and \(0\le a,b'<m\), respectively, the
triangle and turn products on the same family are
\[
 \begin{aligned}
 Q_{ij}(\xi_\kappa)
   &=\prod_{u=1}^{t}R_{(\mu_j^0/\mu_i^0,u),+}(\xi_\kappa^{j-i}),\\
 T_{ab'}(\xi_\kappa)
   &=\prod_{u=1}^{t}R_{(\mu_{a+1}^0\mu_{b'+1}^0,u),+}
                                      (\xi_\kappa^{b'+1}).
 \end{aligned}
\]
The local negative weights are \((j-i)\Delta\) and
\(w_{a+1}(v)+w_{b'+1}(v)\). Both lie in \([t,4b]\), since
\(t\le\Delta\), \(\kappa_\tau\le b\), and \((m+1)\Delta\le b\)
for sufficiently large \(x\). The old tame powers \(j-i\) and
\(a+b'+2\) lie strictly between zero and \(q_2>2m\), so the old
support is exactly \(S_M\). The powers \(j-i,b'+1\) lie in
\([1,m]\), hence have the same exact new support \(Q\).
By \eqref{eq:display-rt1}, each of
\(Q_{ij}(\xi_\kappa)^{-1}\) and \(T_{ab'}(\xi_\kappa)^{-1}\)
has valuation at least \(-t e_\ell\). Half-integral powers of \(q_Q\)
in these formulas also have valuation zero, since their
squares are integral powers of \(q_Q\).

\section{Determinant factors and the gluing operator}
\label{sec:period-gluing-operator}

\sectionopening{In this section, we track the determinant and measure factors in
Section~\ref{sec:normalized-period-exchange} through three period identities.
The resulting adjoint is the common period used in
Section~\ref{sec:common-period-injectivity}.}

\subsection{Products of determinants in the period comparison}
\label{sec:marked-stacked-determinant-transfer}

\sectionopening{In this subsection, we check how the product of archimedean determinants transforms in the
exchange of Section~\ref{sec:marked-partner-period}. The resulting
coefficient is used on the three pairs in
Section~\ref{sec:gluing-three-exchanges}.}

Let the kept real space have positive and negative
dimensions \(m,\ell\), the old small space dimensions
\(p,\ell\), and the definite target dimension \(a\).
The ordered positive basis satisfies \(m=p+a\).
At each infinite place \(\tau\) let \(c_\tau\) be
the number of stacked determinants. In the large and
small Fock spaces use the polynomials
\[
 P(\eta)=\prod_{\tau,i}\det[\eta_{\tau i};Z_\tau],
 \qquad
 Q(\eta)=
  \prod_\tau\det(Z^-_\tau)^{c_\tau+2s_0}
       \prod_{\tau,i}\det(\eta_{\tau i}Y_\tau),
 \qquad 2s_0=\dim W-2\ell .
\]
Here \(Z_\tau\) is the \(a\)-by-\(m\) block of
positive columns, \(Y_\tau\) is the \(m\)-by-\(p\)
positive block of the small vector, and \(Z^-_\tau\)
is its \(\ell\)-by-\(\ell\) negative block.
Each \(\eta_{\tau i}\) is a \(p\)-by-\(m\) coefficient matrix. The
bracket stacks its rows above the \(a\) rows of \(Z_\tau\) in the
ordered positive basis. Put
\(\eta_0=[I_p\ 0]\) in that basis, at every place
and for every index.

Apply the measured partner comparison with the fixed
small vector \(Q(\eta_0)\). Let \(v\) be its large
real vector. We show that \(v\) can be replaced by
\(P(\eta)\) for every \(\eta\), with one
scalar \(\kappa\) multiplying the new large polynomial. At this stage
\(\kappa\) may be zero. Its nonvanishing in the applications with
nonzero projected source will follow from injectivity of the common period.
Write \(\alpha,\beta\) for the two determinant characters of the large
kernel. For the kernels of \((W,I_a)\) and \((V',W)\) the pairs of
characters are \((\alpha,\beta)\) and \((\alpha/\beta,\alpha)\),
respectively. At \((g,h)\) and \((x,g)\), both sides are multiplied
by \(\beta(\det g)^{-1}\) relative to the choice \(\beta=1\).
We cancel this common factor throughout the compact group calculation.
For fixed matched finite data and \(g\), the large period
is a continuous functional \(L_g\) in a real Fock moment
seminorm.

The matched support in the first coordinate is primitive, so the
kernel estimate at exponent zero in the integrated
variable bounds its absolute integral by \(C_g p_N(v)\).
The functionals \(L_g\) are invariant under the definite
target real group by right translation in the original
quotient.
The standard covariance under right translation in the
kept group gives \(L_g(\omega(k,u)v)=L_{gk}(v)\).
Write \(S_g\) for the small period. Proposition~\ref{prop:normalized-period-exchange}
gives \(L_{gk}(v)=c_\nu S_{gk}(Q(\eta_0))\), with \(c_\nu\ne0\).

Let \(t(z)\) on the positive kept basis be diagonal
with \(z^{-2}\) on the \(p\) coordinates and \(z^{-1}\)
on the \(a\) coordinates, with \(|z|=1\). Use its
identity extension on the negative block, and put
\(u(z)=z^{-1} I_a\) in the target. In the real representation the target determinant exponent
is \(\ell-c_\tau\). The vacuum contributes
\(\det(u_\tau)^{-\ell}\), hence their product contributes
\(\det(u_\tau)^{-c_\tau}\).
For the operator
\[
             v\longmapsto
             z^{-a\sum_\tau c_\tau}\omega(t(z),u(z))v
\]
these scalar factors cancel. On polynomial vectors
the substitution scales all variables outside the
\(a\) columns by \(z\) and fixes the variables
in those \(a\) columns. Indeed on the positive columns
the substitutions use \(u(z)\) and \(t(z)^{-1}\);
on the negative columns the signed Fock coordinate
uses the conjugate factor of \(u(z)\).
Thus averaging over \(z\) sets the variables outside the last \(a\)
positive columns to zero. At each place write \(Z_{\tau,0}\) for the
remaining \(a\)-by-\(a\) block.

Under the definite target average this square
polynomial is projected to its left determinant
weight \(c_\tau\). This weight space consists of the
multiples of \(P(\eta_0)=\prod_\tau\det(Z_{\tau,0})^{c_\tau}\).
Indeed, a polynomial with that left covariance under
unitary matrices has it under the complexification
\(\prod_\tau\mathrm{GL}_a(\mathbb C)\), by differentiation
and complex linearity. On invertible square matrices
its value is determined by its value at the identity,
and polynomial continuation gives the statement.

Average first over the circle and then over the compact target.
These bounded operators commute on Fock space. The kept action commutes
with the target action, and \(u(z)=z^{-1}I_a\) is central in the target.
Polynomials are dense, and the resulting line is
closed. Hence for the fixed \(v\) these two averages
give \(\kappa P(\eta_0)\) for one scalar \(\kappa\),
independent of all the finite data, the modulus,
the character and \(g\).
The averages can also be evaluated inside the period.

To obtain polynomial approximation in the moments used by
\(L_g\), write \(p_k(f;A)\) for the integral of
\(|f(z)|(1+\|z\|)^k e^{-c\|z\|^2}\) over \(A\subset\mathbb C^{d_{\rm F}}\)
in the fixed Fock coordinates, and \(p_k(f)=p_k(f;\mathbb C^{d_{\rm F}})\).
For \(0<t_0\le t<1\), the dilation \(D_tf(z)=f(tz)\) satisfies
\[
 p_k(D_tf;\|z\|>R)\le
 t_0^{-2d_{\rm F}-k}p_k(f;\|z\|>t_0R).
\]
This is the substitution \(w=tz\), with the Gaussian decreasing
under the substitution. The tail bound and local uniform convergence
give \(p_k(D_tf-f)\to0\) as \(t\uparrow1\).
Circle extraction gives \(p_k(f_j)\le p_k(f)\) for each homogeneous
term \(f_j\). Thus, for fixed \(0<t<1\),
\[
 p_k\left(D_tf-\sum_{j=0}^J t^j f_j\right)
 \le\sum_{j>J}t^jp_k(f_j)
 \le \frac{t^{J+1}}{1-t}\,p_k(f).
\]
Orthogonality of the degrees gives the two limits in Fock norm as well.
The compact averages are bounded in the radial moments because
their substitutions are unitary and their character factors have
modulus one. Applying the bound for \(L_g\) passes the polynomial
identities to the period of \(v\) on the original finite level.

For \(k=(k_\tau)\in\prod_\tau U(m)\) put
\(\eta_{k,\tau}=\eta_0 k_\tau^{\mathsf t}\) and
\(d(k)=\prod_\tau\det(k_\tau)^{c_\tau}\).
On the large side the kept positive block acts by
\(Z\mapsto Z(k_\tau^{-1})^{\mathsf t}\), and
\[
 \det\begin{bmatrix}\eta_0\\ Z(k_\tau^{-1})^{\mathsf t}\end{bmatrix}
   =(\det k_\tau)^{-1}
       \det\begin{bmatrix}\eta_0k_\tau^{\mathsf t}\\ Z\end{bmatrix}.
\]
On the small side \(Y\mapsto k_\tau^{\mathsf t}Y\), the negative block is
fixed, and the vacuum and splitting contribute \(d(k)^{-1}\).
Thus the two real actions on the indicated Fock vectors satisfy
\[
 \omega_{\rm large}(k,1)P(\eta_0)=d(k)^{-1}P(\eta_k),\qquad
 \omega_{\rm small}(1,k)Q(\eta_0)=d(k)^{-1}Q(\eta_k).
\]

Apply the measured comparison at the real point
\(k t(z)\) and use these covariances.
The small basis there is \(z^{-2}\eta_k\),
so homogeneity contributes \(z^{-2p\sum c_\tau}\).
The factor \(\det(t(z))^{-c_\tau}\) contributes
\(z^{(2p+a)\sum c_\tau}\), which cancels the preceding
power and the additional \(z^{-a\sum c_\tau}\).
The determinant powers of \(k_\tau\) cancel between
the covariances. Thus the projected identity is
\[
       \kappa L_g(P(\eta_k))=c_\nu S_g(Q(\eta_k)).
\]
This equality uses no division by \(\kappa\).

It remains to justify the passage from \(\eta_k\)
to arbitrary families of coefficient matrices. At a fixed place
consider the polynomial space in \(m\)-by-\(p\)
matrices \(Y\) with \(F(Yu)=\det(u)^c F(Y)\).
For \(p=0\) this is the constant line.
For \(p>0\), it is a finite-dimensional representation
of \(\mathrm{GL}_m\) under left substitutions and
of \(U(m)\) by restriction. It is irreducible.
Indeed, in every irreducible summand of this
polynomial representation a highest weight line is
fixed by the lower unitriangular subgroup, using the
corresponding Borel. A polynomial fixed by that subgroup
is determined by one scalar.

On the dense set where the
first \(p\) rows of \(Y\) form an invertible block \(A\),
a lower block shear sends \(Y\) to \((A;0)\).
The right covariance then gives
\(F(Y)=\det(A)^c F((I_p;0))\).
Thus the entire unitriangular invariant space has
dimension one, and complete reducibility under \(U(m)\)
gives irreducibility. The finite space of homogeneous
degree \(pc\) is preserved by the action, and complexified
Lie algebra invariance identifies the summands under
\(U(m)\) and \(\mathrm{GL}_m\).
Consequently \(\prod_i\det(\eta_iY)\) is a finite
linear combination of the corresponding polynomial
at the matrices \(\eta_0 k^{\rm t}\). The same holds over all
places by tensor products.

The very same coefficients work for the two-sided
stacked determinant. At a place with \(c_\tau=0\) both products over
\(i\) are the constant one, including on a rank-deficient block.
At a place with \(c_\tau>0\), for any fixed lower matrix \(Z\),
if it has rank smaller than \(a\), the stacked
determinant vanishes for every \(\eta_i\).
Otherwise extend its rows to an invertible square
matrix. Multiplying on the right by its inverse
turns \(Z\) into \([0\ I_a]\).
It follows that
\(\det[\eta_i;Z]=d\det(\eta_i Y_Z)\)
for a fixed \(Y_Z\), consisting of the first \(p\)
columns of that inverse, and a scalar \(d\) depending
on \(Z\) alone. Apply the polynomial linear combination
at \(Y_Z\), and use the same coefficients in each
large evaluation and in the small polynomial. Thus for every
family \(\eta\), every matched pair of finite arrays and every
adelic \(g\),
\[
 \kappa L_g(P(\eta))=c_\nu S_g(Q(\eta)).
\]
Put \(t_{\rm ex}=(\dim I_a-\dim V')/2\) and
\(A_{\rm ex}=R'_0/D_{\rm final}>0\) as in
Proposition~\ref{prop:normalized-period-exchange}. Since \(c_\nu\ne0\),
the coefficient multiplying the large class sum is
\[
 S_g(Q(\eta))=\frac{\kappa}{c_\nu}L_g(P(\eta)),\qquad
 \left|\frac{\kappa}{c_\nu}\right|^2
    =\frac{|\kappa|^2}{A_{\rm ex}M_\nu(t_{\rm ex})}.
\]
For the fixed real data and ordered basis, \(A_{\rm ex}\) and \(\kappa\)
are chosen before the finite modulus, its character and the family
of covectors. These identities include \(\kappa=0\).

\subsection{Constants in the three finite comparisons}
\label{sec:marked-gluing-finite-normalization}

\sectionopening{In this subsection, we compute the finite orbit and stabilizer measures for the three
applications of Section~\ref{sec:marked-partner-period}. This
identifies the determinant character in each Hecke product and fixes
the finite scalar in Section~\ref{sec:gluing-three-exchanges}.}

We use the covector, column and theta character conventions of
the partner comparison. In the three class-sum
applications the matching character is \(\chi^{-1}\) if
\(\chi\) is the determinant character appearing in the
upper datum. Let \(R\) be the original residue ring, with
its full exponent at each selected place, and write
\(X=(A,B)\) and \(Y=(C;D)\) for the matched blocks of
\(k^{-1}\) and \(k\), in that order.

Here \(B,C\) are square. In the order of
Section~\ref{sec:marked-square-finite-extension}, \(X\) consists
of the last covectors of \(k^{-1}\) and \(Y\) of the first columns
of \(k\). They give a split surjection and injection, respectively.
In this order the identities
\[
                 AC+BD=0,\qquad \det C=\det k\,\det B
\]
hold over \(R\), by block multiplication and complementary
minors. For any \(\Psi'\) invariant under the right action
of the integral column basis group, put
\[
 \Phi(A,B)=\chi(\det B)\Psi'(B^{-1}A),
 \qquad
 \Psi(C,D)=\chi(\det C)\Psi'(D),
\]
with value zero when the determinant of the corresponding
square block is a nonunit of \(R\).

When either block has unit
determinant, both do, and \(D=-B^{-1}AC\).
As \(-C\) is invertible, \(\Psi'(D)=\Psi'(B^{-1}A)\).
Thus \(\Phi(X)=\chi^{-1}(\det k)\Psi(Y)\).
Conversely \(B\) invertible makes the upper array split
surjective, and \(C\) invertible makes the lower array split
injective. This proves support and matching over the actual
residue ring, including the chosen prime-square coordinates.

For a definite target \(G=U(V_0)\), give \(G_\infty\) and the
full integral finite compact \(K_{\rm std}\) mass one. If \(F\)
is right invariant under \(G_\infty J\) for a compact open
\(J\subseteq K_{\rm std}\), disintegration over the finite classes gives
\[
 \int_{[G]} F(h)\,dh
  =\frac{1}{[K_{\rm std}:J]}\sum_{[h]_J}
      \frac{F(h)}{|G(F_0)\cap hJh^{-1}|}.
\]
Here each finite representative \(h\) is taken in \(G(\mathbb A_{F_0,f})\).
At \(J=K_{\rm std}\) the index is one. For the stacked determinant,
the real target action contributes
\(\det(u_\tau)^{\ell-c_\tau}\det(u_\tau)^{-\ell}
 \det(u_\tau)^{c_\tau}=1\).
For \(u\in K_{\rm std}\), write \(u_S\) for its component on
the chosen split side. The array above satisfies
\[
 \Phi(u_S^{-1}A,u_S^{-1}B)=\chi(\det u_S)^{-1}\Phi(A,B).
\]
The finite target character contributes \(\chi(\det u_S)\).
The unselected compact factors preserve the vacuum, so the
full integral compact fixes the twisted test.
Hence these class weights apply to the actual kernels.

For rank zero the quotient is a single point with one
automorphism, so its mass is one. The kernel on the zero
Hom space with characters \((\alpha/\beta,\alpha)\) is
\((\alpha/\beta)(\det g)\) on the retained point.
In the first exchange \(\beta=\alpha\), so it equals one.
Write \(\chi_\alpha\) for the residue character of the source
character \(\alpha\). The matched data are \(\chi_\alpha\circ\det\)
and \(1\), since \(X=k^{-1}\) and \(Y\) has zero columns:
\[
           \chi_\alpha(\det X)=\chi_\alpha(\det k)^{-1}.
\]
In Section~\ref{sec:v2-first-exchange-phase} the source character
is \(\alpha=c^{-1}\), so \(\chi_\alpha=\chi^{-1}\) when \(\chi\)
denotes the residue character of \(c\).
The empty small polynomial is one. The normalization of this
application is therefore a unit-mass constant, with no
unrecorded Haar factor.

For the second and third comparisons, the graph matching
above splits the small kernel into its two blocks.
The standard determinant factor is the kernel in the proper
definite class variable, and the remaining factor is the
polynomial kernel in the retained block. Before any trace averaging,
each column datum on the small side has the factor
\(\chi(\det C)\), so it has full column rank on its support
even with one trace character in place of the invariant sum.
The tensor product of the two real kernel vectors and this
column datum is the original small theta kernel on the
direct sum. The column estimate
\eqref{eq:partner-exchange-ordinary-bound} therefore bounds
its absolute integral for each fixed pair of retained points.
In particular this bound applies term by term in the finite
class and trace sums. The preceding measured identity can be
inserted pointwise and integrated against the remaining factor,
retaining each original class weight. This proves the
multiplication rule for the exchange constants.

In the third comparison \(\Psi'\) is the sum of translates
of the trace character by the group over the untwisted
base ring \(R_0=\mathcal O_E/(\pi_M^2)\), inflated to \(R\).
It is invariant under the right basis group over \(R\)
because its reduction permutes those terms.
Each element of \(\mathrm{GL}_{3m}(R_0)\) lifts to the
original standard level at the chosen split places by
lifting entries and testing determinants.

Extend the lift
by its conjugate unitary coordinate on the other split half
and by identity elsewhere. The splitting character of
the retained small kernel has trivial finite tame part,
and thus has value one on the determinant of this
integral move. The full-level auxiliary lift \(F^\circ\) is unchanged by it.
At a fixed retained point \(g'\in U(W')(\mathbb A)\), let \(K_D(s,g')\)
be the trace kernel for one residue translate \(D\in\mathrm{GL}_{3m}(R_0)\).
The fixed-Schwartz-datum estimate preceding
Section~\ref{sec:four-unitary-constant-terms} gives
\( |K_D(s,g')|\le C_{D,g'}\operatorname{ht}_{S}(s)^A\).
The first two exchanges give \eqref{eq:carrier-definite-class-expansion};
the invertible determinant block there gives
\(\sup_s\operatorname{ht}_S(s)^A|F^\circ(s)|<\infty\)
by Section~\ref{sec:marked-period-carrier-peak}.
Thus \(F^\circ(s)K_D(s,g')\) is integrable on the finite-volume
quotient \([U(S)]\). Right translation in this quotient
makes the integral of each translated trace character equal
to the unshifted integral. The finite sum multiplies by
\(\mathcal N=|\mathrm{GL}_{3m}(R_0)|\).
This is the untwisted group order, even if the modulus \(R\)
has acquired new twist places.

With \(e=m/2\), \(z_*=7e+k/2\) and
\(n_j=j+(t-11m-1)/2\), the three defects and period normalization arguments
are
\[
 \begin{array}{c|c|c}
 \text{kept, large, small}&\text{defect}&
       \text{character in the Hecke product}\\ \hline
 I_m,\ I_m,\ 0&e&c\upsilon^e\\
 I_m\perp S,\ I_{3m},\ I_m^-&m&c^{-1}\upsilon^m\\
 I_{3m}\perp W',\ V'',\ S&z_*&(ca_j)^{-1}\upsilon^{3m+z_*}.
 \end{array}
\]
Here \(m,k\) are even, \(W_j\) has rank \(r=k+m-1\),
and \(W'=W_j\perp I_j^-\perp I_{3m-j}\perp I_{6m+1}\)
has rank \(k+10m=3m+2z_*\).
The negative rank of each kept space equals that of its
small member, respectively \(0,m,m\).
For the third character,
\(ca_j=\mu_j^0\vartheta\,
\upsilon^{m+j+\lfloor(t+k-1)/2\rfloor}\);
the exponent left after inversion is
\[
    3m+z_*-m-j-\lfloor(t+k-1)/2\rfloor
          =-j-(t-11m-1)/2=-n_j .
\]
The first conversion to the negative copy of \(I_m\)
takes the conjugate of its scalar and multiplies by the
unit phase \(\lambda_0\chi(\eta)^{-m}\), with \(\eta\)
the norm-minus-one residue unit. Thus the squared
moduli multiply and the factor from trace averaging is
\(\mathcal N^{-2}\). In particular the coefficient of the
third left class sum taking it to the original auxiliary lift
period is \(\gamma_3\gamma_{12}/\mathcal N\), with
square modulus equal to a fixed constant divided by the
product of the three Hecke factors in the table.

\subsection{The rank-zero comparison and the opposite Hermitian form}
\label{sec:v2-first-exchange-phase}

\sectionopening{In this subsection, the first application of Section~\ref{sec:marked-partner-period} has a
rank-zero smaller space. We compute the source-negation factor with
the original residue character before applying the comparison to the
other two pairs.}

Use the modulus \(\mathfrak m=\pi_M^2\mathfrak q_\vartheta\),
the character \(c=\chi_M^j\vartheta\) of weight zero, and its
residue character \(\chi\). Let \(R=O_E/\mathfrak m\).
Write \(m=2e>0\) and \(P=1\). Apply the measured exchange in
Section~\ref{sec:marked-partner-period} with kept space \(I_m\), positive
space \(I_m\), and smaller space of rank zero. Take the
characters on the positive kernel to be
\(\alpha=\beta=c^{-1}\), and its finite function to be
\(\Psi_+(X)=\chi^{-1}(\det X)\). The two finite functions
are matched.

Indeed the block of an inverse matrix
\(k^{-1}\in\mathrm{GL}_m(R)\) is the entire matrix, and
\(\Psi_+(k^{-1})=\chi(\det k)\), the matching factor
for the residue character \(\chi^{-1}\) of \(\alpha\).
The other block is the unique empty column, with value
one. On singular matrices the extended residue character
is zero. In Section~\ref{sec:marked-stacked-determinant-transfer}
take \(p=\ell=0\), \(a=m\), and \(c_\tau=0\) at every infinite place.
Then \(P(\eta)=Q(\eta)=1\). The definite target average sends
the real vector supplied by the exchange to \(\kappa\) times the
Gaussian with polynomial part \(1\), where \(\kappa\) is independent
of the finite data. On the zero Hom space the only lattice vector is
zero and the vacuum has value one. The surviving character is
\((\alpha/\beta)(\det g)=1\), and the zero group has quotient mass
one. The projected identity is therefore
\(\kappa L_g(1)=c_\nu\), with \(\nu=c\upsilon^e\).
Since \(c_\nu\ne0\), also \(\kappa\ne0\). Set
\(\gamma_1=\kappa/c_\nu\). The measured class sum gives
\[
  \gamma_1\sum_{[h]\in[U(I_m)]/K_{\rm std}}
      \frac{\Theta^{I_m,I_m}_{c^{-1},c^{-1}}
                          [\Psi_+,1](g,h)}
           {|\operatorname{Aut}(h)|}=1,\qquad
  |\gamma_1|^2=\frac{C_1^{\rm ex}}{\mathcal L(c\upsilon^e,e)} .
\]
The first equality is the measured exchange, with its
class weights unchanged. The second is its scalar formula,
since \(\alpha^{-1}\upsilon^e=c\upsilon^e\). Here \(C_1^{\rm ex}\)
depends only on the fixed real and unramified data. In particular
the scalar used in the first equality is nonzero. The
point of evaluation \(g\) is arbitrary on the \(I_m\).

For the source-negation step, take the fixed finite
unit \(u\) with \(u\bar u=-1\) and put \(\eta=u\bmod
\mathfrak m\). At split pairs one can choose \(u\) with
coordinates \(1,-1\); at an unramified quadratic place,
surjectivity of the norm on units gives a norm-minus-one unit
as in Section~\ref{sec:compact-floor-volume}. At a place over
\(2\), the two completions are split: \(-d\equiv1\bmod 8\)
makes \(K\otimes\mathbb Q_2=\mathbb Q_2\times\mathbb Q_2\),
and conjugation interchanges the places of \(E\) above its
two factors. The split choice \((1,-1)\) has norm \(-1\)
already in the local integers. Multiplication by the fixed unit
preserves the integral boxes and their measures also at these places.
Use the original such unit and the original point of the
negative space. In chosen coordinates its point matrix
is \(X_{I_m^-}=(1,u)X_{I_m}\); its rational form is the
negative of the original form. On the Hom space, with
the target unchanged, the new point matrix is therefore
\((1,u^{-1})X_{\mathrm{Hom}}\). Indeed \(x\) is a rational Hom vector precisely when
\(X_TxX_S^{-1}\) is rational in the chosen bases. Replacing \(X_S\)
by \((1,u)X_S\) therefore identifies the new rational summation set
with the old one by \(\phi(x)=(1,u^{-1})x\).

On the coefficient realization of the oscillator data,
put \(Jf(x)=\overline{f(\phi x)}\). This transports the
Heisenberg action to the opposite source, because
\(\beta_{\rm old}(\phi x,\phi y)
=\overline{\beta_{\rm new}(x,y)}\).
At infinity the sign matrix in the bicharacter changes
sign and \(\phi\) is the identity. At finite places the
norm of \(u^{-1}\) is \(-1\). The Gaussian and standard
integral box are preserved, so this is a vacuum-preserving
transport with the same Haar normalization. The standard
holomorphic coordinate of each entry has been conjugated,
so \(P\) becomes its coefficientwise conjugate. For a
residue datum \(\Psi\), reduction of \(\phi\) multiplies
every entry by \(\eta^{-1}\) at the selected places. The
integral support is preserved since \(u\) is a unit, and
the transported residue function is
\[
       \Psi^J(X)=\overline{\Psi(\eta^{-1}X)} .
\]

Both functionals are invariant under rational translations and
modulations in corresponding Schr\"odinger polarizations.
The uniqueness argument of Section~\ref{sec:overlap-fourier-inversion}
applies to each: localizing at a rational point in a finite box, use
only modulations whose finite component is constant on that box.
Their real linear phases span the dual real space, so the localized
functional is evaluation, and rational translations identify its
multipliers at all rational points. Thus the two functionals are
proportional under the transport. Here \(J\) matches the rational Heisenberg
operators and their vacuum translates. The two normalized
theta functionals have the same squared mean, one, on the
vacuum translates over their additive adelic quotients.
The induced quotient isomorphism preserves probability
Haar measure. Consequently for a fixed \(|\lambda_0|=1\),
\[
       \overline{\theta_{\rm old}(f)}
            =\lambda_0\theta_{\rm new}(Jf).
\]
The scalar depends on the two spaces and their
theta normalizations alone.

To track the splitting characters, write
\(A(x)=hxg^{-1}\). Its matrix stays the same on passing
to the opposite source. At a real place the multiplier is
the conjugate of the negative principal minor of \(A\).
Jacobi's complementary-minor identity and
\(A^{-1}=J_{\rm sign}A^*J_{\rm sign}\) give
\(\overline{\Delta_{\infty,\rm old}(A)}
=\det A\,\Delta_{\infty,\rm new}(A)\).
At the finite places let \(L\) be the original standard
lattice and choose \(b\) with \(bL=L\cap A^{-1}L\).

The finite multiplier is
\(\chi_f(\det b)\), where
\(\chi_f=(\upsilon_E)_f|\cdot|_f^{-1/2}\) on
\(\mathbb A_{E,f}^\times\).
The matrix \(Ab\) represents \(L\cap AL\), so
\[
 \Delta_f(A^{-1})=\chi_f(\det A)\Delta_f(A).
\]
Since \(A^{-1}=A^*\), the determinant ideal defining
\(\Delta_f(A^{-1})\) is the conjugate of the one defining
\(\Delta_f(A)\), by local Smith form.
The character \(\chi_f\) is unramified on units and
takes complex-conjugate values on conjugate ideles.
Thus \(\Delta_f(A^{-1})=\overline{\Delta_f(A)}\).
Since \(\det A\) has idele norm one, the total multiplier
under the conjugating transport acquires
\(\upsilon(\det A)\).

In the oscillator integral this is
the only scalar change: the coordinate substitution
intertwines the matrices, conjugates the Heisenberg
bicharacter, and preserves the Gaussian and Haar measure.
As \(\det A=(\det h)^m(\det g)^{-m}\), the original
pair character \(a(\det h)b(\det g)^{-1}\) becomes the
pair character with entries \(a^{-1}\upsilon^m\) and
\(b^{-1}\upsilon^m\).

For \(a=b=c^{-1}\) and \(\Psi_+=\chi^{-1}\circ\det\)
this gives exactly
\[
 \Psi_+^J(X)=\chi(\eta)^{-m}\chi(\det X).
\]
The same formula holds at a singular matrix since both
sides vanish. The level and class stabilizers do not
change, and the class weights are rational. Conjugating
the first measured identity gives
\[
 \overline{\gamma_1}\lambda_0\chi(\eta)^{-m}
 \sum_{[h]\in[U(I_m)]/K_{\rm std}}
  \frac{\Theta^{I_m^-,I_m}_{c\upsilon^m,c\upsilon^m}
                       [\chi\circ\det,1](g^-,h)}
       {|\operatorname{Aut}(h)|}=1.
\]
In particular the same \(\mathcal L(c\upsilon^e,e)\)
controls the square modulus. The factor
\(\lambda_0\) is one fixed complex phase and
\(\chi(\eta)^{-m}\) is the residue phase at the actual
modulus. This proves the normalization used when
inserting the first class sum under the second exchange.

\subsection{Three period identities for the gluing operator}
\label{sec:gluing-three-exchanges}

\sectionopening{In this subsection, we apply the normalized partner comparison to the three pairs defining
the operator of Section~\ref{sec:original-comparison-tests}. The real
factors of Section~\ref{sec:marked-stacked-determinant-transfer} and
the finite factors of
Section~\ref{sec:marked-gluing-finite-normalization} identify the
resulting common period.}

These identities use the partner comparison on the three stated
pairs. In particular the scalar formula is used with its original
vacuum normalization at the top and after all plane removals.
Use the same \(j\), original \(W_j,V'',V^\circ\), and
choice set \(\mathcal I\) as in
Section~\ref{sec:original-comparison-tests}. Put
\(I_a=(E^a,1)\), \(S=I_{2m}\oplus I_m^-\) and
\[
 W'=W_j\oplus I_j^-\oplus I_{3m-j}\oplus I_{6m+1}.
\]
The space underlying \(V''\) is the definite space
with these indices. Set \(e=m/2\) and \(z_*=7e+k/2\), an integer
since \(r\) is odd and \(m\) is even.
Then \(\dim V''=3m+2z_*\). In these comparisons put
\(\upsilon=\upsilon_E|_{T_1(\mathbb A_{F_0})}\), where \(\upsilon_E\)
is the splitting Hecke character of Section~\ref{sec:splitting-level-characters}.
Let \(c=\chi_M^j\vartheta\) be the weight-zero character and \(a_j\)
the twist-free character on the larger pair. Write
\(\mathfrak m=\pi_M^2\mathfrak q_\vartheta\) for the modulus
in the twisted intermediate kernels. The finite character of
\(c\) is the one used in their determinant data and is nontrivial
at every place of the prescribed support of \(\mathfrak m\).

Section~\ref{sec:v2-first-exchange-phase} gives the
measured first exchange on the two copies of \(I_m\).
Its weighted class sum is normalized to one by
\(\gamma_1\). Its squared modulus is a constant divided by
\(\mathcal L(c\upsilon^e,e)\), where
\(\mathcal L(\nu,z)=M_\nu(z)\) is the normalizing factor of
Proposition~\ref{prop:normalized-period-exchange} on the split
support of \(\mathfrak m\). Conjugating the first group to \(I_m^-\) gives the
scalar \(\overline{\gamma_1}\lambda_0\chi(\eta)^{-m}\)
of that section, with the prescribed reduction
and \(|\lambda_0|=1\).

The next exchange uses
\[
          (I_m\oplus S, I_{3m}; I_m^-)
\]
and the source character \(c\upsilon^m\).
Let \(\gamma_2\) be the coefficient multiplying the large class sum
in this second exchange. Its normalization gives
\[
 |\gamma_2|^2=C_2^{\rm ex}/\mathcal L(c^{-1}\upsilon^m,m).
\]
Insert the first constant class sum under this second exchange
and multiply the kernels for the direct sum \(I_m\oplus S\).
The remaining kernel on \((I_m^-,S)\), integrated against the
constant function on \([U(I_m^-)]\), is the auxiliary lift \(F^\circ\).

Its datum is the original vacuum and the polynomial on the
negative block is the \(2m\)-th power of the determinant minor.
For fixed finite and real tests, the convention of
Section~\ref{sec:marked-square-analytic-pairing} gives
\[
 \Theta^{S,I_{3m}}_{\upsilon^m,1}[\chi\circ\det,1](x,h)
  =c(\det x)\Theta^{S,I_{3m}}_{\upsilon^m,c}
                                  [\chi\circ\det,1](x,h).
\]
The small kernel on \((I_m^-,S)\) has the same factor
\(c(\det x)\) from its first character \(c\upsilon^m\).
Cancel this factor and substitute the conjugated first class sum to obtain
\[
 \gamma_{12}=\gamma_2\overline{\gamma_1}\,
                  \lambda_0\chi(\eta)^{-m}.
\]
Here \(\lambda_0\) is the fixed phase of the opposite-form
identification, while \(\chi(\eta)\) is evaluated at the actual
modulus and character. Both have absolute value one. Hence
\[
 |\gamma_{12}|^2
   =\frac{C_{12}^{\rm ex}}{\mathcal L(c^{-1}\upsilon^m,m)
                  \mathcal L(c\upsilon^e,e)} .
\]
More explicitly, write \(\mathcal C_a\) for the full-level finite
classes of \(U(I_a)\), and \(e_a(h)=|\operatorname{Aut}(h)|\).
Give their representatives real component one. Put
\[
 A_c(h)=\sum_{g\in\mathcal C_m}
   \frac{\Theta^{I_m,I_{3m}}_{c,1}[1,1](g,h)}{e_m(g)},
\]
where the first kernel is \(K_0\) of \eqref{eq:display-i5}, with
unit modulus and the integral vacuum. The determinant character
therefore gives \(A_c(h)=c(\det h)A_1(h)\).
The upper residue test on \(I_m\oplus S\) with \(\Psi'=1\)
is \(\chi(\det B)\), where \(B\) is the square block on \(S\).
The polynomial is one. Factoring this tensor and canceling
\(c(\det s)\) gives the pointwise identity
\begin{equation}\label{eq:carrier-definite-class-expansion}
 F^\circ(s)=\gamma_{12}\sum_{h\in\mathcal C_{3m}}
     \frac{A_c(h)}{e_{3m}(h)}
     \Theta^{S,I_{3m}}_{\upsilon^m,c}
                         [\chi\circ\det,1](s,h).
\end{equation}
For \(k\) in the standard finite compact of \(U(I_{3m})\),
the first kernel gives \(A_c(hk)=c(\det k)A_c(h)\).
In the second kernel the residue substitution
\(B\mapsto\bar k^{-1}B\) contributes \(c(\det k)^{-1}\),
since \(\upsilon^m\) is unramified there. Hence the product
has the asserted full level in the class variable.

For \(i=(\Lambda,h^\flat,\sigma)\), where \(\sigma\) is the chosen tuple of short covectors, use the third
exchange on
\[
           (I_{3m}\oplus W',V'';S).
\]
Put \(e_j^{(1)}=j+\lfloor(t+k-1)/2\rfloor\). The first-block
character and the source character are
\[
 b''_{j,\vartheta}=\mu_j^0\vartheta\,\upsilon^{e_j^{(1)}},
 \qquad ca_j=b''_{j,\vartheta}\upsilon^m .
\]
Thus the third kernel has characters \((ca_j,b''_{j,\vartheta})\).
The finite matching pairs the trace character of \(\Lambda\)
with the unnormalized invariant sum of the trace character on the
ratio of the two square blocks. The left datum is
\(\chi(\det x'')\psi_\Lambda^{\rm inv}((x'')^{-1}x')\).

\subsubsection{Matching the real factors}
\label{sec:gluing-real-matched-factors}

\sectionopening{In this subsection, we form the common real polynomial for the period comparisons of Section~\ref{sec:gluing-three-exchanges}.}

For the third comparison the small space \(S\) has \(2m\) positive
columns and \(m\) negative columns. Write \(r_+=k+j-1\)
and \(d_+=\dim V''\).
At \(\tau\) take the following ordered \(2m\)-row matrices
\(\eta_{\tau i}\) in Section~\ref{sec:marked-stacked-determinant-transfer},
extended by zero on the initial \(I_{3m}\) block. For every positive
column of height \(h<r_+\) in the diagram, number the rows by
\(0\le a<2m\). Row \(a\) is covector \(h+a\) of \(W_j\) when
\(h+a<r_+\), and covector \(a\) of \(I_{3m-j}\) otherwise.
The latter index exists since \(3m-j\ge2m\).
For each of the \(\Delta'=(\Delta-t-1)/2\) columns of full positive
height, use the ordered \(2m\) short rows of \(W^\flat\) prescribed
by \(\sigma\) at that place and column. Thus, with the counts of
Section~\ref{sec:peak-reference-coefficient},
\[
 c_{\rm stack,\tau}
    =\sum_{h<r_+}\mu^+_\tau(h)+\Delta'
    =w_j(\tau)-e_j^{(1)} .
\]
Order the negative covectors by \(0\le a<m\), taking row \(a\) of
\(I_j^-\) for \(a<j\) and row \(k+a-1\) of \(W_j\) for \(a\ge j\).
They exhaust the negative block, including when \(j=m\).
Its determinant exponent in the small polynomial is
\[
 c_{\rm stack,\tau}+\dim(I_{3m}\perp W')-2m
   =c_{\rm stack,\tau}+d_++m .
\]
In these orders the small polynomial is exactly \(P_jQ_\sigma\) of
Section~\ref{sec:peak-reference-coefficient}, independent of the
variables of the initial \(I_{3m}\). A fixed alternative order of
the negative rows contributes its permutation sign raised to this
exponent at each place, a scalar independent of the finite data.
The paired large polynomial stacks the same \(\eta_{\tau i}\) over
the positive variable rows in the indicated order. The third exchange scalar
\(\gamma_3\) has squared modulus \(C_3^{\rm ex}\) divided by
\(\mathcal L((ca_j)^{-1}\upsilon^{3m+z_*},z_*)\).
Summing the first block using the identity for the auxiliary lift and then
reducing the trace-character data modulo the untwisted
\(\pi_M^2\) gives
\[
  \gamma_3\gamma_{12}\,\Upsilon_i(g)=\mathcal N F_i(g),
 \qquad
 \mathcal N=|\mathrm{GL}_{3m}(\mathcal O_E/(\pi_M^2))|.
\]
Indeed the square block \(C\) of the matched column array belongs
to the \((S,I_{3m})\) kernel with characters \((\upsilon^m,c)\);
its finite factor is \(\chi(\det C)\) and its real polynomial is one.
The remaining \((S,W')\) kernel has characters \((\upsilon^m,a_j)\)
and the polynomial just computed. Substituting
\(\mu_j^0=\chi_0(\upsilon^{-\Delta}\chi_M)^j\) and
\(c=\chi_M^j\vartheta\) gives
\(a_j=\chi_0\upsilon^{m+e_j^{(1)}-j\Delta}\),
independent of \(\vartheta\). The two source determinant factors
multiply to \((ca_j)(\det s)^{-1}\).
Its trace datum is inflated from \(R_0=\mathcal O_E/(\pi_M^2)\);
at additional twist places the remaining factor is the integral vacuum.
For one translate indexed by \(D\in\mathrm{GL}_{3m}(R_0)\), write
\(K_D(s,g')\) for this remaining trace kernel at a fixed point
\(g'\) of \(U(W')\), and write \(T_h(s)\) for the
\((S,I_{3m})\) kernel in \eqref{eq:carrier-definite-class-expansion}.
That finite class expansion gives
\[
 \int_{[U(S)]}|F^\circ(s)K_D(s,g')|\,ds
 \le |\gamma_{12}|\sum_{h\in\mathcal C_{3m}}
    \frac{|A_c(h)|}{e_{3m}(h)}
    \int_{[U(S)]}|T_h(s)K_D(s,g')|\,ds .
\]
The product in each inner integral is the small theta kernel on
\((S,I_{3m}\perp W')\), up to the fixed unit phase of the
product normalization. Its column datum contains
\(\chi(\det C)\) on the square block at every selected place,
so all its supported columns have full rank over the local field.
Equation~\eqref{eq:partner-exchange-ordinary-bound} in the
conjugate transpose model, with height exponent zero, bounds
the inner integral on the finite-volume quotient. Thus right
translation in this integral against the full-level auxiliary
lift gives the same value for every \(D\).

Their unnormalized sum
has \(\mathcal N\) terms. The integral defining \(F_i\) uses
\([U(S)]\) with the original quotient measure.
Equivalently
\[
 c'_\vartheta=\gamma_3\gamma_{12}/\mathcal N,\qquad
 |c'_\vartheta|^2
  =\frac{C_3^{\rm ex} C_{12}^{\rm ex}/\mathcal N^2}
          {H_j(\vartheta)} .
\]
Here the product in the denominator is
\[
 H_j(\vartheta)=\mathcal L(c\upsilon^e,e)
  \,\mathcal L(c^{-1}\upsilon^m,m)
  \,\mathcal L((ca_j)^{-1}\upsilon^{3m+z_*},z_*).
\]
To match these factors with the signed critical values of
Section~\ref{sec:common-critical-twist}, write
\(n_j=j+(t-11m-1)/2\), \(\mu_j^0\) for the untwisted lift
character, and
\[
 \nu_+=\upsilon^{-m}\chi_M^j,\qquad
 \nu_-=\upsilon^{-e}\chi_M^{-j},\qquad
 \nu_*=\mu_j^0\upsilon^{n_j}.
\]
Then \(c^{-1}\upsilon^m=(\nu_+\vartheta)^{-1}\),
\(c\upsilon^e=(\nu_-\vartheta^{-1})^{-1}\), and
\((ca_j)^{-1}\upsilon^{3m+z_*}=(\nu_*\vartheta)^{-1}\).

\indent{\it The first-factor pairing.}
Write \(R=\mathcal O_E/\mathfrak m\),
\(R_0=\mathcal O_E/(\pi_M^2)\), and let
\(\phi:R\longrightarrow R_0\) be reduction. Let \(T\) be the
coordinate set of \(V''\), \(I\subset T\) the coordinates of \(W_j\),
and \(J_I:R^I\to R^T\) their coordinate inclusion. Write
\(w\) for the covectors \(\Lambda\) stacked over \(\pi_M1_T\) in
the residue coordinates of \(\Xi_{\vartheta,i}\). The top
covector projection \(u\) satisfies \(uw=\Lambda\). Define
\[
 H=\{A\in\mathrm{GL}_T(R):Ae_b=e_b\ (b\in T\setminus I),\
                           w\phi(A)=w\}.
\]
Here \(\phi\) is applied entry by entry. The inverse image of \(H\)
under reduction from the full integral subgroup \(K\) on \(V''\)
is the prescribed level \(J_\Lambda\). Reduction onto \(\mathrm{GL}_T(R)\)
is surjective, since each place in the modulus is split in the
chosen coordinates and reduction of the integral general linear
group onto each local residue ring is surjective.

Let \(Y,Z\) be the two residue blocks of the third kernel, with
\(Z\) square on \(T\) and the columns of \(Y\) indexed by \(I_{3m}\).
The invariant trace function is the unnormalized sum
\[
 B(Y)=\sum_{D\in\mathrm{GL}_{3m}(R_0)}
          \psi\big(\operatorname{tr}(\Lambda\phi(Y)D)\big).
\]
For \(A\in H\), \(\Lambda\phi(A)=uw\phi(A)=\Lambda\), so
\(B(AY)=B(Y)\). The ratio datum \(\Phi(Y,Z)\) is
\(\chi(\det Z)B(Z^{-1}Y)\) for invertible \(Z\), and vanishes
otherwise. Put \(\Phi_H=1_{Z\in H}\Phi\). For \(\sigma=ZJ_I\) set
\[
 C_w(\sigma)=
   \begin{cases}
       \chi(\det(J_I^{\rm t}\sigma)),&
           w\phi(\sigma)=wJ_I,\\
       0,&\text{otherwise}.
   \end{cases}
\]
The character of the residue ring is extended by zero on
nonunits. When the columns of \(Z\) off \(I\) are those of the
identity, \(Z\) is block triangular in \(I\sqcup I^c\) and
\(\det Z=\det(J_I^{\rm t} ZJ_I)\). Also
\(w\phi(Z)=w\) is equivalent to \(w\phi(ZJ_I)=wJ_I\).
On \(H\), \(B(Z^{-1}Y)=B(Y)\). Thus
\[
 \Phi_H(Y,Z)=C_w(ZJ_I)B(Y)1_{Z e_b=e_b\ (b\in I^c)}.       \tag{*}\label{eq:gluing-determinant-test-polynomial}
\]

For \(A\in\mathrm{GL}_T(R)\), define
\[
 (T_A\Psi)(Y,Z)=\chi(\det A)\Psi(A^{-1}Y,A^{-1}Z).
\]
By the finite Weil action, right translation of the kernel
by \(k\in K\) with residue \(A\) replaces \(\Psi\) by \(T_A\Psi\).
The determinant factors give \(T_A\Phi=\Phi\). Hence
\(T_A\Phi_H=1_{Z\in AH}\Phi\) and
\[
 \Phi=\sum_{A H\in\mathrm{GL}_T(R)/H}T_A\Phi_H .
\]
For a full-level class representative \(h\), put
\(\Gamma_h=U(V'')(F_0)\cap hKh^{-1}\). The \(J_\Lambda\)-classes
above \(h\) are the orbits of \(\Gamma_h\) on \(K/J_\Lambda\)
under the action through \(h^{-1}\Gamma_hh\).

For the restricted kernel \(f(hk)\), rational and right
\(J_\Lambda\)-invariance give
\[
 \frac{1}{|\Gamma_h|}\sum_{kJ_\Lambda\in K/J_\Lambda} f(hk)
 =
 \sum_{[k]\in(h^{-1}\Gamma_hh)\backslash K/J_\Lambda}
 \frac{f(hk)}
 {|\Gamma_h\cap hkJ_\Lambda k^{-1}h^{-1}|}.
\]
At a real place \(u\in U(V'')\) left multiplies the bottom rows
of each stacked determinant.
The polynomial contributes \(\det(u)^{c_\tau}\);
the vacuum and splitting give \(\det(u)^{-c_\tau}\).
Their product is one in the full kernel. This proves equality of the full-level class sum of
\(\Phi\) and the \(J_\Lambda\)-class sum of \(\Phi_H\), with
the original class weights. More explicitly, if the Haar measure has
\(\operatorname{vol}(K)=1\), the integral of the single
\(J_\Lambda\)-invariant kernel \(f\) is its weighted
\(J_\Lambda\)-class sum divided by \([K:J_\Lambda]\).
The full kernel is the sum of all \([K:J_\Lambda]\) translates
of this restricted kernel. Integrating that sum restores its
weighted class sum exactly, also when the stabilizer orders
vary between the full-level classes.

Use the ordered block decomposition of the kept space
\[
 I_{3m}\oplus W'= W_j\oplus I_j^-\oplus W^{\mathrm{add}},
\]
where \(W^{\mathrm{add}}\) includes \(I_{3m}\). The equality
is by permutation of the coordinate parameters. In \eqref{eq:gluing-determinant-test-polynomial}
the columns on \(W_j\) carry \(C_w\), those on \(I_j^-\)
carry the inclusion indicator there, and the remaining columns
carry the added-block datum of
Section~\ref{sec:marked-first-factor-bounds}. Indeed the complementary
identity columns split into the negative pad and the fixed positive
columns on \(W^{\mathrm{add}}\), while \(B\) depends only on the
\(I_{3m}\) columns. For the real polynomial, first expand
along the selected coordinate rows. By Section~\ref{sec:bottom-descent-real-type},
\(\mu_\tau^+(0)=0\) and \(\mu_\tau^+(h)=0\) for \(1\le h<k-1\),
since \(\varsigma_{h-1}-\varsigma_h=1\) in that range.
Hence at every contributing height \(h<r_+\),
\[
  1\le r_+-h\le j\le m .
\]
The \(2m\) selected rows of that factor therefore remove exactly the
columns of \(W_j\) indexed by \(h,\ldots,r_+-1\).
The remaining selected rows remove columns of the added positive space.
Expanding in the first \(h\) columns of \(W_j\) gives its height-\(h\)
minors. In each of the \(\Delta'\) full-height factors the short
rows are supported on \(W^\flat\), so expansion in \(W_j\)
uses all \(r_+\) of its positive columns. Thus the Laplace expansion
has precisely the positive Weyl diagram of \(W_j\). On the same
decomposition the real polynomial is
\[
 P(x_{W_j},x_{\rm add})
   =\sum_\lambda L_\lambda(x_{W_j})\zeta_\lambda(x_{\rm add}).
\]
Here \(L_\lambda\) are the positive Weyl coordinate
polynomials, and \(\zeta_\lambda\) are the coordinates used in \eqref{eq:display-i5},
obtained by expanding the stacked determinants in the \(W_j\)
columns with coefficients in the added variables. The polynomial
is constant on the negative pad.

We use the tensor product of the fixed theta functionals on the
standard signed lines, retaining the fixed definite vacuum
normalization on the positive lines. For a Hermitian space \(B\) with point
matrix \(X_B\), let \(Y_B=X_B^{-1}E^{\dim B}\) be its rational vectors
in these coordinates. Landherr's theorem
\cite[pp.~246--247]{LandherrHermitian} applies as in the
plane-removal argument: the finite forms are isometric and
the signatures agree at every real place. If \(J_B,J_{\rm std}\)
are the rational Gram matrices, take a rational isometry
\(P^*J_B P=J_{\rm std}\) and write \(X_{\rm std}\) for the
standard signed point matrix. The point equations imply
\[
 g_B=X_B^{-1}P X_{\rm std},\qquad
 g_B(X_{\rm std}^{-1}E^{\dim B})=Y_B .
\]
Here \(g_B\) is an adelic isometry in the normalized coordinates.
Use the vacuum-normalized Weil action of
Lemma~\ref{lem:weil-vacuum-composition} before the determinant
characters of the kernels are inserted. With \(\theta_{\rm std}\)
the product functional just specified, put
\[
         \theta_B(f)=\theta_{\rm std}(U_{g_B^{-1}}f).
\]
Another rational isometry \(P'\) gives
\(g_B^{-1}g'_B=X_{\rm std}^{-1}(P^{-1}P')X_{\rm std}\), the action
of a rational unitary element of the standard space.
Section~\ref{sec:marked-rational-theta} proves that its action fixes
\(\theta_{\rm std}\). For a block sum take the direct sum of these
maps. The real determinant factors multiply, the finite intersections
are direct sums, and the vacuum is a tensor product, so the block
Weil action is the tensor product of the two actions. Thus
\[
 \theta_{B_1\perp B_2}(f_1\otimes f_2)
       =\theta_{B_1}(f_1)\theta_{B_2}(f_2).
\]
The same identity applies after the indicated permutation of coordinates.
To compare with any previously fixed normalized theta functional
\(\theta_B^{\rm orig}\) on the same model, rational Heisenberg
invariance gives \(\theta_B^{\rm orig}=\epsilon_B\theta_B\).
The symplectic map \(g_B\) carries the standard rational displacement
lattice onto the one for \(B\), and carries the probability Haar on
one quotient to that on the other. The mean-square identity of
Section~\ref{sec:overlap-fourier-inversion} on a unit vector then
gives \(|\epsilon_B|=1\). With previously fixed phases the product
identity acquires exactly \(\epsilon_{B_1\perp B_2}/
(\epsilon_{B_1}\epsilon_{B_2})\). This scalar depends only on the
fixed spaces and their normalizations, so it multiplies the same
comparison constant for every residue array, modulus and twist.
Its squared modulus is one, including in the scalar comparison
at the trivial twist.
In the indicated product convention, applying the identity twice
gives the product of the three kernels. In the Schrödinger model the integral
vacuum and the residue functions factor on the three Hom spaces.
In the real Fock model the vacuum and the displayed polynomial
factor in the same way. The rational summation set is a Cartesian
product, and the finite class sums and Schwartz theta series
converge absolutely.

In the order \(W_j,I_j^-,W^{\mathrm{add}}\), the three character
pairs are
\[
 (\upsilon^{m-j},b''_{j,\vartheta}),\qquad
 (\upsilon^j,b''_{j,\vartheta}),\qquad
 (b''_{j,\vartheta},b''_{j,\vartheta}).
\]
Their first characters multiply to \(ca_j\).
The three first points are \(g\),
the identity of \(I_j^-\), and \(\operatorname{diag}(1,h^\flat,h')\)
in the added block, with the same second point \(h\) in each.
Multiply the added kernel by \(A_c(h')\) and sum over \(h'\)
with the weights of \eqref{eq:carrier-definite-class-expansion}.
The coefficient of \(L_\lambda\) equals the negative-pad kernel
at \(h\) times this weighted sum with \(\zeta_\lambda\). By \eqref{eq:display-i5} these are the coordinates
of \(\Upsilon_{\vartheta,i}(h)\) before the scale and conjugation.
Thus \(\Upsilon_i(g)\) is the \(J_\Lambda\)-class sum of the
\(W_j,V''\) kernel with residue function \(C_w\) and polynomial
vector \(\Upsilon_{\vartheta,i}(h)\).

For the class sum, let \(k_+\in J_\Lambda\) be supported at a new
split member \(w\). Write its residue in \(I\sqcup I^c\) as
\(A=\left(\begin{smallmatrix}D_w&0\\ E_w&1\end{smallmatrix}\right)\).
Then \(J_I^{\mathsf t}A^{-1}=D_w^{-1}J_I^{\mathsf t}\),
\(\det A=\det D_w\), and
\[
 C_w(A^{-1}\sigma)=\chi_w(\det D_w)^{-1}C_w(\sigma).
\]
The \(W_j,V''\) kernel is multiplied by this scalar under
\(h\mapsto hk_+\), since its first splitting character on integral
units is trivial. The added and negative-pad inclusion conditions
are preserved, and the trace datum is unchanged. Their combined
target-leg factor on the residue determinant is \(\chi_w(\det D_w)\).
Thus the two factors of the class summand have inverse level
characters. The conjugated factor
\(\overline{c''\Upsilon_{\vartheta,i}(h)}\) transforms by
\(\chi_w(\det D_w)^{-1}\), the level character of the original
target lift under \(\operatorname{diag}(k_+,1)\). At an old
place above \(\pi_M\), \(k_+\equiv1\bmod\pi_M\) and the residue
characters are all trivial.

The conjugate of the first factor is
\(c''\Upsilon_{\vartheta,i}(h)\) by \eqref{eq:display-i5}. Linearity in the
polynomial, followed by the third exchange, gives
\[
 \Xi_{\vartheta,i}(g)=c''\Upsilon_i(g),\qquad
 c'_\vartheta\Xi_{\vartheta,i}(g)=c''F_i(g).
\]

\subsubsection{Ranks and levels in the three comparisons}

\sectionopening{In this subsection, the applications of Section~\ref{sec:marked-partner-period} require
simultaneous inequalities at every real place and the stated finite
levels on each pair. We verify these conditions for the characters
just used.}

Let \(c=\chi_M^j\vartheta\) and let \(\alpha\) denote
the actual source character in one class-sum comparison.
If the kept rank is \(N\), its number of negative
signs is \(\ell\), and the number of real stacked
determinants is \(c_\tau\), then
\[
 \operatorname{wt}_\tau(\alpha)=\ell-c_\tau,\qquad
 \nu=\alpha^{-1}\upsilon^{N/2},\qquad
 \mu_\tau=N/2-\ell+c_\tau .
\]
This is the same \(\mu_\tau\) used in the real
slice proof and in the normalization. It gives
\(N\le 2\mu_\tau+2\ell\). The comparison needs
\(N\ge2\ell+2\), and the small member must have exactly
\(\ell\) negative signs since the large target is definite.
In order, the small spaces are \(0\), \(I_m^-\), and
\(S=I_{2m}\oplus I_m^-\). The third small quotient is therefore
the one for the indefinite group \(U(2m,m)\), with the Haar
measure of Section~\ref{sec:marked-partner-period}.

For the first comparison the ranks \((N,a,b)\) are
\((m,m,0)\), \(\ell=0\), \(\alpha=c^{-1}\) and
\(c_\tau=0\). The inequality is \(m\ge2\).
For the second they are \((4m,3m,m)\), \(\ell=m\),
\(\alpha=c\upsilon^m\), and again \(c_\tau=0\).
For the third they are
\((6m+2z_*,3m+2z_*,3m)\), \(\ell=m\), and
\(\alpha=ca_j\). The \(m-j\) negative coordinates of
\(W_j\) and the \(j\) negative pad coordinates give
the kept negative rank. The real number of stacked
determinants is explicitly
\[
 c_\tau=w_j(\tau)-j-\lfloor(t+k-1)/2\rfloor ,
\]
the cardinality of the set of determinant
indices.

Since \(ca_j\) has weight
\(m+j+\lfloor(t+k-1)/2\rfloor-w_j(\tau)\)
this is exactly the preceding weight condition.
The parameters give \(m,k\) even, so \(z_*=7m/2+k/2\)
and all three ranks and defects are integers as written.
At each old distinguished place the residue character of \(\chi_M\)
is \(\omega_w\), of order \(q_2\). Thus \(\omega_w^{\pm j}\) has order
\(q_2/\gcd(q_2,j)>1\) for \(1\le j\le m<q_2\). At an additional twist place
the character is nontrivial by exact twist support.

The fixed norm-one splitting and fixed \(a_j\) have
trivial finite tame parts. Thus the same semilocal
matching proof applies to the full modulus
\(\pi_M^2\mathfrak q_\vartheta\) in all three cases.
At an old place the determinant character factors through
\(k_w^\times\), while the test may distinguish classes modulo
\(\varpi_w^2\). Proposition~\ref{prop:normalized-period-exchange}
applies to these locally constant matched tests by the integral
matching and prime-power inflation in
Section~\ref{sec:radical-local-line-terms}. The reversed
matrix Gauss sum in Section~\ref{sec:marked-square-finite-extension}
is taken on its squarefree residue ring.
The polynomial transfer uses a fixed real vector before
the modulus, and it applies to the same determinant
indices for every finite test.

\subsubsection{The common period as the full gluing adjoint}
\label{sec:v2-measured-gluing-adjoint}

\sectionopening{In this subsection, with the three comparisons established, we pair their common period
against arbitrary tuples in the target of
Section~\ref{sec:original-comparison-tests}. The identity identifies
the adjoint used in the injectivity argument of
Section~\ref{sec:common-period-injectivity}.}

Let \(y\) be in the original projected part of \(W_j\) and
let \(\zeta=(\zeta_{i,c})\) be a tuple in the original
second-factor Weyl spaces. Write \(K^\circ_{\zeta_i}(g)\)
for the sum over the original \(V^\circ\)-classes
\(c\) of the second-block kernel at \((g,(1,c))\), with
the polynomial in the negative coordinates associated to \(\overline{\zeta_{i,c}}\).
This kernel uses the untwisted \(W_j\)-character, the
indicated power \(\upsilon_{U(1)}^{m/2-j}\), modulus
\(\pi_M^2\), and the original inclusion indicator.
The class weights and real test of that block are the
ones in the gluing map.

Write \(\theta_\vartheta(g)\) for the vector-valued kernel in
the definite target with its original class pairing. In the
second-linear convention the lift is
\(u_{j,\vartheta}y=\int_{[U(W_j)]}y(g)\theta_\vartheta(g)\,d\mu(g)\)
and the normalized descent is
\(\mathcal D_\vartheta v(g)=\langle\theta_\vartheta(g),v\rangle\).
Consequently the integrand for
\(\langle\zeta,G_\vartheta^*u_{j,\vartheta}y\rangle\) is
\(y(g)\overline{\mathcal D_\vartheta(G_\vartheta\zeta)(g)}\).
For each \((i,c,h)\) in \(G_\vartheta\zeta\), the sum in the
point form cancels the reciprocal stabilizer weight at the glued
target point, by its defining evaluation identity. The first-block
class retains the reciprocal stabilizer weight at \(J_\Lambda\),
which is one here. The second-block level is free, so the
remaining sum over \(c\) is the plain sum for the arbitrary
components \(\zeta_{i,c}\). The product formula
for the original glued lattice expands this integrand
into the finite sum over \(i\) and the congruence group
\(\bar J^{\,1}=\{u\bmod\pi_M^2:u\equiv1\bmod\pi_M\}\)
of translates by the fixed lifts \(k(u)\) of
\[
    M^{l_1}c_{L_i}\,y(g)
                      \Xi_{\vartheta,i}(g) K^\circ_{\zeta_i}(g).
\]
Let \(\delta_i\) be the determinant of \(L_i\) on the old split
members and let \(\chi_{W,r}\) be the fixed split Weil
character of \(W_j\) at rank \(r\). The finite action on the
matrix variable of the glued kernel is
\[
 f(X)\longmapsto c_{L_i}f(L_i^{-1}X),\qquad
 c_{L_i}=\chi_{W,r}(\delta_i)|\delta_i|^{-r/2}.
\]
For a chosen \(w\mid\pi_M\), let \(s_{i,w}\) be the rank of
\(\Lambda\bmod w\). The row lattice of
\(\binom{\Lambda}{\pi_M I_{d_+}}\) has \(s_{i,w}\) elementary
divisors of valuation zero and \(d_+-s_{i,w}\) of valuation one.
Thus \(v_w(\det B_+)=s_{i,w}\) for \(B_+=\pi_M D^{-1}\).
The other diagonal block \(J_z+\pi_M(1-J_z)\) has \(d_--r\)
elementary divisors of valuation one. Hence, with \(v_w(\varpi_w)=1\),
\(v_w(\delta_{i,w})=s_{i,w}+d_--r\).
The target character of Section~\ref{sec:gluing-block-period-vanishing}
is \(a=\upsilon^{m-j}\) on both blocks. In the convention of
\eqref{eq:display-f2}, base change contributes
\(\upsilon_E^{2(m-j)}\), as in Section~\ref{sec:exceptional-jet-hecke}.
The rank-\(r\) matrix variable contributes the split Weil factor
\(\upsilon_E^{-r}\) from Lemma~\ref{lem:weil-vacuum-composition}. Thus
\(a_r=\chi_{W,r}=\upsilon_E^{2(m-j)-r}\) is unramified at the old
split places, and the local factor of \(c_{L_i}\) is
\[
 \bigl(a_{r,w}(\varpi_w)q_w^{r/2}\bigr)^{s_{i,w}+d_--r}.
\]
This expression depends only on the original residue choice and the
two block ranks.
Put \(C_i=M^{l_1}c_{L_i}\), the scalar of this summand in the
common period. Both the character and \(L_i\) are fixed. At the twist
places \(L_i\) and the lifts \(k(u)\) are identity.
At an old split place both \(L_i\) and \(A=\pi_M L_i^{-1}\) are
integral by Section~\ref{sec:marked-comparison-lattice}.
For \(k_+\in J_\Lambda\), let \(D_{\rm old}\) denote the row basis
matrix of \(\binom{\Lambda}{\pi_M I_{d_+}}\).
The two old level conditions give
\(D_{\rm old}(k_+-1)\equiv0\bmod\pi_M^2\). Therefore
\[
 L_i^{-1}\operatorname{diag}(k_+,1)L_i-1
   =\pi_M^{-1} A\operatorname{diag}(k_+-1,0)L_i
   \equiv 0\bmod\pi_M.
\]
This transports the old principal level to the glued point. The condition
\(L_i^{-1}X\equiv I_{\rm incl}\bmod\pi_M\) implies \(X\) is
integral. On integral \(X\) it is equivalent to
\[
 AX\equiv\pi_M I_{\rm incl}\pmod{\pi_M^2}.
\]
The same section splits this condition into the two residue blocks.
Let \(S_+\) denote the first-block matrix after the right
translation by \(k(u)\), with inclusion matrix \(J_I\).
The lower rows of \(w\) give
\(\pi_M S_+\equiv\pi_M J_I\pmod{\pi_M^2}\).
Hence \(J_I^{\mathsf t}S_+\equiv1\pmod{\pi_M}\) at the old
places and every old tame determinant multiplier is one there.
At each new place \(k(u)=1\), so the corresponding factor in
\(\chi(\det(J_I^{\mathsf t}S_+))\) is exactly the prescribed
new determinant weight on the original inclusion block. Thus
the first block supplies the entire finite determinant factor.
The second block retains its inclusion modulo \(\pi_M^2\)
and character \(b_j^\circ=\upsilon^{m/2-j}\). Since
\(e_j^{(1)}+m/2-j=(r+t)/2\), the product on the original target is
\[
 b''_{j,\vartheta}b_j^\circ=\mu_j^0\vartheta\,\upsilon^{(r+t)/2}.
\]
Write \(A_i(g)=\Xi_{\vartheta,i}(g)K^\circ_{\zeta_i}(g)\).
For any left test for which these integrals converge, right
invariance of the quotient measure gives
\[
 \sum_{u\in\bar J^{\,1}}\int y(g) A_i(gk(u))\,d\mu(g)
  =\int \Big(\sum_{u\in\bar J^{\,1}}y(gk(u)^{-1})\Big)
                     A_i(g)\,d\mu(g).
\]
Every \(y\) in the projected source is fixed by the old principal
level, which contains the \(k(u)\). The parenthesized sum is
therefore \(|\bar J^{\,1}|y(g)\) on this source.
The integrals converge absolutely: each class sum is
finite and its theta kernels have moderate growth for
the fixed Schwartz vectors, while \(y\) decreases
rapidly on Siegel sets. For the moderate growth use
continuity of the lattice sum in the Schwartz seminorms
and the polynomial bounds on these seminorms under the
real oscillator action. The powers can be absorbed by
arbitrarily high rapid-decay powers of \(y\).
It follows that
\begin{equation}
 \begin{aligned}
 \langle\zeta,G_\vartheta^*u_{j,\vartheta} y\rangle
  =\sum_i M^{l_1}c_{L_i}|\bar J^{\,1}|
       \int_{[U(W_j)]} y(g)\Xi_{\vartheta,i}(g)
                             K^\circ_{\zeta_i}(g)\,d\mu(g).
 \end{aligned}
\label{eq:unscaled-gluing-pairing}
\end{equation}
At an old place the elements of \(\bar J^{\,1}\) are
\(1+\varpi_w A\) with \(A\in M_r(k_w)\). Thus \(|\bar J^{\,1}|=\prod_{w\in S_M}q_w^{r^2}\).
The earlier trace sum has cardinality
\(\mathcal N=\prod_{w\in S_M}q_w^{(3m)^2}|\mathrm{GL}_{3m}(k_w)|\),
by reduction of invertible matrices modulo \(\varpi_w^2\).
It entered \(c'_\vartheta=\gamma_3\gamma_{12}/\mathcal N\).
Substituting \(c'_\vartheta\Xi_{\vartheta,i}=c'' F_i\) now gives
\[
 c'_\vartheta\langle\zeta,G_\vartheta^*u_{j,\vartheta}y\rangle
 =
 \sum_i M^{l_1}c_{L_i}|\bar J^{\,1}|c''
      \int_{[U(W_j)]} y(g)F_i(g)K^\circ_{\zeta_i}(g)\,d\mu(g).
\]
The normalized descent is conjugated in this pairing, so the first
factor \(\overline{c''\Upsilon_{\vartheta,i}}\) contributes \(c''\),
and the conjugate descent scalar contributes \(c_{L_i}\).
The right side is independent of \(\vartheta\). At the trivial twist
it is \(\langle\zeta,Py\rangle\) for
\(P=c'_1G_1^*u_{j,1}\). Equality for all tuples gives the
operator identity by nondegeneracy.

\begin{proposition}[Common period and the gluing adjoint]
\leavevmode\label{prop:common-period-adjoint}
On the original projected source, at each allowed twist,
\[
                 c'_\vartheta G_\vartheta^*
                       u_{j,\vartheta}=P .
\]
The period \(P\) is independent of \(\vartheta\).
\end{proposition}

In particular \(\ker u_{j,\vartheta}\subseteq\ker P\).
For the proof of injectivity in Section~\ref{sec:common-period-injectivity},
the integral expression for \(\langle\zeta,Py\rangle\) uses the
multipliers \(M^{l_1}c_{L_i}|\bar J^{\,1}|c''\), which are nonzero by
the formula for \(c_{L_i}\) and \eqref{eq:display-i5}.

\section{Coefficient reproduction on the descent space}
\label{sec:period-reproduction}

\sectionopening{In this section, we establish real coefficient separation and the
inner-product identities for the original periods of
Section~\ref{sec:gluing-three-exchanges}. These give the coefficient
tests used in the proof of injectivity in
Section~\ref{sec:common-period-injectivity}.}

All measures, classes, Fock polynomials and finite lattices
in the following argument belong to the original specified space.
The sharp inner-product calculation below includes an
arbitrary rapidly decreasing left test and the compact-source
auxiliary lift. \hyperref[sec:marked-normalized-sparse-overlap]{The adelic overlap calculation}
gives the polarized Fock overlap on the same residue ring,
including the unit coefficient at zero displacement.

The original form in the projected space is read in its holomorphic coefficient
by the untwisted sharp injection and its adjoint.

\subsection{Separation by the circle mean}
\label{sec:marked-period-circle-separation}

\sectionopening{In this subsection, we prove a reproducing identity by integrating
over a circle, then use it to estimate the reference coefficient at a
peak. Section~\ref{sec:peak-period-injectivity} applies the estimate to the common period.}

Fix the original signed real product
\(H=\prod_{\tau\in\Phi}U(J_\tau)\) of rank \(n\), with the
product Haar \(dh\), and let \(f=|\Phi|\).
Use the same injective order \(o\) on the basis, with the
positive block first. The weights \(u,u'\) decrease weakly
within each sign block. At every place the difference of
any negative entry and any positive entry of each weight
is at least \(300n_*\), where
\[
 1\le n\le t,\qquad n_*=2ft^3,\qquad |u_{\tau,i}|\le5b.
\]
For either weight \(v\) define \(m_v\) as the product over
the real places of the minor quotients.

In the
order \(o\), its local factor is
\[
 m_{v_\tau}(h)=\det(h)^{v_{\tau,n}}
       \prod_{a=1}^{n-1}\Delta_a(h)^{v_{\tau,a}-v_{\tau,a+1}},
 \qquad \Delta_a(h)=\det(h_{[1,a],[1,a]}).
 \tag{CM1}\label{eq:circle-minor-coefficient}
\]
The only possible negative exponent on a proper minor is
at the end of the positive block. That minor never
vanishes on \(U(J_\tau)\). Formula \eqref{eq:circle-minor-coefficient} is also the
definition by quotient of two dominant minor products:
add \(2\sum_i|v_{\tau,i}|\) uniformly to the positive
block in the numerator, and divide by the minor product
for that shift. In a compact sign group use the same
formula with no negative proper-minor exponent.

\indent{\it An elementary majorant.}
For a mixed signature write
\(h=\left(\begin{smallmatrix}A&B\\ C&D\end{smallmatrix}\right)\)
in the positive-first order. The equations \(h^*Jh=J=hJh^*\)
imply that the singular values of \(A,D\) are at least one,
\(\|B\|=\|C\|\), and
\(\|A\|^2=\|D\|^2=1+\|C\|^2\).
Consequently \(\|h\|\le\|A\|+\|C\|\le2|\det A|\).
Every positive-initial segment minor has absolute value at most
\(|\det A|\), by the singular-value bound for an exterior
power of the compression of \(A\).
For a longer initial segment, the complementary minor identity
and \(h^{-1}=Jh^*J\), \(|\det h|=1\), express its absolute
value as that of a trailing minor of \(D\).
The same singular-value bound applies, and
\(|\det D|=|\det A|\) by the complementary identity at
the end of the positive block. The nonnegative exponents
and the one cross-block exponent in \eqref{eq:circle-minor-coefficient} therefore give
\[
 |m_{v_\tau}(h)|\le |\det A|^{v_{\tau,1}-v_{\tau,n}}
               \le \min\{1,(2/\|h\|)^{300n_*}\}.
 \tag{CM2}\label{eq:circle-coefficient-decay}
\]
In compact signature \(\|h\|=1\) and every proper minor
has absolute value at most one, giving the same statement.
This argument also proves analyticity locally around each
group point by \eqref{eq:circle-minor-coefficient}.

Let \(I_\tau=\int |m_{u_\tau}|^2\,dh_\tau\).
Take \(L=2\lfloor5b\rfloor+n\), \(\beta=n+L+1\), and
\(\delta_1=(200\beta^5)^{-1}\).
For \(\|h-1\|\le\delta_1\), the logarithm series in each
principal compression gives \(|\log\Delta_a(h)|\le2n\delta_1\).
The sum of absolute exponents in \eqref{eq:circle-minor-coefficient} is at most
\(2n\max_i|u_{\tau,i}|\). Hence \(|m_{u_\tau}(h)|\ge1/2\)
on that ball. A maximal \(\delta_1/Y\)-separated set
in the group ball \(\|h\|\le Y\), \(Y\ge1\), has at most
\((2000\beta^6 Y^2)^{2n^2}\) points, by packing disjoint
ambient matrix balls in dimension \(2n^2\).
The corresponding right translates of the identity
ball cover the group ball, since
\(\|h^{-1}\|=\|h\|\).

In particular ball volumes have
polynomial growth of order \(4n^2\).
With \eqref{eq:circle-coefficient-decay} this first proves \(0<I_\tau<\infty\).
Then the small-ball estimate gives the sharper bound
\[
 dh_\tau\{\|h\|\le Y\}\le4I_\tau(2000\beta^6 Y^2)^{2n^2}.
\]
Put \(c=10n^3\). Integrating the last bound against
\(cY^{-c-1}\,dY\) gives
\[
 \int (2/\|h\|)^c\,dh_\tau
 \le 8I_\tau(2000\beta^6)^{2n^2}2^c
 \le I_\tau 2^{13t^3}b^{14t^2}.
\]
The first inequality uses \(c/(c-4n^2)\le2\).
The second uses the original size condition
\(2000(n+L+1)^6\le b^7\) and \(n\le t\).
No formula for the numerical formal degree is needed
for this estimate.

Write \(X(h)=\max_\tau\|h_\tau\|\) and
\(I=\int_H|m_u|^2\,dh=\prod_\tau I_\tau>0\).
For each \(h\), use \(6n_*\) extra powers in
\eqref{eq:circle-coefficient-decay} at any place attaining \(X(h)\),
and \(c\) powers at the others. Since \(c+6n_*\le300n_*\),
\[
 X(h)^{6n_*}|m_{u'}(h)|
 \le 2^{6n_*}\prod_\tau(2/\|h_\tau\|)^c .
\]
This bound is independent of the choice of maximizing place. Integrating gives
\[
 \int_H |m_{u'}(h)|X(h)^{6n_*}\,dh
 \le 2^{6n_*}(2^{13t^3}b^{14t^2})^f I
 \le b^{n_*}I.
 \tag{CM3}\label{eq:circle-weighted-integral}
\]
In the last step use \(b\ge2^{25}\), \(t\ge14\), and
\(n_*=2ft^3\).

\indent{\it The circle mean and reproduction.}
Let a function \(A\) on \(H\) be holomorphic of weight
\(v\) in the sense of the model. Its analytic
extension to neighborhoods of \(hR\), with \(R\)
invertible upper triangular in order \(o\), transforms
on \(hR\) by the product of diagonal characters
\(\chi_v(R)=\prod_{\tau,i}R_{\tau,ii}^{v_{\tau,i}}\).
Use the circle with probability Haar \(d\theta\), the
simultaneous rotation \(r_s=\operatorname{diag}(s^{o(i)})\)
at every real place, and let \(\chi_v(s)=\chi_v(r_s)\).
Then
\[
 \int_{\mathbb T} A(r_s h)\chi_v(s)^{-1}\,d\theta(s)
                         =A(1)m_v(h).
 \tag{CM4}\label{eq:circle-reproducing-average}
\]
Here is a direct verification. For \(h\) close to the identity
write each \(h_\tau=L_\tau R_\tau\) with \(L_\tau\) lower
unitriangular and \(R_\tau\) upper triangular.

For \(|s|=1\),
\[
 r_s h_\tau=(r_sL_\tau r_s^{-1})(r_sR_\tau).
\]
All strictly lower entries of the first factor have positive
powers \(s^{o(i)-o(j)}\). The factor extends polynomially
to the closed disc, is identity at its center, and stays
in the analytic neighborhood of the identity for \(h\)
sufficiently close to it. Write \(L_s=r_s Lr_s^{-1}\) and
\(\widetilde A\) for the holomorphic extension near the identity.
For \(|s|=1\), the expression \(L_s=(r_s h)(r_s R)^{-1}\)
has first factor in \(H\) and second factor upper triangular.
The character rule gives
\[
 \chi_v(s)^{-1}A(r_s h)=\chi_v(R)\widetilde A(L_s).
\]
The mean value property gives \(\chi_v(R)A(1)\), which equals \(m_v(h)A(1)\)
by the leading-minor formulas for triangular elimination.
For arbitrary \(h\) join identity to \(h\) by a real analytic
path in the connected product \(H\).

Such a path is obtained from the compact and positive
factors of the polar decomposition.
The two sides of \eqref{eq:circle-reproducing-average} are real analytic along it. For the
left side the local analytic extension is uniform in the
circle variable by its compactness. Analytic continuation
along the path proves \eqref{eq:circle-reproducing-average}.

If \(A\overline{m_v}\) is integrable, apply \eqref{eq:circle-reproducing-average} after the
left translations \(h\mapsto r_sh\), which preserve the
original Haar. The principal minors in \eqref{eq:circle-minor-coefficient} transform
by their diagonal weights, so \(m_v(r_sh)=\chi_v(s)m_v(h)\).
Joint integrability follows from the same left
translations and the probability circle. Hence
\[
 \int_H A(h)\overline{m_v(h)}\,dh
    =A(1)\int_H |m_v(h)|^2\,dh.
 \tag{CM5}\label{eq:circle-reproducing-pairing}
\]

\indent{\it The peak estimate.}
Let \(F\) be continuous and bounded on the adelic
group, with \(h\mapsto F(x(h,1))\) holomorphic of weight
\(u\) for every \(x\).
For \(F\ne0\), call \(z\) a peak if \(|F(z)|>\|F\|_\infty/2\).
Such a point exists by the definition of the supremum, since
\(0<\|F\|_\infty<\infty\). Fix one such \(z\).
Suppose \(K_{\rm ref}\) is bounded on \(H\), and at every
peak \(x\) some \(K_x\) of holomorphic weight \(u-u'\)
satisfies
\[
 |K_x(h)-K_{\rm ref}(h)|\le\delta X(h)^{6n_*},
 \qquad
 \int_H\overline{F(x(h,1))}K_x(h)m_{u'}(h)\,dh=0.
 \tag{CM6}\label{eq:circle-approximation-hypotheses}
\]
\begin{lemma}[The circle-mean peak estimate]\label{lem:circle-mean-peak-estimate}
Under the holomorphic weight, peak, and approximation conditions
above, the reference function satisfies
\[
                       |K_{\rm ref}(1)|\le(4b^{n_*}+1)\delta.
\]
\end{lemma}

\begin{proof}

Every \(z_s=z(r_s,1)\) is again a peak: its value differs
from \(F(z)\) by the unit character \(\chi_u(s)\).
Put \(K_*=K_z\) and
\[
 W_s(h)=\overline{F(z(r_sh,1))\chi_u(s)^{-1}}.
\]
Equation \eqref{eq:circle-reproducing-average} at \(z\) gives
\(\int W_s(h)d\theta(s)=\overline{F(z)m_u(h)}\).

Also \(|W_s(h)|\le2|F(z)|\), and \eqref{eq:circle-approximation-hypotheses} eliminates the
integral with \(K_{z_s}\) in place of \(K_*\).
Thus \eqref{eq:circle-weighted-integral} bounds each fixed slice by
\[
 \left|\int_H W_s(h)K_*(h)m_{u'}(h)\,dh\right|
 \le4\delta |F(z)| b^{n_*}I .
\]
Only \(K_*\) occurs in the joint integral. It is jointly
measurable and integrable by \eqref{eq:circle-weighted-integral}, since \(K_*\) is
bounded by a constant times \(X^{6n_*}\).
No regularity in the choice of \(K_{z_s}\) as a function
of \(s\) is needed. Average the slice bound and use
Fubini and \eqref{eq:circle-reproducing-pairing} on the product \(A=K_*m_{u'}\), which is
holomorphic of weight \(u\). Since \(m_{u'}(1)=1\), this gives
\[
 I|F(z)|\,|K_*(1)|\le4\delta |F(z)|b^{n_*}I .
\]
Consequently,
\[
                 |K_{\rm ref}(1)|\le(4b^{n_*}+1)\delta .
 \tag{CM7}\label{eq:circle-peak-bound}
\]
Indeed, cancellation gives \(|K_*(1)|\le4\delta b^{n_*}\), and
\eqref{eq:circle-approximation-hypotheses} at \(h=1\) bounds
\(|K_{\rm ref}(1)-K_*(1)|\) by \(\delta\), since \(X(1)=1\).
\end{proof}

\subsection{The inner product formula at the principal levels}
\label{sec:marked-period-sharp-level}

\sectionopening{In this subsection, we unfold the inner product with an arbitrary rapidly decreasing left
test on the original levels. The calculation extends the use of the
Gram identity in Section~\ref{sec:marked-first-arrow-rallis} to the
left tests needed for injectivity.}

Let \(S,T\) be the original Hermitian spaces of ranks \(n,s\),
with \(s>2n\), and retain their finite unitary bases. Write
\(G=U(S)\), \(H=U(T)\), and use the original quotient measures.
On the finite group normalize the integral compact at each
place of \(F_0\) to have mass one. On additive compact quotients
in the Weil calculation use probability Haar. These conventions
leave the factor \(V_T=\operatorname{vol}([H])\) when passing
from the probability average to the original target pairing.
The argument applies to the pairs
\[
 (n,s)=(r,t-9m-1),\quad (3m,6m+1),\quad (m,3m).
\]
In the first pair \(s-2n=3r-9m-1>0\) for sufficiently
large \(x\), since \(t=5r\), \(r\ge x^5\) and \(m<2\sqrt x\).
In the second and third pairs \(s-2n=1\) and \(m\),
respectively. The first two pairs use the original point datum
at \(\mathfrak m=(\pi_M^2)\). The third uses the unit-modulus
vacuum and has compact real source.

\indent\textit{The sharp datum and its mass.}
At the chosen member \(w\) of each chosen split pair the
Schr\"odinger variable \(X\) is an \(s\)-by-\(n\) matrix.
The point datum is the indicator of \(X\equiv\iota\pmod{\mathfrak m}\)
for the inclusion \(\iota\), with vacuum at the other
places. Write \(R=\mathcal O_E/\mathfrak m\). On the split matrix slice,
the unit vacuum is the indicator of the integral box of additive
volume one. Its residue classes have equal volume, so the inner
product of two tests inflated from \(M_{s,n}(R)\) is their average
there. The squared norm of the point
datum is \(|R|^{-sn}\). Let \(J\) be the original principal
level in the standard compact \(K^0\) at this modulus.
The integral action on the slice is
\[
              \Psi(X)\longmapsto b(\det g)^{-1}\Psi(Xg).
\]
Here \(b\) is the source-actuating character of the lift.
Its tame part is trivial away from the prescribed support, and
on \(J\) it evaluates to one. In the two point tests its
powers are, respectively, \(b=\upsilon^{m/2-j}\) and \(b=1\).
Both characters are trivial on integral determinants at every finite
place, since \(\upsilon_f\) is unramified. By
Proposition~\ref{prop:marked-level-primes}, \(M=p_0p_1\) and
\(N\notin\{p_0,p_1\}\). The chosen members over \(p_0,p_1\) have
exponent two in \(\pi_M^2\). The opposite members are determined by
duality. Every place over \(N\) has the unit vacuum and the full
integral compact in these reproductions. The additional covector
condition at \(\pi_N\) in Section~\ref{sec:ramified-source} is on
the earlier definite source of rank \(k\).

Let \(g\) be integral off the pairs and suppose it gives
nonzero overlap. At a chosen split place dividing \(\mathfrak m\)
the source acts in the Schr\"odinger model by \(X\mapsto Xg\),
times a nonzero determinant scalar. The congruences to the
inclusion matrix give integral left inverses of both \(X\) and \(Xg\). Applying them to \(Xg\) and \(X\),
respectively, makes \(g\) and \(g^{-1}\) integral at the
chosen member of every pair. At the opposite member
this follows by inverse transpose. The congruences
\(X\equiv Xg\equiv\iota\), followed by a left inverse of
\(\iota\) in the residue ring, give \(g\equiv1\).

Conversely \(J\) fixes the test. The character on the groups therefore gives no factor on its coefficient support.
The coefficient on the points integral off the pairs
is exactly \(|R|^{-sn}\mathbf1_J\).
Its integral is
\[
            c_{\mathfrak m}
              =\frac{|R|^{-sn}}{|\operatorname{GL}_n(R)|}.
\]
At a selected place \(w\) of exact modulus exponent two the
local mass is
\[
       \frac{q_w^{-2sn}}{q_w^{n^2}|\operatorname{GL}_n(k_w)|}.
\]
Indeed reduction from exponent two to exponent one is
surjective with kernel \(1+\varpi_w M_n(k_w)\).
The split group over \(F_0\) is one copy of
\(\operatorname{GL}_n(E_w)\). The opposite matrix is
determined by the inverse transpose and does not give a
second volume. At unit modulus \(J=K^0\), the finite
vacuum is a unit vector, and \(c_{(1)}=1\).

\subsubsection{The adelic doubling identity}
\label{sec:reproduction-doubling}\label{sec:marked-graph-doubling}

\sectionopening{In this subsection, we identify the doubled kernel in the original quotient measures after the mass calculation of Section~\ref{sec:marked-period-sharp-level}.}

Let \(v\) be the original real Fock vector and \(f=v\otimes F\)
the datum with the finite part just specified. Denote the
rational Gram matrices by \(J_S,J_T\) and the adelic
point matrices by \(X_S,X_T\). Put \(u=(1,u_f)\) with
\(u_f\overline{u_f}=-1\), so the opposite space has
\(J_{S^-}=-J_S\) and \(X_{S^-}=uX_S\).
At infinity the coordinate transport is the identity
and the alternating form changes sign. The conjugate
oscillator datum there is transported by complex conjugation.
For an oscillator vector \(f\) write
\(D_f(X)=\theta(\mathcal R(X)f)\) on the adelic displacement
space. Let \(Jf'(Y)=\overline{f'(u^{-1}Y)}\) be the
antilinear transport on oscillator data. The probability-mean
argument of Section~\ref{sec:v2-first-exchange-phase} for this pair gives
\(\overline{\theta_{\rm old}(f')}=\lambda_{S,T}\theta_{\rm new}(Jf')\)
with \(|\lambda_{S,T}|=1\). Set \(C=\lambda_{S,T} J\). Moving a Heisenberg
translation across \(J\) gives
\[
 D^{\rm new}_{C f'}(Y)
 =\lambda_{S,T}\theta_{\rm new}
     \bigl(J\mathcal R_{\rm old}(u^{-1}Y)f'\bigr)
 =\overline{D_{f'}(u^{-1}Y)}.
\]
Use the tensor-product normalization of the doubled theta
functional from Section~\ref{sec:gluing-real-matched-factors}.
In a direct sum of rational Schr\"odinger polarizations, summation
over the product of the two rational Lagrangians gives the
displaced value
\[
     D_f(X)\overline{D_{f'}(u^{-1}Y)}
\]
on \(f\otimes C f'\). If the doubled functional has the previously
fixed phase \(\epsilon\) relative to this product, its theta
integral and each graph zero term, obtained by averaging the same
displaced functional, both acquire \(\epsilon\). Dividing the
Siegel--Weil identity by this common phase gives the same pairing.
Write \(a(\det h)b(\det g)^{-1}\) for the pair character
and put \(\upsilon_1=\upsilon_E|_{T_1(\mathbb A_{F_0})}\).
The determinant of \(X\mapsto hXg^{-1}\) on the Hom space is
\((\det h)^n(\det g)^{-s}\). The multiplier calculation of
Section~\ref{sec:v2-first-exchange-phase} gives the opposite pair
\((a^-,b^-)=(a^{-1}\upsilon_1^n,b^{-1}\upsilon_1^s)\).
On the diagonal \(H\) the product of target twists
\(\upsilon_1(\det h)^n\) cancels the mixed-model Levi scalar
\(\upsilon_1(\det h)^{-n}\) of Section~\ref{sec:marked-doubled-action}.

In the split polarization the rational Lagrangian
\[
 L_\delta=\{(\delta^{-1}y,-y):y\in S\},\qquad \delta\in G(F_0),
\]
is parameterized by \(y\mapsto(\delta^{-1}y,-y)\).
For a matrix \(Z\)
in the adelic zero term set
\[
 M=X_T^{-1}J_T^{-1}Z J_S X_S,\qquad
 \widehat\delta=X_S^{-1}\delta X_S.
\]
The adjoint of this parametrization gives the blocks
\((M\widehat\delta,uM)\), since
\((\delta^{-1})^*J_S=J_S\delta\) and
\((-1)^*(-J_S)=J_S\). The second tensor factor
therefore has argument \(M\) in its original coordinates.
The change \(Z\mapsto M\) maps the standard rational
matrices bijectively onto
\[
 \{X_T^{-1} W X_S:W\in\operatorname{Mat}_{s\times n}(E)\},
\]
because \(J_T^{-1}\) and \(J_S\) are rational and
invertible. The induced isomorphism of additive quotients
preserves probability Haar. Rational Weil covariance gives
\(D_f(M\widehat\delta)=D_{\omega_b(\delta)f}(M)\), since
\(b\) is trivial on rational determinants. Under a rational
translation of \(M\), the two displaced values have the same
Heisenberg multiplier. Their product below is therefore a function
on the additive quotient. The probability mean of
Section~\ref{sec:overlap-fourier-inversion} gives its value:
\[
 \int D_{\omega_b(\delta)f}(M)\overline{D_{f'}(M)}\,dM
      =\langle f',\omega_b(\delta)f\rangle .
\]
This is the graph zero term. Translating the two data separately by
\(g,g'\) and using unitarity gives the coefficient
\[
 \langle f',\omega_b(g'^{-1}\delta g)f\rangle .
\]

Only these graphs contribute. If \(L\) is a Lagrangian and
its second projection is injective, it is an isomorphism
since \(\dim L=n\); the isotropy equation then makes the
first projection an isometry. Otherwise \(L\) contains
\((c,0)\) with \(c\ne0\). In the compact-source case this is
impossible by definiteness. For either point test choose a split place
\(w\mid\mathfrak m\), which exists since \(M>1\).
The zero-term functional is invariant under the Heisenberg translations
dual to this vector, acting on the first tensor factor alone.

Indeed, for \(z\in T(\mathbb A)\), the rank-one map
\((v,v')\mapsto z\,h_S(c,v)\) vanishes on \(L\) by isotropy.
It belongs to the isotropic displacement space averaged in the
zero term, so translating by it preserves that mean. It acts
only on the first tensor factor. Taking \(z\) in the local
dual coordinate at the split pair gives a pure modulation on
the chosen Schr\"odinger slice. After moving the source point
through the Heisenberg action these modulations act
by the multipliers
\[
              X\longmapsto
                     \psi_w(\operatorname{tr}(ZXc_w)),
\]
where \(Z\) varies over the dual local vector space and
\(c_w\ne0\) is the coordinate in the chosen basis of \(c\).

An invariant functional kills a compactly supported,
locally constant test supported where \(Xc_w\ne0\):
partition its support into finitely many compact open
pieces on each of which one such multiplier is a constant
different from one, and apply invariance on each piece.
Every matrix in the original point support has an integral
left inverse. After any source translate at \(w\), the support
is obtained by invertible right multiplication, so every matrix
there still has column rank \(n\) over \(E_w\). Hence
\(Xc_w\ne0\) on that support for every nonzero \(c_w\).
This proves vanishing for the translated datum at each pair
\((g,g')\), and the graph terms exhaust the doubled expansion.

The doubled action is identified in Section~\ref{sec:marked-doubled-action}.
Apply the Siegel--Weil formula
\cite[Theorem~2, p.~237]{YamanaSW} to the split polarization
of rank \(n\) of \(S\oplus S^{-}\), with \(T\) as the
Hermitian space of dimension \(s\). The doubled splitting
uses \(\upsilon_E^{\,s}\), whose restriction to
\(\mathbb A_{F_0}^\times\) is \(\varepsilon_{E/F_0}^{\,s}\).
Yamana's division algebra is \(E\), and his Witt index
of \(T\) is at most \(s/2\). Hence the convergent clause
\(s-\operatorname{Witt}(T)>n\) applies with coefficient one
and \(s_0=(s-n)/2>n/2\).

Let \(\Phi=f\otimes C f'\),
\(d(g,g')=\operatorname{diag}(g,(g')^-)\), and let \(\Omega\) be
the doubled action with target twist \(\upsilon^n\) and trivial
source twist. Direct-sum covariance and the opposite character
formula above give the separately translated tensor
\[
 \omega_b(g)f\otimes C\omega_b(g')f'
  =b(\det g)^{-1}b^-(\det g')^{-1}\,\Omega(d(g,g'))\Phi .
\]
The scalar has absolute value one. It multiplies the theta integral
and every Siegel section by the same factor. Use the transported
Fock polynomial identity of Section~\ref{sec:marked-doubled-action}
for the original tensor. Its section at \(s_0\) satisfies
\[
 f_{\Omega(d(g,g'))\Phi}^{(s_0)}(a)
       =f_{\Phi}^{(s_0)}(a\,d(g,g')),
 \qquad f_{\Phi}^{(s_0)}(a)=(\omega(a)\Phi)(0).
\]
Thus the formula for the fixed polynomial datum applies at all
the translated points. The strict-range formula uses target
quotient measure of mass one
\cite[\S52, Theorem~5 and the normalization following (41)]{WeilSiegel}.
Multiplying by \(V_T\) converts it to the original Haar measure
and gives
\[
 \int_{[H]}\overline{\theta_{f'}(g',h)}\,
                       \theta_f(g,h)\,dh
  =V_T\sum_{\delta\in G(F_0)}
       \langle f',\omega_b(g'^{-1}\delta g)f\rangle.
\]
Here \(dh=V_T\,dh_{\mathrm{prob}}\) on \([H]\), with the
probability convention of \cite[\S1.3, p.~1007]{YamanaSingular}.
For the two source integrations we need absolute bounds
in the original quotient measures. In the split
Schr\"odinger model and for a doubled datum \(\Phi\), set
\[
 J(\Phi)=\int_{[H]}\sum_{X\in T^n(F_0)}
                                   |\Phi(h^{-1}X)|\,dh .
\]
Here the rational points of \(T^n\) are taken over \(F_0\)
as for the restriction of scalars. For a fixed finite
factor this is finite on the real Schwartz space.

Indeed the absolute value of the real factor is bounded
by a nonnegative Schwartz function, by taking smooth
annular cutoffs with rapidly decreasing coefficients.
Use that majorant and the absolute value of the finite
factor in the convergent theta integral
\cite[\S51, Proposition~8, p.~75]{WeilSiegel}. The target acts
by substitution in the split model of
Section~\ref{sec:marked-doubled-action}. The involution
of \(E\) is of the second kind, so Weil's type~(I)
parameter is \(1/2\) by \cite[\S\S26--27]{WeilSiegel}.
His convergence condition is therefore
\(s-\operatorname{Witt}(T)>n\), as verified above.

The absolute integral \(J\) is bounded by a continuous
real Schwartz seminorm with the finite factor fixed.
To see this, choose a measurable quotient section
and take finite partial sums of the rational evaluations,
with measurable multipliers of modulus at most one.
After integration over the section these are continuous
linear functionals, each bounded by the real supremum
norm times a finite constant. Their common pointwise
bound is \(J\). Uniform boundedness on the Schwartz
Fr\'echet space gives a common continuous seminorm.
Taking the supremum over partial sums and multipliers
proves the claim. In particular \(J\) bounds the absolute
product of the original theta kernels integrated over
\([H]\).

At the identity of the doubled group the Iwasawa formula bounds
the absolute graph sum by
\[
 E_{\rm sph}(1,s_0)\sup_{k}|(\omega(k)\Phi)(0)|.
\]
Here the supremum is over the maximal compact and the spherical
Eisenstein series converges absolutely since \(s_0>n/2\).
At \(g,g'\) apply this same bound to the translated doubled
datum at the identity.
That supremum is a continuous real Schwartz seminorm
up to a fixed constant. The real compact acts continuously
on Schwartz data and the finite compact has finitely
many translates of the fixed finite factor.
The absolute graph sum is rationally invariant in the
two source variables by reindexing. The corresponding
absolute product on the target side is invariant by
rational theta invariance. Thus we may reduce both
source representatives as in
Section~\ref{sec:adjacent-unfolding}. Their finite
translates range over finite sets.

Multiplication and differentiation operators are
conjugated by the real Weil action into linear
combinations with coefficients from the symplectic
matrix and its inverse. By unitarity their products
give polynomial bounds for the translated weighted
Sobolev norms, which generate the Schwartz topology.
The two seminorm bounds for the doubled translate
are therefore at most
\(C\mathcal H(g)^\nu\mathcal H(g')^\nu\) for the
original reduction heights. Set
\(k_{f',f}(a)=\langle f',\omega_b(a)f\rangle\).
For rapid \(y,y'\), choose decay powers at least \(\nu\) in each
variable. With the original quotient measures this gives
\[
 \int_{[G]^2}|y(g)y'(g')|
 \left(
    \int_{[H]} |\theta_{f'}(g',h)\theta_f(g,h)|\,dh
     +\sum_{\delta\in G(F_0)}|k_{f',f}(g'^{-1}\delta g)|
 \right)\,dg\,dg' <\infty .
\]
Indeed the integrand on \([G]^2\) is bounded and the quotient has
finite volume. The first term permits interchange of the target
and source integrals. The second permits termwise integration
of the rational graph expansion. This argument also applies to
the strict-range doubled first arrow.

For rapid \(y',y\), multiply by
\(\overline{y'(g')}y(g)\) and unfold the rational graph sum.
For each fixed \(g'\), the sum of the translated source
quotient integrals is the \(G(\mathbb A)\)-integral in
\(a=g'^{-1}\delta g\). Indeed the translates by \(G(F_0)\) of a
measurable fundamental domain tile \(G(\mathbb A)\) up to a null set,
and left translation by \(g'^{-1}\) preserves Haar measure.
This uses counting measure on \(G(F_0)\) in the source quotient
formula. Left automorphy gives \(y(g)=y(g'a)\) on that graph term.
Thus
\[
 \langle L_{f'}y',L_fy\rangle_T
    = V_T\int_{G(\mathbb A)}\langle f',\omega_b(a)f\rangle
                      \langle y',R(a)y\rangle_S\,da .
\]
In particular the first test need not be fixed by \(J\).

\subsubsection{Vacuum coefficients away from the chosen places}
\label{sec:reproduction-vacuum-coefficients}

\sectionopening{In this subsection, we evaluate the local factors in Section~\ref{sec:reproduction-doubling}.}

Let \(A\) be the union of conjugation orbits outside the
support of \(\mathfrak m\bar{\mathfrak m}\).
At \(A\) put \(\Lambda=\widehat{\mathcal O}_E^{\,n}\) in
the original finite basis. For \(\gamma\) supported on \(A\)
put \(\Delta(\gamma)=(\upsilon_E)_f(\det B)|\det B|^{-1/2}\),
where \(B\Lambda=\Lambda\cap\gamma^{-1}\Lambda\).
This expression is independent of the integral basis of
the intersection since \((\upsilon_E)_f\) is unramified there.
The index formula is
\[
 I(\gamma)=|\Delta(\gamma)|^2
     =[\Lambda:\Lambda\cap\gamma^{-1}\Lambda]
     =[\gamma\Lambda:\Lambda\cap\gamma\Lambda].
\]
In the normalized finite lattice model the uncharactered
Weil operator for an isometry \(a\) of the displacement
space is
\[
 \omega(a)u=\Delta_{\rm disp}(a)
       \int u(a^{-1}z)\mathcal R(z)v_{\rm vac}\,dz,
\]
where \(z\) is a finite Heisenberg displacement,
\(\mathcal R(z)\) its operator, and the vacuum lattice
has additive measure one.

In the original finite displacement coordinates the vacuum lattice is
\(D=\alpha_{\rm W}\widehat{\mathcal O}_E^{\,sn}\).
The inverse different identity of Section~\ref{sec:period-finite-coordinates}
makes \(D\) self-dual. Write \(z=\alpha_{\rm W}u\) and
\(c_{\rm tr}=(1-j_{\rm W}^{-1})/2\). The phase on \(D\) is
\[
 \eta(z)=\psi_f\!\left(\tfrac12c_{\rm tr} j_{\rm W}
                                  \sum_\nu \bar z_\nu z_\nu\right)
        =\psi_0\!\left(\tfrac12c_{\rm tr}\sum_\nu\bar u_\nu u_\nu\right),
 \qquad \psi_0(t)=\psi_f(j_{\rm W}\bar\alpha_{\rm W}\alpha_{\rm W}t).
\]
Here the finite unitary bases give the identity Hermitian matrix
on the Hom coordinates. On \(D\cap aD\) the product
\(\eta(a^{-1}z)\overline{\eta(z)}\) equals one by isometry,
including at places above two. The integrand vanishes off this
intersection. Since \(a\) preserves additive Haar and \(D\) has mass
one, its integral is
\[
 \operatorname{vol}(D\cap aD)
     =[D:D\cap a^{-1}D]^{-1}=|\Delta_{\rm disp}(a)|^{-2}.
\]
It follows that
\[
 \langle v_{\rm vac},\omega(a)v_{\rm vac}\rangle
       =\Delta_{\rm disp}(a)|\Delta_{\rm disp}(a)|^{-2}
       =\overline{\Delta_{\rm disp}(a)}^{-1}.
\]
On \(\operatorname{Hom}(S,T)\) the action of \(\gamma\)
is \(X\mapsto X\gamma^{-1}\). The \(s\) covectors give
\(\Delta_{\rm disp}(a)=\Delta(\gamma^{-1})^s
                      =\overline{\Delta(\gamma)}^s\).
For the last equality, the inverse of a unitary matrix
is its conjugate transpose. Transposition preserves
the elementary divisors, and conjugation takes the
idele splitting multiplier to its complex conjugate
in the fixed unramified normalization. The character
\(b(\det\gamma)^{-1}\) therefore gives the exact coefficient
\[
       a(\gamma)=\bigl(b(\det\gamma)\Delta(\gamma)^s\bigr)^{-1},
       \qquad |a(\gamma)|=I(\gamma)^{-s/2}.
\]
It multiplies over disjoint conjugation orbits and is
bi-invariant by the standard integral group on \(A\).
Explicitly, if the elementary divisors at \(w\) are
\(\varpi_w^{e_i}\), the local contribution to
\([\Lambda:\Lambda\cap\gamma^{-1}\Lambda]\) is
\(q_w^{\sum_i\max(0,-e_i)}\). At a split pair the opposite
exponents are the negatives of these, so the pair contributes
\(q_w^{\sum_i |e_i|}\). At an inert place the multiset is its own
negative by unitarity, giving \(q_w^{\frac12\sum_i |e_i|}\).
By Section~\ref{sec:different-additive-character-scales}, the
relative extension has no other finite local type.

\subsubsection{Summability and injectivity}
\label{sec:reproduction-summability}

\sectionopening{We use Kneser's method of neighboring lattices
to bound the coefficients of Section~\ref{sec:reproduction-vacuum-coefficients}.
For the construction of Hermitian neighbors at primes where the lattice is
unimodular, see \cite[Propositions~5.3--5.4]{GreenbergVoightLattices}.}

At a split pair a lattice at \(w\) determines the
opposite lattice by duality. If the moved lattice is
not contained in the standard one at \(w\), take a vector
\(x\) in it with \(\varpi_w x\) primitive in the standard
lattice. Enlarge the standard lattice by \(x\) and use
its dual at \(\bar w\). The intersection with the moved lattice
grows by a factor \(q_w\) at \(w\), since the new quotient is
generated by the class of \(x\) with annihilator \((\varpi_w)\).
At \(\bar w\) every element of the old intersection pairs
integrally with \(x\), by duality with the moved lattice at \(w\).
The opposite intersection is therefore unchanged. If containment is
strict in the other direction, choose the vector at
\(\bar w\). Each enlargement at a chosen member is
determined by the line of \(\varpi_wx\) in \(k_w^n\),
giving at most \((q_w^n-1)/(q_w-1)\) steps for that member.

At a nonsplit finite place the extension \(E_w/F_{0,v}\)
is unramified. Take the uniformizer \(\varpi_w\) in \(F_{0,v}\).
Both lattices are self-dual, so containment implies equality. If
they differ, scale a vector in the moved lattice outside the
standard one by a power of \(\varpi_w\) to obtain \(x\) with \(y=\varpi_w x\) primitive
in the standard lattice \(\Lambda_w\), and put
\[
 \Lambda_0=\{z\in\Lambda_w:\langle y,z\rangle\in
                                          \varpi_w\mathcal O_{E_w}\}.
\]
The lattice \(\Lambda'=\Lambda_0+\mathcal O_{E_w}x\) is integral:
\(\langle x,x\rangle\) is integral in the moved lattice and
\(\langle x,\Lambda_0\rangle\subseteq\mathcal O_{E_w}\) by definition.
The primitive vector \(y\) defines a nonzero residue functional, and
\(y=\varpi_w x\) belongs to \(\Lambda_0\). Hence
\([\Lambda_w:\Lambda_0]=[\Lambda':\Lambda_0]=q_w\).
The integral lattice \(\Lambda'\) has the covolume of the self-dual
lattice \(\Lambda_w\), so it too is self-dual.
Write \(D_w\) for the moved lattice and
\(D_0=D_w\cap\Lambda_w\). Integrality of \(D_w\) gives
\(D_0\subset\Lambda_0\), so \(D_w\cap\Lambda_0=D_0\).
Reduction gives an isomorphism
\[
 (D_w\cap\Lambda')/D_0\longrightarrow\Lambda'/\Lambda_0.
\]
Indeed the kernel is zero and \(x\in D_w\cap\Lambda'\) maps to a
generator. Thus the index is exactly \(q_w\).

For a fixed residue line choose a representative of
\(\bar y\); choices of \(x\) then differ by \(z\in\Lambda_w\).
The integrality of the two Hermitian squares requires
\(\operatorname{Tr}_{k_w/k_v}\langle\bar y,\bar z\rangle=0\).
This leaves at most \(\sqrt{q_w}\) classes modulo
\(\Lambda_0\) for each residue line. These neighboring
self-dual lattices admit integral orthonormal bases.
One way to see this last assertion is to find a vector
of unit norm by nondegeneracy of the residual Hermitian
form, normalize it by surjectivity of the unramified
norm on units, and continue on its unimodular orthogonal
complement. The residual norm polynomial is not
identically zero on a nondegenerate Hermitian space:
if every self-pairing vanished, polarizing after
multiplication by all elements of \(k_w\) and using
the nondegenerate trace would make every pairing zero.
This argument also applies in residue characteristic two.

At each orbit, fix for each neighboring lattice an isometry
\(\rho\) carrying it to the standard lattice. Extend it by the
identity off the orbit. At split places the component at the
opposite member is the inverse conjugate transpose.
The intersection
computations at the nonsplit place and the two split places give
\[
 I(\rho\gamma)=q_w^{-1}I(\gamma),\qquad
 |a(\gamma)|=q_w^{-s/2}|a(\rho\gamma)|.
\]
Indeed \(I\) is the index in the moved lattice, and \(\rho\)
preserves indices. For a union \(B\) of whole orbits, take
representatives of \(G_B/K_B\) supported on \(B\).
An index-one moved lattice is contained in the standard lattice,
and dualizing gives equality, hence the identity coset. For
every other coset choose one reducing step at an orbit where
the lattices differ. The positive integral index strictly decreases.
If the successive steps are \(\rho_1,\ldots,\rho_h\), then
\(\gamma K_B=\rho_1^{-1}\cdots\rho_h^{-1}K_B\).
Thus the chosen words distinguish the cosets. The total weight of one-step
choices at \(w\) is bounded, since \(s\ge2n+1\), by
\[
 a_w\le
 \begin{cases}
       \dfrac{q_w^n-1}{q_w-1}q_w^{-s/2}
             \le\dfrac{1}{(q_w-1)\sqrt{q_w}},&w\text{ split},\\[5pt]
       \sqrt{q_w}\,\dfrac{q_w^n-1}{q_w-1}q_w^{-s/2}
             \le\dfrac{1}{q_w-1},&w\text{ nonsplit}.
 \end{cases}
\]
Every \(q_w\ge4\), since \(E\) contains both a primitive
third and a primitive fourth root of unity. A split
orbit has two members of equal norm, so its total is
at most \(1/3\), as is the weight at a nonsplit orbit.
The weights are summable. At split places they are
\(O(q_w^{-3/2})\), and at nonsplit places
\(q_w\ge p^2\) for the rational prime \(p\) underneath.
There are at most \([E:\mathbb Q]\) places over \(p\).

For any union \(B\) of conjugation orbits of weight
\(\sum_{w\in B}a_w\le1/3\), sum the weights over the
descending words. The identity contributes one and
the entire remaining mass is at most
\[
       \sum_{h\ge1}\Bigl(\sum_{w\in B}a_w\Bigr)^h
                     \le\frac12.
\]
This proves absolute convergence and injectivity of
\(\mathcal V_B=\sum_{\gamma\in G_B/K_B}a(\gamma)R(\gamma)\)
on bounded functions right invariant under \(K_B\).
Indeed \(\|\mathcal V_B u-u\|_\infty\le\|u\|_\infty/2\).
Assign to every finite conjugation orbit the numerical upper bound
in the preceding display, adding the two bounds at a split pair.
These majorants depend only on the field and ranks. Choose finitely many orbits
so that the sum of the majorants over all remaining orbits is at most
\(1/3\). Partition the full set of orbits into these single orbits and
the tail, then intersect each part with \(A\). Each resulting piece
\(B\) has weight at most \(1/3\), and \(G_A/K_A\) is the product
of the quotients \(G_B/K_B\). Representatives on
disjoint supports commute. The sum of the absolute weights on
this product is the product of the convergent factor sums
by Tonelli, and multiplicativity gives the same
factorization for \(\mathcal V_A\).
Disjoint supports preserve the other factors' levels, and
bi-invariance of the coefficient preserves each factor's own level.
On the Banach space of bounded functions at the prescribed level,
\(a(1)=1\) and the preceding estimate gives
\[
 \|1-\mathcal V_B\|\le\tfrac12,\qquad
 \mathcal V_B^{-1}=\sum_{d\ge0}(1-\mathcal V_B)^d,\qquad
 \|\mathcal V_B^{-1}\|\le2 .
\]
An empty piece gives the identity. If the original partition has \(h\)
parts, the product is invertible with
\(\|\mathcal V_A^{-1}\|\le 2^h\) for every allowed choice of the
omitted places at these ranks and field.

\subsubsection{Real convolution and the original measures}
\label{sec:reproduction-real-measures}

\sectionopening{In this subsection, we combine the absolutely convergent integrals of Section~\ref{sec:reproduction-summability} with the real projectors.}

Take \(f'=f\) and \(y\) of level \(J\), and put
\(k_\infty(h)=\langle v,\omega_b(h)v\rangle\).
For fixed real \(h\), \(\langle y',R(h\gamma)y\rangle\)
is right \(J\)-invariant in the finite variable, even for a
test \(y'\) of a different level.
On each coset \(\gamma K_A\) the compact at \(A\) has mass
one. Write \(J_{\rm sel}\) for the product of the principal
subgroups at the selected split pairs, so \(J=K_A J_{\rm sel}\).
For the finite part \(F\) of the original tensor, multiplicativity
of the local Weil coefficients gives, for \(k\in K_A\) and an element
\(a_{\rm sel}\) supported at the selected pairs,
\[
 \langle F,\omega_b(\gamma k a_{\rm sel})F\rangle
  =a(\gamma)|R|^{-sn}\mathbf1_{J_{\rm sel}}(a_{\rm sel}).
\]
At all other places where \(\gamma\) is integral the vacuum
coefficient is one. The formula also applies at unit modulus, with the
empty product at the selected pairs. As \(y\) is fixed by \(J\),
the \(K_A\)-integral has mass one and the scaled
\(J_{\rm sel}\)-integral has mass
\(|R|^{-sn}\operatorname{vol}(J_{\rm sel})=c_{\mathfrak m}\).
The sum and real integral of absolute values after taking the
source pairing are bounded by
\[
 c_{\mathfrak m}\|y'\|_2\|y\|_2\,
 \|k_\infty\|_1\sum_{\gamma\in G_A/K_A}|a(\gamma)|,
\]
by unitarity of right translation. The series of absolute
values converges by Section~\ref{sec:reproduction-summability}.
The real \(L^1\) bound is proved below. Fubini therefore gives
\[
 \langle Ly',Ly\rangle_T
    =V_Tc_{\mathfrak m}
          \left\langle y',\mathcal V_A
             \left(\int_{G(\mathbb R)}
                          k_\infty(h)R(h)y\,dh\right)\right\rangle_S.
\]
Here \(L=L_f\). In each holomorphic test let \(M_u=m_u\) be the
minor quotient of \eqref{eq:circle-minor-coefficient}, with \(u\) the
source weight of \(v\) for \(\omega_b\). The coefficient
\(h\mapsto\langle v,\omega_b(h)v\rangle\) is holomorphic of
weight \(u\): in the highest minor, the strictly upper
triangular operators within a sign block annihilate its
flag minors, and the mixed second derivatives vanish because
only one sign of source variable occurs in each of the two
definite-target tests. The diagonal gives the weight \(u\).
The extension to the complex upper triangular group is the
one proved in Section~\ref{sec:marked-adjoint-good-holomorphic}.
Unitarity and the circle action on the first vector give
\[
 k_\infty(r_s h)=\chi_u(s)\,k_\infty(h).
\]
With the pairing linear in its second variable, the circle average
in \eqref{eq:circle-reproducing-average} therefore yields
\(k_\infty(h)=k_\infty(1)m_u(h)=\|v\|^2M_u(h)\).
The opposite-sign gap gives \(M_u\in L^1\). On the compact real group of the auxiliary lift
source the coefficient is also integrable. Absolute
summability of the finite series and integrability on the
real group justify interchanging them.

If \(Ly=0\), the formula shows that
\(z=\mathcal V_A\int k_\infty(h)R(h)y\,dh\) is orthogonal
to every rapid test. Each right translate of \(y\) is such
a test. The function \(z\) is continuous and bounded by
absolute convergence and has finite \(L^2\) norm on the
finite-volume quotient. Expressing the second occurrence
of \(z\) in \(\langle z,z\rangle_S\) by this absolutely
convergent series and integral of translates gives zero.
Thus \(z=0\). The real convolution is right integral-level
invariant on \(A\), so the injectivity above gives
\[
           \int_{G(\mathbb R)}k_\infty(h)y(xh)\,dh=0
               \qquad\text{for every }x.
\]
This is the real convolution implication for both tableau reproductions and the auxiliary lift.

\subsection{The adjoint of a single lift and its holomorphic type}
\label{sec:marked-adjoint-good-holomorphic}

\sectionopening{In this subsection, we apply the Gram identity of
Section~\ref{sec:marked-first-arrow-rallis} to the original lift and
its adjoint. Tracking the complex conjugation of the real coefficient
gives the weight of the projected source used in
Section~\ref{sec:common-period-injectivity}.}

Fix \(j\) and the untwisted descent pair
\(W=W_j\), \(V_j=I_{r+t}\). Write \(d_{\rm tgt}=r+t\),
\(Y^{\rm pr}\) for the subspace of \(L^2([U(W)])\) fixed by the old
principal level and the prescribed real-line projector, and
\(Y=Y_j=e_0Y^{\rm pr}\). Every type summand retained by \(e_0\) is
cuspidal by the constant-term argument of
Section~\ref{sec:marked-comparison-period-vanishing}, so \(Y\) is the
original projected source.
Write \(L:Y^{\rm pr}\to\mathcal A(V_j)\) for the single
lift into the original algebraic form space and
\(L_Y=L|_Y=u_{j,1}\). The metric on \(Y^{\rm pr}\) is its
source Petersson metric, and that on \(\mathcal A(V_j)\)
is the target class metric.

The first-arrow formula of
Section~\ref{sec:marked-first-arrow-rallis} on this pair gives,
for \(x,y\in Y\),
\[
       \langle L_Y x,L_Yy\rangle_{\rm class}
           =A_j\langle x,R_j^W y\rangle_W .
\]
Here \(A_j\ne0\) is the nonzero real edge constant with
its positive class and sharp-level factors. Its real datum
is the joint harmonic of the untwisted pair.
The first-arrow calculation establishes the nonvanishing
before the comparison period is constructed.
Exhaust the vacuum series \(R_j^W\) by finite unions
of double cosets of the standard compact away from the old level.
Each partial Hecke sum preserves cuspidality and that compact level.
It commutes with the old residue projector and the real line
projector, so it sends \(Y\) to itself. Since right translation is
unitary on the Petersson space, the norm of an omitted sum on \(y\)
is at most \(\|y\|_2\sum_{\gamma\ {\rm omitted}}|a(\gamma)|\).
The absolute sum converges by Section~\ref{sec:reproduction-summability}.
Thus \(R_j^W y\) belongs to the finite-dimensional, closed space
\(Y\). The finite vacuum estimate gives injectivity of \(R_j^W\).
Thus \(L_Yy=0\) implies \(R_j^Wy=0\) by applying the
displayed identity with \(x=R_j^W y\), and then \(y=0\).
This deduction uses only the two spaces of the single lift.

Let \(\mathcal K_g\in\mathcal A(V_j)\) denote the vector
theta kernel at \(g\in U(W)(\mathbb A)\).
For a Weyl coordinate vector \(\zeta\), its value at the
finite target point \(c\) is characterized by
\[
 \langle \zeta,\mathcal K_g(c)\rangle_{\rm W}
   =\theta_{P_{\bar\zeta},F_{\rm bl}}(g,(1,c)).
\]
Here \(P_\xi\) is the mixed Cartan polynomial,
\(F_{\rm bl}\) is the original untwisted block datum,
and the characters are
\(\upsilon^{m-j}\) on the target determinant and
\(\mu_j^0\upsilon^{d_{\rm tgt}/2}\) on the source determinant
in the inverse-character convention of the Weil action.

In particular \(P_{\bar\zeta}\) is linear in \(\bar\zeta\).
The right reduction actions on the block test and the
Weyl action on the polynomial give the stated level and
rational covariance of \(\mathcal K_g\).
Explicitly the simultaneous integral source and target
lifts preserve the inclusion reduction, and the target
rational change is the polynomial's Weyl substitution.
Moving the source principal level alone fixes the
residue block, and the tame character is trivial on
its determinant.

Let \(C\) run over the finite target classes at the
original level and \(c_C\) be their representatives.
The integral and the adjoint are
\[
 \begin{split}
 Ly&=\int_{[U(W)]}y(g)\mathcal K_g\,dg,\\
 (L^*\Xi)(g)
   &=\sum_C\frac{\langle \mathcal K_g(c_C),\Xi(c_C)
                    \rangle_{\rm W}}
                      {|\operatorname{Aut}(c_C)|}\\
   &=\sum_C\frac{\overline{
           \theta_{P_{\overline{\Xi(c_C)}},F_{\rm bl}}(g,(1,c_C))}}
                      {|\operatorname{Aut}(c_C)|}.
 \end{split}
\]
Indeed the scalar expression in the last two lines is the
normalized descent at the untwisted block datum.
Equation~\eqref{eq:original-block-descent-decay} applies at the
same old inclusion place to each \(c_C\) and each of its
polynomial coordinates. Their finite sum is rapid and belongs
to \(L^2([U(W)])\). The compact target's highest-line projection cuts out its irreducible
source partner as in Section~\ref{sec:bottom-descent-real-type}.
The joint harmonic lies on the highest line of its lowest source
compact type, which occurs once. Its source line projector therefore
fixes the harmonic. By Cartan equivariance, each
\(P_{\bar\zeta}\) lies in the span of compact target translates of
this harmonic. Those translates commute with the source projector.
Write \(p_S=\int e_S(k)\omega(k,1)\,dk\) for this oscillator
projector. At fixed \(g,c_C,F_{\rm bl}\) the weighted Fock estimate
of Section~\ref{sec:four-unitary-constant-terms} and the coefficient
moment there allow theta evaluation through this integral, as in
Section~\ref{sec:four-unitary-real-projections}. Thus
\[
 R(e_S)\theta_{P_{\bar\zeta},F_{\rm bl}}
  =\theta_{p_S P_{\bar\zeta},F_{\rm bl}}
  =\theta_{P_{\bar\zeta},F_{\rm bl}} .
\]
Conjugating gives \(R(\overline{e_S})\overline{\theta}
=\overline{\theta}\), the original source line in the adjoint. Here
\(\overline{e_S(k)}=e_S(k^{-1})\) in the unitary normalization.
Together with the principal-level
covariance above, this places the class sum in \(Y^{\rm pr}\).
At each of the finitely many target classes the kernel coordinates
are in \(L^2([U(W)])\) by the same decay estimate.
Thus the defining integral for \(L\) converges absolutely
for every \(y\in Y^{\rm pr}\), by Cauchy--Schwarz.
Pairing it with \(\Xi\) and moving the finite class sum
inside the integral proves the adjoint formula.
This pairing has no target level index, because the
original metric is the class sum with reciprocal
stabilizer weights.

It suffices to use the restriction to the projected space. The orthogonal
projection from \(Y^{\rm pr}\) to \(Y\) is the original
unitary finite projector \(e_0\), so \(L_Y^*=e_0 L^*\).
The injective map \(L_Y\) has invertible \(L_Y^*L_Y\).
Thus for every \(y\in Y\) there is \(w\in Y\) with
\[
                      y=e_0L^*Lw.
\]
Conjugating \(L^*Lw\) yields a finite combination of the
scalar kernels with their chosen Fock polynomials.
The operator obtained by conjugating \(e_0\) is a finite
linear combination of finite right translations. It
preserves the real coefficient law at every adelic
left translate. The real weight can therefore be
checked on the scalar kernels of this descent.

At \(\tau\) let \(n(c)\) be the numbering in the
descent datum, counting the positive coordinates in
their original order and the negative coordinates in
reverse. Use the order
\[
 o(c)=
   \begin{cases} -1-n(c),&c\text{ positive},\\
                    n(c),&c\text{ negative}.
   \end{cases}
\]
At \(\tau\), write \(\mu^\pm_\tau(h)=k^\pm(h)\) for the two column
multiplicities of Section~\ref{sec:bottom-descent-real-type},
and put \(\gamma_\tau=\kappa_{j,\tau}-d_{\rm tgt}/2\).
In the highest polynomial each height-\(h\) minor
uses the source columns numbered \(0,\ldots,h-1\) in
its sign block. The source diagonal action contributes
minus the number of such columns containing \(c\) for
the positive variables, and plus the corresponding
number for the negative variables. The real vacuum adds
\(d_{\rm tgt}\) on negative source columns, one for each definite
target space, and zero on positive ones. The inverse
determinant character contributes the common shift
\(-\mathrm{wt}(\mu_j^0\upsilon^{d_{\rm tgt}/2})_\tau=\gamma_\tau\).

The resulting weight is
\[
 u_\tau(c)=
   \begin{cases}
       \displaystyle\gamma_\tau-
              \sum_{n(c)<h<d_{\rm tgt}}\mu^+_\tau(h),
                                 &c\text{ positive},\\[3pt]
       \displaystyle\gamma_\tau+d_{\rm tgt}+
              \sum_{n(c)<h<d_{\rm tgt}}\mu^-_\tau(h),
                                 &c\text{ negative}.
   \end{cases}
\]
The compact highest weight calculation on the same
minors identifies it with the lowest compact weight of
the negated doubled source parameter. Equivalently
writing \(\lambda_{W,\tau}\) for the doubled source parameter, the formula \(2u_\tau=-\lambda_{W,\tau}+
2\rho_{\rm n}-2\rho_{\rm c}\) gives the same entries:
successive differences within each sign block are the
two displayed staircase differences and the vacuum and
character shifts give their initial constants.

The strictly upper triangular root operators for \(o\)
in one sign block kill the joint highest polynomial
by the minor rule. In different sign blocks an upper
triangular operator is in the positive-to-negative
position of this order. Its oscillator action is a sum
of mixed second derivatives in a fixed target space. The
positive and negative minor diagrams use disjoint
target space sets, so every summand is zero. Thus every
strictly upper triangular operator kills the highest
polynomial and the diagonal acts by \(u_\tau\).

For completeness, this infinitesimal statement gives
the holomorphic law of the actual matrix coefficients.
Let \(G_\tau=U(J_\tau)\) and let \(B_\tau\) be the
complex invertible upper triangular group for \(o\).
Their product \(G_\tau B_\tau\) is open, since
\(\mathfrak g_\tau+\mathfrak b_\tau=\mathfrak{gl}_{r}(\mathbb C)\)
as real spaces. The intersection \(G_\tau\cap B_\tau\)
is the diagonal unitary torus, by elimination in the
ordered form with all diagonal entries nonzero.

For the highest Fock vector \(v\), extend the real
orbit map on this product by
\(\rho(gb)v=\chi_{u_\tau}(b)\rho(g)v\).
The torus weight makes the definition independent of
the decomposition. Its real differential at \(1\),
on a tangent \(Z=Y+X\) with
\(Y\in\mathfrak g_\tau\), \(X\in\mathfrak b_\tau\), is
\[
 d\rho(Y)v+d\chi_{u_\tau}(X)v
                  =d\rho_{\mathbb C}(Z)v ,
\]
which is complex linear. The equality uses upper
triangular annihilation and the diagonal weight.
Left translation by \(G_\tau\) and right translation
by \(B_\tau\) transport this calculation throughout
the open product.

Fix an adelic point \(x\), a target class \(C\) and the
finite datum. The resulting theta functional \(T_{x,C}\) is continuous
on the real Schwartz space by the rational lattice-sum estimate.
The real oscillator orbit of the polynomial times its Gaussian is
smooth in this topology. Indeed, near the identity the symplectic
action factors into linear substitutions, multiplication by quadratic
phases and their Fourier conjugates. Parameter differentiation acts
continuously between Schwartz seminorms, and translation gives the
same conclusion near each real group element. Consequently the scalar
function on the product of the open sets \(G_\tau B_\tau\) given by
\[
       gb\longmapsto\chi_u(b)T_{x,C}(\rho(g)v)
\]
is smooth by local smooth sections of multiplication. Its differential
is complex linear by the preceding identity, hence it is holomorphic. The target compact
group commutes with the source action, and its translates
of the target highest vector span the irreducible Weyl
space. The mixed Cartan polynomial map intertwines this
target action, with conjugation of target coordinates in
the displayed descent formula. It follows that the
holomorphic law holds for each prescribed polynomial in the
descent. Taking the finite class sum and applying the finite right
translations in the conjugate projector gives the law for
\(h\mapsto\overline{y(x(h,1))}\) at every adelic \(x\).

Finally, put \(p_-=m-j\), \(d_-=t-9m-1\). The
complementary kernel weight is
\[
 u^\circ_\tau(c)=
  \begin{cases}
    j-m/2,&c\text{ positive},\\
    \displaystyle j-m/2+d_-+
       \sum_{n(c)<h<d_-}\mu^-_\tau(h),&c\text{ negative}.
  \end{cases}
\]
For \(p_->0\), reading the stem at either side
of its \(j\)-th gap gives the last negative column count
\[
          \mu^-_\tau(p_-)
               =\Delta/2-\lfloor t/2\rfloor-1\ge300n_* .
\]
In fact the two stem entries are
\(\Delta/2-\kappa_{j,\tau}\) and
\(-\Delta/2-\kappa_{j,\tau}\), while the truncating
endpoint for the negative count is
\(-\lfloor t/2\rfloor-\kappa_{j,\tau}\).
With \(t=5r\), \(n_*=2ft^3\), \(\Delta=2r^8\),
\(f\le x^{12}\) and \(r\ge x^5\), one has
\[
 \frac{300n_*+t}{\Delta/2}
 =\frac{75000f}{r^5}+\frac{5}{r^7}
 \le 75000x^{-13}+5x^{-35}<1
\]
for sufficiently large \(x\). Since \(\lfloor t/2\rfloor+1\le t\),
this proves the asserted gap. Every negative numbering is smaller than \(p_-<d_-\),
so the count is included at every negative index.
The two weights are dominant on each sign block
because the counts are nonnegative and the positive
numbering is reversed in \(o\). They both have the
required gap \(300n_*\) between a negative and a positive
entry. When \(p_-=0\) this gap condition is empty.
The doubled source parameter has entries of modulus
at most \(2\kappa_{m+1,\tau}\); the signed sum defining
\(\rho_{\rm n}-\rho_{\rm c}\) contributes at most \(r\)
in each index. Thus
\[
                  |u_\tau(c)|\le2\kappa_{m+1,\tau}+r
                                                     \le 5b .
\]
These are the weights and estimates used in the two
peak and circle calculations.

\section{Injectivity of the common period}
\label{sec:common-period-injectivity}

\sectionopening{In this section, we use the reproducing coefficients of
Section~\ref{sec:period-reproduction} to show that the common period of
Section~\ref{sec:gluing-three-exchanges} detects every nonzero vector
in the distinguished source subspace. The circle mean at a peak bounds
the nonzero reference coefficient by the approximation error.}

\subsection{Reproduction and separation at a peak}
\label{sec:reproduction-and-separation}

\sectionopening{In this subsection, we apply the period identity of Proposition~\ref{prop:common-period-adjoint} at a peak of the descent test.}

Fix the \(j\)-th descent pair, put \(W=W_j\), and
let \(V^\circ\) be the positive definite complementary
space of rank \(t-9m-1\). All integrals on the adelic
groups below use the original quotient measures.
Put \(\mathfrak m=(\pi_M^2)\).
For a class \(c\) in \(U(V^\circ)\) at this level take
its prescribed representative \(c_f\).
For a tuple \(\zeta=(\zeta_c)_c\) of Weyl vectors of the
negative diagram, write \(K^\circ_\zeta(g)\) for the
sum over \(c\) of the kernel of \((W,V^\circ)\) at
\((g,c_f)\). Its finite test is the indicator of the
inclusion of \(W\) modulo \(\mathfrak m\).
Its polynomial is the evaluation in the negative coordinates of
\(\overline{\zeta_c}\) for the original numbering of
\(W\). Its splitting characters are \(\upsilon^{m-j}\) and \(\upsilon^{m/2-j}\).
This pairing is conjugate-linear in \(\zeta\).

Write \(F_i\) for the twist-free
first-factor functions of Section~\ref{sec:gluing-three-exchanges}, with the choices specified
below, and \(P\) for the common period operator.
The seesaw identity gives
\[
 \langle\eta,P y\rangle=
 \sum_i C_i^{\mathrm{per}}\int_{[U(W)]} y(g) F_i(g)
                        K^\circ_{\eta_i}(g)\,d\bar g,
 \qquad C_i^{\mathrm{per}}=C_i|\bar J^1|c_{\mathrm{first}} .
\]
Here \(C_i=M^{l_1}c_{L_i}\) is the scalar from
Section~\ref{sec:v2-measured-gluing-adjoint}, so \(C_i^{\mathrm{per}}\ne0\). Here \(c_{\mathrm{first}}\) is the
nonzero multiplier \(c''\) chosen
for the integral first-factor form, and the divided-lattice
Weil scalar \(c_{L_i}\) of Section~\ref{sec:v2-measured-gluing-adjoint} is nonzero.
The separate \(c'_\vartheta\) is the scalar whose
nonvanishing follows from the injectivity of \(P\) below.
Taking \(\eta\) supported in one choice yields
\[
  Py=0\quad\Longrightarrow\quad
   \int_{[U(W)]} y(g)F_i(g)K^\circ_\zeta(g)\,d\bar g=0
                   \quad\text{for every }i,\zeta .
 \tag{CP1}\label{eq:display-cp1}
\]
The class stabilizers at \(J(\pi_M^2)\) are trivial by
Lemma~\ref{lem:marked-source-freeness} at an old odd split place.
Thus a single class component at this original level may be prescribed
arbitrarily and has reciprocal stabilizer weight one.

\subsubsection{Real reproduction on the complementary block}\label{sec:peak-real-reproduction}

\sectionopening{The sharp-level inner product formula of
Section~\ref{sec:reproduction-real-measures} expresses convolution
against the complementary real coefficient in terms of the original
class tests. Section~\ref{sec:peak-period-injectivity} applies the
circle mean to the resulting identity at a peak.}

Let \(M_{u^\circ}\) be the holomorphic quotient of
principal minors for the complementary kernel weight
\(u^\circ\) on \(W\). The parameter inequalities give
\(2r<\dim V^\circ\) and the required holomorphic gap
for \(M_{u^\circ}\in L^1\).
The complementary block reproducing formula says that
for each adelic \(x\), there are \(a_{c,b}(x)\) for which
\[
 \int_{U(W)(\mathbb R)} A(xh)M_{u^\circ}(h)\,dh
  =\sum_{c,b} a_{c,b}(x)
          \int_{[U(W)]} A(g) K^\circ_{\zeta(c,b)}(g)\,d\bar g .
 \tag{CP2}\label{eq:display-cp2}
\]
Here \(A\) can be any rapidly decreasing automorphic
function of level \(\mathfrak m\).
Let \(\mathbf b(X)\) be the column of negative tableau
polynomials, \(G_-\) their Weyl Gram matrix in the coordinates of
Section~\ref{sec:marked-comparison-period-vanishing}, and \(e_b\)
the \(b\)-th coordinate column. The test polynomial at \(c\) is
\(\mathbf b(X)^{\mathsf t}G_-\overline{\zeta_c}\). Thus take
\[
 \zeta(c,b)_{c'}=\mathbf1_{\{c'=c\}}\,
                       \overline{G_-^{-1}e_b}.
\]
This is a vector in the complex tuple space at the original level.
Its test is the \(b\)-th tableau polynomial since
\(G_-\overline{\zeta(c,b)_c}=e_b\). When \(j=m\), the diagram is
empty and the sole polynomial is \(1\). The definite
Fock translates of the highest-weight vector in the negative Weyl space
lie in this tableau span at each class. For completeness,
the reproducing formula follows by viewing the left
side and the finitely many right-side integrals as
linear functionals on the rapid space. Vanishing of
the latter implies vanishing of the highest-vector
theta lift, by the tableau span and the class
description at \(\mathfrak m\). Apply
Section~\ref{sec:reproduction-real-measures} with \(y=A\) and the
complementary highest datum \(f\). Since \(L_fA=0\) and
\(k_\infty=\|f_\infty\|^2 M_{u^\circ}\) with \(f_\infty\ne0\),
the convolution on the left vanishes at every \(x\). For each fixed
\(x\) its linear functional therefore factors through the image of
the finite class-evaluation map. Extending the induced functional
on that image to the full coordinate space gives the coefficients
\(a_{c,b}(x)\) in \eqref{eq:display-cp2}.

Let \(y\) belong to the projected space and suppose \(Py=0\).
Each \(yF_i\) is rapidly decreasing of level \(\mathfrak m\).
By \eqref{eq:display-cp1} its class pairings in
\eqref{eq:display-cp2} vanish. Thus for any constants \(a_i\),
\[
 \int_{U(W)(\mathbb R)}
     y(xh)\Bigl(\sum_i a_iF_i(xh)\Bigr)
                   M_{u^\circ}(h)\,dh=0 .
 \tag{CP3}\label{eq:display-cp3}
\]
Indeed \(y\) is fixed by the old principal level, which
contains \(J_W(\pi_M^2)\). In the expression for \(F_i\) below,
an element of \(J_W(\pi_M^2)\) acts on the target block of the
split residue matrix by the identity modulo \(\pi_M^2\). It
therefore fixes every frequency test \(\psi_\Lambda\).
The target determinant character \(\upsilon^m\) is trivial on
integral units, and elsewhere the integral compact fixes the
vacuum. This proves the asserted level for every \(i\).
Write \(H_S,H_W\) for the two reduction heights. For each fixed
class and polynomial in \(F_i\), the estimate of
Section~\ref{sec:four-unitary-constant-terms} bounds the kernel at
\((s;(g,1,c_{h^\flat,f}))\) by
\(B_i H_S(s)^{A_i}H_W(g)^{A_i}\). Indeed the height of the padded
target is the maximum of the heights of its orthogonal summands.
The fixed definite class contributes a constant. Thus
\[
 |F_i(g)|\le B_i H_W(g)^{A_i}\operatorname{vol}([U(S)])
                   \sup_s H_S(s)^{A_i}|\phi(s)|.
\]
The supremum is finite by Section~\ref{sec:marked-period-carrier-peak}.
For the projected \(y=e_0L^*Lw\) of
Section~\ref{sec:marked-adjoint-good-holomorphic},
\eqref{eq:original-block-descent-decay} gives rapid decay of
\(L^*Lw\). The finite integral translations in \(e_0\) preserve the
moved lattice and its height. Consequently \(yF_i\) is rapid for
each \(i\), as required in \eqref{eq:display-cp3}.

\subsubsection{The finite first-factor choices}\label{sec:peak-finite-first-factor}

\sectionopening{In this subsection, we choose the finite test for the real reproduction of Section~\ref{sec:peak-real-reproduction}.}

Let \(S\) have signature \((2m,m)\),
\(\widetilde W=W\oplus E_-^j\oplus E_+^{3m-j}\),
and \(W^\flat=E_+^{6m+1}\), with the prescribed integral lattices. The auxiliary form \(\phi\) is the theta lift
of \(1\) from \(E_-^m\) to \(S\), with trivial
finite modulus, finite vacuum, splitting characters
\((\upsilon^m,1)\), and
\(\prod_\tau\det(y^-_\tau)^{2m}\) as Fock polynomial.
It is rapidly decreasing of full finite level and
holomorphic of weight \((0;3m)\) on \(S\).
Its norm can be computed on the constant source.
Write \(G_0=U(E_-^m)\) and \(v\) for the Fock vector of its
determinant polynomial. The determinant factors cancel on
\(G_0(F_\infty)\), as computed in
Section~\ref{sec:marked-period-carrier-peak}, so the real coefficient is
\[
 \|v\|^2=\left(\prod_{a=1}^{m}\frac{(2m+a-1)!}{(a-1)!}\right)^f .
\]
Indeed the Gaussian coordinates outside the determinant have total
mass one. In each remaining \(m\)-by-\(m\) matrix of independent
standard complex Gaussians, Gram--Schmidt on the columns writes
\(|\det Y^-|^2\) as the product of squared residual lengths.
If the span of the preceding columns has complex codimension \(a\),
the conditional squared residual length is the sum of \(a\) independent exponentials of
mean one. Its \(2m\)-th moment is \((2m+a-1)!/(a-1)!\).
Successive conditional expectations give the formula at each place.
In each finite block \(B\) of
Section~\ref{sec:reproduction-summability}, right translation fixes
\(1\). Thus \(\mathcal V_B1=\lambda_B1\), where
\[
 \lambda_B=\int_{G_B} a(\gamma)\,d\gamma,\qquad
 |\lambda_B-1|\le\tfrac12 .
\]
The vacuum coefficient satisfies
\(a(\gamma^{-1})=\overline{a(\gamma)}\) by unitarity.
Inversion preserves Haar on \(G_B\), hence
\(\lambda_B\) is real and \(\lambda_B\ge\tfrac12\).
The finitely many blocks include the tail. At unit modulus
\(c_{(1)}=1\), and the sharp inner-product formula gives
\[
 \|\phi\|^2
  =\operatorname{vol}([U(S)])\operatorname{vol}([G_0])
    \operatorname{vol}(G_0(F_\infty))\,
    \|v\|^2\prod_B\lambda_B>0 .
\]

Use the index \(i=(\Lambda,h^\flat,\sigma)\) of
Section~\ref{sec:gluing-real-matched-factors}, with frequency
\(\Lambda\) modulo \(\mathfrak m\), class \(h^\flat\) of \(W^\flat\)
at that level and row choices \(\sigma\) from the short basis of
that class. Write \(c_{h^\flat,f}\) for its fixed finite representative.
Define \(P_j\) on \((S,\widetilde W)\) as the product of the full negative
minor and the positive minors, with exponents read
from the positive diagram of the untwisted pair. The chosen covectors
and the minors are specified in Section~\ref{sec:peak-reference-coefficient}.
Let \(Q_\sigma\) be the paired-determinant polynomial
of the prescribed short covectors in the positive letters
of \((S,W^\flat)\).
Then \(F_i(g)\) is the integral on \([U(S)]\) of
\(\phi\) times the theta kernel of
\((S,\widetilde W\oplus W^\flat)\), at
\((g,1,c_{h^\flat,f})\) in the second variable,
with polynomial \(P_jQ_\sigma\), test
\(\psi_\Lambda(B)=\psi_{\mathfrak m}(\operatorname{tr}(B\Lambda))\), and characters
\((\upsilon^m,\chi_{\widetilde W,j})\). Here
\[
 \chi_{\widetilde W,j}=\chi_0\upsilon^{\,m+e^{(1)}_j-j\Delta},
 \qquad e^{(1)}_j=j+\lfloor(t+k-1)/2\rfloor .
\]
Thus the second character has trivial tame part and weight
\(m+e^{(1)}_j-\kappa_{j,\tau}\) at \(\tau\).
On \(R=\mathcal O_E/\mathfrak m\) we use the additive character
\[
 \psi_{\mathfrak m}(x\bmod\mathfrak m)
 =\psi_f(j_{\rm W}\bar\alpha_{\rm W}\alpha_{\rm W}\, e_{\mathfrak m}s_0^{-1}x),
 \qquad x\in\widehat{\mathcal O}_E ,
\]
where \(s_0\) generates \(\mathfrak m\) in the finite adeles and
\(e_{\mathfrak m}\) is the idempotent of its chosen split places
as in Section~\ref{sec:overlap-finite-character}.
By Section~\ref{sec:period-finite-coordinates},
\(j_{\rm W}\bar\alpha_{\rm W}\alpha_{\rm W}\widehat{\mathcal O}_E
=\mathfrak D_{E/\mathbf Q}^{-1}\widehat{\mathcal O}_E\).
Thus at each selected split member the annihilator of the integers
under \(\psi_f(j_{\rm W}\bar\alpha_{\rm W}\alpha_{\rm W}\,\cdot)\)
is precisely the integers. After division by \(s_0\) its annihilator
on the integers is \(\mathfrak m\), giving a perfect pairing
on the full residue ring at exponent two.

Let \(\theta_{\Psi,P_j}\) be the \((S,\widetilde W)\)
kernel with any function \(\Psi\) modulo
\(\mathfrak m\) and the same two characters.
For \((S,W^\flat)\) write \(X_{\tau,+}\) for the
\((6m+1)\)-by-\(2m\) matrix of positive source columns, and \(e_+\)
for the restriction of the chosen inclusion to those columns.
Take the Fock polynomial
\[
  Q_0(X)=\prod_\tau
   \det\!\left(((X_{\tau,+})_{e_+(i),j})_{1\le i,j\le2m}\right)^{\Delta'},
   \qquad \Delta'=(\Delta-t-1)/2 .
\]
It takes the value \(1\) when the indicated square blocks are
identities. A raising operation within the positive block
repeats a column. The negative columns do not occur, so the
raising operators in the negative block and the mixed second
derivatives also annihilate \(Q_0\). Each positive column has
degree \(\Delta'\) and contributes weight \(-\Delta'\) under
inverse source substitution. The vacuum contributes \(6m+1\)
on each negative column, and the source determinant character is
trivial. Thus the Fock vector is holomorphic of weight
\(u_1=(-\Delta';6m+1)\). The rank inequality is
\(2\dim S=6m<6m+1=\dim W^\flat\).
The classwise span of its compact target translates is
verified below.
Apply this reproduction to the rapidly decreasing
function
\[
       x\longmapsto \phi(x)\theta_{\Psi,P_j}(x,(g,1)).
\]
Here \(g\) is fixed. The same Schwartz seminorm bound gives
polynomial growth in \(x\), which the rapid decay of \(\phi\)
absorbs. At an old split place \(J_S(\pi_M^2)\) acts on the
matrix by right multiplication by the identity modulo \(\pi_M^2\).
Hence it fixes every \(\Psi\). Its determinant has trivial
\(\chi_{\widetilde W,j}\)-factor, and the finite vacuum is
integral-compact invariant elsewhere. These source actions
commute with the fixed target point, and \(\phi\) has full finite
level. This gives precisely the source level of flat reproduction.
The resulting kernels on \((S,W^\flat)\) use
the indicator \(\delta_{\mathrm{incl}}\) of the
inclusion modulo \(\mathfrak m\) and
polynomials \(Q_\sigma\).

The product of the two kernels is the kernel of the
orthogonal sum, with test
\(\Psi(B_1)\delta_{\mathrm{incl}}(B_2)\) and polynomial
\(P_jQ_\sigma\). Finite Fourier inversion writes
that product test as a linear combination of
the \(\psi_\Lambda\) on the
whole residue matrix space, with its original
primitive additive character. Therefore, for each
\(x_\phi\in U(S)(\mathbb A)\) and each \(\Psi\), the function
of \(g\)
\[
 \int_{U(S)(\mathbb R)}
 \phi(x_\phi h')\,\theta_{\Psi,P_j}(x_\phi h',(g,1))
                        M_{u_1}(h')\,dh'
 \tag{CP4}\label{eq:display-cp4}
\]
belongs to the span of the actual \(F_i\).
The reproduction coefficients depend on \(x_\phi\)
and the finite Fourier coefficients on \(\Psi\).
All quotient measures and class representatives are
as above.

The class test used in both reproductions can be checked
before taking any finite Fourier expansion. Write
\(K_T=K_T(\pi_M^2)\) and let \(c_f\) be the chosen
representative of a target class. Every adelic target
point \(a\) can be written
\[
 a=\gamma(1,c_f)(1,k)(h,1),\qquad
 k\in K_T,\quad h=(\gamma^{-1}a)_\infty,
\]
using the finite class equality
\(a_f=\gamma_f c_f k\). For a Fock vector \(v\) fix the
inclusion indicator and denote the value of its
lift at \((1,c_f)\) by \(L_c(v)\).

The lift at \(a\) is
\[
                  L_c\bigl(\omega_\infty(1,h)v\bigr).
\]
Indeed rational left invariance removes \(\gamma\).
Covariance puts the two right translations into the Weil
datum. At the chosen members over \(M\), left multiplication by
\(k_w^{-1}\equiv1\bmod\pi_M^2\) fixes the inclusion residue. Away
from these pairs \(k\) preserves the vacuum. The character on the
target determinant is \(\upsilon^{m-j}\) for the complementary
block and \(\upsilon^m\) for the flat block. Both are trivial on
the determinants of these integral matrices. The real translation
therefore acts on the Fock vector alone.
In the complementary block, writing the translated
polynomial in the negative coordinates in the tableau basis uses the
same \(c_f\), and the tuple supported at \(c\) must take
the conjugated Weyl vector as in \eqref{eq:display-cp2}.

In the reproduction of the auxiliary form the fixed matrix of short covectors
\(A_c\) is the transpose of the chosen rational short
basis of the lattice at this same \(c_f\). It is invertible
at every infinite embedding. For a list \(\sigma\) of \(2m\)
row indices \(I_{\tau,\ell}\) at each \(\tau\) and
\(1\le\ell\le\Delta'\), the corresponding polynomial is
\[
 Q_\sigma(X)=\prod_\tau\prod_{\ell=1}^{\Delta'}
    \det\big(\tau(A_c)_{I_{\tau,\ell},\,*}\,X_{\tau,+}\big).
\]
At \(\tau\) the compact target translates the selected rows of
\(Q_0\) by a linear substitution, up to its determinant scalar.
Express the substituted rows in the basis of rows of \(\tau(A_c)\).
Multilinearity of each determinant, followed by multiplication over
\(\ell\) and \(\tau\), puts every translate in the span of these
\(Q_\sigma\). Hence
vanishing of all finite class tests kills the lift of the
highest vector at every \(a\) in either case.

\subsubsection{Approximation using a peak of the auxiliary form}\label{sec:peak-kernel-approximation}

\sectionopening{In this subsection, we bound the error after applying the test of Section~\ref{sec:peak-finite-first-factor}.}

Put \(n_*=2ft^3\) and \(b\) equal to the large weight scale.
Fix a peak \(x_\phi\) of the nonzero bounded auxiliary form \(\phi\),
so \(0<\|\phi\|_\infty<2|\phi(x_\phi)|\). At this point the
first minimum is at least \((5b\sqrt b)^{-1/2}\)
by the \hyperref[sec:marked-period-peak-minima]{bounds for the first minima}.
If \(\lambda_1(W,x)\ge(5b^{24})^{-1/2}\), the lattice at
\((x,1)\) is the orthogonal sum of the lattice at \(x\) and the two
definite identity lattices. For a nonzero vector in a definite
identity lattice, its signed Hermitian square is a nonzero integral
element of \(F_0\). At least one conjugate has absolute value at
least one, so the Euclidean length is at least one.
Thus adjoining the two summands preserves the first minimum cut off
at one, and
\[
 \lambda_1(S,x_\phi)\lambda_1(\widetilde W,(x,1))
 \ge (5b^{51/4})^{-1}\ge b^{-13}.
\]

We write \((x,1)\) for the block embedding of the full adelic
point \(x\) with identity on the added definite summands.
Let \(c_P(h',h)\) be the real Fock coefficient of
\(P_j\) on \((S,\widetilde W)\), with \(h\) embedded
by the identity on the adjoined summands. On this pair and
these points the sparse kernel estimate gives a
test \(\Psi_x\) modulo \(\mathfrak m\) with
\[
 \left|
   \theta_{\Psi_x,P_j}(x_\phi h',(xh,1))-c_P(h',h)
 \right|
 \le b^{-20n_*} |h'|^{6n_*}|h|^{6n_*}|c_P(1,1)|,
 \tag{CP5}\label{eq:display-cp5}
\]
where the absolute values of real points denote
the maximum of the operator norms.
Here \(\Psi_x\) is the probability average of the additive
characters modulo \(\mathfrak m\) weighted by the conjugates
of their original kernel values at \((x_\phi,(x,1))\).
Section~\ref{sec:overlap-sparse-error} proves the overlap for this
average, with zero term \(c_P(h',h)\).
In the overlap sum the surviving rational differences lie in the
kernel of reduction on the moved finite lattice from
Section~\ref{sec:peak-hom-lattice}. Multiplication by
\(q_1/\pi_M^2\) takes these differences into the integral Hom lattice,
since the only possible denominators of \(\alpha_{\rm W}\) divide
\(D_{\rm rel}\mid q_1\) by Section~\ref{sec:period-finite-coordinates}.
Apply the bound there
and the error estimate of Section~\ref{sec:overlap-sparse-error}
with the minima above, \(q=2n_*\), complex dimension
\(N=f(3m)(r+3m)\le ft^2<q\), and the homogeneous polynomial
\(P_j\) of degree \(d_P\le b\sqrt b\) constructed below.
The block extension of \(h\) has the same maximum operator
norm, since \(|h|\ge1\). Thus the resulting bound is
\eqref{eq:display-cp5} for every \(h',h\) with the same \(\Psi_x\).

Define one reference function
\[
 F_*(h)=\int_{U(S)(\mathbb R)}
       \phi(x_\phi h')c_P(h',h)M_{u_1}(h')\,dh'.
\]
It is bounded by the unitary Fock coefficient bound
and the integrability of \(M_{u_1}\).
\subsubsection{The nonzero reference coefficient}\label{sec:peak-reference-coefficient}

\sectionopening{In this subsection, we use the approximation of Section~\ref{sec:peak-kernel-approximation} to construct the comparison function.}

In the numbering of \(W_j\), put \(r_+=k+j-1\).
The negative covectors used in \(\widetilde W\) are indexed by
\(0\le a<m\): for \(a<j\) take the covector \(a\) of the adjoined negative space, and for \(a\ge j\) take covector \(k+a-1\) of \(W_j\).
The latter covectors lie in the negative block, since
\(k+a-1\ge r_+\), and their indices run through the original
negative covectors. The two disjoint pieces give \(m\) distinct covectors.
For each height \(0\le h<r_+\), take the \(2m\) positive covectors
indexed by \(0\le a<2m\): use covector \(h+a\) of \(W_j\) if
\(h+a<r_+\), and covector \(a\) of the adjoined positive space otherwise.
The positive space has dimension \(3m-j\ge2m\). These covectors are distinct,
and both choices have positive sign.

Let \(Y^-_\tau\) be the matrix on the selected negative covectors
and all \(m\) negative columns of \(S\), and let \(Y_{\tau,h}\)
be the matrix on the height-\(h\) positive covectors and all \(2m\)
positive columns. Write \(n^\pm_\tau(h)\) for the number of columns of height \(h\)
in the positive and negative diagrams, respectively. Put
\[
 \begin{aligned}
  \nu_\tau&=\sum_{h<r_+}n^+_\tau(h),\\
  \mu_\tau&=\nu_\tau+\Delta',\\
  d_1&=r+9m+1,\\
  c_\tau&=\mu_\tau+d_1+m.
 \end{aligned}
\]
This is precisely the polynomial
\[
 P_j=\prod_\tau\det(Y^-_\tau)^{c_\tau}
       \prod_{h<r_+}\det(Y_{\tau,h})^{n^+_\tau(h)} .
\]
Each determinant has independent variables in all its positions
and is nonzero. The polynomial ring is a domain, so their product
is nonzero even when different minors share some variables.
Every variable in the product pairs identical signs on the target and source coordinates.
It thus has positive oscillator sign, and every contraction
involving two opposite oscillator signs kills \(P_j\).

The degree at \(\tau\) is \(mc_\tau+2m\nu_\tau\).
The determinant count on the original positive diagram satisfies
\[
 \mu_\tau+e^{(1)}_j
              =\kappa_{j,\tau}.
\]
Indeed the number of columns at height \(h\) is a successive
difference of the shifted parameter
\(\lambda_{\tau,h-1}+h\) for \(h\ge1\), with the
numbering of Section~\ref{sec:real-marked-filling}, so
\[
 \lambda_{\tau,i}=\begin{cases}
 k/2-1-i, & i<k-1,\\
 \Delta/2-\kappa_{i+2-k,\tau}, & i\ge k-1 .
 \end{cases}
\]
The counts strictly below \(r_+\) telescope to
\(\nu_\tau=k/2-(k+j-1+\Delta/2-\kappa_{j,\tau})\).
Adding the \(\Delta'\) columns at full positive height gives
\(\mu_\tau=\kappa_{j,\tau}-e^{(1)}_j\). The admissible bound
\(d_1+m\le e^{(1)}_j\) now gives
\(c_\tau\le\kappa_{j,\tau}\le2b\).
In particular \(2c_\tau\le5b\).
Since \(\nu_\tau\le c_\tau\), the total degree is at most
\(10fmb\le b\sqrt b\) on the chosen scale.
The second character on the padded pair has weight
\(e_\tau=m+e^{(1)}_j-\kappa_{j,\tau}
=m-\mu_\tau\), and the first has weight \(m\).

Here is the matrix coefficient calculation at one place.
For \(u\in U(2m,m)\), write \(A=\det u_{++}\),
\(D_-=\det u_{--}\), and \(\delta=\det u\).
Write \(u=k_1a k_2\), \(k_i\in U(2m)\times U(m)\).
The Cartan factor \(a\) has \(m\) two-dimensional boosts with
parameters \(s_l\ge0\) and \(m\) fixed positive directions.
Put \(C=\prod_{l=1}^m\cosh s_l\). Both block determinants of \(a\)
are \(C\), and \(\det a=1\).

In signed Bargmann coordinates a boost on one target row has kernel
\[
 \frac{1}{\cosh s}\exp\!\left(
 \varepsilon\tanh s\,(z_1z_2-\bar w_1\bar w_2)
 +\frac{z_1\bar w_1+z_2\bar w_2}{\cosh s}\right).
\]
Here the two variables have opposite oscillator signs and
\(\varepsilon=\pm1\) depends on the choice of sign for the off-diagonal
entry. Indeed the displayed kernel at \(s=0\) is the reproducing
kernel, and its \(s\)-derivative is
\(\varepsilon(z_1z_2-\partial_{z_1}\partial_{z_2})\) applied to it,
the infinitesimal oscillator action of the boost. The Gaussian is
in the Bargmann measure in this formula.
The polynomial \(P_j\) uses only positive oscillator variables.
In pairing the kernel with \(P_j\) at both ends, every term
containing \(z_1z_2\) or \(\bar w_1\bar w_2\) therefore vanishes
by monomial orthogonality. The remaining exponential substitutes
the scalar \((\cosh s_l)^{-1}\) on both participating columns.
Each positive determinant uses every positive source column once,
and each negative determinant uses every negative source column
once. Their substitution scalar is \(C^{-\nu_\tau-c_\tau}\).
The vacuum contributes \(C^{-1}\) for each of the \(r+3m\)
target rows. Thus the local normalized coefficient at \(a\) is
\(C^{-r-3m-\nu_\tau-c_\tau}\).

For \(k=(k_+,k_-)\) put \(\eta_\pm=\det(k_\pm)\).
The polynomial contributes \(\eta_+^{-\nu_\tau}\eta_-^{-c_\tau}\).
The vacuum contributes \(\eta_+^m\) from the negative target rows
and \(\eta_-^{r+2m}\) from the positive ones.
The source twist gives the compact character
\[
 \chi_P(k)=\eta_+^{m-e_\tau-\nu_\tau}
           \eta_-^{r+2m-e_\tau-c_\tau}.
\]
The inner product is linear in its second variable, so the factor
at \(k_1ak_2\) is \(\chi_P(k_1)\chi_P(k_2)\) times the factor at \(a\).
Also \(A=C\det(k_{1,+}k_{2,+})\) and
\(D_-=C\det(k_{1,-}k_{2,-})\). In particular neither minor vanishes.

It follows that the local multiplier of \(c_P(1,1)\) is
\[
   \delta^{r+3m-e_\tau}
           A^{-r-2m-\nu_\tau}D_-^{-m-c_\tau}.
\]
For a weight constant on the two blocks one has
\(m_{(a;b)}(u)=\delta^b A^{a-b}\).
The complementary identity \(A=\delta\overline{D_-}\) and
\(|\delta|=1\) therefore give, on multiplying over all places,
\[
\begin{gathered}
  c_P(u,1)=c_P(1,1)M_{u_2}(u)\overline{M_{u_3}(u)},\\
  u_2=(\Delta';\mu_\tau+r+2m),\quad
  u_3=(0;\mu_\tau+r+11m+1).
\end{gathered}
\]
At the identity, write \(P_j(z)=\sum_\alpha p_\alpha z^\alpha\)
in all signed polynomial coordinates. For \(N\) complex coordinates the
measure is \(\pi^{-N}dz\) and the two Gaussian factors give
\(e^{-\|z\|^2}\). Monomial orthogonality gives
\[
              c_P(1,1)=\|v_{P_j}\|^2=\sum_\alpha
                   \left(\prod_i\alpha_i!\right)|p_\alpha|^2>0 .
\]
Repeated variables in different minors are included in the single
coefficient \(p_\alpha\) of their common monomial.
Unitarity gives \(|c_P(u,h)|\le c_P(1,1)\) on both real groups.

The weight of the auxiliary form and the two factors satisfy
\((0;3m)+u_1+u_2=u_3\), with
\(u_1=(-\Delta';6m+1)\).
The weights \(u_1,u_2,u_3\) are dominant within blocks.
Since \(\Delta/2\ge300n_*+t\) and
\(\Delta'=\Delta/2-(t+1)/2\), the cross-block gaps of
\(u_1,u_3\) are at least \(300n_*\).
Also \(|u_{3,i}|\le c_\tau+m\le2c_\tau\le5b\).
The minor-quotient moment bound thus applies with
\((u,u')=(u_3,u_1)\) on \(S\) in the original product Haar.
In particular \(J_3=\int |M_{u_3}|^2\) is finite and positive.
In particular \(\mu_\tau=\nu_\tau+\Delta'\ge\Delta'\), so the
gap of \(u_3\) follows from \(\Delta'+r+11m+1\ge300n_*\).
Apply the holomorphic reproducing formula to the product
\(h\mapsto\phi(x_\phi h)M_{u_1}(h)M_{u_2}(h)\), of weight \(u_3\).
Its pairing with \(\overline{M_{u_3}}\) is integrable by the
unitary bound on \(c_P\), boundedness of \(\phi\), and the moment
bound. Since \(M_{u_i}(1)=1\), it gives
\[
          F_*(1)=J_3\,\phi(x_\phi)\,c_P(1,1)\ne0 .
\]
All three factors are nonzero at the fixed \(x_\phi\). This reference
value is independent of the subsequent choice of \(x\) and \(\Psi_x\).

The weight comparison gives
\(\int |M_{u_1}(h')||h'|^{6n_*}\,dh'\le b^{n_*}J_3\).
Fix \(x\) and its single test \(\Psi_x\). For every \(h\),
\eqref{eq:display-cp5} and unitarity bound the absolute integrand
of \eqref{eq:display-cp4} at \(g=xh\) by
\[
 \|\phi\|_\infty c_P(1,1)|M_{u_1}(h')|
     \bigl(1+b^{-20n_*}|h'|^{6n_*}|h|^{6n_*}\bigr).
\]
This is integrable in \(h'\), uniformly for \(h\) in compact sets.
The term \(1\) bounds the integrand defining \(F_*(h)\).
Integrating the difference with the same test for every \(h\)
gives \eqref{eq:display-cp6}, since
\[
 b^{-20n_*}\|\phi\|_\infty |c_P(1,1)|\, b^{n_*}J_3
 \le 2b^{-19n_*}|F_*(1)|.
\]

\begin{lemma}[Uniform approximation with a fixed auxiliary form]
\label{lem:carrier-uniform-approximation}
For each adelic \(x\) with \(\lambda_1(W,x)\ge(5b^{24})^{-1/2}\),
there are coefficients \(a_i(x)\), independent of \(h\), such that
\[
 \left|\sum_i a_i(x)F_i(xh)-F_*(h)\right|
    \le 2 b^{-19n_*}|F_*(1)|\,|h|^{6n_*}
                   \quad(h\in U(W)_\infty).
 \tag{CP6}\label{eq:display-cp6}
\]
The reference function \(F_*\) is fixed and \(F_*(1)\ne0\). The test
and the coefficients may depend on \(x\).
\end{lemma}

\subsubsection{Coefficients and weights of the approximants}\label{sec:peak-coefficients-weights}

\sectionopening{In this subsection, we check the weight and source of the approximants in Lemma~\ref{lem:carrier-uniform-approximation}.}

In the reproduction of the auxiliary form the kernel on \((S,W^\flat)\)
has characters \((\upsilon^m,1)\), while that on
\((S,\widetilde W)\) has \((\upsilon^m,\chi_{\widetilde W,j})\).
The product is the kernel on the target sum with
\((\upsilon^m,\chi_{\widetilde W,j})\). In fact the first
character acts separately on each target determinant, and
the two factors on the common source determinant multiply.
Let \(a_r(x_\phi)\) be the coefficients of flat reproduction
at the fixed point for the auxiliary lift, with \(r\) running over the
original class and short-covector choices. These coefficients
work for all rapidly decreasing functions of the given level. With \(R\)
the residue ring, let
\[
 f_\Psi(\Lambda)
   = \frac{1}{\#\operatorname{Mat}_{d_1\times3m}(R)}
       \sum_{B\in\operatorname{Mat}_{d_1\times3m}(R)}
          \Psi(B_1)\,\delta_{\mathrm{incl}}(B_2)
               \overline{\psi_\Lambda(B)} .
\]
Here \(B_1\) has \(r+3m\) rows from \(\widetilde W\), and
\(B_2\) has \(6m+1\) rows from \(W^\flat\). Thus \(B\) has
\(d_1=r+9m+1\) rows and \(3m\) columns. Primitivity of
\(\psi_{\mathfrak m}\) identifies its additive dual with the full space
of \(3m\)-by-\(d_1\) matrices \(\Lambda\) through
\(\operatorname{tr}(B\Lambda)=\sum_{a,b}B_{ab}\Lambda_{ba}\).
These are exactly the gluing indices of
Section~\ref{sec:marked-comparison-lattice}.
Fourier inversion and the product identity give the complex
coefficients \(a_r(x_\phi)f_\Psi(\Lambda)\) of \(F_{(\Lambda,r)}\)
in \eqref{eq:display-cp4}. Once \(\Psi=\Psi_x\) is fixed these
coefficients serve for every \(h\) on the real translate \(xh\).
Indeed the period on the left of \eqref{eq:display-cp3} is
\(\mathbf C\)-linear in the first-factor test, so its vanishing
on all \(F_i\) gives vanishing on this combination at the chosen
peak. Put \(N_1=|R|^{3m(r+3m)}\) and write \(K_t\) for the
\((S,\widetilde W)\) kernel of frequency \(t\) at
\((x_\phi,(x,1))\). At \(h'=h=1\), the definition of
\(\Psi_x=N_1^{-1}\sum_t\overline{K_t}\psi_t\) and
\eqref{eq:display-cp5} give
\[
 N_1^{-1}\sum_t |K_t|^2\le(1+b^{-20n_*})c_P(1,1).
\]
Write \(\Lambda=(t,s)\) in the two column blocks and put
\(N_2=|R|^{3m(6m+1)}\). Orthogonality in the first block and
evaluation at the inclusion in the second give
\[
 f_{\Psi_x}(t,s)
    =\frac{\overline{K_t}\,\overline{\psi_s(\mathrm{incl})}}
            {N_1N_2}.
\]
All \(N_2\) values of \(s\) occur in the original index set.
Summation and Cauchy--Schwarz give
\[
 \sum_\Lambda |f_{\Psi_x}(\Lambda)|
   =N_1^{-1}\sum_t|K_t|
   \le \bigl((1+b^{-20n_*})c_P(1,1)\bigr)^{1/2}.
\]
This bound holds at every \(x\) satisfying the minimum condition
of Lemma~\ref{lem:carrier-uniform-approximation}. The coefficient
bounds at arithmetic embeddings in
Section~\ref{sec:compression-trivial-norm} are applied to the
fixed gluing operator \(G_1\).

The scalar \(c''\) of Section~\ref{sec:marked-first-factor-bounds}
is nonzero. The untwisted identity
\(F_i=(c'_1/c'')\Xi_{1,i}\) therefore writes each \(F_i\) as a
finite combination of \((W_j,V'')\) kernels in the original
order, with polynomial in the positive coordinates of the
diagram with column counts \(n^+_\tau(h)\). The target is positive
definite of rank \(d_1\) and its splitting weight on the
source is \(e^{(1)}_j-\kappa_{j,\tau}
=-\mu_\tau\). The highest-weight transformation law is linear in the function.
Multiplying \(\Xi_{1,i}\) by \(c'_1/c''\) therefore preserves this law
also when \(c'_1=0\). The polynomial in the positive coordinates gives
the holomorphic weight with entries
\[
    \mu_\tau-\sum_{h>n(a)}
                     n^+_\tau(h)
              \quad\text{on positive }a,
    \qquad
    \mu_\tau+d_1\quad\text{on negative }a .
\]
Here the positive entries come from the staircase minors
and the negative entries from the vacuum and splitting
weights. The definite-group translates of this highest
polynomial span the Weyl space in the positive coordinates, so the same
weight applies to every coefficient in that class sum. Put \(d_-=t-9m-1\). The complementary weight at
\(\tau\) has entries
\[
 \begin{cases}
 j-m/2,&\text{on positive }a,\\
 j-m/2+d_-+
 \displaystyle\sum_{n(a)<h<d_-}
                  n^-_\tau(h),&\text{on negative }a.
 \end{cases}
\]

The full staircase shift is \(\kappa_{j,\tau}-(r+t)/2\).
The numbering \(n(a)\) is that of the bottom descent. Since
\(k=r+1-m\) and \(t+k-1\) is even, the first-factor exponent is
\[
 e^{(1)}_j=(r+t)/2+j-m/2,\qquad
 \mu_\tau=\kappa_{j,\tau}-e^{(1)}_j .
\]
The rank terms add because \(d_1+d_-=r+t\), and the positive
columns of height at least \(d_1\) and the negative columns
of height at least \(d_-\) have count zero. Thus the displayed weights add to \(u\) on the actual covector.
For each fixed peak \(x\), the coefficients \(a_i(x)\) in
\eqref{eq:display-cp6} are independent of \(h\), so
\(K_x(h)=\sum_i a_i(x)F_i(xh)\) is holomorphic of weight
\(u-u^\circ\) in the order of the circle estimate. Its reference
function \(F_*\) is bounded and has the fixed value
\(J_3\phi(x_\phi)c_P(1,1)\ne0\) at the identity.

\subsubsection{Injectivity of the common period}\label{sec:peak-period-injectivity}

\sectionopening{In this subsection, we combine Lemmas~\ref{lem:circle-mean-peak-estimate} and~\ref{lem:carrier-uniform-approximation} to prove injectivity.}

Suppose \(0\ne y\in\ker P\). The rapid form \(\bar y\) is bounded,
so choose \(z\) with \(|\bar y(z)|>\|\bar y\|_\infty/2\). The compact
circle acts by its unitary weight, hence each \(z(r_s,1)\) is also
a peak. At any peak \(x\) of \(\bar y\), \eqref{eq:display-pm5} gives
the first minimum required by \eqref{eq:display-cp6}. At \(h=1\) the chosen
coefficients satisfy
\[
 \left|\sum_i a_i(x)F_i(x)\right|
 \ge (1-2b^{-19n_*})|F_*(1)|>0 .
\]
In particular the approximant is nonzero at this peak before using
any orthogonality consequence of \(Py=0\).
Use the same coefficients in
\eqref{eq:display-cp3}. The conjugate weight and the complementary
kernel weight have the required dominance and gap, and their
difference is the weight of each \(F_i(x\,\cdot)\).
Hence Lemma~\ref{lem:circle-mean-peak-estimate}, applied at \(z\), gives
\[
  |F_*(1)|
    \le (4b^{n_*}+1)\,2 b^{-19n_*}|F_*(1)|
    < |F_*(1)| ,
\]
since \(b^{n_*}>10\).
This contradiction proves injectivity of \(P\).
The first orbital moment in Section~\ref{sec:marked-span-proof} gives
\(D>0\) for sufficiently large \(x\), and
\eqref{eq:marked-full-spherical-trace} at the identity gives
\(\dim Y_j=D\) for the projected automorphic space on \(W_j\).
Evaluating \(c'_\vartheta G_\vartheta^*u_{j,\vartheta}=P\) on any
nonzero vector of \(Y_j\) now proves \(c'_\vartheta\ne0\)
for every auxiliary twist \(\vartheta\) at this fixed level.
The inclusion \(\ker u_{j,\vartheta}\subseteq\ker P\) then proves the
injectivity of each original lift \(u_{j,\vartheta}\) on \(Y_j\).
The same \(P\) works for all these twists. The external modulus
\(q_0\) and the quadratic discriminant \(-d\) are fixed
before the large-\(x\) estimates.

\subsection{First minima at peaks of automorphic forms}
\label{sec:marked-period-peak-minima}

\sectionopening{In this subsection, we use transvections at the
prescribed level to bound the first minimum at points satisfying the
half-supremum bound in Section~\ref{sec:peak-period-injectivity}.
The estimate separates the nonzero Hom-lattice terms in
Section~\ref{sec:peak-kernel-approximation}.}

Let \(V\) have the chosen adelic basis, with Hermitian
matrix \(J=\operatorname{diag}(\epsilon_i)\) at infinity,
and order the positive coordinates before the negative
ones. Use the Euclidean norm
\(\|v\|^2=\sum_\tau\sum_i |v_{\tau,i}|^2\) in the
archimedean coordinates of a rational lattice vector.
The first minimum \(\lambda_1(V,g)\) is cut off at \(1\).
Suppose that \(F\) is a bounded automorphic function,
adelically holomorphic of weight \(u_\tau\) in the
following sense: at each adelic translate it extends
along \(h B\), with \(h\in U(J)\), \(B\) upper triangular,
and with value \(\prod_i B_{ii}^{u_{\tau,i}}\) times
the value at \(h\).
Assume
\[
 |u_{\tau,i}|\le A,\quad A\ge2,\qquad
 u_{\tau,i}+1\le u_{\tau,k}
                    \quad(\epsilon_i=1,\ \epsilon_k=-1).
\]
Let \(p=-\bar p\) be an algebraic integer with
\(0<\operatorname{Im}\tau(p)\le A'\) at each chosen
infinite place, \(A'\ge1\), and assume \(F\) is invariant
under the principal level \(K(p)\).
If \(C=\sup |F|>0\) and \(|F(g)|\ge C/2\), then
\[
                   \lambda_1(V,g)\ge (AA')^{-1/2}.
 \tag{PM1}\label{eq:display-pm1}
\]

Suppose instead that an original nonzero rational
lattice vector \(v\) in the basis moved by \(g\) has
length \(<(AA')^{-1/2}\). Its Hermitian square is
integral, because the adelic coordinates at every
finite place are integral and the finite
Hermitian matrices are integral. Its absolute value
at every archimedean embedding is at most \(\|v\|^2<1\).
The integral norm therefore forces the square to
vanish. If the form is definite at any infinite place, this
already contradicts \(v\ne0\). This includes the case \(j=m\) of
the projected source below.

In the remaining case the moved coordinates \(v_\tau\) at every
real place are nonzero and isotropic. The transvection
\[
       t=I-pvv^*J
\]
in the moved adelic basis is unitary. Its conjugate
back to the original basis is the rational
transvection defined by the same rational vector
and the rational Hermitian form of \(V\). Its finite
component is in \(K(p)\), since the moved coordinates
of \(v\) are integral at every finite place.
Writing \(t=(t_\infty,t_f)\), rational left
invariance and right \(K(p)\)-invariance give
\[
       F(g(t_\infty h,1))
        =F\bigl((gtg^{-1})\,g(h,1)(1,t_f^{-1})\bigr)
        =F(g(h,1)).
 \tag{PM2}\label{eq:display-pm2}
\]
This last equality relates the transvection to the
original automorphic function at \(g\).

Here is the triangular factorization and bound for the
complexified transvection. At a distinguished place \(\tau\)
write \(\tau(p)=i a_\tau\), \(a_\tau>0\), and put
\(\eta_\tau=a_\tau\|v_\tau\|^2>0\).
Then \(A\sum_\tau \eta_\tau\le1\), so
\(\eta_\tau\le1/2\). For \(z\) in the upper half-plane
set
\[
  X_\tau(z)=I-d_\tau(z)v_\tau v_\tau^* J,\qquad
  d_\tau(z)=i a_\tau(z-i),\qquad h_z=\operatorname{Im}z .
\]
The isotropy gives
\(X_\tau^*JX_\tau=J-2\operatorname{Re}(d_\tau)
                          Jv_\tau v_\tau^*J\).
For an initial segment of the ordered coordinates
of length \(j\), the determinant of its Gram block,
divided by the determinant of the corresponding
block of \(J\), is
\[
 p_j=1-2\operatorname{Re}(d_\tau)
                         \sum_{i\le j}\epsilon_i|v_{\tau,i}|^2 .
\]
The partial sums in this expression increase to
\(\|v_\tau\|^2/2\) on the positive block and then return
to \(0\).

The \(p_j\) are positive. They start and end at
\(1\), with value \(p_*=1+\eta_\tau(h_z-1)\) at the end
of the positive block. Thus Hermitian elimination
without pivoting gives \(X_\tau=h_\tau B_\tau\)
with \(h_\tau^*Jh_\tau=J\) and \(B_\tau\) upper
triangular. Indeed the successive elimination pivots
have values \(\epsilon_i p_i/p_{i-1}\), and
normalizing their absolute values gives
\(|B_{\tau,ii}|=(p_i/p_{i-1})^{1/2}\).

Write \(L_i=\log(p_i/p_{i-1})\) and
\(q_\tau=1+\eta_\tau h_z\). If \(h_z\ge1\), the
\(L_i\) are nonnegative on the positive block and
nonpositive on the negative block. Their total is
zero and their positive-block sum is \(\log p_*\).
The weight gap gives
\(\sum_i u_{\tau,i}L_i\le-\log p_*\).
Moreover \(p_*=q_\tau-\eta_\tau\ge(1-\eta_\tau)q_\tau\)
and \(-\log(1-\eta_\tau)\le2\eta_\tau\).
If \(0<h_z<1\), the signs of the \(L_i\) reverse,
so
\[
  \tfrac12\sum_i u_{\tau,i}L_i
  \le -A\log p_*
  \le 2A\eta_\tau(1-h_z)
  \le 2A\eta_\tau-\tfrac12\log q_\tau .
\]
The last inequality follows from
\(\log q_\tau\le\eta_\tau h_z\).
Both cases imply
\[
       \left|\prod_i B_{\tau,ii}^{u_{\tau,i}}\right|
          \le \frac{\exp(2A\eta_\tau)}
                         {\sqrt{1+\eta_\tau\operatorname{Im}z}} .
 \tag{PM3}\label{eq:display-pm3}
\]

These factorizations give one holomorphic function along the
upper half-plane. If \(hB=h'B'\), then \(h^{-1}h'=BB'^{-1}\)
is both upper triangular and \(J\)-unitary. Such a matrix is
diagonal: it preserves the successive nondegenerate coordinate
spaces and hence their orthogonal complements. Write
\(\chi_u(B)=\prod_{\tau,i}B_{\tau,ii}^{u_{\tau,i}}\).
The weight rule then gives
\[
 F(g(h,1))\chi_u(B)=F(g(h',1))\chi_u(B').
\]
Near \(X_0=h_0B_0\) choose continuous factors \(h(X),B(X)\)
by Hermitian elimination with positive diagonal moduli and a
fixed diagonal phase. Then \(XB_0^{-1}=h(X)(B(X)B_0^{-1})\)
lies in the holomorphic germ neighborhood of \(h_0\).
Right translation of that germ by \(B_0^{-1}\), multiplied by
\(\chi_u(B_0)\), agrees with the factorized value by the weight
rule. The local compositions along \(z\mapsto(X_\tau(z))_\tau\)
therefore give one holomorphic function \(f(z)\).
Equation \eqref{eq:display-pm3} and \(A\sum_\tau\eta_\tau\le1\)
give \(\sup |f|\le C e^2\), and the same bound with
the denominators in \eqref{eq:display-pm3} gives
\(f(z)\to0\) uniformly as \(\operatorname{Im}z\to\infty\).
The identity \(X(z+1)=t_\infty X(z)\) follows because
\((v_\tau v_\tau^*J)^2=0\). The triangular
character and \eqref{eq:display-pm2} give \(f(z+1)=f(z)\).
In the variable \(q=e^{2\pi i z}\), the function is
holomorphic on \(0<|q|<1\), extends by
zero at the origin, and is bounded by \(C e^2\).
Schwarz's lemma at \(z=i\) gives
\(|F(g)|=|f(i)|\le C e^{2-2\pi}<C/2\), contradicting
the peak condition. This proves \eqref{eq:display-pm1} in the moved
lattice and the original level.

\subsection{The auxiliary automorphic form and its first minimum}
\label{sec:marked-period-carrier-peak}

\sectionopening{In this subsection, we apply the first two exchanges of
Section~\ref{sec:gluing-three-exchanges} to verify the hypotheses of
\eqref{eq:display-pm1} for the auxiliary form and for the descent test.}

Put \(S=I_{2m}\perp I_m^-\), and write \(K_0=U(S)(\widehat{\mathcal O}_F)\)
for the stabilizer of the stated standard finite lattice.
In \eqref{eq:carrier-definite-class-expansion} one may fix
\(\vartheta=1\) and \(c=\chi_M^j\). At either old chosen place the
residue character \(\omega_w^j\) of \(\chi_M^j\) is nontrivial, since
its order is \(q_2/\gcd(q_2,j)\) for \(1\le j\le m<q_2\).
Thus \(\phi=F^\circ\) is a fixed finite sum of the displayed theta
kernels at the original representatives \(h_a\), with coefficients
\(\gamma_{12} A_c(h_a)/e_{3m}(h_a)\) independent of the point. For \(k\in K_0\), the integral
Weil action on the finite datum is
\[
               c(\det k)^{-1}\chi(\det \bar k)=1 ,
\]
where \(\bar k\) denotes residue reduction. Thus every summand and
\(\phi\) itself are right \(K_0\)-invariant. In particular they have
every principal level contained in \(K_0\).
The determinant character extended by zero on nonunits is supported on
invertible residue matrices.

Fix one split member of the nonempty support of \(\mathfrak m\).
The square residue matrix there is invertible on the support of the
determinant test. In the lower-kernel estimate of
Section~\ref{sec:four-unitary-constant-terms}, \(S\) is the varying
source and \(h_a\) is the fixed definite target point. Thus for a short
nonzero source vector \(v\), one has \(Yv_w\ne0\) over \(E_w\)
for every local matrix \(Y\) in the support. This supplies the
nontrivial modulation in the local additive pairing.
Decomposing the finite determinant test into its residue point
indicators and applying that estimate gives, for every \(D\ge0\),
\[
 \sup_s\operatorname{ht}_S(s)^D
   \left|\Theta^{S,I_{3m}}_{\upsilon^m,c}[\chi\circ\det,1](s,h_a)\right|
       \le C_{a,D}<\infty .
\]
Here \(c\) and the finitely many \(h_a\) are fixed once for this
presentation of \(\phi\). Their weighted sum is rapidly decreasing.

Continuity and
rational left invariance follow from the theta sum.
At infinity the Fock polynomial \(1\) is killed by the holomorphic
lowering operators. Its diagonal weight on \(S\), with weight-zero
second splitting character and positive definite rank \(3m\) on the
other side, is \(0\) on the positive block and \(3m\) on the negative
block. For the order numbering the positive block by \(-1,-2,\ldots\)
and the negative block by \(0,1,\ldots\), the upper triangular
root operators in a single block also annihilate the constant polynomial.
The Schwartz-space orbit calculation of
Section~\ref{sec:marked-adjoint-good-holomorphic} therefore gives its
triangular continuation of weight \((0;3m)\). At fixed \((x,h_a)\)
the translated theta functional is continuous on the real Schwartz
space. Applying it to this orbit map gives the continuation for
each adelic left translate. Every summand uses the same order and
weight, so \(\phi\) has these as well.

To prove nonvanishing use the original presentation of \(\phi\) as
the lift of \(1\) from \(I_m^-\) to \(S\), with characters
\((\upsilon^m,1)\). The finite datum is the integral vacuum at every
place. For the sharp convolution formula take the unit modulus,
empty prescribed support, and full integral source level. The finite
sharp condition at the empty support is vacuous, the vacuum fixes
the integral level, its mass there is positive, and the second
character is trivial. In the sharp inner-product formula the finite factor is now the
vacuum series on all finite places. Its neighbour count applies at
each conjugation orbit, so the same absolute convergence and
injectivity argument applies. The rank inequality is \(3m>2m\), and the constant source
is rapidly decreasing since \(I_m^-\) is definite at infinity and
its adelic quotient is compact.

Let \(v\) be the original Fock vector of
\(P=\prod_\tau\det(Y^-_\tau)^{2m}\). This polynomial is nonzero.
For \(u=(u_\tau)\) on the compact negative-definite source,
put \(\delta_\tau=\det u_\tau\). The \(2m\) positive target
coordinates give the vacuum factor \(\delta_\tau^{2m}\)
from the negative Hom block. On the negative target coordinates
the Fock substitution is right multiplication by
\((u_\tau^{-1})^{\mathsf T}\), so that
\(\det(Y_\tau^-)^{2m}\) acquires \(\delta_\tau^{-2m}\).
The source splitting character is trivial and the target is
evaluated at the identity. Consequently
\[
 \omega(u,1)v=v,\qquad
 k(u)=\langle v,\omega(u,1)v\rangle=\|v\|^2 .
\]
If the auxiliary lift vanished, the sharp formula and the injectivity of its
finite vacuum series would give
\[
       0=(k*1)(1)=\|v\|^2\operatorname{vol}(U(I_m^-)(F_\infty))>0 .
\]
The Haar measure is that of the original real convolution. This proves
nonvanishing on the same adelic realization as in the identity for the auxiliary lift.
In particular \(\phi\) is the single nonzero full-level holomorphic
form to which the peak minimum is applied, common to the later source space and
class tests.

Apply this first to the auxiliary form \(\phi\) of signature
\((2m,m)\). Its nonzero rapid full-level form is
holomorphic of weight \((0;3m)\). Full level
contains \(K(\sqrt{-d})\); the given CM embeddings
read \(\sqrt{-d}\) with positive imaginary part
\(\sqrt d\). The admissible scale gives \(3m\le5b\)
and \(d\le b\). Therefore \(A=5b\), \(A'=\sqrt b\)
give, at the chosen peak \(x_\phi\),
\[
                \lambda_1(S,x_\phi)\ge(5b\sqrt b)^{-1/2}.
 \tag{PM4}\label{eq:display-pm4}
\]

For a nonzero form in the projected space \(y\) on \(W_j\), use
\(F=\bar y\) on the original \(W_j\). Its selected finite type is cuspidal and the form is rapidly decreasing.
Its original principal level is \(J_M=K(\pi_M)\). At each chosen
split member \(w\mid\pi_M\), one has
\(v_w(M\sqrt{-d})=1\), since \(M\) is the product of the two
distinct old primes, both prime to \(d\).
Consequently \(K(M\sqrt{-d})\subset J_M\). This containment applies
to \(y\) and to \(\bar y\) at every peak. The exponent-two level
\(K(\pi_M^2)\) enters the reproduction of the products \(yF_i\)
in \eqref{eq:display-cp2}.

\hyperref[sec:marked-adjoint-good-holomorphic]{The adjoint calculation} gives
this holomorphic weight on the original projected space. Its
proof uses only the Gram operator on the projected source of the single lift and
the adjoint. The opposite-sign gap is at least
\(300n_*\) by the last negative column count, with the
definite case having no opposite signs.

Now \(p=M\sqrt{-d}\) is purely imaginary in the
required sense and
\[
       0<\operatorname{Im}\tau(p)=M\sqrt d
                       \le b^{22}\sqrt b\le b^{23}
\]
by the original level-prime interval and \(d\le b\).
Thus \eqref{eq:display-pm1} with \(A=5b,\ A'=b^{23}\) gives
\[
        \lambda_1(W_j,g)\ge(5b^{24})^{-1/2}
                         \quad\text{at every peak of }\bar y.
 \tag{PM5}\label{eq:display-pm5}
\]
The product of \eqref{eq:display-pm4} with \eqref{eq:display-pm5}, after adjoining the fixed definite summands to \(W_j\), is at least
\((5b^{51/4})^{-1}\ge b^{-13}\) at the chosen scale.
Each fixed definite summand at the identity has first minimum
at least \(1\), so the orthogonal sum has the same truncated
first minimum as \(W_j\). At the rotated peak \(z_s=z(r_s,1)\)
the finite lattice is the one at \(z\). The real coordinate
at \(\tau\) is multiplied by \(r_{s,\tau}^{-1}\), which preserves
its Euclidean norm. Hence \(\lambda_1(W_j,z_s)=\lambda_1(W_j,z)\)
throughout the circle of Lemma~\ref{lem:circle-mean-peak-estimate}.

\subsection{Finite overlaps and Gaussian decay}
\label{sec:peak-hom-lattice}

\sectionopening{In this subsection, the error in Section~\ref{sec:common-period-injectivity} is
summed over nonzero maps between two moved lattices. A first-minimum
estimate for this Hom lattice gives a summable Gaussian bound in the
original normalizations.}

The comparison of the first minimum on the Hom space
with the two first minima uses the original self-dual
Hermitian lattices. Let \(L\) be the source lattice in
its Euclidean real coordinates, of dimension
\(D_{\mathbb R}=2f\operatorname{rk}S\). Take \(N_c=dq_1\), so
\(E\subseteq\mathbb Q(\zeta_{N_c})\). The integer \(N_c\) belongs to the
different of \(E\). Indeed it belongs to the different
of the cyclotomic field by differentiating \(X^{N_c}-1\)
at a primitive root and using the cyclotomic factor.
For \(C=\mathbb Q(\zeta_{N_c})\), transitivity gives
\(\mathfrak D_C=\mathfrak D_{C/E}\mathfrak D_E\mathcal O_C
\subseteq\mathfrak D_E\mathcal O_C\).
Contraction to \(E\) gives \(N_c\in\mathfrak D_E\). The Euclidean dual of \(L\) is
\(2J_\infty\mathfrak D_E^{-1}L\). The trace of the
Hermitian pairing is twice the real Euclidean pairing
with the sign matrix \(J_\infty\), and the finite
Hermitian lattice is self-dual. Consequently the
Euclidean dual of \(2L/N_c\) is contained in
\(J_\infty L\), so its nonzero vectors have length at
least the cut-off first minimum \(\lambda_S\).

Banaszczyk's transference theorem \cite[Theorem~2.2, p.~632]{BanaszczykTransference}
gives covering radius
at most \(D_{\mathbb R}/(2\lambda_S)\), hence at most
\(D_{\mathbb R}^2/\lambda_S\), for \(2L/N_c\).
The vectors of a lattice with covering radius at most
\(R>0\) and length at most \(3R\) span its real space.
Otherwise a vector of length \(2R\) in their orthogonal
complement is farther than \(R\) from all those vectors,
and a longer lattice vector is also farther than \(R\).

Thus the source lattice contains an \(E\)-basis with
each norm at most \(6N_c(f\operatorname{rk}S)^2/\lambda_S\).
A nonzero integral Hom vector sends at least one of these
basis vectors to a nonzero vector of the target lattice.
If \(\lambda_T\) is the target minimum, the operator
norm bound by the Euclidean coordinate norm, even with
the weaker factor \(t^2\), gives
\[
 \lambda_{\operatorname{Hom}(S,T)}
      \ge \frac{\lambda_S\lambda_T}{6N_c(ft)^2t^2}.
\]
At the fixed adelic pair \((g,h)\), use the basis matrices
\(X_S,X_T\) of Section~\ref{sec:reproduction-doubling}. Write
\(\mathcal O_f=\widehat{\mathcal O}_E\). The rational lattices are
\[
 L_S=S(E)\cap(X_Sg)_f\mathcal O_f^{\operatorname{rk}S},\quad
 L_T=T(E)\cap(X_Th)_f\mathcal O_f^{\operatorname{rk}T}.
\]
The Euclidean realizations use \((X_Sg)_\infty^{-1}\) and
\((X_Th)_\infty^{-1}\). A rational map \(Z\colon S\to T\) has
coordinates \(z(Z)=(X_Th)^{-1}Z(X_Sg)\).
For a surviving difference, the integral
box and the zero residue condition together give
\[
 \alpha_{\rm W}^{-1}z(Z)_f=\pi_M^2 V_f,\qquad
 V_f\in\operatorname{Mat}_{\operatorname{rk}T\times\operatorname{rk}S}(\mathcal O_f).
\]
Thus \(w=(q_1/\pi_M^2)Z\) carries \(L_S\) into \(L_T\), since
its finite matrix \(z(w)_f=(q_1\alpha_{\rm W})V_f\) is integral
by Section~\ref{sec:period-finite-coordinates}. Write
\[
 \|w\|^2=\sum_{\tau\in\Phi}\sum_{i,j}|z(w)_{\tau,ij}|^2 .
\]
Since \(\pi_M\bar\pi_M=M\), every embedding \(\tau\) satisfies
\(|\tau(\pi_M^2)|=M\). The same calculation applies to the
difference of any two surviving maps, for every choice of
\(g_f,h_f\). In the chosen
CM type, \(j_{\rm W}=-1/\sqrt{-d}\) from
Section~\ref{sec:period-finite-coordinates} gives the real Fock scale
\((2\pi/\sqrt d)^{1/2}\). Therefore
\[
 \|z_{\mathbb R}\|
 =\left(\frac{2\pi}{\sqrt d}\right)^{1/2}\frac{M}{q_1}\|w\|
 \ge
 \left(\frac{2\pi}{\sqrt d}\right)^{1/2}
 \frac{M\lambda_S\lambda_T}{6q_1N_c(ft)^2t^2}
 \ge b^6 .
\]
The last inequality uses \(M\ge b^{20}\),
\(\lambda_S\lambda_T\ge b^{-13}\), and
\(6q_1N_c(ft)^2t^2\le (2\pi/\sqrt d)^{1/2}b\) for sufficiently
large \(x\). The latter follows from polynomial growth in \(x\)
at fixed \(q_0,d\) and the lower bound on \(b\) in
Part~\ref{part:parameters}.

\subsection{The normalized overlap on the adelic Hom space}
\label{sec:marked-normalized-sparse-overlap}

\sectionopening{In this subsection, we calculate the overlap of two adelic tests as a product of a
residue-ring factor and a polarized real Fock factor. This supplies
the coefficient estimates used to compare the leading term and the
error in Section~\ref{sec:common-period-injectivity}.}

Use the original pair \(S,T\), with
\(T=\widetilde W\) for the sparse test, the modulus
\(\mathfrak m=(\pi_M^2)\), and the Weil parameter
\((\psi_f,j_{\rm W},\alpha_{\rm W})\).
Put \(R=\mathcal O_E/\mathfrak m\).
Let \(D=\operatorname{rk}_E\operatorname{Hom}_E(S,T)\).
An adelic displacement \(z\) is read at the finite places
as \(a(z)=\alpha_{\rm W}^{-1}z_f\) and at infinity as
\(z_{\mathbb R}=(\sqrt{2\pi\,\operatorname{Im}
  \tau(j_{\rm W})}\,z_\tau)_\tau\).
Let \(\mathcal R(z)\) be the Heisenberg operator and
\(v_{\rm vac}\) the finite unit vacuum.

\subsubsection{The probability mean and Fourier inversion}
\label{sec:overlap-fourier-inversion}

\sectionopening{In this subsection, we express the overlap by Fourier inversion in the normalization used to bound the error in Section~\ref{sec:overlap-sparse-error}.}

In the rational symplectic displacement space choose a
rational Schr\"odinger polarization \(X\oplus X^\vee\)
over \(F_0\). Give \(X(\mathbb A)\) the Haar measure for which
\(X(F_0)\) has covolume one. The annihilator of \(X(F_0)\)
is \(X^\vee(F_0)\). Poisson summation therefore gives covolume
one on the dual quotient for the Fourier-dual Haar measure.
The theta functional in this realization
is rational summation, up to a unit scalar.

Indeed rational modulation invariance of a functional
forces its support to be on \(X(F_0)\), by adelic
self-duality. At a rational point \(\xi\) fix a finite box
\(\xi_f+L_f\), where \(L_f\) is a compact open additive subgroup.
The rational positions with finite coordinates in the box form a
lattice coset at infinity. Localize at \(\xi_\infty\) in a ball
meeting this coset only there. Rational modulations from
\(X^\vee(F_0)\cap L_f^\perp\) are constant on the finite box,
and their global values at \(\xi\) are one. This intersection is a
full lattice in the real dual under the trace pairing. Select
modulations with real linear phases \(\ell_j\) forming a basis of
that dual. In coordinates \(t\) centered at \(\xi_\infty\), the
localized distribution is killed by each \(e^{i\ell_j(t)}-1\),
hence by each \(\ell_j(t)\) since the quotient by \(\ell_j(t)\)
is smooth and nonzero at zero. Every smooth germ vanishing at zero
is a sum of the \(\ell_j\) times smooth germs, by integration along
radial segments. Thus the distribution acts by evaluation at the
point, including on finite tests since they are locally constant.
Rational translation invariance gives the same multiplier
at every rational point. Denote this common scalar by \(C\).
For a Schwartz vector \(f\) the translate is
\[
 \theta(\mathcal R(u,v)f)
  =C\sum_{\xi\in X(F_0)}
                 f(\xi+u)\psi\bigl(v(\xi+u/2)\bigr).
\]
Multiply by the conjugate of the same expression for \(f'\).
The two half-shift factors cancel. Integration over
\(v\in X^\vee(\mathbb A)/X^\vee(F_0)\) kills each term
whose two rational indices differ. Integration over
\(u\in X(\mathbb A)/X(F_0)\) then unfolds the common
rational index. Consequently, with \(dz\) the probability
Haar on the displacement quotient,
\[
 \int \theta(\mathcal R(z)f)
                 \overline{\theta(\mathcal R(z)f')}\,dz
          =|C|^2\int_{X(\mathbb A)}f(u)\overline{f'(u)}\,du
          =|C|^2\langle f',f\rangle .
\]
On compact fundamental sets the finite support confines the rational
indices to finitely many lattice cosets. The real Schwartz bounds
give absolute uniform convergence there. Taking \(f=f'\) equal to the
unit vacuum gives \(|C|^2=1\) by the prescribed mean-square normalization.

Let \(A_{f,f'}(z)\) denote the integrand. Rational
Heisenberg invariance and cancellation of the common unit
phase make it periodic. For a rational displacement
\(\gamma=(a,b)\) and \(z=(u,v)\), composing the displayed
Schr\"odinger actions gives
\[
 \mathcal R(z)\mathcal R(\gamma)
 =\psi\bigl(b(u)-v(a)\bigr)\mathcal R(\gamma)\mathcal R(z),
 \qquad \psi\bigl(b(u)-v(a)\bigr)=\beta(z,\gamma)^2 .
\]
Applying the last mean identity to
\((\mathcal R(\gamma)f,f')\), rational invariance gives
\[
 \int \beta(z,\gamma)^2 A_{f,f'}(z)\,dz
                     =\langle f',\mathcal R(\gamma)f\rangle .
\]
The inverse characters \(\beta(z,\gamma)^{-2}\) are the
dual of the quotient.

The series of these Fourier
coefficients is absolutely convergent. To see this choose
an open finite additive compact group fixing both \(f,f'\)
under its Heisenberg displacements. The product is invariant
under that group because the remaining commutator factors
have modulus one and cancel. The quotient by this group
is a real torus, by additive approximation. The product
descends to a smooth function there, and integration by
parts bounds its Fourier coefficients by every inverse
power of the frequency norm. Fourier inversion at zero
therefore yields
\[
 \theta(f)\overline{\theta(f')}
        =\sum_{\gamma\ {\rm rational}}
                        \langle f',\mathcal R(\gamma)f\rangle .
\]
In particular the zero-frequency coefficient is
\(\langle f',f\rangle\) with its original scale.

\subsubsection{The finite character and its translate}
\label{sec:overlap-finite-character}

\sectionopening{In this subsection, we compute the finite component of the transform in Section~\ref{sec:overlap-fourier-inversion}.}

Write \(s_0\) for the finite generator of \(\mathfrak m\)
in this paragraph and \(e_{\mathfrak m}\) for the idempotent
of its chosen split places. Define
\(\psi_0(x)=\psi_f(j_{\rm W}\bar\alpha_{\rm W}\alpha_{\rm W} x)\).
The primitive character on \(R\) is the reduction
of \(\psi_0(e_{\mathfrak m}s_0^{-1}x)\).
For a residue frequency matrix \(t\), let
\(\widetilde t\) be an integral lift of its transpose,
indexed by the Hom coordinates. Its normalized coordinate
is
\[
 y_t=\overline{e_{\mathfrak m}s_0^{-1}\widetilde t},
\]
and its actual adelic displacement is \(Y_t=(0,\alpha_{\rm W} y_t)\).
The finite part of the datum of
\(\psi_t(X)=\psi_{\mathfrak m}(\operatorname{tr}(Xt))\)
is \(\mathcal R(Y_t)v_{\rm vac}\). In fact the projection
onto the integral Fourier box fixes \(v_{\rm vac}\).
The translation above has support on the conjugate member
of each chosen split pair, so in the selected
Schr\"odinger box it is multiplication by the indicated
character. The dot product with \(\widetilde t\) is the
trace of \(Xt\). This also shows independence from the
integral lift of \(t\).

For an arbitrary adelic displacement \(z\), move
\(\mathcal R(Y_t)\) across the inner product. Its
commutator with \(\mathcal R(z_f)\) contributes
\(\beta(z_f,Y_t)^2\), with
\[
 \beta(z_f,Y_t)^2
   =\psi_f\!\left(j_{\rm W}\sum_\nu \bar\alpha_{\rm W}\,
                  e_{\mathfrak m}s_0^{-1}\widetilde t_\nu z_{f,\nu}
           \right).
\]
When \(a(z)\) is integral this is
\(\psi_t(a(z)\bmod\mathfrak m)\), by
\(\bar\alpha_{\rm W} z_f=\bar\alpha_{\rm W}\alpha_{\rm W} a(z)\).
The other finite factor is the vacuum coefficient
\[
 \mathbf1_{a(z)\in\widehat{\mathcal O}_E^{\,D}}\,
  \psi_f\!\left(\frac{c_{\rm tr}j_{\rm W}}2
              \sum_\nu\overline{z_{f,\nu}}z_{f,\nu}\right),
         \qquad c_{\rm tr}=(1-j_{\rm W}^{-1})/2 .
\]
Indeed the normalized translation fixes the vacuum on the
integral lattice, so the unnormalized coefficient there is the
inverse of its quadratic phase. The displayed phase equals its
inverse on this support: its square is
\(\psi_0(c_{\rm tr}\sum_\nu\overline{a(z)_\nu}a(z)_\nu)=1\)
by integrality and self-duality. This includes places of residue
characteristic two. For real Fock vectors \(v,w\) the real factor is
\(\langle v,\mathcal R_\infty(z_{\mathbb R})w\rangle\).
Finite character orthogonality is for the perfect trace
pairing of the two rectangular matrix spaces over \(R\).
The probability average over all \(t\) therefore restricts
the product of these factors to
\(a(z)\in\widehat{\mathcal O}_E^{\,D}\) and
\(a(z)\equiv0\bmod\mathfrak m\).
At \(z=0\), the finite factor, all commutators, and the
normalized average of the constant function are one.

\subsubsection{The sparse test and its error}
\label{sec:overlap-sparse-error}

\sectionopening{In this subsection, we apply Section~\ref{sec:overlap-finite-character} to the test used in the comparison.}

Fix the original two adelic points \(g,h\) and let
\(p\) be the induced isometry of \(\operatorname{Hom}_E(S,T)\)
over the adeles, given by \(X\mapsto hXg^{-1}\).
Let \(v\) be the Fock vector for the polynomial \(P\).
Write \(K_t(g,h)\) for its kernel with finite test
\(\psi_t\). Use the single finite test
\[
 \Psi(X)=\frac{1}{|M_{\operatorname{rk}S,\operatorname{rk}T}(R)|}
               \sum_t\overline{K_t(g,h)}\psi_t(X).
\]
At \(g(u,1),h(u',1)\) with real \(u,u'\), put
\(w=\omega_\infty^{a,b}(u,u')v\), including the determinant
weights of the original characters \(a,b\).
Linearity in the finite test writes its kernel as the same
average of the product of the theta value at \(w\) with the
conjugate theta value at \(v\), both moved by \(p\).
For the residue matrix space \(D_R\), write \(K_{[B]}\) for the
kernel of the indicator of the single residue \(B\), and let
\(a=(g,h)\), \(a'=(g(u,1),h(u',1))\). Finite Fourier orthogonality gives
\[
 \frac1{|D_R|}\sum_t K_t(a)=K_{[0]}(a),\qquad
 \frac1{|D_R|}\sum_t\overline{K_t(a)}K_t(a')
   =\sum_{B\in D_R}\overline{K_{[B]}(a)}K_{[B]}(a').
\]
Indeed \(K_t=\sum_B\psi_t(B)K_{[B]}\) with the same indicator of the
integral matrix box in each summand. At \(a'=a\) the second identity
is the squared norm of this vector of residue values.

The starting pair contributes the same unit scalar to both
theta values and cancels in the product. Move \(\omega(p)\)
through the Heisenberg coefficient by unitarity and covariance.
For a rational difference \(Z\colon S\to T\), put
\(z_0=X_T^{-1}ZX_S\). Then
\[
 \omega(p)^*\mathcal R(z_0)\omega(p)=\mathcal R(p^{-1}z_0),
 \qquad p^{-1}z_0=h^{-1}z_0 g=(X_Th)^{-1}Z(X_Sg)=z(Z).
\]
Use the finite reduction and scaled real coordinate of this same
\(z(Z)\) at \((g,h)\) as in Section~\ref{sec:peak-hom-lattice}.
The additional real translations \(u,u'\) act on the vector \(w\).
The Fourier identity just
proved applies for each \(t\). The resulting average is thus
a sum over precisely the rational differences that satisfy the
displayed finite support and congruence, with unit-modulus
finite phases. Its zero difference equals
\(\langle v,\omega_\infty^{a,b}(u,u')v\rangle\).
This is the original Fock matrix coefficient with mass one.

For the nonzero differences let \(N\) be the number of complex
Fock variables. The actual \(P=P_j\) is homogeneous of degree
\(d=\sum_\tau(mc_\tau+2m\nu_\tau)\le b\sqrt b\)
by Section~\ref{sec:peak-reference-coefficient}.
The coherent expression of the normalized real Weil operator
is
\[
       \omega_\infty(\sigma)v=
            \Delta_\infty(\sigma)\int v(x)
                    \mathcal R_\infty(\sigma x)v_{\rm vac,\infty}\,dx .
\]
Here \(dx=\pi^{-N}d\operatorname{Leb}(x)\). To check the
scalar, let \(G\) denote the unit vacuum. For weighted Fock vectors
one has \(\langle\mathcal R_\infty(x)G,v\rangle=v(x)\).
The coherent resolution of the identity follows by integrating
against the Fock reproducing kernel with this measure.
For this polynomial times the Gaussian, \(v\in L^1(dx)\),
hence the vector integral below converges in Hilbert norm for
every unit Schwartz vector \(G'\). The mixed identity is
\[
 \int\langle\mathcal R_\infty(x)G,v\rangle
              \mathcal R_\infty(x)G'\,dx=\langle G,G'\rangle v .
\]
Indeed in a Schr\"odinger realization Fourier inversion in the
modulation variable first identifies the two position arguments,
and integration in the translation variable gives \(\langle G,G'\rangle\).
The scalar in the measure is fixed by \(G'=G\).
Take \(G'=\omega_\infty(\sigma)^{-1}G\) and apply
\(\omega_\infty(\sigma)\), using Heisenberg covariance. This gives
the displayed formula with
\(\Delta_\infty(\sigma)=1/\overline{\langle G,\omega_\infty(\sigma)G\rangle}\).
Taking the constant polynomial in the boost calculation of
Section~\ref{sec:peak-reference-coefficient} gives one inverse
hyperbolic cosine for each boost. Compact factors give the inverse
determinant on their negative-sign block. Thus in the untwisted
unitary normalization the vacuum coefficient is the inverse
determinant of the negative-sign block of \(\sigma\). Each block has
operator norm at most \(\|\sigma\|\), so
\(|\Delta_\infty(\sigma)|\le\|\sigma\|^N\).

Inserting the coherent integral into
\(\langle v,\mathcal R_\infty(z_{\mathbb R})\omega_\infty(\sigma)v\rangle\)
gives the translate \(\mathcal R_\infty(z_{\mathbb R}+\sigma x)G\)
times the unit Heisenberg phase. Its pairing with \(v\) has modulus
\(|v(z_{\mathbb R}+\sigma x)|\), so the displaced coefficient is bounded by
\[
       \|\sigma\|^N\int|v(x)|\,|v(\sigma x+z_{\mathbb R})|\,dx .
\]
The determinant-character factors are units.
If \(q>N\), put \(A_0=d+N\) and \(c^2=1+q/A_0\).
For \(G(x)=|c^dv(x/c)|^2\), homogeneity gives
\[
 \int G\le e^q\|v\|^2,\qquad
 |v(x)|=G(x)^{1/2}
                   e^{-q\|x\|^2/(2(A_0+q))}.
\]
The first inequality uses the dilation Jacobian \(c^{2N}\).
Here the norm on displacements is the Euclidean norm after the
real scaling above. Passage to the holomorphic coordinates of the
signed Fock model conjugates the coordinates with negative oscillator
sign. This is a real orthogonal map, so it preserves the Gaussian,
\(dx\), and the operator norm of the induced Hom action. The latter
preserves real volume by symplectic invariance and has
\(\|\sigma\|\ge1\). The inequality
\[
 \|x\|^2+\|\sigma x+z_{\mathbb R}\|^2
                       \ge \|z_{\mathbb R}\|^2/(2\|\sigma\|^2)
\]
follows from the triangle inequality and the two-square bound.

Apply \(2ab\le a^2+b^2\) to the two square roots of \(G\).
For \(z_{\mathbb R}\ne0\), the remaining Gaussian is at most
\(e^{-q t_z}\), where
\[
 t_z=\frac{\|z_{\mathbb R}\|^2}{4(A_0+q)\|\sigma\|^2}.
\]
Since \(e^{-q t_z}\le e^{-q}t_z^{-q}\) and \(\sigma\) preserves
volume, integration cancels the factor \(e^q\). This gives
\[
  \|\sigma\|^{N+2q}(4(d+N+q))^q
                    \|v\|^2\|z_{\mathbb R}\|^{-2q}.
\]
On Hom, \(\sigma(X)=u'Xu^{-1}\) and
\(\|\sigma\|\le |u'|\,|u|\), since the operator norm of a
signed unitary matrix equals that of its inverse.

On the reduction kernel the original multiplier
\(q_1/\pi_M^2\) sends every surviving difference into the
Hom lattice. When the two given first minima have product
at least \(b^{-13}\), the Hom-lattice calculation gives
separation \(\mathcal D=b^6\) for their Fock coordinates,
including distance at least \(\mathcal D\) from zero.
Packing balls of radius \(\mathcal D/2\) bounds the count
in the annulus of radii \(2^i\mathcal D\) and
\(2^{i+1}\mathcal D\) by \((5\cdot 2^i)^{2N}\).
Summing the geometric series and using \(q\ge N+1\) gives
\[
        \sum_{z_{\mathbb R}\ne0}\|z_{\mathbb R}\|^{-2q}
                     \le (5/\mathcal D)^{2q}.
\]
Here \(N=f(3m)(r+3m)\) since
\(\dim S=3m\) and \(\dim\widetilde W=r+3m\).
Both dimensions are at most \(t\) for sufficiently large \(x\).
Take \(q=2n_*\), so \(N\le ft^2<q\).
The error in the sparse test is at most
\[
 (|u|\,|u'|)^{6n_*}\|v\|^2
                   \left(\frac{100(d+2q)}{b^{12}}\right)^q
 \le b^{-20n_*}(|u|\,|u'|)^{6n_*}\|v\|^2.
\]
Indeed \(100d\le b^2/2\) for \(b\ge200^2\) and
\(200q=400n_*\le b^2/2\) for sufficiently large \(x\) with the
external modulus, quadratic discriminant and \(\epsilon\) fixed.
This uses the polynomial growth of \(n_*\) and
\(\log b=x^\epsilon+O(1)\).
The same finite test works for all \(u,u'\) in this estimate.
At \(u=u'=1\) the estimate gives
\[
 \frac1{|D_R|}\sum_t |K_t(g,h)|^2
       \ge (1-b^{-20n_*})\|v\|^2 .
\]
This is positive for the prescribed nonzero Fock vector.
At every real translate the error uses the same fixed norm
\(\|v\|^2=c_P(1,1)\).

\section{Compression, separation, and reduction at cusps}
\label{sec:comparison-compression-separation}

\sectionopening{In this section, we compress the comparison operator of
Section~\ref{sec:comparison-operator-first-lifts} and establish separation
on the initial unitary group. The cusp estimates justify the decay needed
for the deeper periods of Section~\ref{sec:deeper-periods}.}

\subsection{Compression of the comparison operator}
\label{sec:original-compression-bounds}

\sectionopening{In this subsection, we restrict the comparison of
Section~\ref{sec:original-comparison-tests} to the image of the
injective lift and project with respect to the Petersson form. We
bound the determinant of this operator using the original norm on each
factor.}

Fix one original descent pair and write \(u_\eta=u_{j,\eta}\)
for the injective lift on the projected source \(W\).
The identity test in \eqref{eq:marked-full-spherical-trace}
gives \(\dim W=D\), so \(Z_\eta=\operatorname{im}u_\eta\)
has dimension \(D\) for every allowed twist.
Let \(\Pi_\eta\) be the orthogonal projection onto
\(Z_\eta\) for the target class pairing. It is also the Petersson
projection, since at a fixed level and Weyl type the two pairings
differ by the positive level and compact Schur factor of
Section~\ref{sec:first-lift-sharp-masses}. Use the original gluing operator \(G_\eta\),
its original tuple space, and the operator
\(U_\eta=n_Ah_-^2G_\eta G_\eta^*\) of
Section~\ref{sec:original-comparison-tests}. Put
\[
 \Delta_\eta=\det(\Pi_\eta U_\eta|_{Z_\eta}),\qquad
 Z_j(\eta)=\det(u_\eta^*u_\eta).
\]

\subsubsection{The period and the Hecke factor}\label{sec:compression-period-hecke}

\sectionopening{In this subsection, we use Proposition~\ref{prop:common-period-adjoint} in the compressed pairing.}

\begin{lemma}[Determinant of the compressed comparison]\label{lem:marked-compression-determinant}
Set \(\rho_\eta=H_j(\eta)/H_j(1)\) for the allowed twists. If \(D>0\),
then \(\rho_\eta>0\) and there is a fixed positive definite operator
\(U^\circ\) on the source such that
\[
 u_\eta^*U_\eta u_\eta=\rho_\eta U^\circ,\qquad
 Z_j(\eta)\Delta_\eta=\det(U^\circ)\rho_\eta^D,\qquad
 \frac{Z_j(1)}{Z_j(\eta)}=\frac{\Delta_\eta}{\Delta_1}\rho_\eta^{-D}.
\]
If \(D=0\), all the determinants and ratios of determinants are \(1\).
\end{lemma}

\begin{proof}

For \(D>0\), let \(P\) be the injective common-period map and
\(c'_\eta\ne0\) the scalar in
\(c'_\eta G_\eta^*u_\eta=P\).
The three exchanges of Section~\ref{sec:gluing-three-exchanges} give
\[
 |c'_\eta|^2=C_{\rm ex}/H_j(\eta),
\]
where \(C_{\rm ex}\ne0\) is independent of \(\eta\), and
\(H_j\) is the product of Hecke values including
the local conductor factors. In particular
\[
 \rho_\eta=\frac{H_j(\eta)}{H_j(1)}
     =\frac{|c'_1|^2}{|c'_\eta|^2}>0,\qquad
 u_\eta^*U_\eta u_\eta=\rho_\eta U^\circ,\qquad
 U^\circ=\frac{n_Ah_-^2}{|c'_1|^2}P^*P>0 .
\]
Indeed the pairing of the middle operator against
\(y',y\in W\) is \(n_A h_-^2\langle G_\eta^*u_\eta y',
G_\eta^*u_\eta y\rangle\); using the two scalar factors
gives the formula. Here \(u_\eta^{-1}:Z_\eta\to W\) inverts the lift
onto its image, and \(u_\eta^*(1-\Pi_\eta)=0\) by orthogonality.
Hence
\[
 u_\eta^*U_\eta u_\eta
    =(u_\eta^*u_\eta)
        \bigl(u_\eta^{-1}\Pi_\eta U_\eta u_\eta\bigr).
\]
Taking determinants yields
\[
 Z_j(\eta)\Delta_\eta=\det(U^\circ)\rho_\eta^D,\qquad
 \frac{Z_j(1)}{Z_j(\eta)}
       =\frac{\Delta_\eta}{\Delta_1}\rho_\eta^{-D}.
\]
For the covectors \(f_a\) of
\eqref{eq:original-comparison-test-covector}, transported to the tuple
space of \(G_\eta\) as in Section~\ref{sec:original-comparison-fock},
put \(z=u_\eta y\). Orthogonal projection and the rank-one sum give
\[
 \sum_a|\langle\Pi_\eta G_\eta f_a,z\rangle|^2
   =h_-^2\|G_\eta^*z\|^2
   =h_-^2|c'_\eta|^{-2}\|Py\|^2.
\]
Thus the projected vectors span \(Z_\eta\), and positivity gives
the nonzero determinants used above. In dimension zero the determinant
is the empty product.
\end{proof}

\subsubsection{The projection at complex embeddings}\label{sec:compression-complex-projection}

\sectionopening{In this subsection, we bound the conjugate projection for the operator of Section~\ref{sec:compression-period-hecke}.}

At each complex automorphism \(\sigma\), pull back the positive
Gram pairing of the coefficientwise images of the form values
by \(\sigma^{-1}\). Use the original rational class weights
and Weyl Gram matrices. This pairing is semilinear in its first
argument for \(\sigma^{-1}c\sigma\), where \(c\) is ordinary
complex conjugation. Let \(J_0\) be the original target level and
\(J\) the rational Weyl Gram matrix. For two finite-level forms take a
common compact open \(K'\subseteq J_0\) fixing both, reducing to the
kernels of their finite level characters if necessary. The function
\(b_\sigma(g)=\langle\sigma f(g),\sigma h(g)\rangle_J\) is locally
constant on the finite adelic quotient, by the rational isometry
below. The pairing on the coordinatewise images, before pullback, is
\[
 \frac{1}{\operatorname{vol}(J_0)}\int b_\sigma(g)\,dg
 =\frac{\operatorname{vol}(K')}{\operatorname{vol}(J_0)}
   \sum_{[g]_{K'}}\frac{b_\sigma(g)}{|\operatorname{Aut}_{K'}(g)|}.
\]
For a rational point \(\gamma\), each matrix \(\sigma\tau(\gamma)\)
is unitary, since every embedding of the CM field commutes with
conjugation on it. Thus the rational Weyl action preserves the
summands under changes of representatives. For two forms
of the original level character \(\chi\), the equality
\(\overline{\sigma\chi(j)}\,\sigma\chi(j)=1\) gives \(J_0\)-invariance.
For a \(J_0\)-class represented by \(x\), write
\[
 \Gamma_x=x^{-1}U(I_n)(F_0)x\cap J_0.
\]
The finer classes are the \(\Gamma_x\)-orbits on \(J_0/K'\).
If \((\Gamma_x)_y\) is the stabilizer of a representative \(y\)
of such an orbit, then
\[
 \frac{1}{[J_0:K']}\sum_{y\in\Gamma_x\backslash J_0/K'}
                   \frac1{|(\Gamma_x)_y|}
    =\frac1{|\Gamma_x|}.
\]
This gives the original weight at level \(J_0\).

For \(R(a)f(g)=f(ga)\), right invariance of \(dg\) gives
\(R(a)^*=R(a^{-1})\) on the images with this fixed normalization.
In either finite sum one may further intersect the levels of the
translated test forms. If a translation carries a tame character
factor \(u\), its adjoint carries \(\overline{\sigma(u)}=\sigma(u^{-1})\),
since \(u\) is a root of unity. Pulling back these finite identities
by \(\sigma^{-1}\) gives the adjunction for the conjugate pairing.
Write \(\|\cdot\|_\sigma\) for its positive norm on the images.

The original finite group acts by right translations and
finite-order character factors, so its action is unitary
at each embedding. Put \(\tau_\sigma=\sigma^{-1}c\sigma\). For
each ordinary irreducible character \(\chi\) of the finite group,
\(\tau_\sigma(\chi(g^{-1}))=\chi(g)\): the character values are
cyclotomic, and complex conjugation commutes with every automorphism
of a cyclotomic field. Thus taking adjoints in the finite character
sum for \(e_0\) both inverts \(g\) and conjugates its coefficient,
giving \(e_0^{*_{\sigma}}=e_0\).
By \cite[Corollary~4.12(c)]{MV} there is a rational prime
\(p\equiv1\bmod dq_1\) away from the finite level and its conjugate.
Every place \(w\mid p\) is split over \(F_0\), since
\(E\subseteq\mathbf Q(\zeta_{dq_1})\). Choose one such \(w\).
Let \(T_{w,i}\) be its degree-\(i\) elementary spherical
double coset, with \(T_{w,0}=1\). In the notation of
Proposition~\ref{prop:clearing-on-original-tests} put
\[
 A_w=K_w(X_w),\qquad
 K_w(X)=\sum_{i=0}^{n}(-X)^i q^{\binom{i}{2}}T_{w,i},\qquad
 X_w=(q^{t'-1}\nu_w)^{-1},\quad q=\mathrm Nw .
\]
Here \(n=r'+t'\) is the target rank. On a spherical system with
monic Hecke polynomial \(P_x\), this operator has eigenvalue
\(X_w^nP_x(X_w^{-1})\).

Let \(K=\mathrm{GL}_n(\mathcal O_w)\),
\(d_i=\operatorname{diag}(\varpi_w I_i,I_{n-i})\), and
\(z_w=\varpi_w I_n\). With \(K\) of mass one, inversion interchanges
the operator of \(Kd_iK\) on the original \(K\)-fixed space with
that of \(Kd_i^{-1}K\) in the adjoint. Permuting the diagonal gives
\[
 T_{w,i}^{*_{\sigma}}=T_{\bar w,i}=\rho(z_w^{-1})T_{w,n-i}
                    \qquad (0\le i\le n).
\]
Here \(\rho\) is right translation on the split group. The middle
identity uses the inverse transpose at the conjugate place; changing
the conjugate uniformizer by a unit preserves the double coset.
Consequently the adjoint of \(A_w\) for this pairing is
\(K_{\bar w}(\tau_\sigma(X_w))\), since all the
coefficients \((-1)^i q^{\binom{i}{2}}\) are rational.
These two operators commute in the spherical algebra
of the split pair. Normality gives
\(\ker A_w^*=\ker A_w\), so
\(\operatorname{im}A_w\perp_\sigma\ker A_w\).

The factor \(1-\nu_w q^{t'-1}X\) in the split lift identity of
Section~\ref{sec:normalized-trace-multiplicities} places every
lift generator at the trivial twist in \(\ker A_w\). Thus
\(Z_1\subseteq\operatorname{im}e_0\cap\ker A_w\). For the same subspace
\[
 F_w=Z_1^\perp\cap\operatorname{im}e_0\cap\ker A_w
\]
used in Proposition~\ref{prop:clearing-on-original-tests} and
\eqref{eq:block-period-compression-orthogonality}, the orthogonals
in the original metric satisfy
\(F_w^\perp=Z_1+\ker e_0+\operatorname{im}A_w\).
The last two summands are orthogonal to \(Z_1\) for every
conjugate pairing, including the original one. Hence if
\(v\in F_w^\perp\), its decomposition has the form
\[
 v=z+t,\qquad z=\Pi_1v\in Z_1,\qquad
 t\in\ker e_0+\operatorname{im} A_w.
\]
The same two vectors are orthogonal for every conjugate pairing, so
\[
 \|v\|_\sigma^2=\|\Pi_1v\|_\sigma^2+\|t\|_\sigma^2.
\]
Equation~\eqref{eq:block-period-compression-orthogonality} with
\(\vartheta=1\) gives \(U_1z\in F_w^\perp\) for \(z\in Z_1\).

\subsubsection{The norm at the trivial twist}\label{sec:compression-trivial-norm}

\sectionopening{In this subsection, we combine the projection bound of Section~\ref{sec:compression-complex-projection} with the original gluing sum.}

Write \(A_{\rm glue}=fmbB^2\), \(N_{\rm glue}\) for the number of gluing choices times the number of second-factor classes,
and \(l_1\) for the denominator exponent of the splitting. The original gluing sum has the factor \(M^{l_1}\).
At the trivial twist, the first-block level \(J_\Lambda\),
the two principal levels in the choice of a tuple coordinate,
and the target level all contain \(K(\pi_M^2)\) in their
respective ranks, by Section~\ref{sec:marked-comparison-lattice}.
Since \(\pi_M\bar\pi_M=M\), they contain \(K(M^2)\).
For any of these levels \(J\) of rank \(n\le r+t\),
the lattice count in the proof of Lemma~\ref{lem:terminal-class-count} gives
\[
 \begin{aligned}
 h(J)&\le h(K_{\rm int})[K_{\rm int}:K(M^2)]
       \le h(K_{\rm int})(M^2)^{2fn^2},\\
 \log h(J)&\ll(2fn)^3\log x+4fn^2\log M<b .
 \end{aligned}
\]
Thus each of these class sets has at most \(e^b\) elements.
For a glued point \(z\), let \(\Gamma_z\) be its target stabilizer.
On its single class, \(\delta_z(v)\) is the sum of \(|\Gamma_z|\)
unitary transports of \(v\), including their level characters.
The squared class norm has weight \(1/|\Gamma_z|\).
The transported matrices remain unitary at every \(\sigma\), so
\[
 \|\delta_z(v)\|_\sigma^2\le |\Gamma_z|\,\|v\|_\sigma^2
                               \le n_A\|v\|_\sigma^2 .
\] In the exterior tensor model of
Section~\ref{sec:marked-comparison-lattice}, the two block
inclusions, signed complement maps and determinant twists
preserve the Weyl forms. The normalized Cartan projection is the
identity on the included tensor product, as verified there using
the highest tensor in the prescribed orders. Write \(T\) for this
Cartan map. Its coordinates in the integral exterior models are
defined over \(\mathbf Z\), so at every embedding
\[
 \|\sigma T(y\otimes v)\|^2
     =\|\sigma y\|^2\,\|\sigma v\|^2 .
\]

By \eqref{eq:display-i5}, \(\|\phi''_{1,i}(h)\|_\sigma^2
\le e^{A_{\rm glue}}\) on every original first factor. Apply the triangle
inequality to the sum over \(h\) on each tuple coordinate. For the
sum over coordinates the free second-factor classes give
\(\sum_{i,c}\|\zeta_{i,c}\|_\sigma\le
\sqrt{N_{\rm glue}}\,\|\zeta\|_\sigma\) by Cauchy--Schwarz. Hence
\[
 \|G_1\zeta\|_\sigma^2
     \le M^{2l_1}e^{2b}N_{\rm glue} n_A
                     e^{A_{\rm glue}}\|\zeta\|_\sigma^2 .
\]
There are at most \(e^b\) choices in each of the two
principal-level class sets defining \(N_{\rm glue}\). The residue
matrix has \(3m d_+\) entries, since \(\dim_E S=3m\) and
\(\dim_E V''=d_+\). Its residue ring has order at most \(M^{4f}\).
At each of the \(f\) embeddings the index specifies
\(\Delta'=(\Delta-t-1)/2\) ordered tuples of \(2m\) covectors
from the chosen basis of \(W^\flat\), whose rank is \(6m+1\).
Consequently
\[
 N_{\rm glue}\le e^{2b}M^{4f d_+(3m)}(6m+1)^{2mf\Delta'} .
\]

The column estimate for the
denominator exponent is
\[
 l_1\le 3fm(b+a/x_1+\Delta+m\Delta)+ft^2\Delta .
\]
Here \(B=\log b\), \(\log M\le22B\), \(a/x_1\le b\), and
\(\log n_A\le f(rt)^2\log M\). Since \(f,r,t,m,\Delta\) have
polynomial bounds in \(x\), there is a fixed exponent \(C\) with
\[
 2l_1\log M+2b+\log(N_{\rm glue}n_A)
       \le264fmbB+4b+O(x^C B)=o(A_{\rm glue}).
\]
The same
bounds give \(n_Ah_-^2\le e^{A_{\rm glue}}\), since the minus diagram
has at most \(ft^2\Delta\) boxes and its hook factor is
bounded by \(L_-^{3L_-}\) for \(L_-\ge2\) bounding that count.

For precision about adjunction, use the full spaces of class
values with their rational block Gram matrices. The tuple
space of \(G_1\) is the full tuple-coordinate space, since
its stabilizers are trivial. The image of \(G_1\) lies in
the original target, and its ambient adjoint restricts
there to the adjoint. With \(C_S,C_T\)
the rational Gram matrices, the adjoint matrix is
\(C_S^{-1}c(G_1)^{\mathsf t}C_T\). Therefore
\(\sigma(G_1^*)\) is the ordinary Gram adjoint of
\((c\sigma c)(G_1)\) in the ambient coordinates. On
conjugated class values, the adjoint norm bound follows
from the uniform bound for \(G_1\) just proved on its full
tuple domain. Thus
\[
     \|U_1z\|_\sigma\le n_Ah_-^2e^{2A_{\rm glue}}\|z\|_\sigma
                         \le e^{3A_{\rm glue}}\|z\|_\sigma .
\]
An eigenvector of \(\Pi_1U_1|_{Z_1}\) with eigenvalue
\(\lambda\) now gives \(|\sigma\lambda|\le e^{3A_{\rm glue}}\).
Multiplying over the eigenvalues with algebraic multiplicity
gives \(\log|\sigma\Delta_1|\le3DA_{\rm glue}\).
This uniform bound implies \(\Delta_1\in\overline{\mathbf Q}\) by the
transcendence-basis argument of Section~\ref{sec:gram-product-formula}.

\smallskip
\phantomsection\label{sec:simultaneous-original-clearing}
\indent\textit{One operator on the full projected lift space.}
Choose a basis \(v_1,\ldots,v_s\) from the lift generators
of \(Z\). For each \(i\), construct the first-order operator
\(D_i\) by using the type and original system of \(v_i\).
Every generator \(v_j\) has its own fixed type parameter. Its
projection to the type of \(v_i\) is either zero, if the
types differ, or \(v_j\) itself. In the latter case the
spherical system remains the same after inclusion in the
raised level used for \(D_i\). Thus
\(D_iv_j=d_{ij}v_j\), with the bound \eqref{eq:display-rb1} on \(v_\ell(d_{jj})\).
Put
\[
 \mathcal T_\ell=\{y\text{ integral at the chosen embedding}:
                              y\perp F_w\text{ for every }w\}
\]
in the original target at this twist. Equivalently,
\(\mathcal T_\ell=L_\kappa\cap(\sum_w F_w)^\perp\), where
\(L_\kappa\) is the transported integral lattice and the sum ranges
over the split places outside the original level and its conjugate.
Each \(F_w\) uses the full projector \(e_0\) and \(Z^\perp\).
Thus the test space is the same for each individual type and system
selected in constructing \(D_i\). Proposition~\ref{prop:clearing-on-original-tests}
gives \(D_i\pi_Zy\) integral for each \(y\in\mathcal T_\ell\).

Write \(\iota\) for that embedding and, for each \(j\), select
\(a_j\ne0\) from the numbers \(\iota(d_{ij})\) of maximum
absolute value. Reduce the integral ratios
\(\iota(d_{ij})/a_j\) in the residue field of
\(\overline{\mathbf Q}_\ell\). Each resulting linear form
\(\sum_i X_i\overline{\iota(d_{ij})/a_j}\) is nonzero.
Choose a finite subfield \(k'\subset\overline{\mathbf F}_\ell\)
containing the reduced coefficients with \(q=|k'|>s\). The \(s\)
kernels contain at most \(s q^{s-1}<q^s\) points in \((k')^s\).
For a tuple in the complement, take \(b_i=0\) at its zero
coordinates and the Teichm\"uller lifts at its other coordinates.
Each nonzero \(b_i\) satisfies \(b_i^{q-1}=1\), so every \(b_i\)
is integral and belongs to \(\iota(\overline{\mathbf Q})\). The operator
\(D_Z=\sum_i\iota^{-1}(b_i)D_i\) still gives an integral form
\(D_Z\pi_Zy\) for every \(y\in\mathcal T_\ell\). On \(v_j\) its
eigenvalue has absolute value \(|a_j|_\ell\), at least
\(|\iota(d_{jj})|_\ell\). Each original off-level Hecke operator
preserves every cofinite spherical isotypic space and commutes with
each old type projector, since their local supports are disjoint.
The trace coefficients are scalar on those isotypic spaces.
Thus the restrictions \(D_i\), and hence \(D_Z\), commute with the
off-level Hecke operators on \(Z\), including on joint summands of
dimension greater than one. The determinant bound follows by
multiplying the eigenvalue bounds. If \(Z=0\), use the
identity on the zero space, with determinant one.

\subsubsection{Compression on the original local lattice}\label{sec:compression-local-lattice}

\sectionopening{In this subsection, we apply the integral clearing of Proposition~\ref{prop:clearing-on-original-tests} to the local determinant.}

Fix \(\kappa\mid\ell\), normalize \(v_\kappa(\ell)=1\),
and take the selected twist \(\eta\) whose auxiliary modulus has
norm prime to \(\ell\).
Use the clearing operator \(\Lambda=D_Z\) just constructed
on \(Z=Z_\eta\). The selected-twist bounds in
Section~\ref{sec:common-critical-twist}
give \(v_\kappa(\det\Lambda)\le Dr C_\ell\).
More explicitly, \eqref{eq:display-rb1} and
Proposition~\ref{prop:clearing-on-original-tests} give
\(v_\kappa(\det D_Z)\le D B_{\rm tr}\) on this image.
The selected twist bounds the Selmer length for every generator in the
basis of \(Z\), so \(\mathrm{sel}_Z,\mathrm{ram}\le e_\ell\). Using
\(t'\le t\), \(r'=r\), and the definition of \(C_\ell\) in
Section~\ref{sec:critical-integral-pairings} give
\(B_{\rm tr}+v_\ell(C_G)\le r C_\ell\).
Clearing the projected vectors uses \(D_Z\) itself. The extra
factor \(C_G\) there clears the denominators of a scalar class
pairing. The test for \(D_Z\) must be integral at \(\kappa\)
and perpendicular to every \(F_w\).
The operator \(U_\eta\) preserves integrality at \(\kappa\).
The vanishing in Section~\ref{sec:marked-comparison-period-vanishing}
gives \(U_\eta Z_\eta\subseteq\bigcap_w F_w^\perp\).
Consequently \(\Lambda\Pi_\eta U_\eta|_{Z_\eta}\)
preserves integral vectors in \(Z_\eta\).

For a nonzero eigenvector \(F\) of the composite, take one
representative \(x\) of each original finite class. Write \(x_\ell\)
for its components at all the places above \(\ell\).
The Cayley-chart weak approximation of
Section~\ref{sec:original-comparison-fock} gives a rational
unitary element \(\gamma\) such that \(\gamma_\ell x_\ell\) is in
the standard maximal compact there. Indeed
\(K_{{\rm std},\ell}x_\ell^{-1}\) is a nonempty local open set.
For the Weyl representation \(\rho_{\rm W}\), rational equivariance gives
\[
 \rho_{\rm W}((\gamma_\ell x_\ell)^{-1})F(\gamma x)
       =\rho_{\rm W}(x_\ell^{-1})F(x).
\]
The matrix on the left preserves the unit Weyl lattice in both
directions. Thus the ordinary coordinates at \(\gamma x\) test
the transported lattice norm at \(x\). Right level changes also
preserve it since they act integrally with unit character.
The maximum over these finitely many coordinates is positive.
Divide \(F\) by a coordinate attaining this maximum.
The eigenvalue then has nonnegative valuation, because the
composite preserves integrality and the normalized vector
has a unit coordinate. Taking the product of the
eigenvalues gives
\[
 v_\kappa(\Delta_\eta)\ge-v_\kappa(\det\Lambda)
            \ge-Dr C_\ell,\qquad
 \log|\Delta_\eta|_\kappa\le Dr C_\ell\log\ell .
\]
In particular the finite logarithmic allowance is
\(\beta_{j,\ell}=C_\ell\log\ell\), or the larger
\(B_\ell\log\ell\) used for the paired rank-one test.
Since \(B_\ell=C_\ell+v_\ell(h_1)\) and \(h_1\) is an integer,
in either case \(\beta_{j,\ell}\le B_\ell\log\ell=J_\ell\).
Thus \eqref{eq:display-rt3} applies also to
\(\sum_{\ell\in U}\beta_{j,\ell}\) for every finite \(U\).
These bounds only use the clearing operator on the same
selected twist. The fixed multiplier for the ratio
\(\Delta_\eta/\Delta_1\) is \(\Delta_1\), which does not depend on \(\eta\) and serves as one
height multiplier at every coefficient embedding.

\subsection{Separation on the initial unitary group}
\label{sec:marked-deep-separation}

\sectionopening{In this subsection, we show that the inclusion at the initial deep level separates the
chosen automorphic coefficients. The reproducing coefficient of
Section~\ref{sec:real-reproducing-coefficient} gives the local test,
with the scale selected in Section~\ref{sec:marked-deep-cusp}.}

Write \(V=Y_0\), \(Z=V\perp\langle1\rangle\), and
\(\iota(h)=h\oplus1\) in the original rational basis. In particular,
\(\iota^{-1}(K_Z(a))=K_V(a)\) at each principal modulus \(a\).
The grid numbers the first \(r\) entries of \(Z\) exactly
as for \(V\). Its last entry has index \(r\) and sign
\((-1)^{r+1-k}=(-1)^m=+1\). Thus the real group inclusion is
the block inclusion of this rational orthogonal sum. Set
\(a_0=\pi_M^{b_0}\), \(b_0=\lfloor b\rfloor\), and denote the
original principal source level by \(J_M=K_V(\pi_M)\).
Let \(W\) be the space of square-integrable forms at level \(J_M\)
reproduced by the prescribed real coefficient projector. Define \(W_Z\)
using the initial true grid projector and level \(K_Z(a_0)\).
Every element of each space has a bounded continuous representative,
by the truncation argument below.
We will construct one \(\omega\in U(Z)(F_\infty)\) such that
\[
 y\in W,\quad
 \int_{[U(V)]}\overline{y(h)}\,\psi(\iota(h)(\omega,1))\,dh=0
 \quad(\psi\in W_Z)
 \quad\Longrightarrow\quad y=0 .
 \tag{DS0}\label{eq:deep-separation-condition}
\]
Here \(F\) is the maximal real field of \(E\). We use the original quotient Haar measures.

\indent{\it Weights and the archimedean scale.}
At the infinite place \(v\) indexed by \(\tau\), put
\[
 w_v=\kappa_\tau+\Delta/2,\qquad \Delta=2r^8,
\]
with the integer weight \(\kappa_\tau\) of Section~\ref{sec:notation-fields}.
Let \(l,l_Z\) denote the numbered doubled grid parameters of
\(V,Z\). The actual form types of \(W,W_Z\) use \(l,-l_Z\),
respectively.
In ranks \(r+z\), \(z=0,1\), the numbered doubled entries of the
grid are \(2\lambda_i+z\). On \(0\le i<k-1\),
\(\lambda_i=k/2-1-i\); on the tail \(i=k-1+j\),
\(\lambda_i=-w_v-(j+1/2)\Delta\).
Both initial blocks have sign \(+\), and the tails alternate
starting with \(+\). The correction from doubled parameter to
doubled lowest-weight entry is the signed count of later entries
minus earlier entries, signs taken relative to the entry in question.

On the initial block this equals \(k-2-2i+z\), since the tail has
signed excess \(z\). Thus the lowest-weight entry there is zero.
The correction everywhere is bounded in absolute value by \(r+1\).
Consequently the lowest-weight tails of the unnegated grid parameters
\(l,l_Z\) are strictly decreasing in
\([-3w_v,-w_v]\). Indeed \(w_v\ge\Delta^2\),
\(m\Delta= o(w_v)\), and \(\Delta>r+1\).
The doubled entries have modulus at most \(L_{\rm sep}=\lceil20b\rceil\),
since \(w_v\le2b\) and \((m+1)\Delta=o(b)\).
The opposite-sign gaps for the undoubled parameters are at least
\(\Delta\) before subtracting \(500(r+1)^3\). The same
gap bounds hold after negating a parameter.

Put \(A_{\rm sep}=r^4\), \(D_{\rm sep}=(2r+1)^2\), \(H=20D_{\rm sep}\),
\(\beta=(\sqrt M)^{\lceil x\rceil}\), \(B_d=\beta^{20A_{\rm sep}(f+1)}\),
and let \(\kappa=8fr\operatorname{rd}(E)+48r^6f^2N_D+1\),
where \(f=[F:\mathbb Q]\) and \(N_D=q_1d\) lies in the different.
The parameter bounds \(r\asymp x^5\), \(m\le2\sqrt x\),
\(b^{20}\le M\le b^{22}\), and \(\log b=x^\epsilon+O(1)\) give
\[
 \begin{gathered}
 A_{\rm sep}>8(r+1)^2,\quad A_{\rm sep}\ge2f(r+1)+2,\quad H<A_{\rm sep},\quad
 \beta\ge\max\{10^{12},f,2r+1+L_{\rm sep}\},\\ \beta\ge\max\{2000(r+L_{\rm sep}+1)^6,M\kappa^2\},\\
 (r+1)^4(6w_v)^{2m+2}\le\beta,\qquad
 10^{30}(D_{\rm sep}\log\beta)^3=O(x^{33+3\epsilon})\le r^7,\\
 |\tau(a_0)|_{\mathbb C}=(\sqrt M)^{b_0}\ge B_d^{3r+2}.
 \end{gathered}
 \tag{DS1}\label{eq:deep-reproduction-scale}
\]
Here \(|\cdot|_{\mathbb C}\) is the Euclidean modulus in the moved
lattice. The last equality follows from \(\pi_M\bar\pi_M=M\).
For the first line use the polynomial bounds on \(f,q_1d,
\operatorname{rd}(E)\) of the field and
\(\beta\ge b^{10\lceil x\rceil}\); in particular the admissible
bound on \(f\) gives \(2f(r+1)+2\le r^4\).
For the basis bound, \(r+1,6w_v\le\sqrt M\) and
\(2m+6\le\lceil x\rceil\). For the logarithmic comparison,
\(D_{\rm sep}\log\beta=O(x^{11+\epsilon})\) and
\(33+3\epsilon<35\), since \(0<\epsilon<1/100\).
Thus \(10^{30}(D_{\rm sep}\log\beta)^3=o(r^7)\). Finally
\(b_0\ge20A_{\rm sep}(f+1)(3r+2)\lceil x\rceil\) since its right side
is polynomial in \(x\). The excess gap is at least \(r^7\),
and in particular at least \(4A_{\rm sep}\) in the smaller group and \(2A_{\rm sep}\)
in the larger group.

\indent{\it The real point in the local period integral.}
Let \(\mu_i\), \(0\le i\le m\), be the decreasing tail of the
lowest weight for \(l_Z\), and \(\nu_j\), \(0\le j<m\), the
increasing tail for \(-l\). Then \(\mu_i\in[-3w_v,-w_v]\)
and \(\nu_j\in[w_v,3w_v]\). Write
\[
 p(t)=\prod_i(t-\mu_i),\qquad a(t)=\prod_j(t-\nu_j),\qquad
 [R]_i=\frac{R(\mu_i)}{\sqrt{|p'(\mu_i)a(\mu_i)|}} .
\]
In the tail sign models define
\[
 v_*=[a],\qquad
 u_j=\sqrt{|p(\nu_j)/a'(\nu_j)|}\,[a(t)/(t-\nu_j)].
\]
Since \(m\) is even, \(p'(\mu_i)a(\mu_i)\) has sign
\((-1)^i\). The identity
\(\sum_i R(\mu_i)/p'(\mu_i)=[t^m]R(t)\) for \(\deg R\le m\)
gives \((v_*,v_*)=1\), \((v_*,u_j)=0\), and
\((u_j,u_k)=0\) if \(j\ne k\).
The signed square norm of \([a(t)/(t-\nu_j)]\) is
\(-a'(\nu_j)/p(\nu_j)\): sum the residues of
\(a(t)/((t-\nu_j)^2p(t))\).

Its sign is \((-1)^j\).
Multiplication by \(t\) on the polynomial vectors gives
\(\mu u_j=\nu_j u_j+\sqrt{|p(\nu_j)/a'(\nu_j)|}\,v_*\).
Hence the isometry with columns \(u_j\) compresses \(\mu\)
to \(\nu\). Extend this isometry by the identity on the common
zero-weight positive block \(X\), and append \(v_*\) last to
obtain \(w\in U(Z_v)\). In these coordinates,
\[
 |(v_*)_i|^2=\frac{|a(\mu_i)|}{|p'(\mu_i)|},\qquad
 |(u_j)_i|^2=
 \frac{|p(\nu_j)a(\mu_i)|}
 {|a'(\nu_j)p'(\mu_i)|\,|\mu_i-\nu_j|^2}.
\]
The distinct integral nodes in each tail give
\(|a'(\nu_j)|,|p'(\mu_i)|\ge1\), and the cross distances are in
\([2w_v,6w_v]\). Hence every entry of \(w\) has modulus at most
\((6w_v)^{m+1}\). For \(\omega_v=w^{-1}\), the
equality \(\|w^{-1}\|=\|w\|\) and \eqref{eq:deep-reproduction-scale} give
\(\|\omega_v\|^2\le\beta\).

Put \(g=w(h\oplus1)w^{-1}\), \(h\in U(V_v)\).
Then \(\operatorname{tr}(g-1)=\operatorname{tr}(h-1)\), and
the compression gives
\(\operatorname{tr}(\mu(g-1))=\operatorname{tr}(\nu(h-1))\),
where both weight matrices are extended by zero on \(X\).
The two embeddings contain the diagonal \(U(X)\).
Let \(\theta=c_{l_Z}(g)\overline{c_{-l}(h)}\).
To check the conjugations, negate a numbered parameter \(p\) while
keeping the signed coordinate basis. This reverses its parameter order,
so \(\rho,\rho_K\), the minor exponent \(a=\lambda+\rho\), and
the lowest compact weight \(\mu=a-2\rho_K\) change sign.
For \(t=k_0Z\) in polar coordinates the compact numerator in the
lowest-type integral of Section~\ref{sec:real-coefficient-decay} is
\(m_\mu(k_0k)\overline{m_\mu(k)}\).
Complementary minors give \(m_{-\mu}(k)=\overline{m_\mu(k)}\)
on the compact group, so its squared-norm normalization is the same
for both parameters. For the positive matrix \(D=k^*Zk\),
\[
 m_{-a}(D)=m_a(D^{-1})=m_a(JDJ)=m_a(D)>0 .
\]
The real denominator is consequently identical. Thus
\(c_{-p}(t)=\overline{c_p(t)}\), and the integrand in these coordinates is
\[
 \theta=c_l(h)\overline{c_{-l_Z}(\omega_v^{-1}\iota(h)\omega_v)} .
\]

We apply the one-factor coefficient computations of
Sections~\ref{sec:real-reproducing-coefficient} and
\ref{sec:real-stationary-phase}. Here are the two-factor
localization details. Let \(E_g,E_h\) be the opposite-sign blocks,
\(v_0^2=\|E_g\|_{\mathrm{HS}}^2+\|E_h\|_{\mathrm{HS}}^2\),
\(T=D_{\rm sep}\log\beta\), and \(\eta=10^{-13}/T\).
The signed row identity, weighted by \(-\mu_i\) and \(\nu_i\),
bounds the sum of squared row norms of \(g-1\) off \(X\)
by \(6v_0^2\). Explicitly,
\[
 2\Re(g-1)_{ii}
    =2\sum_j|(E_g)_{ij}|^2-\sum_j|(g-1)_{ij}|^2
\]
and likewise for \(h\); add the weighted identities, use the
equality of weighted traces, discard the nonnegative terms from
the second squared norm, and divide by \(w_v\).

Apply the same argument to the inverse pair for the columns.
If \(v_0\le\eta\), the polar decomposition \(g_{XX}=up\)
gives \(\|p-1\|_{\mathrm{HS}}\le6\eta^2\le\eta\), by the signed
Gram identity on the \(X\)-columns. Left multiply by the diagonal
\(u^*\). Both highest vectors are fixed by \(U(X)\): the
weight-minor thresholds on this zero-weight block cancel its
determinant contribution. Thus \(\theta\) is unchanged.
The moved first squared distance to \(1\) is at most
\(13\eta^2\).

The unweighted row identity and trace equality
give the second squared distance at most \(15\eta^2\).
Both distances are at most \(4\eta\).
The diagonal \(u^*\) commutes with \(w\), since \(w\) fixes \(X\)
and preserves its orthogonal complement. Thus after moving both
arguments we still have \(g=w(h\oplus1)w^{-1}\). The upper-left
block of \(w^{-1}\mu w\) on \(V\) is \(\nu\). For every \(j\ge1\),
\[
 \operatorname{tr}\!\bigl(\mu(g-I)^j\bigr)
   =\operatorname{tr}\!\bigl((w^{-1}\mu w)((h-I)^j\oplus0)\bigr)
   =\operatorname{tr}\!\bigl(\nu(h-I)^j\bigr).
\]
The logarithm series is absolutely convergent at both moved arguments.
These weighted traces of the logarithms are purely imaginary, since
the diagonal entries of the signed unitary Lie algebra are imaginary.
The quadratic approximation in Section~\ref{sec:real-stationary-phase}
has error at most \(a/10+e^{-100T}\) from a number \(a\in[0,1]\)
after removing its imaginary log phase in each factor.

The conjugation on the second coefficient cancels the two phases.
Writing \(s=e^{-100T}\), the product of the two estimates, with real main
terms \(a,b\in[0,1]\), gives
\[
 |\theta-ab|\le0.21ab+1.1s(a+b)+s^2 .
\]
Thus \(|\theta|-2\Re\theta\le-ab+3|\theta-ab|\le12s\). If \(v_0\ge\eta\), the two one-factor decay
estimates give
\[
 |\theta|^4\le
 [(1+\|E_g\|_{\mathrm{HS}}^2)(1+\|E_h\|_{\mathrm{HS}}^2)]^{-r^7}
 \le (1+\eta^2)^{-r^7}\le e^{-400T}.
\]
For this region \(|\theta|-2\Re\theta\le3|\theta|\)
gives the same cone bound.
With \(J(h)=g\oplus h\) and \(Q=(2/\|J(h)\|)^{9D_{\rm sep}}\),
one-factor norm decay at the larger of the two norms gives
\(|\theta|\le Q^2\). Taking the geometric mean of the two bounds
on the positive part, and applying the one-factor small-ball
estimate at both matrices, gives
\[
 |\theta|\le2\Re\theta+6\beta^{-50D_{\rm sep}}Q,\qquad
 |\theta|\ge\tfrac14\quad
 \bigl(\|J(h)-1\|\le(200\beta^5)^{-1}\bigr).
\]
The volume of \(\|J(h)\|\le R\) is at most
\((3R^2/\delta)^{2D_{\rm sep}}\) times that of \(\|J(h)-1\|\le\delta\)
for \(R\ge1\), \(0<\delta\le1\).

Indeed pack \(\delta/R\)-separated matrices, and use
\(\|J(h)^{-1}\|=\|J(h)\|\) to convert each covering ball to
a translate of the \(\delta\)-ball. Let \(V_b,V_s\) be the
volumes for \(R=\beta+1\), \(\delta=(200\beta^5)^{-1}\).
Then \(V_b\le\beta^{16D_{\rm sep}}V_s\). The same packing with
\(\delta=1\), integrated in \(R\) against
\(9D_{\rm sep}\,2^{9D_{\rm sep}}R^{-9D_{\rm sep}-1}\,dR\), gives
\(\int Q\le16^{9D_{\rm sep}}V_b\). Thus
\[
 V_s/4\le 2\left|\int\theta\right|
                +6\beta^{-50D_{\rm sep}}16^{9D_{\rm sep}}V_b,\qquad
 V_b\le16\beta^{16D_{\rm sep}}\left|\int\theta\right|.
\]
For the second inequality, multiply the first by
\(4\beta^{16D_{\rm sep}}\) and absorb the resulting error into \(V_b/2\).

The \(h\)-ball of radius \(1\) about \(1\) lies in the large
ball since \(\|w\|^2\le\beta\). Under its original Haar,
write its volume as \(v_{1,v}\). The resulting point in the original first-path embedding satisfies
\[
 \left|\int_{U(V_v)}c_l(h)
    \overline{c_{-l_Z}(\omega_v^{-1}\iota(h)\omega_v)}\,dh\right|
          \ge v_{1,v}\beta^{-H}.
 \tag{DS2}\label{eq:display-ds2}
\]
The product \(\omega=(\omega_v)_v\) depends only on the arithmetic
datum and the real types. It is fixed before choosing \(y\) or a test in \(W_Z\).

\subsection{Quantitative reduction near the cusp}
\label{sec:marked-deep-cusp}

\sectionopening{In this subsection, we quantify the size of a cusp neighborhood in terms of the first
minimum of the moved integral lattice. This permits the localization
argument of Section~\ref{sec:marked-deep-separation} with a uniform
cutoff.}

Let \(W\) be \(V\) moved to the point under consideration.
Write \(n=\dim_E W\), \(f=[E^+:\mathbb Q]\),
\(D_E\) for the absolute discriminant, and
\(\mathrm{rd}_E=|D_E|^{1/(2f)}\).
Use the rational lattice \(\Lambda\) of vectors with
integral finite coordinates in the moved adelic basis,
and Euclidean lengths in this basis at infinity.
Let \(N_D>0\) be the chosen integer in the different,
and set
\[
  \kappa=8fn\,\mathrm{rd}_E+48n^6f^2N_D+1.
\]
For \(0<\delta\le1\) and a positive integer \(M\), suppose the first
minimum cut off at one satisfies
\(\lambda(W)<[\delta/(M\kappa^2)]^n\).
We construct \(0\ne A\subset W\) rational isotropic with
\[
 N_A(\mathbb A_{E^+})
  \subset N_A(E^+)B_\infty(\delta)J_W(M)
 \tag{CF1}\label{eq:display-cf1}
\]
in this basis. Here \(B_\infty(\delta)\) is the product
of the operator-norm balls about the identity, and
\(J_W(M)\) is the original finite standard principal
level \(M\) in these coordinates.

First let \(y_1,\dots,y_n\) be successively shortest
lattice vectors independent over \(E\), and let their
lengths be \(\alpha_1\le\cdots\le\alpha_n\).
These minima exist since the ordinary real lattice is
discrete and spans its ambient space. For the ordinary
real successive minima \(\mu_j\), one has
\(\alpha_i\le\mu_{2f(i-1)+1}\), by the real dimension
\(2f(i-1)\) of the archimedean extension of the span of
the first \(i-1\) vectors. Hence
\[
 \prod_i\alpha_i^{2f}\le\prod_{j=1}^{2fn}\mu_j.
\]
The Euclidean covolume is \(2^{-fn}|D_E|^{n/2}\).

In fact the covolume in the
archimedean basis multiplies the standard covolume by
the inverse of the adelic absolute value of the basis
determinant, and the latter is one by the Hermitian basis identity and the product formula for \(\det J\).
Compute the covolume from an integral basis of \(O_E\),
the local indices relative to \(O_E^n\), and the real
Jacobian of the inverse infinity basis.
Minkowski's second theorem for the Euclidean ball and
the contained cube of side \(2/\sqrt{2fn}\) give
\[
 \prod_i\alpha_i^{2f}\le(2fn\,\mathrm{rd}_E)^{fn}.
\]
In particular \(\prod_i\alpha_i\le\kappa^n\).

Put \(\theta=\delta/(M\kappa^2)\), and choose \(0\ne z\in\Lambda\)
of length less than \(\theta^n\).
The Hermitian pairing of two vectors of \(\Lambda\)
is integral in \(E\) by the finite basis identity.
If their two lengths have product less than one, the
pairing has absolute value less than one at every
embedding, so is zero. Nondegeneracy of the pairing
against the basis \(y_i\) therefore gives
\(1\le \|z\|\alpha_n\).

If \(n=1\), then \(\|z\|\alpha_n<\theta\kappa\le1\),
which is impossible. Thus \(n\ge2\).
It follows that
\[
 \prod_{i<n}\alpha_i
 \le \|z\|\prod_i\alpha_i <(\theta\kappa)^n .
\]
Pair the \(n-1\) lengths by \(i\leftrightarrow n-i\).
At least one pair satisfies
\(\alpha_i\alpha_{n-i}<\theta\kappa\), since otherwise
the square of the last product would be at least
\((\theta\kappa)^{n-1}\). Choose its index \(i\le n-i\).
Then the first \(i\) vectors span an isotropic
subspace \(A\), orthogonal to the first \(n-i\) vectors.
The latter span \(A^\perp\) by dimension. Write
\(a=\alpha_i\), \(b'=\alpha_{n-i}\), so
\(a\le b'\) and \(a b'\le\delta/(M\kappa)\).
The bases consist of rational vectors in the lattice at the point.

For the radical reduction let \(X,Y\) be the column
matrices of these rational bases of \(A,A^\perp\), and
let \(x,y\) be their adelic coordinates. Write \(J\)
for the standard adelic Hermitian form matrix in the
transported basis. For any matrix \(z_0\) with columns in
\(A^\perp\), in this basis define
\[
 {\cal R}(z_0)=1+2\bigl(
     z_0x^*-xz_0^*-xz_0^*Jz_0x^*\bigr)J .
\]
It is in the radical and satisfies
\[
 {\cal R}(w_0){\cal R}(z_0)
     ={\cal R}(w_0+z_0-xw_0^*Jz_0).
 \tag{CF2}\label{eq:display-cf2}
\]
The product follows on expanding and using
\(x^*Jx=x^*Jz_0=x^*Jw_0=0\).

For surjectivity, use a rational Witt decomposition
\(A\oplus W_1\oplus A'\) with \(X\) the first basis and
the dual basis in \(A'\). Adjoint in the original form
becomes adjoint in the Witt Gram matrix in this calculation.
The radical matrices in Witt coordinates have middle-right
block \(C\), top-middle block \(-C^*J_1\), top-right block
\(D\) with \(D+D^*+C^*J_1C=0\), and identity diagonal.
For \(z_0=(Z_1,Z_2,0)\) the corresponding blocks of
\({\cal R}(z_0)\) have
\(C=2Z_2\) and \(D=2(Z_1-Z_1^*)-\frac12 C^*J_1C\).

Take \(Z_2=C/2\) and
\(Z_1=(D+\frac12 C^*J_1C)/4\).
Thus every \(u\in N_A\) is \({\cal R}(z_0)\), and since
the columns of \(Y\) span \(A^\perp\) one may write
\(z_0=M y c\) for an adelic matrix \(c\).
These identities are in the original basis after
transporting back from the Witt coordinates.

There is a rounding function \(r\) on adelic scalars
such that \(c-r(c)\) is finite integral and has
archimedean absolute values at most \(4f^2N_D\).
To see this, approximate the finite adele by \(t\in E\)
modulo \(\widehat O_E\) by the Chinese remainder theorem.
At infinity use the lattice \(\Lambda_0=(2/N_D)O_E\) with the
real Euclidean inner product. Its dual is
\(N_D\mathfrak D_E^{-1}\subset O_E\): the Euclidean bilinear
form is half the trace pairing with conjugation.

Its nonzero
vectors have length at least one by the product formula.
For any Euclidean lattice in dimension \(d\), dual minimum
at least one implies covering radius at most \(d^2\).
Indeed Minkowski's convex body theorem applied to the lattice
and its dual, using reciprocal covolumes and an inscribed cube,
gives a nonzero lattice vector of length at most \(2d\).
Use a shortest vector \(a\), which is primitive.

Orthogonal
projection to \(a^\perp\) is a lattice since \(a\) extends to a
lattice basis, and its dual is the intersection of the original
dual with \(a^\perp\). By induction the projected remainder
has length at most \((d-1)^2\). Lift the approximant and round
its coefficient along \(a\), adding at most \(\|a\|/2\le d\)
to the error. This proves the bound.
Apply it to the archimedean residue divided by \(N_D\) and
multiply back by \(N_D\). The approximant then lies in \(2O_E\),
so it preserves finite integrality. The archimedean remainder
in each coordinate has modulus at most \(4f^2N_D\).
Apply \(r\) entrywise, so the matrix remainder has
operator norm at most \(4n^2 f^2N_D\) at each real factor.

Choose \(d=-M(y r(c))^*J y c\) in adelic coordinates.
The entries of \(r(c)\) and \(r(d)\) in the expression
\(w_0=M(xr(d)-yr(c))\) are embedded rational
coefficients; thus \({\cal R}(w_0)\) corresponds to
\(\gamma\in N_A(E^+)\).
Here \(d\) may be nonrational. The rational coefficient is \(r(d)\).
In fact \(x^*Jy=0\) gives
\[
 w_0^*Jz_0=-M^2r(c)^*y^*Jyc=Md,\qquad
 w_0+z_0-xw_0^*Jz_0=M\{y(c-r(c))-x(d-r(d))\}.
\]
Thus the central cross term in
\eqref{eq:display-cf2} contributes only its rounded remainder. Substitution yields
\[
 \gamma u=1+2(MT-(MT)^*-(MT)^*J(MT))J,\qquad
 T=(y(c-r(c))-x(d-r(d)))x^* .
\]
The finite entries of \(T\) are integral. The operator
norms of \(x,y\) at infinity are at most \(n^2a,n^2b'\),
so
\[
 \|T_v\|\le4n^6f^2N_D(b'+a)a
        \le8n^6f^2N_Da b'\le\delta/(6M).
\]
Since every \(J_v\) is a sign involution,
\(\|(\gamma u)_v-I\|\le 2\delta/3+\delta^2/18\le\delta\).
The finite formula gives \(\gamma u\equiv I\bmod M\).
This proves \eqref{eq:display-cf1} in the original basis at level \(M\).

\subsection{Vanishing of the constant term under the reproducing kernel}\label{sec:deep-constant-term}

\sectionopening{In this subsection, the reproducing kernel of
Section~\ref{sec:real-reproducing-coefficient} annihilates the
relevant constant term. We prove this pointwise and apply it on the
neighborhood just constructed in Section~\ref{sec:marked-deep-cusp}.}

Now let \(e\) be the small-type reproducing convolution
kernel on the real group \(G_\infty\) of \(V\).
It is integrable, idempotent, and has an integrable
first norm moment, by the prescribed large gap.
Write \(Tf(g)=\int e(u)f(gu)\,du\).
We prove the local vanishing for bounded continuous
functions invariant under the product of the radicals attached
to nonzero isotropic subspaces \(A_v\) of the real factors.
This applies also to their real conjugates. The proof uses
the idempotence and the indicated moment.

For the product of real groups, let \(\mathcal D\)
consist of pairs \((a,c)\) with \(c>0\) such that each
\(a_v\) multiplies \(A_v\) by \(c\) and acts as the
identity on \(A_v^\perp/A_v\). Its kernel over \(c=1\)
is \(N_{A,\infty}\).
Each real factor of \(A\) is nonzero.
Fix one place \(v_0\) and a vector \(0\ne x_{v_0}\in A_{v_0}\).
Put \(H(g)=\log\|g_{v_0}^{-1}x_{v_0}\|\).
Then \(H(ag)=H(g)-\log c\) for \((a,c)\in\mathcal D\)
and
\[
 |H(gu)-H(g)|\le \|u_{v_0}\|.
\]
The latter follows by submultiplicativity in both
directions and \(\|u_{v_0}^{-1}\|=\|u_{v_0}\|\).

For each \(g\) and each \(\epsilon>0\), take at each
place the Euclidean orthogonal projection \(P\) onto
\(g_v^{-1}A_v\), and put \(Q=J_vP J_v\).
The two projections have orthogonal images by
isotropy. The operator
\(b_v=I+\epsilon P-\epsilon Q/(1+\epsilon)\)
is unitary for \(J_v\), acts by \(c=1+\epsilon\)
on \(g_v^{-1}A_v\), and acts trivially on its
orthogonal quotient for \(J_v\).
Its distance from \(I\) is at most \(\epsilon\);
hence \((g b g^{-1},c)\in\mathcal D\).
Translation is continuous on \(L^1(G_\infty)\).
Choose \(\epsilon\) so small that
\(\|e(b^{-1}\cdot)-e\|_1<1/2\)
throughout this neighborhood of \(I\).

Let \(\mathcal B\) consist of the bounded continuous
functions fixed by \(T\) and left invariant under
\(N_{A,\infty}\), with the supremum norm.
For a fixed \((a,c)\in\mathcal D\) with \(c=1+\epsilon\),
left translation \(Uf(g)=f(ag)\) preserves
\(\mathcal B\). Given \(g\), the element \(b\)
just constructed satisfies
\(g b g^{-1}a^{-1}\in N_{A,\infty}\), hence
\(f(ag)=f(gb)\). Reproduction and the \(L^1\) translation
bound give \(\|I-U\|\le1/2\) on \(\mathcal B\).
Put \(\tau=\log c\), take a large positive integer \(k\),
and set \(\chi(g)=\exp(-i\pi H(g)/(k\tau))\) and
\(Sf=T(\chi f)\).

The multiplier is invariant under
\(N_{A,\infty}\), \(T\) commutes with left translations
and \(T^2=T\), so \(S\) preserves \(\mathcal B\).
The shift in \(H\) gives \(US=e^{i\pi/k}SU\).
It follows that \(US^k=-S^kU\); the operator identity
\(2S^k=(I-U)S^k+S^k(I-U)\) and \(\|I-U\|\le1/2\)
then give \(S^k=0\).
On the other hand \(Sf=0\) implies
\[
 \|f\|_\infty\le
 \frac{\pi}{k\tau}
    \left(\int |e(u)|\,\|u_{v_0}\|\,du\right)\|f\|_\infty,
\]
by subtracting \(\chi(g)Tf(g)\) in the integral.
The moment is finite and the coefficient is less than
one for large \(k\). Thus \(S\) is injective, and
\(S^k=0\) proves \(\mathcal B=\{0\}\).
In particular \(T\) kills every bounded continuous
left-\(N_{A,\infty}\)-invariant function.

In original coordinates let \(N\) be the radical of the rational
isotropic subspace just found at \(g\). In the transported coordinates
the two groups are \(N^g=g^{-1}N(\mathbb A_{E^+})g\) and
\(\Gamma_N^g=g^{-1}N(E^+)g\). For \(h\in G_\infty\) take the mean
of \(F(gu(h,1))\) over \(\Gamma_N^g\backslash N^g\).
It is well defined because for \(\gamma\in N(E^+)\) the argument at
\(g^{-1}\gamma g u\) is \(\gamma gu(h,1)\).
The quotient is compact by \eqref{eq:display-cf1}: its representatives
lie in the compact intersection of the radical with
\(B_\infty(\delta)J_W(M)\).
Write this probability mean as \(F_g(h)\). For
\(n\in g_\infty^{-1}N(E^+_\infty)g_\infty\), the change
\(u\mapsto u(n,1)\) on \(\Gamma_N^g\backslash N^g\) gives
\(F_g(nh)=F_g(h)\). It is bounded and continuous.
Since \(\|F\|_\infty\|e\|_1<\infty\), Fubini and reproduction
of \(F\) give \(TF_g=F_g\). Apply the real argument to the
isotropic spaces \(g_v^{-1}A_v\) in the moved real bases to get
\(F_g=0\). Since \(F\) is automorphic and right
invariant under \(J_W(M)\), \eqref{eq:display-cf1} makes every radical
value differ from \(F(g)\) by no more than the
oscillation on the real \(\delta\)-ball.
Averaging against the probability normalization
of the same quotient Haar measure gives
\[
  |F(g)|\le\sup_{\|h-1\|\le\delta}|F(g(h,1))-F(g)|.
\]
This is the required bound at the cusp.

\subsection{Removing the smoothing at the period}
\label{sec:deep-removal-of-smoothing}

\sectionopening{In this subsection, we combine the constant term vanishing of Section~\ref{sec:deep-constant-term} with the quantitative reduction of Section~\ref{sec:marked-deep-cusp}.}

Throughout this subsection write \(A=A_{\rm sep}=r^4\) and
\(L=L_{\rm sep}=\lceil20b\rceil\) as in \eqref{eq:deep-reproduction-scale}.
Let \(G_V,G_Z\) and \(\Gamma_V,\Gamma_Z\) be the adelic and
rational unitary groups in the given bases. Write \(\jmath=\iota\),
\(c_v=c_l\) and \(c_v'=c_{-l_Z}\) at place \(v\). Put
\(e(u)=\prod_v d_v\overline{c_v(u_v)}\), with \(d_v\) the formal degree.
Assume the orthogonality in \eqref{eq:deep-separation-condition}. For an adelic point \(a\)
of \(G_V\) and an infinite point
\(\alpha\), put
\[
 \begin{gathered}
 \Phi_{\alpha,\omega}(g)=\prod_v\overline{c'_v(\alpha_v^{-1}g_v\omega_v)},\\
 \psi_\alpha=\Phi_{\alpha,\omega}\circ\jmath,\qquad
 (R_\alpha y)(a)=\int\psi_\alpha(k)y(a(k,1))\,dk .
 \end{gathered}
\]
For clarity, fix the center \(q=\iota(a)(\alpha,1)\) in \(G_Z\).
The periodization of the conjugate coefficient kernel times
\(\mathbf1_{K_Z(a_0)}\), at arguments \(q^{-1}\gamma g\), is
the conjugate of a bounded form of the required type and level,
divided by the fixed product of formal degrees for normalized
coefficients.

To verify this, restrict the supported terms
\(a_\gamma=q^{-1}\gamma g\) to \(\|a_{\gamma,v}\|\le R_v\)
at every real place, \(R_v\ge1\). For a fixed one of these
terms \(a_{\gamma_*}\), the matrices
\(a_\gamma a_{\gamma_*}^{-1}=q^{-1}\gamma\gamma_*^{-1}q\)
have integral finite coordinates in the center basis.
They are distinct points of the fixed moved rational matrix
lattice and have real norms at most \(R_v^2\).
Packing its disjoint balls of fixed radius less than half its
first minimum gives at most
\(C_q\prod_v R_v^{4(r+1)^2}\) such points. Indeed the balls
fit in the product of local balls whose radii are \(R_v^2\)
plus the fixed packing radius.

Dyadic summation converges uniformly in \(g\), since the
coefficient decay exponent \(2A\) exceeds \(4(r+1)^2\).
It converges locally uniformly term by term by increasing
the real radii by a fixed factor on a compact neighborhood
and restricting the finite neighborhood to \(gK_Z(a_0)\).

Reindexing gives automorphy and right level invariance.
For the type, conjugate and replace \(\gamma\) by
\(\gamma^{-1}\). The conjugate uses arguments \(g^{-1}\gamma q\).
The idempotent kernel \(e_Z=\prod_v d'_v\overline{c'_v}\) gives
\[
 \int e_Z(u)e_Z(u^{-1}g^{-1}\gamma q)\,du
       =e_Z(g^{-1}\gamma q)
\]
at the real places. Boundedness of the absolute sum and
integrability of \(e_Z\) justify convolution term by term.
Thus the conjugate expression is fixed by right reproduction.
Inserting \(g=\iota(h)(\omega,1)\) in this test in \(W_Z\)
and conjugating the vanishing hypothesis of \eqref{eq:deep-separation-condition} gives
zero for the period against \(y(h)\). Its kernel is exactly
\(\Phi_{\alpha,\omega}((\iota(a)^{-1}\gamma\iota(h))_\infty)\)
times the finite level indicator.

For a fixed vector \(y=R(e)y\) in \(L^2\), let
\(y_j=y\mathbf1_{\{|y|\le j\}}\). Right convolution satisfies
\[
 \|R(e)y_j\|_\infty\le j\|e\|_1,\qquad
 \|R(e)y_j-y\|_2\le\|e\|_1\|y_j-y\|_2\longrightarrow0 .
\]
Each projected truncation is continuous by continuity of translation
in \(L^1\), retains the finite level, and is fixed by \(R(e)\).
Section~\ref{sec:paired-lifts-projected-tests} gives cuspidality.
At the fixed level, compact type and infinitesimal character, the
finiteness argument of Section~\ref{sec:four-unitary-real-projections}
places all the projected truncations in one finite-dimensional cusp
space. It is closed in \(L^2\), so it contains \(y\). For this
bounded continuous representative, \(R(e)y=y\) everywhere. Both
sides are continuous and Haar measure has full support on the quotient.
This applies to the initial true grid projector at its level as well.

Take a global bound \(|y|\le S\) for this representative and set
\(\lambda(a)\) equal to the first minimum, cut off at
\(1\), in the lattice of \(V\) moved to \(a\).
If \(B_d^2\le B_d^{3r+2}\lambda(a)\), unfold the period by
the cosets \(\Gamma_Z/\Gamma_V\) of the embedded
rational subgroup (equivalently by the orbits of right
multiplication). Put \(e=e_{\rm last}\),
\(K=K_Z(a_0)\), and \(b_\zeta=\jmath(a_f)^{-1}\zeta_f\).
The finite support condition for the coset \(\zeta\Gamma_V\) is
\[
 \exists h_f\in G_{V,f}:\ b_\zeta\jmath(h_f)\in K
       \quad\Longleftrightarrow\quad b_\zeta e\in K e .
\]
Indeed \(\jmath(h_f)e=e\). Conversely, if \(b_\zeta e=k e\)
with \(k\in K\), then \(k^{-1}b_\zeta\) fixes \(e\) and preserves
\(e^\perp\), so it lies in \(\jmath(G_{V,f})\). Its inverse gives
a suitable \(h_f\). This orbit condition implies the column
congruence used below at the chosen split places. The coset itself
is determined by the rational column \(\zeta e\).

Equal columns mean \(\zeta^{-1}\zeta'\) fixes the
last unit vector and belongs to \(\Gamma_V\).
Write \(c(\zeta)\) for this column in the Euclidean real
coordinates moved to \(\jmath(a)\).
The finite coordinates of every supported column are
\(e_{\rm last}\) modulo \(a_0\). Hence the difference
of two supported columns, divided by \(a_0\), is in the
moved rational lattice of \(Z\). More explicitly, if
\(L_{V,f}\) is the original finite lattice (including its
chosen basis at the places over \(N\)), then the rational vectors
with finite coordinates in \(a_f L_{V,f}\) form \(\Lambda_V(a)\).
At \(\iota(a)\) this lattice for \(Z\) is
\(\Lambda_V(a)\oplus O_E e_{\rm last}\). At infinity the sum
is orthogonal in the Euclidean coordinates moved to \(\iota(a)\).
Each nonzero integer of \(E\) has modulus at least one at some
embedding, so the minimum of this sum is at least \(\lambda(a)\).
Multiplication by \(a_0\) scales each archimedean length by
\(M^{b_0/2}\ge B_d^{3r+2}\).

Consequently
\[
 \|c(\zeta)-c(\zeta')\|\ge B_d^{3r+2}\lambda(a)\ge B_d^2
          \quad(\zeta\Gamma_V\ne\zeta'\Gamma_V),
 \qquad
 \|c(1)\|=\sqrt{f}\le B_d .
 \tag{D3}\label{eq:display-d3}
\]
The last column does not vary under multiplication
on the right by \(\jmath(h)\).
At any archimedean point of a contributing term,
write \(g_v\) for the local group element before the two
translations. Its last column is \(c(\zeta)_v\).
Since \(\|g_v\|\ge1\) and \(\|c(\zeta)\|^2
=\sum_v\|c(\zeta)_v\|^2\), one has
\[
 \prod_v\|g_v\|\ge\max_v\|g_v\|
 \ge \|c(\zeta)\|/\sqrt f\ge \|c(\zeta)\|/B_d .
\]

For clarity, the \(L^1\) slice estimate follows from
decay and a volume calculation. For \(R\ge1\) and a
fixed ambient group element \(q\), take a maximal
\(1/R\)-separated set in the subset of matrices
\(\jmath(h)q\) of operator norm at most \(R\).
Its size is at most \((1+2R^2)^{2(r+1)^2}\) by packing
in the real vector space of complex matrices, using
the operator norm. Two such matrices \(a,b\) at distance
at most \(1/R\) satisfy \(\|ab^{-1}-1\|\le1\) because
\(\|b^{-1}\|=\|b\|\le R\).

Thus the measure of the
original subset of \(h\)'s is at most this same number
times \(v_{1,v}\), by right translation in the subgroup.
Sum over \(2^j\le\|\jmath(h)q\|<2^{j+1}\), \(j\ge0\).
With \(d_{\rm pack}=2(r+1)^2\) the sum for the weight
\((2/\|\jmath(h)q\|)^A\) is at most
\[
 v_{1,v}\,2^A9^{d_{\rm pack}}\sum_{j\ge0}2^{-j(A-2d_{\rm pack})}
        \le v_{1,v}\,16^A,
\]
using \(A>4d_{\rm pack}\). The corresponding bound on the
left-translated slice follows by inversion of the
unitary group and invariance of its Haar measure.

Write
\[
  C_{\alpha,\omega}=\prod_v
       [2(\|\alpha_v\|\|\omega_v\|)^2]^A16^A v_{1,v} .
\]
Here and in the following estimates the unit-ball volume in each
real-place factor is its own \(v_{1,v}\). The local large-type bound
\(\lvert c'_v(g)\rvert\le(2/\|g\|)^{2A}\)
and submultiplicativity give one power of degree \(A\)
in the slice for integration as above, and one power
for \(\|c(\zeta)\|/B_d\). For the latter, use the product
of the operator norms before applying the translations
\(\alpha^{-1}\) and \(\omega\).

For a supported coset, the set of finite \(h\) in its
term is \(h_0K_V(a_0)\). Indeed, after choosing one
\(h_0\) in the support, \(b_\zeta\jmath(h_f)\in K\) is equivalent to
\(\jmath(h_0^{-1}h_f)\in K\), and \(\jmath^{-1}(K)=K_V(a_0)\).
Replacing \(\zeta\) by \(\zeta\jmath(\eta)\), \(\eta\in\Gamma_V\),
replaces this coset by \(\eta_f^{-1}h_0K_V(a_0)\), of the same measure.
Thus its total absolute integral is at most
\[
 \operatorname{vol}(K_V(a_0))\,C_{\alpha,\omega}S
                    (B_d/\|c(\zeta)\|)^A
              \qquad(\zeta\Gamma_V\ne\Gamma_V).
\]
This retains the volume in the unfolded main term.

For the supported cosets in \eqref{eq:display-d3}, all nonbase norms and pairwise distances
are at least \(R_0=B_d^2/2\). In dimension \(2f(r+1)\)
the count in
\(2^j R_0\le\|c(\zeta)\|<2^{j+1}R_0\) is at most
\((5\cdot2^j)^{2f(r+1)}\), by packing balls of radius
\(R_0/2\). With \(d_c=2f(r+1)\) and \(A-d_c\ge2\), this gives
\[
 \sum_{\substack{\zeta\Gamma_V\ne\Gamma_V\\ b_\zeta e\in Ke}}
                (B_d/\|c(\zeta)\|)^A
 \le (2/B_d)^A\,\frac{5^{d_c}}{1-2^{d_c-A}}
 \le (10/B_d)^A\le 10/B_d .
\]
For the absolute unfolding, let \(T_\zeta(h)\) be the term
at \(\zeta\), including \(y(h)\) and the finite indicator.
Then \(T_{\zeta\jmath(\eta)}(h)=T_\zeta(\eta h)\) for
\(\eta\in\Gamma_V\) by automorphy of \(y\). With one representative
per supported orbit, Tonelli and counting measure give
\[
 \int_{[U(V)]}\sum_{\gamma\in\Gamma_Z}|T_\gamma(h)|\,dh
 =\sum_{\zeta\Gamma_V\ {\rm supported}}
            \int_{G_V}|T_\zeta(h)|\,dh .
\]
The nonbase sum is finite by the shell estimate and the base term
by the same real slice estimate with the finite compact indicator.
This proves the signed unfolding with multiplicity one.
For the base coset use \(h=ak\); the finite indicator is
then that of \(K_V(a_0)\) and \(y\) is invariant under
this group since \(\pi_M\) divides \(a_0\).
The base coset contributes precisely
\(\operatorname{vol}(K_V(a_0))\,R_\alpha y(a)\).
For each supported nonbase coset choose a finite solution
\(h_{\zeta,f}\) to its support condition. After dividing by the same
finite volume its exact contribution is
\[
 I_{\zeta,\alpha}(a)=\int_{G_{V,\infty}}
    y(h,h_{\zeta,f})\,
    \Phi_{\alpha,\omega}
     \bigl(\jmath(a_\infty)^{-1}\zeta_\infty\jmath(h)\bigr)\,dh .
\]
Indeed every finite argument of this coset is in
\(h_{\zeta,f}K_V(a_0)\), and \(K_V(a_0)\) fixes \(y\).
Replacing \(\zeta\) by \(\zeta\jmath(\eta)\) changes \(h_{\zeta,f}\)
to \(\eta_f^{-1}h_{\zeta,f}\) up to this level. The substitution
\(h\mapsto\eta_\infty^{-1}h\) then gives the same contribution by
rational left invariance. Hence the vanishing period is exactly
\[
 R_\alpha y(a)+\sum_{\zeta\Gamma_V\ne\Gamma_V\ {\rm supported}}
                              I_{\zeta,\alpha}(a)=0 .
\]
The preceding absolute bounds give
\[
      |R_\alpha y(a)|\le C_{\alpha,\omega} S(10/B_d)
       \qquad\text{if }B_d^2\le B_d^{3r+2}\lambda(a).
 \tag{D4}\label{eq:display-d4}
\]

\indent{\it Removing the smoothing.}

Put
\(\mathfrak M(u)=\prod_v d_v(2/\|u_v\|)^{2A}\).
The integral of each real factor is at most
\(8[2000(r+L+1)^6]^{2r^2}2^{2A}\).
The Casimir calculation in
Section~\ref{sec:real-reproducing-coefficient} gives a
translation Lipschitz bound \(100(r+L+1)^5\) for each
small-type coefficient. In detail, for a unit lowest-type
vector and any real Lie algebra element \(Z\) it gives
\(\|d\pi(Z)v\|\le\sqrt{3r}(r+L+1)\|Z\|_{\rm HS}\).
For \(\|u-1\|\le1/2\) apply this to \(Z=\log u\)
and integrate along \(\exp(tZ)\), using
\(\|Z\|_{\rm HS}\le2\sqrt r\|u-1\|\).

For larger distances the coefficient difference is at most
\(2\). The bound is uniform in the other argument by unitarity.
Now if \(K_0=2000(r+L+1)^6\), the small-type
coefficient has absolute value at least \(1/2\)
in the ball \(\|u-I\|\le 4/K_0\) by its
Lipschitz bound \(100(r+L+1)^5\).
The volume of that ball is at most \(4/d_v\) by
the formal degree identity.
Packing the ambient matrices as above with gap
\((4/K_0)/R\) gives
\(\mu\{\|u\|\le R\}\le (K_0R^2)^{2r^2}(4/d_v)\)
for \(R\ge1\), since \(1+K_0 R^2/2\le K_0R^2\).
Layer integration gives the claimed bound:
\(2A/(2A-4r^2)\le2\).
In particular \eqref{eq:deep-reproduction-scale} gives
\[
   \int\mathfrak M(u)\,du
   \le\frac{B_d}{16\beta\,[\beta^H(2\beta)^{3A}]^f}.
 \tag{D5}\label{eq:display-d5}
\]
For this use
\(1+2r^2+2A+H+6A\le20A\), \(16\le\beta\), and
the product measure of the real factors.

Suppose \(\lambda(a)\ge(B_d^{-2}/(M\kappa^2))^{r}\).
Since \(M\kappa^2\le\beta\le B_d\), the hypothesis in
\eqref{eq:display-d4} follows from
\[
 \lambda(a)\ge B_d^{-3r},
 \qquad B_d^{3r+2}\lambda(a)\ge B_d^2.
\]
For a real \(u\),
change variables on the subgroup to obtain
\[
 (R_\omega y)(a(u,1))=(R_{\jmath(u)\omega}y)(a).
\]
Average this equality against \(e(u)\,du\). Put
\[
 T_e y(a)=\int e(u)y(a(u,1))\,du,\qquad R_k y(a)=y(a(k,1)).
\]
For a bounded form the two convolutions
in \(T_eR_kT_e\) have absolute integral at most
\(\|y\|_\infty\|e\|_1^2\). Their kernel at \(z\) is
\[
 \int e(u)e(k^{-1}u^{-1}z)\,du
 =e(z)\prod_v c_v(k_v)
 =e(z)\overline{e(k)}/e(1)
\]
by the coefficient identity of Section~\ref{sec:real-reproducing-coefficient}
with the formal degrees for these same local Haar measures. Since
\(T_ey=y\), Fubini gives
\[
 \begin{aligned}
 \int e(u)(R_\omega y)(a(u,1))\,du
  &=y(a)\int\psi_\omega(k)\overline{e(k)}/e(1)\,dk
     =\nu\,y(a),\\
 |\nu|&\ge\beta^{-fH}\prod_v v_{1,v}.
 \end{aligned}
\]
Indeed the integral in \(u,k\) is bounded absolutely by
\(\|y\|_\infty\|e\|_1\|\psi_\omega\|_1\), where the slice estimate
above gives \(\psi_\omega\in L^1\). The scalar is exactly
\(\prod_v\int c_l(k)\,\overline{c_{-l_Z}(\omega_v^{-1}\iota(k)\omega_v)}\,dk\)
from \eqref{eq:display-ds2}.
At \(\alpha=\jmath(u)\omega\),
\[
 |e(u)|\,C_{\jmath(u)\omega,\omega}
       \le\mathfrak M(u)(2\beta)^{3Af}\prod_v v_{1,v}.
\]
Indeed \(\|\iota(u_v)\omega_v\|\,\|\omega_v\|
\le\|u_v\|\beta\). After division by the local factor of
\(\mathfrak M(u)\,v_{1,v}\), the product on the left is at most
\[
 (2/\|u_v\|)^{2A}
       [32\beta^2\|u_v\|^2]^A=(128\beta^2)^A
       \le(2\beta)^{3A},
\]
by the \(4A\)-decay and \(\beta\ge16\).
Consequently the averaged unfolding at this \(a\) is
\[
 \nu\,y(a)=-\sum_{\zeta\Gamma_V\ne\Gamma_V\ {\rm supported}}
     \int e(u)\, I_{\zeta,\jmath(u)\omega}(a)\,du .
\]
The sum of the absolute integrals is at most
\((10S/B_d)\int |e(u)|C_{\jmath(u)\omega,\omega}\,du\)
by the coset bounds. This is finite by \eqref{eq:display-d5},
justifying the exchange. Dividing by the lower bound for \(|\nu|\)
gives \(|y(a)|\le 10S/(16\beta)\le S/2\).

If instead \(\lambda(a)<(B_d^{-2}/(M\kappa^2))^{r}\),
use the quantitative inclusion of Section~\ref{sec:marked-deep-cusp} at radius
\(B_d^{-2}\) and level \(K_V(M)\) on \(V\).
It gives a rational isotropic subspace whose unipotent
radical, read in the basis at \(a\), is covered by its
rational translates of the product of the real ball
of radius \(B_d^{-2}\) and \(K_V(M)\).
In the moved coordinates, for each \(u\in N^a\) choose
\(\gamma\in\Gamma_N^a\) such that
\(\gamma u=(u'_\infty,k_f)\), with \(\|u'_\infty-I\|\le B_d^{-2}\)
at each real place and \(k_f\in K_V(M)\). Then
\(a\gamma a^{-1}\) is rational. Since \(K_V(M)\subset J_M\),
\[
 y(au)=y(a\gamma u)=y\bigl(a(u'_\infty,1)\bigr).
\]
For the level inclusion, at a chosen place over either old prime
\(M=\pi_M\bar\pi_M\)
has valuation one and \(J_M=K_V(\pi_M)\). Away from the old pairs
\(J_M\) is the prescribed integral compact. The larger modulus
\(\pi_M^{b_0}\) in \eqref{eq:display-d3} is imposed on the target
column. By the constant-term vanishing of Section~\ref{sec:deep-constant-term},
\[
 |y(a)|\le \sup_{\|u-1\|\le B_d^{-2}}
                    |y(a(u,1))-y(a)|.
\]
Telescoping the product of the small-type coefficients
at \(u^{-1}z\) and \(z\) bounds its difference by
\[
 100(r+L+1)^5\sum_v\|u_v-I\|
 \le 100 f(r+L+1)^5 B_d^{-2}\le\beta^2 B_d^{-2},
\]
since each factor has modulus at most one. With
\(P(z)=\prod_v(2/\|z_v\|)^{2A}\), the same difference
is at most \(P(u^{-1}z)^2+P(z)^2\) by the \(4A\)-decay.
Use \(\min(U,V)\le\sqrt{UV}\) and
\(\sqrt{P_1^2+P_2^2}\le P_1+P_2\) for nonnegative \(P_i\).
Including the formal degrees gives
\[
 |e(u^{-1}z)-e(z)|\le(\beta/B_d)
                    [\mathfrak M(u^{-1}z)+\mathfrak M(z)] .
\]
Applying reproduction of \(y\) and left invariance
of Haar measure, then \eqref{eq:display-d5}, bounds the oscillation by
\(2S(\beta/B_d)\int\mathfrak M\le S/8\).
This proves \(|y(a)|\le S/2\) for all adelic points \(a\).
Repeated halving of any initial bound for \(y\) gives
\(y=0\).

Compose the deep restriction at this \(\omega\) with the
weighted transpose of the true deep lifts from the definite
end of the grid. Section~\ref{sec:deep-measured-square} deduces
\(\Theta^\flat(\mathcal H_m)=\mathcal H_0=W_Z\) from
\eqref{eq:display-tp4} and injectivity of the forward composition.
Thus for each nonzero \(y\in W\) the separating test in
\eqref{eq:deep-separation-condition} is \(\Theta^\flat\eta\)
for some \(\eta\in\mathcal H_m\). The level \(K_Z(a_0)\) and tuple
\(\omega=(\omega_v)_{v\mid\infty}\) are fixed, and the choice of
\(\eta\) may depend on \(y\).

\section{Theta compositions at deeper principal levels}
\label{sec:deeper-periods}

\sectionopening{In this section, we transfer the separation of
Section~\ref{sec:marked-deep-separation} along two longer theta
compositions. The paired period of Section~\ref{sec:marked-square-analytic-pairing}
compares them after identifying the opposite Hermitian forms and gives
the terminal pairings used in Section~\ref{sec:chain-determinants-completion}.}

Put \(\mathfrak b=(\pi_M^{b_0})\) and \(R_b=\mathcal O_E/\mathfrak b\).
Use the trivial twist in the true triangle, with this same ring
at every node. Its rank in column \(n\) is
\(d_n=r+1+n(t-1)\). The identifications \(j_i\) of its base copies
come from passage to the opposite Hermitian form, with the original
finite bases, quotient Haar measures, and principal levels.
Use forward triangle forms of the conjugate type to the true-grid
forms \(\mathcal H_j\) of level \(K_{Z_j}(\mathfrak b)\). In
particular \(\mathcal H_0\) consists of the tested forms \(W_Z\)
of Section~\ref{sec:marked-deep-separation}. Write \(a_{i,n}\) for
the forward point lift and \(w_n=|\mathrm{GL}_{d_n}(R_b)|\).

\subsection{The first-column return}
\label{sec:deep-first-column-return}

\sectionopening{In this subsection, we compute the first return on the level used in Section~\ref{sec:marked-deep-separation}.}

At the first arrow put \(u=d_0=r+1\), \(v=d_1=r+t\) and
\(X_0=\binom{I_u}{0}\) over \(R_b\). The characters at the
two ranks are \(\mu^0_{i+1}\upsilon^{\lceil u/2\rceil}\)
and \(\mu^0_{i+1}\upsilon^{\lceil v/2\rceil}\).
They are unramified off the chosen split pairs. The finite
datum is the indicator of \(X_0\) at the chosen places, with
the integral vacuum elsewhere. At a chosen split place suppose
both \(X\) and \(Xg\) reduce to \(X_0\). Their first \(u\) rows
give square matrices congruent to \(I_u\), hence invertible over
the local integers. Composing the projection onto those rows with
the inverse square matrix gives an integral left inverse for each.

If \(L_X,L_{Xg}\) are these left inverses, then
\(g=L_X(Xg)\) and \(g^{-1}=L_{Xg}X\) are integral there.
The overlap is then equivalent to \(g\equiv I_u\bmod\mathfrak b\),
and the character of \(g\) on the overlap is one.
Give the integral matrix box and the integral source compact
\(K_Z^\circ\) mass one. The squared norm of the residue indicator is
\(|R_b|^{-uv}\), and its source stabilizer has mass
\[
 \operatorname{vol} K_Z(\mathfrak b)=w_0^{-1},\qquad
 w_0=[K_Z^\circ:K_Z(\mathfrak b)]
   =\prod_{w\mid\pi_M}|\mathrm{GL}_u(\mathcal O_{E,w}/
                                    \mathfrak b\mathcal O_{E,w})| .
\]
Reduction is surjective at each chosen split factor, since every
matrix with unit determinant modulo the maximal ideal lifts to an
integral invertible matrix. Unitarity determines the conjugate factor.
The separate test space \(W\) in
Section~\ref{sec:marked-deep-separation} has level \(K_V(\pi_M)\).
The return uses the conjugate of \(W_Z\) at level
\(K_Z(\mathfrak b)\).

Write \(A_i=a_{i,0}j_i\). The open-graph doubling calculation of
Section~\ref{sec:marked-first-arrow-rallis} applies at these ranks.
The doubled source has Siegel rank \(u=r+1\). The target has
rank \(v=r+t=6r\), so \(v/2-u=2r-1>0\) and
\(v-\iota(T)\ge v/2>u\). Weil's convergent formula applies to the
separately translated Schwartz--Bruhat factors
\cite[\S52, Theorem~5, pp.~76--77]{WeilSiegel}. Its parameter is \((v-u)/2\)
and the coefficient is one \cite[Theorem~2]{YamanaSW}.
At each prescribed split place every supported matrix has an
integral left inverse as above, hence full column rank over the local
field. A nongraph orbit would force a rank drop there, so it contributes
zero to the global doubling integral. This gives
\[
 w_0\langle A_i x,A_i y\rangle
   = C_i|\operatorname{Mat}_{v\times u}(R_b)|^{-1}
            \langle x,\mathcal V_i y\rangle,\qquad C_i\ne0 .
 \tag{TP1}\label{eq:display-tp1}
\]
Here \(x,y\) are in the true forward base
\(\overline{\mathcal H_0}\) at level \(K_{Z_0}(\mathfrak b)\).
The projected truncations in Section~\ref{sec:deep-removal-of-smoothing}
supply bounded cuspidal representatives throughout this base before
the doubling. These representatives are rapid by
Section~\ref{sec:paired-lifts-projected-tests}. At the fixed finite
test, Section~\ref{sec:marked-graph-doubling} bounds the absolute
theta-product integral and the absolute graph sum separately by
\(C H(g)^B H(h)^B\). Multiplication by \(|x(g)y(h)|\) and
integration over the two finite-volume source quotients justify the
same rearrangements as in Section~\ref{sec:first-lift-doubling-series}.
Thus the identity applies to the full deep base.

The scalar \(C_i\) is the original target quotient volume times
the integral of the real joint-harmonic coefficient against
\(\overline{e_i(h)}/e_i(1)\), with \(e_i\) the source projector.
We verify the real harmonic at this edge. Fix an infinite
place, let \(L\) be the doubled forward stem
\(2\lambda_j+1\) of Section~\ref{sec:marked-deep-separation}.
The doubled target parameter \(L'\) adds to this stem
the block \(-2w_{i+1}+t-2-2b\), \(0\le b<t-1\),
with positive signs. Here
\(w_j=\kappa_\tau+(j+1/2)\Delta\) in the indexing
of the grid and \(w=w_{i+1}\). The source and target
splitting weights are \(u/2-w,v/2-w\).
Both cross shifts are \(2w\). At tail stem index \(k-1+a\)
the common shifted entry is
\[
 \xi_a=\Delta-2w_{a+1}+1+2w
       =\Delta+1-2(a-i)\Delta .
\]
The ratio of the two signature signs there is
\((-1)^{\min(1,(a-i)_+)}=\operatorname{sgn}\xi_a\).
The head has positive shifted entries and unchanged
signs. The \(v-u=t-1\) inserted entries have positive
signature and shifted values
\(t-2,t-4,\ldots,2-t\). The common shifted entries
have modulus at least \(\Delta-1\). The gap between
opposite signature signs is therefore at least
\((\Delta-t-2)/2>500(u+v)^3+uv\) in the undoubled
parameters.

Use \(s_a,s'_b\) for the source and target signs.
In the oscillator the doubled source parameter is
\(-L\) and the target one is \(L'\).
After subtracting the vacuum exponents, the lowest
compact relative weights are
\[
 \begin{aligned}
 2R^{\rm wt}_b&=\xi_b-s'_b\left(\sum_{b'}s'_{b'}\operatorname{sgn}
                (\xi_b-\xi_{b'})+\sum_{a'}s_{a'}\right),\\
 2C^{\rm wt}_a&=-\xi_a+s_a\left(\sum_{a'}s_{a'}\operatorname{sgn}
                (\xi_a-\xi_{a'})+\sum_{b'}s'_{b'}\right).
 \end{aligned}
\]
Indeed apply \(\lambda+\rho-2\rho_K\) to the two
parameters and subtract the respective vacuum weights
\(-w+s'_b\sum_a s_a/2,\ w-s_a\sum_b s'_b/2\).
At matched indices, the contributions of an index
with the same sign of \(\xi\) cancel using
\(s'_a=s_a\operatorname{sgn}\xi_a\). For the opposite
sign they cancel using
\(\operatorname{sgn}(\xi_a-\xi_b)=\operatorname{sgn}\xi_a\).
The inserted-block contributions also cancel at
matched indices. At an inserted index the matched
terms cancel, and the signed sum of the other
inserted positions equals its shifted value.
It follows that \(R^{\rm wt}_a=-C^{\rm wt}_a\) at matched indices,
\(R^{\rm wt}_b=0\) at inserted indices, and
\(\operatorname{sgn} R^{\rm wt}_a=\operatorname{sgn}\xi_a\).

For a retained target entry, the formula gives
\( |2R^{\rm wt}_b-\xi_b|\le u+v-1<\Delta-1\).
In either group's decreasing parameter order write \(L_i^*\) for
the doubled entries and \(\mu_i\) for the relative compact weights.
For indices \(i<j\) of the same signature, let \(c\) count indices
of that signature in \((i,j]\). The vacuum weight is constant on
that compact block, so
\[
 \mu_i-\mu_j=(L_i^*-L_j^*)/2+(j-i)-2c\ge 2(j-i-c)\ge0 .
\]
Here the successive doubled entries have gaps at least two.
In each sign quadrant take the product over \(h\ge1\) of the
minors on matched indices of absolute weight at least \(h\),
omitting empty factors. On the standard coordinates, weights
at least \(h\) occupy an initial segment. On the dual coordinates,
weights at most \(-h\) occupy a final segment, which becomes initial
in dual order. Since on matched pairs \(R^{\rm wt}=-C^{\rm wt}\) and
\(s'_a=s_a\operatorname{sgn}\xi_a\), the indices matched to each such
segment lie in one compact block of the other group. A compact raising
operator kills the minor or duplicates an index in it.
The row and column weights are \(R^{\rm wt}_b,C^{\rm wt}_a\).
At the matrix with ones on the matched entries and zeros elsewhere
each minor is \(\pm1\), so the product is nonzero. Each source or
target index is used with one oscillator sign. Thus all mixed
derivatives of opposite signs with a common index vanish.

The joint harmonic space contains the indicated pair of compact types
with multiplicity one by \cite[Theorem~2.2(1) and Lemma~2.5]{PaulFirstOccurrence}.
Thus this nonzero highest polynomial spans its highest line. Both ranks
\(u=r+1\) and \(v=6r\) are even. The doubled shifted stem entries are
odd, as are the inserted entries \(t-2,t-4,\ldots,2-t\), so both
parameters lie in the genuine lattice \(\mathbb Z+(u+v+1)/2\).
The signature calculation and gap above give the discrete-series
correspondence of \cite[Theorem~2.7(b), p.~641]{PaulFirstOccurrence}.
Thus the positive Bargmann pairing and the two compact
averages used for the first-arrow constant apply to
this polynomial. They use the source matrix unit of
parameter \(-L/2\) and the contragredient scalar-form
coefficient. At the real place \(\tau\), write \(P_{i,\tau}\)
for the polynomial, \(e_{i,\tau}\) for the source matrix unit,
and \(\omega_\tau\) for the Weil action. The preceding positive
Bargmann integral and compact Schur averages give
\[
 Z_{i,\tau}=
 \int_{U(Z_0)(F_{0,\tau})}
   \langle P_{i,\tau},\omega_\tau(h,1)P_{i,\tau}\rangle
   \frac{\overline{e_{i,\tau}(h)}}{e_{i,\tau}(1)}\,dh>0 .
\]
In projection notation this is \(\|p_{i,\tau}P_{i,\tau}\|^2/e_{i,\tau}(1)\),
where \(p_{i,\tau}=\int\overline{e_{i,\tau}(h)}\,\omega_\tau(h,1)\,dh\)
and \(e_{i,\tau}(1)\) is the positive formal degree for the original
real Haar measure. The product over real places is
\[
 C_i=\operatorname{vol}([U(T)])\prod_{\tau\mid\infty} Z_{i,\tau}>0 .
\]
The passage to the tensor realization in
Section~\ref{sec:four-unitary-real-projections} dualizes the Hom source
before applying the common-entry rule of Paul's theorem.
The operator \(\mathcal V_i\) has the original vacuum weights away from the
prescribed support. The local lattice count of
Section~\ref{sec:marked-first-arrow-rallis}, with source rank
\(u\), gives \(O(q_w^u)\) neighbor choices of weight
\(q_w^{-v/2}\) at each place.

Here \(v/2-u\ge r\) for large \(x\).
At a split pair the number of moves is at most
\(2(1+\cdots+q_w^{u-1})\), and at a nonsplit place it is at most
\(\sqrt{q_w}(1+\cdots+q_w^{u-1})\). Each is at most \(4q_w^u\),
uniformly in \(u\). The total one-step weight therefore satisfies
\[
 \varepsilon\le4[E:\mathbb Q]\sum_{n\ge2}n^{-r}
             \ll [E:\mathbb Q]2^{-r}=o(1).
\]
The polynomial bound for \([E:\mathbb Q]\) in
Section~\ref{sec:datum} gives \(\varepsilon<1/3\) for sufficiently
large \(x\), uniformly in the allowed levels. Each right translate
is unitary on the ambient Petersson \(L^2\), so the absolute sum over
nonidentity words gives
\[
 \|\mathcal V_i-I\|_{2\to2}\le\frac{\varepsilon}{1-\varepsilon}<\tfrac12
\]
on the prescribed level. The series preserves that level by spherical
bi-invariance off the prescribed support. Its summands commute with the
real projector, hence so does their operator-norm limit.
By \eqref{eq:display-tp1}, the adjoint on the fixed-type subspace
at \(K_{Z_0}(\mathfrak b)\) satisfies
\[
 A_i^*A_i=\frac{C_i}{w_0|R_b|^{uv}}\mathcal V_i .
\]
The scalar is positive, and the preceding operator bound gives
\[
 \|A_i y\|^2
   \ge \frac{C_i}{2w_0|R_b|^{uv}}\|y\|^2
\]
on this subspace. Thus \(A_i\) is injective.

The first descent \(d_{i,0}\) is the adjoint return to the next
opposite base copy. To realize it with the original moduli, let
\(P_i\) be the upward joint harmonic and \(P_i^\flat\) the
opposite real edge harmonic. Conjugation of the Weil theta
kernel carries its residue datum \(\Psi\) to
\[
 \Psi_0(X)=\overline{\Psi(u_R^{-1}X)},\qquad
 u_R=u_0\bmod\mathfrak b .
\]
Here \(u_0\in\widehat{\mathcal O}_E^\times\) is the fixed unit
satisfying \(u_0\bar u_0=-1\) in the matrix reversal \eqref{eq:display-mr}.
The scalar \(c_i^{\rm op}\) in this conjugation is fixed
by the normalized theta functionals of the two opposite
bases.

The quotient Haar is carried to the same Haar.
Replace the conjugate polynomial by \(P_i^\flat\)
in the tested period using the fixed real edge convolution.
If its constant is \(k_i'\), the identity with the resulting
descent \(B_0\) is
\(\langle x,A_i y\rangle=\overline{c_i^{\rm op}}\,k_i'
\langle B_0x,j_{i+1}y\rangle\).
We use inner products linear in their second argument.
Thus the datum \(c_i^{\rm op}\overline{k_i'}\Psi_0\) gives
the exact adjoint. This identity uses the opposite-basis
conjugation and the real edge convolution of
Section~\ref{sec:marked-square-analytic-pairing}.
The conjugation and real replacement take place before
any residue expansion, so the initial adjoint datum is a
function on \(R_b\).

Let \(y^\flat\) be the return on the base. Then
\(\langle x,y^\flat\rangle=\langle A_i x,A_i y\rangle\).
Let \(k_b\in\mathrm{GL}_u(R_b)\) be the reduction of the integral
base element \(k\). At the selected split places lift
\(g_b=\operatorname{diag}(k_b,I_{v-u})\) to the integral target.
Then \(g_b^{-1}X_0k_b=X_0\) in the chosen Schr\"odinger coordinates,
and the determinant character contributes a scalar of modulus one.
Consequently the source and target translations intertwine
\(A_i\) up to that unit scalar. Pairing with all \(x\) shows
that \(A_i^*A_i\) commutes with the integral base action.
In particular \(y^\flat\) has the deep level when \(y\) does.
If \(y^\flat=0\), take \(x=y\) and use \eqref{eq:display-tp1} to conclude \(y=0\).

\subsection{The measured square and transposition}
\label{sec:deep-measured-square}

\sectionopening{In this subsection, we use the local return of Section~\ref{sec:deep-first-column-return} in the measured square.}

At the square \((i,n)\), the kept ranks are \(d_{n+1},d_{n+1}\)
and the other ranks are \(d_n,d_{n+2}\). The exchange defect
is \(t-1\). The stem sign with excess parameter \(a_*\) at
\((o,j)\) is
\((-1)^{o+a_*+\min(j,(a_*-o)_+)}\).
Passing from the kept pair to the other pair leaves both signs
unchanged if \(a_*\le i+n+1\). Otherwise the two changed signs
were opposite, so the negative count is preserved. The common
negative count \(q_-\) satisfies
\(q_-\le2(r+1)<d_{n+1}\). The quotient character is
\[
 \nu=\mu^0_{i+n+2}/\mu^0_{i+1}
        =(\upsilon^{-\Delta}\chi_M)^{n+1}.
\]
Its inverse has positive weight \((n+1)\Delta\), and
\(d_{n+1}\le(n+1)\Delta\).

At the chosen split places its
inverse tame part is the inflation of \(\chi_M^{-(n+1)}\),
nontrivial by \(1\le n+1\le m\), and elsewhere it is
unramified. The forward and opposite left-pair factors at
ranks \(p,q\), \(p+q=2d_{n+1}\), have product
\[
 \mu^0_{i+n+2}\upsilon^{\lceil p/2\rceil}
       (\mu^0_{i+1})^{-1}\upsilon^{\lfloor q/2\rfloor}
           =\nu\upsilon^{d_{n+1}} .
\]
The first adds the center \(-w_{i+n+2}\). For the left pair,
pass to the opposite form and add the center \(-w_{i+1}\).
Conjugating both parameters replaces the old determinant
exponent at rank \(h\) by \(h\) minus that exponent,
giving the floor on the left pair.

The positive-weight character in
Section~\ref{sec:marked-partner-period} is \(\nu^{-1}\).
Its prescribed splitting character on the kept space is
\((\nu^{-1})^{-1}\upsilon^{d_{n+1}}=\nu\upsilon^{d_{n+1}}\),
the product just computed. This verifies the remaining character
condition for the measured period on this square. Its ratio is \(k_{i,n}/c_{i,n}\),
where \(c_{i,n}\ne0\) is the measured partner scalar and
\(k_{i,n}\) comes from the real pairing of the four harmonic
vectors. The real replacement is the one in
Section~\ref{sec:marked-square-analytic-pairing}.

We prove that the point extends over \(R_b\).
In the decomposition \(b_1,a\) of the full matrix, invertible
blocks \(Y_1,X_2\) in the first \(b_1\) columns of a matrix and
the last \(a\) rows of its inverse give
\[
 X_2^{-1}X_1=-Y_2Y_1^{-1},\qquad
 \det k=\det Y_1/\det X_2 .
\]
These follow by block elimination over each finite local ring.
Jacobi's complementary minor identity makes invertibility
of either square block equivalent to invertibility of the other.
For fixed \(s\), use the determinant-character weighted
indicators of \(X_2^{-1}X_1=s=-Y_2Y_1^{-1}\) in the
partner period. Covariance cancels the character on the two
invertible blocks. Thus the respective point values are
multiplied by \(|\mathrm{GL}_a(R_b)|\) and
\(|\mathrm{GL}_{b_1}(R_b)|\). Here only invertible blocks
contribute. The argument uses finite groups over \(R_b\)
itself and requires no Gauss sum with a primitive
multiplicative character at the deep power.

For the actual square let \(a_*=d_n\), \(d_+=d_{n+1}\),
and \(b_*=t-1\). In the order \(d_+,a_*,b_*\) the two
arguments are
\[
 Y=\begin{pmatrix}-\eta\\ I_{a_*}\\0\end{pmatrix},
 \qquad
 X=\begin{pmatrix} I_{d_+}&\eta&s\\0&0&I_{b_*}\end{pmatrix}.
\]
They are obtained by choosing the common graph parameter
\(\binom{\eta}{0}\) in the order \(a_*,(d_+,b_*)\) and
applying the determinant-one shear with off-diagonal block \(s\).
The point extension therefore uses the new descent datum
\[
 p^+\begin{pmatrix}\eta&s\\0&I_{b_*}\end{pmatrix}
      =p(-\eta)|R_b|^{-d_+b_*},
 \qquad p^+=0\text{ on the other matrices}.
\]
Indeed the shear fixes the tested covector value by covariance.
Average over \(s\) and sum over \(\eta\) with
weight \(p(-\eta)\). Starting with the adjoint datum above
and applying this recursion gives descents \(d_{i,n}\)
with the measured identity
\[
 w_n a_{i+1,n}d_{i,n}
       =(k_{i,n}/c_{i,n})w_{n+2}d_{i,n+1}a_{i,n+1}.
 \tag{TP2}\label{eq:display-tp2}
\]
To test this equality, take a source form and a target test of
the indicated real types at a common finite level fixing the
residue data. The displayed \(Y\) has an identity minor of order
\(a_*\), and the first and last column blocks of \(X\) give a triangular minor of order \(d_++b_*\) and determinant one.
Thus the primitive coefficient estimate of
Section~\ref{sec:marked-square-analytic-pairing} applies to each
completed pair, before the finite sums in \(s,\eta\). The two
external forms are rapid, and the individual lifts with the
prescribed harmonics are rapid at their respective finite levels
by that same section. The operator identity follows by testing the
difference of the two resulting forms against itself.

\indent{\it Injectivity and transposition.}
Write \(X_{i,h}\) for \(h\) weighted forward arrows in sequence \(i\)
and \(X_{i,0}=I\). Iterating \eqref{eq:display-tp2} from left to right gives,
for some scalar \(\alpha_{i,h}\),
\[
 X_{i+1,h}d_{i,0}a_{i,0}
                  =\alpha_{i,h}d_{i,h}X_{i,h+1}.
 \tag{TP3}\label{eq:display-tp3}
\]
The initial scalar is \(\alpha_{i,0}=w_0^{-1}\).
Multiplying on the left by \(w_h a_{i+1,h}\) and applying
\eqref{eq:display-tp2} gives
\[
 \alpha_{i,h+1}=\alpha_{i,h}\frac{k_{i,h}}{c_{i,h}}
                    \frac{w_{h+2}}{w_{h+1}},
 \qquad X_{i,h+2}=w_{h+1}a_{i,h+1}X_{i,h+1}.
\]
If \(X_{i,h+1}j_i y=0\), \eqref{eq:display-tp3} makes
\(X_{i+1,h}j_{i+1}y^\flat=0\). Induction on the chain length,
starting with \(X_{i,0}=I\)
and using the injective return, proves \(y=0\).
In particular the full true forward chain \(F=X_{0,m}\)
is injective on the original deep base. This implication
uses \eqref{eq:display-tp2} in its displayed direction.
By Section~\ref{sec:marked-span-proof}, the projected space in
\(W\) has positive dimension. Thus
\eqref{eq:deep-separation-condition} gives \(W_Z\ne0\).
For any nonzero forward base vector \(y\), the return
\(y^\flat=A_i^*A_i y\) is nonzero. By injectivity at length
\(h\), the left side of \eqref{eq:display-tp3} is nonzero on
\(j_i y\). Therefore \(\alpha_{i,h}\ne0\) for every admissible
\(i,h\). Its recurrence then gives \(k_{i,h}\ne0\) for each square,
using \(\alpha_{i,h+1}\ne0\). The argument takes place at the trivial
twist and the level \(\mathfrak b\).

On every forward edge of \(F=X_{0,m}\), the residue point is the
specified inclusion matrix. Its left and right actions by the two
principal subgroups at \(\mathfrak b\) fix that matrix. Their determinant
character is one, since it factors through the prescribed tame
reduction. Thus each forward and transposed arrow has the
asserted principal level on the target. At an edge let \(f_1,f_2\)
be the partial lifts from the two ends in the bilinear pairing.
Their prescribed types are cuspidal by
Section~\ref{sec:marked-square-analytic-pairing}. At a definite end
the quotient is compact. For the fixed edge datum its kernel is
bounded by \(C\operatorname{ht}(g)^B\operatorname{ht}(h)^B\)
for some \(B\), by the two-variable theta estimate of that section.
The product with \(|f_1(g)f_2(h)|\) is therefore integrable on the
two quotients. Fubini moves this edge to the other side using the
same two quotient measures. At the next edge take the already
evaluated partial lifts as the two tests and apply the same
two-variable bound. Pullback from a coarser finite level leaves the
automorphic function and quotient integral unchanged. This iteration
puts the same weight \(w_j\) on its transposed edge and gives
\[
       \int_{[Z_m]} (F y)\phi
          =\int_{[Z_0]} y(\Theta^\flat\phi).
 \tag{TP4}\label{eq:display-tp4}
\]
The transpose here is bilinear and has the conjugate
endpoint types. If \(z\in\mathcal H_0\) is perpendicular
to its image, take \(y=\bar z\) and
\(\phi=\overline{F\bar z}\in\mathcal H_m\).
Then \eqref{eq:display-tp4} gives \(\|F\bar z\|^2=0\). Injectivity gives \(z=0\).
The fixed-type fixed-level spaces are finite dimensional,
so \(\Theta^\flat:\mathcal H_m\to\mathcal H_0\) is onto.

\subsection{The first-path array}
\label{sec:deep-first-path-array}

\sectionopening{In this subsection, we write the finite array for the transposed lift of Section~\ref{sec:deep-measured-square}.}

Use the point \(\omega\) from Section~\ref{sec:marked-deep-separation}.
Let \(R_0(\xi)=\mathcal O_E/(\pi_M Q(\xi))\) and
\(R_\xi=\mathcal O_E/(\mathfrak b Q(\xi))\).
The false first column uses its original twist and \(R_0(\xi)\).
The true first column uses the trivial twist and \(R_b\).
Both can be read on \(R_\xi\), since the reduction maps
\(R_\xi\to R_0(\xi),R_b\) are onto. Indeed, with
\(R_Q=\mathcal O_E/Q(\xi)\), the Chinese remainder theorem
gives \(R_\xi=R_b\times R_Q\). On the true half the map is
projection onto \(R_b\). On the false half it is reduction
of each old prime power to its residue field and the identity
at the new places, so its kernel is nilpotent. Pullback of an
edge datum from either quotient is the same adelic datum:
its value on an integral matrix depends on exactly the indicated
class, and it has zero value off the integral box at the chosen places. At the extra places on the true side the pulled-back
value on the integral box is constant, hence is the vacuum.

For a tuple \(S\) of matrices on the two first columns put
\[
 P_0(S)=\prod_{\substack{z\in\{0,1\}\\0\le i<m}}
           \mathbf 1_{\{\bar S_{z,i}=\iota_{z,i}\}},
\]
where \(\iota_{z,i}\) is the original inclusion matrix, and the
bar denotes reduction to the ring of that column.
Put \(D_{z,i}\) for the group order at the source of that
column edge over \(R_0(\xi)\) or \(R_b\), respectively.
Let \(T_{z,i}(S_{z,i})\) be the unweighted transposed
point lift at \(S_{z,i}\) over \(R_\xi\).

At an old place each integral residue box for the original
column is the disjoint union of its lifts at the deeper exponent.
At a new chosen place of the true column the original vector is
the vacuum. In its split Schr\"odinger coordinate this is the
indicator of the entire integral matrix lattice, the disjoint
union of \emph{all} residue boxes there. The conjugate coordinate
belongs to the same split model with its original self-dual
measure. Thus the pullback along either reduction gives the same
finite vector in the original model. By linearity, summing the
point indicators with their original values gives the weighted arrow as
\[
 D_{z,i}\sum_{S_{z,i}}\mathbf 1_{\{\bar S_{z,i}=\iota_{z,i}\}}
                         T_{z,i}(S_{z,i}).
\]
The matrices in a tuple are independent choices on disjoint
edges. Multiplying the two chains and using bilinearity of
the diagonal period at \(\omega\) gives
\[
 \left(\prod_{z,i<m}D_{z,i}\right)
                        \sum_S P_0(S)\mathcal P_\omega(S;y',\eta).
 \tag{TP5}\label{eq:display-tp5}
\]
Here \(\mathcal P_\omega(S;y',\eta)\) is the original diagonal
pairing of the two unweighted transposed point chains, at the
same real \(\omega\), of the false-terminal dual form \(y'\) and
true-terminal form \(\eta\in\mathcal H_m\).

Write \(G_{z,i}\) for the invertible group in the node's rank
over \(R_\xi\). Set \(L_{0,i}=\{1\}\) and
\(L_{1,i}=\ker(G_{1,i}\to\mathrm{GL}_{d_i}(R_b))\). Use the
finite set
\[
                  \mathcal C=\prod_{z,i<m}G_{z,i}/L_{z,i}.
\]
Choose representatives, with identity at the terminal nodes.
In the acting family \(g(c)\) replace the representative at
the true bottom node by the block sum of the false bottom
representative and \(1\). The unused coset coordinate is still
counted in \(\mathcal C\).

For every node \((z,a,b)\) let \(g_{z,a,b}\) be an invertible matrix
of its rank over \(R_\xi\). The family is admissible if
\[
 g_{0,0,m}=1,\qquad g_{1,0,m}\equiv1\pmod{\mathfrak b},\qquad
 g_{1,j,j}=\operatorname{diag}(g_{0,j,j},1)\quad(0\le j\le m)
\]
in the prescribed reindexing at each diagonal node. The other node
matrices are unrestricted. On the first path the only common
diagonal node is \(j=0\), which explains the replacement in
\(g(c)\). At the definite ends these conditions fix the false endpoint
modulo \(\mathfrak b Q(\xi)\) and the true endpoint modulo
\(\mathfrak b\). For an edge \(e:p\to q\), put
\((S^g)_e=g_q^{-1}S_e g_p\). Let \(\chi(g)\) be the product
of the edge characters at \(\det(g_q)/\det(g_p)\), using
the characters of their respective columns. The finite Weil
action on the original boxes gives
\[
 \mathcal P_\omega(S^g;y',\eta)
                  =\chi(g)\mathcal P_\omega(S;y',\eta)
\]
for all admissible families \(g\). Indeed choose integral
lifts of these invertible matrices at each node. Moving a lift
through the transposed edge contributes its determinant
character and replaces the edge matrix by \(g_q^{-1}S_eg_p\).

Every interior group translation cancels between the adjacent
edges. The bottom translations match by the block relation,
commute with \(\omega\), and disappear by right invariance
of the diagonal quotient measure. At the two ends the
translations fix the endpoint forms by their stated levels.
All these equalities hold at the given finite points.

Set
\[
                         P_1(S)=\sum_{c\in\mathcal C}
                         \chi(g(c))P_0(S^{g(c)}).
\]
This array has the required covariance
\(\chi(g)P_1(S^g)=P_1(S)\). To see this, use
\((S^g)^h=S^{gh}\) and send \(c\) to \(gc\). The representative
of \(gc\) differs from \(g\,g(c)\) by the node-level subgroups
on the right. These latter groups fix \(P_0\) and have trivial
character, since they reduce to identity in the ring of their
column. At the true bottom node the discrepancy is identity
by the block relation and \(L_{0,0}=\{1\}\).

This also verifies
independence of representative choices. By reindexing \(S\)
in each summand and using the preceding period covariance,
\[
       \sum_S P_1(S)\mathcal P_\omega(S;y',\eta)
             =|\mathcal C|\sum_S P_0(S)\mathcal P_\omega(S;y',\eta).
\]
Each entry contains at most \(|\mathcal C|\le E_\xi^{2m}\)
roots of unity, where \(E_\xi\) bounds the number of residue
matrices at the allowed grid ranks. If an entry is nonzero,
some translate has nonzero indicator.

Fix such a translate \(g\). On a false edge \(e:p\to q\), let
\(\pi\) project onto the first \(\operatorname{rk}(p)\) coordinates
and put \(C=\pi g_q^{-1}S_e g_p\). Then \(C\equiv I\) over
\(R_0(\xi)\). The reduction \(R_\xi\to R_0(\xi)\) has nilpotent
kernel, so \(C\) is invertible over \(R_\xi\). In the original
coordinates the left inverse is
\[
          (g_p C^{-1}\pi g_q^{-1})S_e=I.
\]

Let \(n_i=r+it\) be the rank at the false first column's
node \(i\). Its factor in \(|\mathcal C|\) is
\(|\mathrm{GL}_{n_i}(R_\xi)|\), while the true factor is the
order over \(R_b\), by onto reduction. The order over
\(R_0(\xi)\) divided by the order over \(R_\xi\) is
\((|R_0(\xi)|/|R_\xi|)^{n_i^2}\).
Indeed a matrix is invertible if and only if its image is,
and each matrix has a fiber of size
\((|R_\xi|/|R_0(\xi)|)^{n_i^2}\).
Thus the ratio of the product in \eqref{eq:display-tp5} to \(|\mathcal C|\)
is the fixed positive number
\[
             A_{\rm col}=
              \prod_{i<m}\left(\frac{\mathrm{Nm}(\pi_M)}
                                  {\mathrm{Nm}(\mathfrak b)}\right)^{n_i^2}.
 \tag{TP6}\label{eq:display-tp6}
\]

For \(\eta\in\mathcal H_m\), take \(r_0\eta\) to be the Riesz
representative on the base \(W\) of the restriction
\(\bigl(\Theta^\flat\eta\bigr)((h\oplus1)(\omega,1))\).
In the original Petersson metric, linear in the second variable,
this means that for \(v\in W\),
\[
 \langle v,r_0\eta\rangle_W
    =\int_{[U(V)]}\overline{v(h)}\,
       (\Theta^\flat\eta)\bigl(\iota(h)(\omega,1)\bigr)\,dh .
\]
The integrand is absolutely integrable since \(v\) is rapid and
\(\Theta^\flat\eta\in W_Z\) is bounded. This defines a linear
map of \(\eta\), with no choice involving \(\xi\).
Take a test \(v\) of the real type of \(W\) and level
\(K_V(\pi_MQ(\xi))\), and pair with its conjugate in the
bilinear integral. Its average at \(K_V(\pi_M)\) belongs
to \(W\), since the finite average preserves the real type.
For \(k\in K_V(\pi_M)\) write \(k=k_bk'\), where \(k_b\)
is supported at the original places and their conjugate
places, and \(k'\) is the identity there. Because the twist
support is disjoint from these places, \(k_b\) fixes every
test of level \(K_V(\pi_MQ(\xi))\).

The component \(k'\)
lies in \(K_V(\mathfrak b)\), so its block embedding fixes
the true function of deep level, and it commutes with \(\omega\).
The restriction paired with \(v\) therefore equals the
restriction paired with its average. Apply the Riesz
identity in \(W\) to that average, and remove the average
from the pairing with \(r_0\eta\), which has the base level.
If \(v\in W\) is perpendicular to \(\operatorname{im}r_0\), the
displayed integral vanishes for every \(\Theta^\flat\eta\).
By \eqref{eq:display-tp4} these tests fill \(W_Z\), so
\eqref{eq:deep-separation-condition} gives \(v=0\).
Finite dimensionality gives \(\operatorname{im}r_0=W\).
To check the real type of the bilinear transpose, take one false
edge with joint harmonic \(P\) and the oscillator line projectors
\(p_S,p_T\) of Section~\ref{sec:four-unitary-real-projections}.
The full source isotypic projection \(Q_S\) is central in both
commuting von Neumann algebras by \cite[Theorem~6.1, p.~550]{HoweDuality}.
Its range is \(V_S\widehat\otimes V_T\), with \(V_T\) irreducible,
as proved there. Hence it is the identity on every copy of \(V_T\):
a target intertwiner between two copies commutes with \(Q_S\).
The matrix-coefficient projector \(p_T\) has range in this
target isotypic component. Since it commutes with the source action,
\(p_TP\) still has the source compact highest weight of \(P\).
That highest line in the lowest compact type of \(V_S\) occurs once.
Consequently
\[
                  p_Sp_TP=p_TP=p_Tp_SP .
\]
Moving the terminal test projector onto the kernel in the bilinear
pairing therefore reproduces the conjugate source line. The
weighted coefficient and theta bounds of
Section~\ref{sec:four-unitary-constant-terms} justify this
operation on each edge as in the forward direction. Iteration
gives the conjugate of a test of the required real type.
The source principal subgroup fixes the residue inclusion and
has trivial determinant character on each edge. Thus the transpose
has level \(K_V(\pi_MQ(\xi))\).
Applying \(r_0\)'s identity and \eqref{eq:display-tp5}--\eqref{eq:display-tp6} yields
\[
 \int_{[Y_m]} y'(h)(T_\xi r_0\eta)(h)\,dh
          =A_{\rm col}\sum_S P_1(S)\mathcal P_\omega(S;y',\eta).
 \tag{TP7}\label{eq:terminal-transpose-restriction}
\]

\subsection{Transport across the consecutive squares}
\label{sec:deep-transport-arrays}

\sectionopening{In this subsection, we transport the arrays of Section~\ref{sec:deep-first-path-array} through the squares.}

Parameterize each half path from its definite endpoint by right
and down steps toward a common diagonal node \((j,j)\),
recording the true half with its coordinates transposed.
At a square fix all residue matrices on unchanged edges.

In the ordered direct sum of the two kept node spaces the
old two edge blocks form a column \(Y\) and the new blocks
a block row \(X\). The number of columns of \(Y\) is the rank of
the removed node, and the number of rows of \(X\) is the rank of the
added node. At an off-diagonal square there
is no line variable. At a diagonal square write \(v\) for
the removed column of the line and \(v'\) for the new covector.
In the prescribed covector and column coordinates define
\[
 P^{\rm F}(Y,v')=\sum_v \psi(\operatorname{tr}(v'v))P_{\rm old}(Y,v).
 \tag{TP8}\label{eq:display-tp8}
\]
The sum for an off-diagonal turn has the one zero line
coordinate and is the old value. The character \(\psi\)
is the one in the line reversal below.

Here is the matched finite array for the partner
comparison of Section~\ref{sec:marked-partner-period}.
For a completion \(k\) in the general linear group of the
kept direct sum, write \(Y\) for its first columns and
\(X\) for the complementary rows of \(k^{-1}\). Define
\[
 P_{\rm new}(X,v')=\chi(\det k)P^{\rm F}(Y,v')
 \tag{TP9}\label{eq:display-tp9}
\]
and set it equal to zero for matrices \(X\) admitting no completion.
This is independent of the completion: another completion
with the same last rows of its inverse has the first columns \(Yg\),
where \(g\) is invertible on \(\ker X\), and its determinant
is multiplied by \(\det g\).

Indeed the change of basis of
the two completions has block form
\(\begin{psmallmatrix}g&*\\0&I\end{psmallmatrix}\).
The old path covariance at the removed node gives
\(P^{\rm F}(Yg,v')=\chi(\det g)^{-1}P^{\rm F}(Y,v')\);
that node acts on \(Y\) on the right and fixes the line
coordinate and the retained edges. Call an edge vertical when it changes the second index
on the false half or the first index on the true half, in
the coordinates of \(n_z(a,b)\) below. On the first path all
edges have the form \((z,0,i)\to(z,0,i+1)\), so precisely the
false edges are vertical. The left inverse just constructed over
\(R_\xi\) starts the induction, also for \(\xi=1\).
At every turn at least one old block is vertical.
If \eqref{eq:display-tp8} is nonzero, some old summand is
nonzero, hence that block has a left inverse. If the old
edge has an extra line column, restrict its left inverse
to the retained columns and project to their coordinates.
Over each local factor of \(R_\xi\), complete the residue
columns of \(Y\) to a basis and lift the extra columns.
The resulting determinant is a unit in that factor. The
Chinese remainder theorem assembles these completions into
one \(k\), allowing a different pivot minor at each place.

The node covariance of \eqref{eq:display-tp9} follows by the same change
of completions. Explicitly, the retained-node block
element \(G\) and the new-node element \(h\) change
\(X\) to \(h^{-1}XG\). One may replace \(k\) by
\(G^{-1}k\,\operatorname{diag}(1,h)\), so that \(Y\)
changes to \(G^{-1}Y\) and the determinant acquires
\(\det h/\det G\). At a diagonal turn the action on the
old column and new covector of the line is inverse at their
common subdiagonal node.

Changing \(v\) in \eqref{eq:display-tp8} therefore
preserves \(\operatorname{tr}(v'v)\).
In detail, let \(\chi_1,\chi_2\) be the characters of
the two old edges, with target determinants \(d_1,d_2\).
The removed node determinant is the same on its line
extension. Thus \(\chi=(\chi_1\chi_2)^{-1}\).
If \(C_w\) is the product from the retained edges, the Fourier
value after the move with the removed element equal to the
identity is multiplied by \([C_w\chi_1(d_1)\chi_2(d_2)]^{-1}\).

Write \(d_h=\det h\). The completion contributes the factor
\(\chi(d_h/(d_1d_2))\). Opposite edges have equal
characters, so the new edge-character product is
\(C_w\chi_2(d_h/d_1)\chi_1(d_h/d_2)\), which cancels
the product of these two factors. Thus the new array
has the required covariance. Each new vertical block is opposite
a split injective old block. Over each residue field, if the
projection of \(Y\) into one retained summand is injective,
then its image intersects the other retained summand
trivially.

Since \(\ker X=\operatorname{im}Y\), \(X\) on
this other summand is injective. Lifting the resulting
minor over the finite local factors makes this a split
injection over \(R_\xi\). If the new edge also contains
the covector \(v'\), extend the left inverse by zero on
that coordinate. The retained edges keep their support conditions. Each new entry has at most as many
root-of-unity terms as the old bound times the number of
choices of \(v\), which is at most \(V_\xi\).

For half \(z=0,1\), the rank at \((a,b)\) is
\(n_z(a,b)=r+z+(t+z)a+(t-z)b\).
Put \(d=n_z(a,b)\). The two retained ranks are
\(d+t-z,d+t+z\); the two integrated ranks are \(d,d+2t\).
Thus both sums are \(N=2(d+t)\).
In the stem sign formula of Section~\ref{sec:four-unitary-real-projections}
the offset is zero. The positive excesses run from \(1\) to
\(r+z-k=m-1+z\), each once. In column \(b\) a stem entry
of excess \(h\) is negative exactly when \(h>b\) and \(h-b\)
is odd. All added entries are positive. With \(s=m-b+z\),
each diagonal therefore has
\[
 q=\lfloor s/2\rfloor+\lfloor(s-1)/2\rfloor=m-b-1+z
\]
negative signs, also for \(b=m-1\). Since \(q\le m<d+t\),
the lower bound \(2q+2\le N\) holds at every place. The product character before inversion is
\(\mu^0_{a+1}\mu^0_{b+1}\xi^h\),
where \(h=b+1\) on the false half and \(h=a+1\) on the
true half. Its inverse has weight at least \(\Delta\ge N/2\).
The old tame character is \(\chi_M^{a+b+2}\) with
\(2\le a+b+2\le2m\), and the new one is \(\xi^h\)
with \(1\le h\le m\). Each is nontrivial at its chosen places, including after inflation to the deep old powers.
At ranks \(p,q'\) the two edge characters have splitting
exponents \(\lfloor(p+1-z)/2\rfloor\) and \(\lfloor(q'+z)/2\rfloor\),
respectively. Their sum is \(N/2\), since \(p+q'=N\) is even.
The twist belongs to the vertical edge in each orientation. The true orientation of the square \((i,j)\)
uses \((a,b)=(j,i)\), so both use \(h=j+1\). The same
square calculation applies at \((j,j)\) before the line flip.

Fix a square \(s=(i,j)\) and the residue matrices on the unchanged edges.
In either half's recorded coordinates the retained nodes are
\(W_1=(a,b+1)\) and \(W_2=(a+1,b)\). The old corner passes through
\(V'=(a,b)\), and the new corner through \(V=(a+1,b+1)\).
Let \(L_B\) be the lift with point matrix \(B\), and
\({}^{\mathsf t}L_B\) its transpose for the bilinear period.
Equation~\eqref{eq:display-tp9} gives matched arrays supported
on split surjections and split injections at the chosen places.
For the incoming partial lift \(f\) on \(W_1\) and a linear
functional \(\lambda\) on the output at \(W_2\), the matched
operator identity of Section~\ref{sec:four-unitary-real-projections},
obtained from \eqref{eq:marked-square-paired-period}, gives
\[
\begin{aligned}
 &\sum_Y P^{\rm F}(Y,v')\,
          \lambda\bigl(L_{Y_2}{}^{\mathsf t}L_{Y_1}f\bigr)\\
 &\qquad=\frac{k_s}{c_{{\rm sq},s}(\xi)}
        \sum_X P_{\rm new}(X,v')\,
          \lambda\bigl({}^{\mathsf t}L_{X_2}L_{X_1}f\bigr).
\end{aligned}
\]
Indeed, the two weighted outputs at a common finite level have the
prescribed cuspidal real type. Testing their difference against itself
in the Petersson pairing proves equality as forms. The measured
partner scalar satisfies \(c_{{\rm sq},s}(\xi)\ne0\). The real coefficient \(k_s\)
is fixed before the twist and the forms; it need not be inverted.

Off the diagonal, \(\lambda\) applies the remaining arrows from
\(W_2\), then pairs with the unchanged half at the fixed real point.
The unchanged true half is the second argument when \(i<j\); the
unchanged false half is the first when \(i>j\). The latter form
already has the conjugate type. The diagonal integral is absolutely
convergent since the smaller-group form is in \(L^1\) and the
translated larger-group form is bounded. All four point lifts have
their original quotient measures. Thus \(c_s=c_{{\rm sq},s}\) off
the diagonal. On the diagonal the false block uses the same square
identity with \(v'\) fixed. The square rank defect is \(t\); the square law is
evaluated at the same determinant shifts and common
Haar as in the comparison. Put
\[
 S_\xi=S_M\cup\operatorname{supp}Q(\xi),\qquad
 \nu_{ij,\xi}=(\mu^0_{i+1}\mu^0_{j+1}\xi^{j+1})^{-1},\qquad
 T_{ij}(\xi)=\prod_{k=1}^t
       \frac{|L^{S_\xi}(k,\nu_{ij,\xi})|}{|L^{S_M}(k,\nu_{ij,1})|}.
\]
The incomplete base-change series omit each indicated place and
its conjugate, as in Section~\ref{sec:marked-partner-period}.
Taking the ratio \(M_{\nu_{ij,\xi}}(t)/M_{\nu_{ij,1}}(t)\) gives
\[
 \frac{|c_{{\rm sq},s}(\xi)|^2}{|c_{{\rm sq},s}(1)|^2}
  =\mathrm{Nm}(Q(\xi))^{t^2}\,T_{ij}(\xi)^2.
 \tag{TP10}\label{eq:display-tp10}
\]
The factors from the old residue cardinalities cancel.
The appearance of \(Q(\xi)\) only to \(t^2\) uses that
the conductor factor is a product over the support,
and the new support is squarefree. The locally constant
matched comparison in Section~\ref{sec:radical-local-line-terms}
applies to the pullbacks of the arrays at the old exponent \(b_0\).
Thus \eqref{eq:display-tp10} holds for the arrays over \(R_\xi\).

At the diagonal \(j\), the two intermediate spaces are
identified in their chosen coordinates. Give this space
rank \(n_j=r+(2j+1)t\). The old and new lines have signs
\(-\epsilon\) and \(\epsilon\), \(\epsilon=(-1)^{j+1}\).
Choose the finite idele \(t_\xi\) generating the grid modulus with
component one outside the selected halves, including at their conjugates.
With \(\alpha_f=\alpha_{\rm W}\) the fixed finite Weil scale, the old and new finite
line points are \(u_\xi^j,u_\xi^{j+1}\) for
\(u_\xi=\overline{t_\xi\alpha_f}/(t_\xi\alpha_f)\).
In the original sign coordinates the polynomial reversal is
\(X_{v,i,1}\mapsto -\epsilon s_i X_{v,1,i}\), where \(s_i\)
are the intermediate signs.

In fact the Hermitian adjoint
transposes and conjugates, while reversing the line changes
its sign and thus the conjugation in the holomorphic Fock
coordinate. The Gaussian is unchanged. Let \(u_0\) be the
fixed finite norm-minus-one unit and \(c_0=u_0^\epsilon\).
On the self-dual coordinates use
\[
 \psi(\operatorname{red} z)=
   \psi_f\!\left(j_{\rm W}\,\alpha_f\overline{\alpha_f}
                         (-\bar c_0)z/t_\xi\right).
\]
The integral box shifts by \(v/t_\xi\) and the line dilation
changes it to \((X-v)/t_\xi\) integral.

The squared modulus of the split Weil multiplier of this dilation
is \(|R_\xi|^{n_j}\). The conjugate-place
Heisenberg displacement gives the multiplier
\(\psi(\operatorname{tr}(v'v))\) on the reversed integral box.
Before changing the determinant character on the reversed line,
write \(C_{\rm rev}(\xi)\) for the scalar of this unnormalized
finite Fourier sum. Its squared modulus is \(|R_\xi|^{-n_j}\).
The rational Hermitian transpose identifies the rational summation spaces.
In the finite adelic bases write the displacement as \(z=\alpha_f X\).
The adjoint change in these bases is
\[
  -c_0\frac{\alpha_f}{\overline{\alpha_f}}z^*
        =\alpha_f(-c_0 X^*).
\]
Thus \(X\mapsto-c_0X^*\) preserves the integral matrix box in normalized
coordinates, including at the primes of the relative different. The
factor \(\overline{\alpha_f}/\alpha_f\) links the two line frames in the
group argument. The dilation \(\bar t_\xi/t_\xi\) acts as
\(t_\xi^{-1}\) on the selected halves and as its inverse conjugate
on their partners. Their product is \(u_\xi\). The transport of the normalized theta functional has unit
modulus. For the determinant phases put \(u=\upsilon|_{U(1)}\)
and set
\[
 H_j(a)=\mu^0_{j+1}u^{\lceil a/2\rceil},\qquad
 V_j(a)=\mu^0_{j+1}\xi^{j+1}u^{\lfloor a/2\rfloor}.
\]
The ranks \(h_j=r+2jt\) and \(h_{j+1}\) are odd, so
\(H_j(h_j+1)=H_j(h_j)\) and
\(V_j(h_{j+1}+1)=u V_j(h_{j+1})\).
With \(k=u_\xi^j\) the determinant factors on the two
products are
\[
 \frac{H_j(h_j)(\det g_{\rm mid})}
      {H_j(n_j)(k\det g_{\rm old})},
 \qquad
 \frac{V_j(n_j)(ku_\xi\det g_{\rm new})}
      {(u V_j(h_{j+1}))(\det g_{\rm mid})}.
\]

They are exactly the true-edge factors by the parity
identities. Reversing the column replaces its character
\(H_j(n_j)\) by \(H_j(n_j)^{-1}u^{n_j}\).
The other character is already \(u\). Thus the
ratio to the required line kernel on the new side is
\[
 \frac{H_j(n_j)^{-1}u^{n_j}}{V_j(n_j)}(ku_\xi)
   = \bigl((\mu^0_{j+1})^{-2}\xi^{-(j+1)}\bigr)(ku_\xi).
\]
Denote this value by \(m_{j,\xi}\) and put
\(c_{\rm line}(\xi)=C_{\rm rev}(\xi)m_{j,\xi}\).
Write \(K^-_{v,p}(h)\) for the old line kernel from \(k\) to
the intermediate point \(h\), and \(K^+_{v',\tau p}(h)\)
for the reversed kernel from \(h\) to \(ku_\xi\), with the
prescribed characters on both edges. Then
\[
 K^-_{v,p}(h)=c_{\rm line}(\xi)
   \sum_{v'\in R_\xi^{1\times n_j}}
       \psi(\operatorname{tr}(v'v))\,K^+_{v',\tau p}(h).
\]
The factor \(m_{j,\xi}\) has complex absolute value one and is
independent of \(h,v,v'\) and \(p\). It is included in
\(c_{\rm line}(\xi)\). The coefficients of the finite sum in
\eqref{eq:display-tp8} are the values of \(\psi\).
The permutation of the block and line coordinates preserves
the quotient Haar measure. In the split Fourier coordinates the self-dual additive
measure gives the integral box volume one. Its pushforward to
\(R_\xi^{n_j}\) gives each residue class mass \(|R_\xi|^{-n_j}\).
In the displayed finite pairing the annihilator of the integral box is
\(t_\xi\widehat{\mathcal O}_E\) by self-duality and the
unit \(-\bar c_0\). In particular the finite pairing is perfect:
for columns \(v,w\in R_\xi^{n_j}\),
\[
 \sum_{v'\in(R_\xi^{n_j})^\vee}
       \psi\bigl(v'(v-w)\bigr)=|R_\xi|^{n_j}\mathbf 1_{v=w}.
\]
Indeed a nonzero coordinate of \(v-w\) generates a nonzero
ideal on which the primitive character is nontrivial, so the
corresponding sum vanishes by translation. Let \(F_{v',v}=\psi(v'v)\)
in a cyclotomic field containing the additive and residue
determinant character values. The identity reads
\(c(F)^{\mathsf T}F=|R_\xi|^{n_j}I\), where \(c\) is complex
conjugation. At a coefficient embedding \(\iota\) the two factors
are read via \(\iota\circ c\) and \(\iota\), respectively.
For a residue function \(D\), integration of \(D(v)\psi(v'v)\) over the integral box gives
\(|R_\xi|^{-n_j}\sum_v D(v)\psi(v'v)\). The sum in
\eqref{eq:display-tp8} uses counting measure. By the displayed
orthogonality its unitary normalization for the same measures on
both residue spaces is \(|R_\xi|^{-n_j/2}\sum_v D(v)\psi(v'v)\).
The unitary reversal supplies the reciprocal squared factor
\(|c_{\rm line}(\xi)|^2=|R_\xi|^{-n_j}\) in every residue characteristic.
Thus
\[
 \frac{|c_{\rm line}(1)|^2}{|c_{\rm line}(\xi)|^2}
   =\left(\frac{|R_\xi|}{|R_1|}\right)^{n_j}
   =\mathrm{Nm}(Q(\xi))^{n_j}.
\]
The common old factor cancels before taking coefficient valuations,
so the ratio is a unit at each coefficient place of residue
characteristic prime to \(q_\xi\).

For the measured diagonal exchange apply the two seesaws on
the intermediate groups. On the old side the true-edge datum
is the product of its original false block indicator and its
line-column indicator. On the new side the factorization
is into the false block and the line-covector indicator.
The factorizations have the actual line points just computed.

Both seesaws use the original quotient measures. Write
\(x=(\xi,y',\eta)\) for the entire triple, \(z\) for the
fixed real point of the intermediate group, and \(p\) for
a complex polynomial in the archimedean line variables. Denote by
\(\mathcal C_{z,p}(x)\) and \(\mathcal B_{z,p}(x)\) the old
column and new covector periods with this line weight.
Let \(\widetilde S,S\) be the true and false intermediate spaces of
rank \(n_j\), and \(s_j:U(\widetilde S)\to U(S)\) their fixed
coordinate reindexing. For each residue tuple let \(f_+\) be the
form from the true path and \(f_-\) the conjugate-type form from
the false path, including its two corner lifts. The term paired
against the corresponding old or new array is
\[
 \int_{[U(\widetilde S)]} K(g)\,f_+(g)\,
                  f_-\bigl(s_j(g)(z,1)\bigr)\,dg .
\]
Here \(K(g)=K^-_{v,p}(s_j(g))\) in the column period and
\(K(g)=K^+_{v',p}(s_j(g))\) in the covector period. The line point
is fixed at \(u_\xi^j\) or \(u_\xi^{j+1}\), respectively. In both
integrals the Haar measure is the original quotient measure, and
the real translation \((z,1)\) occurs only in the false form.
Sum over every residue tuple on the corresponding path, including
its line coordinate.

Suppose that the new path period at every real point is
a scalar multiple of a function \(\Xi(x)\), with the same
scalar for all \(x\) at that point. The new seesaw implies,
for each \(z,p\), that
\(\mathcal B_{z,p}(x)=\beta(z,p)\Xi(x)\) for all \(x\).
Here are the details of the passage to real tests. Fix
the line polynomial \(p\) first. For a polynomial \(A\)
on the false block, the product \(A p\), in the reindexed
orthogonal blocks, gives the whole-edge kernel whose
product factorization is the false block kernel of \(A\)
times the line kernel of \(p\). The residue point datum
on the whole edge is the product of the two point
indicators on those blocks. Substituting this identity
gives a new whole-edge lift with real datum \(A p\).
Let \(G_+\) be the real true group at \((j+1,j+1)\), and
write \(\mathcal P^{\rm new}_k(x)\) for the new path period
with the true form translated on the right by \(k\in G_+\).
The partial source form on \(\widetilde S\) is fixed by the
contragredient source-line projector. Apply
Section~\ref{sec:four-unitary-real-projections} to this source
and the whole true edge. If \(P_+\) is its distinguished harmonic,
there is \(\alpha_{A p}\in C_c^\infty(G_+)\) such that
\(p_{\widetilde S}(A p)=\omega(\alpha_{A p})p_{\widetilde S}P_+\).
Consequently the factored double integral, summed over the new
array, is
\[
                 \int_{G_+}\alpha_{A p}(k)
                     \mathcal P^{\rm new}_k(x)\,dk .
\]
Indeed, at the block point the true lift is evaluated at
\(\mathrm{diag}(h,u_\xi^{j+1})(k,1)\), while the false form
is still evaluated at \(h\). Real convolution leaves the finite
block and line points fixed. The compactly supported function
depends only on the real datum \(A p\). If \(\Xi(x_0)\ne0\), the proportionality
coefficient at each averaging point is the period there
at \(x_0\) divided by \(\Xi(x_0)\). It is integrable against
the smoothing function by the same fixed-triple bound.
If \(\Xi=0\), the average vanishes. In either case the
integrated expression belongs to \(\mathbb C\Xi\) for every
\(A\).

To pass to the actual real translate of the joint
harmonic, regard each summed double integral as a
continuous linear functional of its real Fock vector.
Fix \(x,z,p\) and one residue summand. Write \(f_1(g),f_2(h)\)
for the forms supplied by the two remaining paths. Successive arrows
with the distinguished polynomial have the prescribed cuspidal
target type by Section~\ref{sec:marked-square-analytic-pairing}; at
a definite endpoint the quotient at fixed level is compact.
The two-variable Schwartz bound for the kernel in that section, followed by
the coherent-state estimate of Section~\ref{sec:four-unitary-constant-terms},
majorizes the absolute integrand by
\[
 C p_N(v)\,|f_1(g)f_2(h)|
        \operatorname{ht}(g)^B\operatorname{ht}(h)^B
\]
for some \(C,B,N\). Indeed adjoining the fixed line gives an orthogonal
direct sum of lattice norms, and a fixed real translate changes the
minimum by bounded factors. The fixed line kernel has polynomial
growth by the same kernel estimate at its fixed polynomial, absorbed
in \(B\). The bound is uniform on the compact support of the smoothing
function. The two forms are rapidly decreasing on their quotients,
so the height-weighted product is integrable. There are only
finitely many residue summands for fixed \(x\).
If \(\Xi(x_0)\ne0\), a family of such integrals \(F_v(x)\)
is in \(\mathbb C\Xi\) precisely when
\[
       \Xi(x_0)F_v(x)-\Xi(x)F_v(x_0)=0\qquad\text{for all }x.
\]
Each of these equations is continuous in the real
vector \(v\). If \(\Xi=0\), use the equations \(F_v(x)=0\).
On the false block take \(v\) to be its distinguished harmonic
translated by \(z\) on the intermediate group. Covariance gives
the false form at \(s_j(g)(z,1)\), leaving the line weight \(p\)
at \(s_j(g)\). This translated polynomial times a Gaussian has
finite moments of all orders: the moved Gaussian has a positive
definite real part in its decaying exponent. The Taylor-kernel
approximation of Section~\ref{sec:four-unitary-real-projections}
therefore converges to \(v\) in every \(p_N\). Tensoring with
the fixed line polynomial times its Gaussian preserves this
convergence, since the Gaussian factors and
\((1+|(w_1,w_2)|)^N\le(1+|w_1|)^N(1+|w_2|)^N\).
Taking the limit in the continuous equations gives
\(\mathcal B_{z,p}\in\mathbb C\Xi\). In particular
\(\beta(z,p)\) is independent of every component of \(x\).

The line Fourier law gives \eqref{eq:display-tp8}. For each fixed new
finite covector and each tuple on the unchanged edges,
apply \eqref{eq:display-tp9} and the measured partner identity to the
false block. The archimedean coefficient \(k_s\) is the fixed
one of that square and the line reversal \(p\mapsto\tau p\)
depends on its signed real coordinates only. The measured
square denominator \(c_{\rm sq}(\xi)=c_{{\rm sq},(j,j)}(\xi)\)
and the Fourier scalar \(c_{\rm line}(\xi)\) are independent of the finite
covector and the unchanged edges. Summing the resulting
identities gives
\[
 \mathcal C_{z,p}(x)
   =\frac{k_s c_{\rm line}(\xi)}{c_{\rm sq}(\xi)}
                         \mathcal B_{z,\tau p}(x)
   =k_s\beta(z,\tau p)\frac{\Xi(x)}{c_s(\xi)},\qquad
 c_s(\xi)=\frac{c_{\rm sq}(\xi)}{c_{\rm line}(\xi)} .
\]
Since \(|c_{\rm line}(\xi)|^2=|R_\xi|^{-n_j}>0\) and
\(|R_\xi|/|R_1|=\mathrm{Nm}(Q(\xi))\),
\eqref{eq:display-tp10} gives the full diagonal ratio
\[
 \frac{|c_{(j,j)}(\xi)|^2}{|c_{(j,j)}(1)|^2}
      =\mathrm{Nm}(Q(\xi))^{t^2+n_j}\,T_{jj}(\xi)^2 .
\]

For the old seesaw, products \(A p\) in the false block
and old line variables span all the whole-edge polynomials.
Each product has the factorization into the two block point
indicators and the two real polynomials on those blocks.
This time smooth \(A\) on the false edge with the false
partial source at the smaller diagonal node. Its target is the
intermediate group \(G_-=U(S)_\infty\). The same projected
edge identity and Fubini give the tested old integral as
\[
                \int_{G_-}\alpha_A(k)\mathcal C_{k,p}(x)\,dk ,
                 \qquad \alpha_A\in C_c^\infty(G_-).
\]
Here the translated false form occurs at \(s_j(g)(k,1)\), and the
old line kernel is still evaluated at \(s_j(g)\) with its finite
point \(u_\xi^j\). The averaging preserves the common
proportionality to \(\Xi(x)/c_s(\xi)\). The same continuous
equations, now on the old side and with this function of
\(x\), extend the assertion to the actual translated
whole-edge harmonic. This proves that at every fixed
real point the old path period is proportional to
\(\Xi(x)/c_s(\xi)\) by one scalar for all triples.
The off-diagonal turn has this property directly by its
tested square identity.

\subsection{Composition and terminal measure}
\label{sec:deep-terminal-measure}

\sectionopening{In this subsection, we compose the transports of Section~\ref{sec:deep-transport-arrays} and retain the original terminal measure.}

For a pair of paths ending on the diagonal at \(j\), take
the cells in the \(j\)-by-\(j\) diagonal square, together
with the cells to the left of each down step in both halves,
reflecting the true half. A down step at horizontal
coordinate \(a\) across sequence \(b\) contributes the cells
\((h,b)\), \(0\le h<a\), in its half's recorded coordinates.
At an off-diagonal turn the down step moves one column to
the right, so the union contains one additional cell.

At a diagonal turn
the two strips of length \(j\) next to the \(j\)-by-\(j\)
square are already included, and the enlarged diagonal
square adds just \((j,j)\). Initially the union is empty.
A half with a down step before a right step contains an
adjacent down--right pair and therefore admits the turn. If neither half has one, all
right steps precede the down steps. If \(j<m\), both last
steps are down steps and give a diagonal turn. If \(j=m\),
the last path has been reached. Thus an arbitrary successive
choice of available corners crosses each square exactly
once, ending at the nodes \((z,j,m)\) on each half.

Fix one such order of the \(m^2\) turns for all twists.
Start from \(P_1\) at each twist and choose the arrays by
\eqref{eq:display-tp8}--\eqref{eq:display-tp9} in that order,
before choosing either terminal form. An off-diagonal
turn only multiplies a selected entry by a character value.
The diagonal turn at \(j\) sums over \(|R_\xi|^{n_j}\) columns.
Thus the number of roots of unity in each terminal entry is at most
\[
       E_\xi^{2m}\prod_{j=0}^{m-1}|R_\xi|^{n_j}
                     \le E_\xi^{2m} V_\xi^{m^2}.
\]
The backward implications compose by multiplying the
\(c_s(\xi)\). If \(c_\xi=\prod_s c_s(\xi)\), then \(c_\xi\ne0\) and
\[
 \left|\frac{c_\xi}{c_1}\right|^2
    = \mathrm{Nm}(Q(\xi))^{m^2t^2+m^2t+mr}
                       \prod_{i,j<m}T_{ij}(\xi)^2 .
 \tag{TP11}\label{eq:display-tp11}
\]
Here the sum of the line dimensions is
\(\sum_{j=0}^{m-1}(r+(2j+1)t)=mr+m^2t\).
At the trivial twist all critical ratios are equal to one.
This follows directly from the nonzero Hecke denominators
and the definition of the normalized ratios.

For the last path let \(I_{z,j}\) be the original finite
level index at node \((z,j,m)\). Every such index other
than \(I_{0,0}\) is the general linear group order at
\(\mathfrak b\), because the level there is the kernel of
\(\mathrm{GL}(R_b)\times\mathrm{GL}(R_Q)\to
\mathrm{GL}(R_b)\) in the corresponding rank.
At \((0,0,m)\) the level is the identity over \(R_\xi\)
and the rank is \(r+mt\), so \(I_{0,0}=|\mathrm{GL}_{r+mt}(R_\xi)|\).
The last nodes are the original definite models.

On each
edge with fixed polynomial degree \(D_e\), the scalar
analytic image of a source algebraic form is the analytic
realization of its jet lift times
\[
                        \frac{\beta_e}{s_e I_{z,j}} .
 \tag{TP12}\label{eq:display-tp12}
\]
In detail, the original jet matrix applied to a source
vector is paired with the highest target covector.
At this last-path edge the jet is taken at
\(s_{\rm jet}=p_Fa_d/\omega_0\), as in
Section~\ref{sec:completion-terminal-jets}. For each polynomial \(Q\)
of the edge's homogeneous degree \(D_e\),
\[
 Q(\partial_z)\Theta(s_{\rm jet}z)\big|_{z=0}
       =s_{\rm jet}^{D_e}\,Q(\partial_z)\Theta(z)\big|_{z=0}.
\]
Thus conversion to the normalized analytic joint harmonic gives
\(a_e s_{\rm jet}^{-D_e}\), with \(a_e\) the fixed harmonic scalar.

The original
compact real integral, using the unit highest source
vector, gives \(1/\dim V_{\rm src}\) by Schur orthogonality.
The finite class integral of the algebraic lift gives
the source level mass one, whereas in the analytic
measure its mass is \(I_{z,j}^{-1}\).
Finally the algebraic datum multiplies the same matrix
indicator by \(s_e=s_{\rm vac}^{\dim(\mathrm{Hom})}\).
These are precisely the original source, target and
polynomial conventions in the definite jet lift.

Thus \(\beta_e=a_e s_{\rm jet}^{-D_e}/\dim V_{\rm src}\)
and \(s_e\) are fixed independently of the residue
point and the twist. Iterating \eqref{eq:display-tp12} gives, after \(j\)
edges, the algebraic chain's analytic realization times
\(b_{z,j}I_{z,j}/I_{z,0}\), with
\[
              b_{z,0}=1,\qquad
              b_{z,j+1}=\frac{\beta_e b_{z,j}}{s_e I_{z,j+1}} .
\]
The \(b_{z,j}\) are fixed, because the last index in this
recursion has \(j+1>0\).
Here the integrals for an individual tuple are Haar integrals,
with each node level assigned mass one. A single point-indicator
tuple need only be fixed by a finer level. We take the finite array
sum before replacing these integrals by sums over the displayed
levels.

At the final diagonal the true form's finite point is the
block point \(\mathrm{diag}(h,u_\xi^m)\) in its given reindexing.
The extra line has sign \(+1\) since \(m\) is even.
The two nodes here have ranks \(N=r+2mt\) and \(N+1\).
We construct the bilinear branching pairing for their own weights.
Fix one infinite place, suppress it from the notation and write
\(w_i=\kappa_\tau+i\Delta+\Delta/2\) and
\(s_i=2\lambda_{k-1+i}=\Delta-2w_{i+1}\) for \(0\le i\le m\),
extending the same arithmetic progression for its final term.
The doubled parameters to interlace are the false parameter before
duality and the negative of the true grid parameter.

For an integer threshold \(b'\), write
\(\phi(y)=\mathbf1_{b'<y}\) and
\(g_n(y)=\max(0,\min(n,\lfloor(y-b')/2\rfloor))\).
A descending block \(d,d-2,\ldots,d-2(n-1)\) contributes
\(g_n(d+1)\) entries above the threshold. With \(i\) pairs of
blocks included, the counts on the false and true sides are respectively
\[
\begin{aligned}
 C_i^-&=\sum_{a<k-1+i}\phi(2\lambda_a)
       +\sum_{j=1}^i\{g_t(2w_j+t)+g_t(-2w_j+t)\},\\
 C_i^+&=\sum_{a<k+i}\phi(-2\lambda_a-1)
       +\sum_{j=1}^i\{g_{t+1}(-2w_j+t+1)
                                   +g_{t-1}(2w_j+t-1)\}.
\end{aligned}
\]
These are the base entries and the two progressions inserted at each
step of the diagonal grid, with lengths \(t,t\) on the false side
and \(t+1,t-1\) on the true side. We claim that
\[
 \phi(s_i)\le C_i^+-C_i^-\le\phi(-s_i-1).
\]
For \(i=0\), the central false progression is
\(k-2,k-4,\ldots,2-k\). The central true progression consists of the
same entries diminished by one, together with \(-s_0-1\) above it.
Their counts differ by at most one, and all central entries exceed
\(s_0\), which gives both endpoints.
For the induction use the following block count. The difference between
the counts of the centered progressions of lengths \(n+1,n\) at
center \(c\) is
\(g_{n+1}(c+n+1)-g_n(c+n)\). Their entries alternate, so this is
zero or one, equals one below \(c-n\), and vanishes at or above \(c+n\).
At \(c=-2w_{i+1}, n=t\), this difference lies between
\(\phi(s_{i+1})\) and \(\phi(s_i)\). At
\(c=2w_{i+1}, n=t-1\) it lies between
\(\phi(-s_i-1)\) and \(\phi(-s_{i+1}-1)\).
Indeed \(w_{i+2}-w_{i+1}=\Delta\) and \(t+1\le\Delta/2\).
In passing from \(i\) to \(i+1\), add the first difference,
subtract the second, and add \(\phi(-s_{i+1}-1)-\phi(s_i)\)
from the two new base entries. This proves the claim.
At \(i=m\), since \(r=k-1+m\), the counts of the actual terminal
parameters differ by zero or one at every threshold.

Let the descending highest weights of the larger module and the smaller
primal module be \(\Lambda_1,\ldots,\Lambda_{N+1}\) and
\(\nu_1,\ldots,\nu_N\), respectively. The doubled parameters
\(L_i^+\) and \(L_i^-\) of the two lists have opposite parity.
The base entries of the true list are odd and those of the
false list are even. Their inserted blocks have those same
parities since \(t\) is odd. Thus the count inequality gives
\[
 L_i^+\ge L_i^-+1,\qquad L_i^-\ge L_{i+1}^++1.
\]
Subtracting the two rank shifts gives
\(\Lambda_i=(L_i^+-N-2+2i)/2\) and
\(\nu_i=(L_i^--N-1+2i)/2\), hence
\(\Lambda_i\ge\nu_i\ge\Lambda_{i+1}\).
Put \(e=\Lambda_{N+1}\) and \(e_0=\nu_N\).
Of the \(\Lambda_i-\Lambda_{i+1}\) upper columns at height \(i\),
retain \(\nu_i-\Lambda_{i+1}\) and contract
\(\Lambda_i-\nu_i\) against the last unit vector.
For \(i<N\), the retained columns at \(i\) and contracted columns at
\(i+1\) give exactly \(\nu_i-\nu_{i+1}\) lower columns.
At \(i=N\), the retained full columns contribute \(e_0-e\) to the
lower determinant power. The extra line acts with exponent
\[
 L_\tau=e+\sum_{i=1}^N(\Lambda_i-\nu_i).
\]
This also follows by subtracting the two central weights.
The sum of the highest weights at each rank is half the sum of
its doubled entries. The false primal list therefore has sum of
highest weights \(\sum_{a<r}\lambda_a\); the true list has sum
\(-\sum_{a<r+1}\lambda_a-(r+1)/2-2\sum_{i=1}^m w_i\).
The central head has sum zero and
\(\lambda_{k-1+i}=\Delta/2-w_{i+1}\). Their difference in the
order true minus false gives
\[
 L_\tau=w_{m+1}-(2m+1)\Delta/2-(r+1)/2
       =\kappa_\tau+\Delta-(r+1)/2 .
\]

Let \(E_+,E_-\) be the integral inclusions of the upper and
lower Weyl lattices into their exterior tensor models. Let
\(C\) be the tensor of the retained-column projections and the
last-covector contractions with the above column counts. Reorder
columns and absorb the retained full lower columns into the
lower determinant line. Laplace expansion gives one last-line
factor per contraction. With \(G_-=E_-^{\mathsf t}E_-\) and
\(c_-\) the Fock square of the full lower highest polynomial,
take
\[
 b_{\rm grid}=c_-G_-^{-1} E_-^{\mathsf t} C E_+ .
\]
The Fock clearing calculation of
Section~\ref{sec:integral-terminal-branching}, applied to the
full lower diagram, gives \(c_-G_-^{-1}\) as the Fock matrix
on the integral dual polynomials. Its entries are integers.
The other factors in this formula are integral, and the Weyl
projection and the column operations commute with the smaller group.
The map is nonzero. Pairing the contracted evaluation tensor of a
generic upper flag with the unit lower highest tensor gives a
product of nonzero minors. A retained column uses its leading
rows, and a contracted column uses its shortened leading rows
and the last row.
The number of lower boxes after removing determinant powers is at most
the upper count \(H\) summed over all places, since
\(\nu_i-e_0\le\Lambda_i-e\). The minor and factorial bounds give
\(c_-\le H^{3H}\), with value one for \(H=0\). Before multiplication
by \(c_-\) the map is the composition of the exterior operations and
an orthogonal Weyl projection. Thus \(\|b_{\rm grid}\|\le H^{3H}\),
and \(|L_\tau|\le H+|e|\).

Compose \(b_{\rm grid}\) with the integral invariant duality of the
false node to obtain the integer bilinear matrix
\(B_{\rm grid}=J_{\rm diag} b_{\rm grid}\).
Here complementary coordinate wedges, with the volume sign, identify
the primal representation isometrically with the conjugate of its
dual; the dual tensor Gram matrix composed with this signed
identification gives \(J_{\rm diag}\), as in
Section~\ref{sec:completion-terminal-dual-coordinates} at this rank.
Cauchy--Schwarz bounds the resulting bilinear pairing by
\(H^{3H}\) times the product of the two Weyl norms.
Let \(K\) be the compact real false group, with probability Haar
measure, and \(\iota(k)=\operatorname{diag}(k,1)\) in the prescribed
reindexing. Write \(\rho_-\) for its action on the primal false
module and \(\rho_+\) for the true action. If \(h_0\) is the unit
highest vector of the primal false type and \(\ell_g\) is the true
highest covector after translation by \(g^{-1}\), the real integral
inserts the linear map
\[
 A_g=\int_K \rho_-(k)h_0\,\bigl(\ell_g\circ
            \rho_+(\iota(k)^{-1})\bigr)\,dk .
\]
The false dual coefficient evaluates on \(\rho_-(k)h_0\).
Left invariance gives \(\rho_-(k_0)A_g=A_g\rho_+(\iota(k_0))\)
for every \(k_0\in K\).
The branching line has multiplicity one. After removing the
common determinant power, the Schur character formula on the
union of an \(N\)-variable alphabet and the last variable has
one contribution for each horizontal strip
\cite[Ch.~I, (5.8)--(5.9), pp.~71--72]{MacdonaldHall}.
The strip condition is the interlacing just proved, and the
skew Schur polynomial in one variable is a single monomial.
Thus the average is \(a(g)b_{\rm grid}\). It is independent of
the twist and the finite point and can be zero.

At each
finite point this leaves the bilinear pairing
\(F_-(h)^{\mathsf t}B_{\rm grid} F_+(\mathrm{diag}(h,u_\xi^m))\).
The compact finite measure converts to the class integral
of level mass one with coefficient \(I_{0,m}^{-1}\).
For an individual locally constant summand, take an open
subgroup \(J'\subset J\) of the displayed level that fixes it.
A \(J'\)-coset has mass \([J:J']^{-1}\) in the Haar
normalization giving \(J\) mass one, and mass
\([K_{\rm std}:J']^{-1}\) in the original finite measure.
Dividing the latter mass by the former gives \([K_{\rm std}:J]^{-1}=I_{0,m}^{-1}\), also
when the individual summand has smaller stabilizer than \(J\).
After multiplying by the array, invariance under right
translations at the internal nodes follows by reindexing
their residue matrices.

The determinant phases cancel by \eqref{eq:display-tp9}.
More explicitly, write \(J\) for the product of the levels of
the \(2m+1\) independently integrated class variables, counting the
common terminal variable once. Its two actions at that variable are
\(u\) and \(\operatorname{diag}(u,1)\), the latter commuting with
\(\operatorname{diag}(1,u_\xi^m)\). Write \(F\) for the sum of
all the point-indicator integrands with their array coefficients.
At a change of the common terminal representative by \(k\),
\[
 \mathrm{diag}(hk,u_\xi^m)
       =\mathrm{diag}(h,u_\xi^m)\,\mathrm{diag}(k,1).
\]
Hence the extra-line coordinate stays at its original value.
For fixed outer tests the remaining kernel multipliers are the
node characters in the array covariance, and reindexing its
residue tuples gives \(F(hk)=F(h)\) for \(k\in J\). Choose an
open \(J'\subset J\) fixing the individual summands and let
\(\mathcal C(J')\), \(\mathcal C(J)\) be the respective sets of
class tuples. Each fiber over \(\mathcal C(J)\) has cardinality
\([J:J']\): every component of \(J\) lies in the old principal
subgroup at the split place above \(p_0\), so
Lemma~\ref{lem:marked-source-freeness} applies to each definite group. With \(J\) of mass one, we therefore obtain
\[
  \frac1{[J:J']}\sum_{h\in\mathcal C(J')}F(h)
                  =\sum_{h\in\mathcal C(J)}F(h).
\]
In particular no inverse level index occurs in the final class sum,
including at coefficient primes that divide \([J:J']\). The two outer
levels are respected: the false outer deep principal level is
contained in the original \(K(\pi_MQ(\xi))\), and the true
outer level is \(K(\mathfrak b)\).

There is one common terminal group variable. Its measure gives
one factor \(I_{0,m}^{-1}\). With both chain conversions the scalar is
\[
 I_{0,m}^{-1}\prod_{z=0}^1\frac{b_{z,m}I_{z,m}}{I_{z,0}}
   =\frac{b_{0,m}b_{1,m}}{I_{0,0}}\frac{I_{1,m}}{I_{1,0}} .
\]
Equivalently, before cancelling indices the scalar is
\[
 \frac1{I_{0,m}}\prod_{z=0}^1\prod_{j=0}^{m-1}
        \frac{\beta_{e(z,j)}}{s_{e(z,j)}I_{z,j}},
\]
where \(e(z,j)\) starts at \((z,j,m)\). The common diagonal variable
occurs once. Of these \(2m+1\) measure factors, only \(I_{0,0}\)
depends on the twist. By the same calculation as
\eqref{eq:display-tp6}, \(I_{0,0}\) differs from
\(I_\xi=|\mathrm{GL}_{r+mt}(R_0(\xi))|\) by the fixed
ratio \(I_\xi/I_{0,0}
 =(\mathrm{Nm}(\pi_M)/\mathrm{Nm}(\mathfrak b))^{(r+mt)^2}\).
Consequently the last-period function at any real \(g\)
is a scalar depending only on \(g\) and the datum times
\(\Xi_\xi(y',\eta)/I_\xi\), where \(\Xi_\xi\) is the
final array paired against those algebraic class sums.
Apply the composed backward implication to that function
and evaluate at \(\omega\) in \eqref{eq:terminal-transpose-restriction}. With one scalar \(A_*\)
independent of the twist and both forms, we obtain the following identity.

\Needspace{4\baselineskip}
\begin{proposition}[The measured terminal period]\label{prop:measured-deep-terminal-period}
Let \(\xi\) be a twist for the full chain as in
Section~\ref{sec:notation-theta}, let \(y'\) be a dual test at the
original full principal terminal level, and let \(\eta\in\mathcal H_m\).
With the original paths, real points and quotient measures, one scalar
\(A_*\), independent of \(\xi,y',\eta\), gives
\[
 I_\xi\int_{[Y_m]} y'(h)(T_\xi r_0\eta)(h)\,dh
             = A_*\,\Xi_\xi(y',\eta)/c_\xi .
 \tag{TP13}\label{eq:measured-terminal-period}
\]
Here \(c_\xi\ne0\) and \(I_\xi\) is the index of the original
false terminal test level, computed by surjectivity of split reduction.
\end{proposition}

\section{Chain determinants and completion of the lattice}
\label{sec:chain-determinants-completion}

\sectionopening{In this section, we compute the determinant by composing the individual lifts, using
the terminal periods of Section~\ref{sec:deeper-periods}. The same
calculation gives the auxiliary vectors and the clearing scalar needed in
Section~\ref{sec:gram-product-formula}.}

\subsection{Determinants of the composition}
\label{sec:marked-chain-composition}

\sectionopening{In this subsection, we compute the determinant of the weighted composition by reducing it
to the determinants of the individual lifts. This gives the ratio whose
height is estimated in Section~\ref{sec:marked-chain-ratio-height}.}

Fix \(q_0,d,\epsilon\) as in Section~\ref{sec:datum}, take a
sufficiently large integer \(x\), and fix the resulting admissible
datum. Write \(\kappa=(\ell,\iota)\) for a coefficient embedding. Let \(W\) be the original full
source of level \(K(\pi_M)\), and \(V\subset W\) its projected space of dimension \(D\). Write \(T_\xi=L_\xi\) for the
full weighted chain of Section~\ref{sec:seesaw-normalization}
on \(W\), and \(\Theta_\xi=T_\xi|_V\) for its restriction to the projected space at the allowed \(m\)-fold twist \(\xi\).
The source weights below are \(w_n=A_{\xi,n}\).
The terminal algebraic conversion is that of the same chain.
The target
inner product is the terminal class pairing.
Write
\[
 Z_\xi=\det(\Theta_\xi^*\Theta_\xi),\qquad
 Z_j(\xi^j)=\det(u_{j,\xi^j}^*u_{j,\xi^j}),
\]
where \(u_{j,\xi^j}\) is the \(j\)-th single lift on the
projected space of the specified space \(W_j\), also with its
target class pairing. The same source metrics are
used at all twists. Set
\[
 B=\log b,\qquad J_* =\frac{faB^3}{mk},\qquad
 L=faB^4,\qquad A=mk\,d_{\rm dir}.
\]
Here Section~\ref{sec:common-critical-twist}
gives
\[
 d_{\rm dir}=\sum_\tau
 \bigl(\log(\kappa_\tau!)
      +\kappa_\tau\log(M/\Omega)-S(\kappa_\tau,\tau)\bigr).
\]
Thus \(A=mk\,d_{\rm dir}\) is the same quantity as in
Section~\ref{sec:joint}. Each height bound uses one
multiplier for all coefficient embeddings.

\indent{\it Finite factors in the chain.}
For a twist \(\xi\) put \(q_\xi=\mathrm NQ(\xi)\) and let
\(w_n\) be the order of \(\mathrm{GL}_{r+nt}(\mathcal O_E/(\pi_MQ(\xi)))\).
Let \(a_{i,n}\) be the horizontal arrow and \(d_{i,n}\)
the downward arrow from column \(n+1\) of sequence \(i\) to
column \(n\) of sequence \(i+1\). Denote the composite of the
first \(n\) arrows of sequence \(i\), weighted by \(w_0,\ldots,w_{n-1}\),
by \(X_{i,n}\). Thus \(X_{i,0}=1\) and
\(X_{i,n+1}=w_n a_{i,n}X_{i,n}\) for \(i+n+1\le m\).
Let \(j_i\) be the isometric base identification of sequence \(i\).

\indent{\it Injectivity before dividing by the archimedean scalar.}
The first-source calculation in
Section~\ref{sec:marked-first-arrow-rallis} applies to every
\(y\in W\). In fact \(y\) is invariant under the point
stabilizer at \(\pi_MQ(\xi)\), which is contained in
\(K(\pi_M)\), and is invariant under the integral finite
group off the support. With \(A_i=a_{i,0}j_i\),
Proposition~\ref{prop:first-lift-gram} gives
\[
 w_0\langle A_ix,A_iy\rangle
     =C_i B_\xi^{-1}\langle x,\mathcal V_i y\rangle
 \qquad (x,y\in W).
\]
Put \(R_0=R_0(\xi)=\mathcal O_E/(\pi_MQ(\xi))\) and
\(B_\xi=|R_0|^{r(r+t)}\). The constant \(C_i\) is
nonzero, and \(\mathcal V_i\) is the original convergent spherical operator. Its off-identity operator norm is less
than one, by the neighbor count in that section. On the
projected subspace it is \(R_{i+1}\). In particular the
first arrow is injective on \(W\). Indeed, if \(A_i y=0\), taking
\(x=y\) gives \(\langle y,\mathcal V_i y\rangle=0\), whereas
\(\Re\langle y,\mathcal V_i y\rangle
 \ge(1-\|\mathcal V_i-1\|)\|y\|^2>0\) for \(y\ne0\).

For the descents at the original level pick the first datum by adjointness
on the two isometric copies, and propagate it by the finite
shear formula \eqref{eq:triangle-exact-next-datum} on \(R_0\).
The square gives the equality
\[
 w_n a_{i+1,n}d_{i,n}
      = u_{i,n}w_{n+2}d_{i,n+1}a_{i,n+1},\qquad
 u_{i,n}=k_{i,n}/c_{i,n},\quad c_{i,n}\ne0 .
 \tag{GC1}\label{eq:chain-gram-first-comparison}
\]
Here the measured partner scalar \(c_{i,n}\) is nonzero
from the start. The archimedean pairing factor \(k_{i,n}\) is
allowed to be zero in \eqref{eq:chain-gram-first-comparison}.
At the original triangle modulus \(\pi_MQ(\xi)\), the first-column
calculation in Section~\ref{sec:first-lift-spherical-actions} gives
\(R(u)^{-1}a_{i,0}^*a_{i,0}R(u)=a_{i,0}^*a_{i,0}\) for every
\(u\in K_0\). The isometry \(j_i\) intertwines these right
translations. Hence \(A_i^*A_i\) preserves the \(J_M\)-fixed
space \(W\) and commutes there with \(p=\rho(e_0)\), a finite
linear combination of such translations. The first adjoint identity gives
\[
 B_i=j_{i+1}^{-1}d_{i,0}a_{i,0}j_i=A_i^*A_i
               \qquad\text{on }W .
\]
The spherical series preserves \(W\) as well. Its complete double
cosets are bi-invariant at the unprescribed places and act by the
identity factor at the old and twist places. Testing the first-lift
identity against every \(x\in W\) therefore gives
\[
 B_i=\frac{C_i}{w_0B_\xi}\mathcal V_i\big|_W,\qquad
 B_i\big|_V=\frac{C_i}{w_0B_\xi}R_{i+1},\qquad V=pW .
\]
In particular \(B_i\) is injective on \(W\).
The adjoint here is represented by the original transposed
kernel on the indicated types, with the original Haar.

If \(X_{i,n}\) is the first \(n\) weighted arrows, iterating
\eqref{eq:chain-gram-first-comparison} from \(n=0\) gives
\[
 X_{i+1,n}d_{i,0}a_{i,0}
   = \frac{w_{n+1}}{w_0w_1}\!
       \left(\prod_{s<n}u_{i,s}\right)d_{i,n}X_{i,n+1}.
 \tag{GC2}\label{eq:chain-gram-propagation}
\]
For \(n=0\) this is \(d_{i,0}a_{i,0}=w_0^{-1}d_{i,0}X_{i,1}\).
At the induction step multiply the identity on the left by
\(w_n a_{i+1,n}\). Use \eqref{eq:chain-gram-first-comparison}
and the definition of the weighted composite to obtain
\[
 \begin{aligned}
 w_n a_{i+1,n}d_{i,n}&=u_{i,n}w_{n+2}d_{i,n+1}a_{i,n+1},\\
 w_{n+1}a_{i,n+1}X_{i,n+1}&=X_{i,n+2}.
 \end{aligned}
\]
In particular the implication \(X_{i,n+1}j_i y=0\Rightarrow
X_{i+1,n}j_{i+1}B_i y=0\) involves no inversion of \(k_{i,n}\).
Induct simultaneously over all admissible \(i\) on the chain
length. Injectivity at length \(n\) gives \(B_i y=0\) from the
right-hand vanishing, and the injectivity of \(B_i\) gives \(y=0\). If \(W\ne0\), the left
side of \eqref{eq:chain-gram-propagation} is nonzero on \(j_i y\) for every \(y\ne0\).
Taking \(n=s+1\) shows \(u_{i,s}\ne0\) for every square.
We may therefore use \(\beta'_{i,n}=\overline{c_{i,n}/k_{i,n}}\)
in the Gram computation. When \(D=0\), the empty determinants give the ratio formula directly,
since the individual projected spaces have the same dimension. In the
calculation of the ratio below we can thus assume \(D>0\), so \(W\ne0\).

\indent{\it The terminal level of the chain image.}
For an integral terminal element at level places, let \(g\)
be its reduction, upper block triangular of sizes
\(r,t,\ldots,t\), and \(g_j\) its upper left corner of rank
\(r+jt\). For the coordinate inclusion \(I_j\) from rank \(r+jt\)
to rank \(r+(j+1)t\), the following identities hold over \(R_0\):
\[
 g_{j+1}I_j=I_jg_j,\qquad
 \det(g_{j+1})/\det(g_j)=\det(g^{(j+1)}),
\]
where \(g^{(j+1)}\) denotes the last diagonal block of \(g_{j+1}\).
After lifting each corner integrally by split reduction, the first
identity holds modulo \(\pi_M Q(\xi)\), and the second determines
the tame determinant character on the lifted units.
The simultaneous Schr\"odinger action on the original
point indicator therefore moves the terminal translation
through arrow \(j\), multiplying by the \(j+1\)st lift
character at \(\det(g^{(j+1)})\).

This is the determinant
ratio of the target and source actions, with the same
tame character at the two ranks. Its other local factors
are trivial on integral units. Iterate this equality.
At old places the terminal element is the identity after
reduction, so \(g_0\) fixes the base form. At each new
place \(g_0\) is integral, so it again fixes the base form.
The resulting product of the diagonal-block characters
is the original terminal level character. The passage
from scalar forms to the algebraic compact model intertwines
these translations. Thus the terminal algebraic map with
this character is the very same injective weighted chain.

Write \(I_{\rm par}(\xi)\) for the index of this terminal
level. Its factor at a twist place is the index of the
block upper triangular group. The compact real Schur
average and the finite quotient volume convert the
Petersson norm of the forward chain to its algebraic
class norm by \(d_Y I_{\rm par}(\xi)\).
This differs from the dual-delta test conversion
\(d_Y I_\xi\), where \(I_\xi=|\mathrm{GL}_{r+mt}(R_0)|\).
For that period the dual test has the full principal
level at the twist places. More explicitly, write
\(K_{\rm pr}\subset K_{\rm par}\) for the two terminal levels and
\(G_{\rm pr}(F),G_{\rm par}(F)\) for the two class Gram matrices
of a tuple \(F=(T_\xi v_i)_i\). An old principal factor makes
the rational stabilizer trivial, as proved in
Section~\ref{sec:gram-product-formula}. Hence each
\(K_{\rm par}\)-class splits into exactly
\([K_{\rm par}:K_{\rm pr}]=I_\xi/I_{\rm par}(\xi)\)
principal classes. For \(k\in K_{\rm par}\), every \(F_i(hk)\)
differs from \(F_i(h)\) by the same unitary level character.
Its two factors cancel in each Hermitian matrix entry, giving
\[
 G_{\rm pr}(F)=\frac{I_\xi}{I_{\rm par}(\xi)}G_{\rm par}(F).
\]
The determinant \(Z_\xi\) uses \(G_{\rm par}\) in an orthonormal
source basis. The principal-level delta identity
\eqref{eq:original-dual-delta-readout} tests a single value
\(F_i(h)\), so its integrality assertion applies directly at the
parabolic class representatives, including at primes dividing the
index.

\indent{\it The adjoint square and the Gram scalar.}
The tested adjoint of the upper point arrow is the arrow
with conjugate kernel and opposite determinant
characters, by Fubini on the two original quotients.
Use the original reversed local kernels with Fourier
datum \(\mathcal F(\bar p)\), and let \(\mathcal R_n\)
be right translation by the inverse of the finite central
reversal scalar on rank \(r+nt\). In the notation where
\(D_0,D_1\) are the two descents and \(B_0',B_1'\) their
reversed lifts, the opposite-kernel identities in
Section~\ref{sec:marked-square-finite-extension} read
\[
 \begin{split}
 \langle D_0x',y\rangle
     &=s_0\langle x',B_0'\mathcal R_n y\rangle,\\
 \langle x,B_1'\mathcal R_{n+1}z\rangle
     &=s_1\langle D_1x,z\rangle .
 \end{split}
\]
They are tested identities and the real polynomial
replacement constants in \(s_0,s_1\) may vanish.

Put \(N=r+nt\). The reversed finite square on \(R_0\) gives
\(a_{\rm up}B_0'=\gamma B_1'a_{\rm low}\). Its integrated
nodes \((i+1,n+1)\) and \((i,n+1)\) both have rank \(N+t\).
Thus Proposition~\ref{prop:normalized-period-exchange} applies with
defect zero. The defining products give \(M_\nu(0)=1\), so
\( |c_\nu|^2=R'_0/D_{\rm final}>0\) for this square, with the
same ratio for every admissible tame modulus and character.
The remaining finite factor in \(\gamma\) is
\((\mathcal G|R_0|^{Nt})^{-1}\).
For the lower lift let \(a,b\) be its splitting characters in
the order \(K^{a,b}_{S,T}\), and put \(z=s_{\pi_M Q(\xi)}\).
The source and target ranks are \(N,N+t\). Simultaneous
translation by the finite central scalar \(z\) fixes the Hom variable
and the original residue indicator. The determinant characters give
\[
 K(gz,hz)=\eta K(g,h),\qquad
 \eta=a_{\rm f}(z^{N+t})b_{\rm f}(z^N)^{-1},\qquad |\eta|=1 .
\]
Centrality preserves both right levels. Changing the source variable to
\(g'=gz^{-1}\) in the original quotient Haar gives
\(a_{\rm low}\mathcal R_n=\eta\mathcal R_{n+1}a_{\rm low}\),
where both translations are by \(z^{-1}\). The deeper line period of Section~\ref{sec:deep-transport-arrays}
uses its own generator for \(\mathfrak b Q(\xi)\).
For \(x'=a_{\rm up}^*x\) this gives
\[
 \langle D_0a_{\rm up}^*x,y\rangle
       =s_0\gamma\eta s_1\langle D_1x,a_{\rm low}y\rangle .
 \tag{GC3}\label{eq:chain-terminal-gram}
\]
After removing the fixed real factors, the squared moduli of
\(s_0,\gamma,s_1\) have powers of \(|R_0|\)
\[
        -N(N+t),\qquad -t(t+2N),\qquad (N+t)(N+2t).
\]
The first and last come from \eqref{eq:display-mr}. For the middle
factor the rank defect is zero and the determinant block has size
\(b=t\). The inverse square modulus
\(|\mathcal G|^{-2}|R_0|^{-2Nt}=|R_0|^{-t(t+2N)}\)
comes from Section~\ref{sec:period-finite-determinant-sum}.

Indeed,
\[
 -N(N+t)-t(t+2N)+(N+t)(N+2t)=t^2.
\]
The adjoint square has no critical ratio. Its fixed real
factor may vanish during the iteration. The conjugate kernel represents the adjoint functional
tested on the original cusp types.

In the forward square the integrated ranks differ by \(2t\).
In Proposition~\ref{prop:normalized-period-exchange} the modulus
factor in \(M_\nu(t)\) is \(\prod_{w\in S}q_w^{t^2}\);
at the trivial twist only the old primes occur. The quotient of the
remaining \(L\)-value factors is \(Q_{i+1,i+n+2}(\xi)^2\),
where \(Q_{ij}(\xi)\) is the product of the \(t\) positive critical
ratios for the quotient character of lifts \(i,j\) at \(\xi^{j-i}\).
Thus \(|\beta'_{i,n}|=|\beta'_{i,n}(1)|\,
q_\xi^{t^2/2}Q_{i+1,i+n+2}(\xi)\) when \(W\ne0\). Put
\(\beta_{i,n}=\prod_{s<n}\beta'_{i,s}\) and
\(\epsilon'_{i,n}=s_0\gamma\eta s_1\) as in \eqref{eq:chain-terminal-gram}.

The corresponding adjoint law, starting with the
first-column adjointness, is
\[
 \langle v,a_{i,n}X_{i,n}j_i y\rangle
   =\epsilon_{i,n}\langle d_{i,n}v,X_{i+1,n}j_{i+1}y\rangle,
 \qquad
 \epsilon_{i,n}=\prod_{s<n}\epsilon'_{i,s}.
\]
At the step from \(n\) to \(n+1\), moving
\(a_{i,n+1}\) to the first slot and applying the preceding
length-\(n\) law gives
\[
 w_n\epsilon_{i,n}
 \langle d_{i,n}a_{i,n+1}^*v,X_{i+1,n}j_{i+1}y\rangle.
\]
By \eqref{eq:chain-terminal-gram} this equals
\(\epsilon_{i,n}\epsilon'_{i,n}
 \langle d_{i,n+1}v,X_{i+1,n+1}j_{i+1}y\rangle\),
since \(w_n a_{i+1,n}X_{i+1,n}=X_{i+1,n+1}\).
Let \(\mathcal B_{i,n}(x,y)=
\langle X_{i,n}j_i x,X_{i,n}j_i y\rangle\). Take
\(v=X_{i,n+1}j_i x\) in the adjoint law, retaining the real
weight \(w_n\) in the second argument. By
\eqref{eq:chain-gram-propagation},
\[
 d_{i,n}v=\frac{w_0w_1}{w_{n+1}}\overline{\beta_{i,n}}\,
                X_{i+1,n}j_{i+1}B_i x.
\]
The first argument of the inner product is conjugate-linear, hence
\[
 \mathcal B_{i,n+1}(x,y)=
 \epsilon_{i,n}\beta_{i,n}\frac{w_0w_1w_n}{w_{n+1}}\,
            \mathcal B_{i+1,n}(B_i x,y).
\]
For \(W\ne0\), the left side is positive at \(x=y\ne0\);
in particular the scalar on the right is nonzero whenever it is used.

Carry out the recursion with first argument in \(W\) and second
argument in \(V\). All \(B_i\) preserve \(W\), and all
\(R_j\) preserve \(V\). For \(x\in W\), \(v\in V\),
self-adjointness and the first-lift formula give
\[
 w_0\langle B_i x,v\rangle
   =w_0\langle A_i x,A_i v\rangle
   =C_iB_\xi^{-1}\langle x,R_{i+1}v\rangle .
\]

Consequently the terminal Petersson Gram form is
\[
 g_\xi\frac{w_0w_1^m}{w_m B_\xi^m}
                 \langle x,R_1\cdots R_m y\rangle .
\]
Indeed each substitution of the first-lift formula contributes
\(C_i/(w_0 B_\xi)\) and left-multiplies the remaining operator
by \(R_{i+1}\). The residual weight in the step at \(n\) is
\(w_1w_n/w_{n+1}\), and its product over \(0\le n<m\) is
\(w_0 w_1^m/w_m\). At the turn from sequence \(i\) the squares
are indexed by \(0\le s\le m-i-2\). The map
\((i,s)\mapsto(a,b)=(i+1,i+s+2)\) lists each pair
\(1\le a<b\le m\) once, with twist power \(b-a=s+1\).
For each such square, positivity of the
recursion at lengths \(s+1,s+2\) gives
\(\epsilon_{i,s}\epsilon_{i,s+1}\ne0\), hence
\(\epsilon'_{i,s}\ne0\), at both twists. The three
squared-modulus powers in \eqref{eq:chain-terminal-gram} therefore give
\( |\epsilon'_{i,s}(\xi)/\epsilon'_{i,s}(1)|=q_\xi^{t^2/2}\),
using \(|R_0(\xi)|=|R_0(1)|q_\xi\) by the Chinese remainder theorem.
Multiplication by the corresponding ratio for \(\beta'_{i,s}\)
gives \(q_\xi^{t^2}Q_{i+1,i+s+2}(\xi)\). Thus
\[
 \left|\frac{g_\xi}{g_1}\right|
    =q_\xi^{t^2\binom m2}\prod_{1\le i<j\le m}Q_{ij}(\xi).
\]
The Rallis operators act in the indicated order.
The squared norm at a nonzero \(y\in V\) is positive
by injectivity. Thus the product of archimedean factors in the final Gram form
and its untwisted value are nonzero on \(V\ne0\).

At \(w\in Q(\xi)\) put \(G_n=\mathrm{GL}_{r+nt}(k_w)\),
\(g_n=|G_n|\) and \(g_t'=|\mathrm{GL}_t(k_w)|\).
With \(P_n\subseteq G_n\) the block upper triangular subgroup of
sizes \(r,t,\ldots,t\), its unipotent radical has
\(nrt+\binom n2t^2\) free entries. Hence
\[
 |P_n|=g_0(g_t')^n q_w^{nrt+\binom n2t^2},\qquad
 \delta_w=[G_1:P_1]g_0q_w^{rt}=g_1/g_t'.
\]
It follows that
\[
 [G_m:P_m]
   =q_w^{-mrt-\binom m2t^2}\delta_w^m\frac{g_m}{g_0g_1^m}.
\]
At the fixed chosen places the terminal index is the full general
linear group order.

The Petersson-to-class factor is \(d_YI_{\rm par}(\xi)\), with
\(d_Y\) fixed. Put \(B_w=q_w^{r(r+t)}\) and
\(F_w=\prod_{t<h\le t+r}(1-q_w^{-h})\). At a new place the
node weights, terminal index, square factors and matrix counts give
\[
 \frac{g_0g_1^m}{g_m}[G_m:P_m]q_w^{t^2\binom m2}B_w^{-m}
   =q_w^{-mrt}(\delta_w/B_w)^m
   =q_w^{-mrt}(q_w^{rt}F_w)^m.
\]
Indeed \(\delta_w=g_1/g_t'\) and
\(|\mathrm{GL}_a(k_w)|=q_w^{a^2}\prod_{h=1}^a(1-q_w^{-h})\).
At an old place the
pointwise fixer of the inclusion has index equal to
the number of full-rank matrices, by transitivity.

The sharp mass has the additional inverse source
index \(g_0^{-1}\) (with \(g_0=|\mathrm{GL}_r(k_w)|\) at that
place), so the fixed old factor for each single lift is
\(g_0^{-1}\prod_{t<h\le t+r}(1-q_w^{-h})\).
With \(c_Q(\xi)=\prod_{w\in Q(\xi)}q_w^{rt}F_w\), the Gram
operator of the chain is
\[
                     C_\xi c_Q(\xi)^m R_1\cdots R_m ,
\qquad
 |C_\xi|=C q_\xi^{-mrt}\prod_{1\le i<j\le m}Q_{ij}(\xi),
\qquad C>0 .
\]

The single-lift identities for this comparison are the
Gram operators \(C_j c_Q(\xi) R_j^W\), where
\(C_j\) is independent of the twist. The sharp datum
at the chosen places fixes the inclusion matrix, and
the datum at the twist places has invertible first
block and arbitrary lower block. The Haar index and
the squared residue mass therefore give exactly the
same \(c_Q(\xi)\) above.

It is the same at every
\(\xi^j\) since these powers have the full modulus.
For the Rallis series put \(A\) equal to the complement
of the old and twist supports, including their conjugates.
Write \(K=K_A\). For \(\gamma\in G_A\), let
\(H_\gamma=K\cap\gamma K\gamma^{-1}\) and
\(m_\gamma=[K:H_\gamma]\). Writing \(K=\bigsqcup_{i=1}^{m_\gamma}
k_iH_\gamma\) gives \(K\gamma K=\bigsqcup_i k_i\gamma K\).
For each \(K\)-fixed vector \(y\), integration over these cosets gives
\[
 R\!\left(\frac{\mathbf1_{K\gamma K}}{\operatorname{vol}(K)}\right)y
   =\sum_{i=1}^{m_\gamma}R(k_i\gamma)y
   =m_\gamma P_K R(\gamma)y ,
\]
where \(P_K\) is the probability average over \(K\).
By the bi-invariance of \(u_j\) proved in
Section~\ref{sec:first-lift-doubling-series}, for each finite
\(\mathcal D\subset K_A\backslash G_A/K_A\) the partial series on both
spaces is the action of the same finite Hecke test
\[
 h_{\mathcal D}=
 \sum_{K_A\gamma K_A\in\mathcal D}
 u_j(\gamma)\frac{\mathbf1_{K_A\gamma K_A}}
                  {\operatorname{vol}(K_A)} ,
\]
extended by normalized integral identity factors at the twist
places.

This expression uses every finite place of \(A\),
with the same original bases and vacuum coefficient.
The norm of each normalized double-coset operator is at most
\(m_\gamma\) by unitarity. Thus, on either projected space, the
operator-norm error is bounded by
\[
 \sum_{K\gamma K\notin\mathcal D}m_\gamma|u_j(\gamma)|
      =\sum_{gK:\,KgK\notin\mathcal D}|u_j(g)|.
\]
This is the tail of the absolutely convergent series of
Section~\ref{sec:first-lift-doubling-series}. Hence the same tests
converge to \(R_j\) and \(R_j^W\), respectively. For every positive integer
\(s\), apply \eqref{eq:marked-full-spherical-trace} to the
convolution power \(h_{\mathcal D}^{*s}\), and take the limit.
Continuity of powers and traces in finite dimensions gives
\[
                    \operatorname{tr}(R_j^s)
                         =\operatorname{tr}((R_j^W)^s).
\]
Taking the integral identity test in \eqref{eq:marked-full-spherical-trace}
gives the common dimension \(D\). Newton's identities give the same
characteristic polynomial and hence \(\det R_j=\det R_j^W\).
This limit compares the series in the original complex Hilbert spaces.
All subsequent coefficient embeddings are applied to the fixed height
multipliers or to the scalar ratios after the Rallis determinants cancel.

Set \(W_\xi=q_\xi^{-mrt}\prod_{i<j}Q_{ij}(\xi)\). The two
factorizations just proved give
\[
\begin{aligned}
 Z_\xi &= C^D W_\xi^D c_Q(\xi)^{mD}\prod_j|\det R_j(\xi)|,\\
 \prod_j Z_j(\xi^j) &=
       \left(\prod_j |C_j|\right)^D c_Q(\xi)^{mD}
                                \prod_j|\det R_j^W(\xi)|.
\end{aligned}
\]
The factor \(c_Q(\xi)^{mD}\) is identical on the two sides
because \(Q(\xi^j)=Q(\xi)\) for all \(1\le j\le m\).
At \(w\in Q(\xi_\kappa)\) the numbers \(q_w^h-1\) in
\(F_w\) can be divisible by \(\ell\). The common factor
cancels in the complex equality before taking valuations.
For \(D>0\), \(Z_j(\xi^j)>0\) by injectivity of the single lifts,
and \(Z_\xi>0\) by \eqref{eq:chain-gram-propagation}. The same holds
at the trivial twist. The two factorizations and
\(\det R_j=\det R_j^W\) give
\[
 \frac{Z_\xi}{\prod_j Z_j(\xi^j)}
          =\left(\frac{C}{\prod_j|C_j|}\right)^D W_\xi^D .
\]
Compare this quotient at \(\xi\) and at \(1\), where \(W_1=1\).
\begin{proposition}[The chain determinant ratio]\label{prop:chain-determinant-ratio}
For a twist \(\xi\) for the full chain as in
Section~\ref{sec:notation-theta}, the weighted chain and its
single lifts satisfy
\begin{equation}\label{eq:marked-chain-height-ratio}
 \frac{Z_1}{Z_\xi}
   =\left(q_\xi^{-mrt}\prod_{1\le i<j\le m}Q_{ij}(\xi)\right)^{-D}
                            \prod_{j=1}^m\frac{Z_j(1)}{Z_j(\xi^j)} .
\end{equation}
The formula includes the empty determinant.
\end{proposition}

\subsection{Height of the chain determinant ratio}
\label{sec:marked-chain-ratio-height}

\sectionopening{In this subsection, we apply the critical-value comparisons of Section~\ref{sec:common-critical-twist} to Proposition~\ref{prop:chain-determinant-ratio}.}

Use the family \(\xi_\kappa\) of
Section~\ref{sec:common-critical-twist}
for the finite list of critical inverse ratios and all single
lifts. Its modulus is prime to \(\ell\) and every elementary
critical inverse ratio in the square and turn lists
has height at most \(x_1\). On this same family,
Section~\ref{sec:original-compression-bounds} gives operators
\(U_j(\xi)\) with
\[
 u_{j,\xi^j}^*U_j(\xi)u_{j,\xi^j}
                    =\rho_j(\xi)U_j^\circ
\]
for a fixed positive definite \(U_j^\circ\).
The ratio in Lemma~\ref{lem:marked-compression-determinant} is
\[
 \rho_j(\xi)=\rho_{\xi^j}=\frac{H_j(\xi^j)}{H_j(1)}.
\]
The exchange-product comparison in Section~\ref{sec:critical-exchange-height}
bounds its reciprocal on this family.
Write \(u=u_{j,\xi^j}\) and let \(\Pi\) be orthogonal
projection onto its image.

Set \(\Delta_j(\xi)=
\det(\Pi U_j(\xi)|_{\operatorname{im}u})\).
On the source, conjugate the compression by the
isomorphism \(u\) onto its image. Multiplying the
conjugated operator by \(u^*u\) gives
\(u^*\Pi U_j(\xi)u=u^*U_j(\xi)u\).
Taking determinants gives
\[
        Z_j(\xi^j)\Delta_j(\xi)=C_j'\rho_j(\xi)^D,
 \qquad C_j'=\det U_j^\circ\ne0 .
\]
The required bounds for this comparison are
\[
 \begin{split}
 \log|\sigma\Delta_j(1)|&\le 3DfmbB^2
                                      \quad\text{for all }\sigma,\\
 \log|\Delta_j(\xi_\kappa)|_\kappa
                         &\le Dr\beta_{j,\ell},\qquad
 \sum_{\ell\in U}\beta_{j,\ell}\le J_*\quad(U\text{ finite}),\\
 \operatorname{Height}(\rho_j(\xi_\kappa)^{-1})
                         &\le 10mx_1+k(6J_*-d_{\rm dir}).
 \end{split}
\]
The second line is a logarithmic allowance, independent
of the embedding above \(\ell\). Taking positive parts
of these allowances preserves the finite-sum bound:
apply that bound to the subset of primes with positive
allowance. These are respectively
the untwisted compression estimate, the tested clearing
estimate, and the fixed numerator and inverse ratio
bound on the same regular family.

Dividing the comparison identity by its untwisted
value gives
\[
 \frac{Z_j(1)}{Z_j(\xi_\kappa^j)}
    =\frac{\Delta_j(\xi_\kappa)}{\Delta_j(1)}
                               \rho_j(\xi_\kappa)^{-D}.
\]
For the first ratio use \(\Delta_j(1)\) as its fixed
height multiplier. Multiplication gives the allowance
\[
    D\{3fmbB^2+rJ_*+10mx_1+k(6J_*-d_{\rm dir})\}.
\]
On dividing the four positive terms inside the braces by \(L/m\), using
\(a=xb\), their sum becomes
\[
 \frac{3m^2}{xB^2}+\frac{r}{kB}
    +\frac{10m^2x_1}{f x b B^4}+\frac6B .
\]
It tends to zero under the admissible choices
\(m^2\le4x\), \(r/k\le2\),
\(\log b\sim x^\epsilon\), and
\(\log x_1\le x^{\epsilon/6}\).
This gives height at most \(D(L/m-kd_{\rm dir})\).
For each lift index \(j\), the third measured exchange of
Section~\ref{sec:gluing-three-exchanges} has \(k/2\) points
\(7m/2<u\le(7m+k)/2\) in its upper interval. At each point
Proposition~\ref{prop:main-critical-height} gives the term
\(-d_{\rm dir}\), with multiplicity two in the squared exchange
product of Section~\ref{sec:critical-exchange-height}.
These terms inside the factor \(D\) sum over all lifts to
\[
  2\sum_{j=1}^m\sum_{i=1}^{k/2}(-d_{\rm dir})=-mkd_{\rm dir}=-A.
\]

For \(D>0\) let \(h_j\ne0\) be the fixed multiplier of
Section~\ref{sec:critical-exchange-height}. The determinant
identity at \(1\) has \(Z_j(1)>0\), \(C_j'\ne0\) and
\(\rho_j(1)=1\), so \(\Delta_j(1)\ne0\). A single multiplier
for the entire ratio is
\[
 H=\prod_{j=1}^{m}\bigl(\Delta_j(1)h_j^{D}\bigr).
\]
For \(D=0\) both determinants are empty and the multiplier is one.
Indeed, with \(\xi=\xi_\kappa\),
\[
 H\frac{Z_1}{Z_\xi}
 =q_\xi^{mrtD}\prod_{i<j}Q_{ij}(\xi)^{-D}
   \prod_{j=1}^m
   \left[\Delta_j(\xi)\bigl(h_j\rho_j(\xi)^{-1}\bigr)^D\right].
\]
For \(D>0\) and \(\kappa\mid\ell\), equations
\eqref{eq:display-rt1}--\eqref{eq:display-rt2} and the
signed product bound of Section~\ref{sec:critical-exchange-height} give
\[
 \begin{aligned}
 \frac{1}{D}v_\kappa\!\left(H\frac{Z_1}{Z_{\xi_\kappa}}\right)
 &\ge -mrC_\ell-mk(B_\ell+2e_\ell+2P_\ell)\\
 &\qquad-\bigl(10m^2+t\tbinom m2\bigr)e_\ell .
 \end{aligned}
\]
Here \(Q_{ij}\) uses \(\xi_\kappa^{j-i}\) and \(\rho_j\) uses
\(\xi_\kappa^j\). The conjugated coefficients of the main ratios
are evaluated at \(\bar\kappa\) on the same \(\xi_\kappa^j\),
as in \eqref{eq:display-rt4}--\eqref{eq:display-rt5}.
The right side involves only the datum at \(x\) and the
rational prime \(\ell\).
The archimedean bound for \(H\) uses
\(\log|\sigma H|=\sum_j
 (\log|\sigma\Delta_j(1)|+D\log|\sigma h_j|)\).
For each \(\kappa\mid\ell\), the integer \(q_{\xi_\kappa}\)
is prime to \(\ell\). There are exactly \(t\binom m2\)
elementary inverses in \(W_{\xi_\kappa}^{-1}\).
By \eqref{eq:display-rt1},
\[
 v_\kappa(W_{\xi_\kappa}^{-1})
  =-\sum_{i<j}v_\kappa Q_{ij}(\xi_\kappa)
  \ge-t\binom m2 e_\ell .
\]
Thus \eqref{eq:display-rt3} bounds the height of the triangular
factor by \(t\binom m2 x_1\), with multiplier \(1\), before the
\(D\)-th power. The upper bound is
\[
 Dt\binom m2 x_1+\sum_{j=1}^m D(L/m-kd_{\rm dir})\le D(b+L-A),
 \qquad A=mk\,d_{\rm dir}.
\]
\begin{proposition}[Uniform height of the chain ratio]\label{prop:chain-determinant-height}
Fix \(q_0,d,\epsilon\) as in Section~\ref{sec:datum}. For every
sufficiently large integer \(x\), choose the common
twists \((\xi_\kappa)_\kappa\) as in Section~\ref{sec:common-critical-twist}.
Then
\begin{equation}\label{eq:marked-chain-height-bound}
       \operatorname{Height}\bigl((Z_1/Z_{\xi_\kappa})_\kappa\bigr)
                       \le D(b+L-A).
\end{equation}
\end{proposition}
Here \(m^2t\) grows at most polynomially in \(x\), whereas
\(B=\log b\sim x^\epsilon\) and \(\log x_1\le x^{\epsilon/6}\).
Thus \(2m^2tx_1\le b\) eventually. The cutoff in \(x\) may depend on
\(q_0,d,\epsilon\); it is chosen before any coefficient embedding
or auxiliary twist.

\subsection{Spanning vectors from the terminal period}
\label{sec:marked-chain-completion-bounds}

\sectionopening{In this subsection, we use Proposition~\ref{prop:measured-deep-terminal-period} to bound the norms and denominators of the vectors used in Section~\ref{sec:gram-product-formula}.}

Let \(W\) be the full principal-level source, with
projected subspace \(V\), and \(T_\xi\) its chain.
Let \(\mathcal H\) be the algebraic space at the true outer
node, of level \(K(\pi_M^{b_0})\). Use the map
\(r_0:\mathcal H\to W\), array \(\Psi_\xi\), and scalar
\(c_\xi\) of Section~\ref{sec:deeper-periods}.
The deep transpose and separation prove surjectivity of
\(r_0\), with the same period at all the false
levels \(K(\pi_MQ(\xi))\). The terminal dual delta identity
below gives, for a fixed \(\beta_0\),
\[
 v^{\mathsf t}J_{\rm W}
       \bigl((c_\xi/c_1)T_\xi r_0\eta\bigr)(x)
       =\beta_0\,\Xi_\xi(\delta_x(v),\eta).
\]
The terminal delta and \(\eta\) have the original false and
true outer levels. By \eqref{eq:display-tp11},
\[
 |c_\xi/c_1|=q_\xi^{(m^2t^2+m^2t+mr)/2}
            \prod_{0\le i,j<m} T_{ij}(\xi).
\]
The conductor exponent is integral since \(m\) is even.
Each \(T_{ij}\) consists of the \(t\) critical ratios at power
\(j+1\) on that original square. For \(i\ne j\) the reflected
square \((j,i)\) has the same untwisted character and uses
power \(i+1\).

\smallskip
\subsubsection{Bounds for the two terminal jets}
\label{sec:completion-terminal-jets}

\sectionopening{In this subsection, we apply the jet estimates of Section~\ref{sec:upper-chain-size} to both factors in the measured period.}

Put \(D_0=M^{b_0}\), \(\mathfrak b=(\pi_M)^{b_0}\), and
\(\mathfrak m=\mathfrak b Q(\xi)\).
Use the original edge in either half, with source and target
ranks \(n,n'\), respectively, and fix its residue matrix
\(S\in\operatorname{Mat}_{n'\times n}(\mathcal O_E/\mathfrak m)\).
The indicated point test is the characteristic function
\(\mathbf1_{\{S\}}\) on this matrix ring. On its original
standard definite Hom space, the edge datum \(\Phi_S\) is this
test on the vacuum already scaled by \(s_{\rm vac}^{nn'}\).
Its jet scale is \(s_{\rm jet}\).
The \(\pi_M\)-factor in the earlier forward-chain jet scale
does not enter the last-path kernel.

Write \(z=0,1\) for the false and true halves and put
\[
 p_j=r+z+j(t+z)+m(t-z),\qquad
 \delta_j=(1-z)((j+m)\bmod 2),\qquad 0\le j\le m .
\]
At each \(\tau\), let \(A(c,h)=\{c+h-1-2a:0\le a<h\}\) as a multiset.
With the base entries \(\lambda_a\) of
Section~\ref{sec:marked-deep-separation}, the doubled parameters of the
algebraic path module, in decreasing order, are the entries of
\[
 -\left(\{2\lambda_a+z:0\le a<r+z\}
       \sqcup\bigsqcup_{i=1}^{j} A(2w_i,t+z)
       \sqcup\bigsqcup_{i=1}^{m} A(-2w_i,t-z)\right)+\delta_j .
\]
Indeed \((z,j,m)\) has \(j\) up blocks and \(m\) down blocks
of the respective lengths \(t+z,t-z\) in
Section~\ref{sec:four-unitary-real-projections}; its lowering is
\(\delta_j\). The definite up edge uses the negative of that
lowered parameter. In particular at \(z=j=0\) this module is
the dual of the original outer \(Y_m\)-module, of rank \(r+mt\).
At \(j=m\) the false primal module for the branching pairing has
the opposite parameter, while the true module retains the displayed
parameter, in ranks \(r+2mt,r+2mt+1\), respectively.
All entries have absolute value at most \(P=\lceil20b\rceil\):
\(1\le\kappa_\tau\le b\), \(w_i=\kappa_\tau+(i+1/2)\Delta\), and
\((m+2)\Delta+t+k=o(b)\).

For either sign of the parameter, write its sorted list as
\(L_1>\cdots>L_p\), with rank \(p=p_j\). The highest weight,
determinant exponent and polynomial row lengths are
\[
 h_i=\tfrac12(L_i-p-1+2i),\quad e=h_p=\tfrac12(L_p+p-1),
 \quad h_i-e=\tfrac12(L_i-L_p)-(p-i)\in[0,P].
\]
The entries have a common parity, so each successive difference is at
least two. Since \(p\le 2n_m\le P\), also \(|e|\le P\).

For the displayed path parameter, in passing from \(j\) to \(j+1\)
the old entries acquire the common shift \(\delta_{j+1}-\delta_j\).
The new block lies below them since
\(w_1\gg k+t\) and \(w_{j+1}-w_j=\Delta>t+1\).
Thus the column counts \((L_i-L_{i+1})/2-1\) below height \(p_j\)
are unchanged and the counts strictly between \(p_j\) and \(p_{j+1}\)
vanish. To check the full-height factors, put
\(c_j(p)=\lfloor(p+1-z)/2\rfloor-w_{j+1}\). At \(j=0\),
the last entry is \(-(k-2+z)\) and
\(e_0=m(t+1-z)/2\). At \(j+1\) the last entry is
\(-2w_{j+1}-(t+z-1)+\delta_{j+1}\), giving
\(e_{j+1}=c_j(p_j)\). At the joining height \(p_j\), the difference
of the two adjacent weights is \(e_j-c_j(p_{j+1})\).
Writing \(n=p_j,\ n'=p_{j+1}=n+t+z\), this number of
height-\(n\) minors is
\[
 e_j-c_j(n')=
 \begin{cases}
 w_1-(k+t+2z-1)/2,&j=0,\\
 \Delta-(t+z),&1\le j<m .
 \end{cases}
\]
Both are positive. More explicitly, let \(h_i^{(j)}\) be the
source weights. On the retained indices \(i\le n\) the change of
weight is
\[
 h_i^{(j+1)}-h_i^{(j)}
    =(\delta_{j+1}-\delta_j-t-z)/2=c_j(n)-c_j(n').
\]
On the new indices \(n<i\le n'\), every weight equals \(c_j(n)\).
Thus the rectangular Schur map uses one partition
\(\lambda_i=h_i^{(j)}-c_j(n')\) for \(i\le n\), padded by zeros
in rank \(n'\). The source representation is its Schur module
times \(\det^{c_j(n')}\), and the target representation its
padded Schur module times \(\det^{c_j(n)}\). This gives both
determinant powers on the actual edge.

Every nonzero Schur jet entry is a linear combination of products
with exactly the prescribed number of minors of each height below
\(n'=p_{j+1}\), including the minors of height \(n=p_j\).
Its total degree, summed over all places, therefore satisfies
\[
 H=\sum_{\tau}\sum_{i=1}^{n'}(h_{\tau,i}-e_\tau)
   \le f(2n_m)P\le b_0r^4,
 \qquad n_m=k+m(t+1)=r+mt+1 .
\]
Indeed \(f\le x^{12}\), \(n_m\le x^6\), \(r\ge x^5\) and
\(b_0\ge b/2\) eventually.

The dual has the weights in reverse order and with opposite sign.
At each place put
\(c=h_1-h_p\le P\) and \(\lambda_i=h_i-h_p\).
The dual determinant exponent is \(-h_1\) and its row lengths satisfy
\[
 \lambda_i^\vee=c-\lambda_{p+1-i},\qquad
 \sum_{i=1}^p\lambda_i^\vee=pc-\sum_{i=1}^p\lambda_i\le pP .
\]
Equivalently, each column of height \(d_i\) is replaced by one of
height \(p-d_i\), in reverse order. Summing over the places gives
at most \(f(2n_m)P\le b_0r^4\) boxes in the dual terminal diagram
of Section~\ref{sec:completion-terminal-dual-coordinates}, including
the case \(c=0\).

Fix a prime \(\ell\nmid q_\xi\) and an arbitrary coefficient embedding. The integer \(q_\xi D_0\) belongs to \(\mathfrak m\).
Since \(D_0\mathcal O_E\subseteq\mathfrak b\) and
\((\mathfrak b,Q(\xi))=1\), the Fourier denominator of the point test is
\[
 |\operatorname{Mat}_{n'\times n}(\mathcal O_E/\mathfrak m)|
 = (\operatorname{Nm}(\mathfrak b)q_\xi)^{nn'},\qquad
 \operatorname{Nm}(\mathfrak b)\mid D_0^{2f}.
\]
Its \(\ell\)-part therefore divides \(D_0^{2fnn'}\) for every
\(\ell\nmid q_\xi\), independently of the size of \(q_\xi\).
For \(B_{\mathfrak m}=\operatorname{Mat}_{n'\times n}(\mathcal O_E/\mathfrak m)\)
Fourier inversion gives
\[
  D_0^{2fnn'}\mathbf1_{\{S\}}(X)
   =\frac{D_0^{2fnn'}}{|B_{\mathfrak m}|}
       \sum_{\chi\in\widehat B_{\mathfrak m}}\chi(-S)\chi(X).
\]
The prefactor is integral above \(\ell\) and the character values
are roots of unity. Each \(\chi(X)\) on the vacuum is the modulation
in the dual half of the Heisenberg polarization, up to a root of unity
in the signed convention of Section~\ref{sec:marked-both-norm-jets}.
Thus \(D_0^{2fnn'}\Phi_S\) is an integral combination of these finite
translates of the normalized vacuum. For every
translating vector \(a\) one has \(D_0 a_w\in\mathcal O_{E,w}^{nn'}\)
at each \(w\mid\ell\), and the \(\ell\)-part of the Fourier
denominator divides \(D_0^{2f nn'}\). The residue test is on the
original common modulus \(\mathfrak m=\mathfrak bQ(\xi)\).
In these normalized coordinates the commutator of the translations
by \(x,y\) is \(\psi(y^*x)\): its two half-pairing factors are
\(\psi(y^*x/2)\) and \(\psi(x^*y/2)^{-1}\).
Above \(\ell\nmid q_\xi\) the Fourier character is trivial on
\(D_0\Lambda_\ell\); self-duality therefore gives
\(D_0 x\in\Lambda_\ell\), also above \(2\). The point indicator
has values zero and one at every coefficient embedding.

Let the two finite points be integral at all places above \(\ell\),
and let \(g\) be their induced isometry of the finite Hom space
in the normalized coordinates. For the vacuum lattice \(L\) put
\(H=L\cap g^{-1}L\). The finite vacuum expansion of
Section~\ref{sec:seesaw-integral-tensor} expresses
\(\omega(g)T(x)v_L\) as translated vacua with shifts \(g(x+z)\),
\(z\in L/H\), coefficients differing by roots of unity, and common
coefficient of modulus \([L:H]^{-1/2}\) in the original complex
normalization. At the coefficient embedding we use its actual
value \(\Delta(g)/[L:H]\). Since \(H_\ell=L_\ell\), the index is
prime to \(\ell\) and all the \(z\) may be chosen zero above \(\ell\).
Thus \(D_0g(x+z)\) is integral there whenever \(D_0x\) is.

In the formula for \(\Delta(g)\) the intersection determinant
is a finite idele \(t\) that is a unit above \(\ell\). Its absolute
norm is an \(\ell\)-unit. Its splitting factor \(\upsilon_f(t)\)
is an \(\ell\)-unit as well: raise \(t\) to a ray class exponent and
write \(t^h=a k\) with \(a\in E^\times\) and \(k\) in the ray units.
Then every conjugate of \(a\) is a unit above \(\ell\), and the square
of the archimedean angular character on \(a\) is a monomial in those
conjugates and their inverses. Automorphy gives the assertion.
For each endpoint determinant, which is a norm-one idele integral
above \(\ell\), the same unit assertion is exactly the character
calculation in Section~\ref{sec:critical-integral-pairings}.
Hence the full adelic coefficients are integral at the chosen
embedding, including when the points move the vacuum at other places.
This proves preservation of the standard data of level \(D_0\). At the scale
\(D_0 d\,s_{\rm jet}\), the ordinary integral jet calculation
applies to each Schur matrix entry, since the Schur
polynomials have integral coefficients and degree \(H\).

At the actual scale \(s_{\rm jet}\) homogeneity introduces
\(|D_0d|_\ell^{-H}\). Multiplication by
\(D_0^{2f nn'}(D_0d)^H\) therefore makes the edge matrix
send the source unit lattice into the target unit lattice.
For arbitrary two finite points move the source and target
separately by rational unitary matrices into the open
sets of points integral above \(\ell\). The Cayley chart
gives both weak approximations. Rational equivariance of
the jet changes its argument and its value by the same
weight representations at the two endpoints.

These changes
cancel exactly in the respective transported lattice coordinates. This
proves the bound at every pair of points. In particular
the common integer
\[
 D_{\rm jet}=D_0^{2f(2n_m)^2}(D_0d)^{b_0r^4}
\]
depends only on the original datum. Indeed \(D_0\), \(d\) and the
rank and degree bounds in its exponents precede the choices of
\(S,\xi,\ell\) and the embedding. Thus the same integer works at
both \(\iota\) and \(\iota\circ\mathrm{conj}\) whenever
\(\ell\nmid q_\xi\). Iterating gives the factor \(D_{\rm jet}^m\) on each
half, on the original outer levels and in their transported
Weyl lattices.

For the complex bounds take the trivial twist.
Fourier inversion of the single indicator on the finite
matrix ring of cardinality \(J\) has \(J\) coefficients of
absolute value \(1/J\) at every complex embedding.
Thus the total coefficient norm of the translated integral
vacua is one. The pointwise conjugate vacuum bound for Hom,
of rank \(nn'\), is
\[
                       \exp(C_K(q_1nn'+1)^2).
\]
Apply Lemma~\ref{lem:upper-cm-transport} at the rank of the
Hom space and the Gaussian Schur calculation of
Section~\ref{sec:upper-gaussian-jet} with \(c=L=1\).

In particular the squared scale in
the integral is \(\Omega^{-1}\).
Write \(\lambda_{\tau,i}\ge0\) for the part sizes in
the edge Schur polynomial. They sum to \(H\). The
conjugated Weyl operator norm is at most the exponential
of
\[
 C_K(q_1nn'+1)^2+f nn'\log 4+
 \frac12\sum_{\tau,i}
    \bigl(\lambda_{\tau,i}+\log(\lambda_{\tau,i}!)
            +\lambda_{\tau,i}\log(\Omega^{-1})\bigr).
\]
In fact a matrix coefficient against a covector of
dual Weyl norm one is bounded by the product of the
singular values with the indicated covector exponents.

The factorial estimate with \(L=1\) in the Gaussian
integral leaves the factor \(4^{f nn'}\). This proof
therefore bounds the operator directly between the conjugated
Weyl norms.
Using \(\operatorname{rd}(E)\le q_1d\) and the explicit
formula for \(\Omega\), we have
\(|\log\Omega|\ll_d 1+\log q_1\).
The sum of the factorial logarithms is at most
\(H\log\max(1,H)\). With
\(q_1,f\le x^{12}\), \(n_m\le x^6\), \(r<3x^5\) and
\(\log b\le x\), the exponent is
\(O_d(x^{48}+b x^{21})\), hence at most \(x^{99}b^2\).
This bound is uniform over the original edge points and
holds at every complex embedding.

At the trivial twist each last-path level contains the
principal congruence subgroup of rational level \(D_0\).
Apply the counting argument of Section~\ref{sec:terminal-classes}
at each rank \(n\le2n_m\), with \(D_0\) as level.
There \(u=2fn\ll x^{18}\) and the maximal-level count has
logarithm \(O(u^3\log x)=O(x^{54}\log x)\).
Refining to the path level contributes at most
\(2fn^2\log D_0=O(x^{25}b)\) in the logarithm.
Thus each class count and each unscaled edge norm have
the common bound \(G=\exp(x^{99}b^2)\) used below.

\smallskip
\subsubsection{Local transport and terminal branching}
\label{sec:completion-terminal-branching}

\sectionopening{In this subsection, we use the bounds of Section~\ref{sec:completion-terminal-jets} in the terminal comparison.}

Let \(t_\xi\) generate the finite ideal of \(\mathfrak m\) and
write \(u_\xi=\overline{t_\xi\alpha_{\rm W}}/(t_\xi\alpha_{\rm W})\).
In the original terminal charts the false point is \(g\)
and the true one is \(\operatorname{diag}(g,u_\xi^m)\).
Write \(u'=u_\xi^{-m}\). If \(v',v\) at these points have
lattice norm at most one, their Weyl-coordinate transports \(w',w\) by
\(g^{-1}\) and \(\operatorname{diag}(g^{-1},u')\) have
integral entries. Use the bilinear matrix \(B_{\rm grid}\)
constructed for the two terminal grid nodes in
Section~\ref{sec:deep-terminal-measure}. Its equivariance for the
common diagonal group and the extra line gives
\[
 (w')^{\mathsf t}B_{\rm grid}w
   =\left(\prod_\tau
       \iota_\tau(u')^{L_\tau}\right)
               \iota\bigl((v')^{\mathsf t}B_{\rm grid}v\bigr),
\]
with the weight embeddings read on the corresponding
finite components. The original pairing is bilinear in
these two already-dual Weyl models, as in the displayed
coordinate product.

The scalar \(u_\xi\) and its inverse become integral after
multiplication by \(q_\xi D_0q_1\). For \(u_\xi\), write
this product as
\[
 (q_\xi D_0/t_\xi)\,
              \overline{t_\xi(q_1\alpha_{\rm W})}\,
                      \alpha_{\rm W}^{-1}.
\]
Here \(q_\xi D_0/t_\xi\) is integral since \(q_\xi D_0\in\mathfrak m\);
the choices of the adelic Hermitian basis give
\(q_1\alpha_{\rm W}\) and \(\alpha_{\rm W}^{-1}\) integral. Thus
the three factors are integral on the finite adeles. For the
inverse conjugate this expression, since
\(\bar u_\xi=u_\xi^{-1}\).
Above \(\ell\nmid q_\xi\), either sign of \(L_\tau\) is
therefore cleared by \((D_0q_1)^{m|L_\tau|}\).
Because the left-hand side of the equivariance identity
is integral, the original pairing is cleared by
\[
                         M_0=\prod_\tau(D_0q_1)^{m|L_\tau|}.
\]
Combining with the two half paths gives \(N_{\rm path}=M_0D_{\rm jet}^{2m}\)
on every summand of the original finite expansion.

The true terminal diagram has at most \(b_0r^4\) boxes by
Section~\ref{sec:completion-terminal-jets}. Its determinant exponent
is the last highest-weight entry; the bound \(\lceil20b\rceil\) for
the doubled parameters and \(N+1\le2n_m\le\lceil20b\rceil\)
give \(|e|\le\lceil20b\rceil\) after removing the rank shift.
Hence the exponents of this pairing satisfy
\(|L_\tau|\le b_0r^4+\lceil20b\rceil\).
The rank and box bounds and \(\log M=O(\log b)\) give
\(\log N_{\rm path}\le x^{99}b^2\).
By \eqref{eq:measured-terminal-period}, the original internal integrals are free class sums.
Reindexing the combined array at the internal nodes cancels
the adjacent kernel characters against the array character.
At the linked terminal node the basis changes are \(u\) and
\(\operatorname{diag}(u,1)\).
Each array entry is a sum of roots of unity. Thus the edge
and branching bounds give integrality of
\(N_{\rm path}\Xi_\xi(\delta_x(v),\eta)\) on integral outer data
at all embeddings above \(\ell\nmid q_\xi\).

At the trivial twist, put \(G=\exp(x^{99}b^2)\) and
\(h_0=b_0r^4\). Use the bounds \(G\) for each class
count and each point-indicator edge operator in the
conjugate Weyl metrics, and \(h_0^{3h_0}\) for the
terminal branching pairing. In the residue
ring of the grid modulus, put
\[
 N_b=|\mathcal O_E/(\pi_M)^{b_0}|,\qquad
 E_0=N_b^{(2n_m)^2},\qquad V_0=N_b^{2n_m}.
\]
There are at most \(E_0^{2m}\) edge arrays. The first
path supplies at most \(E_0^{2m}\) root-of-unity terms
per entry and each of the \(m^2\) turns multiplies that
bound by at most \(V_0\).

Hence at most
\(V_0^{m^2}E_0^{4m}G^{2m+1}\) terms remain after summing
over the \(2m+1\) original class variables. Each term
has absolute value at most
\[
       h_0^{3h_0}G^{2m}\|\sigma v\|_{\rm W}\sqrt{B_\eta},
\]
by the two compositions of \(m\) edge operators,
where \(B_\eta\) bounds the squared norm of every
\(\sigma\)-conjugate value of \(\eta\).
The numerical bound
\[
 V_0^{m^2}E_0^{4m}h_0^{3h_0}G^{4m+1}
                                  \le \exp(x^{200}b^2)
\]
now gives the corresponding bound for
\(|\sigma\Xi_1(\delta_x(v),\eta)|\).
Indeed \(f,q_1\le x^{12}\), \(n_m\le x^6\),
\(r\le3x^5\), \(m\le2\sqrt x\) and \(\log M\ll\log b\le x\)
for large \(x\). Thus \(\log N_b\le x^{13}b\).
With \(h_0=b_0r^4\), each term in the logarithm of
the count bound is \(o(x^{200}b^2)\). The same
estimates give \(\log(M_0D_{\rm jet}^{2m})=o(x^{99}b^2)\)
using the bound on \(|L_\tau|\).
This estimate concerns the same path period as the
finite calculation.

\subsubsection{The dual vector on the terminal lattice}
\label{sec:completion-terminal-dual}

\sectionopening{In this subsection, we construct the dual vector for the terminal pairing of Section~\ref{sec:completion-terminal-branching}.}

Write \(K_\xi=K(\pi_MQ(\xi))\) on the original terminal
definite group \(Y_m\), and let \(K_{\rm std}\) denote its maximal
integral group. Put \(I_\xi=[K_{\rm std}:K_\xi]\), and let
\(d_Y\) be the dimension of the terminal Weyl module.
Reduction at the chosen split places is onto, with
kernel \(K_\xi\). Hence \(I_\xi\) is the order of the
general linear group of the false terminal node over
the residue ring of \(\pi_MQ(\xi)\). In particular this is
the terminal divisor appearing in the last-path period.
The measures in this calculation are the original
quotient measures, with the compact real group and
the integral finite group at the terminal space each of
mass one.

Use the dual module of type obtained by reversing and
negating the terminal parameter. Write
\(B(v,z)=v^{\mathsf t}J_{\rm W}z\) for its normalized
invariant bilinear pairing with the terminal module.
Let \(\delta_x(v)\) be the dual form of level \(K_\xi\)
whose value at \(x\) is \(v\), extended by equivariance
on the class of \(x\) and zero elsewhere. The stabilizer is finite by
definiteness. At an old split place above the odd prime \(p_0\) it is
conjugate into \(1+p_0M_{r+mt}(\mathbb Z_{p_0})\), which is torsion-free:
the matrix logarithm is injective there, and \(\log(u^a)=a\log u\)
for every positive integer \(a\). Thus the stabilizer is trivial and
the equivariant extension has the prescribed value at the representative.
Use its conjugate-type real realization in the period.

Let \(w\) be the tensor of unit highest wedges in the flags
ordered by the terminal parameter, and let \(\rho_\infty\) denote
the compact action on this module in the scalar realization. The
negative parameter orders each flag in reverse. With the sign of the
complement map fixed below, the coefficient of the dual test at \(k\)
is \(B(v_h,\rho_\infty(k)w)\), where \(v_h=\delta_x(v)(h)\).
Indeed this equality holds at the identity by the normalization on
the opposed flags, and at \(k\) by invariance of \(B\). The primal
coefficient at \(h\), with \(z_h=(T_\xi y)(h)\), is
\(\langle\rho_\infty(k)w,z_h\rangle_{\rm W}\). For compact Haar
measure of mass one, Schur orthogonality now gives
\[
 \begin{aligned}
 \int B(v_h,\rho_\infty(k)w)\langle\rho_\infty(k)w,z_h\rangle_{\rm W}\,dk
 &= B\!\left(v_h,\int \rho_\infty(k)w
          \langle\rho_\infty(k)w,z_h\rangle_{\rm W}\,dk\right)\\
 &=d_Y^{-1}B(v_h,z_h).
 \end{aligned}
\]
The averaged operator commutes with the irreducible compact action
and has trace \(\|w\|_{\rm W}^2=1\).
The rational left actions cancel in this expression by
invariance of \(B\). Both forms are fixed by \(K_\xi\),
and the finite mass of the class at \(x\) is \(I_\xi^{-1}\)
because its rational stabilizer is trivial. Thus
\begin{equation}
 I_\xi\int_{[U(Y_m)]} \delta_x(v)(h)(T_\xi y)(h)\,dh
               = d_Y^{-1} B(v,(T_\xi y)(x)).
 \label{eq:original-dual-delta-readout}
\end{equation}
This formula pairs the algebraic form with its dual
bilinearly. The dual real type has already absorbed the
complex conjugation.

In Proposition~\ref{prop:measured-deep-terminal-period}, use
\(y'=\delta_x(v)\). Its right stabilizer contains the
false outer deep subgroup \(K(\pi_M^{b_0}Q(\xi))\).
The true form \(\eta\) has level \(K(\pi_M^{b_0})\).
Combining \eqref{eq:measured-terminal-period} with
\eqref{eq:original-dual-delta-readout} gives
\begin{equation}
 B\bigl(v,((c_\xi/c_1)T_\xi r_0\eta)(x)\bigr)
       =\beta_0\,\Xi_\xi(\delta_x(v),\eta),
 \qquad \beta_0=\frac{d_Y A_*}{c_1}.
 \label{eq:original-dual-delta-last}
\end{equation}
The finite index has canceled and the compact dimension
is retained in \(\beta_0\). In particular \(\beta_0\ne0\)
if \(D>0\) and \(T_1\) is injective on \(V\). In that case
take \(y\in V\setminus\{0\}\), lift it by \(r_0\), and
choose a dual vector separating a nonzero value of \(T_1y\).

\subsubsection{Integral dual coordinates and pointwise bounds}
\label{sec:completion-terminal-dual-coordinates}

\sectionopening{In this subsection, we estimate the coordinates of the dual vector of Section~\ref{sec:completion-terminal-dual}.}

At each infinite place let \(d_1\ge\cdots\ge d_c\) be the column
heights of the polynomial diagram in rank \(n\). The dual diagram
has column heights \(n-d_c,\ldots,n-d_1\). In the exterior tensor
model send each coordinate wedge to the complementary wedge with the
sign determined by the unit volume wedge, then reverse the order
of its \(c\) tensor factors. This is an integral isometry between the
two prescribed tensor products, also for \(c=0\). For an invertible basis matrix \(g\),
the complement identity on each exterior power moves
\(g\) to \((g^{-1})^{\mathsf t}\) and supplies one
factor of \(\det g\) per column. Indeed, for a highest
weight \(h_1\ge\cdots\ge h_n\), the polynomial diagram has
\(c=h_1-h_n\) columns and the determinant exponent is \(h_n\).
For the dual, that exponent is \(-h_1\). Hence the total factor
\((\det g)^{c+h_n}\) equals \(\det((g^{-1})^{\mathsf t})^{-h_1}\),
including when either determinant exponent is negative. At each
place, the complement of an initial segment of length \(j\) in the
primal ordered flag is the initial segment of length \(n-j\) in the
reversed flag. Their ordered coordinate wedges agree up to a sign.
Multiply the tensor complement by the product of these fixed signs,
so that it sends the unit highest tensor to the prescribed unit
highest tensor of the dual. This fixes its sign also when the diagram
is empty, in which case the product is one.

Applying this
identity to invertible bases proves that the tensor
isometry restricts to the Weyl spans, since the
evaluation vectors of those bases span. The Weyl
coordinates are integral on these tensors: the
standard-tableau coefficient matrix has an integral left
inverse, obtained from its unitriangular minor
submatrix. The same holds for the inverse complement
map. Denote this unimodular isometry by \(J\) and
the dual Weyl Gram matrix by \(G_\vee\). The pairing
has matrix \(J_{\rm W}=G_\vee J\), so in the
convention linear in the second argument
\[
        B(v,z)=\langle\overline{Jz},v\rangle_{\rm W},
                  \qquad \|\overline{Jz}\|_{\rm W}=\|z\|_{\rm W}.
\]
For the unit highest tensor \(w\) in the real coordinate wedges,
put \(w^\vee=Jw\). Then \(B(w^\vee,w)=1\), as in
\eqref{eq:original-dual-delta-readout}.

The Fock pairings of the integral dual minor basis form
an integral matrix \(F_\vee\) with
\(F_\vee^{\mathsf t}G_\vee=H_\vee I\), where
\(H_\vee>0\) is the Fock square of the dual highest
minor polynomial. Indeed the standard minor polynomials are the
coordinate functions of the evaluation vector in the dual Weyl
module. The Hermitian metric of these coordinate functionals induced
by the Weyl metric has matrix \(G_\vee^{-1}\) in the dual basis,
since \(G_\vee\) is real symmetric. This metric and the Fock metric
are invariant under the compact action on the irreducible polynomial
module, so they are proportional. The highest coordinate functional
is the inner product with the unit highest tensor, by its torus
weight and its value one on that tensor. Its squared norm is one in
the induced metric and \(H_\vee\) in the Fock metric. This gives
\(F_\vee^{\mathsf t}G_\vee=H_\vee I\). Thus
\[
         C_{\rm W}J_{\rm W}=H_\vee I,
             \qquad C_{\rm W}=J^{-1}F_\vee^{\mathsf t}\in M(\mathbb Z).
\]
There are at most \(h_0=b_0r^4\) boxes in the dual diagram.
The determinant expansion bounds the coefficient
\(\ell^1\)-norm of its highest polynomial by \(h_0^{h_0}\),
and each monomial factorial by \(h_0^{h_0}\).
Consequently \(H_\vee\le h_0^{3h_0}\).

If \(x\) is integral above \(\ell\) and \(v\) is an
integer dual coordinate vector, \(\delta_x(v)\) has lattice norm
at most one. Indeed \(x^{-1}\) is integral at those places,
as it is the Hermitian adjoint of \(x\), and the
polynomial representation and determinant units preserve
the unit lattice. Equivariance gives this lattice norm at every
point of the class. At any coefficient automorphism
\(\sigma\), its pointwise Weyl norm is at most
\(\|\sigma v\|_{\rm W}\). Indeed every embedding \(\lambda:E\to\mathbb C\)
satisfies \(\lambda(\bar a)=\overline{\lambda(a)}\) for \(a\in E\),
since \(E\) is CM. Thus \((\sigma\circ\lambda)(\gamma)\) is unitary
for every rational unitary transport \(\gamma\), and its Weyl action
preserves the corresponding rational Gram form. This argument applies
to the transport entries before pairing with the arbitrary vector
\(\sigma v\).

Let \(N_{\rm path}=M_0D_{\rm jet}^{2m}\), the integer fixed above
for every allowed twist and every coefficient prime off its modulus.
Set \(s_0=H_\vee N_{\rm path}/\beta_0\) when \(\beta_0\ne0\)
and \(s_0=1\) otherwise. To see the cancellation in the first case,
put \(F=(c_\xi/c_1)T_\xi(s_0 r_0\eta)\) and take \(v_i\) to be
the transpose of row \(i\) of \(C_{\rm W}\). At an \(x\) integral
above \(\ell\), use \(C_{\rm W}J_{\rm W}=H_\vee I\) in
\eqref{eq:original-dual-delta-last}. This gives the equality
\[
 H_\vee F_i(x)=B(v_i,F(x))=
 H_\vee N_{\rm path}\,\Xi_\xi(\delta_x(v_i),\eta).
\]
The right-hand period multiplied by \(N_{\rm path}\) is integral.
Cancellation in characteristic zero gives the integral coordinate
\(F_i(x)=N_{\rm path}\Xi_\xi(\delta_x(v_i),\eta)\), also when
\(\ell\mid H_\vee\). For \(\beta_0=0\), nondegeneracy in
\eqref{eq:original-dual-delta-last} gives \(F(x)=0\).
In both cases \(s_0\ne0\) is independent of \(\xi\) and \(\ell\).

For arbitrary \(x\), use the Cayley-chart weak approximation of
Section~\ref{sec:original-comparison-fock} to choose
\(\gamma\in U(Y_m)(F_0)\) with
\(\gamma_\ell x_\ell\in K_{{\rm std},\ell}\), where the subscript
\(\ell\) denotes all local components above \(\ell\).
Apply the integral estimate at \(\gamma x\). Writing
\(F=(c_\xi/c_1)T_\xi(s_0r_0\eta)\) and \(\rho_{\rm W}\) for the
terminal Weyl action at each coefficient embedding above \(\ell\),
rational equivariance gives
\[
 \rho_{\rm W}\bigl((\gamma_\ell x_\ell)^{-1}\bigr)F(\gamma x)
     =\rho_{\rm W}(x_\ell^{-1})F(x).
\]
Thus the transported coordinates are integral at the original point.

At the empty twist, if the path bound is
\[
 |\sigma\Xi_1(\delta_x(v),\eta)|
   \le e^{x^{200}b^2}\|\sigma v\|_{\rm W}\sqrt{B_\eta},
\]
then the dual isometry and \eqref{eq:original-dual-delta-last}
give the squared norm bound
\[
     \bigl(H_\vee N_{\rm path}e^{x^{200}b^2}\bigr)^2 B_\eta
                   \le e^{2x^{201}b^2}B_\eta .
\]
The last bound uses \(\log N_{\rm path}\le x^{99}b^2\) and
the estimate for \(H_\vee\) just proved. Indeed the period identity for
the normalized vector has scalar \(s_0\beta_0=H_\vee N_{\rm path}\)
when \(\beta_0\ne0\), hence is fixed by every \(\sigma\).
The identity \(\sigma B(v,F(x))=B(\sigma v,\sigma F(x))\) follows from
the integral pairing matrix. Taking \(\sigma v=\overline{J(\sigma F(x))}\)
in the complex dual space gives the stated norm bound, with the zero
value immediate.
If \(\beta_0=0\), all of these terminal values vanish
by \eqref{eq:original-dual-delta-last}.

\smallskip
\subsubsection{Integral generators at the outer node}
\label{sec:completion-outer-generators}

\sectionopening{In this subsection, we use the dual coordinates of Section~\ref{sec:completion-terminal-dual-coordinates} to bound the original generators.}

At the true outer node \((1,0,m)\) the rank is
\(n_H=r+1+m(t-1)\le 2n_m\) and the level is \(K(\mathfrak b)\).
For each free class \(g\) at that level write
\(L(g)=E^{n_H}\cap g\widehat{\mathcal O}_E^{n_H}\).
The short-basis calculation at the source supplies
\(Y(g)\in\mathrm{GL}_{n_H}(E)\) whose columns lie in \(L(g)\)
and whose entries have absolute value at most
\(R=6(q_1d)(fn_H)^2\) at every complex embedding.
Here the columns form an \(E\)-basis, and may span a
proper \(\mathcal O_E\)-sublattice of \(L(g)\).
Only invertibility over \(E\) and
\(g^{-1}Y(g)\in M_{n_H}(\widehat{\mathcal O}_E)\)
are used. For each actual determinant exponent \(e_\tau\)
apply the same construction at \(g^{e_\tau}\).

At the representative \(g\), multiply the polynomial
Weyl action of \((\iota_\tau Y(g))_\tau\) on the
evaluation vector by
\[
                   \prod_\tau\det\iota_\tau Y(g^{e_\tau}).
\]
After moving into the local lattice norm at \(g\), the determinant
at any coefficient embedding is evaluated on
\[
  \det(g)^{-e_\tau}\det Y(g^{e_\tau})
     =\det\big((g^{e_\tau})^{-1}Y(g^{e_\tau})\big).
\]
The identity also holds for negative \(e_\tau\).
Its right-hand side and the polynomial Weyl action of
\(g^{-1}Y(g)\) are integral. For root-of-unity evaluation
arrays all initial coordinates are integral at every
coefficient embedding.

The rational left action and right level action show that
this is the lattice norm at every class point. Writing such a
point as \(\gamma gk\), with \(k\) in the deep compact
group and \(\gamma\) rational, equivariance cancels the
weight factor of \(\gamma\) against transport by
\((\gamma gk)^{-1}\). The remaining \(k^{-1}\) acts
integrally, with inverse also integral. The free class
is supported at a single class and thus incurs no stabilizer
sum. A corresponding conjugate value has the same norm
as at \(g\), because the rational definite unitary action
remains unitary in the simultaneous conjugated Weyl
coordinates, and the class level has trivial character.

The polynomial Weyl actions of the \(Y(g)\) are invertible,
and all the determinant factors are nonzero, so they preserve
the span of the evaluation vectors. Apply the exterior-power bound in
Lemma~\ref{lem:ramified-short-generators} with rank \(n_H\).
Each column of an evaluation matrix with root-of-unity
entries has norm at most \(\sqrt{n_H}\), and the operator norm of
\(Y(g)\) is at most \(n_H R\). Each exterior tensor of total
degree at most \(b_0r^4\) is bounded by the corresponding power of
\(n_H^{3/2}R\). Each determinant factor is at most
\((\sqrt{n_H}R)^{n_H}\) by Hadamard's inequality. Since
\(R\ge n_H^2\), the squared norm is bounded by
\[
  (n_H^{3/2}R)^{2b_0r^4}(\sqrt{n_H}R)^{2fn_H}
  \le R^{6(fn_H+b_0r^4)}\le \exp(x^{99}b^2).
\] Transport to the analytic form space is
by the same coefficient-type isomorphism at the actual
deep level.

For each free outer class these values span its whole fiber.
In fact put \(H_0=b_0r^4\) and choose any integer \(l>H_0\).
Take one independent matrix of variables for each infinite
place and let every entry range over \(\mu_l\). At this fixed
class the factors \(\det\iota_\tau Y(g^{e_\tau})\) are nonzero
constants for every sign of \(e_\tau\). The variable part consists
of products of minors with total degree at most \(H_0\). A
functional vanishing on the indicated Cartesian grid therefore gives
a polynomial \(P(Y)\) of degree less than \(l\) in each variable.
Writing \(k_0\) for the number of variables, for any multi-index
\(0\le a_i<l\) orthogonality of roots of unity gives
\[
 [Y^a]P=l^{-k_0}\sum_{Y\in\mu_l^{k_0}}P(Y)\prod_{i=1}^{k_0}Y_i^{-a_i}=0 .
\]
This coefficient extraction takes place in \(\mathbb C\).
The integral generators themselves use only polynomial evaluations,
including at singular arrays. The evaluation vectors span the Weyl space
by its minor construction, and the fixed invertible transforms
preserve their span. Taking the union of the generators supported on individual
classes gives a spanning set of \(\mathcal H\); since \(r_0\) is
onto and \(s_0\ne0\), the \(s_0r_0\eta\) span the full source \(W\).

\subsubsection{Projection to the projected space}
\label{sec:completion-good-projection}

\sectionopening{In this subsection, we project the spanning set of Section~\ref{sec:completion-outer-generators} to the projected space, retaining the bounds needed in Section~\ref{sec:marked-determinant-assembly}.}

Use the integral spanning set of \(\mathcal H\) just constructed.
Its forms have integral values in the transported lattices and squared norms at most
\(e^{x^{99}b^2}\) at every complex embedding. Let \(z_\eta=s_0r_0\eta\).
The \(z_\eta\) span \(W\). Section~\ref{sec:completion-terminal-dual-coordinates}
proves that \(\gamma_\xi T_\xi z_\eta\), with \(\gamma_\xi=c_\xi/c_1\),
is integral at every embedding over \(\ell\nmid q_\xi\), and that each
conjugated value of \(T_1z_\eta\) has squared norm at most
\(e^{x^{202}b^2}\). Set
\[
  \mathcal Z=\{\,|\mathcal G|p z_\eta:\eta\in\mathcal S_H\,\}\subseteq V ,
\]
where \(\mathcal G=G_0\times G_1=G\) is the original finite group
and \(p=\rho(e_0)\) projects the principal-level space \(W\) onto
\(V\) as in Section~\ref{sec:notation-automorphic}.

The operator \(|\mathcal G|p:W\to V\) is surjective and is a sum of integral
right translations with coefficients algebraic
integral at every coefficient embedding. Indeed its coefficient at
a group element is the sum of dimension times the
character value of the retained simple types. Its
conjugate absolute value is at most \(|\mathcal G|\),
by the bound of character values by dimension and the
sum-of-squares identity. Apply Lemma~\ref{lem:integral-finite-type-projection}
with the identity translate. Its first-block lifts \(k_m(t)\) give
\[
 (\gamma_\xi T_\xi(|\mathcal G|p z_\eta))(h)
   =\sum_{t\in\mathcal G}|\mathcal G|e_0(t)
                         (\gamma_\xi T_\xi z_\eta)(h k_m(t)).
\]
Each lift and its inverse are integral above \(\ell\); transporting
from \(h k_m(t)\) to \(h\) preserves the unit lattice. Thus the sum
is integral at every coefficient embedding above \(\ell\nmid q_\xi\).
For \(\xi=\xi_\kappa\), apply this at \(\iota\) and
\(\iota\circ\mathrm{conj}\) to the same vector
\(\gamma_\xi T_\xi z\), for every \(z\in\mathcal Z\).
The latter application gives the conjugate coordinates used in the
mixed Gram pairing at \(\kappa\). For any \(\sigma\), the sum
of the absolute values of the conjugated coefficients is bounded by
\(|\mathcal G|^2\); the triangle inequality contributes at most
\(|\mathcal G|^4\) to the squared pointwise bound. The class
count at the original terminal level from Section~\ref{sec:terminal-classes}
and \(|\mathcal G|^4\) multiply the squared
pointwise bound by at most \(e^{x^{99}}\).
The total conjugate squared norm of each projected
form is therefore at most \(e^{2x_1^2b^2}\), by
\(x^{99}+x^{202}b^2\le2x_1^2b^2\).

For \(\xi=\xi_\kappa\) put \(p_{ij,k}=(\mu^0_{i+1}\mu^0_{j+1},k)\)
and
\[
 U_{ij,k}
   =\frac{F_{p_{ij,k},Q(\xi)}(\xi^{j+1})}{F_{p_{ij,k},\varnothing}(1)}.
\]
The support of \(\xi^{j+1}\) is \(Q(\xi)\) for \(0\le j<m\).
The two norm formulas of Section~\ref{sec:critical-signed-systems}
give \(q_\xi^{-1}U_{ij,1}\) for the square of the first factor
of \(T_{ij}(\xi)\) and \(U_{ij,k}\) for each factor \(k\ge2\).
Since \(\mathrm{Nm}(Q(\xi))=q_\xi\), equation~\eqref{eq:display-tp11}
gives
\[
 \gamma_\xi\bar\gamma_\xi
  =q_\xi^{m^2t^2+m^2t+mr-m^2}
       \prod_{0\le i,j<m}\left(U_{ij,1}\prod_{k=2}^t U_{ij,k}^2\right)
  \in\overline{\mathbb Q}\cap\mathbb R_{>0}.
\]
Consequently its valuations at \(\kappa\) and \(\bar\kappa\)
agree for this twist. As \(\ell\nmid q_\xi\),
\eqref{eq:display-rt1} gives
\[
 v_\kappa\bigl((\gamma_\xi\bar\gamma_\xi)^{-1}\bigr)
   =-2\sum_{0\le i,j<m}v_\kappa T_{ij}(\xi)\ge -2m^2t e_\ell .
\]
Hence \eqref{eq:display-rt3} with multiplier \(1\) gives
\[
 \operatorname{Height}\bigl(
      ((\gamma_{\xi_\kappa}\bar\gamma_{\xi_\kappa})^{-1})_\kappa\bigr)
                    \le 2m^2t x_1<b .
\]
Indeed \(m^2\le4x\), \(t<15x^5\) and
\(\log(2m^2t x_1)\le\log120+6\log x+x^{\epsilon/6}
<x^\epsilon<\log b\) for sufficiently large \(x\).

\begin{proposition}[Spanning vectors for the lattice completion]
\label{prop:marked-completion-vectors}
Choose \((\xi_\kappa)_\kappa\) as in Section~\ref{sec:common-critical-twist}.
There is a set \(\mathcal Z\subset V\) with
\[
 \operatorname{span}\mathcal Z=V,\qquad
 \|\sigma\Theta_1 z\|^2\le e^{2x_1^2b^2}
      \quad(z\in\mathcal Z,\ \sigma\in\operatorname{Aut}(\mathbb C)).
\]
For each \(\kappa\), the scalar \(\gamma_{\xi_\kappa}=c_{\xi_\kappa}/c_1\)
is nonzero and \(\gamma_{\xi_\kappa}\Theta_{\xi_\kappa}z\) is integral at
both \(\kappa\) and \(\bar\kappa\) for all \(z\in\mathcal Z\). Moreover,
\[
 \operatorname{Height}\bigl(
           ((\gamma_{\xi_\kappa}\bar\gamma_{\xi_\kappa})^{-1})_\kappa
                    \bigr)\le b .
\]
\end{proposition}

\subsection{Integral derivatives in the two norm formulas}
\label{sec:marked-both-norm-jets}

\sectionopening{In this subsection, we prove the integrality used by the two norm constructions in
Sections~\ref{sec:marked-norm-line}
and~\ref{sec:marked-norm-determinant}. Fourier inversion first reduces
the point indicators to translated vacuum vectors. At coefficient
primes outside the auxiliary support, the same denominator clears each
expansion.}

Fix a coefficient prime \(\ell\nmid q_Q\) and any embedding above it.
For a finite modulus \(\mathfrak m\) with
\(q_QD\in\mathfrak m\), Fourier inversion on the finite matrix
coordinates expresses its integral indicator datum as a combination
of translates, with roots of unity in the numerator. The Fourier
denominator divides
\[
   D^{2f n'n}
\]
in the chosen valuation, since \(\mathcal O_E/\mathfrak m\)
at \(\ell\) is a quotient of \(\mathcal O_E/(D)\) at \(\ell\)
and \([E:\mathbb Q]=2f\). The Fourier dual translations are
integral after multiplication by \(D\) above \(\ell\), by the
self-dual translation pairing. Thus multiplying by \(D^{2fn'n}\)
gives the standard datum of level \(\delta=D\), once the
vacuum has been normalized by \(s_{\rm vac}^{n'n}\). This
argument uses coprimality with \(q_Q\) alone. For finite unitary
points integral invertible above \(\ell\), the lattice index
in Section~\ref{sec:upper-cm-jet-transport} is prime to \(\ell\).
For each norm-one determinant \(u\) of these points, write
\(u^h=a k\), with \(a\) a rational norm-one scalar and \(k\) an
integral norm-one idele, as in that section. For every endpoint
character \(\alpha\) of integral infinity weight,
\[
 \alpha_{\rm f}(u)^h=\alpha_{\rm f}(k)\alpha_\infty(a)^{-1}.
\]
The first factor is a root of unity. Every conjugate of \(a\) is
a unit above \(\ell\), so the infinity weight monomial in the second
factor is an \(\ell\)-adic unit as well. Applied to the splitting
character, the square identity for \(\Delta_g\) in that section now
gives \(v_\ell(\Delta_g)=0\). The Heisenberg phases are roots of unity.
The finite vacuum formula therefore preserves the integral span of
translates with \(Dx\) integral above \(\ell\), also at \(\ell=2\).

Put \(\mathfrak q=\prod_{w\in Q}\mathfrak p_w\), so that
\(q_Q\in\mathfrak q\). In the rank-one test the old column
is prescribed modulo \(\pi_M\). At each new place only its
distinguished entry is prescribed modulo \(\mathfrak p_w\).
Thus this indicator factors through the column space modulo
\(\mathfrak m_1=(\pi_M)\mathfrak q\), including when written on the
larger modulus \((\pi_M^2)\mathfrak q\). On the smaller space it is a sum of distinct
cell indicators with coefficient one. Since \(M=\pi_M\bar\pi_M\), we have
\(q_QM\in\mathfrak m_1\). Apply Fourier inversion at
\(\mathfrak m_1\) with \(D=M\) and \((n',n)=(a_0,1)\).
Since \(\operatorname{Nm}_{E/\mathbb Q}(\pi_M)=M^f\), the full
column space there has cardinality \((M^f q_Q)^{a_0}\).
The column polynomial has integral
coefficients and degree \(u=\sum_v(c_v-z)\). Using the ordinary jet
bound at the integral scale \(Md_Ks_{\rm jet}\) and homogeneity
at the actual scale \(\pi_Ms_{\rm jet}\) shows that the Fourier-cleared
column kernel has lattice norm at most
\(
 |\pi_M/(Md_K)|_\ell^u
\)
at its actual target point. Here the column kernel already includes
the vacuum scalar \(s_{\rm vac}^{a_0}\). Thus
\[
 M^{2fa_0}(Md_K)^u K_Q(g,h)
\]
has lattice norm at most one in the Weyl coordinates transported by \(h\).
The source vector is the constant \(1\) in the trivial
source weight. Rational equivariance and weak approximation
justify the lattice norm at arbitrary \(g,h\). In particular the bound
applies to both coefficient embeddings in the Gram pairing and to the old
level primes and \(2\).

The rank-two determinant point test prescribes the identity
matrix modulo \(\mathfrak m_2=(\pi_M^2)\mathfrak q\).
Its matrix space has cardinality \((M^{2f}q_Q)^4\).
Here \(q_QM^2\in\mathfrak m_2\), so take \(D=M^2\),
\((n',n)=(2,2)\), and \(Q(Z)=\prod_v\det Z_v^{k_v}\),
where \(k_v=c_v-1\). This polynomial has integer coefficients
and degree \(l=2\sum_v k_v\); the datum is the indicator
of the identity residue matrix. The point kernel uses the
unscaled vacuum and the unscaled homogeneous jet. The
normalized kernel is exactly the jet of the scalar-cleared
datum at the scale \(c=M^2d_Ks_{\rm jet}\).
By homogeneity,
\[
 s_p K^{\rm pt}_{Q,1}(h)=Q(\partial)\Theta_{F_h}(cz)\big|_{z=0},
 \qquad s_p=(M^2)^{8f}s_{\rm vac}^4c^l .
\]

Here \(F_h\) is the paired Weil translate at \((1,h)\)
of the identity-residue datum multiplied by
\((M^2)^{8f}s_{\rm vac}^{4}\). Choose \(\gamma\in H(F_0)\)
by weak approximation so that \(\gamma h\) is integral invertible
at all places above \(\ell\). By rational invariance of the adelic
theta functional, left multiplication by \(\gamma\) in the finite
argument acts on the real test by the inverse real action.
At \(\tau=1\) the endpoint character is
\(\alpha_1=\nu\upsilon|_{U(1)}\) of weights \(-k_v\), so its
factor at \(\gamma_\infty^{-1}\) is
\(A_\gamma=\prod_v\Phi_v(\det\gamma)^{k_v}\).
The bare determinant polynomial contributes \(A_\gamma^{-1}\)
by the real action in Section~\ref{sec:marked-norm-determinant}.
Thus \(K^{\rm pt}_{Q,1}(\gamma h)=K^{\rm pt}_{Q,1}(h)\).
The datum \(F_{\gamma h}\) is standard of level \(\delta=M^2\).
The ordinary jet estimate at scale \(M^2d_Ks_{\rm jet}\) gives
\(s_pK^{\rm pt}_{Q,1}(h)\in\mathcal O\).
The multiplier \(s_p\) is independent of the point representative
and of \(Q\).
The statement applies before the squared-modulus coefficient
system is formed and supplies integrality at both embeddings.

\section{Proof of the determinant estimate}
\label{sec:marked-determinant-assembly}

\sectionopening{In this section, we combine the determinant-ratio height of
Proposition~\ref{prop:chain-determinant-height} with the spanning
vectors of Proposition~\ref{prop:marked-completion-vectors}. At each coefficient
place, the same auxiliary twist works in both factors of the terminal
Hermitian pairing. This permits the use of the product formula in
Section~\ref{sec:gram-product-formula} and proves
Proposition~\ref{prop:marked-determinant}.}

\begin{proof}[Proof of Proposition~\ref{prop:marked-determinant}]
The proof applies to the original datum of
Section~\ref{sec:common-filling-assembly} at every sufficiently
large integer \(x\) for fixed \(q_0,d,\epsilon\).
Use the weighted chain \(T_\xi=L_\xi\) of
Section~\ref{sec:seesaw-normalization} and restrict it to \(V\).
The injectivity argument in
Section~\ref{sec:marked-chain-composition} uses the
first-arrow return and \eqref{eq:chain-gram-first-comparison} before division by the real square
factor. The single lifts are injective by
Section~\ref{sec:peak-period-injectivity} and Proposition~\ref{prop:common-period-adjoint}.

On this datum the real scalar in
\eqref{eq:seesaw-measured-chain} is nonzero. The
unscaled diagonal moment of
Proposition~\ref{prop:orbital-diagonal-moment} and the source
codimension bound of Proposition~\ref{prop:ramified-source-codimension}
give
\[
 \|pB P_{C(S)}\|_{\rm HS}^2\ge(1-b^{-8}-5b^{-4})d_X>0 .
\]
Since \(B=(v_G\nu)^{-1/2}\Pi C^{-1}\), some \(\phi\in S\)
has \(p\Pi\phi\ne0\) in \(V\). If \(c_{\rm meas}=0\),
\eqref{eq:seesaw-measured-chain} at the trivial twist gives
\(L_1\Pi\phi=0\). Lemma~\ref{lem:weighted-lower-equivariance}
then gives \(L_1(p\Pi\phi)=0\), contrary to injectivity of the
weighted chain on \(V\). Therefore
\begin{equation}\label{eq:measured-seesaw-scalar-nonzero}
                              c_{\rm meas}\ne0 .
\end{equation}
The diagonal moment and the chain injectivity here concern
the unscaled maps \(B\) and \(L_1\).

For each \(\kappa\) use the twist \(\xi_\kappa\) of
Section~\ref{sec:common-critical-twist}.
Its modulus is prime to the rational prime under \(\kappa\).
It supplies the same signed critical ratios and main Selmer
tests for all the original lift indices. The compression in
Section~\ref{sec:original-compression-bounds} therefore uses
the regular twist and its powers on the original projected subspaces.
The finite spherical trace identity
\eqref{eq:marked-full-spherical-trace} identifies the corresponding
Rallis determinants. Equation~\eqref{eq:marked-chain-height-bound}
gives \(j_Z\le D(b+L-A)\) on this family.

Proposition~\ref{prop:marked-completion-vectors} supplies a spanning
set with
\[
 b_{\rm short}=2x_1^2b^2,\qquad j_\gamma\le b .
\]
At each \(\kappa\), the one scalar \(\gamma_{\xi_\kappa}\) clears these
forms at both \(\kappa\) and \(\bar\kappa\). The height multipliers
are independent of the coefficient place.

Let \(n=\#y\) for an independent tuple as in the proposition.
The images of \(\mathcal Z\) span \(V/\operatorname{span}(y)\);
choose representatives \(z_1,\ldots,z_s\in\mathcal Z\) of a basis
of this quotient. Then \(s=D-n\le4D/b^2\) and \((y,z)\) is a basis
chosen before \(\kappa\). For fixed \(\kappa\), the index hypothesis
allows a clearing matrix for \(y\) at \(\kappa\) and one at
\(\bar\kappa\), both valid for the same \(\xi_\kappa\) since
its modulus is prime to their common residue characteristic. For a tuple of
translated images, Proposition~\ref{prop:translated-norm-and-index} supplies
both clearing matrices for the fixed tuple, each valid for every allowed
twist of level prime to that residue characteristic. The free terminal
classes give integral pairings at \(\kappa\), as proved in
Section~\ref{sec:gram-product-formula}.
Proposition~\ref{prop:mixed-gram-inequality} gives
\[
 \begin{aligned}
 0 &\le c+2\delta+s(2x_1^2b^2+b)+D(b+L-A)\\
   &\le c+2\delta+D(2b+L-A)\\
   &\le c+2\delta+D(3L-A).
 \end{aligned}
\]
\par
\Needspace{9\baselineskip}
Here \(L=faB^4\ge b\), and
\[
 s(2x_1^2b^2+b)\le D(8x_1^2+4/b)\le Db
\]
by \(s\le4D/b^2\) and \(x_1^{10^6}\le b\) for sufficiently large \(x\).
This gives \(D(A-3L)\le c+2\delta\) with a cutoff depending only
on the fixed \(q_0,d,\epsilon\).
\end{proof}

\part{The analytic conclusion}\label{part:analytic}

\section{The analytic criterion}\label{sec:criterion}

\sectionopening{In this section, we give an analytic criterion in terms of a one-sided prime sum. We
choose the auxiliary quadratic twist so that the resulting Dirichlet
\(L\)-functions have no zeros in a prescribed rectangle. Section~\ref{sec:prime-remainder} derives the estimate
used in the Mellin transform argument from
Proposition~\ref{prop:cancel}.}

\subsection{The one-sided Mellin criterion}
\label{sec:one-sided-mellin-argument}

\sectionopening{In this subsection, we use the prime sums of Section~\ref{sec:notation-prime-sums} to rule out zeros by positivity.}

Fix \(8\mid q_0\), a reduced class \(d_0\equiv7\bmod 8\),
and \(1/2<\beta<1\). The required bound is
\(R_{d,q_0}(x)\le x^{\beta-1}\) for all sufficiently large real
\(x\), for each prime \(d\equiv d_0\bmod q_0\). The threshold
may depend on \(d\). Here, with \(\varphi=\varphi(q_0)\),
\[
 \lambda(n)=\mathbf1_{n\equiv1\bmod q_0}\Lambda(n)(1-(n/d)),
 \qquad
 R_{d,q_0}(x)=\sum_{n\le x}\lambda(n)(n^{-1}-x^{-1})
                 -\frac{\log x}{\varphi}-c .
\]
For a character \(\chi\), let \(F_\chi\) be \(L(s,\chi)\)
with its pole removed by multiplication by \(s-1\) when
\(\chi\) is principal.

Dirichlet nonvanishing gives
\(F_\chi(1)\ne0\) for nonprincipal \(\chi\)
\cite[Theorem~4.9]{MV}. In the principal case,
\(F_\chi(1)=\prod_{p\mid\operatorname{mod}(\chi)}(1-p^{-1})\).
Set
\[
 D(s)=\sum_{\chi\bmod q_0}
 \left(\frac{F'_{\chi(\cdot/d)}(s)}{F_{\chi(\cdot/d)}(s)}
       -\frac{F'_\chi(s)}{F_\chi(s)}\right),
 \qquad c=\frac{D(1)-1}{\varphi}.
\]
The same definitions with modulus \(q=1\) give
\(R_{d,1}\) and \(c_{d,1}\). The character set is closed
under conjugation, so \(D(1)\) is real. Since \((d,q_0)=1\), the Chinese remainder theorem identifies the
restriction of each twist to \((\mathbb Z/d\mathbb Z)^\times\)
with the nontrivial quadratic character. Thus the untwisted principal
character contributes the only pole of an unregularized factor at \(1\).
\begin{lemma}[One-sided remainder criterion]\label{lem:positive-mellin-criterion}
Fix \(1/2<\beta<1\), \(H_0>0\), and a prime \(d\) in the prescribed
class. Suppose the twisted factors have no zeros in
\(\Re s\ge\beta,\ |\Im s|\le H_0\) and
\(R_{d,q_0}(x)\le x^{\beta-1}\) for all sufficiently large real \(x\).
Then the untwisted factors \(F_\chi\) have no zeros at
\[
                  \Re s>\beta,\qquad |\Im s|\le H_0 .
\]
\end{lemma}

\begin{proof}

For integers \(N\ge2\), the divisor identity
\(\sum_{n\mid m}\Lambda(n)=\log m\) gives
\[
 \frac N2\sum_{n\le N}\frac{\Lambda(n)}{n}
 \le\sum_{n\le N}\Lambda(n)\lfloor N/n\rfloor=\log(N!).
\]
Since \(0\le\lambda(n)\le2\Lambda(n)\), it follows that
\(R_{d,q_0}(x)=O_{d,q_0}(1+\log x)\).
For \(\sigma=\Re s>1\), character orthogonality and the Euler products give
\[
 \sum_n\lambda(n)n^{-s}=\frac{D(s)+(s-1)^{-1}}{\varphi},
 \qquad
 \int_n^\infty(n^{-1}-x^{-1})x^{-s}\,dx
       =\frac{n^{-s}}{s(s-1)} .
\]
At powers of \(d\) the twisted character vanishes, as required by
\(1-(n/d)\). The kernel is nonnegative for real \(x\ge n\).
The second identity at \(s=\sigma\) and
\(\sum_n\lambda(n)n^{-\sigma}<\infty\) justify the exchange of sum and integral.

The transforms of \(\log x\) and \(1\) are
\((s-1)^{-2}\) and \((s-1)^{-1}\), respectively. Therefore
\begin{equation}\label{eq:remainder-mellin-transform}
 \mathcal M_R(s)=\int_1^\infty R_{d,q_0}(x)x^{-s}\,dx
       =\frac{D(s)-1-s(D(1)-1)}{\varphi s(s-1)}
                   \qquad(\Re s>1).
\end{equation}
The right side is meromorphic on \(\Re s>0\). Its value at the
removable singularity \(s=1\) is
\((D'(1)-D(1)+1)/\varphi\). Each logarithmic derivative has at most a
simple pole at a zero. Thus for \(\Re\rho>0\), \(\rho\ne1\),
\[
 \operatorname{Res}_{s=\rho}\mathcal M_R(s)
 =\frac{\sum_{\chi\bmod q_0}
   (\operatorname{ord}_\rho F_\chi-
    \operatorname{ord}_\rho F_{\chi(\cdot/d)})}
        {\varphi\rho(1-\rho)} .
\]
An Euler factor omitted when passing from a primitive character to its
induced character has the form \(1-\alpha p^{-s}\), \(|\alpha|=1\),
and is nonzero in this half-plane. In the prescribed rectangle every
twisted order vanishes, including on \(|\Im\rho|=H_0\). The numerator
is then the total untwisted multiplicity. At a real
\(\beta<\rho<1\) it gives a strictly positive residue whenever there
is a zero. The Euler products are nonzero on \(\Re s>1\).

For \(x\ge X_0\) sufficiently large put
\(h(x)=x^{\beta-1}-R_{d,q_0}(x)\ge0\).
For \(\Re s>1\), its Mellin transform is
\begin{equation}\label{eq:positive-mellin-tail}
 \int_{X_0}^{\infty}h(x)x^{-s}\,dx
   =\frac{1}{s-\beta}-\mathcal M_R(s)
     -\int_1^{X_0}\bigl(x^{\beta-1}-R_{d,q_0}(x)\bigr)x^{-s}\,dx .
\end{equation}
The last integral is entire because its interval is compact and
\(R_{d,q_0}\) is locally bounded. Denote the tail transform by \(H\)
and its real abscissa of convergence by \(a\in[-\infty,1]\).
Positivity gives absolute convergence and holomorphy on \(\Re s>a\).
By the meromorphic identity theorem, \eqref{eq:positive-mellin-tail}
holds on \(\Re s>\max(a,\beta)\).

Suppose \(a>\beta\). Landau's argument
\cite[Lemma~15.1]{MV} gives a singularity at \(a\). Indeed, if \(H\)
were holomorphic across \(a\), the Taylor series about \(b=a+\delta\)
would reach \(b-3\delta/2\) for sufficiently small \(\delta>0\).
Differentiation is justified by convergence at \(a+\delta/2\).
Since \(X_0\ge1\), Tonelli's theorem gives
\[
 \int_{X_0}^{\infty}h(x)x^{-a+\delta/2}\,dx
 =\sum_{k\ge0}\frac{(3\delta/2)^k}{k!}
       \int_{X_0}^{\infty}h(x)x^{-b}(\log x)^k\,dx
 =\sum_{k\ge0}\frac{(-3\delta/2)^k}{k!}H^{(k)}(b)<\infty ,
\]
contradicting the definition of \(a\). At \(a=1\) the right side
of \eqref{eq:positive-mellin-tail} is regular. If \(\beta<a<1\),
its only possible singularity is a simple pole of strictly negative
residue. Such a pole gives \(H(\sigma)\to-\infty\) as
\(\sigma\downarrow a\), whereas \(H(\sigma)\ge0\).
Thus \(a\le\beta\). On the tail we have
\(|R_{d,q_0}(x)|\le x^{\beta-1}+h(x)\), so \(\mathcal M_R\)
itself converges absolutely and is holomorphic on \(\Re s>\beta\).
Its meromorphic formula agrees with this integral throughout that
half-plane. A zero of an untwisted factor there at height
\(|\Im\rho|\le H_0\) would give a nonzero residue there.
\end{proof}

\subsection{Choosing a twist with no zeros in the rectangle}
\label{sec:twist-zero-free-rectangle}

\sectionopening{In this subsection, we choose the prime needed for the Mellin criterion of Section~\ref{sec:one-sided-mellin-argument}.}

\begin{proposition}[A quadratic twist avoiding a fixed box]\label{prop:zero-free-auxiliary-twist}
Let \(8\mid q_0\), \((d_0,q_0)=1\), and \(d_0\equiv7\pmod8\).
Fix \(1/2<\beta<1\) and \(H_0>0\).
There is a prime \(d\equiv d_0\bmod q_0\) such that
\[
 L(s,\chi(\cdot/d))\ne0
 \quad(\chi\bmod q_0,\ \Re s\ge\beta,\ |\Im s|\le H_0).
\]
\end{proposition}

\begin{proof}

For simultaneous twist selection set
\[
 \mathcal P=\{d\text{ prime}:d\equiv d_0\pmod{q_0}\}.
\]
Fix a character \(\xi\bmod q_0\), and let \(\xi^*\) be its
primitive inducing character. Since \(d_0\) is a reduced class,
\(d\nmid q_0\) for \(d\in\mathcal P\). The two primitive characters in
\(\chi=\xi^*(\cdot/d)\) have coprime conductors, so their product
is primitive with conductor \(Q=\operatorname{cond}(\xi^*)d>1\)
\cite[Lemma~9.3]{MV}.
The Euler factors removed from \(L(s,\chi)\) in passing to the
twist modulo \(q_0d\) are nonzero on \(\Re s>0\).
The Euler product gives nonvanishing for \(\Re s>1\),
so it suffices to prove the assertion on the closed
rectangle \(\beta\le\Re s\le1\), \(|\Im s|\le H_0\).
For \(d\le X\) put \(T_{\rm an}=(q_0X)^{1/2}\ge\sqrt Q\), and let
\(j\in\{0,1\}\) be the parity of \(\chi\). Since \(d\equiv3\pmod4\),
\((-1)^j=-\xi^*(-1)\), which depends only on \(\xi\). For
\(a=\chi,\mu\chi\) set
\[
 \theta_a(y)=y^{1/2}\sum_{n\ge1}a(n)(ny)^j e^{-\pi(ny)^2},
 \qquad
 J_d=\int_{1/T_{\rm an}}^\infty\bigl(|\theta_\chi(y)|+
                         |\theta_{\mu\chi}(y)|\bigr)\frac{dy}{y}.
\]
Here \(J_d=J_{d,\xi}(T_{\rm an})\) depends on both the character
and the cutoff. Put \(\Gamma_j(s)=\pi^{-(s+j)/2}\Gamma((s+j)/2)/2\).
With \(\tau(\chi)\) the Gauss sum, the theta transformation
\cite[Theorem~10.6, (10.16)--(10.17)]{MV} at \(Qy^2\) gives
\[
 \theta_\chi(y)=\varepsilon_\chi\theta_{\bar\chi}(1/(Qy)),
 \qquad \varepsilon_\chi=\frac{\tau(\chi)}{i^j\sqrt Q},
 \qquad |\varepsilon_\chi|=1 .
\]
On \(\Re s>1\), termwise integration with \(t=\pi(ny)^2\)
gives the full Mellin integral \(\Gamma_j(s)L(s,\chi)\).
The theta transformation gives exponential decay at zero and
infinity, so the integral continues this product. Splitting at
\(Q^{-1/2}\) and substituting \(y\mapsto1/(Qy)\) gives
\[
 \Gamma_j(s)L(s,\chi)
   = I_\chi(s)+\varepsilon_\chi Q^{1/2-s}I_{\bar\chi}(1-s),
 \qquad
 I_a(s)=\int_{Q^{-1/2}}^\infty
       \theta_a(y)y^{s-1/2}\frac{dy}{y}.
\]
On \(\Re s=1/2\) the power of \(Q\) has modulus one and
\(|\theta_{\bar\chi}|=|\theta_\chi|\). Since
\(T_{\rm an}\ge\sqrt Q\), each integral is bounded by \(J_d\)
there.
Let \(M(s)\) be the Mellin integral of
\(\theta_{\mu\chi}(y)y^{s-1/2}\) from \(1/T_{\rm an}\) to \(\infty\)
with measure \(dy/y\). Then
\(K(s)=\Gamma_j(s)L(s,\chi)M(s)-\Gamma_j(s)^2\)
satisfies \(|K|\ll 1+J_d^2\) on \(\Re s=1/2\), since
\(|M(s)|\le J_d\) there and \(\Gamma_j\) is bounded on that
line for both parities. On \(\Re s=2\), absolute convergence
and \(\sum_n\mu(n)\chi(n)n^{-s}=1/L(s,\chi)\) identify the
full reciprocal integral with \(\Gamma_j(s)/L(s,\chi)\).
For \(0<y\le1\),
\(\sum_{n\ge1}(ny)^j e^{-\pi(ny)^2}\ll 1/y\), by comparison
with a decreasing Gaussian of smaller exponent. Hence
\(|\theta_{\mu\chi}(y)|\ll y^{-1/2}\), and the removed tail is bounded by
\[
 \int_0^{1/T_{\rm an}}y^{3/2}O(y^{-1/2})\,\frac{dy}{y}=O(T_{\rm an}^{-1}).
\]
The Euler integral for \(\Gamma_j\) and absolute
convergence of the Dirichlet series bound
\(\Gamma_j(s)L(s,\chi)\) uniformly.

Thus \(|K|\ll T_{\rm an}^{-1}\)
there. For fixed \(\chi\), theta reflection gives exponential
decay at both ends of the completed Mellin integral. On
\(1/2\le\Re s\le2\) the integrals defining \(\Gamma_j L\)
and \(M\) have integrable majorants independent of \(\Im s\).
The majorants remain integrable after multiplication by any
fixed power of \(|\log y|\); hence the integrals are holomorphic
in the interior and continuous and bounded on the closed strip.
The Euler integral gives the same strip bound for \(\Gamma_j\).
The two boundary estimates for \(K\) are uniform in the character
and cutoff.
Suppose \(J_d\le T_{\rm an}^\kappa\) with \(\kappa=(\beta-1/2)/8\).

For \(\lambda=(\Re s-1/2)/(3/2)\), the three-lines theorem
gives
\[
  |K(s)|\ll T_{\rm an}^{2\kappa(1-\lambda)-\lambda}
        \le T_{\rm an}^{-5(\beta-1/2)/12}
        \quad(\beta\le \Re s\le1).
\]
In particular \(K=o(1)\) for the fixed range
\(\beta\le\Re s\le1,\ |\Im s|\le H_0\). On a zero in this box
\(K=-\Gamma_j(s)^2\), whose absolute value has a positive
lower bound depending only on \(\beta,H_0\) for the two
parities. Thus, for sufficiently large \(T_{\rm an}\), even
a real zero forces \(J_{d,\xi}(T_{\rm an})>T_{\rm an}^\kappa\).

For \(d\equiv7\bmod8\), quadratic reciprocity gives
\((n/d)=(-d/n_{\rm odd})\), where \(n_{\rm odd}\) is the odd
part of \(n\), with Jacobi symbol equal to one for denominator one.
For each of \(b(n)=\xi^*(n),\,\mu(n)\xi^*(n)\), let
\[
 G_j(t)=t^je^{-\pi t^2},\qquad
 \Theta_b(h,y)=y^{1/2}\sum_{n\ge1}
       b(n)\Big(\frac{-h}{n_{\rm odd}}\Big)G_j(ny).
\]
For \(h=d\in\mathcal P\) these give the two theta sums in
\(J_{d,\xi}\). They converge absolutely and locally uniformly
for \(y>0\), hence are continuous in \(y\). Positivity gives
\[
 \sum_{\substack{d\le X\\d\in\mathcal P}}|\Theta_b(d,y)|^2
       \le\sum_{1\le h\le X}|\Theta_b(h,y)|^2 .
\]
For fixed \(y>0\) absolute Gaussian convergence permits expansion
of the square on the right. At \(u=(nm)_{\rm odd}\) the
character sum is \(\sum_{h\le X}(-h/u)\). It is bounded by
\(X\) for square \(u\). Otherwise a prime divides \(u\) to odd
order, so \((\cdot/u)\) is a nonprincipal character modulo \(u\)
by the Chinese remainder theorem. P\'olya--Vinogradov
\cite[Theorem~9.18]{MV}, which also covers imprimitive characters,
bounds the sum by \(O(\sqrt u\log(2u))\ll_\eta(nm)^{1/2+\eta}\).
For \(n\le A,\ m\le B\), the square case gives
\[
 \begin{aligned}
 \#\{(n,m):(nm)_{\rm odd}\text{ is a square}\}
 &\le \sum_{a,b\ge0}\sum_{k\le\min(A,B)}
                   \frac{\sqrt{AB}}{k\,2^{(a+b)/2}}\\
 &\ll\sqrt{AB}\log(2\min(A,B)).
 \end{aligned}
\]
Here \(n=2^a k h_1^2,\ m=2^b k h_2^2\) with \(k\) odd and squarefree.
For \(y\le1\), split \(ny,my\) into \([0,1]\) and
\((2^{h-1},2^h]\), \(h\ge1\). On the interval with index \(h\ge0\),
\(|G_j|\ll e^{-c4^h}\) for some \(c>0\). With
\(A=2^h/y,\ B=2^{h'}/y\), the contributions before multiplication
by \(y\) are bounded by
\[
 \begin{aligned}
 X\sum_{h,h'\ge0} e^{-c(4^h+4^{h'})}
               \sqrt{AB}\log(2\min(A,B))&\ll X y^{-1}\log(2/y),\\
 \sum_{h,h'\ge0}e^{-c(4^h+4^{h'})}
                         (AB)^{3/2+\eta}&\ll_\eta y^{-3-2\eta}.
 \end{aligned}
\]
Consequently
\[
 \sum_{1\le h\le X}|\Theta_b(h,y)|^2
       \ll_\eta X\log(2/y)+y^{-2-2\eta}\qquad(1/T_{\rm an}\le y\le1).
\]
For \(y>1\), the termwise absolute bound
\(y^{1/2}\sum_{n\ge1}|G_j(ny)|\ll e^{-c' y^2}\)
suffices. Integrating Cauchy--Schwarz in the finite sum over \(h\) gives
\begin{equation}\label{eq:zero-box-mean}
 \sum_{\substack{d\le X\\d\in\mathcal P}}J_d
 \ll_\eta X\!\left(1+\int_{1/T_{\rm an}}^{1}
                    \sqrt{\log(2/y)}\,\frac{dy}{y}\right)
       +\sqrt X\int_{1/T_{\rm an}}^{1}y^{-2-\eta}\,dy
 \ll_{q_0,\eta}T_{\rm an}^{2+\eta}.
\end{equation}
Let \(\mathcal B\subseteq\mathcal P\) consist of the primes for
which at least one prescribed twist has a zero in the fixed box.
For each sufficiently large \(X\), use the same \(T_{\rm an}=\sqrt{q_0 X}\)
in \eqref{eq:zero-box-mean} for every \(\xi\). With \(\eta=\kappa/2\),
\[
 \#\{d\le X:d\in\mathcal B\}
 \le T_{\rm an}^{-\kappa}
       \sum_{\xi\bmod q_0}\sum_{\substack{d\le X\\d\in\mathcal P}}
          J_{d,\xi}(T_{\rm an})
 \ll_{q_0,\beta,H_0} X^{1-\kappa/4}.
\]
Here the outer sum runs over characters. Partial summation gives
\(\sum_{d\in\mathcal B}d^{-1}<\infty\). The reciprocal sum of
primes in the prescribed reduced class diverges
\cite[Corollary~4.12(c)]{MV}, so some prime from the class
avoids all the bad zeros.
\end{proof}

For \(q_0=18480q\) and \(d_0=q_0/2-1\), the hypotheses on the
residue class follow from Lemma~\ref{lem:auxiliary-residue-class}.
The bound for this progression modulus is derived in the next section.

\section{The bound for the prime remainder}\label{sec:prime-remainder}

\sectionopening{In this section, we apply the cancellation estimate of Proposition~\ref{prop:cancel} to
the prime sum of Section~\ref{sec:criterion}. The error from the
auxiliary scale \(x_1\) is estimated with the fixed-modulus prime
number theorem. We then pass from integer to real arguments before
applying the criterion.}

\begin{proposition}[The one-sided prime-remainder bound]\label{prop:real-prime-remainder}
Fix \(q_0\) with \(18480\mid q_0\) and set \(d_0=q_0/2-1\).
For each prime \(d\equiv d_0\pmod{q_0}\) and
\(0<\epsilon<1/100\), there is a
constant \(C_{q_0,d,\epsilon}\) with
\[
 R_{d,q_0}(y)\le C_{q_0,d,\epsilon}y^{-1/2+6\epsilon}
\]
for all sufficiently large real \(y\).
\end{proposition}

\begin{proof}

Write \(B=\log b\), \(u=\#\Phi_a\), \(v=|T|>0\), and
\(\mathcal R\) for the total remainder defined in
Section~\ref{sec:joint}. The exact expansion is
\[
 \frac{\mathcal R}{av}=
 R_{d,q_0}(x)-R_{d,q_0}(x_1)+\frac{u}{v}R_{d,1}(x_1).
\]
The datum has \(u/v\le B\) and \(f=\varphi(q_0)u\).
Thus \(m\mathcal R\le faB^5\) gives
\[
 R_{d,q_0}(x)\le
 \frac{\varphi(q_0)B^6}{m}
 +|R_{d,q_0}(x_1)|+B|R_{d,1}(x_1)|.
\]

For either \(q\in\{q_0,1\}\), put
\[
 S_q(t)=\sum_{n\le t}\lambda_{d,q}(n),\qquad
 E_q(t)=S_q(t)-t/\varphi(q).
\]
The prime number theorem in arithmetic progressions
\cite[Corollary~11.19]{MV} applies at the fixed modulus \(qd\).
Let \(b_a\bmod qd\) be the class \(1\bmod q\), \(a\bmod d\),
and write \(\psi(t;qd,b_a)=\sum_{n\le t,\ n\equiv b_a\ (qd)}\Lambda(n)\).
With \(o_1=1\) and \(o_{q_0}=2\), Lemma~\ref{lem:auxiliary-residue-class}
gives
\[
 S_q(t)=\sum_{a\in(\mathbb Z/d\mathbb Z)^\times}
           (1-\chi_d(a))\psi(t;qd,b_a)
       +\left\lfloor\frac{\log t}{o_q\log d}\right\rfloor\log d .
\]
Indeed \(d^e\equiv1\bmod q\) precisely when \(o_q\mid e\);
these powers have weight one since \(\chi_d(d^e)=0\).
If \(p\mid q_0\), then \(p\) splits in \(K\) by
Lemma~\ref{lem:auxiliary-residue-class}. Hence \(1-\chi_d(p^e)=0\)
for every \(e\ge1\), also in the modulus-one sum.
The reduced-class weights sum to \(d-1\), so their main terms sum
to \((d-1)t/\varphi(qd)=t/\varphi(q)\).
For \(t\ge e^{q_0d}\) both progression moduli satisfy \(qd\le\log t\).
Applying the cited corollary with this exponent one and bounding the
powers of \(d\) gives, for an absolute \(c_1>0\),
\[
 |E_q(t)|\ll (d-1)t e^{-c_1\sqrt{\log t}}+\log t .
\]
Absorbing the remaining bounded interval into the constant gives,
for every \(A>0\),
\[
 E_q(t)=O_{q,d,A}\bigl(t(\log t)^{-A-1}\bigr)\qquad(t\ge2).
\]

For real \(\sigma>1\), the Dirichlet series identity in
Section~\ref{sec:one-sided-mellin-argument} at modulus \(q\) gives
\[
 \int_1^\infty E_q(t)t^{-\sigma-1}\,dt
   =\frac{D_q(\sigma)+(\sigma-1)^{-1}}{\sigma\varphi(q)}
       -\frac{1}{\varphi(q)(\sigma-1)}
   =\frac{D_q(\sigma)-1}{\sigma\varphi(q)} .
\]
The bound on \(E_q\) supplies an integrable majorant at
\(\sigma=1\). Letting \(\sigma\downarrow1\) identifies
\(c_{d,q}=\int_1^\infty E_q(t)t^{-2}\,dt\).
Integrating each term in the finite sum for \(S_q\) up to \(y\)
now gives the exact tail and its bound
\[
 R_{d,q}(y)=-\int_y^\infty E_q(t)t^{-2}\,dt
     =O_{q,d,A}\bigl((\log y)^{-A}\bigr) .
\]

Choose \(0<\epsilon<1/100\) and use a datum at each large
integer \(x\), with \(B\le2x^\epsilon\), \(m\ge\sqrt x\), and
\(\log x_1\ge x^{\epsilon/12}\). Take \(A>24/\epsilon\) in the
fixed-modulus bound. Both remainders at \(x_1\) are at most
\(x^{-2}\) in absolute value for sufficiently large \(x\).
Consequently
\[
 R_{d,q_0}(x)\le (64\varphi(q_0)+3)x^{-1/2+6\epsilon}.
\]
For real \(y\ge2\), the exact tail gives
\[
 R_{d,q_0}(y)-R_{d,q_0}(\lfloor y\rfloor)
    =\int_{\lfloor y\rfloor}^{y}E_{q_0}(t)t^{-2}\,dt
    =O_{q_0,d}(y^{-1}).
\]
Put \(\alpha=-1/2+6\epsilon\). As \(-1<\alpha<0\) and
\(\lfloor y\rfloor\ge y/2\), we have
\(\lfloor y\rfloor^{\alpha}\le2^{-\alpha}y^\alpha\).
The integer bound gives
\[
 R_{d,q_0}(y)\le R_{d,q_0}(\lfloor y\rfloor)+O_{q_0,d}(y^{-1})
                    \le C_{q_0,d,\epsilon}y^{-1/2+6\epsilon}.\qedhere
\]
\end{proof}

\begin{proof}[Proof of Theorem~\ref{thm:rh}]
Let \(\chi\) be a character modulo \(q\). Induce it to
\(\widetilde\chi\) modulo \(q_0=18480q\). The Euler product gives
\[
 L(s,\widetilde\chi)=L(s,\chi)
       \prod_{\substack{p\mid q_0\\p\nmid q}}(1-\chi(p)p^{-s}).
\]
The finite factors are nonzero for \(\Re s>0\).
Suppose \(L(\rho,\chi)=0\) with \(1/2<\Re\rho<1\).
Choose
\[
 \tfrac12<\beta<\Re\rho,\qquad H_0>|\Im\rho|.
\]
Proposition~\ref{prop:zero-free-auxiliary-twist}, with
\(d_0=q_0/2-1\), gives a prime \(d\) for which all the level-\(q_0\)
quadratic twists are zero-free in the prescribed rectangle.
Fix \(d\) and take \(0<\epsilon<1/100\) with
\(6\epsilon<\beta-1/2\). Proposition~\ref{prop:real-prime-remainder}
gives a cutoff \(X_R=X_R(q_0,d,\epsilon)\) and a constant
\(C=C_{q_0,d,\epsilon}\ge1\). Since \(\beta-1/2-6\epsilon>0\),
\[
 x\ge\max\{X_R,C^{1/(\beta-1/2-6\epsilon)}\}
 \quad\Longrightarrow\quad
 R_{d,q_0}(x)\le Cx^{-1/2+6\epsilon}\le x^{\beta-1}.
\]
Hence Lemma~\ref{lem:positive-mellin-criterion} gives
\(L(\rho,\widetilde\chi)\ne0\), a contradiction. For \(\Re s>1\), use the Euler product. At \(s=1\), use
\cite[Theorem~4.9]{MV} and the simple pole for the principal
character. At \(s=1+it\), \(t\ne0\), the nonnegative
identity \(3+4\cos\theta+\cos(2\theta)=2(1+\cos\theta)^2\)
applied termwise to the Euler logarithms gives
\[
 \zeta(\sigma)^3 |L(\sigma+it,\chi)|^4
                      |L(\sigma+2it,\chi^2)|\ge1\qquad(\sigma>1).
\]
Primes dividing \(q\) contribute only the positive zeta term.
As \(\sigma\downarrow1\), the last factor stays bounded
because its only possible pole is at \(1\). A zero of the middle
factor at \(1+it\) would make the product \(O(\sigma-1)\), since
\(\zeta\) has a simple pole. This contradicts the lower bound.

Let \(\psi\) be the primitive character inducing \(\chi\), of
conductor \(q_*\). The extra Euler factors relating their
\(L\)-functions \cite[(10.20)]{MV} can vanish only on
\(\Re s=0\), by taking absolute
values in \(\psi(p)p^{-s}=1\) at a deleted prime.
For \(q_*>1\) put
\[
 \Lambda_\psi(s)=(q_*/\pi)^{(s+j)/2}
                 \Gamma((s+j)/2)L(s,\psi),
 \qquad \psi(-1)=(-1)^j .
\]
The primitive functional equation
\cite[Corollary~10.8, (10.19)]{MV} relates
\(\Lambda_\psi(s)\) to \(\Lambda_{\bar\psi}(1-s)\) by a nonzero
constant. For \(q_*=1\), the inducing function is \(\zeta\) and its
completed function has the corresponding reflection
\cite[Corollary~10.3]{MV}. In the open strip the gamma factors
are finite and nonzero. Thus a zero with \(0<\Re s<1/2\)
would give a zero at \(1-s\) for the primitive conjugate, already
excluded.

For \(q_*>1\) and \(\Re s\le0\), the reflected completed function
is finite and nonzero by the nonvanishing on \(\Re(1-s)\ge1\).
The only zeros of \(L(s,\psi)\) there come from the reciprocal gamma
factor, at \(s=-j-2n\) for \(n\ge0\), each simple. For \(q_*=1\) use
\cite[Corollary~10.4, (10.9)]{MV} in the form
\[
 \zeta(s)=2^s\pi^{s-1}\sin(\pi s/2)\Gamma(1-s)\zeta(1-s).
\]
For \(\Re s\le0\) and \(s\ne0\) this gives precisely the negative even
zeros. At \(s=0\) the sine cancels the pole of \(\zeta(1-s)\) and
gives \(\zeta(0)=-1/2\). For imprimitive characters the deleted Euler
factors account for all additional boundary zeros. For each prime \(p\mid q\) with \(p\nmid q_*\),
the deleted factor contributes a simple zero at each
\(i(\arg\psi(p)+2\pi k)/\log p\), \(k\in\mathbb Z\).
At a nonzero ordinate \(t\), this can occur for at most one prime.
Indeed, \(t\log p/(2\pi)\) is a nonzero rational number because
\(\psi(p)\) is a root of unity, whereas \(\log p/\log p'\) is
irrational for distinct primes by unique factorization.
At the origin the order is
\[
 \operatorname{ord}_{s=0} L(s,\chi)
 =\mathbf1_{\{q_*>1,\ \psi(-1)=1\}}
   +\#\{p\mid q:p\nmid q_*,\ \psi(p)=1\}.\qedhere
\]
\end{proof}

\appendix
\phantomsection
\addcontentsline{toc}{part}{Appendices}
\section{Period rings and integral comparison matrices}
\label{sec:crystalline-coordinate-appendix}

\sectionopening{In this appendix, we construct the coordinate maps and comparison
matrices for the local lattices of Section~\ref{sec:relative-crystalline-colength}.
The approximation and divisibility results give the integral
period bases used in Appendix~\ref{sec:crystalline-basis-appendix}.}

Use the local representations, jump convention and normalized
defect of Section~\ref{sec:relative-crystalline-colength}.
Here \(p\) is the residue characteristic at the coefficient place.

\subsection{Integral coordinates for the period basis}
\label{sec:period-chart-coordinates}

\sectionopening{In this subsection, we construct coordinate maps from \(A_{\inf}\) to the power series
ring \(S\). The construction uses the graded modules for the ideal
\((p,z)\), and it extends to the bounded period ring.
Section~\ref{sec:period-hodge-retraction} uses these maps to preserve
divisibility at the de Rham period.}

Write \(A=A_{\inf}=W(\mathcal O_{\mathbb C_p^\flat})\),
\(S=W(k)[[z]]\subset A\), and \(I=(p,z)A\).
Use the chosen Teichm\"uller coordinate \(z\) and generator
\(E=p+z^e u\) of \(\ker\theta\) from
Section~\ref{sec:coefficient-hodge-tate-intervals}, or
\(E=p-z^e\) in the tame chart constructed below.
Here \(K\) is the finite extension of \(\mathbf Q_p\) defining the chart, \(k\) its
residue field, \(e\) its absolute ramification index, and \(u\in S^\times\).
For the comparison matrix assume that \(T[1/p]\) is crystalline over
\(K\). If \(T\) is defined over the integers \(O\) of a finite
coefficient field, take its underlying \(\mathbf Z_p\)-lattice of rank
\(d=\operatorname{rank}_{\mathbf Z_p}(O)\,\operatorname{rank}_O T\).
The coefficient action is then by \(\mathbf Z_p\)-linear endomorphisms.
The invariant crystalline and de Rham bases below are over
\(K_0=W(k)[1/p]\) and \(K\), respectively.
\begin{lemma}[Graded coordinates of the period ring]\label{lem:period-chart-graded-coordinates}
The ring \(A\) is complete and separated for \(I=(p,z)A\).
For each \(n\ge0\), there is an isomorphism
\[
 \bigoplus_{i=0}^n A/I\xrightarrow{\ \sim\ } I^n/I^{n+1},
 \qquad (\bar a_i)_i\longmapsto
       \left[\sum_{i=0}^n a_i p^i z^{n-i}\right].
\]
\end{lemma}

\begin{proof}

The reduction of \(z\) in the tilt is nonzero and is a nonunit. The ring \(A\) is \(p\)-torsion free and
\(A/p\) is a complete rank-one valuation ring.
Multiplication by \(z\) is injective on \(A/p^m\)
for each \(m\), by induction using \(A/p\).
Also \(A/p^m\) is complete and separated for its
\(z\)-adic topology. To see the induction, use
\(p^{m-1}A/p^m\simeq A/p\), with its induced
topology: if \(z^h a\) is in \(p^{m-1}A\),
then \(a\in p^{m-1}A\) by \(z\)-injectivity modulo
\(p^{m-1}\). Thus the induced topology on the kernel
of \(A/p^m\to A/p^{m-1}\) is its own \(z\)-adic
topology. Thus reduction modulo \(z^h\) preserves the
short exact sequence with these three terms. The transition
maps on the kernels \((A/p)/z^h\) are surjective, so the
inverse limit over \(h\) is still exact. Comparing the two
short exact sequences proves completeness and separatedness
of \(A/p^m\) by induction. Finally
\[
 I^{m+h-1}\subset(p^m,z^h)A\subset I^{\min(m,h)},\qquad
 \varprojlim_{m,h} A/(p^m,z^h)
   =\varprojlim_m A/p^m=A .
\]

For each \(n\), the monomials \(p^i z^{n-i}\),
\(0\le i\le n\), identify \(I^n/I^{n+1}\) with
\((A/I)^{n+1}\). Surjectivity follows by expansion.
For injectivity reduce a relation modulo \(p\)
and cancel \(z^n\) to see that its coefficient of
\(z^n\) lies in \(I\). Subtract that term.
The intersection \(pA\cap I^{n+1}=pI^n\)
follows by reducing a generating expression for
\(I^{n+1}\) modulo \(p\): its pure \(z\)-term has
coefficient divisible by \(p\). Divide the
remaining relation by \(p\) and induct on \(n\).
This proves the asserted graded formula.
\end{proof}

Choose a \(k\)-basis of \(A/I\) containing \(1\),
and lift its vectors to \(A\), using \(1\) for the
distinguished lift. The graded formula shows that
the resulting map \(\bigoplus_i S/(p,z)^n\to A/I^n\)
is an isomorphism for each \(n\).
Taking inverse limits identifies \(A\) with the families
\((s_i)_i\in\prod_i S\) for which, at each \(n\), all but
finitely many \(s_i\) lie in \((p,z)^nS\). The inverse sends
\((s_i)_i\) to the \(I\)-adic sum \(\sum_i b_i s_i\), where
\(b_i\) are the chosen lifts. The projections \(F_i:A\to S\)
are \(S\)-linear and continuous for these adic topologies.
They separate elements and preserve \(p\)- and \(z\)-divisibility.
The coordinate of \(1\) satisfies \(F_0|_S=\operatorname{id}_S\).

For the bounded period ring, realize
\(B^+=A_{\max}[1/p]\) with \(A_{\max}\) the
\(p\)-adic completion of \(A[E/p]\). Equip \(B^+\) with the
ring topology whose neighborhoods of zero are the additive lattices
\(p^n A_{\max}\), \(n\ge0\). All closures and bounded-ring
approximations below refer to this topology.
For any \(S\)-linear \(F:A\to A\), first extend by inverting
\(p\). Since \(E\in S\), the extended map sends
\((E/p)^j a\) to \((E/p)^j F(a)\). For a polynomial relation
multiply by a power of \(p\) clearing the denominators,
use \(S\)-linearity in \(A\), and cancel that power.
The extension preserves \(A[E/p]\), hence extends
continuously to \(A_{\max}\) and \(B^+\).
It commutes with multiplication by the closure
\(S_B\) of \(S[1/p]=S[E/p][1/p]\) in \(B^+\), by continuity.
For the coordinate maps \(F_i:A\to S\) the extension
takes values in \(S_B\). This also identifies \(S_B\) with the
scalars whose multiplication commutes with every \(S\)-linear
endomorphism of \(B^+\) carrying \(A_{\max}\) into itself.
Indeed such endomorphisms are continuous and commute with
\(S[1/p]\), hence with its closure \(S_B\). Conversely, for
any scalar \(b\) in this commutant, the distinguished coordinate
satisfies \(b=b\,\widehat F_0(1)=\widehat F_0(b)\in S_B\).
In particular one denominator suffices at all precisions:
\[
 S_B=\bigl(\overline{S[E/p]}^{\,A_{\max}}\bigr)[1/p].
\]
Indeed, for \(b\in S_B\) take \(m\) with \(p^m b\in A_{\max}\).
Approximate \(p^m b\) by elements of \(A[E/p]\).
The continuous map \(\widehat F_0\) fixes \(p^m b\) and sends these
approximations to \(S[E/p]\). The reverse inclusion follows from
the definition of \(S_B\).
Frobenius preserves \(S_B\): it preserves \(S\), and
\(\varphi(E)/p=1+(z^e)^p\varphi(u)/p\) if
\(E=p+z^e u\). The relation \(z^e/p=(E/p-1)/u\)
puts this expression in \(S[E/p]\). Since \(\varphi(A)\subset A\)
and \(\varphi(p)=p\), Frobenius preserves \(A[E/p]\) and extends
continuously to \(A_{\max}\) and \(B^+\). For any finite
tuple \(y\) over \(B^+\) and any precision \(p^j A_{\max}\),
approximate its entries by polynomials in \(E/p\)
after clearing \(p\)-denominators.

Expand the finitely many coefficients in the lifted basis
\((b_i)_i\) modulo a high power of \((p,z)\). Since
\((p,z)^eA\subset pA_{\max}\), these truncations approximate every
entry, to the required precision, by a finite combination
\(y_0=\sum_{i\in J_0} b_i s_i\) with coefficients in \(S_B\).
Use the union of the supports for all entries. For a finite
\(J\supseteq J_0\) put
\(P_J=\sum_{i\in J}b_i\widehat F_i\). Each \(P_J\) preserves
\(A_{\max}\) and fixes \(y_0\), since
\(\widehat F_i(b_k s)=\delta_{ik}s\) for \(s\in S_B\).
Thus, entrywise,
\[
 y-y_0\in p^j A_{\max}\quad\Longrightarrow\quad
 y-P_Jy=(1-P_J)(y-y_0)\in p^j A_{\max}.
\]
The extended coordinates separate elements of \(B^+\) by
separatedness of its lattice topology.

Suppose \(U\in\operatorname{Mat}_d(B^+)\) is invertible over
\(\operatorname{Frac}(B^+)\). If \(U y_j=E^h x_j\) with
\(y_j\in S_B^d\), an \(A\)-linear relation
among the \(x_j\) is the same as the relation
among the \(y_j\). Applying each extended coordinate
to a relation \(\sum_j a_jy_j=0\) gives
\(\sum_j F_i(a_j)y_j=0\), since \(y_j\in S_B^d\). Thus also
\(\sum_j F_i(a_j)x_j=0\). If the \(x_j\) are \(S\)-independent,
then \(F_i(a_j)=0\) for every \(i,j\). The coordinates separate
elements of \(A\), so all \(a_j=0\). This bounds the cardinality
of the family by \(d\).

Let \(\mathcal P\) be the \(S\)-submodule of \(A^d\)
of vectors \(x\) for which \(Uy=E^h x\) has a solution
\(y\in S_B^d\). It is \((p,z)\)-adically separated inside
\(A^d\), and \(p\)-saturated since \(1/p\in S_B\).
In particular \(\mathcal P\cap pA^d=p\mathcal P\). For
\(x\in\mathcal P\) with \(zx\in p\mathcal P\), reduction in the
domain \(A/p\) gives \(x\in pA^d\), hence \(x\in p\mathcal P\).
Thus \(z\) acts injectively on \(\mathcal P/p\mathcal P\).
Put \(J=(p,z)S\).
In an element of \(p\mathcal P\cap J^{n+1}\mathcal P\),
the coefficient of the pure term \(z^{n+1}\) is zero modulo
\(p\mathcal P\), by this injectivity. Thus
\(p\mathcal P\cap J^{n+1}\mathcal P=pJ^n\mathcal P\). Cancelling
\(p\) in the remaining terms as in the graded calculation for \(A\)
gives the monomial formula for \(J^n\mathcal P/J^{n+1}\mathcal P\)
with coefficient space \(\mathcal P/J\mathcal P\).
Any independent family of residue vectors lifts
to an independent family over \(S\), by taking
the least \((p,z)\)-order in a presumed relation.
So the residue space has dimension at most \(d\).
Lift a basis, and express any element of
\(\mathcal P\) successively modulo \((p,z)^n\)
in those lifts \(x_i\), correcting the coefficients at
step \(n\) by elements of \((p,z)^nS\). The coefficients
converge to \(s_i\in S\). Then
\[
 x-\sum_i s_i x_i\in (p,z)^n\mathcal P\subset I^nA^d
 \qquad(n\ge1).
\]
Separation of \(A^d\) gives equality. Thus the \(x_i\)
form a finite \(S\)-basis of the solution module.

Finally a principal ideal \((\delta)\) with
\(p\nmid\delta\) is closed in \(A\) for the topology
of \((p,z)\). There is \(k>0\) such that
\(\bar\delta\mid\bar z^k\) in the rank-one valuation
ring \(A/p\). Write \(\delta u=z^k+pv\).
Modulo \(p^m\), the geometric sum gives
\[
 \delta u\sum_{i=0}^{m-1}z^{k(m-1-i)}(-pv)^i
   \equiv z^{km}\pmod{p^m A}.
\]
If \(\delta a_j\) is a Cauchy sequence in \(A\) for
\((p,z)\), then modulo \(p^m\) it is Cauchy for \(z\).
Multiplying any difference in \(z^{L+km}(A/p^m)\)
by the displayed multiplier and cancelling the
nonzero divisor \(z^{km}\) puts the difference
of the \(a_j\)'s in \(z^L(A/p^m)\).
This works for every \(L,m\). The \(a_j\) therefore
converge in \(A\), and their product with \(\delta\)
has limit in \((\delta)\).
The topology of \((p,\xi)\) is the same here, since
\(\xi\) is associated to \(E=z^e a+p\) with \(a\)
a unit of the chart.
This supplies the closedness used to pass from
the finite coordinate partial sums to divisibility
in the basis proof.

\subsection{Integral approximation in the fixed period chart}
\label{sec:period-chart-approximation}

\sectionopening{In this subsection, we establish divisibility for the coordinate
maps of Section~\ref{sec:period-chart-coordinates} and use the
Ax--Sen--Tate theorem in the tower of roots of the uniformizer
to approximate fixed tilt vectors. The approximation gives
equality in the prescribed tilt component modulo \(p\). It is used to
prove that the determinant in
Section~\ref{sec:local-integral-period-frame} is nonzero modulo \(p\).}

For the chart \(S\subset A_{\inf}\) of the preceding construction
and each integer \(j\ge0\), put
\(\theta_j=\theta\circ\varphi^{-j}\) and \(E_j=\varphi^j E\).
The map \(\theta_j:A_{\inf}\to\mathcal O_{\mathbb C_p}\)
is surjective, with kernel \((E_j)\), since Frobenius is an
automorphism of \(A\) and \(E\) generates \(\ker\theta\).
Moreover \(\varphi(S)\subset S\), so \(E_j\in S\). For any \(S\)-linear map \(F:A_{\inf}\to S\),
\(\alpha\in S\), \(r\in A_{\inf}\), the implication
\[
 \theta_j(\alpha)\mid\theta_j(r)
 \quad\Longrightarrow\quad
       F(r)-\alpha x\in E_j S\quad\text{for some }x\in S
\]
follows by lifting an integral quotient to \(a\in A_{\inf}\).
Then \(r-\alpha a=E_j b\) and one can take \(x=F(a)\).
This gives the preservation of divisibility by \(p\)
in the residue at the \(j\)-th component for every
chart coordinate map used in the period basis.

Set \(K_0=W(k)[1/p]\) and \(\pi_j=\theta_j(z)\).
The Teichm\"uller coordinate gives
\(\pi_{j+1}^p=\pi_j\), with \(\pi_0\) the chosen uniformizer of \(K\).
If \(H(X)\) is the Eisenstein polynomial of \(\pi_0\)
over \(K_0\), then \(H(X^{p^j})\) is Eisenstein of degree
\(d_j=ep^j\). Thus \(L_j=K_0(\pi_j)\) has degree \(ep^j\) over \(K_0\),
with \(L_0=K\), and has ring of integers
\(W(k)[\pi_j]\). Indeed, in the basis \(1,\ldots,\pi_j^{d_j-1}\)
the nonzero terms \(c_i\pi_j^i\) have valuations
\(v_p(c_i)+i/d_j\) with distinct fractional parts.
Their sum is integral exactly when all \(c_i\in W(k)\).
Now \(\theta_j\) maps \(W(k)\) onto itself by inverse
Frobenius, and every series in \(S\) converges upon
evaluation at \(\pi_j\). Consequently
\[
 R_j=\theta_j(S)=\mathcal O_{L_j},\qquad
 L_\infty=\bigcup_{j\ge0}L_j,\qquad
 \mathcal O_{L_\infty}=\bigcup_{j\ge0}R_j .
\]
For \(\theta_j(r)\ne0\), take
\(n=\lfloor d_j v_p(\theta_j(r))\rfloor\) and \(\alpha=z^n\)
in the divisibility implication. Since \(v_p(\pi_j)=1/d_j\),
\[
 v_p\bigl(\theta_j(F(r))\bigr)
   \ge \frac{\lfloor d_j v_p(\theta_j(r))\rfloor}{d_j}.
\]
If \(\theta_j(r)=0\), then \(r\in E_j A_{\inf}\) and
\(\theta_j(F(r))=0\). In particular all divisibilities by powers
of \(p\) are preserved in the residue. The order retraction in
Section~\ref{sec:period-basis-order-retraction} uses exact vanishing
on the corresponding de Rham grades.
The pointwise fixer of \(S\) fixes each \(R_j\) by
Galois equivariance of \(\theta_j\). Conversely, fixing
\(L_\infty\) fixes \(W(k)\) and every \(\pi_j\), hence
the tilt sequence defining \(z\). Such an element fixes
\(S=W(k)[[z]]\) by continuity.

Let \(y_1,\ldots,y_l\in\mathcal O_{\mathbb C_p^\flat}\)
be fixed by this group, and fix \(m\ge0\). The elements
\(t_i=\theta_m([y_i])\) are integral and fixed.
The algebraic extension \(L_\infty/\mathbf Q_p\) is
Henselian with its induced valuation and has characteristic
zero. The Ax--Sen--Tate fixed field theorem
\cite[theorem, p.~417]{AxZeros} gives
\[
 \mathbb C_p^{\operatorname{Gal}(\overline{\mathbf Q}_p/L_\infty)}
      =\widehat{L_\infty}.
\]
Choose \(x_i\in L_\infty\) with \(v_p(x_i-t_i)>1\) by density.
The \(x_i\) are integral.
Choose \(j_0\ge m\) with \(x_1,\ldots,x_l\in R_{j_0}\).
For every \(j\ge j_0\) the inclusions \(R_{j_0}\subset R_j\)
give \(s_{j,i}\in S\) satisfying
\[
 \theta_j(s_{j,i})=x_i\equiv\theta_m([y_i])\pmod p
 \qquad(1\le i\le l).
\]
Thus at every sufficiently large common level \(j\) the
\(j\)-th tilt components of the \(s_{j,i}\bmod p\) equal
the prescribed \(m\)-th components of all the \(y_i\).

\subsection{Uniform denominator clearing in period coordinates}
\label{sec:period-period-contraction}

\sectionopening{In this subsection, we combine division in \(A_{\max}\) with a Frobenius relation to
obtain a bound for integral denominators in the solution module of
Section~\ref{sec:period-chart-coordinates}. The bound applies
uniformly to all its coordinate maps and is used in
Section~\ref{sec:period-period-saturation}.}

In \(A=A_{\inf}\) choose the cyclotomic generator of \(\ker\theta\)
as in \cite[\S3.3, before Definition~3.22]{BMSIntegral}. Explicitly,
\[
 \tau=[\epsilon^{1/p}]-1,\qquad
 \xi=\sum_{i=0}^{p-1}(1+\tau)^i,\qquad
 \varphi(\tau)=[\epsilon]-1=\mu=\xi\tau .
\]
Thus \(E/\xi\in A^\times\).
Put \(R=A[X]/(pX-\xi)\). If \(pF=(pX-\xi)G\) in \(A[X]\),
then \(\bar\xi\,\bar G=0\) in \((A/p)[X]\).
Since \(\bar\xi\ne0\) in the tilt domain, \(G=pG_1\).
Cancelling \(p\) shows that \(R\) is \(p\)-torsion free.
It therefore embeds in \(A[1/p]\) by \(X\mapsto\xi/p\), since this
map becomes an isomorphism after inverting \(p\).
Its \(p\)-adic completion \(A_{\max}\) is also
\(p\)-torsion free: an element killed by \(p\) has
component modulo \(p^{m+1}\) in \(p^mR/p^{m+1}R\)
for every \(m\), hence is zero. Dividing compatible
components by \(p\) at one higher precision also gives
\(A_{\max}/pA_{\max}=R/pR=(A/(p,\xi))[X]\).
Multiplication by a monic polynomial is injective in a polynomial
ring over any coefficient ring. Hence if \(P\in A_{\max}\) reduces
to a monic polynomial and \(P w\in p^m A_{\max}\), successive
reduction modulo \(p\) and cancellation show \(w\in p^m A_{\max}\).

Write \(\iota:A\to A_{\max}\) for the natural map.
If \(\iota(a)\in p^mA_{\max}\), then \(a\in(p,\xi)^m A\).
The map \(\theta\) extends continuously to \(A_{\max}\) by
\(X\mapsto0\). For \(m>0\) we have \(\theta(a)\in p^m\mathcal O_{\mathbb C_p}\);
lift the quotient under \(\theta\) to \(c\in A\).
Then \(a=\xi w+p^m c\) with \(w\in A\).
Using \(\xi=pX\) and cancelling \(p\) shows
\(Xw\in p^{m-1}A_{\max}\).
The preceding regularity and induction put
\(w\) in \((p,\xi)^{m-1}\).
Conversely, \(\iota((p,\xi)^m A)\subset p^mA_{\max}\) by
\(\xi=pX\). Hence \(\ker\iota\subset\bigcap_m(p,\xi)^m A=0\)
by the separation in Section~\ref{sec:period-chart-coordinates}.
Also \(X\) is a nonzero divisor in \(A_{\max}\) by its monic
reduction and \(p\)-adic separation, so \(\xi=pX\) can be
cancelled in both \(A_{\max}\) and \(B^+\).

Since \(\xi-p\in\tau A\), we have
\(p(X-1)=\tau u\) for an integral \(u\).
For \(c=0,1\), every \(y\in A_{\max}\) can be written
\(a+(X-c)b\) with \(a\in A\) and \(b\in A_{\max}\).
Indeed, reduce \(y\) to \(P(X)\in (A/(p,\xi))[X]\).
Lift \(P(c)\) to \(a_0\in A\) and lift the polynomial
\((P(X)-P(c))/(X-c)\) to \(b_0'\in A[X]\).
Then \(y-a_0-(X-c)b_0'\in pA_{\max}\).
Repeat on the quotient by \(p\); the scalar digits converge
in \(A\) and the polynomial digits converge in \(A_{\max}\)
after multiplication by their successive powers of \(p\).
Taking \(c=1\) and using \(p(X-1)=\tau u\) gives
\(p A_{\max}\subset A+\tau A_{\max}\), and iteration gives
\(p^j A_{\max}\subset A+\tau^j A_{\max}\).
Taking \(c=0\) gives \(p A_{\max}\subset A+\xi A_{\max}\).

Let \(V\subset A_{\max}\) and suppose for each \(v\in V\)
there is \(w\in V\) such that \(\xi^h v=u_0\varphi(w)\)
with \(u_0\in A\) fixed.
Then for every \(m\ge0\) there exists \(a_m\in A\)
with
\[
  p^{h+1}v-a_m\in\tau^{h+1}p^m A_{\max}.
\]
The case \(m=0\) follows from the preceding polynomial
division. For the induction step, if
\(p^{h+1}w=a_m(w)+\tau^{h+1}p^m b_m\), \(b_m\in A_{\max}\), then
\[
 \xi^h p^{h+1}v=u_0\varphi(a_m(w))
  +\xi^{h+1}\tau^{h+1}p^m u_0\varphi(b_m).
\]
Here \(u_0\varphi(a_m(w))\in A\cap\xi^h B^+=\xi^h A\).
Indeed for \(a\in A\cap\xi B^+\) evaluation by \(\theta\) gives
\(a\in\ker(\theta|_A)=\xi A\), and one cancels \(\xi\) and
iterates. Dividing the display by \(\xi^h\) leaves an error in
\(\tau^{h+1}p^{m+1}A_{\max}\), since \(\xi=pX\).
This gives the next integral approximation.
The difference of two approximations at precision
\(m\) lies in \(p^m A_{\max}\cap A=(p,\xi)^m A\).
They converge in \(A\), by
Section~\ref{sec:period-chart-coordinates}.
Their limit in \(A_{\max}\) is \(p^{h+1}v\).
Thus \(p^{h+1}V\) has integral \(A\)-coordinates.

For the period matrix \(U\), apply this to the
coordinates of the solutions \(w\) of
\[
  U y=E^h x,\quad x\in A^d,\qquad
  U\widehat F(y)=E^h w
\]
as \(y\) and the \(S\)-linear maps \(F:A\to A\) vary. There is one power \(p^m\) putting
all such \(w\) in \(A_{\max}\): the inverse of \(U\)
up to \((E\tau)^h\) first bounds a denominator of
\(y\), and integrality of the extended coordinate
maps gives the same bound for \(F(y)\).

Applying
\(U/E^h\) requires one further fixed power.
For the denominator claim, note that
\(\xi\equiv\tau^{p-1}\bmod p\), so
\(\tau^{p-1}=p(X-b_0)\) with \(b_0\in A\).
If \(\mu z\in A_{\max}\) then multiplication
by \(\tau^{p-2}\) gives
\(p^2 X(X-b_0)z\in A_{\max}\).
Multiplication by the monic polynomial \(X(X-b_0)\)
is \(p\)-saturated, by its reduction modulo \(p\).
Clearing an initial power of \(p\) in \(B^+\)
and using this saturation yields \(p^2z\in A_{\max}\).
Iteration gives \(p^{2h}z\in A_{\max}\) if
\(\mu^h z\in A_{\max}\); the same conclusion holds when
\(d^h z\in A_{\max}\) with \(d\mid\mu\) in \(A\).
For the individual factors, \(X\) and \(X-b_0\) are both
monic modulo \(p\). Hence, for \(z\in B^+\),
\[
 \xi z\in A_{\max}\ \text{or}\ \tau z\in A_{\max}
       \quad\Longrightarrow\quad pz\in A_{\max}.
\]
In the first case use \(\xi z=X(pz)\). In the second case
multiply by \(\tau^{p-2}\) to get \((X-b_0)pz\in A_{\max}\).
At \(p=2\) this multiplier is \(1\).

Section~\ref{sec:period-period-saturation} verifies Frobenius descent
for the retracted solutions and applies the contraction with a uniform
exponent.

\subsection{The two period matrices over the bounded ring}
\label{sec:period-period-realization}

\sectionopening{In this subsection, we construct the crystalline comparison matrix and its partial inverse
over \(B^+\). The key scalar assertion is divisibility by \(\mu\) in
the common kernel of the Fontaine evaluations at nonnegative Frobenius
powers. This gives the matrix equations used in
Section~\ref{sec:local-integral-period-frame}.}

For the same chart write \(\mu=\xi\tau\).
In \(A=A_{\inf}\), the common kernel of
\(\theta\varphi^j\), \(j\ge0\), is \(\mu A\).
The individual kernels are generated by the
distinct primes \(\varphi^{-j}\xi\).
For \(0\le i<j\), \(\theta\varphi^j(\varphi^{-i}\xi)=p\),
since \(\theta\varphi^{j-i}([\epsilon^{1/p}])=1\).
On \(\varphi^{-j}\xi\) the same evaluation is zero.
Thus these prime ideals are pairwise distinct. The intersection of any finite initial set of these principal
prime ideals is their product. Its generator is
\[
 \prod_{j=0}^{m-1}\varphi^{-j}(\xi)
      =\mu/\varphi^{-m}(\mu).
\]
The intersection of these products for \(m\ge1\) is \(\mu A\)
by \cite[Lemma~3.23]{BMSIntegral}.
Applying this result to \(\varphi(a)\) gives
the kernel \(\tau A\) for \(\theta\varphi^j\), \(j\ge1\).

For \(y\in B^+\) choose \(c\ge0\) such that
\(r=p^c y\in A_{\max}\).
Section~\ref{sec:period-period-contraction} gives
\[
 pr=a+\tau v,\qquad pv=a_1+\xi w',
 \qquad a,a_1\in A,\quad v,w'\in A_{\max}.
\]
Thus \(p^{c+2}y=pa+\tau a_1+\mu w'\).
Write \(N=c+2\) and \(a'=pa+\tau a_1\in A\).
If \(\theta\varphi^j(y)=0\) for all \(j\ge0\), the same
evaluations kill \(a'\), whence \(a'\in\mu A\) and
\(y\in\mu B^+\). With \(j\ge1\) the same proof gives
\(y\in\tau B^+\). In particular, if
\(\varphi(y)\in\varphi(\tau)B^+=\mu B^+\), then
\(\theta\varphi^j(y)=0\) for \(j\ge1\) and \(y\in\tau B^+\).
If \(\varphi(y)=\mu^h c\) with \(h\ge1\) and \(c\in B^+\),
write \(y=\tau y_1\) by this implication. Cancellation gives
\(\varphi(y_1)=\mu^{h-1}c\). Induction yields
\(y\in\tau^hB^+\). For the single evaluation, \(p^Ny=a'+\mu w'\)
gives \(\theta(a')=0\) if \(\theta(y)=0\). Since
\(\ker(\theta|_A)=\xi A\) and \(\mu=\xi\tau\), this proves
\(\ker(\theta|_{B^+})=\xi B^+=E B^+\).

For comparison with the usual period ring, take the
tilt element \(p^\flat\) from compatible \(p\)-power roots
of \(p\), with \(v_\flat(p^\flat)=1\). Since
\([p^\flat]-p\) also generates \(\ker\theta\)
\cite[Proposition~10.2]{BergerPadic}, one has
\(A[\xi/p]=A[q]\) for \(q=[p^\flat]/p\).
In the Witt coordinates of \(A[1/p]\), Berger's valuation at
\(r_0=(p-1)/p\) is
\[
 V\!\left(\sum_n p^n[x_n],r_0\right)=\inf_n\bigl(n+v_\flat(x_n)\bigr).
\]
This is the valuation defining \(B_{\max}^+\) in \cite[\S13]{BergerPadic}.
Its unit ball is the indicated algebra. Indeed,
a term \(p^{-j}[x]\) of negative index satisfies the
unit-ball condition exactly when \(v_\flat(x)\ge j\),
in which case it belongs to \(([p^\flat]/p)^jA\).
The finitely many negative indices and the integral
tail give the assertion. Completing the unit ball gives \(A_{\max}\)
and inverting \(p\) gives Berger's \(B^+_{\max}=B^+\).
In particular,
\[
 p^mA_{\max}=\{x\in B^+:V(x,r_0)\ge m\}\qquad(m\in\mathbf Z).
\]
For \(P=X\) or \(X-b_0\) from Section~\ref{sec:period-period-contraction},
the reduction of \(P^n\) modulo \(p\) is monic for every \(n\ge1\).
Thus \(0\le nV(P,r_0)<1\) for every \(n\), so \(V(P,r_0)=0\).
The equalities \(\xi=pX\) and \(\tau^{p-1}=p(X-b_0)\) give
\[
 V(E,r_0)=V(\xi,r_0)=1,\qquad
 V(\tau,r_0)=\frac1{p-1},\qquad V(\mu,r_0)=\frac p{p-1}.
\]
Here \(E/\xi\) and its inverse belong to \(A\).
We use the embedding \(B^+\hookrightarrow B_{\rm dR}^+\) of
\cite[Theorem~13.1]{BergerPadic} for the de Rham tests below.
Also \(t=\log(1+\mu)\) and \(t/\mu\) is a unit of \(B^+\).
In fact \(\mu=pX\tau\), so the logarithm quotient
converges and its nonconstant terms lie in \(pA_{\max}\),
using \(j-1-v_p(j)\ge1\) for \(j\ge2\) except \(p=j=2\).
For \(p=2\) one has \(\xi=2+\tau=2X\), whence
\(\mu/2=X\tau=2X(X-1)\in2A_{\max}\).
Hence \(B^+[1/\mu]=B^+[1/t]=B_{\max}\).

Let \(a\ge0\) be an integer. If \(y\in B_{\max}\) satisfies
\(\varphi^j y\in t^{-a}B_{\rm dR}^+\) for every \(j\ge0\),
then \(\mu^a y\in B^+\).
Write \(y=b/\mu^k\) with \(b\in B^+\).
For \(j\ge0\),
\[
 \varphi^j(\mu)=\mu\sum_{i=0}^{p^j-1}(1+\mu)^i,
 \qquad \theta\!\left(\sum_{i=0}^{p^j-1}(1+\mu)^i\right)=p^j .
\]
Thus every \(\varphi^j(\mu)\) has de Rham order one. For either
\(D=A\) or \(D=B^+\) this proves the exact test, for \(n\ge0\),
\[
 \{c\in D:\varphi^j(c)\in\mathrm{Fil}^nB_{\mathrm dR}^+
              \text{ for all }j\ge0\}=\mu^nD .
\]
Indeed the common-kernel calculation in \(D\) gives the first
factor of \(\mu\); after division all orders decrease by one,
and induction gives the remaining factors. If \(k>a\), apply
this test to \(b\) with \(n=k-a\) to obtain \(\mu^a y\in B^+\).

For a crystalline lattice on this chart choose
an invariant matrix \(M\in\mathrm{GL}_d(B_{\max})\)
whose columns are the coordinates of a \(K_0\)-basis of
\(D_{\rm cris}(T[1/p])\) in a \(\mathbf Z_p\)-basis of \(T\),
where \(K_0=W(k)[1/p]\). Crystalline comparison makes \(M\)
invertible and gives \(\rho(g)g(M)=M\) for \(g\in G_K\)
\cite[\S\S4, 13--14]{BergerPadic}.
For \(j\ge0\) put \(C_j=M^{-1}\varphi^j(M)\).
Both \(C_j\) and \(C_j^{-1}\) have entries fixed by \(G_K\), since
Frobenius commutes with \(G_K\) and fixes the entries of \(\rho(g)\).
Thus \(C_j\in\mathrm{GL}_d(K_0)\) by
\cite[\S13, corollary to Theorem~13.1]{BergerPadic}.
The columns of \(M\) form a \(K\)-basis of \(D_{\rm dR}\) by
crystalline comparison.
Choose \(Y\in\mathrm{GL}_d(K)\) so that the columns of
\(MY\) form a filtration-adapted \(K\)-basis \(d_j\) of \(D_{\rm dR}\).
Write the rescaled matrix as
\[
   M=W\,\operatorname{diag}(t^{k_j})\,Y^{-1}
   \quad\text{in }B_{\rm dR},
   \qquad W=(t^{-k_j}d_j)_j\in\operatorname{Mat}_d(B_{\rm dR}^+).
\]
Here \(k_j\) is the filtration index
\(d_j\in\mathrm{Fil}^{k_j}\setminus\mathrm{Fil}^{k_j+1}\).
The Tate twist \(\mathbf Q_p(1)\) has index \(-1\) in this convention.

For every finite \(K/\mathbf Q_p\) one has
\[
   \mathbf C_p^{G_K}=K,\qquad
   \mathbf C_p(n)^{G_K}=0\quad(n\in\mathbf Z\setminus\{0\}).
\]
The first assertion is Tate's Theorem~1. To apply his Theorem~2,
put \(K_a=K(\mu_{p^a})\) and \(K_\infty=\bigcup_a K_a\).
Since \(\mathbf Q_p(\mu_{p^a})/\mathbf Q_p\) is totally ramified,
the absolute residue degrees of the \(K_a\) are bounded by
\([K:\mathbf Q_p]\). Choose \(L=K_a\) where they have stabilized
and \(\chi_{\rm cyc}(G_L)\) is torsion-free, taking \(a\ge2\)
when \(p=2\). The latter image is open in \(\mathbf Z_p^\times\):
each finite \(L/\mathbf Q_p\) gives a subgroup of finite index.
For every \(n\ne0\), the kernels of \(\chi_{\rm cyc}^n\) and
\(\chi_{\rm cyc}\) agree on \(G_L\). Their common fixed field is
\(K_\infty\), totally ramified over \(L\) with group \(\mathbf Z_p\).
Theorem~2 now gives
\(\mathbf C_p(n)^{G_K}\subseteq\mathbf C_p(n)^{G_L}=0\)
\cite[\S3.3, Theorems~1 and~2, p.~176]{TatePDivisible}.

The columns \(w_j=t^{-k_j}d_j\) of \(W\) are integral.
Their residues are independent over \(\mathbf C_p\). Indeed,
normalize one nonzero coefficient in a relation of least support to
\(a_i=1\). Since \(g(w_j)=\chi_{\rm cyc}(g)^{-k_j}w_j\) for the
diagonal action of \(G_K\), minimality gives
\(g(a_j)=\chi_{\rm cyc}(g)^{k_j-k_i}a_j\).
Thus \(a_j\) belongs to \(\mathbf C_p(k_i-k_j)^{G_K}\), so
\(k_j=k_i\) and \(a_j\in K\) for every nonzero coefficient. This would
give \(\sum_j a_jd_j\in\mathrm{Fil}^{k_i+1}\), contrary
to independence in \(\operatorname{gr}^{k_i}D_{\rm dR}\).
Thus \(\det W\) has nonzero residue and \(W^{-1}\) is integral.
Since \(C_j^{-1}Y\) is invertible over \(B_{\rm dR}^+\),
the identity
\[
 \varphi^j(M)C_j^{-1}Y=MY=W\operatorname{diag}(t^{k_i})
\]
gives elementary divisors \(t^{k_i}\) for each \(\varphi^j(M)\).
Thus for jumps in
\([-a,b]\) each \(\varphi^j M\) has entries in
\(t^{-a}B_{\rm dR}^+\) and its inverse has entries in
\(t^{-b}B_{\rm dR}^+\).
The bounded-ring criterion and \(E\tau\sim\mu\)
give the matrices
\[
 U=(E\tau)^a M,\qquad W'=(E\tau)^b M^{-1}
     \quad\text{in }\operatorname{Mat}_d(B^+).
\]
They satisfy \(UW'=(E\tau)^{a+b}I\).

The entries of \(C_1\) and its inverse belong to
\(K_0\subset S_B\).
Also \(\varphi(E)\) is a unit in that bicommutant:
\(E=z^e u+p\) with \(u\) a chart unit, and \(E=pX'\)
in \(A_{\max}\). Thus \(z^e\in pA_{\max}\)
and, more precisely,
\[
 \varphi(E)/p=1+p^{p-1}(z^e/p)^p\varphi(u),\qquad
 z^e/p=(E/p-1)/u\in S[E/p].
\]
The geometric series for \((\varphi(E)/p)^{-1}\) therefore
converges in \(S_B\), where \(p\) is invertible.
Writing \(v_0=\mu/(E\tau)\in A^\times\), we get
\[
 \varphi(U)=v_0^a U\bigl(\varphi(E)^a C_1\bigr)
\]
with an invertible last factor over the bicommutant.

For the Galois equation, \(E/gE\) is a unit
since both elements generate \(\ker\theta\).
The ideal \(\mu A\) is independent of the chosen compatible
primitive roots by \cite[Lemma~3.23]{BMSIntegral}. Applying
\(\varphi^{-1}\) shows that \(g\tau A=\tau A\), so
\(\tau/g\tau\in A^\times\).
The integral representation matrix, scaled by
the unit \((E\tau/g(E\tau))^a\), takes \(g(U)\)
back to \(U\).
Finally the divisibility condition under coordinate
retraction is proved in
Section~\ref{sec:period-hodge-retraction}.
This gives the period-matrix hypotheses of
the integral basis for each lattice in the
modular cup product, on the crystallizing
extension in the common local class.

\subsection{Divisibility under the coordinate retraction}
\label{sec:period-hodge-retraction}

\sectionopening{In this subsection, we prove that the retraction of
Section~\ref{sec:period-chart-coordinates} respects the Hodge
filtration relevant to the equation \(Uy=E^h x\). Approximation of the
change-of-basis coefficients in the bounded ring reduces the assertion
to divisibility by a power of \(E\).}

The corresponding period matrix of height \(h\)
over \(B^+\) is \(U=(E\tau)^aM\).

Use the usual embedding of \(K\) in \(B_{\rm dR}^+\), on which
\(\theta\) is the identity. Each \(u\in K\) is approximable at
every \(\xi\)-adic precision by elements of \(S_B\). First,
\(\theta(S_B)=K\). The chart constants and the
image of its uniformizer generate \(K\), and \(p\) is
invertible in \(S_B\). Conversely, the uncompleted ring
\(S[E/p][1/p]\) has image in \(K\), since \(\theta(E)=0\)
and series in the uniformizer converge in \(K\).
Continuity of \(\theta\) in the bounded-ring \(p\)-adic
topology and completeness of \(K\) give the assertion
after completion.
Start with \(b\in S_B\) having residue \(u\).
For a polynomial \(P\in\mathbb Q_p[X]\) having a
simple root \(u\), lift \(P'(u)^{-1}\) to
\(c\in S_B\) by this same residue map.
If \(b\equiv u\bmod\xi^j\), Taylor expansion gives
\(b-cP(b)\equiv u\bmod\xi^{j+1}\).
This yields the required approximations.

For \(r\in A_{\max}\) write \(r=a+(\xi/p)b\), with \(a\in A\)
and \(b\in A_{\max}\), by polynomial division at zero in
Section~\ref{sec:period-period-contraction}. If \(\theta(r)=0\),
then \(a\in\xi A\), so \(r\in(\xi/p)A_{\max}=(E/p)A_{\max}\).
Since \(E/p\) has de Rham order one, iteration gives
\[
 A_{\max}\cap\mathrm{Fil}^n B_{\rm dR}^+=(E/p)^n A_{\max}
                  \qquad(n\ge0).
\]
After inverting \(p\), the corresponding identities are
\[
 B^+\cap\mathrm{Fil}^n B_{\rm dR}^+=E^nB^+,
 \qquad S_B\cap\mathrm{Fil}^n B_{\rm dR}^+=E^n S_B .
\]
The first identity follows by localization of the integral one.
For the second, apply the \(S_B\)-linear retraction
\(\widehat F_0\) to \(s=E^n y\) with \(s\in S_B\).
For a shifted divisor \(E_j=\varphi^j E\) the integral test is
\[
 a\in E_j^nA\quad\Longleftrightarrow\quad
       \varphi^{-j}(a)\in\mathrm{Fil}^n B_{\rm dR}^+
       \qquad(a\in A,\ j\ge0).
\]
Here one iterates \(\ker(\theta|_A)=E A\) on
\(\varphi^{-j}(a)\) and transports the quotient by \(\varphi^j\) in \(A\).
For \(a\in S\), applying \(F_0\) to an integral quotient gives
\(E_j^n A\cap S=E_j^nS\).

For \(y_i\in B^+\) and \(\beta_i\in K\), let
\(\widehat F\) be the extension of any \(S\)-linear
\(F:A\to A\) constructed in
Section~\ref{sec:period-chart-coordinates}.
If \(\sum_i\beta_i y_i\) is divisible by \(\xi^j\)
in \(B_{\rm dR}^+\), approximate each \(\beta_i\)
modulo \(\xi^j\) by \(b_i\in S_B\).
The first intersection formula then gives
\(\sum_i b_i y_i\in E^jB^+\).
Apply \(\widehat F\), which commutes with \(S_B\)
and multiplication by \(E^j\), and replace \(b_i\)
by \(\beta_i\) again modulo \(\xi^j\).
It follows that
\[
 \operatorname{ord}_\xi\!\left(\sum_i\beta_i\widehat F(y_i)\right)
      \ge \operatorname{ord}_\xi\!\left(\sum_i\beta_i y_i\right).
\]
If the latter sum is zero, this says that the former
is also zero, by separatedness.

Since \(Y^{-1}\in\mathrm{GL}_d(K)\), each of its rows has
coefficients in the field just approximated. There are
\(s_j\in B_{\rm dR}^+\) such that
\((Y^{-1}\widehat F(y))_j=s_j(Y^{-1}y)_j\).
Since the two diagonal factors commute, we obtain
\[
 U\widehat F(y)=W\,\operatorname{diag}(s_j)\,W^{-1}Uy .
\]
Consequently divisibility of \(Uy\) by \(E^h\)
implies the same for \(U\widehat F(y)\), first
in \(B_{\rm dR}^+\) and then in \(B^+\) as above.
Section~\ref{sec:period-period-saturation} applies this retraction
to the bounded-ring solutions.

\subsection{Divisibility of the period matrix against integral covectors}
\label{sec:period-period-saturation}

\sectionopening{In this subsection, we use the denominator bound of
Section~\ref{sec:period-period-contraction} and the retraction of
Section~\ref{sec:period-hodge-retraction} to extend divisibility from
the basis solutions to all bounded-ring solutions.
Section~\ref{sec:local-integral-period-frame} uses this assertion to
prove full rank and then to obtain the integral span.}

Work on the same chart and take the symmetric height bound
\(b\). Put \(h=2b\), \(A=A_{\inf}\), and denote the invariant
matrix \(M\) of Section~\ref{sec:period-period-realization} by
\(P\), so that \(U=(E\tau)^bP\).
For this same \(P\) put
\[
 W=(E\tau)^bP^{-1},\qquad
 v=(\xi/E)^b,\qquad F=\varphi(E)^bP^{-1}\varphi(P).
\]
The comparison in Section~\ref{sec:period-period-realization} gives
\(W\in\operatorname{Mat}_d(B^+)\), \(v\in A^\times\),
\(F\in\mathrm{GL}_d(S_B)\), and
\[
 UW=WU=(E\tau)^h I,\qquad \varphi(U)=vUF.
\]
Write \(\widehat f\) for the extension of an \(S\)-linear
\(f:A\to A\) to \(B^+\). For any solution
\(Uy=E^h x\), \(x\in A^d\), the evaluation divisibility in
Section~\ref{sec:period-hodge-retraction} gives
\[
                         U\widehat f(y)=E^h w.
\]
A single power \(p^N\), depending on the chart and matrix,
puts every such \(w\) in \(A^d\).

Choose \(c_W,c_U\ge0\) clearing the entries of \(W,U\) into
\(A_{\max}\), respectively. Cancelling \(E^h\) in
\((E\tau)^h y=WUy=E^h Wx\) gives \(\tau^h y=Wx\).
The single-factor estimates of
Section~\ref{sec:period-period-contraction} apply iteratively to
\(\tau\) and to \(E\), which is a unit multiple of \(\xi\).
The first gives \(p^{c_W+h}y\in A_{\max}^d\).
Since \(\widehat f\) preserves \(A_{\max}\) and
\(E^h w=U\widehat f(y)\), the second gives
\[
          p^{N_2}w\in A_{\max}^d,\qquad N_2=c_W+c_U+2h
\]
for every solution and every map.

For Frobenius contraction set
\(x'=\varphi^{-1}(vE^h x)\in A^d\).
Apply \(\varphi\) to \(UW x'=(E\tau)^h x'\) and cancel
\(vU\) to get
\[
          F\varphi(Wx')=\varphi(E\tau)^h y.
\]
The inverse of \(F\) is defined over \(B^+\), so
\(\varphi(Wx')\) is divisible there by \(\varphi(\tau)^h\).
The cyclotomic kernel argument of
Section~\ref{sec:period-period-realization} gives
\(y'=\tau^{-h}Wx'\in(B^+)^d\), with \(Uy'=E^h x'\).
Put \(f'=\varphi^{-1} f\varphi\) on \(A\).
It is \(S\)-linear because \(\varphi(S)\subset S\).
Write \(U\widehat{f'}(y')=E^h w'\).
The equality \(\varphi\widehat{f'}=\widehat f\varphi\)
holds on \(A\) and then on \(A[E/p]\) by clearing a power
of \(p\) from each input. Both sides preserve \(A_{\max}\)
and its \(p\)-powers, giving equality on the completion and
on \(B^+\). The \(S_B\)-linearity permits
\(\widehat f\) to pass through \(F\) and \(\varphi(E)^h\).
Since \(F\varphi(y')=\varphi(E)^h y\), we get
\[
 \varphi(E)^h\varphi(w')
   =vUF\widehat f(\varphi(y'))
   =v\varphi(E)^h U\widehat f(y),
 \qquad \varphi(w')=vE^h w .
\]
The entries of \(p^{N_2} w\) thus belong to a single
subset of \(A_{\max}\) closed under the descent relation
\(\xi^h a=u\varphi(a')\), with the fixed unit
\(u=(\xi/E)^h v^{-1}\). The contraction of
Section~\ref{sec:period-period-contraction} yields
\(N=N_2+h+1=c_W+c_U+3h+1\), independent of
\(x,y,f\).

Let \(r=\operatorname{rank}_S\mathcal P\le d\) and write
\(Q\in\operatorname{Mat}_{d\times r}(A)\) for the matrix of an
\(S\)-basis of the integral vectors admitting a solution over
\(S_B\), as in Section~\ref{sec:period-chart-coordinates}. For
every row \(l\in A^{1\times d}\) and \(\delta\in A\) with
\(p\nmid\delta\), use entrywise divisibility for \(lQ\).
If \(x\in A^d\) admits a solution \(y\in(B^+)^d\), then
\[
         \delta\mid lQ\quad\Longrightarrow\quad\delta\mid lx
                         \qquad(Uy=E^h x).
 \tag{1}\label{eq:coordinate-retraction-divisibility}
\]
Choose the topological chart basis \((b_i)_i\).
For \(f=\operatorname{coord}_i\), the entries of
\(\widehat f(y)\) are in \(S_B\).

The uniform exponent gives \(x_i=p^N w_i\in A^d\) and
\[
 U\bigl(p^N\widehat{\operatorname{coord}}_i(y)\bigr)=E^h x_i .
\]
Thus \(x_i=Qs_i\) with \(s_i\in S^r\).
Take \(k=c_U+2h\) to clear the entries of \(U\) and divide
by \(E^h\) in \(A_{\max}\). For any \(m\)
choose one finite subset \(I_m\) of the basis with
\[
 y-\sum_{i\in I_m}b_i\widehat{\operatorname{coord}}_i(y)
             \in p^{m+k}A_{\max}^d.
\]
This is the simultaneous finite approximation of the
finitely many entries.

Set
\[
 t=p^N x-\sum_{i\in I_m}b_i Qs_i\in A^d .
\]
Applying \(p^NU\) to the approximation gives
\(E^h t\in p^{m+2h}A_{\max}^d\). Since \(E\mid\mu\),
the denominator estimate applied to \(p^{-m-2h}t\)
gives \(t\in p^m A_{\max}^d\).
The identity \(p^m A_{\max}\cap A=(p,\xi)^m\)
puts \(lt\) in \((p,\xi)^m\). If \(lQ=\delta g\), then the
multiples \(\delta\sum_{i\in I_m} b_i gs_i\) approximate
\(p^N lx\) to every \((p,\xi)\)-adic order.
By closedness of \(\delta A\) from
Section~\ref{sec:period-chart-coordinates},
\(p^N lx\in\delta A\). Since \(p\) is prime in \(A\)
and \(p\nmid\delta\), repeatedly cancel \(p\) in this
divisibility. This proves \eqref{eq:coordinate-retraction-divisibility}
for every bounded-ring solution on the given chart.

\section{Integral crystalline bases and filtration orders}
\label{sec:crystalline-basis-appendix}

\sectionopening{In this appendix, we use the comparison matrices of
Appendix~\ref{sec:crystalline-coordinate-appendix} to construct an
integral period basis and describe its filtration under the divided
Galois action. The resulting orders are used in the cup calculation
of Appendix~\ref{sec:tame-cup-appendix} and the lifts of
Appendix~\ref{sec:crystalline-lifting-appendix}.}

\subsection{The integral period basis}
\label{sec:local-integral-period-frame}

\sectionopening{In this subsection, we construct an integral basis from the solution module of
Section~\ref{sec:period-chart-coordinates}. The fixed-vector
approximation of Section~\ref{sec:period-chart-approximation} rules
out a kernel modulo \(p\). We then obtain the two partial inverses and
the equations for Frobenius and Galois.}

Use the \(S\)-basis of the solution module from
Section~\ref{sec:period-chart-coordinates}, with its columns
in \(Q\) and solutions in \(Z\). Initially the number of
columns is at most \(d=\operatorname{rank}_{\mathbf Z_p}T\).
For a free \(O\)-lattice this is
\(\operatorname{rank}_{O}(T)[\operatorname{Frac}(O):\mathbf Q_p]\),
as in Section~\ref{sec:period-chart-coordinates}. We take determinants
on the underlying \(\mathbf Z_p\)-lattice. Multiplication by
\(o\in O\) has a matrix \(H_o\) over \(\mathbf Z_p\) and induces a
matrix \(C_o\) over \(K_0\) on the invariant crystalline periods.
Thus \(H_o U=U C_o\), and \(Uy=E^h x\) gives
\(U(C_o y)=E^h H_o x\). Since \(K_0\subset S_B\), the integral
solution module is stable under \(O\).
We have \(UZ=E^h Q\) with \(Z\) over \(S_B\). The divisibility against covectors
of Section~\ref{sec:period-period-saturation} holds before
asserting that \(Q\) is square. If the rank of \(Q\) were
less than \(d\), clearing denominators in its left kernel over
\(\operatorname{Frac}A\) would give \(0\ne l\in A^{1\times d}\) with
\(lQ=0\). Apply \eqref{eq:coordinate-retraction-divisibility} with
\(\delta=z^k\), \(y=We_j\), \(x=\tau^h e_j\). Since \(p\nmid z\), this
gives \(\tau^h l_j\in z^k A\) for every \(j\) and \(k\ge1\).
The intersection of these ideals is zero by \((p,z)\)-adic
separatedness. Since \(A\) is a domain, this contradicts \(l\ne0\).
The covector argument proves the following rank statement.

\begin{lemma}[Full rank of the integral solution basis]\label{lem:period-basis-full-rank}
The solution matrices \(Q,Z\) are square of size
\(d=\operatorname{rank}_{\mathbf Z_p}T\), and \(\det Q\ne0\).
\end{lemma}

If \(Qs=px\) for \(s\in S^d\) and \(x\in A^d\), then
\(U(Zs/p)=E^h x\) and \(Zs/p\in S_B^d\).
It follows that \(x=Qs'\) for \(s'\in S^d\), and
\(s=ps'\) by injectivity of \(Q\). Thus the chart-coordinate
covectors of \(Q\) separate vectors over \(S/p\): if all
coordinate covectors kill a vector modulo \(p\), the
\(z\)-adic coordinate expansion implies \(Qs=0\bmod p\),
hence \(s=0\bmod p\). Over the fraction field of
\(S/p=k[[z]]\), these covectors therefore span the full dual
space. Choose \(d\) rows indexed by pairs \((i_\nu,j_\nu)\) of a
coordinate map and a row of \(Q\), and put
\[
 N_{\nu t}=F_{i_\nu}(Q_{j_\nu t}),\qquad
 (Nu)_\nu=F_{i_\nu}\bigl((Qu)_{j_\nu}\bigr)\quad(u\in S^d),
 \qquad \overline{\det N}\ne0 .
\]

The equation \(\varphi U=vUF\) and the
invertibility of \(\varphi(E)\) in \(S_B\) give
\[
 U\left((E/\varphi E)^h F\varphi Z\right)
                       =E^h v^{-1}\varphi Q .
\]
By the solution basis, \(QB=v^{-1}\varphi Q\) for
some \(B\in M_d(S)\). For nonzero \(a\in A\), the order
\(\operatorname{ord}_p(a)=\max\{j:a\in p^j A\}\) is finite by
separatedness and additive on products since \(A/p\) is a domain.
Frobenius preserves this order and \(v\) is a unit. Taking
determinants gives \(\operatorname{ord}_p(\det B)=0\).
Throughout, \(S/p\hookrightarrow A/p\) by the chart retraction.

Suppose \(\det Q=0\bmod p\). Over the valuation ring
\(A/p\) choose a nonzero kernel vector \(a\) of minimal
support and divide by a coordinate of smallest valuation
to make that coordinate \(1\). If \(g\) fixes \(S\)
pointwise, it belongs to \(G_K\) since
\(K=W(k)[1/p](\theta(z))\), and fixes \(Z\) by continuity on
\(S_B=\overline{S[1/p]}\). Since \(gE=E\), the Galois equation
of Section~\ref{sec:period-period-realization} gives
\[
 U=R_g\,g(U),\qquad
 R_g=(\tau/g\tau)^b\rho(g)\in\mathrm{GL}_d(A).
\]
Applying \(g\) to \(UZ=E^hQ\) and cancelling \(E^h\)
therefore gives \(Q=R_g\,g(Q)\).
Thus \(g(a)\) is again in the kernel. Subtracting \(a\)
gives a kernel vector of strictly smaller support unless
\(g(a)=a\). The normalized vector \(a\) is consequently
fixed by the pointwise fixer of \(S\).

Put \(C_n=B\varphi(B)\cdots\varphi^{n-1}(B)\) and
\(D_0=\overline{\det B}\,\overline{\det N}\).
The Frobenius equation and the determinant product imply
\[
 \bar Q\bar C_n a^{p^n}=0,\qquad
 D_0\,\overline{\det C_n}\mid D_0^{p^n}.
\]
Indeed \(\varphi^n Q=(\prod_{j<n}\varphi^j v)QC_n\), and
\(\overline{\det C_n}=(\overline{\det B})^{(p^n-1)/(p-1)}\) with
\((p^n-1)/(p-1)\le p^n-1\).

Let \(\pi_j:A/p\to\mathcal O_{\mathbb C_p}/p\) be reduction of
\(\theta_j\). Thus \(\pi_{m+n}(x^{p^n})=\pi_m(x)\).
Fix \(m\). By Section~\ref{sec:period-chart-approximation}, each
coordinate of \(a\) admits the required approximation at every
sufficiently large \(j\). Choose a common \(j=m+n\) and
\(s\in S^d\) such that \(\pi_j(\bar s)=\pi_m(a)\). The first
identity above gives \(\pi_j(\overline{QC_ns})=0\). Applying the
coordinate divisibility of Section~\ref{sec:period-chart-approximation}
with \(\alpha=p\) to its \(j_\nu\)-th entry gives
\(\pi_j(\overline{NC_ns})_\nu=0\), since \(C_ns\in S^d\).
One coordinate of \(\pi_j(\bar s)\) is \(1\). Multiplication by the
adjugate over \(\mathcal O_{\mathbb C_p}/p\) gives
\(\pi_j(\overline{\det(NC_n)})=0\). Since
\[
 \overline{\det(NC_n)}\mid D_0\,\overline{\det C_n}
                          \mid D_0^{p^n},
\]
we obtain \(\pi_m(D_0)=\pi_j(D_0^{p^n})=0\). All components vanish,
contradicting \(D_0\ne0\) in the tilt. Hence \(p\nmid\det Q\):
the reduction of \(Q\) has rank \(d\) over \(\operatorname{Frac}(A/p)\).
In particular \(Qt\in pA^d\) implies \(t\in pA^d\) for \(t\in A^d\),
by reduction and the adjugate. Iteration gives
\[
                 QA^d\cap p^s A^d=p^s QA^d\qquad(s\ge0).
\]
Thus for each \(s\ge1\), reduction of the inclusion of the integral
column span identifies \(QA^d/p^sQA^d\) with its image in
\((A/p^s)^d\). In particular,
\[
 QA^d=\bigcap_{s\ge1}(QA^d+p^sA^d).
\]
Indeed the coefficient vector of a residue in this image is
unique modulo \(p^s\) since \(Q\) is injective modulo \(p^s\).
For a vector in the intersection the coefficient vectors are
therefore compatible in \(s\), and lift in the \(p\)-adically
complete ring \(A\). In \(W(\mathbb C_p^\flat)\) the residue
\(\overline{\det Q}\) is a nonzero field element, so \(Q\) is
invertible over this Witt ring and over all its quotients by
\(p^s\). The two assertions use respectively injectivity and
invertibility of the reduced determinant over the indicated rings.

Change the integral basis by \(H\in\mathrm{GL}_d(\mathbf Z_p)\)
and the crystalline period basis by \(C\in\mathrm{GL}_d(K_0)\).
The period matrices and the solution module become
\[
 P_1=H^{-1}PC,\qquad U_1=H^{-1}UC,\qquad
 \mathcal P_1=H^{-1}\mathcal P .
\]
Indeed \(U_1y=E^h x\) is equivalent to \(U(Cy)=E^h Hx\), and
both \(C\) and \(C^{-1}\) preserve \(S_B^d\).
For any \(S\)-basis \(Q_1\) of \(\mathcal P_1\) and its solution
matrix \(Z_1\) there is therefore a \(V\in\mathrm{GL}_d(S)\) with
\[
 Q_1=H^{-1}Q V,\qquad Z_1=C^{-1}Z V .
\]
The formula for \(Z_1\) follows from the injectivity of \(U_1\).
Thus the determinant of the integral solution basis changes by a
unit of \(A\), including when \(O\) is ramified.

For \(\sigma\in G_{\mathbf Q_p}\), transport the representation
to \(G_{\sigma K}\) and the chart to \((\sigma S,\sigma E)\).
Using the same \(\tau\), put \(u=(\tau/\sigma\tau)^b\in A^\times\).
The transported period matrix is \(u\sigma U\), and its solution
module is \(u\sigma\mathcal P\), since \(\sigma S_B\) is the
bounded scalar ring of the new chart. Thus the transported
solution basis \(u\sigma Q\) has determinant
\(u^d\sigma(\det Q)\), again nonzero modulo \(p\).

For \(Uy=E^h x\), \(x\in A^d\), \(y\in(B^+)^d\), apply
\eqref{eq:coordinate-retraction-divisibility} with \(\delta=\det Q\)
and the rows of \(\operatorname{adj}(Q)\). This gives
\(\operatorname{adj}(Q)x=(\det Q)t\) for \(t\in A^d\).
Multiplication by \(Q\) and cancellation give \(x=Qt\).
Conversely \(U(Zt)=E^hQt\), so the integral solution vectors over
\(B^+\) form precisely
\[
 QA^d=A\otimes_S QS^d
    =\{x\in A^d:\operatorname{adj}(Q)x\in(\det Q)A^d\}.
\]
The tensor product embeds via its independent basis columns.
The last equality follows by multiplying the adjugate equation by
\(Q\) and cancelling \(\det Q\). The principal ideal
\((\det Q)A\) is closed in the \((p,z)\)-adic topology by
Section~\ref{sec:period-chart-coordinates}, since \(p\nmid\det Q\).
Its inverse image under the continuous map \(x\mapsto
\operatorname{adj}(Q)x\) is therefore closed in \(A^d\).
With \(y=We_j\) we get
\(QQ'=\tau^h I\) for \(Q'\in M_d(A)\).

Put \(v_0=\mu/(E\tau)\in A^\times\). To construct the
second partial inverse, take \(y\in B^+\) and \(s,k\ge0\).
For some \(c\ge0\), \(a\in A\), and \(b\in A_{\max}\),
\[
                 p^c y=a+p^{c+k}\mu^s b .
\]
Indeed choose \(c_0\) with \(p^{c_0}y\in A_{\max}\).
Approximate by a polynomial in \(\xi/p\) modulo
\(p^{c_0+k+2s}A_{\max}\), then multiply by \(p^{c-c_0}\)
to put all of its evaluated terms in \(A\).
The resulting error is \(p^{c+k+2s}b_0\).
For \(X=\xi/p\), polynomial division in
Section~\ref{sec:period-period-contraction} gives
\(A_{\max}=A+(X-1)A_{\max}=A+XA_{\max}\)
by division at \(1\) and \(0\), respectively.
Iterating and using \(p(X-1)\in\tau A\), \(pX=\xi\), gives
\[
 p^s b_0=a_1+\tau^s b_1,\qquad
 p^s b_1=a_2+\xi^s b,\qquad a_1,a_2\in A .
\]
The remaining error after absorbing integral terms is
\(p^{c+k}\tau^s\xi^s b=p^{c+k}\mu^s b\).

Apply this entrywise to \(U\) with \(k=0,s=h\).
Take one common \(c\) and obtain \(p^cU=A_0+p^c\mu^h V_0\)
with integral \(A_0\) and \(V_0\) over \(A_{\max}\).
Then
\[
 A_0=p^c U(I-v_0^h W V_0),\qquad
 p^c\mu^h V_0=p^c(Ev_0)^h QQ' V_0 .
\]
Taking \(y=p^c E^h(I-v_0^h WV_0)\) gives \(Uy=E^h A_0\).
Hence \(A_0=QA_1\) for \(A_1\in M_d(A)\), and
\[
  Y=p^{-c}A_1+(Ev_0)^hQ'V_0\in M_d(B^+),\qquad U=QY .
\]
Since \(B^+\) is a domain containing \(A\), cancellation of \(Q\) in
\(QYZ=E^hQ\) gives \(YZ=E^h I\).

Finally let \(g\) fix the chart field, so it fixes the
Witt constants and \(\theta(z)\). If \(U=A_g gU\) is the period equation, then
\(UgZ=(gE)^hA_g gQ\), so
\(UgZ=E^hR_g gQ\) with \(R_g=(gE/E)^hA_g\) integral.
The integral span gives \(R_g gQ=QG_g\) for \(G_g\)
integral. Cancel \(U\) in the resulting equality to
obtain \(gZ=ZG_g\). Altogether the basis satisfies
\[
  QQ'=\tau^h I,\qquad UZ=E^hQ,\qquad YZ=E^h I,\qquad
  \varphi Q=vQB,\qquad gZ=ZG_g .
\]

\subsection{The period basis, divided action, and tensor orders}
\label{sec:v2-crystalline-frame-action}

\sectionopening{In this subsection, we use the de Rham evaluations of the
integral basis in Section~\ref{sec:local-integral-period-frame}
to control division of Galois differences by \(z\tau\).
The retraction preserves filtration orders, and the tensor
estimates give the orders in Section~\ref{sec:local-tame-moves}.}

Keep the chart of Sections~\ref{sec:period-chart-coordinates}
and~\ref{sec:period-period-contraction}, and the original
crystalline \(\mathbf Z_p\)-lattice \(T\). Write
\(A=A_{\inf}\), \(\mathrm{Fil}^j=t^jB_{\mathrm dR}^+\),
\(E_a=\varphi^a E\), and let \(K'\) be the local fixed field
of the chart subgroup. We use the symmetric height bound \(b\)
for the Hodge jumps. For an invariant crystalline period
matrix \(P\), use the matrices of
Section~\ref{sec:local-integral-period-frame} with \(h=2b\), so
\[
 U=(E\tau)^bP,\qquad UZ=E^{2b}Q,\qquad
 QQ'=\tau^{2b}I,\qquad YZ=E^{2b}I.
\]
The columns \(n_i\) of \(\tau^{-b}Q\), in the original
lattice basis, satisfy
\[
 E_a^b\lambda\Bigl(\sum_i c_i\varphi^a n_i\Bigr)
  =(\lambda_{\mathrm{lat}}\varphi^aP)\varphi^a(Z)c
 \qquad(a\ge0).
 \tag{1}\label{eq:period-basis-evaluation}
\]
For the integral tests \(c\) has entries in
\(A_{(p)}=(A\setminus pA)^{-1}A\subset B_{\mathrm dR}\).
The formula holds for every \(c\in B_{\mathrm dR}^d\) after
extending the chart-equivariant covector \(\lambda\) by linearity.
The row \(\lambda_{\mathrm{lat}}\) gives its values on the original
lattice basis. The matrix
identity underlying \eqref{eq:period-basis-evaluation} is
\(E^b n=PZ\).
Write \(\rho(g)\) for the matrix on \(T\) and \(H=G_{K'}\).
Equivariance gives \(g(\lambda_{\mathrm{lat}})=\lambda_{\mathrm{lat}}\rho(g)\).
Since \(\rho(g)\) has \(\mathbf Z_p\)-entries,
\(g(\varphi^aP)=\rho(g)^{-1}\varphi^aP\).
Thus \(\lambda_{\mathrm{lat}}\varphi^aP\) has entries in
\(B_{\mathrm dR}^{H}=K'\). For a crystalline functional both
\(\lambda(T)\) and the entries of \(\varphi^aP\) lie in \(B_{\max}\),
so the invariant products belong to \(B_{\max}^{H}=K_0\subset S_B\)
\cite[\S13, corollary to Theorem~13.1]{BergerPadic}.
For each standard covector \(e_j^{\mathsf T}\), define
\((\lambda_{a,j})_{\mathrm{lat}}=e_j^{\mathsf T}(\varphi^a P)^{-1}\).
Its entries are crystalline periods. Since
\(g((\varphi^aP)^{-1})=(\varphi^aP)^{-1}\rho(g)\),
its row defines a chart-equivariant crystalline functional.
Its product with \(\varphi^aP\) is \(e_j^{\mathsf T}\). For
\(x_a=\sum_i c_i\varphi^a n_i\), put
\(L_a(x_a)=(\lambda_{a,j}(x_a))_j\). Substitution in
\eqref{eq:period-basis-evaluation} and \(YZ=E^{2b}I\) give
\[
 \varphi^a(Z)c=E_a^b L_a(x_a),\qquad
 E_a^{2b}c=\varphi^a(Y)\,E_a^b L_a(x_a).
 \tag{2}\label{eq:period-basis-coefficient-orders}
\]
In particular, for any \(B^+\)-submodule \(J\subset B_{\mathrm dR}\),
\(E_a^b L_a(x_a)\in J^d\) implies \(E_a^{2b}c\in J^d\),
since \(\varphi^a(Y)\) has entries in \(B^+\).

The \(A\)-span \(M\) of the \(n_i\) is exactly the part of
\(A[1/\tau]\otimes T\) specified by the de Rham conditions
\[
     \lambda(\varphi^a x)\in \mathrm{Fil}^0 B_{\rm dR}
     \qquad(a\ge1,\ \lambda\text{ a chart-equivariant form}).
\] In one direction \eqref{eq:period-basis-evaluation} proves the integrality, since
\(E_a\) has de Rham order zero for \(a>0\).
In the other direction \(QQ'=\tau^{2b}I\) first gives a
presentation \(x=\tau^{-m}\sum_i c_i n_i\), \(c_i\in A\).
Take \(m\) least. If \(m>0\), for every \(a>0\) the values
of \(\varphi^a(\tau^m x)\) have order at least one, since
\(\theta\varphi^a(\tau)=0\).

Apply \eqref{eq:period-basis-coefficient-orders} at that level to
\(\varphi^a(c_i)\) and use \(\theta(E_a)\ne0\).
Then \(\theta\varphi^a(c_i)=0\) for all \(a>0\). The
cyclotomic intersection of these kernels is \(\tau A\).
This contradicts the least choice of \(m\).
The columns of \(n=\tau^{-b}Q\) are independent by
\(QQ'=\tau^{2b}I\). Since \(nQ'=\tau^b I\) and \(\tau^b n=Q\) is
integral, they give
\[
 \tau^b(A\otimes T)\subset M\subset \tau^{-b}(A\otimes T).
 \tag{3}\label{eq:period-basis-lattice-inclusions}
\]
The positive Frobenius tests defining \(M\) are independent of the
chosen bound. Hence one may take \(b\) to be the maximum absolute
value of the jumps, with \(b=0\) for the zero lattice.
Every element acting on \(T\) and normalizing the chart subgroup
preserves these tests, since it permutes the equivariant forms and
commutes with Frobenius. This applies to the tame generator because
the tame chart field is Galois over \(\mathbf Q_p\).

Put \(N=\sum_i S n_i\) and invert all \(E_a\) to obtain
\(N_{\rm loc}\) and \(M_{\rm loc}\) from \(N\) and \(M\).
The integral change of basis is
\[
                   E^b\varphi(n)=nB,\qquad B\in M_d(S).
 \tag{4}\label{eq:period-basis-frobenius-matrix}
\]
Indeed the scalar in \(\varphi Q=vQB\) is
\(v=(\varphi(\tau)/(E\tau))^b\) by
Section~\ref{sec:period-period-realization}, so
\(E^b\varphi(\tau)^{-b}\varphi Q=\tau^{-b}QB\). Since \(n\) and
\(\varphi(n)\) are invertible over \(W(\mathbb C_p^\flat)\),
by reduction of \(\det Q\) and \(\tau\) in its residue field,
and \(E\ne0\), we have \(\det B\ne0\). Iterating \eqref{eq:period-basis-frobenius-matrix}
gives, with the integral product \(C_a=B\varphi(B)\cdots\varphi^{a-1}(B)\)
from Section~\ref{sec:local-integral-period-frame},
\[
 \varphi^a(n)=n\mathcal B_a,\qquad
 \mathcal B_a=\left(\prod_{j=0}^{a-1}E_j^b\right)^{-1}C_a
       \in\mathrm{GL}_d(\operatorname{Frac}S),\qquad \mathcal B_0=I .
\]
Define the value field
\[
        \mathcal K_a=K'(S_B,\varphi^{-a}S)\subset B_{\mathrm dR}.
\]
By \eqref{eq:period-basis-evaluation}, forms take values in \(\mathcal K_0\) on \(n_i\);
by the identity
\(\varphi^{-a}(n)=n\varphi^{-a}(\mathcal B_a^{-1})\), they take values
in \(\mathcal K_a\) on \(\varphi^{-a}(N_{\rm loc})\).
The inverse Frobenius here acts on Witt coefficients,
before taking their de Rham images.

\subsubsection{Integral divided operators}
\label{sec:period-basis-integral-operators}

\sectionopening{We use the basis and the partial inverse
\(QQ'=\tau^h I\) constructed in
Section~\ref{sec:local-integral-period-frame} to bound the
divided Galois operators.}

For the scalar integral cancellation, let \(y=z\) or
\(y=\tau\). In each case \(p\nmid y\), \(\theta(y)\ne0\),
and \(\xi-pu\in yA\) for some \(u\in A\).
For \(z\) use \(E=p+z^e u'\) and \(E\sim\xi\).
For \(\tau\) use \(\xi\equiv p\bmod\tau\).
The map \(A[X]/(pX-\xi)\to A/(y,p^j)\) defined by
\(X\mapsto u\) extends to the \(p\)-adic completion.
If \(a\in A\cap y B^+\), clear \(p^m\) from its
quotient by \(y\) in \(B^+\) and apply these maps.

It follows that \(p^m a\in yA+p^j A\) for every \(j\).
The ideal \(yA\) is closed for the \(p\)-adic topology:
since \(\bar y\ne0\) in the tilt domain, divisibility of
\(y c\) by \(p^j\) implies divisibility of \(c\) by \(p^j\).
The approximating quotients therefore converge. The same
regularity cancels \(p^m\) and gives \(a\in yA\).
For \(a\in A\) and \(r\ge0\), contract \(E^r a\in yB^+\)
to \(E^r a=yc\) with \(c\in A\). Since \(E\) is prime and
\(\theta(y)\ne0\),
\(E^r\) divides \(c\). This proves
\[
                   E^r a\in y B^+\ \Longrightarrow\ a\in y A.
 \tag{5}\label{eq:period-basis-integral-cancellation}
\]
For a product \(y_1\cdots y_k\) of copies of \(z,\tau\), suppose
\(E^r a\in y_1\cdots y_k B^+\). Equation
\eqref{eq:period-basis-integral-cancellation} gives \(a=y_1 a_1\)
with \(a_1\in A\). Cancelling \(y_1\) in \(B^+\) gives
\(E^r a_1\in y_2\cdots y_k B^+\). Induction yields
\(a\in y_1\cdots y_k A\), including all \(y=\tau^m z^n\).

Write \(\delta=z\tau\) and let \(H\) be the subgroup fixing
the chart field \(K'\), so \(g\in H\) fixes the Witt
constants and \(\theta(z)\).
For \(g\in H\) the chart action satisfies
\[
             gs-s\in \delta E A\quad(s\in S),\qquad
             gy-y\in \delta B^+\quad(y\in S_B).
 \tag{6}\label{eq:period-basis-divided-action}
\]
For the first inclusion \(gz/z=[v]\) is a Teichm\"uller lift
with \(v^\sharp=1\). For every \(j\ge0\),
\(\theta\varphi^j([v]-1)=(v^\sharp)^{p^j}-1=0\).
The common kernel calculation in
Section~\ref{sec:period-period-realization} gives
\([v]-1\in\mu A\), with \(\mu=\varphi(\tau)\).
Thus \(gz-z\in z\mu A=\delta EA\).
For \(s=\sum_n a_n z^n\) with its Witt coefficients fixed,
\[
 gs-s=(gz-z)\sum_{n\ge1}a_n\sum_{j=0}^{n-1}(gz)^j z^{n-1-j}.
\]
The inner sum lies in \(z^{n-1}A\), so the series converges
in the \((p,z)\)-adic topology of \(A\).

For the second inclusion put \(R=A_{\max}\) and
\(\partial_g=(g-1)/(z\tau)\). The action preserves \(R\), since
\(g\) sends the generator \(\xi\) of \(\ker\theta\) to a unit
multiple of \(\xi\). We have \(\partial_g(S)\subset EA\).
Writing \(E=p+z^e u\), the inclusions \(gz/z-1\in\mu A\)
and \(gu-u\in z\mu A\) give \(gE-E\in z^e\mu A\). Hence
\[
 \partial_g(E/p)\in z^{e-1}(\xi/p)A\subset R,\qquad
 \partial_g(vw)=g(v)\partial_g(w)+\partial_g(v)w.
\]
Consequently \(\partial_g(S[E/p])\subset R\).
For \(y\in S_B\cap R\), write \(y=\sum_{m\ge0}p^m H_m\)
with \(H_m\in A[E/p]\) by \(p\)-adic completeness.
The extended coordinate \(\widehat F_0\) from
Section~\ref{sec:period-chart-coordinates} fixes \(y\) and
sends \(H_m\) to \(h_m\in S[E/p]\). Therefore
\[
 gy-y=z\tau\sum_{m\ge0}p^m\partial_g(h_m)\in z\tau R.
\]
The series on the right converges in \(R\). For arbitrary
\(y\in S_B\), apply this to \(p^m y\in S_B\cap R\) for a
sufficiently large \(m\). This gives \(gy-y\in z\tau B^+\).
For \(x\in N\), \eqref{eq:period-basis-evaluation} gives \(h=E^b\lambda(x)\in S_B\)
for each crystalline \(\lambda\).

Equivariance and \eqref{eq:period-basis-divided-action}
give
\[
 E^{2b}\lambda(gx-x)=
 \frac{E^b}{(gE)^b}\bigl(E^b(gh-h)+(E^b-(gE)^b)h\bigr)
 \in \delta B^+ ,
\]
since \(gE/E\) is an integral unit. Write
\(gx-x=\sum_i c_i n_i\), with \(c_i\in A\) by stability of \(M\).
At \(a=0\), \eqref{eq:period-basis-coefficient-orders} yields
\[
 E^{3b}c=Y E^{2b}L_0(gx-x)\in (\delta B^+)^d.
\]
Equation \eqref{eq:period-basis-integral-cancellation} gives
\(c_i\in\delta A\) for each \(i\), hence
\((gx-x)/\delta\in M\). For \(x=y/D\in N_{\rm loc}\),
each factor \(E_a\) of \(D\) satisfies
\(gE_a/E_a=\varphi^a(gE/E)\in A^\times\). Thus
\(gx\in M_{\rm loc}\), and
\[
 \frac{gx-x}{\delta}
   =D^{-1}\left(\frac{gy-y}{\delta}
                      -\frac{gD-D}{\delta}\,gx\right)\in M_{\rm loc}
\]
by \eqref{eq:period-basis-divided-action} applied to \(D\in S\).
Write \(\pi_a=\theta\varphi^{-a}(z)\). At every \(a\ge0\),
\[
 \theta\varphi^{-a}(\delta)=\pi_a(\zeta_{p^{a+1}}-1)\ne0.
\]
At \(p=2,a=0\) this value is \(-2\pi_0\). Thus
\(\varphi^{-a}(\delta)\) is a unit in \(B_{\mathrm dR}^+\)
at every level. Division by it preserves the de Rham orders
of the difference. The coefficient projector commutes with
this scalar and with \(g\), which fixes the Witt constants.

\subsubsection{Retraction at each order}
\label{sec:period-basis-order-retraction}

\sectionopening{In this subsection, we construct a retraction that preserves the
specified orders and depths after applying the divided differences
of Section~\ref{sec:period-basis-integral-operators}.}

Let \(F:A\to S\) be any \(S\)-linear coordinate map,
\(F_a=\varphi^{-a}F\varphi^a:A\to A\). On \(A\) it
commutes with multiplication by \(\varphi^{-a}(s)\) for
every \(s\in S\), since \(F\) is \(S\)-linear. In
particular \(F_a\) is \(S\)-linear by \(\varphi^a S\subset S\).
Extend it to \(B^+\) by applying it to coefficients of
divided polynomials in \(E/p\) and then taking limits.
Each scalar \(\varphi^{-a}(s)\) belongs to \(A\), so
its commutation survives this extension. The extended map
is also \(S_B\)-linear.
The image of \(S[1/p]\) under \(\theta\) is \(K'\),
generated by the Witt constants and \(\theta(z)\). This
field is complete and \(\theta\) is continuous on \(B^+\),
so \(\theta(S_B)=K'\). Successive Newton
approximation lifts each constant of \(K'\) to \(S_B\)
modulo any specified positive de Rham order. For \(c\in K'\)
take its minimal polynomial \(P\in\mathbf Q_p[X]\). Choose
\(b,d\in S_B\) with \(\theta b=c\) and \(\theta d=1/P'(c)\).
If \(b-c\in\mathrm{Fil}^j\) for \(j\ge1\), Taylor expansion gives
\(b-dP(b)-c\in\mathrm{Fil}^{j+1}\).

Consequently, for
\(u_i\in B^+\), \(c_i\in K'\), replacing \(u_i\) by
\(\widehat F_a(u_i)\) can only increase the order of
\(\sum_i c_i u_i\). For positive precision,
approximate the \(c_i\) by elements of \(S_B\),
use \(B^+\cap\mathrm{Fil}^j=E^jB^+\), and apply
\(S_B\)-linearity. Order zero is automatic.
For the values \(v_i\in\mathcal K_a\), clear a common
nonzero denominator and expand each numerator as a
finite sum of \(c\,b\,\varphi^{-a}(s)\) with
\(c\in K', b\in S_B,s\in S\). Applying the preceding
argument to the resulting single finite sum and using both
scalar actions gives
\[
 \operatorname{ord}_t\Bigl(\sum_i\varphi^{-a}(F(r_i))v_i\Bigr)
  \ \ge\ 
 \operatorname{ord}_t\Bigl(\sum_i\varphi^{-a}(r_i)v_i\Bigr),
                   \qquad r_i\in A.
 \tag{7}\label{eq:period-basis-order-retraction}
\]
The clearing factor shifts both orders equally.
Separatedness includes the case of an exact relation.

Take \(F=\operatorname{coord}_0\), which fixes \(1\). Then, inside
\(B_{\mathrm dR}\), for every \(a\ge0\),
\[
                  A\cap\mathcal K_a=\varphi^{-a}S .
\]
For \(b\) in the intersection, use \(v_1=b,v_2=-1\) with
\(r_1=1,r_2=\varphi^a b\in A\) in
\eqref{eq:period-basis-order-retraction}. The right-hand sum is
zero. Thus \(b=\varphi^{-a}F(\varphi^a b)\in\varphi^{-a}S\).
The reverse inclusion follows from \(\varphi^{-a}S\subset A\)
and the definition of \(\mathcal K_a\).

Write \(x=\sum_i r_i n_i\) in \(M\). If all its form
values lie in \(\mathcal K_0\), apply
\eqref{eq:period-basis-order-retraction} with \(a=0\) to the
values \(-\lambda(x),\lambda(n_i)\) and coefficients
\(1,r_i\). The exact relation survives the replacement
\(r_i\mapsto F(r_i)\). Crystalline forms separate vectors
by the inverse period matrix, so \(x=\sum_i F(r_i)n_i\in N\).
For every \(l\ge0\) the same argument gives
\[
 \{x\in\varphi^l M:\lambda(x)\in\mathcal K_0
         \text{ for every chart-equivariant }\lambda\}
       =\sum_i S\varphi^l(n_i).
\]
Indeed the basis values \(\lambda(\varphi^l n_i)\) lie in
\(\mathcal K_0\) by \eqref{eq:period-basis-evaluation};
the coefficients \(\varphi^l(r_i)\) of \(x\) may therefore
be replaced by \(F(\varphi^l(r_i))\in S\). The reverse
inclusion uses \(S\subset\mathcal K_0\) and \(\varphi^l A=A\).
Since \(\varphi^l(S)\subset S\), the right side is
\(S\cdot\varphi^l(N_T)\).

Let \(D\) be a finite product of the \(E_j\), allowing
repetitions and the empty product. On the localized modules define
\[
 R_T\!\left(D^{-1}\sum_i r_i n_i\right)
     =D^{-1}\sum_i F(r_i)n_i,\qquad r_i\in A .
\]
Passing to a common denominator gives the same expression, by
the independence of the \(n_i\) and \(S\)-linearity of \(F\).
Under the changes of basis in Section~\ref{sec:local-integral-period-frame},
the period columns become \(n_1=H^{-1}nV\). For a column of
coefficients \(r\in A^d\), \(F(V^{-1}r)=V^{-1}F(r)\) entrywise.
Hence the new retracted column is
\(D^{-1}n_1 F(V^{-1}r)=H^{-1}D^{-1}nF(r)\).
This is the same vector in the new integral coordinates, so its
orders on all equivariant forms are preserved.
In particular \(R_T(ay)=F(a)y\) for \(a\in A\) and
\(y\in N_{\rm loc}\), and \(R_T\) fixes \(N_{\rm loc}\).
At each fixed denominator \(D\) its integral bound is
\[
 R_T\bigl(p^mD^{-1}M_T\bigr)=p^mD^{-1}N_T\qquad(m\ge0).
\]
The inclusion follows from \(F(p^mA)\subset p^mS\);
the restriction to the right-hand side is the identity.
For a chart-equivariant form \(\lambda\), put
\(v_a(y)=\lambda(\varphi^{-a}y)\). For \(y\in N_{\rm loc}\)
we have \(v_a(y)\in\mathcal K_a\) by
\eqref{eq:period-basis-evaluation} and
\(\varphi^{-a}(n)=n\varphi^{-a}(\mathcal B_a^{-1})\).
If \(D=\prod_j E_j^{d_j}\), then
\(\operatorname{ord}_t(\varphi^{-a}D)=d_a\).
Indeed \(E\) generates \(\operatorname{Fil}^1 B_{\rm dR}^+\);
for \(j>a\) the value \(\theta\varphi^{j-a}(\xi)=p\) is nonzero,
and for \(s=a-j>0\) one has
\[
 \theta\varphi^{-s}(\xi)
   =\frac{\zeta_{p^s}-1}{\zeta_{p^{s+1}}-1}\ne0 .
\]
The distinction between the three cases is preserved by
\(E/\xi\in A^\times\). Write \(x=\sum_i r_i(n_i/D)\). Applying
\eqref{eq:period-basis-order-retraction} to the values
\(v_a(n_i/D)\) and the additional value \(-v\) with coefficient
\(1=F(1)\) gives, for every \(a\ge0\) and \(v\in\mathcal K_a\),
\[
 \operatorname{ord}_t\bigl(v_a(R_Tx)-v\bigr)
       \ge \operatorname{ord}_t\bigl(v_a(x)-v\bigr).
\]
In particular \(v=0\) preserves every lower order bound.
For \(H=G_{K'}\) and \(d\in\mathbf Z\), the depth condition on a
value \(y\in B_{\mathrm dR}\) is equivalent to a congruence
\[
  (g-1)y\in\mathrm{Fil}^d\quad(g\in H)
  \qquad\Longleftrightarrow\qquad
  y-c\in\mathrm{Fil}^d\quad\text{for some }c\in K'.
\]
Indeed, if \(m=\operatorname{ord}_t y<d\), its leading class
in \(\mathrm{gr}^m B_{\mathrm dR}=\mathbf C_p(m)\) is nonzero and
fixed by \(H\). Tate's fixed-point theorems give \(m=0\) and
\(\theta(y)\in K'\). Subtract this constant. Any further
nonzero leading class below \(d\) would have positive degree
and be fixed by \(H\), which is impossible
\cite[\S3.3, Theorems~1 and~2, p.~176]{TatePDivisible}.
The converse follows from the \(H\)-stability of \(\mathrm{Fil}^d\).
Apply the preceding inequality to \(v=c\in K'\subset\mathcal K_a\)
for each level and each form. Thus the same \(R_T x\) retains all
prescribed order and depth bounds simultaneously.

An equivariant lattice map \(f:T\to T'\) carries \(M_T\) to
\(M_{T'}\) and \(N_T\) to \(N_{T'}\) by pulling back the
forms in their defining tests. In the two chart bases, write
\(f(n_i)=\sum_j s_{ji}n'_j\), \(s_{ji}\in S\). Using the same
scalar coordinate \(F\) for both lattices gives
\[
 R_{T'} f(x)
 =D^{-1}\sum_j F\Bigl(\sum_i s_{ji}r_i\Bigr)n'_j
 =D^{-1}\sum_{i,j}s_{ji}F(r_i)n'_j=fR_T(x).
\]
This applies to the lattice projections.
It also gives commutation with their integral spectral idempotents,
which are polynomials in the coefficient action over \(W(k)\).

Let \(\gamma\) be a \(\mathbf Z_p\)-automorphism of the coefficient
ring \(O\). On the same chart put \(T^\gamma=O_\gamma\otimes_O T\),
where \(O_\gamma\) has its usual left action and right action
\(b\cdot a=b\gamma(a)\). The map \(f_\gamma(x)=1\otimes x\) is
Galois-equivariant and \(\mathbf Z_p\)-linear, with
\(f_\gamma(ax)=\gamma(a)f_\gamma(x)\). Extending it by the
identity on the chart scalars gives
\[
 R_{T^\gamma}f_\gamma=f_\gamma R_T,\qquad
 f_\gamma e_i(\zeta)=e_i(\gamma(\zeta))f_\gamma .
\]
Here \(e_i(\zeta)\) is the projector onto the \([i]\)-eigenspace
of multiplication by the coefficient root \(\zeta\), as in
Section~\ref{sec:local-standard-orders}. Its projector polynomial
has coefficients in \(W(k)\), which the extended map fixes, while
it carries multiplication by \(\zeta\) to multiplication by
\(\gamma(\zeta)\).

In particular the coordinate retraction is \(O\)-linear on the lattice
with its given coefficient action. For a finite extension \(O'/O\)
put \(T'=O'\otimes_O T\), on the same chart. An \(O\)-basis of \(O'\)
gives equivariant inclusions and projections exhibiting \(T'\) as
a direct sum of copies of \(T\) over \(O\). Applying the preceding
identities to these maps gives, under the natural maps,
\[
 M_{T'}=O'\otimes_O M_T,\qquad N_{T'}=O'\otimes_O N_T,\qquad
 R_{T'}=1\otimes R_T .
\]
The last identity holds after inverting the \(E_a\) as well.
For \(l\ge0\), let \(P_l=\varphi^l(\tau^b M_T)\subset A\otimes T\).
Since \(\tau^b M_T=QA^d\) and \(\varphi\) is bijective on \(A\),
the saturation of the integral column span in
Section~\ref{sec:local-integral-period-frame} gives
\(P_l\cap p^m(A\otimes T)=p^m P_l\).
If \(\varpi\) is a uniformizer of \(O\), write
\(p=u\varpi^{e_O}\) with \(u\in O^\times\). For \(y\in A\otimes T\)
with \(\varpi^s y\in P_l\), choose \(m\) such that \(e_Om\ge s\).
The \(O\)-action on \(T\) commutes with Frobenius and preserves
\(M_T\), so \(p^m y=u^m\varpi^{e_Om-s}(\varpi^s y)\in P_l\).
Thus \(y\in P_l\), which gives saturation at every power of
\(\varpi\) as well.

For \(\sigma\in G_{\mathbf Q_p}\), transport the lattice and chart as in
Section~\ref{sec:local-integral-period-frame}. Taking
\(F^\sigma=\sigma F\sigma^{-1}:A\to\sigma S\) in the transported
basis \(\sigma n_i\) gives
\[
 R^\sigma(\sigma x)
   =(\sigma D)^{-1}\sum_i F^\sigma(\sigma r_i)\,\sigma n_i
   =\sigma R_T(x),\qquad x=D^{-1}\sum_i r_i n_i.
\]
Indeed \(\sigma E_j=\varphi^j(\sigma E)\), so these are the divisors
in the transported chart. For translates by the subgroup \(H\) with
the original coordinate map fixed, we use the following residue calculation.

For \(v\in\mathcal K_a\cap\mathrm{Fil}^0\), choose
\(r\in A\) and \(m\ge0\) with
\(p^mv\equiv r\bmod\mathrm{Fil}^1\), using the
surjectivity of \(\theta:A\to\mathcal O_{\mathbb C_p}\).
In \eqref{eq:period-basis-order-retraction} use the two values \(p^m v,-1\) and the two
coefficients \(1,\varphi^a(r)\). This replaces \(r\)
by \(\varphi^{-a}F\varphi^a(r)\) modulo
\(\mathrm{Fil}^1\). With \(\pi_a=\theta(\varphi^{-a}z)\)
this gives
\[
 \theta(v)=p^{-m}\theta\varphi^{-a}F(\varphi^a r)
            \in K'(\pi_a).
\]
The integral tower calculation in
Section~\ref{sec:period-chart-approximation} gives the first equality
below. The preceding inclusion gives one direction of the second,
and the residue field contains the residues of \(K'\) and
\(\varphi^{-a}z\), hence contains \(K'(\pi_a)\). Thus
\[
 \theta\varphi^{-a}(S)=W(k)[\pi_a]=\mathcal O_{K'(\pi_a)},\qquad
 \theta\bigl(\mathcal K_a\cap B_{\mathrm dR}^+\bigr)=K'(\pi_a).
\]
More explicitly, let \(g_b\in H\), \(m_b\in\mathbf Z\), and
suppose \(\sum_b m_b g_b(s)\in E_a A\) for all \(s\in S\).
For \(v\in\mathcal K_a\cap\mathrm{Fil}^0\) the residue
calculation above supplies \(s\in S\) with
\(p^m v\equiv\varphi^{-a}s\bmod\mathrm{Fil}^1\) for some \(m\).
Commuting each \(g_b\) with \(\varphi^{-a}\) and reducing gives
\(\sum_b m_b g_b(v)\in\mathrm{Fil}^1\). In particular
\(\pi_0\in K'\), so an integral form value at level zero
of \(x\in N_{\rm loc}\) has \(g\)-difference of order at least
one for every \(g\in H\). This is the extra order needed for
the level-one division by \(E\) below. For the tame generator acting nontrivially on \(\theta(z)\),
Section~\ref{sec:tame-cup-original-periods} proves divisibility
by \(z\) using the annihilating polynomial of its linear period action.

\subsubsection{Tensor pairings and orders}
\label{sec:period-basis-tensor-orders}

\sectionopening{In this subsection, we apply the order estimates of Section~\ref{sec:period-basis-order-retraction} to a pairing.}

Finally let \(b_0:T_1\times T_2\to T_3\) be an equivariant pairing, bilinear over the coefficient ring.
The localized coefficients and all their Frobenius translates
lie in \(A_{(p)}\), since each inverted
period has nonzero reduction modulo \(p\). The inclusions
\[
 W(\mathbb C_p^\flat)\ \longleftarrow\ A_{(p)}
                         \ \longrightarrow\ B_{\rm dR}
\]
give Witt reduction modulo \(p^m\) on the left and the de Rham
value on the right. For the value at level \(a\), apply
\(\varphi^{-a}\) before taking the right inclusion. The Witt
inclusion commutes with Frobenius and the chart Galois action,
and scalar extension preserves the pairing. The coefficient idempotents
for an index \(i\) commute with this base change, and
bilinearity gives
\(b_0(e_i y,x)=e_i b_0(y,x)=b_0(y,e_i x)\). In particular
\(b_0(e_i y,e_{i'}x)=0\) for \(i\ne i'\). At level \(a\) the
relevant index in both factors and in the target is \(j=i^{p^{-a}}\).
A form of type \(u\) on \(T_s\) takes the lattice
projected to \(j\) into \(\mathrm{Fil}^{-u}\).

For each projected factor, choose a \(K'\)-basis of the invariant
de Rham dual adapted to its filtration. Write \(-u_s\) for the
filtration index of \(\lambda_s\), so its type is \(u_s\).
The functionals \(t^{u_s}\lambda_s\) form an
integral de Rham dual basis, as proved in
Section~\ref{sec:local-standard-orders}. On the two parts
indexed by \(j\), expand the pullback of a form \(\lambda\)
of type at most \(v\) as
\[
 \lambda\circ b_0=\sum_{s,h}\alpha_{sh}
                  \lambda_s\otimes\lambda'_h,\qquad \alpha_{sh}\in K'.
\]
Indeed the tensor products give a de Rham basis, and its
coefficients are fixed by the chart subgroup. Bilinearity on
the original lattices makes \(t^v(\lambda\circ b_0)\)
integral; in the rescaled dual bases this means
\(t^{v-u_s-u'_h}\alpha_{sh}\in B_{\mathrm dR}^+\).
A nonzero scalar in \(K'\) has de Rham order zero, so
\(\alpha_{sh}=0\) if \(u_s+u'_h>v\). For two vectors whose
orders are bounded below by \(c_1,c_2\), it follows that
\[
 \operatorname{ord}_t\lambda\bigl(
         \varphi^{-a} b_0(x_1,x_2)\bigr)
 \ \ge\ \min_{u+u'\le v}\bigl(c_1(a,u)+c_2(a,u')\bigr),
\]
where \(u,u'\) run through the jumps at \(j\), and the minimum
of an empty set is \(+\infty\).

For \(x_s\in M_{T_s}\), bilinearity on the original lattices
first gives \(b_0(x_1,x_2)\in A[1/\tau]\otimes T_3\).
For every \(n>0\), expand the pullback of any invariant form
on \(T_3\) in products of invariant dual bases of the two factors.
All coefficients lie in \(K'\); both evaluations on
\(\varphi^n x_s\) lie in \(B_{\mathrm dR}^+\) by the defining
test for \(M_{T_s}\). This proves
\(b_0(x_1,x_2)\in M_{T_3}\). For localized vectors with
denominators \(D_1,D_2\), the product \(D_1D_2\) clears
the paired vector, giving the same conclusion in \(M_{T_3,\rm loc}\).

For a saturated inclusion \(T_0\subset T\) of lattices crystalline
on this chart, restriction of equivariant de Rham forms from
\(T\) to \(T_0\) is surjective. After extending scalars from
\(K'\) to \(B_{\rm dR}\), the comparison isomorphisms identify
this restriction with ordinary dual restriction. Thus it is already
surjective over \(K'\). With \(A_{\rm loc}=A[1/E_a:a\ge0]\),
the positive Frobenius tests and the value-field test give
\[
 \begin{aligned}
 M_{T_0,\mathrm{loc}}
  &= M_{T,\mathrm{loc}}\cap(A_{\rm loc}[1/\tau]\otimes T_0),\\
 N_{T_0,\mathrm{loc}}
  &= N_{T,\mathrm{loc}}\cap(A_{\rm loc}[1/\tau]\otimes T_0).
 \end{aligned}
\]
Indeed, clear a common finite product of the \(E_a\) and apply
those tests before localization. Saturation makes \(T_0\) a direct
summand of the underlying \(\mathbf Z_p\)-module, so the cleared
vector lies in \(A[1/\tau]\otimes T_0\). For the orders at each
individual jump on the kernel of a cocycle extension we use the
type-preserving extension of forms in Section~\ref{sec:local-tame-moves}.

\subsection{Standard orders and projections on the chart}
\label{sec:local-standard-orders}

\sectionopening{In this subsection, we identify the lattices defined by the
filtration in Section~\ref{sec:v2-crystalline-frame-action}
with \(\varphi^l M_T\). Frobenius transports the projection
cokernel along the residue-field embeddings, as used in
Section~\ref{sec:v2-crystalline-lift-ideals}.}

For a lattice \(T\) crystalline on the chart subgroup use the
two bases \(M_T=\sum A n_i\) and \(N_T=\sum S n_i\) just
constructed. Write \(E_a=\varphi^a E\) for \(a\ge0\), and
localize both bases at the \(E_a\). A residue index is a nonzero
residue-field element \(i\). Write \(m=|k^\times|\) and let \(Z_\zeta\) be the endomorphism
of \(T\) given by the chosen coefficient root \(\zeta\), of order
dividing \(m\). Its integral projector for \(i\in k^\times\) is
\[
 e_i=\frac1m\sum_{r=0}^{m-1}[i]^{-r}Z_\zeta^r .
\]
Indeed \(m\in W(k)^\times\), \(Z_\zeta^m=1\), and the distinct
Teichm\"uller roots have unit differences. These projectors commute
with \(R_T\) by \(W(k)\)-linearity and naturality for \(Z_\zeta\).
After extending the de Rham scalars to \(\mathbf C_p\), \(e_i\)
selects exactly the coefficient embeddings \(\gamma\) with
\(\gamma(\zeta)=[i]\), since \(Z_\zeta\) acts on this summand
by \(\gamma(\zeta)\). For a further coefficient embedding
\(\gamma'\) extending \(\gamma\), this eigenvalue is unchanged.
At level \(a\) Frobenius on the Witt coefficients gives
\(\varphi^{-a}(e_i)=e_{i^{p^{-a}}}\).

Let \(\lambda\) be a chart-equivariant de Rham form of type
\(u\) on that latter part, meaning that its restriction to the
original lattice, after projection, takes values in
\(t^{-u}B_{\mathrm dR}^+\). The Hodge filtration is split after
choosing an adapted basis of invariant forms. In particular an
element of the de Rham vector space belongs to the original
\(B_{\mathrm dR}^+\)-lattice precisely when its values on the
basis forms of types \(u\) are in \(t^{-u}B_{\mathrm dR}^+\).

Indeed the reductions of the rescaled forms \(t^u\lambda\)
are independent over \(\mathbf C_p\). In a minimal
dependence, write \(a_i=1\) for the normalized coefficient
and \(u_j\) for the type of \(\lambda_j\). The residue of
\(t^{u_j}\lambda_j\) transforms by \(\chi_{\rm cyc}^{u_j}\).
Uniqueness of the normalized relation gives
\[
 h(a_j)=\chi_{\rm cyc}(h)^{u_i-u_j}a_j,\qquad
 a_j\in\mathbf C_p(u_j-u_i)^H .
\]
By the fixed-point calculation in Section~\ref{sec:period-period-realization},
every nonzero coefficient has \(u_j=u_i\) and belongs to \(K'\).
Lifting this relation with these constants gives a relation in the
corresponding graded space, contrary to the adapted basis.

The rescaled forms therefore give a basis of the
\(B_{\mathrm dR}^+\)-dual: their number is the dimension of
the projected space by de Rham comparison. For each residue index \(j\), let \(\lambda_s\) be this adapted dual
basis and write
\[
 D_j=D_{\rm dR}(T[1/p])_j,\qquad
 u_s=-\max\{k:\lambda_s\in\mathrm{Fil}^k(D_j^\vee)\}.
\]
Here \(D_j^\vee=\operatorname{Hom}_{K'}(D_j,K')\) carries the
decreasing dual filtration.
The \(u_s\) are the jumps of the original representation on that
part. Fix a finite set \(\mathcal J_j\) containing all these types.
For an invariant form \(\lambda=\sum_s a_s\lambda_s\), the
coordinates of \(t^u\lambda\) in the integral dual basis are
\(a_s t^{u-u_s}\). Every nonzero \(a_s\in K'\) has de Rham order
zero. Hence the forms of type at most \(u\) are precisely the
\(K'\)-span of the \(\lambda_s\) with \(u_s\le u\). A vector
\(x\in e_i M_{T,\mathrm{loc}}\) has order bounds \(c\)
when for all \(a\ge0\), \(j=i^{p^{-a}}\), \(u\in\mathcal J_j\),
and all chart-equivariant forms of type at most \(u\) on the
part at \(j\),
\[
 \operatorname{ord}_t\lambda(\varphi^{-a} x)\ge c(a,u).
\]
Equivalently, in the adapted dual basis for \(j=i^{p^{-a}}\)
the exact tests at level \(a\) are
\[
 \operatorname{ord}_t\lambda_s(\varphi^{-a}x)
      \ge \max_{\substack{u\in\mathcal J_j\\u_s\le u}} c(a,u)
                    \qquad(\text{every }s).
\]
Each \(\lambda_s\) occurs in every test indexed by \(u\ge u_s\).
Conversely, expanding any form of type at most \(u\) in the indicated
\(K'\)-span gives its order bound by the ultrametric inequality.
In particular the standard bounds \(0\) and \(-u\) give respectively
\(0\) and \(-u_s\) on \(\lambda_s\). Its intersection with
\(N_{T,\mathrm{loc}}\) is the module of chart vectors with those
orders. The coordinate retraction maps the first module onto
the latter, preserves congruences of the values to constants,
and commutes with an equivariant projection of lattices.

The shifted divisors can be distinguished by their Fontaine values.
For every integer \(j\) we have \(v_p(\theta\varphi^j z)=p^j/e\).
The two terms of \(\varphi^j E=p+(\varphi^j z)^e\varphi^j(u_E)\),
where \(E=p+z^eu_E\) and \(u_E\in S^\times\), have valuations \(1,p^j\).
For \(j\ne0\) these valuations are distinct, so
\[
 v_p(\theta\varphi^j E)=\min\{1,p^j\}.
\]
Since \(E_a\) generates \(\ker(\theta\varphi^{-a}|_A)\) and
\(E\) is a uniformizer of \(B_{\mathrm dR}^+\), this gives
\(\operatorname{ord}_t(\varphi^{-a}E_b)=\mathbf1_{a=b}\).
Also
\[
 \theta\varphi^j\tau=
 \begin{cases}
  0,&j>0,\\
  \zeta_{p^{1-j}}-1,&j\le0.
 \end{cases}
\]
In particular \(\varphi^{-a}(z\tau)\) has de Rham order zero for
all \(a\ge0\), including \(p=2,a=0\) where \(\theta(\tau)=-2\).

Define \(\operatorname{std}_l(a,u)=0\) for \(a<l\), and
\(\operatorname{std}_l(a,u)=-u\) for \(a\ge l\).
The part at \(i\) having these orders is precisely the part
at \(i\) of \(\varphi^l M_T\). To prove the direct inclusion,
\(\varphi^l M_T\) is contained in \(M_{T,\mathrm{loc}}\) by
\(E^b\varphi(N_T)\subset N_T\) and induction. Write a vector
there as \(\varphi^l y\). At \(a<l\), positive Frobenius
integrality of \(y\) gives order at least \(0\). At \(a\ge l\), a power
of \(\tau\) makes \(y\) an original integral vector. That power has order zero under every \(\theta\varphi^{-(a-l)}\),
so the lattice condition gives order at least \(-u\).
Conversely, clear from \(x\) a finite product of the \(E_a\)
and a power of \(\tau\) so that all its original coordinates
are integral. For \(a\ge l\) the adapted-basis test puts every original
coordinate of \(\varphi^{-a}x\) in \(B_{\rm dR}^+\). If an original
coordinate of \(x\) is written \(c/(E_a^hD')\), with \(c\in A\),
\(h>0\), and \(D'\) a product of \(\tau\) and the other shifted divisors,
then \(\theta\varphi^{-a}D'\ne0\). Since \(E\) has de Rham order one,
\(\varphi^{-a}c\in\mathrm{Fil}^h\), in particular \(c\in E_a A\).
Division by \(E_a\) reduces \(h\) by one. Repeating for every
\(a\ge l\) removes precisely the denominators excluded by these tests.

Write the remaining common denominator as
\(\tau^m\prod_{a<l} E_a^{m_a}\), with nonnegative exponents.
Since \(E_a\) is associated to \(\varphi^a\xi\) and
\(\varphi^l\tau=\tau\prod_{a<l}\varphi^a\xi\),
this denominator divides \((\varphi^l\tau)^h\) for
\(h\ge m\) and \(h\ge m_a\) for all \(a<l\).
Thus \(\tau^h\varphi^{-l}x\in A\otimes T\).
Put \(y=\varphi^{-l}x\). To test \(\varphi^n y\) for
\(1\le n\le l\), take \(a=l-n\) in the standard orders:
all basis forms have integral values. If \(n>l\), use
\(x\in M_{T,\mathrm{loc}}\) to choose a product \(D'\) of
the \(E_a\) with \(D'x\in M_T\). The positive Frobenius test
makes every tested value of \(\varphi^{n-l}(D'x)\) integral.
Each factor of \(\varphi^{n-l}D'\) has de Rham order zero, since
\(\theta\varphi^{n-l}(E_a)\ne0\) for every \(a\ge0\).
Dividing gives integral values on \(\varphi^n y\).
Thus \(y\in M_T\), proving the equality.
For a fixed index, the corresponding chart lattice is
\[
 e_i\varphi^lM_T\cap N_{T,\mathrm{loc}}
       =e_i\sum_j S\,\varphi^l(n_j).
\]
Indeed the values of every chart vector lie in \(\mathcal K_0\) at
level zero, so the coefficient extraction of
Section~\ref{sec:period-basis-order-retraction} gives the inclusion
from left to right. The reverse follows from
\eqref{eq:period-basis-frobenius-matrix} and preservation of
\(N_{T,\mathrm{loc}}\) by \(e_i\).

Use one scalar coordinate \(F=F_0:A\to S\), with \(F(1)=1\),
for all lattices on this chart. For a trivial finite free parameter
lattice \(C\) over \(O\), choose a \(\mathbf Z_p\)-basis of \(C\).
At the positive Frobenius powers its coordinate rows are equivariant
de Rham forms. If \(\tau^{-m}a\), \(a\in A\), is one coefficient
of a point of \(M_C\) and \(m>0\), then
\(\theta\varphi^j(a)=0\) for every \(j>0\). The cyclotomic
intersection gives \(a\in\tau A\). Iterating on \(m\) proves
\(M_C=A\otimes_{\mathbf Z_p}C\). The intrinsic chart test at
level zero puts each coordinate of an element of \(N_C\) in
\(\mathcal K_0\). The exact relation in
\eqref{eq:period-basis-order-retraction} fixes each such coefficient
under \(F\), so the coefficients belong to \(S\). Conversely, all
vectors of \(S\otimes_{\mathbf Z_p} C\) pass the test since the values
on constant basis vectors lie in \(K'\). Thus
\(N_C=S\otimes_{\mathbf Z_p} C\), with \(e_iN_C\) generated over
\(S\) by \(e_i(1\otimes c)\), \(c\in C\). On the parameter \(R_C\)
acts coefficientwise by \(F\), also after localization. The assertion applies both to a
single-cocycle parameter \(C=O\) and to the finite free module of
crystalline cocycles in Section~\ref{sec:v2-crystalline-deep-constants}.

For an extension on any of the three sides, write \(\pi:X\to C\)
for its parameter projection. Write \(V_l(i)\)
for the level-\(l\) standard-order submodule of the part at
\(i\) of \(M_{X,\rm loc}\). The multiplier ideal is the
following ideal of \(A\):
\[
 I_l(i)=\{f\in A:f e_i(A\otimes C)\subset\pi(V_l(i))\}.
\]
For \(c\in C\) put \(c_i=e_i(1\otimes c)\). These projected
constants span \(e_i(A\otimes C)\) over \(A\); for the trivial
lattice each lies in its chart module, since the constant
\(1\otimes c\) is Frobenius fixed with integral constant form
values and the projector preserves the chart module. Thus for
\(q=p^\nu\), membership \(f\in I_l(i)+(q)\) is equivalent to
\[
 \exists w\in A\quad\forall c\in C\quad\exists x_c\in V_l(i),
 \qquad \pi(x_c)=(f-qw)c_i.
\]
For \(f\in S\), retraction replaces \(w\) by the single
element \(F(w)\in S\) and \(x_c\) by a lift in
\(V_l(i)\cap N_{X,\rm loc}\). Indeed \(R\) commutes with
\(\pi\), fixes the projected constants and satisfies
\(R(ac_i)=F(a)c_i\).

Witt Frobenius identifies \(V_l(i)=e_i\varphi^l M_X\)
bijectively with \(V_{l+1}(i^p)\) and sends \(c_i\) to \(c_{i^p}\).
Since it fixes \(q\), for each \(c\in C\) we have
\[
 \pi(x_c)=(f-qw)c_i
 \quad\Longleftrightarrow\quad
 \pi(\varphi x_c)=(\varphi(f)-q\varphi(w))c_{i^p}.
\]
The witness \(\varphi(w)\) works simultaneously for all \(c\).
In the reverse direction apply \(\varphi^{-1}\) to the common witness
and every lift over \(A=W(\mathcal O_{\mathbf C_p}^{\flat})\). Thus
\[
 \varphi\bigl(I_l(i)+(q)\bigr)=I_{l+1}(i^p)+(q).
\]
For \(f\in S\) this gives an equivalence of membership for
\(f\) and \(\varphi(f)\), with chart lifts obtained by retraction
at the corresponding level. More explicitly, for \(f\in S\),
\(\pi x=(f-qw)c_i\) and \(a\ge0\), put \(i_a=i^{p^a}\). Then
\[
 \begin{aligned}
 \pi R_X(\varphi^a x)
   &=\bigl(\varphi^a f-q F(\varphi^a w)\bigr)c_{i_a},\\
 \pi\varphi^a(R_X x)
   &=\bigl(\varphi^a f-q\varphi^a(F(w))\bigr)c_{i_a}.
 \end{aligned}
\]
Subtracting the two parameter projections gives
\[
 \pi(R_X\varphi^a-\varphi^a R_X)x
       =q\bigl(\varphi^a(F(w))-F(\varphi^a w)\bigr)c_{i_a}.
\]
Thus both projections reduce to \((\varphi^a f)c_{i_a}\) modulo \(q\).
For exact lifts \(w=0\) they agree integrally at every level.
These identities apply componentwise to the parameter maps of the
two cocycle extensions and the extension of the pairing.
At level \(a\) the same scalar coordinate is used on each lattice
in the form \(F_a=\varphi^{-a}F\varphi^a\) of
\eqref{eq:period-basis-order-retraction}.

Multiplication by \(E_a\) raises
the level-\(a\) order by one and leaves all other orders unchanged,
since \(\theta\varphi^{-j}(E_a)\ne0\) for \(j\ne a\). If
\(\omega\ge0\) and \(u\le\omega\) for every \(u\in\mathcal J_i\),
multiplication by \(E^\omega\) carries level-zero lifts to
level-one lifts: at level zero the new bound is
\(\omega-u\ge0\), and at every positive level the required
bound is still \(-u\).

For \(g\) in the chart subgroup and \(x\in N_{T,\rm loc}\) of
standard order \(l\), Section~\ref{sec:period-basis-integral-operators}
gives \((gx-x)/(z\tau)\in M_{T,\rm loc}\) with the same lower
bounds on its orders. When \(l=1\), we may also divide by \(E\)
while preserving these bounds. To see this, apply at level zero
the residue calculation of Section~\ref{sec:period-basis-order-retraction}.
If \(\lambda(x)\in\mathrm{Fil}^0\), then
\(\theta(\lambda(x))\in K'\); hence
\(\lambda(gx-x)=g\lambda(x)-\lambda(x)\in\mathrm{Fil}^1\).
The other levels are unaffected by
division by \(E\). On the trivial parameter, the same
quotient of a chart scalar \(w\in S\) is integral by
\(gw-w\in z\tau E A\). Thus division and subsequent
retraction give the integral lift-ideal operations with
divisor \(z\tau\) at level zero and \(z\tau E\) at level
one, for the original parameter, index, and modulus.

\section{The tame cup-product estimate}
\label{sec:tame-cup-appendix}

\sectionopening{In this appendix, we apply the filtration orders of
Appendix~\ref{sec:crystalline-basis-appendix} on the tame chart.
The level-zero estimates control the four terms in the cup
calculation and give the bound at the tame level primes in
Section~\ref{sec:relative-crystalline-colength}.}

\subsection{Level-zero steps for tame periods}
\label{sec:local-tame-depth}

\sectionopening{In this subsection, we use the divided difference in
Section~\ref{sec:local-standard-orders} to raise the order of a Galois
difference at \(E=0\). The valuation of the scalar
\((\omega^j-1)/(z\tau)\) determines the power of \(p\) required at
each step. The operator preserves every filtration bound already attained.}

Let \(H=G_{K'}\) for the finite tame chart field \(K'\)
over which the lattices are crystalline. At \(p\mid M\),
Section~\ref{sec:tame-cup-original-periods} takes \(K'=K_0(\pi)\)
with \(K_0/\mathbf Q_p\) unramified of degree \(r!\),
\(e=\operatorname{lcm}_{1\le b\le r}(p^b-1)\), and \(\pi^e=p\).
Write \(E=p-z^e\), \(\Delta_z=z\tau\),
and \(\mu=E\tau v_0\) with \(v_0\in A^\times\). Choose
\(g\in H\) such that \(gz=[\underline\epsilon]z\). Indeed,
\(X^{p^a}-\theta(z)\) is Eisenstein over \(K'\), so at each
finite level such a conjugate of the chosen root exists. The
conditions are nested and \(H\) is compact. Write \(gE=E\omega\),
where \(\omega\in A^\times\) since \(g\) preserves \(\ker\theta\). Then
\[
 \frac{\omega-1}{\Delta_z}
   =-v_0 z^{e-1}\sum_{i=0}^{e-1}[\underline\epsilon]^i .
\]
Its value under \(\theta\) has normalized \(p\)-adic
valuation \((e-1)/e\), since \(e\) is prime to \(p\).
Let \(\delta=1\) for \(e>1\) and \(\delta=0\) for \(e=1\).
For \(0\ne j\in\mathbf Z\) put \(X_j=(\omega^j-1)/\Delta_z\).
This lies in \(A\) by the geometric sum and the invertibility
of \(\omega\). Since \(v_p(\theta z)=1/e\) and
\(v_p(\theta\tau)=1/(p-1)\), we obtain
\[
 v_p(\theta\omega-1)=\frac{p}{p-1},\qquad
 v_p(\theta X_j)=\frac{e-1}{e}+v_p(j).
\]
Indeed, if \(v_p(a)>1/(p-1)\), the unique term of least valuation
in \((1+a)^p-1\) is \(pa\). Taking a power prime to \(p\) preserves
\(v_p(a)\) by the geometric sum, and inversion preserves it since
\(1+a\) is a unit. These calculations also apply at \(p=2\),
where \(v_2(\theta\omega-1)=2\).
Every nonzero step in the tame level-prime application satisfies
\(p\nmid j\) by the arithmetic radius. For these \(j\) the ratio
\(p^\delta/\theta(X_j)\) is integral. Surjectivity of \(\theta\)
therefore gives \(\upsilon_j\in A\) with
\(\upsilon_jX_j\equiv p^\delta\bmod E\).

If \(y\in B_{\rm dR}\) satisfies \((h-1)y\in\mathrm{Fil}^j\)
for all \(h\in H\), then \(y\in K'+\mathrm{Fil}^j\).
Indeed, if its order \(n\) is below \(j\), the nonzero leading
class belongs to \(\mathbf C_p(n)^H\). This space vanishes for
\(n\ne0\), and for \(n=0\) it equals \(K'\). In the latter case
subtract the representative in \(K'\) and apply the same argument
to the difference. The fixed-point calculation for every
integer \(n\) is given in Section~\ref{sec:period-period-realization}.

Let \(x\) be a localized chart point over a constant parameter,
with selected form value \(y=\lambda(x)\) of depth \(j\). Choose
\(k\in K'\) with \(y-k\in\mathrm{Fil}^j\).
By Section~\ref{sec:period-basis-order-retraction},
\(y,E\in\mathcal K_0\) and \(K'\subset\mathcal K_0\).
Thus \(h=E^{-j}(y-k)\)
belongs to \(\mathcal K_0\cap B_{\rm dR}^+\) for either sign of
\(j\). Its residue is in \(K'(\pi_0)=K'\) by the same section.
Consequently \(gh-h\in\mathrm{Fil}^1\) and
\(g(y-k)-\omega^j(y-k)=\omega^jE^j(gh-h)\in\mathrm{Fil}^{j+1}\).
Since \(gk=k\),
\[
 \theta\!\left(E^{-j}(gy-y)/\Delta_z\right)
       =\theta(X_j)\,\theta(h)\qquad (j\ne0).
\]
Consequently
\[
 x\longmapsto R\!\left(p^\delta x-
             \upsilon_j(gx-x)/\Delta_z\right)
\]
moves the values to the next filtration step modulo the
scaled constants. The divided difference and its retraction
preserve the old orders and lie over zero. This follows
from Section~\ref{sec:local-standard-orders}, since \(g\)
fixes the parameter constant and \(\Delta_z\) has order
zero at all places of the chart. At \(j=0\) a chart value
that is integral at \(E\) has residue fixed by the chart
group, so \(p^\delta x\) already advances this step.
Thus the contribution of each step is at most \(\delta\).

Perform the steps in increasing order of \(j\). For a
line whose target depth is \(j+1\), all earlier needed
steps have then been performed, while the operator
preserves previously obtained depths at every line. Indeed
a value congruent to a constant modulo \(\operatorname{Fil}^{h}\)
has Galois difference of order at least \(h\). Division by
\(\Delta_z\) has order zero, and the retraction preserves
the congruence. This advances the orders over the finite union of step intervals.
 If a form value originally
had positive order, the constant in its higher-depth
congruence must be zero. Otherwise the new guaranteed
order is the minimum of \(0\) and the new depth. In
both cases the operator also retains the initially
imposed orders. In the tame application the starting
order is \(-u\), and the final depth at zero is
\(\eta-u\). The final order is at least \(-u\). On the side whose far jumps are negative we have \(u\le\eta\).
Its increased depth is nonnegative, so its final order is at least
\(\max(0,-u)\), in particular at least \(\max(0,-u-\eta)\).

\subsection{The four terms of the crystalline cup calculation}
\label{sec:local-tame-moves}

\sectionopening{In this subsection, we prove the tame cup estimate by comparing the two pairs of
crystalline lifts before and after Frobenius. The order estimates of
Section~\ref{sec:local-tame-depth} place four bilinear terms in the
image of the Witt vectors of the maximal ideal of the tilt. The
Artin--Schreier--Witt equation then produces a continuous cochain
whose differential is the required cup product.}

On both sides choose the resulting points at each index
with the same constant multiplier, multiplying by extra powers of \(p\) to reach the common exponent. Their orders are at least \(-u\) at
every level and their depths are \(\eta-u\) at zero and
\(-u\) above. A difference of two points over the same
parameter belongs to the kernel lattice with those orders.

For strictness of the filtration, take an invariant
form \(\lambda\) of type \(u\) on the kernel and a
crystalline primitive \(b(c)\) for the cocycle of a
basis element \(c\) of the trivial parameter lattice
over \(\mathbf Z_p\). Since
\[
 (g-1)\lambda(b(c))=\lambda(z(c)(g))\in\operatorname{Fil}^{-u},
\]
the invariant calculation in Section~\ref{sec:local-tame-depth}
gives \(k(c)\in K'\) with
\(\lambda(b(c))-k(c)\in\operatorname{Fil}^{-u}\).
Extend \(\lambda\) on \((v,c)\) by adding this last
value to \(\lambda(v)\), and extend linearly over that
parameter basis. The resulting form is equivariant
and has type at most \(u\). The integral localized basis
of the kernel is the intersection with its ambient
subspace, by clearing the \(E_a\) and applying the
positive Frobenius tests. Hence the order bounds on
a difference apply to its original kernel components.

At index \(i\), write \(x,x'\) for the two current points
and \(y,y'\) for the Frobenius translates from the preceding
index. These old points have order at least \(0\) at level zero and
at least \(-u\) at every positive level. Their differences
under Galois have order at least \(1\) at level zero:
for a chart vector in \(N_{\rm loc}\), the values of its
Frobenius translate under \(\theta\) belong to the chart
field, by the change of basis, the invariance of the period
covector, and the residue map on \(S_B\). At higher levels
the old differences inherit the depth at the preceding
index.

Suppose the first side is the negative-far side. Take the
path with \(w=x,w'=x'\). Its first two moves are
\((y,y')\to(x,y')\to(x,x')\).
The first move \(y-x\) has order at least zero at level
zero because \(x\) and \(y\) both do. Use the Frobenius
translate of the preceding third period, of order at least zero
at level zero. The second-side resting difference at
\(y'\) has order at least \(1\) there, giving total at
least \(1\).

At positive levels the two orders sum to at least
\(-u-u'\). In the expansion of the trilinear form with Tate target
from Section~\ref{sec:period-basis-tensor-orders}, a third jump \(v\)
can contribute only when \(u+u'+v\le-1\). If there is no such
triple, the component of the trilinear form is zero. Otherwise
each contributing third jump has standard order
\(-v\ge u+u'+1\). In the second move \(y'-x'\), the moving side
has \(u'\ge-\eta\), so its order is at least \(-u'-\eta\)
at zero. The resting difference at \(x\) has depth
\(\eta-u\), and the current third period has order at least
\(u+u'+1\) on each contributing term. Thus total order is at least \(1\), also at
higher levels using the standard orders.

The last two moves are zero. If the second side has the
negative far jumps, take \(w=y,w'=x'\). The two active moves
\((y,y')\to(y,x')\to(x,x')\) are now the second and third,
respectively. They use the preceding and current third periods in
that order. The same order bounds
apply with the sides interchanged. The pair of cocycle
differences at two current points also passes the order
test with the current third period, since both depths
dominate their respective standard orders.

In each active test the factor \(\Delta_z\) comes from
the Galois difference of the resting point, or one of the
two cocycle differences. Dividing by \(\Delta_z\) preserves
its orders and the integral localized span by the basis
calculation. This applies also to a Galois translate of the difference:
\(g(z)/z\in A^\times\), and \(g(\tau)A=\tau A\), since \(g\)
preserves every \(\ker(\theta\varphi^j|_A)\), \(j>0\), whose
intersection is \(\tau A\). Thus \(g(\Delta_z)/\Delta_z\in A^\times\).
The third extension form is equivariant, so
orders add by the filtered tensor rule. The resulting
Tate value has its standard orders \(1\) at every level.
It is therefore in the Tate basis at level zero by
Section~\ref{sec:local-standard-orders}. The Tate basis
bound says that multiplication by \(\tau\) puts it in
\(A\otimes\mathbf Z_p(1)\). The remaining factor \(z\)
places the reduced pairing in the small Witt ideal.
This verifies all the order tests for the shared
Frobenius cochain calculation. Summing over indices
preserves the tests because the equivariant bilinear
pairing is \(O\)-balanced, so distinct spectral idempotents
have zero cross pairing. Since \(\sum_i e_i=1\) and
\(\varphi(e_i)=e_{i^p}\), summing the points at all indices puts
each point on side \(j\) over \(p^{e_j}(1\otimes 1)\), and summing
the translates gives the Frobenius translate of the same point.

\subsubsection{The Frobenius cup calculation}\label{sec:tame-cup-frobenius}

\sectionopening{In this subsection, we use the four terms and depth bounds in Section~\ref{sec:local-tame-moves} to construct the cochain.}

Fix two extensions \(Y\oplus O,D\oplus O\) for \(z_1,z_2\)
and the third period of the pairing. Work over the chart
subgroup \(H\). Write \(x,w\) for first-side points over
\(p^{e_1}\), and \(x',w'\) for second-side points over
\(p^{e_2}\). The two endpoints and intermediate points
give the path
\[
 (\varphi x,\varphi x')\longrightarrow (w,\varphi x')
 \longrightarrow(w,w')\longrightarrow(x,w')
 \longrightarrow(x,x').
\]
In the following computation the order tests for the four
moves mean that their four cross pairings, with the
resting difference as the other argument, lie in the
subgroup
\[
 I=\ker\!\left(W_m(\mathcal O_{\mathbf C_p^\flat})
          \longrightarrow W_m(\overline{\mathbf F}_p)\right)
           \otimes_{\mathbf Z_p}\mathbf Z_p(1)
 \ \subset\ W_m(\mathbf C_p^\flat)\otimes_{\mathbf Z_p}\mathbf Z_p(1).
\]
Membership means that all \(m\) Witt coordinates have positive
valuation in a Tate basis. This subgroup is stable under \(H\).

Fix a \(\mathbf Z_p\)-linear functional \(\lambda:O\to\mathbf Z_p\).
The integral \(O\)-bilinear pairing \(\psi\) is equivariant
modulo \(p^m\), so
\[
 \psi_m=((\lambda\otimes1)\circ\psi)\bmod p^m:
       (Y/p^mY)\times(D/p^mD)\longrightarrow(\mathbf Z/p^m)(1)
\]
is equivariant and bilinear over \(\mathbf Z/p^m\). It satisfies
\(\psi_m(by,d)=\psi_m(y,bd)\) for \(b\in O\).
Let \(B\) be the Witt-bilinear extension of \(p^{e_3}\psi_m\)
on the two kernel components. For points in
\(Y\oplus O\) and \(D\oplus O\), use the first-component
projections \(\operatorname{pr}_Y\) and \(\operatorname{pr}_D\)
of these module splittings. A difference over zero already
belongs to its kernel. The pairing takes values in
\(W_m(\mathbf C_p^\flat)\otimes_{\mathbf Z_p}\mathbf Z_p(1)\).
All localized evaluations are defined there. Indeed
on a chart \(E=p+z^e u\) with \(u\in S^\times\), one has
\[
 E_a\bmod p=\bar z^{ep^a}\overline{\varphi^a(u)}\ne0.
\]
Thus \(E_a\) is a unit in \(W_m(\mathbf C_p^\flat)\), including
on the tame chart \(u=-1\),
and \(\tau\) likewise has nonzero reduction. Thus the map
\(A_{\inf}\to W(\mathbf C_p^\flat)\to W_m(\mathbf C_p^\flat)\)
extends to the localization used for the points. On the underlying
lattices it is scalar extension over \(\mathbf Z_p\). In particular
the functional on Witt coefficients is \(1\otimes\lambda\);
its integrality follows from \(\lambda(O)\subseteq\mathbf Z_p\).
The order estimates place the indicated pairings in \(I\). The third period over
\(p^{e_3}\) gives this multiplier by its extension form.
The smallness order test is applied to points over zero.
The sum of the three orders must be at least \(1\),
the Tate standard order. A difference of Galois
translates has the additional factor \(z\tau\), which
puts the integral Tate value in \(I\): \(\tau\) clears
the Tate basis. Write \(z=[z_0]\) with \(v_\flat(z_0)>0\).
In ordinary Witt coordinates
\[
 z(w_0,\ldots,w_{m-1})=(z_0^{p^r}w_r)_{0\le r<m},
\]
so every coordinate of the product has positive valuation.

Solve
\[
(\varphi-1)y=B(\operatorname{pr}_Y w,\varphi x'-w')
                     +B(\operatorname{pr}_Y x,w'-x').
\]
Frobenius acts on the Witt scalars and fixes a Tate basis.
It commutes with Galois. In each successive quotient for the
filtration by initial zero Witt coordinates, \(\varphi-1\) acts as
\(v\mapsto v^p-v\) on \(\mathbf C_p^\flat\). This map is onto
because the tilt is algebraically closed
\cite[p.~223]{BMSIntegral}. Induction over the \(m\) coordinates
therefore gives \(y\). Write \(a=\operatorname{pr}_Y x\) and
\(b=\operatorname{pr}_D x'\). For \(g\in H\) set
\(A_g=gx-x\), \(D_g=gx'-x'\) and
\(c_{i,g}=p^{e_i}z_i(g)\), so that
\[
 A_g=ga-a+c_{1,g},\qquad D_g=gb-b+c_{2,g}.
\]
Define
\[
 F(g)=B(\operatorname{pr}_Y(gx),gb)
       -g\bigl(B(a,b)-y\bigr)-B(a,D_g)-y
     =B(c_{1,g},gb)-B(a,D_g)+(g-1)y .
\]
Here \(gb\) is the linear action on the kernel of the
second extension. Since \(\varphi\) fixes \(c_{i,g}\) and
commutes with \(g\), expanding \((\varphi-1)F(g)\) and
substituting the equation for \(y\) gives the four terms
\[
\begin{aligned}
 (\varphi-1)F(g)
 &=B(gw-w,g(\varphi x'-w'))\\
 &\quad-B(\varphi x-w,g\varphi x'-\varphi x')\\
 &\quad-B(w-x,gw'-w')+B(gx-x,g(w'-x')) .
\end{aligned}
\]
Each term lies in \(I\) by the indicated path.

For every \(d\) with \((\varphi-1)d\in I\) there is a unique
\(\gamma\in(\mathbf Z/p^m)(1)\) with \(d-\gamma\in I\).
In a Tate basis, filter \(I\) by the number of initial zero Witt
coordinates. On each successive quotient \(\varphi-1\) acts as
\(v\mapsto v^p-v\) on the maximal ideal of
\(\mathcal O_{\mathbf C_p^\flat}\). For \(u\) in this ideal,
Hensel's lemma gives a unique root of \(v^p-v=u\) reducing to
zero. Induction on the finite filtration proves that
\(\varphi-1:I\to I\) is bijective. Subtract from \(d\) the unique
element of \(I\) with the same Frobenius difference. The remainder
lies in the fixed submodule \(W_m(\mathbf F_p)(1)=(\mathbf Z/p^m)(1)\),
since each of its Witt coordinates is fixed by \(v\mapsto v^p\).
Apply this to \(d=F(g)\) and denote the resulting constant by
\(\gamma(g)\).

Use additive notation on the Tate module and the inhomogeneous
differential \((dF)(g,h)=gF(h)-F(gh)+F(g)\). The relations
\(D_{gh}=D_g+gD_h\) and \(c_{1,gh}=c_{1,g}+gc_{1,h}\) give
\[
\begin{aligned}
 (dF)(g,h)
 &= B(c_{1,g},gD_h)-B(c_{1,g},g(hb-b))-B(A_g,gD_h)\\
 &= B(c_{1,g},gc_{2,h})-B(A_g,gD_h).
\end{aligned}
\]
The first term on the last line is
\(p^{e_1+e_2+e_3}(\psi_m)_*(\bar z_1\cup\bar z_2)(g,h)\)
by the definition of \(B\).
The same order test applied to the two cocycle
differences puts the final term in \(I\).
As \(I\) is stable under the Galois action, the identity
modulo \(I\) is exactly
\(p^{e_1+e_2+e_3}(\psi_m)_*(\bar z_1\cup\bar z_2)=d\gamma\)
on the chart subgroup. Since \(F(1)=0\), all \(m\) Witt coordinates
of \(F(g)\) tend to zero as \(g\to1\). To check the action on an
integral tilt element \(v\), let \(\pi_n(v)\in\mathcal O_{\mathbf C_p}/p\)
be its \(n\)-th component. Approximate a lift of this component modulo
\(p\) by an algebraic integer over \(\mathbf Z_p\). Its stabilizer
is open and fixes \(\pi_n(v)\). On this subgroup
\(\pi_n(gv-v)=0\), so \(v_\flat(gv-v)\ge p^n\), with
\(v_\flat(p^\flat)=1\). This proves continuity on integral tilt
elements, and passing to quotients gives it on \(\mathbf C_p^\flat\).
The input cocycles are continuous, the finite Witt operations are
polynomial, and the Tate action modulo \(p^m\) is locally constant.
The maximal ideal of the tilt is open, so \(F(g)\in I\) on a common
neighborhood of \(1\). There \(\gamma(g)=0\) by uniqueness. Choose
an open subgroup in this neighborhood on which \(\bar z_2=0\).
For each \(g\), the
differential identity gives \(\gamma(gh)=\gamma(g)\) for \(h\) in
this subgroup. Hence \(\gamma\) is a continuous cochain.
For a \(\mathbf Z_p\)-basis \((b_i)\) of \(O\) and its dual
coordinate forms \(\lambda_i\), write \(\gamma_{\lambda_i}\) for
these scalar cochains. On the chart subgroup the continuous cochain
\(\Gamma_O=\sum_i b_i\gamma_{\lambda_i}\) satisfies
\[
 d\Gamma_O=p^{e_1+e_2+e_3}
    (\psi\bmod p^m)_*(\bar z_1\cup\bar z_2)
       \quad\text{in }(O/p^mO)(1).
\]
Indeed the Galois action is trivial on the coefficients \(b_i\),
and evaluating the equality at each \(\lambda_i\) gives the
scalar calculation.

\subsubsection{The three extensions at a prescribed level prime}\label{sec:tame-cup-original-periods}

\sectionopening{In this subsection, we apply Section~\ref{sec:tame-cup-frobenius} on the original local group and bound the common exponent.}

For the chosen infinite place \(\tau_\infty\), put
\(p=\ell\mid M\), \(\beta=S_j(\tau_\infty)\), \(W_*=4t\),
\(N=(r+t)r\), \(k_0=r^2(N+1)^2\), and
\(n_0=6(\eta+W_*)\). Both level primes split completely in \(E\)
by Section~\ref{sec:ramified-source}, so \(E_v=\mathbf Q_p\).
Proposition~\ref{prop:marked-level-primes} gives \(p\ge b^{10}\).
Here \(\beta=\kappa_{\tau_\infty}+\delta_j\le b+\delta_j\).
Since \(t=5r\), \(\eta=r^{16}+1\), \(N=6r^2\) and
\(\delta_j\le(2m+1)r^8\), we have
\(\delta_j+n_0+k_0+2=O(r^{16})=o(b)\), uniformly for \(1\le j\le m\).
Thus \(p>\beta+n_0+k_0+2\) for sufficiently large \(x\).
In particular the
dyadic place lies outside \(M\) and is treated in
Appendix~\ref{sec:crystalline-colength-appendix}. Set \(D=L^\dagger\)
for the original summand in the colength estimate. Fix
\(Y\in\mathcal Q_v\), crystalline cocycles \(z_1,z_2\) of
\(Y,D\), and an integral bilinear form \(\psi\) giving one of
the pairings modulo \(p^m\). The third lattice is
\(J=(Y\otimes_O D)^\vee(1)\). Its cocycle is
\[
              z_3(g)=(g\psi-\psi)/p^m .
\]
It is integral by finite-level equivariance, and crystalline
with period primitive \(\psi/p^m\). On these three lattices
\(V_1=Y,V_2=D,V_3=J\), use the split-as-modules extensions
\[
  X_i=V_i\oplus O,\qquad
  g(v,c)=(gv+c\,z_i(g),c).
\]
The constructions are \(O\)-linear in \(c\). For \(X_3\) keep
this parameter even if \(z_3=0\). Over \(\operatorname{Frac}(O)\) the map
\((\lambda,c)\mapsto(\lambda+c\psi/p^m,c)\) identifies \(X_3[1/p]\)
equivariantly with \(J[1/p]\oplus\operatorname{Frac}(O)\), by the
formula for \(z_3\). Thus its rational jump multiset and its linear
period inertia are independent of \(m\) for the fixed \(Y,D\).
On the original lattice the extension form is
\[
 E_3((\lambda,c),y\otimes d)
        =p^m\lambda(y\otimes d)+c\,\psi(y\otimes d).
\]
The integral map \(e_m:X_3\to J\), given by
\(e_m(\lambda,c)=p^m\lambda+c\psi\), is equivariant because
\(p^m z_3(g)=g\psi-\psi\). Its restriction to \(J\) is multiplication
by \(p^m\). For a chart-equivariant de Rham form \(f\) of type at
most \(u\) on \(J\), the form \(p^{-m}f\circ e_m\) on \(X_3\)
extends \(f\) with the same type bound: \(e_m(X_3)\subseteq J\)
and \(p\) is a unit of \(B_{\rm dR}^+\).
Thus the positive Frobenius tests and the order and depth bounds
on kernel vectors apply at every shift with the original form \(f\).
Using the same scalar coordinate map on both lattices,
Section~\ref{sec:period-basis-order-retraction} gives
\(R_J e_m=e_m R_{X_3}\) on the localized modules.
In \(W_m(\mathbf C_p^\flat)\)
all chart denominators used in the three points are units.
Take \(1\in O\) as the third parameter. Modulo \(p^m\) the
\(p^m\lambda\) term vanishes, so a standard lift over
\(p^{k_0}\cdot1\) evaluates to \(p^{k_0}\psi\). When \(z_3=0\),
the Galois-equivariant section \(c\mapsto(0,c)\) gives a standard
lift of \(1\) in the trivial summand, which we multiply by
\(p^{k_0}\). In both cases Frobenius fixes the projected parameter.
Composing with the chosen \(\mathbf Z_p\)-linear functional
therefore gives the same \(B=p^{k_0}\psi_m\) at each move.

Put \(e=\operatorname{lcm}_{1\le a\le r}(p^a-1)\),
\(q=p^{r!}\), and \(K_0=W(\mathbf F_q)[1/p]\).
The unramified field of degree \(\operatorname{lcm}(1,\ldots,r)\) in
Section~\ref{sec:marked-inertia-matrix-span} lies in \(K_0\).
Since \(a\mid r!\) for \(1\le a\le r\), we have \(e\mid q-1\)
and \(\mu_e\subset K_0\). Take \(K'=K_0(\pi)\) with \(\pi^e=p\)
and set \(H=G_{K'}\). The Eisenstein polynomial \(X^e-p\) has
all its roots in \(\mu_e\pi\subset K'\). Since \(K_0/\mathbf Q_p\)
is Galois, so is \(K'/\mathbf Q_p\), and its maximal unramified
subfield is \(K_0\). Thus \(H\) is the pointwise stabilizer in
\(G_{E_v}\) of \(W(\mathbf F_q)\) and \(\pi\). Replacing \(\pi\) by
another root gives the same field since their ratio lies in \(K_0\).
In particular
\[
 [G_{E_v}:H]=[K':K_0][K_0:\mathbf Q_p]=e\,r!,\qquad (p,e\,r!)=1,
\]
since \(p\nmid e\) and \(p>r\). On inertia the character
\[
 I_p\longrightarrow\mu_e,\qquad g\longmapsto g(\pi)/\pi
\]
is surjective with kernel \(I_p\cap H\). Indeed \(I_p\) fixes
\(K_0\) and \(X^e-p\) remains Eisenstein over every finite
unramified extension of \(K_0\).

On the potentially crystalline periods of \(Y\) and \(L\), the
type polynomials of Section~\ref{sec:marked-inertia-matrix-span}
and the definition of \(\mathcal Q_v\) in
Section~\ref{sec:relative-crystalline-colength} give trivial wild
inertia and tame inertia eigenvalues in \(\mu_e\). On
\(D=L^\dagger\) these eigenvalues are inverted. For a basis
\(v_{\rm cyc}\) of \(\mathbf Q_p(1)\) the crystalline period
\(t^{-1}\otimes v_{\rm cyc}\) is Galois invariant, so the Tate
factor contributes the trivial period character. The finite linear
inertia operators are semisimple in characteristic zero and
satisfy \(X^e-1=0\). The kernel of the displayed character in
the procyclic tame quotient consists of the \(e\)-th powers.
Consequently \(I_p\cap H\) acts trivially on the potentially
crystalline periods of both \(Y\) and \(D\).

For \(V=Y[1/p]\) or \(D[1/p]\), take the Galois closure over
\(\mathbf Q_p\) of the composite of \(K'\) with a crystallizing
field of \(V\). The resulting finite Galois extension \(L'/\mathbf
Q_p\) contains \(K'\), and \(V\) remains crystalline over \(L'\).
Let \(L'_0\) be its maximal unramified subfield over \(\mathbf Q_p\).
The inertia subgroup of \(\operatorname{Gal}(L'/K')\) acts trivially
on \(D_{\rm cris,L'}(V)\) by the preceding calculation. Restriction
to \(L'_0\) identifies the quotient by inertia with
\(\operatorname{Gal}(L'_0/K_0)\), since \(K'L'_0/K'\) is the maximal
unramified subextension. The quotient acts semilinearly with
invariants \(D_{\rm cris,K'}(V)\). Galois descent gives
\[
 L'_0\otimes_{K_0}D_{\rm cris,K'}(V)\simeq D_{\rm cris,L'}(V).
\]
Thus \(Y,D\) and their tensor dual \(J\) are crystalline over
\(H\). Each \(z_i|_H\) is a coboundary over the crystalline period
ring and splits the corresponding extension there. Hence all
three \(X_i\) are crystalline over \(H\).
Use the same Kummer coordinate \(z\) with
\(\theta(z)=\pi\), \(E=p-z^e\). Choose an inertia element
\(s\) inducing a generator of the order-\(e\) inertia quotient
of this chart and satisfying \(s(z)=\zeta z\), where \(\zeta\)
has order \(e\). Inertia fixes the unramified chart coefficients.
Here \(\zeta\bmod p\) has that same order.

To choose \(s\), fix compatible roots \(\pi_a\) of
\(\pi\) of degrees \(p^a\), and let \(\zeta_a^{p^a}=\zeta\) in
\(\mu_e\). The polynomial \(X^{ep^a}-p\) is Eisenstein
over every finite unramified field, so inertia takes
\(\pi_a\) to \(\zeta_a\pi_a\). The closed conditions at
successive \(a\) are nested; compactness gives one
\(s\) for all \(a\). Their Teichm\"uller ratios give
\(s(z)=\zeta z\).
The unramified part has trivial different. On adjoining \(\pi\),
the different is generated by \(e\pi^{e-1}\), of normalized
valuation \((e-1)/e\). Its ceiling is the \(\delta\) used above.

Take the first two annihilating polynomials \(q_i\) at \(s\)
over \(\mathbf Z_p\). Each inertia eigenvalue of these types
has order dividing \(p^a-1\) for some \(1\le a\le r\). The
value and its inverse are Teichm\"uller units in an unramified
extension of degree \(a\), hence have monic minimal polynomials
over \(\mathbf Z_p\) of degree at most \(a\). Multiplying these
polynomials over the original upper type gives degree at most
\(N=(r+t)r\). Finite inertia is semisimple in characteristic zero.
This gives the indicated linear-period annihilators, with the
inverse type for \(D\). To pass to the underlying \(\mathbf Q_p\)-period
space, put \(F_O=\operatorname{Frac}(O)\). Over
\(\overline{\mathbf Q}_p\) the scalar algebra \(K_0\otimes_{\mathbf Q_p}F_O\)
splits into factors indexed by pairs \((\tau,\iota)\) of embeddings.
The element \(s\) acts linearly over this algebra. Applying any
\(g\in G_{\mathbf Q_p}\) to the matrix identity on one factor gives
the identity on \((g\tau,g\iota)\), since \(q_i\in\mathbf Z_p[X]\).
For fixed \(g\iota\), crystalline Frobenius transports that identity
to every embedding of \(K_0\): it commutes with \(s\) and acts on
\(K_0\) by a generator of its unramified Galois group. All pairs are
therefore covered, and faithful scalar extension gives the identity
on the original period space.
For the third side take annihilators of \(Y,D\) at \(s^{-1}\).
The characteristic polynomial of multiplication by \(X\otimes X\)
on the tensor product of their two polynomial quotient algebras
is monic of degree at most \(N^2\). By Cayley--Hamilton it kills
the tensor period operator at \(s^{-1}\), hence kills the
dual operator at \(s\).

The Tate factor has trivial linear
inertia on its crystalline period. This is the polynomial
\(q_3\). On \(X_i\) the polynomial \(Q_i=(X-1)q_i\) kills
the period action, by the exact sequence with trivial
quotient. The possible jump sets for the extensions are
the two original intervals with \(0\) adjoined and
\(\{0\}\cup\{-u-u'-1\}\) for the third side.
Here \(u,u'\) range over the two jump sets at the corresponding indices, already
with \(0\) adjoined. If both jumps are low then
\(|u+u'|\le2\eta\). If only one is far then
\(|u+u'|\le\beta+W_*+\eta\), and if both are far, their opposite
signs give \(|u+u'|\le W_*\). Consequently
\[
 |-u-u'-1|\le 1+\max(2\eta,\beta+W_*+\eta)\le\omega,
 \qquad \omega=\beta+n_0 .
\]
The two original jump sets also lie in \([-\omega,\omega]\).
Use this symmetric height for the three integral period bases. All coefficient embeddings
restrict to the same chosen embedding of the local
field \(E_v=\mathbf Q_p\). Thus the union of their jump
sets still has the specified far sign and center on
each side, even if the coefficient ring is ramified.
Choose \(1\in O^\times\) as the coefficient root of unity on all
three extensions. For \(i\in\mathbf F_q^\times\), the integral
projector of Section~\ref{sec:local-standard-orders} acts by the scalar
\[
 e_i=\frac{1}{q-1}\sum_{n=0}^{q-2}[i]^{-n}
 =\begin{cases}1,&i=1,\\0,&i\ne1.\end{cases}
\]
Thus all coefficient embeddings occur at \(i=1\), which Frobenius
fixes. The primitive root \(\zeta\in W(\mathbf F_q)^\times\) of
order \(e\) gives the inertia action \(s(z)=\zeta z\) on the chart
coordinate.

\smallskip
\indent\textit{The exact cokernel calculation on the basis.}
Use the integral crystalline basis on each \(X_i\), and the
coordinate retraction on its localization by all the
\(\varphi^a E\). At level \(l\), the standard orders
are \(0\) at \(a<l\), and \(-u\) at \(a\ge l\), where
\(a\) is Frobenius level and \(u\) the jump. At an index
let \(I_l\subset A_{\inf}\) be the annihilator of the
quotient of the trivial parameter module by the image
of the standard-order module. Put
\(J_l=I_l+(p^{k_0+1})\).
The basis includes \(\tau^\omega(A_{\inf}\otimes X_i)\) in
its standard module for the symmetric height \(\omega\),
where \(\tau=[\underline\epsilon^{1/p}]-1\).

Thus \(\tau^\omega\in J_0\) at every index. Modulo \(p\)
the tilt is a valuation ring with rank-one valuation,
and \(z\) and \(\tau\) both have strictly positive value.
Consequently \(z\) is in the radical of \((p,\tau)\).
There is a least integer \(\nu\ge0\) for which
\(z^\nu\in J_0\) at every index. Multiplying by \(E^\omega\)
passes to level one. Since \(z^e\in(E,p)\), binomial
expansion puts \(z^{\nu+e(\omega+k_0)}\) into \(J_1\).
Indeed a term of the added power has either \(E\)-degree
at least \(\omega\) or \(p\)-degree at least \(k_0+1\).

For level-one lifts use the divided Kummer difference with
divisor \(z\tau E\). For \(j>0\), put \(v=v_p(j)\).
The cyclotomic factorization of \eqref{eq:local-cyclotomic-LI1}
leaves \(v\) factors after dividing by \(z\tau E\), namely
\[
 \frac{(g-1)z^j}{z\tau E}
    =u_j z^{j-1}\prod_{a=1}^{v}\varphi^a E,\qquad
       u_j\in A_{\inf}^{\times}.
\]
An empty product is one. Since \(\varphi^a E=p-z^{ep^a}\),
the product is \(p^v+zv_j\) with \(v_j\in A_{\inf}\).
After multiplying the divided lift by \(u_j^{-1}\) we may
subtract the correction \(v_j z^j\in J_1\).
The chart error modulo \(p^{k_0+1}\) stays integral because
\((g-1)S\subset z\tau E A_{\inf}\).
Thus \(z^j\in J_1\) gives \(p^v z^{j-1}\in J_1\).
Every step down to the nearest multiple
\(p\lfloor n/p\rfloor\) has \(v=0\). Inverse
Frobenius passes from the level-one lift ideal to the
level-zero lift ideal at the preceding index and
divides this last exponent by \(p\). The Frobenius permutation
runs over all indices. Minimality gives
\[
 \nu\le\left\lfloor\frac{\nu+e(\omega+k_0)}{p}\right\rfloor,
 \qquad (p-1)\nu\le e(\omega+k_0)<e(p-1).
\]
In particular \(\nu<e\).

Suppose \(c\ge0\), \(1\le j\le\nu\), and \(p^c z^j\in J_0\).
Choose \(u\in A_{\inf}\) with \(p^c z^j-p^{k_0+1}u\in I_0\)
before choosing the parameter. For each constant parameter \(d\)
take a standard-order point with projection
\((p^c z^j-p^{k_0+1}u)d\). By
Section~\ref{sec:period-basis-order-retraction}, retracting this point
preserves its orders and changes its projection to
\((p^c z^j-p^{k_0+1}F(u))d\). The scalar coordinate
\(F=\operatorname{coord}_0\) is fixed on this chart, so \(F(u)\in S\)
is independent of \(d\). Write \(F(u)=w+zv\), with
\(w\in W(\mathbf F_q)\), \(v\in S\). The inertia element \(s\) fixes
\(w\) and the trivial parameter \(d\). Apply \(Q_i(s)\) to the
retracted point \(x\).
For a crystalline form \(\lambda\), put
\(\lambda_a=s^{-a}\lambda s^a\). In an \(H\)-invariant
period basis \((P_b)_b\), write \(x=\sum_b h_bP_b\). The
values \(\lambda(s^aP_b)\) lie in the unramified constant
field, which \(s\) fixes. Thus
\(\lambda_a(x)=\sum_b h_b\lambda(s^aP_b)\) evaluates the
\(a\)-th power of the linear period operator. The coefficients
\(b_a\) of \(Q_i=\sum_a b_a X^a\) are in \(\mathbf Z_p\) and
commute with the covector. The annihilator identity gives
\(\sum_a b_a\lambda_a(x)=0\) on the lifted point \(x\).
Thus the value on \(Q_i(s)x\) is
\(\sum_a b_a(s^a-1)\lambda_a(x)\).

The retracted standard-order lift belongs to
\(M\cap N_{\rm loc}=N\). For a crystalline covector the formula
\eqref{eq:period-basis-evaluation} at level zero puts
\(E^\omega\lambda_a(x)\) in \(S_B\). Indeed its row on the invariant period basis has entries in
the unramified fixed field. To divide the difference on \(S_B\),
first take \(f\in S\). Since \(s^a(z)=\zeta^a z\) and \(s\) fixes
the Witt constants, \((s^af-f)/z\in S\). Also \(s(E/p)=E/p\), so
the same operator preserves \(S[E/p]\). By
Section~\ref{sec:period-chart-coordinates}, \(S_B\cap A_{\max}\)
is the closure of \(S[E/p]\) in \(A_{\max}\).
Expand any element of that closure as a sum
\(\sum_{n\ge0}p^n h_n\) with \(h_n\in S[E/p]\), by applying
the extended first coordinate to an expansion in \(A[E/p]\).
The divided differences converge term by term in \(A_{\max}\),
as each \((s^ah_n-h_n)/z\) belongs to \(S[E/p]\).
After inverting \(p\), this yields
\(E^\omega\lambda(Q_i(s)x)\in z B^+\). Apply this to the rows
\(L_0\) in \eqref{eq:period-basis-coefficient-orders}.
Since \(s\) preserves \(M\), the basis coordinates of
\(Q_i(s)x\) are integral. Multiplying by the matrix \(Y\) and
using \eqref{eq:period-basis-integral-cancellation} shows that
each is divisible by \(z\) in \(A\). Hence \(Q_i(s)x/z\in M\). Division by \(z\) preserves each
order since \(\theta\varphi^{-a}(z)\ne0\).
Using \(Q_i(1)=0\) and \(s^a(z)=\zeta^a z\), the projection
of \(Q_i(s)x/z\) is
\[
 \left(p^cz^{j-1}Q_i(\zeta^j)
       -p^{k_0+1}\sum_{a=0}^{\deg Q_i} b_a\zeta^a s^a(v)\right)d .
\]
The sum belongs to \(S\) and is the same for every \(d\).
Since these constants generate the parameter module over \(A_{\inf}\),
the multiplier in parentheses belongs to \(I_0\). Therefore
\(p^cz^{j-1}Q_i(\zeta^j)\in J_0\). If \(Q_i(\bar\zeta^j)\ne0\),
then \(Q_i(\zeta^j)\) is a unit in \(W(\mathbf F_q)\), so
\(p^cz^{j-1}\in J_0\). In all cases the Kummer difference
\(g(z)=[\underline\epsilon]z\) also gives
\[
             p^c z^j\in J_0
 \quad\Longrightarrow\quad
             p^{c+v_p(j)+1}z^{j-1}\in J_0 .
\]
Indeed modulo \(z\) one has \(E\equiv p\), and
\(\tau(p^{v_p(j)+1}+zB)=A([\underline\epsilon]^j-1)\)
with \(A,B\in A_{\inf}\). Multiply this identity by
\(p^c z^j\) and divide the lifted Kummer difference
by \(z\tau\); the correction proportional to the
old multiplier lies in the lift ideal.

The monic reduction of \(Q_i\) vanishes on at most
\(\deg q_i+1\) of the distinct \(\bar\zeta^j\) for
\(1\le j\le\nu<e\). Also
\[
 j<e\le\prod_{a=1}^r(p^a-1)<p^{r(r+1)/2}\le p^{r^2}.
\]
Thus \(v_p(j)+1\le r^2\). On the three extensions the total
costs are at most \(r^2\deg Q_i\), respectively, with
\[
 \begin{aligned}
 r^2\deg Q_1,\ r^2\deg Q_2&\le r^2(N+1),\\
 r^2\deg Q_3&\le r^2(N^2+1)\le r^2(N+1)^2=k_0 .
 \end{aligned}
\]
Multiplying each resulting scalar by the unused power of \(p\)
gives \(p^{k_0}\in J_0\)
for each of the three sides. Nakayama applied to the
cyclic image of \((p^{k_0})\) modulo \(I_0\) removes
the modulus \(p^{k_0+1}\). Explicitly, for some \(w\in A\)
we have \(p^{k_0}(1-pw)\in I_0\), and the geometric
series in \(pw\) inverts \(1-pw\) in \(A\). This proves
\[
 p^{k_0} \text{ kills each of the three standard-order
 projection cokernels at each index.}
\]

\smallskip
\indent\textit{Application of the chart depth and cup calculation.}
Use the three cokernel annihilators with \(E\)-exponent
\(0\). The negative-jump exclusion \([-0,0)\) is empty.
There is no nonzero multiple of \(p\) within \(n_0\)
of \(0\) or \(\pm\beta\), since \(p>\beta+n_0\).
The level-zero depth increase from \(-u\) to \(\eta-u\)
uses the integers \(-u\le j<\eta-u\). For \(|u|\le\eta\)
these lie in \([-\eta,2\eta]\). For \(|u|\in[\beta,\beta+W_*]\)
they lie within \(\eta+W_*\) of one of \(\pm\beta\).
Thus each required step lies within \(n_0\) of one of
the three centers and every nonzero step is prime to \(p\).
The contribution at a step \(j\ne0\) is \(v_p(j)+\delta\),
where \(\delta\le1\) is the tame different allowance,
and step \(0\) contributes \(\delta\). Summing over at most
four intervals, each of length \(2n_0+1\), contributes at most
\((8n_0+4)\delta\) per side.

The two points can therefore be taken over the same
\(p^{k_0+(8n_0+4)\delta}\) at every index. A third period over
\(p^{k_0}\) already has the standard orders needed in
Section~\ref{sec:local-tame-moves}. Take \(e_1=e_2=k_0+(8n_0+4)\delta\) and \(e_3=k_0\).
In Section~\ref{sec:tame-cup-frobenius} the two cocycles
are scaled by \(p^{e_1},p^{e_2}\) and the reduced pairing by
\(p^{e_3}\). Their product therefore multiplies the restricted
cup class by \(p^{e_1+e_2+e_3}\), with
\[
 e_1+e_2+e_3
 \le 3k_0+2(8n_0+4)\delta
 \le 3k_0+96(\eta+W_*)+8
 \le 3k_0+150(\eta+W_*+1).
\]
For \(G=G_{E_v}\), put
\(\alpha=p^{e_1+e_2+e_3}[(\psi_m)_*(\bar z_1\cup\bar z_2)]
\in H^2(G,\mu_{p^m})\). The cochain just constructed gives
\(\operatorname{res}_H^G\alpha=0\).
The cocycle for \(\alpha\) and the cochain for its restriction
take values in the finite module \(\mu_{p^m}\). By compactness
and continuity, choose a common open normal subgroup of \(G\),
contained in \(H\), through whose cosets both factor and which
acts trivially on the module. On the resulting finite quotient,
\cite[\href{https://www.jmilne.org/math/CourseNotes/CFTc.pdf\#page=79}{Ch.~II, Proposition~1.30, p.~70}]{MilneCFT}
gives, after inflation,
\[
 [G:H]\alpha=\operatorname{cor}_H^G\operatorname{res}_H^G\alpha=0.
\]
Here \([G:H]=er!\) is prime to \(p\) since \(p>r\) and
\(p\nmid e\). As \(p^m\alpha=0\), this proves \(\alpha=0\).
For the applications at \(p\nmid M\),
Section~\ref{sec:off-level-coefficient-defect} works on
\(H=G=G_{E_v}\) itself. This includes all \(p\mid q_1d\), since
the two primes of \(M\) are \(1\bmod q_2^2\) and \(q_1d\mid q_2^2\).

\begin{proposition}[The tame cup-product exponent]\label{prop:tame-cup-exponent}
Under the tame chart and annihilator hypotheses above, let
\(z_1,z_2\) be the crystalline cocycles of an allowed cup test
defining \(\Delta_v(L)\), with pairing \(\psi\) and depth \(m\).
For each \(\mathbf Z_p\)-linear functional \(\lambda:O\to\mathbf Z_p\),
let \(\psi_m\) be the reduction of \((\lambda\otimes1)\circ\psi\)
modulo \(p^m\). Then
\[
 p^{\,3k_0+150(\eta+W_*+1)}
       [(\psi_m)_*(\bar z_1\cup\bar z_2)]=0
            \quad\text{in }H^2(E_v,\mu_{p^m}).
\]
\end{proposition}

\section{Depth operators and crystalline lifts}
\label{sec:crystalline-lifting-appendix}

\sectionopening{In this appendix, we extend the divided operators from
Appendix~\ref{sec:tame-cup-appendix} to the higher Frobenius levels.
We use the projection ideals and the period basis to construct lifts
with prescribed filtration orders. The construction gives uniform
cup exponents for the colength estimates in
Appendix~\ref{sec:crystalline-colength-appendix}.}

\subsection{Depth operators on the integral chart}
\label{sec:local-offlevel-depth-operators}

\sectionopening{In this subsection, we extend the divided operators of
Section~\ref{sec:local-tame-depth} to \(E_a=\varphi^a(E)\),
including at \(p=2\). Each operator raises the order of the
Galois difference at \(E_a\) by one and preserves the earlier
orders. Its factor of \(p\) enters the annihilator bounds of
Section~\ref{sec:v2-crystalline-lift-ideals}.}

Keep the integral localized basis of
Section~\ref{sec:local-standard-orders}. For the side of any pairing
let \(x\) be a chart point over a constant at index \(j\), and put
\(v_{a,\lambda}(x)=\lambda(\varphi^{-a}x)\).
Depth \(d(a,u)\) means that for each equivariant form \(\lambda\)
on the jump line \(u\) at \(\varphi^{-a}j\),
\[
   v_{a,\lambda}(gx-x)\in\mathrm{Fil}^{d(a,u)}
                    \qquad(g\in H=G_{K'}).
\]
By Section~\ref{sec:local-tame-depth} this is equivalent to
\(v_{a,\lambda}(x)-c\in\mathrm{Fil}^{d(a,u)}\) for some
\(c\in K'\), which may depend on \(a,\lambda,x\).
If the value already has order \(l>0\) and \(d(a,u)>0\), then
\(c\in K'\cap\mathrm{Fil}^1=0\). The new lower bound on the
order of the value is therefore \(d(a,u)\) when \(l>0\) and
\(\min(0,d(a,u))\) otherwise, in addition to its previous bound.

Write \(e=e(K'/\mathbf Q_p)\) for the absolute ramification index
of the chart field. In the off-level application \(K'=E_v\) for
each finite coefficient extension. Put \(D_E=0\) if \(e=1\) and
\(D_E=1+v_p(e)\) otherwise. A depth step from \(n\) to \(n+1\)
at level \(a\) requires
\[
               C(a,n)=v_p(n)+D_E+3a,\qquad v_p(0)=0 .
\]
Here is a construction including the case \(p=2\).
Put \(\pi_b=\theta\varphi^{-b}(z)\). For every integer \(j\)
one can choose \(g_j\in H\) with
\(g_jz=[\underline\epsilon]^jz\). Indeed
\(X^{p^b}-\pi_0\) is Eisenstein over \(K'\). Hence the condition
\(g_j\pi_b=\zeta_{p^b}^{\,j}\pi_b\) cuts out a nonempty closed
subset of \(H\). These subsets are nested in \(b\) by
compatibility of the roots, so compactness gives \(g_j\).
For \(n\ne0\) take \(g=g_{p^a}\). Write \(gE=E\omega\);
then \(gE_a=E_a\varphi^a(\omega)\). Put
\(\Delta_z=z\tau\). The chart action gives
\(\Delta_z\mid\omega-1\), and \(\Delta_z\mid\varphi^a(\Delta_z)\).
Thus \(X=(\varphi^a(\omega)^n-1)/\Delta_z\) is integral,
also for negative \(n\).

Let \(F\) be the Eisenstein minimal polynomial of \(\theta(z)\)
over the unramified coefficient field. Weierstrass division
in \(S\) identifies the kernel of specialization at
\(\theta(z)\) with \(ES\). Thus \(F(z)=E G\) for \(G\in S\).
The reductions of \(F(z)\) and \(E\) both have order \(e\)
in \(z\), so \(G\in S^\times\). The subgroup \(H\) fixes
\(\theta(G)\). Taylor expansion of \(F(gz)\) at \(z\),
followed by division by \(E\) and Fontaine evaluation,
therefore gives
\[
 \theta(\omega-1)\theta(G)
           = F'(\theta z)\,\theta((gz-z)/E).
\]
The derivative has valuation at most \(D_E\). If \(e>1\),
the nonzero derivative terms have distinct valuations
modulo \(\mathbf Z\), and the term of degree \(e-1\) has
valuation \(v_p(e)+(e-1)/e\); for \(e=1\) the derivative is a
unit.

Since \(gz=[\underline\epsilon]^{p^a}z\),
\[
 \frac{gz-z}{E\Delta_z}
   =\frac{\mu}{E\tau}
          \sum_{i=0}^{p^a-1}[\underline\epsilon]^i.
\]
The first factor is a unit of \(A_{\inf}\), and the sum has
Fontaine value \(p^a\). Thus the valuation of the displayed
quotient after applying \(\theta\) is \(a\).
Hence
\[
 v_p(\theta(\omega-1))\le D_E+a+v_p(\theta\Delta_z).
\]
Moreover \(v_p(\theta(\omega-1))\ge v_p(\theta\Delta_z)>
1/(p-1)\), since \(\theta(\tau)=\zeta_p-1\) and
\(\theta(z)\) has positive valuation. If \(v_p(y)>1/(p-1)\),
the linear term \(py\) has uniquely smallest valuation in
\((1+y)^p-1\). Iterating and using the geometric sum for a
prime-to-\(p\) exponent gives
\(v_p((1+y)^n-1)=v_p(y)+v_p(n)\) for \(n>0\).
For \(n=-m<0\), use
\((1+y)^{-m}-1=-((1+y)^m-1)/(1+y)^m\), with unit denominator.
At \(p=2\) this uses \(v_2(\theta\Delta_z)=1+1/e>1\).
Moreover, \(\theta\varphi^{-a}(z)=\pi_a\) and
\(\theta\varphi^{-a}(\tau)=\zeta_{p^{a+1}}-1\), so
\[
 v_p(\theta\Delta_z)-v_p(\theta\varphi^{-a}\Delta_z)
  =(1-p^{-a})\left(\frac1e+\frac1{p-1}\right)\le 2a .
\]
Since \(\theta\varphi^{-a}(\Delta_z)\ne0\), we have
\[
 \begin{aligned}
 v_p(\theta\varphi^{-a} X)
 &=v_p(\theta(\omega)^n-1)-v_p(\theta\varphi^{-a}\Delta_z)\\
 &\le v_p(n)+D_E+a+
       v_p(\theta\Delta_z)-v_p(\theta\varphi^{-a}\Delta_z)
 \le C(a,n).
 \end{aligned}
\]
More precisely, \(G\) and \(\mu/(E\tau)\) are units, so the
derivative identity and the valuation of the geometric sum give
\[
 v_p(\theta\varphi^{-a}X)
 =v_p(F'(\pi_0))+a+v_p(n)
      +(1-p^{-a})\left(\frac1e+\frac1{p-1}\right).
\]
In particular \(\theta\varphi^{-a}X\ne0\), and
\(p^{C(a,n)}/\theta\varphi^{-a}X\in\mathcal O_{\mathbf C_p}\).
Lift this quotient to \(\upsilon\in A_{\inf}\) by
\(\theta\varphi^{-a}\). Since its kernel is \(E_aA_{\inf}\),
we have \(\upsilon X\equiv p^{C(a,n)}\pmod{E_a}\).

At level \(a\) the selected \(g\) fixes \(K'(\pi_a)\), since
\(\theta\varphi^{-a}([\underline\epsilon]^{p^a})=1\).
Write \(v_{a,\lambda}(x)=c+y\), with \(c\in K'\) and
\(y\in\mathrm{Fil}^n\cap\mathcal K_a\).
The element \(\varphi^{-a}(E_a)=E\) has filtration order one.
The residue of \(y/E^n\in\mathcal K_a\cap\mathrm{Fil}^0\)
lies in \(K'(\pi_a)\) by
Section~\ref{sec:period-basis-order-retraction}. Hence
\((g-1)(y/E^n)\in\mathrm{Fil}^1\). Multiplying by \(g(E^n)\)
gives
\[
 (g-1)y\equiv
    \bigl(\varphi^{-a}(g(E_a)/E_a)^n-1\bigr)y
                             \pmod{\operatorname{Fil}^{n+1}}.
\]
Division by \(\varphi^{-a}(\Delta_z)\) has order zero.
Thus the divided difference has value
\(\varphi^{-a}(X)y\) modulo the next grade, for either
sign of \(n\).
Therefore the operator
\[
        p^{C(a,n)}x-R\!\left(\upsilon(gx-x)/\Delta_z\right)
\]
advances the depth. For \(n\ne0\), put
\(y_1=p^{C(a,n)}x-\upsilon(gx-x)/\Delta_z\) and
\(c_1=p^{C(a,n)}c\in K'\). Since \(R(x)=x\), the displayed vector
is \(R(y_1)\). The relative estimate of
Section~\ref{sec:period-basis-order-retraction} gives
\[
 \operatorname{ord}_t\bigl(v_{a,\lambda}(R(y_1))-c_1\bigr)
 \ge \operatorname{ord}_t\bigl(v_{a,\lambda}(y_1)-c_1\bigr)
 \ge n+1.
\]
The selected form thus has Galois differences of order at least
\(n+1\), since \(c_1\) is fixed. The correction projects to zero
by naturality of \(R\), and retains the earlier order bounds on
the other lines by the same retraction inequality.

For \(n=0\) set \(C_0=D_E+3a\) and
\(y=\sum_{j=0}^{p^{C_0}-1}g_jx\). A value of depth zero is
integral since it differs from a constant in \(K'\) by an element
of \(\mathrm{Fil}^0\). Write its residue in \(K'(\pi_a)\) as
\(\sum_{i=0}^{p^a-1} b_i\pi_a^i\), with \(b_i\in K'\).
In \(\mathbf C_p\) its conjugates under the \(g_j\) are obtained
by putting \(\zeta_{p^a}^j\pi_a\) in place of \(\pi_a\).
As \(C_0\ge a\), each embedding occurs \(p^{C_0-a}\) times.
Thus
\[
  \theta(v_{a,\lambda}(y))=
  \sum_{i=0}^{p^a-1} b_i\pi_a^i
      \sum_{j=0}^{p^{C_0}-1}\zeta_{p^a}^{ij}=p^{C_0}b_0\in K'.
\]
Take the coordinate retraction after this sum. Since \(R(x)=x\),
\[
 R(y)=p^{C_0}x+\sum_{j=0}^{p^{C_0}-1}R(g_jx-x).
\]
The retraction preserves congruence of the value to the constant
\(p^{C_0}b_0\) modulo \(\mathrm{Fil}^1\), so all its Galois
differences have depth at least one at this step. Each correction
lies over zero and preserves the old order bounds, as in the
nonzero step (multiply the divided difference by \(\Delta_z\)).
The projection of \(R(y)\) is \(p^{C_0}\) times the original
parameter.

For each \((a,n)\) fix the Kummer elements used above and, for
\(n\ne0\), the lift \(\upsilon\), independently of the chart
vector and of the equivariant form. Each resulting operator is
\(\mathbf Z_p\)-linear, since the Galois action and the coordinate
retraction are \(\mathbf Z_p\)-linear.
For finitely many required depth increases take the union of
their steps \((a,n)\) over the jump lines and perform them in
increasing order of \(n\), in any order among ties. If \(x\) has
orders \(o\) and depth \(d\), every difference \(gx-x\) has orders
at least both \(o\) and \(d\) by equivariance. At every level
\(\varphi^{-a}(z\tau)\) has de Rham order zero. The divided
difference belongs to the localized integral module by
Section~\ref{sec:period-basis-integral-operators}, and the retraction
preserves its orders and its zero projection by
Section~\ref{sec:period-basis-order-retraction}. Thus the correction in either step has orders at least both
\(o\) and \(d\). Its own Galois differences have orders at least \(d\). The resulting point keeps
the old orders and depths on all lines and raises the selected depth
to \(n+1\) on every line where that step is required. If a line is
already deeper, it retains its bound. The projection is multiplied
by \(p^{C(a,n)}\) once for the whole step. Induction gives the sum
of the costs of the distinct steps, independent of the rank.

\subsection{Ideals controlling the crystalline lifts}
\label{sec:v2-crystalline-lift-ideals}

\sectionopening{In this subsection, we bound the annihilator of the cokernel of projection to the trivial
parameter module in Section~\ref{sec:local-standard-orders}. The
divided operators of Section~\ref{sec:local-offlevel-depth-operators}
give a power of \(p\) in the sparse case and an \(E\)-power under the
unramified Hodge bounds.}

We give the calculation of the two types of annihilator used
at a cut. Write \(p=\ell\), \(\mathfrak S=W(\mathbf F_{p^f})
[[z]]\), \(E_a=\varphi^a E\),
\(\tau=[\underline\epsilon^{1/p}]-1\), and \(D_E=0\) if \(e=1\)
and \(\operatorname{ord}_p(e)+1\) otherwise. The coefficient ring
of the chart is a subring of \(A_{\inf}\). Keep the extension
\(0\to Y\to X\xrightarrow{\pi}C\to0\) of one fixed side, in
particular any of the two universal sides or the Tate-dual side.
Localize the period modules \(M_X,N_X\) at all \(E_a\).
By Section~\ref{sec:local-standard-orders} the
\(A_{\inf}\)-module of order \(r\) at \(j\) is
\(e_j\varphi^r M_X\), and its chart vectors form
\(e_j\varphi^r M_X\cap N_{X,\mathrm{loc}}\).
On each tested form of type at most \(u\) at \(j^{p^{-a}}\), the
order required at \(E_a\) is \(\ge0\) for \(a<r\) and \(\ge-u\) for \(a\ge r\).

For a nonnegative integer \(R\) let \(J_r(j;R)\subset A_{\inf}\) be the ideal
of scalars lifting constants at order \(r\), plus \((R)\).
The integer \(R\) in this notation is the modulus. Write
\(R_X,R_C\) for the coordinate retractions.
Concretely, a member \(s\) means that for some \(t\in A_{\inf}\)
and every \(c\in C\), the constant
\((s-Rt)c_j\) is \(\pi(x)\) for some order-\(r\) vector \(x\)
in component \(j\) of the \(A_{\inf}\)-module. The
following properties of the period basis are used here.
The basis spans the order-zero module over \(A_{\inf}\), and
each way between that module and \(A_{\inf}\otimes X\) is
bounded by a power of \(\tau\).

Section~\ref{sec:period-basis-order-retraction} gives
\(R_C\pi=\pi R_X\). Applying this to \(\pi(x)=(s-Rt)c_j\)
with \(s\in\mathfrak S\) gives
\[
 \pi(R_X x)=(s-RF(t))c_j,\qquad F(t)\in\mathfrak S .
\]
The same \(F(t)\) works for every \(c\), and each \(R_X x\)
retains its prescribed orders. Frobenius transports the
standard-order module separately, as in
Section~\ref{sec:local-standard-orders}.

For such a localized chart vector \(y\), the divided difference
\((gy-y)/(z\tau)\) lies in the localized \(A_{\inf}\)-module by
Section~\ref{sec:period-basis-integral-operators}. The divisor
has order zero at every \(E_a\); equivariance therefore preserves
the prescribed lower bounds after division. At order \(1\) we may
also divide by \(E\). At \(a=0\) the residues of the integral
form values lie in \(K'\) by Section~\ref{sec:period-basis-order-retraction}.
The chart subgroup fixes \(K'\), so each value of \(gy-y\) at this
level is in \(\mathrm{Fil}^1\). At \(a>0\) the factor \(E\) has order zero.
These properties apply to each of the three extensions.
They are the properties of the basis used here.
Suppose \(F\in\mathfrak S\cap J_r(j;R)\) and
\(\alpha,\beta,s\in A_{\inf}\) satisfy
\[
 \alpha(gF-F)=D(s+\beta F)
\]
for \(D=z\tau\) at order \(0\), or \(D=z\tau E\) at order
\(1\).

Then \(s\in J_r(j;R)\).
Indeed lift \((F-Rt)c_j\) over \(\mathfrak S\) and apply
\(\alpha(g-1)/D\). As \(g\) fixes \(c_j\), the new lift has
projection \((s+\beta F-R\alpha(gt-t)/D)c_j\).
Here \(t\in\mathfrak S\) by the preceding retraction, so
\(gt-t\in z\tau E A_{\inf}\) by
\eqref{eq:period-basis-divided-action}.
Thus \(s+\beta F\in J_r(j;R)\), and subtracting
\(\beta F\in J_r(j;R)\) proves the assertion.

Choose a chart element \(g\) with
\(g z=[\underline\epsilon]z\).
The cyclotomic factorization for \(n>0\) gives
\[
 \alpha\bigl([\underline\epsilon]^n-1\bigr)
 =\tau E\bigl(p^{\operatorname{ord}_p(n)}+\beta z\bigr).
 \tag{LI1}\label{eq:local-cyclotomic-LI1}
\]
Indeed \( [\underline\epsilon]^n-1\) is a unit times
\(\tau \prod_{a=0}^{\operatorname{ord}_p(n)} E_a\).
For the prime-to-\(p\) part of \(n\) the geometric sum is a
unit. For its \(p\)-part the consecutive cyclotomic factors
give the displayed generators. Each \(E_a\equiv p\bmod z\),
which proves \eqref{eq:local-cyclotomic-LI1}. Applying the difference
operator just established to \(F=p^c z^n\) gives the two
implications
\begin{align*}
 p^c z^n\in J_1(j;R)
  &\ \Longrightarrow\
 p^{c+\operatorname{ord}_p(n)}z^{n-1}\in J_1(j;R),\\
 p^c z^n\in J_0(j;R)
  &\ \Longrightarrow\
 p^{c+\operatorname{ord}_p(n)+1}z^{n-1}\in J_0(j;R).
\end{align*}
The additional factor \(p\) in the second line is
\(E\bmod z\). Put
\(d_{n,m}=\operatorname{ord}_p(n!)-\operatorname{ord}_p(m!)\)
for \(0\le m\le n\). Iterating down from \(n\) to \(m\) increases
the exponent of \(p\) by at most \(d_{n,m}\) at level one and
\(n-m+d_{n,m}\) at level zero.

For the order-one operation on the scalar ideal, take
\(a,c\ge0\) and \(n\ge1\) integral. Use
\(gE=\omega E\), \(z\tau X=\omega^n-1\),
\(\upsilon X\equiv p^{\operatorname{ord}_p(n)+D_E}\pmod E\).
The chart supplies this by choosing \(gz=[\underline\epsilon]z\).
To see the valuation, write \(F_0(z)=E G_0\) with \(F_0\) the
minimal polynomial of the chosen uniformizer over the unramified
field. As in Section~\ref{sec:local-offlevel-depth-operators},
\(G_0\in S^\times\), so \(g\) fixes \(\theta(G_0)\).
Write \(gz-z=Ew\). Taylor expansion and reduction by
Fontaine's map \(\theta\) give
\[
 \theta(\omega-1)\theta(G_0)
      = F_0'(\theta z)\theta(w),\qquad
 \theta(w)=\theta(z\tau)\cdot\text{unit}.
\]
The different bound is \(v_p(F_0'(\theta z))\le D_E\).

For \(e>1\) the nonzero derivative terms have distinct
valuations modulo \(\mathbf Z\); the leading term has valuation
\(v_p(e)+(e-1)/e\). For \(e=1\) the derivative is a unit.
The derivative calculation of
Section~\ref{sec:local-offlevel-depth-operators} applies
also at \(p=2\): \(v_p(\theta(z\tau))>1/(p-1)\)
and \(\omega-1\) is divisible by \(z\tau\). Hence
\(v_p(\theta X)\le D_E+\operatorname{ord}_p(n)\).
Lift a quotient \(p^{D_E+\operatorname{ord}_p(n)}/\theta X\)
by the surjectivity of \(\theta\). Inserting \(\upsilon\) in
the difference operator at order \(1\) trades \(E^n\) for
\(zE^{n-1}\) at contribution \(D_E+\operatorname{ord}_p(n)\).
In detail, if \(F=p^c z^a E^n\), then
\(g z^a/z^a=1+\tau E b\) for some \(b\in A_{\inf}\).
Write \(\upsilon X=p^{D_E+\operatorname{ord}_p(n)}+E k\).
With \(\beta=zk+\upsilon b\omega^n\in A_{\inf}\), we obtain
\[
 \frac{z\upsilon(gF-F)}{z\tau E}
 =p^{c+D_E+\operatorname{ord}_p(n)}z^{a+1}E^{n-1}
   +\beta F .
\]
This is the asserted trade on \(J_1(j;R)\).

Frobenius carries the order-zero ideal at \(j\) to the
order-one ideal at \(\varphi j\). Thus
\(p^c z^{pa}\in J_1(\varphi j;R)\) is equivalent to
\(p^c z^a\in J_0(j;R)\).
No integer between \(p\lfloor n/p\rfloor+1\) and \(n\)
is divisible by \(p\). Thus lowering to that multiple of
\(p\) is free, and subsequent inverse Frobenius gives
\(p^c z^{\lfloor n/p\rfloor}\) in the preceding ideal.

Put \(I=\bigcap_j J_0(j;p^{\kappa+1})\), and let \(b_X\)
bound the absolute jumps of \(X\). For every \(c\in C\), take
\(x_c\in X\) with \(\pi(x_c)=c\).
By \eqref{eq:period-basis-lattice-inclusions} the component at
\(j\) of \(\tau^{b_X}(1\otimes x_c)\) has standard orders and
maps to \(\tau^{b_X}c_j\). Hence
\((\tau^{b_X},p^{\kappa+1})\subseteq I\).
In the rank-one valuation ring \(A_{\inf}/p\), the reduction
\(\bar z\) has positive tilt valuation since \(z\) comes from
compatible roots of a uniformizer. Choose \(m_0\ge1\) with
\(m_0v^\flat(\bar z)\ge b_Xv^\flat(\bar\tau)\).
Then \(z^{m_0}\in(\tau^{b_X},p)\), so
\[
 z^{m_0(\kappa+1)}\in(\tau^{b_X},p^{\kappa+1})\subseteq I .
\]
Let \(a\) be the least nonnegative integer with \(z^a\in I\). If
\(\omega_0\) bounds all jumps above, multiplying the lifts by
\(E^{\omega_0}\) puts them at order \(1\).
The relation \(z^e=uE+vp\) with \(u,v\in\mathfrak S\) then
puts \(z^{a+e(\omega_0+\kappa)}\) in every order-one ideal
modulo \(p^{\kappa+1}\). Indeed, every term of its binomial
expansion contains \(E^{\omega_0}\) or \(p^{\kappa+1}\).
The preceding Frobenius step therefore puts
\(z^{\lfloor(a+e(\omega_0+\kappa))/p\rfloor}\) in every
order-zero ideal. Minimality of \(a\) yields
\[
 (p-1)a\le e(\omega_0+\kappa).
 \tag{LI2}\label{eq:display-li2}
\]

For an unramified chart suppose \(p\ge5\),
\(E=z+p\), and \(\omega_0\le(p-3)(s+t)\).
Set \(\kappa=s+t\). Then
\[
 (p-1)a\le(p-2)\kappa,\qquad
 a+\operatorname{ord}_p(a!)\le\kappa ,
\]
by \((p-1)\operatorname{ord}_p(a!)\le a\).
The order-zero lowering operation puts \(p^{\kappa-m}z^m\)
in \(J_0(j;p^{\kappa+1})\) for \(0\le m\le\kappa\):
when \(m<a\) use the accumulated contribution formula, and when
\(m\ge a\) use ideal closure. In particular the ideal contains
\(p^\kappa\). Modulo \(J_0(j;0)\) this means
\(p^\kappa=p^{\kappa+1}h\). Since \(1-ph\) is a unit,
\(p^\kappa\in J_0(j;0)\), so
\(J_0(j;p^{\kappa+1})=J_0(j;0)\).
Expand \(p^s E^t=p^s(z+p)^t\). Each term is one of those
monomials times its integer binomial coefficient. Thus \(p^s
E^t\) kills the constant cokernel.

For a cut followed by \(J\) indices of jump bound \(\omega_1\),
assume \(8e\le p^{J+1}\).
Use \eqref{eq:display-li2} with the global jump bound \(\omega_0\).
Start with \(z^a\) at the \(J\)-th successor of the cut,
and work backwards through the successors. At each successor
multiply the lifts by \(E^{\omega_1}\) and trade its factors for
\(z\). Put \(d=D_E\omega_1+\operatorname{ord}_p(\omega_1!)\)
and \(R=p^{\kappa+1}\). For each preceding index \(i\) this gives
\[
 p^c z^n\in J_0(\varphi i;R)
 \ \Longrightarrow\ p^{c+d}z^{n+\omega_1}\in J_1(\varphi i;R)
 \ \Longrightarrow\ p^{c+d}z^{\lfloor(n+\omega_1)/p\rfloor}
                         \in J_0(i;R).
\]
The second implication first lowers to the nearest multiple of
\(p\) at level one and then applies inverse Frobenius.
After \(J\) steps the monomial \(p^c z^n\) in the ideal at
the cut satisfies
\[
 c\le J\bigl(D_E\omega_1+\operatorname{ord}_p(\omega_1!)\bigr),
 \qquad n\le a/p^J+\omega_1/(p-1).
\]
Since \(p\ge2\), \eqref{eq:display-li2} gives
\(a/p^J\le 2e\omega_0/p^{J+1}+\kappa/4\).

Removing the remaining \(z^n\) requires
\(n+\operatorname{ord}_p(n!)\le2n\). Therefore \(p^\kappa\)
lies in the ideal modulo \(p^{\kappa+1}\) if \(c+2n\le\kappa\).
As above, \(1-ph\) is a unit, so the same power lies in the
exact ideal. Indeed Legendre's inequality and the preceding bounds give
\[
 c+2n\le JD_E\omega_1+(J+2)\frac{\omega_1}{p-1}
                +\frac{4e\omega_0}{p^{J+1}}+\frac\kappa2 .
\]
Thus it suffices that
\[
 2JD_E\omega_1+(2J+4)\frac{\omega_1}{p-1}
                +\frac{8e\omega_0}{p^{J+1}}\le\kappa .
 \tag{LI3}\label{eq:display-li3}
\]
For the stated cut parameters the global bound is \(2b_0\)
and the bound for a thin next component is \(G_1\) in the cut-size
conditions. They impose
\(D_EG_1\le X_0\), \(G_1/(p-1)\le X_0\), and
\(2eb_0/p^{J+1}\le D_0\).
The left side of \eqref{eq:display-li3} is at most
\((4J+4)X_0+8D_0\), below
\(\kappa=12((J+1)X_0+D_0)+1\).

At a shifted cut the unramified argument applies with
\(s=\kappa\), \(t=s_j\), since the global jump bound is
at most \((p-3)(\kappa+s_j)\). This proves the two
annihilator assertions on all three sides. It also gives
the annihilator with no \(p\)-factor \(E^A\) in the separated-interval
regime, using the jump bound \((p-3)A\).

The basis used for these divisibilities is built in
Sections~\ref{sec:local-integral-period-frame}
and~\ref{sec:local-standard-orders}.
Over the finite crystallizing field \(K'\), whose residue field
is perfect, the underlying \(\mathbf Z_p\)-lattice of a crystalline
representation is the \'etale realization of a prismatic
\(F\)-crystal on \(\operatorname{Spf}\mathcal O_{K'}\), also
when \(p=2\). The chosen generator \(E=p+z^e u\) differs
from the Eisenstein polynomial of \(\theta(z)\) by a unit of
\(\mathfrak S\), as in Section~\ref{sec:local-offlevel-depth-operators}.
Thus \((\mathfrak S,(E))\) is the Breuil--Kisin prism of the chosen
uniformizer. Evaluation gives a finite free \(\mathfrak S\)-module
\(\mathfrak M\) with
\((\varphi^*\mathfrak M)[1/E]\simeq\mathfrak M[1/E]\)
\cite[Examples~2.6(1), 4.3, Theorem~5.6 and \S7]{BhattScholzeCrystalline}.
Evaluation on this prism is fully faithful on crystalline lattices
\cite[Theorem~7.9, p.~31 of the authors' version]{BhattScholzeCrystalline}. For finite free
Breuil--Kisin modules with linearization invertible over
\(\mathfrak S[1/E]\), base change to \(\mathfrak S[1/E]^{\wedge}_p\)
is fully faithful by \cite[Theorem~7.2, p.~28 of the authors' version]{BhattScholzeCrystalline}.
Kisin states this for integral linearizations with cokernel killed by
a power of \(E\) \cite[\S\S2.1.1, 2.1.3 and Proposition~2.1.12]{KisinCrystalline}.
His completed ring is \(\mathfrak S[1/z]^{\wedge}_p\).
Modulo \(p^m\), \(m\ge1\), the nilpotence of \(p\) and
\(E=p+z^e u\), \(u\in\mathfrak S^\times\), give the same
localization upon inverting \(E\) or \(z\). Taking inverse limits
identifies the two completed rings. For a pair of modules, multiplying both linearizations by the same
sufficiently high power of \(E\) gives that formulation and preserves
the Frobenius intertwining equation.

One can identify the integral model of \(T\) in these statements
directly. Let \(N_T=\sum_j\mathfrak S n_j\) be the chart lattice and
\(\Phi:\varphi^*N_T\to N_T[1/E]\) its linearized Frobenius.
If the jumps lie in \([-b,b]\), then
\[
 E^bN_T\subseteq\Phi(\varphi^*N_T)\subseteq E^{-b}N_T.
\]
The right inclusion is \eqref{eq:period-basis-frobenius-matrix}.
For \(x\in N_T\), the adapted dual basis of
Section~\ref{sec:local-standard-orders} gives
\(\operatorname{ord}_t\lambda(x)\ge -u_s\ge-b\) for every basis
form of type \(u_s\). Multiplication by \(E^b\) therefore gives
the level-one standard orders at \(a=0\). At \(a>0\) the
level-zero and level-one bounds coincide and
\(\theta\varphi^{-a}E\ne0\). Thus \(E^bx\) lies in
\(\varphi M_T\cap N_{T,\mathrm{loc}}\), which on summing the
projected equalities in that section is
\(\sum_j\mathfrak S\,\varphi(n_j)=\Phi(\varphi^*N_T)\).
In particular \(E^b\Phi\) has integral cokernel killed by \(E^{2b}\).

Put \(K'_\infty=\bigcup_{a\ge0}K'(\theta(\varphi^{-a}z))\)
and \(W_\flat=W(\mathbb C_p^\flat)\). The group \(G_{K'_\infty}\)
fixes \(z\), the Witt constants, and hence \(S_B\) by continuity.
In the column identity \(E^b n=PZ\) of
Section~\ref{sec:v2-crystalline-frame-action}, \(P\) is invariant
for the diagonal action and \(Z\) has entries in \(S_B\).
Hence \(G_{K'_\infty}\) fixes each \(n_j\) for that action.
The column matrix \(n\) is invertible over \(W_\flat\), so the
inclusion of its columns identifies
\[
 W_\flat\otimes_{\mathfrak S} N_T
      =W_\flat\otimes_{\mathbf Z_p}T
\]
with its Frobenius and \(G_{K'_\infty}\)-action. Its
Frobenius-fixed vectors form the original \(T|_{G_{K'_\infty}}\),
since \(W_\flat^{\varphi=1}=\mathbf Z_p\).
For the convention \(\widetilde\theta=\theta\varphi^{-1}\) on
the perfect prism, use \(\widetilde\iota=\varphi\circ\iota\) for the
chart embedding \(\iota:\mathfrak S\to A_{\inf}\)
\cite[\S5 and Notation~7.1]{BhattScholzeCrystalline}. The map sends
the chart coordinate to \(z^p\) and acts on the Witt constants by
Frobenius, as in \cite[p.~281, before Lemma~4.30]{BMSIntegral}.
Thus \(\widetilde\theta\circ\widetilde\iota=\theta\circ\iota\).
Let \(j:N_T\to W_\flat\otimes_{\mathbf Z_p}T\) be the inclusion
given by the columns \(n_j\) in the chart \(\iota\). In the shifted
chart the \(W_\flat\)-linear identification is
\[
 W_\flat\otimes_{\widetilde\iota,\mathfrak S}N_T
 \longrightarrow W_\flat\otimes_{\mathbf Z_p}T,\qquad
 a\otimes y\longmapsto a\,\varphi(j(y)),
\]
where \(\varphi\) acts on the coefficients of \(T\).
Indeed \(\varphi(j(sy))=\widetilde\iota(s)\varphi(j(y))\) for
\(s\in\mathfrak S\), and its column matrix \(\varphi(n)\) is invertible
over \(W_\flat\). The map intertwines Frobenius since \(j\) does
after inverting \(E\), and intertwines \(G_{K'_\infty}\) since this
group fixes the columns \(n_j\) and commutes with \(\varphi\).
It commutes with the \(O\)-action since the endomorphisms of \(T\)
giving that action are \(\mathbf Z_p\)-linear. This is the integral
identification of \cite[Theorem~4.4, p.~265]{BMSIntegral}.
The equivalence for the completed \'etale ring used in
\cite[Corollary~3.8, proof of Corollary~7.10 and Remark~7.11]{BhattScholzeCrystalline}
therefore identifies the \'etale realization of \(N_T\) with that of
the evaluated prismatic crystal. Theorem~7.2 lifts this identification
and its inverse uniquely to \(\mathfrak S\). The resulting
Breuil--Kisin--Fargues module is
\[
 A_{\inf}\otimes_{\widetilde\iota,\mathfrak S}N_T
     =\varphi^*(A_{\inf}\otimes_{\iota,\mathfrak S}N_T).
\]
Its linearization is invertible after inverting
\(\widetilde\xi=\varphi(\xi)\), since
\(\widetilde\iota(E)A_{\inf}=\widetilde\xi A_{\inf}\).
This is the scalar extension of
\cite[Definition~4.22, p.~276, and Proposition~4.32, p.~281]{BMSIntegral}.
Its integral \'etale realization is described in
\cite[Lemma~4.26, p.~278]{BMSIntegral}.
The numerical bounds \eqref{eq:display-li2}--\eqref{eq:display-li3}
follow from the coordinate divisibilities above.

\subsection{The convergence radius of the period basis}
\label{sec:v2-crystalline-frame-radius}

\sectionopening{In this subsection, we prove a generation bound for positive
Frobenius translates of the basis in
Section~\ref{sec:local-integral-period-frame}. The stronger de
Rham orders at \(E=0\) supply the intermediate lifts in
Section~\ref{sec:v2-crystalline-stops}.}

Use the chart and period comparison of
Sections~\ref{sec:period-chart-coordinates}
and~\ref{sec:period-period-contraction}. Write
\(\Lambda=W(\mathbf F_{p^f})\), \(S=\Lambda[[z]]\),
\(A=A_{\max}\), \(B^+=A[1/p]\), and \(S_B\subset B^+\)
for the closure of \(S[E/p][1/p]\) in the bounded-ring topology.
The extended first coordinate \(\mathcal R:B^+\to S_B\)
is the identity on \(S_B\) and is \(S_B\)-linear.
It restricts to \(\operatorname{coord}_0:A_{\inf}\to S\)
on \(A_{\inf}\). It is compatible with approximation by finite
chart coordinates.

Here are two scalar consequences. If \(\nu\in S\) and
\(\xi-p b_0\in\nu A_{\inf}\) for some \(b_0\in A_{\inf}\),
then for \(r\in A\) and \(M\ge0\) there are \(s\in S\),
\(y\in S_B\) with
\begin{equation}
 p^M\mathcal R(r)=s+\nu^M y.
 \label{eq:v2-radius-approx-coordinate}
\end{equation}
To see this, express \(r\) as a \(p\)-convergent sum of
polynomials in \(\xi/p\) with integral coefficients and
expand in \(\xi/p-b_0\). After multiplication by \(p^M\)
a term of degree \(k<M\) becomes a multiple of
\(p^{M-k}(\xi-pb_0)^k\) in \(A_{\inf}\). For \(k\ge M\) use
\[
 p^M(\xi/p-b_0)^k=(\xi-pb_0)^M(\xi/p-b_0)^{k-M}.
\]
Writing the original expansion as \(\sum_{i\ge0}p^i h_i(\xi/p)\),
with \(h_i\in A_{\inf}[X]\), retains the factor \(p^i\) in each
of the two sums, including after extracting \((\xi-pb_0)^M\)
in the second sum.
Apply \(\mathcal R\). This applies to \(\nu=E\) with
\(b_0=0\). It applies to \(\nu=z\) because
\(E=p+z^e u\) and \(E\) is associated to \(\xi\).
It follows, using \(\nu=z\), that every element of
\(S_B\) is constant in \(\Lambda[1/p]\) modulo \(zS_B\).

For \(g\) fixing \(\Lambda\) and \(\theta(z)\), one has
\begin{equation}
 z\tau\mid gy-y\quad\text{in }B^+,\qquad y\in S_B .
 \label{eq:v2-radius-galois-difference}
\end{equation}
Indeed the bounded integral calculation in
Section~\ref{sec:period-basis-integral-operators} gives
\((g-1)(S_B\cap A_{\max})\subset z\tau A_{\max}\).
Inverting \(p\) gives \eqref{eq:period-basis-divided-action} on \(S_B\).
The cancellation \eqref{eq:period-basis-integral-cancellation} gives the contraction
used below: if \(z^r\tau\mid E^m x\) in \(B^+\)
for \(m,r\ge0\) and \(x\in A_{\inf}\), then
\(z^r\tau\mid x\) in \(A_{\inf}\).

Suppose \(Z\) has entries in \(S_B\),
\(Z' Z=E^h I\), with \(Z'\in M_d(B^+)\), and \(tZ=ZG\) for
\(G\in M_d(A_{\inf})\), where \(t\) fixes \(K'=\theta(S_B)\) and
\(t z=[\underline\epsilon]z\).
Apply \(\mathcal R\) entrywise to \(Z'\), which preserves
\(Z'Z=E^hI\). The constant term of \(p^{-h}\mathcal R(Z')\)
modulo \(z\) gives a matrix \(B_0\) over
\(\Lambda[1/p]\) with
\[
 Y=B_0 Z=I+zW,\qquad W\in M_d(S_B).
\]
Indeed \(p^{-h}E^h\equiv1\bmod z\), and multiplication
by \(Z\) preserves the discarded ideal. Moreover
\(E^h(G-I)=Z'(tZ-Z)\), hence \(G-I\) is divisible by
\(z\tau\) over \(A_{\inf}\), by
\eqref{eq:v2-radius-galois-difference} and contraction.
The relation \(tY=YG\) holds since \(t B_0=B_0\).

For clarity the successive approximation proving the
radius starts with the period coefficients \(Y_i\). Suppose
\[
 A_n=q_n Y_i-z^{n+1}W_n\in A_{\inf}^d,\qquad
 W_n\in S_B^d,\quad q_n\in\mathbf N .
\]
The initial values satisfy \(q_0=1,\ W_0=W_i,\ A_0=I_i\).
The period term cancels from \(A_nG-tA_n\), leaving
\[
 z^{n+1}\bigl(([\underline\epsilon]^{n+1}-1)W_n
 +[\underline\epsilon]^{n+1}(tW_n-W_n)-W_n(G-I)\bigr).
\]
The last two terms in parentheses are divisible by \(z\tau\).
For \(v_n=v_p(n+1)\), the cyclotomic factorization
\eqref{eq:local-cyclotomic-LI1} and \(E=p\bmod z\) give a multiplier
\(\alpha\in A_{\inf}\) with
\[
 \alpha([\underline\epsilon]^{n+1}-1)
           =\tau(p^{v_n+1}+z u),\qquad u\in A_{\inf}.
\]
The numerator \(\alpha(A_nG-tA_n)\) has entries in
\(A_{\inf}\). Its divisibility by \(z^{n+1}\tau\) in
\(B^+\) therefore puts the quotient in \(A_{\inf}^d\) by contraction. In \(B^+\) it equals
\(p^{v_n+1}W_n+zV_n\). Applying \(\mathcal R\) produces
\(s_n\in S^d\) and \(W_{n+1}\in S_B^d\) such that
\[
 p^{v_n+1}W_n=s_n+zW_{n+1}.
\]
Consequently
\(A_{n+1}=p^{v_n+1}A_n+z^{n+1}s_n\) and
\(q_{n+1}=p^{v_n+1}q_n\) obey the same recurrence.

Take \(N=p^l\ge3e\). Since \(z^e/p=(E/p-1)u^{-1}\in S[E/p]\),
the value \(X=z^N/p^2=pz^{N-3e}(z^e/p)^3\) is in
\(pS[E/p]\subset pA\). Thus evaluation
\(\mathrm{ev}:S[[X]]\to A\) is defined by \(p\)-adic
convergence. Put \(a_n=n-v_p(n!)+2\). Applying \(\varphi^l\)
to the recurrence for the period coefficients gives
\[
 p^{a_n}\varphi^l W_n
   =p^{a_{n+1}-2}\varphi^l s_n
      + X\,p^{a_{n+1}}\varphi^l W_{n+1}.
\]
Indeed \(v_p((n+1)!)\le n+1\) and
\(a_{n+1}+v_n+1=a_n+2\). In particular at \(p=2\),
\(a_{n+1}-2\) is the sum of the binary digits of \(n+1\),
so the coefficients still lie in \(S\).
Put \(c_n=p^{a_{n+1}-2}\varphi^l s_n\in S^d\).
For each \(m\ge1\), finite iteration gives
\[
 p^2\varphi^l W_0-\sum_{n=0}^{m-1}X^n c_n
       =X^m p^{a_m}\varphi^l W_m\in X^m(B^+)^d .
\]
The tail of \(\sum_{n\ge0}X^nc_n\) belongs to \(X^m A^d\),
since \(X\in pA\). Thus the difference between the full sum and
\(p^2\varphi^l W_0\) lies in \(X^m(B^+)^d\) for every \(m\).

To see that \(B^+\) is \(z\)-adically separated, write
\(\xi=u_0E\) with \(u_0\in A_{\inf}^\times\) and set
\(X_0=\xi/p\). In \(B^+\) the element \(z^e\) is a unit
multiple of \(X_0-u_0\), since \(E=p+z^e u\).
If a nonzero element were divisible by every power
of \(z\), rescale by a power of \(p\) to obtain an
element of \(A_{\max}\setminus pA_{\max}\) with
the same divisibility. Each power of \(X_0-u_0\)
has monic reduction in \((A_{\inf}/(p,\xi))[X_0]\).
By the \(p\)-saturation proved in
Section~\ref{sec:period-period-contraction}, division
by that power takes place in \(A_{\max}\). Reduction
would make a nonzero polynomial divisible by monic
polynomials of arbitrarily high degree. This proves
separatedness.

Applying this construction to every row of \(W\) gives a matrix \(F\) over \(S[[X]]\)
and
\[
 \varphi^l(B_0)\varphi^l(Z)=\mathrm{ev}(I+XF).
\]
Also \(\varphi^l(E)=p\,\mathrm{ev}(\varepsilon)\) with
\(\varepsilon=1+pX\varphi^l((E-p)/z)\).
The series \(I+XF\) and \(\varepsilon\) are invertible over
\(S[[X]]\), with constant terms \(I\) and \(1\).
Thus for any \(b\ge0\),
\[
 \varphi^l(Z)=\varphi^l(E)^b\,C\,\mathrm{ev}(V),\qquad
 C=(p^b\varphi^l(B_0))^{-1},\qquad
 V=\varepsilon^{-b}(I+XF).
\]
The matrix \(C\) is fixed by every \(g\) fixing \(\Lambda\).
Its inverse exists because \(\varphi^l(B_0)\varphi^l(Z)\)
is invertible over \(B^+\).

For integers \(A_0,\delta\ge0\) with \(e A_0\le N\delta\),
approximate each formal series
\(H\in S[[X]]\) by its truncation of degree \(<\delta\).
After multiplying by \(p^{2\delta}\), the truncated part
evaluates into \(S\), and the tail is \(z^{N\delta}\)
times the evaluation of another series, so lies in
\(p^{A_0}A\). In fact \(\mathrm{ev}(H)\in S_B\), since \(X\in pS[E/p]\)
and the partial sums converge in the \(p\)-adic completion of
that ring. Write the tail after clearing the denominator as
\(p^{A_0}r\), with \(r\in A\). This tail is fixed by
\(\mathcal R\), as both the full sum and the truncated polynomial
are in \(S_B\). Hence \(p^{A_0}r=p^{A_0}\mathcal R(r)\).
Applying \eqref{eq:v2-radius-approx-coordinate} with
\(M=A_0\) and \(\nu=E\) gives
\[
  p^{2\delta}\mathrm{ev}(H)\equiv s\pmod{E^{A_0}B^+},
  \qquad s\in S.
\]
Apply this entrywise to \(V^{-1}\) and \(V\), obtaining
matrices \(M,M'\in M_d(S)\) congruent modulo \(E^{A_0}B^+\) to
\(p^{2\delta}\operatorname{ev}(V)^{-1}\) and
\(p^{2\delta}\operatorname{ev}(V)\), respectively. In that quotient,
\[
 \mathrm{ev}(V)M=p^{2\delta}I,\qquad
 \mathrm{ev}(V)(p^{4\delta}I-MM')=0.
\]
The columns of \(M\) and \(p^{4\delta}I-MM'\) have images
\(p^{2\delta}e_i\) and \(0\), respectively, under
\(\operatorname{ev}(V)\) in this quotient, where \(e_i\) is the
corresponding coordinate vector. For every \(y\in S^d\),
\[
 p^{4\delta}y=(p^{4\delta}I-MM')y+M(M'y),
\]
so these columns span a module containing \(p^{4\delta}S^d\).

For the period basis of a lattice, take \(b\) to be the
nonnegative uniform jump bound used in choosing the invariant
period matrix \(P\). The basis construction in
Section~\ref{sec:local-integral-period-frame} supplies \(Z'Z=E^{2b}I\)
and \(tZ=ZG_t\) with \(G_t\) integral. Choose the element \(t\)
with \(tz=[\underline\epsilon]z\) as in
Section~\ref{sec:local-offlevel-depth-operators}.
The basis vectors \(n_\nu=\tau^{-b}Q_\nu\) satisfy the
formula for every equivariant de Rham form \(\lambda\)
and column \(y\) of basis coordinates:
\[
 \lambda\!\left(\sum_\nu y_\nu\varphi^l n_\nu\right)
 =\bigl(\lambda_{\mathrm{lat}}\varphi^l P\bigr)
       \,\frac{\varphi^l Z}{\varphi^l(E)^b}\,y .
\]
This follows from \(UZ=E^{2b}Q\) and
\(U=(E\tau)^b P\), in the given integral basis.
Write \(b_\lambda=\lambda_{\mathrm{lat}}\varphi^l(P)C\).
Its entries belong to \(K'\subset B_{\mathrm{dR}}^+\) by
\eqref{eq:period-basis-evaluation} and the invariance of \(C\).
For a column \(y\) with \(\mathrm{ev}(V)y=r+E^{A_0}w\),
\(r\in\mathbb Z^d\), \(w\in(B^+)^d\), we obtain
\[
  \lambda\!\left(\sum_\nu y_\nu\varphi^l n_\nu\right)
  = b_\lambda r+E^{A_0}b_\lambda w
  \equiv b_\lambda r\in K'
  \pmod{\mathrm{Fil}^{A_0}B_{\mathrm{dR}}}.
\]
Here \(\mathrm{Fil}^{A_0}B_{\mathrm{dR}}=E^{A_0}B_{\mathrm{dR}}^+\),
since the image of \(E\) has de Rham order one. The span statement therefore gives the precise
\(p^{4\delta}\) radius property at level \(l\) used in
the long-gap coefficient extraction. For an index \(j\), apply
\(e_j\) to the resulting span. It preserves the translated chart
module by Section~\ref{sec:local-standard-orders}. For a generator
\(w\), the value \(\lambda(e_jw)=(\lambda\circ e_j)(w)\) is again
constant modulo \(\mathrm{Fil}^{A_0}\), since the chart group fixes
the Witt scalars of \(e_j\) and its coefficient action is equivariant.
Thus the projected generators span a module containing
\(p^{4\delta}e_j\sum_\nu S\varphi^l(n_\nu)\).

\subsection{Crystalline lifts and termination of the construction}
\label{sec:v2-crystalline-stops}

\sectionopening{In this subsection, we apply the annihilator estimates of
Section~\ref{sec:v2-crystalline-lift-ideals} to construct lifts of
constants at selected coefficient summands. Between successive
summands we transport the lifts by Frobenius. The radius estimate from
Section~\ref{sec:v2-crystalline-frame-radius} supplies the additional
depth needed for the longer intervals.}

Use the period basis and its depth operators on each
of the three sides of the pairing. At an index \(j\) let
\(\beta_j\) be the far start, let \(\eta\ge1\) be the low radius,
let \(V\) be the width of the separated intervals, and fix
\(W_0\ge V+4\eta\), \(J\ge1\), \(n_0\), \(b_0\), \(G_1\),
\(X_0\) and \(D_0\). The numerical requirements, verified
for the two schemes below, are
\[
\begin{gathered}
 \beta_j+s_j+W_0\le 2b_0,\quad
 s_j+W_0\le pn_0,\quad
 G_1\le(p-1)n_0,\quad 2eb_0\le p^{J+1}D_0,\\
 n_0+\lfloor\log_p(2b_0)\rfloor
 +(J+1+D_E)W_0+D_EG_1\le X_0,\\
 \kappa=12((J+1)X_0+D_0)+1 .
\end{gathered}
\]
At a cut we apply one of the two estimates in
Section~\ref{sec:v2-crystalline-lift-ideals}. Either its shift
\(s_j=0\) and all three sides at the next \(J\) indices have upper
jump bound \(G_1\), with
\(8e\le p^{J+1}\), or the chart is unramified, \(p\ge5\)
and \(s_j<\beta_j\) with all jumps bounded above by
\((p-3)(\kappa+s_j)\).
In either case \(p^\kappa E^{s_j}\) kills all three
standard projection cokernels. The third-side jumps are \(0\) and \(-u-v-1\), for a jump
\(u\) from the first side and \(v\) from the second side.
If both are in their low intervals then \(|u+v+1|\le2\eta+1\).
If one is far and one low then \(|u+v+1|\le\beta_j+V+\eta+1\).
If both are far, their opposite signs give \(|u+v+1|\le V+1\).
Since \(W_0\ge V+4\eta\) and \(\eta\ge1\), all three sides have
absolute jump bound \(\beta_j+W_0\) at index \(j\).

At a cut \(j\) let \(l\le2J+1\) be the distance to the
preceding cut and \(s_j\) its shift. The annihilator
\(p^\kappa E^{s_j}\) gives a point over \(p^\kappa c_j\)
by lifting the multiple of the constant, applying the coordinate
retraction, and dividing by \(E^{s_j}\). The retraction fixes the
projected multiple. The resulting orders and depths are at least
\(-u-s_j\mathbf1_{a=0}\). Increase the depth to
\[
 D(a,u)=-u+s_j\mathbf1_{a=0}+\eta\mathbf1_{a\le l}.
\]
Call the improved orders \(C_j(a,u)\), obtained from
\(-u-s_j\mathbf1_{a=0}\) and \(D(a,u)\) by the constant-term
rule of Section~\ref{sec:local-offlevel-depth-operators}. Explicitly, if
\(I_j(a,u)=-u-s_j\mathbf1_{a=0}\), then
\[
 C_j(a,u)=
 \begin{cases}
 D(a,u),&I_j(a,u)>0,\\
 \min(0,D(a,u)),&I_j(a,u)\le0 .
 \end{cases}
\]

At the cut \(j\) and level \(0\le a\le l\), put
\(i=F^{-a}j\), where \(F(i)=i^p\), and \(r_0=s_j\mathbf1_{a=0}\).
On either moving side a filtration step \(n\mapsto n+1\)
on a line of jump \(u\) has
\[
  -u-r_0\le n<-u+r_0+\eta.
\]
For \(|u|\le\eta\) this gives \(|n|\le r_0+2\eta\).
For \(u\in\epsilon_i[\beta_i,\beta_i+V]\), with \(\epsilon_i\)
the sign of the far interval at \(i\), it gives
\(|n+\epsilon_i\beta_i|\le r_0+\eta+V\).
Thus the two integer intervals have centers \(0,-\epsilon_i\beta_i\)
and common radius \(R=r_0+W_0\le pn_0\), since
\(W_0\ge V+4\eta\). Their absolute endpoint bound is
\(\beta_j+s_j+W_0\le2b_0\) at \(a=0\) and
\(\beta_i+W_0\le2b_0\) at \(a>0\).
There are no steps for \(a>l\), and \(l\le2J+1\).
If \(r_0>0\), then \(a=0\) and \(D_E=0\).
The sum of \(v_p(n)\) for nonzero \(n\) in either interval
of radius \(R\) is at most
\(2R/(p-1)+\lfloor\log_p(2b_0)\rfloor\), where
\(2b_0\) bounds all absolute endpoints.

Indeed, for either such interval \(I\) this sum equals
\[
 \sum_{q=1}^{\lfloor\log_p(2b_0)\rfloor}
    \#\{n\in I\cap\mathbb Z:n\ne0,\ p^q\mid n\}.
\]
Each inner count is at most \(2R/p^q+1\). Thus the two intervals
contribute at most \(4R/(p-1)+2\lfloor\log_p(2b_0)\rfloor\)
from the valuations. Since \(R\le pn_0\) and \(p\ge2\), the first
term is at most \(8n_0\). If \(r_0>0\), the constant part
\(D_E+3a\) of the cost is zero. Thus both intervals on one side at one level have total
exponent at most
\[
 8 n_0+2\lfloor\log_p(2b_0)\rfloor
       +30(J+1+D_E)W_0\ \le\ 30 X_0 .
\]
For the last term use \(4W_0+2\le5W_0\) and
\(D_E+3a\le6(J+1+D_E)\) when the shift is zero.
The size conditions
\[
 n_0+\lfloor\log_p(2b_0)\rfloor+(J+1+D_E)W_0
          +D_E G_1\le X_0,\qquad
 12((J+1)X_0+D_0)<\kappa
\]
give \(\kappa+4D_0+(2J+2)30X_0<25\kappa\).
Hence on each moving side the standard points may be taken
over \(p^{25\kappa}c_j\), using any unused powers of \(p\) as scalars.
The Hodge intervals determine the cuts and their distances on
the finite Frobenius cycle. Fix the two parameters and choose
one such point over each scaled parameter at every cut. Each
segment uses these chosen points at its two ends.
The Tate-dual side uses its own projection annihilator as in the
order criterion below.
At \(a=0\), the absence of a jump in \([-s_j,-\eta)\)
shows that \(C_j(0,u)\ge-u-\eta\):
for \(u<-s_j\), the initial order is positive and improves
to \(D(0,u)\); otherwise \(u\ge-\eta\), the required
order is nonpositive and is bounded by
\(\min(D(0,u),0)\). At all other levels \(C_j(a,u)\ge-u\).

\subsubsection{Termination at a plain cut}
\label{sec:lifting-plain-cut}

\sectionopening{In this subsection, we use the ideal bounds of Section~\ref{sec:v2-crystalline-lift-ideals} for a plain cut.}

For the stop at a plain cut, choose a sign such that
for every index at distance \(a<l\) the far interval
has that sign or lies inside the low interval. The side
with negative separated interval for the chosen sign has
all its jumps \(v\le\eta\) at those levels; the moving
side has all \(u\ge-\eta\). Use the former side's standard
point as the stop. Its depth remains \(D(a,v)=-v+s_j\mathbf1_{a=0}
+\eta\mathbf1_{a\le l}\), and it has orders at least
\[
 O_l(a,v)=\begin{cases}
 \max(0,-v-\eta),&a<l,\\
 -v,&a\ge l .
 \end{cases}
\]
At level \(0\), if \(v\ge-\eta\), then \(v\le\eta\) and
\(D(0,v)=-v+s_j+\eta\ge0\), so the improved order is
nonnegative. If \(v<-\eta\) and \(s_j=0\), the initial
order is \(-v\ge1\). If \(v<-\eta\) and \(s_j>0\), then
\(v\le-\beta_j\) and the shifted cut has \(s_j<\beta_j\).
Thus \(-v-s_j\ge1\) also in this case. In both cases
the improved order is \(D(0,v)\).

At the other levels the asserted orders follow from the depth increase.
Let \(x_0\) be a standard point at the preceding cut with constant
parameter, and take \(g\in G_{K'}\). By
Section~\ref{sec:period-basis-integral-operators},
\[
 w_g=(gx_0-x_0)/(z\tau)\in M_{T,\mathrm{loc}}\cap\ker\pi.
\]
Here \(\pi\) is the projection to the parameter. We have
\[
 (g-1)\varphi^l x_0=\varphi^l(z\tau)\,\varphi^l w_g .
\]
At level \(a<l\), put \(n=l-a>0\). Choose a finite product
\(D=\prod_b E_b^{m_b}\) such that \(Dw_g\in M_T\). For every
chart-equivariant form \(\lambda\),
\[
 \lambda(\varphi^n w_g)=
 \frac{\lambda(\varphi^n(Dw_g))}{\prod_b E_{b+n}^{m_b}}
 \in B_{\mathrm dR}^+ .
\]
The positive Frobenius test makes the numerator integral, and
\(\theta(E_{b+n})\ne0\) since \(b+n>0\).
Further, \(\varphi^n\tau=[\underline\epsilon]^{p^{n-1}}-1\)
has Fontaine value zero. Thus the difference has order at least one.
At \(a\ge l\), \(\theta\varphi^{-(a-l)}(z\tau)\ne0\);
the factor has order zero and the preceding depth at \(a-l\)
is preserved.
The order of the difference of two points over the same
constant is bounded below by the minimum of their two
orders. For the move from the preceding cut to the present
point on the moving side, that minimum is at least
\(\min(0,-u-\eta)\) at level \(0\),
\(\min(0,-u)\) at \(0<a<l\), \(-u-\eta\) at \(a=l\),
and \(-u\) above. At \(a=0\) we have \(u\ge-\eta\), so adding \(D(0,v)\)
gives at least \(s_j-u-v\). For \(0<a<l\), use
\(\min(0,-u)\ge-u-\eta\). At \(a=l\) the extra \(\eta\)
in \(D(l,v)\) cancels the possible loss \(\eta\) in the order.
The resulting sum is at least \(-u-v\) for every \(a>0\). Thus
the annihilator of the Tate-dual side at the current cut
satisfies the order criterion below.

For a long gap, the preceding cut has shift \(0\).
A pivot \(j_*\) occurs \(L=J+1\) steps after it, with
the rest of the indices at distance \(<l\) of opposite
sign or bounded by \(\eta\). Let \(g_0\) bound the positive
jumps of the resting side at the \(L\) indices up to the
pivot and let \(E^{s_0}\) kill the pivot's cokernel.
The gap imposes
\(N=s_0+g_0\sum_{0\le a<L}p^a<\beta_{j_*}\).

On the unramified chart take \(A=2b_0\), with
\(N<A\le p^L D_0\). Start with the point \(y\) at \(j_*\)
over \(c_{j_*}\) obtained from \(E^{s_0}\). Multiply it by
\[
 f(E)=E^{s_0}\prod_{0\le a<L}E_a^{g_0},\qquad
 E_a=(E-p)^{p^a}+p .
\]
This is monic of degree \(N\) in \(E\). Since \(E_a\) has order one at its own level and zero at the
others, the orders of \(f(E)y\) are at least \(g_0-u\) at
\(0\le a<L\) and \(-u\) at \(a\ge L\). The former are
nonnegative by the upper bound at these \(L\) indices.
Write \(T\) for the lattice of the extension on the resting
side. Since \(y\in N_{T,\rm loc}\) and \(f(E)\in\mathfrak S\),
every equivariant form has value in \(\mathcal K_0\) on
\(f(E)y\). Its standard orders place it in \(\varphi^L M_T\) by
Section~\ref{sec:local-standard-orders}. The value-field test of
Section~\ref{sec:period-basis-order-retraction} therefore puts
\(f(E)y\) in \(\sum_i\mathfrak S\varphi^L(n_i)\).
Since \(A\le p^L D_0\) and \(3\le p^L\) on this unramified chart,
Section~\ref{sec:v2-crystalline-frame-radius} expresses
\(p^{4D_0}f(E)y\) as an \(\mathfrak S\)-linear combination
of vectors of this module whose form values are fixed modulo
\(\mathrm{Fil}^A\). Project this equality to the summand at
\(j_*\). The projector fixes its left side, preserves the module
and preserves the filtration condition by pullback of forms.
Since \(E=z+p\), expand the \(\mathfrak S\)-coefficients in
powers of \(E\) over \(\Lambda=W(\mathbf F_{p^f})\).
Indeed evaluation at \(z=-p\) maps \(\mathfrak S\) onto
\(\Lambda\) with kernel \((E)\). Repeated division supplies
the first \(N+1\) coefficients and a remainder in
\(E^{N+1}\mathfrak S\). The subgroup fixes those coefficients, giving
\[
 p^{4D_0}f(E)y=\sum_{n=0}^{N}E^n\xi_n+E^{N+1}y_1
\]
Each \(\xi_n\) and \(y_1\) belongs to the level-\(L\) standard
module, since its coefficients after division by \(E\) still lie in
\(\mathfrak S\). At level zero the value of every tested form
\(\lambda\) on each \(\xi_n\) is fixed modulo \(\mathrm{Fil}^A\). By Section~\ref{sec:local-tame-depth}
this value is congruent modulo \(E^A B_{\mathrm dR}^+\) to a
constant in \(K'\).

Fix a form on a line of jump \(u<-N\) at the pivot. The value
of the left side has order at least \(g_0-u\ge N+1\).
After removing the first \(k\) terms and dividing by \(E^k\),
the residual value has order at least \(N+1-k\). For \(k\le N\),
its reduction modulo \(E\) is the constant of \(\xi_k\), because
the later coefficients and the remainder have order at least zero
at the pivot. The constant lies in
\(K'\cap\mathrm{Fil}^1B_{\mathrm dR}=0\), so \(\xi_k\) has order
at least \(A\) on this form. Since \(A\ge N+1\), subtraction
of its value and division by \(E\) give the induction at \(k+1\).
On the trivial quotient write \(f_n=[X^n]f(X)\) and
\(w_n=\pi(\xi_n)-p^{4D_0}f_nc_{j_*}\). Since \(\pi(y)=c_{j_*}\),
\[
 \sum_{n=0}^N E^n w_n=-E^{N+1}\pi(y_1).
\]
The right side has order at least \(N+1\) at level zero, and each
coordinate of \(w_n\) is congruent to an element of \(K'\) modulo
\(\mathrm{Fil}^A\). The same induction gives order at least \(A\) for
each \(w_n\). Since \(f_N=1\), this yields
\(\pi(\xi_N)=p^{4D_0}c_{j_*}+w\). In a
\(\mathbf Z_p\)-basis of the trivial parameter lattice \(C\), all
coordinates of \(w\) belong to \(S\): the coefficient \(\xi_N\)
is an \(S\)-linear combination of \(\varphi^L n_i\), and
\(\pi(\varphi^L n_i)=\varphi^L\pi(n_i)\in S\otimes C\) by
\(\pi(N_T)\subset N_C=S\otimes C\) and \(\varphi^L(S)\subset S\).
The quotient calculation gives each coordinate de Rham order at
least \(A\). On this unramified chart \(\ker(\theta|_S)=(z+p)S=ES\).
Divide successively by \(E\) in \(S\), decreasing the de Rham
order by one each time. Hence \(w/E^A\in S\otimes C\).
It is supported at \(j_*\) and has the standard orders at level
\(L\), since the only jump of \(C\) is zero.

Over \(A_{\inf}\), \(\varphi^L\) carries the preceding
summand and its level-zero standard module to the pivot summand
and its level-\(L\) standard module. Projection to the trivial
quotient commutes with \(\varphi^L\), which fixes \(p^\kappa\).
The trivial quotient has only jump zero, so its standard module at
level \(L\) is \(e_{j_*}(A_{\inf}\otimes C)\). The projected parameters
\(e_{j_*}(1\otimes c)\), \(c\in C\), generate this module over
\(A_{\inf}\). Thus the transported multiplier \(p^\kappa\) lifts
every \(p^\kappa\)-multiple of this module, in particular
\(p^\kappa w/E^A\). Apply the coordinate retraction of
Section~\ref{sec:v2-crystalline-lift-ideals} to a lift. It
preserves the standard orders and fixes \(w/E^A\) in the
localized \(\mathfrak S\)-module of the trivial quotient.
This gives a point \(x_1\) of the same orders with
\(\pi(x_1)=p^\kappa w/E^A\).
Put \(q_*=p^\kappa\xi_N-E^A x_1\). This point lies over the
exact constant \(p^{\kappa+4D_0}c_{j_*}\), has depth \(A\) at
the pivot, and has order at least \(A\) there for \(u<-N\).
Both terms lie in \(N_{T,\mathrm{loc}}\cap\varphi^L M_T\).
Put \(t=l-L=l-J-1\) and take \(x=\varphi^t q_*\in\varphi^l M_T\)
at the current cut. Its pivot is at level \(t\).

For \(a<t\) write \(h=t-a>0\). For each \(g\) fixing the chart field,
Section~\ref{sec:period-basis-integral-operators} gives
\(y_g=(gq_*-q_*)/(z\tau)\in M_{T,\mathrm{loc}}\).
Since \(g\) commutes with Frobenius, every tested form \(\lambda\) satisfies
\[
 \lambda\!\left(\varphi^{-a}(gx-x)\right)
   =\varphi^h(z\tau)\,\lambda(\varphi^h y_g).
\]
Clear a product \(D\) of shifted divisors \(E_j\) so that
\(Dy_g\in M_T\). The positive Frobenius test gives
\(\lambda(\varphi^h(Dy_g))\in B_{\mathrm dR}^+\), and each
factor of \(\varphi^h D\) is \(E_{j+h}\), a de Rham unit by
Section~\ref{sec:local-standard-orders}. Thus
\(\lambda(\varphi^h y_g)\in B_{\mathrm dR}^+\).
The scalar \(\varphi^h z\) has order zero. Moreover
\[
 \varphi^h\tau=(1+\mu)^{p^{h-1}}-1
    \equiv p^{h-1}\mu\pmod{\mu^2 B_{\mathrm dR}^+}.
\]
Here \(\mu=\xi\tau\) has order one, and \(p^{h-1}\) is a
de Rham unit. Hence \(\varphi^h(z\tau)\) has order one and
the depth of \(x\) at every \(a<t\) is at least one.

\subsubsection{Raising depth before termination}
\label{sec:lifting-depth-ranges}

\sectionopening{In this subsection, we apply the depth step to the remaining indices before checking the two schemes in Section~\ref{sec:v2-crystalline-stop-hypotheses}.}

At the levels \(t<a<l\) raise depth from \(0\) to \(1\)
by the step \(n=0\). The chart is unramified, so its cost is
\(C(a,0)=3a\). Since \(l\le2J+1\) and \(l-t-1=J\),
\[
 \sum_{a=t+1}^{l-1}3a=\frac{3J(l+t)}2
 \le 3J(2J+1)<\kappa .
\]
The strict inequality follows from
\(X_0\ge (J+1)W_0\), \(W_0\ge4\eta\ge4\) and the definition of
\(\kappa\). The multiplier exponent is therefore at most
\(\kappa+4D_0+\sum_{a=t+1}^{l-1}3a<3\kappa\). Scale the
point to the common exponent \(25\kappa\). The resting side
has \(v\ge-\eta\) off the pivot below \(l\). At the pivot
either \(v\ge-\eta\) or \(v<-N\). In the latter case the order
is at least \(A\), and \(-v\le\beta_{j_*}+V\le2b_0=A\).
Thus the point has order bounds \(O_l(a,v)\) and
depth bounds \(1\) below \(l\) off the pivot, \(A\) at the pivot,
and \(-v\) at \(a\ge l\).

At a long gap use the preceding cut's Tate-dual
annihilator. Off the pivot the moving side has \(u\le\eta\)
below \(l\), so \(C_j(a,u)\ge0\) there. Its old point
carried from the preceding cut has order \(\ge0\) below
\(l\). At \(a=t\) the order of their difference is at
least \(\min(0,\eta-u)\), and \(u<\eta+A\).
At \(a\ge l\) the order is at least \(-u\); the shift of
the old cut is zero. Adding the depth of the stop yields
at least \(1\) below \(l\) and \(-u-v\) thereafter, as
required for the move.

Here is the order criterion used in both cases. The dual jumps
\(-u-v-1\) can lie in both far intervals. At the cut \(l_1\)
steps before the current index use the Tate-dual projection
annihilator \(p^\kappa E^s\). Lift the multiple \(p^\kappa E^s\)
of the constant, divide by \(E^s\), and carry the resulting
point over \(p^\kappa\) by Frobenius. Denote its order bound
on a dual jump \(w\) by \(P_a(w)\). This is zero for \(a<l_1\)
and \(-w-s\mathbf1_{a=l_1}\) for \(a\ge l_1\), hence is
nonincreasing in \(w\). The Tate order \(1\) follows if the mover's
order \(o_a(u)\) and the resting point's depth \(d_a(v)\) satisfy
\[
 o_a(u)+d_a(v)\ge 1-P_a(-u-v-1)
  =\begin{cases}
    1,&a<l_1,\\
    s-u-v,&a=l_1,\\
    -u-v,&a>l_1 .
   \end{cases}
\]
A tensor covector of type at most \(r\) is a sum of products of
covectors with \(u+v\le r\), by the tensor filtration. Then
\(P_a(-r-1)\ge P_a(-u-v-1)\), so the tensor
has order at least \(1-P_a(-r-1)\) on weight \(r\). In the final
pairing a contributing dual weight \(w\) satisfies \(w+r\le-1\),
the Tate jump. Hence \(P_a(w)\ge P_a(-r-1)\), giving total order
at least one. The bounds hold for sums of terms since each
filtration step is a subspace. In the Galois difference of the resting point
divide by \(z\tau\) before taking this contraction.

This
keeps the orders by the chosen basis. The contraction has
\(\tau\) times its value in the integral Tate lattice.
Returning the removed \(z\tau\) places the value in
\(W(\mathfrak m^\flat)\), since \(z\) lies in that ideal.
Reduction modulo the pairing exponent kills the component
of the Tate-dual extension over zero, so the remaining
value is precisely \(p^\kappa\) times the desired reduced
pairing. For a move on the second side the test pairs the first difference
with the Galois translate of the move. Equivariance preserves its
orders, and conjugating the group element gives the same depth
bound. This proves the required pairing condition for either side.

The moves to and from a stop also satisfy this criterion.
The move from the old cut to the stop uses \(l_1=l\)
and that cut's shift \(s_{\rm old}\). Below \(l\) its
order is \(\ge0\) and the other old point's depth is
\(\ge1\). At \(l\) the lower bounds are \(-u-\eta\) and
\(-v+s_{\rm old}+\eta\), and above \(l\) their sum is
at least \(-u-v\). The move from the stop to the new
point uses \(l_1=0\) and \(s=s_j\).

At \(0\) its order
is at least \(-u-\eta\) and the depth is at least
\(-v+s_j+\eta\); above \(0\) their sum is at least
\(-u-v\). Finally the pair of cocycle differences of
two standard points transported from a cut satisfies
the criterion with that cut's annihilator. Their depths
are at least \(1\) below the transport length, at the
transport length the order of each increases by the shift of the cut, and
above it each is at least its standard order.
The cochain calculation now applies to either cut scheme.

\subsubsection{Termination in the two filtration schemes}
\label{sec:v2-crystalline-stop-hypotheses}

\sectionopening{We check the hypotheses of Section~\ref{sec:v2-crystalline-stops} for
the two distributions of Hodge jumps in
Section~\ref{sec:relative-crystalline-colength}. In the first case the
large jumps occur sparsely in the Frobenius cycle. In the second case
we use the two possible magnitudes of the centers of their intervals.}

Work in the Frobenius cycle of one original chart. Let \(B_0,G_0\)
be the larger and smaller starting points of the separated intervals,
of width \(V\). Index a segment by
\(0,1,\ldots,l\) in the Frobenius direction, with the preceding
cut at \(0\) and the current cut at \(l\). An index at either
end is counted in this direction even if the segment wraps
around the cycle.

\indent\textit{Sparse cuts.}
Here \(\eta_0\ge1\) and \(G_0+V\le\eta_0\). Extend the
Frobenius cycle periodically to \(\mathbf Z\), and let \(T\) be
the positions with the larger start \(B_0\). These positions are
thick and their complement is thin. The sparsity hypothesis is
\[
 \{k,k+a,k+a+b\}\not\subseteq T
       \qquad(k\in\mathbf Z,\quad 1\le a,b\le J).
\]
The three positions may repeat residues on the cycle. The signed
ratios in Section~\ref{sec:coefficient-hodge-tate-intervals} exclude
these triples as well.
Assume either that the chart is unramified,
\(p\ge7\) and \(10p^J\eta_0\le B_0\), or that the signs are constant
and \(8e\le p^{J+1}\). In both cases require
\[
 2e(B_0+5\eta_0)\le p^{J+1}D_0,\qquad
 20(J+1+D_E)\eta_0+\lfloor\log_p(2(B_0+5\eta_0))\rfloor\le X_0.
\]
Take \(s_j=0\) at every sparse cut. A cut is an index whose
next \(J\) indices are thin. If a thick
index \(w\) fails the cut condition, some thick index \(w'\)
occurs \(a\) steps after it with \(1\le a\le J\). The \(J\)
indices immediately before \(w\) must all be thin, since any
thick one would give a forbidden triple ending at \(w'\).
Therefore \(w-(J+1)\) is a cut.

For any index \(i\), either \(i-(J+1)\) is a cut or there is a
thick \(w\) in the interval from \(i-J\) to \(i-1\). In the
second case one of \(w\) and \(w-(J+1)\) is a cut. Thus the
nearest strictly preceding cut has \(l\le2J+1\). Within its
segment, every thick intermediate index occurs exactly \(J+1\)
steps after that cut. It occurs at least \(J+1\) steps later by
the definition of cut. If it occurred more than \(J+1\) steps
later, the preceding argument would put another cut strictly
inside the segment.

Suppose all thick indices in \(1,\ldots,l\) have the same sign.
The plain stop applies since every thin separated interval is contained in
the low interval. This includes the constant-sign alternative.
Otherwise there are exactly two thick indices in \(1,\ldots,l\).
Take the earlier one, at step \(L=J+1\), as the pivot. The
only other thick index is the current cut, whose sign is
opposite to the pivot sign.

The opposite signs force the unramified alternative, so
\(p\ge7\) and \(10p^J\eta_0\le B_0\). Both cut shifts are zero. The low radius satisfies
\(\eta_0\ge1\). For the
resting side take \(g_0=\eta_0\), including at the pivot, where
its separated interval is negative. The integer
\[
 s_0=\left\lfloor\frac{B_0+10\eta_0}{p-3}\right\rfloor+1
\]
is an \(E\)-annihilator exponent for the resting side at the
pivot. Indeed, at every index its jumps are at most
\(\max(\eta_0,B_0+V)\le B_0+\eta_0\le(p-3)s_0\).
This is the global upper bound required by the unramified
lift-ideal estimate. The long-stop
degree satisfies
\[
 N=s_0+\eta_0\sum_{a=0}^Jp^a
     \le (B_0+10\eta_0)/4+1+2p^J\eta_0 < B_0 .
\]
For the last inequality, \(10p^J\eta_0\le B_0\),
\(p^J\ge7\), and \(B_0\ge70\). Set \(b_0=B_0+5\eta_0\).
The radius parameter in the long stop is
\(A=2b_0\le p^L D_0\), by the sparse depth hypothesis with
\(e=1\). In particular \(N<A\). This proves the geometric
and numerical hypotheses of the long stop at every sparse cut.

For the sparse lift ideals and the step bound take the auxiliary
width \(W_0=5\eta_0\), and take \(G_1=n_0=10\eta_0\).
The \(J\) successors of a cut satisfy \(\beta+W_0\le G_1\).
Also \(V+4\eta_0\le W_0\), \(W_0\le pn_0\),
\(G_1\le(p-1)n_0\), and \(\beta+W_0\le2b_0\).
The depth hypothesis gives \(2eb_0\le p^{J+1}D_0\).
The ramification inequality \(8e\le p^{J+1}\)
is assumed in the constant-sign case and follows
from \(e=1\), \(p\ge7\), \(J\ge1\) in the other.
Finally
\[
 n_0+\lfloor\log_p(2b_0)\rfloor
       +(J+1+D_E)W_0+D_EG_1
 \le 20(J+1+D_E)\eta_0+\lfloor\log_p(2b_0)\rfloor\le X_0 .
\]
Thus the lift and step estimates apply at every sparse cut,
where all \(J\) succeeding indices are thin.

\indent\textit{Two far starts.}
Assume an unramified chart with \(p\ge7\), \(\eta\ge1\),
\(0<G_0\le B_0<Kp^2\), \(2G_0<Kp\), and
\(\widetilde W=V+4\eta\le K\). Put
\(q=\lfloor B_0/p\rfloor\) and \(s_0=q+12K\).
At every index, independently of its shift and sign,
\[
 (p-3)s_0-(\beta_j+\widetilde W)
 \ge -3q+12K(p-3)-(p-1)-\widetilde W>0.
\]
Indeed \(B_0\le pq+p-1\), \(q<Kp\),
\(\widetilde W\le K\), \(K\ge1\) and \(p\ge7\).
This bounds the jumps on all three sides wherever the
unramified lift-ideal estimate is used. Let the shift at an index be
\(\max(0,q-7K)\) if its successor is thick, and zero
otherwise. An index is a cut when its shift is smaller
than its far start. Every thick index is a cut, since
its shift is at most \(q<B_0\). Every index with thin next component is a cut, since its shift is zero. Therefore
if the immediate predecessor is no cut, it is thin and
the second predecessor is a cut. Thus \(l\in\{1,2\}\).

For \(l=1\), or \(l=2\) with equal signs on the indices at
steps \(1,2\), use the plain stop. When \(l=1\), the preceding
cut may be thick with the opposite sign. Only the current index
enters the sign condition at \(a<l\). At \(a=l\) the two jump
indices are both \(F^{-l}j\), as in
Section~\ref{sec:lifting-plain-cut}. Otherwise \(l=2\), the
index at step \(1\) is thin and the current index is
thick. The preceding cut has shift zero. The no-cut
condition at the thin index gives \(G_0+7K\le q\), since
\(G_0>0\). The pivot is the current index. Take
\(L=2\) (the long-stop parameter is \(J=1\)),
and \(g_0=G_0+\widetilde W\). The global jump bound
gives the annihilator \(E^{s_0}\) at the pivot.
The bound on \(g_0\)
uses the thin start at step \(1\) and the negative separated interval of the resting side at the pivot. Since
\(G_0+\widetilde W\le q-6K\),
\[
 \begin{aligned}
 N
   &=s_0+g_0(1+p)\\
   &\le (p+2)q-(6p-6)K
     \le B_0+2q-(6p-6)K <B_0 .
 \end{aligned}
\]
With \(b_0=B_0+\widetilde W\), one has \(N<2b_0\le4Kp^2\),
so the radius parameter \(A=2b_0\) meets the long-stop
bound with \(D_0=4K\).

For the step bound use \(n_0=K,\ G_1=G_0+\widetilde W\),
\(J=1,\ D_0=4K\), and \(X_0=3K+\lfloor\log_p(2b_0)\rfloor\).
We have \(G_1\le(p-1)K\) by \(2G_0<Kp\). Also
\(s+\widetilde W\le pK\) by \(q<pK\) and the two cases of
\(\max(0,q-7K)\), and \(\beta+s+\widetilde W\le2b_0\).
Because \(D_E=0\), the step-bound follows from
\(K+2\widetilde W\le3K\). For a cut with a thin next component
use the lift-ideal estimate for a thin next component with shift zero.
For a cut with thick successor (even if its shift happens to be zero),
\(s_0-s\le19K<\kappa\), so the unramified lift-ideal estimate
gives \(p^\kappa E^s\). All cuts and stops use the original
signs and coefficient intervals of the original components.

\smallskip
\indent\textit{The total exponent.}
For each index \(j\), let \(l(j)\) be the distance to the
nearest preceding cut strictly before it, and let
\(d(j)=0\) at a cut and \(d(j)=l(j)\) elsewhere. In a cycle of
length \(h_c\) with a single cut \(j\), one has \(l(j)=h_c\).
The preceding point at this cut is consequently transported by
\(\varphi^{h_c}\), even though its residue index returns to \(j\).
Transport
the two selected points at the last cut by \(\varphi^{d(j)}\).
The points chosen at the cuts have the improved orders \(C_j\)
and the lower orders needed at a plain cut.
Write \(Fj=j^p\). The identity \(l(j)=1+d(F^{-1}j)\)
gives, for either family \((x_i)\) of points at the cuts,
\[
 \varphi\!\left(\sum_j\varphi^{d(j)}x_{F^{-d(j)}j}\right)
   =\sum_j\varphi^{l(j)}x_{F^{-l(j)}j}.
\]
Reindex the left sum by the permutation \(F\). On either moving
side, for \(a=d(j)\) or \(a=l(j)\), the projection to its parameter
lattice sends the transported point to
\[
 \pi\bigl(\varphi^a x_{F^{-a}j}\bigr)
   =p^{25\kappa}\varphi^a\bigl(e_{F^{-a}j}(1\otimes c)\bigr)
   =p^{25\kappa}e_j(1\otimes c).
\]
Here \(\varphi\) acts on the scalar ring and fixes \(1\otimes c\).
Summing over \(j\) gives \(p^{25\kappa}(1\otimes c)\).
Off a cut the transported points agree since \(d(j)=l(j)\).
At a cut the two points have equal projections even when the
preceding cut has the same residue index. Write \((y,y')\) for the
transported pair and \((x,x')\) for the chosen pair there.
In every summand use the same order of changes,
\[
 (y,y')\longrightarrow(w,y')\longrightarrow(w,w')
 \longrightarrow(x,w')\longrightarrow(x,x').
\]
If the first side has the stopping point \(q\), take
\((w,w')=(q,x')\); the fourth change is zero. If the second
side has the stopping point \(q'\), take \((w,w')=(y,q')\);
the first change is zero. In each case the other three changes
satisfy the order inequalities just proved. Thus each position
in this sequence can be summed across coefficient summands.

For the cocycle-difference test, take the third
period from the cut supplying the endpoints and transport it
by \(\varphi^{d(j)}\) from \(j^{p^{-d(j)}}\). At a cut this
is the current third period. Off a cut it is the transported
preceding period. Write \(d=d(j)\) and let \(s\) be the shift of
that cut. Below level \(d\), each transported Galois difference
has order at least one by the factor \(\varphi^{d-a}(z\tau)\)
at level \(a<d\). At \(a=d\) the two differences have orders at
least \(-u+s+\eta\) and \(-v+s+\eta\), whose sum is at least
\(s-u-v\). At \(a>d\) their sum is at least \(-u-v\).
These are the bounds \(1-P_a(-u-v-1)\) for the third period
transported through these same \(d\) steps.
At the fixed pairing depth \(m\), the
third period is over \(p^\kappa e_j\), where \(e_j\)
is the projector to index \(j\). The reduction of its
extension form is \(p^\kappa\psi(e_jy,y')\), followed
by the same \(\mathbf Z_p\)-linear functional on \(O\)
for every test. On index \(j\) this is \(p^\kappa\psi(y,y')\).

The cross pairings of distinct indices vanish: balance
for the coefficient root of unity multiplies such a
pairing by the difference of two Teichm\"uller lifts,
which is a Witt unit. Thus the sum over all indices
satisfies the four move tests and the two-cocycle difference test of
Section~\ref{sec:tame-cup-frobenius}. Including the two moving sides
and the Tate-dual period gives the factor
\[
              (p^{25\kappa})^2p^\kappa=p^{51\kappa}.
\]

In the sparse regime \(\kappa=12((J+1)X_0+D_0)+1\),
so \(51\kappa\le 10^3((J+1)X_0+D_0)\).
In the two-start regime \(\kappa=12(2X_0+4K)+1\), where
\(X_0=3K+\lfloor\log_p(2b_0)\rfloor\). Consequently
\[
 51\kappa\le 10^4\bigl(K+\lfloor\log_p(2b_0)\rfloor\bigr).
\]
Here \(K\ge1\).

\begin{proposition}[Uniform cup annihilation in the two filtration schemes]
\label{prop:crystalline-filtration-cup-exponent}
Use the Kummer chart of \(E_v\) specified in
Section~\ref{sec:coefficient-hodge-tate-intervals}, with working
group \(H=G_{E_v}\). Suppose the allowed partner \(Y\) and the
fixed Tate-dual lattice \(L^\dagger\) are crystalline over
\(E_v\). Let \(z\in Z_f(Y)\), \(z'\in Z_f(L^\dagger)\),
with the pairing \(\psi\) and depth \(m\) of
Section~\ref{sec:v2-crystalline-stops}. For any \(\mathbf Z_p\)-linear
functional \(\lambda:O\to\mathbf Z_p\), put
\(\psi_m=((\lambda\otimes1)\circ\psi)\bmod p^m\). Use one of the two filtration schemes under the interval and cut
hypotheses of Section~\ref{sec:v2-crystalline-stop-hypotheses}.
Then
\[
 p^{51\kappa}[(\psi_m)_*(\bar z\cup\bar z')]=0
       \quad\text{in }H^2(E_v,\mu_{p^m}),\qquad
 51\kappa\le
 \begin{cases}
 10^3((J+1)X_0+D_0),&\text{sparse regime},\\
 10^4(K+\lfloor\log_p(2b_0)\rfloor),&\text{two-start regime}.
 \end{cases}
\]
\end{proposition}

\section{Local colength away from the original level}
\label{sec:crystalline-colength-appendix}

\sectionopening{In this appendix, we apply the lifts of
Appendix~\ref{sec:crystalline-lifting-appendix} to the local
coefficient lattices. We handle separated Hodge--Tate weights and the
cases of exact integrality, then combine the bounds at the primes
away from the original level in
Section~\ref{sec:off-level-coefficient-defect}.}

\subsection{Colength for separated Hodge--Tate weights}
\label{sec:v2-crystalline-deep-constants}

\sectionopening{In this subsection, we treat the unramified case in which the nonzero multiples of \(p\)
avoid the intervals of filtration orders used by the divided
operators. Using the \(E\)-power annihilator of
Section~\ref{sec:v2-crystalline-lift-ideals}, we construct one
submodule of crystalline cocycles on which all the requisite cup
products vanish. Its colength gives the local estimate.}

Use an unramified chart with residue degree \(f\) and \(p\ge5\).
At each residue-field index \(j\), the jumps of the two sides lie in
\([-\eta,\eta]\) or in their respective far intervals
\(\pm\varepsilon_j[\beta_j,\beta_j+W]\), with opposite signs as in
Section~\ref{sec:coefficient-hodge-tate-intervals}. Here
\(\eta\ge1\) and \(W\ge0\) are integers. Let \(A\ge0\) be an integer with
\[
 \beta_j+W+\eta\le(p-3)A,\qquad
 6(A+\eta+W)<\beta_j
\]
at every index \(j\). Require that no nonzero multiple of
\(p\) lie within \(6(A+\eta+W)\) of \(0,\pm\beta_j\).
For jumps \(u,v\) of the two sides at the same index,
one of them is at least \(-\eta\), since their separated
intervals have opposite signs. Both are at least
\(-(\beta_j+W)\), as \(\beta_j>\eta\). It follows that
\[
     -u-v-1\le \beta_j+W+\eta-1\le(p-3)A .
\]
The trivial quotient contributes only the jump \(0\).
On each original side the jumps are at most \(\beta_j+W\).
Thus the unramified lift-ideal estimate in
Section~\ref{sec:v2-crystalline-lift-ideals} applies to all three
extensions, with annihilators \(p^s E^{A-s}\) for \(0\le s\le A\).
Write \(R(t)\) for the orders \(R(t)(a,u)=-u+t\mathbf1_{a=0}\).

Let \(L\) be the coefficient lattice of the fixed side, so
\(C=Z_f^1(L)\) is its module of crystalline cocycles.
Use the single extension on \(L\oplus C\), with trivial
action on \(C\) and with
\[
             g(y,c)=(gy+c(g),c).
\]
The cocycle equation gives the action law. The module \(C\)
is finite free: integral \(H^1\) is finitely generated,
the coboundaries have finite rank, and \(C\) is torsion free
over \(O\). A crystalline primitive for each element of an
\(O\)-basis of \(C\) splits this extension over \(B_{\rm cris}\).
In the de Rham grading it adds copies of jump zero.
The nonzero-jump estimates in the extension thus use the
multiplicities of \(L\); the depth estimates over zero use
the kernel \(L\) itself.

Let \(\rho_j\) be the projection of \(e_jN_{L\oplus C,\mathrm{loc}}\)
to its trivial quotient at index \(j\), and put
\(\iota_j(c_0)=e_j(1\otimes c_0)\) for \(c_0\in C\).
Define \(\mathcal D_j(t,c)\subseteq\iota_j(C)\) by
\[
 \iota_j(c_0)\in\mathcal D_j(t,c)
 \quad\Longleftrightarrow\quad
 \exists x\in e_j N_{L\oplus C,\mathrm{loc}},\quad
 \rho_j(x)=\iota_j(c_0),\quad x\text{ has orders }R(t)
                    \text{ and depth }c .
\]
The retraction of Section~\ref{sec:period-basis-order-retraction}
fixes each such constant and preserves the orders of the form values.
It also preserves depth, since a value of depth \(d\) is congruent
modulo \(\mathrm{Fil}^d\) to a scalar in the chart field.
Thus \(\mathcal D_j(t,c)\) is the image under \(\rho_j\) of these
points with the specified orders and depth. All preimages under
\(\iota_j\) are taken in the same \(C\), before choosing the
partner, its cocycle, the finite level or the pairing.
For each element of the intersection the lifts at every index
therefore project to the components of that element itself,
including when it is divisible by \(p\) in \(C\).
For a partner or a Tate-dual period in an individual cup test
one may use the extension by \(O\) defined by that single
cocycle. Let \(n_j,r_j\) bound \(f\) times the low-jump and
all-jump multiplicities of \(L\) at index \(j\).

We first prove the integral linear algebra for one order
step. Suppose a system of orders \(c\) is being raised
at \((a,u)\), and \(c(a,v)>c(a,u)\) for \(v<u\) at the same
index and level. Take \(f m_u+1\) vectors of orders \(c\)
in the localized chart module, where \(m_u\) is the multiplicity
of the jump at the index \(\varphi^{-a}j\) read at level \(a\) from
the current index \(j\). Consider the leading coefficient covectors on a basis
of the \(m_u\)-dimensional graded space of de Rham forms
at \(u\). Lower lines make no contribution by the strict
inequality on orders. Over \(C_p\) choose a maximal-rank
minor with largest absolute determinant among minors of
that rank. Cramer's rule then writes every covector in terms of
its \(d\le m_u\) selected covectors with coefficients in
\(\mathcal O_{C_p}\). For each coefficient \(\alpha\) choose
\(\widetilde\alpha\in A_{\inf}\) with
\(\theta(\widetilde\alpha)=\alpha\). Use
\(\varphi^a\widetilde\alpha\) in the relation at \(E_a\), since
\((\theta\varphi^{-a})(\varphi^a\widetilde\alpha)=\alpha\).
The difference is
a vector whose order is raised at \((a,u)\).

Write \(M(c)\) for the submodule of chart vectors of
orders \(c\). Let \(x_j\) denote the original vectors
and \(x_{s(k)}\) the \(d\) selected ones. The retraction
fixes these vectors. It preserves the raised order of their
difference and sends a product \(b x_{s(k)}\), \(b\in A_{\inf}\),
to \(F(b)x_{s(k)}\), with \(F(b)\in\mathfrak S\). Hence
\[
 x_j-\sum_{k=1}^{d}b_{jk}x_{s(k)}
      \in M(c+\mathbf1_{(a,u)}),\qquad b_{jk}\in\mathfrak S .
\]
Modulo \(E\), the \(d\)-tuples of coefficients are in
\(W(\mathbf F_{p^f})^d\), of \(\mathbf Z_p\)-rank
\(fd\le f m_u\). Divide a nonzero integral relation by
a coefficient of smallest valuation. The resulting
coefficients \(a_j\) include \(1\) and satisfy
\[
 \sum_j a_j b_{jk}\in E\mathfrak S,\qquad
 \sum_j a_jx_j\in
    M(c+\mathbf1_{a=0})+M(c+\mathbf1_{(a,u)}).
\]
Here \(E\) raises exactly the level-zero orders. This argument is also
valid for the actual extension lattice. Its multiplicity
at any nonzero \(u\) is that of the sublattice \(L\),
since the parameter quotient is trivial.

\subsubsection{Uniform comparison of the depth steps}
\label{sec:crystalline-depth-comparison}

\sectionopening{In this subsection, we compare successive depths using the lift construction of Section~\ref{sec:v2-crystalline-stops}.}

The depth version compares \(\mathcal D_j(0,c)\) with the
constants for the stronger depth \(c^+\) at \((a,u)\).
Assume all current depth orders are \(\ge-u\), so the order
module of vectors over zero with those orders is contained
in the no-pole order module. Denote this kernel order module
by \(M\). Let \(X(c)\) consist of the pole-free points of depth
\(c\) whose projection under \(\rho_j\) belongs to \(\iota_j(C)\).
Then \(\rho_j(X(c))=\mathcal D_j(0,c)\) by definition.
Every vector of \(M\) lies in \(X(c)\): its order bounds imply
the corresponding bounds on its Galois differences. Its projection is
zero. Consequently \(\rho_j\) induces a surjection
\[
 Q:=X(c)/(X(c^+)+M)\ \longrightarrow \mathcal D_j(0,c)/\mathcal D_j(0,c^+).
\]
Put \(k=C(a,c(a,u))\). The one-step operator sends \(x\in X(c)\)
to a point of \(X(c^+)\) congruent to \(p^kx\) modulo \(M\);
in particular \(p^k Q=0\).

Choose \(f m_u+1\) classes in \(Q[p]\) and representatives \(x_j\).
By the definition of \(Q\), write \(p x_j=x'_j+m_j\) with
\(x'_j\in X(c^+)\) and \(m_j\in M\).
The inclusion \(L\subset L\oplus C\) is saturated since \(C\)
is free. The localized kernel is its intersection with the ambient
chart lattice by Section~\ref{sec:period-basis-tensor-orders}.
Every form of type at most \(u\) on this kernel extends to one of
the same type bound on the extension by Section~\ref{sec:local-tame-moves}.
Thus each \(m_j\) lies in the kernel chart lattice with all the
specified orders. Applying the preceding relation there, with the
multiplicity \(m_u\) of \(L\), gives a primitive combination \(x\)
of the \(x_j\) for which \(p x-w\) has the raised depth
and \(w\) is over zero with all level-zero orders raised
by \(1\).

Suppose further that no multiple of \(p\) lies strictly
between \(c(0,v)\) and \(c(0,v)+p^a\) at any line \(v\).
All the depth operations that raise \(w\) from the first
level-zero order increment to \(p^a\) increments are
free. Use their common composite \(B\) on \(x,w,px-w\).
By Section~\ref{sec:local-offlevel-depth-operators}, its factors
are \(\mathbf Z_p\)-linear, each has scalar multiplier \(1\), and
its correction on the original system lies in \(M\). Thus
\(B(px-w)=pBx-Bw\) and \(Bx-x\in M\). Since \(w\) lies over
zero and has the increased orders, it has at least those
depths as well.

The free steps and improvement of orders
give \(B w=E^{p^a}w_1\) with \(w_1\in M\).
In checking this improvement, either
\(c(0,v)+1>0\) or \(c(0,v)+p^a\le0\), because \(0\) is absent
from the open interval. In the former case the positive initial
order forces the improved order to equal or exceed the target.
In the latter case the target is nonpositive, so congruence to a
constant to that depth already gives the target order. This includes
the endpoint \(c(0,v)+p^a=0\).
Use \(E^{p^a}=E_a+pF\) with \(F\in\mathfrak S\). Then
\[
 p(Bx-Fw_1)=B(px-w)+E_a w_1.
\]
The first term has the raised depth because the free
operator preserves depths already acquired. The second
term has that depth because \(w_1\) has orders \(c\)
and the multiplication raises every order at \(a\) by
\(1\).

At each period prime \(E_b\), the scalar \(p\)
has order zero, so \(Bx-Fw_1\) has the raised depth as well.
It differs from \(x\) by an element of \(M\).
It remains pole-free, since \(M\) is contained in the no-pole module.
The resulting primitive relation is a nontrivial relation in \(Q[p]\),
so \(\dim_{\mathbf F_p}Q[p]\le f m_u\). Multiplication by
\(p^{i-1}\) injects \(Q[p^i]/Q[p^{i-1}]\) into \(Q[p]\) for
\(1\le i\le k\). Hence
\[
 \operatorname{length}_{\mathbf Z_p}
 \bigl(\mathcal D_j(0,c)/\mathcal D_j(0,c^+)\bigr)
 \le \operatorname{length}_{\mathbf Z_p}Q\le k f m_u .
\]

One may iterate this estimate toward a target depth that
is antitone in the jumps at the level being raised.
At each step take the least jump \(u\) whose current
depth is short of the target \(d\). For every smaller jump \(v<u\),
\[
 c(a,u)<d(a,u)\le d(a,v)\le c(a,v).
\]
Thus the primitive relation applies even when two target depths
are equal. Write \(i=\varphi^{-a}j\) and let \(m_{i,u}\) be the
multiplicity at \(u\) there. With the other levels fixed and the
open-interval condition valid at each step, the total loss in length is at most
\[
 f\sum_{u\in\mathcal J_i}m_{i,u}
   \sum_{c(a,u)\le n<d(a,u)}\bigl(v_p(n)+3a\bigr),
 \qquad v_p(0)=0 .
\]
Each increase removes one pair \((u,n)\) from this sum. Equal step
values on different jumps occur separately with their respective
multiplicities. For level \(0\) the open intervals of length \(p^0=1\) are
always empty of integers. For level \(1\) we retain
the same depth at level \(0\) on every line, a multiple
of \(p\), so the intervals of length \(p\) contain no
multiple in their interiors.

Fix an index \(j\). For each integer \(0\le s\le A\), the
unramified lift-ideal bound with \(t=A-s\) and \(\kappa=s+t=A\)
gives the annihilator \(p^sE^{A-s}\) on all three extensions.
On either original side let \(d=e_j(1\otimes c)\) for an element
\(c\) of its trivial parameter lattice. Choose a standard-order
lift of \(p^sE^{A-s}d\). Its coordinate retraction preserves these
orders and fixes the projection, since \(d\) is a chart vector
and \(E\in S\). Dividing by \(E^{A-s}\) gives a localized chart
point over \(p^sd\) with orders \(R(-(A-s))\). Indeed \(E\) has
order one at level zero and order zero at every positive level.

Raise its depth to \(R(2A+\eta)\). Every level-zero step is within
\(6(A+\eta+W)\) of \(0\) or \(\pm\beta_j\) and costs no factor
of \(p\), so the projection stays \(p^sd\). At \(s=0\) the gap
\([-A,-\eta)\) gives improved orders at least \(R(-\eta)\) over
\(d\). At \(s=A\) they are at least \(R(0)\) over \(p^Ad\).
The depth \(R(A)\) for the tensor-dual extension is verified for
its own jumps in Section~\ref{sec:crystalline-pole-removal}.

\subsubsection{Pole removal and the final integral lift}
\label{sec:crystalline-pole-removal}

\sectionopening{In this subsection, we apply the comparison of Section~\ref{sec:crystalline-depth-comparison} at each pole order.}

Remove the pole one step at a time. For \(t=-\eta,\ldots,-1\),
compare \(\mathcal D_j(t,R(2A+\eta))\) with
\(\mathcal D_j(t+1,R(2A+\eta))\).
Their quotient is killed by \(p^A\), since \(p^A\iota_j(C)\)
already has lifts of order \(R(0)\) and the prescribed depth.
For a jump \(v<t\), the old order \(t-v>0\) and the retained
depth \(2A+\eta-v>0\) force the approximating constant to be zero.
Thus the value has order at least \(2A+\eta-v\ge t+1-v\).
For \(v>t\), its order is at least \(\min(0,2A+\eta-v)\ge t+1-v\),
since \(v\ge t+1\) and \(2A+\eta\ge t+1\). If \(m_t=0\), this
covers every basis form, so the two modules of constants coincide.
If \(m_t>0\), take \(f m_t+1\) representatives and apply the
single-order relation at \(t\), using
\(R(t)(0,v)=t-v>0=R(t)(0,t)\) for \(v<t\).
Both terms in the resulting primitive combination have order at
least one at \((0,t)\). The combination retains its depth by
\(\mathbf Z_p\)-linearity, so the bounds at \(v\ne t\) apply to it
as well. It therefore lies in the smaller module.
If \(Q_t\) denotes this quotient of constants, any
\(f m_t+1\) classes in \(Q_t/pQ_t\) therefore satisfy a nontrivial
relation. Thus \(\dim_{\mathbf F_p} Q_t/pQ_t\le f m_t\), and
each layer \(p^iQ_t/p^{i+1}Q_t\) is a quotient of \(Q_t/pQ_t\).
Since \(p^A Q_t=0\), its length is at most \(A f m_t\).
The trivial parameter quotient adds no jump at \(t<0\).
All these jumps are low, so pole removal has total
\(\mathbf Z_p\)-colength at most
\[
       A f\sum_{-\eta\le t<0} m_t\le A n_j .
\]

At an index where the separated interval of the fixed side is
positive, retain the depth \(R(2A+\eta)\).
At a negative-interval index seek the depths
\[
\begin{aligned}
 d_j(0,u)&=\beta_j+W+A+\eta,\\
 d_j(1,u)&=W+1-u+2A\mathbf1_{|u|\le\eta},\\
 d_j(a,u)&=-u\qquad(a\ge2).
\end{aligned}
\]
At level zero first increase the depth to the least
multiple \(T\) of \(p\) at least \(\beta_j+W+A+\eta\).
First replace each larger old requirement by \(T\), enlarging the
module of admissible constants. The old and new requirements are
both at least \(-u\), since \(-u\le\beta_j+W\) here. Where the
old requirement is smaller, the exact step range is
\(2A+\eta-u\le y<T\); otherwise there is no step.
Intersecting the resulting lattice with the original one costs
at most the length bound for these increases.
For \(|u|\le\eta\) the steps lie in \([0,T)\), with
\(C(0,0)=0\). The choice of \(T\) gives
\(T\le\beta_j+W+A+\eta+p-1\le(p-2)A+p-1\).
By Legendre's formula,
\[
 \sum_{2A+\eta-u\le y<T} C(0,y)
 \le v_p((T-1)!)\le (T-1)/(p-1)\le A+1 .
\]
For a far negative jump, a multiple
of \(p\) below \(T\) is at most \(T-p\). If such a multiple \(y\) occurs among the steps, then
\[
 \beta_j+2A+\eta\le -u+2A+\eta\le y\le T-p
 <\beta_j+W+A+\eta .
\]
The exclusion of nonzero multiples near \(\beta_j\) makes
these steps free.
By the depth version of the order estimate the contribution in
length is at most \((A+1)n_j\).

Let \(i=\varphi^{-1}j\) and keep \(T\) at level zero.
As \(p\mid T\), the open interval \((T,T+p)\) contains no multiple
of \(p\). For each jump \(u\) at \(i\) the first-level steps are
\[
       -u\le n<-u+W+1+2A\mathbf1_{|u|\le\eta}.
\]
For a low jump, \(|n|\le\eta+W+2A\). For a far jump,
\(n\) is within \(2W\) of one of \(\pm\beta_i\).
The exclusion at \(i\) therefore makes each step zero or prime to
\(p\), of cost \(3\) at level one. The target is antitone in \(u\)
on each interval and across the gap since \(\beta_i-\eta>2A\).
The total first-level contribution is at most
\(3(W+1)r_i+6A n_i\).
Finally the depth at level zero may be lowered from
\(T\) to \(d_j(0,u)\) without any loss.
The contribution at \(j\) is therefore at most
\[
 (2A+1)n_j+3(W+1)r_{\varphi^{-1}j}+6A n_{\varphi^{-1}j}.
\]
Let \(D_j\subseteq\iota_j(C)\) be the resulting lattice of
constants. Set \(V=\bigcap_j\iota_j^{-1}(D_j)\). The map
\[
 C/V\hookrightarrow\bigoplus_j\iota_j(C)/D_j,\qquad
 c_0+V\longmapsto\bigl(\iota_j(c_0)+D_j\bigr)_j
\]
is injective by the definition of \(V\). Since Frobenius
permutes the indices, summing gives
\[
 \operatorname{length}_{\mathbf Z_p}(C/V)
 \le (8A+1)\sum_j n_j+3(W+1)\sum_j r_j
 \le 11\sum_j\bigl((A+1)n_j+(W+1)r_j\bigr).
\]
Thus \(V\) has full \(\mathbf Z_p\)-rank in \(C\). The multiplicities
in this bound are those of the fixed lattice \(L\). The same
\(V\) is used for every allowed partner, including at a cut
where the partner point stays fixed during the middle change.

For any allowed partner and pairing, the constructed
points give safe moves, with a cut at every index.
On the partner take points over every constant of depth
\(R(2A+\eta)\), starting from \(E^A\) and using the free
steps. Their orders are the improvement of \(R(-A)\)
by that depth. In particular their orders at level zero
are at least \(-u-\eta\). Suppose the fixed side has a positive
separated interval at \(j\). Write \(u\) for a partner jump and \(v\)
for a fixed-side jump at the index read at the level in question.
Put \(i=\varphi^{-1}j\), and denote the chosen partner and fixed-side
points at \(k\) by \(x_k,y_k\), respectively. The path is
\[
 (\varphi x_i,\varphi y_i)\longrightarrow(x_j,\varphi y_i)
                         \longrightarrow(x_j,y_j).
\]
These are the first two changes of the four-move path.

For the first change use the Tate-dual point obtained from the
\(E^A\)-annihilator at \(i\), transported once. The required
order-plus-depth bounds of Section~\ref{sec:lifting-depth-ranges}
are \(1\) at level zero, \(A-u-v\) at level one and \(-u-v\)
above. At zero, \(u\le\eta\) and the improved order of \(x_j\)
is nonnegative. The transported partner point also has nonnegative
order, and \(\varphi y_i\) has depth at least one by the transported
\(z\tau\) divided difference. At level one the partner difference
has order at least \(\min(-u,-u-\eta)=-u-\eta\).
At \(i\) the old fixed depth at level zero is \(2A+\eta-v\)
for the positive sign and \(\beta_i+W+A+\eta\) for the negative sign.
Both are at least \(A+\eta-v\), since \(-v\le\beta_i+W\).
At levels \(a\ge2\) the partner order is at least \(-u\)
and the fixed depth at least \(-v\).

For the second change use the Tate-dual point from the \(E^A\)-annihilator
at \(j\). The required sum is \(A-u-v\) at zero and \(-u-v\) above.
At zero, \(v\ge-\eta\), so the fixed-side order and current partner
depth satisfy
\[
       \min(0,-v)+2A+\eta-u\ \ge\ A-u-v .
\]
At every positive level both fixed-side points have order at least
\(-v\), including at level one by pole removal at \(i\). The current
partner has depth at least \(-u\), giving the required sum.

Suppose now that the fixed side has a negative separated interval
at \(j\). With \(x_k\) the partner point and \(y_k\) the fixed-side
point, change the fixed side and then the partner along
\[
 (\varphi x_i,\varphi y_i)\longrightarrow(\varphi x_i,y_j)
                         \longrightarrow(x_j,y_j),\qquad
                         i=\varphi^{-1}j .
\]
They occupy the second and third positions of the four-move path.
Continue to write \(u\) for a partner jump and \(v\) for a fixed-side
jump at the index read at level \(a\).

For the initial change use the once-transported Tate-dual lift from
the \(E^A\)-annihilator at \(i\). At level zero the depth
\(d_j(0,v)=\beta_j+W+A+\eta>0\) puts each tested form value of \(y_j\)
in \(K'+\mathrm{Fil}^1\subset\mathrm{Fil}^0\). The transported
\(\varphi y_i\) also has nonnegative order, and \(\varphi x_i\)
has depth at least one by the transported \(z\tau\) divided
difference. At level one both fixed points have order at least
\(-v\), using pole removal at \(i\) for \(\varphi y_i\).
The resting depth there is \(2A+\eta-u\) for either sign at \(i\).
Thus the sum is at least \(A-u-v\). At \(a\ge2\) the fixed
difference has order at least \(-v\) and the resting depth is
at least \(-u\). These are the three required bounds of
Section~\ref{sec:lifting-depth-ranges}, also when \(i=j\).

The partner now moves with \(y_j\) at rest. On the fixed side,
\(-v\le\beta_j+W\) at level zero, so \(d_j(0,v)\ge A+\eta-v\);
at level one, \(d_j(1,v)\ge-v\). Its pole-free orders dominate
\(R(-\eta)\), and the partner has orders at least \(R(-\eta)\)
and depth \(R(2A+\eta)\).
For this move use a Tate-dual point at the predecessor
\(i=\varphi^{-1}j\), of depth \(R(A)\) and orders given by the
improvement of \(R(-A)\) by \(R(A)\). Its jumps are \(0\) and
\(-u-v-1\) for jumps \(u,v\) at \(i\). If both are low or both
are far, the dual jump has absolute value at most \(2\eta+W+1\)
because the two far intervals have opposite signs. If only one is far,
the dual jump is within \(\eta+W+1\) of one of \(\pm\beta_i\).
Thus the exclusion at \(i\) supplies the free steps. After Frobenius
transport from \(i\) to \(j\) the level-zero dual order is at least zero. For the moving order write \(o(u)\) for the improvement of
\(-A-u\) by \(2A+\eta-u\). If \(u<-A\), then \(o(u)>0\).
If \(u\ge-A\), then \(o(u)=\min(0,2A+\eta-u)\).
Since \(u\le\beta_j+W\), it follows in both cases that
\(\min(0,o(u))+\beta_j+W+A+\eta>0\). At levels \(>1\) the three
orders sum to the standard Tate order \(1\).

At level \(1\), both jumps are read at \(\varphi^{-1}j\), where
the two far intervals have opposite signs and the same start.
Write
\[
 E_q(w)=
 \begin{cases}q-w,&w<-A,\\ \min(q-w,0),&w\ge-A.\end{cases}
\]
The sum of the three orders there is at least
\[
 E_A(-u-v-1)+\min(-u,E_{2A+\eta}(u))
        + W+1-v+2A\mathbf1_{|v|\le\eta}.
\]
This is positive. For \(|v|\le\eta\), use
\(E_A(-u-v-1)\ge u+v+1-A\) and the moving order
\(\ge-u-A\), giving \(W+2\).
If \(v\) is far and \(u\) is low, the gap gives
\(E_A(-u-v-1)\ge A+u+v+1\), again giving \(W+2\).
If both are far, then \(|u+v|\le W\), the moving order
equals \(-u\), and \(E_A(-u-v-1)\ge u+v-W\).
The sum is at least \(1\). The safe-move criterion, with
the dual point in place of its standard lift, now applies.
Fix a partner cocycle \(z\) and \(c\in V\). Let \(e_j\) denote
the coefficient idempotent at \(j\). The partner point \(x_j\) and
the fixed-side point \(y_j\) lie over \(e_j(1\otimes1)\) and
\(e_j(1\otimes c)\). Put \(x=\sum_jx_j=(a,1)\) and
\(y=\sum_jy_j=(b,c)\). With \(i(j)=\varphi^{-1}j\),
\[
     \sum_j\varphi x_{i(j)}=\varphi x,\qquad
     \sum_j\varphi y_{i(j)}=\varphi y,
     \qquad \varphi(e_{i(j)})=e_j .
\]
For the intermediate first point use \(w_j=x_j\) at a positive
fixed sign and \(w_j=\varphi x_{i(j)}\) at a negative one;
the intermediate second point is \(w'_j=y_j\) for both signs.
Their sums lie over the same parameters as \(x,y\). The four
changes therefore add in the slots of
Section~\ref{sec:tame-cup-frobenius}. In a second-side change
the second argument is the Galois translate of the difference.
The group \(H\) fixes the idempotents and preserves the order
bounds. Differences over a common parameter have the asserted
kernel orders by Section~\ref{sec:local-tame-moves}.

There is also the endpoint test for every \(g,h\in H\).
At \(j\), pair \(gx_j-x_j\) with \(g(hy_j-y_j)\).
At level zero the sum of their depths is at least
\[
          (2A+\eta-u)+(A+\eta-v)\ge A-u-v
\]
for either fixed sign. At every positive level the sum is at
least \(-u-v\). Thus the current Tate-dual lift from the
\(E^A\)-annihilator gives total order at least one at every
level. Dividing the first difference by \(z\tau\) and returning
this factor after the contraction puts the value in the small
Witt ideal, as in Section~\ref{sec:lifting-depth-ranges}.

At depth \(m\) the integral pairing \(\psi\) gives the third
cocycle \((g\psi-\psi)/p^m\). Each third point used at \(j\)
projects to \(e_j(1\otimes1)\), also after transport from
\(i(j)\). Its extension form therefore reduces to \(\psi\bmod p^m\)
on the \(j\)-components by Section~\ref{sec:tame-cup-original-periods}.
Balance over \(O\) gives
\(\psi(e_j v,e_k v')=0\) for \(j\ne k\) after extension to
Witt scalars. Hence the four changes and the endpoint test still
hold after summation and composition with any
\(\mathbf Z_p\)-linear \(\lambda:O\to\mathbf Z_p\).
The two kernel differences are \(ga-a+z(g)\) and \(gb-b+c(g)\).
Section~\ref{sec:tame-cup-frobenius} with \(e_1=e_2=e_3=0\) now
gives a continuous coboundary for
\(((\lambda\otimes1)\psi\bmod p^m)_*(\bar z\cup\bar c)\)
on \(H=G_{E_v}\). This holds for every \(c\in V\), partner, depth
and test pairing. The annihilator is \(O\)-linear:
multiplying a cocycle by a coefficient scalar composes each
integral coefficient functional in the test with multiplication
by that scalar. This allows the constructed \(\mathbf Z_p\)
sublattice to bound the colength of the \(O\)-annihilator.

Write \(K_O=\operatorname{Frac}(O)\) and let \(d_j\) count its
coefficient embeddings at index \(j\). These are counted in the
de Rham grading over \(\mathbb C_p\), so each embedding occurs
at one index. Restriction of scalars in the order steps is
over \(\mathbf Z_p\): the multiplicity at a jump is the sum of
the graded multiplicities of those embeddings, and
\(\mathfrak S/(E)\) has rank \(f\) over \(\mathbf Z_p\).
Thus \(n_j=f d_j\mu\), \(r_j=f d_j\rho\) suffice for the
respective bounds per embedding \(\mu,\rho\), with
\(\sum_j d_j=[K_O:\mathbf Q_p]\). For the coefficient field of
Section~\ref{sec:relative-crystalline-colength}, the root \(\zeta_v\)
generates the unramified subfield of degree \(f=f_v\). Each of its
\(f\) embeddings extends to \([K_O:\mathbf Q_p]/f\) embeddings of
\(K_O\). Hence \(d_j\) equals this quotient on the Frobenius orbit
and vanishes elsewhere. These counts are over \(\mathbf C_p\).
Let \(\mathcal A\subset C\) be the annihilator on the integral
\(O\)-lattice. Since \(V\subseteq\mathcal A\), the
map \(C/V\to C/\mathcal A\) is a surjection over \(\mathbf Z_p\).
The residue degree of \(O\) converts \(O\)-length to
\(\mathbf Z_p\)-length, and \(v_O(p)\) is the ramification index.
Thus
\[
 \frac{\operatorname{length}_O(C/\mathcal A)}{v_O(p)}
 =\frac{\operatorname{length}_{\mathbf Z_p}(C/\mathcal A)}
        {[K_O:\mathbf Q_p]}
 \le\frac{\operatorname{length}_{\mathbf Z_p}(C/V)}
          {[K_O:\mathbf Q_p]}
 \le 11f\bigl(\mu(A+1)+\rho(W+1)\bigr).
\]

This estimate persists under finite extensions \(O'/O\) of
coefficient DVRs with the prescribed jump supports. Set
\(V_O=O\cdot V\subset C\), which is still in the cup annihilator.
Choose an \(O\)-basis \((b_i)\) of \(O'\), with coordinate maps
\(h_i:O'\to O\). Applying these maps to a cocycle and to its
crystalline primitive gives
\[
 Z_f^1(O'\otimes_O L)=\bigoplus_i b_i Z_f^1(L)=O'\otimes_O C.
\]
An allowed partner \(Y'\) over \(O'\) restricts to an allowed
\(O\)-lattice: at each coefficient embedding the graded spaces are
the direct sum over its extensions, and the class allows arbitrary
multiplicities on the prescribed support. For an integral
\(O'\)-pairing \(\psi'\) equivariant modulo \(p^m\), restrict the
fixed input to \(L\). Each \(h_i\circ\psi'\) is an integral
\(O\)-pairing equivariant modulo the same \(p^m\).
For any \(\lambda':O'\to\mathbf Z_p\) linear over \(\mathbf Z_p\),
\[
 \lambda'(x)=\sum_i \lambda_i(h_i(x)),\qquad
 \lambda_i(s)=\lambda'(b_i s).
\]
Write \(\mathcal A_O\) and \(\mathcal A_{O'}\) for the full
annihilators in \(C\) and \(O'\otimes_O C\). The coordinate
identity proves that every \(c\in\mathcal A_O\) passes the
tests against \(Y'\). Absorbing a scalar of \(O'\) into
\(\lambda'\) gives \(O'\otimes_O\mathcal A_O\subseteq\mathcal A_{O'}\).
Conversely, write \(c'=\sum_i b_i c_i\in\mathcal A_{O'}\).
For an allowed \(O\)-partner \(Y\), cocycle \(z\), integral pairing
\(\psi\) equivariant modulo \(p^m\), and functional \(\lambda\),
take \(O'\otimes_OY,\,1\otimes z,\,1\otimes\psi\) in the extended
class. Each extended embedding has the jumps prescribed at its
restriction. Composing the pairing with \(\lambda\circ h_i\)
gives exactly the original cup against \(c_i\) at that depth.
It vanishes, so \(c_i\in\mathcal A_O\) for every \(i\). Thus
\[
 \mathcal A_{O'}=O'\otimes_O\mathcal A_O,\qquad
 (O'\otimes_O C)/\mathcal A_{O'}=O'\otimes_O(C/\mathcal A_O).
\]
In particular \(O'\otimes_O V_O\) lies in this annihilator.
If \(e'=e(O'/O)\), then
\[
 \operatorname{length}_{O'}\bigl(O'\otimes_O(C/V_O)\bigr)
      =e'\operatorname{length}_O(C/V_O),\qquad
 v_{O'}(p)=e'v_O(p),
\]
so this lattice has the same normalized colength. The same
identities with \(\mathcal A_O\) in place of \(V_O\) give equality
of the normalized defects for the fixed lattice and its scalar
extension.

\subsection{Exact integrality in two special cases}
\label{sec:v2-crystalline-zero-loss}

\sectionopening{In this subsection, in two special cases the filtration bounds in
Section~\ref{sec:v2-crystalline-deep-constants} let us use the full
crystalline cocycle module. We verify the order conditions for the
coboundary calculation of Section~\ref{sec:tame-cup-frobenius} with
lifts over the unscaled constants.}

Keep the unramified chart and separated-interval
hypotheses of Section~\ref{sec:v2-crystalline-deep-constants}. At each index
\(i\) the two jump sets fit opposite signed intervals of far
start \(\beta_i\). Put \(R=6(A+\eta+W)\). The required
arithmetic hypothesis is that no nonzero multiple of \(p\) lies
within \(R\) of \(0\) or \(\pm\beta_i\). The three
upper jump bounds \((p-3)A\) give annihilators \(E^A\)
for both side cokernels and the Tate-dual cokernel.
Use the integral lift-ideal estimate of
Section~\ref{sec:v2-crystalline-lift-ideals}.
Write \(P_s\) for the order function
\(-u-s\mathbf1_{a=0}\), where \(a\) is Frobenius
level, and put
\[
 \begin{aligned}
 D_{s,h,h_1}(a,u)
   &=-u+(s+h)\mathbf1_{a=0}+h_1\mathbf1_{a=1},\\
 O_h(a,u)&=\begin{cases}\max(0,-u-h),&a=0,\\ -u,&a>0.\end{cases}
 \end{aligned}
\]
A point of order \(P_s\) has a pole of order at most \(s\)
at level zero. All points below lie over the given constant
in the given index, with no \(p\)-multiplier.

If neither side has a jump in \([-A,0)\), lift using \(E^A\)
and divide by \(E^A\). The initial orders and depths are
\(P_A\). Raise the depth to \(D_{A,\eta,0}\).
The steps lie at level zero in
\([-u-A,A+\eta-u)\). A low jump gives a step within
\(A+2\eta\) of zero and a far jump gives one within
\(A+\eta+W\) of its opposite far start. Both are within
\(R\). The contribution is zero by the unramified depth calculation.
If \(u<-A\), the initial order \(-u-A>0\) eliminates its
constant term and the new order is at least \(A+\eta-u\).
Otherwise \(u\ge0\), and the improved order is at least
\(\min(A+\eta-u,0)\ge-u\). Thus the new points have
order \(P_0\).

At each index let the side with the negative separated interval
rest. Its jumps \(v\le\eta\) imply that its just-constructed
point has order \(O_\eta\): for \(v<-\eta\) use \(P_0\),
and otherwise use the nonnegative depth
\(A+\eta-v\) to rule out a pole. The moving side has
\(u\ge-\eta\) at level zero and both its current and
preceding points have order \(P_0\).
The difference \(\varphi x_{i^{p^{-1}}}-x_i\) has
orders at least \(\min(0,-u)\) at level zero and
\(-u\) at levels \(\ge1\).
Adding the resting depth gives at least \(A-u-v\)
at zero and \(-u-v\) at the other levels. These are
exactly the order bounds required for the Tate-dual annihilator
\(E^A\) at the current index. At index \(i\), put \(i_-=i^{p^{-1}}\). The path from
\((\varphi x_{i_-},\varphi y_{i_-})\) to \((x_i,y_i)\) passes through
\((x_i,\varphi y_{i_-})\) when the first side has negative far sign
and through \((\varphi x_{i_-},y_i)\) when the second side has
negative far sign.

Now let the second side have jump set in \([-\eta,0]\),
and put \(h=A+\eta\). Its lift ideal contains (1), by the
small-jump estimate with maximum jump \(0\).
Give the points on the first side orders \(P_h\) and
depth \(D_{A,h,0}\) by starting at \(P_A\) and
raising the depth for free. Give the points on the second
side orders \(P_0\) and the same depth by starting at
\(P_0\). The steps in either case are within
\(A+h+\eta+W\le R\) of the three centers. For the first
side the improved orders dominate \(P_h\) because
\(h\ge A\): if \(-u-A>0\) the order becomes
\(A+h-u\), and otherwise the improved order
\(\min(0,A+h-u)\) is at least \(-u-h\).
For the second side \(u\le0\), so the same constant-term
rule retains \(P_0\).

Where the first side has negative separated interval, its point
has stop orders \(O_h\) by the constant-term rule and
\(u\le\eta\le h\). The moving second side has pole-free
preceding and current points, with \(v\ge-h\).
The sum of the moved order and resting depth is at least
\(A-u-v\) at level zero and \(-u-v\) above.

Where the first side has positive separated interval, transport the
preceding point of the second side by \(\varphi\).
Its depth is initially \(1\) at level zero,
\(A+h-v\) at level one, and \(-v\) above.
Its orders are \(0\) at zero and \(-v\) above.
Raise the depth at zero to \(A+h-v\) by free steps;
these lie between \(1\) and \(A+h+\eta\).
The result has order at least \(O_h\), since
\(v\ge-\eta\ge-h\). Its depth dominates
\(D_{A,h,h}\), including the extra \(h\) at level one.
The moving first side has \(u\ge-h\) at level zero.
The difference between its preceding translated point
and its present point has orders bounded below by
\[
            -u-h\quad(a=0,1),\qquad -u\quad(a>1).
\]
The resting depth therefore gives \(A-u-v\) at zero
and \(-u-v\) at all other levels. This explains the
additional first-level depth required when the moving
point has a pole at the preceding index.
Let \(i\) precede \(j\) under Frobenius, write \(x_i,y_i\) and
\(x_j,y_j\) for the chosen first- and second-side points, and
write \(\widehat y_j\) for the level-zero improvement of
\(\varphi y_i\). At a positive far index on the first side the moves are
\[
 (\varphi x_i,\varphi y_i)\ \longrightarrow\
 (\varphi x_i,\widehat y_j)\ \longrightarrow\
 (x_j,\widehat y_j)\ \longrightarrow\ (x_j,y_j).
\]
In the middle move \(\widehat y_j\) retains the transported depth
\(A+h-v\) at level one. At a negative far index on the first side
use the two moves through \((x_j,\varphi y_i)\).

In the two cases put \(h_0=\eta\) and \(h_0=h\),
respectively. Every index is a cut of shift \(A\) and the
transport length is one. Write \(s,t\) for a moving and a
resting jump. For the move from the preceding point to the stop,
both moving points have nonnegative order at zero and the
transported resting point has depth at least one. At level one
the moving order is at least \(-s-h_0\), and the transported
resting depth at least \(A+h_0-t\). Above level one their sum
is at least \(-s-t\). These are the bounds for the transported
third point from the preceding cut. For a move from the stop
to the current point, the bounds at zero are \(-s-h_0\) and
\(A+h_0-t\); at positive levels their sum is at least \(-s-t\).
Use the third point at the current cut for this move.

At the endpoints both depths at zero are at least
\(A+h_0\) minus the respective jump, and both positive-level
depths dominate minus that jump. The pair of endpoint cocycle
differences thus passes the test against the current third
point, also after translating the second difference by an element
of \(H\). The constants on all three sides are unscaled.
The four-move calculation gives an actual cup coboundary
on the chart. In the off-level application below its group already
equals \(G_{E_v}\), so the cochain gives the required coboundary on
the original group. This yields zero normalized
defect.

\subsection{The bound for the original coefficient lattices}
\label{sec:off-level-coefficient-defect}

\sectionopening{In this subsection, we apply the preceding estimates to the coefficient primes away from
the original level. The bounds from
Sections~\ref{sec:v2-crystalline-stops},
\ref{sec:v2-crystalline-deep-constants}
and~\ref{sec:v2-crystalline-zero-loss} yield the required normalized
colength on the local class of
Section~\ref{sec:relative-crystalline-colength}.}

Use the original local class in
Section~\ref{sec:relative-crystalline-colength}, with
\(p=\ell\nmid M\), fixed local dual \(L^\dagger\), and
\(d_v=[E_v:\mathbf Q_\ell]=e_v^{\rm ram}f_v\). Write
\(\rho=\operatorname{rank}_O L\); write \(c_*=1\) for its
main summand and \(c_*=t\) for each additional character.
Only the actual jumps of \(L^\dagger\) are used on the
fixed side. The test-side jump set is the full original
upper set. The Kummer chart is that of \(E_v\), with
working group \(G_{E_v}\) itself. Take the common coefficient DVR
large enough to contain the residue-field root \(\zeta_v\) of
Section~\ref{sec:coefficient-hodge-tate-intervals}. Fix the lift and
its descent pair before varying the upper lattice. Denote the lift
number by \(j_0\), and the variable residue-field index by \(j\).
Then \(\delta_{j_0}=j_0\Delta+h-2t\) is fixed, while the start
\(\beta_j\) can vary along the Frobenius cycle.
At a residue parameter \(i\) write \(U_i,V_i\) for the unions of
the partner and fixed-side jump sets over coefficient embeddings
with residue index \(i\), with \(0\) adjoined to both. Put \(T_1=Y\),
\(T_2=L^\dagger\), \(T_3=(Y\otimes_O L^\dagger)^\dagger\).
Off \(M\), the first two are crystalline over \(E_v\) by
Section~\ref{sec:relative-crystalline-colength}, and the tensor
product and Tate dual give crystallinity of \(T_3\).
At depth \(\nu\), write \(c_1,c_2\) for the first two cocycles and
\(c_3(g)=(g\psi-\psi)/p^\nu\). The last is integral by equivariance
modulo \(p^\nu\). Each has a primitive \(b_i\) over \(B_{\rm cris}\).
For \(i=1,2\), apply \(\varphi\otimes1\) to a \(B_{\max}\)-primitive.
Frobenius commutes with \(g\), fixes \(1\otimes c_i(g)\), and
\(\varphi(B_{\max})\subset B_{\rm cris}\). For \(i=3\), take
\(b_3=\psi/p^\nu\). Under \(g(v,s)=(gv+s c_i(g),s)\) the vector
\((-b_i,1)\) is an invariant lift of \(1\), so the rational
representation of each extension is crystalline. For the extension
by the entire cocycle module in Section~\ref{sec:v2-crystalline-deep-constants},
choose primitives on an \(O\)-basis and extend \(O\)-linearly.
The parameter quotient contributes only the jump \(0\).
The jump set for the third extension at \(i\) lies in
\(\{0\}\cup\{-u-v-1:u\in U_i,\ v\in V_i\}\), since
an \(O\)-bilinear pairing pairs only equal coefficient embeddings.

\indent\textit{Small coefficient primes.}
Suppose \(\ell\le Q_*\). Apply the sparse case of
Proposition~\ref{prop:crystalline-filtration-cup-exponent} with
\[
 \eta_0=2G,\quad B_0=b+\delta_{j_0},\quad G_0=G+\delta_{j_0},\quad
 J=\lfloor\log_\ell h^3\rfloor,\quad D_0=2K_\ell,\quad X_0=K_\ell .
\]
The actual far width is \(W_*\). The thin intervals fit
inside \([-\eta_0,\eta_0]\), since \(G_0+W_*\le2G\).
The thick-index sparsity for successive gaps at most
\(J\) is the one proved for these coefficient parameters.
Here \(Q_*=kt=(r+1-m)5r=R_{\rm prog}\), so the signed-ratio
exclusion applies throughout this range of coefficient primes.
If \(\ell\nmid q_1d\), the chart is unramified, \(\ell>7\), and
\(10\ell^J\eta_0=20\ell^JG\le20h^3G<B_0\). On a segment with
opposite thick signs use the pivot at step \(J+1\) of
Section~\ref{sec:v2-crystalline-stop-hypotheses}.
If \(\ell\mid q_1d\), the thick signs are constant at the fixed
completion and \(8 e_v^{\rm ram}<h^3<\ell^{J+1}\). The thin
intervals of either sign lie in \([-\eta_0,\eta_0]\), including
at \(\ell=d\). On every segment with a common thick sign choose the
side with negative far interval there as the resting side, and choose
either side if the segment is entirely thin. On a segment of length
\(l\), the resting jumps satisfy \(v\le\eta_0\) and the moving jumps
satisfy \(u\ge-\eta_0\) at levels \(0\le a<l\). At \(a=l\) the
plain stop uses the depth transported from the preceding cut.
All sparse cut shifts are zero, so these segments have the plain stop
of Section~\ref{sec:lifting-plain-cut}.
At the dyadic prime,
\(\lambda_2\ge v_2((2b)!)\ge b\), by the \(b\) even
integers in the factorial, so \(K_2=\lfloor
2\lambda_2/Q_*\rfloor\ge1\) for large \(x\).
At a prime dividing \(q_1\), use the splitting of \(K\) and put
\(a_\ell=v_\ell(q_1)\). The roots of order prime to \(\ell\) give
an unramified extension, and \(\Phi_{\ell^{a_\ell}}(1+X)\)
is Eisenstein. Thus \(e_v^{\rm ram}=\varphi(\ell^{a_\ell})\)
and \(D_E=1+v_\ell(e_v^{\rm ram})=a_\ell\). Here \(e_v^{\rm ram}>1\)
also at \(\ell=2\), since \(16\mid q_0\). For the fixed progression
modulus \(q_0\), the definition \(q_1=3L_0\,2^s q_0\) gives
\[
 a_2=s+v_2(q_0),\qquad
 a_3=\lfloor\log_3s\rfloor+5+v_3(q_0),\qquad
 a_\ell=v_\ell(q_0)\quad(\ell\ge5,\ \ell\mid q_0).
\]
At \(\ell=d\) the cyclotomic part is unramified and the quadratic
part gives \(e_v^{\rm ram}=2\), \(D_E=1\). In the dyadic depth
calculation of Section~\ref{sec:local-offlevel-depth-operators},
\(v_2(\theta\Delta_z)=1+1/e_v^{\rm ram}>1\).
The plain stop applies at all these ramified primes. At \(\ell=2\),
let \(j\) be a thick position on the periodically extended cycle.
Section~\ref{sec:coefficient-hodge-tate-intervals} excludes another
thick position at positive distance less than \(L_0\) in either
direction, and \(J<L_0\). Therefore \(j\) and \(j-J-1\) are cuts.
Every position strictly between them has \(j\) among its next \(J\)
positions. The preceding cut is exactly \(J+1\) steps before \(j\),
and both shifts are zero. If \(j\) is a thin cut, the positions
\(j,\ldots,j+J-1\) are thin, so \(j-1\) is a cut and the preceding
distance is one. This includes every position on a cycle without
thick positions. The plain-cut argument in both cases uses the
divisor \(z\tau\) on the same localized kernel.

The sparse depth inequality follows from
\[
 B_0+5\eta_0<2B_0,\quad
 2B_0<h^2K_\ell,\quad e_v^{\rm ram}\le h,\quad h^3<\ell^{J+1}.
\]
These give \(2e_v^{\rm ram}(B_0+5\eta_0)
<\ell^{J+1}(2K_\ell)\).
For the other sparse condition,
\(J+1+D_E\le h\) and
\(\lfloor\log_\ell(2(B_0+5\eta_0))\rfloor <m\le h\).
Indeed \(h^3\), \(e_v^{\rm ram}\), and \(2(B_0+5\eta_0)\)
are at most \(4b\), and \(\log_\ell(4b)\le\log(4b)/\log2<m\).
Thus \(J\le m-1\), \(D_E\le m\) and \(J+1+D_E\le2m\le h\).
The second sparse condition is
\[
 20(J+1+D_E)\eta_0
   +\lfloor\log_\ell(2(B_0+5\eta_0))\rfloor
 \le 40hG+h\le41hG<K_\ell .
\]
Here \(\eta_0=2G\), \(G\ge1\), and the last inequality follows from
\(K_\ell>2B_0/h^2\ge 2b/h^2\) and \(b\ge99G(h^3+1)\). The uniform exponent is
at most \(10^3 K_\ell(J+3)\le10^4 K_\ell m\).
At \(\ell=2\), the parameters above give \(\kappa=12(J+3)K_2+1\).
For each allowed partner, pairing, functional and finite depth,
Proposition~\ref{prop:crystalline-filtration-cup-exponent} gives
the same annihilating power \(2^{51\kappa}\). The chosen points may
depend on the test. By bilinearity and the definition in
Section~\ref{sec:relative-crystalline-colength}, the cup annihilator
therefore contains \(2^{51\kappa}Z_f(L^\dagger)\).

\indent\textit{Large coefficient primes with large loss.}
All rational primes ramified in \(E\) divide \(q_1d\).
Since \(q_1/q_0=3L_0\cdot 2^s\) has prime divisors only \(2,3\), and
\(6\mid q_0\), this set is contained in the prime support of the
fixed integer \(q_0d\). The cutoff \(Q_*=kt\) tends to infinity
with \(x\). For sufficiently large \(x\), every \(\ell>Q_*\)
therefore gives an unramified chart and satisfies \(\ell\ge7\).
Suppose \(\ell>Q_*\) and \(\lambda_\ell\ge h^3\). Use the two-start
case of Proposition~\ref{prop:crystalline-filtration-cup-exponent}
with \(B_0=b+\delta_{j_0}\), \(G_0=G+\delta_{j_0}\).
Here \(\lambda_\ell\ge h^3>Q_*\), so
\(K_\ell=\lfloor2\lambda_\ell/Q_*\rfloor\ge\lambda_\ell/Q_*\).
Since \(\lambda_\ell\ge1\) and \(\ell(\lambda_\ell+1)>2b\),
we have \(\ell\lambda_\ell\ge\ell(\lambda_\ell+1)/2>b\).
The required sizes follow from
\[
 \begin{aligned}
 K_\ell\ell^2-b&>
    (\ell\lambda_\ell/Q_*)(\ell-Q_*)\ge h^3>\delta_{j_0},\\
 K_\ell\ell&>b/Q_*>2(G+\delta_{j_0}),\\
 K_\ell&\ge h^3/Q_*>W_*+4\eta .
 \end{aligned}
\]
The argument of the logarithm in the two-start estimate
is at most \(4b\). Since
\(2\log(4b)<m\) and \(\log\ell>1/2\), the integer
\(\lfloor\log_\ell(2(B_0+W_*+4\eta))\rfloor\)
is at most \(m-1\). With \(K_\ell\ge1\), the resulting
uniform exponent is at most \(10^4 K_\ell m\).

\indent\textit{Large coefficient primes with small loss.}
For \(\ell>Q_*\) and \(\lambda_\ell<h^3\), put
\(A=2(\lambda_\ell+1)\). The missing-multiple radius
\(6(A+\eta+W_*)\) is exactly the one proved above.
It is smaller than \(G\), so in particular smaller
than both far starts. Also
\[
 \beta_j+W_*+\eta<2b<\ell A/2\le(\ell-3)A.
\]
The last inequality uses \(\ell\ge7\). For any jumps \(u,v\)
at the same index, the opposite signs of the separated intervals
give the estimate of Section~\ref{sec:v2-crystalline-deep-constants}
for the Tate-dual jump,
\[
 -u-v-1\le\beta_j+W_*+\eta-1<(\ell-3)A .
\]
Thus the required upper bound on the jumps holds on all
three sides of the cup pairing. For \(0\le a<t'\), the dual of the
Tate character lattice \(O(-a)\) is \(O(a+1)\), of jump
\(-a-1\in[-\eta,0]\). Section~\ref{sec:v2-crystalline-zero-loss}
therefore applies to its full crystalline cocycle module.

Equation~\eqref{eq:fixed-character-cup-annihilator} gives zero
defect for the full original sum of characters. For the main summand, if
\(\lambda_\ell<\lfloor h/5\rfloor\), then \(A<h/2\)
for the scales. Every negative jump on either
side is at most \(-h/2\), so the same section gives
zero defect by the gap below zero. In particular
\(K_\ell=0\) implies \(2\lambda_\ell<Q_*\) and
\(\lambda_\ell<\lfloor h/5\rfloor\) by the scale inequalities,
so both defects vanish in that case.

It remains to treat a main summand with
\(\lfloor h/5\rfloor\le\lambda_\ell<h^3\). Use the common submodule from
Section~\ref{sec:v2-crystalline-deep-constants}, with
at most \(\mu=\lfloor\rho m/(kc_*)\rfloor\) low jumps
per embedding. Its normalized colength is at most
\[
 11f_v\bigl(\mu(A+1)+\rho(W_*+1)\bigr)
  \ \le\ 11d_v\bigl(6\mu\lambda_\ell+\rho K_\ell\bigr).
\]
Here \(A+1\le6\lambda_\ell\) and \(W_*+1\le K_\ell\).
For the second inequality use
\(\lambda_\ell\ge\lfloor h/5\rfloor\) and
\(99(W_*+1)Q_*\le h\) to see
\(2\lambda_\ell\ge Q_*(W_*+1)\).
The bound comes from the colength of this single
submodule, whose cup vanishes unscaled for every
partner and every pairing.

On the two branches with a uniform exponent the
rank bound for the crystalline cocycle lattice gives
normalized defect at most \(2d_v\rho\cdot10^4K_\ell m\).
In the remaining branch the bound is
\(11d_v(6\mu\lambda_\ell+\rho K_\ell)\).
Both estimates are bounded by the expression below, since \(m\ge1\).

\begin{proposition}[Off-level crystalline colength]
\label{prop:crystalline-off-level-colength}
Fix \(q_0,d,\epsilon\) and take the arithmetic datum at a
sufficiently large scale \(x\) as in Section~\ref{sec:datum}.
For every \(\ell\nmid M\) and every \(v\mid\ell\), fix the
lift, descent pair, common coefficient ring and local class as above.
Let \(L\) be a main or
additional character lattice in this class, and let
\(\mu=\lfloor\rho m/(kc_*)\rfloor\). Then the normalized defect
of the cup annihilator in \(Z_f(L^\dagger)\) satisfies
\[
 \frac{\Delta_v(L)}{v_{\rm ord}(\ell)}
       \le 10^5 d_v(\rho K_\ell m+\mu\lambda_\ell).
\]
For \(\ell>Q_*\), the defect vanishes for an additional character if
\(\lambda_\ell<h^3\), and for a main summand if
\(\lambda_\ell<\lfloor h/5\rfloor\).
\end{proposition}
The comparison with \(d_v\rho F_\ell/c_*\) in
Section~\ref{sec:coefficient-local-defects} now applies.

\phantomsection

\begin{samepage}
\bigskip
\noindent Department of Mathematics, Purdue University, West Lafayette, IN, 47907, USA

\smallskip
\noindent\href{mailto:liu2053@purdue.edu}{\texttt{liu2053@purdue.edu}}
\end{samepage}
\end{document}